\pdfoutput=1

\documentclass[12pt,oneside]{book}

\usepackage[a4paper,lmargin=40mm,rmargin=25mm,tmargin=25mm,bmargin=25mm,headsep=2.5mm,headheight=15pt]{geometry}
\usepackage[T1]{fontenc}
\usepackage[utf8]{inputenc}
\usepackage{microtype}
\usepackage{setspace}
\usepackage{fancyhdr}
\usepackage{titlesec}
\usepackage{chngcntr}

\newcommand{\kaustthesistitle}{From Imperfect Information to Controlled Steps: Stochastic Optimization under Heavy-Tailed Noise and Compression}
\newcommand{\kaustauthor}{Abdurakhmon Sadiev}
\newcommand{\kaustapprovaldate}{October 2026}
\newcommand{\kaustcommitteechair}{Peter Richtárik}
\newcommand{\kaustcommitteemembers}{%
  Francesco~Orabona, Rolf~Krause, Arkadi~Nemirovski, Guanghui~Lan}
\providecommand{\ThesisIncludeAcknowledgments}{1}

\titleformat{\part}[display]
  {\normalfont\bfseries\centering\fontsize{14}{16}\selectfont}
  {\partname\ \thepart}{0.5em}{}
\titlespacing*{\part}{0pt}{20pt}{2\baselineskip}
\titleformat{\chapter}[display]
  {\normalfont\bfseries\fontsize{14}{16}\selectfont}
  {\chaptername\ \thechapter}{0.5em}{}
\titlespacing*{\chapter}{0pt}{20pt}{2\baselineskip}
\titleformat{name=\chapter,numberless}[display]
  {\normalfont\bfseries\centering\fontsize{14}{16}\selectfont}
  {}{0pt}{}
\titleformat{\section}
  {\normalfont\bfseries\fontsize{14}{16}\selectfont}
  {\thesection}{0.75em}{}
\titleformat{\subsection}
  {\normalfont\bfseries\fontsize{14}{16}\selectfont}
  {\thesubsection}{0.75em}{}
\titleformat{\subsubsection}
  {\normalfont\bfseries\fontsize{14}{16}\selectfont}
  {\thesubsubsection}{0.75em}{}

\fancypagestyle{plain}{%
  \fancyhf{}
  \fancyhead[C]{\thepage}
}
\counterwithout{footnote}{chapter}

\usepackage{amsmath,amssymb,amsfonts,amsthm,mathtools,bm,mathrsfs}
\usepackage{bbm}
\usepackage{nicefrac}
\usepackage{accents}
\usepackage{aligned-overset}
\usepackage{derivative}
\usepackage{suffix}
\usepackage{thmtools}
\usepackage{thm-restate}
\DeclareMathAlphabet{\mathsfit}{\encodingdefault}{\sfdefault}{m}{sl}
\SetMathAlphabet{\mathsfit}{bold}{\encodingdefault}{\sfdefault}{bx}{n}

\usepackage{graphicx}
\usepackage{float}
\usepackage{booktabs,array,tabularx,tabulary,multirow,makecell}
\newcolumntype{C}[1]{>{\centering\arraybackslash}p{#1}}
\newcolumntype{R}[1]{>{\raggedright\arraybackslash}p{#1}}
\usepackage{threeparttable,colortbl,hhline,diagbox,adjustbox}
\usepackage{caption}
\usepackage{subcaption}
\usepackage{placeins}
\usepackage{wrapfig}
\usepackage{listings}
\usepackage{tikz}
\usetikzlibrary{arrows.meta,calc,positioning}
\usepackage{tikz-cd}
\usepackage{pgfplots}
\usepackage{tcolorbox}
\usepackage{enumitem}
\usepackage{algorithm}
\usepackage{algpseudocode}
\let\thesisAlgorithmicFor\For
\let\thesisAlgorithmicEndFor\EndFor
\let\thesisAlgorithmicIf\If
\let\thesisAlgorithmicEndIf\EndIf
\let\thesisAlgorithmicElsIf\ElsIf
\providecommand{\STATE}{\State}
\providecommand{\FOR}{\For}
\providecommand{\ENDFOR}{\EndFor}
\providecommand{\IF}{\If}
\providecommand{\ENDIF}{\EndIf}
\providecommand{\ELSIF}{\ElsIf}
\providecommand{\REQUIRE}{\Require}
\providecommand{\ENSURE}{\Ensure}
\AtBeginEnvironment{algorithmic}{%
  \let\For\thesisAlgorithmicFor
  \let\EndFor\thesisAlgorithmicEndFor
  \let\If\thesisAlgorithmicIf
  \let\EndIf\thesisAlgorithmicEndIf
  \let\ElsIf\thesisAlgorithmicElsIf
}
\usepackage{pifont}
\usepackage{xspace}
\usepackage[dvipsnames]{xcolor}
\definecolor{BrickRed}{RGB}{143,20,2}
\definecolor{niceblue}{rgb}{0.0,0.19,0.56}
\providecommand{\ThesisProofreadingColours}{0}
\ifnum\ThesisProofreadingColours=1
  \definecolor{navy}{RGB}{0,0,128}
  \definecolor{darkred}{RGB}{132,0,0}
\else
  \colorlet{navy}{black}
  \colorlet{darkred}{black}
\fi
\definecolor{unitdirred}{RGB}{132,0,0}
\definecolor{mydarkgreen}{RGB}{39,130,67}
\definecolor{linen}{HTML}{FAF0E6}
\definecolor{steponepurple}{HTML}{4B136F}
\providecommand{\StepOneReviewMode}{0}
\ifnum\StepOneReviewMode=1
  \colorlet{steponeink}{steponepurple}
\else
  \colorlet{steponeink}{black}
\fi
\newcommand{\steponechange}[1]{{\color{steponeink}#1}}
\definecolor{steponeonepurple}{HTML}{4B136F}
\providecommand{\StepOneOneReviewMode}{0}
\ifnum\StepOneOneReviewMode=1
  \colorlet{steponeoneink}{steponeonepurple}
  \DeclareRobustCommand{\stepfive}[1]{{\color{red}#1}}
\else
  \colorlet{steponeoneink}{black}
  \DeclareRobustCommand{\stepfive}[1]{#1}
\fi
\DeclareRobustCommand{\steponeonechange}[1]{{\color{steponeoneink}#1}}
\definecolor{stepfivemetapurple}{HTML}{4B136F}
\providecommand{\StepFiveMetaReviewMode}{0}
\ifnum\StepFiveMetaReviewMode=1
  \colorlet{stepfivemetaink}{stepfivemetapurple}
\else
  \colorlet{stepfivemetaink}{black}
\fi
\DeclareRobustCommand{\stepfivemetachange}[1]{{%
  \ifnum\StepFiveMetaReviewMode=1
    \colorlet{niceblue}{stepfivemetaink}%
    \colorlet{algnamecolor}{stepfivemetaink}%
    \color{stepfivemetaink}%
  \fi
  #1}}
\usepackage{xurl}
\usepackage{csquotes}

\usepackage[
  backend=biber,
  style=authoryear,
  citestyle=authoryear,
  natbib=true,
  maxcitenames=2,
  mincitenames=1,
  maxbibnames=99,
  uniquelist=false,
  uniquename=init,
  sorting=nyt
]{biblatex}
\DeclareDelimFormat{nameyeardelim}{\addcomma\space}
\DeclareDelimFormat{multicitedelim}{\addsemicolon\space}
\newcommand{\tabcite}[1]{\textcite{#1}}
\usepackage[
  colorlinks=true,
  hypertexnames=false,
  linkcolor=niceblue,
  citecolor=niceblue,
  urlcolor=niceblue,
  filecolor=niceblue,
  menucolor=niceblue,
  runcolor=niceblue
]{hyperref}
\hypersetup{
  pdftitle={\kaustthesistitle},
  pdfauthor={\kaustauthor}
}

\DeclareFieldFormat{citehyperref}{%
  \DeclareFieldAlias{bibhyperref}{noformat}%
  \bibhyperref{#1}}

\DeclareFieldFormat{textcitehyperref}{%
  \DeclareFieldAlias{bibhyperref}{noformat}%
  \bibhyperref{%
    #1%
    \ifbool{cbx:parens}
      {\bibcloseparen\global\boolfalse{cbx:parens}}
      {}}}

\savebibmacro{cite}
\savebibmacro{textcite}

\renewbibmacro*{cite}{%
  \printtext[citehyperref]{%
    \restorebibmacro{cite}%
    \usebibmacro{cite}}}

\renewbibmacro*{textcite}{%
  \ifboolexpr{ ( not test {\iffieldundef{prenote}} and test {\ifnumequal{\value{citecount}}{1}} )
               or ( not test {\iffieldundef{postnote}} and test {\ifnumequal{\value{citecount}}{\value{citetotal}}} ) }
    {\DeclareFieldAlias{textcitehyperref}{noformat}}
    {}%
  \printtext[textcitehyperref]{%
    \restorebibmacro{textcite}%
    \usebibmacro{textcite}}}
\defbibheading{kaustreferences}[REFERENCES]{%
  \chapter*{\MakeUppercase{#1}}%
  \phantomsection
  \addcontentsline{toc}{chapter}{\textbf{\MakeUppercase{#1}}}%
}
\AtBeginBibliography{\singlespacing\setlength{\bibitemsep}{12pt}}

\providecommand{\backrefalt}[4]{}
\providecommand{\bibsection}{}
\providecommand{\eq}{}
\providecommand{\st}{}

\providecommand{\Circled}[1]{\textcircled{\scriptsize #1}}
\providecommand{\say}[1]{``#1''}

\newenvironment{paperabstract}
  {\begin{quote}\noindent\textbf{Abstract.}\ }{\end{quote}}

\theoremstyle{plain}
\newtheorem{theorem}{Theorem}[chapter]
\newtheorem{lemma}[theorem]{Lemma}

\newtheorem{corollary}[theorem]{Corollary}
\newtheorem{assumption}[theorem]{Assumption}
\theoremstyle{definition}
\newtheorem{definition}[theorem]{Definition}

\theoremstyle{remark}
\newtheorem{remark}[theorem]{Remark}

\NewDocumentCommand{\Exp}{e{_} o g}{%
  \mathbb{E}%
  \IfValueT{#1}{_{#1}}%
  \IfValueTF{#2}{\left[#2\right]}{%
    \IfValueT{#3}{\left[#3\right]}%
  }%
}
\newcommand{\Prob}{\mathbb{P}}

\newcommand{\R}{\mathbb{R}}
\newcommand{\N}{\mathbb{N}}
\newcommand{\cO}{\mathcal{O}}

\newcommand{\cD}{\mathcal{D}}
\newcommand{\cF}{\mathcal{F}}
\newcommand{\cL}{\mathcal{L}}
\newcommand{\cU}{\mathcal{U}}
\newcommand{\sign}{\operatorname{sign}}

\newcommand{\eqdef}{\overset{\mathrm{def}}{=}}
\WithSuffix\newcommand\eqdef*{&\overset{\mathrm{def}}{=}}

\newcommand{\norm}[1]{\left\lVert #1 \right\rVert}

\newcommand{\normtwo}[1]{\left\lVert #1 \right\rVert_{2}}

\newcommand{\norminf}[1]{\left\lVert #1 \right\rVert_{\infty}}
\newcommand{\normop}[1]{\left\lVert #1 \right\rVert_{\mathrm{op}}}

\newcommand{\normstar}[1]{\left\lVert #1 \right\rVert_{\star}}
\newcommand{\sqnorm}[1]{\left\lVert #1 \right\rVert^{2}}

\newcommand{\abs}[1]{\left\lvert #1 \right\rvert}
\newcommand{\inner}[2]{\left\langle #1, #2 \right\rangle}

\newcommand{\circledOne}{\text{\ding{172}}}
\newcommand{\circledTwo}{\text{\ding{173}}}
\newcommand{\circledThree}{\text{\ding{174}}}
\newcommand{\circledFour}{\text{\ding{175}}}
\newcommand{\circledFive}{\text{\ding{176}}}
\newcommand{\circledSix}{\text{\ding{177}}}
\newcommand{\cmark}{{\color{PineGreen}\ding{51}}}
\newcommand{\xmark}{{\color{BrickRed}\ding{55}}}
\newcommand{\circledSeven}{\text{\ding{178}}}
\newcommand{\circledEight}{\text{\ding{179}}}
\newcommand{\circledNine}{\text{\ding{180}}}
\newcommand{\circledTen}{\text{\ding{181}}}
\newcommand{\circledEleven}{\text{\ding{182}}}
\newcommand{\circledTwelve}{\text{\ding{183}}}
\renewcommand{\Circled}[1]{%
  {\color{navy}%
  \ifcase#1\relax
  \or\ifmmode\circledOne\else\ding{172}\fi
  \or\ifmmode\circledTwo\else\ding{173}\fi
  \or\ifmmode\circledThree\else\ding{174}\fi
  \or\ifmmode\circledFour\else\ding{175}\fi
  \or\ifmmode\circledFive\else\ding{176}\fi
  \or\ifmmode\circledSix\else\ding{177}\fi
  \or\ifmmode\circledSeven\else\ding{178}\fi
  \or\ifmmode\circledEight\else\ding{179}\fi
  \or\ifmmode\circledNine\else\ding{180}\fi
  \or\ifmmode\circledTen\else\ding{181}\fi
  \or\ifmmode\circledEleven\else\ding{182}\fi
  \or\ifmmode\circledTwelve\else\ding{183}\fi
  \else\textcircled{\scriptsize #1}%
  \fi}}

\newcommand{\unitdir}[1]{{\color{unitdirred}\frac{#1}{\norm{#1}}}}

\colorlet{algnamecolor}{BrickRed}
\colorlet{thesisAlgNameRed}{BrickRed}
\DeclareRobustCommand{\algname}[1]{{\color{thesisAlgNameRed}\small\sf#1}\xspace}
\newcommand{\algnamesmall}[1]{{\color{thesisAlgNameRed}\scriptsize\sf#1}\xspace}
\newcommand{\algnametiny}[1]{{\color{thesisAlgNameRed}\tiny\sf#1}\xspace}

\newif\ifshowoutrev
\showoutrevfalse
\ifshowoutrev
  \newcommand{\outrev}[1]{{\color{blue}#1}}
\else
  \newcommand{\outrev}[1]{#1}
\fi

\newcommand{\thesisnoformulascaling}[1]{%
  \ifmmode
    \PackageError{thesis}{#1 used inside a formula}{Formulas must stay at full
      size. Break the formula over several lines instead of scaling it.}%
  \fi}
\AtBeginDocument{%
  \let\thesisorigresizebox\resizebox
  \renewcommand{\resizebox}{\thesisnoformulascaling{\string\resizebox}\thesisorigresizebox}%
  \let\thesisorigscalebox\scalebox
  \renewcommand{\scalebox}{\thesisnoformulascaling{\string\scalebox}\thesisorigscalebox}%
  \let\thesisorigadjustbox\adjustbox
  \renewcommand{\adjustbox}{\thesisnoformulascaling{\string\adjustbox}\thesisorigadjustbox}%
}

\graphicspath{%
  {figures/ICML2023_High_probability_bounds_bounded_alpha_moment/arXiv/}%
  {figures/ICML2023_High_probability_bounds_bounded_alpha_moment/}%
  {figures/ICML2024_Prox_Clip/ICML2024/}%
  {figures/ICML2024_Prox_Clip/}%
  {figures/ICML2024_Prox_Clip/plots/}%
  {figures/Hessian_clipping-2/}%
  {figures/Hessian_clipping-2/plots/}%
  {figures/NSGD_with_heavy_tail-5/}%
  {figures/NSGD_with_heavy_tail-5/plots/}%
  {figures/ef21_generalized_smooth/NeurIPS_2025/}%
  {figures/ef21_generalized_smooth/}%
  {figures/2025_EFSkip-3/AISTAT2026/}%
  {figures/2025_EFSkip-3/}%
  {figures/2025_EFSkip-3/AISTAT2026/plots/}%
  {figures/MVR_second_order_momentum-2/neurips2026/}%
  {figures/MVR_second_order_momentum-2/}%
  {figures/MVR_second_order_momentum-2/iclr2026/}%
  {figures/MVR_second_order_momentum-2/iclr2026/plots/}%
}

\providecommand{\StepEightAReviewMode}{0}%
\definecolor{stepEightAReviewPurple}{HTML}{4B136F}%
\newcommand{\stepEightALegacyColor}[1]{%
  \ifnum\StepEightAReviewMode=1
    \colorlet{#1}{stepEightAReviewPurple}%
  \else
    \colorlet{#1}{black}%
  \fi
}%
\newcommand{\stepEightAAppendixColors}[1]{%
  \ifnum#1=2 \stepEightALegacyColor{algnamecolor}\fi
  \ifnum#1=3 \stepEightALegacyColor{algnamecolor}\fi
  \ifnum#1=4 \stepEightALegacyColor{blue}\fi
  \ifnum#1=5
    \stepEightALegacyColor{navy}%
    \stepEightALegacyColor{darkred}%
    \stepEightALegacyColor{red}%
    \stepEightALegacyColor{algnamecolor}%
  \fi
  \ifnum#1=6
    \stepEightALegacyColor{darkred}%
    \stepEightALegacyColor{algnamecolor}%
  \fi
  \ifnum#1=7 \stepEightALegacyColor{algnamecolor}\fi
  \ifnum#1=8 \stepEightALegacyColor{red}\fi
}%
\providecommand{\StepEightBReviewMode}{0}
\definecolor{stepEightBReviewPurple}{HTML}{4B136F}
\newcommand{\stepEightBchange}[1]{%
  {\ifnum\StepEightBReviewMode=1
    \color{stepEightBReviewPurple}%
    \colorlet{algnamecolor}{stepEightBReviewPurple}%
  \fi #1}%
}
\newcommand{\stepEightBdeletion}{%
  \ifnum\StepEightBReviewMode=1
    \rlap{\kern.15em\smash{\textcolor{stepEightBReviewPurple}{\rule[-.2ex]{.7pt}{1.4ex}}}}%
  \fi
}
\colorlet{thesisAlgorithmInk}{algnamecolor}
\providecommand{\StepEightCReviewMode}{0}
\definecolor{stepEightCReviewPurple}{HTML}{4B136F}
\newcommand{\stepEightCAlgorithmNames}{%
  \colorlet{algnamecolor}{thesisAlgorithmInk}%
  \ifnum\StepEightCReviewMode=1
    \colorlet{algnamecolor}{stepEightCReviewPurple}%
  \fi
}
\newcommand{\stepEightCchange}[1]{{%
  \ifnum\StepEightCReviewMode=1\color{stepEightCReviewPurple}\fi#1}}
\newcommand{\stepEightCTableAlgname}[1]{{%
  \stepEightCAlgorithmNames\colorlet{BrickRed}{algnamecolor}%
  \let\xspace\relax\algnamesmall{#1}}\xspace}
\DeclareRobustCommand{\normalgname}[1]{{%
  \let\xspace\relax
  \algname{\mbox{\ensuremath{\lVert}#1\ensuremath{\rVert}}}}\xspace}
\DeclareRobustCommand{\normEFname}[1]{{%
  \stepEightCAlgorithmNames\let\xspace\relax
  \normalgname{EF21\if\relax\detokenize{#1}\relax\else-#1\fi}}\xspace}
\newcommand{\kaustpreliminary}[1]{%
  \chapter*{\MakeUppercase{#1}}%
  \phantomsection
  \addcontentsline{toc}{chapter}{\textbf{\MakeUppercase{#1}}}%
}

\newcommand{\thesischaptersourcenote}[2]{%
  \begin{quote}
    \small\noindent
    This chapter is based on the joint work \enquote{#1} \parencite{#2}.
    The material has been adapted to the notation and organization of this
    dissertation.
  \end{quote}%
}

\providecommand{\StepSevenReviewMode}{0}
\definecolor{stepsevenpurple}{HTML}{4B136F}
\newcommand{\stepsevenink}{%
  \ifnum\StepSevenReviewMode=1
    \color{stepsevenpurple}%
    \colorlet{algnamecolor}{stepsevenpurple}%
    \colorlet{niceblue}{stepsevenpurple}%
  \fi
}
\newcommand{\stepsevenchange}[1]{{\stepsevenink #1}}
\newenvironment{stepsevenchanged}{\begingroup\stepsevenink}{\par\endgroup}

\begin{document}

\begin{titlepage}
  \centering
  \begin{doublespace}
    {\bfseries\fontsize{14}{16}\selectfont \stepsevenchange{\kaustthesistitle}\par}
    \vspace{0.25\baselineskip}
    Dissertation by\par
    \kaustauthor\par
    \vspace{0pt plus 1.7fill}
    In Partial Fulfillment of the Requirements\par
    For the Degree of\par
    Doctor of Philosophy\par
    \vfill
    King Abdullah University of Science and Technology\par
    Thuwal, Kingdom of Saudi Arabia\par
  \end{doublespace}
  \vfill
  \begin{singlespace}
    \copyright\ \stepEightCchange{\kaustapprovaldate}\\
    \kaustauthor\\
    All rights reserved%
  \end{singlespace}
\end{titlepage}

\pagenumbering{arabic}
\setcounter{page}{2}

\kaustpreliminary{Examination Committee Page}
\vspace*{3\baselineskip}
\begin{onehalfspace}
  \raggedright
  The dissertation of \kaustauthor\ is approved by the examination committee\par
  \vspace{6\baselineskip}
  Committee Chairperson: \kaustcommitteechair\par
  Committee Members: \kaustcommitteemembers\par
\end{onehalfspace}

\kaustpreliminary{Abstract}
\begin{center}
  \stepsevenchange{\kaustthesistitle}\\[0.5\baselineskip]
  \kaustauthor
\end{center}
\begin{stepsevenchanged}
\noindent
Stochastic optimization relies on information whose statistical quality, availability, and computational cost can vary substantially. Heavy-tailed observations challenge concentration, compression distorts communicated vectors, and generalized smoothness changes how far an algorithm can move. This dissertation studies these difficulties through two linked design decisions: how to construct an informative search vector and how to convert it into a controlled step.

We first develop high-probability methods under bounded central $p$-th moments, $p\in(1,2]$, permitting infinite variance when $p<2$. Clipped gradient, accelerated, extragradient, and gradient descent--ascent methods provide guarantees for stochastic minimization and variational inequalities, with smoothness and moment conditions imposed only on regions containing the iterates with high probability. For composite and distributed problems, shifted and innovation clipping preserve the deterministic fixed point that direct clipping can disturb. The distributed guarantees also retain linear speedup in the stated stochastic terms.

We then characterize the value of richer stochastic information through lower complexity bounds and corresponding algorithms. The bounds quantify when same-sample gradient differences and stochastic Hessian information improve heavy-tailed nonconvex optimization. Normalized momentum variance reduction and Hessian-corrected momentum attain matching or near-matching in-expectation bounds in the stated oracle regimes. Separately clipping gradients and their correction terms yields batch-free high-probability methods with polylogarithmic confidence dependence.

For compressed distributed optimization, we establish convergence of normalized error-feedback methods under generalized smoothness without a uniform gradient-dissimilarity assumption. Under standard smoothness, we develop parameter-agnostic variants with five momentum estimators whose decreasing stepsizes require neither smoothness constants nor the initial optimality gap. Finally, we extend the estimator hierarchy to linear minimization oracles over arbitrary norm balls. The analysis bounds the minimum expected dual gradient norm along the trajectory, recovering the corresponding Euclidean convergence orders under estimator-specific generalized smoothness.

Together, these results connect the statistical quality of stochastic information, its reconstruction under compression, and the geometry of the resulting update. They identify how controlled transformations and steps can improve convergence, and where lower bounds delimit the gains available from richer information.
\end{stepsevenchanged}

\ifnum\ThesisIncludeAcknowledgments=1
  \kaustpreliminary{Acknowledgments}

First and foremost, I would like to express my deepest gratitude to my PhD advisor, Professor Peter Richtárik. From the beginning of my internship and throughout more than four years of doctoral studies at KAUST, Peter has been a constant source of guidance, encouragement, and genuine care. Our many discussions, not only about mathematics and research but also about life, have meant a great deal to me. He has been much more than a scientific advisor: a mentor who inspired me to explore new ideas, encouraged me to take on challenging problems, and cared about my growth both as a researcher and as a person. Perhaps most importantly, I have always felt that Peter believed in me, and his belief has helped me believe more in myself. For everything he has done for me throughout these years, I will always be grateful.

My first research advisor, Professor Alexander Gasnikov, played a fundamental role in the beginning of my scientific career. He introduced me to interesting research problems, guided my early steps, and opened the door to the optimization research community. Through him, I had the opportunity to meet and collaborate with many talented researchers, including Alexander Beznosikov, Eduard Gorbunov, and Dmitry Kovalev, who have each contributed to my scientific development. I also greatly value the warm relationship that Professor Gasnikov and I have maintained to this day.

Eduard Gorbunov has been another important influence on my development as a researcher. Although we had known each other before, our scientific collaboration began during my internship at KAUST and continued throughout my PhD. Our joint research forms the basis of the first two technical chapters of this dissertation. Working with Eduard taught me the importance of dedication, perseverance, and hard work in scientific research. To me, he is a remarkable example of how mathematical talent, combined with exceptional diligence, can help solve challenging problems. I have learned a great deal from our collaboration and have always enjoyed our conversations, both about research and beyond.

I owe special thanks to Dmitry Kovalev, who has had a profound influence on the way I think about mathematics and conduct research. Working with him helped me discover not only the beauty of optimization, but also the intellectual pleasure of developing mathematical proofs and understanding the techniques behind them. He was the first person who truly showed me how research in optimization is done: how to approach challenging problems, search for the right ideas, and turn mathematical intuition into rigorous arguments. Through our discussions and joint work, Dmitry generously shared his mathematical insights, problem-solving approaches, and research techniques, many of which I have adopted in my own work. He played an essential role in shaping my mathematical taste and in developing my appreciation for elegant proofs and creative approaches to optimization. I am deeply grateful for everything I learned from him, and his influence continues to shape the way I approach research today.

Working with Laurent Condat over the past two years has been another particularly rewarding experience. Our collaboration and stimulating discussions on primal–dual optimization and splitting methods have allowed me to explore mathematical ideas that I find especially interesting. In many ways, this work has further developed the mathematical taste and appreciation for elegant optimization methods that I first gained through working with Dmitry Kovalev. I particularly value our continuing discussions and exchange of ideas.

I also wish to thank Professor Martin Takáč for the opportunity to spend two months as a research intern at MBZUAI. Although our collaboration was relatively short, his advice and guidance were valuable to me, and I thoroughly enjoyed working with him.

I am sincerely thankful to the members of my dissertation committee, Professors Francesco Orabona, Rolf Krause, Arkadi Nemirovski, and Guanghui Lan, for the time and attention they devoted to reading and evaluating my dissertation. Their thoughtful questions, constructive feedback, and participation in my PhD defense made this important milestone especially meaningful. It was a privilege to present and discuss my research with them.

My doctoral years were also shaped by the wonderful people I met in the Optimization and Machine Learning Lab at KAUST. I wish to acknowledge its past and present members and visiting researchers, including Adrien Fradin, Ahmad Rammal, Alexander Tyurin, Alina Abdikarimova, Ammar Mahran, Andrei Panferov, Artavazd Maranjyan, Artem Riabinin, Avetik Karagulyan, Dymitr Lubczyk, Elnur Gasanov, Georg Meinhardt, Hanmin Li, Hussein Rammal, Igor Klimczak, Igor Sokolov, Ivan Ilin, Kai Yi, Kaja Gruntkowska, Kirill Acharya, Konstantin Burlachenko, Lukang Sun, Mher Safaryan, Michał Grudzień, Omar Shaikh Omar, Philip Zmushko, Samuel Horváth, Sarit Khirirat, Simone Maria Giancola, Slavomír Hanzely, Tigran Melkonyan, Timur Kharisov, Wojciech Anyszka, Yassine Maziane, Yury Demidovich, and Zhirayr Tovmasyan. The many scientific discussions, exchanges of ideas, and conversations extending far beyond research made our laboratory a stimulating and welcoming environment. Being surrounded by people from different countries, cultures, and academic backgrounds broadened my perspective, both scientifically and personally. The collaborative atmosphere, friendships, and sense of community made my time at KAUST truly special.

I would also like to thank Marina Danilova, Gauthier Gidel, and Pavel Dvurechensky for their contributions to the joint research presented in Chapters 2 and 3 of this dissertation. I am also grateful to all my other coauthors for their contributions to our joint publications.

Beyond my academic mentors and collaborators, I am deeply indebted to the friends who have made this journey meaningful, both before and during my PhD.

My friendship with Ilyas Fatkhullin goes back to our second year at MIPT. I fondly remember the countless hours we spent studying together, laughing, and discussing mathematics and research. In many ways, it was our friendship and our shared experience of attending Alexander Gasnikov's lectures that first inspired me to pursue optimization. Seeing Ilyas's enthusiasm for the subject made me want to explore it myself, and what began as a shared interest eventually became the direction of my scientific career. Our common research interests and many conversations over the years have shaped my thinking in important ways. It means a great deal to me that we remain close and continue to exchange thoughts not only about science, but also about life.

Alexander (Sasha) Beznosikov has been both a close friend and an important mentor throughout my scientific journey. He has generously shared his experience, advised me on research directions, and helped me navigate the world of optimization. Beyond our scientific collaboration, his understanding, encouragement, and friendship have been especially valuable.

A very special place in my memories of KAUST belongs to Grigory (Grisha) Malinovsky and Anastasia (Nastya) Malinovskaya. Their warmth, kindness, and constant support made their home feel like a second home to me. I will always cherish the many evenings we spent together, our conversations, and the moments playing with their daughter, Masha. Being welcomed into such a warm family atmosphere brought me a sense of comfort and belonging while living far from my own family. Our conversations about mathematics, life, and philosophy, as well as Grisha's encouragement and belief in me, have meant more than I can easily express.

Egor Shulgin has been one of my closest friends during my years at KAUST. We have shared countless conversations about research, life, health, and sports, along with many memorable experiences traveling together. Our scientific discussions have always been stimulating, but what I value most is the friendship that extends far beyond our academic interests. His support, advice, and companionship have made these years richer, and I am glad that our friendship continues beyond KAUST.

I also want to warmly acknowledge Alexander (Sasha) Gaponov. Although our friendship developed relatively recently, we have shared a remarkable amount in a short time. From conversations about art and life to traveling together, exchanging jokes, and discovering new perspectives, our friendship has been a wonderful part of my final years at KAUST. Our discussions have introduced me to ideas and ways of looking at things that I might otherwise never have considered.

Above all, my deepest love and gratitude belong to my mother. Her unconditional love, care, and unwavering support have accompanied me throughout my life. She raised me, gave me the opportunity to pursue an education, and played an essential role in shaping the person I am today. Although reaching this point has required years of my own effort, I know that much of this journey would not have been possible without her. Her belief in me and everything she has done for my future have provided the foundation upon which I have been able to build my life.

My beloved sister, Laylo, has also been a source of love, encouragement, and support throughout the years. I feel fortunate to have her in my life and deeply cherish the bond we share.

To my mother and sister, I want to say that I love you both dearly. Thank you for being such an important part of my life and for helping me become who I am today. This achievement is not mine alone, and I am deeply grateful to share it with you.

\fi

\begin{onehalfspacing}
  \renewcommand{\contentsname}{TABLE OF CONTENTS}
  \tableofcontents
\end{onehalfspacing}

\cleardoublepage
\kaustpreliminary{List of Abbreviations}
\begin{singlespace}
  \begin{description}[style=nextline,leftmargin=42mm,labelwidth=36mm]
  \item[CC] Communication Compression
  \item[EF] Error Feedback
  \item[\steponeonechange{EF21}] \steponeonechange{Error Feedback 21}
  \item[EG] Extragradient
  \item[FOSP] First-Order Stationary Point
  \item[GD] Gradient Descent
  \item[GS] Generalized Smoothness
  \item[HM] Hessian-Corrected Momentum
  \item[IGT] Implicit Gradient Transport
  \item[LMO] Linear Minimization Oracle
  \item[MSGS] Mean-Squared Generalized Smoothness
  \item[MSS] Mean-Squared Smoothness
  \item[MVR] Momentum Variance Reduction
  \item[NSGD] Normalized Stochastic Gradient Descent
  \item[$p$-BCM] Bounded Central $p$-th Moment
  \item[$q$-WAS] $q$-Weak Average Smoothness
  \item[SEG] Stochastic Extragradient
  \item[SGD] Stochastic Gradient Descent
  \item[SGDA] Stochastic Gradient Descent--Ascent
  \item[SOGS] Second-Order Generalized Smoothness
  \item[SSTM] Stochastic Similar Triangles Method
  \item[VI] Variational Inequality
\end{description}

\end{singlespace}

\begin{onehalfspacing}
  \cleardoublepage
  \renewcommand{\listfigurename}{\centering\textbf{LIST OF FIGURES}}
  \phantomsection
  \addcontentsline{toc}{chapter}{\textbf{LIST OF FIGURES}}
  \listoffigures

  \cleardoublepage
  \renewcommand{\listtablename}{\centering\textbf{LIST OF TABLES}}
  \phantomsection
  \addcontentsline{toc}{chapter}{\textbf{LIST OF TABLES}}
  \listoftables
\end{onehalfspacing}

\chapter{Introduction}
\label{ch:chapter_1}
\providecommand{\FigureTablePolishReviewMode}{0}
\ifnum\FigureTablePolishReviewMode=1
  \definecolor{figuretablepolishink}{HTML}{4B136F}
\else
  \colorlet{figuretablepolishink}{black}
\fi
\DeclareRobustCommand{\figuretablepolishchange}[1]{{%
  \leavevmode
  \ifnum\FigureTablePolishReviewMode=1
    \colorlet{niceblue}{figuretablepolishink}%
    \colorlet{algnamecolor}{figuretablepolishink}%
    \colorlet{BrickRed}{figuretablepolishink}%
    \colorlet{PineGreen}{figuretablepolishink}%
    \color{figuretablepolishink}%
  \fi
  #1}}
\providecommand{\StepSixReviewMode}{0}
\ifnum\StepSixReviewMode=1
  \definecolor{stepsixink}{HTML}{4B136F}
\else
  \colorlet{stepsixink}{black}
\fi
\DeclareRobustCommand{\stepsixchange}[1]{{%
  \leavevmode
  \ifnum\StepSixReviewMode=1
    \colorlet{niceblue}{stepsixink}%
    \colorlet{algnamecolor}{stepsixink}%
    \color{stepsixink}%
  \fi
  #1}}
\label{ch1:introduction_chapter_1}

\section{Motivation and Scope}
\label{ch1:sec:intro_motivation_chapter_1}

Stochastic optimization studies problems in which an objective function is accessed through random estimates rather than evaluated or differentiated exactly at every iteration. A canonical formulation is
\begin{equation}
    \min_{x\in\R^d}
    \left\{
        f(x)
        \eqdef
        \Exp_{\xi\sim\cD}\bigl[f(x;\xi)\bigr]
    \right\},
    \label{ch1:eq:intro_stochastic_problem_chapter_1}
\end{equation}
where $x\in\R^d$ is the decision variable and $\xi\sim\cD$ is a data sample or stochastic-oracle realization. A standard stochastic-gradient method generates iterates of the form
\begin{equation}
    x^{k+1}
    =
    x^k-\gamma_k g(x^k;\xi^k),
    \qquad \gamma_k>0,
    \label{ch1:eq:intro_standard_sgd_chapter_1}
\end{equation}
where $g(x^k;\xi^k)$ is a stochastic gradient estimate. This iteration originates in stochastic approximation and remains a basic mechanism for large-scale optimization \parencite{robbins1951stochastic,nemirovski2009robust,lan2020first,bottou2018optimizationmethodsML}.

Training a neural network from sampled data is a principal example. Here $x\in\R^d$ denotes the vector of trainable model parameters, while $x^k$ denotes its value at iteration $k$. The random variable $\xi^k$ represents a sampled data point or mini-batch, and $g(x^k;\xi^k)$ is the corresponding stochastic-gradient estimate. When the model or dataset is too large for one machine, such estimates may be computed by multiple workers and communicated to a coordinating server. When the messages are compressed, the method must account both for the statistical quality of the local estimates and for the communication cost, compression distortion, and reconstruction accuracy \parencite{bottou2018optimizationmethodsML}.

Beyond minimization, Chapters~\ref{ch:chapter_2}--\ref{ch:chapter_3} study stochastic variational inequalities. A convenient composite formulation is to find $x^\star\in\R^d$ such that
\begin{equation}
    \inner{F(x^\star)}{x-x^\star}
    +\Psi(x)-\Psi(x^\star)
    \geq 0
    \qquad
    \text{for all }x\in\R^d,
    \label{ch1:eq:intro_composite_vi_chapter_1}
\end{equation}
where $F(x)=\Exp_{\xi\sim\cD}[F(x;\xi)]$ and $\Psi$ is a proper closed convex function. This formulation contains the first-order optimality conditions of composite convex minimization and also covers monotone inclusions, equilibrium problems, and convex--concave saddle-point problems \parencite{facchinei2003finite,ryu2021large}. Min--max formulations are central in game-theoretic models and machine learning; generative adversarial networks are a prominent example, and their training dynamics have been studied explicitly from a VI perspective \parencite{goodfellow2014generative,gidel2019variational}.

The update~\eqref{ch1:eq:intro_standard_sgd_chapter_1} couples two decisions. The orientation of $g(x^k;\xi^k)$ determines where the method moves, while its norm contributes directly to the displacement length. This coupling is effective when the stochastic estimate is sufficiently accurate and a conventional local model provides a trustworthy scale. The chapters of this dissertation study settings in which additional mechanisms are needed: the stochastic oracle may satisfy only weak moment bounds; compression may remove information from communicated messages; local curvature may grow with the gradient norm; theoretically valid schedules may depend on unavailable quantities; and richer gradient-difference or second-order observations may be available but must be used efficiently.

\paragraph{Weak moments and high-probability guarantees.}
Several empirical studies report strongly non-Gaussian or heavy-tailed stochastic-gradient behavior in the investigated deep neural networks, attention models, and reinforcement-learning problems \parencite{simsekli2019tail,zhang2020adaptive,garg2021proximal,battash2024revisiting}. These observations are model- and task-dependent: finite empirical samples cannot prove that a population variance is infinite, and heavy-tailedness does not by itself imply infinite variance. They nevertheless motivate oracle models that do not require a uniformly bounded second moment.

A recurring assumption in Chapters~\ref{ch:chapter_2}--\ref{ch:chapter_5} is access to an unbiased stochastic oracle satisfying a bounded central $p$-th moment condition:
\begin{equation}
    \Exp_{\xi\sim\cD}\bigl[g(x;\xi)\bigr] = \nabla f(x),
    \quad
    \Exp_{\xi\sim\cD}\bigl[\lVert g(x;\xi)-\nabla f(x)\rVert^p\bigr] \leq \sigma^p,
    \quad
    p\in(1,2].
    \label{ch1:eq:intro_p_bcm_chapter_1}
\end{equation}
For $p=2$, the second inequality is the usual uniform bounded-variance condition. For $p<2$, it controls only a lower-order central moment and does not require the variance to be finite. Analogous unbiased central-moment conditions are imposed on stochastic operators, and Chapters~\ref{ch:chapter_4}--\ref{ch:chapter_5} introduce corresponding assumptions for stochastic Hessian information. We use \emph{weak-moment noise} for the full range $p\in(1,2]$ and reserve \emph{infinite-variance regime} for $p<2$.

The literature also contains stronger noncentral moment assumptions, for example
\begin{equation}
    \Exp_{\xi\sim\cD}\left[\norm{g(x;\xi)}^p\right]\leq G^p.
    \label{ch1:eq:intro_p_bm_chapter_1}
\end{equation}
Condition~\eqref{ch1:eq:intro_p_bm_chapter_1} is not interchangeable with~\eqref{ch1:eq:intro_p_bcm_chapter_1}; in particular, it also restricts the magnitude of the population gradient. This distinction matters in comparisons with methods analyzed under bounded noncentral moments, including the normalized-momentum analysis of \textcite{cutkosky2021high}.

Expectation and high-probability guarantees answer different questions. For a performance measure $\mathcal{M}$, an expectation result controls $\Exp[\mathcal{M}(\widehat{x})]$, whereas a high-probability result requires
\begin{equation}
    \Prob
    \left\{
        \mathcal{M}(\widehat{x})\leq\varepsilon
    \right\}
    \geq
    1-\beta,
    \qquad
    \beta\in(0,1).
    \label{ch1:eq:intro_high_probability_chapter_1}
\end{equation}
The latter directly controls the probability of an unsuccessful run. Markov's inequality can convert an expectation estimate into a tail bound, but it typically introduces polynomial dependence on $1/\beta$. Direct high-probability analyses seek, under the same performance criterion, to retain the target-accuracy dependence of the relevant in-expectation benchmark while replacing that polynomial confidence penalty by logarithmic or polylogarithmic dependence on $1/\beta$, whenever the oracle and algorithmic structure permit it \parencite{nazin2019algorithms,gorbunov2020stochastic,cutkosky2021high,sadiev2023high,nguyen2023improved,fatkhullin2025can}.

\outrev{The output convention is part of the guarantee. A uniformly sampled iterate corresponds, in expectation, to the arithmetic trajectory average of the same performance measure, while a stepsize-weighted random iterate corresponds to a weighted trajectory average. Best-iterate and terminal-iterate criteria are different and cannot be substituted for these averages without an additional argument. The technical chapters therefore retain their theorem-specific output rules.}

\outrev{The technical obstruction is not only that an individual oracle response can be large. Natural optimization proofs generate martingale terms whose increments depend on the stochastic trajectory, for example
\(
    \inner{g(x^k;\xi^k)-\nabla f(x^k)}{x^k-x^\star}.
\)
On an unbounded domain, neither the distance to a solution nor the resulting increment is a priori bounded. Classical analyses therefore impose light-tail, bounded-domain, or bounded-gradient conditions. More recent approaches control the trajectory, transform the stochastic observation, or construct concentration arguments directly around the optimization recursion.}

\outrev{\stepsixchange{Such concentration arguments can yield high-probability guarantees under lighter-tailed models. \textcite{pmlr-v202-liu23aa} exploit optimization-specific martingale structure under sub-Gaussian noise, while \textcite{madden2024high} use self-normalized concentration for norm-sub-Weibull noise in vanilla nonconvex \algname{SGD}. These tail assumptions differ from the bounded-$p$-moment setting of Chapters~\ref{ch:chapter_2}--\ref{ch:chapter_5}; Section~\ref{ch1:subsec:intro_heavy_tail_chapter_1} gives the detailed comparison.}}

\outrev{Weak moments do not imply a criterion-independent failure of \algname{SGD}. On an unbounded quadratic with infinite-variance noise, even one nonadaptive \algname{SGD} step can make a second-moment quantity such as $\Exp{\norm{\nabla f(x^1)}^2}$ infinite. This shows that a second-moment criterion may be incompatible with a $p<2$ oracle; it does not show that every meaningful notion of convergence fails. \textcite{fatkhullin2025can} establish positive guarantees for vanilla or projected \algname{SGD} under particular combinations of domain, regularity, output, and moment assumptions. Conversely, in a convex heavy-tailed model, they show that a broad cone class of first-order schemes with predetermined coefficients necessarily incurs polynomial dependence on $1/\beta$ in high probability. Thus, whether \algname{SGD} ``handles heavy tails'' depends on the criterion, oracle, domain, stepsize class, and probability mode. Chapters~\ref{ch:chapter_2}--\ref{ch:chapter_5} focus on regimes in which robust transformations or richer stochastic information improve the available guarantee or attain a sharper complexity benchmark.}

\paragraph{Nonlinear transformations: robustness, bias, and structural preservation.}
The central issue is not only which nonlinear transformation is used, but also where it is applied. A bounded nonlinear map can limit the influence of extreme oracle realizations, but it generally changes the expectation of the transformed vector. In composite, distributed, and operator problems, applying the transformation to a component that does not vanish at a solution can also shift the fixed point. Robustification must therefore control extreme observations while preserving the equation that characterizes the target solution.

Clipping, normalization, sign mappings, quantization, and related nonlinear gradient maps are instances of a broader stochastic-approximation design principle \parencite{polyak1979adaptive,jakovetic2023nonlinear,pmlr-v258-armacki25a}. In neural-network training, clipping was popularized as a remedy for exploding gradients \parencite{pascanu2013difficulty}; in weak-moment optimization, it is often combined with momentum or normalization \parencite{NEURIPS2021_26901deb,NEURIPS2023_4c454d34,pmlr-v202-sadiev23a}. The statistical tradeoff is bounded influence versus transformation bias: the threshold or nonlinear map controls extreme observations but can distort the mean, and this distortion can persist even for small stepsizes \parencite{koloskova2023revisiting,vatankhah2025differentially}. Recent analyses further clarify when normalization, clipped adaptive methods, or unclipped algorithms can be preferable \parencite{hübler2025gradientclippingnormalizationheavy,pmlr-v267-chezhegov25a,liu2025nonconvex}.

Structural placement is therefore as important as the map itself. Chapter~\ref{ch:chapter_2} transforms stochastic gradients and operators directly. Chapter~\ref{ch:chapter_3} clips shifted innovations or differences because individual components need not vanish at the global solution. Chapters~\ref{ch:chapter_4}--\ref{ch:chapter_5} separately robustify raw-gradient and correction components, including same-sample differences and Hessian--vector corrections. Clipping is one important nonlinear transformation, but it is neither the only one nor universally necessary; its role depends on the oracle assumptions, problem structure, probability mode, and controlled output.

\paragraph{Richer stochastic information and complexity limits.}
Clipping and normalization control how available information influences an update. A different route to faster convergence is to construct a more accurate search vector. Momentum variance reduction (\algname{MVR}) uses correlated same-sample gradient differences to improve the search-vector estimate \parencite{cutkosky2019momentum}. Hessian-corrected momentum (\algname{HM}) instead uses stochastic Hessian--vector products to transport or correct a previous estimate as the iterate moves \parencite{tran2022better,salehkaleybar2022momentum,sadievsecond}. In the heavy-tailed setting, these two information models underlie the lower- and upper-bound comparisons in Chapter~\ref{ch:chapter_5} \parencite{fradin2025tight}.

\outrev{Implicit gradient transport (\algname{IGT}) was introduced by \textcite{arnold2019reducing} as a way to compensate for the staleness of reused gradients without explicitly evaluating Hessians. \textcite{pmlr-v119-cutkosky20b} subsequently incorporated an \algname{IGT}-style extrapolation into normalized \algname{SGD} and analyzed the resulting momentum construction for second-order-smooth nonconvex objectives. Ordinary normalized momentum and \algname{IGT} should therefore not be treated as interchangeable.}

The black-box oracle-complexity viewpoint originates in the work of \textcite{NemiYudin}. Chapters~\ref{ch:chapter_4}--\ref{ch:chapter_5} place such lower bounds at the center of the analysis. Oracle complexity asks how many stochastic first- or second-order observations are necessary, in the worst case over a specified function, oracle, output, and algorithm class, to reach a prescribed stationarity criterion. Deterministic lower bounds and stochastic first- and second-order lower bounds provide the baseline for this comparison \parencite{carmon2020lower,arjevani2020second,arjevani2022lower}. The new lower bounds quantify how the attainable complexity changes with the moment indices of stochastic gradients and Hessians and with assumptions governing same-sample gradient differences. Matching or near-matching algorithms then identify when normalization, variance reduction, and stochastic curvature extract the benefit permitted by the oracle.

\paragraph{Communication constraints.}
Chapters~\ref{ch:chapter_6}--\ref{ch:chapter_7} analyze the uniform finite-average problem
\begin{equation}
    f(x)
    =
    \frac{1}{n}
    \sum_{i=1}^{n} f_i(x),
    \qquad
    f_i(x)
    =
    \Exp_{\xi_i\sim\cD_i}
    \bigl[f_i(x;\xi_i)\bigr].
    \label{ch1:eq:intro_distributed_objective_chapter_1}
\end{equation}
Here $n$ is the number of workers, $f_i$ is the local population objective, and $\cD_i$ is the local data distribution. The distributions $\cD_i$ need not coincide, so the model permits statistical heterogeneity. This is distinct from the stochastic-sampling assumptions used in the analysis: independence or conditional independence of fresh samples is stated separately in each technical chapter.

\outrev{The formal model is a synchronous star topology with a fixed, fully participating worker population. Workers construct local gradient-type information and communicate contractively compressed worker-to-server messages; the server aggregates reconstructed vectors before applying a normalized update. This is also a common server-based abstraction of federated learning, but Chapters~\ref{ch:chapter_6}--\ref{ch:chapter_7} isolate the interaction among uplink compression, error-feedback reconstruction, momentum, and the server step rather than attempting to model every cross-device systems feature.}

\paragraph{Generalized smoothness and controlled steps.}
\outrev{Chapters~\ref{ch:chapter_6} and~\ref{ch:chapter_8} study objectives for which one global quadratic scale does not adequately describe the local curvature. Chapter~\ref{ch:chapter_6} uses the differentiable symmetric $(L_0,L_1)$-smoothness condition of \textcite{chen2023generalized} on the local objectives. Chapter~\ref{ch:chapter_8} uses arbitrary-norm, estimator-specific extensions governing the population gradient, same-sample stochastic-gradient differences, and the Hessian. Normalization and centered-ball LMOs are natural in these settings because the magnitude of a search vector need not itself prescribe a safe displacement. The detailed literature and the distinct chapter-specific assumptions are reviewed in Sections~\ref{ch1:subsec:intro_communication_chapter_1}--\ref{ch1:subsec:intro_lmo_chapter_1}; they are not collapsed into one blanket generalized-smoothness model.}

Normalization separates the orientation of a search vector from the length of the update and can support schedules that avoid selected problem constants, but it does not make an algorithm tuning-free. In Chapter~\ref{ch:chapter_7}, parameter agnosticism means specifically that the prescribed decreasing stepsizes do not require advance knowledge of the smoothness constant or the initial suboptimality gap. The methods retain structural parameters and a user-selected base scale. Chapter~\ref{ch:chapter_8} extends the bounded-step mechanism from Euclidean normalization to arbitrary-norm linear minimization oracles, allowing the displacement geometry to reflect vector or matrix structure.

These considerations lead to the central question of the dissertation:
\begin{quote}
    \emph{How should an algorithm construct an informative search vector, convert it into a stable step, and exploit the available information as efficiently as possible?}
\end{quote}
\begingroup\widowpenalty=10000
A stable step refers here to an update whose magnitude or geometry is explicitly controlled, so that an extreme stochastic estimate, a large gradient, a compression error, or locally growing curvature does not directly produce an uncontrolled displacement. Efficiency is assessed through first- and second-order oracle complexity, confidence dependence, and communication complexity, with lower bounds used to identify the regimes in which the available information is exploited optimally.\par
\endgroup

\stepsixchange{The preceding difficulties do not call for one universal algorithm. They can nevertheless be compared through two design decisions: how available information is combined into a search vector, and how that vector determines a controlled step. The next section develops this common lens before turning to the four research directions.}

\section{An Algorithmic Lens: Information, Search Vectors, and Steps}
\label{ch1:sec:intro_lens_chapter_1}

\stepsixchange{To formalize these two decisions, at iteration} $k$, let $\mathcal{H}_k$ denote the information available before the next update. This history may contain current and past stochastic oracle responses, momentum variables, same-sample gradient differences, Hessian--vector products, communicated messages, and error-feedback memories. We write
\begin{equation}
    m^k
    =
    \mathsf{E}_k(\mathcal{H}_k),
    \qquad
    s^k
    \in
    \mathsf{U}_k(x^k,m^k),
    \qquad
    x^{k+1}=x^k+s^k.
    \label{ch1:eq:intro_estimator_update_chapter_1}
\end{equation}
Here, $m^k$ is the constructed search vector and $s^k$ is the resulting displacement. The map $\mathsf{E}_k$ describes how available information is combined, while $\mathsf{U}_k$ describes how the constructed vector changes the iterate.

The decomposition~\eqref{ch1:eq:intro_estimator_update_chapter_1} is a conceptual comparison device rather than a unique factorization. For example, clipping can be interpreted as part of the estimator, by setting $m^k$ equal to a clipped stochastic estimate, or as part of the step rule applied to an unclipped estimate. The first interpretation is useful when asking which stochastic quantity should be clipped; the second is useful when comparing clipping with normalization and norm-constrained steps. Similarly, the accelerated and \algname{SEG} algorithms in Chapters~\ref{ch:chapter_2}--\ref{ch:chapter_3} use several applications of basic update primitives and should not be forced into a single one-stage iteration.

A broad class of displacement rules can be represented by the following regularized local subproblem with an optional norm constraint:
\begin{equation}
    s^k
    \in
    \operatorname*{argmin}_{\norm{s}\leq r_k}
    \left\{
        \inner{m^k}{s}
        +
        \frac{1}{2}\inner{s}{B^k s}
        +
        \Psi(x^k+s)
    \right\},
    \qquad
    x^{k+1}=x^k+s^k,
    \label{ch1:eq:intro_local_model_chapter_1}
\end{equation}
where $B^k\succeq0$, $\Psi:\R^d\to\R\cup\{+\infty\}$ is a proper closed convex function when a composite term is present, and
\begin{equation}
    r_k\in(0,+\infty],
    \qquad
    \left\{s\in\R^d:\norm{s}\leq+\infty\right\}
    =
    \R^d.
    \label{ch1:eq:intro_infinite_radius_chapter_1}
\end{equation}
Thus, $r_k=+\infty$ makes the feasible region equal to the whole space, whereas a finite value imposes a hard bound on the displacement. The quadratic term provides a metric or soft regularization, and $\Psi$ preserves nonsmooth regularization or constraints.

A related idealized viewpoint is provided by ball-proximal, or Broximal, methods, which minimize the objective itself over a ball centered at the current iterate \parencite{gruntkowska2025broximal}. Their non-Euclidean extension studies this construction in general norm geometries \parencite{gruntkowska2025noneuclidean}. These methods are conceptually close to the optional-radius viewpoint used here, but differ from~\eqref{ch1:eq:intro_local_model_chapter_1}: the present framework minimizes a regularized local model over a displacement set rather than exactly minimizing the original objective over that set.

Related first-order trust-region formulations recover clipped and normalized stochastic-gradient updates under generalized smoothness \parencite{xie2024trust}. Non-Euclidean trust-region formulations similarly recover normalized, sign-type, and orthogonalized matrix directions \parencite{kovalev2025understanding}, while norm-constrained LMOs provide a direct framework for centered-ball displacement rules \parencite{pethick2025training}. Equation~\eqref{ch1:eq:intro_local_model_chapter_1} is used here as a broader \stepsixchange{dissertation-level} comparison device: the optional quadratic term and composite term allow it to represent the different outer update mechanisms appearing across the dissertation.

The matrix $B^k$ should not be confused with the stochastic Hessian information used later in the \stepsixchange{dissertation}. In Chapters~\ref{ch:chapter_4}, \ref{ch:chapter_5}, \ref{ch:chapter_7}, and~\ref{ch:chapter_8}, Hessian--vector products primarily improve the construction of $m^k$ through \algname{HM}. They are not generally inserted as $B^k$ in the outer subproblem~\eqref{ch1:eq:intro_local_model_chapter_1}.

\paragraph{Gradient-type and proximal steps.}
Let $r_k=+\infty$, choose Euclidean geometry, and set $B^k=\gamma_k^{-1}I$ with $\gamma_k>0$. If $\Psi\equiv0$, then completing the square gives
\begin{align}
    s^k
    &=
    \operatorname*{argmin}_{s\in\R^d}
    \left\{
        \inner{m^k}{s}
        +
        \frac{1}{2\gamma_k}\norm{s}^2
    \right\}
    \notag\\
    &=
    \operatorname*{argmin}_{s\in\R^d}
    \frac{1}{2\gamma_k}
    \norm{s+\gamma_km^k}^2
    =
    -\gamma_km^k.
    \label{ch1:eq:intro_gradient_special_case_chapter_1}
\end{align}
Consequently, $x^{k+1}=x^k-\gamma_km^k$. The local model is quadratic, while its minimizer is linear in the constructed search vector.

With a nonzero composite term $\Psi$, the same specialization gives
\begin{equation}
\begin{aligned}
    x^{k+1}
    &=
    \operatorname{prox}_{\gamma_k\Psi}\bigl(x^k-\gamma_km^k\bigr),
    \\
    \operatorname{prox}_{\gamma\Psi}(z)
    &\eqdef
    \operatorname*{argmin}_{x\in\R^d}
    \left\{
        \Psi(x)+\frac{1}{2\gamma}\norm{x-z}^2
    \right\}.
\end{aligned}
    \label{ch1:eq:intro_prox_special_case_chapter_1}
\end{equation}
This is the relevant outer primitive for the composite methods in Chapter~\ref{ch:chapter_3} \parencite{parikh2014proximal,beck2017first}.

\paragraph{Clipping as a bounded quadratic step.}
For $v\in\R^d$ and $\lambda>0$, the notation used throughout the dissertation is
\begin{equation}
    \operatorname{clip}(v,\lambda)
    \eqdef
    \begin{cases}
        \displaystyle
        \min\left\{1,\frac{\lambda}{\norm{v}}\right\}v,
        & v\neq0,\\[2mm]
        0,
        & v=0.
    \end{cases}
    \label{ch1:eq:intro_clip_definition_chapter_1}
\end{equation}
Let $\mathbb{B}_2=\{u\in\R^d:\norm{u}\leq1\}$. Then $\operatorname{clip}(v,\lambda)$ is the Euclidean projection of $v$ onto $\lambda\mathbb{B}_2$. In~\eqref{ch1:eq:intro_local_model_chapter_1}, set $\Psi\equiv0$, $B^k=\gamma_k^{-1}I$, and $r_k=\gamma_k\lambda_k$. Completing the square yields
\begin{align}
    s^k
    &\in
    \operatorname*{argmin}_{\norm{s}\leq\gamma_k\lambda_k}
    \left\{
        \inner{m^k}{s}
        +
        \frac{1}{2\gamma_k}\norm{s}^2
    \right\}
    \notag\\
    &=
    \operatorname{proj}_{\gamma_k\lambda_k\mathbb{B}_2}
    \left(-\gamma_km^k\right)
    =
    -\gamma_k\operatorname{clip}(m^k,\lambda_k).
    \label{ch1:eq:intro_clipping_special_case_chapter_1}
\end{align}
Thus, clipping produces the usual gradient-type displacement when $\norm{m^k}\leq\lambda_k$ and a normalized displacement of length $\gamma_k\lambda_k$ when $\norm{m^k}>\lambda_k$. The threshold $\lambda_k$ bounds the clipped search vector, while $r_k=\gamma_k\lambda_k$ bounds the resulting displacement.

\paragraph{Normalization and centered norm-ball LMOs.}
Set $B^k=0$, $\Psi\equiv0$, and $r_k<+\infty$. Then
\begin{equation}
    s^k
    \in
    \operatorname*{argmin}_{\norm{s}\leq r_k}
    \inner{m^k}{s},
    \qquad
    \min_{\norm{s}\leq r_k}
    \inner{m^k}{s}
    =
    -r_k\normstar{m^k}.
    \label{ch1:eq:intro_lmo_special_case_chapter_1}
\end{equation}
Here $\normstar{\cdot}$ is the dual norm. In Euclidean geometry and for $m^k\neq0$,
\begin{equation}
    s^k
    =
    -r_k\frac{m^k}{\norm{m^k}}.
    \label{ch1:eq:intro_normalized_special_case_chapter_1}
\end{equation}
This is the normalized displacement used in Chapters~\ref{ch:chapter_4}--\ref{ch:chapter_7}. Chapter~\ref{ch:chapter_8} replaces the Euclidean ball by an arbitrary norm ball. The resulting subproblem is a linear minimization oracle over a centered set of displacements. It should be distinguished from the classical \algname{Frank--Wolfe} update, which minimizes a linear model over the original feasible set and then moves toward the selected feasible point \parencite{frank1956algorithm,jaggi2013revisiting}.

If $B^k=0$, $\Psi\equiv0$, and $r_k=+\infty$, the local model is unbounded below whenever $m^k\neq0$. Hence the finite radius is essential for the normalized/LMO specialization, whereas the positive-definite quadratic regularizer makes the whole-space gradient and proximal specializations well posed.

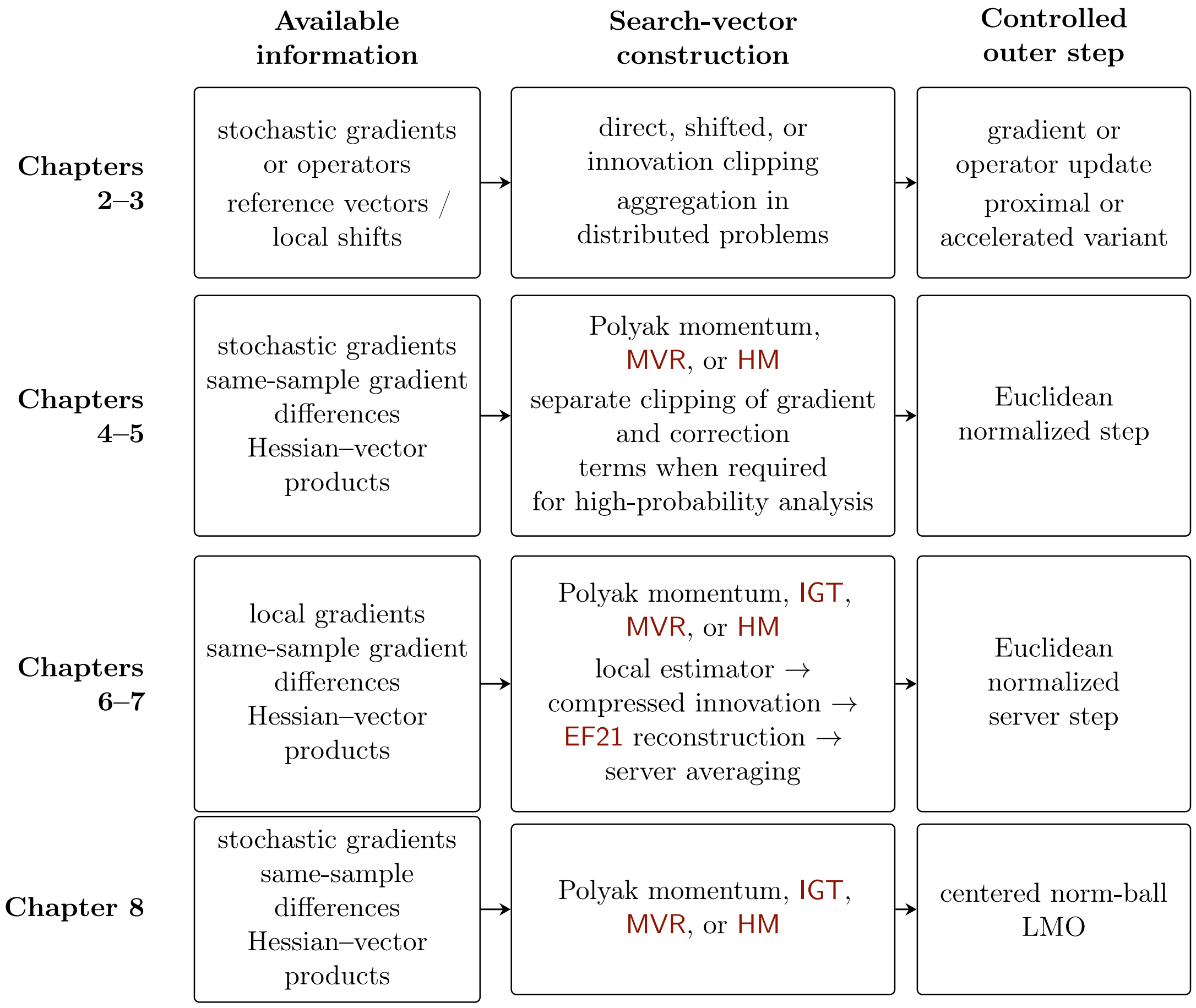
\begin{figure}[H]
    \centering
    {\color{steponeink}
    \stepEightCAlgorithmNames
    \resizebox{\textwidth}{!}{\begin{tikzpicture}[
    x=1pt,
    y=1pt,
    stepone box/.style={
        draw=steponeink,
        text=steponeink,
        rounded corners=2pt,
        line width=0.55pt,
        align=center,
        inner xsep=4pt,
        inner ysep=5pt,
        font=\small
    },
    stepone heading/.style={
        text=steponeink,
        align=center,
        font=\small\bfseries
    },
    stepone lane/.style={
        text=steponeink,
        align=right,
        font=\small\bfseries
    },
    stepone arrow/.style={
        draw=steponeink,
        line width=0.65pt,
        -{Stealth[length=4.5pt,width=5pt]}
    }
]
    \def\laneLabelWidth{56pt}
    \def\infoWidth{106pt}
    \def\constructWidth{145pt}
    \def\stepWidth{101pt}
    \def\columnGap{13pt}

    \node[stepone heading, text width=\infoWidth] (head-info) at (130,0)
        {Available\\information};
    \node[stepone heading, text width=\constructWidth] (head-construct) at (276,0)
        {Search-vector construction};
    \node[stepone heading, text width=\stepWidth] (head-step) at (416,0)
        {Controlled outer step};

    \node[stepone lane, text width=\laneLabelWidth, anchor=east] (lane23) at (57,-58)
        {Chapters\\2--3};
    \node[stepone box, text width=\infoWidth, minimum height=76pt] (info23) at (130,-58)
        {stochastic gradients\\or operators\\[2pt]\figuretablepolishchange{reference vectors /}\\\figuretablepolishchange{local shifts}};
    \node[stepone box, text width=\constructWidth, minimum height=76pt] (construct23) at (276,-58)
        {direct, shifted, or\\innovation clipping\\[2pt]aggregation in\\distributed problems};
    \node[stepone box, text width=\stepWidth, minimum height=76pt] (step23) at (416,-58)
        {gradient or\\operator update\\[2pt]proximal or\\accelerated variant};

    \node[stepone lane, text width=\laneLabelWidth, anchor=east] (lane45) at (57,-151)
        {Chapters\\4--5};
    \node[stepone box, text width=\infoWidth, minimum height=96pt] (info45) at (130,-151)
        {stochastic gradients\\same-sample gradient\\differences\\Hessian--vector products};
    \node[stepone box, text width=\constructWidth, minimum height=96pt] (construct45) at (276,-151)
        {\mbox{Polyak momentum},\\\algname{MVR}, or \algname{HM}\\[2pt]separate clipping of gradient\\and correction terms when required\\for high-probability analysis};
    \node[stepone box, text width=\stepWidth, minimum height=96pt] (step45) at (416,-151)
        {Euclidean\\normalized step};

    \node[stepone lane, text width=\laneLabelWidth, anchor=east] (lane67) at (57,-258)
        {Chapters\\6--7};
    \node[stepone box, text width=\infoWidth, minimum height=102pt] (info67) at (130,-258)
        {local gradients\\same-sample gradient\\differences\\Hessian--vector products};
    \node[stepone box, text width=\constructWidth, minimum height=102pt] (construct67) at (276,-258)
        {\mbox{Polyak momentum}, \algname{IGT},\\\algname{MVR}, or \algname{HM}\\[3pt]local estimator $\to$\\compressed innovation $\to$\\\algname{EF21} reconstruction $\to$\\server averaging};
    \node[stepone box, text width=\stepWidth, minimum height=102pt] (step67) at (416,-258)
        {Euclidean normalized\\server step};

    \node[stepone lane, text width=\laneLabelWidth, anchor=east] (lane8) at (57,-348)
        {Chapter 8};
    \node[stepone box, text width=\infoWidth, minimum height=68pt] (info8) at (130,-348)
        {stochastic gradients\\same-sample differences\\Hessian--vector products};
    \node[stepone box, text width=\constructWidth, minimum height=68pt] (construct8) at (276,-348)
        {\mbox{Polyak momentum}, \algname{IGT},\\\algname{MVR}, or \algname{HM}};
    \node[stepone box, text width=\stepWidth, minimum height=68pt] (step8) at (416,-348)
        {centered norm-ball\\LMO};

    \foreach \row in {23,45,67,8}{
        \draw[stepone arrow] (info\row.east) -- (construct\row.west);
        \draw[stepone arrow] (construct\row.east) -- (step\row.west);
    }
\end{tikzpicture}}
    \captionsetup{labelfont={color=steponeoneink},textfont={color=steponeoneink}}
    \caption[Information--search-vector--step structure across Chapters~2--8]{Information, search-vector construction, and controlled steps across Chapters~2--8. Each row shows the information used by a chapter group, how it is converted into a search vector, and the corresponding outer step. Individual chapters use only the relevant subset of the displayed information.}
    \label{fig:intro-information-search-step}
    \par\smallskip
    {\footnotesize\color{steponeoneink}\textbf{Abbreviations.}
    \algname{MVR}: momentum variance reduction;
    \algname{HM}: Hessian-corrected momentum;
    \algname{IGT}: implicit gradient transport;
    \algname{EF21}: Error Feedback 21;
    LMO: linear minimization oracle.\par}
    }
\end{figure}

\steponechange{\stepsixchange{Figure~\ref{fig:intro-information-search-step} uses the decomposition in~\eqref{ch1:eq:intro_estimator_update_chapter_1} to compare algorithmic choices across the chapters: similar information can support different estimators, and related estimators can be paired with different outer steps.}} \steponeonechange{Moment and smoothness assumptions are not shown because they characterize the analytical setting rather than the information used directly to construct an update.}

\steponechange{Chapters~4--5 and~8 use closely related estimator families but different outer geometries. Chapters~6--7 insert an additional communication layer between the local estimator and the outer step: workers send compressed innovations, \algname{EF21} updates the local reconstructions, and the server averages them before the normalized update.}

\section{Research Directions and Main Contributions}
\label{ch1:sec:intro_directions_chapter_1}

The dissertation develops four research directions. The first establishes high-probability robustness under weak moments. The second asks how richer first- and second-order information changes the attainable complexity and develops matching lower and upper bounds. The third transfers normalized and momentum-based constructions to communication-efficient distributed optimization. The fourth extends improved momentum estimators to arbitrary-norm LMO steps. Each direction is described below through its motivation, the closest prior theory, the remaining gap, and the contribution of the corresponding chapters.

\subsection{High-Probability Optimization under Heavy-Tailed Noise}
\label{ch1:subsec:intro_heavy_tail_chapter_1}

The high-probability literature is best understood along two independent axes: the tail assumptions imposed on the stochastic oracle and the output criterion being controlled. Classical stochastic mirror descent and accelerated stochastic approximation obtain expectation and high-confidence guarantees under light-tail assumptions, typically with domain or gradient conditions that make the relevant martingale increments controllable \parencite{nemirovski2009robust,lan2012optimal,lan2020first}. Under finite central second moments but without exponential tails, robust truncation and clipping recover sub-Gaussian-type confidence in more restricted settings \parencite{nazin2019algorithms,gorbunov2020stochastic}.

\outrev{Two later developments are especially useful for separating the role of concentration from the role of robustification. \textcite{pmlr-v202-liu23aa} derive their martingale sequence directly from the optimization proof and exploit the recursive structure of the resulting process; this ``white-box'' analysis gives high-probability guarantees under sub-Gaussian noise without replacing the initial distance by a global domain diameter and extends to adaptive stepsizes. \textcite{madden2024high} instead establish a self-normalized concentration inequality for norm-sub-Weibull martingale differences and use it to analyze vanilla nonconvex \algname{SGD}. Neither result is a bounded-$p$-moment heavy-tail theorem: sub-Gaussian, sub-exponential, sub-Weibull, finite-variance, and finite-$p$-moment assumptions describe different tail classes \stepfive{and should not be merged under one informal label.}}

\outrev{Adaptive stepsizes create an additional dependence because the stepsize is constructed from past stochastic gradients and is therefore correlated with the update direction. \stepfivemetachange{\textcite{li2020high} and \textcite{kavis2022high} analyze adaptive methods under light-tail assumptions, while the white-box framework of \textcite{pmlr-v202-liu23aa} also accommodates adaptive stepsizes under its sub-Gaussian model. Under bounded lower-order moments, \textcite{pmlr-v267-chezhegov25a} show that clipping improves Adam-Norm and AdaGrad-Norm in the heavy-tailed regimes they consider.}}

The lower-order moment regime requires additional care. \textcite{cutkosky2021high} combine clipping, momentum, and normalization for nonconvex high-probability optimization, but use a bounded noncentral $p$-th moment that is stronger than the $p$-BCM model in~\eqref{ch1:eq:intro_p_bcm_chapter_1}. Chapters~\ref{ch:chapter_2}--\ref{ch:chapter_3} instead work with unbiased stochastic gradients or operators satisfying bounded central $p$-th moments. Subsequent and concurrent work improves horizon dependence, removes some prior parameter requirements, and clarifies the transition from clipping to normalization \parencite{nguyen2023improved,hübler2025gradientclippingnormalizationheavy}. A complementary line analyzes broader bounded nonlinear mappings under symmetry or related distributional structure rather than the same central-moment model \parencite{jakovetic2023nonlinear,pmlr-v258-armacki25a}. \stepfive{These branches should not be conflated: their assumptions, outputs, and attainable guarantees differ.}

Variational inequalities bring an additional algorithmic lineage. Korpelevich's deterministic extragradient method (\algname{EG}) was followed by Mirror--Prox and dual-extrapolation schemes for monotone Lipschitz operators \parencite{dem1972numerical,korpelevich1976extragradient,nemirovski2004prox,nesterov2007dual}. Stochastic extragradient (\algname{SEG}) and stochastic Mirror--Prox extend this viewpoint to stochastic operators \parencite{juditsky2011solving,mishchenko2020revisiting,gorbunov2022stochastic}. In Euclidean projection geometry, \algname{SEG} can be viewed as a specialization of stochastic Mirror--Prox; stochastic Mirror--Prox is the broader non-Euclidean proximal framework. Related analyses establish last-iterate convergence for optimistic-gradient methods and high-probability guarantees for clipped stochastic VI methods under heavy-tailed noise \parencite{gorbunov2022lastiterate,gorbunov2022clipped}. Distributed stochastic saddle-point problems have also been studied through oracle lower bounds and near-optimal centralized and decentralized schemes \parencite{beznosikov2020distributed}. The first part of this dissertation asks how high-probability clipping mechanisms can be combined with these optimization and operator structures under the weaker oracle model~\eqref{ch1:eq:intro_p_bcm_chapter_1}.

Chapter~\ref{ch:chapter_2} studies stochastic minimization and stochastic variational inequalities when neither a globally bounded variance nor a globally bounded objective gradient is assumed. For minimization, it covers smooth nonconvex, Polyak--\L{}ojasiewicz, convex, strongly convex, and quasi-strongly convex regimes. For VIs, it covers Lipschitz and star-cocoercive operators under monotone or quasi-strongly monotone structure. The chapter develops \algname{clipped-SGD}, accelerated clipped methods, \algname{clipped-SGDA}, and \algname{clipped-SEG} \parencite{sadiev2023high}. The output and criterion are regime-specific: the nonconvex result controls a realized trajectory average of squared gradient norms; monotone VI results use restricted gaps evaluated at averaged iterates or averaged operator norms; convex results use averaged iterates; and accelerated or strongly structured regimes use the stated terminal criteria. The chapter's detailed tables state separately which complexities recover known optimal orders up to logarithmic factors and which regimes lack a matching lower bound.

A distinctive feature is the local scope of the assumptions. The proof constructs a high-probability event on which the iterates remain in a prescribed neighborhood or level set. Smoothness and moment conditions are required only on that region rather than globally. This replaces a bounded-domain assumption by a trajectory-control argument, while also explaining why the prescribed stepsizes and clipping levels depend on the selected confidence level.

\outrev{The chapter also contains a complementary negative result for ordinary \algname{SGD}. Even under bounded variance, it constructs a smooth strongly convex problem for which controlling the terminal squared distance with probability at least $1-\beta$ requires inverse-power dependence on $\beta$. This theorem does not prove that every linear or nonadaptive first-order method has the same limitation, but it gives a concrete reason for introducing a nonlinear transformation when logarithmic confidence dependence is part of the target guarantee \parencite{sadiev2023high}.}

Chapter~\ref{ch:chapter_3} begins from a limitation of direct clipping. In a composite problem, the solution generally balances a smooth term and a nonsmooth term; the smooth gradient need not vanish at the solution. In a distributed problem, local gradients may remain nonzero even when their average vanishes. Clipping complete local gradients can therefore shift the deterministic fixed point and prevent exact convergence. Classical stochastic composite optimization already provides proximal and accelerated methods under standard oracle assumptions \parencite{lan2012optimal}; the new question is how to preserve that structure under weak-moment high-probability robustification.

The chapter develops shifted and innovation-clipping methods for composite stochastic minimization, distributed optimization, and corresponding VIs \parencite{gorbunov2023high}. Rather than clipping a complete local gradient indiscriminately, the methods clip a stochastic difference relative to a reference vector or an innovation that vanishes at the relevant target structure. The resulting proximal, accelerated, \algname{SGDA}, and \algname{SEG} methods retain bounded increments for concentration while preserving the deterministic fixed point. The distributed theorems obtain linear speedup in the stated stochastic terms. Accuracy dependence, output rule, and parameter requirements are reported separately by regime in the chapter's theorem tables \stepfive{rather than summarized by one universal optimality claim.}

The relation between Chapters~\ref{ch:chapter_2} and~\ref{ch:chapter_3} is therefore precise. Chapter~\ref{ch:chapter_2} establishes that bounded nonlinear processing can provide broad high-probability control under weak moments; Chapter~\ref{ch:chapter_3} shows that the transformed quantity must be chosen to preserve composite, distributed, or operator structure. \stepfive{These results do not imply that clipping is necessary for every notion of convergence. The positive and negative results of \textcite{fatkhullin2025can} show why that conclusion would be too broad.} The contribution of Chapters~\ref{ch:chapter_2}--\ref{ch:chapter_3} is the breadth of their nonasymptotic high-probability guarantees and the mild confidence dependence achieved across structured minimization and VI settings.

\subsection{Complexity Limits and Richer Stochastic Information}
\label{ch1:subsec:intro_complexity_chapter_1}

Once robust convergence is available, a second question becomes unavoidable: how much improvement is permitted by the information supplied to the algorithm? A convergence theorem for one method cannot answer this question. The comparison requires a lower bound over a specified function, oracle, output, and algorithm class, followed by an upper bound under the same information model.

For smooth nonconvex optimization, deterministic lower bounds characterize the first-order cost of finding approximate stationary points \parencite{carmon2020lower}. Stochastic lower bounds add the effect of oracle noise, and second-order lower bounds ask whether stochastic Hessians or Hessian--vector products improve sample complexity \parencite{arjevani2020second,arjevani2022lower}. Under weak moments, the answer depends on the moment index and on whether the oracle supplies independent gradients, correlated same-sample differences, or stochastic Hessian information.

\outrev{Under unbiased stochastic gradients satisfying $p$-BCM, \textcite{zhang2020adaptive} establish the first-order heavy-tailed benchmark
\[
    \Exp{\norm{\nabla f(x_{\mathrm{out}})}}
    =
    \Omega\!\left(K^{-\frac{p-1}{3p-2}}\right),
\]
or, suppressing the remaining problem parameters,
\[
    K
    =
    \Omega\!\left(\varepsilon^{-\frac{3p-2}{p-1}}\right)
\]
to obtain $\Exp{\norm{\nabla f(x_{\mathrm{out}})}}\leq\varepsilon$. This is an in-expectation lower bound for the non-squared gradient norm over a broad stochastic first-order algorithm class under standard smoothness. \stepfive{It should not be transferred without proof to a high-probability criterion, a squared-gradient criterion, a trajectory-average output, or a different oracle class.} It provides the standard first-order heavy-tailed baseline against which Chapters~\ref{ch:chapter_4}--\ref{ch:chapter_5} measure the value of same-sample and second-order information.}

Chapter~\ref{ch:chapter_4} develops lower and upper bounds for stochastic second-order optimization under bounded central moments for gradient and Hessian noise \parencite{sadiev2025second}. Its lower bound is proved for zero-respecting stochastic higher-order algorithms. In the leading stochastic term, the second-order oracle improves the target-accuracy exponent by one power of $1/\varepsilon$ relative to the corresponding first-order heavy-tailed lower bound. The construction also shows that derivatives above order two do not improve this leading term within the stated oracle class.

The expectation upper bound uses normalized \algname{HM}. Stochastic Hessian--vector products transport the previous gradient estimate; they are not used to solve a Newton subproblem. The theorem controls the average expected gradient norm over the generated trajectory. Its complexity matches the lower-bound dependence on the target accuracy. If, in addition,
\[
    L\leq \sigma_h
    \qquad\text{and}\qquad
    \Delta_0(L+\sigma_h)\geq \sigma\varepsilon,
\]
then its leading stochastic term matches the corresponding lower-bound term up to universal constants. The method does not require Lipschitz continuity of the Hessian.

Chapter~\ref{ch:chapter_4} also introduces clipping of both the stochastic gradient and the Hessian--vector correction. The resulting batch-free method gives a high-probability bound for the trajectory-average gradient norm with polylogarithmic dependence on $1/\beta$. Its dependence on $\varepsilon$ is near-optimal, but the dependence on $\Delta_0$, $L$, $\sigma$, and $\sigma_h$ does not match the lower bound. The method also requires several tuning parameters and is analyzed for a fixed confidence level on which those parameters depend. The chapter states these gaps explicitly.

\outrev{Chapter~\ref{ch:chapter_5} broadens the information model by distinguishing the noise of one stochastic gradient from the regularity of a same-sample stochastic-gradient difference. For $q\in[1,2]$, $q$-weak average smoothness ($q$-WAS) requires a finite constant $\bar L\geq0$ such that
\[
    \Exp_{\xi\sim\cD}
    \left[\norm{\nabla f(x;\xi)-\nabla f(y;\xi)}^q\right]
    \leq
    \bar L^q\norm{x-y}^q
    \qquad\text{for all }x,y\in\R^d.
\]
By contrast, $(q,\delta)$-similarity controls the centered same-sample difference:
\[
    \Exp_{\xi\sim\cD}
    \left[
        \norm{
            [\nabla f(x;\xi)-\nabla f(y;\xi)]
            -[\nabla f(x)-\nabla f(y)]
        }^q
    \right]
    \leq
    \delta^q\norm{x-y}^q.
\]
The $p$-BCM assumption controls noise at a fixed point; $q$-WAS and $(q,\delta)$-similarity control how two gradients evaluated with the same sample change as the iterate moves. For $q=2$, these conditions recover familiar mean-squared-smoothness and similarity structures used in stochastic variance reduction \parencite{arjevani2020second,fradin2025tight}.}

Chapter~\ref{ch:chapter_5} permits different moment indices for stochastic gradients, same-sample differences, and Hessian information, and studies both expectation and high-probability guarantees \parencite{fradin2025tight}.

For known moment indices, the expectation results in Theorems~\ref{ch5:thm:nsgd-mvr-convergence-analysis_chapter_5}--\ref{ch5:thm:nsgd-mvr-convergence-analysis-2_chapter_5} sample $x_{\mathrm{out}}$ uniformly from the generated iterates and control $\Exp{\norm{\nabla f(x_{\mathrm{out}})}}$. Under $q$-weak average smoothness, \algname{NSGD-MVR} matches the corresponding first-order lower bound in expectation. Under $(q,\delta)$-similarity, \algname{NSGD-MVR} attains the variance-reduced branch of the upper bound. The complete minimum in the lower bound of Theorem~\ref{ch5:thm:lower-bound-sofo-mean-squared-smoothness-2_chapter_5} is matched, up to constant factors, by taking the better guarantee among the appropriate normalized methods: \algname{NSGD-Mom} (or normalized minibatch \algname{SGD} for the corresponding first branch) and \algname{NSGD-MVR} for the variance-reduced branch. These are a collection of methods whose guarantees recover different branches of the minimum, not a new combined algorithm.

When the moment indices are unknown, Theorems~\ref{ch5:thm:nsgd-mvr-convergence-analysis-unknown-p-q_chapter_5}--\ref{ch5:thm:nsgd-mvr-convergence-analysis-unknown-p-q-2_chapter_5} instead control the arithmetic mean of expected gradient norms,
\[
    \frac{1}{K}\sum_{k=0}^{K-1}\Exp{\norm{\nabla f(x^k)}}.
\]
The \algname{HM} methods in the chapter give near-matching upper bounds for the corresponding stochastic second-order regimes.

For high-probability convergence, Chapter~\ref{ch:chapter_5} develops double-clipped methods that bound the gradient component and the gradient-difference or Hessian-correction component separately. The principal result, Theorem~\ref{ch5:thm:clipped-nsgd-mvr-convergence-analysis_chapter_5}, controls the realized trajectory average
\[
    \frac{1}{K}\sum_{k=0}^{K-1}\norm{\nabla f(x^k)}
\]
with probability at least $1-\beta$. When the gradient and difference moment indices coincide, its dependence on $\varepsilon$ matches the first-order lower bound under the stated $q$-WAS model. When $p<q$, the expectation result can be faster than the available high-probability result; whether this gap is fundamental remains open. The chapter therefore provides a nearly complete, but not universal, complexity characterization.

\outrev{\stepfivemetachange{The lower bounds of \textcite{fatkhullin2025can} are algorithm-specific: they characterize vanilla \algname{SGD} under their stated regularity, output, and stepsize classes. In contrast, \textcite{zhang2020adaptive} establish the broader first-order $p$-BCM oracle-complexity benchmark used in the comparisons of Chapters~\ref{ch:chapter_4}--\ref{ch:chapter_5}.}}

Chapters~\ref{ch:chapter_4}--\ref{ch:chapter_5} serve a different role from the local-model identities of Section~\ref{ch1:sec:intro_lens_chapter_1}. Their lower bounds do not follow from the estimator--step decomposition. \stepsixchange{The next research direction concerns a different restriction on information: what can be communicated to construct the search vector.}

\subsection{Communication Compression, Generalized Smoothness, and Error-Feedback Reconstruction}
\label{ch1:subsec:intro_communication_chapter_1}

\outrev{The third research direction concerns information that is expensive to transmit rather than statistically heavy-tailed. Section~\ref{ch1:sec:intro_motivation_chapter_1} already fixes the model: a synchronous star topology, a fixed fully participating population, and contractive worker-to-server compression. The scientific issue is how the server can construct a reliable direction from messages that intentionally discard information.}

\outrev{Compression methods differ in how their approximation error is controlled. Unbiased compressors satisfy $\Exp{\mathcal C(v)}=v$. Contractive compressors control the mean-square distortion and may be biased; a representative unconditional condition is
\[
    \Exp{\norm{\mathcal C(v)-v}^2}
    \leq
    (1-\alpha)\norm{v}^2,
    \qquad \alpha\in(0,1].
\]
The technical chapters state the required history-conditional version explicitly. Quantization and sparsification are representative mechanisms, but unbiasedness is an expectation identity rather than a synonym for random coordinate selection, and biased compression is broader than Top-$k$ \parencite{alistarh2017qsgd,stich2018sparsified,beznosikov2023biased}.}

\outrev{A direct compressed-gradient update can fail under contractive compression and heterogeneous local objectives. Classical error feedback stores the residual between an intended update and the message actually transmitted, then reinjects this residual in later rounds \parencite{seide20141,karimireddy2019error}. \algname{EF21} uses a different state: each worker maintains an estimator of its local gradient and communicates a compressed innovation relative to that estimator \parencite{richtarik2021ef21}. This distinction matters because the two mechanisms track different fixed-point quantities.}

\outrev{\stepsixchange{Two complementary approaches broaden this picture: changing the compressor and generalizing the reconstruction rule.}
\stepfivemetachange{\textcite{induced} convert a contractive compressor into an
induced unbiased compressor, so that methods designed for unbiased
compression apply without any residual-feedback construction.
\textcite{fatkhullin2025ef21} subsequently extend \algname{EF21} to stochastic
gradients, variance reduction, partial participation, bidirectional
compression, momentum, and proximal optimization. \textcite{richtarik20223pc}
introduce the \algname{3PC} framework, which generalizes the underlying
communication rule and recovers \algname{EF21} as a special case, while
\algname{EF-BV} treats biased contractive and unbiased random compressors
within a common two-parameter bias--variance framework, recovering both
\algname{EF21} and the compressed-gradient-difference method
\parencite{condat2022ef,mishchenko2019distributed}.}}

\outrev{Related stateful approximation mechanisms reach beyond first-order
gradient tracking. \stepfivemetachange{\algname{FedNL} learns compressed local
Hessian information without relying on classical error feedback
\parencite{safaryan2022fednl}, and \algname{RandProx} replaces a dual proximal
update by a randomized one, using variance reduction to retain exact
convergence in a primal--dual framework \parencite{condat2023randprox}.}}

\outrev{For stochastic nonconvex optimization, \algname{EF21-SGDM} and \algname{EF21-MVR} combine \algname{EF21} reconstruction with Polyak momentum and momentum variance reduction. The former is the direct algorithmic predecessor normalized in Chapter~\ref{ch:chapter_6}, and both are direct predecessors of the estimator hierarchy in Chapter~\ref{ch:chapter_7} \parencite{fatkhullin2023momentum}. \algname{EControl} provides a neighboring two-parameter error-control mechanism for stochastic heterogeneous optimization with contractive compressors \parencite{gao2024econtrol}. \stepsixchange{To study reconstruction together with normalized steps, Chapter~\ref{ch:chapter_6} also accounts for local curvature that grows with the gradient.}}

\outrev{\stepfivemetachange{The generalized-smoothness model used in Chapter~\ref{ch:chapter_6} has developed through several complementary lines of work. \textcite{zhang2020gradient} introduced the original Hessian-growth $(L_0,L_1)$ condition, while \textcite{chen2023generalized} developed the differentiable symmetric formulation adopted in Chapter~\ref{ch:chapter_6}. \textcite{crawshaw2022robustness} studied coordinate-wise generalized smoothness for \algname{SignSGD}- and Adam-type methods, and \textcite{hubler2024parameter} developed parameter-agnostic normalized momentum together with lower bounds under relaxed smoothness. Concurrent work of \textcite{gorbunov2024methods} studied clipping, Polyak stepsizes, stochastic extensions, acceleration, and adaptive methods, while \textcite{vankov2024optimizing} developed structural inequalities, normalized methods, nonconvex guarantees, and a distinct accelerated construction.}}

Chapter~\ref{ch:chapter_6} assumes that each local objective satisfies the symmetric $(L_0,L_1)$ condition
\begin{equation}
    \norm{\nabla f_i(x)-\nabla f_i(y)}
    \leq
    \left(
        L_0
        +L_1\sup_{\theta\in[0,1]}
        \norm{\nabla f_i(\theta x+(1-\theta)y)}
    \right)
    \norm{x-y}.
    \label{ch1:eq:intro_generalized_smoothness_chapter_1}
\end{equation}
The contribution is not another centralized method under this condition. It shows that \algname{EF21} reconstruction from compressed innovations remains compatible with a normalized server step when the local curvature scale may grow with the gradient norm \parencite{khirirat2024errorfeedbackl0l1smoothnessnormalization}. For deterministic gradients, the theorem controls an analytical best iterate and gives the $\cO(K^{-1/2})$ order known for \algname{EF21} under ordinary smoothness. The stochastic momentum variant controls the same type of criterion at order $\cO(K^{-1/4})$, uses one fresh stochastic gradient per worker and iteration, and does not impose a uniform gradient-dissimilarity bound. \stepfive{These are upper-bound comparisons, not a minimax lower-bound result for the complete compressed generalized-smoothness model.}

Chapter~\ref{ch:chapter_7} keeps the \algname{EF21} reconstruction and normalized server step conceptually fixed while changing the local estimator \parencite{sadiev2025improved}. Its five variants use Polyak momentum, \algname{IGT}, \algname{MVR}, and the \algname{HM} variants \algname{EF21-RHM} and \algname{EF21-HM}. The decreasing stepsizes do not require the smoothness constant $L$ or the initial gap $f(x^0)-f^{\inf}$, although the algorithms retain structural parameters and a user-selected base scale. Their theorems control the expected gradient norm of a stepsize-weighted random iterate rather than the terminal iterate. Up to logarithmic factors, the guarantees recover the convergence orders of corresponding problem-dependent \algname{EF21} analyses and centralized normalized methods using analogous estimators: $K^{-1/4}$ for Polyak momentum, $K^{-2/7}$ for \algname{IGT}, and $K^{-1/3}$ for \algname{MVR} and the \algname{HM} variants. \stepfive{This is a comparison with established upper-bound orders, not a separate minimax lower-bound claim for the complete compressed distributed model.}

The formal communication scope remains limited to a central server, a fixed worker population, and the stated contractive uplink compressor model. Partial participation, multiple local epochs, downlink compression, privacy, asynchrony, client churn, and peer-to-peer network topologies are outside the simultaneous scope of Chapters~\ref{ch:chapter_6}--\ref{ch:chapter_7}.

\subsection{A Geometry-Aware Extension: LMO-Based Momentum Methods}
\label{ch1:subsec:intro_lmo_chapter_1}

The final research direction extends normalized steps beyond Euclidean geometry. Classical \algname{Frank--Wolfe} methods minimize a linear model over the original feasible set \parencite{frank1956algorithm,jaggi2013revisiting}. Recent deep-learning optimizers use a related primitive in a different role: an LMO is applied to a centered norm ball of displacements. The selected norm then determines the geometry of the step.

Different norm balls recover familiar updates. The Euclidean ball yields normalization, the $\ell_\infty$ ball yields a sign-type direction, and the spectral-norm ball yields the polar factor of a matrix search direction. The latter is the exact-oracle direction associated with idealized \algname{Muon} updates; practical \algname{Muon} implementations approximate it by a small number of Newton--Schulz iterations \parencite{pmlr-v80-bernstein18a,sun2023momentum,jordan2024muon}. This recent use in matrix-aware language-model optimizers has renewed interest in centered norm-ball LMOs.

\textcite{pethick2025training} formulate a broad LMO framework for deep-learning optimization, while \textcite{kovalev2025understanding} relates normalized, sign, and orthogonalized updates to non-Euclidean trust-region subproblems. Recent native-geometry heavy-tail analysis further shows that smoothness, stochastic error, and stationarity should ideally be formulated in the selected norm and its dual, since a preliminary reduction to Euclidean geometry can introduce dimension-dependent norm-equivalence losses \parencite{hubler2026free}.

\begingroup
\color{steponeoneink}
\ifnum\StepOneOneReviewMode=1 \colorlet{algnamecolor}{steponeoneink}\fi
\stepsixchange{Chapter~\ref{ch:chapter_8} asks which improvements in the search-vector estimator transfer from Euclidean normalization to arbitrary-norm LMOs. This transfer is studied under three estimator-specific generalized-smoothness conditions in the selected norm.} They control, respectively, changes of the population gradient, same-sample stochastic-gradient differences in mean square, and changes of the Hessian. We refer to them as generalized smoothness, mean-squared generalized smoothness, and second-order generalized smoothness. For $x,y\in\R^d$ and $\theta\in[0,1]$, write
$w_\theta(x,y)\eqdef\theta x+(1-\theta)y$.

For constants $L_0,L_1\geq0$, the generalized-smoothness condition is
\begin{equation}
\begin{aligned}
\normstar{\nabla f(x)-\nabla f(y)}
&\leq
\left(
L_0+L_1\sup_{\theta\in[0,1]}
\normstar{\nabla f\!\left(w_\theta(x,y)\right)}
\right)
\norm{x-y}.
\end{aligned}
\label{eq:intro-ch8-generalized-smoothness}
\end{equation}
This is the arbitrary-norm symmetric generalized-smoothness condition used by all Chapter~\ref{ch:chapter_8} variants.

For distinct constants $\cL_0,\cL_1\geq0$, mean-squared generalized smoothness is
\begin{equation}
\begin{aligned}
\Exp_{\xi\sim\cD}{
\normstar{\nabla f(x;\xi)-\nabla f(y;\xi)}^2
}
&\leq
\Exp_{\xi\sim\cD}{
\left(
\cL_0+\cL_1\sup_{\theta\in[0,1]}
\normstar{\nabla f\!\left(w_\theta(x,y);\xi\right)}
\right)^2
}
\\[-0.25ex]
&\quad{}\times\norm{x-y}^2.
\end{aligned}
\label{eq:intro-ch8-mean-squared-generalized-smoothness}
\end{equation}
The same sample $\xi$ is used at $x$ and $y$, so the condition controls same-sample stochastic-gradient differences rather than only the noise of a stochastic gradient at one point. This is an arbitrary-norm analogue of the expected symmetric generalized-smoothness condition of \textcite{chen2023generalized}, specialized to $\alpha=1$.

For constants $M_0,M_1\geq0$, second-order generalized smoothness is
\begin{equation}
\begin{aligned}
\normstar{\nabla^2 f(x)-\nabla^2 f(y)}
&\leq
\left(
M_0+M_1\sup_{\theta\in[0,1]}
\normstar{\nabla f\!\left(w_\theta(x,y)\right)}
\right)
\norm{x-y}.
\end{aligned}
\label{eq:intro-ch8-second-order-generalized-smoothness}
\end{equation}
For a Hessian, the displayed norm is the induced operator norm associated with the selected primal norm and its dual, consistently with Chapter~\ref{ch:chapter_8}. We call this condition second-order generalized smoothness because it controls changes of the Hessian; it is distinct from a moment bound on stochastic Hessian noise.

\begingroup\clubpenalty=10000
The assumptions are estimator-dependent. Polyak momentum uses generalized smoothness. \algname{MVR} additionally uses mean-squared generalized smoothness. \algname{IGT} and the direct-Hessian \algname{HM} construction additionally use second-order generalized smoothness. The randomized-interpolation \algname{HM} construction attains its $\cO(K^{-1/3})$ order without the second-order generalized-smoothness assumption. Complete theorem assumptions and convergence results are given in Chapter~\ref{ch:chapter_8} and summarized in Table~\ref{ch8:tab:TheoreticalComparisons_chapter_8}.\par
\endgroup
\endgroup

\stepsixchange{Under these estimator-specific assumptions, Chapter~\ref{ch:chapter_8} establishes convergence for the corresponding LMO-based momentum methods} \parencite{khirirat2025better}. The theorems control
\[
    \min_{0\leq k\leq K}\Exp{\normstar{\nabla f(x^k)}},
\]
where $\normstar{\cdot}$ is the dual norm. The prior Polyak-momentum baseline has order $\cO(K^{-1/4})$; \algname{IGT} attains $\cO(K^{-2/7})$; and \algname{MVR} and the \algname{HM} variants attain $\cO(K^{-1/3})$. These rates match the corresponding Euclidean estimator orders under standard smoothness. \stepfive{The chapter does not establish a matching lower bound for its full arbitrary-norm relaxed-smoothness class, so the $\cO(K^{-1/3})$ result should not be described as minimax optimal for that complete setting.}

\stepsixchange{The connection to Chapters~\ref{ch:chapter_4}, \ref{ch:chapter_5}, and~\ref{ch:chapter_7} is the shared hierarchy of \algname{IGT}, \algname{MVR}, and \algname{HM} estimators, now paired with a common LMO-based outer step.} The experiments include multilayer perceptrons and NanoGPT training on TinyShakespeare and WikiText-2, providing an initial empirical comparison among momentum constructions and LMO geometries.

Finite-radius LMOs bound the displacement independently of the search-vector magnitude and are invariant to positive rescaling of that vector. These properties make them compatible with weak-moment robustness and parameter-agnostic design, but do not provide those guarantees automatically. \stepEightBchange{Chapter~\ref{ch:chapter_8} is an in-expectation exact-LMO theory; high-probability heavy-tailed analysis and parameter-agnostic schedules remain outside its present scope. Practical inexact realizations of the LMO, including Newton--Schulz approximations in matrix settings, are likewise not analyzed in this dissertation.}

\section{Thesis Map, Publication Scope, and Organization}
\label{ch1:sec:intro_map_chapter_1}

The four research directions form a progression rather than a collection of unrelated algorithms. Chapters~\ref{ch:chapter_2}--\ref{ch:chapter_3} develop high-probability robustification under weak moments and identify the correct object to clip. Chapters~\ref{ch:chapter_4}--\ref{ch:chapter_5} move from robust upper bounds to a minimax account of first- and second-order information. Chapters~\ref{ch:chapter_6}--\ref{ch:chapter_7} transfer normalized search-vector constructions to communication compression and then remove important parameter dependencies while comparing richer estimators. Chapter~\ref{ch:chapter_8} changes the geometry of the outer step and transports the estimator hierarchy to arbitrary-norm LMOs.

\stepsixchange{Where Figure~\ref{fig:intro-information-search-step} compares algorithmic structure, Table~\ref{ch1:tab:intro_chapter_map_chapter_1} records the scope and theorem-covered outputs of each chapter.} Exact rates, oracle models, and comparison scopes are stated in Sections~1.3.1--1.3.4 and in the theorem tables of the individual chapters, because guarantees formulated for different outputs, oracle classes, or algorithm classes are not interchangeable.

\begin{table}[htbp]
\footnotesize
\let\algname\stepEightCTableAlgname
\setlength{\tabcolsep}{2pt}
\renewcommand{\arraystretch}{1.12}
\caption[Scope and main contributions of Chapters~2--8]
        {Scope, theorem outputs, and main contributions of Chapters~2--8.
        \figuretablepolishchange{Mode records whether the principal guarantees
        are in expectation or high probability; LB, CC, and GS mark lower bounds,
        communication compression, and generalized smoothness; Geometry and
        Output describe the outer step and the theorem-covered criterion.}}
\label{ch1:tab:intro_chapter_map_chapter_1}
\begin{tabularx}{\textwidth}{@{}c C{0.90cm} C{0.50cm} C{0.50cm} C{0.50cm} R{1.50cm} R{1.55cm} X@{}}
\toprule
Ch. & Mode & LB & CC & GS & Geometry & Output & Setting and main contribution\\
\midrule
2 & HP & \xmark & \xmark & \xmark & Euclidean & Average or terminal, by regime
  & Clipped \algnamesmall{SGD}, \algnamesmall{SSTM}, \algnamesmall{SGDA}, and \algnamesmall{SEG} methods under bounded lower-order
    moments. Smoothness and moment conditions are imposed only on a neighborhood
    in which the iterates are shown to remain, rather than globally.\\
\addlinespace[3pt]
3 & HP & \xmark & \xmark & \xmark & Euclidean & Average or terminal, by regime
  & Clipping a full local gradient can move the deterministic fixed point in
    composite and distributed problems. Clipping a shift or an innovation instead
    leaves that fixed point intact, and the distributed theorems obtain linear
    speedup in the stochastic terms.\\
\addlinespace[3pt]
4 & Exp.\ \& HP & \cmark & \xmark & \xmark & Euclidean & Trajectory average
  & \figuretablepolishchange{Second-order information improves the leading
    heavy-tailed lower-bound dependence by one power of $1/\varepsilon$, and
    normalized \algname{HM} attains this improvement. Clipping the gradient and
    Hessian--vector correction gives a batch-free high-probability method.}\\
\addlinespace[3pt]
5 & Exp.\ \& HP & \cmark & \xmark & \xmark & Euclidean & Random iterate or trajectory average
  & \figuretablepolishchange{New lower bounds quantify the value of same-sample
    gradient differences and Hessian information. \algname{NSGD-MVR} is optimal
    under $q$-WAS in the stated regime, while double clipping provides
    high-probability guarantees.}\\
\addlinespace[8pt]
6 & Exp. & \xmark & \cmark & \cmark & Euclidean & Best iterate
  & \figuretablepolishchange{Combines \algname{EF21} reconstruction with normalized
    server steps under generalized smoothness, without a uniform
    gradient-dissimilarity bound.}\\
\addlinespace[3pt]
7 & Exp. & \xmark & \cmark & \xmark & Euclidean & Weighted random iterate
  & \figuretablepolishchange{Combines normalized \algname{EF21} with five momentum
    estimators and parameter-agnostic decreasing stepsizes that avoid the
    smoothness constant and initial gap.}\\
\addlinespace[8pt]
8 & Exp. & \xmark & \xmark & \cmark & \figuretablepolishchange{Arbitrary} & Best iterate in the dual norm
  & \figuretablepolishchange{Extends \algname{IGT}, \algname{MVR}, and \algname{HM}
    to arbitrary-norm LMO steps under generalized smoothness, preserving the
    corresponding Euclidean convergence orders.}\\
\bottomrule
\end{tabularx}

\vspace{5pt}
\noindent\begin{minipage}{\textwidth}
\footnotesize
\textbf{Note.} \figuretablepolishchange{Exp. denotes an in-expectation guarantee
and HP a guarantee holding with probability at least $1-\beta$. A \cmark{} means
that the feature is part of the chapter's main theory, while an \xmark{} means
that it is not. ``Arbitrary'' includes Euclidean geometry as a special case.
Chapter~2 additionally contains an algorithm-specific negative result for
\algname{SGD}, which is not counted as an oracle-complexity LB. Chapter~6 uses
the symmetric $(L_0,L_1)$ condition. For Chapter~8, GS includes the
estimator-specific mean-squared and second-order generalized-smoothness
conditions described in Section~\ref{ch1:subsec:intro_lmo_chapter_1}.}
\end{minipage}
\end{table}

\FloatBarrier

\paragraph{Publications underlying the dissertation.}
Chapters~\ref{ch:chapter_2}--\ref{ch:chapter_8} are based, respectively, on the following seven coauthored works: high-probability stochastic optimization and VIs under unbounded variance \parencite{sadiev2023high}; high-probability composite and distributed optimization and VIs \parencite{gorbunov2023high}; stochastic second-order optimization under heavy-tailed noise \parencite{sadiev2025second}; tight lower bounds and optimal heavy-tailed nonconvex algorithms \parencite{fradin2025tight}; normalized error feedback under generalized smoothness \parencite{khirirat2024errorfeedbackl0l1smoothnessnormalization}; parameter-agnostic error feedback with richer momentum \parencite{sadiev2025improved}; and improved arbitrary-norm LMO momentum methods with second-order information \parencite{khirirat2025better}. Each technical chapter preserves the scientific content and detailed theorem statements of its source paper while using the \stepsixchange{dissertation-wide} notation conventions, with supplementary proofs collected in the corresponding appendices.

\paragraph{Research trajectory and publications outside the dissertation.}
The seven works listed above form the technical core of this dissertation, but they emerged from a broader research trajectory that began before the doctoral studies and continued through several parallel projects during the PhD. The grouping below follows the period in which the research was initiated or substantially developed, rather than the date on which its final archival version appeared. This distinction is relevant because several projects begun during the BS or MS studies were published only at a later stage.

\emph{Min--max, saddle-point, and derivative-free optimization.}
The earliest stage of this research, developed during the BS and the beginning of the MS studies, concerned optimization with restricted or inexact oracle information. This line included gradient-free and zeroth-order methods for stochastic saddle-point problems, mixed-oracle algorithms for min--min and min--max problems, and structured nonconvex--concave saddle-point formulations \parencite{beznosikov2020gradient,sadiev2021zeroth,gladin2021solving,alkousa2022approach}. These works introduced themes that later reappeared in the dissertation: operator and saddle-point formulations, restrictions on the information available to an algorithm, and the question of which convergence guarantees remain possible under partial or inexact local information.

\emph{Decentralized, personalized, and communication-efficient optimization.}
The next stage shifted the emphasis from restrictions on a single oracle to restrictions created by a network of communicating agents. The corresponding works studied decentralized personalized federated optimization, decentralized stochastic variational inequalities, communication acceleration through local computation and inexact proximal operations, adaptive compression, and random reshuffling with compressed gradient differences \parencite{borodich2021decentralized,sadiev2022decentralized,kovalev2022optimal,sadiev2022communication,makarenko2023adaptive,sadiev2022federatedoptimizationalgorithmsrandom}. This line established two principles that remain visible in Chapters~\ref{ch:chapter_3}, \ref{ch:chapter_6}, and~\ref{ch:chapter_7}: communicated information should preserve the structure needed by the global method, and improved upper bounds should be compared with the relevant lower complexity limits whenever possible.

\emph{Stochastic estimators, recursion, and update geometry.}
A parallel line of work investigated how the construction of a stochastic direction interacts with the update applied to that direction. Preconditioned stochastic-gradient methods and adaptive implicit recursive estimators examined the roles of preconditioning, adaptivity, implicitness, and variance reduction \parencite{sadiev2024stochastic_gradient,shi2023ai}. These projects help motivate the estimator--step decomposition used in Section~\ref{ch1:sec:intro_lens_chapter_1} and the later emphasis on momentum, \algname{HM}, normalization, and geometry-aware updates.

The doctoral research brought these earlier themes together around the question of how optimization algorithms should construct useful updates from stochastic, partial, or communication-constrained information. \stepsixchange{The technical chapters develop these themes through the four research directions above, while the parallel projects below pursue adjacent problem classes.}

\emph{Parallel research during the PhD.}
Several additional doctoral projects are not included as technical chapters because they address adjacent but distinct problem classes. One line develops proximal-point, primal--dual splitting, and monotone-inclusion methods, including accelerated splitting schemes and momentum variance reduction for implicit updates \parencite{condat2026nesterov,sadiev2026unified,karagulyan2024spam,richtarik2024unified,sadiev2024stochastic_proximal}. A second line studies optimization under structural or system-level restrictions, including asynchronous LMO updates, randomized low-rank adaptation, and differentially private block-coordinate methods \parencite{sadiev2026ringmaster,sokolov2025bernoulli,maranjyan2024differentially}. These works broaden the overall research program toward nonsmooth structure, delayed information, privacy, low-rank models, and partial-coordinate access, but they are omitted as full chapters in order to preserve a focused dissertation narrative and a manageable manuscript length.

\paragraph{Organization of the dissertation.}
Part~I contains Chapters~\ref{ch:chapter_2}--\ref{ch:chapter_3} and develops high-probability methods for stochastic minimization and variational inequalities under heavy-tailed noise. Part~II contains Chapters~\ref{ch:chapter_4}--\ref{ch:chapter_5} and studies nonconvex first- and second-order oracle complexity, lower bounds, normalized variance reduction, and clipping of gradient and Hessian information. Part~III contains Chapters~\ref{ch:chapter_6}--\ref{ch:chapter_7} and develops normalized communication-efficient error-feedback methods under generalized smoothness, including parameter-agnostic variants with five momentum constructions. Part~IV contains Chapter~\ref{ch:chapter_8} and extends advanced momentum estimators to arbitrary-norm LMO-based optimization. Chapter~\ref{ch:chapter_9} concludes the dissertation and discusses open directions suggested by the results and limitations identified throughout the technical chapters.

\newcommand{\thesischtwotheHalgorithm}{\arabic{algorithm}}
\definecolor{bgcolor}{rgb}{0.8,1,1}
\definecolor{bgcolor2}{rgb}{0.8,1,0.8}
\definecolor{niceblue}{rgb}{0.0,0.19,0.56}
\definecolor{PineGreen}{RGB}{0,110,51}
\definecolor{BrickRed}{RGB}{143,20,2}
\newcommand{\thesischtwocmark}{{\color{PineGreen}\ding{51}}}%
\newcommand{\thesischtwoxmark}{{\color{BrickRed}\ding{55}}}%
\newcommand{\thesischtwoExp}{\Exp}
\newcommand{\thesischtwoProb}{\Prob}
\newcommand{\thesischtwooldstuff}[1]{ {\small \color{blue} #1}}
\newcommand{\thesischtworedstuff}[1]{ {\small \color{red} #1}}
\newcolumntype{Y}{>{\centering\arraybackslash}X}
\newcommand\thesischtwofigscale{0.25}
\newcommand\thesischtwoevofigscale{0.21}
\newcommand{\thesischtwosqueeze}{} %
\newcommand{\thesischtwoalgnamex}[1]{\algname{#1}}
\newcommand\thesischtwotagthis{\addtocounter{equation}{1}\tag{\theequation}}
\newcommand{\thesischtwoargmin}{\mathop{\arg\!\min}}
\newcommand{\thesischtwoSp}{\mathrm{Sp}}
\newcommand{\thesischtwoId}{\mathrm{Id}}
\newcommand{\thesischtwoTr}{\mathrm{Tr}}
\newcommand{\thesischtwocA}{{\cal A}}
\newcommand{\thesischtwocB}{{\cal B}}
\newcommand{\thesischtwocC}{{\cal C}}
\newcommand{\thesischtwocD}{{\cal D}}
\newcommand{\thesischtwocE}{{\cal E}}
\newcommand{\thesischtwocF}{{\cal F}}
\newcommand{\thesischtwocG}{{\cal G}}
\newcommand{\thesischtwocH}{{\cal H}}
\newcommand{\thesischtwocJ}{{\cal J}}
\newcommand{\thesischtwocK}{{\cal K}}
\newcommand{\thesischtwocL}{{\cal L}}
\newcommand{\thesischtwocM}{{\cal M}}
\newcommand{\thesischtwocN}{{\cal N}}
\newcommand{\thesischtwocP}{{\cal P}}
\newcommand{\thesischtwocQ}{{\cal Q}}
\newcommand{\thesischtwocR}{{\cal R}}
\newcommand{\thesischtwocS}{{\cal S}}
\newcommand{\thesischtwocT}{{\cal T}}
\newcommand{\thesischtwocU}{{\cal U}}
\newcommand{\thesischtwocV}{{\cal V}}
\newcommand{\thesischtwocX}{{\cal X}}
\newcommand{\thesischtwocY}{{\cal Y}}
\newcommand{\thesischtwocW}{{\cal W}}
\newcommand{\thesischtwocZ}{{\cal Z}}
\newcommand{\thesischtwoVar}{\mathrm{Var}}
\newcommand{\thesischtwomA}{{\bf A}}
\newcommand{\thesischtwomB}{{\bf B}}
\newcommand{\thesischtwomC}{{\bf C}}
\newcommand{\thesischtwomE}{{\bf E}}
\newcommand{\thesischtwomF}{{\bf F}}
\newcommand{\thesischtwomG}{{\bf G}}
\newcommand{\thesischtwomH}{{\bf H}}
\newcommand{\thesischtwomI}{{\bf I}}
\newcommand{\thesischtwomJ}{{\bf J}}
\newcommand{\thesischtwomK}{{\bf K}}
\newcommand{\thesischtwomL}{{\bf L}}
\newcommand{\thesischtwomM}{{\bf M}}
\newcommand{\thesischtwomN}{{\bf N}}
\newcommand{\thesischtwomO}{{\bf O}}
\newcommand{\thesischtwomP}{{\bf P}}
\newcommand{\thesischtwomQ}{{\bf Q}}
\newcommand{\thesischtwomR}{{\bf R}}
\newcommand{\thesischtwomS}{{\bf S}}
\newcommand{\thesischtwomT}{{\bf T}}
\newcommand{\thesischtwomU}{{\bf U}}
\newcommand{\thesischtwomV}{{\bf V}}
\newcommand{\thesischtwomW}{{\bf W}}
\newcommand{\thesischtwomX}{{\bf X}}
\newcommand{\thesischtwomY}{{\bf Y}}
\newcommand{\thesischtwomZ}{{\bf Z}}
\newcommand{\thesischtwosign}{\mathrm{sign}}
\newcommand{\thesischtwocnorm}{x}
\newcommand{\thesischtwoEE}{\Exp}
\newcommand{\thesischtwoPP}{\Prob}
\newcommand{\thesischtwoCC}{\mathbb{C}}
\newcommand{\thesischtwoVV}{\mathbf{V}}
\newcommand{\thesischtwoprox}{\mathop{\mathrm{prox}}\nolimits}
\newcommand{\thesischtwoproxR}{\thesischtwoprox_{\gamma R}}
\newcommand{\thesischtwoproxkR}{\thesischtwoprox_{\gamma^k R}}
\newcommand{\thesischtwomean}{\overline}
\newcommand{\thesischtwosumin}{\sum_{i=1}^n}
\newcommand{\thesischtwotx}{\widetilde{x}}
\newcommand{\thesischtwotX}{\widetilde{X}}
\newcommand{\thesischtwoty}{\widetilde{y}}
\newcommand{\thesischtwotF}{\widetilde{F}}
\newcommand{\thesischtwotnabla}{\widetilde{\nabla}}
\newcommand{\thesischtwohx}{\widehat{x}}
\newcommand{\thesischtwohy}{\widehat{y}}
\def\Bxi{\boldsymbol{\xi}}
\def\clip{\texttt{clip}}
\def\avg{\texttt{avg}}
\def\gap{\texttt{Gap}}
\ifcsname thesischtwodhatheight\endcsname\else\newlength{\thesischtwodhatheight}\fi
\newcommand{\thesischtwodoublehat}[1]{%
    \settoheight{\thesischtwodhatheight}{\ensuremath{\hat{#1}}}%
    \addtolength{\thesischtwodhatheight}{-0.35ex}%
    \hat{\vphantom{\rule{1pt}{\thesischtwodhatheight}}%
    \smash{\hat{#1}}}}
\usetikzlibrary{automata, positioning, arrows, shapes, fit, calc, intersections}
\usepgfplotslibrary{statistics}
\def\la{\langle}
\def\ra{\rangle}

\newcommand{\thesischthreetheHalgorithm}{\arabic{algorithm}}
\newcommand{\thesischthreecmark}{{\color{PineGreen}\ding{51}}}%
\newcommand{\thesischthreexmark}{{\color{BrickRed}\ding{55}}}%
\newcommand{\thesischthreeExp}{\Exp}
\newcommand{\thesischthreeProb}{\Prob}
\newcommand{\thesischthreeoldstuff}[1]{ {\small \color{blue} #1}}
\newcommand{\thesischthreeredstuff}[1]{ {\small \color{red} #1}}
\newcommand\thesischthreefigscale{0.25}
\newcommand\thesischthreeevofigscale{0.21}
\newcommand{\thesischthreesqueeze}{} %
\newcommand{\thesischthreealgnamex}[1]{\algname{#1}}
\newcommand\thesischthreetagthis{\addtocounter{equation}{1}\tag{\theequation}}
\newcommand{\thesischthreeargmin}{\mathop{\arg\!\min}}
\newcommand{\thesischthreeSp}{\mathrm{Sp}}
\newcommand{\thesischthreeId}{\mathrm{Id}}
\newcommand{\thesischthreeTr}{\mathrm{Tr}}
\newcommand{\thesischthreecA}{{\cal A}}
\newcommand{\thesischthreecB}{{\cal B}}
\newcommand{\thesischthreecC}{{\cal C}}
\newcommand{\thesischthreecD}{{\cal D}}
\newcommand{\thesischthreecE}{{\cal E}}
\newcommand{\thesischthreecF}{{\cal F}}
\newcommand{\thesischthreecG}{{\cal G}}
\newcommand{\thesischthreecH}{{\cal H}}
\newcommand{\thesischthreecJ}{{\cal J}}
\newcommand{\thesischthreecK}{{\cal K}}
\newcommand{\thesischthreecL}{{\cal L}}
\newcommand{\thesischthreecM}{{\cal M}}
\newcommand{\thesischthreecN}{{\cal N}}
\newcommand{\thesischthreecP}{{\cal P}}
\newcommand{\thesischthreecQ}{{\cal Q}}
\newcommand{\thesischthreecR}{{\cal R}}
\newcommand{\thesischthreecS}{{\cal S}}
\newcommand{\thesischthreecT}{{\cal T}}
\newcommand{\thesischthreecU}{{\cal U}}
\newcommand{\thesischthreecV}{{\cal V}}
\newcommand{\thesischthreecX}{{\cal X}}
\newcommand{\thesischthreecY}{{\cal Y}}
\newcommand{\thesischthreecW}{{\cal W}}
\newcommand{\thesischthreecZ}{{\cal Z}}
\newcommand{\thesischthreeVar}{\mathrm{Var}}
\newcommand{\thesischthreemA}{{\bf A}}
\newcommand{\thesischthreemB}{{\bf B}}
\newcommand{\thesischthreemC}{{\bf C}}
\newcommand{\thesischthreemE}{{\bf E}}
\newcommand{\thesischthreemF}{{\bf F}}
\newcommand{\thesischthreemG}{{\bf G}}
\newcommand{\thesischthreemH}{{\bf H}}
\newcommand{\thesischthreemI}{{\bf I}}
\newcommand{\thesischthreemJ}{{\bf J}}
\newcommand{\thesischthreemK}{{\bf K}}
\newcommand{\thesischthreemL}{{\bf L}}
\newcommand{\thesischthreemM}{{\bf M}}
\newcommand{\thesischthreemN}{{\bf N}}
\newcommand{\thesischthreemO}{{\bf O}}
\newcommand{\thesischthreemP}{{\bf P}}
\newcommand{\thesischthreemQ}{{\bf Q}}
\newcommand{\thesischthreemR}{{\bf R}}
\newcommand{\thesischthreemS}{{\bf S}}
\newcommand{\thesischthreemT}{{\bf T}}
\newcommand{\thesischthreemU}{{\bf U}}
\newcommand{\thesischthreemV}{{\bf V}}
\newcommand{\thesischthreemW}{{\bf W}}
\newcommand{\thesischthreemX}{{\bf X}}
\newcommand{\thesischthreemY}{{\bf Y}}
\newcommand{\thesischthreemZ}{{\bf Z}}
\newcommand{\thesischthreesign}{\mathrm{sign}}
\newcommand{\thesischthreecnorm}{x}
\newcommand{\thesischthreeEE}{\Exp}
\newcommand{\thesischthreeSD}{\mathbb{S}}
\newcommand{\thesischthreeCC}{\mathbb{C}}
\newcommand{\thesischthreeVV}{\mathbf{V}}
\newcommand{\thesischthreeprox}{\mathop{\mathrm{prox}}\nolimits}
\newcommand{\thesischthreeproj}{\mathop{\mathrm{proj}}\nolimits}
\newcommand{\thesischthreeproxR}{\thesischthreeprox_{\gamma R}}
\newcommand{\thesischthreeproxkR}{\thesischthreeprox_{\gamma^k R}}
\newcommand{\thesischthreeproxkPsi}{\thesischthreeprox_{\gamma \Psi}}
\newcommand{\thesischthreemean}{\overline}
\newcommand{\thesischthreesumin}{\sum_{i=1}^n}
\newcommand{\thesischthreePP}{\Prob}
\newcommand{\thesischthreetx}{\widetilde{x}}
\newcommand{\thesischthreetX}{\widetilde{X}}
\newcommand{\thesischthreetg}{\widetilde{g}}
\newcommand{\thesischthreety}{\widetilde{y}}
\newcommand{\thesischthreetF}{\widetilde{F}}
\newcommand{\thesischthreehx}{\widehat{x}}
\newcommand{\thesischthreehy}{\widehat{y}}
\newcommand{\thesischthreehg}{\widehat{g}}
\newcommand{\thesischthreethh}{\widetilde{h}}
\newcommand{\thesischthreehh}{\widehat{h}}
\newcommand{\thesischthreetnabla}{\widetilde{\nabla}}
\ifcsname thesischthreedhatheight\endcsname\else\newlength{\thesischthreedhatheight}\fi
\newcommand{\thesischthreedoublehat}[1]{%
    \settoheight{\thesischthreedhatheight}{\ensuremath{\hat{#1}}}%
    \addtolength{\thesischthreedhatheight}{-0.35ex}%
    \hat{\vphantom{\rule{1pt}{\thesischthreedhatheight}}%
    \smash{\hat{#1}}}}
\def\pd#1{{\color{black}#1}}

\definecolor{darkgreen}{rgb}{0.0, 0.5, 0.0}
\definecolor{darkblue}{rgb}{0.0, 0.0, 0.55}
\newcommand{\thesischfourspecialcell}[2][c]{%
  \begin{tabular}[#1]{@{}c@{}}#2\end{tabular}}
\newcommand{\thesischfourcmark}{{\color{PineGreen}\ding{51}}}%
\newcommand{\thesischfourxmark}{{\color{BrickRed}\ding{55}}}%
\newcommand{\thesischfourargmin}{\mathop{\arg\!\min}}
\newcommand{\thesischfourargmax}{\mathop{\arg\!\max}}
\newcommand{\thesischfourEndProof}{\begin{flushright}$\square$\end{flushright}}
\newcommand{\thesischfourcA}{{\cal A}}
\newcommand{\thesischfourcB}{{\cal B}}
\newcommand{\thesischfourcC}{{\cal C}}
\newcommand{\thesischfourcD}{{\cal D}}
\newcommand{\thesischfourcE}{{\cal E}}
\newcommand{\thesischfourcF}{{\cal F}}
\newcommand{\thesischfourcG}{{\cal G}}
\newcommand{\thesischfourcH}{{\cal H}}
\newcommand{\thesischfourcJ}{{\cal J}}
\newcommand{\thesischfourcK}{{\cal K}}
\newcommand{\thesischfourcL}{{\cal L}}
\newcommand{\thesischfourcM}{{\cal M}}
\newcommand{\thesischfourcN}{{\cal N}}
\newcommand{\thesischfourcP}{{\cal P}}
\newcommand{\thesischfourcQ}{{\cal Q}}
\newcommand{\thesischfourcR}{{\cal R}}
\newcommand{\thesischfourcS}{{\cal S}}
\newcommand{\thesischfourcT}{{\cal T}}
\newcommand{\thesischfourcU}{{\cal U}}
\newcommand{\thesischfourcV}{{\cal V}}
\newcommand{\thesischfourcX}{{\cal X}}
\newcommand{\thesischfourcY}{{\cal Y}}
\newcommand{\thesischfourcW}{{\cal W}}
\newcommand{\thesischfourcZ}{{\cal Z}}
\newcommand{\thesischfourVar}{\mathrm{Var}}
\newcommand{\thesischfourmA}{{\bf A}}
\newcommand{\thesischfourmB}{{\bf B}}
\newcommand{\thesischfourmC}{{\bf C}}
\newcommand{\thesischfourmE}{{\bf E}}
\newcommand{\thesischfourmF}{{\bf F}}
\newcommand{\thesischfourmG}{{\bf G}}
\newcommand{\thesischfourmH}{{\bf H}}
\newcommand{\thesischfourmI}{{\bf I}}
\newcommand{\thesischfourmJ}{{\bf J}}
\newcommand{\thesischfourmK}{{\bf K}}
\newcommand{\thesischfourmL}{{\bf L}}
\newcommand{\thesischfourmM}{{\bf M}}
\newcommand{\thesischfourmN}{{\bf N}}
\newcommand{\thesischfourmO}{{\bf O}}
\newcommand{\thesischfourmP}{{\bf P}}
\newcommand{\thesischfourmQ}{{\bf Q}}
\newcommand{\thesischfourmR}{{\bf R}}
\newcommand{\thesischfourmS}{{\bf S}}
\newcommand{\thesischfourmT}{{\bf T}}
\newcommand{\thesischfourmU}{{\bf U}}
\newcommand{\thesischfourmV}{{\bf V}}
\newcommand{\thesischfourmW}{{\bf W}}
\newcommand{\thesischfourmX}{{\bf X}}
\newcommand{\thesischfourmY}{{\bf Y}}
\newcommand{\thesischfourmZ}{{\bf Z}}
\newcommand{\thesischfoursign}{\mathrm{sign}}
\newcommand{\thesischfourcnorm}{\omega}
\newcommand{\thesischfourEE}{\Exp}
\newcommand{\thesischfourVV}{\mathbf{V}}
\def\p{\mathfrank{p}}
\newcommand{\thesischfourE}{\Exp}
\def\W{\mathcal W}
\def\C{\mathcal C}
\def\X{\mathcal X}
\def\Y{\mathcal Y}
\def\M{\mathcal M_+^1}
\def\thesischfourE{\Exp}
\def\thesischfourEE{\Exp}
\def\PP{\Prob}
\def\D{\mathbb D}
\def\U{\mathcal{U}}
\def\V{\mathcal{V}}
\def\e{\varepsilon}
\def\vp{f}
\def\y{\mathbf{y}}
\def\a{\mathbf{a}}
\def\x{\mathbf{x}}
\def\w{\mathbf{w}}
\def\Rbf{\mathbf{R}}
\def\one{{\mathbf 1}}
\def\tf{\tilde{f}}
\def\tL{\tilde{L}}
\def\tg{\tilde{g}}
\def\tnabla{\tilde{\nabla}}
\def\tx{\tilde{x}}
\def\ty{\tilde{y}}
\def\Blm{\boldsymbol{\lambda}}
\def\Bzeta{\boldsymbol{\zeta}}
\def\Beta{\boldsymbol{\eta}}
\def\BBlm{\boldsymbol{\overline{\lambda}}}
\def\BBzeta{\boldsymbol{\overline{\zeta}}}
\def\BBeta{\boldsymbol{\overline{\eta}}}
\def\Bx{\boldsymbol{x}}
\def\p{\mathtt{p}}
\def\lm{\lambda}
\def\eig{\lambda_{max}(W)}
\def\bld{\boldsymbol}
\def\ag#1{{\color{red}#1}}
\def\xs{x^{\star}}
\newcommand{\thesischfourExp}[2][]{
	\ifthenelse{\equal{#1}{}}{\thesischfourE \left[ #2 \right]}{\thesischfourE \left[ #2 \, \middle| \, #1 \right]}
}
\newcommand{\thesischfourphit}{\phi_t}
\newcommand{\thesischfourpsit}{\psi_t}
\newcommand{\thesischfourtSum}{\sum_{\ii=\fin}^\li}
\newcommand{\thesischfourExpc}[2]{\Exp\left[#1 \,\middle| \, #2\right]}
\newcommand{\thesischfourtinN}{{t \in \N}}
\newcommand{\thesischfourtinNp}{{t \geq 1}}
\newcommand{\thesischfourAa}{\eta \thesischfoursm}
\newcommand{\thesischfournnFxt}{\norm{\nFxt}}
\newcommand{\thesischfourood}{{ \nicefrac 1 \delta }}
\newcommand{\thesischfourlogood}{\log \pare \thesischfourood}
\newcommand{\thesischfourEp}{E_p}
\newcommand{\thesischfoursm}{L}
\newcommand{\thesischfourlsmooth}{(\thesischfoursm-Smoothness)}
\newcommand{\thesischfourlb}{(LB)}
\newcommand{\thesischfourpBCM}{(p-BCM)}
\newcommand{\thesischfournsgdm}{N-SGDM}
\newcommand{\thesischfourxfin}{\iterationx \fin}
\newcommand{\thesischfourrb}[1]{\left(#1\right)}
\providecommand{\polylog}{\operatorname{polylog}}
\definecolor{linen}{HTML}{FAF0E6} %

\newcommand{\thesischfivecommentalgo}[1]{\textcolor{blue}{#1}}%
\definecolor{bgcolor}{rgb}{0.76,0.88,0.50}
\definecolor{bgcolor0}{rgb}{0.93,0.99,1}
\definecolor{bgcolor1}{rgb}{0.8,1,1}
\definecolor{bgcolor3}{rgb}{0.50,0.90,0.50}
\definecolor{mydarkgreen}{RGB}{39,130,67}
\definecolor{mydarkorange}{RGB}{236,147,14}
\definecolor{mydarkred}{RGB}{192,47,25}
\definecolor{ruby}{RGB}{155,17,30}
\definecolor{chili}{RGB}{191,0,0}
\definecolor{sangria}{RGB}{146,0,10}
\definecolor{burgundy}{RGB}{128,0,32}
\ifnum\ThesisProofreadingColours=1 \definecolor{darkred}{RGB}{132,0,0}\fi
\definecolor{cherry}{RGB}{192,0,0}
\newcommand{\thesischfivegreen}{\color{mydarkgreen}}
\newcommand{\thesischfiveorange}{\color{mydarkorange}}
\newcommand{\thesischfiveblue}{\color{blue}}
\newcommand{\thesischfivered}{\color{cherry}}
\definecolor{linen}{HTML}{FAF0E6}
\definecolor{mydarkblue}{RGB}{20,20,192}
\newcommand{\thesischfivemyblue}{\color{mydarkblue}}
\newcommand{\thesischfiveniceblue}{\color{niceblue}}
\newcommand{\thesischfivecmark}{{\color{PineGreen}\ding{51}}}%
\newcommand{\thesischfivexmark}{{\color{BrickRed}\ding{55}}}%
\newcommand{\thesischfivedataname}[1]{{\tt\footnotesize\color{blue}#1}\xspace}
\newcommand{\thesischfivealgnamebf}[1]{\algname{\textbf{#1}}}
\newcommand{\thesischfiveind}{\perp\!\!\!\!\perp}
\newcommand{\thesischfiveTVdist}[1]{\left\| #1 \right\|_{\textnormal{TV}}}
\newcommand{\thesischfivelin}[1]{\left\langle #1\right\rangle} %
\newcommand{\thesischfiveE}[1]{\Exp\left[#1\right]}
\newcommand{\thesischfiveMed}[1]{\mathrm{Med}\left[#1\right]}
\newcommand{\thesischfiveVar}[1]{\mathrm{Var}\left(#1\right)}
\newcommand{\thesischfiveExp}[1]{\Exp\left[#1\right]}
\newcommand{\thesischfiveExpSub}[2]{\Exp_{#1}\left[#2\right]}
\newcommand{\thesischfiveExpCond}[2]{\Exp\left[\left.#1\,\right\vert\,#2\right]}
\newcommand{\thesischfiveExpSubCond}[3]{\Exp_{#1}\left[\left.#2\,\right\vert\,#3\right]}
\newcommand{\thesischfiveProba}[1]{\Prob\left(#1\right)} %
\newcommand{\thesischfiveProbCond}[2]{\Prob\left(#1\middle\vert#2\right)}
\newcommand{\thesischfiveABC}{\algname{ABC}}
\newcommand{\thesischfiveSGD}{\algname{SGD}}
\newcommand{\thesischfivePAGE}{\algname{PAGE}}
\newcommand{\thesischfiveVanSGD}{\algname{Vanilla SGD}}
\newcommand{\thesischfiveCompr}[1]{\mathcal{C}\left(#1\right)}
\newcommand{\thesischfiveRandK}{\algname{RandK}}
\newcommand{\thesischfiveRandC}[1]{\algname{Rand#1}}
\newcommand{\thesischfiveAlgZer}{\mathcal{A}_{\normalfont\text{rz}}}
\newcommand{\thesischfiveFuncLB}{\mathcal{F}_{\Delta, L}}
\newcommand{\thesischfivecA}{\mathcal{A}}
\newcommand{\thesischfivecB}{\mathcal{B}}
\newcommand{\thesischfivecC}{\mathcal{C}}
\newcommand{\thesischfivecD}{\mathcal{D}}
\newcommand{\thesischfivecE}{\mathcal{E}}
\newcommand{\thesischfivecF}{\mathcal{F}}
\newcommand{\thesischfivecH}{\mathcal{H}}
\newcommand{\thesischfivecL}{\mathcal{L}}
\newcommand{\thesischfivecT}{\mathcal{T}}
\newcommand{\thesischfivecP}{\mathcal{P}}
\newcommand{\thesischfivecQ}{\mathcal{Q}}
\newcommand{\thesischfivecS}{\mathcal{S}}
\newcommand{\thesischfivecW}{\mathcal{W}}
\newcommand{\thesischfivecY}{\mathcal{Y}}
\newcommand{\thesischfivecZ}{\mathcal{Z}}
\newcommand{\thesischfivecJ}{\mathcal{J}}
\newcommand{\thesischfivecN}{\mathcal{N}}
\newcommand{\thesischfivemA}{\mathbf{A}}
\newcommand{\thesischfivemX}{\mathbf{X}}
\newcommand{\thesischfivemM}{\mathbf{M}}
\newcommand{\thesischfivemN}{\mathbf{N}}
\newcommand{\thesischfivemI}{\mathbf{I}}
\newcommand{\thesischfivemJ}{\mathbf{J}}
\newcommand{\thesischfivemL}{\mathbf{L}}
\newcommand{\thesischfivemD}{\mathbf{D}}
\newcommand{\thesischfivemO}{\mathbf{O}}
\newcommand{\thesischfiveaddeq}{\addtocounter{equation}{1}}%
\newcommand\thesischfivenumberthis{\thesischfiveaddeq\tag{\theequation}}%
\NewEnviron{thesischfivefitdisplay}{%
  \begin{equation*}
    \hbox{$\displaystyle\BODY$}%
  \end{equation*}}
\NewEnviron{thesischfivefitalignat}[1]{%
  \begin{equation*}
    \hbox{$\displaystyle
      \begin{alignedat}{#1}\BODY\end{alignedat}$}%
  \end{equation*}}
\NewEnviron{thesischfivefitalignatnum}[1]{%
  \begin{equation}
    \hbox{$\displaystyle
      \begin{alignedat}{#1}\BODY\end{alignedat}$}%
  \end{equation}}
\makeatletter
\DeclareRobustCommand\widecheck[1]{{\mathpalette\@widecheck{#1}}}
\makeatother
\makeatletter
\def\@widecheck#1#2{%
    \setbox\z@\hbox{\m@th$#1#2$}%
    \setbox\tw@\hbox{\m@th$#1%
       \widehat{%
          \vrule\@width\z@\@height\ht\z@
          \vrule\@height\z@\@width\wd\z@}$}%
    \dp\tw@-\ht\z@
    \@tempdima\ht\z@ \advance\@tempdima2\ht\tw@ \divide\@tempdima\thr@@
    \setbox\tw@\hbox{%
       \raise\@tempdima\hbox{\scalebox{1}[-1]{\lower\@tempdima\box
\tw@}}}%
    {\ooalign{\box\tw@ \cr \box\z@}}}
\makeatother
\newcommand{\thesischfivervline}{\hspace*{-\arraycolsep}\vline\hspace*{-\arraycolsep}}
\newcommand{\thesischfiveoversetref}[2]{\overset{{\scriptscriptstyle\mathrm{#1}\!~\text{#2}}}}
\newcommand{\thesischfiveoversetlab}[1]{\overset{\text{#1}}}
\ifcsname c@thesischfiverelctr\endcsname\else\newcounter{thesischfiverelctr}\fi %
\ifcsname c@thesischfiverelgroup\endcsname\else\newcounter{thesischfiverelgroup}\fi %
\renewcommand*\thethesischfiverelctr{\alph{thesischfiverelctr}} %
\newcommand*{\thesischfiveRoot}{{\normalfont\text{root}}}
\AtBeginDocument{\let\originallabel\label} %
\makeatletter
\newcommand\thesischfiveoversetrel[1]{%
    \ifnum\value{thesischfiverelctr}=0\relax%
      \stepcounter{thesischfiverelgroup}%
    \fi%
    \refstepcounter{thesischfiverelctr}%
    \phantomsection%
    \originallabel{rel:#1@\thethesischfiverelgroup}%
    \overset{\scriptscriptstyle\text{(\alph{thesischfiverelctr})}}%
}
\makeatother
\makeatletter
\newcommand\thesischfiverelref[1]{%
    \@ifundefined{r@rel:#1@\arabic{thesischfiverelgroup}}{%
        \eqref{}%
    }{%
        \eqref{rel:#1@\arabic{thesischfiverelgroup}}%
    }%
}
\makeatother
\newcommand{\thesischfivereqdef}{=:}
\newcommand{\thesischfiveMod}[1]{\ (\mathrm{mod}\ #1)}
\makeatletter
\newcommand{\thesischfivevast}{\bBigg@{4}}
\makeatother
\newcommand*{\thesischfivesketchproofname}{Sketch of Proof}
\newenvironment{sketchproof}[1][\thesischfivesketchproofname]{
  \begin{proof}[#1] } {
  \end{proof}
}
\newcommand{\thesischfiveprog}{\mathrm{prog}\,}
\newcommand{\thesischfiveconv}{\mathrm{Conv}\,}
\newcommand{\thesischfivecconv}{\overline{\mathrm{Conv}}\,}
\newcommand{\thesischfiveextr}{\mathrm{Extr}\,}
\newcommand{\thesischfiveepi}{\mathrm{epi}\,}
\newcommand{\thesischfivehypo}{\mathrm{hypo}\,}
\newcommand{\thesischfiveinte}{\mathrm{int}\,}
\newcommand{\thesischfiveexte}{\mathrm{ext}\,}
\newcommand{\thesischfiverank}{\mathrm{rank}\,}
\newcommand{\thesischfiveInt}[2]{\left\{ #1, \ldots, #2 \right\}}
\newcommand{\thesischfivesupp}{\mathrm{supp}\,}
\newcommand*{\thesischfivelbd}{\lambda}%
\newcommand*{\thesischfiveeps}{\varepsilon}%
\newcommand*{\thesischfivezeros}{\mathop{\mathscr{Z}}}%
\newcommand*{\thesischfivedom}{\mathop{\mathrm{dom}}}%
\newcommand*{\thesischfivemul}{\mathop{\mathrm{m}}}%
\newcommand*{\thesischfiveres}{\mathop{\mathrm{Res}}}%
\newcommand*{\thesischfiveFloor}[1]{\left\lfloor #1 \right\rfloor}%
\newcommand*{\thesischfiveCeil}[1]{\left\lceil #1 \right\rceil}%
\newcommand*{\thesischfiveval}{\nu}%
\newcommand*{\thesischfivepgcd}{\mathop{\mathrm{PGCD}}}%
\newcommand*{\thesischfiveppcm}{\mathop{\mathrm{PPCM}}}%
\newcommand*{\thesischfiveprimepi}{\mathop{\pi}}%
\newcommand*{\thesischfivemat}{\mathop{\mathrm{mat}}\limits}%
\newcommand*{\thesischfivecoord}{\mathop{\mathrm{Coord}}\limits}%
\newcommand*{\thesischfiveT}[1]{{#1}^{\top}}%
\newcommand*{\thesischfiveKer}{\mathop{\mathrm{Ker}}}%
\newcommand*{\thesischfiveImg}{\mathop{\mathrm{Im}}}%
\newcommand*{\thesischfiverg}{\mathop{\mathrm{rg}}}%
\newcommand*{\thesischfiveId}{\mathop{\mathrm{Id}}}%
\newcommand*{\thesischfiveTr}{\mathop{\mathrm{Tr}}}%
\newcommand*{\thesischfiveCom}{\mathop{\mathrm{Com}}}%
\newcommand*{\thesischfiveiddots}{\rotatebox[origin = c]{90}{$\ddots$}}%
\newcommand*{\thesischfiveGL}{\mathop{\mathrm{GL}}}%
\newcommand*{\thesischfiveSL}{\mathop{\mathrm{SL}}}%
\newcommand*{\thesischfiveSO}{\mathop{\mathrm{SO}}}%
\newcommand*{\thesischfiveSto}{\mathop{\mathrm{ST}}}%
\newcommand*{\thesischfiveBSto}{\mathop{\mathrm{BST}}}%
\newcommand*{\thesischfiveSp}{\mathop{\mathrm{Sp}}}%
\newcommand*{\thesischfiveVect}[2][]{\mathop{\mathrm{Vect}}_{#1}\left(#2\right)}%
\newcommand*{\thesischfiveintervalle}[4]{\mathopen{#1}#2\mathclose{},#3\mathclose{#4}}%
\newcommand*{\thesischfiveintff}[2]{\thesischfiveintervalle{[}{#1}{#2}{]}}%
\newcommand*{\thesischfiveintof}[2]{\thesischfiveintervalle{(}{#1}{#2}{]}}%
\newcommand*{\thesischfiveintfo}[2]{\thesischfiveintervalle{[}{#1}{#2}{)}}%
\newcommand*{\thesischfiveintoo}[2]{\thesischfiveintervalle{(}{#1}{#2}{)}}%
\newcommand*{\thesischfiveens}[1]{\left\{#1\right\}}%
\newcommand*{\thesischfiveenstq}[2]{\left\{#1\,:\,#2\right\}}%
\newcommand*{\thesischfiveCard}[1]{\left|#1\right|}%
\newcommand*{\thesischfivePartie}{\mathop{\mathcal{P}}}%
\newcommand*{\thesischfiveZ}{\mathbb{Z}}%
\newcommand*{\thesischfiveQ}{\mathbb{Q}}%
\newcommand*{\thesischfiveK}{\mathbb{K}}%
\newcommand*{\thesischfiveF}{\mathbb{F}}%
\newcommand*{\thesischfiveRelB}{\mathop{\mathscr{R}}}%
\newcommand*{\thesischfiveCl}[1]{\overline{#1}}%
\newcommand*{\thesischfivegr}{\mathop{\mathrm{gr}}}%
\newcommand*{\thesischfiveord}{\mathop{\mathrm{ord}}}%
\newcommand*{\thesischfiveOrb}{\mathop{\mathrm{Orb}}}%
\newcommand*{\thesischfiveCar}{\mathop{\mathrm{Car}}}%
\newcommand*{\thesischfiveAlg}{\mathop{\mathrm{Alg}}} %
\newcommand{\thesischfivedist}{\mathop{\mathrm{d}}}%
\newcommand{\thesischfivediam}{\delta}
\newcommand{\thesischfivebouleO}{\mathop{\mathrm{B}}}%
\newcommand{\thesischfivebouleF}{\overline{\mathop{\mathrm{B}}}}%
\newcommand{\thesischfivesphere}{\mathop{\mathrm{S}}}%
\newcommand{\thesischfiveadh}[1]{\overline{#1}}%
\newcommand{\thesischfiveext}{\mathop{\mathrm{ext}}}%
\newcommand{\thesischfivefr}{\mathop{\mathrm{fr}}}%
\newcommand*{\thesischfiveConv}{\mathop{\mathrm{Conv}}}%
\newcommand{\thesischfiveEpi}{\mathop{\mathrm{Epi}}}%
\newcommand{\thesischfiveortho}[1]{{#1}^{\bot}}%
\newcommand{\thesischfiveorthosum}{\mathop{\overset{\bot}{\oplus}}}%
\newcommand{\thesischfiveorthoproj}{\mathop{\mathrm{p}}^{\bot}}%
\newcommand{\thesischfiveps}[2]{\left\langle #1, #2 \right\rangle}%
\newcommand*{\thesischfiveSym}{\mathop{\mathfrak{S}}}%
\newcommand*{\thesischfivesig}{\mathop{\varepsilon}}%
\newcommand*{\thesischfiveInv}{\mathop{\mathrm{Inv}}}%
\newcommand*{\thesischfiveconj}[1]{\overline{#1}}%
\newcommand{\thesischfiveargsinh}{\mathop{\mathrm{argsinh}}}%
\newcommand{\thesischfiveargcosh}{\mathop{\mathrm{argcosh}}}%
\newcommand{\thesischfiveargtanh}{\mathop{\mathrm{argtanh}}}%
\newcommand*{\thesischfiveargmax}{\mathop{\mathrm{arg\,max}}}%
\newcommand*{\thesischfiveargmin}{\mathop{\mathrm{arg\,min}}}%
\newcommand*{\thesischfiveJac}{\mathrm{Jac}}%
\newcommand*{\thesischfivediag}{\mathrm{diag}}
\newcommand*{\thesischfiveproba}[1]{\Prob\left( #1 \right)}%
\newcommand*{\thesischfiveprobac}[2]{\Prob\left( #1 \, \vert \, #2 \right)}%
\newcommand*{\thesischfiveesp}{\Exp}%
\newcommand*{\thesischfivevar}{\mathop{\mathbb{V}}}%
\newcommand*{\thesischfivecov}{\mathop{\mathrm{cov}}}%
\newcommand*{\thesischfivecorr}{\mathop{\rho}}%
\newcommand*{\thesischfiveBer}{\mathrm{Ber}}
\newcommand{\thesischfiveclip}{{\normalfont\texttt{clip}}}
\newcommand{\thesischfivesign}{\mathop{\mathrm{sign}}}
\definecolor{table-green}{rgb}{0., 0.43, 0.2}
\newcommand{\thesischfivelya}{H}
\ifcsname thesischfivetest\endcsname\else\newsavebox{\thesischfivetest}\fi %

\definecolor{darkscarlet}{rgb}{0.34, 0.01, 0.1}
\definecolor{yaleblue}{rgb}{0.06, 0.3, 0.57}
\definecolor{darkpowderblue}{rgb}{0.0, 0.2, 0.6}
\definecolor{midnightblue}{HTML}{0059b3}
\definecolor{noonblue}{HTML}{e5eef7}
\definecolor{chromered}{HTML}{f14233}
\definecolor{olivedrab}{HTML}{6b8e23}
\newcommand{\thesischsixouralg}{\algname{Fed-$\alpha$-NormEC}}
\newcommand{\thesischsixDPSGDalg}{\algname{DP-SGD}}
\newcommand{\thesischsixprioralg}{\algname{$\alpha$-NormEC}}
\newcommand{\thesischsixg}{{\color{red}g}}
\newcommand{\thesischsixh}{{\color{blue}h}}
\newcommand{\thesischsixlin}[1]{\left\langle #1\right\rangle} %
\newcommand{\thesischsixcA}{\mathcal{A}}
\newcommand{\thesischsixcB}{\mathcal{B}}
\newcommand{\thesischsixcC}{\mathcal{C}}
\newcommand{\thesischsixcD}{\mathcal{D}}
\newcommand{\thesischsixcH}{\mathcal{H}}
\newcommand{\thesischsixcL}{\mathcal{L}}
\newcommand{\thesischsixcQ}{\mathcal{Q}}
\newcommand{\thesischsixcS}{\mathcal{S}}
\newcommand{\thesischsixcW}{\mathcal{W}}
\newcommand{\thesischsixcY}{\mathcal{Y}}
\newcommand{\thesischsixcZ}{\mathcal{Z}}
\newcommand{\thesischsixmI}{\mathbf{I}}
\newcommand{\thesischsixhotidea}{{\color{red}\bf HOT IDEA: }}
\newcommand{\thesischsixdone}{{\color{blue}\bf DONE: }}
\newcommand{\thesischsixdel}[1]{}
\let\la=\langle
\let\ra=\rangle
\newcommand{\thesischsixU}{\mathbb{U}} %
\newcommand{\thesischsixPermComp}{\thesischsixU\mathbb{P}} %
\newcommand{\thesischsixst}{\;:\;} %
\newcommand{\thesischsixProb}{\Prob} %
\newcommand{\thesischsixProbCond}[2]{\thesischsixProb\left(\left.#1\right\vert#2\right)}
\newcommand{\thesischsixnnz}[1]{{\color{red}\|#1\|_0}}
\providecommand{\card}{\operatorname{card}}         %
\providecommand{\diam}{\operatorname{diam}}        %
\providecommand{\vol}{\operatorname{vol}}               %
\newcommand{\thesischsixExp}[1]{\Exp\left[#1\right]}
\newcommand{\thesischsixExpCond}[2]{\Exp\left[\left.#1\right\vert#2\right]}
\newcommand{\thesischsixExpSub}[2]{\Exp_{#1}\left[#2\right]}
\providecommand{\Cov}{\operatorname{Cov}}         %
\providecommand{\Var}{\operatorname{Var}}         %
\providecommand{\Corr}{\operatorname{Corr}}       %
\NewEnviron{thesischsixfiteqnarray}{%
  \begin{equation*}
    \hbox{$\displaystyle
      \begin{alignedat}{2}\BODY\end{alignedat}$}%
  \end{equation*}%
}
\NewEnviron{thesischsixfitalign}{%
  \begin{equation*}
    \hbox{$\displaystyle
      \begin{aligned}\BODY\end{aligned}$}%
  \end{equation*}%
}
\providecommand{\signum}{\operatorname{sign}}     %
\providecommand{\dom}{\operatorname{dom}}         %
\providecommand{\epi}{\operatorname{epi}}         %
\providecommand{\Ker}{\operatorname{null}}        %
\providecommand{\nullspace}{\operatorname{null}}  %
\providecommand{\range}{\operatorname{range}}     %
\providecommand{\Image}{\operatorname{Im}}        %
\providecommand{\interior}{\operatorname{int}}    %
\providecommand{\ri}{\operatorname{rint}}         %
\providecommand{\rint}{\operatorname{rint}}       %
\providecommand{\bdry}{\operatorname{bdry}}       %
\providecommand{\cl}{\operatorname{cl}}           %
\providecommand{\linspan}{\operatorname{span}}
\providecommand{\linspace}{\operatorname{linspace}}
\providecommand{\cone}{\operatorname{cone}}
\providecommand{\tr}{\operatorname{tr}}           %
\providecommand{\rank}{\operatorname{rank}}       %
\providecommand{\conv}{\operatorname{conv}}       %
\providecommand{\Diag}{\operatorname{Diag}}       %
\providecommand{\diag}{\operatorname{diag}}       %
\providecommand{\Arg}{\operatorname{Arg}}         %
\newcommand{\thesischsixSample}{\mathcal{S}}
\newcommand{\thesischsixB}{\mathbb{B}}
\tcolorboxenvironment{algorithms}{
	colback=gray!10,
	boxrule=0pt,
	boxsep=1pt,
	left=2pt,right=2pt,top=2pt,bottom=2pt,
	oversize=2pt,
	sharp corners,
	before skip=\topsep,
	after skip=\topsep,
}
\newcommand{\thesischsixforarxiv}[1]{{\color{blue}#1}}
\newcommand{\thesischsixfigleft}{{\em (Left)}}
\newcommand{\thesischsixfigcenter}{{\em (Center)}}
\newcommand{\thesischsixfigright}{{\em (Right)}}
\newcommand{\thesischsixfigtop}{{\em (Top)}}
\newcommand{\thesischsixfigbottom}{{\em (Bottom)}}
\newcommand{\thesischsixcaptiona}{{\em (a)}}
\newcommand{\thesischsixcaptionb}{{\em (b)}}
\newcommand{\thesischsixcaptionc}{{\em (c)}}
\newcommand{\thesischsixcaptiond}{{\em (d)}}
\newcommand{\thesischsixnewterm}[1]{{\bf #1}}
\def\figref#1{figure~\ref{#1}}
\def\Figref#1{Figure~\ref{#1}}
\def\twofigref#1#2{figures \ref{#1} and \ref{#2}}
\def\quadfigref#1#2#3#4{figures \ref{#1}, \ref{#2}, \ref{#3} and \ref{#4}}
\def\secref#1{section~\ref{#1}}
\def\Secref#1{Section~\ref{#1}}
\def\twosecrefs#1#2{sections \ref{#1} and \ref{#2}}
\def\secrefs#1#2#3{sections \ref{#1}, \ref{#2} and \ref{#3}}
\def\Eqref#1{Equation~\ref{#1}}
\def\plaineqref#1{\ref{#1}}
\def\chapref#1{chapter~\ref{#1}}
\def\Chapref#1{Chapter~\ref{#1}}
\def\rangechapref#1#2{chapters\ref{#1}--\ref{#2}}
\def\algref#1{algorithm~\ref{#1}}
\def\Algref#1{Algorithm~\ref{#1}}
\def\twoalgref#1#2{algorithms \ref{#1} and \ref{#2}}
\def\Twoalgref#1#2{Algorithms \ref{#1} and \ref{#2}}
\def\partref#1{part~\ref{#1}}
\def\Partref#1{Part~\ref{#1}}
\def\twopartref#1#2{parts \ref{#1} and \ref{#2}}
\def\ceil#1{\lceil #1 \rceil}
\def\floor#1{\lfloor #1 \rfloor}
\def\1{\bm{1}}
\newcommand{\thesischsixtrain}{\mathcal{D}}
\newcommand{\thesischsixvalid}{\mathcal{D_{\mathrm{valid}}}}
\newcommand{\thesischsixtest}{\mathcal{D_{\mathrm{test}}}}
\def\eps{{\epsilon}}
\def\reta{{\textnormal{$\eta$}}}
\def\rb{{\textnormal{b}}}
\def\rc{{\textnormal{c}}}
\def\rd{{\textnormal{d}}}
\def\re{{\textnormal{e}}}
\def\rf{{\textnormal{f}}}
\def\rg{{\textnormal{g}}}
\def\rh{{\textnormal{h}}}
\def\rj{{\textnormal{j}}}
\def\rk{{\textnormal{k}}}
\def\rl{{\textnormal{l}}}
\def\rn{{\textnormal{n}}}
\def\ro{{\textnormal{o}}}
\def\rp{{\textnormal{p}}}
\def\rq{{\textnormal{q}}}
\def\rr{{\textnormal{r}}}
\def\rs{{\textnormal{s}}}
\def\rt{{\textnormal{t}}}
\def\ru{{\textnormal{u}}}
\def\rv{{\textnormal{v}}}
\def\rw{{\textnormal{w}}}
\def\rx{{\textnormal{x}}}
\def\ry{{\textnormal{y}}}
\def\rz{{\textnormal{z}}}
\def\rvepsilon{{\mathbf{\epsilon}}}
\def\rvtheta{{\mathbf{\theta}}}
\def\rva{{\mathbf{a}}}
\def\rvb{{\mathbf{b}}}
\def\rvc{{\mathbf{c}}}
\def\rvd{{\mathbf{d}}}
\def\rve{{\mathbf{e}}}
\def\rvf{{\mathbf{f}}}
\def\rvg{{\mathbf{g}}}
\def\rvh{{\mathbf{h}}}
\def\rvu{{\mathbf{i}}}
\def\rvj{{\mathbf{j}}}
\def\rvk{{\mathbf{k}}}
\def\rvl{{\mathbf{l}}}
\def\rvm{{\mathbf{m}}}
\def\rvn{{\mathbf{n}}}
\def\rvo{{\mathbf{o}}}
\def\rvp{{\mathbf{p}}}
\def\rvq{{\mathbf{q}}}
\def\rvr{{\mathbf{r}}}
\def\rvs{{\mathbf{s}}}
\def\rvt{{\mathbf{t}}}
\def\rvu{{\mathbf{u}}}
\def\rvv{{\mathbf{v}}}
\def\rvw{{\mathbf{w}}}
\def\rvx{{\mathbf{x}}}
\def\rvy{{\mathbf{y}}}
\def\rvz{{\mathbf{z}}}
\def\erva{{\textnormal{a}}}
\def\ervb{{\textnormal{b}}}
\def\ervc{{\textnormal{c}}}
\def\ervd{{\textnormal{d}}}
\def\erve{{\textnormal{e}}}
\def\ervf{{\textnormal{f}}}
\def\ervg{{\textnormal{g}}}
\def\ervh{{\textnormal{h}}}
\def\ervi{{\textnormal{i}}}
\def\ervj{{\textnormal{j}}}
\def\ervk{{\textnormal{k}}}
\def\ervl{{\textnormal{l}}}
\def\ervm{{\textnormal{m}}}
\def\ervn{{\textnormal{n}}}
\def\ervo{{\textnormal{o}}}
\def\ervp{{\textnormal{p}}}
\def\ervq{{\textnormal{q}}}
\def\ervr{{\textnormal{r}}}
\def\ervs{{\textnormal{s}}}
\def\ervt{{\textnormal{t}}}
\def\ervu{{\textnormal{u}}}
\def\ervv{{\textnormal{v}}}
\def\ervw{{\textnormal{w}}}
\def\ervx{{\textnormal{x}}}
\def\ervy{{\textnormal{y}}}
\def\ervz{{\textnormal{z}}}
\def\rmA{{\mathbf{A}}}
\def\rmB{{\mathbf{B}}}
\def\rmC{{\mathbf{C}}}
\def\rmD{{\mathbf{D}}}
\def\rmE{{\mathbf{E}}}
\def\rmF{{\mathbf{F}}}
\def\rmG{{\mathbf{G}}}
\def\rmH{{\mathbf{H}}}
\def\rmI{{\mathbf{I}}}
\def\rmJ{{\mathbf{J}}}
\def\rmK{{\mathbf{K}}}
\def\rmL{{\mathbf{L}}}
\def\rmM{{\mathbf{M}}}
\def\rmN{{\mathbf{N}}}
\def\rmO{{\mathbf{O}}}
\def\rmP{{\mathbf{P}}}
\def\rmQ{{\mathbf{Q}}}
\def\rmR{{\mathbf{R}}}
\def\rmS{{\mathbf{S}}}
\def\rmT{{\mathbf{T}}}
\def\rmU{{\mathbf{U}}}
\def\rmV{{\mathbf{V}}}
\def\rmW{{\mathbf{W}}}
\def\rmX{{\mathbf{X}}}
\def\rmY{{\mathbf{Y}}}
\def\rmZ{{\mathbf{Z}}}
\def\ermA{{\textnormal{A}}}
\def\ermB{{\textnormal{B}}}
\def\ermC{{\textnormal{C}}}
\def\ermD{{\textnormal{D}}}
\def\ermE{{\textnormal{E}}}
\def\ermF{{\textnormal{F}}}
\def\ermG{{\textnormal{G}}}
\def\ermH{{\textnormal{H}}}
\def\ermI{{\textnormal{I}}}
\def\ermJ{{\textnormal{J}}}
\def\ermK{{\textnormal{K}}}
\def\ermL{{\textnormal{L}}}
\def\ermM{{\textnormal{M}}}
\def\ermN{{\textnormal{N}}}
\def\ermO{{\textnormal{O}}}
\def\ermP{{\textnormal{P}}}
\def\ermQ{{\textnormal{Q}}}
\def\ermR{{\textnormal{R}}}
\def\ermS{{\textnormal{S}}}
\def\ermT{{\textnormal{T}}}
\def\ermU{{\textnormal{U}}}
\def\ermV{{\textnormal{V}}}
\def\ermW{{\textnormal{W}}}
\def\ermX{{\textnormal{X}}}
\def\ermY{{\textnormal{Y}}}
\def\ermZ{{\textnormal{Z}}}
\def\vzero{{\bm{0}}}
\def\vone{{\bm{1}}}
\def\vmu{{\bm{\mu}}}
\def\vtheta{{\bm{\theta}}}
\def\va{{\bm{a}}}
\def\vb{{\bm{b}}}
\def\vc{{\bm{c}}}
\def\vd{{\bm{d}}}
\def\ve{{\bm{e}}}
\def\vf{{\bm{f}}}
\def\vg{{\bm{g}}}
\def\vh{{\bm{h}}}
\def\vi{{\bm{i}}}
\def\vj{{\bm{j}}}
\def\vk{{\bm{k}}}
\def\vl{{\bm{l}}}
\def\vm{{\bm{m}}}
\def\vn{{\bm{n}}}
\def\vo{{\bm{o}}}
\def\vp{{\bm{p}}}
\def\vq{{\bm{q}}}
\def\vr{{\bm{r}}}
\def\vs{{\bm{s}}}
\def\vt{{\bm{t}}}
\def\vu{{\bm{u}}}
\def\vv{{\bm{v}}}
\def\vw{{\bm{w}}}
\def\vx{{\bm{x}}}
\def\vy{{\bm{y}}}
\def\vz{{\bm{z}}}
\def\evalpha{{\alpha}}
\def\evbeta{{\beta}}
\def\evepsilon{{\epsilon}}
\def\evlambda{{\lambda}}
\def\evomega{{\omega}}
\def\evmu{{\mu}}
\def\evpsi{{\psi}}
\def\evsigma{{\sigma}}
\def\evtheta{{\theta}}
\def\eva{{a}}
\def\evb{{b}}
\def\evc{{c}}
\def\evd{{d}}
\def\eve{{e}}
\def\evf{{f}}
\def\evg{{g}}
\def\evh{{h}}
\def\evi{{i}}
\def\evj{{j}}
\def\evk{{k}}
\def\evl{{l}}
\def\evm{{m}}
\def\evn{{n}}
\def\evo{{o}}
\def\evp{{p}}
\def\evq{{q}}
\def\evr{{r}}
\def\evs{{s}}
\def\evt{{t}}
\def\evu{{u}}
\def\evv{{v}}
\def\evw{{w}}
\def\evx{{x}}
\def\evy{{y}}
\def\evz{{z}}
\def\mA{{\bm{A}}}
\def\mB{{\bm{B}}}
\def\mC{{\bm{C}}}
\def\mD{{\bm{D}}}
\def\mE{{\bm{E}}}
\def\mF{{\bm{F}}}
\def\mG{{\bm{G}}}
\def\mH{{\bm{H}}}
\def\mJ{{\bm{J}}}
\def\mK{{\bm{K}}}
\def\mM{{\bm{M}}}
\def\mN{{\bm{N}}}
\def\mO{{\bm{O}}}
\def\mP{{\bm{P}}}
\def\mQ{{\bm{Q}}}
\def\mR{{\bm{R}}}
\def\mS{{\bm{S}}}
\def\mT{{\bm{T}}}
\def\mU{{\bm{U}}}
\def\mV{{\bm{V}}}
\def\mW{{\bm{W}}}
\def\mX{{\bm{X}}}
\def\mY{{\bm{Y}}}
\def\mZ{{\bm{Z}}}
\def\mBeta{{\bm{\beta}}}
\def\mPhi{{\bm{\Phi}}}
\def\mLambda{{\bm{\Lambda}}}
\def\mSigma{{\bm{\Sigma}}}
\newcommand{\thesischsixtens}[1]{\bm{\mathsfit{#1}}}
\def\tA{{\thesischsixtens{A}}}
\def\tB{{\thesischsixtens{B}}}
\def\tC{{\thesischsixtens{C}}}
\def\tD{{\thesischsixtens{D}}}
\def\tE{{\thesischsixtens{E}}}
\def\tF{{\thesischsixtens{F}}}
\def\tG{{\thesischsixtens{G}}}
\def\tH{{\thesischsixtens{H}}}
\def\tI{{\thesischsixtens{I}}}
\def\tJ{{\thesischsixtens{J}}}
\def\tK{{\thesischsixtens{K}}}
\def\tL{{\thesischsixtens{L}}}
\def\tM{{\thesischsixtens{M}}}
\def\tN{{\thesischsixtens{N}}}
\def\tO{{\thesischsixtens{O}}}
\def\tP{{\thesischsixtens{P}}}
\def\tQ{{\thesischsixtens{Q}}}
\def\tR{{\thesischsixtens{R}}}
\def\tS{{\thesischsixtens{S}}}
\def\tT{{\thesischsixtens{T}}}
\def\tU{{\thesischsixtens{U}}}
\def\tV{{\thesischsixtens{V}}}
\def\tW{{\thesischsixtens{W}}}
\def\tX{{\thesischsixtens{X}}}
\def\tY{{\thesischsixtens{Y}}}
\def\tZ{{\thesischsixtens{Z}}}
\def\gA{{\mathcal{A}}}
\def\gB{{\mathcal{B}}}
\def\gC{{\mathcal{C}}}
\def\gD{{\mathcal{D}}}
\def\gE{{\mathcal{E}}}
\def\gF{{\mathcal{F}}}
\def\gG{{\mathcal{G}}}
\def\gH{{\mathcal{H}}}
\def\gI{{\mathcal{I}}}
\def\gJ{{\mathcal{J}}}
\def\gK{{\mathcal{K}}}
\def\gL{{\mathcal{L}}}
\def\gM{{\mathcal{M}}}
\def\gN{{\mathcal{N}}}
\def\gO{{\mathcal{O}}}
\def\gP{{\mathcal{P}}}
\def\gQ{{\mathcal{Q}}}
\def\gR{{\mathcal{R}}}
\def\gS{{\mathcal{S}}}
\def\gT{{\mathcal{T}}}
\def\gU{{\mathcal{U}}}
\def\gV{{\mathcal{V}}}
\def\gW{{\mathcal{W}}}
\def\gX{{\mathcal{X}}}
\def\gY{{\mathcal{Y}}}
\def\gZ{{\mathcal{Z}}}
\def\sA{{\mathbb{A}}}
\def\sB{{\mathbb{B}}}
\def\sC{{\mathbb{C}}}
\def\sD{{\mathbb{D}}}
\def\sF{{\mathbb{F}}}
\def\sG{{\mathbb{G}}}
\def\sH{{\mathbb{H}}}
\def\sI{{\mathbb{I}}}
\def\sJ{{\mathbb{J}}}
\def\sK{{\mathbb{K}}}
\def\sL{{\mathbb{L}}}
\def\sM{{\mathbb{M}}}
\def\sN{{\mathbb{N}}}
\def\sO{{\mathbb{O}}}
\def\sP{{\mathbb{P}}}
\def\sQ{{\mathbb{Q}}}
\def\sR{{\mathbb{R}}}
\def\sS{{\mathbb{S}}}
\def\sT{{\mathbb{T}}}
\def\sU{{\mathbb{U}}}
\def\sV{{\mathbb{V}}}
\def\sW{{\mathbb{W}}}
\def\sX{{\mathbb{X}}}
\def\sY{{\mathbb{Y}}}
\def\sZ{{\mathbb{Z}}}
\def\emLambda{{\Lambda}}
\def\emA{{A}}
\def\emB{{B}}
\def\emC{{C}}
\def\emD{{D}}
\def\emE{{E}}
\def\emF{{F}}
\def\emG{{G}}
\def\emH{{H}}
\def\emI{{I}}
\def\emJ{{J}}
\def\emK{{K}}
\def\emL{{L}}
\def\emM{{M}}
\def\emN{{N}}
\def\emO{{O}}
\def\emP{{P}}
\def\emQ{{Q}}
\def\emR{{R}}
\def\emS{{S}}
\def\emT{{T}}
\def\emU{{U}}
\def\emV{{V}}
\def\emW{{W}}
\def\emX{{X}}
\def\emY{{Y}}
\def\emZ{{Z}}
\def\emSigma{{\Sigma}}
\newcommand{\thesischsixetens}[1]{\mathsfit{#1}}
\def\etLambda{{\thesischsixetens{\Lambda}}}
\def\etA{{\thesischsixetens{A}}}
\def\etB{{\thesischsixetens{B}}}
\def\etC{{\thesischsixetens{C}}}
\def\etD{{\thesischsixetens{D}}}
\def\etE{{\thesischsixetens{E}}}
\def\etF{{\thesischsixetens{F}}}
\def\etG{{\thesischsixetens{G}}}
\def\etH{{\thesischsixetens{H}}}
\def\etI{{\thesischsixetens{I}}}
\def\etJ{{\thesischsixetens{J}}}
\def\etK{{\thesischsixetens{K}}}
\def\etL{{\thesischsixetens{L}}}
\def\etM{{\thesischsixetens{M}}}
\def\etN{{\thesischsixetens{N}}}
\def\etO{{\thesischsixetens{O}}}
\def\etP{{\thesischsixetens{P}}}
\def\etQ{{\thesischsixetens{Q}}}
\def\etR{{\thesischsixetens{R}}}
\def\etS{{\thesischsixetens{S}}}
\def\etT{{\thesischsixetens{T}}}
\def\etU{{\thesischsixetens{U}}}
\def\etV{{\thesischsixetens{V}}}
\def\etW{{\thesischsixetens{W}}}
\def\etX{{\thesischsixetens{X}}}
\def\etY{{\thesischsixetens{Y}}}
\def\etZ{{\thesischsixetens{Z}}}
\newcommand{\thesischsixpdata}{p_{\rm{data}}}
\newcommand{\thesischsixptrain}{\hat{p}_{\rm{data}}}
\newcommand{\thesischsixPtrain}{\hat{P}_{\rm{data}}}
\newcommand{\thesischsixpmodel}{p_{\rm{model}}}
\newcommand{\thesischsixPmodel}{P_{\rm{model}}}
\newcommand{\thesischsixptildemodel}{\tilde{p}_{\rm{model}}}
\newcommand{\thesischsixpencode}{p_{\rm{encoder}}}
\newcommand{\thesischsixpdecode}{p_{\rm{decoder}}}
\newcommand{\thesischsixprecons}{p_{\rm{reconstruct}}}
\newcommand{\thesischsixlaplace}{\mathrm{Laplace}} %
\newcommand{\thesischsixE}{\Exp}
\newcommand{\thesischsixLs}{\mathcal{L}}
\newcommand{\thesischsixemp}{\tilde{p}}
\newcommand{\thesischsixlr}{\alpha}
\newcommand{\thesischsixreg}{\lambda}
\newcommand{\thesischsixrect}{\mathrm{rectifier}}
\newcommand{\thesischsixsoftmax}{\mathrm{softmax}}
\newcommand{\thesischsixsigmoid}{\sigma}
\newcommand{\thesischsixsoftplus}{\zeta}
\newcommand{\thesischsixKL}{D_{\mathrm{KL}}}
\newcommand{\thesischsixstandarderror}{\mathrm{SE}}
\newcommand{\thesischsixnormlzero}{L^0}
\newcommand{\thesischsixnormlone}{L^1}
\newcommand{\thesischsixnormltwo}{L^2}
\newcommand{\thesischsixnormlp}{L^p}
\newcommand{\thesischsixnormmax}{L^\infty}
\newcommand{\thesischsixparents}{Pa} %
\providecommand{\argmax}{\operatorname*{arg\,max}}
\providecommand{\argmin}{\operatorname*{arg\,min}}
\providecommand{\sign}{\operatorname{sign}}
\providecommand{\Tr}{\operatorname{Tr}}
\let\ab\allowbreak

\lstdefinestyle{mystyle}{
    language=Python,
    basicstyle=\ttfamily\small,
    keywordstyle=\color{blue},
    stringstyle=\color{red},
    commentstyle=\color{gray},
    numbers=left,
    numberstyle=\tiny,
    stepnumber=1,
    breaklines=true,
    tabsize=4
}
\lstset{style=mystyle}
\providecommand{\prox}{\operatorname*{prox}}
\providecommand{\brox}{\operatorname*{brox}}
\newcommand{\thesischsevenProd}[2]{\prod_{#1}^{#2}(1-\eta_\tau)}
\newcommand{\thesischsevenz}{\zeta}
\newcommand{\thesischsevenlt}{\left}
\newcommand{\thesischsevenal}{\alpha}
\newcommand{\thesischsevenp}{\partial}
\newcommand{\thesischsevenD}{\Delta}
\newcommand{\thesischsevennfr}{\nicefrac}
\newcommand{\thesischsevenfr}{\frac}
\newcommand{\thesischsevendfr}{\dfrac}
\newcommand{\thesischsevenmbf}{\mathbf}
\newcommand{\thesischsevenmds}{\mathds}
\newcommand{\thesischsevenbb}{\mathbb}
\newcommand{\thesischsevenwtd}{\widetilde}
\newcommand{\thesischsevenwhat}{\widehat}
\newcommand{\thesischsevenavg}{\overline}
\newcommand{\thesischsevenwtc}{{\thesischsevenwtd c}}
\newcommand{\thesischsevenwS}{\thesischsevenwtd{S}}
\newcommand{\thesischsevenwA}{\thesischsevenwtd{A}}
\newcommand{\thesischsevensigmt}{\sigma^{t}}
\newcommand{\thesischsevensigmtpo}{\sigma^{t+1}}
\newcommand{\thesischseventtheta}{\tilde{\theta}}
\newcommand{\thesischsevenwgit}{\wg_i^t}
\newcommand{\thesischsevenwgitpo}{\wg_i^{t+1}}
\newcommand{\thesischsevenwgix}{\wg_i(x)}
\newcommand{\thesischsevenwgiy}{\wg_i(y)}
\newcommand{\thesischsevenwgik}{\wg_i(x^k)}
\newcommand{\thesischsevenwgikpo}{\wg_i(x^{k+1})}
\newcommand{\thesischsevenhL}{\hat{L}}
\newcommand{\thesischsevenhG}{\hat{G}}
\newcommand{\thesischsevenhAi}{\hat{A}_i}
\newcommand{\thesischsevenhA}{\hat{A}}
\newcommand{\thesischsevenhB}{\hat{B}}
\newcommand{\thesischsevenhC}{\hat{C}}
\newcommand{\thesischsevenhD}{\hat{D}}
\newcommand{\thesischsevenpA}{A^{\prime}}
\newcommand{\thesischsevenpB}{B^{\prime}}
\newcommand{\thesischsevenpC}{C^{\prime}}
\newcommand{\thesischsevenwpA}{\thesischsevenwA^{\prime}}
\newcommand{\thesischsevenwpB}{\wB^{\prime}}
\newcommand{\thesischsevenwpC}{\wC^{\prime}}
\newcommand{\thesischsevenppA}{A^{\prime \prime}}
\newcommand{\thesischsevenppB}{B^{\prime \prime}}
\newcommand{\thesischsevenppC}{C^{\prime \prime}}
\newcommand{\thesischsevenhg}{\hat{g}}
\newcommand{\thesischsevenhatf}{\hat{f}}
\newcommand{\thesischsevenhgi}{\hat{g}_i}
\newcommand{\thesischsevenhgix}{\hat{g}_i(x)}
\newcommand{\thesischsevenhgiy}{\hat{g}_i(y)}
\newcommand{\thesischsevenhgik}{\hat{g}_i(x^k)}
\newcommand{\thesischsevenhgikpo}{\hat{g}_i(x^{k+1})}
\newcommand{\thesischsevenhgit}{\hat{g}_i(x^t)}
\newcommand{\thesischsevenhgitpo}{\hat{g}_i(x^{t+1})}
\newcommand{\thesischsevenhgu}{\hat{g}^u}
\newcommand{\thesischsevenhs}{\hat{s}}
\newcommand{\thesischseveneik}{e_i^t}
\newcommand{\thesischsevengik}{g_i^t}
\newcommand{\thesischsevengikpo}{g_i^{t+1}}
\newcommand{\thesischseveneikpo}{e_i^{t+1}}
\newcommand{\thesischsevenvik}{v_i^t}
\newcommand{\thesischsevenvikpo}{v_i^{t+1}}
\newcommand{\thesischseveneit}{e_i^t}
\newcommand{\thesischseveneitpo}{e_i^{t+1}}
\newcommand{\thesischsevenLa}{\Lambda}
\newcommand{\thesischsevenlam}{\lambda}
\newcommand{\thesischsevenlamt}{\lambda_t}
\newcommand{\thesischsevenlamtpo}{\lambda_{t+1}}
\newcommand{\thesischsevennut}{\nu_t}
\newcommand{\thesischsevennutpo}{\nu_{t+1}}
\newcommand{\thesischsevenlami}{\thesischsevenlam_i}
\newcommand{\thesischsevenopn}{\operatorname}
\newcommand{\thesischsevenbM}{\thesischsevenmbf M}
\newcommand{\thesischsevenbS}{\thesischsevenmbf S}
\newcommand{\thesischsevenbQ}{\thesischsevenmbf Q}
\newcommand{\thesischsevenbG}{\thesischsevenmbf G}
\newcommand{\thesischsevenbI}{\thesischsevenmbf I}
\newcommand{\thesischsevenbH}{\thesischsevenmbf H}
\newcommand{\thesischsevenbD}{\thesischsevenmbf D}
\newcommand{\thesischsevenbA}{\thesischsevenmbf A}
\newcommand{\thesischsevenbC}{\thesischsevenmbf C}
\newcommand{\thesischsevenbP}{\thesischsevenmbf P}
\newcommand{\thesischsevenbB}{\thesischsevenmbf B}
\newcommand{\thesischsevenbbB}{\thesischsevenbb B}
\newcommand{\thesischsevenbU}{\thesischsevenbb U}
\newcommand{\thesischsevenbLa}{\thesischsevenmbf {\thesischsevenLa}}
\newcommand{\thesischsevenbE}{\thesischsevenmbf E}
\newcommand{\thesischsevenhS}{\hat{S}}
\newcommand{\thesischsevenbL}{\mathbf{L}}
\newcommand{\thesischsevenstepsize}{\gamma}
\newcommand{\thesischsevenlin}[1]{\left\langle #1\right\rangle} %
\newcommand{\thesischsevenlinf}[1]{\left\langle #1\right\rangle_{\operatorname{F}}} %
\newcommand{\thesischsevencCw}{\thesischsevencC_{\omega}}
\newcommand{\thesischsevencCd}{\thesischsevencC_{\delta}}
\newcommand{\thesischsevenCtk}[1]{\thesischsevencC_{Tk}\rb{#1}}
\newcommand{\thesischsevengixk}{g_{i}(x_k)}
\newcommand{\thesischsevenbN}{\mathbf{N}}
\newcommand{\thesischsevenbx}{\bar{x}}
\newcommand{\thesischsevenby}{\bar{y}}
\newcommand{\thesischsevenbarA}{\overline{A}}
\newcommand{\thesischsevenbarB}{\overline{B}}
\newcommand{\thesischsevenbarL}{\overline{L}}
\newcommand{\thesischsevenbarC}{\overline{C}}
\newcommand{\thesischsevenbarD}{\overline{D}}
\newcommand{\thesischsevenbara}{\overline{a}}
\newcommand{\thesischseveneqtext}[1]{\ensuremath{\stackrel{\text{#1}}{=}}}
\newcommand{\thesischsevenletext}[1]{\ensuremath{\stackrel{\text{#1}}{\le}}}
\newcommand{\thesischsevengetext}[1]{\ensuremath{\stackrel{\text{#1}}{\ge}}}
\newcommand{\thesischsevenltext}[1]{\ensuremath{\stackrel{\text{#1}}{<}}}
\newcommand{\thesischsevengtext}[1]{\ensuremath{\stackrel{\text{#1}}{>}}}
\newcommand{\thesischsevensucctext}[1]{\ensuremath{\stackrel{\text{#1}}{\succ}}}
\newcommand{\thesischsevenpreceqtext}[1]{\ensuremath{\stackrel{\text{#1}}{\preceq}}}
\newcommand{\thesischsevengnorm}[1]{\left\| #1 \right\|^{1 + \gamma }}      %
\newcommand{\thesischsevenanorm}[1]{\left\| #1 \right\|^{\thesischsevenal}}      %
\newcommand{\thesischsevenbnorm}[1]{\left\| #1 \right\|^{\beta}}      %
\newcommand{\thesischsevenahnorm}[1]{\left\| #1 \right\|^{\thesischsevenal/2}}      %
\newcommand{\thesischsevenbbone}[1]{\mathbb{1}}      %
\newcommand{\thesischseveng}{\gamma}
\newcommand{\thesischsevennf}[1]{\nabla f\rb{#1}}
\newcommand{\thesischsevenpf}[1]{\thesischsevenp f\rb{#1}}
\newcommand{\thesischsevenpfi}[1]{\thesischsevenp f_i\rb{#1}}
\newcommand{\thesischsevenhotidea}{{\color{red}\bf HOT IDEA: }}
\newcommand{\thesischsevendone}{{\color{blue}\bf DONE: }}
\newcommand{\thesischsevenDt}{{\Delta}^t}
\newcommand{\thesischsevenDtpo}{{\Delta}^{t+1}}
\newcommand{\thesischsevenrmXt}{{\mathrm X}^t}
\newcommand{\thesischsevenrmYt}{{\mathrm Y}^t}
\newcommand{\thesischsevenrmZt}{{\mathrm Z}^t}
\newcommand{\thesischsevenrmXtpo}{{\mathrm X}^{t+1}}
\newcommand{\thesischsevenrmYtpo}{{\mathrm Y}^{t+1}}
\newcommand{\thesischsevenrmZtpo}{{\mathrm Z}^{t+1}}
\newcommand{\thesischsevendel}[1]{}
\newcommand{\thesischsevenfinf}{f^\mathrm{inf}}
\newcommand{\thesischsevenfiinf}{f_i^\mathrm{inf}}
\newcommand{\thesischsevenii}{{}^{(i)}}
\newcommand{\thesischsevenProb}{\Prob}
\newcommand{\thesischsevennc}[2]{{\color{red}\|#1\|_{(#2)}}}
\newcommand{\thesischsevenncs}[2]{\|#1\|^2_{(#2)}}
\newcommand{\thesischsevenncc}[2]{{\color{red}\|#1\|^*_{(#2)}}}
\newcommand{\thesischsevenls}[1]{{\color{red} \mathcal S(#1)}}
\newcommand{\thesischsevenRw}[2]{\mathcal R_{#1}(#2)}
\newcommand{\thesischsevenRws}[2]{\mathcal R^2_{#1}(#2)}
\newcommand{\thesischsevennbp}[2]{\|#1\|_{(#2)}}   %
\newcommand{\thesischsevennbd}[2]{\|#1\|_{(#2)}^*} %
\newcommand{\thesischsevenlf}{\mathcal L}
\newcommand{\thesischsevenU}{U}
\newcommand{\thesischsevenN}{N}
\newcommand{\thesischsevenmLi}{{\color{red}m^{(i)}}}
\newcommand{\thesischsevengLi}{{\color{red}f{(i)}}}
\newcommand{\thesischsevenTLi}[1]{{\color{blue}T^{(#1)}}}
\newcommand{\thesischsevenLip}{L}
\newcommand{\thesischsevenRc}[1]{{\color{red}  \mathbf{RC}_{(#1)}}}
\newcommand{\thesischsevenNRCDM}{{\color{red}NRCDM}\  }
\newcommand{\thesischsevennnz}[1]{{\color{red}\|#1\|_0}}
\providecommand{\MVEE}{\operatorname{MVEE}}       %
\newcommand{\thesischsevenExp}[1]{\Exp\left[#1\right]}
\newcommand{\thesischsevenExpCond}[2]{\Exp\left[#1 | #2\right]}
\newcommand{\thesischsevenExpu}[2]{\Exp_{#1}\left[#2\right]}
\newcommand{\thesischsevenExpSub}[2]{\Exp_{#1}\left[#2\right]}
\newcommand{\thesischsevenlamax}[1]{\thesischsevenlam_{\max}\left(#1\right)}
\newcommand{\thesischsevenlamin}[1]{\thesischsevenlam_{\min}\left(#1\right)}
\newcommand{\thesischsevencb}[1]{\thesischsevenlt\{#1\rt\}}
\newcommand {\thesischsevenrbl}[1]{(#1)} %
\renewcommand{\thesischsevennf}[1]{\nabla f\thesischsevenrbl{#1}}
\newcommand{\thesischsevennfi}[1]{\nabla f_i\thesischsevenrbl{#1}}
\newcommand{\thesischsevennfj}[1]{\nabla f_j\thesischsevenrbl{#1}}
\newcommand{\thesischsevenhf}[1]{\nabla^2 f\thesischsevenrbl{#1}}
\newcommand{\thesischsevenhfi}[1]{\nabla^2 f_i\thesischsevenrbl{#1}}
\newcommand{\thesischsevenvtp}{v^{t-1}}
\newcommand{\thesischsevenvtpo}{v^{t+1}}
\newcommand{\thesischsevenyijt}{y_{ij}^t}
\newcommand{\thesischsevenyilt}{y_{il}^t}
\newcommand{\thesischsevenyijtpo}{y_{ij}^{t+1}}
\newcommand{\thesischsevenyiltpo}{y_{il}^{t+1}}
\newcommand{\thesischsevenvitp}{v_i^{t-1}}
\newcommand{\thesischsevenvit}{v_i^{t}}
\newcommand{\thesischsevenvitpo}{v_i^{t+1}}
\newcommand{\thesischsevengit}{g_i^t}
\newcommand{\thesischsevengitpo}{g_i^{t+1}}
\newcommand{\thesischsevengitpm}{g_i^{t-1}}
\newcommand{\thesischsevendgit}{\thesischsevenvit}
\newcommand{\thesischsevendgitpo}{\thesischsevenvitpo}
\newcommand{\thesischsevendgitpm}{\thesischsevengitpm}
\newcommand{\thesischsevenxk}{x^{k}}
\newcommand{\thesischsevenxt}{x^{t}}
\newcommand{\thesischsevensit}{s_{i}^{t}}
\newcommand{\thesischsevenyt}{y^{t}}
\newcommand{\thesischsevenxtpo}{x^{t+1}}
\newcommand{\thesischsevenytpo}{y^{t+1}}
\newcommand{\thesischsevenwt}{W^{t}}
\newcommand{\thesischsevenws}{W^{*}}
\newcommand{\thesischsevenwtpo}{W^{t+1}}
\newcommand{\thesischsevenetpo}{e^{t+1}}
\newcommand{\thesischsevenxitpo}{\xi^{t+1}}
\newcommand{\thesischsevenxit}{\xi^{t}}
\newcommand{\thesischsevenxiit}{\xi_{i}^{t}}
\newcommand{\thesischsevenxtmo}{x^{t-1}}
\newcommand{\thesischsevenwtdx}{\thesischsevenwtd{x}}
\newcommand{\thesischsevenwtdxk}{\thesischsevenwtd{x}^k}
\newcommand{\thesischsevenwtdxkpo}{\thesischsevenwtd{x}^{k+1}}
\newcommand{\thesischsevenomtpo}{\omega^{t+1}}
\newcommand{\thesischsevenomt}{\omega^{t}}
\newcommand{\thesischsevenomi}{\omega_i}
\newcommand{\thesischsevenom}{\omega}
\newcommand{\thesischsevenxkpo}{x^{k+1}}
\newcommand{\thesischsevenfxk}{f(\thesischsevenxk)}
\newcommand{\thesischsevengxk}{g(\thesischsevenxk)}
\newcommand{\thesischsevenfxkpo}{f(\thesischsevenxkpo)}
\newcommand{\thesischsevendz}{\delta_0}
\newcommand{\thesischsevenxs}{x^{*}}
\newcommand{\thesischsevenfxt}{f(\thesischsevenxt)}
\newcommand{\thesischsevenfxs}{f(x^*)}
\newcommand{\thesischsevenfxtpo}{f(\thesischsevenxtpo)}
\newcommand{\thesischsevenhfxk}{\thesischsevenhf{\thesischsevenxk}}
\newcommand{\thesischsevennfxk}{\thesischsevennf{\thesischsevenxk}}
\newcommand{\thesischsevennfxkpo}{\thesischsevennf{x^{k+1}}}
\newcommand{\thesischsevennfxt}{\thesischsevennf{\thesischsevenxt}}
\newcommand{\thesischsevennfxtpo}{\thesischsevennf{\thesischsevenxtpo}}
\newcommand{\thesischsevennfomt}{\thesischsevennf{\thesischsevenomt}}
\newcommand{\thesischsevennfomtpo}{\thesischsevennf{\thesischsevenomtpo}}
\newcommand{\thesischsevennfxs}{\thesischsevennf{\thesischsevenxs}}
\newcommand{\thesischsevennfx}{\thesischsevennf{x}}
\newcommand{\thesischsevennfy}{\thesischsevennf{y}}
\newcommand{\thesischsevenhfx}{\thesischsevenhf{x}}
\newcommand{\thesischsevenhfy}{\thesischsevenhf{y}}
\newcommand{\thesischsevennfxijk}{\nabla f_{\xi_{ij}^k} (x^k)}
\newcommand{\thesischsevennfxijkpo}{\nabla f_{\xi_{ij}^k}(x^{k+1})}
\newcommand{\thesischsevennfixk}{\thesischsevennfi{\thesischsevenxk}}
\newcommand{\thesischsevennfixt}{\thesischsevennfi{\thesischsevenxt}}
\newcommand{\thesischsevennfix}{\thesischsevennfi{x}}
\newcommand{\thesischsevennfiy}{\thesischsevennfi{y}}
\newcommand{\thesischsevennfij}[1]{\nabla f_{ij}\thesischsevenrbl{#1}}
\newcommand{\thesischsevennfil}[1]{\nabla f_{il}\thesischsevenrbl{#1}}
\newcommand{\thesischsevennfijxt}{\thesischsevennfij{\thesischsevenxt}}
\newcommand{\thesischsevennfijxtpo}{\thesischsevennfij{\thesischsevenxtpo}}
\newcommand{\thesischsevennfjxtpo}{\thesischsevennfj{\thesischsevenxtpo}}
\newcommand{\thesischsevennfjxt}{\thesischsevennfj{\thesischsevenxt}}
\newcommand{\thesischsevennfjomtpo}{\thesischsevennfj{\thesischsevenomtpo}}
\newcommand{\thesischsevennfijomiit}{\thesischsevennfij{\thesischsevenomiit}}
\newcommand{\thesischsevennfiomiit}{\thesischsevennfi{\thesischsevenomiit}}
\newcommand{\thesischsevennfijomiitpo}{\thesischsevennfij{\thesischsevenomiitpo}}
\newcommand{\thesischsevennfilomiit}{\thesischsevennfil{\thesischsevenomiit}}
\newcommand{\thesischsevennfilomiitpo}{\thesischsevennfil{\thesischsevenomiitpo}}
\newcommand{\thesischsevennfilxt}{\thesischsevennfil{\thesischsevenxt}}
\newcommand{\thesischsevennfilxtpo}{\thesischsevennfil{\thesischsevenxtpo}}
\newcommand{\thesischsevennfijx}{\nabla f_{ij}(x)}
\newcommand{\thesischsevennfijy}{\nabla f_{ij}(y)}
\newcommand{\thesischsevendgt}{\Delta g^{t}}
\newcommand{\thesischsevendgtpo}{\Delta g^{t+1}}
\newcommand{\thesischsevengt}{ g^{t}}
\newcommand{\thesischsevengtpo}{ g^{t+1}}
\newcommand{\thesischsevenopern}{\operatorname}
\newcommand{\thesischsevenomiitpo}{w_i^{t+1}}
\newcommand{\thesischsevenomijt}{w_{ij}^{t}}
\newcommand{\thesischsevenomilt}{w_{il}^{t}}
\newcommand{\thesischsevenomijtpo}{w_{ij}^{t+1}}
\newcommand{\thesischsevenomiltpo}{w_{il}^{t+1}}
\newcommand{\thesischsevenomiit}{w_i^{t}}
\newcommand{\thesischsevennfiomtpo}{\thesischsevennfi{\omitpo}}
\newcommand{\thesischsevennfitpoxtpo}{\nabla f_{i_{t+1}}(x^{t+1})}
\newcommand{\thesischsevennfixs}{\thesischsevennfi{\thesischsevenxs}}
\newcommand{\thesischsevennfjxk}{\thesischsevennfj{\thesischsevenxk}}
\newcommand{\thesischsevennfjxs}{\thesischsevennfj{\thesischsevenxs}}
\newcommand{\thesischsevennfixkpo}{\thesischsevennfi{x^{k+1}}}
\newcommand{\thesischsevennfixtpo}{\thesischsevennfi{x^{t+1}}}
\newcommand{\thesischsevennfitxt}{\nabla f_{i_{t}}(x^{t})}
\newcommand{\thesischsevennfitomt}{\nabla f_{i_{t}}(\thesischsevenomt)}
\newcommand{\thesischsevennfitpoomtpo}{\nabla f_{i_{t+1}}(\thesischsevenomtpo)}
\newcommand{\thesischsevennfjxkpo}{\thesischsevennfj{x^{k+1}}}
\newcommand{\thesischsevenhik}{h^{k}_{i}}
\newcommand{\thesischsevenhikpo}{h^{k+1}_{i}}
\newcommand{\thesischsevenhit}{h^{t}_{i}}
\newcommand{\thesischsevenhitpo}{h^{t+1}_{i}}
\newcommand{\thesischsevenali}{\thesischsevenal_{i}}
\newcommand{\thesischsevenmik}{m^{k}_{i}}
\newcommand{\thesischsevenhk}{h^{k}}
\newcommand{\thesischsevennmo}{\thesischsevenfr{1}{n}}
\newcommand{\thesischsevennmt}{\thesischsevenfr{1}{n^2}}
\newcommand{\thesischseventnmo}{\thesischsevenfr{2}{n}}
\newcommand{\thesischseventnmt}{\thesischsevenfr{2}{n^2}}
\newcommand{\thesischsevenub}{\underbrace}
\newcommand {\thesischsevenmaxin}{\max_{i=1,\dots, n}}
\newcommand {\thesischsevenmaxi}{\max_{i}}
\newcommand {\thesischsevenminin}{\min_{i=1,\dots, n}}
\newcommand {\thesischsevenmini}{\min_{i}}
\newcommand {\thesischsevensumin}{\sum_{i=1}^{n}}
\newcommand {\thesischsevensumtt}{\sum_{t=0}^{T}}
\newcommand {\thesischsevensumlmi}{\thesischsevenfr{1}{m_i}\sum_{l=1}^{m_i}}
\newcommand {\thesischsevensuminn}{\thesischsevenfr{1}{n}\sum_{i=1}^{n}}
\newcommand {\thesischsevensumjm}{\sum_{j=1}^{m}}
\newcommand {\thesischsevensumjmi}{\sum_{j=1}^{m_i}}
\newcommand {\thesischsevensumjmm}{\thesischsevenfr{1}{m}\sum_{j=1}^{m}}
\newcommand {\thesischsevensumim}{\sum_{i=1}^{m}}
\newcommand {\thesischsevensumimm}{\thesischsevenfr{1}{m}\sum_{i=1}^{m}}
\newcommand {\thesischsevensumjnn}{\thesischsevenfr{1}{n}\sum_{j=1}^{n}}
\newcommand {\thesischsevensumjmmi}{\thesischsevenfr{1}{m_i}\sum_{j=1}^{m_i}}
\newcommand {\thesischsevensumjn}{\sum_{j=1}^{n}}
\newcommand {\thesischsevensumis}{\sum_{i \in S}}
\newcommand {\thesischsevenwit}{w_{i}^{t}}
\newcommand {\thesischsevenwik}{w_{i}^{k}}
\newcommand {\thesischsevenwjk}{w_{j}^{k}}
\newcommand {\thesischsevenwk}{w^{k}}
\newcommand {\thesischsevenwitpo}{w_{i}^{t+1}}
\newcommand {\thesischsevenwikpo}{w_{i}^{k+1}}
\newcommand {\thesischsevenwjkpo}{w_{j}^{k+1}}
\newcommand {\thesischsevensumid}{\sum_{i=1}^{d}}
\newcommand {\thesischsevenghx}{\hat{g}(x)}
\newcommand {\thesischsevengh}{\hat{g}}
\newcommand {\thesischsevenhx}{\hat{x}}
\newcommand {\thesischsevengih}{\hat_i}
\newcommand {\thesischsevengihx}{\hat{g_i}(x)}
\newcommand{\thesischsevencmark}{{\color{PineGreen}\ding{51}}}%
\newcommand{\thesischsevenxmark}{{\color{BrickRed}\ding{55}}}%
\providecommand{\tr}{\operatorname{Tr}}
\newcommand{\thesischseventrace}[1]{\tr\left(#1\right)}
\newcommand{\thesischsevenrbtop}[1]{\left(#1\right)^{\top}}
\newcommand{\thesischsevencA}{{\cal A}}
\newcommand{\thesischsevencB}{{\cal B}}
\newcommand{\thesischsevencC}{{\cal C}}
\newcommand{\thesischsevencD}{{\cal D}}
\newcommand{\thesischsevencE}{{\cal E}}
\newcommand{\thesischsevencF}{{\cal F}}
\newcommand{\thesischsevencG}{{\cal G}}
\newcommand{\thesischsevencH}{{\cal H}}
\newcommand{\thesischsevencJ}{{\cal J}}
\newcommand{\thesischsevencK}{{\cal K}}
\newcommand{\thesischsevencL}{{\cal L}}
\newcommand{\thesischsevencM}{{\cal M}}
\newcommand{\thesischsevencN}{{\cal N}}
\newcommand{\thesischsevencP}{{\cal P}}
\newcommand{\thesischsevencQ}{{\cal Q}}
\newcommand{\thesischsevencR}{{\cal R}}
\newcommand{\thesischsevencS}{{\cal S}}
\newcommand{\thesischsevencT}{{\cal T}}
\newcommand{\thesischsevencU}{{\cal U}}
\newcommand{\thesischsevencV}{{\cal V}}
\newcommand{\thesischsevencX}{{\cal X}}
\newcommand{\thesischsevencY}{{\cal Y}}
\newcommand{\thesischsevencW}{{\cal W}}
\newcommand{\thesischsevencZ}{{\cal Z}}
\newcommand*{\thesischsevenqedb}{\hfill\ensuremath{\square}}   %
\newcommand{\thesischseventopic}[1]{\vspace{2mm}\noindent{{\bf #1:}}}
\newcommand{\thesischsevenhead}[1]{\noindent{{\bf #1:}}}
\definecolor{mydarkblue}{RGB}{39,47,180}
\definecolor{mydarkpurple}{RGB}{48,25,52}
\newcommand{\thesischsevengreen}{\color{mydarkgreen}}
\newcommand{\thesischsevenred}{\color{mydarkred}}
\newcommand{\thesischsevenblue}{\color{mydarkblue}}
\newcommand{\thesischsevenpurple}{\color{mydarkpurple}}
\newcommand{\thesischsevenzerosarah}{\algname{ZeroSARAH}}
\newcommand{\thesischsevenclip}{\operatorname{clip}}
\definecolor{lightblue}{RGB}{224, 240, 255} %
\newenvironment{proofof}[1]{\smallskip\noindent{\bf \em Proof of #1.}}{\par}
\newcommand{\thesischsevenrA}{{\thesischsevenred A}}
\newcommand{\thesischsevenrB}{{\thesischsevenred B}}
\newcommand{\thesischsevenxp}{x_p^*}
\newcommand{\thesischsevendist}{\mathrm{dist}}
\providecommand{\proj}{\operatorname*{proj}}
\newcommand\thesischsevenpbr[1]{\left \{ #1 \right \} }
\newcommand\thesischsevenbr[1]{\left ( #1 \right )}
\definecolor{LightCyan}{rgb}{0.88,1,1}
\newcommand{\thesischsevenfigleft}{{\em (Left)}}
\newcommand{\thesischsevenfigcenter}{{\em (Center)}}
\newcommand{\thesischsevenfigright}{{\em (Right)}}
\newcommand{\thesischsevenfigtop}{{\em (Top)}}
\newcommand{\thesischsevenfigbottom}{{\em (Bottom)}}
\newcommand{\thesischsevencaptiona}{{\em (a)}}
\newcommand{\thesischsevencaptionb}{{\em (b)}}
\newcommand{\thesischsevencaptionc}{{\em (c)}}
\newcommand{\thesischsevencaptiond}{{\em (d)}}
\newcommand{\thesischsevennewterm}[1]{{\bf #1}}
\newcommand{\thesischseventrain}{\mathcal{D}}
\newcommand{\thesischsevenvalid}{\mathcal{D_{\mathrm{valid}}}}
\newcommand{\thesischseventest}{\mathcal{D_{\mathrm{test}}}}
\def\ri{{\textnormal{i}}}
\def\mI{{\bm{I}}}
\def\mL{{\bm{L}}}
\newcommand{\thesischseventens}[1]{\bm{\mathsfit{#1}}}
\def\tA{{\thesischseventens{A}}}
\def\tB{{\thesischseventens{B}}}
\def\tC{{\thesischseventens{C}}}
\def\tD{{\thesischseventens{D}}}
\def\tE{{\thesischseventens{E}}}
\def\tF{{\thesischseventens{F}}}
\def\tG{{\thesischseventens{G}}}
\def\tH{{\thesischseventens{H}}}
\def\tI{{\thesischseventens{I}}}
\def\tJ{{\thesischseventens{J}}}
\def\tK{{\thesischseventens{K}}}
\def\tL{{\thesischseventens{L}}}
\def\tM{{\thesischseventens{M}}}
\def\tN{{\thesischseventens{N}}}
\def\tO{{\thesischseventens{O}}}
\def\tP{{\thesischseventens{P}}}
\def\tQ{{\thesischseventens{Q}}}
\def\tR{{\thesischseventens{R}}}
\def\tS{{\thesischseventens{S}}}
\def\tT{{\thesischseventens{T}}}
\def\tU{{\thesischseventens{U}}}
\def\tV{{\thesischseventens{V}}}
\def\tW{{\thesischseventens{W}}}
\def\tX{{\thesischseventens{X}}}
\def\tY{{\thesischseventens{Y}}}
\def\tZ{{\thesischseventens{Z}}}
\newcommand{\thesischsevenetens}[1]{\mathsfit{#1}}
\def\etLambda{{\thesischsevenetens{\Lambda}}}
\def\etA{{\thesischsevenetens{A}}}
\def\etB{{\thesischsevenetens{B}}}
\def\etC{{\thesischsevenetens{C}}}
\def\etD{{\thesischsevenetens{D}}}
\def\etE{{\thesischsevenetens{E}}}
\def\etF{{\thesischsevenetens{F}}}
\def\etG{{\thesischsevenetens{G}}}
\def\etH{{\thesischsevenetens{H}}}
\def\etI{{\thesischsevenetens{I}}}
\def\etJ{{\thesischsevenetens{J}}}
\def\etK{{\thesischsevenetens{K}}}
\def\etL{{\thesischsevenetens{L}}}
\def\etM{{\thesischsevenetens{M}}}
\def\etN{{\thesischsevenetens{N}}}
\def\etO{{\thesischsevenetens{O}}}
\def\etP{{\thesischsevenetens{P}}}
\def\etQ{{\thesischsevenetens{Q}}}
\def\etR{{\thesischsevenetens{R}}}
\def\etS{{\thesischsevenetens{S}}}
\def\etT{{\thesischsevenetens{T}}}
\def\etU{{\thesischsevenetens{U}}}
\def\etV{{\thesischsevenetens{V}}}
\def\etW{{\thesischsevenetens{W}}}
\def\etX{{\thesischsevenetens{X}}}
\def\etY{{\thesischsevenetens{Y}}}
\def\etZ{{\thesischsevenetens{Z}}}
\newcommand{\thesischsevenpdata}{p_{\rm{data}}}
\newcommand{\thesischsevenptrain}{\hat{p}_{\rm{data}}}
\newcommand{\thesischsevenPtrain}{\hat{P}_{\rm{data}}}
\newcommand{\thesischsevenpmodel}{p_{\rm{model}}}
\newcommand{\thesischsevenPmodel}{P_{\rm{model}}}
\newcommand{\thesischsevenptildemodel}{\tilde{p}_{\rm{model}}}
\newcommand{\thesischsevenpencode}{p_{\rm{encoder}}}
\newcommand{\thesischsevenpdecode}{p_{\rm{decoder}}}
\newcommand{\thesischsevenprecons}{p_{\rm{reconstruct}}}
\newcommand{\thesischsevenlaplace}{\mathrm{Laplace}} %
\newcommand{\thesischsevenLs}{\mathcal{L}}
\newcommand{\thesischsevenemp}{\tilde{p}}
\newcommand{\thesischsevenlr}{\alpha}
\newcommand{\thesischsevenreg}{\lambda}
\newcommand{\thesischsevenrect}{\mathrm{rectifier}}
\newcommand{\thesischsevensoftmax}{\mathrm{softmax}}
\newcommand{\thesischsevensigmoid}{\sigma}
\newcommand{\thesischsevensoftplus}{\zeta}
\newcommand{\thesischsevenKL}{D_{\mathrm{KL}}}
\newcommand{\thesischsevenstandarderror}{\mathrm{SE}}
\newcommand{\thesischsevennormlzero}{L^0}
\newcommand{\thesischsevennormlone}{L^1}
\newcommand{\thesischsevennormltwo}{L^2}
\newcommand{\thesischsevennormlp}{L^p}
\newcommand{\thesischsevennormmax}{L^\infty}
\newcommand{\thesischsevenparents}{Pa} %

\newcommand{\thesischeightg}{{\color{red}g}}
\newcommand{\thesischeighth}{{\color{blue}h}}
\newcommand{\thesischeightClip}[2]{{\rm Clip}_{#1}\left( #2 \right)}
\newcommand{\thesischeightNormalize}[2]{{\rm Normalize}_{#1}\left( #2 \right)}
\newcommand{\thesischeightubar}[1]{\underaccent{\bar}{#1}}
\newcommand{\thesischeightlin}[1]{\left\langle #1\right\rangle} %
\newcommand{\thesischeightcA}{\mathcal{A}}
\newcommand{\thesischeightcB}{\mathcal{B}}
\newcommand{\thesischeightcC}{\mathcal{C}}
\newcommand{\thesischeightcD}{\mathcal{D}}
\newcommand{\thesischeightcF}{\mathcal{F}}
\newcommand{\thesischeightcH}{\mathcal{H}}
\newcommand{\thesischeightcL}{\mathcal{L}}
\newcommand{\thesischeightcT}{\mathcal{T}}
\newcommand{\thesischeightcN}{\mathcal{N}}
\newcommand{\thesischeightcQ}{\mathcal{Q}}
\newcommand{\thesischeightcS}{\mathcal{S}}
\newcommand{\thesischeightcU}{\mathcal{U}}
\newcommand{\thesischeightcW}{\mathcal{W}}
\newcommand{\thesischeightcY}{\mathcal{Y}}
\newcommand{\thesischeightcZ}{\mathcal{Z}}
\newcommand{\thesischeightmI}{\mathbf{I}}
\newcommand{\thesischeightdotparagraph}{\paragraph{$\bullet$}}
\newcommand{\thesischeighthotidea}{{\color{red}\bf HOT IDEA: }}
\newcommand{\thesischeightdone}{{\color{blue}\bf DONE: }}
\newcommand{\thesischeightdel}[1]{}
\definecolor{junglegreen}{rgb}{0.16, 0.67, 0.53}
\definecolor{lasallegreen}{rgb}{0.03, 0.47, 0.19}
\newcommand{\thesischeightU}{\mathbb{U}} %
\newcommand{\thesischeightPermComp}{\thesischeightU\mathbb{P}} %
\newcommand{\thesischeightst}{\;:\;} %
\newcommand\thesischeightbr[1]{\left ( #1 \right )}
\newcommand\thesischeightpbr[1]{\left \{ #1 \right \} }
\newcommand{\thesischeightProb}{\Prob} %
\newcommand{\thesischeightProbCond}[2]{\thesischeightProb\left(\left.#1\right\vert#2\right)}
\newcommand{\thesischeightnnz}[1]{{\color{red}\|#1\|_0}}
\newcommand{\thesischeightExp}[1]{\Exp\left[#1\right]}
\newcommand{\thesischeightExpCond}[2]{\Exp\left[\left.#1\right\vert#2\right]}
\newcommand{\thesischeightExpSub}[2]{\Exp_{#1}\left[#2\right]}
\newcommand{\thesischeightSample}{\mathcal{S}}
\newcommand{\thesischeightB}{\mathbb{B}}
\definecolor{red}{rgb}{1.0, 0.01, 0.24}
\definecolor{orange}{rgb}{1.0, 0.43, 0.29}
\newcommand{\thesischeightfigleft}{{\em (Left)}}
\newcommand{\thesischeightfigcenter}{{\em (Center)}}
\newcommand{\thesischeightfigright}{{\em (Right)}}
\newcommand{\thesischeightfigtop}{{\em (Top)}}
\newcommand{\thesischeightfigbottom}{{\em (Bottom)}}
\newcommand{\thesischeightcaptiona}{{\em (a)}}
\newcommand{\thesischeightcaptionb}{{\em (b)}}
\newcommand{\thesischeightcaptionc}{{\em (c)}}
\newcommand{\thesischeightcaptiond}{{\em (d)}}
\newcommand{\thesischeightnewterm}[1]{{\bf #1}}
\newcommand{\thesischeighttrain}{\mathcal{D}}
\newcommand{\thesischeightvalid}{\mathcal{D_{\mathrm{valid}}}}
\newcommand{\thesischeighttest}{\mathcal{D_{\mathrm{test}}}}
\newcommand{\thesischeighttens}[1]{\bm{\mathsfit{#1}}}
\def\tA{{\thesischeighttens{A}}}
\def\tB{{\thesischeighttens{B}}}
\def\tC{{\thesischeighttens{C}}}
\def\tD{{\thesischeighttens{D}}}
\def\tE{{\thesischeighttens{E}}}
\def\tF{{\thesischeighttens{F}}}
\def\tG{{\thesischeighttens{G}}}
\def\tH{{\thesischeighttens{H}}}
\def\tI{{\thesischeighttens{I}}}
\def\tJ{{\thesischeighttens{J}}}
\def\tK{{\thesischeighttens{K}}}
\def\tL{{\thesischeighttens{L}}}
\def\tM{{\thesischeighttens{M}}}
\def\tN{{\thesischeighttens{N}}}
\def\tO{{\thesischeighttens{O}}}
\def\tP{{\thesischeighttens{P}}}
\def\tQ{{\thesischeighttens{Q}}}
\def\tR{{\thesischeighttens{R}}}
\def\tS{{\thesischeighttens{S}}}
\def\tT{{\thesischeighttens{T}}}
\def\tU{{\thesischeighttens{U}}}
\def\tV{{\thesischeighttens{V}}}
\def\tW{{\thesischeighttens{W}}}
\def\tX{{\thesischeighttens{X}}}
\def\tY{{\thesischeighttens{Y}}}
\def\tZ{{\thesischeighttens{Z}}}
\newcommand{\thesischeightetens}[1]{\mathsfit{#1}}
\def\etLambda{{\thesischeightetens{\Lambda}}}
\def\etA{{\thesischeightetens{A}}}
\def\etB{{\thesischeightetens{B}}}
\def\etC{{\thesischeightetens{C}}}
\def\etD{{\thesischeightetens{D}}}
\def\etE{{\thesischeightetens{E}}}
\def\etF{{\thesischeightetens{F}}}
\def\etG{{\thesischeightetens{G}}}
\def\etH{{\thesischeightetens{H}}}
\def\etI{{\thesischeightetens{I}}}
\def\etJ{{\thesischeightetens{J}}}
\def\etK{{\thesischeightetens{K}}}
\def\etL{{\thesischeightetens{L}}}
\def\etM{{\thesischeightetens{M}}}
\def\etN{{\thesischeightetens{N}}}
\def\etO{{\thesischeightetens{O}}}
\def\etP{{\thesischeightetens{P}}}
\def\etQ{{\thesischeightetens{Q}}}
\def\etR{{\thesischeightetens{R}}}
\def\etS{{\thesischeightetens{S}}}
\def\etT{{\thesischeightetens{T}}}
\def\etU{{\thesischeightetens{U}}}
\def\etV{{\thesischeightetens{V}}}
\def\etW{{\thesischeightetens{W}}}
\def\etX{{\thesischeightetens{X}}}
\def\etY{{\thesischeightetens{Y}}}
\def\etZ{{\thesischeightetens{Z}}}
\newcommand{\thesischeightpdata}{p_{\rm{data}}}
\newcommand{\thesischeightptrain}{\hat{p}_{\rm{data}}}
\newcommand{\thesischeightPtrain}{\hat{P}_{\rm{data}}}
\newcommand{\thesischeightpmodel}{p_{\rm{model}}}
\newcommand{\thesischeightPmodel}{P_{\rm{model}}}
\newcommand{\thesischeightptildemodel}{\tilde{p}_{\rm{model}}}
\newcommand{\thesischeightpencode}{p_{\rm{encoder}}}
\newcommand{\thesischeightpdecode}{p_{\rm{decoder}}}
\newcommand{\thesischeightprecons}{p_{\rm{reconstruct}}}
\newcommand{\thesischeightlaplace}{\mathrm{Laplace}} %
\newcommand{\thesischeightE}{\Exp}
\newcommand{\thesischeightLs}{\mathcal{L}}
\newcommand{\thesischeightemp}{\tilde{p}}
\newcommand{\thesischeightlr}{\alpha}
\newcommand{\thesischeightreg}{\lambda}
\newcommand{\thesischeightrect}{\mathrm{rectifier}}
\newcommand{\thesischeightsoftmax}{\mathrm{softmax}}
\newcommand{\thesischeightsigmoid}{\sigma}
\newcommand{\thesischeightsoftplus}{\zeta}
\newcommand{\thesischeightKL}{D_{\mathrm{KL}}}
\newcommand{\thesischeightstandarderror}{\mathrm{SE}}
\newcommand{\thesischeightnormlzero}{L^0}
\newcommand{\thesischeightnormlone}{L^1}
\newcommand{\thesischeightnormltwo}{L^2}
\newcommand{\thesischeightnormlp}{L^p}
\newcommand{\thesischeightnormmax}{L^\infty}
\newcommand{\thesischeightparents}{Pa} %

\newcommand{\thesischaptersetup}[1]{\csname thesischaptersetup\romannumeral#1\endcsname}
\newcommand{\thesischaptersetupii}{%
  \renewcommand*{\backrefalt}[4]{%
      \ifcase ##1 \footnotesize{(Not cited.)}%
      \or        \footnotesize{(Cited on page~##2)}%
      \else      \footnotesize{(Cited on pages~##2)}%
      \fi}
  \renewcommand{\algorithmicrequire}{\textbf{Input:}}
  \renewcommand{\algorithmicensure}{\textbf{Output:}}
  \renewcommand{\Re}{\mathrm{Re}}
  \renewcommand{\Im}{\mathrm{Im}}
}
\newcommand{\thesischaptersetupiii}{%
  \renewcommand{\algorithmicrequire}{\textbf{Input:}}
  \renewcommand{\algorithmicensure}{\textbf{Output:}}
  \renewcommand{\Re}{\mathrm{Re}}
  \renewcommand{\Im}{\mathrm{Im}}
}
\newcommand{\thesischaptersetupiv}{%
  \renewcommand{\bibname}{References}
  \renewcommand{\bibsection}{\section*{\bibname}}
}
\newcommand{\thesischaptersetupv}{%
  \definecolor{blue}{RGB}{0,0,255}
  \everydisplay\expandafter{\the\everydisplay\setcounter{thesischfiverelctr}{0}}%
  \def\<{\left\langle}
  \def\>{\right\rangle}
  \def\({\left(}
  \def\){\right)}
  \renewcommand*{\deg}{\mathop{\mathrm{deg}}}%
  \renewcommand*{\det}{\mathop{\mathrm{det}}}%
  \renewcommand*{\O}{\mathop{\mathrm{O}}}%
  \renewcommand*{\dim}{\mathop{\mathrm{dim}}}%
  \renewcommand*{\P}{\Prob}%
  \renewcommand*{\H}{\mathbb{H}}%
  \renewcommand*{\L}{\mathbb{L}}%
  \renewcommand*{\S}{\mathbb{S}}%
  \renewcommand*{\Bar}{\mathop{\mathrm{Bar}}}%
  \renewcommand*{\Re}{\mathop{\mathrm{Re}}}%
  \renewcommand*{\Im}{\mathop{\mathrm{Im}}}%
  \renewcommand*{\lim}{\mathop{\mathrm{lim}}\limits}%
  \renewcommand*{\o}{\mathop{\mathrm{o}}\limits}%
  \renewcommand*{\O}{\mathop{\mathrm{O}}\limits}%
  \renewcommand*{\eq}{\mathop{\sim}\limits}%
}
\newcommand{\thesischaptersetupvi}{%
  \renewcommand*{\backrefalt}[4]{}
  \def\eqref##1{(\ref{##1})}
}
\newcommand{\thesischaptersetupvii}{%
  \renewcommand{\bibname}{References}
  \renewcommand{\bibsection}{\subsubsection*{\bibname}}
  \restylefloat{table}
  \def\<{\left\langle}
  \def\>{\right\rangle}
  \def\({\left(}
  \def\){\right)}
  \renewcommand{\sb}[1]{\left[##1\right]}
  \renewcommand{\st}{s^{t}}
  \floatstyle{plaintop}
  \restylefloat{table}
  \def\eqref##1{equation~\ref{##1}}
}
\newcommand{\thesischaptersetupviii}{%
  \def\eqref##1{(\ref{##1})}
  \definecolor{blue}{named}{black}
}

\part{Foundational Bounds and Distributed Methods under Heavy-Tailed Noise}

\chapter{High-Probability Bounds for Stochastic Optimization and Variational Inequalities: the Case of Unbounded Variance}
\label{ch:chapter_2}
\thesischaptersourcenote
  {High-Probability Bounds for Stochastic Optimization and Variational Inequalities: The Case of Unbounded Variance}
  {sadiev2023high}
  \begingroup
  \thesischaptersetup{2}%
  
\begin{paperabstract}
During recent years the interest of optimization and machine learning communities in high-probability convergence of stochastic optimization methods has been growing. One of the main reasons for this is that high-probability complexity bounds are more accurate and less studied than in-expectation ones. However, state-of-the-art high-probability non-asymptotic convergence results are derived under strong assumptions such as the boundedness of the gradient noise variance or of the objective's gradient itself. This chapter develops several algorithms with high-probability convergence results under less restrictive assumptions. In particular, we derive new high-probability convergence results under the assumption that the gradient/operator noise has bounded central ${p}$-th moment for ${p} \in (1,2]$ in the following setups: (i) smooth non-convex / Polyak-{\L}ojasiewicz / convex / strongly convex / quasi-strongly convex  minimization problems, (ii) Lipschitz / star-cocoercive and monotone / quasi-strongly monotone variational inequalities. These results justify the usage of the considered methods for solving problems that do not fit standard functional classes studied in stochastic optimization.
\end{paperabstract}

\section{Introduction}\label{ch2:sec:intro_chapter_2}

Training of machine learning models is usually performed via stochastic first-order optimization methods, e.g., Stochastic Gradient Descent (\algname{SGD}) \parencite{robbins1951stochastic}
\begin{equation}
    x^{k+1} = x^k - \gamma \nabla f(x^k;\xi^k), \label{ch2:eq:SGD_update_chapter_2}
\end{equation}
where $\nabla f(x^k;\xi^k)$ represents the stochastic gradient of the objective/loss function $f$ at point $x^k$. Despite numerous empirical studies and observations validating the good performance of such methods, it is also important for the field to understand their theoretical convergence properties, e.g., under what assumptions a method converges and what the rate is. However, since the methods of interest are stochastic, one needs to specify what type of convergence is considered before moving on to further questions.

Typically, the convergence of the stochastic methods is studied only in expectation, i.e., for some performance metric\footnote{Examples of performance metrics for minimization of function $f$: $\thesischtwocP(x) = f(x) - f(x^*)$, $\thesischtwocP(x) = \|\nabla f(x)\|^2$, $\thesischtwocP(x) = \|x - x^*\|^2$, where $x^* \in \arg\min_{x\in\R^d} f(x)$.} $\thesischtwocP(x)$, upper bounds are derived for the number of iterations $K$ needed to achieve $\Exp[\thesischtwocP(x^K)] \leq \varepsilon$, where $x^K$ is the output of the method after $K$ steps, $\varepsilon$ is an optimization error, and $\Exp[\cdot]$ is the full expectation. These bounds can be ``blind'' to some important properties like light-/heavy-tailedness of the noise distribution and, as a result, such guarantees do not accurately describe the methods' convergence in practice \parencite{gorbunov2020stochastic}. In contrast, high-probability convergence guarantees are more sensitive to the noise distribution and thus are more accurate. Such results provide upper bounds for the number of iterations $K$ needed to achieve $\Prob\{\thesischtwocP(x^K) \leq \varepsilon\} \geq 1 - \beta$ for some confidence level $\beta \in (0,1]$, where $\Prob\{\cdot\}$ denotes some probability measure determined by a setup.

With the ultimate goal of bridging the theory and practice of stochastic methods, recent works on high-probability convergence guarantees \parencite{nazin2019algorithms, davis2021low, gorbunov2020stochastic, gorbunov2021near, gorbunov2022clipped, cutkosky2021high} focus on relaxing the assumptions under which these guarantees are derived. \emph{This chapter complements that line of work in two main aspects: across a broad range of settings, we derive new high-probability results allowing the variance of the noise and the gradient of the objective to be unbounded.}

\subsection{Technical Preliminaries}
Before presenting the main results, we introduce the problem classes and the necessary preliminaries. In particular, we consider stochastic unconstrained optimization problems 
\begin{equation}
    \min\limits_{x\in\R^d} \left\{f(x) = \Exp_{\xi\sim\thesischtwocD}\left[f(x;\xi)\right]\right\}, \label{ch2:eq:min_problem_chapter_2}
\end{equation}
where $\xi$ is a random variable with distribution $\thesischtwocD$. Such problems often arise in machine learning, where $f(x;\xi)$ represents the loss function on the data sample $\xi$ \parencite{shalev2014understanding}.

Another class considered in this chapter is unconstrained variational inequality problems (VIP), i.e., non-linear equations \parencite{harker1990finite, ryu2021large}:
\begin{equation}
    \text{find } x^*\in \R^d \text{ such that } F(x^*) = 0, \label{ch2:eq:VIP_chapter_2}
\end{equation}
where $F(x) = \Exp_{\xi\sim \thesischtwocD}[F(x;\xi)]$. These problems arise in adversarial/game formulations of machine learning tasks \parencite{goodfellow2014generative, gidel2019variational}.

\textbf{Notation.} We use standard notation: $\|x\| = \sqrt{\langle x, x\rangle}$ denotes the standard Euclidean norm in $\R^d$, $\Exp_{\xi}[\cdot]$ denotes an expectation w.r.t.\ the randomness coming from random variable $\xi$, $B_{R}(x) = \{y \in \R^d\mid \|y-x\| \leq R\}$ is a ball with center at $x$ and radius $R$. We define restricted gap-function as $\gap_{R}(x) = \max_{y \in B_{R}(x^*)}\langle F(y), x - y \rangle$ -- a standard convergence criterion for monotone VIP \parencite{nesterov2007dual}. Finally, $\cO(\cdot)$ hides numerical factors and $\widetilde\cO(\cdot)$ hides poly-logarithmic and numerical factors.

\textbf{Assumptions on a subset.} Although we consider unconstrained problems, our analysis does not require any assumptions to hold on the whole space. For our purposes, it is sufficient to introduce all assumptions only on some subset of $\R^d$, since we prove that the considered methods do not leave some ball around the solution or some level-set of the objective function with high probability. This allows us to consider quite large classes of problems.

\textbf{Stochastic oracle.} We assume that at given point $x$ we have an access to the unbiased stochastic oracle returning $\nabla f(x;\xi)$ or $F(x;\xi)$ that satisfy the following conditions.
\begin{assumption}\label{ch2:as:bounded_alpha_moment_chapter_2}
    We assume that there exist some set $Q \subseteq \R^d$ and values $\sigma \geq 0$, ${p} \in (1,2]$ such that for all $x \in Q$
    \begin{itemize}
        \item[(i)] for problem \eqref{ch2:eq:min_problem_chapter_2} $\Exp_{\xi\sim \thesischtwocD}[\nabla f(x;\xi)] = \nabla f(x)$ and
    \begin{equation}
        \Exp_{\xi\sim \thesischtwocD}[\left\|\nabla f(x;\xi) - \nabla f(x)\right\|^{p}] \leq \sigma^{p}, \label{ch2:eq:bounded_alpha_moment_gradient_chapter_2}
    \end{equation}
        \item[(ii)] for problem \eqref{ch2:eq:VIP_chapter_2} $\Exp_{\xi\sim \thesischtwocD}[F(x;\xi)] = F(x)$ and
    \begin{equation}
        \Exp_{\xi\sim \thesischtwocD}[\left\|F(x;\xi) - F(x)\right\|^{p}] \leq \sigma^{p}. \label{ch2:eq:bounded_alpha_moment_operator_chapter_2}
    \end{equation}
    \end{itemize}
\end{assumption}
When ${p} = 2$, the above assumption recovers the standard uniformly bounded variance assumption \parencite{nemirovski2009robust, ghadimi2012optimal, ghadimi2013stochastic}. However, Assumption~\ref{ch2:as:bounded_alpha_moment_chapter_2} allows the variance of the estimator to be \emph{unbounded} when ${p} \in (1,2)$, i.e., the noise can follow some heavy-tailed distribution. For example, the distribution of the gradient noise in the training of large attention models resembles L\'evy ${p}$-stable distribution with ${p} < 2$ \parencite{zhang2020adaptive}.

\textbf{Assumptions on $f$.} We start with a very mild assumption since without it, problem~\eqref{ch2:eq:min_problem_chapter_2} does not make sense.
\begin{assumption}\label{ch2:as:lower_boundedness_chapter_2}
    We assume that there exist some set $Q \subseteq \R^d$ such that $f$ is uniformly lower-bounded on $Q$: $f^{\inf} = \inf_{x \in Q}f(x) > -\infty$.
\end{assumption}

Moreover, when working with minimization problems \eqref{ch2:eq:min_problem_chapter_2}, we always assume smoothness of $f$.
\begin{assumption}\label{ch2:as:L_smoothness_chapter_2}
    We assume that there exist some set $Q \subseteq \R^d$ and constant $L > 0$ such that for all $x, y \in Q$
    \begin{eqnarray}
        \|\nabla f(x) - \nabla f(y)\| &\leq& L\|x - y\|, \label{ch2:eq:L_smoothness_chapter_2}\\
        \|\nabla f(x)\|^2 &\leq& 2L\left(f(x) - f^{\inf}\right), \label{ch2:eq:L_smoothness_cor_2_chapter_2}
    \end{eqnarray}
    where $f^{\inf} = \inf_{x \in Q}f(x) > -\infty$.
\end{assumption}
We notice here that \eqref{ch2:eq:L_smoothness_cor_2_chapter_2} follows from \eqref{ch2:eq:L_smoothness_chapter_2} for $Q = \R^d$, but in the general case, the implication is slightly more involved (see the details in Appendix~\ref{ch2:appendix:useful_facts_chapter_2}).

In each particular special case, we also make \emph{one of the following assumptions} about the structured non-convexity of the objective function. The previous two assumptions hold for a very broad class of functions. The next assumption -- Polyak-{\L}ojasiewicz condition \parencite{polyak1963gradient, lojasiewicz1963topological} -- narrows the class of non-convex functions.
\begin{assumption}\label{ch2:as:PL_chapter_2}
    We assume that there exist some set $Q \subseteq \R^d$ and constant $\mu > 0$ such that $f$ satisfies Polyak-{\L}ojasiewicz (P{\L}) condition/inequality on $Q$, i.e., for all $x \in Q$ and $x^* = \arg\min_{x\in \R^d} f(x)$
    \begin{equation}
        \|\nabla f(x)\|^2 \geq 2\mu\left(f(x) - f(x^*)\right). \label{ch2:eq:PL_chapter_2}
    \end{equation}
\end{assumption}
When function $f$ is $\mu$-strongly convex, it satisfies P{\L} condition. However, {P{\L}} inequality can hold even for non-convex functions. Some analogs of this assumption have been observed for over-parameterized models \parencite{liu2022loss}.

We also consider another relaxation of convexity.
\begin{assumption}\label{ch2:as:QSC_chapter_2}
    We assume that there exist some set $Q \subseteq \R^d$ and constant $\mu \geq 0$ such that $f$ is $\mu$-quasi-strongly convex, i.e., for all $x \in Q$ and $x^* = \arg\min_{x\in \R^d} f(x)$
    \begin{equation}
        f(x^*) \geq f(x) + \langle \nabla f(x), x^* - x \rangle + \frac{\mu}{2}\|x - x^*\|^2. \label{ch2:eq:QSC_chapter_2}
    \end{equation}
\end{assumption}
As P{\L} condition, this assumption holds for any $\mu$-strongly convex function but does not imply convexity. Nevertheless, for the above two assumptions, some standard deterministic methods such as Gradient Descent (\algname{GD}) converge linearly; see \textcite{necoara2019linear} for more details and examples.

In the analysis of the accelerated method, we also need standard (strong) convexity.
\begin{assumption}\label{ch2:as:str_cvx_chapter_2}
    We assume that there exist some set $Q \subseteq \R^d$ and constant $\mu \geq 0$ such that $f$ is $\mu$-strongly convex, i.e., for all $x, y \in Q$
    \begin{equation}
        f(y) \geq f(x) + \langle \nabla f(x), y - x \rangle + \frac{\mu}{2}\|y - x\|^2. \label{ch2:eq:str_cvx_chapter_2}
    \end{equation}
    When $\mu = 0$ function $f$ is called convex.
\end{assumption}

\textbf{Assumptions on $F$.} In the context of solving \eqref{ch2:eq:VIP_chapter_2}, we assume Lipschitzness of $F$ -- a standard assumption for VIP.

\begin{assumption}\label{ch2:as:L_Lip_chapter_2}
    We assume that there exist some set $Q \subseteq \R^d$ and constant $L > 0$ such that for all $x, y \in Q$
    \begin{eqnarray}
        \|F(x) - F(y)\| \leq L\|x - y\|, \label{ch2:eq:L_Lip_chapter_2}
    \end{eqnarray}
\end{assumption}

Similarly to the case of minimization problems, we make \emph{one or two of the following assumptions} about the structured non-monotonicity of the operator $F$. The first assumption we consider is the standard monotonicity.

\begin{assumption}\label{ch2:as:monotonicity_chapter_2}
    We assume that there exist some set $Q \subseteq \R^d$ such that $F$ is monotone on $Q$, i.e., for all $x, y \in Q$
    \begin{equation}
        \langle F(x) - F(y), x - y \rangle \geq  0.\label{ch2:eq:monotonicity_chapter_2}
    \end{equation}
\end{assumption}
Monotonicity can be seen as an analog of convexity for VIP. When \eqref{ch2:eq:monotonicity_chapter_2} holds with $\mu\|x-y\|^2$ in the r.h.s. instead of just 0, operator $F$ is called $\mu$-strongly monotone.

Next, we consider quasi-strong monotonicity \parencite{mertikopoulos2019learning, song2020optimistic, loizou2021stochastic} -- a relaxation of strong monotonicity. There exist examples of non-monotone problems such that the assumption below holds \parencite[Appendix A.6]{loizou2021stochastic}.
\begin{assumption}\label{ch2:as:QSM_chapter_2}
    We assume that there exist some set $Q \subseteq \R^d$ and constant $\mu > 0$ such that $F$ is $\mu$-quasi strongly monotone on $Q$, i.e., for all $x \in Q$ and $x^*$ such that $F(x^*) = 0$ we have
    \begin{equation}
        \langle F(x), x - x^* \rangle \geq  \mu \|x - x^*\|^2.\label{ch2:eq:QSM_chapter_2}
    \end{equation}
\end{assumption}

Another structured non-monotonicity assumption considered in this chapter is star-cocoercivity.
\begin{assumption}\label{ch2:as:star-cocoercivity_chapter_2}
    We assume that there exist some set $Q \subseteq \R^d$ and constant $\ell > 0$ such that $F$ is star-cocoercive on $Q$, i.e., for all $x \in Q$ and $x^*$ such that $F(x^*) = 0$
    \begin{equation}
        \|F(x)\|^2 \leq \ell\langle F(x), x - x^* \rangle.\label{ch2:eq:star-cocoercivity_chapter_2}
    \end{equation}
\end{assumption}
This assumption can be seen as a relaxation of the standard cocoercivity: $\|F(x) - F(y)\|^2 \leq \ell \langle F(x) - F(y), x - y \rangle$. However, unlike cocoercivity, star-cocoercivity implies neither monotonicity nor Lipschitzness of operator $F$ \parencite[Appendix A.6]{loizou2021stochastic}.

\subsection{Closely Related Work and Contributions}

\begin{table*}[t]
    \centering
    \scriptsize
    \caption[High-probability complexity results for stochastic optimization]{\scriptsize Summary of known and new high-probability complexity results for solving smooth problem \eqref{ch2:eq:min_problem_chapter_2}. Column ``Setup'' indicates the assumptions made in addition to Assumptions~\ref{ch2:as:bounded_alpha_moment_chapter_2} and \ref{ch2:as:L_smoothness_chapter_2}. All assumptions are made only on some ball around the solution with radius $\sim R \geq \|x^0 - x^*\|$ (unless the opposite is indicated). By the complexity we mean the number of stochastic oracle calls needed for a method to guarantee that $\Prob\{\text{Metric} \leq \varepsilon\} \geq 1 - \beta$ for some $\varepsilon> 0$, $\beta \in (0,1]$ and ``Metric'' is taken from the corresponding column. For simplicity, we omit numerical and logarithmic factors in the complexity bounds. Column ``${p}$'' shows the allowed values of ${p}$, ``UD?'' shows whether the analysis works on unbounded domains, and ``UG?'' indicates whether the analysis works without assuming boundedness of the gradient. Notation: $L$ = Lipschitz constant; $D$ = diameter of the domain in the result of \textcite{nazin2019algorithms}; $\sigma$ = parameter from Assumption~\ref{ch2:as:bounded_alpha_moment_chapter_2}; $R$ = any upper bound on $\|x^0 - x^*\|$; $\mu$ = (quasi-)strong convexity/Polyak-{\L}ojasiewicz parameter; $\Delta$ = any upper bound on $f(x^0) - f^{\inf}$; $G$ = parameter such that $\Exp_{\xi\sim\thesischtwocD}\left[\|\nabla f(x;\xi)\|^{p}\right] \leq G^{p}$, as assumed by \textcite{cutkosky2021high}. The results established in this chapter are highlighted in blue.
    }
    \label{ch2:tab:comparison_of_rates_minimization_chapter_2}
    \begin{threeparttable}
        \resizebox{0.95\textwidth}{!}{%
        \begin{tabular}{|c|c c c c c c c|}
        \hline
Setup & Method & Citation & Metric & Complexity & ${p}$ & UD? & UG?\\
        \hline\hline
        \multirow{8}{*}{\makecell{As.~\ref{ch2:as:str_cvx_chapter_2}\\ ($\mu = 0$)}} & \algname{RSMD} & \tabcite{nazin2019algorithms}\tnote{\color{blue}(1)} & $f(\overline{x}^K) - f(x^*)$ & $\max\left\{\frac{LD^2}{\varepsilon}, \frac{\sigma^2 D^2}{\varepsilon^2}\right\}$ & {\color{BrickRed}$2$} & \thesischtwoxmark & \thesischtwocmark\\
        & \algname{clipped-SGD} & \makecell{\tabcite{gorbunov2020stochastic}\\\tabcite{gorbunov2021near}} & $f(\overline{x}^K) - f(x^*)$ & $\max\left\{\frac{LR^2}{\varepsilon}, \frac{\sigma^2 R^2}{\varepsilon^2}\right\}$ & {\color{BrickRed}$2$} & \thesischtwocmark & \thesischtwocmark \\
        & \algname{clipped-SSTM} & \makecell{\tabcite{gorbunov2020stochastic}\\\tabcite{gorbunov2021near}} & $f(y^K) - f(x^*)$ & $\max\left\{\sqrt{\frac{LR^2}{\varepsilon}}, \frac{\sigma^2 R^2}{\varepsilon^2}\right\}$ & {\color{BrickRed}$2$} & \thesischtwocmark & \thesischtwocmark \\
        & \cellcolor{linen}\algname{clipped-SGD} & \cellcolor{linen}Theorems~\ref{ch2:thm:clipped_SGD_main_theorem_chapter_2} \& \ref{ch2:thm:clipped_SGD_convex_case_appendix_chapter_2} & \cellcolor{linen}$f(\overline{x}^K) - f(x^*)$ &\cellcolor{linen} $\max\left\{\frac{LR^2}{\varepsilon}, \left(\frac{\sigma R}{\varepsilon}\right)^{\frac{{p}}{{p}-1}}\right\}$ &\cellcolor{linen} {\color{PineGreen}$(1,2]$} &\cellcolor{linen} \thesischtwocmark &\cellcolor{linen} \thesischtwocmark \\
        & \cellcolor{linen} \algname{clipped-SSTM} &\cellcolor{linen}Theorems~\ref{ch2:thm:clipped_SSTM_main_theorem_chapter_2} \& \ref{ch2:thm:clipped_SSTM_convex_case_appendix_chapter_2} &\cellcolor{linen} $f(y^K) - f(x^*)$ &\cellcolor{linen} $\max\left\{\sqrt{\frac{LR^2}{\varepsilon}}, \left(\frac{\sigma R}{\varepsilon}\right)^{\frac{{p}}{{p}-1}}\right\}$ & \cellcolor{linen}{\color{PineGreen}$(1,2]$} &\cellcolor{linen} \thesischtwocmark &\cellcolor{linen} \thesischtwocmark \\
        \hline
        \multirow{8}{*}{\makecell{As.~\ref{ch2:as:str_cvx_chapter_2}\\ ($\mu > 0$)}} & \algname{restarted-RSMD} & \tabcite{nazin2019algorithms}\tnote{\color{blue}(1)} & $f(\overline{x}^K) - f(x^*)$ & $\max\left\{\frac{L}{\mu}, \frac{\sigma^2}{\mu\varepsilon}\right\}$ & {\color{BrickRed}$2$} & \thesischtwoxmark & \thesischtwocmark\\
        & \algname{proxBoost} & \tabcite{davis2021low}\tnote{\color{blue}(1)} & $f(\overline{x}^K) - f(x^*)$ & $\max\left\{\sqrt{\frac{L}{\mu}}, \frac{\sigma^2}{\mu\varepsilon}\right\}$\tnote{\color{blue}(2)} & {\color{BrickRed}$2$} & \thesischtwocmark & \thesischtwocmark\\
        & \algname{R-clipped-SGD} & \makecell{\tabcite{gorbunov2020stochastic}\\\tabcite{gorbunov2021near}} & $f(\overline{x}^K) - f(x^*)$ & $\max\left\{\frac{L}{\mu}, \frac{\sigma^2}{\mu\varepsilon}\right\}$ & {\color{BrickRed}$2$} & \thesischtwocmark & \thesischtwocmark \\
        & \algname{R-clipped-SSTM} & \makecell{\tabcite{gorbunov2020stochastic}\\\tabcite{gorbunov2021near}} & $f(y^K) - f(x^*)$ & $\max\left\{\sqrt{\frac{L}{\mu}}, \frac{\sigma^2}{\mu\varepsilon}\right\}$ & {\color{BrickRed}$2$} & \thesischtwocmark & \thesischtwocmark \\
        & \cellcolor{linen}\algname{R-clipped-SSTM} & \cellcolor{linen}Theorems~\ref{ch2:thm:clipped_SSTM_main_theorem_chapter_2} \& \ref{ch2:thm:R_clipped_SSTM_main_theorem_appendix_chapter_2} & \cellcolor{linen}$f(y^K) - f(x^*)$ & \cellcolor{linen}$\max\left\{\sqrt{\frac{L}{\mu}}, \left(\frac{\sigma^2}{\mu\varepsilon}\right)^{\frac{{p}}{2({p}-1)}}\right\}$ &\cellcolor{linen} {\color{PineGreen}$(1,2]$} & \cellcolor{linen}\thesischtwocmark &\cellcolor{linen} \thesischtwocmark \\
        \hline
        \makecell{As.~\ref{ch2:as:QSC_chapter_2}\\ ($\mu > 0$)} & \cellcolor{linen}\algname{clipped-SGD} & \cellcolor{linen}Theorems~\ref{ch2:thm:clipped_SGD_main_theorem_chapter_2} \& \ref{ch2:thm:main_result_QSC_SGD_chapter_2} &\cellcolor{linen} $\|x^K - x^*\|^2$ &\cellcolor{linen} $\max\left\{\frac{L}{\mu}, \left(\frac{\sigma^2}{\mu^2\varepsilon}\right)^{\frac{{p}}{2({p}-1)}}\right\}$ &\cellcolor{linen} {\color{PineGreen}$(1,2]$} &\cellcolor{linen} \thesischtwocmark &\cellcolor{linen} \thesischtwocmark \\
        \hline
        \multirow{5.5}{*}{As.~\ref{ch2:as:lower_boundedness_chapter_2}} & \algname{MSGD} & \tabcite{li2020high}\tnote{\color{blue}(1)} & $\frac{1}{K+1}\sum\limits_{k=0}^K\|\nabla f(x^k)\|^2$ & $\max\left\{\frac{L^2 \Delta^2}{\varepsilon}, \frac{\sigma^4}{\varepsilon^2}\right\}$ & {\color{BrickRed}\thesischtwoxmark}\tnote{\color{blue}(3)} & \thesischtwocmark & \thesischtwocmark\\
        & \algname{clipped-NMSGD} & \tabcite{cutkosky2021high}\tnote{\color{blue}(1)} & $\left(\frac{1}{K+1}\sum\limits_{k=0}^K\|\nabla f(x^k)\|\right)^2$ \tnote{\color{blue}(4)} & $\left(\frac{G^2}{\varepsilon}\right)^{\frac{3{p}-2}{2{p}-2}}$ & {\color{PineGreen}$(1,2]$} & \thesischtwocmark & \thesischtwoxmark\\
        &\cellcolor{linen} \algname{clipped-SGD} &\cellcolor{linen} Theorems~\ref{ch2:thm:clipped_SGD_main_theorem_chapter_2} \& \ref{ch2:thm:clipped_SGD_non_convex_main_chapter_2} \tnote{\color{blue}(5)}&\cellcolor{linen} $\frac{1}{K+1}\sum\limits_{k=0}^K\|\nabla f(x^k)\|^2$ &\cellcolor{linen} $\max\left\{\frac{L \Delta}{\varepsilon}, \left(\frac{\sqrt{L\Delta} \sigma}{\varepsilon}\right)^{\frac{{p}}{{p}-1}}\right\}$ &\cellcolor{linen} {\color{PineGreen}$(1,2]$} &\cellcolor{linen} \thesischtwocmark & \cellcolor{linen}\thesischtwocmark \\
        \hline
        As.~\ref{ch2:as:PL_chapter_2} & \cellcolor{linen}\algname{clipped-SGD} & \cellcolor{linen}Theorems~\ref{ch2:thm:clipped_SGD_main_theorem_chapter_2} \& \ref{ch2:thm:clipped_SGD_PL_main_chapter_2} \tnote{\color{blue}(5)}& \cellcolor{linen}$f(x^K) - f(x^*)$ & \cellcolor{linen}$\max\left\{\frac{L}{\mu}, \left(\frac{L\sigma^2}{\mu^2\varepsilon}\right)^{\frac{{p}}{2({p}-1)}}\right\}$ &\cellcolor{linen} {\color{PineGreen}$(1,2]$} &\cellcolor{linen} \thesischtwocmark &\cellcolor{linen} \thesischtwocmark \\
        \hline
    \end{tabular}}
    \par\begin{minipage}{0.95\textwidth}
        \scriptsize\raggedright\setlength{\parindent}{0pt}
        {\color{blue}(1)} All assumptions are made on the whole domain.\par
        {\color{blue}(2)} Complexity has extra logarithmic factor of $\ln(\nicefrac{L}{\mu})$.\par
        {\color{blue}(3)} \textcite{li2020high} assume that the noise is sub-Gaussian: $\Exp\left[\exp\left(\nicefrac{\|\nabla f(x;\xi) - \nabla f(x)\|^2}{\sigma^2}\right)\right] \leq \exp(1)$ for all $x$ from the domain.\par
        {\color{blue}(4)} We notice that $\left(\frac{1}{K+1}\sum_{k=0}^K\|\nabla f(x^k)\|\right)^2 \leq \frac{1}{K+1}\sum_{k=0}^K\|\nabla f(x^k)\|^2$ and in the worst case the left-hand side is $K+1$ times smaller than the right-hand side.\par
        {\color{blue}(5)} All assumptions are made on the level set $Q = \{x \in \R^d\mid \exists y\in \R^d:\; f(y) \leq f^{\inf} + 2\Delta \text{ and } \|x-y\|\leq \nicefrac{\sqrt{\Delta}}{20\sqrt{L}}\}$.
    \end{minipage}\par
    \end{threeparttable}
\end{table*}

\begin{table*}[t]
    \centering
    \scriptsize
    \caption[High-probability complexity results for stochastic variational inequalities]{\scriptsize Summary of known and new high-probability complexity results for solving \eqref{ch2:eq:VIP_chapter_2}. Column ``Setup'' indicates the assumptions made in addition to Assumption~\ref{ch2:as:bounded_alpha_moment_chapter_2}. All assumptions are made only on some ball around the solution with radius $\sim R \geq \|x^0 - x^*\|$ (unless the opposite is indicated). By the complexity we mean the number of stochastic oracle calls needed for a method to guarantee that $\Prob\{\text{Metric} \leq \varepsilon\} \geq 1 - \beta$ for some $\varepsilon> 0$, $\beta \in (0,1]$ and ``Metric'' is taken from the corresponding column. For simplicity, we omit numerical and logarithmic factors in the complexity bounds. Column ``${p}$'' shows the allowed values of ${p}$, ``UD?'' shows whether the analysis works on unbounded domains, and ``UG?'' indicates whether the analysis works without assuming boundedness of the gradient. Notation: $\thesischtwotx^K_{\text{avg}} = \frac{1}{K+1}\sum_{k=0}^K \thesischtwotx^k$ (for \algname{clipped-SEG}), $x^K_{\avg} = \frac{1}{K+1}\sum_{k=0}^K x^k$ (for \algname{clipped-SGDA}); $L$ = Lipschitz constant; $D$ = diameter of the domain used by \textcite{juditsky2011solving}; $\gap_{D}(x) = \max_{y \in \thesischtwocX}\langle F(y), x - y\rangle$, where $\thesischtwocX$ is a bounded domain with diameter $D$ where the problem is defined, as used by \textcite{juditsky2011solving}; $D$ = diameter of the domain in the result of \textcite{juditsky2011solving}; $\sigma$ = parameter from Assumption~\ref{ch2:as:bounded_alpha_moment_chapter_2}; $R$ = any upper bound on $\|x^0 - x^*\|$; $\mu$ = quasi-strong monotonicity parameter; $\ell$ = star-cocoercivity parameter. The results established in this chapter are highlighted in blue.}
    \label{ch2:tab:comparison_of_rates_VIP_chapter_2}
    \begin{threeparttable}
        \resizebox{0.95\textwidth}{!}{%
        \begin{tabular}{|c|c c c c c c c|}
        \hline
Setup & Method & Citation & Metric & Complexity & ${p}$ & UD? & UG?\\
        \hline\hline
        \multirow{4}{*}{As.~\ref{ch2:as:L_Lip_chapter_2} \& \ref{ch2:as:monotonicity_chapter_2}} & \algname{Mirror-Prox} & \tabcite{juditsky2011solving}\tnote{\color{blue}(1)} & $\gap_{D}(\thesischtwotx_{\avg}^K)$ & $\max\left\{\frac{LD^2}{\varepsilon}, \frac{\sigma^2 D^2}{\varepsilon^2}\right\}$ & \thesischtwoxmark\tnote{\color{blue}(2)} & \thesischtwoxmark & \thesischtwocmark\\
        & \algname{clipped-SEG} & \tabcite{gorbunov2022clipped} & $\gap_{R}(\thesischtwotx_{\avg}^K)$ & $\max\left\{\frac{LR^2}{\varepsilon}, \frac{\sigma^2 R^2}{\varepsilon^2}\right\}$ & {\color{BrickRed}$2$} & \thesischtwocmark & \thesischtwocmark \\
        & \cellcolor{linen}\algname{clipped-SEG} & \cellcolor{linen}Theorems~\ref{ch2:thm:clipped_SEG_main_theorem_chapter_2} \& \ref{ch2:thm:main_result_gap_SEG_chapter_2} & \cellcolor{linen}$\gap_{R}(\thesischtwotx_{\avg}^K)$ &\cellcolor{linen} $\max\left\{\frac{LR^2}{\varepsilon}, \left(\frac{\sigma R}{\varepsilon}\right)^{\frac{{p}}{{p}-1}}\right\}$ &\cellcolor{linen} {\color{PineGreen}$(1,2]$} &\cellcolor{linen} \thesischtwocmark &\cellcolor{linen} \thesischtwocmark \\
        \hline
        \multirow{2.5}{*}{As.~\ref{ch2:as:L_Lip_chapter_2} \& \ref{ch2:as:QSM_chapter_2}} & \algname{clipped-SEG} & \tabcite{gorbunov2022clipped} & $\|x^k - x^*\|^2$ & $\max\left\{\frac{L}{\mu}, \frac{\sigma^2}{\mu^2\varepsilon}\right\}$ & {\color{BrickRed}$2$} & \thesischtwocmark & \thesischtwocmark\\
        & \cellcolor{linen}\algname{clipped-SEG} & \cellcolor{linen}Theorems~\ref{ch2:thm:clipped_SEG_main_theorem_chapter_2} \& \ref{ch2:thm:main_result_str_mon_SEG_chapter_2}& \cellcolor{linen}$\|x^k - x^*\|^2$ &\cellcolor{linen} $\max\left\{\frac{L}{\mu}, \left(\frac{\sigma^2}{\mu^2\varepsilon}\right)^{\frac{{p}}{2({p}-1})}\right\}$ &\cellcolor{linen} {\color{PineGreen}$(1,2]$} &\cellcolor{linen} \thesischtwocmark &\cellcolor{linen} \thesischtwocmark \\
        \hline\hline
        \multirow{2.5}{*}{As.~\ref{ch2:as:monotonicity_chapter_2} \& \ref{ch2:as:star-cocoercivity_chapter_2}} & \algname{clipped-SGDA} & \tabcite{gorbunov2022clipped} & $\gap_{R}(x_{\avg}^K)$ & $\max\left\{\frac{\ell R^2}{\varepsilon}, \frac{\sigma^2R^2}{\varepsilon^2}\right\}$ & {\color{BrickRed}$2$} & \thesischtwocmark & \thesischtwocmark\\
        & \cellcolor{linen}\algname{clipped-SGDA} & \cellcolor{linen}Theorems~\ref{ch2:thm:clipped_SGDA_main_theorem_chapter_2} \& \ref{ch2:thm:main_result_gap_SGDA_chapter_2} & \cellcolor{linen}$\gap_{R}(x_{\avg}^K)$ &\cellcolor{linen} $\max\left\{\frac{\ell R^2}{\varepsilon}, \left(\frac{\sigma R}{\varepsilon}\right)^{\frac{{p}}{{p}-1}}\right\}$ &\cellcolor{linen} {\color{PineGreen}$(1,2]$} &\cellcolor{linen} \thesischtwocmark &\cellcolor{linen} \thesischtwocmark \\
        \hline
        \multirow{3}{*}{As.~\ref{ch2:as:star-cocoercivity_chapter_2}} & \algname{clipped-SGDA} & \tabcite{gorbunov2022clipped} & $\frac{1}{K+1}\sum\limits_{k=0}^{K}\|F(x^k)\|^2$ & $\max\left\{\frac{\ell^2 R^2}{\varepsilon}, \frac{\ell^2\sigma^2R^2}{\varepsilon^2}\right\}$ & {\color{BrickRed}$2$} & \thesischtwocmark & \thesischtwocmark\\
        & \cellcolor{linen}\algname{clipped-SGDA} & \cellcolor{linen}Theorems~\ref{ch2:thm:clipped_SGDA_main_theorem_chapter_2} \& \ref{ch2:thm:main_result_non_mon_SGDA_chapter_2} & \cellcolor{linen}$\frac{1}{K+1}\sum\limits_{k=0}^{K}\|F(x^k)\|^2$ &\cellcolor{linen} $\max\left\{\frac{\ell^2 R^2}{\varepsilon}, \left(\frac{\ell \sigma R}{\varepsilon}\right)^{\frac{{p}}{{p}-1}}\right\}$ &\cellcolor{linen} {\color{PineGreen}$(1,2]$} &\cellcolor{linen} \thesischtwocmark &\cellcolor{linen} \thesischtwocmark \\
        \hline
        \multirow{2.5}{*}{As.~\ref{ch2:as:QSM_chapter_2} \& \ref{ch2:as:star-cocoercivity_chapter_2}} & \algname{clipped-SGDA} & \tabcite{gorbunov2022clipped} & $\|x^K - x^*\|^2$ & $\max\left\{\frac{\ell}{\mu}, \frac{\sigma^2}{\mu^2\varepsilon}\right\}$ & {\color{BrickRed}$2$} & \thesischtwocmark & \thesischtwocmark\\
        & \cellcolor{linen}\algname{clipped-SGDA} & \cellcolor{linen}Theorems~\ref{ch2:thm:clipped_SGDA_main_theorem_chapter_2} \& \ref{ch2:thm:main_result_str_mon_SGDA_chapter_2} & \cellcolor{linen}$\|x^K - x^*\|^2$ &\cellcolor{linen} $\max\left\{\frac{\ell}{\mu}, \left(\frac{\sigma^2}{\mu^2\varepsilon}\right)^{\frac{{p}}{2({p}-1})}\right\}$ &\cellcolor{linen} {\color{PineGreen}$(1,2]$} &\cellcolor{linen} \thesischtwocmark &\cellcolor{linen} \thesischtwocmark \\
        \hline
    \end{tabular}}
    \par\begin{minipage}{0.95\textwidth}
        \scriptsize\raggedright\setlength{\parindent}{0pt}
        {\color{blue}(1)} All assumptions are made on the whole domain.\par
        {\color{blue}(2)} \textcite{juditsky2011solving} assume that the noise is sub-Gaussian: $\Exp\left[\exp\left(\nicefrac{\|F(x;\xi) - F(x)\|^2}{\sigma^2}\right)\right] \leq \exp(1)$ for all $x$ from the domain.
    \end{minipage}\par
    \end{threeparttable}
\end{table*}

This subsection reviews closely related work and describes the contributions of this chapter. Additional related work is discussed in Appendix~\ref{ch2:appendix:additional_related_chapter_2}.

\textbf{Convex optimization and monotone VIPs.} Classical high-probability results for (strongly) convex minimization \parencite{nemirovski2009robust, ghadimi2012optimal} and monotone VIP \parencite{juditsky2011solving} are derived under the so-called light-tails assumption, meaning that the noise in the stochastic gradients/operators is assumed to be sub-Gaussian: $\Exp_{\xi \sim \thesischtwocD}[\exp(\nicefrac{\|\nabla f(x;\xi) - \nabla f(x)\|^2}{\sigma^2})] \leq \exp(1)$ or $\Exp_{\xi \sim \thesischtwocD}[\exp(\nicefrac{\|F(x;\xi) - F(x)\|^2}{\sigma^2})] \leq \exp(1)$. In these settings, optimal (up to logarithmic factors) rates of convergence are derived in the mentioned papers.

The first high-probability results with logarithmic dependence\footnote{Note that from in-expectation convergence guarantee, one can always get a high-probability one using Markov's inequality. For example, under bounded variance, smoothness, and strong convexity assumptions \algname{SGD} achieves $\Exp\left[\|x^k - x^*\|^2\right] \leq \varepsilon$ after $k = \widetilde\cO(\max\{\nicefrac{L}{\mu}, \nicefrac{\sigma^2}{\mu\varepsilon}\})$ iterations. Therefore, taking $k$ such that $\Exp\left[\|x^k - x^*\|^2\right] \leq \varepsilon\beta$ we get from Markov's inequality that $\Prob\{\|x^k - x^*\|^2 > \varepsilon\} \leq \beta$. However, in this case, we get bound $k = \widetilde\cO(\max\{\nicefrac{L}{\mu}, \nicefrac{\sigma^2}{\mu\varepsilon\beta}\})$, having undesirable inverse-power dependence on $\beta$.} on $\nicefrac{1}{\beta}$ under just bounded variance assumption are given by \textcite{nazin2019algorithms}, where the authors show non-accelerated rates of convergence for a version of Mirror Descent with a special truncation operator for smooth convex and strongly convex problems defined on the bounded sets. Then, \textcite{davis2021low} derive accelerated rates in the strongly convex case using robust distance estimation techniques. \textcite{gorbunov2020stochastic, gorbunov2021near} propose an accelerated method with clipping for unconstrained (strongly) convex problems with Lipschitz / H\"older continuous gradients and derive the first high-probability results for \algname{clipped-SGD}. In the context of VIP, \textcite{gorbunov2022clipped} derive the first high-probability results for the stochastic methods for solving VIP under bounded variance assumption and different assumptions on structured non-monotonicity.

However, there are no high-probability results (with logarithmic dependence on the confidence level) for smooth (strongly) convex minimization problems and Lipschitz VIP without imposing bounded variance assumption. Only recently, \textcite{zhang2022parameter} derived optimal regret-bounds under Assumption~\ref{ch2:as:bounded_alpha_moment_chapter_2} in the convex case with bounded gradients on $\R^d$. However, the bounded gradients assumption is quite restrictive when assumed on the whole space. Thus, a noticeable gap in the stochastic optimization literature remains.

\textbf{\emph{Contribution.}} We obtain new high-probability convergence results under Assumption~\ref{ch2:as:bounded_alpha_moment_chapter_2} for smooth convex minimization problems and Lipschitz VIP; see the summary in Tables~\ref{ch2:tab:comparison_of_rates_minimization_chapter_2} and \ref{ch2:tab:comparison_of_rates_VIP_chapter_2}. In particular, for Clipped Stochastic Similar Triangles Method (\algname{clipped-SSTM}) \parencite{gorbunov2020stochastic} and its restarted version, we derive high-probability convergence results  for smooth convex and strongly convex problems. The high-probability complexity in the strongly convex case matches (up to logarithmic factors) the known in-expectation lower bound \parencite{zhang2020adaptive} and deterministic lower bound \parencite{NemiYudin}. In other words, we derive the first optimal high-probability complexity results for smooth strongly convex optimization. Noticeably, the derived results have clear separation between accelerated part and stochastic part that emphasizes a potential of \algname{clipped-SSTM} for efficient parallelization. Next, we derive high-probability results for \algname{clipped-SGD} for smooth star-convex and quasi-strongly convex objectives under Assumption~\ref{ch2:as:bounded_alpha_moment_chapter_2}. Finally, under the same assumption, we prove the high-probability convergence of Clipped Stochastic Extragradient (\algname{clipped-SEG}) \parencite{korpelevich1976extragradient, juditsky2011solving, gorbunov2022clipped} for Lipschitz monotone and quasi-strongly monotone VIP and also obtain high-probability results for Clipped Stochastic Gradient Descent-Ascent \algname{(clipped-SGDA)} for star-cocoercive and monotone / quasi-strongly monotone VIP. In the special case of ${p} = 2$, our analysis recovers state-of-the-art high-probability results under the bounded-variance assumption.

\textbf{Non-convex optimization.} Under the light-tails and smoothness assumption \textcite{li2020high} derive high-probability convergence rates to the first-order stationary point for \algname{SGD}. These rates match the known in-expectation guarantees for \algname{SGD} and are optimal up to logarithmic factors \parencite{arjevani2022lower}. Recently, \textcite{cutkosky2021high} derived the first high-probability results for non-convex optimization under Assumption~\ref{ch2:as:bounded_alpha_moment_chapter_2} for a version of \algname{SGD} with gradient clipping and normalization of the momentum. The results are obtained for the non-standard metric -- $\frac{1}{K+1}\sum_{k=0}^K\|\nabla f(x^k)\|$ -- and match the in-expectation lower bound for the expected (non-squared) norm of the gradient established by \textcite{zhang2020adaptive}. However, \textcite{cutkosky2021high} make an additional assumption that the norm of the gradient is bounded\footnote{More precisely, instead of Assumption~\ref{ch2:as:bounded_alpha_moment_chapter_2}, \textcite{cutkosky2021high} assume $\Exp_{\xi\sim}\left[\|\nabla f(x;\xi)\|^{p}\right] \leq G^{p}$ for some $G > 0$. This assumption implies Assumption~\ref{ch2:as:bounded_alpha_moment_chapter_2} and boundedness of $\|\nabla f(x)\|$.} on $\R^d$, which is quite restrictive.

\textbf{\emph{Contribution.}} We derive the first high-probability result with logarithmic dependence on the confidence level for finding first-order stationary points of smooth (possibly, non-convex) functions without bounded gradients assumption. The result is derived for simple \algname{clipped-SGD}. Moreover, we extend the analysis to the functions satisfying Polyak-{\L}ojasiewicz condition; see Table~\ref{ch2:tab:comparison_of_rates_minimization_chapter_2} for the summary.

\textbf{Gradient clipping} received a lot of attention in the machine learning community due to its successful empirical applications in the training of deep neural networks \parencite{pascanu2013difficulty, goodfellow2016deep}. The clipping operator is defined as $\clip(x,\lambda) = \min\left\{1 , \nicefrac{\lambda}{\|x\|}\right\}x$ ($\clip(x,\lambda) = 0$, when $x = 0$). 
From the theoretical perspective, gradient clipping is used for multiple different purposes: to handle structured non-smoothness in the objective function \parencite{zhang2020gradient}, to robustify aggregation \parencite{karimireddy2021learning} and to provide privacy guarantees \parencite{abadi2016deep} in the distributed training. Moreover, as we already mentioned before, gradient clipping is used to handle heavy-tailed noise (satisfying Assumption~\ref{ch2:as:bounded_alpha_moment_chapter_2}) in the stochastic gradients \parencite{zhang2020adaptive} and, in particular, to derive better high-probability guarantees under bounded variance assumption \parencite{nazin2019algorithms, gorbunov2020stochastic}. However, there are no results showing the necessity of modifying standard methods like \algname{SGD} and its accelerated variants to achieve high-probability convergence with logarithmic dependence on the confidence level under bounded variance assumption.

\textbf{\emph{Contribution.}} We construct an example of a strongly convex smooth problem and stochastic oracle with bounded variance such that to achieve $\Prob\{\|x^k - x^*\|^2 > \varepsilon\} \leq \beta$ \algname{SGD} requires $\Omega\left(\nicefrac{\stepEightBchange{\sigma}}{\mu\sqrt{\varepsilon\beta}}\right)$ iterations, i.e., the algorithm has inverse-power dependence on the confidence level. This justifies the importance of using some non-linearity such as gradient clipping to achieve logarithmic dependence on the confidence level even in the bounded variance case.

\section{Failure of Standard \algname{SGD}}
It is known that \algname{SGD} $x^{k+1} = x^k - \gamma \nabla f(x^k;\xi^k)$ can diverge in expectation, when Assumption~\ref{ch2:as:bounded_alpha_moment_chapter_2} is satisfied with ${p} < 2$ \parencite[Remark 1]{zhang2020adaptive}. However, it does converge in expectation when ${p} = 2$, i.e., when the variance is bounded. In contrast, there are no high-probability convergence results for \algname{SGD} having logarithmic dependence on $\nicefrac{1}{\beta}$. The next theorem establishes the impossibility of deriving such high-probability results.

\begin{theorem}\label{ch2:thm:SGD_high-prob_conv_chapter_2}
    For any $\varepsilon > 0$ and sufficiently small $\beta \in (0,1)$ there exists a problem \eqref{ch2:eq:min_problem_chapter_2} such that Assumptions~\ref{ch2:as:bounded_alpha_moment_chapter_2}, \ref{ch2:as:L_smoothness_chapter_2}, and \ref{ch2:as:str_cvx_chapter_2} hold with $Q = \R^d$, ${p} = 2$, $0 < \mu \leq L$ and for the iterates produced by \algname{SGD} with any stepsize $\gamma > 0$
    \begin{equation}
        \Prob\left\{\|x^k - x^*\|^2 \geq \varepsilon\right\} \leq \beta \;\; \Longrightarrow\;\; k = \Omega\left(\frac{\sigma}{\mu\sqrt{\varepsilon\beta}}\right). \notag
    \end{equation}
\end{theorem}

The proof is deferred to Appendix~\ref{ch2:appendix:failure_of_SGD_chapter_2}. We believe that similar examples can be constructed for any stochastic first-order methods having linear dependence on the stochastic gradients in their update rules. Thus, Theorem~\ref{ch2:thm:SGD_high-prob_conv_chapter_2} motivates the use of non-linear operators such as gradient clipping in stochastic methods to achieve logarithmic dependence on the confidence level in the high-probability bounds.

\section{Main Results for Minimization Problems}

\subsection{\algname{SGD} with Clipping}
We start with \algname{clipped-SGD}:
\begin{gather}
    x^{k+1} = x^k - \gamma \cdot \clip\left(\nabla f(x^k;\xi^k), \lambda_k\right), \label{ch2:eq:clipped_SGD_update_chapter_2}
\end{gather}
where $\xi^k$ is sampled from $\thesischtwocD_k$ independently from previous steps. We emphasize here and below that distribution of the noise is allowed to be dependent on $k$: we require just independence of $\xi^k$ from the previous steps. Our main convergence results for \algname{clipped-SGD} are summarized in the following theorem.

\begin{theorem}[Convergence of \algname{clipped-SGD}]\label{ch2:thm:clipped_SGD_main_theorem_chapter_2}
{\color{white}}
Let $K\geq 0$ and $\beta \in (0,1]$ be such that $A = \ln \frac{4(K+1)}{\beta} \geq 1$.

\noindent\textbf{Case 1.} Let Assumptions~\ref{ch2:as:bounded_alpha_moment_chapter_2}, \ref{ch2:as:lower_boundedness_chapter_2}, and \ref{ch2:as:L_smoothness_chapter_2} hold, with
\begin{align*}
Q &= \bigl\{x \in \R^d \bigm| \exists y \in \R^d:\ f(y) \leq f^{\inf} + 2\Delta,\\[-0.25\baselineskip]
  &\hspace{39mm} \|x-y\| \leq \nicefrac{\sqrt{\Delta}}{20\sqrt{L}}\bigr\},
\qquad \Delta \geq f(x^0) - f^{\inf},\\
0 < \gamma &\leq \cO\left(\min\left\{\nicefrac{1}{LA},
\nicefrac{\sqrt{\Delta}}{\sigma \sqrt{L}K^{\nicefrac{1}{{p}}}A^{\nicefrac{({p}-1)}{{p}}}}\right\}\right),
\qquad \lambda_k = \lambda = \Theta\left(\nicefrac{\sqrt{\Delta}}{\sqrt{L}\gamma A}\right).
\end{align*}

\noindent\textbf{Case 2.} Let Assumptions~\ref{ch2:as:bounded_alpha_moment_chapter_2}, \ref{ch2:as:L_smoothness_chapter_2}, and \ref{ch2:as:PL_chapter_2} hold, with $Q$ and $\Delta$ as in Case~1, and choose
\begin{align*}
0 < \gamma &= \cO\left(\min\left\{\nicefrac{1}{LA}, \nicefrac{\ln(B_K)}{\mu(K+1)}\right\}\right),\\
B_K &= \Theta\left(\max\left\{2,
\nicefrac{(K+1)^{\nicefrac{2({p}-1)}{{p}}}\mu^2 \Delta}
{L\sigma^2 A^{\nicefrac{2({p}-1)}{{p}}}\ln^2(B_K)}\right\}\right),\\
\lambda_k &= \Theta\left(\nicefrac{\exp(-\gamma\mu(1+\nicefrac{k}{2}))\sqrt{\Delta}}
{\sqrt{L}\gamma A}\right).
\end{align*}

\noindent\textbf{Case 3.} Let Assumptions~\ref{ch2:as:bounded_alpha_moment_chapter_2}, \ref{ch2:as:L_smoothness_chapter_2}, and \ref{ch2:as:str_cvx_chapter_2} hold with $\mu = 0$, $Q = B_{3R}(x^*)$, and $R \geq \|x^0-x^*\|$. Choose
\begin{align*}
0 < \gamma &\leq \cO\left(\min\left\{\nicefrac{1}{LA},
\nicefrac{R}{\sigma K^{\nicefrac{1}{{p}}}A^{\nicefrac{({p}-1)}{{p}}}}\right\}\right),
\qquad \lambda_k = \lambda = \Theta\left(\nicefrac{R}{\gamma A}\right).
\end{align*}

\noindent\textbf{Case 4.} Let Assumptions~\ref{ch2:as:bounded_alpha_moment_chapter_2}, \ref{ch2:as:L_smoothness_chapter_2}, and \ref{ch2:as:QSC_chapter_2} hold with $\mu > 0$, $Q = B_{3R}(x^*)$, and $R \geq \|x^0-x^*\|$. Choose
\begin{align*}
0 < \gamma &= \cO\left(\min\left\{\nicefrac{1}{LA}, \nicefrac{\ln(B_K)}{\mu(K+1)}\right\}\right),\\
B_K &= \Theta\left(\max\left\{2,
\nicefrac{(K+1)^{\nicefrac{2({p}-1)}{{p}}}\mu^2 R^2}
{\sigma^2 A^{\nicefrac{2({p}-1)}{{p}}}\ln^2(B_K)}\right\}\right),\\
\lambda_k &= \Theta\left(\nicefrac{\exp(-\gamma\mu(1+\nicefrac{k}{2}))R}{\gamma A}\right).
\end{align*}

With probability at least $1-\beta$, the respective guarantees are
\begin{align*}
\textbf{Case 1:}\quad &\frac{1}{K+1}\sum_{k=0}^K \|\nabla f(x^k)\|^2 \leq \varepsilon,\\
\textbf{Case 2:}\quad &f(x^K)-f(x^*) \leq \varepsilon,\\
\textbf{Case 3:}\quad &f(\bar{x}^K)-f(x^*) \leq \varepsilon,
\qquad \bar{x}^K = \frac{1}{K+1}\sum_{k=0}^K x^k,\\
\textbf{Case 4:}\quad &\|x^K-x^*\|^2 \leq \varepsilon.
\end{align*}
For these guarantees, \algname{clipped-SGD} requires
\begin{align}
    \text{{\bf \emph{Case 1}}:}&\quad \widetilde\cO\left(\max\left\{\frac{L \Delta}{\varepsilon}, \left(\frac{\sqrt{L\Delta} \sigma}{\varepsilon}\right)^{\frac{{p}}{{p}-1}}\right\}\right) \label{ch2:eq:SGD_main_result_non_cvx_chapter_2}\\
   \text{{\bf \emph{Case 2}}:}&\quad \widetilde\cO\left(\max\left\{\frac{L}{\mu}, \left(\frac{L\sigma^2}{\mu^2\varepsilon}\right)^{\frac{{p}}{2({p}-1)}}\right\}\right) \label{ch2:eq:SGD_main_result_PL_chapter_2}\\
   \text{{\bf \emph{Case 3}}:}&\quad \widetilde\cO\left(\max\left\{\frac{LR^2}{\varepsilon}, \left(\frac{\sigma R}{\varepsilon}\right)^{\frac{{p}}{{p}-1}}\right\}\right) \label{ch2:eq:SGD_main_result_cvx_chapter_2}\\
    \text{{\bf \emph{Case 4}}:}&\quad \widetilde\cO\left(\max\left\{\frac{L}{\mu}, \left(\frac{\sigma^2}{\mu^2\varepsilon}\right)^{\frac{{p}}{2({p}-1)}}\right\}\right) \label{ch2:eq:SGD_main_result_QSC_chapter_2}
\end{align}
oracle calls.
\end{theorem}

The complete formulation of the result and full proofs are deferred to Appendix~\ref{ch2:appendix:clipped_SGD_chapter_2}. As one can see from Table~\ref{ch2:tab:comparison_of_rates_minimization_chapter_2}, for ${p} = 2$ the derived complexity bounds match the best-known ones for \algname{clipped-SGD} in the setups where it was analyzed. Next, we emphasize that the second term under the maximum in \eqref{ch2:eq:SGD_main_result_QSC_chapter_2} (quasi-strongly convex functions) is optimal up to logarithmic factors \parencite{zhang2020adaptive}. In the convex case, there are no lower bounds, but we conjecture that the second term in \eqref{ch2:eq:SGD_main_result_cvx_chapter_2} is optimal (up to logarithms) in this case as well.

Next, in the case of P{\L}-functions, we are not aware of any high-probability convergence results in the literature. In the special case of ${p} = 2$, the derived complexity bound \eqref{ch2:eq:SGD_main_result_PL_chapter_2} matches the best-known in-expectation complexity bound for \algname{SGD} \parencite{karimi2016linear, khaled2020better} and the first term coincides (up to logarithms) with the lower bound for deterministic first-order methods in this setup \parencite{yue2022lower}.

Finally, in the non-convex case, bound \eqref{ch2:eq:SGD_main_result_non_cvx_chapter_2} is the first high-probability result under Assumption~\ref{ch2:as:bounded_alpha_moment_chapter_2} without the additional assumption of the boundedness of the gradients. For ${p} = 2$ it matches (up to logarithms) in-expectation lower bound \parencite{arjevani2022lower}. However, when ${p} < 2$, bound \eqref{ch2:eq:SGD_main_result_non_cvx_chapter_2} is inferior to the existing one $\widetilde\cO\left(\left(\nicefrac{G^2}{\varepsilon}\right)^{\nicefrac{(3{p}-2)}{2({p}-1)}}\right)$ by \textcite{cutkosky2021high}, which relies on the stronger assumption that $\Exp_{\xi\sim \thesischtwocD}\left[\|\nabla f(x;\xi)\|^{p}\right] \leq G^{p}$ for some $G > 0$ and all $x\in\R^d$, and also do not match the lower bound by \textcite{zhang2020adaptive} derived for $\Exp\left[\|\nabla f(x^k)\|\right]$, where $x^k$ is the output of the stochastic first-order method. It is also worth mentioning that \textcite{cutkosky2021high} use a different performance metric: $\hat\thesischtwocP_K = \left(\frac{1}{K+1}\sum_{k=0}^K\|\nabla f(x^k)\|\right)^2$. This metric is always smaller than $\thesischtwocP_K =\frac{1}{K+1}\sum_{k=0}^K\|\nabla f(x^k)\|^2$, which we use in our result. In the worst case, $\thesischtwocP_K$ can be $K+1$ times larger than $\hat\thesischtwocP_K$. Moreover, the lower bound of \textcite{zhang2020adaptive} is derived for $\Exp\left[\|\nabla f(x^k)\|\right]$ that is also always smaller than the standard quantity of interest $\Exp\left[\|\nabla f(x^k)\|^2\right]$. Therefore, the question of optimality of the bound \eqref{ch2:eq:SGD_main_result_non_cvx_chapter_2} remains open for ${p} < 2$. Moreover, it will also be interesting to modify our analysis in this case to derive a better bound for metric $\hat\thesischtwocP_K$ than \eqref{ch2:eq:SGD_main_result_non_cvx_chapter_2}.

\subsection{Acceleration}

Next, we focus on the accelerated version of \algname{clipped-SGD} called Clipped Stochastic Similar Triangles Method \algname{clipped-SSTM} \parencite{gorbunov2020stochastic}. The method constructs three sequences of points $\{x^k\}_{k \geq 0}$, $\{y^k\}_{k\geq 0}$, $\{z^k\}_{k \geq 0}$ satisfying the following update rules: $x^0 = y^0 = z^0$ and
\begin{align}
    x^{k+1} &= \frac{A_k y^k + a_{k+1} z^k}{A_{k+1}}, \label{ch2:eq:clipped_SSTM_x_update_chapter_2}\\
    z^{k+1} &= z^k - a_{k+1} \cdot \clip\left(\nabla f(x^{k+1};\xi^k), \lambda_k\right), \label{ch2:eq:clipped_SSTM_z_update_chapter_2}\\
    y^{k+1} &= \frac{A_k y^k + a_{k+1} z^{k+1}}{A_{k+1}}, \label{ch2:eq:clipped_SSTM_y_update_chapter_2}
\end{align}
where $A_0 = a_0 = 0$, $a_{k+1} = \frac{k+2}{2aL}$, $A_{k+1} = A_k + a_{k+1}$, and  $\xi^k$ is sampled from $\thesischtwocD_k$ independently from previous steps. Our main convergence result for \algname{clipped-SSTM} is given in the following theorem.

\begin{theorem}[Convergence of \algname{clipped-SSTM}]\label{ch2:thm:clipped_SSTM_main_theorem_chapter_2} Let Assumptions~\ref{ch2:as:bounded_alpha_moment_chapter_2}, \ref{ch2:as:L_smoothness_chapter_2}, \ref{ch2:as:str_cvx_chapter_2} with $\mu = 0$ hold for $Q = B_{3R}(x^*)$, $R \geq \|x^0 - x^*\|$ and
\[
    a = \Theta\left(\max\left\{A^2, \frac{\sigma K^{(p+1)/p}A^{(p-1)/p}}{LR}\right\}\right), \qquad \lambda_k = \Theta\left(\frac{R}{a_{k+1}A}\right),
\]
where $A = \ln(4K/\beta)$, $\beta \in (0,1]$ are such that $A \geq 1$. Then to guarantee $f(y^K) - f(x^*) \leq \varepsilon$ with probability $\geq 1 - \beta$ \algname{clipped-SSTM} requires
\begin{equation}
    \widetilde\cO\left(\max\left\{\sqrt{\frac{LR^2}{\varepsilon}}, \left(\frac{\sigma R}{\varepsilon}\right)^{\frac{{p}}{{p}-1}}\right\}\right) \quad \text{oracle calls.} \label{ch2:eq:clipped_SSTM_main_result_chapter_2}
\end{equation}
Moreover, with probability $\geq 1-\beta$ the iterates of \algname{clipped-SSTM} stay in the ball $B_{2R}(x^*)$: $\{x^k\}_{k=0}^{K+1}, \{y^k\}_{k=0}^{K}, \{z^k\}_{k=0}^{K} \subseteq B_{2R}(x^*)$.
\end{theorem}

The derived high-probability bound matches (see the proof in Appendix~\ref{ch2:appendix:clipped_SSTM_chapter_2}) the best-known one in the case of ${p} = 2$. For ${p} < 2$ there are no lower bounds in the convex case. However, the first term in \eqref{ch2:eq:clipped_SSTM_main_result_chapter_2} is optimal and matches the deterministic lower bound in the convex case \parencite{NemiYudin}. The second term is the same as in the bound for \algname{clipped-SGD} \eqref{ch2:eq:SGD_main_result_cvx_chapter_2} and we conjecture that it cannot be improved.

In the strongly convex case, we consider \algname{clipped-SSTM} with restarts (\algname{R-clipped-SSTM}). This method consists of $\tau$ stages. On the $t$-th stage \algname{R-clipped-SSTM} runs \algname{clipped-SSTM} for $K_0$ iterations from the starting point $\hat x^t$, which is the output from the previous stage ($\hat x^t = x^0$), and defines the obtained point as $\hat x^{t+1}$; see Algorithm~\ref{ch2:alg:R-clipped-SSTM_chapter_2} in Appendix~\ref{ch2:appendix:R_clipped_SSTM_chapter_2}. For this procedure we have the following result.

\begin{theorem}[Convergence of \algname{R-clipped-SSTM}]\label{ch2:thm:R_clipped_SSTM_main_theorem_chapter_2}
    Let Assumptions~\ref{ch2:as:bounded_alpha_moment_chapter_2}, \ref{ch2:as:L_smoothness_chapter_2}, \ref{ch2:as:str_cvx_chapter_2} with $\mu > 0$ hold for $Q = B_{3R}(x^*)$, $R \geq \|x^0 - x^*\|$ and \algname{R-clipped-SSTM} runs \algname{clipped-SSTM} $\tau = \lceil \log_2(\nicefrac{\mu R^2}{2\varepsilon}) \rceil$ times. Let $K_t = \widetilde\Theta(\max\{\sqrt{\nicefrac{LR_{t-1}^2}{\varepsilon_t}}, (\nicefrac{\sigma R_{t-1}}{\varepsilon_t})^{\nicefrac{{p}}{({p}-1)}}\})$, $a_t = \widetilde\Theta(\max\{1, \nicefrac{\sigma K_t^{\nicefrac{{p}+1}{{p}}}}{LR_{t-1}}\})$, $\lambda_k^t = \widetilde{\Theta}(\nicefrac{R_{t-1}}{a_{k+1}^t})$ be the parameters for the stage $t$ of \algname{R-clipped-SSTM}, where $R_{t-1} = 2^{-\nicefrac{(t-1)}{2}} R$, $\varepsilon_t = \nicefrac{\mu R_{t-1}^2}{4}$, $\ln \frac{4\tau K_t}{\beta} \geq 1$ for all $t=1,\ldots,\tau$ and some $\beta \in (0,1]$. Then to guarantee $f(\hat x^\tau) - f(x^*) \leq \varepsilon$ with probability $\geq 1 - \beta$ \algname{R-clipped-SSTM} requires
\begin{equation}
    \widetilde\cO\left(\max\left\{\sqrt{\frac{L}{\mu}}, \left(\frac{\sigma^2}{\mu\varepsilon}\right)^{\frac{{p}}{2({p}-1)}}\right\}\right) \quad \text{oracle calls.} \label{ch2:eq:R_clipped_SSTM_main_result_chapter_2}
\end{equation}
Moreover, with probability $\geq 1-\beta$ the iterates of \algname{R-clipped-SSTM} at stage $t$ stay in the ball $B_{2R_{t-1}}(x^*)$.
\end{theorem}

The obtained complexity bound (see the proof in Appendix~\ref{ch2:appendix:R_clipped_SSTM_chapter_2}) is the first optimal (up to logarithms) high-probability complexity bound under Assumption~\ref{ch2:as:bounded_alpha_moment_chapter_2} for the smooth strongly convex problems. Indeed, the first term cannot be improved in view of the deterministic lower bound by \textcite{NemiYudin}, and the second term is optimal due to \textcite{zhang2020adaptive}.

\section{Main Results for Variational Inequalities}

\subsection{Clipped Stochastic Extragradient}
For (quasi strongly) monotone VIPs we consider Clipped Stochastic Extragradient method (\algname{clipped-SEG}):
\begin{align}
    \thesischtwotx^k &= x^k - \gamma \cdot \clip(F(x^k;\xi_1^k), \lambda_k), \label{ch2:eq:clipped_SEG_extr_step_chapter_2}\\
    x^{k+1} &= x^k - \gamma \cdot \clip(F(\thesischtwotx^k;\xi_2^k), \lambda_k), \label{ch2:eq:clipped_SEG_upd_step_chapter_2}
\end{align}
where $\xi_1^k$, $\xi_2^k$ are sampled from $\thesischtwocD_k$ independently from previous steps. Our main convergence results for \algname{clipped-SEG} are summarized below.

\begin{theorem}[Convergence of \algname{clipped-SEG}]\label{ch2:thm:clipped_SEG_main_theorem_chapter_2}
\textbf{Case 1.} Let Assumptions~\ref{ch2:as:bounded_alpha_moment_chapter_2}, \ref{ch2:as:L_smoothness_chapter_2}, \ref{ch2:as:monotonicity_chapter_2} hold for $Q = B_{4R}(x^*)$ and
\[
    0 < \gamma = \cO\left(\min\left\{\frac{1}{LA}, \frac{R}{K^{1/p}\sigma A^{(p-1)/p}}\right\}\right), \qquad \lambda_k = \lambda = \Theta\left(\frac{R}{\gamma A}\right),
\]
where $A = \ln(6(K+1)/\beta) \geq 1$, $\beta \in (0,1]$.\par\noindent
\textbf{Case 2.} Let Assumptions~\ref{ch2:as:bounded_alpha_moment_chapter_2}, \ref{ch2:as:L_smoothness_chapter_2}, \ref{ch2:as:QSM_chapter_2} with $\mu > 0$ hold for $Q = B_{3R}(x^*)$ and
\begin{gather*}
    0 < \gamma =  \cO\left(\min\left\{\frac{1}{LA}, \frac{\ln(B_K)}{\mu(K+1)}\right\}\right), \qquad \lambda_k = \Theta\left(\frac{\exp(-\gamma\mu(1+k/2))R}{\gamma A}\right),\\
    B_K = \Theta\left(\max\left\{2, \frac{(K+1)^{2(p-1)/p}\mu^2 R^2}{\sigma^2 A^{2(p-1)/p}\ln^2(B_K)}\right\}\right),
\end{gather*}
where $A = \ln(6(K+1)/\beta)$, $\beta \in (0,1]$ are such that $A \geq 1$.\par\noindent
Then to guarantee $\gap_{R}(\thesischtwotx_{\text{avg}}^K) \leq \varepsilon$ in \textbf{Case 1} with $\thesischtwotx_{\text{avg}}^K = \frac{1}{K+1}\sum_{k=0}^K \thesischtwotx^k$, $\|x^K - x^*\|^2 \leq \varepsilon$ in \textbf{Case 2} with probability $\geq 1 - \beta$ \algname{clipped-SEG} requires
\begin{align}
    \text{{\bf \emph{Case 1}}:}&\quad \widetilde\cO\left(\max\left\{\frac{LR^2}{\varepsilon}, \left(\frac{\sigma R}{\varepsilon}\right)^{\frac{{p}}{{p}-1}}\right\}\right) \label{ch2:eq:SEG_main_result_monotone_chapter_2}\\
    \text{{\bf \emph{Case 2}}:}&\quad \widetilde\cO\left(\max\left\{\frac{L}{\mu}, \left(\frac{\sigma^2}{\mu^2\varepsilon}\right)^{\frac{{p}}{2({p}-1)}}\right\}\right) \label{ch2:eq:SEG_main_result_QSM_chapter_2}
\end{align}
oracle calls.
\end{theorem}

The proofs are deferred to Appendix~\ref{ch2:appendix:SEG_chapter_2}. When ${p} =2$, the above bounds recover state-of-the-art high-probability bounds for monotone and quasi-strongly monotone Lipschitz VIP \parencite{gorbunov2022clipped}. For the case of ${p} < 2$ \eqref{ch2:eq:SEG_main_result_monotone_chapter_2} and \eqref{ch2:eq:SEG_main_result_QSM_chapter_2} are the first high-probability results for the mentioned classes. Next, the first terms in these complexity bounds are optimal (up to logarithms) due to the lower bounds for the deterministic methods \parencite{ouyang2021lower, zhang2022lower}. The second term in \eqref{ch2:eq:SEG_main_result_QSM_chapter_2} is also optimal (up to logarithms) due to the lower bounds for stochastic strongly convex minimization \parencite{zhang2020adaptive}. Similarly to the convex case in minimization, we conjecture that the second term in \eqref{ch2:eq:SEG_main_result_monotone_chapter_2} cannot be improved in the monotone case as well.

\subsection{Clipped Stochastic Gradient Descent-Ascent}

In the star-cocoercive case, we focus on Clipped Stochastic Gradient Descent-Ascent (\algname{clipped-SGDA}):
\begin{equation}
    x^{k+1} = x^k - \gamma \cdot \clip(F(x^k;\xi^k), \lambda_k), \label{ch2:eq:clipped_SGDA_update_chapter_2}
\end{equation}
where $\xi^k$ is sampled from $\thesischtwocD_k$ independently from previous steps. For this method we derive the following convergence guarantees.
\begin{theorem}[Convergence of \algname{clipped-SGDA}]\label{ch2:thm:clipped_SGDA_main_theorem_chapter_2}
\noindent\textbf{Case 1.} Let Assumptions~\ref{ch2:as:bounded_alpha_moment_chapter_2}, \ref{ch2:as:star-cocoercivity_chapter_2}, and \ref{ch2:as:monotonicity_chapter_2} hold for $Q=B_{2R}(x^*)$. Choose
\begin{align*}
0 < \gamma &= \cO\left(\min\left\{\nicefrac{1}{\ell A},
\nicefrac{R}{K^{\nicefrac{1}{{p}}}\sigma A^{\nicefrac{({p}-1)}{{p}}}}\right\}\right),
\qquad \lambda_k = \lambda = \Theta\left(\nicefrac{R}{\gamma A}\right),
\end{align*}
where $\stepEightBchange{A = \ln \frac{6(K+1)}{\beta},}\; \beta \in (0,1]$ and $A \geq 1$.

\noindent\textbf{Case 2.} Let Assumptions~\ref{ch2:as:bounded_alpha_moment_chapter_2} and \ref{ch2:as:star-cocoercivity_chapter_2} hold for $Q=B_{2R}(x^*)$. Choose
\begin{align*}
0 < \gamma &= \cO\left(\min\left\{\nicefrac{1}{\ell A},
\nicefrac{R}{K^{\nicefrac{1}{{p}}}\sigma A^{\nicefrac{({p}-1)}{{p}}}}\right\}\right),
\qquad \lambda_k = \lambda = \Theta\left(\nicefrac{R}{\gamma A}\right),
\end{align*}
where $A = \ln \frac{4(K+1)}{\beta}$, $\beta \in (0,1]$, and $A \geq 1$.

\noindent\textbf{Case 3.} Let Assumptions~\ref{ch2:as:bounded_alpha_moment_chapter_2}, \ref{ch2:as:star-cocoercivity_chapter_2}, and \ref{ch2:as:QSM_chapter_2} hold with $\mu > 0$ and $Q=B_{2R}(x^*)$. Choose
\begin{align*}
0 < \gamma &= \cO\left(\min\left\{\nicefrac{1}{\ell A}, \nicefrac{\ln(B_K)}{\mu(K+1)}\right\}\right),\\
B_K &= \Theta\left(\max\left\{2,
\nicefrac{(K+1)^{\nicefrac{2({p}-1)}{{p}}}\mu^2 R^2}
{\sigma^2 A^{\nicefrac{2({p}-1)}{{p}}}\ln^2(B_K)}\right\}\right),\\
\lambda_k &= \Theta\left(\nicefrac{\exp(-\gamma\mu(1+\nicefrac{k}{2}))R}{\gamma A}\right),
\end{align*}
where $A = \ln \frac{4(K+1)}{\beta}$, $\beta \in (0,1]$, and $A \geq 1$.

With probability at least $1-\beta$, the respective guarantees are
\begin{align*}
\textbf{Case 1:}\quad &\gap_{R}(\thesischtwotx_{\text{avg}}^K) \leq \varepsilon,
\qquad \thesischtwotx_{\text{avg}}^K = \frac{1}{K+1}\sum_{k=0}^K \thesischtwotx^k,\\
\textbf{Case 2:}\quad &\frac{1}{K+1}\sum_{k=0}^K \|F(x^k)\|^2 \leq \ell\varepsilon,\\
\textbf{Case 3:}\quad &\|x^K-x^*\|^2 \leq \varepsilon.
\end{align*}
For these guarantees, \algname{clipped-SGDA} requires
\begin{align}
    \text{{\bf \emph{Case 1}} and {\bf \emph{2}}:}&\quad \widetilde\cO\left(\max\left\{\frac{\ell R^2}{\varepsilon}, \left(\frac{\sigma R}{\varepsilon}\right)^{\frac{{p}}{{p}-1}}\right\}\right) \label{ch2:eq:SGDA_main_result_monotone_and_star_cocoercive_chapter_2} \\
    \text{{\bf \emph{Case 3}}:}&\quad \widetilde\cO\left(\max\left\{\frac{\ell}{\mu}, \left(\frac{\sigma^2}{\mu^2\varepsilon}\right)^{\frac{{p}}{2({p}-1)}}\right\}\right) \label{ch2:eq:SGDA_main_result_QSM_chapter_2}
\end{align}
oracle calls.
\end{theorem}

One can find the proofs in Appendix~\ref{ch2:appendix:SGDA_chapter_2}. The derived high-probability results generalize the existing state-of-the-art results from the case of ${p} = 2$ \parencite{gorbunov2022clipped} to the case of ${p} < 2$.

\section{Key Lemma and Intuition Behind the Proofs}

The proofs of the results in this chapter follow a similar pattern. To illustrate the main idea, we consider the analysis of \algname{clipped-SGD} in the non-convex case. Mimicking the proof of deterministic \algname{GD} we derive the following inequality:
\begin{equation}
    \gamma\sum\limits_{k=0}^K\|\nabla f(x^k)\|^2 \lessapprox \Delta_0 - \Delta_{K+1} - \gamma \sum\limits_{k=0}^K \langle \nabla f(x^k), \theta_k \rangle + L\gamma^2 \sum\limits_{k=0}^{K-1} \|\theta_k\|^2.\label{ch2:eq:induction_inequality_example_chapter_2}
\end{equation}
where $\Delta_k = f(x^k) - f^{\inf}$ and $\theta_k = \thesischtwotnabla f(x^{k};\xi^k) - \nabla f(x^{k})$. In other words, we separate the deterministic part of the method from its stochastic part. To obtain the result of Theorem~\ref{ch2:thm:clipped_SGD_main_theorem_chapter_2} (Case 1) it remains to upper bound with high-probability the sums from the second line of the formula above. We do it with the help of Bernstein's inequality (Lemma~\ref{ch2:lem:Bernstein_ineq_chapter_2}). However, it requires several preliminary steps. In particular, Bernstein's inequality needs the random variables to be bounded. The magnitudes of summands depend on $\nabla f(x^k)$ that can be arbitrarily large due to the stochasticity in $x^k$. However, \eqref{ch2:eq:induction_inequality_example_chapter_2} allows to bound $\Delta_{K+1}$ inductively and, using smoothness, to bound $\|\nabla f(x^{K+1})\|$. Secondly, Bernstein's inequality requires knowing the bounds on the bias and variance of the clipped stochastic estimator. For such purposes, we derive the following result.
\begin{lemma}\label{ch2:lem:bias_and_variance_clip_chapter_2}
    Let $X$ be a random vector in $\R^d$ and $\thesischtwotX = \clip(X,\lambda)$. Then, $\|\thesischtwotX - \Exp[\thesischtwotX]\| \leq 2\lambda$. 
    Moreover, if for some $\sigma \geq 0$ and ${p} \in (1,2]$ we have $\Exp[X] = x\in\R^d$, $\Exp[\|X - x\|^{p}] \leq \sigma^{p}$, 
    and $\|x\| \leq \nicefrac{\lambda}{2}$, then
    \begin{eqnarray}
        \left\|\Exp[\thesischtwotX] - x\right\| &\leq& \frac{2^{p}\sigma^{p}}{\lambda^{{p}-1}}, \label{ch2:eq:bias_X_main_chapter_2}\\
        \Exp\left[\left\|\thesischtwotX - x\right\|^2\right] &\leq& 18 \lambda^{2-{p}} \sigma^{p}, \label{ch2:eq:distortion_X_main_chapter_2}\\
        \Exp\left[\left\|\thesischtwotX - \Exp[\thesischtwotX]\right\|^2\right] &\leq& 18 \lambda^{2-{p}} \sigma^{p}. \label{ch2:eq:variance_X_main_chapter_2}
    \end{eqnarray}
\end{lemma}

This lemma can be useful on its own for analyses involving clipping operators. Moreover, our high-probability analysis does not rely on the choice of clipping explicitly: in the proofs, we use only $\|\thesischtwotX\| \leq \lambda$ and inequalities \eqref{ch2:eq:bias_X_main_chapter_2}-\eqref{ch2:eq:variance_X_main_chapter_2}. Therefore, the results hold for the methods considered in this chapter with any other non-linearity $\phi_\lambda(x)$ (not necessarily clipping), if it satisfies the conditions from the above lemma for $\thesischtwotX = \phi_\lambda(X)$.

\section{Conclusion and Future Work}

This chapter developed high-probability methods for stochastic minimization and variational inequalities when the oracle noise has only a bounded central \(p\)-th moment, \(p\in(1,2]\), and may therefore have unbounded variance. Clipping converts rare, extreme oracle errors into controlled increments, enabling logarithmic dependence on the failure probability for \algname{clipped-SGD}, accelerated \algname{clipped-SSTM}, \algname{clipped-SEG}, and \algname{clipped-SGDA}. The analysis applies to several smooth optimization and structured variational-inequality classes, with assumptions imposed locally in many regimes.

The guarantees refer to different outputs and should not be read as one universal last-iterate statement. The nonconvex result controls the realized average of squared gradient norms; convex results use an averaged iterate, while the accelerated convex method reports \(y^K\). Polyak--{\L}ojasiewicz, quasi-strongly convex, and quasi-strongly monotone regimes provide terminal-iterate criteria. Monotone variational inequalities instead use restricted gaps or averaged operator norms. The restarted accelerated strongly convex complexity is optimal up to logarithmic factors, and several deterministic or quasi-strongly structured terms match the relevant lower bounds. Other cases, including the \(p<2\) nonconvex metric comparison and stochastic terms for convex or monotone problems, do not yet have matching lower bounds.

The prescribed stepsizes and clipping levels still require horizon, confidence, noise, smoothness, and radius information. Extending the framework to H\"older-continuous gradients, other robust nonlinearities, and less parameter-dependent schedules remains open. Chapter~\ref{ch:chapter_3} next addresses composite and distributed problems, where clipping a complete local gradient can destroy the fixed-point structure and shifted innovation clipping becomes necessary.

\newpage
  \endgroup

\chapter{High-Probability Convergence for Composite and Distributed Stochastic Minimization and Variational Inequalities with Heavy-Tailed Noise}
\label{ch:chapter_3}
\thesischaptersourcenote
  {High-Probability Convergence for Composite and Distributed Stochastic Minimization and Variational Inequalities with Heavy-Tailed Noise}
  {gorbunov2023high}
  \begingroup
  \thesischaptersetup{3}%
  
\begin{paperabstract}
High-probability analysis of stochastic first-order optimization methods under mild assumptions on the noise has been gaining a lot of attention in recent years. Typically, gradient clipping is one of the key algorithmic ingredients to derive good high-probability guarantees when the noise is heavy-tailed. However, if implemented na\"ively, clipping can spoil the convergence of the popular methods for composite and distributed optimization (\algname{Prox-SGD}/\algname{Parallel SGD}) even in the absence of any noise. Due to this reason, many works on high-probability analysis consider only unconstrained non-distributed problems, and the existing results for composite/distributed problems do not include some important special cases (like strongly convex problems) and are not optimal. To address this issue, we propose new stochastic methods for composite and distributed optimization based on the clipping of stochastic gradient differences and prove tight high-probability convergence results (including nearly optimal ones) for the new methods. In addition, we also develop new methods for composite and distributed variational inequalities and analyze the high-probability convergence of these methods.
\end{paperabstract}

\section{Introduction}
Many recent works on stochastic optimization have the ultimate goal of bridging the theory and practice in machine learning. This is mostly reflected in the attempts at the theoretical analysis of optimization methods under weaker assumptions than the standard ones. Moreover, as the motivating example of \textcite{gorbunov2020stochastic} shows, some phenomena cannot be explained using classical in-expectation convergence analysis, which results in the growing interest in more accurate ways to the analysis of stochastic methods, for example, \emph{high-probability convergence analysis}.

However, despite the significant attention to this topic \parencite{nazin2019algorithms, davis2021low, gorbunov2020stochastic, gorbunov2022clipped, cutkosky2021high, sadiev2023high, nguyen2023high, liu2023stochastic, liu2023high}, several important directions remain unexplored. In particular, all mentioned works either consider unconstrained problems or consider general composite/constrained minimization/variational inequality problems but have some noticeable limitations, such as bounded domain assumption, extra logarithmic factors \pd{in the complexity bounds}, not optimal (not accelerated) \pd{convergence} rates, or no analysis of (quasi-) strongly convex (monotone) case. The importance of composite/constrained formulations for the machine learning community can be justified in many ways. For example, composite optimization and distributed optimization have a lot of similarities, i.e., one can view a distributed optimization problem as a special composite optimization problem \parencite{parikh2014proximal}. Due to the large sizes of modern machine learning models and datasets, many important problems can be solved in a reasonable time only via distributed methods. Moreover, some problems have a distributed nature, e.g., data determining the problem can be (privately) stored on multiple devices connected with a central server, which is a classical setup for Federated Learning \parencite{konevcny2016federated, kairouz2021advances}. Next, composite formulations are very useful for handling different regularizations popular in machine learning and statistics \parencite{zou2005regularization, shalev2014understanding, beck2017first}. Finally, variational inequalities are usually considered with constraints as well.

The discrepancy between the importance of composite/constrained and distributed formulations and the lack of high-probability convergence results in this setup can be partially explained as follows. State-of-the-art high-probability convergence results are derived for the algorithms that use \emph{gradient clipping} \parencite{pascanu2013difficulty}, i.e., the clipping operator defined as $\clip(x,\lambda) = \min\{1, \nicefrac{\lambda}{\|x\|}\}x$ for $x\neq 0$ and $\clip(0,\lambda) = 0$ with some clipping level $\lambda > 0$ is applied to the stochastic gradients. If $\lambda$ is too small, then na\"ive Proximal Gradient Descent with gradient clipping is not a fixed point method, i.e., the method escapes the solution even if it is initialized there (see a technical explanation in Section~\ref{ch3:sec:the_issue_and_star_method_chapter_3}). This fact implies that one has either to increase the clipping level or to decrease the stepsize to converge to the exact solution asymptotically; the latter approach leads to a slower convergence rate. On the other hand, even in the unconstrained case, the existing results with acceleration/linear convergence are derived for the methods using decreasing clipping level \parencite{gorbunov2020stochastic,sadiev2023high}. Therefore, new algorithms and analyses are required to handle this issue.

This chapter addresses the gap by proposing new stochastic methods for composite and distributed problems via the clipping of \emph{gradient differences} that converge to zero with high probability. This allows us to achieve the desirable acceleration and linear convergence. Before presenting the main contributions, we introduce the problem settings formally.

\subsection{Setup}

\paragraph{Notation.} The standard Euclidean norm of vector $x \in \R^d$ is denoted as $\|x\| = \sqrt{\langle x,x \rangle}$. $B_R(x) = \{y\in \R^d\mid \|y-x\| \leq R\}$ is the ball centered at $x$ with radius $R$. Bregman divergence w.r.t.\ function $f$ is denoted as $D_f(x,y) \eqdef f(x) - f(y) - \langle \nabla f(y), x - y \rangle$. In $\cO(\cdot)$, we omit the numerical factors, and in $\widetilde\cO(\cdot)$, we omit numerical and logarithmic factors. For natural $n\geq 1$ the set $\{1,2,\ldots,n\}$ is denoted as $[n]$. We use $\Exp_{\xi}[\cdot]$ to denote the expectation w.r.t.\ the randomness coming from $\xi$.

\paragraph{Considered problems.} This chapter focuses on stochastic composite minimization problems:
\begin{equation}
    \min\limits_{x\in\R^d} \left\{\Phi(x) = f(x) + \Psi (x)\right\}, \label{ch3:eq:min_problem_chapter_3}
\end{equation}
where $f(x) = \Exp_{\xi\sim\thesischthreecD}[f(x;\xi)]$ is a differentiable function satisfying some properties to be defined later and $\Psi(x)$ is a proper, closed, convex function (composite/regularization term). The examples of problem \eqref{ch3:eq:min_problem_chapter_3} arise in various applications, e.g., machine learning \parencite{shalev2014understanding}, signal processing \parencite{combettes2011proximal}, image processing \parencite{luke2020proximal}. We also consider variational inequality problems, see Appendix~\ref{ch3:appendix:VIs_chapter_3}.

The distributed version of \eqref{ch3:eq:min_problem_chapter_3} has the following structure:
\begin{equation}
    f(x) = \frac{1}{n}\sum\limits_{i=1}^n \left\{f_i(x) = \Exp_{\xi_i\sim \thesischthreecD_i}[f_i(x;\xi_i)]\right\}. \label{ch3:eq:distributed_problem_chapter_3}
\end{equation}
In this case, there are $n$ workers connected in a centralized way with some parameter server; worker $i$ can query some noisy information (stochastic gradients/estimates) about $f_i$.

\paragraph{In-expectation and high-probability convergence.} In-expectation convergence guarantees provide the upper bounds on the number of iterations/oracle calls $\hat K = \hat K(\varepsilon)$ for a method needed to find point $x^{\hat K}$ such that $\Exp[\thesischthreecC(x^{\hat K})] \leq \varepsilon$ for given convergence criterion $\thesischthreecC(x)$ (e.g., $\thesischthreecC(x)$ can be $f(x) - f(x^*)$, $\|x - x^*\|^2$, $\|\nabla f(x)\|^2$) and given accuracy $\varepsilon > 0$. High-probability convergence guarantees give the upper bounds on the  number of iterations/oracle calls $K = K(\varepsilon, \beta)$ for a method needed to find point $x$ such that $\Prob\{\thesischthreecC(x^K) \leq \varepsilon\} \geq 1 - \beta$, where $\beta \in (0,1)$ is a confidence level. It is worth noting that Markov's inequality implies $\Prob\{\thesischthreecC(x^K) > \varepsilon\} < \nicefrac{\Exp[\thesischthreecC(x^K)]}{\varepsilon}$, meaning that it is sufficient to take $K = \hat K(\beta\varepsilon) = \hat K$: $\Prob\{\thesischthreecC(x^{\hat K}) > \varepsilon\} < \nicefrac{\Exp[\thesischthreecC(x^{\hat K})]}{\varepsilon} \leq \beta$. However, this typically leads to the polynomial dependence on $\nicefrac{1}{\beta}$ that significantly spoils the complexity of the method when $\beta$ is small. Therefore, we focus on the high-probability convergence guarantees that depend on $\nicefrac{1}{\beta}$ poly-logarithmically. Moreover, such high-probability results are more sensitive to the noise distribution (and, thus, more accurate) than in-expectation ones \parencite{gorbunov2020stochastic, sadiev2023high}.

\paragraph{Proximal operator.} We assume that function $\Psi(x)$ has a relatively simple structure such that one can efficiently compute \emph{proximal operator}: $\thesischthreeprox_{\gamma\Psi}(x) = \arg\min_{y\in \R^d}\{\gamma\Psi(y) + \tfrac{1}{2}\|y - x\|^2\}$. For the properties of the proximal operator and examples of functions $\Psi(x)$ such that $\thesischthreeprox_{\gamma\Psi}(x)$ can be easily computed, we refer the reader to \textcite{beck2017first}.

\paragraph{Bounded central ${p}$-th moment.} We consider the situation when $f_i$ and $F_i$ are accessible through the stochastic oracle calls. The stochastic estimates satisfy the following assumption.\footnote{Following \textcite{sadiev2023high}, we consider all assumptions only on some bounded set $Q \subseteq \R^d$; the diameter of $Q$ depends on the starting point. We emphasize that we do not assume boundedness of the domain of the original problem. Instead, we prove via induction that the iterates of the considered methods stay in some ball around the solution with high probability (see the details in Section~\ref{ch3:sec:proofs_structure_chapter_3}). Thus, it is sufficient for us to assume everything just on this ball, though our analysis remains unchanged if we introduce all assumptions on the whole domain.}
\begin{assumption}\label{ch3:as:bounded_alpha_moment_chapter_3}
    There exist some set $Q \subseteq \R^d$ and values $\sigma \geq 0$, ${p} \in (1,2]$ such that for all $x \in Q$ we have $\Exp_{\xi_i\sim \thesischthreecD_i}[\nabla f_i(x;\xi_i)] = \nabla f_i(x)$ and
    \begin{equation}
        \Exp_{\xi_i\sim \thesischthreecD_i}[\left\|\nabla f_i(x;\xi_i) - \nabla f_i(x)\right\|^{p}] \leq \sigma^{p}. \label{ch3:eq:bounded_alpha_moment_gradient_chapter_3}
    \end{equation}
\end{assumption}
For ${p} = 2$, Assumption~\ref{ch3:as:bounded_alpha_moment_chapter_3} reduces to the bounded variance assumption, and for ${p} \in (1,2)$ variance of the stochastic estimator can be unbounded, e.g., the noise can have L\'evy ${p}$-stable distribution \parencite{zhang2020adaptive}, which is heavy-tailed.

\paragraph{Assumptions on $f_i$.} We assume that functions $\{f_i\}_{i\in [n]}$ are $L$-smooth.
\begin{assumption}\label{ch3:as:L_smoothness_chapter_3}
    We assume that there exist some set $Q \subseteq \R^d$ and constant $L > 0$ such that for all $x, y \in Q, i\in [n]$ and for all $x^* \in \arg\min_{x\in \R^d} \Phi(x)$
    \begin{eqnarray}
        \|\nabla f_i(x) - \nabla f_i(y)\| &\leq& L\|x - y\|, \label{ch3:eq:L_smoothness_chapter_3}\\
        \|\nabla f_i(x) - \nabla f_i(x^*)\|^2 &\leq& 2LD_{f_i}(x, x^*). \label{ch3:eq:L_smoothness_cor_2_chapter_3}
    \end{eqnarray}
\end{assumption}
As noted by \textcite[Appendix B]{sadiev2023high}, \eqref{ch3:eq:L_smoothness_cor_2_chapter_3} is satisfied on the set $Q \neq \R^d$ if \eqref{ch3:eq:L_smoothness_chapter_3} holds on a slightly larger set in the case of $\Psi \equiv 0$, $n = 1$ (unconstrained single-node case). For simplicity, we assume that both \eqref{ch3:eq:L_smoothness_chapter_3} and \eqref{ch3:eq:L_smoothness_cor_2_chapter_3} hold on $Q$. This is always the case for $L$-smooth functions on $Q = \R^d$ when $\Psi \equiv 0$, $n = 1$. In a more general situation, condition \eqref{ch3:eq:L_smoothness_cor_2_chapter_3} can be viewed as an assumption on the structured non-convexity of $\{f_i\}_{i\in [n]}$. Finally, if $\{f_i\}_{i\in [n]}$ are convex and $L$-smooth on the whole domain of the problem \eqref{ch3:eq:min_problem_chapter_3}, then Assumption~\ref{ch3:as:L_smoothness_chapter_3} holds.

Next, for each particular result about the convergence of methods for \eqref{ch3:eq:min_problem_chapter_3}, we make one of the following assumptions.
\begin{assumption}\label{ch3:as:str_cvx_chapter_3}
    There exist some set $Q \subseteq \R^d$ and constant $\mu \geq 0$ such that $f$ is $\mu$-strongly convex: $\forall \, x, y \in Q$
    \begin{equation}
        f(y) \geq f(x) + \langle \nabla f(x), y - x \rangle + \frac{\mu}{2}\|y - x\|^2. \label{ch3:eq:str_cvx_chapter_3}
    \end{equation}
    When $\mu = 0$, function $f$ is called convex on $Q$.
\end{assumption}
This is a standard assumption for optimization literature \parencite{nesterov2018lectures}. We also consider a relaxation of strong convexity.
\begin{assumption}\label{ch3:as:QSC_chapter_3}
    There exist some set $Q \subseteq \R^d$ and constant $\mu \geq 0$ such that $f_1,\ldots,f_n$ are $(\mu, x^*)$-quasi-strongly convex for all $x^* \in \arg\min_{x\in \R^d} \Phi(x)$: $\forall \, x \in Q,\, i\in[n]$
    \begin{equation}
        f_i(x^*) \geq f_i(x) + \langle \nabla f_i(x), x^* - x \rangle + \frac{\mu}{2}\|x - x^*\|^2. \label{ch3:eq:QSC_chapter_3}
    \end{equation}
\end{assumption}
Condition \eqref{ch3:eq:QSC_chapter_3} is weaker than \eqref{ch3:eq:str_cvx_chapter_3} and holds even for some non-convex functions \parencite{necoara2019linear}.

\subsection{Our Contributions}

\begin{table*}[t]
    \centering
    \scriptsize
    \caption[High-probability complexity results for composite stochastic minimization]{\scriptsize Summary of known and new high-probability complexity results for solving (non-) composite (non-) distributed smooth optimization problem \eqref{ch3:eq:min_problem_chapter_3}. Column ``Setup'' indicates the assumptions made in addition to Assumptions~\ref{ch3:as:bounded_alpha_moment_chapter_3} and \ref{ch3:as:L_smoothness_chapter_3}. All assumptions are made only on some ball around the solution with radius $\sim R \geq \|x^0 - x^*\|$. Complexity is the number of stochastic oracle calls (per worker) needed for a method to guarantee that $\Prob\{\text{Metric} \leq \varepsilon\} \geq 1 - \beta$ for some $\varepsilon> 0$, $\beta \in (0,1]$ and ``Metric'' is taken from the corresponding column. Numerical and logarithmic factors are omitted for simplicity. Column ``C?'' shows whether the problem \eqref{ch3:eq:min_problem_chapter_3} is composite, ``D?'' indicates whether the problem \eqref{ch3:eq:min_problem_chapter_3} is distributed. Notation: $L$ = Lipschitz constant; $\sigma$ = parameter from Assumption~\ref{ch3:as:bounded_alpha_moment_chapter_3}; $R$ = any upper bound on $\|x^0 - x^*\|$; $\zeta_*$ = any upper bound on $\sqrt{\frac{1}{n}\sum_{i=1}^n \|\nabla f_i(x^*)\|^2}$; $\widehat{R}^2 = R\left(3R + L^{-1}(2\eta \sigma + \|\nabla f(x^0)\|)\right)$ for some $\eta > 0$ in the result of \textcite{nguyen2023improved}; one can show that $\widehat{R}^2 = \Theta(R^2 + \nicefrac{R\zeta_*}{L})$ when $n=1$, see the discussion after Theorem~\ref{ch3:thm:D_prox_clipped_SGD_main_cvx_chapter_3}; %
    $\mu$ = (quasi-)strong convexity parameter. The results established in this chapter are highlighted in blue.
    }
    \label{ch3:tab:comparison_of_rates_minimization_chapter_3}
    \begin{threeparttable}
        \resizebox{0.95\textwidth}{!}{%
        \begin{tabular}{|c|c c c c c|}
        \hline
Setup & Method & Metric & Complexity & C? & D?\\
        \hline\hline
        \multirow{12}{*}{\makecell{As.~\ref{ch3:as:str_cvx_chapter_3}\\ ($\mu = 0$)}} &\begin{tabular}{c}
             \algname{clipped-SGD} \\ \tabcite{sadiev2023high}
        \end{tabular}  & $f(\overline{x}^K) - f(x^*)$ & $\max\left\{\frac{LR^2}{\varepsilon}, \left(\frac{\sigma R}{\varepsilon}\right)^{\frac{{p}}{{p}-1}}\right\}$ & \thesischthreexmark & \thesischthreexmark \\
        & \begin{tabular}{c}
             \algname{clipped-SSTM} \\ \tabcite{sadiev2023high}
        \end{tabular} & $f(y^K) - f(x^*)$ & $\max\left\{\sqrt{\frac{LR^2}{\varepsilon}}, \left(\frac{\sigma R}{\varepsilon}\right)^{\frac{{p}}{{p}-1}}\right\}$ &\thesischthreexmark &\thesischthreexmark\\
        & \begin{tabular}{c}
             \algname{Clipped-SMD}\tnote{\color{blue}(1),(2)} \\ \tabcite{nguyen2023improved}
        \end{tabular} &$\Phi(\overline{x}^K) - \Phi(x^*)$ &$\max\left\{\frac{L \widehat R^2}{\varepsilon}, \left(\frac{\sigma R}{\varepsilon}\right)^{\frac{{p}}{{p}-1}}\right\}$ &\thesischthreecmark &\thesischthreexmark \\
        &\begin{tabular}{c}
             \algname{Clipped-ASMD}\tnote{\color{blue}(1)} \\ \tabcite{nguyen2023improved}
        \end{tabular} &$\Phi(y^K) - \Phi(x^*)$ &$\max\left\{\sqrt{\frac{L R^2}{\varepsilon}}, \left(\frac{\sigma R}{\varepsilon}\right)^{\frac{{p}}{{p}-1}}\right\}$ &\thesischthreecmark\thesischthreexmark \tnote{\color{blue}(3)}  &\thesischthreexmark \\
        & \cellcolor{linen} \begin{tabular}{c}
        
             \algname{DProx-clipped-SGD-shift} \\ Theorem~\ref{ch3:thm:D_prox_clipped_SGD_main_cvx_chapter_3}
        \end{tabular}  & \cellcolor{linen}$\Phi(\overline{x}^K) - \Phi(x^*)$ & \cellcolor{linen}$\max\left\{\frac{LR^2}{\varepsilon}, \frac{R\zeta_*}{\sqrt{n}\varepsilon}, \frac{1}{n}\left(\frac{\sigma R}{\varepsilon}\right)^{\frac{{p}}{{p}-1}}\right\}$ &\cellcolor{linen} \thesischthreecmark &\cellcolor{linen} \thesischthreecmark \\
        & \cellcolor{linen} \begin{tabular}{c}

             \algname{DProx-clipped-SSTM-shift} \\ Theorem~\ref{ch3:thm:D_prox_clipped_SSSTM_main_chapter_3}
        \end{tabular} &\cellcolor{linen} $\Phi(y^K) - \Phi(x^*)$ &\cellcolor{linen} $\max\left\{\sqrt{\frac{LR^2}{\varepsilon}}, \sqrt{\frac{R\zeta_*}{\sqrt{n}\varepsilon}}, \frac{1}{n}\left(\frac{\sigma R}{\varepsilon}\right)^{\frac{{p}}{{p}-1}}\right\}$ &\cellcolor{linen} \thesischthreecmark &\cellcolor{linen} \thesischthreecmark \\
        \hline
        \multirow{3}{*}{\makecell{As.~\ref{ch3:as:QSC_chapter_3}\\ ($\mu > 0$)}} &  \begin{tabular}{c}
             \algname{clipped-SGD} \\ \tabcite{sadiev2023high}
        \end{tabular} & $\|x^K - x^*\|^2$ & $\max\left\{\frac{L}{\mu}, \left(\frac{\sigma^2}{\mu^2\varepsilon}\right)^{\frac{{p}}{2({p}-1)}}\right\}$ &\thesischthreexmark
        &\thesischthreexmark \\
        & \cellcolor{linen} \begin{tabular}{c}
             
             \algname{DProx-clipped-SGD-shift} \\ Theorem~\ref{ch3:thm:D_prox_clipped_SGD_main_chapter_3}
        \end{tabular} &\cellcolor{linen} $\|x^K - x^*\|^2$ &\cellcolor{linen} $\max\left\{\frac{L}{\mu}, \stepEightBchange{\frac{\zeta_*}{\sqrt{n}\mu R},} \frac{1}{n}\left(\frac{\sigma^2}{\mu^2\varepsilon}\right)^{\frac{{p}}{2({p}-1)}}\right\}$ &\cellcolor{linen}\thesischthreecmark &\cellcolor{linen}\thesischthreecmark \\
        \hline
    \end{tabular}}
    \par\begin{minipage}{0.95\textwidth}
        \scriptsize\raggedright\setlength{\parindent}{0pt}
        {\color{blue}(1)} All assumptions are made on the whole domain.\par
        {\color{blue}(2)} The authors additionally assume that for a chosen point $\widehat x$ from the domain and for $\eta > 0$ one can compute an estimate $\widehat g$ such that $\Prob\lbrace \| \widehat g - \nabla f(\hat x) \| > \eta \sigma \rbrace \leq \varepsilon$. Such an estimate can be found using geometric median of $\cO(\ln \varepsilon^{-1})$ samples \parencite{minsker2015geometric}.\par
        {\color{blue}(3)} The authors assume that $\nabla f(x^*) = 0$, which is not true for general composite optimization.
    \end{minipage}\par
    \end{threeparttable}
\end{table*}

$\bullet$ \textbf{Methods with clipping of gradient differences for distributed composite minimization.} We develop two stochastic methods for composite minimization problems -- Proximal Clipped \algname{SGD} with shifts (\algname{Prox-clipped-SGD-shift}) and Proximal Clipped Similar Triangles Method with shifts (\algname{Prox-clipped-SSTM-shift}). Instead of clipping stochastic gradients, these methods clip the difference between the stochastic gradients and the shifts that are updated on the fly. This trick allows us to use decreasing clipping levels, and, as a result, we derive the first accelerated high-probability convergence rates and tight high-probability convergence rates for the non-accelerated method in the quasi-strongly convex case. We also generalize the proposed methods to the distributed case (\algname{DProx-clipped-SGD-shift} and \algname{DProx-clipped-SSTM-shift}) and prove that they benefit from parallelization. To the best of our knowledge, our results are the first showing linear speed-up under Assumption~\ref{ch3:as:bounded_alpha_moment_chapter_3}.

$\bullet$ \textbf{Methods with clipping of gradient differences for distributed composite VIPs.} We also apply the proposed trick to the methods for variational inequalities. In particular, we propose \algname{DProx-clipped-SGDA-shifts} and \algname{DProx-clipped-SEG-shifts} and rigorously analyze their high-probability convergence. As in the minimization case, the proposed methods have provable benefits from parallelization.

$\bullet$ \textbf{Tight convergence rates.} As a separate contribution, we highlight the tightness of our analysis: in the known special cases ($\Psi \equiv 0$ and/or $n=1$), the derived complexity bounds either recover or outperform previously known ones (see Table~\ref{ch3:tab:comparison_of_rates_minimization_chapter_3} and also Table~\ref{ch3:tab:comparison_of_rates_VIP_chapter_3} in the appendix). Moreover, in certain regimes, the results have optimal (up to logarithms) dependencies on $\varepsilon$. This is achieved under quite general assumptions.

\subsection{Closely Related Work}

We discuss closely related work here and defer additional discussion to Appendix~\ref{ch3:appendix:extra_related_chapter_3}.

\paragraph{High-probability bounds for unconstrained convex problems.} Standard high-probability convergence results are obtained under the so-called light-tails assumption (sub-Gaussian noise) \parencite{nemirovski2009robust, juditsky2011solving, ghadimi2012optimal}. The first work addressing this limitation is \textcite{nazin2019algorithms}, where the authors derive the first high-probability complexity bounds for the case of minimization on a bounded set under bounded variance assumption. In the unconstrained case, these results are extended and accelerated by \textcite{gorbunov2020stochastic} for smooth convex and strongly convex minimization problems. \textcite{gorbunov2021near} tightens them and generalizes to the case of problems with H\"older-continuous gradients and \textcite{gorbunov2022clipped} derives high-probability convergence rates in the case of VIPs. \textcite{sadiev2023high} relaxes the assumption of bounded variance to Assumption~\ref{ch3:as:bounded_alpha_moment_chapter_3} for all problem classes mentioned above, and the results under the same assumption are also derived for \algname{clipped-SGD} (without acceleration) by \textcite{nguyen2023high} in the convex and non-convex cases.

\paragraph{High-probability bounds for composite convex problems.} \textcite{nazin2019algorithms} propose a truncated version of Mirror Descent for convex and strongly convex composite problems and prove non-accelerated rates of convergence under bounded variance and \emph{bounded domain} assumptions. Accelerated results under bounded variance assumption for strongly convex composite problems are proven by \textcite{davis2021low}, who propose an approach based on robust distance estimation. Since this approach requires solving some auxiliary problem at each iteration of the method, the complexity bound from \textcite{davis2021low} contains extra logarithmic factors independent of the confidence level. Finally, in their very recent work, \textcite{nguyen2023improved} prove high-probability convergence for Clipped Stochastic Mirror Descent (\algname{Clipped-SMD}) for \emph{convex} composite problems. Moreover, the authors also propose Accelerated \algname{Clipped-SMD} (\algname{Clipped-ASMD}) and show that the algorithm is indeed accelerated \emph{but only under the additional assumption that $\nabla f(x^*) = 0$}.

\section{Main Results}\label{ch3:sec:the_issue_and_star_method_chapter_3}

In this section, we consider problem \eqref{ch3:eq:min_problem_chapter_3} and methods for it.

\paragraph{Failure of the na\"ive approach.}For simplicity, consider a non-stochastic case with strongly convex $f(x)$, $n = 1$. The standard deterministic first-order method for solving problems like \eqref{ch3:eq:min_problem_chapter_3} is Proximal Gradient Descent (\algname{Prox-GD}) \parencite{combettes2011proximal,nesterov2013gradient}: $x^{k+1} = \thesischthreeprox_{\gamma\Psi}(x^k - \gamma \nabla f(x^k))$. 
Due to the good interplay between the structure of the problem, properties of the proximal operator, and the structure of the method, \algname{Prox-GD} has the same (linear) convergence rate as \algname{GD} for minimization of $f(x)$. One of the key reasons for that is that any solution $x^*$ of problem \eqref{ch3:eq:min_problem_chapter_3} satisfies $x^* = \thesischthreeprox_{\gamma \Psi}(x^* - \gamma \nabla f(x^*))$, 
i.e., the solutions of \eqref{ch3:eq:min_problem_chapter_3} are fixed points of \algname{Prox-GD} (and vice versa), which is equivalent to $-\nabla f(x^*) \in \partial \Psi(x^*)$, where $\partial \Psi(x^*)$ is a subdifferential of $\Psi$ at $x^*$. However, if we apply gradient clipping to \algname{Prox-GD} na\"ively
\begin{equation}
    x^{k+1} = \thesischthreeprox_{\gamma\Psi}\left(x^k - \gamma \clip(\nabla f(x^k), \lambda)\right), \label{ch3:eq:prox_clip_gd_naive_chapter_3}
\end{equation}
then the method loses a fixed point property if $\|\nabla f(x^*)\| > \lambda$, because in this case, $-\clip(\nabla f(x^*), \lambda)$ does not necessarily belong to $\partial \Psi(x^*)$ and $x^* \neq \thesischthreeprox_{\gamma \Psi}(x^* - \gamma \clip(\nabla f(x^*), \lambda))$ in general. Therefore, for such $\lambda$, one has to decrease the stepsize $\gamma$ to achieve any accuracy of the solution. This approach slows down the convergence making it sublinear even without any stochasticity in the gradients. To avoid this issue, it is necessary to set $\lambda$ large enough. This strategy works in the deterministic case but becomes problematic for a stochastic version of \eqref{ch3:eq:prox_clip_gd_naive_chapter_3}:
\begin{equation}
    x^{k+1} = \thesischthreeprox_{\gamma\Psi}\left(x^k - \gamma \clip(\nabla f(x^k;\xi^k), \lambda_k)\right), \label{ch3:eq:prox_clip_sgd_naive_chapter_3}
\end{equation}
where $\xi^k$ is sampled independently from previous iterations. The problem comes from the fact the existing analysis in the unconstrained case (which is a special case of the composite case) requires taking decreasing $\lambda_k$ \parencite{gorbunov2021near, sadiev2023high} that contradicts the requirement that clipping level has to be large enough. Therefore, more fundamental algorithmic changes are needed.

\paragraph{Non-implementable solution.} Let us reformulate the issue: (i) to handle the heavy-tailed noise, we want to use decreasing clipping level $\lambda_k$, (ii) but the method should also converge linearly without the noise, i.e., when $\nabla f(x^k;\xi^k) = \Exp_{\xi^k}[\nabla f(x^k;\xi^k)] = \nabla f(x^k)$. In other words, \emph{the expectation of the vector that is clipped in the method should converge to zero with the same rate as $\lambda_k$}. The method should converge, i.e., with high probability, we should have $\nabla f(x^k) \to \nabla f(x^*)$. These observations lead us to the following purely theoretical algorithm that we call \algname{Prox-clipped-SGD-star}\footnote{The idea behind and the name of this method is inspired by \algname{SGD-star} proposed by \textcite{gorbunov2020unified, hanzely2019one}.}:
\begin{gather}
    x^{k+1} = \thesischthreeproxkPsi \left(x^k - \gamma \thesischthreetg^k\right),\;\; \text{where} \label{ch3:eq:prox_clipped_SGD_star_chapter_3}\\
    \thesischthreetg^k = \nabla f(x^*) + \clip\left(\nabla f(x^k;\xi^k) - \nabla f(x^*), \lambda_k\right). 
\end{gather}
The method is non-implementable since $\nabla f(x^*)$ is unknown in advance. Nevertheless, as we explain in the next subsection, the method is useful in designing and analyzing implementable versions. The following theorem gives the complexity of \algname{Prox-clipped-SGD-star}.

\begin{theorem}\label{ch3:thm:prox_clipped_SGD_main_chapter_3}
    Let $n=1$ and Assumptions~\ref{ch3:as:bounded_alpha_moment_chapter_3}, \ref{ch3:as:L_smoothness_chapter_3}, and \ref{ch3:as:QSC_chapter_3} with $\mu > 0$ hold for $Q = B_{2R}(x^*)$, $R \geq \|x^0 - x^*\|$, for some\footnote{In all of our results, one can use any solution $x^*$, e.g., one can take $x^*$ being a projection of $x^*$ on the solution set.} $x^* \in \arg\min_{x\in \R^d} \Phi(x)$. Assume that $K \geq 1$, $\beta \in (0,1)$, $A = \ln \tfrac{4(K+1)}{\beta}$, 
    \begin{gather*}
        0 < \gamma = \cO\left(\min\left\{\frac{1}{LA}, \frac{\ln(B_K)}{\mu(K+1)}\right\}\right),\\
        B_K = \Theta\left(\max\left\{2, \frac{(K+1)^{\nicefrac{2({p}-1)}{{p}}}\mu^2 R^2}{\sigma^2 A^{\nicefrac{2({p}-1)}{{p}}}\ln^2(B_K)}\right\}\right),\\
        \lambda_k = \Theta\left(\frac{\exp(-\gamma\mu(1+\nicefrac{k}{2}))R}{\gamma A}\right).
    \end{gather*}
    Then to guarantee $\|x^K - x^*\|^2 \leq \varepsilon$ with probability $\geq 1 - \beta$ \algname{Prox-clipped-SGD-star} requires
    \begin{equation}
        \widetilde\cO\left(\max\left\{\frac{L}{\mu}, \left(\frac{\sigma^2}{\mu^2\varepsilon}\right)^{\frac{{p}}{2({p}-1)}}\right\}\right) \label{ch3:eq:SGD_star_main_result_QSC_chapter_3}
    \end{equation}
    iterations/oracle calls.
\end{theorem}
\begin{proof}[Sketch of the proof]
    Following \textcite{gorbunov2020stochastic, sadiev2023high}, we prove by induction\footnote{We use the induction to apply Bernstein's inequality for the estimation of the sums appearing due to the stochasticity of the gradients. We refer to Section~\ref{ch3:sec:proofs_structure_chapter_3} for the details.} that $\|x^k - x^*\|^2 \leq 2\exp(-\gamma\mu k)R^2$ with high probability. This and $L$-smoothness imply that $\|\nabla f(x^k) - \nabla f(x^*)\| \sim \exp(-\nicefrac{\gamma\mu k}{2})$ and $\|\nabla f(x^k) - \nabla f(x^*)\| \leq \nicefrac{\lambda_k}{2}$ with high probability. These facts allow us to properly clip the heavy-tailed noise without sacrificing the convergence rate. See the complete formulation of Theorem~\ref{ch3:thm:prox_clipped_SGD_main_chapter_3} and the full proof in Appendix~\ref{ch3:appendix:prox_clipped_SGD_star_chapter_3}.
\end{proof}

The above complexity bound for \algname{Prox-clipped-SGD-star} coincides with the known one for \algname{clipped-SGD} for the unconstrained problems under the same assumptions \parencite{sadiev2023high} -- similarly as the complexity of \algname{Prox-GD} coincides with the complexity of \algname{GD} for unconstrained smooth problems.

\begin{remark}[On the logarithmic factors.]
    In our results, we do not focus on the logarithmic factors to keep the proofs simpler. Therefore, the main text omits logarithmic dependencies, while the appendix provides them. We believe it is possible to improve the logarithmic factors in our complexity bounds following the proof technique of \textcite{nguyen2023improved}.
\end{remark}

\begin{remark}[Dependence of the parameters on $R$.]
    Similarly to many previous works \parencite{gorbunov2020stochastic, gorbunov2022clipped, sadiev2023high, nguyen2023high, nguyen2023improved}, in all our results, the choice of the stepsize and clipping level depends on some upper bound\footnote{One can avoid this via considering problems defined on bounded domain \parencite{nazin2019algorithms} (then one can take $R$ equal to the diameter of the domain) or assuming that the gradient is bounded \parencite{cutkosky2021high} (in this case, one can at least remove the dependence on $R$ from the clipping level). However, both options have clear limitations.} $R$ for $\|x^0 - x^*\|$. Deriving matching high-probability results for parameter-free methods is an important open problem beyond the scope of this chapter. In practice, one can run the method with $R = \widetilde R$ for a sufficiently large number of steps (prescribed by theory). If the target accuracy is not achieved, then one needs to take $R = 2\widetilde R$ and repeat the procedure. Then, with a high probability after the logarithmic number of restarts, either parameter $R$ will be chosen appropriately, or the method will reach the desired accuracy.
\end{remark}

\paragraph{\algname{Prox-clipped-SGD-shift}.} As mentioned before, the key limitation of \algname{Prox-clipped-SGD-star} is that it explicitly uses shift $\nabla f(x^*)$, which is not known in advance. Therefore, guided by the literature on variance reduction and communication compression \parencite{gorbunov2020unified, gower2020variance, mishchenko2019distributed}, it is natural to approximate $\nabla f(x^*)$ via shifts $h^k$. This leads us to a new method called \algname{Prox-clipped-SGD-shift}: as before $x^{k+1} = \thesischthreeproxkPsi \left(x^k - \gamma \thesischthreetg^k\right)$ but now
\begin{gather}
    \thesischthreetg^k = h^k + \hat\Delta^k,\;\; h^{k+1} = h^k + \nu \hat \Delta^k, \label{ch3:eq:prox_clipped_SGD_shift_shifts_chapter_3}\\
    \hat\Delta^k = \clip\left(\nabla f(x^k;\xi^k) - h^k, \lambda_k\right), 
\end{gather}
where $\nu > 0$ is a stepsize for learning shifts. Similar shifts are proposed by \textcite{mishchenko2019distributed} in the context of distributed optimization with communication compression. Despite the similarities between compression operators and clipping, the latter does not have a constant contraction factor, making it non-trivial to extend the techniques from the literature on communication compression to the methods with clipping; see, e.g., \textcite[Section 5]{khirirat2023clip21}. Moreover, the proofs of \textcite{mishchenko2019distributed} rely on the unbiasedness of the compression operator. In contrast, gradient clipping introduces bias, and we consider high-probability analysis that differs significantly from the in-expectation one considered by \textcite{mishchenko2019distributed}.

Since \algname{Prox-clipped-SGD-shift} is a special case of its distributed variant, we continue our discussion with the distributed version of the method.

\paragraph{Distributed \algname{Prox-clipped-SGD-shift}.} We propose a generalization of \algname{Prox-clipped-SGD-shift} to the distributed case \eqref{ch3:eq:distributed_problem_chapter_3} called Distributed \algname{Prox-clipped-SGD-shift} (\algname{DProx-clipped-SGD-shift}):
\begin{gather}
    x^{k+1} = \thesischthreeproxkPsi \left(x^k - \gamma \thesischthreetg^k\right),\;\; \text{where}\label{ch3:eq:Dprox_clipped_SGD_shift_chapter_3}\\
    \thesischthreetg^k = \frac{1}{n}\sum\limits_{i=1}^n \thesischthreetg_i^k,\;\; \thesischthreetg_i^k = h_i^k + \hat\Delta_i^k,\;\; h_i^{k+1} = h_i^k + \nu \hat \Delta_i^k, \notag\\
    \hat\Delta_i^k = \clip\left(\nabla f_i(x^k;\xi_i^k) - h_i^k, \lambda_k\right), \label{ch3:eq:Dprox_clipped_SGD_shift_shifts_chapter_3}
\end{gather}
where $\xi_1^k,\ldots, \xi_n^k$ are sampled independently from each other and previous steps. In this method, worker $i$ updates the shift $h_i^k$ and sends clipped vector $\hat \Delta_i^k$ to the server. Since $\thesischthreetg^k = h^k + \tfrac{1}{n}\sum_{i=1}^n\hat \Delta_i^k$ and $h^{k+1} = h^k + \tfrac{\nu}{n}\sum_{i=1}^n\hat \Delta_i^k$, where $h^k = \frac{1}{n}\sum_{i=1}^n h_i^k$, workers do not need to send $h_i^k$ to the server for $k > 0$. We notice that even when $\Psi \equiv 0$, i.e., the problem is unconstrained, individual gradients $\{\nabla f_i(x^*)\}_{i\in [n]}$ of the clients' function at the solution of problem \eqref{ch3:eq:min_problem_chapter_3} are not necessarily zero, though their sum equals zero. However, if applied without any shifts to the local (stochastic) gradients, then, similarly to the case of non-distributed \algname{Prox-GD} \eqref{ch3:eq:prox_clip_gd_naive_chapter_3}, the clipping operation also breaks the fixed point property, since $\frac{1}{n}\sum_{i=1}^n \clip(\nabla f_i(x^*), \lambda) \neq 0$ for small values of $\lambda$. This highlights the importance of the shifts for distributed unconstrained case.

For the proposed method, we derive the following result. 

\begin{theorem}[Convergence of \algname{DProx-clipped-SGD-shift}: quasi-strongly convex case]\label{ch3:thm:D_prox_clipped_SGD_main_chapter_3}
    Let $K \geq 1$, $\beta \in (0,1)$, $A = \ln \frac{48n(K+1)}{\beta}$. Let Assumptions~\ref{ch3:as:bounded_alpha_moment_chapter_3}, \ref{ch3:as:L_smoothness_chapter_3}, and \ref{ch3:as:QSC_chapter_3} with $\mu > 0$ hold for $Q = B_{3n\sqrt{2}R}(x^*)$, where $R \geq \|x^0 - x^*\|^2.$ 
    Assume that $\zeta_* \geq \sqrt{\frac{1}{n}\sum_{i=1}^n \|\nabla f_i(x^*)\|^2}$, $\nu = \Theta(\nicefrac{1}{A})$,
    \begin{gather*}
        0 < \gamma = \cO\left(\min\left\{\frac{1}{LA}, \frac{\sqrt{n}R}{A\zeta_*}, \frac{\ln(B_K)}{\mu(K+1)}\right\}\right),\\
        B_K = \Theta\left(\max\left\{2, \frac{(K+1)^{\nicefrac{2({p}-1)}{{p}}}\mu^2 n^{\nicefrac{2({p}-1)}{{p}}} R^2}{\sigma^2 A^{\nicefrac{2({p}-1)}{{p}}}\ln^2(B_K)}\right\}\right),\\
        \lambda_k = \Theta\left(\frac{n\exp(-\gamma\mu(1+\nicefrac{k}{2}))R}{\gamma A}\right).
    \end{gather*}
    Then to guarantee $\|x^K - x^*\|^2 \leq \varepsilon$ with probability $\geq 1 - \beta$ \algname{DProx-clipped-SGD-shift} requires
    \begin{align}
        \widetilde\cO\left(\max\left\{\frac{L}{\mu}, \frac{\zeta_*}{\sqrt{n}\mu R}, \frac{1}{n}\left(\frac{\sigma^2}{\mu^2\varepsilon}\right)^{\frac{{p}}{2({p}-1)}}\right\}\right) \label{ch3:eq:d_prox_clipped_SGD_shift_main_result_QSC_chapter_3}
    \end{align}
    iterations/oracle calls per worker.
\end{theorem}
\begin{proof}[Sketch of the proof]
    The proof follows similar steps to the proof of Theorem~\ref{ch3:thm:prox_clipped_SGD_main_chapter_3} up the change of the Lyapunov function: by induction, we prove that $V_k \leq 2\exp(-\gamma\mu k)V$ with high probability, where $V_k = \|x^k - x^*\|^2 + \frac{C^2\gamma^2 A^2}{n}\sum_{i=1}^n \|h_i^k - \nabla f_i(x^*)\|^2$. The choice of the Lyapunov function reflects the importance of the ``quality'' of shifts $\{h_i^k\}_{i\in[n]}$, i.e., their proximity to $\{\nabla f_i(x^*)\}_{i\in[n]}$. Moreover, we increase the clipping level $n$ times to balance the bias and variance of $\thesischthreetg^k$; see Appendix~\ref{ch3:appendix:aux_results_chapter_3}. This allows us to reduce the last term in the complexity bound $n$ times. See the complete formulation of Theorem~\ref{ch3:thm:D_prox_clipped_SGD_main_chapter_3} and the full proof in Appendix~\ref{ch3:appendix:dprox_clipped_SGD_shift_chapter_3}.
\end{proof}

The next theorem gives the result in the convex case.

\begin{theorem}[Convergence of \algname{DProx-clipped-SGD-shift}: convex case]\label{ch3:thm:D_prox_clipped_SGD_main_cvx_chapter_3}
    Let $K \geq 1$, $\beta \in (0,1)$, $A = \ln \frac{48n(K+1)}{\beta}$. Let Assumptions~\ref{ch3:as:bounded_alpha_moment_chapter_3}, \ref{ch3:as:L_smoothness_chapter_3},  and \ref{ch3:as:str_cvx_chapter_3} with $\mu = 0$ hold for $Q = B_{\sqrt{2R}}(x^*)$, where $R \geq \|x^0 - x^*\|$. Assume that $\nu = 0$, $\zeta_* \geq \sqrt{\frac{1}{n}\sum_{i=1}^n \|\nabla f_i(x^*)\|^2}$,
    \begin{gather*}
        0 < \gamma = \cO\left(\min\left\{\frac{1}{LA}, \frac{\sqrt{n}R}{A\zeta_*}, \frac{n^{\nicefrac{({p}-1)}{{p}}}R}{\sigma K^{\nicefrac{1}{{p}}}A^{\nicefrac{({p}-1)}{{p}}}}\right\}\right),\\
        \lambda_k = \lambda = \Theta\left(\frac{nR}{\gamma A}\right).
    \end{gather*}
    Then to guarantee $\Phi(\bar{x}^K) - \Phi(x^*) \leq \varepsilon$ for $\bar{x}^K = \frac{1}{K+1}\sum_{k=0}^K x^k$ with probability $\geq 1 - \beta$ \algname{DProx-clipped-SGD-shift} requires
    \begin{align}
        \widetilde\cO\left(\max\left\{\frac{LR^2}{\varepsilon}, \frac{R\zeta_*}{\sqrt{n}\varepsilon}, \frac{1}{n}\left(\frac{\sigma R}{\varepsilon}\right)^{\frac{{p}}{{p}-1}}\right\}\right) \label{ch3:eq:d_prox_clipped_SGD_shift_main_result_cvx_chapter_3}
    \end{align}
    iterations/oracle calls per worker.
\end{theorem}

Discussion of the results for \algname{DProx-clipped-SGD-shift}. Up to the difference between $V$ and $\|x^0 - x^*\|^2$, in the single-node case, the derived results coincide with ones known for \algname{clipped-SGD} in the unconstrained case \parencite{sadiev2023high}. In the composite non-distributed case ($n=1$), the result of Theorem~\ref{ch3:thm:D_prox_clipped_SGD_main_chapter_3} is the first known of its type, and Theorem~\ref{ch3:thm:D_prox_clipped_SGD_main_cvx_chapter_3} recovers (up to logarithmic factors) the result of \textcite{nguyen2023improved} for a version of Stochastic Mirror Descent with gradient clipping (\algname{Clipped-SMD}), see Table~\ref{ch3:tab:comparison_of_rates_minimization_chapter_3}. Indeed, parameter $\widehat{R}^2 = R\left(3R + L^{-1}(2\eta \sigma + \|\nabla f(x^0)\|)\right)$ for some $\eta > 0$ from the result by \textcite{nguyen2023improved} equals $\Theta(R^2 + \nicefrac{R\zeta_*}{L})$, when $\eta$ is sufficiently small (otherwise $\widehat R$ can be worse than $\Theta(R^2 + \nicefrac{R\zeta_*}{L})$), which can be seen from the following inequalities following smoothness: $\|\nabla f(x^0)\| \leq \|\nabla f(x^*)\| + \|\nabla f(x^0) - \nabla f(x^*)\| \leq \|\nabla f(x^*)\| + L\|x^0 - x^*\|$ and $\|\nabla f(x^*)\| \leq \|\nabla f(x^0)\| + \|\nabla f(x^0) - \nabla f(x^*)\| \leq \|\nabla f(x^0)\| + L\|x^0 - x^*\|$. Since this chapter does not focus on logarithmic factors, we do not show them in the main text and provide the complete expressions in the appendix. \textcite{nguyen2023improved} have better dependencies on the parameters under logarithms than our results. In particular, the logarithmic factor in front of the last term in \eqref{ch3:eq:d_prox_clipped_SGD_shift_main_result_cvx_chapter_3} is $\ln\left(\beta^{-1}\left(\nicefrac{\sqrt{2}\sigma R}{\varepsilon}\right)^{\nicefrac{{p}}{({p}-1)}}\right)$ (see \eqref{ch3:eq:prox_clipped_SGD_convex_case_complexity_chapter_3} from Theorem~\ref{ch3:thm:prox_clipped_SGD_convex_case_chapter_3}), which matches the logarithmic factors from the result obtained by \textcite{sadiev2023high}. In the same setting with $n=1$, \textcite{nguyen2023improved} obtained a high-probability complexity result with $\ln(\beta^{-1})$ logarithmic factor, which is smaller than the one we have. However, this discrepancy can be relatively small in many situations. For example, if $\sigma \sim 10^3$ (high variance), $R \sim 10^3$ (starting point is far from the solution), ${p} = \nicefrac{3}{2}$, $\varepsilon \sim 10^{-6}$, $\beta \sim 10^{-6}$ (very accurate solution and small failure probability), then $\ln\left(\beta^{-1}\left(\nicefrac{\sqrt{2}\sigma R}{\varepsilon}\right)^{\nicefrac{{p}}{({p}-1)}}\right) \approx 98$ and $\ln(\beta^{-1}) \approx 14$, i.e., our logarithmic factor is just $\sim 7$ times larger in this case. We conjecture that adjusting the proof technique of \textcite{nguyen2023improved} one can improve the logarithmic factors in our results as well.

It is worth mentioning that shifts are not needed in the convex case because the method does not have fast enough convergence, which makes it work with a constant clipping level, i.e., the method in the convex case requires less tight gradient estimates and is more robust to the bias than in strongly convex. In the quasi-strongly convex case, the shifts' stepsize is chosen as $\nu \sim \Theta(\nicefrac{1}{A})$ and it does not explicitly affect the rate since $\gamma\mu = \Theta(\nicefrac{1}{A})$, see the details in Section~\ref{ch3:sec:proofs_structure_chapter_3} and Appendix~\ref{ch3:appendix:dprox_clipped_SGD_shift_chapter_3}. Moreover, as we show later, shifts are important for achieving acceleration in the convex case as well. A similar tradeoff between the needed accuracy of the estimate and the speed of convergence of the method is observed in the literature on the convergence of methods with inexact oracle \parencite{devolder2013exactness, devolder2014first}.

Next, as expected for a distributed method, the terms in the complexity bounds related to the noise improve with the growth of $n$. More precisely, the terms depending on the noise level $\sigma$ are proportional to $\nicefrac{1}{n}$, i.e., our results show so-called linear speed-up in the complexity -- a desirable feature for a stochastic distributed method. This aspect highlights the benefits of parallelization. To the best of our knowledge, the results for the distributed methods proposed in this chapter are the only existing ones under Assumption~\ref{ch3:as:bounded_alpha_moment_chapter_3} (even if we take into account the in-expectation convergence results). In the special case of ${p} = 2$, our results match (up to logarithmic factors) the state-of-the-art ones of \textcite{gorbunov2021near} since parallelization with linear speed-up follows for free under the bounded variance assumption, if the clipping is applied after averaging as it should be in the parallelized version of the methods of \textcite{gorbunov2021near} to keep the analysis of \textcite{gorbunov2021near} unchanged. Indeed, when $\{\nabla f_i(x;\xi_i)\}_{i\in [n]}$ are independent stochastic gradients satisfying Assumption~\ref{ch3:as:bounded_alpha_moment_chapter_3} with parameters $\sigma > 0$ and ${p} =2$, then $\frac{1}{n}\sum_{i\in[n]}\nabla f_i(x;\xi_i)$ also satisfies Assumption~\ref{ch3:as:bounded_alpha_moment_chapter_3} with parameters $\nicefrac{\sigma}{\sqrt{n}}$ and ${p} =2$. However, when ${p} < 2$ achieving linear speed-up is not that straightforward. If $\{\nabla f_i(x;\xi_i)\}_{i\in [n]}$ are independent stochastic gradients satisfying Assumption~\ref{ch3:as:bounded_alpha_moment_chapter_3} with parameters $\sigma > 0$ and ${p} < 2$, then \textcite[Lemma 7]{wang2021convergence} give a weaker guarantee: $\frac{1}{n}\sum_{i\in[n]}\nabla f_i(x;\xi_i)$ satisfies Assumption~\ref{ch3:as:bounded_alpha_moment_chapter_3} with parameters $\frac{2^{2-{p}}d^{\frac{1}{{p}}-\frac{1}{2}}\sigma}{n^{\frac{{p}-1}{{p}}}}$, which is \emph{dimension dependent}, and the same ${p}$. Therefore, if one applies this result to the known ones of \textcite{sadiev2023high, nguyen2023improved}, then the resulting complexity will have an extra factor of $d^{\frac{1}{{p}-1} - \frac{{p}}{2({p}-1)}}$ in the term that depends on $\sigma$. For large-scale or even medium-scale heavy-tailed problems, this factor can be huge, e.g., when $d = 1000$ and ${p} = \frac{7}{6}$, this factor is $1000^{6 - \frac{7}{3}} > 1000^3 = 10^9$.

To avoid these issues, we apply gradient clipping on the workers and then average clipped vectors, not vice versa. This is also partially motivated by the popularity of gradient clipping for ensuring differential privacy guarantees \parencite{abadi2016deep, chen2020understanding} in Federated Learning. Therefore, the proposed distributed methods can be useful for differential privacy as well, though this chapter does not study that aspect.

\begin{remark}[Dependence of the parameters on $\zeta_*$.]
    In our results, $\gamma$ and $\{\lambda_k\}_{k\geq 0}$ depend on $\zeta_*$, which is an upper bound for
    \[
        \sqrt{\frac{1}{n}\sum_{i=1}^n \|\nabla f_i(x^*)\|^2}.
    \]
    Since $f_i$ are $L$-smooth for all $i\in [n]$, we have
    \[
        \sqrt{\frac{1}{n}\sum_{i=1}^n \|\nabla f_i(x^*)\|^2}
        \leq G_0 + \sqrt{\frac{1}{n}\sum_{i=1}^n \|\nabla f_i(x^*) - \nabla f_i(x^0)\|^2}
        \leq G_0 + L\|x^0 - x^*\|,
    \]
    where
    \[
        G_0 = \sqrt{\frac{1}{n}\sum_{i=1}^n \|\nabla f_i(x^0)\|^2}.
    \]
    As suggested by \textcite{nazin2019algorithms, nguyen2023improved}, one can estimate $\|\nabla f_i(x^0)\|$ with probability at least $1 - \beta$ using the procedure of \textcite{minsker2015geometric}, which requires $\cO(\ln(n\beta^{-1}))$ samples on each client $i\in [n]$. Therefore, to estimate $\zeta_*$, it is sufficient to know $R$.
\end{remark}

\paragraph{Acceleration.} Next, we propose a distributed version of clipped Stochastic Similar Triangles Method \parencite{gorbunov2020stochastic, gasnikov2016universal} for composite problems (\algname{DProx-clipped-SSTM-shift}): $x^0 = y^0 = z^0$, $A_0 = a_0 = 0$, $a_{k+1} = \frac{k+2}{2aL}$, $A_{k+1} = A_k + a_{k+1}$ and
\begin{gather}
    x^{k+1} = \frac{A_k y^k + {a}_{k+1}z^k}{A_{k+1}},\\
    z^{k+1} = \thesischthreeprox_{a_{k+1}\Psi}\left(z^k - a_{k+1}\thesischthreetg(x^{k+1})\right), \label{ch3:eq:D_prox_clipped_SSTM_shift_x_z_chapter_3}\\
    y^{k+1} = \frac{A_k y^k + {a}_{k+1}z^{k+1}}{A_{k+1}} \label{ch3:eq:D_prox_clipped_SSTM_shift_y_chapter_3}
\end{gather}
with
\begin{gather}
    \thesischthreetg(x^{k+1}) = \frac{1}{n}\sum\limits_{i=1}^n \thesischthreetg_i(x^{k+1}),\;\; \thesischthreetg_i(x^{k+1}) = h_i^k + \hat\Delta_i^k, \label{ch3:eq:D_prox_clipped_SSTM_shift_estimator_chapter_3}\\
    h_i^{k+1} = h_i^k + \nu_k \hat \Delta_i^k,\\
    \hat\Delta_i^k = \clip\left(\nabla f_i(x^{k+1};\xi_i^k) - h_i^k, \lambda_k\right), \label{ch3:eq:D_prox_clipped_SSTM_shift_estimator_2_chapter_3}
\end{gather}
where $\xi_1^k,\ldots, \xi_n^k$ are sampled independently from each other and previous steps. For the proposed method, we derive the following result.

\begin{theorem}[Convergence of \algname{DProx-clipped-SSTM-shift}]\label{ch3:thm:D_prox_clipped_SSSTM_main_chapter_3}
    Let Assumptions~\ref{ch3:as:bounded_alpha_moment_chapter_3}, \ref{ch3:as:L_smoothness_chapter_3},  and \ref{ch3:as:str_cvx_chapter_3} with $\mu = 0$ hold for $Q = B_{5\sqrt{2}nR}(x^*)$, where $R \geq \|x^0 - x^*\|^2$. Let $\zeta_* \geq \sqrt{\frac{1}{n}\sum_{i=1}^n \|\nabla f_i(x^*)\|^2}$, $C = \Theta(\nicefrac{A}{\sqrt{n}})$, $K_0 = \Theta(A^2)$, where $K \geq 1$, $\beta \in (0,1)$, $A = \ln \tfrac{10nK}{\beta}$. Let
    \begin{gather*}
        \nu_k = \begin{cases}
            \frac{2k+5}{(k+3)^2},& \text{if } k \geq K_0,\\ \frac{(k+2)^2}{C^2(K_0+2)^2n},& \text{if } k < K_0,
        \end{cases},\;\; \lambda_k = \Theta\left(\frac{n R}{{a}_{k+1}A}\right),\\
        a = \Theta\left(\max\left\{2, \frac{A^4}{n}, \frac{A^3\zeta_*}{L\sqrt{n}R}, \frac{\sigma K^{\nicefrac{({p}+1)}{{p}}}A^{\nicefrac{({p}-1)}{{p}}}}{LRn^{\nicefrac{{p}-1}{{p}}}}\right\}\right).
    \end{gather*}
    Then to guarantee $\Phi(y^K) - \Phi(x^*) \leq \varepsilon$ with probability $\geq 1 - \beta$ \algname{DProx-clipped-SSTM-shift} requires
    \begin{align}
        \widetilde\cO\left(\max\left\{\sqrt{\frac{LR^2}{\varepsilon}}, \sqrt{\frac{R\zeta_*}{\sqrt{n}\varepsilon}}, \frac{1}{n}\left(\frac{\sigma R}{\varepsilon}\right)^{\frac{{p}}{{p}-1}}\right\}\right) \label{ch3:eq:d_prox_clipped_SSTM_shift_main_result_cvx_chapter_3}
    \end{align}
    iterations/oracle calls per worker.
\end{theorem}
\begin{proof}[Sketch of the proof]
    The proof resembles the one for \algname{clipped-SSTM} of \textcite{sadiev2023high} but has some noticeable differences. In addition to handling the extra technical challenges appearing due to the composite structure (e.g., one cannot apply some useful formulas like $z^k - z^{k+1} = a_{k+1}\thesischthreetg(x^{k+1})$ that hold in the unconstrained case), we use a non-standard potential function $M_k$ defined as $M_k = \|z^k - x^*\|^2 + \frac{C^2{a}_{K_0+1}^2}{n}\sum_{i=1}^n \|h_i^k - \nabla f_i(x^*)\|^2$ for $k \leq K_0$ and $M_k = \|z^k - x^*\|^2 + \tfrac{C^2a_{k+1}^2}{n}\sum_{i=1}^n \|h_i^k - \nabla f_i(x^*)\|^2$ for $k > K_0$. See more details in Appendix~\ref{ch3:appendix:dprox_clipped_SSTM_shift_chapter_3}. %
\end{proof}

When $n=1$, the derived result has optimal dependence on $\varepsilon$ (up to logarithmic factors) \parencite{NemiYudin, zhang2020adaptive}. In contrast to the result of \textcite{nguyen2023improved}, we do not assume that $\nabla f(x^*) = 0$. Moreover, as \algname{DProx-clipped-SGD-shift}, \algname{DProx-clipped-SSTM-shift} benefits from parallelization since the second term in \eqref{ch3:eq:d_prox_clipped_SSTM_shift_main_result_cvx_chapter_3} is proportional $\nicefrac{1}{n}$. When $n$ is sufficiently large, the effect of acceleration can become significant even for large $\sigma$. In Appendix~\ref{ch3:appendix:SSTM_str_cvx_chapter_3}, we also provide the convergence results for the restarted version of \algname{DProx-clipped-SSTM-shift} assuming additionally that $f$ is strongly convex and one can compute starting shifts $h_i^0$ as $\nabla f_i(x^0)$.

\section{Structure of the Proofs\label{ch3:sec:proofs_structure_chapter_3}}

\begin{table*}[t]
    \centering
    \scriptsize
    \caption[High-probability complexity results for distributed variational inequalities]{\scriptsize Summary of known and new high-probability complexity results for solving (non-) composite (non-) distributed variational inequality problem \eqref{ch3:eq:main_problem_chapter_3}. Column ``Setup'' indicates the assumptions made in addition to Assumptions~\ref{ch3:as:bounded_alpha_moment_chapter_3}. In the results of \textcite{sadiev2023high}, all assumptions are made only on some ball around the solution with radius $\sim R \geq \|x^0 - x^*\|$; in Theorems~\ref{ch3:thm:D_prox_clipped_SGDA_main_chapter_3} and \ref{ch3:thm:D_prox_clipped_SEG_main_chapter_3}, the radius is $\sim\sqrt{V}$. Complexity is the number of stochastic oracle calls(per worker) needed for a method to guarantee that $\Prob\{\text{Metric} \leq \varepsilon\} \geq 1 - \beta$ for some $\varepsilon> 0$, $\beta \in (0,1]$ and ``Metric'' is taken from the corresponding column. Numerical and logarithmic factors are omitted for simplicity. Column ``C?'' shows whether the problem \eqref{ch3:eq:min_problem_chapter_3} is composite, ``D?'' indicates whether the problem \eqref{ch3:eq:min_problem_chapter_3} is distributed. Notation: $\thesischthreetx^K_{\text{avg}} = \frac{1}{K+1}\sum_{k=0}^K \thesischthreetx^k$ (for \algname{SEG}-type methods), $x^K_{\avg} = \frac{1}{K+1}\sum_{k=0}^K x^k$ (for \algname{SGDA}-type methods); $L$ = Lipschitz constant; $\sigma$ = parameter from Assumption~\ref{ch3:as:bounded_alpha_moment_chapter_3}; $R$ = any upper bound on $\|x^0 - x^*\|$ in the results of \textcite{sadiev2023high}; $V$ = any upper bound on $\|x^0 - x^*\|^2 + \frac{409600\gamma^2\ln^2\frac{48n(K+1)}{\beta}}{n^2}\sum_{i=1}^n\|F_i(x^*)\|^2$ for the results established in this chapter; $\mu$ = quasi-strong monotonicity parameter; $\ell$ = star-cocoercivity parameter. The results established in this chapter are highlighted in blue.}
    \label{ch3:tab:comparison_of_rates_VIP_chapter_3}
    \begin{threeparttable}
        \resizebox{0.95\textwidth}{!}{%
        \begin{tabular}{|c|c c c c c|}
        \hline
Setup & Method  & Metric & Complexity & C? & D?\\
        \hline\hline
        \multirow{3}{*}{As.~\ref{ch3:as:L_Lip_chapter_3} \& \ref{ch3:as:monotonicity_chapter_3}} & \begin{tabular}{c}
            \algname{clipped-SEG} \\ \tabcite{sadiev2023high}
        \end{tabular} & $\gap_{R}(\thesischthreetx_{\avg}^K)$ & $\max\left\{\frac{LR^2}{\varepsilon}, \left(\frac{\sigma R}{\varepsilon}\right)^{\frac{{p}}{{p}-1}}\right\}$ & \thesischthreexmark & \thesischthreexmark \\
        & \cellcolor{linen}\begin{tabular}{c}
            \algname{DProx-clipped-SEG-shift}\\ Theorem~\ref{ch3:thm:D_prox_clipped_SEG_monotone_appendix_chapter_3}
        \end{tabular} & \cellcolor{linen}$\gap_{\sqrt{V}}(\thesischthreetx_{\avg}^K)$ &\cellcolor{linen} $\max\left\{\frac{LV}{\varepsilon}, \frac{1}{n}\left(\frac{\sigma \sqrt{V}}{\varepsilon}\right)^{\frac{{p}}{{p}-1}}\right\}$ &\cellcolor{linen} \thesischthreecmark &\cellcolor{linen} \thesischthreecmark \\
        \hline
        \multirow{3}{*}{As.~\ref{ch3:as:L_Lip_chapter_3} \& \ref{ch3:as:str_monotonicity_chapter_3}} & \begin{tabular}{c}
            \algname{clipped-SEG} \\ \tabcite{sadiev2023high}
        \end{tabular} & $\|x^k - x^*\|^2$ & $\max\left\{\frac{L}{\mu}, \left(\frac{\sigma^2}{\mu^2\varepsilon}\right)^{\frac{{p}}{2({p}-1})}\right\}$ & \thesischthreexmark & \thesischthreexmark\\
        & \cellcolor{linen}\begin{tabular}{c}
            \algname{DProx-clipped-SEG-shift}\\ Theorem~\ref{ch3:thm:D_prox_clipped_SEG_QSM_appendix_chapter_3}
        \end{tabular} & \cellcolor{linen}$\|x^k - x^*\|^2$ &\cellcolor{linen} $\max\left\{\frac{L}{\mu}, \frac{1}{n}\left(\frac{\sigma^2}{\mu^2\varepsilon}\right)^{\frac{{p}}{2({p}-1})}\right\}$ &\cellcolor{linen} \thesischthreecmark &\cellcolor{linen} \thesischthreecmark \\
        \hline\hline
        \multirow{3}{*}{As.~\ref{ch3:as:monotonicity_chapter_3} \& \ref{ch3:as:star_cocoercivity_chapter_3} \& \ref{ch3:as:cocoercivity_chapter_3}} & \begin{tabular}{c}
            \algname{clipped-SGDA} \\ \tabcite{sadiev2023high}
        \end{tabular} & $\gap_{R}(x_{\avg}^K)$ & $\max\left\{\frac{\ell R^2}{\varepsilon}, \left(\frac{\sigma R}{\varepsilon}\right)^{\frac{{p}}{{p}-1}}\right\}$ & \thesischthreexmark & \thesischthreexmark\\
        & \cellcolor{linen}\begin{tabular}{c}
            \algname{DProx-clipped-SGDA-shift}\\ Theorem~\ref{ch3:thm:main_result_gap_prox_clipped_SGDA_chapter_3}
        \end{tabular} & \cellcolor{linen}$\gap_{\sqrt{V}}(x_{\avg}^K)$ &\cellcolor{linen} $\max\left\{\frac{\ell V}{\varepsilon}, \frac{1}{n}\left(\frac{\sigma \sqrt{V}}{\varepsilon}\right)^{\frac{{p}}{{p}-1}}\right\}$ &\cellcolor{linen} \thesischthreecmark &\cellcolor{linen} \thesischthreecmark \\
        \hline
        \multirow{3}{*}{As.~\ref{ch3:as:str_monotonicity_chapter_3} \& \ref{ch3:as:star_cocoercivity_chapter_3}} & \begin{tabular}{c}
            \algname{clipped-SGDA} \\ \tabcite{sadiev2023high}
        \end{tabular} & $\|x^K - x^*\|^2$ & $\max\left\{\frac{\ell}{\mu}, \left(\frac{\sigma^2}{\mu^2\varepsilon}\right)^{\frac{{p}}{2({p}-1})}\right\}$ & \thesischthreexmark & \thesischthreexmark\\
        & \cellcolor{linen}\begin{tabular}{c}
            \algname{DProx-clipped-SGDA-shift}\\ Theorem~\ref{ch3:thm:main_result_prox_clipped_SGDA_2_chapter_3}
        \end{tabular} & \cellcolor{linen}$\|x^K - x^*\|^2$ &\cellcolor{linen} $\max\left\{\frac{\ell}{\mu}, \frac{1}{n}\left(\frac{\sigma^2}{\mu^2\varepsilon}\right)^{\frac{{p}}{2({p}-1})}\right\}$ &\cellcolor{linen} \thesischthreecmark &\cellcolor{linen} \thesischthreecmark \\
        \hline
    \end{tabular}}
    \end{threeparttable}
\end{table*}

In this section, we elaborate on the structure of the proofs of our results and highlight additional challenges appearing due to the presence of the composite term and distributed nature of the methods. The proof of each result consists of two parts: optimization/descent lemma and the analysis of the sums appearing due to the stochasticity and biasedness of the updates (due to the clipping). In the first part, we usually follow some standard analysis of corresponding deterministic method without clipping and separate the stochastic part from the deterministic one (though for \algname{DProx-clipped-SSTM-shift} we use quite non-standard Lyapunov function, which can be interesting on its own). For example, in the analysis\footnote{In the appendix, we analyze this case in the generality of variational inequalities. Here we provide a simplified version for minimization.} of \algname{DProx-clipped-SGD-shift} under Assumption~\ref{ch3:as:QSC_chapter_3}, we prove the following inequality:
\begin{align}
        &V_{K+1} \leq (1 - \gamma \mu)^{K+1}V_0 \notag\\
        &\;\; + \frac{2\gamma}{n} \sum\limits_{k=0}^K\sum\limits_{i=1}^n (1-\gamma\mu)^{K-k}\langle s_k, \omega_{i,k} \rangle\notag\\
        &\;\; + \gamma^2\sum\limits_{k=0}^K(1-\gamma \mu)^{K-k}\left(\|\omega_k\|^2 + \frac{1}{n^2}\sum\limits_{i=1}^n  \|\omega_{i,k}\|^2\right), \notag
\end{align}
where $V_k = \|x^k - x^*\|^2 + \frac{C^2\gamma^2 A^2}{n}\sum_{i=1}^n \|h_i^k - \nabla f_i(x^*)\|^2$ for some numerical constant $C> 0$, $s_k = x^k - x^* - \gamma (\nabla f(x^k)- \nabla f(x^*))$ and vectors $\omega_{i,k} = \nabla f_i(x^k) - \thesischthreetg_i^k$ represent the discrepancy between the full gradients and their estimates. Moreover, to use this inequality for some $K \geq 0$ we need to show that $\{x^k\}_{k=0}^K$ belong to the set where the assumptions hold (in this particular case, to $B_{3n\sqrt{2}R}(x^*)$) with high probability. We do it always by induction. More precisely, we prove that $\Prob\{E_k\} \geq 1 - \nicefrac{k\beta}{(K+1)}$ for the probability event $E_k$ defined as follows: inequalities $V_t \leq 4\exp(-\gamma\mu t)R^2$ and $\left\|\frac{\gamma}{n}\sum_{i=1}^{r-1}\omega_{i,t-1}^u\right\| \leq \exp(-\nicefrac{\gamma\mu(t-1)}{2})\sqrt{\nicefrac{R^2}{2}}$ hold for $t = 0,1,\ldots,k$ and $r = 1, 2, \ldots, n$ simultaneously, where $\omega_{i,t}^u = \Exp_{\xi_i^t}[\thesischthreetg_i^t] - \thesischthreetg_i^t$ and $\Exp_{\xi_i^t}[\cdot]$ denotes an expectation w.r.t.\ $\xi_i^t$. To prove this, we use Bernstein inequality for martingale difference (see Lemma~\ref{ch3:lem:Bernstein_ineq_chapter_3}). However, to apply Bernstein inequality we need to circumvent multiple technical difficulties related to the estimation of the norm of the clipped vector (that involves derivations related to the shifts $\{h_i^k\}_{i\in [n]}$), proper choice of the clipping level to control the bias and variance and achieve desired linear speed-up (see Lemma~\ref{ch3:lem:bias_and_variance_clip_distributed_chapter_3} and the following discussion). Moreover, when $n > 1$ (distributed case), we also need to apply additional induction over clients to estimate sums like $\circledSix$ from \eqref{ch3:eq:SGDA_2_prox_clipped_12345_bound_chapter_3}.

\section{Extensions for Variational Inequalities}

In addition to the minimization problems, we also consider stochastic composite variational inequality problems (VIPs):
\begin{equation}
    \text{find } x^* \in \R^d \text{ such that } \langle F(x^*), x-x^* \rangle + \Psi(x) - \Psi(x^*) \geq 0. \label{ch3:eq:main_problem_chapter_3}
\end{equation}
where the assumptions on operator $F(x) = \Exp_{\xi\sim \thesischthreecD}[F(x;\xi)]: \R^d \to \R^d$ will be specified later and, as in the case of minimization, $\Psi(x)$ is a proper, closed, convex function. When $f(x)$ is convex problem \eqref{ch3:eq:min_problem_chapter_3} is a special case of \eqref{ch3:eq:main_problem_chapter_3} with $F(x) = \nabla f(x)$. For the examples of problems of type \eqref{ch3:eq:main_problem_chapter_3}, we refer to \textcite{alacaoglu2022stochastic, beznosikov2022stochastic}.

The distributed version of \eqref{ch3:eq:main_problem_chapter_3} has the following structure:
\begin{equation}
    F(x) = \frac{1}{n}\sum\limits_{i=1}^n \left\{F_i(x) = \Exp_{\xi_i\sim \thesischthreecD_i}[F_i(x;\xi_i)]\right\}. \label{ch3:eq:distributed_problem_VI_chapter_3}
\end{equation}
In this case, there are $n$ workers connected in a centralized way with some parameter server; worker $i$ can query some noisy information (stochastic gradients/estimates) about $F_i$.

Similarly to the minimization case, we consider the situation when $F_i$ is accessible through the stochastic oracle calls. The stochastic estimates satisfy the following assumption.
\begin{assumption}\label{ch3:as:bounded_alpha_moment_VI_chapter_3}
    There exist some set $Q \subseteq \R^d$ and values $\sigma \geq 0$, ${p} \in (1,2]$ such that for all $x \in Q$
    we have $\Exp_{\xi_i\sim \thesischthreecD_i}[F_i(x;\xi_i)] = F_i(x)$ and
    \begin{equation}
        \Exp_{\xi_i\sim \thesischthreecD_i}[\left\|F_i(x;\xi_i) - F_i(x)\right\|^{p}] \leq \sigma^{p}. \label{ch3:eq:bounded_alpha_moment_operator_chapter_3}
    \end{equation}
\end{assumption}

Under this and additional assumptions on $\{F_i\}_{i=1}^n$, we obtain several convergence results for the distributed versions\footnote{\algname{GDA}-type methods are analogs of \algname{GD} iteration $x^{k+1} = x^k - \gamma F(x^k)$. Such methods converge under (star-)cocoercivity but fail for general monotone problems. To address this issue, versions of Extragradient-type iteration \parencite{korpelevich1976extragradient} $x^{k+1} = x^k - \gamma F(x^k - \gamma F(x^k))$ are used; see further details in Appendix~\ref{ch3:appendix:VIs_chapter_3}.} of \algname{clipped-SGDA} and \algname{clipped-SEG} developed by \textcite{gorbunov2022clipped}. The results are summarized in Table~\ref{ch3:tab:comparison_of_rates_VIP_chapter_3}. As one can see from the table, our results generalize the results of \textcite{sadiev2023high} to the composite distributed case and achieve linear speedup. Further details on the assumptions and our main results for distributed composite variational inequalities are deferred to Appendix~\ref{ch3:appendix:VIs_chapter_3}.

\section{Conclusion and Future Work}

This chapter extended high-probability optimization under bounded central \(p\)-th moments to composite objectives and to a centralized distributed model with \(n\) workers. Directly clipping a local stochastic gradient can displace the fixed point because an individual component need not vanish at a solution. The proposed proximal methods instead clip stochastic innovations relative to shifts that are learned online. This construction preserves the target balance while allowing decreasing clipping levels, acceleration, and analogous updates for composite variational inequalities.

The theorem outputs depend on the structural regime. Convex minimization is assessed at an averaged iterate, and the accelerated method at \(y^K\); quasi-strongly convex minimization uses terminal squared distance. For variational inequalities, the monotone results control a restricted gap at averaged extragradient or gradient-descent--ascent iterates, whereas quasi-strong monotonicity yields terminal squared distance. Complexity is measured per worker, and the heavy-tailed stochastic terms exhibit the stated linear speedup in \(n\). Tightness or optimality holds only in the documented special cases and up to logarithmic factors; the chapter does not establish a universal lower bound for every composite and distributed regime.

The schedules require an upper radius and, in relevant results, the heterogeneity quantity \(\zeta_*\); the logarithmic factors are not optimized. The communication model also assumes a fixed, fully participating worker population connected to a central server. Although worker-side clipping is compatible with mechanisms used in private learning, no differential-privacy guarantee is proved here. Chapter~\ref{ch:chapter_4} turns to a different source of information: stochastic Hessian--vector products, and asks through matching lower and upper bounds when second-order information improves heavy-tailed nonconvex optimization.

\newpage
  \endgroup

\part{Advanced Non-Convex and Second-Order Optimization in Heavy-Tailed Settings}

\chapter{Second-order Optimization under Heavy-Tailed Noise: Hessian Clipping and Sample Complexity Limits}
\label{ch:chapter_4}
\thesischaptersourcenote
  {Second-Order Optimization under Heavy-Tailed Noise: Hessian Clipping and Sample Complexity Limits}
  {sadiev2025second}
  \begingroup
  \thesischaptersetup{4}%
  
\begin{paperabstract}
  Heavy-tailed noise is pervasive in modern machine learning applications, arising from data heterogeneity, outliers, and non-stationary stochastic environments. While second-order methods can significantly accelerate convergence in light-tailed or bounded-noise settings, such algorithms are often brittle and lack guarantees under heavy-tailed noise---precisely the regimes where robustness is most critical. This chapter takes a first step toward a theoretical understanding of second-order optimization under heavy-tailed noise. We consider a setting where stochastic gradients and Hessians have only bounded $p$-th moments, for some $p\in (1,2]$, and establish tight lower bounds on the sample complexity of {any zero-respecting stochastic $q$-th-order algorithm}. We then develop a variant of normalized stochastic gradient descent that leverages second-order information and provably matches these lower bounds. To address the instability caused by large deviations, we introduce a novel algorithm based on gradient and Hessian clipping, and prove high-probability upper bounds that nearly match the fundamental limits. The results provide a comprehensive sample-complexity characterization for second-order optimization under heavy-tailed noise and position Hessian clipping as a robust strategy for second-order algorithm design in heavy-tailed regimes.
\end{paperabstract}

\section{Introduction}

We consider the stochastic optimization problem
\begin{eqnarray}\label{ch4:eq:main_problem_chapter_4}
\min_{x \in \mathbb{R}^d} f(x), \qquad f(x) := \Exp_{\xi \sim {\cD}}[f(x;\xi)],
\end{eqnarray}
with $f$ a smooth (but potentially nonconvex) objective, and $\xi$ a random variable drawn from an unknown distribution. While first-order methods are widely used in practice due to their simplicity and scalability, the application of second-order methods remains limited in stochastic settings. This is in stark contrast to their strong theoretical properties: second-order methods such as Newton’s method \parencite{bennett1916newton,kantorovich1948functional}, cubic regularization \parencite{nesterov2006cubic,nesterov2008accelerating}, their variants for finite-sum minimization \parencite{kovalev2019stochastic}, quasi-Newton \parencite{dennis1977quasi}, and adaptive second-order algorithms \parencite{doikov2024super} can offer provably faster convergence rates. Furthermore, complexity-theoretic results demonstrate that second-order methods can outperform first-order ones in noiseless (or finite sum) settings, often with  moderate additional computational cost \parencite{agarwal2018lower,arjevani2019oracle,arjevani2020second} only.

\paragraph{Second-order methods are highly sensitive to noise.}
First-order stochastic optimization (FOSO) methods--such as stochastic gradient descent (SGD) and its adaptive variants--are well understood, both algorithmically and from a complexity--theoretic perspective. These methods enjoy optimal convergence guarantees under a broad class of noise models, including bounded variance and infinite variance regimes \parencite{gower2019sgd,khaled2020better,yang2023two_sides,Wang_SGD_inf_var_21,fatkhullin2025can}. In contrast, the development of second-order stochastic optimization (SOSO) methods has been significantly hampered by their extreme sensitivity to noise, particularly in Hessian estimation. Existing second-order methods suffer from both practical instability and overly restrictive theoretical assumptions. For instance, stochastic extensions of cubic regularization \parencite{ghadimi2017second,tripuraneni2018stochastic} and its momentum variants \parencite{chayti2024improving} require bounded noise assumptions, which are rarely satisfied in modern large-scale learning tasks. Other frameworks--such as trust-region methods \parencite{arjevani2020second}, recursive momentum schemes \parencite{BetterSGDUsingSOM_Tran_2021}, and extrapolation-based approaches \parencite{antonakopoulos2022extra,agafonov2023advancing}--relax this to bounded variance, but still fall short of handling more realistic heavy-tailed noise distributions. To the best of our knowledge, no existing SOSO method can operate reliably under noise conditions with unbounded variance, highlighting a fundamental gap in the current landscape of stochastic optimization.

\paragraph{Heavy-tailed noise in gradients and Hessians.} In recent years, the first-order optimization literature has increasingly moved beyond the bounded variance assumption, adopting more general noise models that better reflect empirical observations in modern applications. A particularly influential framework is the \emph{bounded central moment} condition ($p$-BCM), which assumes
\[
\Exp\left[\|\nabla f(x;\xi) - \nabla f(x)\|^p\right] \le \sigma^p, \qquad \text{for some } p \in (1, 2],
\]
thereby allowing for heavy-tailed noise with unbounded variance. This model aligns well with empirical findings in deep learning and reinforcement learning (RL), where heavy-tailed gradient noise is often observed \parencite{garg2021proximal,zhang2020adaptive,ahn2023linear,simsekli2019tail,battash2024revisiting}. Algorithms designed for robustness under $p$-BCM noise, such as those employing gradient clipping or normalization, have demonstrated both strong empirical performance and rigorous convergence guarantees \parencite{zhang2020adaptive,Hubler2024clip_to_norm}. These developments suggest that heavy-tailed noise models are not only practically relevant, but also theoretically tractable, providing a principled foundation for algorithm design. Given that second-order methods are even more sensitive to noise than first-order ones, this naturally motivates the extension of such robustness principles to the second-order setting. This chapter takes that step and aims to
\begin{center}
\textit{develop a comprehensive theory of second-order stochastic optimization under heavy-tailed noise.}
\end{center}
Specifically, we assume access to unbiased stochastic gradients and Hessian-vector products, each satisfying a $p$-BCM-type condition:
\[
\Exp\left[{\normop{\nabla^2 f(x;\xi) - \nabla^2 f(x)}^p}\right] \le \sigma_h^p, \qquad \text{for some } p \in (1, 2].
\]
Our goal is to characterize the fundamental performance limits and develop practical algorithms for minimizing $f(x)$ in this setting, with a focus on obtaining guarantees for finding points with small gradient norm $\|\nabla f(x)\| \leq \varepsilon$, either in expectation or with high probability.\footnote{It would be interesting to extend our results to finding second-order stationary points \parencite{nesterov2006cubic,arjevani2020second}, but we leave this investigation for future work.}

\paragraph{From gradient to Hessian clipping.} Gradient clipping has become a standard tool in modern machine learning, particularly due to its empirical success in stabilizing training under heavy-tailed noise and ill-conditioned objectives \parencite{pascanu2013difficulty,schulman-et-al17ppo}. Theoretically, it has been shown to offer robustness under relaxed smoothness and moment assumptions \parencite{polyak1979adaptive,jakovetic2023nonlinear}, and enable high-probability guarantees with logarithmic (sub-Gaussian) dependence on the failure probability. The latter property of gradient clipping is in stark contrast to the classical linear algorithms, e.g., SGD, Momentum SGD, for which high probability lower bounds were recently established \parencite{fatkhullin2025can}. These properties are well-documented across both convex \parencite{nazin2019algorithms,gorbunov2020stochastic,davis2021low,gorbunov2023high,liu2023stochastic,gorbunov2024high,puchkin2024breaking,NonlinearHP2023Armacki} and nonconvex \parencite{sadiev2023high,nguyen2023improved,cutkosky2021high,Hubler2024clip_to_norm} regimes. Crucially, high-probability results are both theoretically and practically appealing as they guarantee the performance of individual runs rather than for the average-case behavior.

Despite these advances, the robustness enabled by gradient clipping remains largely confined to first-order methods. In the second-order setting, where both gradients and Hessians may be corrupted by heavy-tailed noise, comparable algorithmic tools are conspicuously lacking. The core difficulty is not just technical but conceptual: existing SOSO algorithms do not include mechanisms to suppress the influence of extreme noise in Hessian estimates. This leads us to a fundamental question for the design of robust second-order stochastic methods:
\begin{center}
\textit{Can we develop second-order algorithms with high-probability convergence guarantees under simultaneous heavy-tailed noise in both gradients and Hessians?}
\end{center}

\begin{figure}[ht]
    \begin{minipage}[t]{0.48\linewidth}
        \centering
        \includegraphics[width=\linewidth]{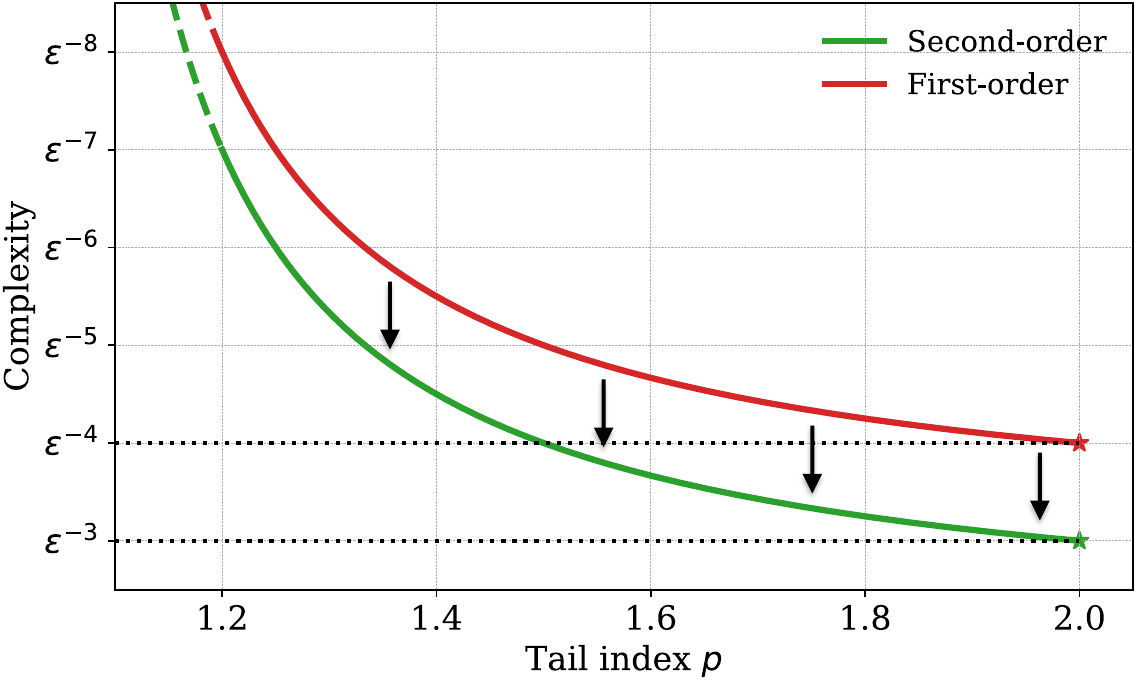}
        \caption[Sample complexity comparison for FOSO and SOSO depending on the tail index $p$.]{Sample complexity comparison for FOSO and SOSO depending on the tail index $p$. Each line corresponds to the leading term in the sample complexity for each class of algorithms. These leading terms match in upper and lower bounds, so this characterization is exact. We establish the characterization along the entire green line complementing  prior work \parencite{arjevani2023lower} for $p=2$. 
        }
        \label{ch4:fig:complexities_chapter_4}
    \end{minipage}%
    \hfill
    \begin{minipage}[t]{0.47\linewidth}
        \centering
        \includegraphics[width=\linewidth]{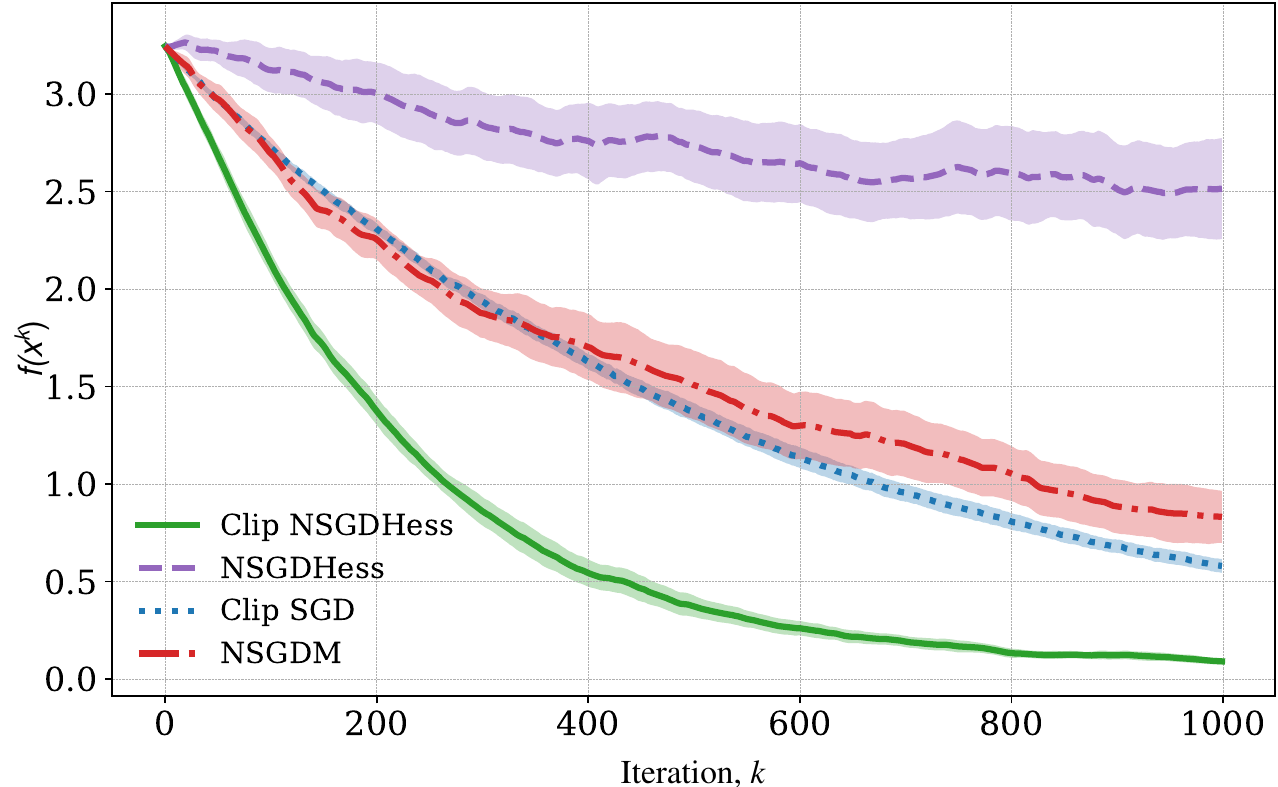}
        \caption{Performance of algorithms on a simple problem, $f(x) = 0.5 \norm{x}^2$, $d = 10$ with synthetic noise generated from a two-sided Pareto distribution with tail index $p=1.1$. We observe that algorithms without clipping, \algname{NSGDM} and \algname{NSGDHess}, suffer significantly from noise. This motivates our more in-depth study involving gradient and Hessian clipping for high probability convergence.}
        \label{ch4:fig:compare-algs_chapter_4}
    \end{minipage}
\end{figure}

\paragraph{Contributions:}
\begin{itemize}
    \item \textbf{Tight lower bound.} We establish a minimax lower bound on the sample complexity of {any zero-respecting stochastic $q$-th-order algorithm} under $p$-BCM noise model:
\(
\Omega\left( \frac{\Delta \sigma_h}{\varepsilon^2} \left( \frac{\sigma}{\varepsilon} \right)^{\frac{1}{p - 1}} \right),
\)
where $\Delta$ denotes the initial suboptimality, and $\sigma$, $\sigma_h$ denote the scale of noise in gradients and Hessians, respectively. This bound improves over the best known complexity for first-order methods \parencite{zhang2020adaptive,Hubler2024clip_to_norm} by a factor of $1/\varepsilon$ uniformly for all $p\in (1,2]$, demonstrating that second-order methods can yield provable advantages even in the heavy-tailed regime; see Figure~\ref{ch4:fig:complexities_chapter_4}. Moreover, our result implies that increasing the order of the algorithm (e.g., using higher-order derivatives) cannot improve the rate under the same noise assumptions significantly.

\item \textbf{Optimal algorithm.} We develop a second-order stochastic algorithm that achieves the lower bound (up to constants) in the regime $L \leq \sigma_h$, where $L$ denotes the Lipschitz constant of the gradient. The algorithm attains the sample complexity 
\(
\mathcal{O}\left( \frac{\Delta (L + \sigma_h)}{\varepsilon^2} \left( \frac{\sigma}{\varepsilon} \right)^{\frac{1}{p - 1}} \right),
\)
and is, to the best of our knowledge, the first second-order method that provably works under $p$-BCM noise. Notably, our method does not require second-order smoothness of the objective. This result holds even in the classical setting $p = 2$, highlighting the broader implications of the approach.

\item \textbf{High-probability convergence via Hessian clipping.} We propose a clipped variant of our algorithm that incorporates both gradient and \emph{Hessian clipping}--the latter introduced here for the first time. We show that this variant achieves near-optimal sample complexity with high probability, incurring only a poly-logarithmic overhead in the failure probability. This significantly extends the robustness of SOSO methods, enabling strong guarantees not only in expectation but also with high confidence.

\end{itemize}

Our results are primarily theoretical and validated through synthetic experiments with controllable heavy-tailed noise. While one component of our algorithm (in its unclipped form) has previously been applied in reinforcement learning settings \parencite{Fatkhullin_SPGM_FND_2023}, its broader practical deployment remains an open challenge. In particular, several questions must be addressed to make these methods viable in real-world applications: How can we reduce hyperparameter tuning overhead? Are there default configurations that are robust across problem classes? Is clipping truly necessary, or can it be systematically replaced by normalization techniques, as recent work suggests in the first-order context \parencite{Hubler2024clip_to_norm}? We view our contributions as a foundational step towards a principled understanding of second-order methods under realistic, heavy-tailed noise conditions, and hope they spur further theoretical and algorithmic advances in this direction.

\section{Notation and Assumptions}
We use the common notation for natural numbers $\mathbb{N} = \left\{0, 1, \dots\right\}$, and let $[n] = \left\{1, 2, \dots, n\right\}$. Throughout this chapter, $d\in\mathbb{N}_{\geq1}$ denotes the dimension of the optimization problem \eqref{ch4:eq:main_problem_chapter_4}. We use ${\cO(\cdot), \Omega(\cdot)}$ as sample-complexity notation that preserves the dominating term in the target accuracy $\varepsilon>0$ along with the multiplicative constants such as initial suboptimality $\Delta$, smoothness constants $L_r$, $r \geq 1$, the variance of higher-order stochastic oracles $\sigma_r$, $r\geq 1$. When we write ${\bar\cO(\cdot)}$, we suppress the dependence on all parameters except for accuracy $\varepsilon$ in the dominating term.
As in \textcite{carmon2020lower}, for any $q$-th order tensor $T \in \R^{\otimes^q d}$, we define the support of $T$ as $\operatorname{supp}\{T\} \eqdef \left\{i \in [d]~|~ T_i \neq 0\right\}$, where $T_i$ is the $(q-1)$-th order tensor, given by $[T_i]_{j_1,\ldots,j_{q-1}} =T_{i,j_1,\ldots,j_{q-1}} $.

We make the following assumptions throughout the chapter.  
\begin{assumption}[Lower Boundedness]\label{ch4:assum:lower_bounded_chapter_4}
	The objective function $f$ is lower bounded by $f^{\inf} > -\infty$.
\end{assumption}

\begin{assumption}[{$L$-smoothness}]\label{ch4:assum:l_smooth_chapter_4}
	The objective function $f$ is $L$-smooth, i.e.,\ $f$ is differentiable and for all $x,y\in\R^d$, we have $\norm{\nabla f(x) -\nabla f(y)} \leq L \norm{x-y}$.
\end{assumption}

\begin{assumption}[$p$-BCM]\label{ch4:ass:p_BCM_chapter_4}
    We have access to stochastic gradients $\nabla f(x;\xi)$ such that $\Exp\left[\nabla f(x;\xi)\right] = \nabla f(x)$. Moreover, there exist $p \in (1, 2]$ and $\sigma > 0$ such that 
    $$
    \Exp\left[\norm{\nabla f(x;\xi) - \nabla f(x)}^{p}\right] \leq \sigma^p  \qquad \text{for all } x \in \R^d .
    $$
\end{assumption}

\begin{assumption}[$p$-BCM for Hessian]\label{ch4:ass:p_BCM_hess_chapter_4}
    The function $f$ is twice differentiable. For simplicity, we assume access to the full stochastic Hessian estimator $\nabla^2 f(x;\xi)$ satisfying $\Exp\left[\nabla^2 f(x;\xi)\right] = \nabla^2 f(x)$.\footnote{The algorithms developed in this chapter only multiply the stochastic Hessian estimator by a vector and therefore can be implemented using Hessian-vector products.} For the same $p\in(1,2]$ as in Assumption~\ref{ch4:ass:p_BCM_chapter_4}, there exists $\sigma_h>0$ such that
    $$
    \Exp\left[{\normop{\nabla^2 f(x;\xi) - \nabla^2 f(x)}^{p}}\right] \leq \sigma_h^p  \qquad \text{for all } x \in \R^d .
    $$
\end{assumption}
The lower-bound oracle model introduced below is more general: it provides stochastic estimators of all derivatives of orders $r\in[q]$, so the resulting lower bounds cover first-, second-, and higher-order stochastic algorithms.

\section{Lower Complexity Bounds for SOSO under $p$-BCM}
\label{ch4:sec:lower_bounds_chapter_4}

To establish lower bounds, we build on a technique developed in a series of works \parencite{arjevani2023lower, carmon2020lower, arjevani2020second}, which introduced a \textit{worst-case} non-convex function exhibiting the zero-chain property: starting from the origin, each algorithmic iteration reveals information about a single coordinate only. Leveraging this framework, \textcite{zhang2020adaptive} derived the first lower bounds for FOSO under heavy-tailed noise in the smooth nonconvex setting, using an unbiased stochastic-gradient oracle satisfying a bounded central $p$-th moment condition ($p$-BCM), as in Assumption~\ref{ch4:ass:p_BCM_chapter_4}.

In this chapter, we assume access to a \( q \)-th order stochastic oracle, meaning that each query yields stochastic approximations of derivatives of orders \( 1 \) to \( q \) at a given point \( x \). This oracle model was formally introduced by \textcite{arjevani2020second}. Building on their proof technique, we derive new lower bounds under the assumption that each stochastic derivative has a bounded \( p \)-th central moment.

\textbf{Function class.} We consider the class of \( q \)-times differentiable functions (\( q \geq 1 \)) satisfying two key properties: (i) functional boundedness, i.e., \( f(0) - f^{\inf} \leq \Delta \); and (ii) Lipschitz continuity of each derivative up to order \( q \), i.e., for all \( x, y \in \R^d \) and \( r \in [q] \), 
\[
{\normop{\nabla^{r}f(x) - \nabla^{r}f(y)}} \leq L_{r}\norm{x - y}.
\]
We denote this function class by \( {\cF}(\Delta, L_{1:q}) \), where \( L_{1:q} = \{L_1, L_2, \dots, L_q\} \). Note that the constant \( L \) from Assumption~\ref{ch4:assum:l_smooth_chapter_4} corresponds to \( L_1 \) in this notation.

\textbf{Oracle class.} We now define the oracle model. For a fixed function \( f \in {\cF}(\Delta, L_{1:q}) \), a stochastic \( q \)-th order oracle is defined as follows: for any \( x \in \R^d \) and \( \xi \sim {\cD} \),
\begin{equation*}
    \text{O}^q_f(x;\xi) \eqdef \left( f(x;\xi), \nabla f(x;\xi), \nabla^2 f(x;\xi), \dots, \nabla^q f(x;\xi) \right),
\end{equation*}
where each component satisfies the unbiasedness conditions:
\begin{align*}
\Exp_{\xi\sim {\cD}}[f(x;\xi)] &= f(x),\\
\Exp_{\xi\sim {\cD}}[\nabla^r f(x;\xi)] &= \nabla^r f(x), \quad \forall~r \in [q].
\end{align*} We define the oracle class \( \cO_q(f, \sigma_{1:q}) \) as the set of all such stochastic oracles for which the \( p \)-th moment of the estimation error (with \( p \in (1,2] \)) satisfies
\begin{equation*}
    \Exp_{\xi\sim {\cD}}\left[{\normop{\nabla^r f(x;\xi) - \nabla^r f(x)}^p}\right] \leq \sigma_r^p, \quad \forall~r \in [q],~~ \forall~x \in \R^d.
\end{equation*}
For \( r = 1 \), this recovers \( \sigma_1 = \sigma \) from Assumption~\ref{ch4:ass:p_BCM_chapter_4}, and for \( r = 2 \), we similarly have \( \sigma_2 = \sigma_h \) from Assumption~\ref{ch4:ass:p_BCM_hess_chapter_4}.

\textbf{Algorithm class.} To prove lower bounds, we will work with a class of \textit{zero-respecting} algorithms. The formal definition is provided below. The reason why this class of algorithms is interesting for us is based on its property: at each iteration of such an algorithm, at most one new coordinate can be  revealed. 

\begin{definition}
    We call a stochastic $q$-th order algorithm $A$ zero-respecting if, for any function $f$ and any $q$-th order oracle $\text{O}^q_f$, the iterates $\{x^k\}_{k\in\mathbb{N}}$ generated by $A$ via queries to $\text{O}^q_f$ satisfy
    \[
        \operatorname{supp}\{x^k\} \subseteq \bigcup_{i<k}\operatorname{supp}\{\text{O}^q_f(x^i;\xi^i)\},
        \qquad \forall k\in\mathbb{N},
    \]
    with probability one with respect to the randomness of the algorithm and the realizations of $\{\xi^{k}\}_{k \in \mathbb{N}}$.
\end{definition}

\begin{theorem}
\label{ch4:thm:lower_bounds_chapter_4}
Let $q\in\mathbb{N}_{\geq1}$, $p\in(1,2]$, $\Delta>0$,
$L_{1:q}\eqdef(L_1,\dots,L_q)\in\mathbb{R}_{++}^q$, and
$\sigma_{1:q}\eqdef(\sigma_1,\dots,\sigma_q)\in\mathbb{R}_{++}^q$. Define
\[
\mathfrak M_q\eqdef
\min\left\{
\min_{r\in\{2,\dots,q\}}
\left(\frac{\sigma_r}{\sigma_1}\right)^{\frac1{r-1}},
\min_{r'\in[q]}
\left(\frac{L_{r'}}{\varepsilon}\right)^{\frac1{r'}}
\right\},
\]
where the minimum over the empty set is understood to be $+\infty$ when
$q=1$. There exist constants $c_q,C_q>0$, depending only on $q$ and
the numerical constants of the zero-chain construction, such that if
$0<\varepsilon\leq c_q\sigma_1$ and
$\Delta\mathfrak M_q/\varepsilon\geq C_q$, then there exist
$f\in{\cF}(\Delta,L_{1:q})$ and a corresponding noisy oracle
$\text{O}^q_f\in\cO_q(f,\sigma_{1:q})$ such that, for any
zero-respecting stochastic $q$-th-order algorithm, the number of oracle
queries required to find an $\varepsilon$-stationary point with constant
probability is lower bounded by
\begin{equation}
    \Omega(1)\cdot\frac{\Delta}{\varepsilon}
    \left(\frac{\sigma_1}{\varepsilon}\right)^{\frac{p}{p-1}}
    \mathfrak M_q.
\end{equation}
Moreover, this lower bound is realized by a construction of dimension
\[
\Theta\left(\frac{\Delta}{\varepsilon}\mathfrak M_q\right)
\].
\end{theorem}
The proof of Theorem~\ref{ch4:thm:lower_bounds_chapter_4} is deferred to Appendix~\ref{ch4:app:missing_proofs_lower_bounds_chapter_4}. In Theorem~\ref{ch4:thm:lower_bounds_chapter_4} we provide lower bounds for any $q \in \mathbb{N}_{\geq 1}$, which corresponds to the order of an oracle used in an optimization algorithm. 
When the oracle order satisfies $q \ge 2$, the algorithm gains access to stochastic curvature information, enabling more accurate prediction of the function or gradient behavior between consecutive iterates. This richer information can potentially lead to tighter worst-case complexity limits compared to first-order methods.
Since the main focus of this chapter is the second-order setting, we now discuss the complexity result from Theorem~\ref{ch4:thm:lower_bounds_chapter_4} for $q = 2$. For second-order methods, $q=2$, under Assumptions~\ref{ch4:assum:l_smooth_chapter_4},~\ref{ch4:ass:p_BCM_chapter_4},~\ref{ch4:ass:p_BCM_hess_chapter_4} and, additionally, (possibly) assuming the Hessian is Lipschitz continuous with constant $L_h$, we have 
\begin{equation}
\label{ch4:eq:lower_bound_second_order_chapter_4}
\Omega(1)\cdot \min\left\{
\begin{aligned}
 \frac{\Delta \sigma_h}{\varepsilon^2}\left(\frac{\sigma}{\varepsilon}\right)^{\frac{1}{p-1}},\frac{\Delta L \sigma}{\varepsilon^{3}} \left(\frac{\sigma}{\varepsilon}\right)^{\frac{1}{p-1}},\frac{\Delta \sqrt{L_h}\sigma}{\varepsilon^{\nicefrac{5}{2}}} \left(\frac{\sigma}{\varepsilon}\right)^{\frac{1}{p-1}}
\end{aligned}
\right\},
\end{equation}
where the parameters satisfy $\sigma = \sigma_1$, $\sigma_h = \sigma_2$, $L = L_1$, $L_2 = L_h$. We will now discuss each term separately, in different regimes.

\textbf{First-order optimization with Lipschitz gradient.} When we do not have access to the second-order information and the objective function has Lipschitz continuous gradient only, but not any higher derivative, our lower bounds \eqref{ch4:eq:lower_bound_second_order_chapter_4} reduces to the second term:
\begin{equation*}
    \Omega\left(\frac{\Delta L \sigma}{\varepsilon^{3}} \left(\frac{\sigma}{\varepsilon}\right)^{\frac{1}{p-1}}\right).
\end{equation*}
This bound exactly matches (up to a numerical constant factor) the previously known lower bounds for first-order stochastic optimization under $p$-BCM \parencite{zhang2020adaptive}, and the corresponding upper bound achieved by \algname{normalized SGD} \parencite{Hubler2024clip_to_norm}. Several other algorithms, like \algname{Clip-SGD} and \algname{NSGD with momentum and clipping}, also nearly match this complexity lower bound, up to additional polylogarithmic terms in $1/\varepsilon$ and polynomial terms in ${\beta}$, $L$ and $\sigma$ \parencite{zhang2020adaptive,cutkosky2021high,nguyen2023improved}.

\paragraph{First-order optimization under higher-order smoothness.} The second case worth discussing is when we (still) do not have access to stochastic Hessian, but we know that the objective function has Lipschitz continuous second-order derivatives. Then our lower bound \eqref{ch4:eq:lower_bound_second_order_chapter_4} becomes 
\begin{equation*}
    \Omega\left(\min\left\{\frac{\Delta L \sigma}{\varepsilon^{3}} \left(\frac{\sigma}{\varepsilon}\right)^{\frac{1}{p-1}}, \frac{\Delta \sqrt{L_h}\sigma}{\varepsilon^{\nicefrac{5}{2}}} \left(\frac{\sigma}{\varepsilon}\right)^{\frac{1}{p-1}} \right\}\right).
\end{equation*}
The last term in the above bound nearly matches with the complexity of the method called \algname{NIGT} with clipping, by \textcite{cutkosky2021high}, in terms of its dependence on $1/\varepsilon$, up to logarithmic factors. Unfortunately, there is still a small discrepancy with the upper bound in the dependence on other parameters and with the fact that \textcite{cutkosky2021high} use a slightly stronger (non-central) assumption $\Exp\left[\norm{\nabla f(x;\xi)}^p\right] \leq G^p$ instead of our $p$-BCM. We believe our bound is tight and can be achieved with a more careful analysis of \algname{NIGT} type algorithm under the $p$-BCM assumption.

\paragraph{Second-order optimization with only Lipschitz gradient.}
Finally, we wish to discuss an important setting, which is rarely discussed in the literature on second-order optimization. Higher-order smoothness can be challenging to verify in practice. Moreover, even when it is possible, the estimates of $L_h$ can be too large to be useful. Thus we pay attention to the case when we have access to {a} stochastic Hessian, but the second derivative can be non-continuous, i.e., $L_h \rightarrow +\infty$. Our lower bound \eqref{ch4:eq:lower_bound_second_order_chapter_4} provides new insights on this interesting setting, simplifying to
\begin{equation*}
    \Omega\left( \min\left\{\frac{\Delta \sigma_h}{\varepsilon^2}\left(\frac{\sigma}{\varepsilon}\right)^{\frac{1}{p-1}}, \frac{\Delta L \sigma}{\varepsilon^{3}} \left(\frac{\sigma}{\varepsilon}\right)^{\frac{1}{p-1}}\right\}\right).
\end{equation*}
As we can see, the first term corresponding to the second-order information has better dependence on {the accuracy} $\varepsilon$, and can be potentially much smaller than the second term. Moreover, this complexity is not impacted with the constants corresponding to higher-order smoothness, which potentially gives us the possibility of designing second-order methods with complexity that does not depend on $L_h$ and higher order smoothness constants. However, this is merely a lower bound, which does not give us any algorithmic solution yet. The question we investigate in the next section is 
\begin{center}
\textit{Can we design an algorithm handling heavy-tailed noise with complexity ${\cO\left(\frac{\Delta \sigma_h}{\varepsilon^2}\left(\frac{\sigma}{\varepsilon}\right)^{\frac{1}{p-1}}\right)}$? }
\end{center}

\section{Near-optimal Second-Order Method under $p$-BCM}
\label{ch4:sec:NSGD_SOM_conv_chapter_4}

In the last decade, {second-order stochastic optimization (SOSO)} has been extensively studied under stronger noise assumptions. We will first review the existing approaches, particularly focusing on the development in the non-convex setup. First, \textcite{tripuraneni2018stochastic} proposed and analyzed the \algname{Stochastic Cubic Newton (SCN)} method, achieving $\bar\cO(\varepsilon^{-\nicefrac{7}{2}})$ sample complexity. They require a strong bounded noise assumption, use large mini-batches and importantly, this complexity is not order-optimal. Later \textcite{arjevani2020second} established lower bounds for SOSO for finding first-order stationary points and proposed two algorithms: \algname{SGD with HVP-RVR} and a \algname{subsampled cubic-regularized trust-region method with HVP-RVR}. Both algorithms achieve $\bar\cO(\varepsilon^{-3})$  
sample complexities, which is minimax optimal when the variances of stochastic gradient and Hessian are bounded. Besides strong noise assumptions, their algorithms require large Hessian batch sizes, since they use the \algname{HVP-RVR} subroutine to update the momentum term. To overcome the large batch-size requirement, \textcite{BetterSGDUsingSOM_Tran_2021} designed a simpler algorithm, \algname{SGDHess}, with Hessian-corrected momentum (and optional normalization). Thanks to the Hessian corrected momentum term, these methods also achieve $\bar\cO(\varepsilon^{-3})$ sample complexity with any batch size $\geq 1$. Recently, \textcite{chayti2024improving} proposed a variant of \algname{SCN} with momentum without a large batch requirement. However, their analysis requires the bounded noise assumption and only achieves $\bar\cO(\varepsilon^{-\nicefrac{7}{2}})$ sample complexity.\footnote{However, in the noiseless case their complexity has better $\bar \cO(\varepsilon^{-\nicefrac{3}{2}})$ iteration complexity compared to $\bar \cO(\varepsilon^{-2})$ for \algname{SGDHess}.}

To summarize, {all the above-mentioned works} require strong noise assumptions: bounded noise and bounded variance. Moreover, their sample complexities also depend on second-order smoothness $L_h$, which can be potentially infinite, as discussed in the previous section. We refer to Table~\ref{ch4:tab:summary_chapter_4} for the summary of existing results. 

\begin{table}[t]
    \centering
    \caption{Summary of sample complexities of stochastic second-order (SOSO) methods for finding an $\varepsilon$-stationary point, {with high probability or in expectation}; that is, the number of stochastic gradient and Hessian evaluations needed to find $\bar x$ with $\norm{\nabla f(\bar x)} \leq \varepsilon.$ The column ``$p$'' indicates the range of moments in Assumptions~\ref{ch4:ass:p_BCM_chapter_4}, \ref{ch4:ass:p_BCM_hess_chapter_4} for which the result holds, that is when $p=2$, the corresponding result holds only under bounded variance. When $p=\infty$ it means that at least for stochastic gradient or stochastic Hessian bounded noise assumption is required, which corresponds to all moments being bounded. We are not aware of any prior works for SOSO for $p<2.$ The column ``\textbf{HP?}'' denotes whether the high probability guarantee with polylogarithmic dependence on the inverse of failure probability $1/\beta$ is available. }
    \label{ch4:tab:summary_chapter_4}
 \begin{threeparttable}
\resizebox{\textwidth}{!}{%
\begin{tabular}{cccc}

\bf Algorithm & \bf Sample Complexity &  $p$  & \bf HP?\\ \hline \hline 
\begin{tabular}{c} \algname{SCN} \\ \tabcite{tripuraneni2018stochastic} \end{tabular} & $\frac{\Delta\sqrt{L_h}\sigma^2}{\varepsilon^{\nicefrac{7}{2}}} $ & {\color{BrickRed}$\infty$} & \cmark$^{(1)}$ \\ 
\begin{tabular}{c} \algname{SGD with HVP-RVR} \\ \tabcite{arjevani2020second} \end{tabular} & $\frac{\Delta\sigma\sigma_h}{\varepsilon^{3}} + \frac{\Delta \sqrt{L_h}\sigma}{\varepsilon^{\nicefrac{5}{2}}}+ \frac{\Delta L}{\varepsilon^2}$ & {\color{BrickRed}$2$} & \xmark\\ 
\begin{tabular}{c} \algname{SN$^{{\color{niceblue}(2)}}$ with HVP-RVR} \\ \tabcite{arjevani2020second} \end{tabular} &  $\frac{\Delta\sigma\sigma_h}{\varepsilon^{3}} + \frac{\Delta \sqrt{L_h}\sigma}{\varepsilon^{\nicefrac{5}{2}}}+ \frac{\Delta \sigma_h}{\varepsilon^2}$  & {\color{BrickRed}$2$} & \xmark\\
\begin{tabular}{c} \algname{NSGDHess} \\ \tabcite{BetterSGDUsingSOM_Tran_2021} \end{tabular} &  $\frac{\sigma^3}{\varepsilon^{3}}+ \frac{\Delta\sigma\sigma_h}{\varepsilon^{3}} + \frac{\Delta \sqrt{L_h}\sigma}{\varepsilon^{\nicefrac{5}{2}}}+ \frac{\Delta \sigma_h}{\varepsilon^2}$  & {\color{BrickRed}$2$}& \xmark \\
\begin{tabular}{c} \algname{SCN with IT-HB} \\ \tabcite{chayti2024improving} \end{tabular} &  $\frac{\Delta\sqrt{L_h}}{\varepsilon^{\nicefrac{3}{2}}} + \frac{\Delta L_h^{\nicefrac{1}{4}}\sigma_h^{\nicefrac{1}{2}}}{\varepsilon^{\nicefrac{7}{4}}}+ \frac{\Delta \sqrt{L}\sigma^2}{\varepsilon^{\nicefrac{7}{2}}}$  & {\color{BrickRed}$\infty$}& \xmark\\
\hline\hline
 \cellcolor{linen} \begin{tabular}{c} Lower Bounds \\   Theorem~\ref{ch4:thm:lower_bounds_chapter_4} \end{tabular}
& \cellcolor{linen}  $\min\left\{\frac{\Delta \sigma_h}{\varepsilon^2}, \frac{\Delta L \sigma}{\varepsilon^{3}} , \frac{\Delta \sqrt{L_h}\sigma}{\varepsilon^{\nicefrac{5}{2}}} \right\} \left(\frac{\sigma}{\varepsilon}\right)^{\frac{1}{p-1}}$ &\cellcolor{linen} {\color{PineGreen}$(1,2]$} & \cellcolor{linen}\\
\hline
 \cellcolor{linen} \begin{tabular}{c} \algname{NSGDHess} \\   Theorem~\ref{ch4:thm:nsgdm_hess_chapter_4} \end{tabular}
& \cellcolor{linen} $\frac{\Delta (L +\sigma_h)}{\varepsilon^2} + \left(\frac{\Delta (L + \sigma_h)}{\varepsilon^2}  + \frac{\sigma}{\varepsilon}\right)\left(\frac{\sigma}{\varepsilon}\right)^{\frac{1}{p-1}}$  &\cellcolor{linen} {\color{PineGreen}$(1,2]$} & \cellcolor{linen} \xmark \\
 \cellcolor{linen} \begin{tabular}{c} \algname{Clip NSGDHess}  \\  Theorem~\ref{ch4:thm:main_theorem_hp_chapter_4} \end{tabular} & \cellcolor{linen}$\left(\frac{\sqrt{\Delta_0(L+\sigma_h)}+ \sigma}{\varepsilon}\right)^{2 +\frac{1}{p-1}}$ & \cellcolor{linen} {\color{PineGreen}$(1,2]$} & \cellcolor{linen} \cmark \\ 
\end{tabular}
}
\begin{tablenotes}
        \scriptsize \item [{(1)}]
        \textcite{tripuraneni2018stochastic} provide analysis under stronger  assumptions, which implies the noise has light tails.  
        \item [{\color{niceblue}(2)}] \algnametiny{SN} = Subsampled Newton.
    \end{tablenotes}
    \end{threeparttable}
\end{table}

\begin{algorithm}[ht!]
\caption{\algname{NSGDHess} (Normalized SGD with Hessian correction)}
\label{ch4:alg:NSGD_SOM_chapter_4}
	\begin{algorithmic}[1]
		\State \textbf{Input:} Starting point $x^0 \in \R^d$, a vector $g^0 \in \R^d$, a stepsize $\gamma > 0$, momentum parameter $\eta\in(0,1]$, and the number of iterations $K$.
		\State $x^1 = x^0 - \gamma\frac{g^0}{\norm{g^0}}$\footnotemark
		\For{ $k = 1, 2, \ldots, K-1$}
		\State Sample ${q^k} \sim {\cU}\left([0,1]\right)$
		\State $\hat{x}^k = {q^k} x^k + (1-{q^k})x^{k-1}$
		\State Sample $\xi^{k}, \hat{\xi}^{k} \sim {\cD}$ independently
		\State $g^k = (1-\eta)\left( g^{k-1} + \nabla^2 f(\hat{x}^k; \hat{\xi}^{k})(x^k - x^{k-1})\right) + \eta \nabla f(x^k;\xi^{k})$
		\State $x^{k+1} = x^k -\gamma \frac{g^k}{\norm{g^k}}$
		\EndFor\stepEightBdeletion
	\end{algorithmic}
\end{algorithm}
\footnotetext{Throughout this chapter, we use the convention $\frac{g^k}{\norm{g^k}}=0$ whenever $g^k=0$.}

To design an optimal second-order algorithm without these limitations, we take an inspiration from the recent developments in reinforcement learning \parencite{salehkaleybar-et-al22,Fatkhullin_SPGM_FND_2023}. Their algorithms, \algname{SHARP} and \algname{(N)HAR-PG}, are variants of the \algname{SGDHess} method \parencite{BetterSGDUsingSOM_Tran_2021} with an optional normalization step. We adapt one variant of their method to our general stochastic optimization setting and present it in Algorithm~\ref{ch4:alg:NSGD_SOM_chapter_4}. This algorithm computes a random interpolation point $\hat{x}^k$ uniformly distributed between consecutive iterates $x^{k-1}$ and $x^k$. The method then evaluates a stochastic gradient and a stochastic Hessian-vector product using independent samples $\xi^{k}$ and $\hat{\xi}^{k}$. Finally, a recursive momentum is constructed before applying the normalization step.\footnote{Normalization step for this method is important for our analysis under $p$-BCM. However, it can be removed when $p=2$, or when additional clipping is used, as in the next section.} We are now ready to state the convergence guarantee for Algorithm~\ref{ch4:alg:NSGD_SOM_chapter_4}.

\begin{theorem}
\label{ch4:thm:nsgdm_hess_chapter_4}
{
    Suppose Assumptions~\ref{ch4:assum:lower_bounded_chapter_4},~\ref{ch4:assum:l_smooth_chapter_4},~\ref{ch4:ass:p_BCM_chapter_4}, and~\ref{ch4:ass:p_BCM_hess_chapter_4} hold, and define
    \[
        \Delta_0 \eqdef f(x^0)-f^{\inf},
        \qquad
        \mathcal E_0\geq
        \Exp\left[\norm{g^0-\nabla f(x^0)}\right].
    \]
    For $K\geq1$, choose
    \[
        \eta=\min\{1,\eta_{\mathrm{eff}}\},
        \qquad
        \gamma=
        \sqrt{\frac{\Delta_0\eta^{1/p}}{(L+\sigma_h)K}},
    \]
    where
    \[
        \eta_{\mathrm{eff}}
        \eqdef \max\left\{
            \left(\frac{\mathcal E_0}{K\sigma}\right)^{\frac{p}{2p-1}},
            \left(\frac{\Delta_0(L+\sigma_h)}{K\sigma^2}\right)^{\frac{p}{2p-1}}
        \right\}.
    \]
    Then the iterates of Algorithm~\ref{ch4:alg:NSGD_SOM_chapter_4}
    satisfy
    \begin{align*}
        \frac1K\sum_{k=0}^{K-1}
        \Exp\left[\norm{\nabla f(x^k)}\right]
        =\cO\Bigg(
            &\sqrt{\frac{\Delta_0(L+\sigma_h)}{K}}
            +\sigma
            \left(
                \frac{\Delta_0(L+\sigma_h)}{K\sigma^2}
            \right)^{\frac{p-1}{2p-1}}\\
            &+\frac{\mathcal E_0}{K}
            +\sigma
            \left(
                \frac{\mathcal E_0}{K\sigma}
            \right)^{\frac{p-1}{2p-1}}
        \Bigg).
    \end{align*}
    In particular, take
    \[
        g^0=\frac1{B_{\mathrm{init}}}
        \sum_{j=1}^{B_{\mathrm{init}}}\nabla f(x^0;\xi_j^0),
        \qquad
        B_{\mathrm{init}}
        =
        \max\left\{
            1,
            \left(\frac{\sigma}{\varepsilon}\right)^{\frac{p}{p-1}}
        \right\}.
    \]
    Then $\mathcal E_0\leq
    2\sigma/B_{\mathrm{init}}^{(p-1)/p}$, and the expected average
    gradient norm is at most $\varepsilon$ with total sample
    complexity
    \begin{equation}
    \label{ch4:eq:convergence_rate_thm_chapter_4}
        \cO\left(
            \frac{\Delta_0(L+\sigma_h)}{\varepsilon^2}
            +
            \frac{\Delta_0(L+\sigma_h)}{\varepsilon^2}
                \left(\frac{\sigma}{\varepsilon}\right)^{\frac{1}{p-1}}
            +
            \frac{\sigma}{\varepsilon}
                \left(\frac{\sigma}{\varepsilon}\right)^{\frac{1}{p-1}}
        \right).
    \end{equation}
}
\end{theorem}
The proof is deferred to Appendix~\ref{ch4:app:missing_proofs_NSGD_SOM_chapter_4}. We also investigate different choices of $g^{0}$ and their influence on the total sample complexity rate in Appendix~\ref{ch4:app:missing_proofs_NSGD_SOM_chapter_4}. We find that  while the initial batch size to estimate $g^{0}$ is not necessary for convergence, it  helps to slightly improve the total sample complexity.

\paragraph{Discussion:} 
Comparing this result with Theorem~\ref{ch4:thm:lower_bounds_chapter_4}, we observe that our upper bound \eqref{ch4:eq:convergence_rate_thm_chapter_4} matches our lower bound--the first term in \eqref{ch4:eq:lower_bound_second_order_chapter_4} in terms of the target accuracy $\varepsilon$. Moreover, when $\Delta(L+\sigma_h) \geq \sigma \varepsilon$, we have a tight upper bound, which exactly matches the lower bound from the previous section in the leading term (up to a numerical constant). Moreover, since our upper bound does not depend on constant $L_h$, it explains our observations in Section~\ref{ch4:sec:lower_bounds_chapter_4} regarding the limit case $L_h \rightarrow \infty, $ and answers the raised question affirmatively. We refer to Table~\ref{ch4:tab:summary_chapter_4} for more technical comparison with prior work.

\section{Hessian Clipping for High-probability Convergence}
\label{ch4:sec:high_probability_analysis_chapter_4}

In Section~\ref{ch4:sec:NSGD_SOM_conv_chapter_4} we provided an in-expectation guarantee for Algorithm~\ref{ch4:alg:NSGD_SOM_chapter_4}. While in-expectation guarantees are useful when a method can be run multiple times to analyze its average-case behavior, {they can be impractical when computational or time constraints permit only a single run}. Thus, our main goal in this section is to conduct a high-probability analysis of Algorithm~\ref{ch4:alg:NSGD_SOM_chapter_4} or its modification. 

Since we work in the unbounded variance case, which corresponds to heavy-tailed noise, it is important to ensure that the analyzed method is robust to such noise. Following prior work \parencite{gorbunov2020stochastic,cutkosky2021high, sadiev2023high, nguyen2023improved}, we propose a new algorithm, called \algname{Clip NSGDHess}. It is formally presented as Algorithm~\ref{ch4:alg:NSGD_SOM_clipped_chapter_4}; compared to Algorithm~\ref{ch4:alg:NSGD_SOM_chapter_4}, it also incorporates gradient clipping and Hessian clipping. Formally, a clipping operator (clipping for short) is defined as 
\begin{equation}\label{ch4:eq:clipping_def_chapter_4}
    \texttt{clip}(v,\lambda)=
    \begin{cases}
        \min\left\{1,\dfrac{\lambda}{\norm{v}}\right\}v,
            & v\neq0,\\[1mm]
        0,  & v=0,
    \end{cases}
    \qquad v\in\R^d,
\end{equation}
where $\lambda >0$ is called the clipping level/threshold. 

\begin{algorithm}[ht!]
\caption{\algname{Clip NSGDHess} (Normalized SGD with Hessian correction and clipping)}
\label{ch4:alg:NSGD_SOM_clipped_chapter_4}
	\begin{algorithmic}[1]
		\State \textbf{Input:} Starting point $x^0 \in \R^d$, a stepsize $\gamma > 0$, momentum parameter $\eta\in(0,1]$, clipping levels $\lambda > 0$, ${\bar\lambda_h} > 0$, and the number of iterations $K$.
		\State $x^1 = x^0$ and $g^0 = 0$
		\For{ $k = 1, 2, \ldots, K-1$}
		\State Sample ${q^k} \sim {\cU}\left([0,1]\right)$
		\State $\hat{x}^k = {q^k} \, x^k + (1-{q^k}) \, x^{k-1}$
		\State Sample $\xi^{k}, \hat{\xi}^{k} \sim {\cD}$ independently
		\State $\begin{aligned}[t]
			g^k ={}&(1-\eta)\left(g^{k-1} + \gamma\,\texttt{clip}\left(\gamma^{-1}\nabla^2 f(\hat{x}^k;\hat{\xi}^{k})(x^k-x^{k-1}),{\bar\lambda_h}\right)\right)\\
			&\quad+\eta\,\texttt{clip}\left(\nabla f(x^k;\xi^{k}),\lambda\right)
		\end{aligned}$
		\State $x^{k+1} = x^k -\gamma \frac{g^k}{\norm{g^k}}$
		\EndFor\stepEightBdeletion
	\end{algorithmic}
\end{algorithm}

As in Algorithm~\ref{ch4:alg:NSGD_SOM_chapter_4}, we use normalization in the gradient update (line 8), which allows the next point $x^{k+1}$ to stay in the ball with center at $x^k$ of radius $\gamma$. The main difference between Algorithms~\ref{ch4:alg:NSGD_SOM_chapter_4} and~\ref{ch4:alg:NSGD_SOM_clipped_chapter_4} is in the momentum term. To enhance robustness, we do not only clip the gradient, but also the Hessian-vector product. As we mentioned above, gradient clipping is a common approach, but Hessian clipping is new, and it allows us to provide the first high-probability guarantees for second-order stochastic optimization without light-tail noise assumptions. {It is worth drawing attention} to the term
\begin{equation}\label{ch4:eq:hessian_clipping_term_chapter_4}
\gamma \, \texttt{clip}\left(\gamma^{-1} \nabla^2 f(\hat{x}^k; \hat{\xi}^{k})(x^k - x^{k-1}), {\bar\lambda_h} \right).
\end{equation}
Observe that clipping the entire Hessian based on its operator norm is very costly since it typically requires $\cO(d^2)$ arithmetic operations (e.g., using power iteration \parencite{golub2013matrix} or Lanczos algorithm \parencite{lanczos1950iteration,saad2011numerical}). 
Instead of clipping the stochastic Hessian, we clip the Hessian-vector product, which is computationally tractable, and only requires $\cO(d)$ operations. Thus, the computation of this clipped Hessian-vector product can be easily implemented via backpropagation and subsequent ``vector'' clipping \eqref{ch4:eq:clipping_def_chapter_4}. Another reason why we prefer to use this form of clipping comes from the analysis; we want to preserve the following useful property: 
\begin{eqnarray}\label{ch4:eq:expectation_of_hv_chapter_4}
    \Exp_{{q^k}, \hat{\xi}^{k}}\left[\nabla^2 f(\hat{x}^k; \hat{\xi}^{k})(x^k - x^{k-1})\right] &=& \Exp_{{q^k}}\left[\Exp_{ \hat{\xi}^{k}}\left[\nabla^2 f(\hat{x}^k; \hat{\xi}^{k})\right](x^k - x^{k-1})\right] \notag \\
    &=& \Exp_{{q^k}}\left[\nabla^2 f(\hat{x}^k)(x^k - x^{k-1})\right] \notag \\
    &=&\int^1_0 \nabla^2 f(q x^k + (1-q)x^{k-1})(x^k - x^{k-1})\,dq \notag \\
    &=& \nabla f(x^k) - \nabla f(x^{k-1}).
\end{eqnarray}
This would not hold were we to use direct Hessian matrix clipping. Instead, we prove that the expectation of \eqref{ch4:eq:hessian_clipping_term_chapter_4} will approximately be equal to \eqref{ch4:eq:expectation_of_hv_chapter_4} when ${\bar\lambda_h}$ is selected properly. By smoothness and the gradient step, we have $\|\nabla f(x^{k}) - \nabla f(x^{k-1})\| \leq L\|x^{k}-x^{k-1}\|= L\gamma$, meaning the expectation is bounded. 
To provide the analysis, we need to rewrite the clipped Hessian-vector product as follows 
$$\gamma \, \texttt{clip}\left(\gamma^{-1} \nabla^2 f(\hat{x}^k; \hat{\xi}^{k})(x^k - x^{k-1}), {\bar\lambda_h} \right) = \texttt{clip}\left( \nabla^2 f(\hat{x}^k; \hat{\xi}^{k})(x^k - x^{k-1}), {\gamma\bar\lambda_h} \right),$$
where we denote ${\lambda_h = \gamma\bar\lambda_h}$. The former formulation in the above formula is useful for implementation, since as we will see, $\lambda_{h} \sim \cO(\gamma\eta^{-\nicefrac{1}{p}})$ and $\lambda \sim \cO(\eta^{-\nicefrac{1}{p}})$. Thus, $\lambda$ and ${\bar\lambda_h}$ are of the same order $\sim\cO(\eta^{-\nicefrac{1}{p}}) = \cO(K^{\frac{1}{2p-1}})$, which makes it easier to tune these two parameters together in practice. The latter formulation will be used in the high-probability analysis. We are now ready to state the main theorem of this section.

\begin{theorem}
\label{ch4:thm:main_theorem_hp_chapter_4}
{
Let Assumptions~\ref{ch4:assum:lower_bounded_chapter_4},~\ref{ch4:assum:l_smooth_chapter_4},~\ref{ch4:ass:p_BCM_chapter_4}, and~\ref{ch4:ass:p_BCM_hess_chapter_4} hold, and define $\Delta_0 \eqdef f(x^{0}) - f^{\inf}$. Let $\beta \in (0,1]$ and $K \geq 1$ satisfy $\log(8K/\beta) \geq 1$.
Suppose Algorithm~\ref{ch4:alg:NSGD_SOM_clipped_chapter_4} is executed with
\[
    \begin{aligned}
    \eta &= K^{-\frac{p}{2p-1}},\qquad
    \lambda = \max\left\{4\sqrt{L\Delta_0}, \frac{\sigma}{\eta^{1/p}}\right\},\\
    \lambda_h &= \frac{2\gamma(L+\sigma_h)}{\eta^{1/p}},
    \qquad
    \bar\lambda_h=\frac{\lambda_h}{\gamma}
    =\frac{2(L+\sigma_h)}{\eta^{1/p}}.
    \end{aligned}
\]
and stepsize
\[
\begin{aligned}
\gamma = \min\Bigg\{&
\frac{1}{2}\sqrt{\frac{\Delta_0}{LK}},\,
\frac{\eta}{12}\sqrt{\frac{\Delta_0}{L}},\,
\frac{1}{1408\eta K\log\frac{8K}{\beta}}\sqrt{\frac{\Delta_0}{L}},\\
&
\frac{\Delta_0}{352\sigma\eta^{\frac{p-1}{p}}K\log\frac{8K}{\beta}},\,
\sqrt{\frac{\Delta_0\eta^{1/p}}
{968(L+\sigma_h)K\log\frac{8K}{\beta}}}
\Bigg\}.
\end{aligned}
\]
Then, with probability at least $1-\beta$, the generated iterates satisfy
\[
    \frac{1}{K} \sum_{k=0}^{K-1} \norm{\nabla f(x^k)}
    \leq \frac{2\Delta_0}{\gamma K}.
\]
Consequently,
\begin{equation*}
    \frac{1}{K} \sum_{k=0}^{K-1} \norm{\nabla f(x^{k})}
    = \mathcal{O}\left(
        \frac{\sqrt{\Delta_0(L+\sigma_h)}+\sigma}
             {K^{\frac{p-1}{2p-1}}}
        \log\left(\frac{K}{\beta}\right)
    \right)
\end{equation*}
with high probability.
}
\end{theorem}
{The complete proof of Theorem~\ref{ch4:thm:main_theorem_hp_chapter_4} can be found in Appendix~\ref{ch4:app:missing_proofs_hp_chapter_4}. Disregarding problem-dependent constants and logarithmic factors, the momentum parameter is selected as $\eta = K^{-p/(2p-1)}$, and the active stochastic terms in the stepsize have the same polynomial scaling, $\gamma = \widetilde{\cO}(K^{-p/(2p-1)})$. The clipping levels $\lambda$ and $\bar\lambda_h$ both scale as $\cO(K^{1/(2p-1)})$.}

\begin{corollary}
    \label{ch4:cor:high_probability_chapter_4}
    In the setting of Theorem~\ref{ch4:thm:main_theorem_hp_chapter_4}, Algorithm~\ref{ch4:alg:NSGD_SOM_clipped_chapter_4} ensures that
    \[
        \frac{1}{K}\sum_{k=0}^{K-1}\|\nabla f(x^{k})\| \leq \varepsilon,
    \]
    with probability at least $1-\beta$. To achieve this, the algorithm requires at most
    \[
        \cO\!\left(
            \left[
                \frac{\sqrt{\Delta_0(L+\sigma_h)}+\sigma}{\varepsilon}
                \log\!\left(\frac{1}{\beta\varepsilon}\right)
            \right]^{\frac{2p-1}{p-1}}
        \right)
        \quad \text{first/second-order oracle calls.}
    \]
\end{corollary}

According to the results of Theorem~\ref{ch4:thm:main_theorem_hp_chapter_4} and Corollary~\ref{ch4:cor:high_probability_chapter_4},  we have shown that Algorithm~\ref{ch4:alg:NSGD_SOM_clipped_chapter_4} has a near-optimal convergence rate in the dependence on $\varepsilon$. Moreover, the complexity has a logarithmic dependence on $\nicefrac{1}{\beta}$, which is better than a polynomial dependence that can be achieved by directly translating the in-expectation result of Theorem~\ref{ch4:thm:nsgdm_hess_chapter_4} to a high probability one using Markov's inequality. Unfortunately, the dependence on $\Delta_0, L, \sigma, \sigma_h$ does not match the  lower bound \eqref{ch4:eq:lower_bound_second_order_chapter_4}. 

\paragraph{Comparison with \textcite{UnboundedClippedNSGDM2023Liu}.} Recently, and under additional strong assumptions,  \textcite{UnboundedClippedNSGDM2023Liu} analyzed a momentum-based, variance-reduced first-order method, and obtained rates similar to those in our Corollary~\ref{ch4:cor:high_probability_chapter_4}. However, they assume individual smoothness (i.e., $\|\nabla f(x;\xi) - \nabla f(y;\xi)\| \leq \ell \|x - y\|$ for all $x, y\in \R^d$ and all realizations of the random variable $\xi$). However, individual smoothness is a very strong assumption which implies boundedness of the stochastic Hessian, i.e., $\|\nabla^2 f(x;\xi)\| \leq \ell$ a.s., for all $x\in \R^d$. We also note that under individual smoothness,  fast $p$-independent rates can be obtained, e.g., \textcite{lei2019stochastic} achieve $\bar \cO(\varepsilon^{-4})$ complexity for SGD under  individual smoothness.

\section{Conclusion and Future Work}

This chapter studied smooth nonconvex optimization with heavy-tailed stochastic gradients and higher-order oracle information. Lower bounds for zero-respecting stochastic \(q\)-th-order algorithms quantify the limitations of this information model. On the algorithmic side, \algname{NSGDHess} combines normalized steps with Hessian-corrected recursive momentum, while \algname{Clip-NSGDHess} clips both gradient and Hessian--vector contributions to obtain a batch-free high-probability method. The algorithms use Hessian--vector products rather than requiring a full Hessian solve.

The expectation and high-probability conclusions are distinct. The expectation theorem controls the trajectory average \(K^{-1}\sum_{k=0}^{K-1}\Exp[\norm{\nabla f(x^k)}]\) and matches the relevant lower-bound dependence in the stated dominant second-order regime. The clipped result controls the realized average \(K^{-1}\sum_{k=0}^{K-1}\norm{\nabla f(x^k)}\) with probability at least \(1-\beta\). Its dependence on accuracy is near-optimal up to logarithmic factors, but the dependencies on the initial gap, smoothness, and gradient and Hessian noise levels do not fully match the lower bound. Neither statement is a terminal-iterate guarantee, and the lower bound itself concerns the specified zero-respecting oracle class and theorem conditions.

Several limitations remain. The high-probability method introduces additional tuned parameters, depends on the selected confidence level, and does not recover the best deterministic second-order complexity. Hessian--vector evaluations and two clipping thresholds also create practical overhead. Closing these gaps may require different concentration tools, adaptive tuning, or preconditioned updates. Chapter~\ref{ch:chapter_5} refines the comparison by separating gradient and Hessian tail indices, considering weaker similarity and mean-squared smoothness structures, and developing first-order variance-reduced and double-clipped methods alongside sharper lower bounds.

\newpage
  \endgroup

\chapter{Tight Lower Bounds and Optimal Algorithms for Stochastic Nonconvex Optimization with Heavy-Tailed Noise}
\label{ch:chapter_5}
\thesischaptersourcenote
  {Tight Lower Bounds and Optimal Algorithms for Stochastic Nonconvex Optimization with Heavy-Tailed Noise}
  {fradin2025tight}
  \begingroup
  \thesischaptersetup{5}%
  
\begin{paperabstract}
We study stochastic nonconvex optimization under heavy-tailed noise. In this setting, the stochastic gradients only have bounded $p$-th central moment ($p$-BCM) for some $p \in \thesischfiveintof{1}{2}$. Building on the foundational work of \textcite{arjevani2022lower} in stochastic optimization, we establish tight sample complexity lower bounds for\stepEightBdeletion{} first-order methods under \emph{relaxed} mean-squared smoothness ($q$-WAS) and $\delta$-similarity ($(q, \delta)$-S) assumptions, allowing any exponent $q \in [1,2]$ instead of the standard $q = 2$. These results substantially broaden the scope of existing lower bounds. To complement them, we show that Normalized Stochastic Gradient Descent with Momentum Variance Reduction (\algname{NSGD-MVR}), a known algorithm, matches these bounds in expectation. Beyond expectation guarantees, we introduce a new algorithm, Double-Clipped \algname{NSGD-MVR}, which allows the derivation of high-probability convergence rates under weaker assumptions than in previous works. Finally, for second-order methods with stochastic Hessians satisfying bounded $q$-th central moment assumptions for some exponent $q \in [1, 2]$ (allowing $q \neq p$), we establish sharper lower bounds than previous works while improving over \textcite{sadiev2025second} (where only $p = q$ is considered) and yielding stronger convergence exponents. Together, these results provide a nearly complete complexity characterization of stochastic nonconvex optimization in heavy-tailed regimes.
\end{paperabstract}

\section{Introduction}
We consider the stochastic optimization problem
    \[ \min_{x \in \R^d} f(x), \qquad f(x) \eqdef \thesischfiveExpSub{\xi \sim \thesischfivecD}{f(x;\xi)}, \thesischfivenumberthis\label{ch5:eq:optim-problem_chapter_5} \]
where $d\geq 1$ is the dimension, $f \colon \R^d \to \R$ is a smooth possibly nonconvex objective function, and $\xi$ is a random variable drawn from an unknown distribution $\thesischfivecD$. In the nonconvex setting, our goal is to find an $\thesischfiveeps$--stationary point~\parencite{nesterov2018lectures}; that is, a vector $\bar{x}$ such that $\Exp[\norm{\nabla f(\bar{x})}] \le \thesischfiveeps$. Problems of the form~\eqref{ch5:eq:optim-problem_chapter_5} are pervasive in machine learning where $f\left(x;\xi\right)$ denotes the loss of a model with weights $x$ on a data sample $\xi \sim \mathcal{D}$, and $\thesischfivecD$ is the distribution of the training samples from the dataset~\parencite{bottou2018optimizationmethodsML}. While Gradient Descent (\algname{GD}) is known to achieve the optimal rate $\cO(\thesischfiveeps^{-2})$ for finding an $\thesischfiveeps$--stationary point \parencite{carmon2020lower}, it requires access to exact gradients $\nabla f(\cdot)$ which is infeasible in practice. Therefore, reliance on noisy gradients has become a gold standard approach, giving rise to stochastic gradient methods like Stochastic Gradient Descent (\algname{SGD})~\parencite{MR42668}. Under the standard $\sigma_1^2$--bounded variance assumption, \algname{SGD} is provably optimal~\parencite{ghadimi2013stochastic}, matching the lower bound $\Omega(\nicefrac{L_1 \Delta}{\thesischfiveeps^2} + \nicefrac{L_1 \Delta \sigma_1^2}{\thesischfiveeps^4})$~\parencite{arjevani2022lower} for $L_1$--smooth functions with %
$\smash{\Delta \eqdef f(x^0) - f^{\inf}}$.

However, empirical observations in modern machine learning, e.g., image classification~\parencite{simsekli2019heavytailedtheorystochasticgradient,pmlr-v97-simsekli19a,pmlr-v238-battash24a}, large language models~\parencite{10.5555/3495724.3497014,ahn2024linear} and reinforcement learning~\parencite{pmlr-v139-garg21b}, have shed light on the importance of \emph{heavy-tailed} noise as a more realistic setting than the standard bounded variance assumption. That is, noisy gradients only have bounded $p$--th central moment, $p \in \thesischfiveintof{1}{2}$. Such a relaxed moment assumption introduces challenges on both the algorithmic and theoretical sides, as \algname{SGD} may fail to converge when $p < 2$. In the heavy-tailed regime,~\textcite{10.5555/3495724.3497014} establish a lower bound of $\Omega(K^{-\nicefrac{(p - 1)}{(3 p - 2)}})$ on the convergence rate of first-order methods, where $K$ denotes the number of iterations. Yet, their analysis is confined to the standard stochastic setting, where algorithms access smooth (potentially nonconvex) functions via unbiased stochastic gradients. {Crucially, they do not address more structured assumptions on the gradient noise, such as \textit{mean-squared smoothness} or \textit{$\delta$--similarity}, which play a central role in the design and analysis of \emph{variance-reduced} methods and remain largely unexplored in the heavy-tailed regime.}

This chapter aims to: $(1)$ \textbf{develop lower complexity bounds for first-order variance-reduced methods in the heavy-tailed regime}, thereby closing an important open question from~\textcite{pmlr-v195-liu23c}, while advancing beyond the results of~\textcite{sun2024gradientnormalizationprovablybenefits}; $(2)$ \textbf{extend existing high-probability analysis to these new assumptions}, and finally, for second-order methods; $(3)$ \textbf{characterize how the tail indices $p$ (of gradients) and $q$ (of Hessians) jointly affect convergence in the heavy-tailed regime}.

Beyond the development of lower bounds, a natural question arises:

{\centering\textit{Does there exist an algorithm that provably matches, in expectation, the lower bounds developed in $(1)$?}\par}
Here, we answer this question affirmatively by revisiting a classical technique: gradient normalization. Remarkably, our results show that this simple method remains effective even in the heavy-tailed regime.

\subsection{Our Contributions}

\begin{table*}[t]
    \centering
    \caption[Sample complexities of stochastic methods for finding $\varepsilon$-stationary points]{Sample complexities of stochastic methods for finding an $\varepsilon$-stationary point (in expectation or with high probability). Stochastic gradients satisfy $p$-BCM with $p \in (1,2]$ (Assumption~\ref{ch5:ass:p-bounded-central-moment-gradient_chapter_5}). The column “$q$” specifies the relevant moment assumptions: Assumption~\ref{ch5:ass:mean-squared-smoothness_chapter_5} ($q$-WAS), Assumption~\ref{ch5:ass:mean-squared-smoothness-2_chapter_5} ($(q,\delta)$-S), or Assumption~\ref{ch5:ass:q-bounded-central-moment-hessian_chapter_5} ($q$-BCM); $q=\infty$ corresponds to bounded noise. The column “\textbf{HP?}” indicates whether a high-probability guarantee with only polylogarithmic dependence on $\nicefrac{1}{\beta}$ is available. Lower bounds are stated \emph{up to constant factors}. For the algorithm names: \algname{NSGD}: Normalized \algname{SGD}, \algname{Acc}: Accelerated, \algname{SFOM}: Stochastic First-Order Method, \algname{RM}: Recursive Momentum, \algname{MVR}: Momentum Variance Reduction, \algname{Hess}: Hessian-corrected momentum, \algname{Clip} / \algname{D-Clip}: Clipping / Double-Clipped.}
    \label{ch5:tab:summary_chapter_5}
\begin{minipage}{\textwidth}
\resizebox{\textwidth}{!}{%
\begin{tabular}{c|cccc}

\bf Setup& \bf Algorithm & \bf Sample Complexity & $q$& \bf HP?\\ \hline \hline  
    \multirow{10}{*}{\makecell{$q$-WAS}} &\begin{tabular}{c}  \algname{AccNSGD}  \\  \tabcite{pmlr-v195-liu23c} \end{tabular} & $\left(\frac{\sqrt{\bar L\Delta}+ \sigma_1}{\varepsilon}\right)^{2 +\frac{1}{p-1}}$ & {\color{BrickRed} $\ \ \ \ \infty$}$^{\color{blue}(1)}$ &  \ \thesischfivecmark \\ 
    &\begin{tabular}{c}  \algname{NSFOM with RM}  \\  \tabcite{he2025complexitynormalizedstochasticfirstorder} \end{tabular} & $\left(\frac{\Delta +\sigma_1^p+L_1+L_1^p+\bar{L}^p}{\varepsilon}\right)^{2 +\frac{1}{p-1}}$ $^{\color{blue}(2)}$ & {\color{BrickRed} $\ p=q$}&  \ \thesischfivexmark \\ 
 &
  \cellcolor{linen} \begin{tabular}{c} \algname{NSGD-MVR}  \\  Theorem~\ref{ch5:thm:nsgd-mvr-convergence-analysis_chapter_5} \end{tabular} & \cellcolor{linen}$\frac{\bar L\Delta}{\varepsilon^2}+\frac{\bar L\Delta}{\varepsilon^2}\left(\frac{\sigma_1}{\varepsilon}\right)^{\frac{p}{q(p-1)}}+ \left(\frac{\sigma_1}{\varepsilon}\right)^{\frac{p}{p-1}} $ & \cellcolor{linen} {\color{PineGreen}$(1,2]$} & \cellcolor{linen} \thesischfivexmark \\ 
 &\cellcolor{linen} \begin{tabular}{c} \algname{D-Clip-NSGD-MVR}  \\  Theorem~\ref{ch5:thm:clipped-nsgd-mvr-convergence-analysis_chapter_5} \end{tabular} & \cellcolor{linen}$\left(\frac{\sqrt{\bar L\Delta}+ \sigma_1}{\varepsilon}\right)^{2 +\max\left\{\frac{1}{p-1}, \frac{p}{q(p-1)}\right\}}$ & \cellcolor{linen} {\color{PineGreen}$(1,2]$} & \cellcolor{linen} \thesischfivecmark \\ \cline{2-5}
 & \cellcolor{linen} \begin{tabular}{c} Lower Bound \\ Theorem~\ref{ch5:thm:lower-bound-sofo-mean-squared-smoothness_chapter_5} \end{tabular} & \cellcolor{linen}$\frac{\bar L\Delta}{\varepsilon^2}+\frac{\bar L\Delta}{\varepsilon^2}\left(\frac{\sigma_1}{\varepsilon}\right)^{\frac{p}{q(p-1)}}+ \left(\frac{\sigma_1}{\varepsilon}\right)^{\frac{p}{p-1}} $ & \cellcolor{linen} {\color{PineGreen}$(1,2]$} & \cellcolor{linen} - \\ 
 \hline
 \multirow{6}{*}{\makecell{$(q,\delta)$-S}} 
  &
 \cellcolor{linen} \begin{tabular}{c} \algname{NSGD-MVR}  \\  Theorem~\ref{ch5:thm:nsgd-mvr-convergence-analysis-2_chapter_5} \end{tabular} & \cellcolor{linen}$\frac{(L_1+\delta)\Delta}{\varepsilon^2}+\frac{\delta\Delta}{\varepsilon^2}\left(\frac{\sigma_1}{\varepsilon}\right)^{\frac{p}{q(p-1)}}+ \left(\frac{\sigma_1}{\varepsilon}\right)^{\frac{p}{p-1}} $ & \cellcolor{linen} {\color{PineGreen}$(1,2]$} & \cellcolor{linen} \thesischfivexmark \\ 
 &\cellcolor{linen} \begin{tabular}{c} \algname{D-Clip-NSGD-MVR}  \\  Theorem~\ref{ch5:thm:clipped-nsgd-mvr-convergence-analysis-2_chapter_5} \end{tabular} & \cellcolor{linen}$\left(\frac{\sqrt{L_1\Delta}+ \sigma_1}{\varepsilon}\right)^{2 +\max\left\{\frac{1}{p-1}, \frac{p}{q(p-1)}\right\}}+ \left(\frac{\sqrt{\delta\Delta}}{\varepsilon}\right)^{2 + \frac{p}{q(p-1)}}$ & \cellcolor{linen} {\color{PineGreen}$(1,2]$} & \cellcolor{linen} \thesischfivecmark \\ \cline{2-5}
    & \cellcolor{linen} \begin{tabular}{c} Lower Bound \\ Theorem~\ref{ch5:thm:lower-bound-sofo-mean-squared-smoothness-2_chapter_5} \end{tabular} & \cellcolor{linen}$\min\left\{ \frac{L_1 \Delta}{\thesischfiveeps^2} + \frac{L_1 \Delta}{\thesischfiveeps^2} \left( \frac{\sigma_1}{\thesischfiveeps} \right)^{\frac{p}{p - 1}}, \frac{(L_1+\delta)\Delta}{\varepsilon^2}+\frac{\delta\Delta}{\varepsilon^2}\left(\frac{\sigma_1}{\varepsilon}\right)^{\frac{p}{q(p-1)}}+ \left(\frac{\sigma_1}{\varepsilon}\right)^{\frac{p}{p-1}}\right\}$ & \cellcolor{linen} {\color{PineGreen}$(1,2]$} & \cellcolor{linen} - \\ 
    \hline
    \multirow{10}{*}{\makecell{$q$-BCM}} &\begin{tabular}{c}  \algname{NSGDHess}  \\  \tabcite{sadiev2025second} \end{tabular} & $\frac{(L_1 +\sigma_2)\Delta}{\varepsilon^2} + \left(\frac{(L_1 + \sigma_2)\Delta }{\varepsilon^2}  + \frac{\sigma_1}{\varepsilon}\right)\left(\frac{\sigma_1}{\varepsilon}\right)^{\frac{1}{p-1}} $ & {\color{BrickRed} $\ p=q$}&  \ \thesischfivexmark \\ 
    &\begin{tabular}{c}  \algname{{Clip NSGDHess}}  \\  \tabcite{sadiev2025second} \end{tabular} & $\left(\frac{\sqrt{(L_1+\sigma_2)\Delta}+ \sigma_1}{\varepsilon}\right)^{2 +\frac{1}{p-1}}$ & {\color{BrickRed} $\ p=q$}&  \ \thesischfivecmark \\ 
 &
  \cellcolor{linen} \begin{tabular}{c} \algname{NSGDHess}  \\  Theorem~\ref{ch5:thm:nsgd-mvr-hess-convergence-analysis_chapter_5} \end{tabular} & \cellcolor{linen}$\frac{\sigma_2\Delta}{\varepsilon^2}\left(\frac{\sigma_1}{\varepsilon}\right)^{\frac{p}{q(p-1)}}+ \left(\frac{\sigma_1}{\varepsilon}\right)^{\frac{p}{p-1}} +\frac{\sqrt{L_2} \Delta \sigma_1^{\nicefrac{1}{4}}}{\thesischfiveeps^{\nicefrac{7}{4}}} \left( \frac{\sigma_1}{\thesischfiveeps} \right)^{\frac{1}{4 (p - 1)}}$ $^{{\color{blue}(3)}}$ & \cellcolor{linen} {\color{PineGreen}$(1,2]$} & \cellcolor{linen} \thesischfivexmark \\ 
 &\cellcolor{linen} \begin{tabular}{c} \algname{{D-Clip-NSGDHess}}  \\  Theorem~\ref{ch5:thm:clipped-nsgd-mvr-convergence-analysis-hess_chapter_5} \end{tabular} & \cellcolor{linen}$\left(\frac{\sqrt{L_1\Delta}+ \sigma_1}{\varepsilon}\right)^{2 +\max\left\{\frac{1}{p-1}, \frac{p}{q(p-1)}\right\}}+ \left(\frac{\sqrt{\sigma_2\Delta}}{\varepsilon}\right)^{2 + \frac{p}{q(p-1)}}$ $^{{\color{blue}(3)}}$ & \cellcolor{linen} {\color{PineGreen}$(1,2]$} & \cellcolor{linen} \thesischfivecmark \\ \cline{2-5}
 & \cellcolor{linen} \begin{tabular}{c} Lower Bound \\ Theorem~\ref{ch5:thm:lower-bound-sofo-bounded-central-moments_chapter_5} \end{tabular} & \cellcolor{linen}${\min\left\{\frac{\sigma_2\Delta }{\varepsilon^2}\left(\frac{\sigma_1}{\varepsilon}\right)^{\frac{p}{q(p-1)}}, \frac{L_1\Delta}{\varepsilon^{2}}\left(\frac{\sigma_1}{\varepsilon}\right)^{\frac{p}{p-1}} , \frac{\sqrt{L_2}\Delta}{\varepsilon^{\nicefrac{3}{2}}}\left(\frac{\sigma_1}{\varepsilon}\right)^{\frac{p}{p-1}} \right\}^{(3)}}$ & \cellcolor{linen} {\color{PineGreen}$(1,2]$} & \cellcolor{linen} - \\ \hline
\end{tabular}
}
{\scriptsize
\par\noindent{\color{blue}(1)} \textcite{pmlr-v195-liu23c} analyze stronger assumptions that imply bounded stochastic Hessians.
\par\noindent{\color{blue}(2)} \textcite{he2025complexitynormalizedstochasticfirstorder} establish a near-optimal dependence on $\varepsilon$ when $p=q$.
\par\noindent{\color{blue}(3)} Only the stochastic part of the complexity is displayed.
}
\end{minipage}
\end{table*}

We revisit the analysis of Normalized \algname{SGD} (\algname{NSGD}) with Momentum Variance Reduction (\algname{MVR}) in the heavy-tailed regime. Our contributions are:

$\bullet$ \textbf{Tight lower bounds.} {We obtain novel and tight lower bounds on the sample complexity of\stepEightBdeletion{} first-order method\stepEightBchange{s}} in the $p$--BCM noise model, under \emph{relaxed} mean-squared smoothness and $\delta$-similarity assumptions, allowing any exponent $q \in \thesischfiveintff{1}{2}$ instead of $q = 2$. Our bound, $\smash{\Omega(\thesischfiveeps^{-\nicefrac{(p (2 q + 1) - 2 q)}{q (p - 1)}})}$, recovers the known $\smash{\Omega(\thesischfiveeps^{-3})}$ when $p = q = 2$ and achieves an improvement over the minimax complexity for first-order methods~\parencite{10.5555/3495724.3497014,sun2024gradientnormalizationprovablybenefits,hübler2025gradientclippingnormalizationheavy}, thereby demonstrating the effectiveness of variance reduction in the heavy-tailed regime. Moreover, our result \emph{continuously interpolates} between variance-reduced and non-variance-reduced regimes, as in the limit $q \to 1$ we recover the known $\Omega(\thesischfiveeps^{-\nicefrac{(3 p - 2)}{(p - 1)}})$ complexity. 
\emph{Our proofs are based on a generalization of~\textcite[Lemma~10]{arjevani2022lower} to the $p$--BCM noise assumption, which resolves an open question from~\textcite{pmlr-v195-liu23c}}.

$\bullet$ \textbf{Optimal first-order methods.} Revisiting \algname{NSGD-MVR}, we show that, with a suitable choice of parameters, the method matches the lower bounds (up to constant factors) with respect to all parameters, thereby establishing the optimality of our analysis and the robustness of gradient normalization in the context of heavy-tailed noise. %
We further provide convergence rates with \emph{unknown} tail indices $p$ and $q$, recovering the best known rate of $\cO(\thesischfiveeps^{-\nicefrac{2 p}{(p - 1)}})$~\parencite{liu2025nonconvex}.

$\bullet$ \textbf{High-probability upper bounds.} We propose a clipped variant of \algname{NSGD-MVR} leading to a \emph{new} algorithm: Double-Clipped \algname{NSGD-MVR} (\algname{D-Clip-NSGD-MVR}), which provably achieves high-probability convergence rates under weaker assumptions than in previous works \parencite{pmlr-v195-liu23c}.

$\bullet$ \textbf{Extension to second-order methods.} Lastly, for second-order methods using stochastic Hessians with bounded $q$-th central moment, for some exponent $q \in \thesischfiveintff{1}{2}$, we establish sharper lower bounds than previous works~\parencite{pmlr-v125-arjevani20a,sadiev2025second}, extend previous analysis to the more general setting where $p \neq q$ (unlike~\textcite{sadiev2025second}), and derive upper bounds with stronger convergence exponents, thereby offering a more complete characterization of the limits of stochastic second-order optimization under heavy-tailed noise. Overall, our results reveal a striking parallel in the rates of variance-reduced methods and second-order optimization.

\subsection{Related Work}

{\noindent\textbf{Lower bounds in the heavy-tailed regime.}\par
\noindent In nonconvex optimization,~\textcite{10.5555/3495724.3497014,liu2025nonconvex}
establish, under $p$--BCM noise and standard $L_1$--smoothness, the lower bound
\[
    \Omega\!\left(
        \frac{L_1\Delta}{\thesischfiveeps^2}
        +
        \frac{L_1\Delta}{\thesischfiveeps^2}
        \left(\frac{\sigma_1}{\thesischfiveeps}\right)^{\frac{p}{p-1}}
    \right).
\]
For $p=2$, this recovers the well-known $\Omega(\thesischfiveeps^{-4})$ stochastic lower bound.}
Another line of work establishes algorithm-specific lower bounds for \algname{NSGD}~\parencite{10.5555/3666122.3669370} and, more recently, \algname{SGD}~\parencite{fatkhullin2025sgdhandleheavytailednoise} under their respective oracle, regularity, output, and stepsize classes; these results characterize the limitations of the analyzed algorithms and do not establish a criterion-independent failure of \algname{SGD} under heavy-tailed noise. Beyond first-order methods,~\textcite{sadiev2025second} obtain near-optimal lower bounds for all second-order methods, extending the work of~\textcite{pmlr-v125-arjevani20a} to $p$--BCM noise.

\textbf{Convergence under heavy-tailed noise:} A substantial body of work investigates upper bounds for stochastic optimization with heavy-tailed noise, both in expectation and high probability. For smooth and nonconvex objectives, clipped or normalized variants of \algname{SGD} achieve the rate $\cO(\thesischfiveeps^{-\nicefrac{(3 p - 2)}{(p - 1)}})$ up to logarithmic factors~\parencite{10.5555/3495724.3497014,NEURIPS2021_26901deb,pmlr-v195-liu23c,pmlr-v202-sadiev23a}. Among these,~\textcite{10.5555/3495724.3497014} obtained the sharpest in-expectation guarantees, later shown to be tight by \textcite{hübler2025gradientclippingnormalizationheavy} via \algname{minibatch-NSGD}, while~\textcite{NEURIPS2023_4c454d34} established improved high-probability guarantees, avoiding extra $\cO(\log K)$ factors. Later works extended known high-probability analyses to settings with higher-order smoothness~\parencite{sadiev2025second} and to more general nonlinear \algname{SGD}-type methods~\parencite{DBLP:journals/corr/abs-2410-13954,pmlr-v258-armacki25a}. In parallel,~\textcite{sun2024gradientnormalizationprovablybenefits} studied \algname{clip-SGD} under stronger smoothness assumptions, obtaining a rate of $\cO(\thesischfiveeps^{-\nicefrac{(2 p - 1)}{(p - 1)}})$, while~\textcite{liu2025nonconvex} derived in-expectation bounds for \algname{NSGD-Mom} under a general $(\sigma_0,\sigma_1)$--affine $p$--BCM model, albeit without matching lower bounds.

\textbf{Gradient clipping and normalization:} Gradient clipping and normalization are two closely related techniques that have become central in modern optimization. Gradient clipping, originally popularized to stabilize training across various machine learning applications~\parencite{10.5555/3042817.3043083,schulman2017proximalpolicyoptimizationalgorithms}, has been analyzed extensively, providing robustness under relaxed moment assumptions~\parencite{PolyTsy79,jakovetic2023nonlinear}, high-order smoothness~\parencite{sadiev2025second}, and enabling high-probability convergence guarantees with only logarithmic dependence on the failure probability. These results hold in both convex~\parencite{Nazin2019,10.5555/3495724.3496985,JMLR:v22:20-821,10.5555/3692070.3692710,liu2023stochasticnonsmoothconvexoptimization,Gorbunov2024,pmlr-v238-puchkin24a,pmlr-v258-armacki25a,armacki2025optimalhighprobabilityconvergencenonlinear} and nonconvex settings~\parencite{10.5555/3495724.3497014,NEURIPS2021_26901deb,pmlr-v202-sadiev23a,NEURIPS2023_4c454d34,pmlr-v195-liu23c,sadiev2025second}. Importantly, high-probability guarantees are valuable both theoretically and practically, as they capture the behavior of individual runs rather than merely characterizing the average-case performance.%

Normalized gradient methods date back to Nesterov's pioneering work~\parencite{nesterov1984minimization} and were later extended to the smooth and stochastic regimes~\parencite{kiwiel2001,hazan2015beyond,NIPS2017_ce5140df,pmlr-v89-nacson19b}. In deep learning, normalization addresses exploding and vanishing gradients~\parencite{you2017largebatchtrainingconvolutional,you2020Large}, though rigorous nonconvex guarantees were first obtained by~\textcite{pmlr-v119-cutkosky20b}, showing how Polyak's momentum ensures convergence without large batches. Later works explored \algname{NSGD}'s strengths and limits, including saddle-point escape~\parencite{levy2016powernormalizationfasterevasion}, lower bounds~\parencite{10.5555/3666122.3669370} and parameter-agnostic convergence~\parencite{pmlr-v238-hubler24a}. Recently, extensions to heavy-tailed noise have been investigated~\parencite{NEURIPS2021_26901deb,pmlr-v195-liu23c,sun2024gradientnormalizationprovablybenefits}, though typically under stronger smoothness assumptions. 

{\centering\emph{Despite strong empirical performance and complementarity with clipping, the theory of gradient normalization under heavy-tailed noise remains incomplete, motivating our contributions}.\par}

\textbf{Variance reduction for stochastic optimization:} Variance reduction originated as a tool to accelerate convergence in convex finite-sum optimization, with seminal works such as~\textcite{10.5555/2999325.2999432,10.5555/2999611.2999647,10.5555/2567709.2502598,10.5555/3042817.3043024,10.5555/2968826.2969010} introducing various algorithms, e.g., \algname{SAG}, \algname{SVRG}, \algname{SAGA}. Building on these foundations and the key advances of~\textcite{10.1145/3055399.3055448}, a sequence of works 
\textcite{lan2019unified,zhou2019direct,song2020variance,kovalev2020don} developed algorithms attaining near-optimal or optimal rates under various regimes. In nonconvex optimization, variance reduction can also improve convergence: in the general stochastic setting~\eqref{ch5:eq:optim-problem_chapter_5}, several works~\parencite{SPIDER,cutkosky2019momentum,tran2019hybrid,liu2020optimal,PAGE} established an $\cO(\thesischfiveeps^{-3})$ convergence rate in expectation, which improves upon the classical $\Theta(\thesischfiveeps^{-4})$ of \algname{SGD} and matches the minimax lower bound $\Omega(\thesischfiveeps^{-3})$ under mean-squared smoothness~\parencite{arjevani2022lower}. More recently, there has been growing interest in extending variance reduction techniques to the heavy-tailed regime~\parencite{pmlr-v195-liu23c,sun2024gradientnormalizationprovablybenefits,he2025complexitynormalizedstochasticfirstorder}, where normalization and clipping have emerged as effective mechanisms to guarantee robustness under weaker moment assumptions. 

{\centering\emph{This chapter characterizes the fundamental limits of first-order \say{variance-reduced} methods in the heavy-tailed regime and provides optimal methods that attain these limits in the stated regimes}.\par}

\section{Notation and Assumptions}

{We review the basic notation and assumptions used in this chapter; further details appear in Section~\ref{ch5:appdx-sec:notation_chapter_5}.}

\textbf{Notation.} For integer $n > 0$, $[n] \eqdef \thesischfiveens{1, 2, \ldots, n}$. We let $d \ge 1$ be the dimension, $\thesischfiveps{\cdot}{\cdot}$ the standard dot product on $\R^d$, $\norm{\cdot}$ the $\ell^2$--norm and $\normop{\cdot}$ the canonical spectral/operator norm. Here, {$\nabla f(\cdot;\cdot)$ and $\nabla^2 f(\cdot;\cdot)$} denote the stochastic gradient and Hessian oracles; $a \wedge b \eqdef \min\thesischfiveens{a, b}$ and $a \vee b \eqdef \max\thesischfiveens{a, b}$. We use the standard $\cO(\cdot)$ and $\Omega(\cdot)$ for complexity notation. To avoid confusion, we use subscript $1$ for gradient parameters (e.g., $L_1$, $\sigma_1$) and $2$ for Hessian ones.

The analysis in this chapter uses the following assumptions.
\begin{assumption}[Lower Boundedness]\label{ch5:ass:lower-boundedness_chapter_5}
    The objective $f$ is lower bounded: $f^{\inf} \eqdef \inf_{x \in \R^d} f(x) > -\infty$. %
\end{assumption}
We then let $\Delta \eqdef f(x^0) - f^{\inf}$ be the initial suboptimality where $x^0$ is the starting point.%

\begin{assumption}[$L_1$--Lipschitz Gradients]\label{ch5:ass:L-lipschitz-gradients_chapter_5}
    The objective $f$ is differentiable over $\R^d$ and  its gradient is $L_1$--Lipschitz for some $L_1 \ge 0$, i.e., for all $x, y \in \R^d$,
        \[ \norm{\nabla f(x) - \nabla f(y)} \le L_1 \norm{x - y}. \]
\end{assumption}

\begin{assumption}[$p$--BCM for Gradients]\label{ch5:ass:p-bounded-central-moment-gradient_chapter_5}
    We have access to unbiased stochastic gradients $\nabla f(x;\xi)$, which have $p$-bounded central moment for some $p \in \thesischfiveintof{1}{2}$, 
    i.e., there exists {$\sigma_1 \ge 0$} such that for all $x \in \R^d$,
    \[
        \thesischfiveExpSub{\xi \sim \thesischfivecD}{\nabla f(x;\xi)} = \nabla f(x)
        \quad\text{and}\quad
        \thesischfiveExpSub{\xi \sim \thesischfivecD}{\norm{\nabla f(x;\xi)-\nabla f(x)}^p} \le \sigma_1^p.
    \]
\end{assumption}

\begin{assumption}[$L_2$--Lipschitz Hessians]
\label{ch5:ass:L-lipschitz-hessians_chapter_5}
    The objective $f$ is twice continuously differentiable over $\R^d$ 
    and we have access to stochastic Hessian-vector products $\nabla^2 f(x;\xi) \cdot v$ for any $v \in \R^d$. Moreover, the Hessian of $f$ is $L_2$--Lipschitz, i.e., for all $x, y \in \R^d$,
        \[ \normop{\nabla^2 f(x) - \nabla^2 f(y)} \le L_2 \norm{x - y}. \] 
\end{assumption}

\begin{assumption}[$q$--BCM for Hessians]\label{ch5:ass:q-bounded-central-moment-hessian_chapter_5}
    The stochastic Hessians $\nabla^2 f(x;\xi)$ are unbiased and have $q$-bounded central moment for some $q \in \thesischfiveintff{1}{2}$,
    i.e., there exists {$\sigma_2 \ge 0$} such that for all $x \in \R^d$
    \[
        \thesischfiveExpSub{\xi \sim \thesischfivecD}{\nabla^2 f(x;\xi)} = \nabla^2 f(x)
        \quad\text{and}\quad
        \thesischfiveExpSub{\xi \sim \thesischfivecD}{\normop{\nabla^2 f(x;\xi)-\nabla^2 f(x)}^q} \le \sigma_2^q.
    \]
\end{assumption}

\begin{assumption}[{$q$-Weak Average Smoothness~\textcite{he2025complexitynormalizedstochasticfirstorder}}]\label{ch5:ass:mean-squared-smoothness_chapter_5}
    For some $q \in [1,2]$, there exists a finite constant $\bar{L} \ge 0$ such that, for all $x, y \in \R^d$ we have
        \[ \thesischfiveExpSub{\xi \sim \thesischfivecD}{\norm{\nabla f(x;\xi) - \nabla f(y;\xi)}^q} \le \bar{L}^q \norm{x - y}^q. \]
\end{assumption}

The next assumption has been introduced in~\textcite{pmlr-v125-arjevani20a} in the special case $q = 2$.
\begin{assumption}[$(q, \delta)$--Similarity]\label{ch5:ass:mean-squared-smoothness-2_chapter_5}
    For some exponent $q \in [1,2]$, there exists a finite constant $\delta \ge 0$ such that, for all $x, y \in \R^d$ we have
    \[
        \thesischfiveExpSub{\xi \sim \thesischfivecD}{\norm{[\nabla f(x;\xi)-\nabla f(y;\xi)]-[\nabla f(x)-\nabla f(y)]}^q}
        \le \delta^q \norm{x-y}^q.
    \]
\end{assumption}
\vspace*{-0.25cm}
{
Assumption~\ref{ch5:ass:mean-squared-smoothness-2_chapter_5} together with Assumption~\ref{ch5:ass:L-lipschitz-gradients_chapter_5} implies Assumption~\ref{ch5:ass:mean-squared-smoothness_chapter_5} with
\[
    \bar L^q \le 2^{q-1}\bigl(L_1^q+\delta^q\bigr).
\]
Conversely, Assumption~\ref{ch5:ass:mean-squared-smoothness_chapter_5} implies $L_1\le\bar L$ by Jensen's inequality and implies Assumption~\ref{ch5:ass:mean-squared-smoothness-2_chapter_5} with $\delta\le 2\bar L$.  Thus, for $q\in[1,2]$ the two descriptions are equivalent up to universal factors.  In the special case $q=2$, the variance identity sharpens the latter relation to $\delta\le\bar L$.
}
\vspace*{-0.075cm}
\section{Lower Complexity Bounds}
\vspace*{-0.075cm}

{
To establish our lower bounds, we analyze the distributional complexity of finding an $\thesischfiveeps$-stationary point. We focus on the broad class of \textit{zero-respecting} algorithms, denoted by $\mathcal{A}_{\texttt{zr}}$, which only query coordinates previously discovered to have non-zero gradients~\parencite{carmon2020lower}. To keep the query-batch arity distinct from the iteration horizon $K$, we denote it by $b$. At each oracle round, an algorithm $\texttt{A} \in \mathcal{A}_{\texttt{zr}}$ makes $b\ge1$ simultaneous multi-point queries evaluated using the same fresh seed $\xi$. If this seed reveals no new coordinate, none of the $b$ simultaneous queries can reveal one. Consequently, $b$ does not enter the final lower bounds. The MVR algorithms use $b=2$, whereas the Hessian-corrected algorithm uses $b=1$.
}

Formally, following~\textcite{arjevani2022lower}, we define the minimax sample complexity as:
\[
\begin{aligned}
{\mathfrak{m}^{\texttt{zr}}_{\thesischfiveeps}(b, \Delta\ldots)}
&\eqdef \sup_{f \in \mathcal{F}(\Delta\ldots)}
\sup_{(\Prob_{\xi},\texttt{O}_f) \in \mathcal{O}(f\ldots)}
{\inf_{\texttt{A} \in \mathcal{A}_{\texttt{zr}}(b)}} \\
&\quad \inf \thesischfiveenstq{N \ge 1}{\thesischfiveExpSub{\texttt{O}_f,\texttt{A}}{\norm{\nabla f\left(\mathbf{x}_{{\normalfont\texttt{A}}[\normalfont\texttt{O}_f]}^{(N)}\right)}} \le \thesischfiveeps}.
\end{aligned}
\]
{where the expectation is taken over the randomness in the oracle $\texttt{O}_f$ and in the algorithm $\texttt{A} \in \mathcal{A}_{\texttt{zr}}$, if any.}

To derive these bounds, we rely on a technique developed by~\textcite{carmon2020lower,arjevani2022lower}, introducing a \say{worst-case} nonconvex function $f_N \colon \R^N \to \R$ that is hard for any zero-respecting algorithm to optimize. Our main theoretical novelty here is~Lemma~\ref{ch5:appdx-lem:lower-bound-global-stichastic-model-mean-squared-smoothness-p-BCM_chapter_5}, an important result lying at the core of all our lower bounds. By generalizing~\textcite[Lemma~10]{arjevani2022lower} to the $p$--BCM assumption, this lemma yields a tight improvement from $\Omega(\nicefrac{\sigma_1^2}{\thesischfiveeps^2})$ to $\Omega((\nicefrac{\sigma_1}{\thesischfiveeps})^{\nicefrac{p}{(p - 1)}})$. This handles a critical edge case in~Theorem~\ref{ch5:thm:lower-bound-sofo-mean-squared-smoothness_chapter_5}, Theorem~\ref{ch5:thm:lower-bound-sofo-mean-squared-smoothness-2_chapter_5}, and Theorem~\ref{ch5:thm:lower-bound-sofo-bounded-central-moments_chapter_5} where the parameters would otherwise force the dimension $N$ of the \say{hard instance} $f_N$ to be too small. 

While the key complexity definitions are provided above, we defer the formal definitions of the function and oracle classes ($\mathcal{F}, \mathcal{O}$), the exact construction of $f_N$, and further motivations behind our generalization to~Remark~\ref{ch5:appdx-rem:lower-bound-global-stochastic-model-motivations_chapter_5} and~Definition~\ref{ch5:appdx-sec:lower-bound-def_chapter_5}. The full proofs for the following three lower bounds can be found in~\ref{ch5:appdx-subsec:proof-lower-bound-sofo-mean-squared-smoothness_chapter_5}, \ref{ch5:appdx-subsec:proof-lower-bound-sofo-mean-squared-smoothness-2_chapter_5}, and \ref{ch5:appdx-subsec:proof-lower-bound-sofo-bounded-central-moments_chapter_5}.

\vspace*{-0.08cm}
\subsection{\texorpdfstring{{Under $q$-Weak Average Smoothness}}{Under q-Weak Average Smoothness}}
\vspace*{-0.08cm}
We first establish a lower bound on the oracle complexity of any zero-respecting algorithm under~Assumption~\ref{ch5:ass:mean-squared-smoothness_chapter_5}, an extension of the mean-squared smoothness (MSS) assumption. This result generalizes one of the main contributions of~\textcite{arjevani2022lower}, who obtained optimal lower bound in the $p = q = 2$ setting.
\vspace*{-0.4cm}
\begin{theorem}\label{ch5:thm:lower-bound-sofo-mean-squared-smoothness_chapter_5}
    {Given $\Delta,\bar L>0$, $\sigma_1\ge0$, and $0<\thesischfiveeps\le c_1\sqrt{\bar L\Delta}$ for a universal constant $c_1>0$, the minimax sample complexity over functions and oracles satisfying Assumptions~\ref{ch5:ass:lower-boundedness_chapter_5}, \ref{ch5:ass:mean-squared-smoothness_chapter_5}, and \ref{ch5:ass:p-bounded-central-moment-gradient_chapter_5} obeys}
    \[
        {\mathfrak{m}^{\normalfont\texttt{zr}}_{\thesischfiveeps}\left(b,\Delta,\bar{L}^q,\sigma_1^p\right)}
        \ge \Omega(1) \cdot \left(\left(\frac{\sigma_1}{\thesischfiveeps}\right)^{\frac{p}{p-1}} + \frac{\bar{L}\Delta}{\thesischfiveeps^2} + \frac{\bar{L}\Delta}{\thesischfiveeps^2}\left(\frac{\sigma_1}{\thesischfiveeps}\right)^{\frac{p}{q(p-1)}}\right).
    \]
\end{theorem}
In the special case $p = q = 2$, our result recovers the optimal complexity lower bound $\Omega(\nicefrac{\sigma_1^2}{\thesischfiveeps^2} + \nicefrac{\bar{L} \Delta}{\thesischfiveeps^2} + \nicefrac{\bar{L} \Delta \sigma_1}{\thesischfiveeps^3})$ derived by~\textcite{arjevani2022lower}. We further establish the tightness of this bound in~Theorem~\ref{ch5:thm:nsgd-mvr-convergence-analysis_chapter_5}, where we prove a matching upper bound.
\vspace*{-0.075cm}
\subsection{\texorpdfstring{{Under $(q,\delta)$--Similarity}}{Under (q, delta)-Similarity}}
\vspace*{-0.075cm}
The previous bound captures limits of first-order methods under~Assumption~\ref{ch5:ass:mean-squared-smoothness_chapter_5}, but this assumption couples smoothness and noise together into $\bar{L}$. Assuming separately $(q, \delta)$–similarity with $L_1$-smoothness decouples these effects, yielding sharper lower bounds.
\begin{theorem}\label{ch5:thm:lower-bound-sofo-mean-squared-smoothness-2_chapter_5}
    {Given $\Delta,L_1,\delta>0$, $\sigma_1\ge0$, and $0<\thesischfiveeps\le c_1\sqrt{L_1\Delta}$ for a universal constant $c_1>0$, the minimax sample complexity over functions and oracles satisfying Assumptions~\ref{ch5:ass:lower-boundedness_chapter_5}, \ref{ch5:ass:L-lipschitz-gradients_chapter_5}, \ref{ch5:ass:p-bounded-central-moment-gradient_chapter_5}, and \ref{ch5:ass:mean-squared-smoothness-2_chapter_5} obeys}
    \begin{align*}
        {\mathfrak{m}^{\normalfont\texttt{zr}}_{\thesischfiveeps}\left(b,\Delta,L_1,\delta^q,\sigma_1^p\right)}
        &\ge \Omega(1) \cdot \min\left\{\frac{L_1\Delta}{\thesischfiveeps^2} + \frac{L_1\Delta}{\thesischfiveeps^2}\left(\frac{\sigma_1}{\thesischfiveeps}\right)^{\frac{p}{p-1}},\right.\\[-0.2\baselineskip]
        &\hspace{24mm}\left.\left(\frac{\sigma_1}{\thesischfiveeps}\right)^{\frac{p}{p-1}} + \frac{(L_1+\delta)\Delta}{\thesischfiveeps^2} + \frac{\delta\Delta}{\thesischfiveeps^2}\left(\frac{\sigma_1}{\thesischfiveeps}\right)^{\frac{p}{q(p-1)}}\right\}.
    \end{align*}
\end{theorem}
Expressing the rate in terms of $L_1$ and $\delta$ is more precise than using $\bar{L}$ alone; in particular, when $\delta \ll \bar{L}$, we show in~Theorem~\ref{ch5:thm:nsgd-mvr-convergence-analysis-2_chapter_5} that \algname{NSGD-MVR} can achieve improved complexity, even matching the second term of the above $\min\thesischfiveens{\ldots}$.

\subsection{Under Bounded Central Moments}

Compared to the previous results, which focused on first-order methods under $q$-weak average smoothness and $(q, \delta)$–similarity, we consider here the setting where the oracle's noise is controlled only via bounded central moments, allowing distinct exponents $p$ and $q$ for the stochastic gradients and Hessians.
\begin{theorem}\label{ch5:thm:lower-bound-sofo-bounded-central-moments_chapter_5}
    {Let $\Delta,L_1,L_2>0$ and $\sigma_1,\sigma_2\ge0$. For a universal constant $c_1>0$, assume
    \[
        0<\thesischfiveeps\le
        c_1\min\!\left\{\sqrt{L_1\Delta},L_2^{1/3}\Delta^{2/3}\right\}.
    \]
    Then the minimax sample complexity over functions and oracles satisfying Assumptions~\ref{ch5:ass:lower-boundedness_chapter_5}, \ref{ch5:ass:L-lipschitz-gradients_chapter_5}, \ref{ch5:ass:p-bounded-central-moment-gradient_chapter_5}, \ref{ch5:ass:L-lipschitz-hessians_chapter_5}, and \ref{ch5:ass:q-bounded-central-moment-hessian_chapter_5} obeys}
    \begin{align*}
        \MoveEqLeft {\mathfrak{m}^{\normalfont\texttt{zr}}_{\thesischfiveeps}
        \left(b,\Delta,L_1,L_2,\sigma_1^p,\sigma_2^q\right)}\\
        &\ge \Omega(1)\min\Bigg\{
        \frac{L_1\Delta}{\thesischfiveeps^2}
        \left[1+\left(\frac{\sigma_1}{\thesischfiveeps}\right)^{p/(p-1)}\right],\\
        &\qquad
        \frac{L_2^{1/2}\Delta}{\thesischfiveeps^{3/2}}
        \left[1+\left(\frac{\sigma_1}{\thesischfiveeps}\right)^{p/(p-1)}\right],\\
        &\qquad
        \min\left\{\frac{L_1\Delta}{\thesischfiveeps^2},
        \frac{L_2^{1/2}\Delta}{\thesischfiveeps^{3/2}}\right\}
        +\left(\frac{\sigma_1}{\thesischfiveeps}\right)^{p/(p-1)}
        +\frac{\sigma_2\Delta}{\thesischfiveeps^2}
        \left[1+\left(\frac{\sigma_1}{\thesischfiveeps}\right)^{p/(q(p-1))}\right]
        \Bigg\}.
    \end{align*}
\end{theorem}
\begin{remark}
    In the noiseless setting, i.e., $\sigma_1 = \sigma_2 = 0$, the complexity of~Theorem~\ref{ch5:thm:lower-bound-sofo-bounded-central-moments_chapter_5} reduces to
        \[ \min\thesischfiveens{\frac{L_1 \Delta}{\thesischfiveeps^2}, \frac{L_2^{\nicefrac{1}{2}} \Delta}{\thesischfiveeps^{\nicefrac{3}{2}}}}, \thesischfivenumberthis\label{ch5:b26a6d07-2917-471c-a3c4-5c9b1b33f5cc_chapter_5} \]
    which is matched by the combination of gradient descent (\algname{GD}) and cubic regularized Newton method~\parencite{Nesterov2006}, hence is optimal. This is slightly better than the $\smash{\cO(L_1^{\nicefrac{1}{2}} L_2^{\nicefrac{1}{4}} \Delta \thesischfiveeps^{-\nicefrac{7}{4}})}$ bound %
    achieved in~\textcite{10.5555/3305381.3305449}, since
        \[ \min\thesischfiveens{\frac{L_1 \Delta}{\thesischfiveeps^2}, \frac{L_2^{\nicefrac{1}{2}} \Delta}{\thesischfiveeps^{\nicefrac{3}{2}}}} \le \sqrt{\frac{L_1 \Delta}{\thesischfiveeps^2} \cdot \frac{L_2^{\nicefrac{1}{2}} \Delta}{\thesischfiveeps^{\nicefrac{3}{2}}}} = \frac{L_1^{\nicefrac{1}{2}} L_2^{\nicefrac{1}{4}} \Delta}{\thesischfiveeps^{\nicefrac{7}{4}}}, \]
    but achieving~\eqref{ch5:b26a6d07-2917-471c-a3c4-5c9b1b33f5cc_chapter_5} requires full Hessian access in order to get rid of the dependency in $L_1$~\parencite{pmlr-v125-arjevani20a,Carmon2021}.
\end{remark}

\section{Optimal Method under First-Order BCM}\label{ch5:sec:optimal-method-nsgd-mvr_chapter_5}

To establish the optimal rate, we consider the well-known Normalized \algname{SGD} algorithm with the momentum variance reduction (\texttt{MVR}) technique~\parencite{cutkosky2019momentum}, and present its pseudo-code in Algorithm~\ref{ch5:algo:nsgd-mvr_chapter_5}.

{Throughout this chapter, we write normalized vectors as
$v/\norm{v}$.}\footnote{{Whenever $v=0$, the normalized quotient
$v/\norm{v}$ is defined to be the zero vector. This convention makes every
normalized update well defined.}}

    \begin{algorithm}[htbp]
\caption{\algname{NSGD-MVR} (Normalized SGD with \texttt{MVR})}
\label{ch5:algo:nsgd-mvr_chapter_5}
\begin{algorithmic}[1]
\Require $x^0 \in \R^d$, the starting point;              $K > 0$, the number of iterations;              $g^0 \in \R^d$, an initial vector;              $\gamma > 0$, the stepsize;              $\eta \in \thesischfiveintof{0}{1}$, the momentum parameter for \texttt{MVR};
\State {$x^1 = x^0 - \gamma \frac{g^0}{\norm{g^0}}$}
\For{$k = 1, 2, \ldots, K - 1$}
\State \textit{Apply \texttt{MVR}.}
\State {$g^{k} =(1 - \eta)\!\left(g^{k-1}+\nabla f(x^{k};\xi^{k})-\nabla f(x^{k-1};\xi^{k})\right)+\eta\nabla f(x^{k};\xi^{k})$}
\State \textit{Do one descent step.}
\State {$x^{k+1} = x^{k} - \gamma \frac{g^{k}}{\norm{g^{k}}}$}
\EndFor
\Ensure {$x_{\mathrm{out}}\sim\mathcal U\{x^0,\ldots,x^{K-1}\}$}
\end{algorithmic}
\end{algorithm}%

\subsection{Convergence Analysis}

\subsubsection{Case of Known Tail Indices $p$ and $q$}

\begin{theorem}\label{ch5:thm:nsgd-mvr-convergence-analysis_chapter_5}
    Under~Assumption~\ref{ch5:ass:lower-boundedness_chapter_5}, Assumption~\ref{ch5:ass:p-bounded-central-moment-gradient_chapter_5}, and Assumption~\ref{ch5:ass:mean-squared-smoothness_chapter_5}, {choose, for a sufficiently large universal constant $C>0$,}
    {
    \[
        K=\left\lceil C\left(1+\left(\frac{\sigma_1}{\thesischfiveeps}\right)^{\frac{p}{p-1}}
        +\frac{\bar L\Delta}{\thesischfiveeps^2}
        +\frac{\bar L\Delta}{\thesischfiveeps^2}
        \left(\frac{\sigma_1}{\thesischfiveeps}\right)^{\frac{p}{q(p-1)}}\right)\right\rceil .
    \]}
    {Let the initial gradient estimate $g^{0}$ be given by}
        \[ g^{0} = \frac{1}{B_{\textnormal{init}}} \sum_{j = 0}^{B_{\textnormal{init}} - 1} \nabla f\left( x^{0}; \xi_{j}^{0} \right), \]
    where {$B_{\textnormal{init}} = \left\lceil\max\thesischfiveens{1, \left( \frac{2\sigma_1}{\thesischfiveeps} \right)^{\frac{p}{p - 1}}}\right\rceil$}, let the stepsize $\gamma = \sqrt{\frac{\Delta \eta^{\nicefrac{1}{q}}}{\bar{L} K}}$\footnote{If $\bar L=0$, the assumptions imply $\nabla f\equiv0$, so the conclusion is immediate; the displayed stepsize is used for $\bar L>0$.}, and let the momentum parameter $\eta = \min\thesischfiveens{1, \eta_{\textnormal{eff}}}$, where
        \[ \eta_{\textnormal{eff}} = \max\thesischfiveens{\left( \frac{\thesischfiveeps}{\sigma_1 K} \right)^{\frac{p}{2p - 1}}, \left( \frac{\bar{L} \Delta}{\sigma_1^2 K} \right)^{\frac{pq}{p (2 q + 1) - 2 q}}}. \thesischfivenumberthis\label{ch5:2c0dbb02-323f-4955-a4f2-8f3b73a0f89d_chapter_5} \]
    {When $\sigma_1=0$, we use the deterministic convention $\eta=1$ and $B_{\textnormal{init}}=1$.}
    {The output $x_{\mathrm{out}}$ is sampled uniformly from $\{x^0,\ldots,x^{K-1}\}$. Then Algorithm~\ref{ch5:algo:nsgd-mvr_chapter_5} guarantees $\Exp{\norm{\nabla f(x_{\mathrm{out}})}}\le\thesischfiveeps$ with total sample complexity}
        \[ \cO\left( {1+} \left( \frac{\sigma_1}{\thesischfiveeps} \right)^{\frac{p}{p - 1}} + \frac{\bar{L} \Delta}{\thesischfiveeps^2} + \frac{\bar{L} \Delta}{\thesischfiveeps^2} \left( \frac{\sigma_1}{\thesischfiveeps} \right)^{\frac{p}{q (p - 1)}} \right). \]
\end{theorem}
The proof is deferred to~\ref{ch5:appdx:proof-nsgd-mvr_chapter_5}.

Using~Assumption~\ref{ch5:ass:mean-squared-smoothness-2_chapter_5} combined with~Assumption~\ref{ch5:ass:L-lipschitz-gradients_chapter_5}, we obtain the following refined result.
\begin{theorem}\label{ch5:thm:nsgd-mvr-convergence-analysis-2_chapter_5}
    Under~Assumption~\ref{ch5:ass:lower-boundedness_chapter_5}, Assumption~\ref{ch5:ass:L-lipschitz-gradients_chapter_5}, Assumption~\ref{ch5:ass:p-bounded-central-moment-gradient_chapter_5}, and Assumption~\ref{ch5:ass:mean-squared-smoothness-2_chapter_5}, {choose, for a sufficiently large universal constant $C>0$,}
    {
    \[
        K=\left\lceil C\left(1+\left(\frac{\sigma_1}{\thesischfiveeps}\right)^{\frac{p}{p-1}}
        +\frac{(L_1+\delta)\Delta}{\thesischfiveeps^2}
        +\frac{\delta\Delta}{\thesischfiveeps^2}
        \left(\frac{\sigma_1}{\thesischfiveeps}\right)^{\frac{p}{q(p-1)}}\right)\right\rceil .
    \]}
    {Let the initial gradient estimate $g^{0}$ be given by}
        \[ g^{0} = \frac{1}{B_{\textnormal{init}}} \sum_{j = 0}^{B_{\textnormal{init}} - 1} \nabla f\left( x^{0}; \xi_{j}^{0} \right), \]
    with {$B_{\textnormal{init}} = \left\lceil\max\thesischfiveens{1, \left( \frac{2\sigma_1}{\thesischfiveeps} \right)^{\frac{p}{p - 1}}}\right\rceil$}, stepsize $\gamma = \min\thesischfiveens{\sqrt{\frac{\Delta}{L_1 K}}, \sqrt{\frac{\Delta \eta^{\nicefrac{1}{q}}}{\delta K}}}$\footnote{A stepsize term with zero denominator is interpreted as $+\infty$ and omitted from the minimum. If $L_1=\delta=0$, then $\nabla f\equiv0$ and the conclusion is immediate.}, and momentum parameter $\eta = \min\thesischfiveens{1, \eta_{\textnormal{eff}}}$, where
        \[ \eta_{\textnormal{eff}} = \max\thesischfiveens{\left( \frac{\thesischfiveeps}{\sigma_1 K} \right)^{\frac{p}{2p - 1}}, \left( \frac{\delta \Delta}{\sigma_1^2 K} \right)^{\frac{pq}{p (2 q + 1) - 2 q}}}. \thesischfivenumberthis\label{ch5:2c0dbb02-323f-4955-a4f2-8f3b73a0f89d-2_chapter_5} \]
    {When $\sigma_1=0$, we use the deterministic convention $\eta=1$ and $B_{\textnormal{init}}=1$.}
    {For the uniformly sampled output $x_{\mathrm{out}}$, Algorithm~\ref{ch5:algo:nsgd-mvr_chapter_5} guarantees $\Exp{\norm{\nabla f(x_{\mathrm{out}})}}\le\thesischfiveeps$ with total sample complexity}
        \[ \cO\left( {1+} \left( \frac{\sigma_1}{\thesischfiveeps} \right)^{\frac{p}{p - 1}} + \frac{(L_1 + \delta) \Delta}{\thesischfiveeps^2} + \frac{\delta \Delta}{\thesischfiveeps^2} \left( \frac{\sigma_1}{\thesischfiveeps} \right)^{\frac{p}{q (p - 1)}} \right). \]
\end{theorem}
The proof is deferred to~\ref{ch5:appdx:proof-nsgd-mvr_chapter_5}. Essentially, the constant $L_1 + \delta$ is of the same order as the constant $\bar{L}$ from~Assumption~\ref{ch5:ass:mean-squared-smoothness_chapter_5}.

By~\textcite[Corollary~18]{hübler2025gradientclippingnormalizationheavy_nips_workshop}, this implies that the lower bound from~Theorem~\ref{ch5:thm:lower-bound-sofo-mean-squared-smoothness-2_chapter_5} is \emph{tight}: it is matched, up to some constant factors, by the combination of \algname{NSGD-Mom} (Normalized \algname{SGD} with Momentum) with \algname{NSGD-MVR} (Algorithm~\ref{ch5:algo:nsgd-mvr_chapter_5}). Additionally, it is worth noting that Normalized Minibatch \algname{SGD}, like \algname{NSGD-Mom}, also attains the first term in the $\min\thesischfiveens{\ldots}$ of~Theorem~\ref{ch5:thm:lower-bound-sofo-mean-squared-smoothness-2_chapter_5}~\parencite[Corollary~3]{hübler2025gradientclippingnormalizationheavy}.

\subsubsection{Case of Unknown Tail Indices $p$ and $q$}

\begin{theorem}\label{ch5:thm:nsgd-mvr-convergence-analysis-unknown-p-q_chapter_5}
    Under~Assumption~\ref{ch5:ass:lower-boundedness_chapter_5}, Assumption~\ref{ch5:ass:p-bounded-central-moment-gradient_chapter_5}, and Assumption~\ref{ch5:ass:mean-squared-smoothness_chapter_5}, 
    assume 
    $g^{0} = \nabla f(x^{0};\xi^{0})$, %
    $\gamma = \sqrt{\frac{\Delta \eta}{\bar{L} K}}$\footnote{If $\bar L=0$, the assumptions imply $\nabla f\equiv0$, so the conclusion is immediate; the displayed stepsize is used for $\bar L>0$.}, 
    $\eta = K^{-\frac{1}{2}} \in \thesischfiveintof{0}{1}$.
    {Algorithm~\ref{ch5:algo:nsgd-mvr_chapter_5} satisfies}
        \[ \frac{1}{K} \sum_{k=0}^{K - 1} \thesischfiveE{\norm{\nabla f(x^{k})}} = \cO\left( \frac{\sigma_1}{K^{\frac{p - 1}{2 p}}} + \frac{\sqrt{\bar{L} \Delta}}{K^{\nicefrac{1}{4}}} \right). \thesischfivenumberthis\label{ch5:12c01754-28b6-4bc2-9ab4-c6f380187237_chapter_5} \]
\end{theorem}

Our~Theorem~\ref{ch5:thm:nsgd-mvr-convergence-analysis-unknown-p-q_chapter_5} makes the pessimistic assumption $q = 1$, under which we recover the usual $\cO(\smash{K^{-\frac{p - 1}{2 p}}})$ rate with unknown tail index $p$~\parencite{hübler2025gradientclippingnormalizationheavy,liu2025nonconvex}. The bound~\eqref{ch5:12c01754-28b6-4bc2-9ab4-c6f380187237_chapter_5} can be slightly refined under~Assumption~\ref{ch5:ass:L-lipschitz-gradients_chapter_5} and Assumption~\ref{ch5:ass:mean-squared-smoothness-2_chapter_5}.

\begin{theorem}\label{ch5:thm:nsgd-mvr-convergence-analysis-unknown-p-q-2_chapter_5}
    Under~Assumption~\ref{ch5:ass:lower-boundedness_chapter_5}, Assumption~\ref{ch5:ass:L-lipschitz-gradients_chapter_5}, Assumption~\ref{ch5:ass:p-bounded-central-moment-gradient_chapter_5}, and Assumption~\ref{ch5:ass:mean-squared-smoothness-2_chapter_5}, 
    assume $g^{0} = \nabla f(x^{0};\xi^{0})$ and choose the stepsize\footnote{A stepsize term with zero denominator is interpreted as $+\infty$ and omitted from the minimum. If $L_1=\delta=0$, then $\nabla f\equiv0$ and the conclusion is immediate.}
    \[
        \gamma = \min\thesischfiveens{\sqrt{\frac{\Delta}{L_1K}}, \sqrt{\frac{\Delta\eta}{\delta K}}}.
    \]
    Let $\eta = K^{-\frac{1}{2}} \in \thesischfiveintof{0}{1}$.
    {Algorithm~\ref{ch5:algo:nsgd-mvr_chapter_5} satisfies}
        \[ \frac{1}{K} \sum_{k=0}^{K - 1} \thesischfiveE{\norm{\nabla f(x^{k})}} = \cO\left( \frac{\sigma_1}{K^{\frac{p - 1}{2 p}}} + \frac{\sqrt{\delta \Delta}}{K^{\nicefrac{1}{4}}} + \frac{\sqrt{L_1 \Delta}}{K^{\nicefrac{1}{2}}} \right). \]
\end{theorem}
The proofs are deferred to~\ref{ch5:appdx:proof-nsgd-mvr_chapter_5}.

\subsection{Discussion of the Obtained Results}

\stepEightBchange{Under $q$-weak average smoothness, Theorem~\ref{ch5:thm:nsgd-mvr-convergence-analysis_chapter_5} shows that \algname{NSGD-MVR} matches the lower bound of Theorem~\ref{ch5:thm:lower-bound-sofo-mean-squared-smoothness_chapter_5} up to constant factors. Under $(q,\delta)$-similarity, Theorem~\ref{ch5:thm:nsgd-mvr-convergence-analysis-2_chapter_5} shows that \algname{NSGD-MVR} attains the second, variance-reduced branch of the minimum in Theorem~\ref{ch5:thm:lower-bound-sofo-mean-squared-smoothness-2_chapter_5}. The first branch is attained by \algname{NSGD-Mom} according to~\textcite[Corollary~18]{hübler2025gradientclippingnormalizationheavy}; normalized minibatch \algname{SGD} provides an alternative guarantee for this branch~\parencite[Corollary~3]{hübler2025gradientclippingnormalizationheavy}. Therefore, selecting the better guarantee among the appropriate normalized methods matches the full lower-bound minimum in Theorem~\ref{ch5:thm:lower-bound-sofo-mean-squared-smoothness-2_chapter_5} up to constant factors.} For the standard case $q=2$, the stochastic terms are
\[
    \frac{\bar L\Delta}{\varepsilon^2}\left(\frac{\sigma_1}{\varepsilon}\right)^{\nicefrac{p}{2(p-1)}}
    \quad\text{and}\quad
    \frac{\delta\Delta}{\varepsilon^2}\left(\frac{\sigma_1}{\varepsilon}\right)^{\nicefrac{p}{2(p-1)}},
\]
under Assumption~\ref{ch5:ass:mean-squared-smoothness_chapter_5} ($q$-WAS) and Assumption~\ref{ch5:ass:mean-squared-smoothness-2_chapter_5} ($(q,\delta)$-S), respectively. Compared to \textcite{pmlr-v195-liu23c}, this yields faster rates under weaker assumptions.

In the case $p=q$, these stochastic terms coincide with those for \algname{NSGDHess} in \textcite{sadiev2025second} under Assumption~\ref{ch5:ass:q-bounded-central-moment-hessian_chapter_5}:
$\tfrac{\sigma_2 \Delta}{\varepsilon^2}\!\left(\tfrac{\sigma_1}{\varepsilon}\right)^{\nicefrac{1}{(p-1)}}.$
{
When the stochastic losses are twice differentiable and differentiation may be interchanged with expectation and the local limit, $(q,\delta)$--similarity yields the corresponding Hessian-noise control. Without these additional regularity conditions, the implication to Assumption~\ref{ch5:ass:q-bounded-central-moment-hessian_chapter_5} should not be used. For empirical risk minimization with $p=q=2$, the standard regular setting gives the familiar equivalence~\parencite{pmlr-v125-arjevani20a}.
}

These insights motivate the study of a more general framework for second-order stochastic optimization.

\section{Near-Optimal Method under BCM Gradients and BCM Hessians}\label{ch5:sec:near-optimal-method-nsgd-mvr-hess_chapter_5}

Unlike Algorithm~\ref{ch5:algo:nsgd-mvr_chapter_5}, we consider a method that employs Hessian–vector products instead of differences of stochastic gradients. This modification leads to the Hessian-corrected momentum (\texttt{Hess}) technique~\parencite{10.5555/3600270.3600525,salehkaleybar-et-al22}.
    \begin{algorithm}[htbp]
\caption{\algname{{NSGDHess}} (Normalized SGD with Hessian-corrected Momentum)}
\label{ch5:algo:nsgd-mvr-hess_chapter_5}
\begin{algorithmic}[1]
\Require $x^0 \in \R^d$, the starting point;              $K > 0$, the number of iterations;              $g^0 \in \R^d$, an initial vector;              $\gamma > 0$, the stepsize;              $\eta \in \thesischfiveintof{0}{1}$, the momentum parameter for \texttt{Hess};
\State {$x^1 = x^0 - \gamma \frac{g^0}{\norm{g^0}}$}
\For{$k = 1, 2, \ldots, K - 1$}
\State {Sample $u^k \sim \mathcal{U}\left(\thesischfiveintff{0}{1}\right)$}
\State {$\hat{x}^k = u^k x^{k} + (1 - u^k) x^{k-1}$}
\State \textit{$\xi^{k}, \hat{\xi}^{k} \sim \thesischfivecD$ are independent.}
\State \textit{Apply Hessian-vector products}
\State $g^{k} = (1 - \eta) \left( g^{k-1} + \nabla^2 f(\hat{x}^k; \hat{\xi}^{k})(x^{k} - x^{k-1}) \right) + \eta \nabla f\left(x^{k};\xi^{k}\right)$
\State \textit{Do one descent step.}
\State {$x^{k+1} = x^{k} - \gamma \frac{g^{k}}{\norm{g^{k}}}$}
\EndFor
\Ensure {$x_{\mathrm{out}}\sim\mathcal U\{x^0,\ldots,x^{K-1}\}$}
\end{algorithmic}
\end{algorithm}%

\begin{remark}
    {Algorithm~\ref{ch5:algo:nsgd-mvr-hess_chapter_5} was previously considered by~\textcite{sadiev2025second}; here, we refine its analysis and improve on their Theorem~2. It is zero-respecting in the sense of Definition~\ref{ch5:appdx-def:zero-especting-algo_chapter_5}, because
    \[
        \thesischfivesupp(\hat{x}^k)
        \subseteq
        \thesischfivesupp(x^{k-1})\cup\thesischfivesupp(x^k),
    \]
    where $\thesischfivesupp(x)\eqdef\thesischfiveenstq{i\in[d]}{x_i\neq0}$ denotes the support of $x\in\R^d$.}
    
    Additionally, the use of $\hat{x}^k$ taken uniformly at random on the line $\thesischfiveintff{x^{k-1}}{x^{k}}$ allows to use~Lemma~\ref{ch5:appdx-technical-lem:app-von-Bahr-and-Essen_chapter_5}, which provides a better overall sample complexity than the one derived in~\textcite{10.5555/3600270.3600525}.
\end{remark}

\begin{theorem}\label{ch5:thm:nsgd-mvr-hess-convergence-analysis_chapter_5}
    Under~Assumption~\ref{ch5:ass:lower-boundedness_chapter_5}, Assumption~\ref{ch5:ass:L-lipschitz-gradients_chapter_5}, Assumption~\ref{ch5:ass:L-lipschitz-hessians_chapter_5}, Assumption~\ref{ch5:ass:p-bounded-central-moment-gradient_chapter_5}, and Assumption~\ref{ch5:ass:q-bounded-central-moment-hessian_chapter_5}, {choose, for a sufficiently large universal constant $C>0$,}
    {
    \[
    \begin{aligned}
        K=\Bigg\lceil C\Bigg(&1+\left(\frac{\sigma_1}{\thesischfiveeps}\right)^{\frac{p}{p-1}}
        +\frac{\Delta}{\thesischfiveeps}
        \left(\frac{L_1+\sigma_2}{\thesischfiveeps}
        +\sqrt{\frac{L_2}{\thesischfiveeps}}\right)\\
        &+\frac{\Delta\sigma_2}{\thesischfiveeps^2}
        \left(\frac{\sigma_1}{\thesischfiveeps}\right)^{\frac{p}{q(p-1)}}
        +\frac{L_2^{1/2}\Delta\sigma_1^{1/4}}{\thesischfiveeps^{7/4}}
        \left(\frac{\sigma_1}{\thesischfiveeps}\right)^{\frac{1}{4(p-1)}}\Bigg)\Bigg\rceil .
    \end{aligned}
    \]}
    {Let the initial gradient estimate $g^{0}$ be given by}
        \[ g^{0} = \frac{1}{B_{\textnormal{init}}} \sum_{j = 0}^{B_{\textnormal{init}} - 1} \nabla f\left( x^{0}; \xi_{j}^{0} \right), \]
    where {$B_{\textnormal{init}} = \left\lceil\max\thesischfiveens{1, \left( \frac{2\sigma_1}{\thesischfiveeps} \right)^{\frac{p}{p - 1}}}\right\rceil$}, let the stepsize\footnote{A stepsize term with zero denominator is interpreted as $+\infty$ and omitted from the minimum.}
        \[ \gamma = \min\thesischfiveens{\sqrt{\frac{\Delta}{L_1 K}}, \sqrt{\frac{\Delta \eta^{\nicefrac{1}{q}}}{\sigma_2 K}}, \sqrt[3]{\frac{\Delta \eta^{\nicefrac{1}{2}}}{L_2 K}}}, \]
    {the momentum parameter $\eta = \min\thesischfiveens{1, \eta_{\textnormal{eff}}}$, where}
        \begin{alignat*}{2}
            &{\eta_{\textnormal{eff}} =}\max\left\{\left( \frac{\thesischfiveeps}{\sigma_1 K} \right)^{\frac{p}{2p - 1}}\!, \left( \frac{\sigma_2 \Delta}{\sigma_1^2 K} \right)^{\frac{pq}{p (2 q + 1) - 2 q}}\!, 
                \left( \frac{L_2^{\nicefrac{1}{2}} \Delta}{\sigma_1^{\nicefrac{3}{2}} K} \right)^{\frac{4 p}{7 p - 6}}\right\}.
        \end{alignat*}
    {When $\sigma_1=0$, we use the deterministic convention $\eta=1$ and $B_{\textnormal{init}}=1$.}
    {For the uniformly sampled output $x_{\mathrm{out}}$, Algorithm~\ref{ch5:algo:nsgd-mvr-hess_chapter_5} guarantees $\Exp{\norm{\nabla f(x_{\mathrm{out}})}}\le\thesischfiveeps$ with total sample complexity}
    \[
    \begin{aligned}
        \cO\Bigg(&1+\left(\frac{\sigma_1}{\thesischfiveeps}\right)^{\frac{p}{p-1}}
        +\frac{\Delta}{\thesischfiveeps}\left(\frac{L_1+\sigma_2}{\thesischfiveeps}
        +\sqrt{\frac{L_2}{\thesischfiveeps}}\right)\\
        &+\frac{\Delta\sigma_2}{\thesischfiveeps^2}
        \left(\frac{\sigma_1}{\thesischfiveeps}\right)^{\frac{p}{q(p-1)}}
        +\frac{L_2^{1/2}\Delta\sigma_1^{1/4}}{\thesischfiveeps^{7/4}}
        \left(\frac{\sigma_1}{\thesischfiveeps}\right)^{\frac{1}{4(p-1)}}\Bigg).
    \end{aligned}
    \]
\end{theorem}
The proof is deferred to~\ref{ch5:appdx:proof-nsgd-mvr-hess_chapter_5}.

\section{Gradient Clipping for High-Probability Convergence}\label{ch5:sec:clipped-nsgd-mvr_chapter_5}

In this section, we conduct a high-probability analysis under heavy-tailed noise. Since we work under Assumption~\ref{ch5:ass:p-bounded-central-moment-gradient_chapter_5}, we incorporate gradient clipping in Algorithm~\ref{ch5:algo:nsgd-mvr_chapter_5}; the pseudo-code of the new method is provided in Algorithm~\ref{ch5:algo:clipped-nsgd-mvr_chapter_5}. {Formally, for every $v\in\R^d$,
\[
\thesischfiveclip(v,\thesischfivelbd)\eqdef
\begin{cases}
\min\thesischfiveens{1,\thesischfivelbd/\norm{v}}\,v,&v\neq0,\\
0,&v=0,
\end{cases}
\]}
where $\thesischfivelbd > 0$ is called clipping level/threshold. The proofs are deferred to~\ref{ch5:appdx:proof-clipped-nsgd-mvr_chapter_5}.

    \begin{algorithm}[htbp]
\caption{\algname{D-\thesischfiveclip-NSGD-MVR} (Double Clipped Normalized SGD with \texttt{MVR})}
\label{ch5:algo:clipped-nsgd-mvr_chapter_5}
\begin{algorithmic}[1]
\Require $x^0 \in \R^d$, the starting point;              $K > 0$, the number of iterations;              $g^0 \in \R^d$, an initial vector;              $\gamma > 0$, the stepsize;              $\eta \in \thesischfiveintof{0}{1}$, the momentum parameter for \texttt{MVR};              $\thesischfivelbd_1, \thesischfivelbd_2 > 0$, the clipping thresholds;
\State {$x^1 = x^0 - \gamma \frac{g^0}{\norm{g^0}}$}
\For{$k = 1, 2, \ldots, K - 1$}
\State \textit{Apply \texttt{MVR} with \thesischfiveclip~operator.}
\State $\begin{aligned}[t]
                    g^{k} ={}&(1 - \eta) \bigl(g^{k-1} + \thesischfiveclip(\nabla f\left(x^{k};\xi^{k}\right) - \nabla f\left(x^{k-1};\xi^{k}\right), \thesischfivelbd_1)\bigr) \\
                    &\quad+ \eta \, \thesischfiveclip(\nabla f\left(x^{k};\xi^{k}\right), \thesischfivelbd_2)
                \end{aligned}$
\State \textit{Do one descent step.}
\State {$x^{k+1} = x^{k} - \gamma \frac{g^{k}}{\norm{g^{k}}}$}
\EndFor\stepEightBdeletion
\end{algorithmic}
\end{algorithm}%

Algorithm~\ref{ch5:algo:clipped-nsgd-mvr_chapter_5} modifies Algorithm~\ref{ch5:algo:nsgd-mvr_chapter_5} by applying clipping not only to the stochastic gradient — to control heavy-tailed noise — but also to the gradient difference term. This additional clipping step is the primary feature distinguishing our method from Accelerated NSGD with clipping and momentum \parencite{pmlr-v195-liu23c}. Crucially, this modification enables us to establish high-probability guarantees under assumptions weaker than those in \textcite{pmlr-v195-liu23c}. While their analysis relies on individual smoothness (i.e., assuming $f(\cdot;\xi)$ is $\bar L$-smooth almost surely), our analysis holds under $q$-WAS (Assumption~\ref{ch5:ass:mean-squared-smoothness_chapter_5}). Individual smoothness is a significantly stronger condition; notably, under that assumption alone, \textcite{lei2019stochastic} showed it is possible to attain an $\cO(\varepsilon^{-4})$ rate independent of the heavy-tail index $p$.

\begin{theorem}\label{ch5:thm:clipped-nsgd-mvr-convergence-analysis_chapter_5}
    Under~Assumption~\ref{ch5:ass:lower-boundedness_chapter_5}, Assumption~\ref{ch5:ass:p-bounded-central-moment-gradient_chapter_5}, and Assumption~\ref{ch5:ass:mean-squared-smoothness_chapter_5}, let $K \ge 1$ and $\beta \in \thesischfiveintof{0}{1}$ be such that $\log \frac{8 K}{\beta} \ge 1$. Let $x^{0}\in\mathbb{R}^d$ and define $\Delta_0 \eqdef f(x^{0}) - f^{\inf}$. 
    Suppose that Algorithm~\ref{ch5:algo:clipped-nsgd-mvr_chapter_5} is run with $g^{0} = 0$, momentum parameter $\eta = \max\{K^{-\frac{p}{2 p - 1}}, K^{-\frac{p q}{p (2 q + 1) - 2 q}}\}$, clipping thresholds $\thesischfivelbd_1 = 2 \gamma \bar{L} \eta^{-\frac{1}{q}}$ and $\thesischfivelbd_2 = \max\{4 \sqrt{\bar{L} \Delta_0}, \sigma_1 \eta^{-\frac{1}{p}}\}$, and {stepsize
    \[
    \begin{aligned}
    \gamma=\min\Bigg\{&
        \sqrt{\frac{\Delta_0}{2\bar L K}},\
        \frac{\eta}{8}\sqrt{\frac{\Delta_0}{2\bar L}},\
        \frac{\sqrt{\Delta_0/\bar L}}{704\,\eta K\log(8K/\beta)},\\
        &\frac{\Delta_0}{176\,\sigma_1\eta^{(p-1)/p}K\log(8K/\beta)},\
        \sqrt{\frac{\Delta_0\eta^{1/q}}{736\,\bar L K\log(8K/\beta)}}
    \Bigg\}.
    \end{aligned}
    \]}
    {Any term whose denominator contains a zero value of $\bar L$ or $\sigma_1$ is interpreted as $+\infty$ and hence does not constrain the minimum.}
    
    {Then, with probability at least $1-\beta$, the iterates generated by Algorithm~\ref{ch5:algo:clipped-nsgd-mvr_chapter_5} satisfy}
        \[ \frac{1}{K} \sum_{k=0}^{K - 1} \norm{\nabla f(x^{k})} \le \frac{2 \Delta_0}{\gamma K}, \]
    and, by our choice of parameters, the gradient norm converges with high probability at the rate
        \[ \frac{1}{K} \sum_{k=0}^{K - 1} \norm{\nabla f(x^{k})} = \cO\left( \left( \frac{\sqrt{\bar{L} \Delta_0} + \sigma_1}{K^{\frac{p - 1}{2 p - 1} \wedge \frac{q (p - 1)}{p (2 q + 1) - 2 q}}} \right)\log \frac{K}{\beta} \right). \]
    
\end{theorem}

The detailed proof is provided in~\ref{ch5:appdx:proof-clipped-nsgd-mvr_chapter_5}.

Notably, in the case where $p=q$, we establish the same {sample complexity $\widetilde \cO\left(\varepsilon^{-\frac{2p-1}{p-1}}\right)$} as \textcite{pmlr-v195-liu23c} did for Algorithm~2, but under weaker assumptions. Moreover, our lower bounds (see Theorem~\ref{ch5:thm:lower-bound-sofo-mean-squared-smoothness_chapter_5}) indicate that this rate is optimal in terms of the dependence on $\varepsilon$. 

Another interesting observation is that when $p\leq q$, the high-probability rate remains identical to the $p=q$ case. In contrast, for in-expectation guarantees, we observe a faster rate (see Theorem~\ref{ch5:thm:nsgd-mvr-convergence-analysis_chapter_5}). Determining whether this gap is fundamental or a limitation of the current proof technique remains an open question, and calls for a tighter high-probability analysis for these methods.

Furthermore, \ref{ch5:appdx:proof-clipped-nsgd-mvr_chapter_5} establishes high-probability convergence guarantees for Algorithm~\ref{ch5:algo:clipped-nsgd-mvr_chapter_5} under Assumption~\ref{ch5:ass:mean-squared-smoothness-2_chapter_5}, and for Algorithm~\ref{ch5:algo:clipped-nsgd-mvr-hess_chapter_5} under Assumption~\ref{ch5:ass:L-lipschitz-hessians_chapter_5} and Assumption~\ref{ch5:ass:q-bounded-central-moment-hessian_chapter_5}.

\textbf{Interpretation of Clipping Thresholds and Practical Implementation.} 
While deriving optimal high-probability bounds in heavy-tailed settings remains an open challenge, our double-clipping mechanism provides a significant theoretical advantage: it establishes convergence guarantees without assuming bounded noise---a restrictive assumption often required by prior methods lacking this specific mechanism \parencite[e.g.][]{pmlr-v195-liu23c}. 

Beyond the theoretical guarantees, our analysis offers a crucial insight for practical implementation. Under Assumption~\ref{ch5:ass:p-bounded-central-moment-gradient_chapter_5} (heavy-tailed noise, exponent $p$) and Assumption~\ref{ch5:ass:mean-squared-smoothness_chapter_5} and Assumption~\ref{ch5:ass:mean-squared-smoothness-2_chapter_5} ($q$-WAS and $(q,\delta)$-S, exponent $q$), our theoretical clipping thresholds scale as $\lambda_1 \sim \mathcal{O}(\gamma \eta^{-1/q})$ for the gradient differences, and $\lambda_2 \sim \mathcal{O}(\eta^{-1/p})$ for the raw stochastic gradients. In practice, we can rewrite the gradient difference clipping step by factoring out the stepsize $\gamma$:
\[
    {
    \hbox{$\displaystyle
    \thesischfiveclip\!\left(\nabla f(x^{k};\xi^{k})-\nabla f(x^{k-1};\xi^{k}),\lambda_1\right)
    =
    \gamma\,\thesischfiveclip\!\left(
        \frac{\nabla f(x^{k};\xi^{k})-\nabla f(x^{k-1};\xi^{k})}{\gamma},
        \bar{\lambda}_1
    \right)$}.}
\]
where we define the rescaled threshold $\bar{\lambda}_1 := \gamma^{-1} \lambda_1$. This reveals that $\bar{\lambda}_1$ scales as $\mathcal{O}(\eta^{-1/q})$. {Consequently, when $p=q$, $\bar\lambda_1$ and $\lambda_2$ have the same exponent in $\eta$, namely $\mathcal O(\eta^{-1/p})$. Their multiplicative constants still depend on different problem parameters, so setting $\bar\lambda_1=\lambda_2$ is a practical one-parameter heuristic rather than a consequence required by the theorem.}

\section{Conclusion and Future Work}

This chapter characterized stochastic nonconvex optimization under heavy-tailed gradient information with tail index $p$, together with same-sample gradient differences or stochastic Hessian information governed by index $q$. The lower bounds distinguish what can be achieved by first-order variance-reduced and higher-order zero-respecting algorithms. Normalized momentum variance reduction exploits correlated gradient differences, while \algname{NSGDHess} transports the search vector with Hessian--vector information. Double clipping separately controls the innovation and raw stochastic-gradient terms when a high-probability statement is required.

The output criterion is regime-dependent. The known-index in-expectation results for \algname{NSGD-MVR} and \algname{NSGDHess} use a uniformly sampled iterate from the generated trajectory; unknown-index results control the arithmetic mean of expected gradient norms. Under $q$-weak average smoothness, the first-order bounds match the corresponding lower bound in the stated oracle class. Under $(q,\delta)$-similarity, the full lower-bound minimum is attained by combining the appropriate normalized methods rather than by attributing every branch to \algname{NSGD-MVR} alone. The second-order results are near-optimal in their stated regimes. For \algname{D-Clip-NSGD-MVR}, the theorem controls a realized trajectory-average gradient norm with probability at least $1-\beta$; when $p=q$, its accuracy exponent matches the lower bound up to logarithmic factors.

For general $p$ and $q$, a gap remains between the in-expectation and high-probability rates. The methods also rely on the stated same-sample oracle structures and parameter-dependent momentum, stepsize, and clipping choices; none of the results is a universal terminal-iterate guarantee. Chapter~\ref{ch:chapter_6} changes the principal bottleneck from heavy-tailed oracle information to lossy communication, combining normalization with error-feedback reconstruction under generalized smoothness.

\clearpage%
  \endgroup

\part{Communication-Efficient Distributed Learning}

\chapter{Error Feedback under $(L_0, L_1)$-Smoothness: Normalization and Momentum}
\label{ch:chapter_6}
\thesischaptersourcenote
  {Error Feedback under $(L_0,L_1)$-Smoothness: Normalization and Momentum}
  {khirirat2024errorfeedbackl0l1smoothnessnormalization}
  \begingroup
  \thesischaptersetup{6}%
  
\begin{paperabstract}
		We provide the first proof of convergence for normalized error feedback algorithms across a wide range of machine learning problems. 
		Despite their popularity and efficiency in training deep neural networks, traditional analyses of error feedback algorithms rely on the smoothness assumption that does not capture the properties of objective functions in these  problems. 
		Rather, these problems have recently been shown to satisfy generalized smoothness assumptions, and the theoretical understanding of error feedback algorithms under these assumptions remains largely unexplored. 
		{Moreover, to the best of our knowledge, all existing analyses under generalized smoothness either (i) focus on single-node settings or (ii) impose restrictive assumptions in distributed settings, such as bounded data heterogeneity or an almost-sure bound on stochastic-gradient noise.}
		{In this chapter, we propose distributed error feedback algorithms that utilize normalization to achieve the $\cO(1/\sqrt{K})$ convergence rate for nonconvex problems under generalized smoothness. Our analyses apply to distributed settings without a bounded-heterogeneity assumption, and enable stepsize tuning that is independent of problem parameters.}
		Additionally, we provide strong convergence guarantees of normalized error feedback algorithms for stochastic settings. 
		Finally, we show that due to their larger allowable stepsizes, our new normalized error feedback algorithms outperform their non-normalized counterparts on various tasks, including the minimization of polynomial functions, logistic regression, and ResNet-20 training. 
	\end{paperabstract}

	\section{Introduction}
	\label{ch6:sec:intro_chapter_6}

	Machine learning models achieve impressive prediction and classification power by employing sophisticated architectures, comprising vast numbers of model parameters, and requiring training on massive datasets.
	Distributed training has emerged as an important approach, where multiple machines with their own local training data collaborate to train a model efficiently within a reasonable time.
	Many optimization algorithms can be easily adapted for distributed training frameworks. For example, stochastic gradient descent (\algname{SGD}) can be modified into distributed stochastic gradient descent within a data parallelism framework, and into federated averaging algorithms~\parencite{mcmahan2017communication} in a federated learning framework.
	However, the communication overhead of running these distributed algorithms poses a significant barrier to scaling up to large models. For example, training the VGG-16 model~\parencite{simonyan2014very} using distributed stochastic gradient descent involves communicating $138.34$ million parameters, thus consuming over $500$MB of storage and posing an unmanageable burden on the communication network between machines.

	One approach to mitigate the communication burden is to apply compression.
	In this approach, the information, such as gradients or model parameters, is compressed using sparsifiers or quantizers to be transmitted with much lower communicated bits between machines. 
	However, while this reduces communication overhead, too coarse compression often brings substantial challenges in maintaining high training performance due to information loss, and in extreme cases, it may potentially lead to divergence. 
	Therefore, error feedback mechanisms have been developed to improve the convergence performance of compression algorithms, while ensuring high communication efficiency. 
	Examples of error feedback mechanisms include \algname{EF14}~\parencite{seide20141,stich2018sparsified,alistarh2018convergence,wu2018error,gorbunov2020linearly}, \algname{EF21}~\parencite{richtarik2021ef21,fatkhullin2021ef21}, \algname{EF21-SGDM} ~\parencite{fatkhullin2024momentum}, \algname{EF21-P}~\parencite{gruntkowska2023ef21}, and \algname{EControl}~\parencite{gao2023econtrol}. 
	Several studies developing error feedback algorithms often assume  the smoothness of an objective function, i.e., 
	its gradient is Lipschitz continuous.

	However, many modern learning problems, such as distributionally robust optimization~\parencite{jin2021non} and deep neural network training, are often non-smooth. 
	For instance, the gradient of the loss computed for deep neural networks, such as LSTM~\parencite{zhang2019gradient}, ResNet20~\parencite{zhang2019gradient}, and transformer models~\parencite{crawshaw2022robustness}, is not Lipschitz continuous.
	These empirical findings highlight the need for a new smoothness assumption.
	One such assumption is $(L_0,L_1)$-smoothness, originally introduced by~\textcite{zhang2019gradient}, for twice differentiable functions, and later extended to differentiable functions by~\textcite{chen2023generalized}.

	To solve generalized smooth problems, clipping and normalization have been widely utilized in first-order algorithms. 
	Gradient descent with gradient clipping was initially shown by~\textcite{zhang2019gradient} to achieve lower iteration complexity, i.e., fewer iterations needed to attain a target solution accuracy,  than classical gradient descent.
	Subsequent works have further refined the convergence theory of clipped gradient descent~\parencite{koloskova2023revisiting}, and improved its convergence performance by employing momentum updates~\parencite{zhang2020improved}, variance reduction techniques~\parencite{reisizadeh2023variance}, and adaptive step sizes~\parencite{wang2024provable,li2024convergence,takezawa2024polyak}.
	Similar convergence results have been obtained for gradient descent using normalization~\parencite{zhao2021convergence}, and its momentum variants~\parencite{hubler2024parameter}, including generalized \algname{SignSGD}~\parencite{crawshaw2022robustness}.
	However, these first-order algorithms have mostly been explored in training on a single machine. 
	To the best of our knowledge, distributed algorithms under generalized smoothness have been investigated in only a few works, e.g., by \textcite{crawshaw2024federated,liu2022communication}.
	Nonetheless, these works rely on assumptions limiting families of optimization problems, including data heterogeneity, almost sure variance bounds, and symmetric noise distributions around the mean assumptions. 
	Furthermore, these first-order algorithms under generalized smoothness do not incorporate compression techniques to improve communication efficiency. 
	These aspects motivate us to develop \emph{distributed communication-efficient algorithms for solving nonconvex generalized smooth problems}.

	\subsection{Contributions}
	
		{In this chapter, we develop distributed error-feedback methods for nonconvex optimization under generalized smoothness. Our contributions are summarized below.}
	
	\textbf{$\bullet$ Importance of normalization.}
	Just as gradient clipping is crucial for gradient descent, we empirically demonstrate that normalization stabilizes the convergence of error feedback algorithms for minimizing nonconvex generalized smooth functions. {In this chapter, we introduce a normalized variant of \algname{EF21}, the error feedback algorithm of \textcite{richtarik2021ef21}, and establish its convergence for nonconvex generalized-smooth problems.}
	In a single-node setting, this new method, which we call \normEFname{GD}, or more compactly as \normEFname{}, provides larger stepsize, and faster convergence rate than its non-normalized counterpart \algname{EF21} for minimizing simple nonconvex polynomial functions that satisfy generalized smoothness, as shown by Figure~\ref{ch6:fig:exp1_chapter_6}.

	\begin{figure}[h]
		\centering
			\includegraphics[width=0.3\textwidth]{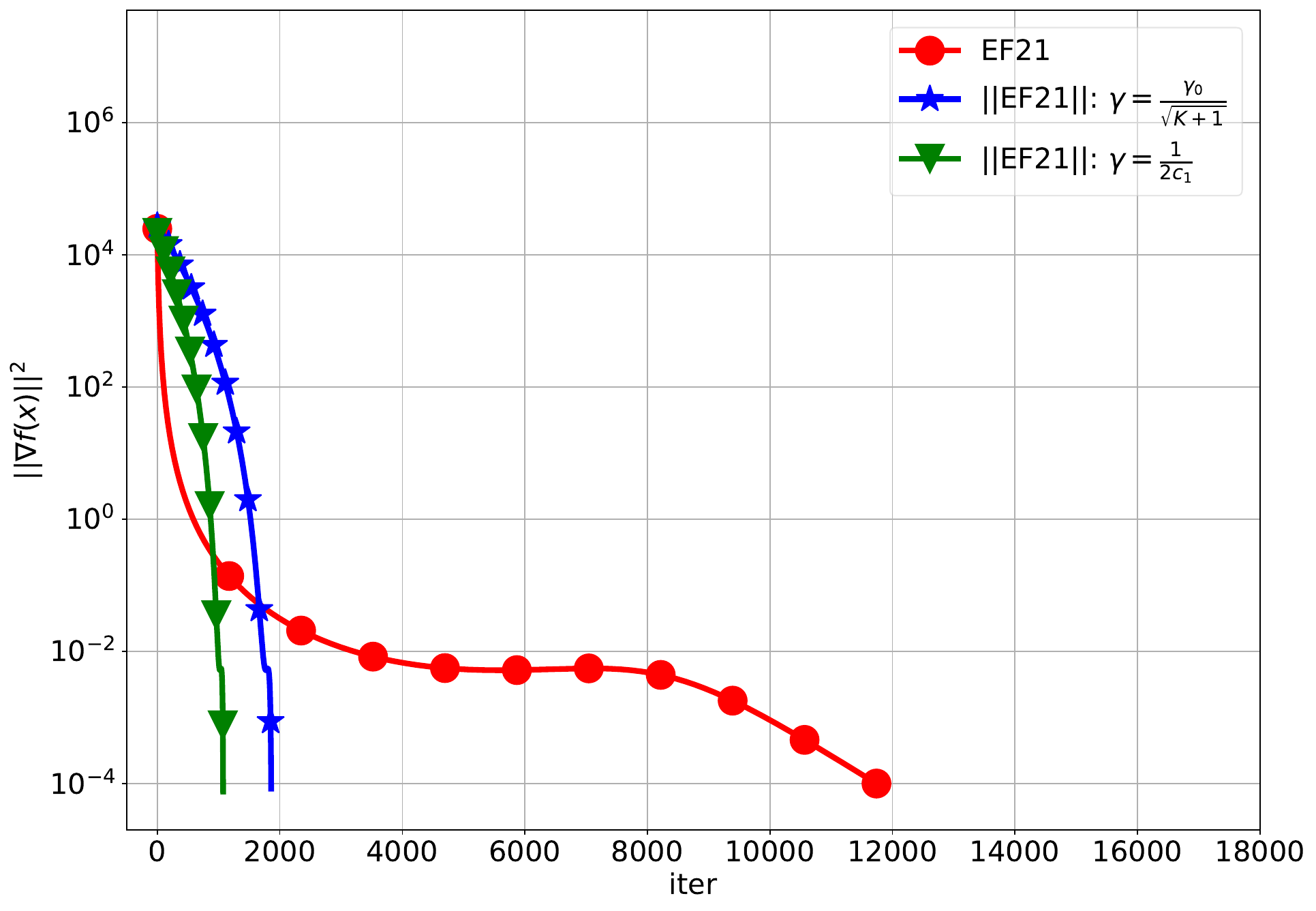}
		\hfill
			\includegraphics[width=0.3\textwidth]{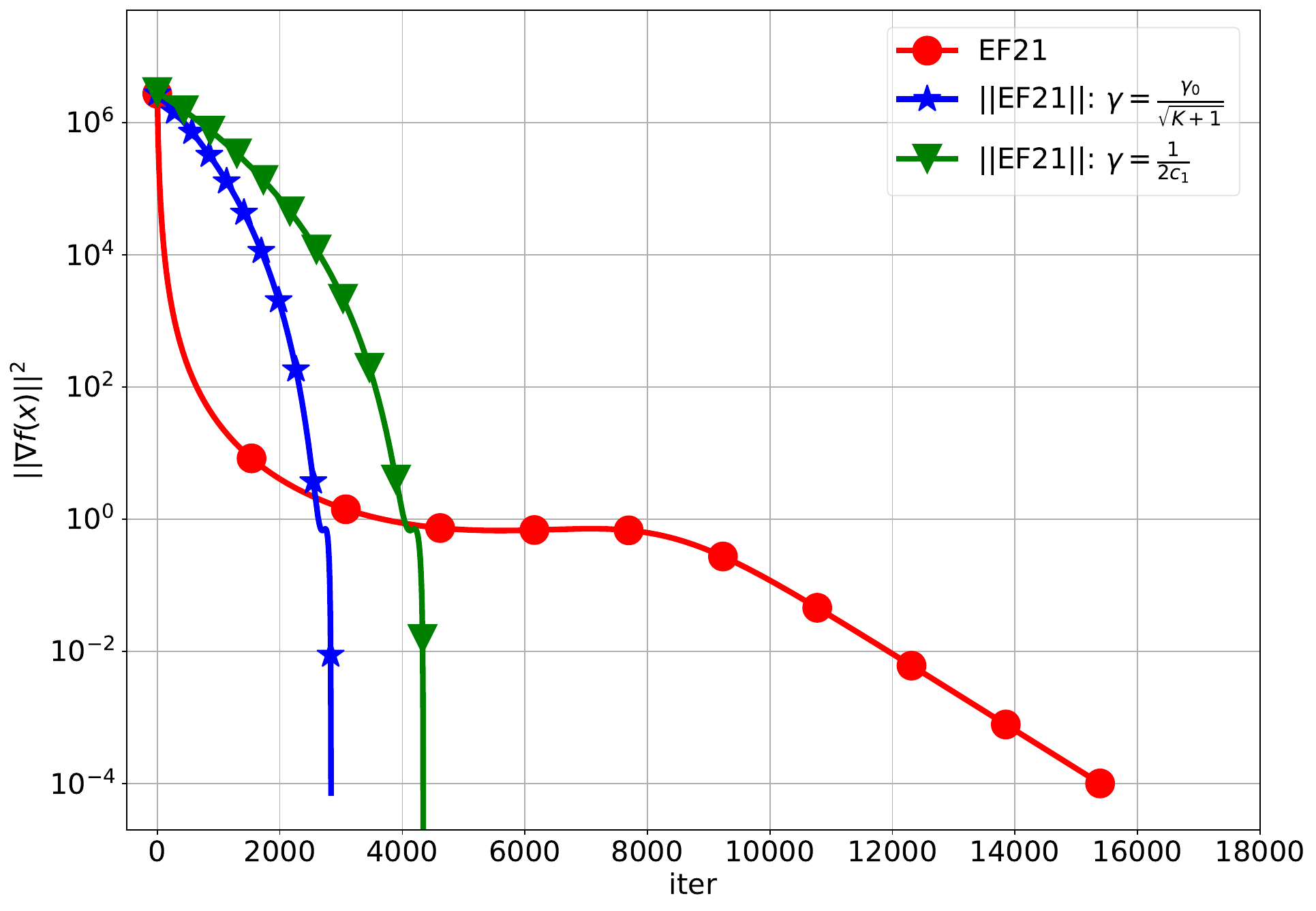}
		\hfill
			\includegraphics[width=0.3\textwidth]{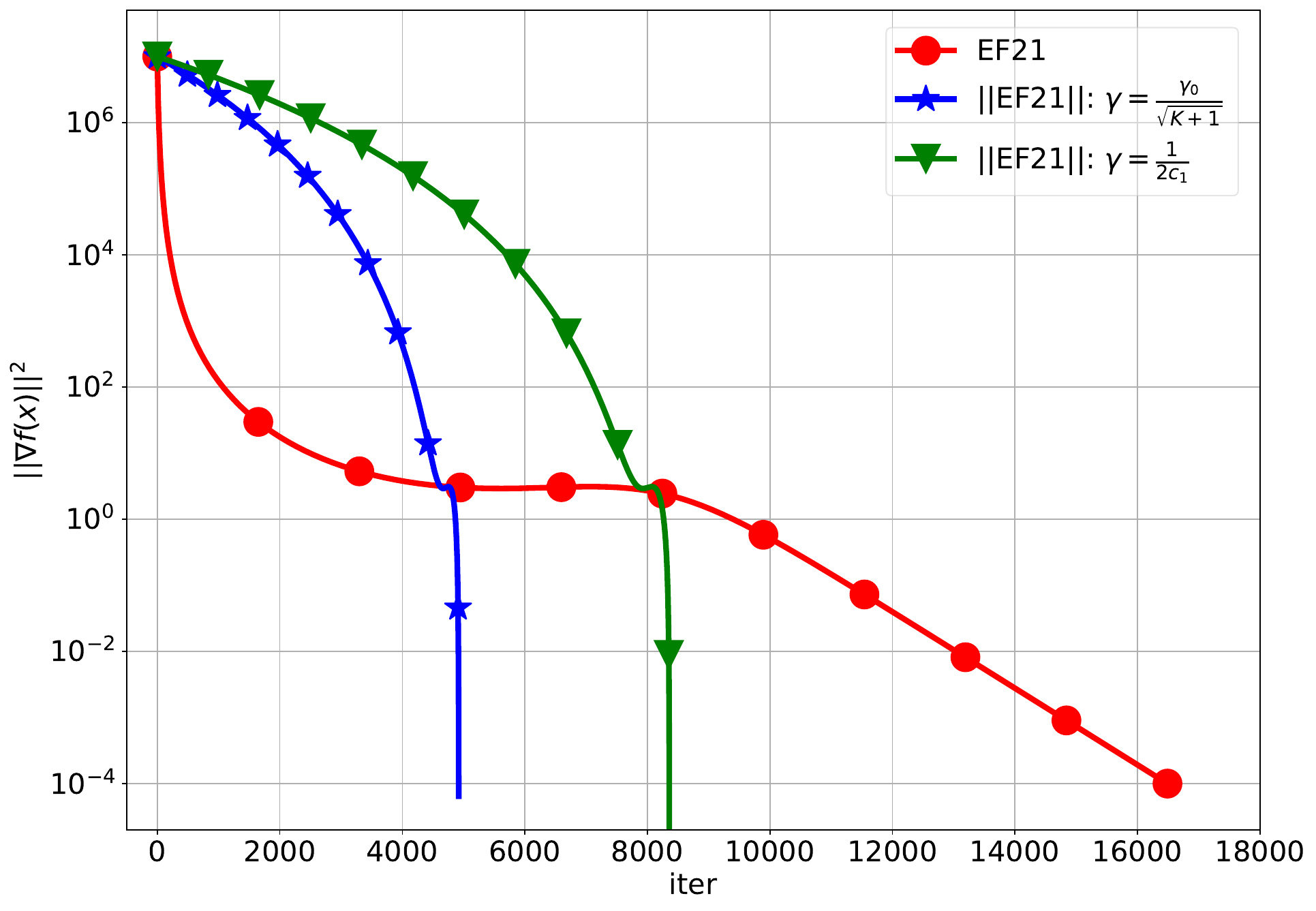}
		\caption{The minimization of polynomial functions using \algname{EF21} with $\gamma = \frac{1}{L + L \sqrt{\frac{\beta}{\theta}}}$, and \normEFname{}  with $\gamma = \frac{\hat\gamma}{\sqrt{K+1}}$, $\hat\gamma = 1$ (blue line) and $\gamma = \frac{1}{2c_1}$ (green line). Here, we ran both algorithms for (1) $L_0 = 4$, $L_1 = 1$, and $K=2,000$ (left), (2) $L_0 = 4$, $L_1 = 4$, and $K=5,000$ (middle), and (3) $L_0 = 4$, $L_1 = 8$, and $K=16,000$ (right). }
		\label{ch6:fig:exp1_chapter_6}
	\end{figure}

	\textbf{$\bullet$ Convergence of normalized error feedback algorithms.}
	We establish an $\mathcal{O}(1/\sqrt{K})$ convergence rate in the gradient norm for \normEFname{} on nonconvex generalized smooth problems. 
	\normEFname{} achieves the same  rate as  \algname{EF21} under $L$-smoothness by \textcite{richtarik2021ef21}. Our results are derived under standard assumptions, i.e., generalized smoothness and the existence of lower bounds on the objective function, and are applicable in distributed settings regardless of any data heterogeneity degree, unlike the results by \textcite{crawshaw2024federated,liu2022communication}. Additionally, our stepsize rules for \normEFname{} ensure convergence without requiring knowledge of the generalized smoothness constants $L_0$ or $L_1$, in contrast to \textcite{richtarik2021ef21}, where the stepsize depends on the smoothness constant $L$ (which is often inaccessible).

	\textbf{$\bullet$ Extension to stochastic settings.}
	Furthermore, we propose a variant of \algname{EF21-SGDM}, an error feedback algorithm with momentum updates by~\textcite{fatkhullin2024momentum}, that employs normalization for solving nonconvex, stochastic optimization under generalized smoothness. 
	Specifically, we prove that \normEFname{SGDM}  with suitable stepsize choices attains the same $\cO(1/K^{1/4})$ convergence rate in the gradient norm as  \algname{EF21-SGDM}.

	\textbf{$\bullet$ Numerical evaluation.} 
	We implemented \normEFname{} using the stepsize rules derived from our theory, and compared its performance against  \algname{EF21}. Both algorithms were evaluated on three learning tasks: minimizing nonconvex polynomial functions, solving logistic regression with a nonconvex regularizer, and training ResNet-20 on the CIFAR-10 dataset. Thanks to its larger stepsizes, \normEFname{} outperforms \algname{EF21}, in terms of both convergence speed and solution accuracy across these tasks.

	\begin{table}[h]
		
		\caption{Comparisons of complexities and assumptions between known and our results for \algname{EF21} variants. The  complexity is defined by the iteration count $K$ required by the algorithms to attain $\underset{k=0,1,\ldots,K}{\min}\Exp\!\left[\norm{\nabla f(x^k)}\right] \leq \varepsilon$. $(L_0,L_1)$-smoothness refers to generalized smoothness in Assumption~\ref{ch6:assum:LzeroLoneSmooth_chapter_6}. The variance bound in expectation is defined in Assumption~\ref{ch6:assum:bounded_variance_chapter_6}.}
		\label{ch6:sample-table_chapter_6}
\begin{center}
			\resizebox{0.95\textwidth}{!}{%
				\begin{tabular}{ccccc}
					\bf Methods & \bf Complexity & \bf Smoothness & \bf Variance bound & \bf Normalization \\ \hline 
					\begin{tabular}{c} \algname{EF21} \\ \tabcite{richtarik2021ef21} \end{tabular} & $\cO(1/\varepsilon^2)$ & $L$ & No & No \\ 
					\begin{tabular}{c} \algname{EF21-SGDM}  \\ \tabcite{fatkhullin2024momentum} \end{tabular} &  $\cO( {1}/{\varepsilon^4})$  & $L$ & expectation & No \\
					\rowcolor{linen} \begin{tabular}{c} \normEFname{} \\  {\bf NEW} (Alg.~\ref{ch6:alg:normalized_ef21_chapter_6}) \end{tabular}
					&  $\cO(1/\varepsilon^2)$ & $(L_0,L_1)$ & No & Yes \\
					\rowcolor{linen} \begin{tabular}{c} \normEFname{SGDM}  \\ {\bf NEW} (Alg.~\ref{ch6:alg:normalized_ef21_sgdm_chapter_6}) \end{tabular} & $\cO(1/\varepsilon^4)$ & $(L_0,L_1)$ & Expectation & Yes \\ \hline 
				\end{tabular}
			}
		\end{center}

	\end{table}

	\section{Related Work}
	\label{ch6:sec:related_work_chapter_6}

	\textbf{Error feedback.}
	Error feedback mechanisms have been utilized in various algorithms with communication compression, leading to significant improvements in solution accuracy, while reducing communication. 
	As the first version of these mechanisms, \algname{EF14} was introduced by~\textcite{seide20141}, and later analyzed for first-order algorithms in both single-node \parencite{stich2018sparsified,karimireddy2019error} and distributed settings \parencite{alistarh2018convergence,wu2018error,tang2019doublesqueeze,basu2019qsparse,gorbunov2020linearly,li2020acceleration,qian2021error,tang2021errorcompensatedx}. 
	Next, \algname{EF21} is another error feedback variant proposed by \textcite{richtarik2021ef21}, which offers strong convergence guarantees for distributed gradient algorithms with any contractive compressors, without requiring bounded gradient norm or bounded data heterogeneity assumptions. \algname{EF21} can also be adapted for  stochastic optimization through sufficiently large mini-batches \parencite{fatkhullin2021ef21} or momentum updates \parencite{fatkhullin2024momentum}. More recently, \algname{EControl} was developed by \textcite{gao2023econtrol}  to guarantee provably superior complexity results for distributed stochastic optimization compared to prior error feedback mechanisms. %
	To the best of our knowledge, these existing works on error feedback have focused solely on optimization under traditional \( L \)-smoothness. {In this chapter, we introduce a normalized variant of the \algname{EF21} methods \parencite{richtarik2021ef21} for solving nonconvex generalized-smooth problems.}
	In particular, we prove that  \normEFname{} under generalized smoothness achieves the same $\cO(1/\sqrt{K})$ rate as \algname{EF21} under traditional smoothness, and demonstrate in experiments that \normEFname{} permits larger step sizes, and thus attains faster convergence than  \algname{EF21}.

	\textbf{Non-smoothness assumptions.}
	Empirical findings suggest that the traditional smoothness used for analyzing optimization algorithms does not capture the properties of objective functions in many machine learning problems, especially deep neural network training problems. 
	This motivates researchers to consider different assumptions to replace this traditional smoothness condition. 
	First introduced by~\textcite{zhang2019gradient}, the $(L_0,L_1)$-smoothness condition on a twice differentiable function $f(x)$ is defined by $\norm{ \nabla^2 f(x)} \leq L_0+L_1 \norm{\nabla f(x)}$ for $x\in\R^d$. This $(L_0,L_1)$-smoothness has been extended to
	differentiable functions without assuming the existence of the Hessian.  
	For instance, the smoothness with a differentiable function $\ell(x)$~\parencite{li2024convex}, and symmetric generalized smoothness~\parencite{chen2023generalized} cover the $(L_0,L_1)$-smoothness when the Hessian exists, and includes many important machine learning problems, such as phase retrieval problems~\parencite{chen2023generalized}, and distributionally robust optimization~\parencite{levy2020large}. Other classes of non-smoothness assumptions, which are not related to the generalized smoothness but capture other optimization problems, include H{\"o}lder's continuity of the gradient \parencite{devolder2014first},  the relative smoothness~\parencite{bauschke2017descent}, and the polynomial growth of the gradient norm~\parencite{mai2021stability}.
	{In this chapter, we impose generalized smoothness to establish convergence of \normEFname{} in deterministic and stochastic optimization.}

	\textbf{Gradient clipping and normalization.} Clipping and normalization are commonly employed in gradient-based methods for solving generalized smooth problems. Clipped (stochastic) gradient descent has been studied for both nonconvex and convex problems under  $(L_0,L_1)$-smoothness conditions by \textcite{zhang2019gradient,koloskova2023revisiting}. Extensions to clipped gradient algorithms have been proposed, including momentum updates \parencite{zhang2020improved}, variance reduction methods \parencite{reisizadeh2023variance}, and adaptive step sizes \parencite{wang2024provable,li2024convergence,takezawa2024polyak,gorbunov2024methods}. Comparable complexities have been achieved for normalized gradient descent \parencite{zhao2021convergence}, and its momentum-based variants \parencite{hubler2024parameter}, including  \algname{SignSGD} \parencite{crawshaw2022robustness} and its variance-reduction variants~\parencite{jiang2024efficient}.
	Convergence properties of gradient-based algorithms have also been explored under more generalized forms of non-uniform smoothness, extending beyond the $(L_0,L_1)$-smoothness by \textcite{zhang2019gradient} to cover a wider range of optimization problems. For example, variants of (stochastic) gradient descent have been analyzed under $\alpha$-symmetric generalized smoothness by \textcite{chen2023generalized}, and under $\ell$-smoothness involving certain differentiable functions $\ell(\cdot)$ by \textcite{li2024convex,li2024convergence}. 
	However, the majority of these analyses focus on the single-node setting. To the best of our knowledge, only a limited number of works, such as those by \textcite{crawshaw2024federated,liu2022communication}, have examined federated averaging algorithms for nonconvex problems under generalized smoothness. These works, however, often rely on restrictive assumptions, including data heterogeneity, almost sure variance bounds, and symmetric noise distributions centered around their means.
	{In this chapter, we develop distributed error feedback algorithms that avoid the restrictive assumptions mentioned above and rely on standard conditions on the objective functions, compressors, and stochastic-oracle noise.}

	\section{Preliminaries}

	  \textbf{Notation.}
	We use \([n]\) to denote the set \(\{1, 2, \ldots, n\}\), and \(\Exp[u]\) to represent the expectation of a random variable \(u\). Additionally, \(\norm{\cdot}\) indicates the Euclidean norm for vectors or the spectral norm for matrices, and \(\norm{\cdot}_1\) is the $\ell_1$-norm for vectors, while \(\langle x, y \rangle\) denotes the inner product between \(x, y \in \mathbb{R}^d\). Lastly, for a square matrix \(A \in \mathbb{R}^{d \times d}\), \(\lambda_{\min}(A)\) refers to its minimum eigenvalue, and \(I \in \mathbb{R}^{d \times d}\) is the identity matrix.

    \textbf{Problem Formulation.}
	We focus on the following distributed optimization problem: 
	{
	\begin{equation}\label{ch6:eqn:Problem_chapter_6}
		\min_{x\in\R^d}\left\{f(x) \eqdef \frac{1}{n}\sum_{i=1}^n f_i(x)\right\}.
	\end{equation}
	}
	where $n$ refers to the number of clients, and $f_i(x)$ is the loss of a model parameterized by vector $x\in \R^d$ over its local data $\thesischsixcD_i$ owned by client $i\in [n]$.

	\textbf{Assumptions.}
	To facilitate our convergence analysis, we make standard assumptions on objective functions and compression operators.

	\begin{assumption}[Lower Boundedness of $f$]\label{ch6:assum:lowerbound_whole_f_chapter_6}
		The function $f$ is bounded from below, i.e.,
		 $$f^{\inf} = \inf_{x\in\R^d} f(x)>-\infty.$$ 
	\end{assumption}

	\begin{assumption}[Lower Boundedness of $f_i$]\label{ch6:assum:lowerbound_f_i_chapter_6}
		For each $i\in [n]$, the function $f_i$ is bounded from below, i.e., $$f^{\inf}_{i} \eqdef \inf_{x\in\R^d} f_i(x)>-\infty.$$
	\end{assumption}
	
	Assumptions~\ref{ch6:assum:lowerbound_whole_f_chapter_6} and~\ref{ch6:assum:lowerbound_f_i_chapter_6} are standard for analyzing optimization algorithms for unconstrained problems.

	\begin{assumption}[Generalized Smoothness of $f_i$]\label{ch6:assum:LzeroLoneSmooth_chapter_6}
		{A function $f_i(x)$ is symmetrically generalized smooth if there exist $L_0>0$ and $L_1\geq 0$ such that, for $u_\theta = \theta x + (1-\theta)y$ and all $x,y\in\R^d$,}
		\begin{eqnarray}\label{ch6:eqn:LzeroLoneSmooth_chapter_6}
			\norm{ \nabla f_i(x) - \nabla f_i(y) } \leq \left( L_0+L_1 \sup_{\theta\in[0,1]}\norm{\nabla f_i(u_\theta)} \right) \norm{x-y}.
		\end{eqnarray}	
	\end{assumption}

	Assumption~\ref{ch6:assum:LzeroLoneSmooth_chapter_6} refers to symmetric generalized smoothness  by~\textcite{chen2023generalized}, which covers asymmetric generalized smoothness~\parencite{koloskova2023revisiting,chen2023generalized}, and the original $(L_0,L_1)$-smoothness by~\textcite{zhang2019gradient}.
	Moreover, Assumption~\ref{ch6:assum:LzeroLoneSmooth_chapter_6} covers the functions with unbounded classical smoothness constant, e.g., exponential function. 
	Additionally, Assumption~\ref{ch6:assum:LzeroLoneSmooth_chapter_6} with $L_1=0$ reduces to the traditional $L_0$-smoothness~\parencite{nesterov2018lectures,beck2017first}, under which the convergence of optimization algorithms has been extensively studied.

	\begin{assumption}[Contractive Compressor]\label{ch6:assum:contractive_comp_chapter_6}
		{Let $\gH_k$ contain all randomness revealed before the compression performed at iteration $k$. An operator $\thesischsixcC^k:\R^d\rightarrow\R^d$ is an $\alpha$-contractive compressor if there exists $\alpha\in(0,1]$ such that, for every $k\geq0$ and every $\gH_k$-measurable vector $v\in\R^d$,}
		{
		\begin{equation} \label{ch6:eqn:contractive_comp_chapter_6}
			\Exp\!\left[\left.\sqnorm{\thesischsixcC^k(v)-v}\,\right|\,\gH_k\right] \leq (1-\alpha) \sqnorm{v}.
		\end{equation}
		}
	\end{assumption}
	
	Furthermore, compressors defined by Assumption~\ref{ch6:assum:contractive_comp_chapter_6} cover top-$k$ sparsifiers~\parencite{alistarh2018convergence,stich2018sparsified}, low-rank approximation~\parencite{vogels2019powersgd,safaryan2021fednl}, and various other compressors described by~\textcite{safaryan2022uncertainty,JMLR:v24:21-1548,demidovich2023guide}.
	
	\begin{assumption}[Bounded Variance]\label{ch6:assum:bounded_variance_chapter_6}
		{Let $\gF_k$ denote the history immediately before the stochastic samples at iteration $k$ are drawn. Each client draws a fresh sample $\xi_i^k\sim\thesischsixcD_i$, conditionally independently across clients given $\gF_k$. For every $\gF_k$-measurable $x\in\R^d$,}
			{
			\begin{equation}\label{ch6:eqn:bounded_variance_chapter_6}
				\begin{aligned}
				\Exp\!\left[\left.\nabla f_i(x;\xi_i^k)\,\right|\,\gF_k\right]
				&= \nabla f_i(x),\\
				\Exp\!\left[\left.
				\sqnorm{\nabla f_i(x;\xi_i^k)-\nabla f_i(x)}
				\,\right|\,\gF_k\right]
				&\leq \sigma^2.
				\end{aligned}
			\end{equation}
			}
	\end{assumption}

	{Assumption~\ref{ch6:assum:bounded_variance_chapter_6} is standard in stochastic optimization~\parencite{nemirovski2009robust,ghadimi2012optimal,ghadimi2013stochastic}. It controls only the noise of each local stochastic oracle and does not impose a uniform bound on the dissimilarity between local and global gradients.}

		{Throughout this chapter, every normalized quotient is
		written explicitly as $v/\norm{v}$.}\footnote{At $v=0$,
		the quotient $v/\norm{v}$ is defined to be the zero vector.}
		
		\begin{algorithm}
		\caption{Normalized Error Feedback (\normEFname{})}
		\label{ch6:alg:normalized_ef21_chapter_6}
		\begin{algorithmic}[1]
			\STATE {\textbf{Input:} Stepsizes $\gamma_k>0$ for $k\geq0$; starting point $x^0$ and memories $g_i^{-1}\in\R^d$, $i\in[n]$; set $g^{-1}=\frac1n\sum_{i=1}^n g_i^{-1}$; and $\alpha$-contractive compressors $\thesischsixcC^k$.}
			\FOR{each iteration $k = 0, 1, \dots, K$}
			\FOR{each client $i = 1, 2, \dots, n$ in parallel}
			\STATE Compute local gradient $\nabla f_i(x^k)$
			\STATE Transmit $\Delta_i^k = \thesischsixcC^k(\nabla f_i(x^k) - g_i^{k-1})$
			\STATE Update  $g_i^k = g_i^{k-1} +\Delta_i^k $
			\ENDFOR
			\STATE {Central server computes $g^{k}=g^{k-1}+\frac{1}{n}\sum_{i=1}^n\Delta_i^k=\frac{1}{n}\sum_{i=1}^n g_i^k$}
			\STATE {Central server updates $x^{k+1}=x^k-\gamma_k\frac{g^k}{\norm{g^k}}$}
			\ENDFOR\stepEightBdeletion
		\end{algorithmic}
		\end{algorithm}

	\section{Normalized Error Feedback (\normEFname{})}
	\label{ch6:sec:ef21_chapter_6}

	For nonconvex deterministic optimization under generalized smoothness,  we  develop a distributed error feedback algorithm. One challenge is that the generalized smoothness parameter scales with the gradient norm $\norm{\nabla f(x^k)}$. To resolve this issue, we apply gradient normalization to the algorithms. 
	In particular, we consider \normEFname{}, the normalized version of \algname{EF21}~\parencite{richtarik2021ef21} that updates the next iterates $x^{k+1}$ using the \normEFname{} update. 
	The full description of \normEFname{} can be found in Algorithm~\ref{ch6:alg:normalized_ef21_chapter_6}.

	Our new method \normEFname{}, just like \algname{EF21}~\parencite{richtarik2021ef21} under traditional smoothness, enjoys the $\cO(1/\sqrt{K})$ convergence in the gradient norm under generalized smoothness, as shown below.
	
	\begin{theorem}[Convergence of \normEFname{}]\label{ch6:thm:ef21_chapter_6}
		Consider Problem~\eqref{ch6:eqn:Problem_chapter_6}, where Assumption~\ref{ch6:assum:lowerbound_whole_f_chapter_6} (lower bound on $f$), Assumption~\ref{ch6:assum:lowerbound_f_i_chapter_6} (lower bound on $f_i$), Assumption~\ref{ch6:assum:LzeroLoneSmooth_chapter_6} (generalized smoothness of $f_i$), and Assumption~\ref{ch6:assum:contractive_comp_chapter_6} (contractive compressor) hold. {Assume in addition that $L_1>0$.}
		Then, the iterates $\{x^k\}$ generated by \normEFname{} (Algorithm~\ref{ch6:alg:normalized_ef21_chapter_6}) with 
		{
		\begin{equation*}
			\gamma_k\equiv\gamma_K=\frac{\gamma_0}{\sqrt{K+1}}
		\end{equation*}
		}
		{for $K\geq0$ and $\gamma_0>0$ satisfy}
		{
		\begin{equation*}
			\underset{k=0,1,\ldots,K}{\min} \Exp\!\left[\norm{\nabla f(x^k)}\right]
			\leq  \frac{\Exp{V^0}\exp\!\left(8c_1L_1\exp(L_1\gamma_K)\gamma_0^2\right)}{\gamma_0\sqrt{K+1}}
			+B\frac{\gamma_0\exp(L_1\gamma_K)}{\sqrt{K+1}}.
		\end{equation*}
		}
		{where
		\begin{align*}
			V^k &\eqdef f(x^k)-f^{\inf}
			+\frac{2\gamma_K}{1-\sqrt{1-\alpha}}\,
			\frac{1}{n}\sum_{i=1}^n \norm{\nabla f_i(x^k)-g_i^k},\\
			B &\eqdef 2c_0+\frac{8L_1c_1}{n}
			\sum_{i=1}^n(f^{\inf}-f_i^{\inf}),\\
			c_i&\eqdef\left(\frac{1}{2}
			+2\frac{\sqrt{1-\alpha}}{1-\sqrt{1-\alpha}}\right)L_i,
			\quad i\in\{0,1\}.
		\end{align*}}
	\end{theorem}

	Theorem~\ref{ch6:thm:ef21_chapter_6} establishes the $\cO(1/\sqrt{K})$ convergence in the expectation of gradient norms for \normEFname{} on nonconvex deterministic problems under generalized smoothness.
	This rate is the same as Theorem 1 of \textcite{richtarik2021ef21} for \algname{EF21} under traditional smoothness, and does not depend on data heterogeneity conditions in contrast to \textcite{crawshaw2024federated,liu2022communication}. 
	{The stepsize depends on an arbitrary positive base constant $\gamma_0$ and the horizon $K$, but does not require the smoothness constants $L_0,L_1$, unlike the rule in \textcite{richtarik2021ef21}.}
	{For $L_1>0$, choosing $\gamma_0 = 1/(8c L_1)$ gives the following specialization of Theorem~\ref{ch6:thm:ef21_chapter_6}:}
	\begin{eqnarray*}
		\underset{k=0,1,\ldots,K}{\min} \Exp\!\left[\norm{\nabla f(x^k)}\right]
		&\leq&  \frac{32c L_1 V^0 +L_0/L_1 + 2L_1 \delta^{\inf} }{\sqrt{K+1}}, 
	\end{eqnarray*}
	where $c=\frac{1}{2} + 2\frac{\sqrt{1-\alpha}}{1-\sqrt{1-\alpha}}$, and $\delta^{\inf} = \frac{1}{n}\sum_{i=1}^n (f^{\inf} -f_i^{\inf} )$.

	\textbf{Comparisons between \normEFname{} and \algname{EF21} under traditional smoothness.}
	{For the separate $L$-smooth specialization ($L_1=0$), a direct application of the classical smoothness descent inequality shows that \normEFname{} achieves the same $\cO(1/\sqrt{K})$ rate in the expected gradient norm as \algname{EF21} analyzed by~\textcite{richtarik2021ef21}, with a displayed constant at most twice as large. This specialization is proved directly because the generalized-smoothness proof of Theorem~\ref{ch6:thm:ef21_chapter_6} divides by $L_1$.}
		We refer to the detailed derivation and discussion in Appendix~\ref{ch6:sec:diss_thm_ef21_chapter_6}.

	In the following section, we demonstrate how to integrate normalization into \algname{EF21-SGDM}~\parencite{fatkhullin2024momentum}, an error feedback algorithm that allows each node to compute its local stochastic gradient, for solving nonconvex stochastic problems.

	\begin{algorithm}
		\caption{Normalized Error Feedback with Stochastic Gradients \& Momentum (\normEFname{SGDM})}
		\label{ch6:alg:normalized_ef21_sgdm_chapter_6}
		\begin{algorithmic}[1]
			\STATE {\textbf{Input:} Stepsizes $\gamma_k>0$ for $k\geq0$ and momentum parameters $\eta_k\in(0,1]$ for $k\geq1$; starting point $x^0$ and memories $g_i^{-1}=0$, with $g^{-1}=0$; $\alpha$-contractive compressors $\thesischsixcC^k$.}
			\STATE {Each client draws $\xi_i^0$ independently and sets $v_i^0=\nabla f_i(x^0;\xi_i^0)$.}
			\FOR{each iteration $k = 0, 1, \dots, K$}
			\FOR{each client $i = 1, 2, \dots, n$ in parallel}
			\IF{$k\geq1$}
			\STATE {Draw $\xi_i^k$ and set $v_i^k=(1-\eta_k)v_i^{k-1}+\eta_k\nabla f_i(x^k;\xi_i^k)$.}
			\ENDIF
			\STATE Transmit $\Delta_i^k = \thesischsixcC^k( v_i^k - g_i^{k-1})$
			\STATE Update  $g_i^k = g_i^{k-1} +\Delta_i^k $
			\ENDFOR
			\STATE {Central server computes $g^{k}=g^{k-1}+\frac{1}{n}\sum_{i=1}^n\Delta_i^k=\frac{1}{n}\sum_{i=1}^n g_i^k$}
			\STATE {Central server updates $x^{k+1}=x^k-\gamma_k\frac{g^k}{\norm{g^k}}$}
			\ENDFOR\stepEightBdeletion
		\end{algorithmic}
		\end{algorithm}

	\section{Normalized Error Feedback with Stochastic Gradients \& Momentum (\normEFname{SGDM})}
	\label{ch6:sec:ef21_sgdm_chapter_6}

	{Having established convergence of \normEFname{} for deterministic optimization, we next develop a distributed error feedback algorithm that incorporates stochastic gradients and normalization under generalized smoothness.}
		In particular, we focus on \normEFname{SGDM}  (Algorithm~\ref{ch6:alg:normalized_ef21_sgdm_chapter_6}), the normalized version of \algname{EF21-SGDM} due to \textcite{fatkhullin2024momentum}. We also note that  \normEFname{SGDM}  recovers many optimization algorithms of interest in the special cases. For instance, it  reduces to   
	\begin{itemize}
		\item  normalized version of \algname{EF21}~\parencite{richtarik2021ef21}, which we call \normEFname{}, when we let $\eta_k=1$ and $\nabla f_i(x^k;\xi_i^k)=\nabla f_i(x^k)$, 
		\item normalized version of \algname{EF21-SGD}~\parencite{fatkhullin2021ef21}, which we call \normEFname{SGD}, when we let $\eta_k=1$, and 
		\item  normalized version of \algname{SGDM}~\parencite{cutkosky2020momentum}, which we call \normalgname{SGDM}\footnote{This method is also known as \algname{NSGD-M}.}, when we let $\eta_k=1-\beta_k$ and $\thesischsixcC^k(\cdot)$ is the identity compressor/mapping. 
	\end{itemize}

	In the next theorem, we demonstrate that \normEFname{SGDM}  attains the same \(\mathcal{O}(1/K^{1/4})\) convergence rate as both \algname{EF21-SGDM}  and \normalgname{SGDM}.

	\begin{theorem}[Convergence of \normEFname{SGDM}]\label{ch6:thm:ef21_sgdm_chapter_6}
		Consider Problem~\eqref{ch6:eqn:Problem_chapter_6}, where Assumption~\ref{ch6:assum:lowerbound_whole_f_chapter_6} (lower bound on $f$), Assumption~\ref{ch6:assum:lowerbound_f_i_chapter_6} (lower bound on $f_i$), Assumption~\ref{ch6:assum:LzeroLoneSmooth_chapter_6} (generalized smoothness of $f_i$), Assumption~\ref{ch6:assum:contractive_comp_chapter_6} (contractive compressor), and Assumption~\ref{ch6:assum:bounded_variance_chapter_6} (bounded variance) hold.
		{Suppose $L_1>0$, $g_i^{-1}=0$ for all $i$, and}
		{
		\begin{eqnarray*}
			\gamma_k\equiv \gamma_K &=& \frac{\gamma_0}{(K+1)^{\nicefrac{3}{4}}},\text{ with } 0 < \gamma_0 \leq \frac{1}{16L_1}\min\left\{(K+1)^{\nicefrac{1}{2}}C_\alpha, 1\right\}, \quad \text{and} \\
			\quad \eta_k \equiv \eta &=& \frac{1}{(K+1)^{\nicefrac{1}{2}}},
		\end{eqnarray*}
		}
		where $C_\alpha \eqdef 1 - \sqrt{1-\alpha}$, then the iterates $\{x^k\}$ generated by \normEFname{SGDM}  (Algorithm~\ref{ch6:alg:normalized_ef21_sgdm_chapter_6})
		satisfy for $K \geq 0$ 
			{
			\begin{equation*}
			\begin{aligned}
				\min_{0\leq k\leq K}\Exp\!\left[\norm{\nabla f(x^k)}\right]
				\leq {}&
				\gO\!\left(
				\frac{\frac{\delta^0}{\gamma_0}+\frac{\sigma}{\sqrt n}
				+\gamma_0\!\left(L_0+L_1^2\delta^{\inf}\right)}
				{(K+1)^{1/4}}\right)\\
				&+\gO\!\left(\frac{\sqrt{1-\alpha}}{\alpha}
				\left[
				\frac{L_1(\delta^0+\delta^{\inf})+\frac{L_0}{L_1}+\sigma}{K+1}
				+\frac{\sigma}{(K+1)^{1/2}}\right.\right.\\
				&\hspace{47mm}\left.\left.
				+\frac{\gamma_0\!\left(L_0+L_1^2\delta^{\inf}\right)}
				{(K+1)^{3/4}}
				\right]\right).
			\end{aligned}
			\end{equation*}}
		where $\delta^0 \eqdef f(x^0) -f^{\inf}$, and $\delta^{\inf} \eqdef \frac{1}{n}\sum^n_{i=1}(f^{\inf} - f^{\inf}_i)$.
	\end{theorem}
	
	From Theorem~\ref{ch6:thm:ef21_sgdm_chapter_6},  \normEFname{SGDM}  under generalized smoothness achieves the  $\mathcal{O}(1/K^{1/4})$ convergence rate in the expectation of gradient norms. 
	This rate is the same as that of \algname{EF21-SGDM}, previously analyzed under traditional smoothness by \textcite[Theorem 3]{fatkhullin2024momentum}. 
	{The result does not assume a uniform bound on data heterogeneity. Algorithm~\ref{ch6:alg:normalized_ef21_sgdm_chapter_6} uses one fresh stochastic-gradient sample per client and iteration, so no growing mini-batch is required. The quantity $\delta^{\inf}$ records the gap between the available local and global lower bounds; it is not a bounded-gradient-dissimilarity assumption.}
	We also notice that the stepsize \(\gamma_0\) for \normEFname{SGDM}, unlike in the case of \normEFname{}, depends on the generalized smoothness constant \(L_1\), and the compression parameter \(\alpha\). However, the considered choice of stepsizes is agnostic to $\sigma$ and $L_0$.

	Furthermore, Theorem~\ref{ch6:thm:ef21_sgdm_chapter_6} with \(\alpha = 1\) (i.e., \(\mathcal{C}^k\) is the identity compressor) implies the convergence bound of the distributed version of normalized SGD with momentum (\normalgname{SGDM}) \parencite{cutkosky2020momentum} using $\beta=1-\eta$: 
	{
	\begin{eqnarray}\label{ch6:eqn:result_ef21_sgdm_no_comp_chapter_6}
		\underset{k=0,1,\ldots,K}{\min} \Exp\!\left[\norm{\nabla f(x^k)}\right]
		&\leq& \gO\left(\frac{\nicefrac{(f(x^0)-f^{\inf})}{\gamma_0}+\nicefrac{\sigma}{\sqrt{n}} +\gamma_0 L_0 + \gamma_0 L_1^2\delta^{\inf} }{(K+1)^{\nicefrac{1}{4}}}\right).
	\end{eqnarray} 
	For the single-node \algname{SGDM}, where \(n=1\) and
	\(\delta^{\inf}=0\), choosing
	$\gamma_0=\Theta(\nicefrac{1}{L_1})$ gives
	$\gO\!\left((L_1(f(x^0)-f^{\inf})+\sigma+L_0/L_1)/(K+1)^{1/4}\right)$,
	which matches the rate obtained by
	\textcite[Corollary~3]{hubler2024parameter}. Unlike that single-node
	result, our result also covers multiple nodes. In this case,
	$\sigma/\sqrt n$ displays the variance reduction across clients, while
	$\gamma_0L_1^2\delta^{\inf}$ records the gap between the local and global
	lower bounds.}

	\textbf{Novel proof techniques for \normEFname{} and \normEFname{SGDM} under generalized smoothness.} Our analysis demonstrates that \normEFname{}
	 achieves the convergence rate under generalized smoothness equivalent to \algname{EF21} under traditional smoothness. However, our proof techniques differ significantly from prior work. 
	 We employ different Lyapunov functions. For \normEFname{}, we use 
     $V^k := f(x^k)-f^{\inf} + \frac{A}{n}\sum_{i=1}^n \norm{\nabla f_i(x^k)-g_i^k}$, 
	     {in contrast to \textcite{richtarik2021ef21}, which uses $V^k := f(x^k)-f^{\inf} + \frac{B}{n}\sum_{i=1}^n \sqnorm{\nabla f_i(x^k)-g_i^k}$.}
     For \normEFname{SGDM}, we use $V^k := f(x^k)-f^{\inf} + \frac{C}{n}\sum_{i=1}^n \norm{v_i^k - g_i^k} + \frac{D}{n}\sum_{i=1}^n \norm{v_i^k - \nabla f_i(x^k)}$, unlike \textcite{fatkhullin2024momentum} that uses $V^k := f(x^k)-f^{\inf} + \frac{E}{n}\sum_{i=1}^n \sqnorm{v_i^k - g_i^k} + \frac{F}{n}\sum_{i=1}^n \sqnorm{v_i^k - \nabla f_i(x^k)}$.  
	These new Lyapunov functions necessitate the Lyapunov-based convergence analysis, distinct from standard techniques for error feedback methods. 
	Our analysis leverages Lemma~\ref{ch6:lemma:bound_avg_norm_grad_multinodeef21_chapter_6} to handle generalized smoothness.
	For \normEFname{}, we rely on Lemma~\ref{ch6:lemma:conv_v2_chapter_6}. 
	For \normEFname{SGDM}, we derive a new upper-bound on $\Exp\!\left[\norm{v^k-\nabla f(x^k)}\right]$, unlike \textcite{fatkhullin2024momentum} to show the $\sqrt{n}$-speedup for the term proportional to $\sigma$, and utilize non-uniform weights to obtain convergence in the gradient norm.

	\section{Experiments}

	In this section, we evaluate the performance of \normEFname{}, and compare it against \algname{EF21}~\parencite{richtarik2021ef21}. We test these algorithms for three nonconvex, generalized smooth problems: the problem of minimizing  polynomial functions, the logistic regression problem with a nonconvex regularization term over synthetic  and benchmark datasets from LIBSVM~\parencite{chang2011libsvm}, and the training of the ResNet-20~\parencite{7780459} model over the CIFAR10~\parencite{Krizhevsky2009LearningML} dataset\footnote{We implemented \algnamesmall{EF21} and  \normEFname{} on training the ResNet-20 model by using PyTorch. Our source codes can be found in \href{https://anonymous.4open.science/r/error-feedback-generalized-smoothness-9DDC}{the link to error-feedback-generalized-smoothness-paper}.}. {To avoid confusing the sparsity budget with the iteration index $k$, we denote it by $s$: all experiments use a top-$s$ sparsifier, which is an $s/d$-contractive compressor.}

	\subsection{Logistic Regression with a Nonconvex Regularizer}

	First, we consider a logistic regression problem with a nonconvex regularizer, i.e., Problem~\eqref{ch6:eqn:Problem_chapter_6} with 
    \begin{align*}
    f_i(x)=\log(1+\exp(-b_i a_i^T x)) + \lambda \sum_{j=1}^d \frac{x_j^2}{1+x_j^2},
    \end{align*}
	where $a_i \in \R^d$ is the $i^{\rm th}$ feature vector of data matrix $A\in\R^{n\times d}$ with its class label $b_i \in \{-1,1\}$, and $\lambda>0$ is a regularization parameter. Here, $f(x)$ is nonconvex, and $L$-smooth with $L=\sqnorm{A}/(4n) +2\lambda$. Also, each $f_i(x)$ is $\hat L_i$-smooth with $\hat L_i = \sqnorm{a_i}/4+2\lambda$, and generalized smooth with $L_0=2\lambda+\lambda \sqrt{d} \max_i \norm{a_i}$ and $L_1 = \max_i \norm{a_i}$. The derivations of smoothness parameters can be found in Appendix~\ref{ch6:app:smooth_logistic_chapter_6}.

	In this experiment, we initialized \( x^0 \in \R^d \), where each coordinate was drawn from a standard normal distribution \(\mathcal{N}(0,1)\), and set \(\lambda = 0.1\). %
	Here, the condition $\lambda >  \lambda_{\min}\left(A^{\top} A\right)/(2n)$  ensures that $f(x)$ is nonconvex.
	We ran  \normEFname{} and \algname{EF21} on the following datasets: (1) two from LIBSVM~\parencite{chang2011libsvm}: \texttt{Breast Cancer} (\(n = 683\), \(d = 10\), and scaled to $[-1,1]$), and \texttt{a1a} (\(n = 1605\), \(d = 123\)); and (2) a synthetically generated dataset (\(n = 20\), \(d = 10\)), where the data matrix \( A \in \R^{n \times d} \) had entries drawn from \(\mathcal{N}(0,1)\), and the class label \( b_i \) was set to either \(-1\) or \(1\) with equal probability.
	For \algname{EF21}, we selected the stepsize $\gamma_k  = {1}/{\left( L + \tilde L \sqrt{{\beta}/{\theta}} \right)}$ with $\tilde L = \sqrt{ \sum_{i=1}^n  \hat L_i^2/n}$, $\theta=1-\sqrt{1-\alpha}$, and $\beta={(1-\alpha)}/{(1-\sqrt{1-\alpha})}$,  given by \textcite[Theorem 1]{richtarik2021ef21}.
	For \normEFname{}, we chose  $\gamma_k= \gamma/\sqrt{K+1}$ with $\gamma>0$ from Theorem~\ref{ch6:thm:ef21_chapter_6}, by setting $\gamma_0=1$, $K = 100$ for the generated data and \texttt{Breast Cancer}, and $K=400$ for \texttt{a1a}.
	We choose $\gamma_0=1$, because \normEFname{} with $\gamma_0 \in [1,10]$ converges faster than that with small values of $\gamma_0$ (e.g. $0.1$), when we run the algorithm on a single node ($n=1$) for minimizing polynomial function and solving logistic regression.
	We determine $K$ as the smallest number of iterations required to achieve the desired accuracy by performing a grid search with a stepsize of $50$.

	Figure~\ref{ch6:fig:exp2_chapter_6} shows that \normEFname{} reaches smaller gradient norms than \algname{EF21} in fewer iterations on all three datasets. The experiment compares the two methods under the stepsizes given by their theory: \algname{EF21} uses the stepsize of \textcite[Theorem 1]{richtarik2021ef21}, and \normEFname{} uses $\gamma_k=\gamma_0/\sqrt{K+1}$ from Theorem~\ref{ch6:thm:ef21_chapter_6}, whose free parameters $\gamma_0$ and $K$ were set as described above. Neither method was tuned beyond these theoretical rules, so the comparison reflects the larger stepsize allowed by Theorem~\ref{ch6:thm:ef21_chapter_6}, not a comparison of separately tuned methods.

	\begin{figure}[h]
		\centering
			\includegraphics[width=0.3\textwidth]{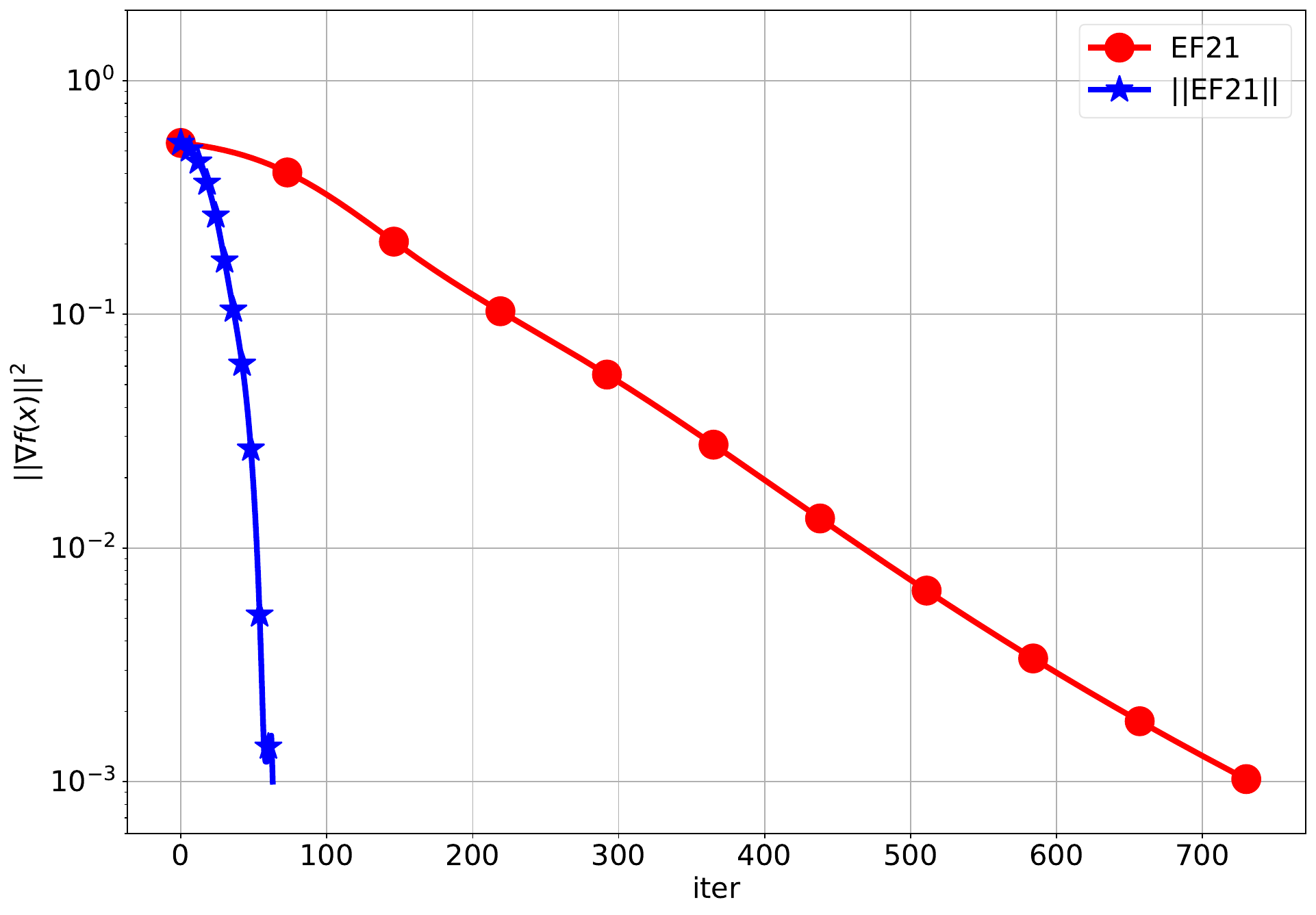}
		\hfill
			\includegraphics[width=0.3\textwidth]{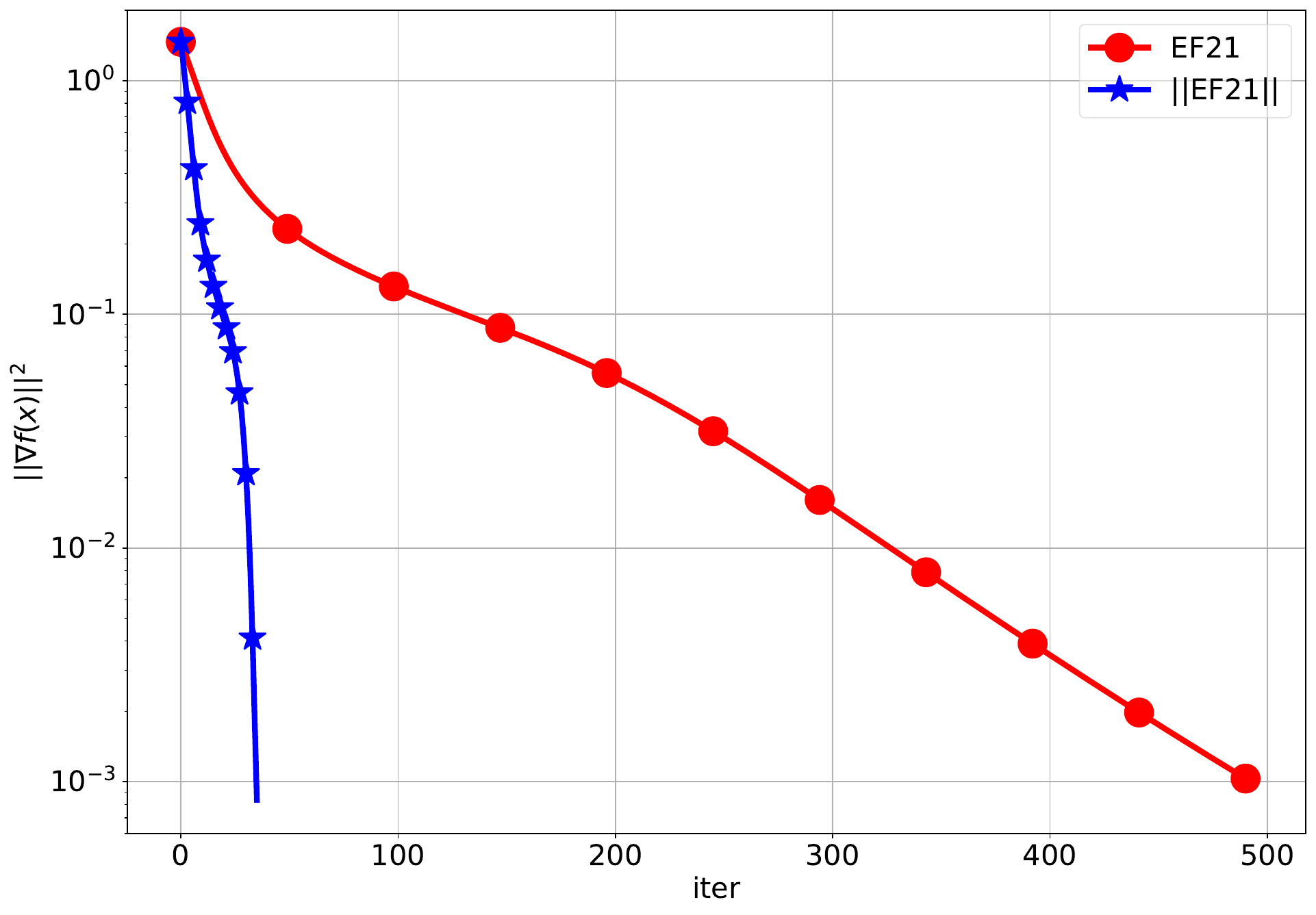}
		\hfill
			\includegraphics[width=0.3\textwidth]{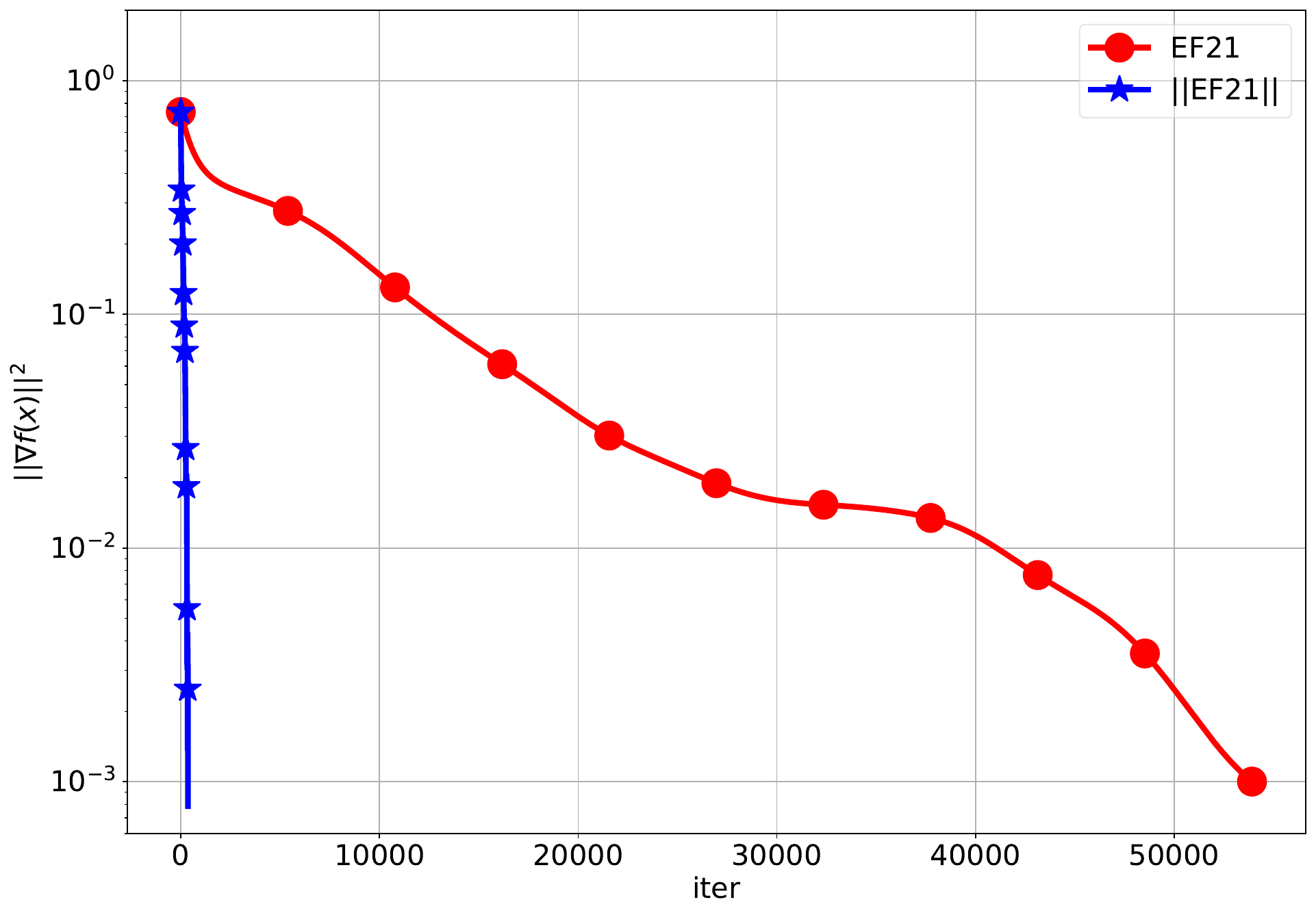}
			\caption{Logistic regression with a nonconvex regularizer using normalized \normEFname{} and  \algname{EF21}. We reported $\sqnorm{\nabla f(x^k)}$ with respect to iteration count $k$. We used the constant stepsize $\gamma = \frac{1}{L + \tilde L \sqrt{\frac{\beta}{\theta}}}$ for \algname{EF21}, and $\gamma = \frac{\hat \gamma}{\sqrt{K+1}}$, $\hat\gamma = 1$ for \normEFname{}. Here, $K=100$ for our generated data (left), and \texttt{Breast Cancer} (middle), while $K=400$ for \texttt{a1a}  (right). }
			\label{ch6:fig:exp2_chapter_6}
		\end{figure}

		\subsection{ResNet20 Training Over CIFAR-10}\label{ch6:subsec:resnet_20_chapter_6}

		Next, we trained the ResNet20~\parencite{7780459} model on the CIFAR-10~\parencite{Krizhevsky2009LearningML} dataset, which was demonstrated  empirically by \textcite{zhang2019gradient}  to satisfy the $(L_0,L_1)$-smoothness condition. 
		{In these experiments, we used a top-$s$ compressor over $50,000$ training images, with evaluation on $10,000$ test images.}
		The dataset was evenly distributed  among $5$ clients, each using a mini-batch size of $128$. Both algorithms were run for $100$ epochs with a constant stepsize $\gamma=5$. They use the same random seed, data split and initialization, and each configuration was run once; run-to-run variability was therefore not assessed. Here, one epoch refers to a full pass through the entire dataset processed by all clients.

		From Figure%
		~\ref{ch6:fig:plot_exp_DL_chapter_6}, {under the same constant stepsize and a top-$s$ sparsifier with $s=0.1d$}, \normEFname{} outperforms \algname{EF21}, in terms of  convergence speed (in gradient norms and losses) and accuracy, relative to the number of bits communicated from each client to the server. Specifically, \normEFname{} achieved accuracy gains of up to $10$\% over \algname{EF21}.

		\begin{figure}[h]
			\centering
			\includegraphics[width=\textwidth]{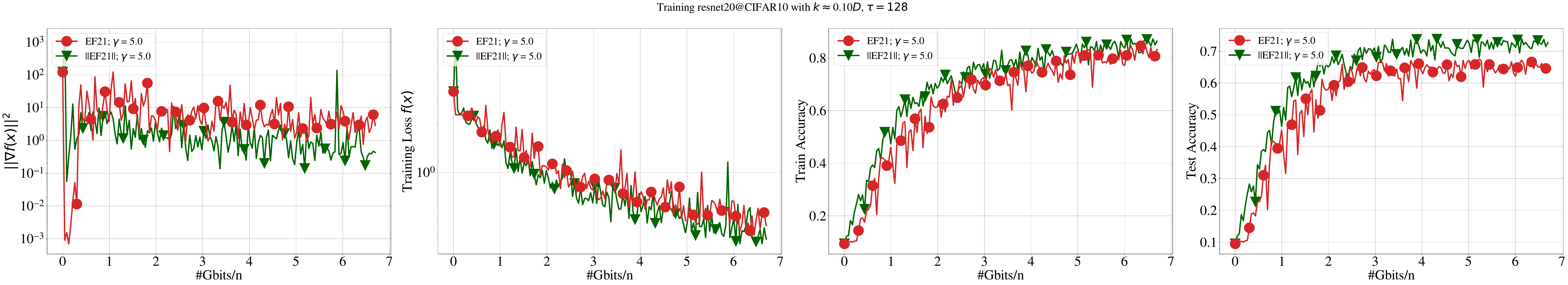}
			\caption{{ResNet20 training on CIFAR-10 using \algname{EF21} and \normEFname{} with the same stepsize $\gamma=5$ and a top-$s$ sparsifier with $s=0.1d$.}}
			\label{ch6:fig:plot_exp_DL_chapter_6}
		\end{figure}

		\section{Conclusion and Future Work}

		This chapter considered synchronous distributed nonconvex optimization with a fixed worker population, generalized smooth local objectives, and contractive, possibly biased compression. Normalization controls the server displacement, while \algname{EF21} reconstruction returns information omitted by earlier compressed messages. The stochastic variant additionally uses Polyak momentum to stabilize the local gradient estimates.

		The deterministic \normEFname{} theorem controls the analytical best-iterate criterion \(\min_{0\leq k\leq K}\Exp[\norm{\nabla f(x^k)}]\) at order \(\mathcal{O}(K^{-1/2})\). The stochastic \normEFname{SGDM} theorem controls the same in-expectation best-iterate criterion at order \(\mathcal{O}(K^{-1/4})\), together with faster compression-dependent terms. These are neither high-probability nor terminal-iterate guarantees. The rates match the established orders for the corresponding \algname{EF21} methods under classical smoothness, but this chapter does not prove a matching lower bound. The analysis also avoids a uniform gradient-dissimilarity assumption; the quantity \(\delta^{\inf}\) records the gap between available local and global lower bounds instead.

		The deterministic horizon-based rule does not require \(L_0\) or \(L_1\), whereas the stochastic base stepsize remains constrained by \(L_1\) and the compressor parameter. The experiments indicate, rather than prove, that normalization can permit useful larger stepsizes. Extensions to decreasing or adaptive schedules, partial participation, asynchrony, and broader federated models remain open. Chapter~\ref{ch:chapter_7} develops decreasing parameter-agnostic schedules and compares a wider family of first- and second-order momentum estimators within normalized error feedback.

  \endgroup

\chapter{Improved Convergence in Parameter-Agnostic Error Feedback through Momentum}
\label{ch:chapter_7}
\thesischaptersourcenote
  {Improved Convergence in Parameter-Agnostic Error Feedback through Momentum}
  {sadiev2025improved}
  \begingroup
  \thesischaptersetup{7}%
  
\begin{paperabstract}
Communication compression is essential for scalable distributed training of modern machine learning models, but it often degrades convergence due to the noise it introduces. 
Error Feedback (EF) mechanisms are widely adopted to mitigate this issue of distributed compression algorithms. 
Despite their popularity and training efficiency, existing distributed EF algorithms often require prior knowledge of problem parameters (e.g., smoothness constants) to fine-tune stepsizes. 
This limits their practical applicability especially in large-scale  neural network training. 
In this chapter, we study normalized error feedback algorithms that combine EF with normalized updates, various momentum variants, and parameter-agnostic, time-varying stepsizes.
Here, parameter-agnostic means that the schedules do not require prior knowledge of the smoothness constants or the initial optimality gap; the methods still contain structural parameters and user-selected settings and are not tuning-free.
We analyze the convergence of these algorithms for minimizing smooth functions, and establish parameter-agnostic complexity bounds that are close to the best-known bounds with carefully-tuned problem-dependent stepsizes. 
Specifically, for the expected gradient norm at a stepsize-weighted random iterate, normalized \algname{EF21} methods achieve rates of $\widetilde{\cO}(1/K^{1/4})$ for Polyak's heavy-ball momentum, $\widetilde{\cO}(1/K^{2/7})$ for Implicit Gradient Transport (IGT), and $\widetilde{\cO}(1/K^{1/3})$ for STORM and Hessian-corrected momentum, where logarithmic factors are suppressed.
Our results hold with decreasing stepsizes and small mini-batches.
Finally, our empirical experiments are consistent with the theoretical predictions.

\end{paperabstract}

\section{Introduction}

Distributed optimization has become essential for efficiently training modern machine learning models on large-scale datasets. 
This shift is driven by the increasing size of models—such as deep neural networks with billions of parameters—and the increasing volume of training data~\parencite{Kolesnikov2019BigT,LMFewShot}. In distributed settings, multiple clients collaborate in parallel, sharing local information (e.g., stochastic gradients) with a central server to jointly minimize the average of their objective functions, thereby keeping training computationally feasible. However, a major challenge in this paradigm is the communication bottleneck, which becomes particularly severe with large models. For instance, transmitting the VGG-16 model~\parencite{simonyan2014very}, which contains 138.34 million parameters, requires over 500 MB of data per exchange—placing a significant strain on network resources during distributed stochastic gradient descent (SGD), which aims to find an $\varepsilon$-approximate critical point $x^k$ such that $\Exp{\norm{\nabla f(x^k)}} \leq \varepsilon$.

A popular strategy to reduce communication overhead is  compression, where clients apply compression operators to their local gradients before sending them to the server. Contractive but potentially biased compressors—including many widely used compressors  from sparsification~\parencite{alistarh2018convergence,StitchSparseSGD} to quantization~\parencite{wen2017terngrad,alistarh2017qsgd,beznosikov2023biased,horvoth2022natural}—have been shown to have favorable theoretical properties~\parencite{StitchSparseSGD,karimireddy2019error,EF21,fatkhullin2025ef21}, and to outperform randomized (unbiased) compressors in practice for distributed gradient-based algorithms~\parencite{lin2017deep,sun2019sparse,vogels2019powersgd}. However, naively aggregating these biased compressed gradients from the clients in general does not converge~\parencite{khirirat2020compressed} or even  diverges~\parencite{karimireddy2019error,beznosikov2023biased,bao2025efskip}.

To improve the convergence stability of distributed compression methods, error feedback (EF) mechanisms have been proposed. A widely studied EF variant, \algname{EF14}~\parencite{Seide20141bitSG}, has been extensively explored in both centralized~\parencite{StitchSparseSGD,karimireddy2019error} and distributed optimization settings~\parencite{alistarh2018convergence,tang2019doublesqueeze,tang2019texttt,gorbunov2020linearly,khirirat2020compressed,qian2021error}. However, many of these works rely on the assumption that the norms of the stochastic gradients are uniformly bounded, which restricts the class of objective functions the EF algorithms can solve.
Under this assumption, distributed EF algorithms have been shown to minimize smooth nonconvex functions with a convergence rate of $\cO(1/K^{1/3})$ in the gradient norm, where $K$ denotes the total number of iterations~\parencite{koloskova2019decentralized}.

To further improve theoretical convergence guarantees, many novel EF variants have been developed, including \algname{EControl}~\parencite{gao2024econtrol} and \algname{EF21}~\parencite{EF21}. Notably, \algname{EF21} achieves a convergence rate of $\cO(1/\sqrt{K})$ in the gradient norm for deterministic optimization problems—matching the rate of classical gradient descent. This method has also been extended to broader settings. 
One extension is \algname{EF21-P}~\parencite{EF21-P}, which adapts \algname{EF21} for federated optimization, where model parameters (rather than gradients) are exchanged between clients and the server. 
Another extension is \algname{EF21-SGDM}~\parencite{fatkhullin2023momentum,khirirat2024errorfeedbackl0l1smoothnessnormalization}, which incorporates Polyak momentum into the original \algname{EF21} algorithms to handle stochastic optimization, where clients compute local stochastic gradients. \algname{EF21-SGDM} achieves the convergence rate of $\cO(1/{K}^{1/4})$. 
To further accelerate the convergence, \algname{EF21-MVR}~\parencite{fatkhullin2023momentum} integrates the momentum-with-variance-reduction (MVR) update from STORM~\parencite{cutkosky2019momentum}, and enjoys the improved convergence rate of $\cO(1/K^{1/3})$.

The smoothness assumption of objective functions plays a key role in fine-tuning stepsizes  to ensure significant convergence improvements in distributed EF algorithms. However, stepsize rules that depend on smoothness are often impractical to implement. For example, estimating smoothness constants is typically infeasible, especially in deep neural network training. This motivates the need for \emph{parameter-agnostic} stepsize rules that do not depend on the smoothness constants. Crucially, it is important to design these stepsize rules so that they still achieve \emph{near-optimal} convergence, implying the performance of distributed EF algorithms with parameter-agnostic stepsizes  almost matches that of EF algorithms with optimally tuned, problem-dependent stepsizes.

One  approach for incorporating parameter-agnostic stepsizes into gradient-based algorithms—while still achieving near-optimal convergence—is through \emph{normalization}. For example, normalized stochastic momentum methods~\parencite{cutkosky2020momentum} have been shown by~\textcite{hubler2024parameter} to attain a near-optimal convergence rate of $\tilde{\cO}(1/{K}^{1/4})$ using parameter-agnostic, decreasing stepsize rules. However, to the best of our knowledge, this normalization approach has been so far limited to centralized algorithms.

\section{Contributions}

We summarize our  key contributions as follows:

\renewcommand\labelitemi{\ding{117}}
\begin{itemize}
    \item \textbf{Distributed \algname{EF21} algorithms with parameter-agnostic stepsizes and five momentum variants.} 
    In~Section~\ref{ch7:sec:ef21_momentum_chapter_7}, 
    we propose distributed \algname{EF21} algorithms that exploit normalization and momentum for solving stochastic optimization. Specifically, our algorithms 
    employ normalization, which enables stepsize schedules that do not require prior knowledge of the smoothness constant $L$ or the initial optimality gap $f(x^0) - f^{\inf}$. These schedules retain structural parameters and user-selected settings and therefore are not tuning-free.
    Furthermore, our algorithms leverage {\color{black} five variants of momentum updates widely adopted in centralized stochastic gradient algorithms: (1) Polyak momentum, (2) IGT momentum~\parencite{cutkosky2020momentum}, (3) MVR momentum~\parencite{cutkosky2019momentum}, (4) two second-order momentum schemes by~\textcite{salehkaleybar2022momentum,tran2022better}.

    } 
    
    \item \textbf{Near-optimal convergence for non-convex, smooth functions.} 
    In~Section~\ref{ch7:sec:convergence_chapter_7}, we prove that our proposed algorithms, using parameter-agnostic and decreasing stepsizes, achieve near-optimal convergence rates for the expected gradient norm at a stepsize-weighted random iterate when minimizing non-convex smooth functions.
    Our results match—up to logarithmic factors—the convergence guarantees of momentum-based \algname{EF21} algorithms from prior work, which typically rely on problem-dependent stepsizes.
     They also align with the rates achieved by centralized stochastic methods using analogous momentum variants. 
     Furthermore, our algorithms only require a batch size of $\cO\thesischsevenbr{1}$, in contrast to existing EF algorithms such as \textcite{fatkhullin2023momentum}, which often rely on problem-dependent stepsizes and larger mini-batch sizes.
    A summary of theoretical comparisons between our results and existing analyses in provided in~Table~\ref{ch7:table:complexities-ef_chapter_7}.

    \begin{table*}
	\centering
		\caption[Comparison of distributed error-feedback methods]{A comparison of distributed error feedback methods using contractive compressors for stochastic optimization under data heterogeneity. In the table, \textbf{PA} indicates whether the methods use parameter-agnostic stepsizes, \textbf{Stepsize} and \textbf{MP} list compact forms of the theoretical stepsize and momentum-parameter selections, \textbf{SO} denotes the use of second-order information of the functions, and \textbf{Complexity} refers to the sample complexity required to make the expected gradient norm of the theorem-covered output at most $\varepsilon$ for some $\varepsilon>0$. For the methods proposed in this chapter, that output is the stepsize-weighted random iterate $\tilde{x}^K$, not necessarily the final iterate $x^K$.\textsuperscript{\color{blue}(1)} } 
	\label{ch7:table:complexities-ef_chapter_7}
	\begin{threeparttable}
    \resizebox{\linewidth}{!}{%
		\begin{tabular}{cccccc}
			{\bf Method} & \bf PA & \bf Stepsize & \bf MP & \bf SO & \bf \shortstack{Complexity}\\			
			\hline
			\makecell[c]{\algname{EF14}\\ \tabcite{Seide20141bitSG}} & \thesischsevenxmark & -- & -- & \thesischsevenxmark 
			& $\cO\thesischsevenbr{\varepsilon^{-4}}$
			\\ 
			\hline
			\makecell[c]{\algname{Choco-SGD}\\ \tabcite{koloskova2020decentralized}} & \thesischsevenxmark & $\gamma\asymp K^{-1/2}$\textsuperscript{\color{blue}(2)} & -- & \thesischsevenxmark
			& $\cO\thesischsevenbr{\varepsilon^{-4}}$\\ 
			
			\hline
			\makecell[c]{\algname{EF21-SGD}\\ \tabcite{fatkhullin2025ef21}} & \thesischsevenxmark & $\gamma\asymp\bigl(L+\frac{\widetilde L}{\alpha}\bigr)^{-1}$\textsuperscript{\color{blue}(3)} & -- & \thesischsevenxmark
			& $\cO\thesischsevenbr{\varepsilon^{-4}}$\\

			\hline
			\makecell[c]{\algname{EF21-SGDM}\\ \tabcite{fatkhullin2023momentum}}
			& \thesischsevenxmark & $\gamma\lesssim\min\{\frac{\alpha}{\widetilde L},\frac{\eta}{L}\}$ & $\eta\asymp K^{-1/2}$ & \thesischsevenxmark & $\cO\thesischsevenbr{\varepsilon^{-4}}$\\ 

            \hline
			\makecell[c]{\algname{EF21-MVR}\\ \tabcite{fatkhullin2023momentum}}
			& \thesischsevenxmark & \makecell[c]{$\gamma\lesssim\min\{\frac{\alpha}{\widetilde L},\frac{\sqrt{\alpha}}{\widetilde\ell},$\\ $\frac{\sqrt{n\eta}}{\widetilde\ell}\}$} & $\eta\asymp K^{-2/3}$ & \thesischsevenxmark & $\cO\thesischsevenbr{\varepsilon^{-3}}$\\ 
            
			\hline
			\makecell[c]{\normEFname{SGDM}\\ \tabcite{khirirat2024errorfeedbackl0l1smoothnessnormalization}\textsuperscript{\color{blue}(4)}}
			& \thesischsevenxmark & $\gamma_k=\dfrac{\gamma}{(K+1)^{3/4}}$ & $\eta_k=\dfrac{1}{(K+1)^{1/2}}$ & \thesischsevenxmark & $\cO\thesischsevenbr{\varepsilon^{-4}}$\\ 
            
			\hline
			\cellcolor{linen}\begin{tabular}{c} \normEFname{SGDM}\\ Theorem~\ref{ch7:thm:sgdm_chapter_7} \end{tabular} & \cellcolor{linen}\thesischsevencmark & \cellcolor{linen}$\gamma_k=\gamma_0\left(\frac{2}{k+2}\right)^{3/4}$ & \cellcolor{linen}$\eta_k=\left(\frac{2}{k+2}\right)^{1/2}$ & \cellcolor{linen}\thesischsevenxmark
			& \cellcolor{linen}$\Tilde{\cO}\thesischsevenbr{\varepsilon^{-4}}$ \\ 
			\hline
			\cellcolor{linen}\begin{tabular}{c} \normEFname{IGT}\\ Theorem~\ref{ch7:thm:IGT_chapter_7} \end{tabular} & \cellcolor{linen}\thesischsevencmark & \cellcolor{linen}$\gamma_k=\gamma_0\left(\frac{2}{k+2}\right)^{5/7}$ & \cellcolor{linen}$\eta_k=\left(\frac{2}{k+2}\right)^{4/7}$ & \cellcolor{linen}\thesischsevenxmark
			& \cellcolor{linen}$\Tilde{\cO}\thesischsevenbr{\varepsilon^{-\nicefrac{7}{2}}}$ \\ 

            \hline
			\cellcolor{linen}\begin{tabular}{c} \normEFname{RHM}\\ Theorem~\ref{ch7:thm:RHM_chapter_7} \end{tabular} & \cellcolor{linen}\thesischsevencmark & \cellcolor{linen}$\gamma_k=\gamma_0\left(\frac{2}{k+2}\right)^{2/3}$ & \cellcolor{linen}$\eta_k=\left(\frac{2}{k+2}\right)^{2/3}$ & \cellcolor{linen}\thesischsevencmark
			& \cellcolor{linen}$\Tilde{\cO}\thesischsevenbr{\varepsilon^{-3}}$ \\

            \hline
			\cellcolor{linen}\begin{tabular}{c} \normEFname{HM}\\ Theorem~\ref{ch7:thm:HM_chapter_7} \end{tabular} & \cellcolor{linen}\thesischsevencmark & \cellcolor{linen}$\gamma_k=\gamma_0\left(\frac{2}{k+2}\right)^{2/3}$ & \cellcolor{linen}$\eta_k=\left(\frac{2}{k+2}\right)^{2/3}$ & \cellcolor{linen}\thesischsevencmark
			& \cellcolor{linen}$\Tilde{\cO}\thesischsevenbr{\varepsilon^{-3}}$ \\

            \hline
			\cellcolor{linen}\begin{tabular}{c} \normEFname{MVR}\\ Theorem~\ref{ch7:thm:MVR_chapter_7} \end{tabular} & \cellcolor{linen}\thesischsevencmark & \cellcolor{linen}$\gamma_k=\gamma_0\left(\frac{2}{k+2}\right)^{2/3}$ & \cellcolor{linen}$\eta_k=\left(\frac{2}{k+2}\right)^{2/3}$ & \cellcolor{linen}\thesischsevenxmark
			& \cellcolor{linen}$\Tilde{\cO}\thesischsevenbr{\varepsilon^{-3}}$ \\ 
			\hline
		\end{tabular}
        }
	\end{threeparttable}
		\begin{minipage}{\textwidth}
		{\scriptsize
		\begin{list}{}{\setlength{\leftmargin}{1.5em}\setlength{\labelwidth}{1em}\setlength{\labelsep}{0.5em}\setlength{\itemsep}{0pt}\setlength{\parsep}{0pt}}
		\item [{\color{blue}(1)}] In all displayed scalings, quantities other than $K$ are fixed and universal numerical constants are omitted. Here $\alpha\in(0,1]$ is the contractivity parameter of the compressor, $L$ is the smoothness constant of the global objective, and $\widetilde L:=\bigl(\frac{1}{n}\sum_{i=1}^n L_i^2\bigr)^{1/2}$ is the root-mean-square smoothness of the local objectives. For EF21-MVR, $\widetilde\ell:=\bigl(\frac{1}{n}\sum_{i=1}^n \ell_i^2\bigr)^{1/2}$ is the root-mean-square individual stochastic smoothness constant; this is a stronger assumption than ordinary smoothness. The exact selections also depend on the initial gap, noise level, number of workers, and, where applicable, the initial batch size. When each iteration uses a constant number of stochastic-gradient evaluations, sample and iteration complexity agree up to constants.
		\item [{\color{blue}(2)}] To match the thesis notation, $\gamma$ denotes the CHOCO-SGD learning rate (called $\eta$ in the cited paper). Its analysis uses $T$ for the time horizon; throughout this thesis we write $K:=T$. The method additionally uses a separate network-consensus parameter, which is omitted here.
		\item [{\color{blue}(3)}] The displayed EF21-SGD choice is the constant-stepsize scaling after the auxiliary parameters in the original analysis are fixed. Its bounded-variance guarantee additionally uses an accuracy-dependent mini-batch size.
		\item [{\color{blue}(4)}] Theorem~\ref{ch6:thm:ef21_sgdm_chapter_6} establishes this result under generalized $(L_0,L_1)$-smoothness. Its parameter selection requires the time horizon $K$ in advance; the additional upper bound on $\gamma$ is omitted from the table.
		\end{list}
		}
		\end{minipage}
\end{table*}

    \item \textbf{Numerical evaluation.} 
    In~Section~\ref{ch7:sec:exp_chapter_7}, we benchmark our proposed algorithms for solving the image classification task with the CIFAR-10 dataset using the ResNet-18 model. Among five momentum variants, \normEFname{HM}, which employs second-order momentum, achieves the fastest per-epoch convergence. However, this performance comes at the cost of higher computations per epoch. Notably, \normEFname{IGT} outperforms other error feedback algorithms in terms of the solution accuracy against the wall-clock time. 
\end{itemize}

\section{Related Work}

\paragraph{Error feedback.}
Error feedback (EF) mechanisms have been widely adopted to enhance the solution accuracy of gradient-based algorithms that employ communication compression. The first of these EF mechanisms, \algname{EF14}, was introduced by~\textcite{Seide20141bitSG} and later analyzed in both centralized~\parencite{StitchSparseSGD,karimireddy2019error} and distributed settings~\parencite{alistarh2018convergence,tang2019doublesqueeze,tang2019texttt,gorbunov2020linearly,khirirat2020compressed,qian2021error}. 
Another novel EF mechanism, \algname{EF21}, was proposed by~\textcite{EF21} and offers an improved convergence rate of $\cO(1/\sqrt{K})$ in the gradient norm, compared to the $\cO(1/K^{1/3})$ rate of earlier EF algorithms like~\textcite{koloskova2019decentralized}. 
More importantly, \algname{EF21} achieves this rate without requiring restrictive conditions, such as uniformly bounded gradient norms or bounded data heterogeneity.
Furthermore, \algname{EF21} has been extended to broader problem settings. For stochastic optimization, it can be adapted by using large mini-batches~\parencite{fatkhullin2025ef21} or by employing Polyak momentum, resulting in \algname{EF21-SGDM} and \algname{EF21-MVR}~\parencite{fatkhullin2023momentum}. In the context of federated optimization, where model parameters are exchanged instead of gradients, its variant called \algname{EF21-P}~\parencite{EF21-P} has been proposed. 
More recently, \algname{EControl}~\parencite{gao2024econtrol} was introduced to provide even stronger convergence guarantees for distributed stochastic optimization, thus advancing the development of provably efficient, distributed EF methods.

\paragraph{Normalization.}
Normalization has been commonly used to stabilize the training of randomized gradient-based algorithms for solving problems under relaxed smoothness conditions or in the presence of heavy-tailed noise.  
For minimizing relaxed smooth functions, normalization enables  stochastic momentum methods~\parencite{zhao2021convergence,hubler2024parameter} in the centralized setting and \algname{EF21-SGDM}~\parencite{fatkhullin2023momentum,khirirat2024errorfeedbackl0l1smoothnessnormalization} in the distributed setting to converge. 
In the heavy-tailed noise setting, normalization has been shown to ensure the convergence of stochastic methods in high-probability guarantees, including stochastic momentum methods as demonstrated by~\textcite{cutkosky2021high,hubler2024gradient}.
Furthermore, normalization allows for the use of parameter-agnostic stepsize rules, as shown by~\textcite{fatkhullin2023momentum}. 
However, these results are limited to the centralized setting.

\paragraph{Stochastic momentum algorithms.}
Stochastic momentum algorithms are widely used and studied for minimizing smooth objective functions. These algorithms are inspired by Polyak’s heavy-ball momentum~\parencite{polyak1964speeding}, which achieves accelerated linear convergence compared to classical gradient descent when applied to twice continuously differentiable, strongly convex, and smooth functions~\parencite{ghadimi2015global}. Several works~\parencite{yan2018unified,yu2019linear,liu2020improved,cutkosky2020momentum,hubler2024gradient} have shown that stochastic momentum algorithms can achieve a convergence rate of $\cO(1/K^{1/4})$, matching that of classical stochastic gradient descent (SGD).
To further accelerate convergence, recent research has proposed various momentum variants that modify the gradient estimators used in momentum updates. One notable approach is extrapolated momentum, introduced by~\textcite{cutkosky2020momentum}, which achieves a convergence rate of $\cO(1/K^{2/7})$. This has been further improved by algorithms such as \algname{STORM}~\parencite{cutkosky2019momentum}, \algname{MARS}~\parencite{yuan2024mars}, and two second-order momentum algorithms~\parencite{salehkaleybar2022momentum,tran2022better}, all of which attain a rate of $\cO(1/K^{1/3})$. This rate is known to be optimal for minimizing smooth nonconvex functions under mild conditions~\parencite{arjevani2023lower}.
However, these algorithms have primarily been developed and analyzed in the context of centralized optimization. As a result, their  applicability to the distributed optimization setting for broader training applications remains limited.

\paragraph{Parameter-agnostic algorithms.}
Other approaches for parameter-agnostic stepsizes, in addition to normalization, include adaptive stochastic methods, such as backtracking line search~\parencite{armijo1966minimization}, AdaGrad~\parencite{duchi2011adaptive}, AdaGrad-Norm~\parencite{streeter2010less}, and ADAM~\parencite{kingma2017adammethodstochasticoptimization}.

\section{Preliminaries}

\subsection{Notation}

We denote the expectation of a random variable $u$ by $\Exp[u]$. For any vectors $x, y \in \R^d$, $\inner{x}{y}$ refers to their inner product, and $\norm{x} = \sqrt{\inner{x}{x}}$ denotes the Euclidean norm of the vector $x$. For a matrix $A \in \R^{m \times n}$, $\norm{A}_{\mathrm{op}}$ denotes its spectral norm, i.e., the largest singular value of $A$. The notation $\cO(h(x))$ implies that a function $f(x)$ satisfies $f(x) \leq c \cdot h(x)$ for some constant $c > 0$, while $\tilde{\cO}(h(x))$ hides both constant and logarithmic factors. 
Finally, we use $\underset{x\in\R^d}{\inf} f(x)$ to denote the infimum of a function $f \colon \R^d \to \R$, and $\underset{x\in\R^d}{\min} f(x)$ to denote its minimum, when it exists.

\subsection{Problem Formulation}
Consider a distributed stochastic optimization problem: 
\begin{eqnarray}\label{ch7:eqn:problem_chapter_7}
 \underset{x\in\R^d}{\min} \quad f(x):=\frac{1}{n}\sum_{i=1}^n f_i(x), 
\end{eqnarray}
where $f_i(x)= \thesischsevenExpSub{\xi_i\sim \thesischsevencD_i}{f_i(x;\xi_i)}$, and  $f_i(x;\xi_i)$ is a possibly nonconvex function  parameterized by the vector $x\in\R^d$ on the random variable $\xi_i$ drawn from the data distribution $\thesischsevencD_i$ known to client $i$. 
Problem~(\ref{ch7:eqn:problem_chapter_7}) often appears in supervised machine learning applications~\parencite{friedman2009elements}.

\subsection{Assumptions}

To facilitate our  analysis, we impose standard assumptions on compression operators and objective functions.

First, we assume an $\alpha$-contractive  compression operator, which covers many popular biased compressors of interest, including TopK~\parencite{alistarh2018convergence} and RandK sparsifiers~\parencite{beznosikov2023biased}.
\begin{assumption}[Contractive compression]\label{ch7:assum:contract_comp_chapter_7}
	A biased but possibly randomized compressor $\thesischsevencC:\R^d\rightarrow\R^d$ is $\alpha$-contractive with its sample $\xi_i \sim \thesischsevencD_i$ if there exists $\alpha \in (0,1]$ such that 
\begin{equation*}
	\Exp\!\left[\sqnorm{\thesischsevencC(v)-v}\right]
	\leq (1-\alpha)\sqnorm{v}, \quad \forall v \in \R^d.
\end{equation*}
\end{assumption}

Second, we introduce commonly used conditions on the objective function, including the existence of a finite infimum, and Lipschitz continuity of the gradient and Hessian of component and stochastic functions.
\begin{assumption}\label{ch7:assum:inf_chapter_7}
The function $f:\R^d\rightarrow\R$ is bounded from below, i.e., $f^{\inf}= \underset{x\in\R^d}{\inf} f(x) > -\infty$.
\end{assumption}

\begin{assumption}\label{ch7:assum:L_stoc_comp_fcn_chapter_7}
The stochastic component function $f_i(x;\xi_i)$ has the $L_{\text{ms},i}$-Lipschitz continuous gradient, if there exists $L_{\text{ms},i} >0$ such that for all $x,y\in\R^d$, 
\begin{eqnarray*}
    \thesischsevenExpSub{\xi_i}{\sqnorm{\nabla f_i(x;\xi_i) - \nabla f_i(y;\xi_i)}} \leq L_{\text{ms},i}^2 \sqnorm{x-y}.
\end{eqnarray*}
\end{assumption}

\begin{assumption}\label{ch7:assum:L_comp_fcn_chapter_7}
The component function $f_i:\R^d\rightarrow\R$ has the $L_i$-Lipschitz continuous gradient, if there exists $L_i >0$ such that for all $x,y\in\R^d$, 
\begin{eqnarray*}
    \norm{\nabla f_i(x) - \nabla f_i(y)} \leq L_i \norm{x-y}.
\end{eqnarray*}
\end{assumption}

\begin{assumption}\label{ch7:assum:L_fcn_chapter_7}
The function $f:\R^d\rightarrow\R$ has the $L$-Lipschitz continuous gradient, if there exists $L >0$ such that for all $ x,y\in\R^d$, 
\begin{eqnarray*}
    \norm{\nabla f(x) - \nabla f(y)} \leq L \norm{x-y}.
\end{eqnarray*}
\end{assumption}

\begin{assumption}\label{ch7:assum:L_Hessian_comp_fcn_chapter_7}
The component function $f_i:\R^d\rightarrow\R$ has the $L_{h,i}$-Lipschitz continuous Hessian, if there exists $L_{h,i} >0$ such that for all $x,y\in\R^d$, 
\begin{eqnarray*}
    \norm{\nabla^2 f_i(x) - \nabla^2 f_i(y)}_{\mathrm{op}} \leq L_{h,i} \norm{x-y}.
\end{eqnarray*}
\end{assumption}

\begin{assumption}\label{ch7:assum:L_Hessian_fcn_chapter_7}
The  function $f:\R^d\rightarrow\R$ has the $L_{h}$-Lipschitz continuous Hessian, if there exists $L_{h} >0$ such that for all $x,y\in\R^d$, 
\begin{eqnarray*}
    \norm{\nabla^2 f(x) - \nabla^2 f(y)}_{\mathrm{op}} \leq L_{h} \norm{x-y}.
\end{eqnarray*}
\end{assumption}

Note that Assumption~\ref{ch7:assum:L_stoc_comp_fcn_chapter_7} implies Assumption~\ref{ch7:assum:L_comp_fcn_chapter_7} with $L_i = L_{\text{ms},i}$, while  Assumption~\ref{ch7:assum:L_comp_fcn_chapter_7} implies Assumption~\ref{ch7:assum:L_fcn_chapter_7} with $L = \frac{1}{n}\sum_{i=1}^n L_i$. Furthermore, Assumption~\ref{ch7:assum:L_Hessian_comp_fcn_chapter_7} implies Assumption~\ref{ch7:assum:L_Hessian_fcn_chapter_7} with $L_h = \frac{1}{n}\sum_{i=1}^n L_{h,i}$.

Third, we assume that each client can query its stochastic oracle to obtain its local noisy gradient and Hessian that satisfy the unbiased and variance-bounded condition, which is commonly used for analyzing stochastic gradient methods~\parencite{cutkosky2019momentum,cutkosky2020momentum,tran2022better,hubler2024parameter,hubler2024gradient}. 
\begin{assumption}\label{ch7:assum:stoc_g_h_chapter_7}
	At iteration $k$, let $\mathcal{F}_k$ denote the history available before the auxiliary randomness and stochastic-oracle samples used in the momentum update are drawn. Conditional on this history, the samples are fresh and independent across clients and iterations. For \normEFname{RHM}, the master draws one scalar $q_k\sim\thesischsevencU(0,1)$, independently of $\mathcal{F}_k$, and shares it with all clients. Conditional on $\mathcal{F}_k$ and $q_k$, each client draws independent samples $\hat{\xi}_i^{k+1}$ and $\xi_i^{k+1}$ for the Hessian and gradient evaluations, respectively. Let $\mathcal{H}_k=\mathcal{F}_k$ for the other variants and $\mathcal{H}_k=\sigma(\mathcal{F}_k,q_k)$ for \normEFname{RHM}.

	The local stochastic gradient $\nabla f_i(x;\xi_i)$ at client $i$ is an unbiased estimator of $\nabla f_i(x)$ with bounded variance. For every query point $x\in\R^d$ measurable with respect to the history preceding its oracle call, 
	\begin{equation*}
	\begin{aligned}
		\Exp\!\left[\nabla f_i(x;\xi_i)\mid\mathcal{H}_k\right] &= \nabla f_i(x),\\
		\Exp\!\left[\sqnorm{\nabla f_i(x;\xi_i) - \nabla f_i(x)}\mid\mathcal{H}_k\right] &\leq \sigma_g^2.
	\end{aligned}
	\end{equation*}
	Furthermore, 
	the local stochastic Hessian $\nabla^2 f_i(x;\xi_i)$ at client $i$ is an unbiased estimator of $\nabla^2 f_i(x)$ with bounded variance:
	\begin{equation*}
	\begin{aligned}
		\Exp\!\left[\nabla^2 f_i(x;\xi_i)\mid\mathcal{H}_k\right] &= \nabla^2 f_i(x),\\
		\Exp\!\left[\norm{\nabla^2 f_i(x;\xi_i) - \nabla^2 f_i(x)}_{\mathrm{op}}^2\mid\mathcal{H}_k\right] &\leq \sigma_h^2.
	\end{aligned}
	\end{equation*}
\end{assumption}

\section{\algname{EF21} Methods with Five Momentum Updates}
\label{ch7:sec:ef21_momentum_chapter_7}

To solve Problem~(\ref{ch7:eqn:problem_chapter_7}), we consider \algname{EF21} algorithms (Algorithm~\ref{ch7:alg:EF_momentum_chapter_7}) that leverage normalized descent updates and momentum variants. 

At each iteration $k = 0,1, \ldots, K-1$, key updating rules for~Algorithm~\ref{ch7:alg:EF_momentum_chapter_7} are described as follows: 
\begin{itemize}
    \item Each client $i=1,2,\ldots,n$ computes stochastic gradients $\nabla f_i(\cdot;\xi_i^{k+1})$, and updates its local momentum vector $v_i^{k+1}$ depending on the choice of momentum.
    \item Each client $i=1,2,\ldots,n$ compresses $v_i^{k+1}-g_i^k$, transmits $c_i^{k+1} = \mathcal{C}_i^{k+1}(v_i^{k+1}-g_i^k)$, and updates its local memory vector as $g_i^{k+1} = g_i^k + c_i^{k+1}$.
    \item The server receives $c_i^{k+1}$ from every client, updates $g^{k+1} = g^k + \frac{1}{n}\sum_{i=1}^n c_i^{k+1}$, and computes $x^{k+1} = x^k -\gamma_k \cdot \nicefrac{g^k}{\norm{g^k}}$.
\end{itemize}

Depending on the choice for updating the momentum,~Algorithm~\ref{ch7:alg:EF_momentum_chapter_7} is referred to as:
\begin{itemize}
    \item \normEFname{SGDM}~\parencite{khirirat2024errorfeedbackl0l1smoothnessnormalization} when using the Polyak momentum, 
    \item \normEFname{IGT} when using the Implicit Gradient Transport (IGT) momentum~\parencite{cutkosky2020momentum}, 
    \item \normEFname{RHM} when using the second-order momentum variant by~\textcite{salehkaleybar2022momentum}, \item \normEFname{HM} when using the second-order momentum variant by~\textcite{tran2022better}, and 
    \item \normEFname{MVR} when using the momentum variant from \algname{STORM}~\parencite{cutkosky2019momentum}.
\end{itemize}
Here, \normEFname{SGDM} and  \normEFname{MVR} are the normalized versions of \algname{EF21-SGDM} and \algname{EF21-MVR}, respectively, as analyzed by~\textcite{fatkhullin2023momentum}.
Moreover, by setting the compression operator $\thesischsevencC_i^k$ to the identity and choosing $n=1$, 
our algorithms reduce to the corresponding normalized stochastic momentum methods in the centralized setting: 
(1) \normEFname{SGDM} becomes normalized stochastic momentum methods~\parencite{cutkosky2020momentum,hubler2024parameter}, (2) \normEFname{IGT} becomes normalized SGD with IGT momentum~\parencite{cutkosky2020momentum}, (3) \normEFname{RHM} becomes normalized second-order momentum methods analyzed by~\textcite{salehkaleybar2022momentum}, (4) \normEFname{HM} becomes second-order momentum methods proposed by~\textcite{tran2022better}, and (5) \normEFname{MVR} becomes STORM~\parencite{cutkosky2019momentum}. 
Finally, notice that~Algorithm~\ref{ch7:alg:EF_momentum_chapter_7} with $v_i^{k+1} = \nabla f_i(x^{k+1})$ recovers \normEFname{} analyzed by~\textcite{khirirat2024errorfeedbackl0l1smoothnessnormalization}.

\begin{algorithm}[t]
\caption{\algname{EF21} with Normalized Updates and Momentum Variants}
\label{ch7:alg:EF_momentum_chapter_7}
\begin{algorithmic}[1]
\STATE \textbf{Input:} Initial point $x^0, g_i^0,v_i^0 \in \mathbb{R}^d$, stepsize $\gamma_k > 0$, momentum parameter $0<\eta_k \leq 1$, compressor $\mathcal{C}_i^k$, and number of iterations $K$
\FOR{$k = 0, 1, \ldots, K-1$}
    \STATE Master computes $x^{k+1} = x^k - \gamma_k \frac{g^k}{\norm{g^k}}$, and broadcasts $x^{k+1}$ to all the clients
    \STATE \normEFname{RHM} only: Master draws one $q_k\sim\thesischsevencU(0,1)$
    \STATE \quad Set $\hat x^{k+1}=q_kx^{k+1}+(1-q_k)x^k$ and broadcast $\hat x^{k+1}$
    \FOR{Every client $i=1,2,\ldots,n$}
        \vspace{0.3em}
        \STATE \normEFname{SGDM}:
        \STATE  \quad Set $v_i^{k+1}=(1-\eta_k)v_i^k + \eta_k \nabla f_i \left( x^{k+1};\xi_i^{k+1} \right)$

        \vspace{0.3em}
        \STATE \normEFname{IGT}:
        \STATE \quad Set $y^k = x^{k+1} + \frac{1-\eta_k}{\eta_k}(x^{k+1}-x^k)$
        \STATE  \quad Set $v_i^{k+1}=(1-\eta_k)v_i^k + \eta_k \nabla f_i \left( y^k;\xi_i^{k+1} \right)$

        \vspace{0.3em}
        \STATE \normEFname{RHM}:
        \STATE \quad Draw independent samples $\hat{\xi}_i^{k+1}$ and $\xi_i^{k+1}$
        \STATE \quad Set $\tilde v_i^k = v_i^k + \nabla^2 f_i(\hat x^{k+1};\hat{\xi}_i^{k+1})(x^{k+1} - x^k)$
        \STATE  \quad Set $v_i^{k+1}=(1-\eta_k)\tilde v_i^k + \eta_k \nabla f_i(x^{k+1};\xi_i^{k+1})$

        \vspace{0.3em}
        \STATE \normEFname{HM}:
        \STATE \quad Set $\tilde v_i^k = v_i^k + \nabla^2 f_i(x^{k+1};\xi_i^{k+1})(x^{k+1} - x^k)$
        \STATE  \quad Set $v_i^{k+1}=(1-\eta_k)\tilde v_i^k + \eta_k \nabla f_i(x^{k+1};\xi_i^{k+1})$

        \vspace{0.3em}
        \STATE \normEFname{MVR}:
        \STATE \quad Set $\tilde v_i^k = v_i^k + \nabla f_i(x^{k+1};\xi_i^{k+1}) - \nabla f_i(x^k;\xi_i^{k+1})$
        \STATE  \quad Set $v_i^{k+1}=(1-\eta_k)\tilde v_i^k + \eta_k \nabla f_i(x^{k+1};\xi_i^{k+1})$

        \vspace{0.3cm}

        \STATE Transmit $c_i^{k+1} = \thesischsevencC_i^{k+1}(v_i^{k+1}-g_i^k)$
        \STATE Update $g_i^{k+1}=g_i^k + c_i^{k+1}$
    \ENDFOR 
    \STATE Master receives $c_1^{k+1},\ldots,c_n^{k+1}$, and computes $g^{k+1} = g^k + \frac{1}{n}\sum_{i=1}^n c_i^{k+1}$
\ENDFOR
\STATE \textbf{Theoretical output:} Draw $\tilde{x}^K$ from $\{x^0,\ldots,x^{K-1}\}$ with
\STATE \quad $\Pr(\tilde{x}^K=x^k)=\frac{\gamma_k}{\sum_{t=0}^{K-1}\gamma_t}$; implementations may also retain $x^K$
\end{algorithmic}
\end{algorithm}

\section{Convergence Theorems}
\label{ch7:sec:convergence_chapter_7}

To this end, we present the convergence of Algorithm~\ref{ch7:alg:EF_momentum_chapter_7} that use five momentum variants: Polyak momentum, IGT momentum, MVR momentum in \algname{STORM}, and two second-order momentum schemes.  
Throughout this section, $K$ denotes the number of iterations, the guarantees concern the expected gradient norm at the stepsize-weighted random iterate $\tilde{x}^K$, and $\widetilde{\cO}$ suppresses logarithmic factors. Each method uses a constant number of stochastic-oracle evaluations per iteration, so its sample and iteration complexities agree up to method-dependent constants.

\paragraph{\normEFname{SGDM}.}
We begin with the convergence rate of \normEFname{SGDM}.

\begin{theorem}\label{ch7:thm:sgdm_chapter_7}
Consider Problem~(\ref{ch7:eqn:problem_chapter_7}), where Assumptions~\ref{ch7:assum:contract_comp_chapter_7},~\ref{ch7:assum:inf_chapter_7},~\ref{ch7:assum:L_comp_fcn_chapter_7},~\ref{ch7:assum:L_fcn_chapter_7}, and~\ref{ch7:assum:stoc_g_h_chapter_7} hold. 
Let tuning parameters satisfy 
\begin{eqnarray*}
\eta_k = \left( \frac{2}{k+2} \right)^{1/2} \quad \text{and} \quad \gamma_k = \gamma_0 \left( \frac{2}{k+2} \right)^{3/4}
\end{eqnarray*}
with $\gamma_0 >0$. 
Then, the iterates $\{x^k\}$ governed by \algname{EF21-SGDM}
satisfy
\begin{equation*}
\Exp\!\left[\norm{\nabla f(\tilde{x}^K)}\right]
= \widetilde{\cO}\left(\frac{\nicefrac{V_0}{\gamma_0} + \gamma_0\left(L + \nicefrac{\bar L}{\alpha^2}\right) + \sigma_g\left(\nicefrac{1}{\sqrt{n}} + \nicefrac{1}{\alpha^2}\right)}{K^{\nicefrac{1}{4}}}\right).
\end{equation*}
where $\tilde{x}^K$ is randomly chosen from $\{x^0,x^1,\ldots,x^{K-1}\}$ with probability $\nicefrac{\gamma_k}{\sum_{t=0}^{K-1}\gamma_t}$ for $k = 0,1, \ldots, K-1$. Here, $\bar L = \frac{1}{n}\sum_{i=1}^n L_i$.
\end{theorem}

From~Theorem~\ref{ch7:thm:sgdm_chapter_7}, \normEFname{SGDM} achieves an $\widetilde{\cO}(1/K^{1/4})$ rate for the expected gradient norm at the stepsize-weighted random iterate, with logarithmic factors suppressed.
This result holds without requiring the stepsizes to depend on problem-dependent parameters (e.g. smoothness constants), but still yields the convergence almost matching the $\cO(1/K^{1/4})$ convergence bound for \normEFname{SGDM} by~\textcite{khirirat2024errorfeedbackl0l1smoothnessnormalization} and for \algname{EF21-SGDM} analyzed by~\textcite{fatkhullin2023momentum}.
Furthermore, unlike Theorem 3 in \textcite{fatkhullin2023momentum}, our result does not rely on initializing the algorithm with a sufficiently large mini-batch size. 
Finally, by setting the compression operator  $\mathcal{C}^k(\cdot)$ be the identity operator and by letting $n=1$, our result recovers the same 
$\tilde{\cO}(1/K^{1/4})$ convergence rate as that of normalized stochastic momentum methods in the centralized setting by~\textcite{hubler2024parameter}.

\paragraph{\normEFname{IGT}.}
To further improve the convergence performance of \normEFname{SGDM}, we replace Polyak momentum with IGT momentum~\parencite{cutkosky2020momentum} in~Algorithm~\ref{ch7:alg:EF_momentum_chapter_7}.
The resulting algorithm, referred to as \normEFname{IGT}, admits the following convergence guarantee: 
\begin{theorem}\label{ch7:thm:IGT_chapter_7}
Consider Problem~(\ref{ch7:eqn:problem_chapter_7}), where
Assumptions~\ref{ch7:assum:contract_comp_chapter_7}, \ref{ch7:assum:inf_chapter_7},
\ref{ch7:assum:L_comp_fcn_chapter_7}, and \ref{ch7:assum:L_fcn_chapter_7} hold.
Suppose also that Assumptions~\ref{ch7:assum:L_Hessian_comp_fcn_chapter_7},
\ref{ch7:assum:L_Hessian_fcn_chapter_7}, and \ref{ch7:assum:stoc_g_h_chapter_7} hold.
Let tuning parameters satisfy 
\begin{eqnarray*}
\eta_k = \left( \frac{2}{k+2} \right)^{4/7} \quad \text{and} \quad \gamma_k = \gamma_0 \left( \frac{2}{k+2} \right)^{5/7}
\end{eqnarray*}
with $\gamma_0 > 0$. 
Then, the iterates $\{x^k\}$ governed by \algname{EF21-IGT}
satisfy
\begin{equation*}
\Exp\!\left[\norm{\nabla f(\tilde{x}^K)}\right]
= \widetilde{\cO}\left(\frac{\nicefrac{V_0}{\gamma_0} + \sigma_g\left(\nicefrac{1}{\sqrt{n}} + \nicefrac{1}{\alpha^2}\right) + \gamma_0B_1 + \gamma_0^2B_2}{K^{\nicefrac{2}{7}}}\right).
\end{equation*}
where $\tilde{x}^K$ is randomly chosen from $\{x^0,x^1,\ldots,x^{K-1}\}$ with probability $\nicefrac{\gamma_k}{\sum_{t=0}^{K-1}\gamma_t}$ for $k = 0,1, \ldots, K-1$.
Here, $B_1 = L+\nicefrac{\bar L}{\alpha^2}$ and $B_2 = L_h +\nicefrac{\bar L_h}{\alpha^2}$. Here, $\bar L =\frac{1}{n}\sum_{i=1}^n L_i$ and $\bar L_h = \frac{1}{n}\sum_{i=1}^n L_{h,i}$.
\end{theorem}

From~Theorem~\ref{ch7:thm:IGT_chapter_7}, \normEFname{IGT}, which uses decreasing, parameter-agnostic stepsize rules, attains an $\widetilde{\cO}(1/K^{2/7})$ rate for the expected gradient norm at the stepsize-weighted random iterate, with logarithmic factors suppressed. 
When $\mathcal{C}^k(\cdot) = I$ and $n=1$, this result recovers the same $\tilde{\cO}(1/K^{2/7})$ rate as stochastic methods with IGT momentum in the centralized setting--specifically, Theorem 3 of~\textcite{cutkosky2020momentum}, Theorem 6 of~\textcite{cutkosky2021high} in the bounded variance case,  and  Theorem 2 of~\textcite{sun2023momentum}.
Unlike these prior results, however, our method achieves this rate without requiring stepsizes to know smoothness constants.

\paragraph{\normEFname{RHM}.}
Another momentum variant we incorporate into \algname{EF21} methods using momentum (Algorithm~\ref{ch7:alg:EF_momentum_chapter_7}) to enhance their convergence performance is second-order momentum proposed by~\textcite{salehkaleybar2022momentum}. 
We refer to the resulting algorithm as~\normEFname{RHM}, and its convergence guarantee is presented below:  
\begin{theorem}\label{ch7:thm:RHM_chapter_7}
Consider Problem~(\ref{ch7:eqn:problem_chapter_7}), where Assumptions~\ref{ch7:assum:contract_comp_chapter_7},~\ref{ch7:assum:inf_chapter_7},~\ref{ch7:assum:L_comp_fcn_chapter_7},~\ref{ch7:assum:L_fcn_chapter_7}
and~\ref{ch7:assum:stoc_g_h_chapter_7} hold. 
Let tuning parameters satisfy 
\begin{eqnarray*}
\eta_k = \left( \frac{2}{k+2} \right)^{2/3}\quad \text{and} \quad \gamma_k = \gamma_0 \left( \frac{2}{k+2} \right)^{2/3}
\end{eqnarray*}
with $\gamma_0 >0$. 
Then, the iterates $\{x^k\}$ governed by \algname{EF21-RHM}
satisfy
\begin{equation*}
\Exp\!\left[\norm{\nabla f(\tilde{x}^K)}\right]
= \widetilde{\cO}\left(\frac{\nicefrac{V_0}{\gamma_0} + \left(\sigma_g + \gamma_0\sigma_h\right)\left(\nicefrac{1}{\sqrt{n}} + \nicefrac{1}{\alpha^2}\right) + \gamma_0 C}{K^{\nicefrac{1}{3}}}\right).
\end{equation*}
where $\tilde{x}^K$ is randomly chosen from $\{x^0,x^1,\ldots,x^{K-1}\}$ with probability $\nicefrac{\gamma_k}{\sum_{t=0}^{K-1}\gamma_t}$ for $k = 0,1, \ldots, K-1$. Here, $C= L +\nicefrac{\bar L}{\alpha^2}$ and $\bar L = \frac{1}{n}\sum_{i=1}^n L_i$.
\end{theorem}

From~Theorem~\ref{ch7:thm:RHM_chapter_7}, \normEFname{RHM} achieves an $\widetilde{\cO}(1/K^{1/3})$ rate for the expected gradient norm at the stepsize-weighted random iterate, with logarithmic factors suppressed, thus improving upon the rates of \normEFname{SGDM} and \normEFname{IGT}.
When $\mathcal{C}^k(\cdot) = I$ and $n=1$, our result nearly matches the ${\cO}(1/K^{1/3})$ convergence for stochastic second-order momentum methods in the centralized setting  by~\textcite{salehkaleybar2022momentum}, Theorem 2 of \textcite{sadiev2025second} in the bounded variance case.

\paragraph{\normEFname{HM}.}
In addition to the second-order momentum by~\textcite{salehkaleybar2022momentum}, we also incorporate the second-order momentum proposed by~\textcite{tran2022better} into~Algorithm~\ref{ch7:alg:EF_momentum_chapter_7}. 
The resulting algorithm, which we call \normEFname{HM}, achieves the following convergence rate:  
\begin{theorem}\label{ch7:thm:HM_chapter_7}
Consider Problem~(\ref{ch7:eqn:problem_chapter_7}), where Assumptions~\ref{ch7:assum:contract_comp_chapter_7},~\ref{ch7:assum:inf_chapter_7},
\stepEightBchange{~\ref{ch7:assum:L_comp_fcn_chapter_7},}
~\ref{ch7:assum:L_fcn_chapter_7},
~\ref{ch7:assum:L_Hessian_comp_fcn_chapter_7}
and~\ref{ch7:assum:stoc_g_h_chapter_7} hold. 
Let tuning parameters satisfy 
\begin{eqnarray*}
\eta_k = \left( \frac{2}{k+2} \right)^{2/3}  \quad \text{and} \quad \gamma_k = \gamma_0 \left( \frac{2}{k+2} \right)^{2/3}
\end{eqnarray*}
with $\gamma_0 >0$. 
Then, the iterates $\{x^k\}$ governed by \algname{EF21-HM}
satisfy
\begin{equation*}
\Exp\!\left[\norm{\nabla f(\tilde{x}^K)}\right]
= \widetilde{\cO}\left(\frac{\nicefrac{V_0}{\gamma_0} + D_\sigma\left(\nicefrac{1}{\sqrt{n}} + \nicefrac{1}{\alpha^2}\right) + \gamma_0D_1 + \gamma_0^2D_2}{K^{\nicefrac{1}{3}}}\right).
\end{equation*}
where $\tilde{x}^K$ is randomly chosen from $\{x^0,x^1,\ldots,x^{K-1}\}$ with probability $\nicefrac{\gamma_k}{\sum_{t=0}^{K-1}\gamma_t}$ for $k = 0,1, \ldots, K-1$. Here, $D_\sigma=\sigma_g + \gamma_0\sigma_h$, $D_1 =L + \nicefrac{\bar L}{\alpha^2}$,  $D_2 = L_h +\nicefrac{\bar L_h}{\alpha^2}$, $\bar L = \frac{1}{n}\sum_{i=1}^n L_i$, and $\bar L_h = \frac{1}{n}\sum_{i=1}^n L_{h,i}$.
\end{theorem}

From~Theorem~\ref{ch7:thm:HM_chapter_7}, \normEFname{HM} achieves the same $\widetilde{\cO}(1/K^{1/3})$ rate for the expected gradient norm at the stepsize-weighted random iterate, with logarithmic factors suppressed, as \normEFname{RHM}. 
When $\mathcal{C}^k(\cdot) = I$ and $n=1$, our result closely matches the ${\cO}(1/K^{1/3})$ convergence for stochastic second-order momentum methods in the centralized setting  by~\textcite{tran2022better}.

\paragraph{\normEFname{MVR}.}
Finally, we consider Algorithm~\ref{ch7:alg:EF_momentum_chapter_7} with the MVR momentum from \algname{STORM}~\parencite{cutkosky2019momentum}, resulting in \normEFname{MVR}. 
The next theorem presents its parameter-agnostic convergence guarantee.
\begin{theorem}\label{ch7:thm:MVR_chapter_7}
Consider Problem~(\ref{ch7:eqn:problem_chapter_7}), where Assumptions~\ref{ch7:assum:contract_comp_chapter_7},~\ref{ch7:assum:inf_chapter_7},
~\ref{ch7:assum:L_stoc_comp_fcn_chapter_7},
~\ref{ch7:assum:L_fcn_chapter_7}
and~\ref{ch7:assum:stoc_g_h_chapter_7} hold. 
Let tuning parameters satisfy 
\begin{eqnarray*}
\eta_k = \left( \frac{2}{k+2} \right)^{2/3} \quad \text{and} \quad \gamma_k = \gamma_0 \left( \frac{2}{k+2} \right)^{2/3}
\end{eqnarray*}
with $\gamma_0 > 0$. 
Then, the iterates $\{x^k\}$ governed by \algname{EF21-MVR}
satisfy
\begin{equation*}
\Exp\!\left[\norm{\nabla f(\tilde{x}^K)}\right]
= \widetilde{\cO}\left(\frac{\nicefrac{V_0}{\gamma_0} + \sigma_g\left(\nicefrac{1}{\sqrt{n}} + \nicefrac{1}{\alpha^2}\right) + \gamma_0E_1}{K^{\nicefrac{1}{3}}}\right).
\end{equation*}
where $\tilde{x}^K$ is randomly chosen from $\{x^0,x^1,\ldots,x^{K-1}\}$ with probability $\nicefrac{\gamma_k}{\sum_{t=0}^{K-1}\gamma_t}$ for $k = 0,1, \ldots, K-1$. Here, $E_1=L + {\bar L}_{\text{ms}} + \nicefrac{{\bar L}_{\text{ms}}}{\alpha^2}$ and $\bar L_{\text{ms}} = \frac{1}{n}\sum_{i=1}^n L_{\text{ms},i}$.
\end{theorem}

From~Theorem~\ref{ch7:thm:MVR_chapter_7}, \normEFname{MVR} achieves the same $\widetilde{\cO}(1/K^{1/3})$ rate for the expected gradient norm at the stepsize-weighted random iterate, with logarithmic factors suppressed, as \normEFname{RHM} and \normEFname{HM}. 
Unlike \normEFname{RHM} and \normEFname{HM}, which require Hessian-vector evaluations per iteration, 
\normEFname{MVR} needs two stochastic gradient computations. 
Our result nearly matches the ${\cO}(1/K^{1/3})$ convergence for \algname{EF21-STORM/MVR} by~\textcite{fatkhullin2023momentum}, but with notable improvements. 
\normEFname{MVR} does not require a large mini-batch size at initialization, and operates with parameter-agnostic stepsizes without the need to know problem-specific quantities like smoothness constants. 
Furthermore, when $\mathcal{C}^k(\cdot) = I$ and $n=1$, our result closely matches the ${\cO}(1/K^{1/3})$ convergence for \algname{STORM} in the centralized setting  by~\textcite{cutkosky2019momentum}, but without assuming additional restrictive conditions like bounded norms of stochastic gradients.

\section{Numerical Experiments}\label{ch7:sec:exp_chapter_7}
We evaluated the performance of \algname{EF21} methods\footnote{To support reproducibility, we provide the full experimental code in our \href{https://github.com/IgorSokoloff/Parameter-Agnostic-Error-Feedback-through-Momentum}{GitHub repository}.} with momentum variants (Algorithm~\ref{ch7:alg:EF_momentum_chapter_7}) for training the ResNet-18 model~\parencite{he2016deep} (with $d = 11,173,962$ parameters) on the CIFAR-10 dataset~\parencite{krizhevsky2009learning}.
In particular, we benchmark~Algorithm~\ref{ch7:alg:EF_momentum_chapter_7} with five momentum variants against existing distributed error feedback methods: \algname{EF21-SGD}~\parencite{fatkhullin2025ef21}, \algname{EF21-SGDM}~\parencite{fatkhullin2023momentum}, and \algname{EControl}~\parencite{gao2024econtrol}.

In these experiments, we adopted a data distribution strategy inspired by~\textcite{gao2024econtrol}. Specifically, we first held out 10\% of the CIFAR-10 training set as a common validation set. Of the remaining training images, one half was distributed randomly and uniformly among the 10 clients, and the other half was allocated by class label, so that images with the $i$-th label (for $i \in \{0, \dots, 9\}$) were assigned to client $i+1$. All local datasets were then truncated to the size of the smallest one, so that every client holds the same number of images. Test accuracy is measured on the standard CIFAR-10 test set of $10{,}000$ images. This partitioning scheme introduces data heterogeneity, a common characteristic of federated settings. For communication compression, we employed the Top-$K$ sparsifier, retaining 10\% of the coordinates (i.e., $\nicefrac{K}{d} = 0.1$).
All experiments were implemented using PyTorch~\parencite{paszke2019pytorch}, and were performed on a server-grade machine running Ubuntu 18.04 (Linux Kernel v5.4.0). This system was equipped with dual 16-core 3.3 GHz Intel Xeon processors (totaling 32 cores) and four NVIDIA A100 GPUs, each with 40GB of memory.
Each method was trained once, with the fixed random seed $42$; Table~\ref{ch7:main_table:accuracy_results_chapter_7} and the figures report these single runs. Run-to-run variability was therefore not assessed, and small differences between methods should be read with this limitation in mind. We report convergence with respect to the number of epochs, GPU time, and wall-clock time.
The hyperparameter tuning procedure for each method and the definitions of GPU and wall-clock time are given in Appendix~\ref{ch7:sec:restored_tuning} and Appendix~\ref{ch7:sec:restored_timing}, respectively.

\begin{table*}[htbp]
\centering
\caption{Best performance metrics achieved by each method when training ResNet-18 on CIFAR-10. The baselines (top block) and the proposed methods (bottom block) are each listed in increasing order of best validation accuracy. The two methods with the highest validation accuracy are highlighted in \textbf{bold}; their remaining entries are bold for identification only and are not necessarily the best in their columns.}
\label{ch7:main_table:accuracy_results_chapter_7}

\begin{threeparttable}
\resizebox{0.95\textwidth}{!}{%
\begin{tabular}{@{}lccccc@{}}
\toprule
\textbf{Method} & \textbf{\makecell{Best Val.\\Accuracy (\%)}} & \textbf{\makecell{Corr. Test\\Accuracy (\%)}} & \textbf{\makecell{Epoch of\\Best}} & \textbf{\makecell{GPU Time\\to Best}} & \textbf{\makecell{Wall Time\\to Best}} \\ \midrule
\algname{EF21-SGDM}          & $74.66$          & $73.53$          & $76$ & 0h 35m 59s & 0h 35m 59s \\
\algname{EControl}           & $76.38$          & $74.49$          & $79$ & 0h 38m 06s & 0h 38m 06s \\
\algname{EF21-SGD}           & $77.50$          & $75.90$          & $69$ & 0h 30m 35s & 0h 30m 36s \\
\midrule
\normEFname{MVR}  & $79.30$          & $78.77$          & $72$ & 0h 44m 33s & 0h 44m 33s \\
\normEFname{RHM} & $81.48$          & $80.54$          & $66$ & 1h 35m 10s & 1h 35m 10s \\
\normEFname{SGDM} & $82.66$          & $81.70$          & $79$ & 0h 37m 54s & 0h 37m 55s \\
\textbf{\normEFname{IGT}}  & $\mathbf{83.48}$          & $\mathbf{81.77}$          & $\mathbf{70}$ & \textbf{0h 39m 03s} & \textbf{0h 39m 03s} \\
\textbf{\normEFname{HM}}  & $\mathbf{84.32}$ & $\mathbf{83.22}$ & $\mathbf{75}$ & \textbf{1h 37m 55s} & \textbf{1h 37m 56s} \\ \bottomrule
\end{tabular}}

\end{threeparttable}
\end{table*}

\begin{figure*}[t]
    \centering
    \begin{subfigure}{0.32\textwidth}
        \includegraphics[width=\linewidth]{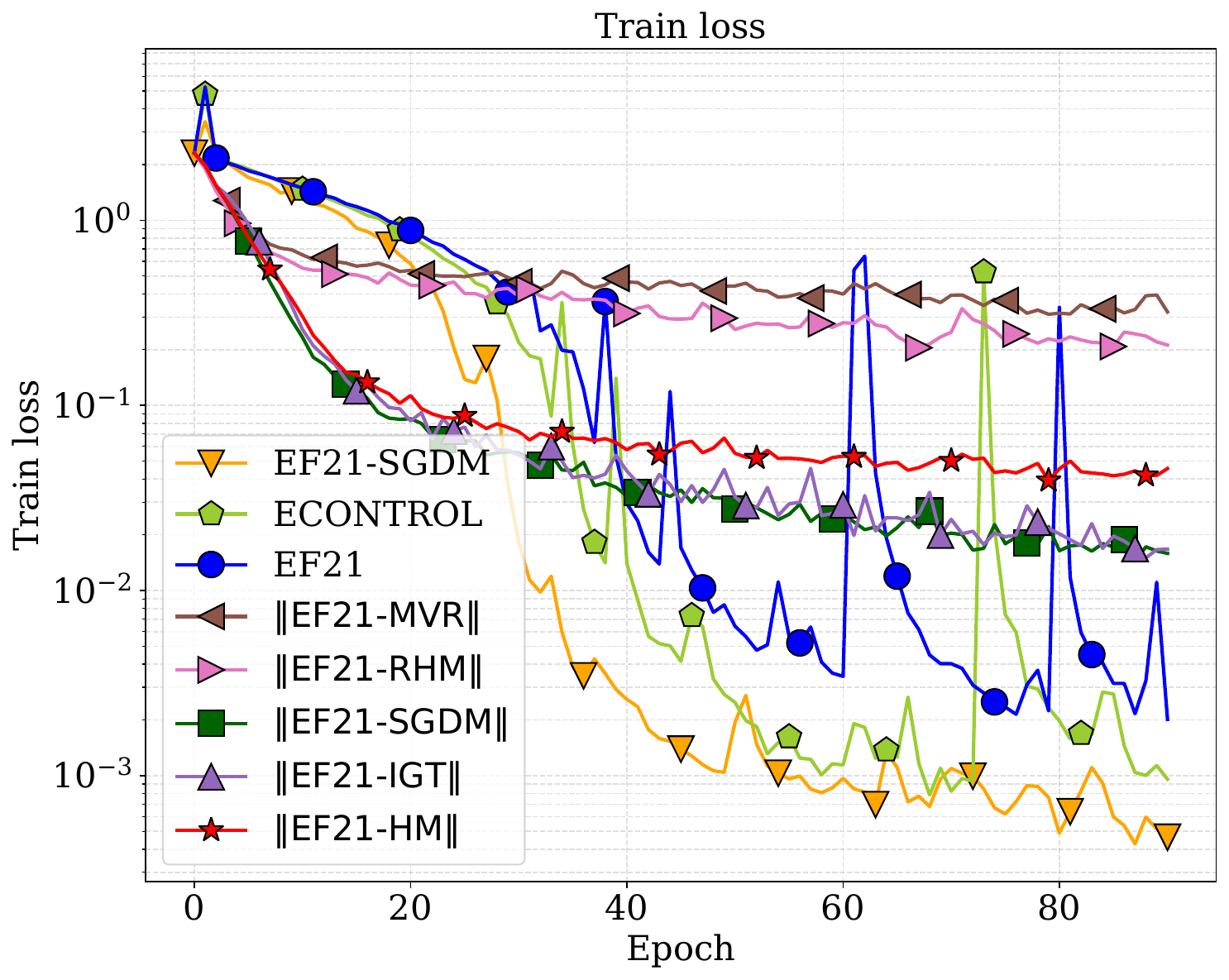}
        \caption{Training Loss}
        \label{ch7:main_fig:train_loss_epoch_chapter_7}
    \end{subfigure}
    \hfill %
    \begin{subfigure}{0.32\textwidth}
        \includegraphics[width=\linewidth]{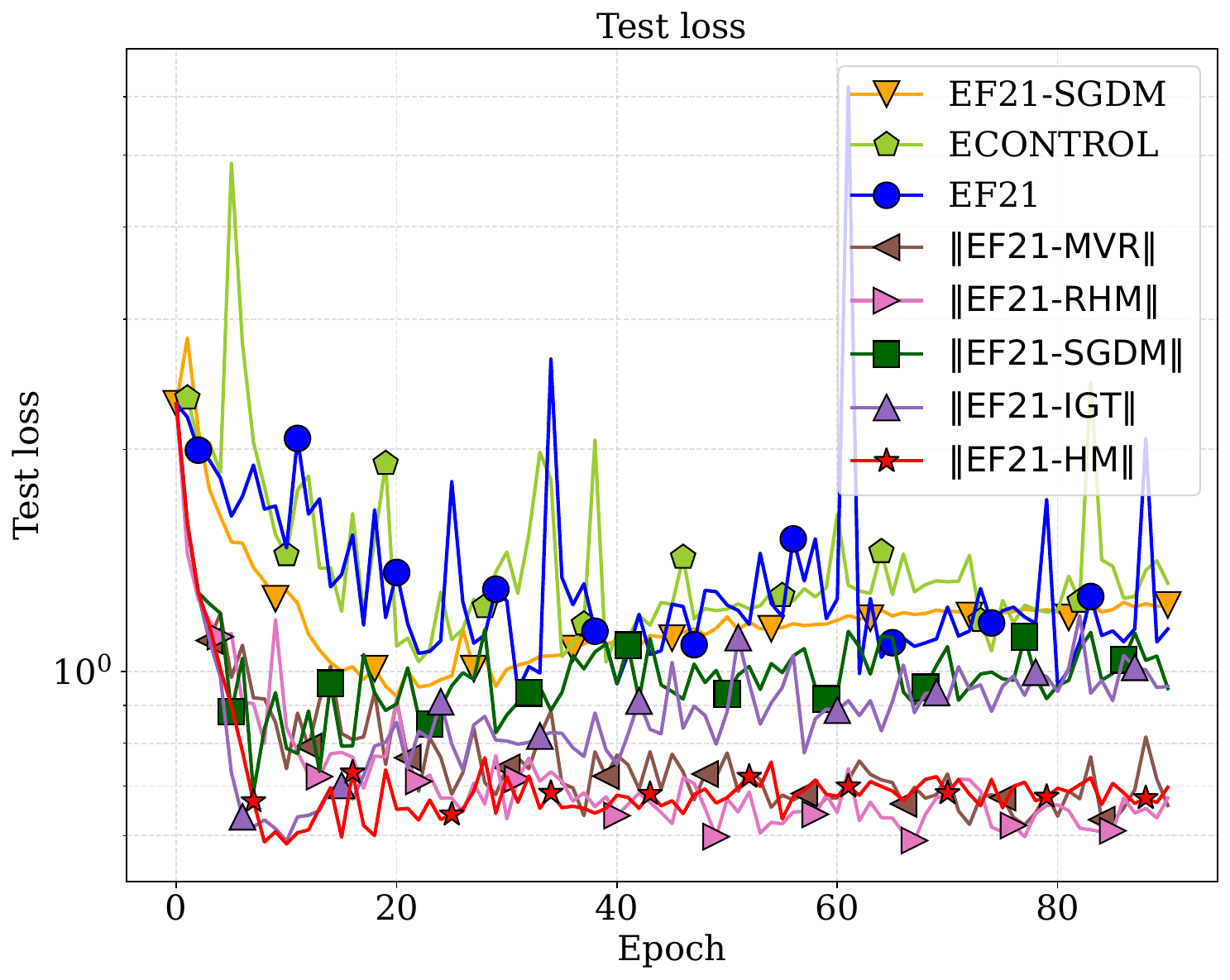}
        \caption{Test Loss}
        \label{ch7:main_fig:test_loss_epoch_chapter_7}
    \end{subfigure}
    \hfill %
    \begin{subfigure}{0.32\textwidth}
        \includegraphics[width=\linewidth]{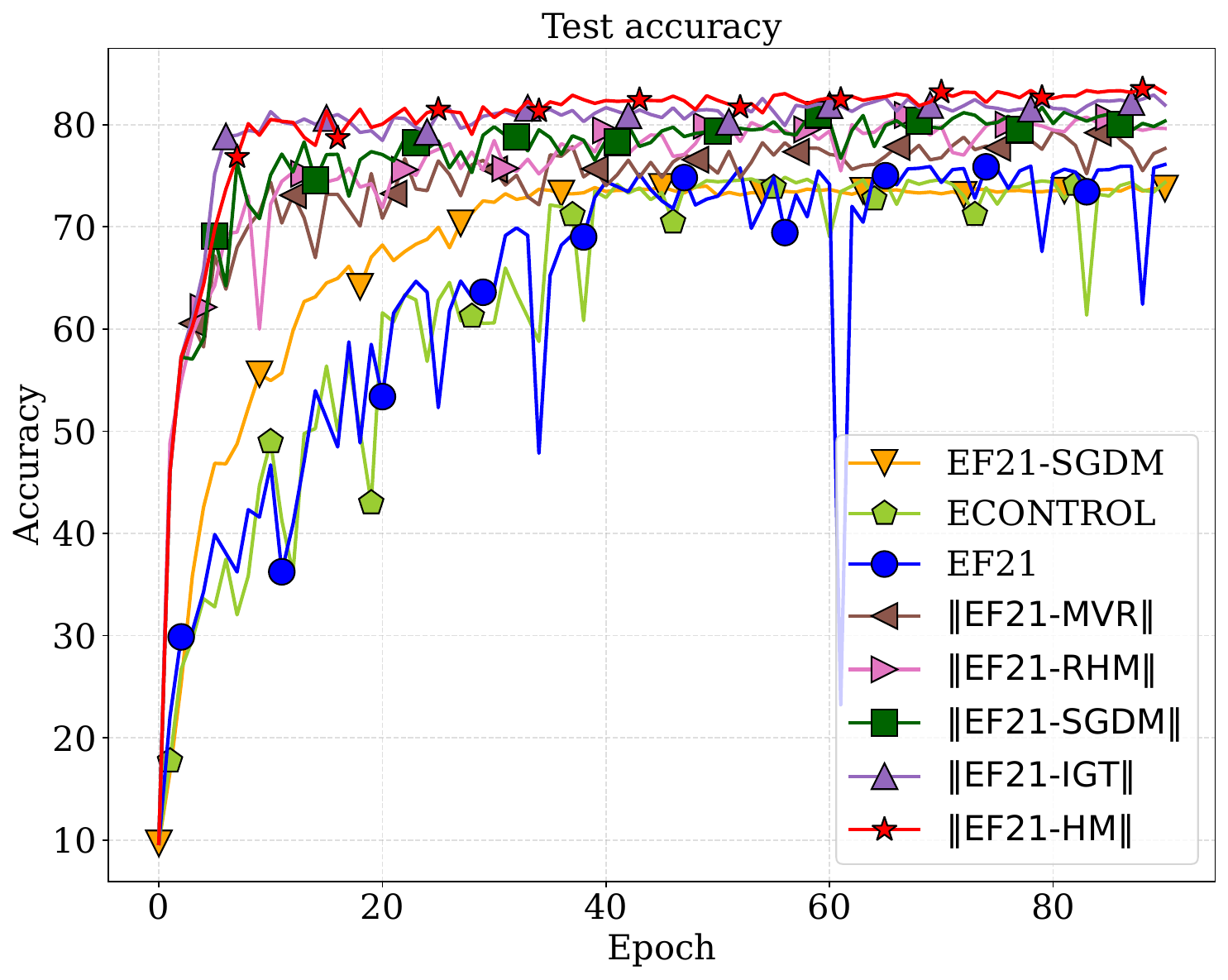}
        \caption{Test Accuracy}
        \label{ch7:main_fig:test_acc_epoch_chapter_7}
    \end{subfigure}
    \caption{Performance comparison of all methods on CIFAR-10 with ResNet-18, plotted as a function of epochs. The proposed momentum variants, particularly \normEFname{HM} and \normEFname{IGT}, show superior sample efficiency.}
    \label{ch7:main_fig:performance_plots_epoch_chapter_7}
\end{figure*}

\textbf{Implementation details.}
Implementing \normEFname{SGDM}, \normEFname{IGT} and \normEFname{MVR} is straightforward. 
\normEFname{IGT} requires one additional vector to memorize the extrapolated point, while \normEFname{MVR}  requires two stochastic gradient evaluations on the same minibatch. 
For \normEFname{HM} and \normEFname{RHM}, we can compute Hessian-vector products (HVPs) efficiently.
This is achieved via the well-known identity
\begin{align*}
    \nabla^2 f(x;\xi)\,v \;=\; \nabla_x\big\langle \nabla_x f(x;\xi),\, v \big\rangle,
\end{align*}
which allows for the computation of the HVP without materializing the full Hessian matrix, which is computationally expensive.

\textbf{$\bullet$}
\normEFname{HM}:
At each iteration $k$ and worker $i$, we form the displacement vector $\Delta^{k+1} = x^{k+1} - x^k$. The stochastic gradient $\nabla f_i(x^{k+1};\xi_i^{k+1})$ and the HVP
\[
h_v^{k+1} \gets \nabla^2 f_i(x^{k+1};\xi_i^{k+1})\,\Delta^{k+1}
\]
are computed at the new point $x^{k+1}$ using the same minibatch $\xi_i^{k+1}$.
The momentum buffer is subsequently updated using this Hessian correction:
\begin{align*}
    v_i^{k+1} \;=\; (1-\eta_k)\bigl(v_i^k + h_v^{k+1}\bigr) \;+\; \eta_k\, \nabla f_i(x^{k+1};\xi_i^{k+1}).
\end{align*}
This approach has a computational cost equivalent to two backpropagations per minibatch.

\textbf{$\bullet$}
\normEFname{RHM}: In this randomized variant, the master draws one scalar $q_k\sim\mathcal{U}(0,1)$ and broadcasts the shared interpolation point
\[
\hat{x}^{k+1} = q_k\,x^{k+1} + (1-q_k)\,x^k .
\]
At each iteration $k$, each worker $i$ defines $\Delta^{k+1} = x^{k+1} - x^k$ and draws two independent minibatches: $\hat{\xi}_i^{k+1}$ for the HVP at $\hat{x}^{k+1}$ and $\xi_i^{k+1}$ for the stochastic gradient at $x^{k+1}$. Thus,
\[
    h_v^{k+1} \;\gets\; \nabla^2 f_i(\hat{x}^{k+1};\hat{\xi}_i^{k+1})\,\Delta^{k+1}
\]
is combined with $\nabla f_i(x^{k+1};\xi_i^{k+1})$ in the RHM update. The overall computational overhead corresponds to three backpropagations per iteration and worker.

\paragraph{Epoch-based  performance.}
We observe that normalization significantly improves the per-epoch convergence performance of  distributed error feedback algorithms, in both validation and test accuracy. 
Among all methods, \normEFname{HM}  achieves the  highest accuracy, outperforming other momentum variants. 
Meanwhile, \normEFname{IGT}, \normEFname{RHM}, and \normEFname{SGDM} also perform competitively, each surpassing the $80\%$ accuracy threshold.

\paragraph{Wall-Clock Time Performance.}
Next, we evaluate the convergence performance with respect to wall-clock time.
\normEFname{IGT} achieves the best overall performance in terms of reaching the highest accuracy within the shortest time, despite all distributed error feedback algorithms being run under the same timing conditions.
While \normEFname{HM} achieves the strongest per-epoch accuracy, it requires roughly twice the wall-clock time to reach its peak performance.
This increased cost is due to its use of Hessian-vector products, which require the equivalent of two backpropagations per minibatch—compared to a single backpropagation used by  \normEFname{IGT} and \normEFname{SGDM}.

\section{Conclusion and Future Work}

This chapter combined normalized \algname{EF21} reconstruction with five momentum estimators for compressed distributed nonconvex optimization. Polyak momentum, implicit gradient transport, momentum variance reduction, and two Hessian-corrected constructions provide different ways to improve the local search vector before error-feedback compression. The schedules are parameter-agnostic in the precise sense that they do not require prior knowledge of smoothness constants or the initial optimality gap. They still retain structural parameters, base constants, compression information, and user-selected settings, and are therefore not tuning-free.

Every convergence theorem controls the expected gradient norm at the stepsize-weighted random iterate $\tilde{x}^K$, selected from the trajectory with probabilities proportional to the stepsizes. The rates are $\widetilde{\mathcal{O}}(K^{-1/4})$ for Polyak momentum, $\widetilde{\mathcal{O}}(K^{-2/7})$ for implicit gradient transport, and $\widetilde{\mathcal{O}}(K^{-1/3})$ for momentum variance reduction and the Hessian-corrected variants. They nearly match, up to logarithmic factors, the best-known problem-dependent rates in the corresponding regimes; no new lower bound or terminal-iterate guarantee is asserted.

The analysis remains in expectation and assumes synchronous full participation, contractive compression, and the stated smoothness and stochastic-oracle conditions. Hessian-corrected variants also incur additional Hessian--vector or backpropagation cost, while the experiments cover a limited training scale. Extensions to adaptive rules, (star-)convex objectives, and less restrictive participation models remain open. Chapter~\ref{ch:chapter_8} next replaces the Euclidean normalized displacement with centered-ball linear minimization oracles, transporting the estimator hierarchy to arbitrary norms without claiming a high-probability heavy-tailed theory.

\clearpage
  \endgroup

\part[Advanced Momentum and LMO Methods]{Advanced Momentum and Frank-Wolfe (LMO) Methods}

\chapter{Better LMO-based Momentum Methods with Second-Order Information}
\label{ch:chapter_8}
\thesischaptersourcenote
  {Better LMO-Based Momentum Methods with Second-Order Information}
  {khirirat2025better}
  \begingroup
  \thesischaptersetup{8}%
  
\begin{paperabstract}
  The use of momentum in stochastic optimization algorithms has shown empirical success across a range of machine learning tasks. 
Recently, a new class of stochastic momentum algorithms has emerged within the Linear Minimization Oracle (LMO) framework, thus leading to state-of-the-art methods, such as \algname{Muon}, \algname{Scion}, and \algname{Gluon}, that effectively solve deep neural network training problems.
However, traditional stochastic momentum methods  offer convergence guarantees no better than the $\mathcal{O}(1/K^{1/4})$ rate. 
While several approaches, such as Hessian-Corrected Momentum (HM), have aimed to improve this rate, their theoretical results are generally restricted to the Euclidean norm setting. 
This limitation hinders their applicability in problems,  where arbitrary norms are often required.
This chapter extends the LMO-based framework by integrating HM and provides convergence guarantees under relaxed smoothness and arbitrary norm settings. 
We establish improved convergence rates of $\mathcal{O}(1/K^{1/3})$ for HM, which achieve a faster rate than traditional momentum. 
Experimental results on {\color{blue}training Multi-Layer Perceptrons (MLPs) and transformer models (NanoGPT)} verify our theoretical observations.
\end{paperabstract}

\section{Introduction}

Stochastic momentum methods have been widely used for  solving data-intensive machine learning problems.
They are employed as the default Stochastic Gradient Descent (SGD) optimizer in open-sourced software libraries for deep neural network training, such as Pytorch~\parencite{paszke2019pytorch}.
Originating from Polyak's heavy-ball method~\parencite{polyak1964some}, stochastic momentum methods have been analyzed extensively in the literature~\parencite{yan2018unified,yu2019linear,liu2020improved,cutkosky2020momentum,hubler2024gradient}.

Recent studies have extended the convergence analysis of stochastic momentum methods to arbitrary norm settings by employing the Linear Minimization Oracle (LMO) framework~\parencite{pethick2025training,kovalev2025understanding,riabinin2025gluon}.
This framework encompasses many state-of-the-art optimization algorithms specifically used in deep neural network training, such as \algname{Muon}~\parencite{jordan2024muon} and \algname{Scion}~\parencite{pethick2025training}.
However, existing analyses of LMO-based momentum algorithms typically rely on standard smoothness assumptions, which limits their applicability to a wider range of machine learning problems.

Many modern learning problems, including distributionally robust optimization~\parencite{chen2023generalized} and deep neural network training, often violate the standard $L$-smoothness assumptions. 
Empirical studies have shown that the Hessian norm  is not uniformly upper-bounded, as required by $L$-smoothness, for neural network models, including LSTM by~\textcite{zhang2019gradient}, ResNet20 by ~\textcite{zhang2019gradient}, and transformers by~\textcite{crawshaw2022robustness}.  
These observations have motivated the development of relaxed smoothness conditions that more accurately capture the behavior of modern learning problems.
Key examples include the $(L_0,L_1)$-smoothness~\parencite{zhang2019gradient,gorbunov2024methods,chen2023generalized,vankov2024optimizing}, the $(L_0,L_1,\alpha)$-smoothness~\parencite{chen2023generalized}, and  the $\ell$-smoothness~\parencite{li2023convex}.

\algname{Gluon}, proposed by \textcite{riabinin2025gluon}, is a layer-wise implementation of LMO-based algorithms using momentum, including \algname{Muon} and \algname{Scion}, for minimizing relaxed smooth functions.
However, its convergence under arbitrary norm settings--analogous to stochastic momentum methods under Euclidean norm settings--does not surpass the $\mathcal{O}(1/K^{1/4})$ rate in terms of the expected gradient norm.
Thus, the current analysis does not provide theoretical advantage from using momentum.

To improve the convergence of Euclidean stochastic algorithms using momentum, many works have focused on modifying the gradient estimators used in the momentum updates. 
\textcite{arnold2019reducing} introduced implicit gradient transportation (IGT) as a novel gradient estimator, and subsequently, \textcite{cutkosky2020momentum} established the $\cO(1/K^{2/7})$ convergence for IGT-based momentum methods. 
This convergence was further improved by recent variants, such as \algname{STORM}~\parencite{cutkosky2019momentum} and \algname{MARS}~\parencite{yuan2024mars}, which employ stochastic gradient differences as gradient estimators. 
Other variants by~\textcite{zhang2020one,salehkaleybar2022momentum,tran2022better}
leverage second-order information to construct more accurate estimators.  
Both variants ensure the improved rate of $\cO(1/K^{1/3})$, which is proven to be optimal for minimizing smooth nonconvex functions under mild  conditions~\parencite{arjevani2023lower}. 
However, these momentum methods have been studied only under the standard smoothness and Euclidean norm setting.

\textbf{Contributions.}
{\color{blue}
To solve modern nonconvex problems, we integrate  two second-order momentum  variants~\parencite{salehkaleybar2022momentum,tran2022better} into LMO-based algorithms~\parencite{pethick2025training,kovalev2025understanding}. 
}
Our methods generalize second-order momentum methods originally analyzed under the Euclidean setting by~\textcite{salehkaleybar2022momentum,tran2022better} to broader norm choices associated with the  problems. 
{\color{blue}In particular,} we establish an $\cO(1/K^{1/3})$ convergence  in the expected gradient norm for LMO-based methods using two second-order momentum variants under relaxed smoothness assumptions in the arbitrary norm setting. This improves upon the previously known $\cO(1/K^{1/4})$ rate for LMO-based methods using Polyak momentum by~\textcite{pethick2025training,pethick2025generalized,riabinin2025gluon,kovalev2025understanding}. Also, our results match the known rates for second-order momentum methods under the standard $L$-smoothness and Euclidean setting~\parencite{salehkaleybar2022momentum,tran2022better,sadiev2025second}. A comparison of theoretical guarantees is summarized in~Table~\ref{ch8:tab:TheoreticalComparisons_chapter_8}.
{\color{blue}Moreover, our analysis can be utilized to derive the convergence of LMO-based methods using extrapolated momentum and momentum variance reduction.
}
Finally, we evaluate LMO-based methods with second-order momentum and other momentum updates on solving logistic regression problems, Multi-Layer Perceptron (MLP) training problems, {\color{blue} Long Short-Term Memory (LSTM) network training problems, and transformer training problems}. 
{\color{blue}
Specifically, our experiments on NanoGPT using the WikiText-2 and TinyShakespeare datasets highlight that our proposed second-order momentum methods attain faster convergence and lower final loss than LMO-based methods using traditional Polyak momentum.
}

\begin{table*}[h!]

\caption[Convergence comparison for LMO-based momentum methods]{Comparisons of convergence for
LMO-based methods that utilize Polyak momentum such as \algname{Gluon}, second-order momentum, extrapolated momentum (implicit gradient transportation), and momentum variance reduction. We use the following abbreviations: {\bf GS} = \steponeonechange{Generalized Smoothness}, {\bf SOGS} = \steponeonechange{Second-Order Generalized Smoothness}, and {\bf MSGS} = \steponeonechange{Mean-Squared Generalized Smoothness}. \steponeonechange{Notation in the comparison columns is normalized to the thesis conventions; the cited results retain their original assumptions.}
}
\label{ch8:tab:TheoreticalComparisons_chapter_8}
\begin{center}
\renewcommand{\arraystretch}{1.0}
\resizebox{0.95\textwidth}{!}{%
\begin{tabular}{cccccc}
\bf Methods & \bf Rate & \bf GS & \bf SOGS & \bf MSGS  & \bf Norm Type \\ \hline
\begin{tabular}{c} \algname{uSCG} \\ \tabcite{pethick2025training} \end{tabular} & $\cO\left( \nicefrac{1}{K^{1/4}} \right)$ & $L$ & -- & -- & Arbitrary   \\ 
\begin{tabular}{c} \algname{GGNC} \\ \tabcite{pethick2025generalized} \end{tabular} & $\cO\left( \nicefrac{1}{K^{1/4}} \right)$ & $(L_0,L_1)$ & -- & -- & Arbitrary   \\ 
\begin{tabular}{c} \algname{Gluon} \\ \tabcite{riabinin2025gluon} \end{tabular} & $\cO\left( \nicefrac{1}{K^{1/4}} \right)$ & $(L_0,L_1)$ & -- & -- & Arbitrary   \\ 
\begin{tabular}{c} \algname{LMO-IGT} \\ \tabcite{kovalev2025understanding} \end{tabular} &  $\cO\left( \nicefrac{1}{K^{2/7}} \right)$  & $L$ & $M$& -- & Arbitrary  \\
\begin{tabular}{c} \algname{LMO-MVR}  \\ \tabcite{kovalev2025non} \end{tabular} &  $\cO\left( \nicefrac{1}{K^{1/3}} \right)$  & $L$ & --& $\steponeonechange{\cL}$  & Arbitrary  \\
\begin{tabular}{c} \algname{LiMuon}\\ \tabcite{huang2025limuon} \end{tabular} &  $\cO\left( \nicefrac{1}{K^{1/3}} \right)$  & $(L_0,L_1)$ & --& $(\steponeonechange{\cL_0},\steponeonechange{\cL_1})$ & Euclidean  \\
\begin{tabular}{c} \algname{SHARP}  \\ \tabcite{salehkaleybar2022momentum} \end{tabular} &  $\cO\left( \nicefrac{1}{K^{1/3}} \right)$  & $L$ & --& --  & Euclidean  \\
\begin{tabular}{c} \algname{N-SGDHess}  \\ \tabcite{tran2022better} \end{tabular} &  $\cO\left( \nicefrac{1}{K^{1/3}} \right)$  & $L$ & $M$& --  & Euclidean  \\ \hline
 \rowcolor{linen} \begin{tabular}{c} \algname{LMO-IGT}  \\ Theorem~\ref{ch8:thm:IGT_relaxed_smoothness_chapter_8} \end{tabular} & $\cO\left( \nicefrac{1}{K^{2/7}} \right)$ & $(L_0,L_1)$ & $(M_0,M_1)$ & -- & Arbitrary   \\
 \rowcolor{linen} \begin{tabular}{c} \algname{LMO-MVR}  \\ Theorem~\ref{ch8:thm:MVR_relaxed_smoothness_chapter_8} \end{tabular} & $\cO\left( \nicefrac{1}{K^{1/3}} \right)$ & $(L_0,L_1)$ & -- & $(\steponeonechange{\cL_0}, \steponeonechange{\cL_1})$ & Arbitrary   \\  \hline 
 \rowcolor{linen} \begin{tabular}{c} \algname{LMO-SOM-V1}  \\ Theorem~\ref{ch8:thm:V1_relaxed_smoothness_chapter_8} \end{tabular}
&  $\cO\left( \nicefrac{1}{K^{1/3}} \right)$ & $(L_0,L_1)$ & --& --  & Arbitrary   \\
 \rowcolor{linen} \begin{tabular}{c} \algname{LMO-SOM-V2}  \\ Theorem~\ref{ch8:thm:V2_relaxed_smoothness_chapter_8} \end{tabular} & $\cO\left( \nicefrac{1}{K^{1/3}} \right)$ & $(L_0,L_1)$ & $(M_0,M_1)$ & -- & Arbitrary   \\
 \hline 
\end{tabular}
}
\end{center}

\end{table*}

\section{Related Work}

\textbf{Stochastic momentum methods.}
Polyak momentum is widely used  to expedite the training process in stochastic optimization methods.
In the Euclidean norm setting, stochastic  methods with Polyak momentum enjoy the $\cO(1/K^{1/4})$ convergence under standard smoothness~\parencite{yan2018unified,yu2019linear,liu2020improved,cutkosky2020momentum}, and recently under relaxed smoothness~\parencite{hubler2024gradient}. 
Moreover, Polyak momentum has been incorporated to enhance the performance of other methods, such as distributed methods with communication compression for communication-efficient learning~\parencite{fatkhullin2023momentum,khirirat2024error}, and LMO-based optimization methods for deep neural network training~\parencite{jordan2024muon,pethick2025training}.
To improve upon the $\cO(1/K^{1/4})$ convergence rate, several novel momentum updates have been proposed. 
For instance, stochastic methods employing IGT momentum (or extrapolated momentum)~\parencite{arnold2019reducing} have been shown by~\textcite{cutkosky2020momentum} to achieve an improved rate of $\cO(1/K^{2/7})$. Furthermore, the methods based on \algname{STORM}~\parencite{cutkosky2019momentum}, \algname{MARS}~\parencite{yuan2024mars}, and Hessian-corrected momentum~\parencite{salehkaleybar2022momentum,tran2022better,zhang2020improved} attain the convergence rate of $\cO(1/K^{1/3})$, which matches the known lower bound established by~\textcite{arjevani2023lower}.
However, the convergence guarantees of these stochastic  methods using novel momentum techniques were restricted to the Euclidean norm setting. 
This chapter incorporates Hessian-corrected momentum into LMO-based optimization methods to guarantee $\cO(1/K^{1/3})$ convergence under the arbitrary-norm setting, improving upon the $\cO(1/K^{1/4})$ convergence achieved by LMO-based methods using Polyak momentum.

\textbf{\algname{Muon}.}
\algname{Muon} optimizers~\parencite{jordan2024muon} update the model parameters (with their inherent matrix structure) with orthogonalized gradient momentum.
Empirical studies, e.g. by~\textcite{jordan2024muon,liu2025muon}, highlight the superior performance of \algname{Muon} over Shampoo~\parencite{gupta2018shampoo}, SOAP~\parencite{vyas2024soap}, and AdamW~\parencite{loshchilov2017decoupled}, for language model training.
From these encouraging observations, the convergence of \algname{Muon} was analyzed initially by~\textcite{li2025note} under the standard smoothness with respect to the Frobenius norm. 
Later, LMO-based optimization methods with Polyak momentum were proposed by~\textcite{pethick2025training} to capture optimizers, including \algname{Muon}, \algname{Scion}, and stochastic Euclidean methods with Polyak momentum.  
Under the arbitrary norm setting, the methods were shown to enjoy the $\cO(1/K^{1/4})$ convergence under the standard smoothness~\parencite{pethick2025training,kovalev2025understanding}, and under the relaxed smoothness~\parencite{pethick2025generalized,riabinin2025gluon}.

\textbf{Relaxed smoothness.}
Training deep neural networks, especially large-scale language models, poses significant challenges due to their highly non-convex functions. To better capture this, a variety of relaxed smoothness conditions on the functions have been proposed. One early contribution by~\textcite{zhang2019gradient} introduced a function class in which the Hessian norm grows linearly with the gradient norm, providing a more flexible alternative to classical smoothness assumptions. Building upon this, several works have introduced broader classes of functions, such as the symmetric $(L_0, L_1, \alpha)$-smoothness condition~\parencite{chen2023generalized}, the $\ell$-smoothness condition~\parencite{li2023convex}, and the $(K_0,K_1,\alpha)$-smoothness condition~\parencite{liu2025theoretical}. These conditions aim to more accurately reflect the optimization landscape encountered in deep learning. Under such relaxed assumptions, the convergence behavior of various gradient-based methods has been extensively studied~\parencite{crawshaw2022robustness, hubler2024gradient, zhang2019gradient, koloskova2023revisiting, chezhegov2025convergence, gorbunov2024methods, vankov2024optimizing, khirirat2024error, faw2023beyond,wang2023convergence,zhao2021convergence,hubler2024parameter}, providing theoretical insights that align more closely with empirical observations in neural network  training.

\section{Problem Formulation}
Throughout this chapter, we consider the following unconstrained stochastic optimization problem
\begin{eqnarray}\label{ch8:eqn:problem_chapter_8}
    \underset{x\in\R^d}{\min} \ f(x) := \thesischeightExpSub{\xi \sim \thesischeightcD}{f(x;\xi)},
\end{eqnarray}
where $f(x;\xi)$ is a differentiable but possibly nonconvex function,  $x\in\R^d$ is the $d$-dimensional vector of model parameters, and $\xi$ is a random vector, which represents a training data sample drawn from an unknown data distribution $\thesischeightcD$.
Often, we  access a stochastic oracle to compute the stochastic gradient and Hessian of $f$ at $x$ with respect to $\xi\sim \thesischeightcD$ denoted by $\nabla f(x;\xi)$ and $\nabla^2 f(x;\xi)$, respectively. 
We  assume that the stochastic gradient and stochastic Hessian satisfy the following unbiased and variance-bounded properties: 
\begin{assumption}\label{ch8:assum:Hessian_boundedvariance_v2_chapter_8}
 The stochastic gradient $\nabla f(x;\xi)$ and the stochastic Hessian $\nabla^2 f(x;\xi)$ are unbiased estimators of $\nabla f(x)$ and $\nabla^2 f(x)$ with bounded variance, i.e., for all $x,w\in\R^d$,
\begin{align*}
& \thesischeightExpSub{\xi}{\nabla f(x;\xi)}= \nabla f(x), \quad  \thesischeightExpSub{\xi}{\normtwo{\nabla f(x;\xi)-\nabla f(x)}^2} \leq \sigma^2_g, \quad  \\  &
\thesischeightExpSub{\xi}{\nabla^2 f(x;\xi)}= \nabla^2 f(x), \quad \text{and} \quad  \thesischeightExpSub{\xi}{\normtwo{(\nabla^2 f(x;\xi)-\nabla^2 f(x))w}^2} \leq \sigma^2_H \normtwo{w}^2. 
\end{align*}
\end{assumption}

Assumption~\ref{ch8:assum:Hessian_boundedvariance_v2_chapter_8} is commonly used for analyzing the convergence of optimization algorithms using stochastic gradients~\parencite{ghadimi2013stochastic,cutkosky2019momentum,gorbunov2020linearly} and stochastic Hessians~\parencite{tran2022better,jiang2024stochastic}.

We further impose the following standard assumptions on objective functions.

\begin{assumption}\label{ch8:assum:lower_bound_chapter_8}
The function $f:\R^d\rightarrow \R$ is bounded from below, i.e., $f^{\inf}=\inf_{x\in\R^d} f(x) > -\infty$.    
\end{assumption}

\begin{assumption}[\steponeonechange{Generalized Smoothness}]\label{ch8:assum:Relaxed_Lipschitz_Grad_chapter_8}
The function $f:\R^d\rightarrow\R$ is differentiable, and its gradient $\nabla f(x)$ is symmetrically $(L_0,L_1)$-Lipschitz continuous with respect to the norm $\norm{\cdot}$, i.e. for all $x,y\in\R^d$, 
\begin{align*}
     \normstar{\nabla f(x) - \nabla f(y)}  \leq (L_0 + L_1 \underset{\theta \in [0,1]}{\sup} \normstar{\nabla f(w_{\theta}(x,y))}) \norm{x-y},
\end{align*}
where $w_{\theta}(x,y)=\theta x + (1-\theta)y$.
\end{assumption}

{\color{black}
\begin{assumption}[\steponeonechange{Mean-Squared Generalized Smoothness}]
\label{ch8:assum:Relaxed_Mean_Squared_Lipschitz_Grad_chapter_8}
    The stochastic gradients of $f$ are mean-squared symmetrically $(\steponeonechange{\cL_0}, \steponeonechange{\cL_1})$-Lipschitz continuous with respect to the norm $\|\cdot\|$, i.e., for all $x, y \in \R^d$,
    \begin{align*}
        \MoveEqLeft \thesischeightExpSub{\xi}{\normstar{\nabla f(x;\xi)-\nabla f(y;\xi)}^2}\\
        &\leq \thesischeightExpSub{\xi}{\left(\steponeonechange{\cL_0}+\steponeonechange{\cL_1}\sup_{\theta \in [0,1]}\normstar{\nabla f(w_{\theta}(x,y);\xi)}\right)^2} \norm{x-y}^2,
    \end{align*}
    where $w_{\theta}(x,y)=\theta x+ (1-\theta)y$.
\end{assumption}
}

We refer to the condition in Assumption~\ref{ch8:assum:Relaxed_Mean_Squared_Lipschitz_Grad_chapter_8} as the \textit{\steponeonechange{mean-squared generalized smoothness assumption}} \parencite{chen2023generalized}. As shown by \textcite{chen2023generalized}, it is a stronger condition than the generalized smoothness in Assumption~\ref{ch8:assum:Relaxed_Lipschitz_Grad_chapter_8}.

\begin{assumption}[\steponeonechange{Second-Order Generalized Smoothness}]\label{ch8:assum:Relaxed_Lipschitz_Hessian_chapter_8}
The function $f:\R^d\rightarrow \R$ is twice differentiable, and its Hessian $\nabla^2 f(x)$ is symmetrically $(M_0,M_1)$-Lipschitz continuous with respect to the norm $\norm{\cdot}$, i.e. for all $x,y\in\R^d$, 
\begin{align*}
    \normstar{\nabla^2 f(x) - \nabla^2 f(y)}   \leq (M_0 + M_1 \underset{\theta \in [0,1]}{\sup} \normstar{ \nabla f( w_{\theta}(x,y))}) \norm{x-y},
\end{align*}
where $w_\theta(x,y) = \theta x + (1-\theta)y$.
\end{assumption}

Assumptions~\ref{ch8:assum:lower_bound_chapter_8},~\ref{ch8:assum:Relaxed_Lipschitz_Grad_chapter_8},~\ref{ch8:assum:Relaxed_Mean_Squared_Lipschitz_Grad_chapter_8} and~\ref{ch8:assum:Relaxed_Lipschitz_Hessian_chapter_8} ensure the lower-bound of the objective function $f$, \steponeonechange{generalized smoothness of its gradient} $\nabla f$,
{\color{black}\steponeonechange{mean-squared generalized smoothness}},
and \steponeonechange{second-order generalized smoothness of its Hessian} $\nabla^2 f$, respectively.

Assumption~\ref{ch8:assum:Relaxed_Lipschitz_Grad_chapter_8} is a generalization of the symmetric $(L_0,L_1)$-Lipschitz continuity of $\nabla f(\cdot)$ with respect to the Euclidean norm by~\textcite{gorbunov2024methods,chen2023generalized}. 
Also, Assumption~\ref{ch8:assum:Relaxed_Lipschitz_Grad_chapter_8}  recovers $L$-Lipschitz continuity with respect to the arbitrary norm by~\textcite{kovalev2025understanding,pethick2025training} when we let $L_0=L$ and $L_1=0$.
{\color{black}
Moreover,~\steponeonechange{the mean-squared generalized-smoothness condition in} Assumption~\ref{ch8:assum:Relaxed_Mean_Squared_Lipschitz_Grad_chapter_8} generalizes the standard mean-squared smoothness assumption proposed by \textcite{arjevani2023lower}.
}
\stepEightBchange{Assumption~\ref{ch8:assum:Relaxed_Lipschitz_Hessian_chapter_8} is related to the $(M_0,M_1)$-type Hessian regularity considered by~\textcite[Assumption~3]{xie2024trust}. Here, the condition is formulated symmetrically in the selected norm and its dual, consistently with the geometry used throughout this chapter. When $M_1=0$ and $M_0=L_H$, Assumption~\ref{ch8:assum:Relaxed_Lipschitz_Hessian_chapter_8} reduces to standard $L_H$-Lipschitz continuity of the Hessian.}

Finally, note that the variance (Assumption~\ref{ch8:assum:Hessian_boundedvariance_v2_chapter_8}) is measured with respect to the Euclidean norm, while \steponeonechange{generalized smoothness} (Assumption~\ref{ch8:assum:Relaxed_Lipschitz_Grad_chapter_8}), {\color{black}\steponeonechange{mean-squared generalized smoothness} (Assumption~\ref{ch8:assum:Relaxed_Mean_Squared_Lipschitz_Grad_chapter_8}),}
and \steponeonechange{second-order generalized smoothness} (Assumption~\ref{ch8:assum:Relaxed_Lipschitz_Hessian_chapter_8}) are formulated with respect to the arbitrary norm and its dual norm. 
Hence, we describe a connection between these norms using the following inequality: 
\begin{align}
{\color{red}\thesischeightubar{\rho}}\normtwo{x} \leq \normstar{x} \leq {\color{red}\bar\rho} \normtwo{x}, \quad \text{and}   \quad 
 {\color{blue}\thesischeightubar{\theta}} \normtwo{x} \leq  \norm{x} \leq {\color{blue}\bar \theta} \normtwo{x}.\label{ch8:eqn:norm_diff_chapter_8}
\end{align}

\subsection{LMO-based  Methods with Momentum}\label{ch8:sec:LMO_chapter_8}

To solve the problem in~\eqref{ch8:eqn:problem_chapter_8} for neural network training tasks, we focus on the LMO-based optimization methods using Polyak momentum~\parencite{pethick2025training,kovalev2025understanding}, which update the iterates $\{x^{k}\}_{k \geq 0}$ as follows: Given $x^{0},m^{0} \in\R^d$, 
\begin{align}
    x^{k+1} &= x^{k} + \text{lmo}(m^{k}), \quad \text{where}\label{ch8:eqn:LMO_momentum_chapter_8}\\ 
   m^{k+1} &= (1-\eta_k)m^{k} + \eta_k \nabla f(x^{k+1};\xi^{k+1}).\label{ch8:eqn:LMO_momentum_m_chapter_8} 
\end{align}
Here,  $\gamma_k>0$ is a stepsize, $\eta_k \in (0,1)$ is a momentum parameter,  $\nabla f(x;\xi)$ is the stochastic gradient, and $\text{lmo}(m^{k}) \eqdef {\text{argmin}}_{ \norm{x} \leq \gamma_k } \inner{m^{k}}{x}$. 
The methods in~\eqref{ch8:eqn:LMO_momentum_chapter_8} encompass a class of stochastic momentum methods under specific norm choices. 
For instance, they recover normalized stochastic momentum methods~\parencite{cutkosky2020momentum} when we let $\norm{\cdot} = \normtwo{\cdot}$, and    sign stochastic momentum methods~\parencite{jiang2025improved} when we let $\norm{\cdot} = \norm{\cdot}_{\infty}$.

The LMO-based methods using Polyak momentum in~\eqref{ch8:eqn:LMO_momentum_chapter_8} were  shown by~\textcite{pethick2025training,kovalev2025understanding,riabinin2025gluon} to attain the $\cO\left( 1/K^{1/4} \right)$ convergence rate in the gradient norm  under the arbitrary norm setting, The obtained rate is the same as that of stochastic momentum methods under the Euclidean norm setting, e.g. by~\textcite{cutkosky2020momentum,jiang2025improved,zhang2020improved,hubler2024parameter}.

To further improve the iteration complexity of the momentum methods in~\eqref{ch8:eqn:LMO_momentum_chapter_8}, we may use three momentum variants: extrapolated momentum, \algname{STORM}, and second-order momentum.
First, extrapolated momentum~\parencite{arnold2019reducing,cutkosky2020momentum}   updates the momentum vector $m^{k}$ according to: 
\begin{eqnarray}\label{ch8:eqn:extrapolated_m_chapter_8}
   m^{k+1} = (1-\eta_k)m^{k} + \eta_k \nabla f\left(y^{k+1};\xi^{k+1}\right),
\end{eqnarray}
where $y^{k+1} =  x^{k+1} + (\nicefrac{(1-\eta_k)}{\eta_k} )(x^{k+1} - x^{k})$. 
LMO-based methods using extrapolated momentum enjoys the $\cO(1/K^{2/7})$ convergence rate under the standard smoothness and arbitrary norm setting, as shown by~\textcite{kovalev2025understanding}, which matches the rate under the Euclidean norm setting by~\textcite[Theorem 3]{cutkosky2020momentum} and by~\textcite[Theorem 5]{cutkosky2021high} in the bounded variance case.

To further improve the convergence of stochastic methods using extrapolated momentum, we may consider \algname{STORM}~\parencite{cutkosky2019momentum}. The key technique is to use the stochastic gradient difference $\nabla f(x^{k+1};\xi^{k+1}) - \nabla f(x^{k};\xi^{k+1})$ as the correction term to expedite the training process. In \algname{STORM}, the momentum vector $m^{k}$ is updated via: 
\begin{align}
m^{k+1} =& (1-\eta_k)(m^{k} + d_{\xi^{k+1}}(x^{k+1},x^{k})) + \eta_k \nabla f\left(x^{k+1};\xi^{k+1}\right),    
\end{align}
where $d_{\xi^{k+1}}(x^{k+1},x^{k}) = \nabla f(x^{k+1};\xi^{k+1}) - \nabla f(x^{k};\xi^{k+1})$. 
LMO-based methods using \algname{STORM} have been recently shown to achieve the $\cO(1/K^{1/3})$ convergence rate under the relaxed smoothness and Frobenius norm setting recently by~\textcite{huang2025limuon}, which recovers the rate under the standard smoothness and Euclidean norm setting by~\textcite{cutkosky2019momentum}.
To better control the gradient correction term in \algname{STORM},~\textcite{yuan2024mars} introduce a scaling parameter $\rho_k \in (0,1)$. 
The resulting algorithm, named \algname{MARS}, modifies \algname{STORM} by adjusting the update rule for the momentum vector $m^{k}$ as follows: 
\begin{align}
m^{k+1} =& (1-\eta_k) \left(m^{k} + {\color{orange} \frac{\rho_k}{1-\eta_k}}d_{\xi^{k+1}}(x^{k+1},x^{k}) \right) + \eta_k \nabla f\left(x^{k+1};\xi^{k+1}\right).
\label{ch8:eq:MARS_chapter_8}
\end{align}
By tuning the scaling factor $\rho_k$ properly, \algname{MARS} can provide finer control over the gradient difference term, thus leading to improved convergence performance. 
Note that \algname{MARS} with $\rho_k =1-\eta_k$ becomes \algname{STORM}.

Finally, other variants that achieve the same rate as \algname{STORM} are second-order momentum~\parencite{salehkaleybar2022momentum,tran2022better,zhang2020one} (or Hessian-corrected momentum). 
They utilize the Hessian-vector product as the correction term, and perform  the following update: 
\begin{align}
    m^{k+1} =& (1-\eta_k)(m^{k} + \nabla^2 f(\hat x^{k+1};\xi^{k+1})(x^{k+1}-x^{k}))  + \eta_k \nabla f(x^{k+1};\xi^{k+1}).\label{ch8:eqn:second_order_momentum_update_chapter_8}
\end{align}
The second-order momentum update in~\eqref{ch8:eqn:second_order_momentum_update_chapter_8} becomes the one proposed by~\textcite{tran2022better} when $\hat x^{k+1} = x^{k+1}$, and the one by~\textcite{salehkaleybar2022momentum} when $\hat x^{k+1} = b_{k+1} x^{k+1} + (1-b_{k+1})x^{k}$ with $b_{k+1} \sim \thesischeightcU(0,1)$.
Both \algname{STORM} and second-order momentum ensure the $\cO(1/K^{1/3})$ convergence, which is faster than the $\cO(1/K^{2/7})$ convergence of extrapolated momentum and the $\cO(1/K^{1/4})$ convergence of Polyak momentum.

\section{LMO-based Methods with Second-order Momentum}

Now, we consider LMO-based methods~\parencite{pethick2025training,kovalev2025understanding} using  second-order momentum~\parencite{salehkaleybar2022momentum,tran2022better,zhang2020one}.
The methods leverage the LMO-based oracle in~\eqref{ch8:eqn:LMO_momentum_chapter_8},  
second-order momentum updates in~\eqref{ch8:eqn:second_order_momentum_update_chapter_8}, and the scaling factor $\rho_k$ introduced in \algname{MARS}~\parencite{yuan2024mars}.
See Algorithm~\ref{ch8:alg:2_order_momentum_chapter_8} for the detailed description of the methods. 

In~Algorithm~\ref{ch8:alg:2_order_momentum_chapter_8}, the second-order momentum incorporates the Hessian-vector product term into its update to accelerate the convergence of the LMO-based momentum methods. 
The Hessian-vector product requires roughly the same computational time as the gradient~\parencite{pearlmutter1994fast}.
This can be achieved by using the automatic differentiation package to compute $\nabla^2 f(x) v$ for $x,v\in\R^d$ by evaluating $\nabla h(x)$, where  $h(x) = \inner{\nabla f(x)}{v}$~\parencite{tran2022better}.
At iteration $k$, the analysis uses one fresh sample $\xi^{k+1}\sim\thesischeightcD$ for both the stochastic gradient and the stochastic Hessian-vector product. For Variant~1, $b_{k+1}\sim\thesischeightcU(0,1)$ is also drawn freshly and independently of the past and of $\xi^{k+1}$. Thus, this chapter analyzes the shared-sample construction; independent gradient and Hessian samples are not covered by the present proof.

\begin{algorithm}
\caption{\algname{LMO-SOM} (LMO-based Optimization Methods with Second-Order Momentum)}
\label{ch8:alg:2_order_momentum_chapter_8}
\begin{algorithmic}[1]
\STATE \textbf{Input:} Tuning parameters $\gamma_k>0$ and $\eta_k , {\color{orange} \rho_k} \in (0,1)$ for $k=0,1,\ldots$; Initial points $x^{0} , m^{0} \in \R^d$; Horizon $K\geq 1$
\FOR{each iteration $k = 0, 1, \dots, K$}
        \STATE Update $x^{k+1} = x^{k} + \text{lmo}(m^{k})$, where $\text{lmo}(m^{k}) \eqdef \underset{ \norm{x} \leq \gamma_k }{\text{argmin}} \inner{m^{k}}{x}$, and draw a fresh sample $\xi^{k+1} \sim \thesischeightcD$
        \IF{Variant 1}
        \STATE Set $\hat x^{k+1}=b_{k+1} x^{k+1} + (1-b_{k+1})x^{k} $, where $b_{k+1} \sim \thesischeightcU(0,1)$
        \STATE Update $m^{k+1} =  (1-\eta_k)\left(m^{k} + {\color{orange} \frac{\rho_k}{1-\eta_k}}\nabla^2 f(\hat x^{k+1};\xi^{k+1}) (x^{k+1} - x^{k}) \right) + \eta_k \nabla f\left(x^{k+1};\xi^{k+1}\right)$
        \ELSIF{Variant 2}
        \STATE Update $m^{k+1} =  (1-\eta_k)\left( m^{k} + {\color{orange} \frac{\rho_k}{1-\eta_k}}\nabla^2 f(x^{k+1};\xi^{k+1}) (x^{k+1} - x^{k}) \right) + \eta_k \nabla f\left(x^{k+1};\xi^{k+1}\right)$
\ENDIF
\ENDFOR
\STATE \textbf{Output:} Generated trajectory $(x^0,x^1,\ldots,x^{K+1})$
\end{algorithmic}
\end{algorithm}

Algorithm~\ref{ch8:alg:2_order_momentum_chapter_8} uses a horizon $K\geq1$ and performs $K+1$ updates, indexed by $k=0,\ldots,K$. The theorems below report the analytical best-iterate criterion $\min_{0\leq k\leq K}\Exp{\normstar{\nabla f(x^k)}}$; this is neither a terminal-iterate guarantee nor an implementable argmin unless the exact gradient norms are available.

To this end, we provide the convergence theorems for~Algorithm~\ref{ch8:alg:2_order_momentum_chapter_8} with $\rho_k=1-\eta_k$.

\begin{theorem}[\algname{LMO-SOM-V1}: Variant 1 of~Algorithm~\ref{ch8:alg:2_order_momentum_chapter_8}]\label{ch8:thm:V1_relaxed_smoothness_chapter_8} Consider the problem of minimizing $f(x)=\thesischeightExpSub{\xi \sim \thesischeightcD}{f(x;\xi)}$. Let $f$ be twice differentiable, let $K\geq1$ be an integer horizon, and let Assumptions~\ref{ch8:assum:Hessian_boundedvariance_v2_chapter_8},~\ref{ch8:assum:lower_bound_chapter_8}, and~\ref{ch8:assum:Relaxed_Lipschitz_Grad_chapter_8} hold. 
Then, the iterates $\{x^{k}\}$ generated by Variant 1 of Algorithm~\ref{ch8:alg:2_order_momentum_chapter_8} with $\rho_k=1-\eta_k$ and with
\begin{equation*}
    \eta_k =  \eta = \frac{1}{(K+1)^{2/3}} \quad \text{and} \quad  \gamma_k= \gamma = \frac{\hat \gamma}{(K+1)^{2/3}},
\end{equation*}
where $\hat \gamma =\frac{1}{80L_1} \left({\color{red}\bar\rho}/{\color{red}\thesischeightubar\rho}\right)^{-1}$ satisfy  
\begin{align*}
 & \min_{k \in \{0,1,\ldots,K\}} \Exp{\normstar{\nabla f(x^{k})}}
\leq  \frac{2 \Delta_0}{\hat \gamma (K+1)^{\nicefrac{1}{3}}} +  \frac{4\left( \Exp{\normstar{r^{0}}} + {\color{red}\bar\rho}\sigma_g\right)}{(K+1)^{\nicefrac{1}{3}}} +  \frac{\hat \gamma L_0 \exp(L_1 \hat\gamma)}{(K+1)^{\nicefrac{2}{3}}}  \\
& \qquad   + 12\frac{{\color{red}\bar\rho}}{{\color{red} \thesischeightubar \rho}} \frac{  \hat\gamma L_0 }{(K+1)^{ \nicefrac{1}{3} }}   \left(\exp( L_1 \hat\gamma) + \sqrt{L_1 \hat \gamma} \exp( L_1 \hat\gamma) + 2\right)  + 8\frac{ {\color{red}\bar\rho}}{{\color{blue}\thesischeightubar\theta}} \frac{ \hat\gamma \sigma_H}{(K+1)^{ \nicefrac{1}{3} }},
\end{align*}
where $\Delta_0=f(x^{0}) - f^{\inf}$ and $r^0=m^0-\nabla f(x^0)$.
\end{theorem}

\begin{theorem}[\algname{LMO-SOM-V2}: Variant 2 of~Algorithm~\ref{ch8:alg:2_order_momentum_chapter_8}]\label{ch8:thm:V2_relaxed_smoothness_chapter_8}
 Consider the problem of minimizing $f(x)=\thesischeightExpSub{\xi \sim \thesischeightcD}{f(x;\xi)}$. Let $f$ be twice differentiable, let $K\geq1$ be an integer horizon, and let Assumptions~\ref{ch8:assum:Hessian_boundedvariance_v2_chapter_8},~\ref{ch8:assum:lower_bound_chapter_8},~\ref{ch8:assum:Relaxed_Lipschitz_Grad_chapter_8}, and~\ref{ch8:assum:Relaxed_Lipschitz_Hessian_chapter_8} hold. 
Then, the iterates $\{x^{k}\}$ generated by Variant 2 of Algorithm~\ref{ch8:alg:2_order_momentum_chapter_8} with $\rho_k=1-\eta_k$ and with 
\begin{equation*}
  \eta_k =  \eta = \frac{1}{(K+1)^{2/3}} \quad \text{and} \quad  \gamma_k= \gamma = \frac{\hat \gamma}{(K+1)^{2/3}},
\end{equation*}
where $\hat \gamma =  \frac{1}{3}\min\left\{ \frac{1}{L_1} , \frac{1}{\sqrt{M_1 } }  \right\}$  satisfy  
\begin{align*}
& \min_{k \in \{0,1,\ldots,K\}} \Exp{\normstar{\nabla f(x^{k})}} 
\leq   \frac{2\Delta_0}{\hat\gamma(K+1)^{\nicefrac{1}{3}}} + \frac{4\left( \Exp{\normstar{r^{0}}} + {\color{red}\bar\rho} \sigma_g\right)}{ (K+1)^{\nicefrac{1}{3}}}  + \frac{\hat\gamma L_0 \exp(L_1 \hat \gamma) }{(K+1)^{\nicefrac{2}{3}}}  \\
& \qquad + \frac{2\hat \gamma^2}{(K+1)^{\nicefrac{2}{3}}}\left(M_0 + \frac{2}{3}M_1 \hat \gamma L_0 \exp(L_1\hat \gamma)\right)   + 4\frac{{\color{red}\bar\rho}}{\color{blue} \thesischeightubar{\theta}} \frac{\hat \gamma \sigma_H }{(K+1)^{\nicefrac{1}{3}}}, 
\end{align*}
where $\Delta_0=f(x^{0}) - f^{\inf}$ and $r^0=m^0-\nabla f(x^0)$.
\end{theorem}

From Theorems~\ref{ch8:thm:V1_relaxed_smoothness_chapter_8} and~\ref{ch8:thm:V2_relaxed_smoothness_chapter_8}, Algorithm~\ref{ch8:alg:2_order_momentum_chapter_8} using two second-order momentum variants with $\rho_k=1-\eta_k$ achieves the $\cO(1/K^{1/3})$ convergence in the expected gradient norm for minimizing relaxed smooth functions under the arbitrary norm setting.
Our methods match the convergence rates of second-order momentum methods established for standard smooth functions under the Euclidean norm setting by~\textcite{salehkaleybar2022momentum,tran2022better} and under bounded variance assumption by~\textcite{sadiev2025second}.
Our results do not require large minibatch sizes for computing stochastic gradients, and do not impose restrictive assumptions, such as  bounded stochastic gradient norms by~\textcite{tran2022better}.
Moreover, our methods converge faster than LMO-based optimization methods using Polyak momentum analyzed by~\textcite{pethick2025training, pethick2025generalized,kovalev2025understanding,riabinin2025gluon}. 
Also, note that Variant 2 of Algorithm~\ref{ch8:alg:2_order_momentum_chapter_8} assumes the symmetric smoothness on the Hessian, which Variant 1 of Algorithm~\ref{ch8:alg:2_order_momentum_chapter_8} does not. 
Moreover, although our theoretical results are established for Algorithm~\ref{ch8:alg:2_order_momentum_chapter_8} with $\rho_k=1-\eta_k$, the analysis can be extended to the case of general $\rho_k \in (0,1]$. In such cases, the convergence bounds include an additional term depending on the choice of $\rho_k$. 
Nonetheless, determining an optimal choice of $\rho_k$ remains unclear. We leave a tight convergence analysis under appropriate $\rho_k$ selection as an open direction for future work.

\textbf{Extension to extrapolated momentum.}
We can leverage our analysis to establish  the convergence of LMO-based methods in~\eqref{ch8:eqn:LMO_momentum_chapter_8}, which incorporate extrapolated momentum in~\eqref{ch8:eqn:extrapolated_m_chapter_8}, for minimizing relaxed smooth functions in the arbitrary norm setting (see~\ref{ch8:app:IGT_relaxed_smoothness_chapter_8}). 
The resulting convergence matches the $\cO(1/K^{2/7})$ rate obtained by \textcite[Corollary 5]{kovalev2025understanding} under traditional smoothness. 
{\color{black}
\begin{theorem}
    Consider the problem of minimizing $f(x)=\thesischeightExpSub{\xi \sim \thesischeightcD}{f(x;\xi)}$. Let $f$ be twice differentiable, and let   Assumptions~\ref{ch8:assum:Hessian_boundedvariance_v2_chapter_8},~\ref{ch8:assum:lower_bound_chapter_8},~\ref{ch8:assum:Relaxed_Lipschitz_Grad_chapter_8}, and~\ref{ch8:assum:Relaxed_Lipschitz_Hessian_chapter_8} hold. 
    Then, the iterates $\{x^{k}\}$ generated by LMO-based methods in~\eqref{ch8:eqn:LMO_momentum_chapter_8} using extrapolated momentum in~\eqref{ch8:eqn:extrapolated_m_chapter_8}  and stepsizes: 
    \begin{equation*}
        \eta_k =  \eta = \frac{1}{(K+1)^{4/7}} \quad \text{and} \quad  \gamma_k= \gamma = \frac{\hat \gamma}{(K+1)^{5/7}},
    \end{equation*}
    where $\hat \gamma = \min\left\{\frac{1}{3L_1}, \frac{1}{3\sqrt{M_1}}\right\}$  satisfy  
    \begin{equation*}
        \min_{k \in \{0,1,\ldots,K\}} \Exp{\normstar{\nabla f(x^{k})}} \leq \cO\left( \frac{1}{(K+1)^{2/7}} \right) + \cO\left( \frac{1}{(K+1)^{5/7}} \right).
    \end{equation*}
\end{theorem}
}

\textbf{Extension to momentum variance reduction.}
Our analysis can also be used to derive the convergence of  
the LMO-based methods in~\eqref{ch8:eqn:LMO_momentum_chapter_8} that incorporate momentum variance reduction (MVR/\algname{STORM}/\algname{MARS}), as defined in~\eqref{ch8:eq:MARS_chapter_8}. We derive the results for minimizing relaxed smooth functions in an arbitrary norm setting (see~\ref{ch8:app:MVR_relaxed_smoothness_chapter_8}), assuming \steponeonechange{mean-squared generalized smoothness} in an arbitrary norm setting (see~Assumption~\ref{ch8:assum:Relaxed_Mean_Squared_Lipschitz_Grad_chapter_8}). 
Notably, our derived convergence rate matches the $\cO(1/K^{1/3})$ rate obtained by \textcite[Theorem 3]{kovalev2025non} under traditional matrix-type smoothness.
 
{\color{black}
\begin{theorem}
    Consider the problem of minimizing $f(x)=\thesischeightExpSub{\xi \sim \thesischeightcD}{f(x;\xi)}$. Suppose Assumptions~\ref{ch8:assum:Hessian_boundedvariance_v2_chapter_8},~\ref{ch8:assum:lower_bound_chapter_8},~\ref{ch8:assum:Relaxed_Lipschitz_Grad_chapter_8}, and~\ref{ch8:assum:Relaxed_Mean_Squared_Lipschitz_Grad_chapter_8} hold.
    Then, the iterates $\{x^{k}\}$ generated by the LMO-based method in~\eqref{ch8:eqn:LMO_momentum_chapter_8} using the momentum variance reduction in~\eqref{ch8:eq:MARS_chapter_8} with  $\rho_k = 1-\eta_k$, and  step sizes
    \begin{equation*}
        \eta_k = \eta = \frac{1}{(K+1)^{2/3}} \quad \text{and} \quad \gamma_k= \gamma = \frac{\hat \gamma}{(K+1)^{2/3}},
    \end{equation*}
    where $\hat \gamma = \frac{1}{3}\min\left\{\frac{1}{L_1}, \frac{1}{11\steponeonechange{\cL_1}}\left(\frac{\color{red}\bar \rho}{\color{red} \thesischeightubar \rho}\right)^{-1}\right\}$, satisfy
    \begin{equation*}
        \min_{k \in \{0,1,\ldots,K\}} \Exp{\normstar{\nabla f(x^{k})}} \leq \cO\left( \frac{1}{(K+1)^{1/3}} \right) + \cO\left( \frac{1}{(K+1)^{2/3}} \right).
    \end{equation*}
\end{theorem}
}

\section{Numerical Experiments}

We conduct numerical experiments to investigate the performance of LMO-based  methods using Polyak momentum, extrapolated momentum~\parencite{cutkosky2020momentum}, and two second-order momentum variants (Algorithm~\ref{ch8:alg:2_order_momentum_chapter_8}). 
Specifically, the methods using two second-order momentum variants with $\rho_k=1-\eta_k$ are referred to as SOM-V1 and SOM-V2, while those with general $\rho_k$ values are denoted by $\rho$-SOM-V1 and $\rho$-SOM-V2. 
We implemented these algorithms using PyTorch~\parencite{paszke2019pytorch}, and benchmarked them for solving
{\color{blue}
problems of training Multi-Layer Perceptrons (MLPs), Long Short-Term Memory (LSTM) networks, and 
transformer training.
We reported our results for logistic regression problems and training problems of LSTM networks in~Section~\ref{ch8:sec:exp_chapter_8}. 
}
For all experiments, each element of the initial point $x^{0}\in \R^d$ was generated from the standard normal distribution $\thesischeightcN(0,1)$, the random seed was fixed, and also the learning rate $\gamma_k$ and the momentum parameter $\eta_k$ were chosen as follows: (1) $\eta_k = \frac{1}{\sqrt{k+1}}$ and $\gamma_k = \frac{\gamma_0}{(k+1)^{3/4}}$ for Polyak momentum, (2) $\eta_k = \frac{1}{(k+1)^{4/7}}$ and $ \gamma_k = \frac{\gamma_0}{(k+1)^{5/7}}$ for extrapolated momentum, and $\eta_k = \frac{1}{(k+1)^{2/3}}$ and $\gamma_k = \frac{\gamma_0}{(k+1)^{2/3}}$ for all second-order momentum variants.

\subsection{MLP}

We evaluated the algorithms for binary classification tasks using the MLP model with two hidden layers over  the \texttt{splice} dataset from the \texttt{libsvm} library~\parencite{chang2011libsvm}. 
The dataset comprises $1000$ training samples and $60$ features. 
We minimize the objective function $f(x) = \thesischeightcL(x) + R(x)$, where $\thesischeightcL(x)$ is the binary cross-entropy loss with logits, and $R(x)$ is the nonconvex Welsch regularizer defined by $R(x) = \lambda \sum_{i=1}^d x_i^2/(1+x_i^2)$ with $\lambda = 0.01$.

From {\color{black}Figure~\ref{ch8:fig:momentum_comparison_main_text_chapter_8}}, LMO-based methods using second-order momentum outperform those using Polyak or extrapolated momentum. In particular, Variant 2 of second-order momentum achieves superior convergence in training loss, compared to other momentum variants.

\begin{figure}[h!]
    \centering
    \includegraphics[width=0.4\textwidth]{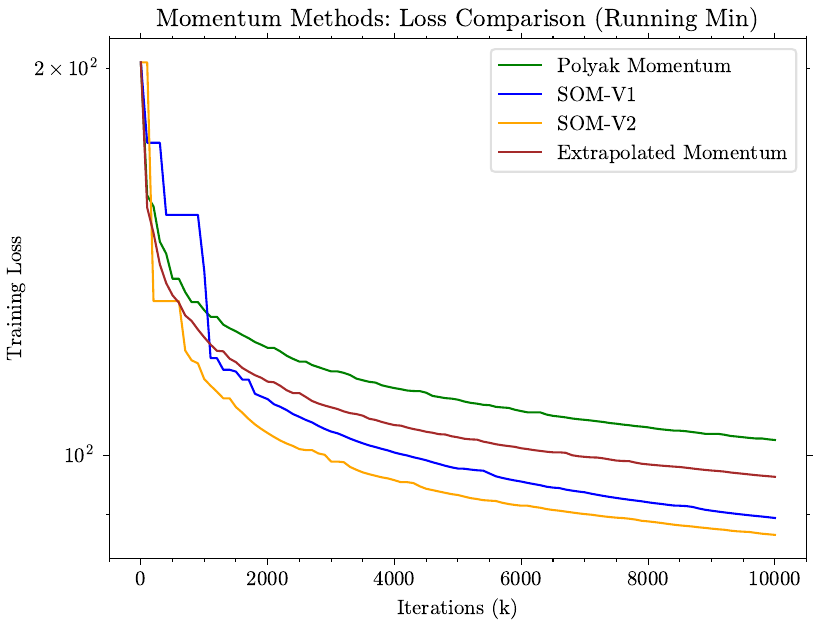}
    \includegraphics[width=0.4\textwidth]{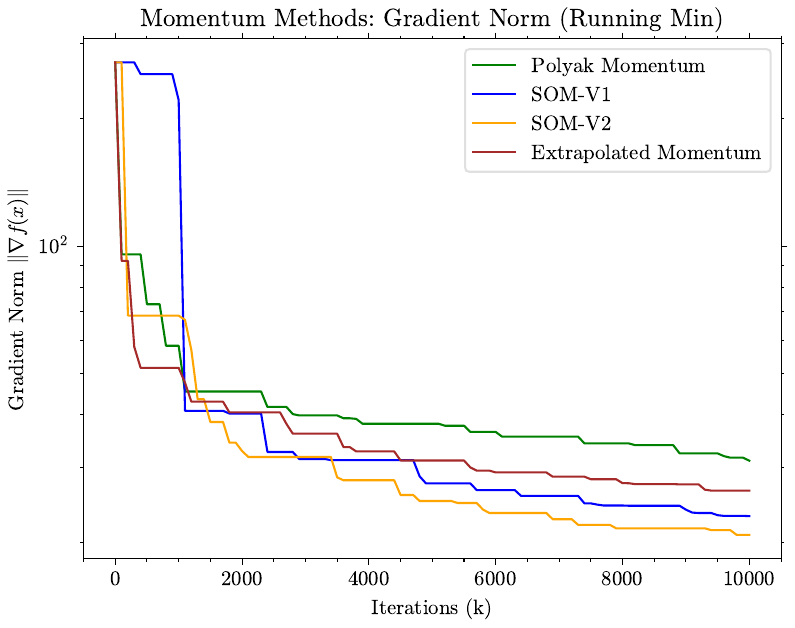}
    \caption{Comparison of training the MLP model on the \texttt{splice} dataset using LMO-based methods with various momentum updates  in the training loss (left) and gradient norm (right).
    }
    \label{ch8:fig:momentum_comparison_main_text_chapter_8}
\end{figure}

{\color{blue}
\subsection{Transformer Training on TinyShakespeare}
We evaluate our proposed methods on a character-level language modeling task. The model is a GPT-2 style Transformer (Nano-GPT) consisting of an embedding layer with a dimension of $128$, two Transformer blocks (each with $4$ attention heads), and a final linear decoder layer. The architecture utilizes standard LayerNorm and GELU activations.
The experiment is conducted on the TinyShakespeare dataset, a standard benchmark containing approximately $1$ million characters. We use a sequence length (block size) of $32$ and a batch size of $16$. We benchmark our proposed \algname{LMO-SOM-V1} and \algname{LMO-SOM-V2} algorithms against Polyak momentum, an LMO method using \algname{STORM}, and an Adam baseline.
For the LMO-based algorithms, we apply the spectral Newton-Schulz orthogonalization oracle to the 2D hidden matrices. To ensure a stable training trajectory and prevent the implicit weight decay effect of the standard LMO update from aggressively eroding the learned parameters, we utilize an effective learning rate of $0.02$ by setting a conservative base step size $\gamma_{base} = 0.0001$ combined with an expanded oracle radius $R = 200.0$. The training is performed over $K = 2000$ iterations using a cosine learning rate schedule with a 10\% linear warmup. For our methods, we apply the damping factor $\rho_k = 0.25(1 - \eta_k)$ to stabilize the highly non-convex landscape of the Transformer.
\begin{figure}[h!]
\centering
\includegraphics[width=0.4\textwidth]{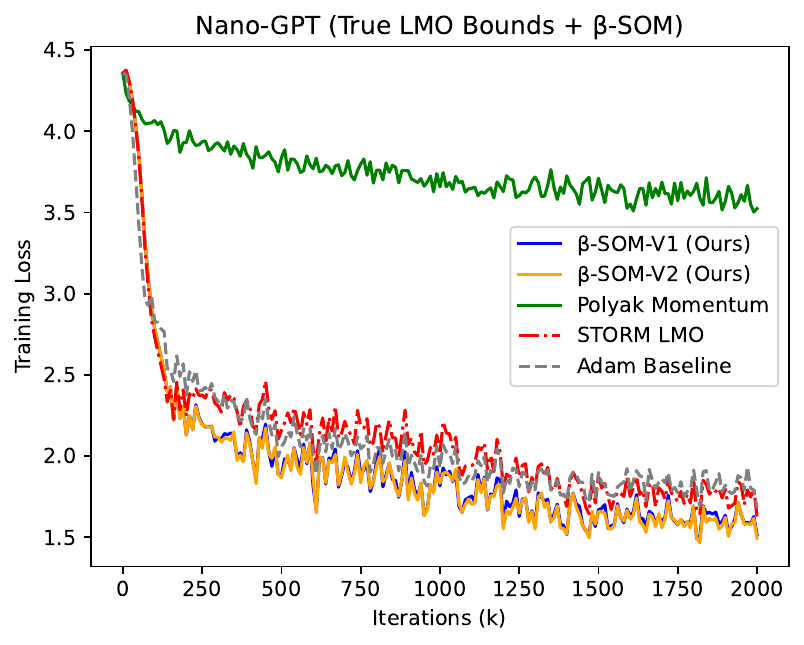}
\includegraphics[width=0.4\textwidth]{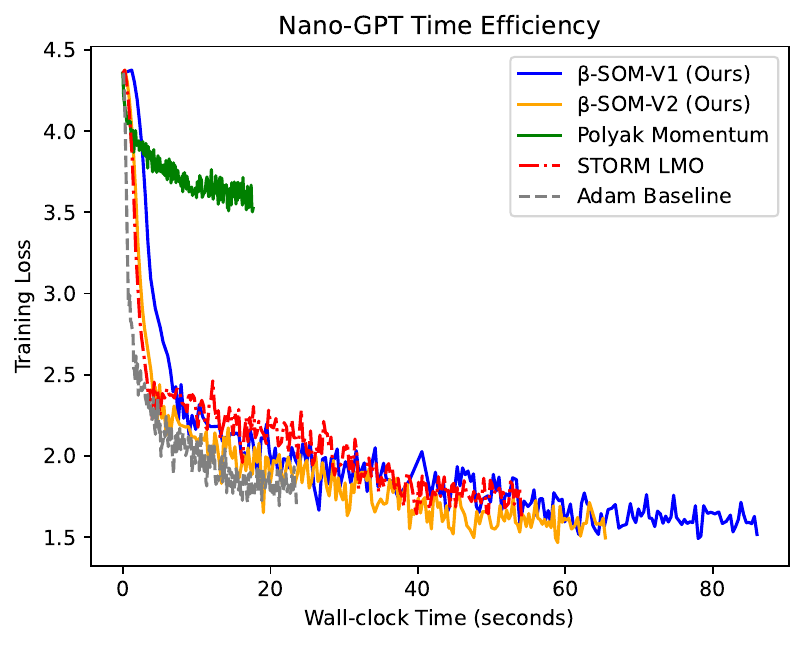}
\caption{Comparison of training the Nano-GPT model on the TinyShakespeare dataset. Left: Training loss vs. iterations. Right: Training loss vs. wall-clock time.}
\label{ch8:fig:nanogpt_shakespeare_chapter_8}
\end{figure}
Figure~\ref{ch8:fig:nanogpt_shakespeare_chapter_8} illustrates the convergence behavior of the evaluated optimizers. In terms of sample efficiency (Figure~\ref{ch8:fig:nanogpt_shakespeare_chapter_8}, left), both \algname{LMO-SOM} variants outperform Polyak momentum and the LMO method using \algname{STORM} in the reported run. They also attain lower loss per iteration than the Adam baseline under this configuration. These observations are consistent with a benefit from the Hessian-vector-product correction, but the experiment does not isolate curvature as the sole causal mechanism.

We also note that the convergence trajectories of \algname{LMO-SOM-V1} and \algname{LMO-SOM-V2} are nearly indistinguishable under the tested setup. Variant~2 is less expensive in the reported implementation because it avoids the additional interpolation step used by Variant~1; this observation does not establish a universal preference across deep-learning tasks.
Finally, on the approximately three-million-parameter model used in this experiment, Adam has a lower wall-clock cost per iteration than the \algname{LMO-SOM} variants. The reported experiments do not determine how the relative overhead of Hessian--vector products and Newton--Schulz iterations scales to billion-parameter models. Establishing performance in that regime requires dedicated large-scale benchmarks and is left for future work.

\subsection{Scaling to WikiText-2}
To evaluate the methods on a larger corpus, we extend the character-level Nano-GPT experiment to the WikiText-2 (WT2) dataset. Unlike TinyShakespeare, WikiText-2 retains original punctuation and capitalization and has a larger vocabulary. The training set contains approximately $2.5$ million characters.
To accommodate the increased complexity, we expand the sequence length to $128$ and the batch size to $32$. The optimization hyperparameters, including the LMO radius $R=200.0$, base step size $\gamma_{base}=0.0001$, and the damping factor $\rho_k$, remain strictly identical to the previous experiment to ensure a fair evaluation of the algorithm's inherent stability. The models are trained over $K = 4000$ iterations.

\begin{figure}[h!]
\centering
\includegraphics[width=0.4\textwidth]{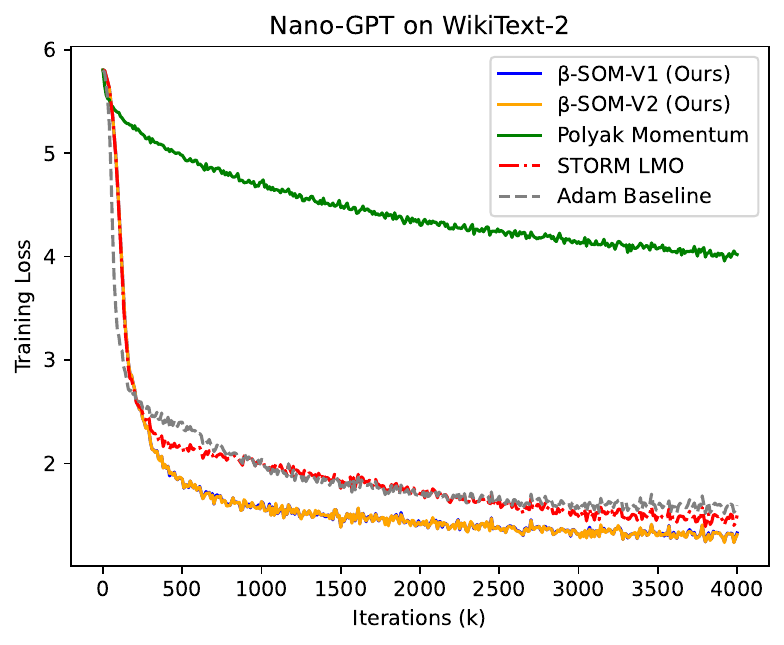}
\includegraphics[width=0.4\textwidth]{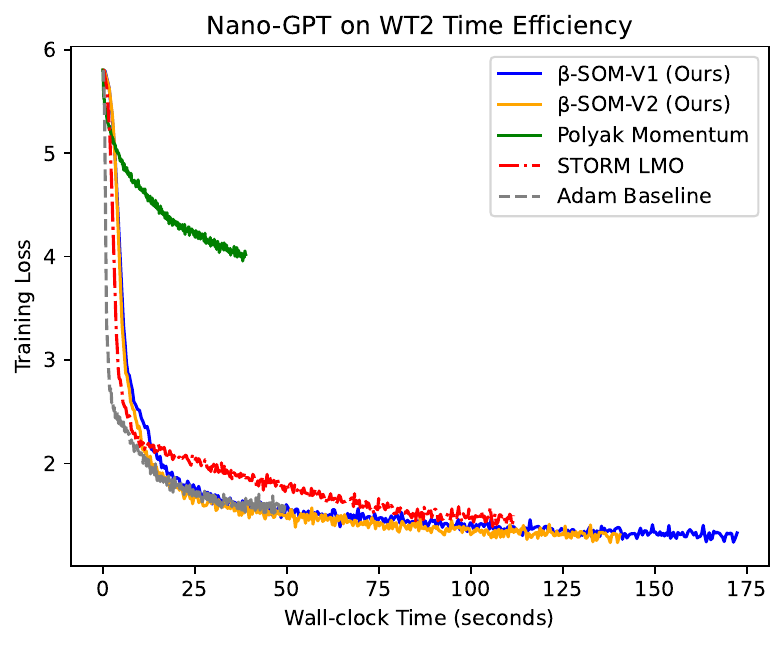}
\caption{Comparison of training the Nano-GPT model on the WikiText-2 dataset. Left: Training loss vs. iterations. Right: Training loss vs. wall-clock time.}
\label{ch8:fig:nanogpt_wt2_chapter_8}
\end{figure}

Under the tested WikiText-2 configuration, the second-order variants outperform the \algname{STORM}-LMO baseline. This observation is consistent with a benefit from curvature information, but the experiment does not isolate curvature as the sole causal mechanism. In the reported runs, the two \algname{LMO-SOM} variants follow nearly indistinguishable trajectories and attain sample efficiency comparable to the Adam baseline.
}

\section{Conclusion and Future Work}

This chapter extended advanced stochastic momentum estimators to nonconvex optimization with centered-ball linear minimization oracles in arbitrary norms. The oracle converts the search vector into a geometry-aware displacement, while implicit gradient transport, shared-sample momentum variance reduction, or Hessian--vector corrections improve the vector supplied to that oracle. In particular, the Hessian information refines the search direction; it is not generally used as the quadratic matrix in an outer local model.

The two \algname{LMO-SOM} variants and the MVR/\algname{STORM}/\algname{MARS} construction attain \(\mathcal{O}(K^{-1/3})\) rates, while the iterative-gradient-transport extension attains \(\mathcal{O}(K^{-2/7})\). These guarantees concern the analytical criterion \(\min_{0\leq k\leq K}\Exp[\normstar{\nabla f(x^k)}]\). They are in expectation, do not identify an implementable argmin without exact gradient norms, and do not guarantee the terminal iterate. The rates improve the arbitrary-norm Polyak-momentum order and match known orders in related standard-smooth or Euclidean settings, but no new lower bound establishes optimality for the arbitrary-norm relaxed-smooth model.

Several gaps separate this framework from a complete theory of practical matrix optimizers. The analysis does not cover high-probability heavy-tailed noise, parameter-agnostic schedules, or independent sampling in place of the stated shared-sample structure. The best choice of the additional momentum parameter remains open. Practical \algname{Muon}-style implementations also use inexact Newton--Schulz or polar computations, whereas the theory assumes an exact centered-ball oracle; finite-precision and inner-iteration errors require separate analysis. The experiments are limited in scale and do not establish a billion-parameter training claim. Chapter~\ref{ch:chapter_9} synthesizes these distinctions and identifies the remaining links among robustness, communication, estimator information, geometry, and implementable outputs.

\newpage

  \endgroup

\chapter{Conclusion and Future Directions}
\label{ch:chapter_9}
\section[Principal Conclusions]{\stepsevenchange{Principal Conclusions}}

\begin{stepsevenchanged}
This dissertation examined how stochastic optimization methods can make progress when the information available to them is heavy-tailed, compressed, or costly to refine. The common viewpoint separates the construction of a search vector from the choice of the displacement applied to the iterate. This distinction connects methods that clip stochastic observations, transport earlier gradient estimates, reconstruct compressed messages, or select steps through a norm-ball linear minimization oracle. This viewpoint makes the role of each algorithmic component explicit. The lower-bound results complement it by distinguishing improvements over existing algorithms from rates that reach the limits of a specified information model.

Chapters~\ref{ch:chapter_2}--\ref{ch:chapter_3} show that high-probability control under weak moments must respect the structure of the underlying problem. Chapter~\ref{ch:chapter_2} develops clipped methods for stochastic minimization and variational inequalities under bounded central $p$-th moments, including regimes that permit infinite variance. A central part of the analysis is to control the stochastic trajectory: the proof constructs an event on which the iterates remain in a prescribed region, allowing smoothness and moment assumptions to be imposed locally. Chapter~\ref{ch:chapter_3} addresses the additional difficulty that a smooth component or local gradient need not vanish at a solution. Direct clipping can then change the balance that defines the deterministic fixed point. Clipping differences relative to suitable reference vectors or learned shifts preserves that balance, while retaining high-probability guarantees, acceleration in the stated regimes, and linear speedup in the distributed stochastic terms. The combined lesson is that robustness depends not only on the transformation, but on the quantity being transformed.

Chapters~\ref{ch:chapter_4}--\ref{ch:chapter_5} determine when additional stochastic information can improve the complexity of finding first-order stationary points. Chapter~\ref{ch:chapter_4} establishes lower bounds for zero-respecting stochastic higher-order algorithms. In the covered noise-dominated regime, second-order information improves the leading dependence on target accuracy by one power of $1/\varepsilon$ relative to the corresponding first-order benchmark; higher derivative orders do not further improve this accuracy exponent. Normalized Hessian-corrected momentum (\algname{HM}) attains the same accuracy dependence in expectation, using Hessian--vector products to transport an estimate rather than solve a Newton subproblem. Separate clipping of gradient and Hessian--vector corrections gives a batch-free high-probability counterpart. Chapter~\ref{ch:chapter_5} separates the moment condition on an individual gradient from the regularity of same-sample differences and permits distinct gradient and Hessian moment indices. Under $q$-weak average smoothness, normalized momentum variance reduction (\algname{MVR}) matches the corresponding lower bound in expectation; under similarity, different normalized methods recover different branches of the lower-bound minimum. These comparisons quantify the value of richer information while locating the remaining gaps in confidence dependence and problem-parameter dependence.

Chapters~\ref{ch:chapter_6}--\ref{ch:chapter_7} address the loss of information through communication compression. Chapter~\ref{ch:chapter_6} combines \algname{EF21} reconstruction from compressed innovations with a normalized server step under symmetric generalized smoothness, without assuming uniformly bounded gradient dissimilarity. Reconstruction and normalization play complementary roles: the former controls the discrepancy between local estimates and communicated information, while the latter controls the server displacement. Under standard smoothness, Chapter~\ref{ch:chapter_7} develops five momentum variants with decreasing stepsizes that do not require the smoothness constants or initial optimality gap, although structural parameters and a user-selected base scale remain. Their guarantees concern a stepsize-weighted random iterate and recover the corresponding established convergence orders up to logarithmic factors. Thus, in the analyzed regimes, communication compression need not remove the convergence-rate improvements obtained from better local estimators.

Chapter~\ref{ch:chapter_8} extends this comparison by changing the geometry of the outer step. Centered norm-ball linear minimization oracles replace Euclidean normalization, while Polyak momentum, implicit gradient transport (\algname{IGT}), \algname{MVR}, and \algname{HM} provide the search vectors. All variants use generalized smoothness; \algname{MVR} additionally uses mean-squared generalized smoothness, while \algname{IGT} and direct-Hessian \algname{HM} use second-order generalized smoothness. It controls the minimum expected dual gradient norm along the trajectory and recovers the corresponding Euclidean rate exponents. These are upper-bound comparisons, not a minimax characterization of the full arbitrary-norm setting. Nevertheless, their recurrence under normalization, compressed reconstruction, and norm-aware steps shows that the benefits of improved estimators extend beyond a single update mechanism.

Taken together, the results identify three interdependent design choices: which quantity to transform, which stochastic information to obtain and retain, and how to control the resulting displacement. Lower bounds provide a complementary assessment of these choices. When upper and lower bounds match, the covered complexity dependence cannot be improved within the same function, oracle, algorithm, probability, and output model. When they do not match, the discrepancy identifies a specific question for further analysis. The contribution is therefore both constructive and diagnostic: algorithms demonstrate how imperfect information can be used, while complexity bounds distinguish attained limits from remaining opportunities.
\end{stepsevenchanged}

\section[Future Directions]{\stepsevenchange{Future Directions}}

\begin{stepsevenchanged}
A first direction is to close the remaining gaps in high-probability complexity for first-order stationarity. The clipped methods developed here achieve polylogarithmic confidence dependence without growing batches, but do not match all problem-parameter dependencies of the available lower bounds. Chapter~\ref{ch:chapter_5} also leaves a gap between expectation and high-probability rates when the moment indices of gradients and richer information differ. It remains to determine whether this gap reflects the analysis or an intrinsic limitation of the corresponding oracle and output model. Sharper concentration arguments could exploit same-sample differences and curvature corrections more effectively, while preserving the weak-moment assumptions and making the required confidence-dependent tuning explicit.

A natural extension of Chapters~\ref{ch:chapter_4}--\ref{ch:chapter_5} is to find approximate second-order stationary points under simultaneous heavy-tailed gradient and Hessian noise. The current methods use Hessian information to improve first-order stationarity; the stronger target is a point satisfying, with high probability,
\[
\|\nabla f(x)\|_2\leq\varepsilon,
\qquad
\lambda_{\min}\!\left(\nabla^2 f(x)\right)\geq-\varepsilon_H,
\]
where $\varepsilon,\varepsilon_H>0$ are the gradient and curvature tolerances. The local-model viewpoint of Section~\ref{ch1:sec:intro_lens_chapter_1} suggests a route: use $m^k$ for a robust gradient estimate, $B^k$ for a robust curvature model, and a finite radius $r_k$ to control the step. To represent negative curvature directly, this would extend the positive-semidefinite model in \eqref{ch1:eq:intro_local_model_chapter_1} to potentially indefinite symmetric matrices, with a trust-region constraint or suitable regularization keeping the subproblem well posed. Robust trust-region and regularized second-order methods are natural candidates. Their analysis would need to distinguish genuine negative curvature from estimation error, control approximate subproblem solutions, and characterize the oracle cost of satisfying both stationarity conditions under weak moments.

Another challenge is to combine high-probability robustness, parameter agnosticism, and non-Euclidean geometry. These properties are developed in different parts of the dissertation, but the current results do not establish them jointly. A natural objective is a method that handles weak-moment noise directly in the chosen primal and dual norms while avoiding prior knowledge of selected smoothness and noise parameters. Such an analysis could clarify which dimension dependencies are intrinsic and which result from converting the problem to Euclidean geometry. It would also revisit the structural question from Chapter~\ref{ch:chapter_3}: whether clipping should act on a stochastic gradient, a momentum innovation, or the displacement selected by the oracle.

The output certified by a convergence theorem raises a separate question. Uniformly sampled and stepsize-weighted random iterates already provide implementable outputs, whereas an analytical best-iterate bound need not identify a point that can be selected from observable information. For the nonconvex regimes studied here, it remains important to understand when trajectory control can yield a terminal-iterate guarantee or a reliable selection procedure with a comparable oracle budget. This requires identifying additional stability or structural conditions, or developing certificates that remain meaningful under heavy-tailed noise. The cost of selecting and validating an output should be included in the complexity rather than treated as an operation available without further observations.

The communication model can also be extended beyond synchronous full participation with a fixed worker population. Partial participation, asynchronous messages, decentralized communication, and multiple local steps raise different questions about how reconstructed information should be maintained. Studying these extensions separately would help isolate the effects of sampling, delay, and local computation on the estimator and its communication error. The appropriate performance measure may also change: communication rounds, transmitted bits, local oracle calls, and elapsed time capture different aspects of efficiency. A useful extension would specify which of these resources is being reduced and how that reduction affects the convergence guarantee.

Finally, empirical evaluation should determine when improved oracle complexity leads to improved computational efficiency. The experiments in Chapter~\ref{ch:chapter_7} already illustrate that stronger per-epoch performance need not imply better wall-clock performance when Hessian--vector products increase the cost of an update. Broader comparisons should therefore account separately for gradient evaluations, Hessian--vector products, compression, communicated bits, synchronization, memory, and tuning effort. Evaluation across different noise levels, data distributions, and computational budgets would help identify when richer estimators justify their additional cost. The aim is not simply to test larger models, but to understand which theoretical advantages remain useful when all major resources are counted.
\end{stepsevenchanged}

\section[Closing Perspective]{\stepsevenchange{Closing Perspective}}

\begin{stepsevenchanged}
The results suggest approaching new optimization problems by first asking what information can be trusted, what must be estimated or reconstructed, and what displacement that information can justify. This perspective shifts attention from individual algorithm names to the interaction between estimation and movement. The next challenge is to combine the mechanisms developed here without losing their statistical guarantees or their computational advantages. Doing so will require analyses that account for the information actually queried, the messages actually communicated, and an output that can be returned and assessed. The path from imperfect information to controlled steps is therefore a continuing design problem: learning how much information is needed for progress, and how to use it without allowing its errors to dictate the update.
\end{stepsevenchanged}

\begingroup
\urlstyle{same}
\setlength{\emergencystretch}{3em}
\sloppy
\printbibliography[heading=kaustreferences,title={REFERENCES}]
\endgroup

\appendix
\chapter*{\centering APPENDICES}
\phantomsection
\addcontentsline{toc}{chapter}{\textbf{APPENDICES}}
\chapter{Supplementary Material and Proofs for Chapter 2}
  \begingroup
  \thesischaptersetup{2}%
  \stepEightAAppendixColors{2}%
  
{\small\tableofcontents}

\clearpage

\section{Additional Related Work}\label{ch2:appendix:additional_related_chapter_2}

In this section, we provide an overview of the existing in-expectation convergence results under Assumption~\ref{ch2:as:bounded_alpha_moment_chapter_2}.

\paragraph{Convex minimization.} The first in-expectation result under Assumption~\ref{ch2:as:bounded_alpha_moment_chapter_2} is given by \textcite{NemiYudin}, who derive\footnote{In this section, we hide in $\cO(\cdot)$ all dependencies except the dependency on $\varepsilon$.} $\cO(\varepsilon^{-\nicefrac{{p}}{({p}-1)}})$ complexity for Mirror Descent applied to the minimization of convex functions with bounded gradients.  This result was recently extended by \textcite{vural2022mirror} to the uniformly convex functions, and matching lower bounds were derived. 
In the strongly convex case, \textcite{zhang2020adaptive} prove $\cO(\varepsilon^{-\nicefrac{{p}}{2({p}-1)}})$ complexity for \algname{clipped-SGD}. However, all these results rely on the boundedness of the gradient. To the best of our knowledge, there are no results for smooth convex problems under Assumption~\ref{ch2:as:bounded_alpha_moment_chapter_2} without assuming that the gradient is bounded even in terms of expectation.

\paragraph{Non-convex minimization.} In the non-convex smooth case, \textcite{zhang2020adaptive} prove $\cO(\varepsilon^{-\nicefrac{(3{p}-2)}{({p}-1)}})$ complexity for \algname{clipped-SGD} to produce a point $x$ such that $\Exp\left[\|\nabla f(x)\|\right] \leq \varepsilon$. In the same work, the authors derive the matching lower bound. However, both upper and lower bounds are derived for $\Exp\left[\|\nabla f(x)\|\right]$ which is smaller than $\sqrt{\Exp\left[\|\nabla f(x)\|^2\right]}$. The later one is stronger and is more standard performance metric for stochastic non-convex optimization. Therefore, the question of deriving lower and matching upper bounds for the standard metric remains open.

\clearpage

\section{Useful Facts}\label{ch2:appendix:useful_facts_chapter_2}

\paragraph{Smoothness.} If $f$ is $L$-smooth on convex set $Q \subseteq \R^d$, then for all $x, y \in Q$ \parencite{nesterov2018lectures}
\begin{equation}
    f(y) \leq f(x) + \langle \nabla f(x), y-x \rangle + \frac{L}{2}\|y - x\|^2. \label{ch2:eq:L_smoothness_cor_1_chapter_2}
\end{equation}
In particular, if $x$ and $y = x - \frac{1}{L}\nabla f(x)$ lie in $Q$, then the above inequality gives
\begin{align*}
f(y) & \leq f(x) - \frac{1}{L}\|\nabla f(x)\|^2 + \frac{1}{2L}\|\nabla f(x)\|^2 \\
& = f(x) - \frac{1}{2L}\|\nabla f(x)\|^2
    \end{align*}
and
\begin{equation*}
    \|\nabla f(x)\|^2 \leq 2L\left(f(x) - f(y)\right) \leq 2L\left(f(x) - f^{\inf}\right)
\end{equation*}
under the assumption that $f^{\inf} = \inf_{x \in Q}f(x) > -\infty$. In other words, \eqref{ch2:eq:L_smoothness_cor_2_chapter_2} holds for any $x \in Q$ such that $(x - \frac{1}{L}\nabla f(x)) \in Q$. For example, if $x^*$ is an optimum of $f$, then $L$-smoothness on $B_{2R}(x^*)$ implies that \eqref{ch2:eq:L_smoothness_cor_2_chapter_2} holds on $B_R(x^*)$: indeed, for any $x \in B_R(x^*)$ we have
\begin{align}
\left\|x - \frac{1}{L}\nabla f(x) - x^* \right\| & \leq \|x - x^*\| + \frac{1}{L}\|\nabla f(x)\| \notag \\
& \overset{\eqref{ch2:eq:L_smoothness_chapter_2}}{\leq} 2\|x - x^*\| \leq 2R. \notag
    \end{align}
This derivation means that, in the worst case, to have \eqref{ch2:eq:L_smoothness_cor_2_chapter_2} on a set $Q$ we need to assume smoothness on a slightly larger set.

\paragraph{Parameters in \algname{clipped-SSTM}.} To analyze \algname{clipped-SSTM} we use the following lemma about its parameters $a_k$ and $A_k$.

\begin{lemma}[Lemma~E.1 from \textcite{gorbunov2020stochastic}]\label{ch2:lem:alpha_k_A_K_lemma_chapter_2}
    Let sequences $\{a_k\}_{k\ge0}$ and $\{A_k\}_{k\ge 0}$ satisfy
    \begin{equation}
        a_{0} = A_0 = 0,\quad A_{k+1} = A_k + a_{k+1},\quad a_{k+1} = \frac{k+2}{2aL}\quad \forall k \geq 0, \label{ch2:eq:alpha_k_A_k_def_chapter_2}
    \end{equation}
    where $a > 0$, $L > 0$. Then for all $k\ge 0$
    \begin{align}
        A_{k+1} & = \frac{(k+1)(k+4)}{4aL}, \label{ch2:eq:A_k+1_explicit_chapter_2}\\
        A_{k+1} & \geq a L a_{k+1}^2. \label{ch2:eq:A_k+1_lower_bound_chapter_2}
    \end{align}
\end{lemma}

\paragraph{Bernstein inequality.} One of the final steps in our proofs is in the proper application of the following lemma known as {\it Bernstein inequality for martingale differences} \parencite{bennett1962probability,dzhaparidze2001bernstein,freedman1975tail}.
\begin{lemma}\label{ch2:lem:Bernstein_ineq_chapter_2}
    Let the sequence of random variables $\{X_i\}_{i\ge 1}$ form a martingale difference sequence, i.e.\ $\Exp\left[X_i\mid X_{i-1},\ldots, X_1\right] = 0$ for all $i \ge 1$. Assume that conditional variances $\sigma_i^2\eqdef\Exp\left[X_i^2\mid X_{i-1},\ldots, X_1\right]$ exist and are bounded and assume also that there exists deterministic constant $c>0$ such that $|X_i| \le c$ almost surely for all $i\ge 1$. Then for all $b > 0$, $H > 0$ and $n\ge 1$
    \begin{equation}
        \Prob\left\{\Big|\sum\limits_{i=1}^nX_i\Big| > b \text{ and } \sum\limits_{i=1}^n\sigma_i^2 \le H\right\} \le 2\exp\left(-\frac{b^2}{2H + 2cb/3}\right).
    \end{equation}
\end{lemma}

\clearpage

\section{Proof of Lemma~\ref{ch2:lem:bias_and_variance_clip_chapter_2}}

\begin{lemma}[Lemma~\ref{ch2:lem:bias_and_variance_clip_chapter_2}]\label{ch2:lem:bias_variance_chapter_2}
    Let $X$ be a random vector in $\R^d$ and $\thesischtwotX = \clip(X,\lambda)$. Then,
    \begin{equation}
        \left\|\thesischtwotX - \Exp[\thesischtwotX]\right\| \leq 2\lambda. \label{ch2:eq:bound_X_chapter_2}
    \end{equation} 
    Moreover, if for some $\sigma \geq 0$ and ${p} \in (1,2]$
    \begin{equation}
        \Exp[X] = x\in\R^d,\quad \Exp[\|X - x\|^{{p}}] \leq \sigma^{{p}} \label{ch2:eq:UBV_X_chapter_2}
    \end{equation}
    and $\|x\| \leq \lambda/2$, then
    \begin{align}
        \left\|\Exp[\thesischtwotX] - x\right\| & \leq \frac{2^{{p}}\sigma^{{p}}}{\lambda^{{p} -1}}, \label{ch2:eq:bias_X_chapter_2}\\
        \Exp\left[\left\|\thesischtwotX - x\right\|^2\right] & \leq 18\lambda^{2-{p}}\sigma^{{p}}, \label{ch2:eq:distortion_X_chapter_2}\\
        \Exp\left[\left\|\thesischtwotX - \Exp[\thesischtwotX]\right\|^2\right] & \leq 18\lambda^{2-{p}}\sigma^{{p}}. \label{ch2:eq:variance_X_chapter_2}
    \end{align}
\end{lemma}
\begin{proof}
    \textbf{Proof of \eqref{ch2:eq:bound_X_chapter_2}:} by definition of a clipping operator, we have 
    \begin{align*}
        \left\|\thesischtwotX -\Exp\left[\thesischtwotX\right]\right\| & \leq \left\|\thesischtwotX\right\| + \left\|\Exp\left[\thesischtwotX\right]\right\|\\
        & = \left\|\texttt{clip}(X, \lambda)\right\| + \left\|\Exp\left[\texttt{clip}(X,\lambda)\right]\right\|\\
        & \leq \left\|\min\left\{1,\frac{\lambda}{\left\|X\right\|}\right\}X\right\| + \Exp\left[\left\|\min\left\{1,\frac{\lambda}{\left\|X\right\|}\right\}X\right\|\right]\\
        & = \min\left\{\left\|X\right\|,\lambda\right\} + \Exp\left[\min\left\{\left\|X\right\|,\lambda\right\}\right]\\
        & \leq \lambda + \lambda = 2\lambda.
    \end{align*}

    \textbf{Proof of \eqref{ch2:eq:bias_X_chapter_2}}: To start the proof, we introduce two indicator random variables. Let 
    \begin{equation}
        \begin{aligned}
            \chi &= \mathbb{I}_{\left\{X: \left\|X\right\| > \lambda\right\}} = \begin{cases} 1, & \text{if } \|X\| > \lambda,\\ 0, & \text{otherwise},\end{cases}\\
            \eta &= \mathbb{I}_{\left\{X: \left\|X-x\right\| > \lambda/2\right\}} = \begin{cases} 1, & \text{if } \|X - x\| > \lambda/2,\\ 0, & \text{otherwise}.\end{cases}
        \end{aligned}
    \end{equation}
    Moreover, since $\|X\| \leq \|x\| + \|X-x\| \leq \lambda/2 + \|X-x\|$ (we use $\|x\| \leq \lambda/2$), we have $\chi \leq \eta$. We are now in a position to show \eqref{ch2:eq:bias_X_chapter_2}. Using that
    \begin{equation*}
        \thesischtwotX = \min\left\{1,\frac{\lambda}{\left\|X\right\|}\right\}X = \chi \frac{\lambda}{\left\|X\right\|} X + (1-\chi)X,
    \end{equation*}
    we obtain 
    \begin{align*}
        \left\|\Exp\left[\thesischtwotX\right] -x\right\| & = \left\|\Exp\left[X + \chi \left(\frac{\lambda}{\left\|X\right\|}-1\right) X\right] -x\right\|\\
        & =\left\|\Exp\left[\chi \left(\frac{\lambda}{\left\|X\right\|}-1\right) X\right] \right\|\\
        & \leq\Exp\left[\Big|\chi \left(\frac{\lambda}{\left\|X\right\|}-1\right)\Big|\left\|X\right\|\right]\\
        & =\Exp\left[\chi \left(1-\frac{\lambda}{\left\|X\right\|}\right)\left\|X\right\|\right].
    \end{align*}
    Since $1-\lambda/(\left\|X\right\|)\in (0,1)$ when $\chi \ne 0$, we derive
    \begin{align*}
        \left\|\Exp\left[\thesischtwotX\right] -x\right\| & \leq \Exp\left[\chi \|X\|\right]\\
        & \overset{\chi\leq \eta}{\leq} \Exp\left[\eta \|X\|\right]\\
        & \leq \Exp\left[\eta \|X-x\| + \eta \|x\|\right]\\
        & \overset{(\ast)}{\leq} \left(\Exp\left[\left\|X-x\right\|^{{p}}\right]\right)^{1/({p})}\left(\Exp\left[\eta^{p/(p - 1)}\right]\right)^{(p - 1)/p} \\
        &\quad + \|x\|\Exp\left[\eta\right]\\
        & \overset{\eqref{ch2:eq:UBV_X_chapter_2}}{\leq}\sigma\left(\Exp\left[\eta^{p/(p-1)}\right]\right)^{(p-1)/p} + \frac{\lambda}{2} \Exp\left[\eta\right],
    \end{align*}
    where in $(\ast)$ we applied Hölder's inequality. By Markov's inequality,
    \begin{align}
    \Exp\left[\eta^{p/(p-1)}\right] & = \Exp\left[\eta\right] = \Prob\left\{\|X-x\|>\frac{\lambda}{2}\right\} \notag \\
    & = \Prob\left\{\|X-x\|^{{p}}>\frac{\lambda^{{p}}}{2^{{p}}}\right\} \label{ch2:eq:bound_eta_chapter_2} \notag\\
        & \leq \frac{2^{{p}}}{\lambda^{{p}}}\Exp\left[\|X-x\|^{{p}}\right]\notag\\
        & \leq \left(\frac{2\sigma}{\lambda}\right)^{{p}}.
    \end{align}
    Thus, in combination with the previous chain of inequalities, we finally have
    \begin{align*}
        \left\|\Exp\left[\thesischtwotX\right] -x\right\| & \leq \sigma \left(\frac{2\sigma}{\lambda}\right)^{{p}-1} + \frac{\lambda}{2}\left(\frac{2\sigma}{\lambda}\right)^{{p}} = \frac{2^{{p}}\sigma^{{p}}}{\lambda^{{p}-1}}.
    \end{align*}

    \textbf{Proof of \eqref{ch2:eq:distortion_X_chapter_2}}: Using $\|\thesischtwotX - x\| \leq \|\thesischtwotX\| +\|x\| \leq \lambda + \frac{\lambda}{2} = \frac{3\lambda}{2}$, we have 
    \begin{align*}
        \Exp\left[\|\thesischtwotX - x\|^2\right] & = \Exp\left[\|\thesischtwotX - x\|^{{p}} \|\thesischtwotX - x\|^{2-{p}}\right]\\
        & \leq \left(\frac{3\lambda}{2}\right)^{2-{p}}\Exp\left[\|\thesischtwotX - x\|^{{p}}\chi + \|\thesischtwotX - x\|^{{p}}(1-\chi)\right]\\
        & = \left(\frac{3\lambda}{2}\right)^{2-{p}}\Exp\left[\chi\left\|\frac{\lambda}{\|X\|}X - x\right\|^{{p}} + \|X - x\|^{{p}}(1-\chi)\right]\\
        & \leq \left(\frac{3\lambda}{2}\right)^{2-{p}}\Exp\left[\chi\left(\left\|\frac{\lambda}{\|X\|}X\right\| + \left\|x\right\|\right)^{{p}} + \|X - x\|^{{p}}(1-\chi)\right]\\
        & \leq \left(\frac{3\lambda}{2}\right)^{2-{p}}\left(\Exp\left[\chi\left(\frac{3\lambda}{2}\right)^{{p}} + \sigma^{{p}}\right]\right),
    \end{align*}
    where in the last inequality we applied $\|x\|\leq \lambda/2$, \eqref{ch2:eq:UBV_X_chapter_2} and $1 - \chi \leq 1$. By \eqref{ch2:eq:bound_eta_chapter_2} and $\chi \leq \eta$, we obtain
    \begin{align*}
        \Exp\left[\|\thesischtwotX - x\|^2\right] & \leq  \frac{9\lambda^2}{4}\left(\frac{2\sigma}{\lambda}\right)^{{p}} + \left(\frac{3\lambda}{2}\right)^{2-{p}}\sigma^{{p}}\\
        & \leq \frac{9\lambda^2}{4}\left(2^{{p}} +\frac{2^{{p}}}{3^{{p}}}\right)\frac{\sigma^{{p}}}{\lambda^{{p}}}\\
        & \leq 18 \lambda^{2-{p}}\sigma^{{p}}.
    \end{align*}

    \textbf{Proof of \eqref{ch2:eq:variance_X_chapter_2}}: Using variance decomposition and \eqref{ch2:eq:distortion_X_chapter_2}, we have 
    \begin{equation*}
        \Exp\left[\left\|\thesischtwotX - \Exp[\thesischtwotX]\right\|^2\right] \leq \Exp\left[\|\thesischtwotX - x\|^2\right] \leq 18 \lambda^{2-{p}}\sigma^{{p}}.
    \end{equation*}
\end{proof}

\clearpage

\section{Proof of Theorem~\ref{ch2:thm:SGD_high-prob_conv_chapter_2}}\label{ch2:appendix:failure_of_SGD_chapter_2}

In this section, we give an example of the problem for which \algname{SGD} without clipping leads to a weak high-probability convergence guarantee even under the strong assumption of bounded variance. Theorem below formally states our result, showing that, in the worst-case, the bound for \algname{SGD} scales worse than that of \algname{clipped-SGD} in terms of the probability $\beta$.

\begin{theorem}\label{ch2:thm:SGD_high-prob_conv_app_chapter_2}
    For any $\varepsilon > 0$, $\beta \in (0,1)$, and \algname{SGD} parameterized by the number of steps $K$ and stepsize $\gamma$, there exists problem \eqref{ch2:eq:min_problem_chapter_2} such that Assumptions~\ref{ch2:as:bounded_alpha_moment_chapter_2}, \ref{ch2:as:L_smoothness_chapter_2}, and \ref{ch2:as:str_cvx_chapter_2} hold with ${p} = 2$, $0 < \mu \leq L$ and for the iterates produced by \algname{SGD} with any stepsize $\gamma > 0$
    \begin{equation}
        \Prob\left\{\|x^K - x^*\|^2 \geq \varepsilon\right\} \leq \beta \;\; \Longrightarrow\;\; K = \Omega\left(\frac{\sigma}{\mu\sqrt{\beta \varepsilon}}\right). \label{ch2:eq:SGD_high-prob_conv_chapter_2}
    \end{equation}
\end{theorem}

\begin{proof}
    To prove the above theorem, we consider the simple one-dimensional problem $f(x) = \mu x^2/2$. It is easy to see that the considered problem is $\mu$-strongly convex, $\mu$-smooth, and has optimum at $x^* = 0$. We construct the noise in an adversarial way with respect to the parameters of the \algname{SGD}. Concretely, the noise depends on the number of iterations $K$, stepsize $\gamma$, target precision $\varepsilon$, the starting point $x^0$, and bound on the variance $\sigma^2$ such that 
    \begin{align*}
        \nabla f(x^k;\xi_k) = \mu x^k + \sigma z_k, 
    \end{align*}
    where 
    \begin{align}
    z_k &=
        \begin{cases}
            0, & \text{if } k < K-1 \text{ or } (1 - \gamma \mu)^K |x^0| > \sqrt{\varepsilon} ,\\ 
            \begin{cases}
                - A,  & \text{with probability } \frac{1}{2A^2}, \\
                 0 & \text{with probability } 1 - \frac{1}{A^2},\\
                  A,  & \text{with probability } \frac{1}{2A^2}, \\
            \end{cases}
            &\text{otherwise}
        \end{cases}
        \notag\\
        &\quad \forall k \in \{0, 1, \ldots, K-1\},
        \label{ch2:eq:sgd_noise_chapter_2}
    \end{align}
    where $A = \max\left\{\frac{2\sqrt{\varepsilon}}{\gamma\sigma} , 1\right\}$. We note that $\Exp\left[z_k\right] = 0$. Therefore, $\Exp\left[\nabla f(x^k;\xi_k)\right] = \mu x^k = \nabla f (x^k).$ Furthermore, 
    \begin{align*}
        \thesischtwoVar\left[z_k\right] = \Exp\left[(z_k)^2\right] \leq \frac{1}{2A^2}  A^2 + \frac{1}{2A^2}  A^2 = 1, 
    \end{align*}
    which implies that Assumption~\ref{ch2:as:bounded_alpha_moment_chapter_2} holds for  ${p} = 2$. We note that our construction depends on the parameters of the algorithm and the target value $\varepsilon$. However, our analysis of the methods with clipping works in such generality. 

    Let us now analyze the properties of the introduced problem. We are interested in the situation when  
    \begin{align*}
        \Prob\left\{\|x^K - x^*\|^2 \geq \varepsilon\right\} \leq \beta
    \end{align*}
    for $\beta \in (0,1)$. We first prove that this implies that $(1 - \gamma\mu)^K |x^0| \leq \sqrt{\varepsilon}$. To do that we proceed by contradiction and assume that 
    \begin{equation}
        \label{ch2:eq:counterexample_proof_1_chapter_2}
        (1 - \gamma\mu)^K |x^0| > \sqrt{\varepsilon}.
    \end{equation}
    By construction, this implies that  $z_k = 0, \; \forall k \in \{0, 1, \ldots, K-1\}$. This, in turn, implies that $x^K = (1-\gamma\mu)^K x^0$, and, further, by \eqref{ch2:eq:counterexample_proof_1_chapter_2} and since $x^*=0$, that
    \begin{align*}
        \Prob\left\{\|x^K - x^*\|^2 \geq \varepsilon\right\} = \Prob\left\{\|x^K\|^2 \geq \varepsilon\right\} = 1.
    \end{align*}
    Thus, the contradiction shows that $(1 - \gamma\mu)^K |x^0| \leq \sqrt{\varepsilon}$, which yields $K \geq \frac{1}{\gamma \mu} \ln\frac{|x^0|}{\sqrt{\varepsilon}}$. Using \eqref{ch2:eq:sgd_noise_chapter_2} with $K \geq \frac{1}{\gamma\mu} \log\frac{|x^0|}{\sqrt{\varepsilon}}$, we obtain
    \begin{align*}
        \|x^K - x^*\|^2 = ((1 - \gamma\mu)^K x^0 + \gamma \sigma z_{K-1})^2.
    \end{align*}
    Furthermore,
    \begin{align*}
        \Prob\left\{\|x^K - x^*\|^2 \geq \varepsilon\right\} &= \Prob\left\{\left|(1 - \gamma\mu)^K  x^0 + \gamma \sigma z_{K-1} \right| \geq \sqrt{\varepsilon}\right\} \\
        &= \Prob\left\{\gamma \sigma z_{K-1} \geq \sqrt{\varepsilon} - (1 - \gamma\mu)^K  x^0 \right\}\\
        &\qquad + \Prob\left\{\gamma \sigma z_{K-1}  \leq -\sqrt{\varepsilon} - (1 - \gamma\mu)^K  x^0\right\} \\
        &\geq \Prob\left\{\gamma \sigma z_{K-1} \geq \sqrt{\varepsilon} + (1 - \gamma\mu)^K  x^0 \right\}\\
        &\qquad + \Prob\left\{\gamma \sigma z_{K-1}  \leq -\sqrt{\varepsilon} - (1 - \gamma\mu)^K  x^0\right\} \\
        &=  \Prob\left\{\left|\gamma \sigma z_{K-1}\right| \geq \sqrt{\varepsilon} + (1 - \gamma\mu)^K  x^0 \right\} \\
        &\geq \Prob\left\{\left|\gamma \sigma z_{K-1}\right| \geq 2\sqrt{\varepsilon} \right\} = \Prob\left\{\left|z_{K-1}\right| \geq \frac{2\sqrt{\varepsilon}}{\gamma \sigma} \right\}.
    \end{align*}
    Now if $\frac{2\sqrt{\varepsilon}}{\gamma \sigma} < 1$ then $A = 1$. Therefore, 
    \begin{align*}
        1 = \Prob\left\{\left|z_{K-1}\right| \geq \frac{2\sqrt{\varepsilon}}{\gamma \sigma} \right\} \leq \Prob\left\{\|x^K - x^*\|^2 \geq \varepsilon\right\} \leq \beta < 1, 
    \end{align*}
    yielding contradiction, which implies that if $\Prob\left\{\|x^K - x^*\|^2 \geq \varepsilon\right\} \leq \beta$ for our constructed problem, then $\frac{2\sqrt{\varepsilon}}{\gamma \sigma} \geq 1$, i.e., $\gamma \leq \frac{2\sqrt{\varepsilon}}{ \sigma}$. For $\gamma \leq \frac{2\sqrt{\varepsilon}}{ \sigma}$, we have 
    \begin{align*}
         \beta \geq  \Prob\left\{\|x^K - x^*\|^2 \geq \varepsilon\right\} \geq \Prob\left\{\left|z_{K-1}\right| \geq \frac{2\sqrt{\varepsilon}}{\gamma \sigma} \right\} = \frac{1}{A^2} =  \frac{\gamma^2 \sigma^2}{4\varepsilon}.
    \end{align*}
    This implies that $\gamma \leq \frac{2\sqrt{\beta \varepsilon}}{\sigma}$. Combining this inequality with $K \geq \frac{1}{\gamma\mu} \log\frac{|x^0|}{\sqrt{\varepsilon}}$ yields
    \begin{align*}
        K \geq \frac{\sigma}{2\mu\sqrt{\beta \varepsilon}}\log\frac{|x^0|}{\sqrt{\varepsilon}}
    \end{align*}
    and concludes the proof.
    \end{proof}

\clearpage

\section{Additional Proofs for \algname{clipped-SGD}}\label{ch2:appendix:clipped_SGD_chapter_2}
In this section, we provide all the missing details and proofs of the results for \algname{clipped-SGD}. For brevity, we will use the following notation: $\thesischtwotnabla f(x^k;\xi^k) = \clip(\nabla f(x^k;\xi^k), \lambda_k)$.

\begin{algorithm}[h]
\caption{Clipped Stochastic Gradient Descent (\algname{clipped-SGD}) \parencite{pascanu2013difficulty}}
\label{ch2:alg:clipped-SGD_chapter_2}   
\begin{algorithmic}[1]
\REQUIRE starting point $x^0$, number of iterations $K$, stepsize $\gamma > 0$, clipping levels $\{\lambda_k\}_{k=0}^{K-1}$.
\FOR{$k=0,\ldots, K-1$}
\STATE Compute $\thesischtwotnabla f(x^{k};\xi^k) = \clip\left(\nabla f(x^{k};\xi^k), \lambda_k\right)$ using a fresh sample $\xi^k \sim \thesischtwocD_k$
\STATE $x^{k+1} = x^k - \gamma \thesischtwotnabla f(x^{k};\xi^k)$
\ENDFOR
\ENSURE $x^K$ 
\end{algorithmic}
\end{algorithm}

\subsection{Non-Convex Functions}

We start the analysis of \algname{clipped-SGD} in the non-convex case with the following lemma that follows the proof of deterministic \algname{GD} and separates the stochasticity from the determinisitc part of the method.

\begin{lemma}\label{ch2:lem:main_opt_lemma_clipped_SGD_non_convex_chapter_2}
    Let Assumptions~\ref{ch2:as:lower_boundedness_chapter_2} and \ref{ch2:as:L_smoothness_chapter_2} hold on
    \begin{equation*}
        Q = \left\{x \in \R^d\mid \exists y\in \R^d:\; f(y) \leq f^{\inf} + 2\Delta \text{ and } \|x-y\|\leq \frac{\sqrt{\Delta}}{20\sqrt{L}}\right\},
    \end{equation*}
    where $\Delta \geq \Delta_0 = f(x^0) - f^{\inf}$, and let stepsize  $\gamma$ satisfy $\gamma < 2/L$. If the segments $[x^k, x^{k+1}]$ lie in $Q$ for all $k = 0,1,\ldots,K-1$, $K \ge 0$, then after $K$ iterations of \algname{clipped-SGD} we have
    \begin{align}
        \MoveEqLeft[1] \gamma\left(1-\frac{L\gamma}{2}\right)\sum\limits_{k=0}^{K-1}\|\nabla f(x^k)\|^2 \notag \\
        & \leq (f(x^0) - f^{\inf}) - (f(x^K) - f^{\inf}) \notag \\
        &\quad - \gamma(1-L\gamma)\sum\limits_{k=0}^{K-1}\langle \nabla f(x^k), \theta_k \rangle \notag\\
        &\quad + \frac{L\gamma^2}{2} \sum\limits_{k=0}^{K-1} \|\theta_k\|^2,\label{ch2:eq:main_opt_lemma_clipped_SGD_non_convex_chapter_2}\\
        \theta_{k} & \eqdef \thesischtwotnabla f(x^{k};\xi^k) - \nabla f(x^{k}).\label{ch2:eq:theta_k_def_clipped_SGD_non_convex_chapter_2}
    \end{align}
\end{lemma}
\begin{proof}
    Using $x^{k+1} = x^k - \gamma \thesischtwotnabla f(x^{k};\xi^k)$ and smoothness of $f$ on the segment $[x^k, x^{k+1}]$ (Assumption~\ref{ch2:as:L_smoothness_chapter_2}) we get that for all $k = 0,1,\ldots, K-1$
\begin{align*}
        f(x^{k+1}) & \leq f(x^k) + \langle \nabla f(x^k), x^{k+1}-x^k\rangle + \frac{L}{2}\|x^{k+1} - x^k\|^2\notag\\
        & = f(x^k) - \gamma\langle \nabla f(x^k), \thesischtwotnabla f(x^{k};\xi^k) \rangle + \frac{L\gamma^2}{2}\|\thesischtwotnabla f(x^{k};\xi^k)\|^2\notag\\ 
        & \overset{\eqref{ch2:eq:theta_k_def_clipped_SGD_non_convex_chapter_2}}{=} f(x^k) - \gamma\|\nabla f(x^k)\|^2 - \gamma \langle\nabla f(x^k), \theta_k\rangle + \frac{L\gamma^2}{2}\|\theta_k\|^2\notag\\ 
        &\quad + \frac{L\gamma^2}{2}\|\nabla f(x^k)\|^2 + L\gamma^2 \langle \nabla f(x^k), \theta_k\rangle\notag\\ 
        & = f(x^k) - \gamma\left(1-\frac{L\gamma}{2}\right)\|\nabla f(x^k)\|^2 \\
        &\quad - \gamma(1-L\gamma)\langle \nabla f(x^k), \theta_k\rangle + \frac{L\gamma^2}{2}\|\theta_k\|^2.
\end{align*}
We rearrange the terms and get
\begin{align*}
        \gamma\left(1-\frac{L\gamma}{2}\right)\|\nabla f(x^k)\|^2 & \leq f(x^{k}) - f(x^{k+1}) \\
        &\quad -  \gamma(1-L\gamma)\langle \nabla f(x^k), \theta_k\rangle + \frac{L\gamma^2}{2}\|\theta_k\|^2.
\end{align*}
Finally, summing up these inequalities for $k=0,\ldots,K-1$, we get
\begin{align*}
        \MoveEqLeft[1] \gamma\left(1-\frac{L\gamma}{2}\right)\sum\limits_{k=0}^{K-1}\|\nabla f(x^k)\|^2 \\
        & \leq \sum\limits_{k=0}^{K-1}\left(f(x^{k}) -  f(x^{k+1})\right) \\
        &\quad -  \gamma(1-L\gamma)\sum\limits_{k=0}^{K-1}\langle \nabla f(x^k), \theta_k\rangle \notag\\
        &\quad + \frac{L\gamma^2}{2}\sum\limits_{k=0}^{K-1}\|\theta_k\|^2\\
        & = (f(x^0) - f^{\inf}) - (f(x^K) - f^{\inf}) \\
        &\quad - \gamma(1-L\gamma)\sum\limits_{k=0}^{K-1}\langle \nabla f(x^k), \theta_k \rangle \notag\\
        &\quad + \frac{L\gamma^2}{2} \sum\limits_{k=0}^{K-1} \|\theta_k\|^2,
\end{align*}
which concludes the proof.
\end{proof}

Using this lemma, we prove the main convergence result for \algname{clipped-SGD} in the non-convex case.
\begin{theorem}\label{ch2:thm:clipped_SGD_non_convex_main_chapter_2}
Let Assumptions~\ref{ch2:as:bounded_alpha_moment_chapter_2}, \ref{ch2:as:lower_boundedness_chapter_2}, \ref{ch2:as:L_smoothness_chapter_2} hold on
\begin{equation*}
    Q = \left\{x \in \R^d\mid \exists y\in \R^d:\; f(y) \leq f^{\inf} + 2\Delta \text{ and } \|x-y\|\leq \frac{\sqrt{\Delta}}{20\sqrt{L}}\right\},
\end{equation*}
where $\Delta \geq \Delta_0 = f(x^0) - f^{\inf}$. Let $K \geq 0$ and $\beta \in (0,1]$ be such that $A = \ln(4(K+1)/\beta) \geq 1$, and choose stepsize
\begin{gather}
        \gamma \leq \min\left\{\frac{1}{80L A}, \; \frac{\sqrt{\Delta}}{27^{1/p}\cdot 20 \sigma \sqrt{L} (K+1)^{1/p}A^{(p-1)/p}}\right\}\label{ch2:eq:clipped_SGD_non_convex_step_size_chapter_2}
\end{gather}
and clipping level
 \begin{gather}
        \lambda_k = \lambda = \frac{\sqrt{\Delta}}{20\sqrt{L}\gamma A}. \label{ch2:eq:clipped_SGD_non_convex_clipping_level_chapter_2}
\end{gather}
Then, after $K+1$ iterations of \algname{clipped-SGD} the iterates with probability at least $1 - \beta$ satisfy
\begin{equation}
       \frac{1}{K+1}\sum\limits_{k=0}^{K}\|\nabla f(x^k)\|^2 \leq \frac{2\Delta}{\gamma\left(1-\frac{L\gamma}{2}\right)(K+1)}. \label{ch2:eq:clipped_SGD_non_convex_case_appendix_chapter_2}
\end{equation}
In particular, when $\gamma$ equals the minimum from \eqref{ch2:eq:clipped_SGD_non_convex_step_size_chapter_2}, then the iterates produced by \algname{clipped-SGD} after $K+1$ iterations with probability at least $1-\beta$ satisfy
    \begin{align}
        \MoveEqLeft \frac{1}{K+1}\sum\limits_{k=0}^{K}\|\nabla f(x^k)\|^2 \notag\\
        &= \cO\biggl(\max\biggl\{\frac{L\Delta \ln(K/\beta)}{K}, \frac{\sqrt{L\Delta} \sigma \ln^{(p-1)/p}(K/\beta)}{K^{(p-1)/p}}\biggr\}\biggr), \label{ch2:eq:clipped_SSTM_non_convex_case_2_appendix_chapter_2}
    \end{align}
    meaning that to achieve $\frac{1}{K+1}\sum_{k=0}^{K}\|\nabla f(x^k)\|^2 \leq \varepsilon$ with probability at least $1 - \beta$ \algname{clipped-SGD} requires
    \begin{align}
        K = \cO\biggl(\max\biggl\{&\frac{L\Delta}{\varepsilon} \ln(L\Delta/(\beta\varepsilon)),\notag\\
        &\left(\frac{\sqrt{L\Delta}\sigma}{\varepsilon}\right)^{p/(p-1)} \ln \left(\frac{1}{\beta} \left(\frac{\sqrt{L\Delta}\sigma}{\varepsilon}\right)^{p/(p-1)}\right) \biggr\}\biggr) \label{ch2:eq:clipped_SSTM_non_convex_case_complexity_appendix_chapter_2}
    \end{align}
    iterations/oracle calls.
\end{theorem}

\begin{proof}
    Let $\Delta_k = f(x^k) - f^{\inf}$ for all $k\geq 0$. %
    Next, our goal is to show by induction that $\Delta_{t} \leq 2\Delta$ with high probability, which allows to apply the result of Lemma~\ref{ch2:lem:main_opt_lemma_clipped_SGD_non_convex_chapter_2} and then use Bernstein's inequality to estimate the stochastic part of the upper-bound. More precisely, for each $k = 0,\ldots, K+1$ we consider probability event $E_k$ defined as follows: inequalities
    \begin{align}
        \frac{L\gamma^2}{2}\sum\limits_{l=0}^{t-1}\|\theta_l\|^2 - \gamma(1-L\gamma)\sum\limits_{l=0}^{t-1}\langle \nabla f(x^l), \theta_l\rangle & \leq \Delta, \label{ch2:eq:clipped_SGD_non_convex_induction_inequality_1_chapter_2}\\
        \Delta_t & \leq 2\Delta \label{ch2:eq:clipped_SGD_non_convex_induction_inequality_2_chapter_2}
    \end{align}
    hold for all $t = 0,1,\ldots, k$ simultaneously. We want to prove via induction that $\Prob\{E_k\} \geq 1 - k\beta/(K+1)$ for all $k = 0,1,\ldots, K+1$. For $k = 0$ the statement is trivial. Assume that the statement is true for some $k = T - 1 \leq K$: $\Prob\{E_{T-1}\} \geq 1 - (T-1)\beta/(K+1)$. One needs to prove that $\Prob\{E_{T}\} \geq 1 - T\beta/(K+1)$. First, we notice that probability event $E_{T-1}$ implies that $\Delta_t \leq 2\Delta$ for all $t = 0,1, \ldots, T-1$, i.e., $x^t \in \{y \in \R^d\mid f(y) \leq f^{\inf} + 2\Delta\}$ for $t = 0,1, \ldots, T-1$. Moreover, due to the choice of clipping level $\lambda$ we have for all $t = 0,1,\ldots, T-1$
    \begin{equation*}
        \|x^{t+1} - x^{t}\| = \gamma \|\thesischtwotnabla f(x^t;\xi^t)\| \leq \gamma \lambda \overset{\eqref{ch2:eq:clipped_SGD_non_convex_clipping_level_chapter_2}}{=} \frac{\sqrt{\Delta}}{20\sqrt{L}A} \leq \frac{\sqrt{\Delta}}{20\sqrt{L}}.
    \end{equation*}

    Therefore, $E_{T-1}$ implies that the segments $[x^t, x^{t+1}]$, $t = 0,1,\ldots, T-1$, lie in $Q$, meaning that the assumptions of Lemma~\ref{ch2:lem:main_opt_lemma_clipped_SGD_non_convex_chapter_2} are satisfied and we have
    \begin{align}
        \MoveEqLeft[1] \hat\gamma\sum\limits_{l=0}^{t-1}\|\nabla f(x^l)\|^2 \notag \\
        &\leq \Delta_0 - \Delta_t
        - \gamma(1-L\gamma)\sum\limits_{l=0}^{t-1}
        \langle \nabla f(x^l), \theta_l \rangle \notag\\
        &\quad + \frac{L\gamma^2}{2}
        \sum\limits_{l=0}^{t-1} \|\theta_l\|^2,
        \label{ch2:eq:clipped_SGD_non_convex_technical_1_chapter_2}\\
        \hat\gamma
        &\eqdef \gamma\left(1-\frac{L \gamma}{2}\right)
        \geq \frac{\gamma}{2} > 0
    \end{align}
    for all $t = 0,1, \ldots, T$ simultaneously, where the inequality for $\hat\gamma$ follows from \eqref{ch2:eq:clipped_SGD_non_convex_step_size_chapter_2}. Since $\Delta_0 - \Delta_t \leq \Delta$, for all $t = 0, 1, \ldots, T-1$ this probability event also implies that
    \begin{align}
        \MoveEqLeft[1] \sum\limits_{l=0}^{t-1}\|\nabla f(x^l)\|^2 \notag \\
        &\overset{\eqref{ch2:eq:clipped_SGD_non_convex_technical_1_chapter_2}}{\leq}
        \frac{1}{\hat\gamma}\Biggl(
        \Delta - \gamma(1-L\gamma)\sum\limits_{l=0}^{t-1}
        \langle \nabla f(x^l), \theta_l \rangle \notag\\
        &\qquad + \frac{L\gamma^2}{2}
        \sum\limits_{l=0}^{t-1} \|\theta_l\|^2
        \Biggr) \notag\\
        &\overset{\eqref{ch2:eq:clipped_SGD_non_convex_induction_inequality_1_chapter_2}}{\leq}
        \frac{2\Delta}{\hat\gamma}.
        \label{ch2:eq:clipped_SGD_non_convex_technical_1_1_chapter_2}
    \end{align}
    Taking into account that $\hat\gamma\sum_{t=0}^{T-1}\|\nabla f(x^t)\|^2 \geq 0$ and $\Delta_0 \leq \Delta$, we also derive that $E_{T-1}$ implies
    \begin{equation}
        \Delta_T \leq \Delta + \frac{L\gamma^2}{2}\sum\limits_{t=0}^{T-1}\|\theta_t\|^2 - \gamma(1-L\gamma)\sum\limits_{t=0}^{T-1}\langle \nabla f(x^t), \theta_t\rangle. \label{ch2:eq:clipped_SGD_non_convex_technical_2_chapter_2}
    \end{equation}
    Next, we define random vectors
    \begin{equation}
        u^{t} = \begin{cases} \nabla f(x^t),& \text{if } \|\nabla f(x^t)\| \leq 2\sqrt{L\Delta},\\ 0,&\text{otherwise}, \end{cases} \notag
    \end{equation}
    for all $t = 0,1, \ldots, T-1$. By definition these random vectors are bounded with probability 1
    \begin{equation}
        \|u^{t}\| \leq 2\sqrt{L\Delta}. \label{ch2:eq:clipped_SGD_non_convex_technical_6_chapter_2}
    \end{equation}
    Moreover, for $t = 0,1, \ldots, T-1$ event $E_{T-1}$ implies
    \begin{equation}
        \|\nabla f(x^{t})\| \overset{\eqref{ch2:eq:L_smoothness_cor_2_chapter_2}}{\leq} \sqrt{2L(f(x^t) - f^{\inf})} = \sqrt{2L\Delta_t} \leq  2\sqrt{L\Delta} \overset{\eqref{ch2:eq:clipped_SGD_non_convex_step_size_chapter_2},\eqref{ch2:eq:clipped_SGD_non_convex_clipping_level_chapter_2}}{\leq} \frac{\lambda}{2},\label{ch2:eq:clipped_SGD_non_convex_technical_4_chapter_2} 
    \end{equation}
    meaning that $E_{T-1}$ implies that $u^{t} = \nabla f(x^t)$ for all $t = 0,1, \ldots, T-1$. Next, we define the unbiased part and the bias of $\theta_{t}$ as $\theta_{t}^u$ and $\theta_{t}^b$, respectively:
    \begin{equation}
        \begin{aligned}
            \theta_{t}^u &= \thesischtwotnabla f(x^{t};\xi^t) - \Exp_{\xi^t}\left[\thesischtwotnabla f(x^{t};\xi^t)\right],\\
            \theta_{t}^b &= \Exp_{\xi^t}\left[\thesischtwotnabla f(x^{t};\xi^t)\right] - \nabla f(x^{t}).
        \end{aligned}\label{ch2:eq:clipped_SGD_non_convex_theta_u_b_chapter_2}
    \end{equation}
    We notice that $\theta_{t} = \theta_{t}^u + \theta_{t}^b$. Using new notation and $\|\theta_t\|^2 \leq 2\|\theta_t^u\|^2 + 2\|\theta_t^b\|^2$, we get that $E_{T-1}$ implies
    \begin{align}
        \Delta_T 
        & \leq \Delta \underbrace{-\gamma(1-L\gamma)\sum\limits_{t=0}^{T-1}\langle \theta_t^u, u^{t}\rangle}_{\circledOne}  \underbrace{-\gamma(1-L\gamma)\sum\limits_{t=0}^{T-1}\langle \theta_t^b, u^{t}\rangle}_{\circledTwo}\notag\\
        &\quad + \underbrace{L\gamma^2\sum\limits_{t=0}^{T-1}\left(\left\|\theta_{t}^u\right\|^2 - \Exp_{\xi^t}\left[\left\|\theta_{t}^u\right\|^2\right]\right)}_{\circledThree}\notag\\
        &\quad + \underbrace{L\gamma^2\sum\limits_{t=0}^{T-1}\Exp_{\xi^t}\left[\left\|\theta_{t}^u\right\|^2\right]}_{\circledFour} + \underbrace{L\gamma^2\sum\limits_{t=0}^{T-1}\left\|\theta_{t}^b\right\|^2}_{\circledFive}. \label{ch2:eq:clipped_SGD_non_convex_technical_7_chapter_2}
    \end{align}
    It remains to derive good enough high-probability upper-bounds for the terms $\circledOne, \circledTwo, \circledThree, \circledFour, \circledFive$, i.e., to finish our inductive proof we need to show that $\circledOne + \circledTwo + \circledThree + \circledFour + \circledFive \leq \Delta$ with high probability. In the subsequent parts of the proof, we will need to use many times the bounds for the norm and second moments of $\theta_{t}^u$ and $\theta_{t}^b$. First, by definition of clipping operator, we have with probability $1$ that
     \begin{equation}
        \|\theta_{t}^u\| \leq 2\lambda. \label{ch2:eq:clipped_SGD_non_convex_norm_theta_u_bound_chapter_2}
    \end{equation}
    Moreover, since $E_{T-1}$ implies that $\|\nabla f(x^{t})\| \leq \lambda/2$ for $t = 0,1, \ldots, T-1$ (see \eqref{ch2:eq:clipped_SGD_non_convex_technical_4_chapter_2}), then, in view of Lemma~\ref{ch2:lem:bias_and_variance_clip_chapter_2}, we have that $E_{T-1}$ implies
    \begin{align}
        \|\theta_{t}^b\| & \leq \frac{2^{p}\sigma^{p}}{\lambda^{{p}-1}}, \label{ch2:eq:clipped_SGD_non_convex_norm_theta_b_bound_chapter_2} \\
        \Exp_{\xi^t}\left[\|\theta_{t}^u\|^2\right] & \leq 18 \lambda^{2-{p}}\sigma^{p}. \label{ch2:eq:clipped_SGD_non_convex_second_moment_theta_u_bound_chapter_2}
    \end{align}

\textbf{Upper bound for $\circledOne$.} By definition of $\theta_{t}^u$, we have $\Exp_{\xi^t}[\theta_{t}^u] = 0$ and
    \begin{equation}
        \Exp_{\xi^t}\left[-\gamma(1-L\gamma)\langle\theta_t^u, u^{t}\rangle\right] = 0. \notag
    \end{equation}
    Next, sum $\circledOne$ has bounded with probability $1$ terms:
    \begin{align}
        |\gamma(1-L\gamma)\left\la \theta_{t}^u, u^{t}\right\ra| &\overset{\eqref{ch2:eq:clipped_SGD_non_convex_step_size_chapter_2}}{\leq} \gamma \|\theta_{t}^u\| \cdot \|u^{t}\| \notag\\
        &\overset{\eqref{ch2:eq:clipped_SGD_non_convex_technical_6_chapter_2},\eqref{ch2:eq:clipped_SGD_non_convex_norm_theta_u_bound_chapter_2}}{\leq} 4\gamma \lambda \sqrt{L \Delta}\overset{\eqref{ch2:eq:clipped_SGD_non_convex_clipping_level_chapter_2}}{=} \frac{\Delta}{5A} \eqdef c. \label{ch2:eq:clipped_SGD_non_convex_technical_8_chapter_2} 
    \end{align}
    The summands also have bounded conditional variances $\sigma_t^2$ (the conditional expectations $\Exp_{\xi^t}$ of their squares):
    \begin{align}
\sigma_t^2 & \leq \Exp_{\xi^t}\left[\gamma^2(1-L\gamma)^2\|\theta_{t}^u\|^2\cdot \|u^{t}\|^2\right] \notag \\
& \overset{\eqref{ch2:eq:clipped_SGD_non_convex_technical_6_chapter_2}}{\leq} 4\gamma^2(1-L\gamma)^2L\Delta \Exp_{\xi^t}\left[\|\theta_{t}^u\|^2\right] \notag \\
& \overset{\eqref{ch2:eq:clipped_SGD_non_convex_step_size_chapter_2}}{\leq} 4\gamma^2L\Delta \Exp_{\xi^t}\left[\|\theta_{t}^u\|^2\right]. \label{ch2:eq:clipped_SGD_non_convex_technical_9_chapter_2}
    \end{align}
    In other words, we showed that $\{-\gamma(1-L\gamma)\left\la \theta_{t}^u, u^{t}\right\ra\}_{t=0}^{T-1}$ is a bounded martingale difference sequence with bounded conditional variances $\{\sigma_t^2\}_{t=0}^{T-1}$. Next, we apply Bernstein's inequality (Lemma~\ref{ch2:lem:Bernstein_ineq_chapter_2}) with $X_t = -\gamma(1-L\gamma)\left\la \theta_{t}^u, u^{t}\right\ra$, parameter $c$ as in \eqref{ch2:eq:clipped_SGD_non_convex_technical_8_chapter_2}, $b = \frac{\Delta}{5}$, $H = \frac{\Delta^2}{150A}$:
    \begin{equation*}
        \Prob\left\{|\circledOne| > \frac{\Delta}{5}\quad \text{and}\quad \sum\limits_{t=0}^{T-1} \sigma_{t}^2 \leq \frac{\Delta^2}{150A}\right\} \leq 2\exp\left(- \frac{b^2}{2H + 2cb/3}\right) = \frac{\beta}{2(K+1)}.
    \end{equation*}
    Equivalently, we have $\Prob\{ E_{\circledOne} \} \geq 1 - \beta/(2(K+1))$ for
    \begin{equation}
        E_{\circledOne} = \left\{ \text{either }  \sum\limits_{t=0}^{T-1} \sigma_{t}^2 > \frac{\Delta^2}{150A} \text{ or } |\circledOne| \leq \frac{\Delta}{5}\right\}. \label{ch2:eq:clipped_SGD_non_convex_sum_1_upper_bound_chapter_2}
    \end{equation}
    In addition, $E_{T-1}$ implies that
    \begin{align}
        \sum\limits_{t=0}^{T-1} \sigma_{t}^2 & \overset{\eqref{ch2:eq:clipped_SGD_non_convex_technical_9_chapter_2}}{\leq} 4\gamma^2L\Delta \sum\limits_{t=0}^{T-1}  \Exp_{\xi^t}\left[\|\theta_{t}^u\|^2\right] \notag \\
        & \overset{\eqref{ch2:eq:clipped_SGD_non_convex_second_moment_theta_u_bound_chapter_2}}{\leq} 72\gamma^2 L \Delta \sigma^{{p}}T \lambda^{2-{p}} \notag\\
        & \overset{\eqref{ch2:eq:clipped_SGD_non_convex_clipping_level_chapter_2}}{=} \frac{9 \cdot 20^{{p}}\sqrt{\Delta}^{4-{p}}\sigma^{p} T\sqrt{L}^{p}\gamma^{{p}}}{50A^{2-{p}}} \overset{\eqref{ch2:eq:clipped_SGD_non_convex_step_size_chapter_2}}{\leq} \frac{\Delta^2}{150A}. \label{ch2:eq:clipped_SGD_non_convex_sum_1_variance_bound_chapter_2}
    \end{align}

\textbf{Upper bound for $\circledTwo$.} From $E_{T-1}$ it follows that
\begin{align}
        \circledTwo & = -\gamma(1-L\gamma)\sum\limits_{t=0}^{T-1}\langle \theta_t^b, u^{t} \rangle \notag \\
        & \overset{\eqref{ch2:eq:clipped_SGD_non_convex_step_size_chapter_2}}{\leq} \gamma\sum\limits_{t=0}^{T-1}\|\theta_{t}^b\|\cdot \|u^{t}\| \overset{\eqref{ch2:eq:clipped_SGD_non_convex_technical_6_chapter_2},\eqref{ch2:eq:clipped_SGD_non_convex_norm_theta_b_bound_chapter_2}}{\leq}  \frac{2\cdot 2^{p} \gamma \sigma^{p} T \sqrt{L \Delta}}{\lambda^{{p}-1}} \notag\\
        & \overset{\eqref{ch2:eq:clipped_SGD_non_convex_clipping_level_chapter_2}}{=} \frac{40^{p}}{10} \cdot \frac{\sigma^{p} T \sqrt{\Delta}^{2-{p}}\sqrt{L}^{p}\gamma^{p}}{A^{1-{p}}}\overset{\eqref{ch2:eq:clipped_SGD_non_convex_step_size_chapter_2}}{\leq} \frac{\Delta}{5}. \label{ch2:eq:clipped_SGD_non_convex_sum_2_upper_bound_chapter_2}
\end{align}

\textbf{Upper bound for $\circledThree$.} First, we have
    \begin{equation}
        \Exp_{\xi^t}\left[L\gamma^2\left(\left\|\theta_{t}^u\right\|^2 - \Exp_{\xi^t}\left[\left\|\theta_{t}^u\right\|^2\right]\right)\right] = 0. \notag
    \end{equation}
    Next, sum $\circledThree$ has bounded with probability $1$ terms:
    \begin{align}
        \MoveEqLeft[1] \left|L\gamma^2\left(\left\|\theta_{t}^u\right\|^2 - \Exp_{\xi^t}\left[\left\|\theta_{t}^u\right\|^2\right]\right)\right| \notag \\
        & \leq L\gamma^2\left( \|\theta_{t}^u\|^2 +   \Exp_{\xi^t}\left[\left\|\theta_{t}^u\right\|^2\right]\right)\notag\\
        & \overset{\eqref{ch2:eq:clipped_SGD_non_convex_norm_theta_u_bound_chapter_2}}{\leq} 8L\gamma^2\lambda^2\overset{\eqref{ch2:eq:clipped_SGD_non_convex_clipping_level_chapter_2}}{=} \frac{\Delta}{50A^2} \le \frac{\Delta}{5A}\eqdef c. \label{ch2:eq:clipped_SGD_non_convex_technical_10_chapter_2}
    \end{align}
    The summands also have bounded conditional variances $\widetilde\sigma_t^2$ (the conditional expectations $\Exp_{\xi^t}$ of their squares):
    \begin{align}
        \widetilde\sigma_t^2 & \overset{\eqref{ch2:eq:clipped_SGD_non_convex_technical_10_chapter_2}}{\leq} \frac{\Delta}{5A} \Exp_{\xi^t}\left[L\gamma^2\left|\left\|\theta_{t}^u\right\|^2 - \Exp_{\xi^t}\left[\left\|\theta_{t}^u\right\|^2\right]\right|\right] \notag \\
        & \leq \frac{2L\gamma^2 \Delta}{5A} \Exp_{\xi^t}\left[\|\theta_{t}^u\|^2\right], \label{ch2:eq:clipped_SGD_non_convex_technical_11_chapter_2}
    \end{align}
    where in \eqref{ch2:eq:clipped_SGD_non_convex_technical_10_chapter_2} we used $A \geq 1$. In other words, we showed that the summands of $\circledThree$ form a bounded martingale difference sequence with bounded conditional variances $\{\widetilde\sigma_t^2\}_{t=0}^{T-1}$. Next, we apply Bernstein's inequality (Lemma~\ref{ch2:lem:Bernstein_ineq_chapter_2}) with $X_t$ equal to the summands of $\circledThree$, parameter $c$ as in \eqref{ch2:eq:clipped_SGD_non_convex_technical_10_chapter_2}, $b = \frac{\Delta}{5}$, $H = \frac{\Delta^2}{150A}$:
    \begin{equation*}
        \Prob\left\{|\circledThree| > \frac{\Delta}{5}\quad \text{and}\quad \sum\limits_{t=0}^{T-1} \widetilde\sigma_{t}^2 \leq \frac{\Delta^2}{150A}\right\} \leq 2\exp\left(- \frac{b^2}{2H + 2cb/3}\right) = \frac{\beta}{2(K+1)}.
    \end{equation*}
    Equivalently, we have $\Prob\{ E_{\circledThree} \} \geq 1 - \beta/(2(K+1))$ for
    \begin{equation}
        E_{\circledThree} = \left\{ \text{either }  \sum\limits_{t=0}^{T-1} \widetilde\sigma_{t}^2 > \frac{\Delta^2}{150A} \text{ or } |\circledThree| \leq \frac{\Delta}{5}\right\}. \label{ch2:eq:clipped_SGD_non_convex_sum_3_upper_bound_chapter_2}
    \end{equation}
    In addition, $E_{T-1}$ implies that
    \begin{align}
        \sum\limits_{t=0}^{T-1} \widetilde\sigma_{t}^2 & \overset{\eqref{ch2:eq:clipped_SGD_non_convex_technical_11_chapter_2}}{\leq} \frac{2L\gamma^2\Delta}{5A} \sum\limits_{t=0}^{T-1}  \Exp_{\xi^t}\left[\|\theta_{t}^u\|^2\right] \notag \\
        & \overset{\eqref{ch2:eq:clipped_SGD_non_convex_second_moment_theta_u_bound_chapter_2}}{\leq}  
        \frac{36L\gamma^2\Delta\lambda^{2-{p}}\sigma^{{p}}T}{5A} \notag\\
        & \overset{\eqref{ch2:eq:clipped_SGD_non_convex_clipping_level_chapter_2}}{=}
        \frac{9\cdot 20^{p}}{500}\cdot\frac{\sigma^{p} T \sqrt{\Delta}^{4-{p}}\sqrt{L}^{{p}}\gamma^{{p}}}{A^{3-{p}}} \overset{\eqref{ch2:eq:clipped_SGD_non_convex_step_size_chapter_2}}{\leq} \frac{\Delta^2}{150A}. \label{ch2:eq:clipped_SGD_non_convex_sum_3_variance_bound_chapter_2}
    \end{align}

\textbf{Upper bound for $\circledFour$.} From $E_{T-1}$ it follows that
\begin{align}
        \circledFour & = L\gamma^2\sum\limits_{t=0}^{T-1}\Exp_{\xi^t}\left[\left\|\theta_{t}^u\right\|^2\right] \overset{\eqref{ch2:eq:clipped_SGD_non_convex_second_moment_theta_u_bound_chapter_2}}{\leq} 18L\gamma^2\lambda^{2-{p}}\sigma^{{p}}T \notag \\
        & \overset{\eqref{ch2:eq:clipped_SGD_non_convex_clipping_level_chapter_2}}{=}
        \frac{9\cdot 20^{{p}}}{200}\cdot\frac{\sqrt{L}^{{p}}\gamma^{{p}}\sigma^{{p}}T\sqrt{\Delta}^{2-{p}}}{A^{2-{p}}} \overset{\eqref{ch2:eq:clipped_SGD_non_convex_step_size_chapter_2}}{\leq} \frac{\Delta}{5}. \label{ch2:eq:clipped_SGD_non_convex_sum_4_upper_bound_chapter_2}
\end{align}

\textbf{Upper bound for $\circledFive$.} From $E_{T-1}$ it follows that
\begin{align}
        \circledFive & = L\gamma^2\sum\limits_{t=0}^{T-1}\left\|\theta_{t}^b\right\|^2 \overset{\eqref{ch2:eq:clipped_SGD_non_convex_norm_theta_b_bound_chapter_2}}{\leq} \frac{4^{p} \sigma^{2{p}}TL\gamma^2}{\lambda^{2({p}-1)}} \notag \\
        & \overset{\eqref{ch2:eq:clipped_SGD_non_convex_clipping_level_chapter_2}}{=}   \frac{1600^{p}}{400}\cdot \frac{\sigma^{2{p}}TL^{p}\gamma^{2{p}}\Delta^{1-{p}}}{A^{2(1-{p})}} \overset{\eqref{ch2:eq:clipped_SGD_non_convex_step_size_chapter_2}}{\leq} \frac{\Delta}{5}. \label{ch2:eq:clipped_SGD_non_convex_sum_5_upper_bound_chapter_2}
\end{align}

Now, we have the upper bounds for  $\circledOne, \circledTwo, \circledThree, \circledFour, \circledFive$. In particular, probability event $E_{T-1}$ implies
\begin{gather*}
        \Delta_T \overset{\eqref{ch2:eq:clipped_SGD_non_convex_technical_7_chapter_2}}{\leq} \Delta + \circledOne + \circledTwo + \circledThree + \circledFour + \circledFive,\\
        \circledTwo \overset{\eqref{ch2:eq:clipped_SGD_non_convex_sum_2_upper_bound_chapter_2}}{\leq} \frac{\Delta}{5},\quad \circledFour \overset{\eqref{ch2:eq:clipped_SGD_non_convex_sum_4_upper_bound_chapter_2}}{\leq} \frac{\Delta}{5}, \quad \circledFive \overset{\eqref{ch2:eq:clipped_SGD_non_convex_sum_5_upper_bound_chapter_2}}{\leq} \frac{\Delta}{5},\\
        \sum\limits_{t=0}^{T-1} \sigma_t^2 \overset{\eqref{ch2:eq:clipped_SGD_non_convex_sum_1_variance_bound_chapter_2}}{\leq} \frac{\Delta^2}{150A},\quad \sum\limits_{t=0}^{T-1} \widetilde\sigma_t^2 \overset{\eqref{ch2:eq:clipped_SGD_non_convex_sum_3_variance_bound_chapter_2}}{\leq} \frac{\Delta^2}{150A}.
\end{gather*}
    Moreover, we also have (see \eqref{ch2:eq:clipped_SGD_non_convex_sum_1_upper_bound_chapter_2}, \eqref{ch2:eq:clipped_SGD_non_convex_sum_3_upper_bound_chapter_2} and our induction assumption)
\begin{equation*}
        \Prob\{E_{T-1}\} \geq 1 - \frac{(T-1)\beta}{K+1},\quad \Prob\{E_{\circledOne}\} \geq 1 - \frac{\beta}{2(K+1)},\quad \Prob\{E_{\circledThree}\} \geq 1 - \frac{\beta}{2(K+1)},
\end{equation*}
    where 
\begin{align*}
        E_{\circledOne} & = \left\{ \text{either} \quad  \sum\limits_{t=0}^{T-1} \sigma_{t}^2 > \frac{\Delta^2}{150A} \quad \text{or}\quad |\circledOne| \leq \frac{\Delta}{5}\right\},\\
        E_{\circledThree} & = \left\{ \text{either} \quad  \sum\limits_{t=0}^{T-1} \widetilde\sigma_{t}^2 > \frac{\Delta^2}{150A} \quad \text{or}\quad |\circledThree| \leq \frac{\Delta}{5}\right\}.
\end{align*}
    Thus, probability event $E_{T-1} \cap E_{\circledOne} \cap E_{\circledThree}$ implies
\begin{align}
        \Delta_T & \leq \Delta + \frac{\Delta}{5} + \frac{\Delta}{5} + \frac{\Delta}{5} + \frac{\Delta}{5} + \frac{\Delta}{5} = 2\Delta, \notag
\end{align}
    which is equivalent to \eqref{ch2:eq:clipped_SGD_non_convex_induction_inequality_1_chapter_2} and \eqref{ch2:eq:clipped_SGD_non_convex_induction_inequality_2_chapter_2} for $t = T$, and 
\begin{align*}
        \Prob\{E_T\}
        &\geq \Prob\left\{E_{T-1} \cap E_{\circledOne} \cap E_{\circledThree}\right\} \\
        &= 1 - \Prob\left\{\overline{E}_{T-1} \cup \overline{E}_{\circledOne} \cup \overline{E}_{\circledThree}\right\} \\
        &\geq 1 - \Prob\{\overline{E}_{T-1}\} - \Prob\{\overline{E}_{\circledOne}\} - \Prob\{\overline{E}_{\circledThree}\} \\
        &\geq 1 - \frac{T\beta}{K+1}.
\end{align*}
    This finishes the inductive part of our proof, i.e., for all $k = 0,1,\ldots,K+1$ we have $\Prob\{E_k\} \geq 1 - k\beta/(K+1)$. In particular, for $k = K+1$ we have that with probability at least $1 - \beta$
\begin{equation*}
        \frac{1}{K+1}\sum\limits_{k=0}^{K}\|\nabla f(x^k)\|^2 \overset{\eqref{ch2:eq:clipped_SGD_non_convex_technical_1_1_chapter_2}}{\leq}\frac{2\Delta}{\hat\gamma(K+1)} = \frac{2\Delta}{\gamma\left(1-\frac{L\gamma}{2}\right)(K+1)}
\end{equation*}
    and $\{x^k\}_{k=0}^{K+1} \subseteq Q$, which follows from \eqref{ch2:eq:clipped_SGD_non_convex_induction_inequality_2_chapter_2}.
    
Finally, if $\gamma$ equals the minimum in \eqref{ch2:eq:clipped_SGD_non_convex_step_size_chapter_2},
then with probability at least $1-\beta$
\begin{align*}
        \frac{1}{K+1}\sum\limits_{k=0}^{K}\|\nabla f(x^k)\|^2 & \leq \frac{2\Delta}{\gamma\left(1-\frac{L\gamma}{2}\right)(K+1)} \leq \frac{4\Delta}{\gamma(K+1)}\\
        & = \max\left\{\frac{320\Delta L A}{K+1}, \frac{80\cdot 27^{1/p}\sqrt{\Delta}\sigma\sqrt{L}(K+1)^{1/p}A^{(p-1)/p}}{K+1}\right\}\notag\\
        & = \cO\biggl(\max\biggl\{\frac{L\Delta \ln(K/\beta)}{K}, \frac{\sqrt{L\Delta} \sigma \ln^{(p-1)/p}(K/\beta)}{K^{(p-1)/p}}\biggr\}\biggr).
\end{align*}
To get $\frac{1}{K+1}\sum\limits_{k=0}^{K}\|\nabla f(x^k)\|^2 \leq \varepsilon$ with probability at least $1-\beta$ it is sufficient to choose $K$ such that both terms in the maximum above are $\cO(\varepsilon)$. This leads to
\begin{align*}
        K = \cO\biggl(\max\biggl\{&\frac{L\Delta}{\varepsilon} \ln(L\Delta/(\beta\varepsilon)),\notag\\
        &\left(\frac{\sqrt{L\Delta}\sigma}{\varepsilon}\right)^{p/(p-1)} \ln \left(\frac{1}{\beta} \left(\frac{\sqrt{L\Delta}\sigma}{\varepsilon}\right)^{p/(p-1)}\right) \biggr\}\biggr)
    \end{align*}
which concludes the proof.
\end{proof}

\subsection{Polyak-{\L}ojasiewicz Functions}

In this subsection, we provide a high-probability analysis of \algname{clipped-SGD} in the case of Polyak-{\L}ojasiewicz functions. As in the non-convex case, we start with the lemma that handles optimization part of the algorithm and separates it from the stochastic one.

\begin{lemma}\label{ch2:lem:main_opt_lemma_clipped_SGD_PL_chapter_2}
    Let Assumptions~\ref{ch2:as:L_smoothness_chapter_2} and \ref{ch2:as:PL_chapter_2} hold on
    \begin{equation*}
        Q = \left\{x \in \R^d\mid \exists y\in \R^d:\; f(y) \leq f^{\inf} + 2\Delta \text{ and } \|x-y\|\leq \frac{\sqrt{\Delta}}{20\sqrt{L}}\right\},
    \end{equation*}
    where $\Delta \geq \Delta_0 = f(x^0) - f^{\inf}$, and let stepsize  $\gamma$ satisfy $\gamma \leq 1/L$. If the segments $[x^k, x^{k+1}]$ lie in $Q$ for all $k = 0,1,\ldots,K$, $K \ge 0$, then after $K+1$ iterations of \algname{clipped-SGD} we have
    \begin{align}
       f(x^{K+1}) -f^{\inf} & \leq (1-\gamma\mu)^{K+1}(f(x^0) - f^{\inf}) \notag \\
       &\quad - \gamma(1-L\gamma)\sum\limits_{k=0}^K(1-\gamma\mu)^{K-k} \langle \nabla f(x^k), \theta_k \rangle \notag\\
        &\quad + \frac{L\gamma^2}{2} \sum\limits_{k=0}^{K} (1-\gamma\mu)^{K-k} \|\theta_k\|^2,\label{ch2:eq:main_opt_lemma_clipped_SGD_PL_chapter_2}
    \end{align}
    where $\theta_k$ is defined in \eqref{ch2:eq:theta_k_def_clipped_SGD_non_convex_chapter_2}.
\end{lemma}
\begin{proof}
    Using $x^{k+1} = x^k - \gamma \thesischtwotnabla f(x^{k};\xi^k)$ and smoothness of $f$ on the segment $[x^k, x^{k+1}]$ (Assumption~\ref{ch2:as:L_smoothness_chapter_2}) we get that for all $k = 0,1,\ldots, K$
\begin{align*}
        f(x^{k+1}) & \leq f(x^k) + \langle \nabla f(x^k), x^{k+1}-x^k\rangle + \frac{L}{2}\|x^{k+1} - x^k\|^2\notag\\
        & = f(x^k) - \gamma\langle \nabla f(x^k), \thesischtwotnabla f(x^{k};\xi^k) \rangle + \frac{L\gamma^2}{2}\|\thesischtwotnabla f(x^{k};\xi^k)\|^2\notag\\
        & \overset{\eqref{ch2:eq:theta_k_def_clipped_SGD_non_convex_chapter_2}}{=} f(x^k) - \gamma\left(1-\frac{L\gamma}{2}\right)\|\nabla f(x^k)\|^2 \\
        &\quad - \gamma(1-L\gamma)\langle \nabla f(x^k), \theta_k\rangle + \frac{L\gamma^2}{2}\|\theta_k\|^2 \notag\\ 
        & \leq
        f(x^k) - \frac{\gamma}{2}\|\nabla f(x^k)\|^2 - \gamma(1-L\gamma)\langle \nabla f(x^k), \theta_k\rangle \\
        &\quad + \frac{L\gamma^2}{2}\|\theta_k\|^2 \notag\\ 
        & \overset{\eqref{ch2:eq:PL_chapter_2}}{\leq}
        f(x^k) - \gamma \mu (f(x^k) - f^{\inf}) - \gamma(1-L\gamma)\langle \nabla f(x^k), \theta_k\rangle \\
        &\quad + \frac{L\gamma^2}{2}\|\theta_k\|^2.
\end{align*}
where we used $\gamma \leq 1/L$ and, in the last step, \eqref{ch2:eq:PL_chapter_2} together with $f(x^*) \leq f^{\inf}$. By rearranging the terms and subtracting $f^{\inf}$, we obtain 
\begin{align*}
        f(x^{k+1}) - f^{\inf} & \leq (1-\gamma\mu)(f(x^k)-f^{\inf})- \gamma(1-L\gamma)\langle \nabla f(x^k), \theta_k\rangle \\
        &\quad + \frac{L\gamma^2}{2}\|\theta_k\|^2.
\end{align*}
Unrolling the recurrence, we obtain \eqref{ch2:eq:main_opt_lemma_clipped_SGD_PL_chapter_2}.
\end{proof}
\begin{theorem}\label{ch2:thm:clipped_SGD_PL_main_chapter_2}
Let Assumptions~\ref{ch2:as:bounded_alpha_moment_chapter_2}, \ref{ch2:as:L_smoothness_chapter_2}, \ref{ch2:as:PL_chapter_2} hold on
\begin{equation*}
    Q = \left\{x \in \R^d\mid \exists y\in \R^d:\; f(y) \leq f^{\inf} + 2\Delta \text{ and } \|x-y\|\leq \frac{\sqrt{\Delta}}{20\sqrt{L}}\right\},
\end{equation*}
where $\Delta \geq \Delta_0 = f(x^0) - f^{\inf}$. Let $K > 0$ and $\beta \in (0,1]$ be such that $A = \ln(4(K+1)/\beta) \geq 1$, choose stepsize
\begin{align}
        0 < \gamma & \leq \min\left\{\frac{1}{500 L A}, \; \frac{\ln(B_K)}{\mu(K+1)}\right\},\label{ch2:eq:clipped_SGD_PL_step_size_chapter_2}\\
        B_K & = \max\left\{2, \frac{(K+1)^{2(p-1)/p}\mu^2\Delta}{264600^{2/p}L\sigma^2A^{2(p-1)/p}\ln^2(B_K)} \right\} \label{ch2:eq:B_K_SGD_PL_1_chapter_2} \\
        & = \Theta\biggl(\max\biggl\{1, \frac{K^{2(p-1)/p}\mu^2 \Delta}{L\sigma^2\ln^{2(p-1)/p}(K/\beta)}\notag\\
        &\qquad\qquad \cdot \frac{1}{\ln^2\left(\max\left\{2, \frac{K^{2(p-1)/p}\mu^2 \Delta}{L\sigma^2\ln^{2(p-1)/p}(K/\beta)} \right\}\right)} \biggr\}\biggr), \label{ch2:eq:B_K_SGD_PL_2_chapter_2}
\end{align}
and clipping level
 \begin{gather}
        \lambda_k = \frac{\exp(-\gamma\mu(1 + k/2))\sqrt{\Delta}}{120\sqrt{L}\gamma A}, \label{ch2:eq:clipped_SGD_PL_clipping_level_chapter_2}
\end{gather}
where $K$ is large enough that the maximum in \eqref{ch2:eq:B_K_SGD_PL_1_chapter_2} is attained at its second term, i.e., $B_K \ln^2(B_K) = (K+1)^{2(p-1)/p}\mu^2\Delta/(264600^{2/p}L\sigma^2A^{2(p-1)/p})$ with $B_K \geq 2$. Then, after $K+1$ iterations of \algname{clipped-SGD} the iterates with probability at least $1 - \beta$ satisfy
\begin{equation}
       f\left(x^{K+1}\right) - f^{\inf} \leq 2 \exp(-\gamma \mu (K+1))\Delta. \label{ch2:eq:clipped_SGD_PL_case_appendix_chapter_2}
\end{equation}
In particular, when $\gamma$ equals the minimum from \eqref{ch2:eq:clipped_SGD_PL_step_size_chapter_2}, then the iterates produced by \algname{clipped-SGD} after $K+1$ iterations with probability at least $1-\beta$ satisfy
    \begin{align}
        \MoveEqLeft f\left( x^{K+1}\right) - f^{\inf}\notag\\
        &= \cO\Biggl(\max\Biggl\{
        \Delta\exp\left(- \frac{\mu K}{L \ln(K/\beta)}\right),
        \frac{L\sigma^2\ln^{2(p-1)/p}(K/\beta)}{K^{2(p-1)/p}\mu^2}\notag\\
        &\qquad\qquad\cdot\ln^2\left(\max\left\{2,
        \frac{K^{2(p-1)/p}\mu^2\Delta}{L\sigma^2\ln^{2(p-1)/p}(K/\beta)}
        \right\}\right)
        \Biggr\}\Biggr).
        \label{ch2:eq:clipped_SGD_PL_case_2_appendix_chapter_2}
    \end{align}
    meaning that \algname{clipped-SGD} requires the following number of iterations/oracle calls to achieve $f\left( x^{K+1}\right) - f^{\inf} \leq \varepsilon$ with probability at least $1 - \beta$:
    \begin{align}
        K &= \cO\biggl(\max\biggl\{\frac{L}{\mu}\ln\left(\frac{\Delta}{\varepsilon}\right)\ln\left(\frac{L}{\mu \beta}\ln\frac{\Delta}{\varepsilon}\right),\notag\\
        &\quad \left(\frac{L\sigma^2}{\mu^2\varepsilon}\right)^{p/(2(p-1))}\ln \left(\frac{1}{\beta} \left(\frac{L\sigma^2}{\mu^2\varepsilon}\right)^{p/(2(p-1))}\right)\ln^{p/(p-1)}\left(B_\varepsilon\right)\biggr\}\biggr), \label{ch2:eq:clipped_SGD_PL_case_complexity_appendix_chapter_2}
    \end{align}
    where
    \begin{equation*}
        B_\varepsilon = \max\left\{2, \frac{\Delta}{\varepsilon \ln \left(\frac{1}{\beta} \left(\frac{L\sigma^2}{\mu^2\varepsilon}\right)^{p/(2(p-1))}\right)}\right\}.
    \end{equation*}
\end{theorem}
\begin{proof}
    As in the previous results, the proof is based on the induction argument and shows that the iterates do not leave some set with high probability. More precisely, for each $k = 0,1,\ldots,K+1$ we consider probability event $E_k$ as follows: inequalities
    \begin{equation}
        \Delta_t \leq 2 \exp(-\gamma\mu t) \Delta \label{ch2:eq:induction_inequality_SGD_PL_chapter_2}
    \end{equation}
    hold for $t = 0,1,\ldots,k$ simultaneously, where $\Delta_t = f(x^t) - f^{\inf}$. We want to prove $\Prob\{E_k\} \geq  1 - k\beta/(K+1)$ for all $k = 0,1,\ldots,K+1$ by induction. The base of the induction is trivial: for $k=0$ we have $\Delta_0 \leq \Delta < 2\Delta$ by definition. Next, assume that for $k = T-1 \leq K$ the statement holds: $\Prob\{E_{T-1}\} \geq  1 - (T-1)\beta/(K+1)$. Given this, we need to prove $\Prob\{E_{T}\} \geq  1 - T\beta/(K+1)$. Since $\Delta_t \leq 2\exp(-\gamma\mu t) \Delta \leq 2\Delta$, event $E_{T-1}$ implies $x^t \in \{y\in \R^d\mid f(y) \leq f^{\inf} + 2\Delta\}$ for $t = 0,1, \ldots, T-1$. Thus, $E_{T-1}$ implies
    \begin{align}
        \|\nabla f(x^t)\| & \overset{\eqref{ch2:eq:L_smoothness_cor_2_chapter_2}}{\leq} \sqrt{2L(f(x^t) - f^{\inf})} \overset{\eqref{ch2:eq:induction_inequality_SGD_PL_chapter_2}}{\leq} 2\sqrt{L\exp(- \gamma\mu t)\Delta} \overset{\eqref{ch2:eq:clipped_SGD_PL_step_size_chapter_2},\eqref{ch2:eq:clipped_SGD_PL_clipping_level_chapter_2}}{\leq} \frac{\lambda_t}{2} \label{ch2:eq:operator_bound_x_t_SGD_PL_chapter_2}
    \end{align}
    for all $t = 0, 1, \ldots, T-1$, where the last inequality uses $\gamma \leq 1/(500LA)$ and $\exp(\gamma\mu) \leq \exp(1/500)$. Moreover, for all $t = 0, 1, \ldots, T-1$
    \begin{equation*}
        \|x^{t+1} - x^{t}\| = \gamma \|\thesischtwotnabla f(x^t;\xi^t)\| \leq \gamma\lambda_t \overset{\eqref{ch2:eq:clipped_SGD_PL_clipping_level_chapter_2}}{\leq} \frac{\sqrt{\Delta}}{20\sqrt{L}},
    \end{equation*}
    meaning that $E_{T-1}$ implies that the segments $[x^t, x^{t+1}]$, $t = 0, 1, \ldots, T-1$, lie in $Q$. Using Lemma~\ref{ch2:lem:main_opt_lemma_clipped_SGD_PL_chapter_2}, $\Delta_0 \leq \Delta$ and $(1 - \gamma\mu)^T \leq \exp(-\gamma\mu T)$, we obtain that $E_{T-1}$ implies
    \begin{align}
        \Delta_T & \leq \exp(-\gamma\mu T)\Delta \notag \\
        &\quad - \gamma (1-L\gamma) \sum\limits_{l=0}^{T-1} (1-\gamma\mu)^{T-1-l} \langle \nabla f(x^l), \theta_l \rangle \notag\\
        &\quad\quad + \frac{L\gamma^2}{2} \sum\limits_{l=0}^{T-1} (1-\gamma\mu)^{T-1-l} \|\theta_l\|^2. \notag
    \end{align}
    To handle the sums above, we introduce a new notation:
    \begin{gather}
        u^{t} = \begin{cases} \nabla f(x^t),& \text{if } \|\nabla f(x^t)\| \leq 2\sqrt{L}\exp(-\gamma\mu t/2)\sqrt{\Delta},\\ 0,& \text{otherwise}, \end{cases} \label{ch2:eq:eta_t_SGD_PL_chapter_2}
    \end{gather}
    for $t = 0, 1, \ldots, T-1$. These vectors are bounded almost surely:
    \begin{equation}
        \|u^{t}\| \leq 2\sqrt{L}\exp(-\gamma\mu t/2)\sqrt{\Delta} \label{ch2:eq:zeta_t_eta_t_bound_SGD_PL_chapter_2} 
    \end{equation}
    for all $t = 0, 1, \ldots, T-1$. In other words, $E_{T-1}$ implies $u^{t} = \nabla f(x^t)$ for all $t = 0,1, \ldots, T-1$, meaning that from $E_{T-1}$ it follows that 
    \begin{align}
        \Delta_T & \leq \exp(-\gamma\mu T)\Delta - \gamma(1-L\gamma) \sum\limits_{l=0}^{T-1} (1-\gamma\mu)^{T-1-l} \langle u^{l}, \theta_l \rangle \notag\\
        &\quad + \frac{L\gamma^2}{2} \sum\limits_{l=0}^{T-1} (1-\gamma\mu)^{T-1-l} \|\theta_l\|^2. \notag
    \end{align}
    To handle the sums appeared on the right-hand side of the previous inequality we consider unbiased and biased parts of $\theta_l$:
    \begin{align}
\theta_{t}^u & = \thesischtwotnabla f(x^{t};\xi^t) - \Exp_{\xi^t}\left[\thesischtwotnabla f(x^{t};\xi^t)\right],\quad \theta_{t}^b \notag \\
& = \Exp_{\xi^t}\left[\thesischtwotnabla f(x^{t};\xi^t)\right] - \nabla f(x^{t}). \label{ch2:eq:theta_unbias_bias_SGD_PL_chapter_2}
    \end{align}
    for all $l = 0,\ldots, T-1$. By definition we have $\theta_l = \theta_l^u + \theta_l^b$ for all $l = 0,\ldots, T-1$. Therefore, $E_{T-1}$ implies
    \begin{align}
        \Delta_T
        &\leq \exp(-\gamma\mu T) \Delta
        \notag\\
        &\quad \underbrace{- \gamma (1-L\gamma) \sum\limits_{l=0}^{T-1}
        (1-\gamma\mu)^{T-1-l} \langle u^{l}, \theta_l^u \rangle}_{\circledOne}
        \notag\\
        &\quad \underbrace{-\gamma(1-L\gamma) \sum\limits_{l=0}^{T-1}
        (1-\gamma\mu)^{T-1-l} \langle u^{l}, \theta_l^b \rangle}_{\circledTwo}
        \notag\\
        &\quad + \underbrace{L\gamma^2 \sum\limits_{l=0}^{T-1}
        (1-\gamma\mu)^{T-1-l}
        \Exp_{\xi^l}\left[\|\theta_l^u\|^2\right]}_{\circledThree}
        \notag\\
        &\quad + \underbrace{L\gamma^2 \sum\limits_{l=0}^{T-1}
        (1-\gamma\mu)^{T-1-l}
        \left(\|\theta_l^u\|^2
        -\Exp_{\xi^l}\left[\|\theta_l^u\|^2\right]\right)}_{\circledFour}
        \notag\\
        &\quad + \underbrace{L\gamma^2 \sum\limits_{l=0}^{T-1}
        (1-\gamma\mu)^{T-1-l} \|\theta_l^b\|^2}_{\circledFive}.
        \label{ch2:eq:SGD_PL_12345_bound_chapter_2}
    \end{align}
    where we also apply inequality $\|a+b\|^2 \leq 2\|a\|^2 + 2\|b\|^2$ holding for all $a,b \in \R^d$ to upper bound $\|\theta_l\|^2$. It remains to derive good enough high-probability upper-bounds for the terms $\circledOne, \circledTwo, \circledThree, \circledFour, \circledFive$, i.e., to finish our inductive proof we need to show that $\circledOne + \circledTwo + \circledThree + \circledFour + \circledFive \leq \exp(-\gamma\mu T) \Delta$ with high probability. In the subsequent parts of the proof, we will need to use many times the bounds for the norm and second moments of $\theta_{t}^u$ and $\theta_{t}^b$. First, by definition of clipping operator, we have with probability $1$ that
    \begin{equation}
        \|\theta_l^u\| \leq 2\lambda_l.\label{ch2:eq:theta_omega_magnitude_SGD_PL_chapter_2}
    \end{equation}
    Moreover, since $E_{T-1}$ implies that $\|\nabla f(x^l)\| \leq \lambda_l/2$ for all $l = 0,1, \ldots, T-1$ (see \eqref{ch2:eq:operator_bound_x_t_SGD_PL_chapter_2}), from Lemma~\ref{ch2:lem:bias_and_variance_clip_chapter_2} we also have that $E_{T-1}$ implies
    \begin{gather}
        \left\|\theta_l^b\right\| \leq \frac{2^{{p}}\sigma^{p}}{\lambda_l^{{p}-1}}, \label{ch2:eq:bias_theta_omega_SGD_PL_chapter_2}\\
        \Exp_{\xi^l}\left[\left\|\theta_l^u\right\|^2\right] \leq 18 \lambda_l^{2-{p}}\sigma^{p}, \label{ch2:eq:distortion_theta_omega_SGD_PL_chapter_2}
    \end{gather}
    for all $l = 0,1, \ldots, T-1$.
    
\paragraph{Upper bound for $\circledOne$.} By definition of $\theta_l^u$, we have $\Exp_{\xi^l}[\theta_{l}^u] = 0$ and
    \begin{equation*}
        \Exp_{\xi^l}\left[-\gamma (1-L\gamma) (1-\gamma\mu)^{T-1-l} \langle u^{l}, \theta_l^u \rangle\right] = 0.
    \end{equation*}
    Next, sum $\circledOne$ has bounded with probability $1$ terms:
    \begin{align}
        &\left|-\gamma (1-L\gamma)(1-\gamma\mu)^{T-1-l}
        \langle u^{l}, \theta_l^u \rangle \right|
        \notag\\
        &\quad \overset{\eqref{ch2:eq:clipped_SGD_PL_step_size_chapter_2}}{\leq}
        \gamma\exp(-\gamma\mu (T - 1 - l)) \|u^{l}\|\cdot \|\theta_l^u\|
        \notag\\
        &\quad \overset{\eqref{ch2:eq:zeta_t_eta_t_bound_SGD_PL_chapter_2},\eqref{ch2:eq:theta_omega_magnitude_SGD_PL_chapter_2}}{\leq}
        4\sqrt{L\Delta}\gamma
        \exp(-\gamma\mu (T - 1 - l/2))\lambda_l
        \notag\\
        &\quad \overset{\eqref{ch2:eq:clipped_SGD_PL_step_size_chapter_2},\eqref{ch2:eq:clipped_SGD_PL_clipping_level_chapter_2}}{\leq}
        \frac{\exp(-\gamma\mu T)\Delta}{5A}
        \eqdef c.
        \label{ch2:eq:SGD_PL_technical_3_1_chapter_2}
    \end{align}
    The summands also have bounded conditional variances
    \begin{equation*}
        \sigma_l^2 \eqdef
        \Exp_{\xi^l}\left[
        \gamma^2 (1-L\gamma)^2(1-\gamma\mu)^{2T-2-2l}
        \langle u^{l}, \theta_l^u \rangle^2\right]:
    \end{equation*}
    \begin{align}
        \sigma_l^2 & \leq \Exp_{\xi^l}\left[\gamma^2(1-L\gamma)^2\exp(-\gamma\mu (2T - 2 - 2l)) \|u^{l}\|^2\cdot \|\theta_l^u\|^2\right]\notag\\
        & \overset{\eqref{ch2:eq:zeta_t_eta_t_bound_SGD_PL_chapter_2}, \eqref{ch2:eq:clipped_SGD_PL_step_size_chapter_2}}{\leq} 4\gamma^2 L\Delta \exp(-\gamma\mu (2T - 2 - l))  \Exp_{\xi^l}\left[\|\theta_l^u\|^2\right]\notag\\
        & \overset{\eqref{ch2:eq:clipped_SGD_PL_step_size_chapter_2}}{\leq} 10\gamma^2L\Delta\exp(-\gamma\mu (2T - l)) \Exp_{\xi^l}\left[\|\theta_l^u\|^2\right]. \label{ch2:eq:SGD_PL_technical_3_2_chapter_2}
    \end{align}
    In other words, we showed that $\{-\gamma(1-L\gamma) (1-\gamma\mu)^{T-1-l} \langle u^{l}, \theta_l^u \rangle\}_{l = 0}^{T-1}$ is a bounded martingale difference sequence with bounded conditional variances $\{\sigma_l^2\}_{l = 0}^{T-1}$. Next, we apply Bernstein's inequality (Lemma~\ref{ch2:lem:Bernstein_ineq_chapter_2}) with $X_l = -\gamma(1-L\gamma) (1-\gamma\mu)^{T-1-l} \langle u^{l}, \theta_l^u \rangle$, parameter $c$ as in \eqref{ch2:eq:SGD_PL_technical_3_1_chapter_2}, $b = \tfrac{1}{5}\exp(-\gamma\mu T) \Delta$, $H = \tfrac{\exp(-2 \gamma\mu T) \Delta^2}{150A}$:
    \begin{align*}
        &\Prob\left\{
        |\circledOne| > \frac{1}{5}\exp(-\gamma\mu T) \Delta
        \ \text{and}\
        \sum\limits_{l=0}^{T-1}\sigma_l^2
        \leq \frac{\exp(- 2\gamma\mu T) \Delta^2}{150A}
        \right\}\\
        &\qquad \leq
        2\exp\left(- \frac{b^2}{2H + 2cb/3}\right)
        = \frac{\beta}{2(K+1)}.
    \end{align*}
    Equivalently, we have $\Prob\{E_{\circledOne}\} \geq 1 - \frac{\beta}{2(K+1)}$ for
    \begin{equation}
        E_{\circledOne} = \left\{\text{either } \sum\limits_{l=0}^{T-1}\sigma_l^2 > \frac{\exp(- 2\gamma\mu T) \Delta^2}{150A} \text{ or } |\circledOne| \leq \frac{1}{5}{\exp(-\gamma\mu T)} \Delta\right\}. \label{ch2:eq:bound_1_SGD_PL_chapter_2}
    \end{equation}
    In addition, $E_{T-1}$ implies that
    \begin{align}
        \MoveEqLeft \sum\limits_{l=0}^{T-1}\sigma_l^2\notag\\
        &\overset{\eqref{ch2:eq:SGD_PL_technical_3_2_chapter_2}}{\leq} 10\gamma^2L\Delta\exp(- 2\gamma\mu T)\sum\limits_{l=0}^{T-1} \frac{\Exp_{\xi^l}\left[\|\theta_l^u\|^2\right]}{\exp(-\gamma\mu l)}\notag\\
        &\overset{\eqref{ch2:eq:distortion_theta_omega_SGD_PL_chapter_2}, T \leq K+1}{\leq} 180\gamma^2L\Delta\exp(-2\gamma\mu T)  \sigma^{p} \sum\limits_{l=0}^{K} \frac{\lambda_l^{2-{p}}}{\exp(-\gamma\mu l)}\notag\\
        &\overset{\eqref{ch2:eq:clipped_SGD_PL_clipping_level_chapter_2}}{=} \frac{180\gamma^{{p}}\sqrt{L}^{{p}}\sqrt{\Delta}^{4-{p}}\exp(-2\gamma\mu T)  \sigma^{p}}{120^{2-{p}} A^{2-{p}}} \notag\\
        &\qquad\cdot\sum\limits_{l=0}^{K} \frac{1}{\exp(-\gamma\mu l)} \cdot \left(\exp(-\gamma\mu(1 + l/2))\right)^{2-{p}} \notag\\
        &= \frac{180\gamma^{{p}}\sqrt{L}^{{p}}\sqrt{\Delta}^{4-{p}}\exp(-2\gamma\mu T)  \sigma^{p}}{120^{2-{p}} A^{2-{p}}} \sum\limits_{l=0}^{K} \exp(\gamma\mu({p}-2)) \cdot \exp\left(\frac{\gamma\mu{p} l}{2}\right) \notag\\
        &\leq \frac{180\gamma^{{p}}\sqrt{L}^{{p}}\sqrt{\Delta}^{4-{p}}\exp(-2\gamma\mu T)  \sigma^{p}(K+1)\exp(\frac{\gamma\mu{p} K}{2})}{120^{2-{p}} A^{2-{p}}} \notag\\
        &\overset{\eqref{ch2:eq:clipped_SGD_PL_step_size_chapter_2}}{\leq} \frac{\exp(-2\gamma\mu T)\Delta^2}{150A}, \label{ch2:eq:bound_3_variances_SGD_PL_chapter_2}
    \end{align}
\begin{equation}
        \gamma^2 L \Delta\sum\limits_{l=0}^{K} \frac{\lambda_{l}^{2-{p}}}{\exp(-\gamma\mu l)} \leq \frac{\gamma^{p} \sqrt{L}^{p} \sqrt{\Delta}^{4-{p}}(K+1)\exp(\frac{\gamma\mu{p} K}{2})}{120^{2-{p}}A^{2-{p}}}. \label{ch2:eq:useful_inequality_on_lambda_SGD_PL_chapter_2}
\end{equation}

\paragraph{Upper bound for $\circledTwo$.} From $E_{T-1}$ it follows that
    \begin{align}
        \circledTwo
        &\overset{\eqref{ch2:eq:clipped_SGD_PL_step_size_chapter_2}}{\leq}
        \gamma \exp(-\gamma\mu (T-1))
        \sum\limits_{l=0}^{T-1} \frac{\|u^{l}\|\cdot \|\theta_l^b\|}{\exp(-\gamma\mu l)}
        \notag\\
        &\overset{\eqref{ch2:eq:zeta_t_eta_t_bound_SGD_PL_chapter_2}, \eqref{ch2:eq:bias_theta_omega_SGD_PL_chapter_2}}{\leq}
        \begin{multlined}[t]
        2^{1+{p}} \gamma \exp(-\gamma\mu (T-1)) \sqrt{L\Delta} \sigma^{p} \\
        \sum\limits_{l=0}^{T-1}
        \frac{1}{\lambda_l^{{p}-1} \exp(-\gamma\mu l/2)}
        \end{multlined}
        \notag\\
        &\overset{\eqref{ch2:eq:clipped_SGD_PL_clipping_level_chapter_2}}{=}
        \begin{multlined}[t]
        2^{1+{p}}\cdot 120^{{p}-1}\sqrt{L}^{p} \sqrt{\Delta}^{2-{p}}
        \exp(-\gamma\mu (T-1)) \gamma^{{p}} \sigma^{p}
        A^{{p}-1} \\
        \sum\limits_{l=0}^{T-1}
        \frac{1}{\exp\left(-\gamma\mu(1 + l/2)\right)^{{p}-1}
        \cdot\exp(-\gamma\mu l/2)}
        \end{multlined}
        \notag\\
        &\overset{T \leq K+1}{\leq}
        \begin{multlined}[t]
        2^{1+{p}}\cdot 120^{{p}-1}\sqrt{L}^{p} \sqrt{\Delta}^{2-{p}}
        \exp(-\gamma\mu (T-{p})) \gamma^{{p}} \sigma^{p}
        A^{{p}-1} \\
        \sum\limits_{l=0}^{K} \exp\left(\frac{\gamma\mu {p} l}{2}\right)
        \end{multlined}
        \notag\\
        &\leq
        \begin{multlined}[t]
        2^{1+{p}}\cdot 120^{{p}-1}\sqrt{L}^{p} \sqrt{\Delta}^{2-{p}}
        \exp(-\gamma\mu (T-{p})) \gamma^{{p}} \sigma^{p}
        A^{{p}-1} \\
        (K+1)\exp\left(\frac{\gamma\mu {p} K}{2}\right)
        \end{multlined}
        \notag\\
        &\overset{\eqref{ch2:eq:clipped_SGD_PL_step_size_chapter_2}}{\leq}
        \frac{1}{5}\exp(-\gamma\mu T) \Delta.
        \label{ch2:eq:bound_4_SGD_PL_chapter_2}
    \end{align}

\paragraph{Upper bound for $\circledThree$.} From $E_{T-1}$ it follows that
    \begin{align}
        \circledThree & = L\gamma^2 \exp(-\gamma\mu (T-1)) \sum\limits_{l=0}^{T-1} \frac{\Exp_{\xi^l}\left[\|\theta_l^u\|^2\right]}{\exp(-\gamma\mu l)} \notag\\
        & \overset{\eqref{ch2:eq:distortion_theta_omega_SGD_PL_chapter_2}}{\leq} 18 L \gamma^2\exp(-\gamma\mu (T-1)) \sigma^{p}\sum\limits_{l=0}^{T-1} \frac{\lambda_l^{2-{p}}}{\exp(-\gamma\mu l)} \notag\\
        & \overset{\eqref{ch2:eq:useful_inequality_on_lambda_SGD_PL_chapter_2}}{\leq}  \frac{18\gamma^{p} \sqrt{L}^{{p}}\sqrt{\Delta}^{2-{p}}\exp(-\gamma\mu (T-1)) \sigma^{p} (K+1)\exp(\frac{\gamma\mu{p} K}{2})}{120^{2-{p}} A^{2-{p}}} \notag\\
        & \overset{\eqref{ch2:eq:clipped_SGD_PL_step_size_chapter_2}}{\leq} \frac{1}{5} \exp(-\gamma\mu T) \Delta. \label{ch2:eq:bound_5_SGD_PL_chapter_2}
    \end{align}

\paragraph{Upper bound for $\circledFour$.} First, we have
    \begin{equation*}
        L\gamma^2 (1-\gamma\mu)^{T-1-l}\Exp_{\xi^l}\left[\|\theta_l^u\|^2 -\Exp_{\xi^l}\left[\|\theta_l^u\|^2\right] \right] = 0.
    \end{equation*}
    Next, sum $\circledFour$ has bounded with probability $1$ terms:
    \begin{align}
        \MoveEqLeft[1] L\gamma^2 (1-\gamma\mu)^{T-1-l}\left| \|\theta_l^u\|^2 -\Exp_{\xi^l}\left[\|\theta_l^u\|^2\right] \right| \notag \\
        & \overset{\eqref{ch2:eq:theta_omega_magnitude_SGD_PL_chapter_2}}{\leq} \frac{8 L \gamma^2 \exp(-\gamma\mu T) \lambda_l^2}{\exp(-\gamma\mu (1+l))}\notag\\
        & \overset{\eqref{ch2:eq:clipped_SGD_PL_clipping_level_chapter_2}}{=} \frac{\exp(-\gamma\mu (T+1))\Delta}{1800A^2}\notag\\
        & \leq \frac{\exp(-\gamma\mu T)\Delta}{5A}\notag\\
        & \eqdef c. \label{ch2:eq:SGD_PL_technical_6_1_chapter_2}
    \end{align}
    The summands also have conditional variances
    \begin{equation*}
        \widehat\sigma_l^2 \eqdef \Exp_{\xi^l}\left[L^2\gamma^4 (1-\gamma\mu)^{2T-2-2l} \left| \|\theta_l^u\|^2 -\Exp_{\xi^l}\left[\|\theta_l^u\|^2\right]  \right|^2\right]
    \end{equation*}
    that are bounded
    \begin{align}
        \widehat\sigma_l^2 & \overset{\eqref{ch2:eq:SGD_PL_technical_6_1_chapter_2}}{\leq} \frac{L\gamma^2\exp(-2\gamma\mu T)\Delta}{5\exp(-\gamma\mu (1+l))A} \Exp_{\xi^l}\left[\left| \|\theta_l^u\|^2 -\Exp_{\xi^l}\left[\|\theta_l^u\|^2\right] \right|\right]\notag\\
        & \leq \frac{2L\gamma^2\exp(-2\gamma\mu T)\Delta}{5\exp(-\gamma\mu (1+l))A} \Exp_{\xi^l}\left[\|\theta_l^u\|^2\right]. \label{ch2:eq:SGD_PL_technical_6_2_chapter_2}
    \end{align}
    In other words, we showed that $\left\{L\gamma^2 (1-\gamma\mu)^{T-1-l}\left( \|\theta_l^u\|^2 -\Exp_{\xi^l}\left[\|\theta_l^u\|^2\right] \right)\right\}_{l = 0}^{T-1}$ is a bounded martingale difference sequence with bounded conditional variances $\{\widehat\sigma_l^2\}_{l = 0}^{T-1}$. Next, we apply Bernstein's inequality (Lemma~\ref{ch2:lem:Bernstein_ineq_chapter_2}) with $X_l = L\gamma^2 (1-\gamma\mu)^{T-1-l}\left( \|\theta_l^u\|^2 -\Exp_{\xi^l}\left[\|\theta_l^u\|^2\right] \right)$, parameter $c$ as in \eqref{ch2:eq:SGD_PL_technical_6_1_chapter_2}, $b = \tfrac{1}{5}\exp(-\gamma\mu T) \Delta$, $H = \tfrac{\exp(-2 \gamma\mu T) \Delta^2}{150A}$:
    \begin{align*}
        \MoveEqLeft \Prob\left\{|\circledFour| > \frac{1}{5}\exp(-\gamma\mu T) \Delta \text{ and } \sum\limits_{l=0}^{T-1}\widehat\sigma_l^2 \leq \frac{\exp(-2\gamma\mu T) \Delta^2}{150A}\right\}\\
        &\leq 2\exp\left(- \frac{b^2}{2H + 2cb/3}\right) = \frac{\beta}{2(K+1)}.
    \end{align*}
    Equivalently, we have $\Prob\{E_{\circledFour}\} \geq 1 - \frac{\beta}{2(K+1)}$ for
    \begin{equation}
        E_{\circledFour} = \left\{\text{either } \sum\limits_{l=0}^{T-1}\widehat\sigma_l^2 > \frac{\exp(-2\gamma\mu T) \Delta^2}{150A} \text{ or } |\circledFour| \leq \frac{1}{5}{\exp(-\gamma\mu T)} \Delta\right\}. \label{ch2:eq:bound_6_SGD_PL_chapter_2}
    \end{equation}
    In addition, $E_{T-1}$ implies that
    \begin{align}
        \sum\limits_{l=0}^{T-1}\widehat\sigma_l^2 & \overset{\eqref{ch2:eq:SGD_PL_technical_6_2_chapter_2}}{\leq} \frac{2L\gamma^2\exp(-\gamma\mu (2T-1))\Delta}{5A} \sum\limits_{l=0}^{T-1} \frac{\Exp_{\xi^l}\left[\|\theta_l^u\|^2 \right]}{\exp(-\gamma\mu l)}\notag\\
        & \overset{\eqref{ch2:eq:distortion_theta_omega_SGD_PL_chapter_2}, T \leq K+1}{\leq} \frac{36L\gamma^2\exp(-\gamma\mu (2T-1)) \Delta \sigma^{p}}{5A} \sum\limits_{l=0}^{K} \frac{\lambda_l^{2-{p}}}{\exp(-\gamma\mu l)}\notag\\
        & \overset{\eqref{ch2:eq:useful_inequality_on_lambda_SGD_PL_chapter_2}}{\leq} \frac{36\sqrt{L}^{p}\gamma^{p}\exp(-\gamma\mu (2T-1)) \sqrt{\Delta}^{4-{p}} \sigma^{p} (K+1)\exp(\frac{\gamma\mu{p} K}{2})}{5\cdot 120^{2-{p}}A^{3-{p}}} \notag\\
        & \overset{\eqref{ch2:eq:clipped_SGD_PL_step_size_chapter_2}}{\leq} \frac{\exp(-2\gamma\mu T)\Delta^2}{150A}. \label{ch2:eq:bound_6_variances_SGD_PL_chapter_2}
\end{align}

\paragraph{Upper bound for $\circledFive$.} From $E_{T-1}$ it follows that
\begin{align}
        \circledFive
        &= L\gamma^2 \sum\limits_{l=0}^{T-1}
        \exp(-\gamma\mu (T-1-l)) \|\theta_l^b\|^2
        \notag\\
        &\overset{\eqref{ch2:eq:bias_theta_omega_SGD_PL_chapter_2}}{\leq}
        \begin{multlined}[t]
        2^{2{p}}L\gamma^2 \exp(-\gamma\mu (T-1)) \sigma^{2{p}} \\
        \sum\limits_{l=0}^{T-1}
        \frac{1}{\lambda_l^{2{p}-2} \exp(-\gamma\mu l)}
        \end{multlined}
        \notag\\
        &\overset{\eqref{ch2:eq:clipped_SGD_PL_clipping_level_chapter_2}, T \leq K+1}{\leq}
        \begin{multlined}[t]
        \frac{2\cdot 2^{2{p}}\cdot 120^{2{p}-2}\gamma^{2{p}}\sqrt{L}^{2{p}}
        \exp(-\gamma\mu T) \sigma^{2{p}}
        A^{2{p}-2}}{\sqrt{\Delta}^{2{p}-2}} \\
        \sum\limits_{l=0}^{K}
        \exp\left(\gamma\mu(2{p}-2)\left(1 + \frac{l}{2}\right)\right)
        \exp(\gamma\mu l)
        \end{multlined}
        \notag\\
        &\leq
        \begin{multlined}[t]
        \frac{4\cdot 2^{2{p}}\cdot 120^{2{p}-2}\gamma^{2{p}}\sqrt{L}^{2{p}}
        \exp(-\gamma\mu T) \sigma^{2{p}}
        A^{2{p}-2}}{\sqrt{\Delta}^{2{p}-2}} \\
        \sum\limits_{l=0}^{K} \exp(\gamma\mu {p} l)
        \end{multlined}
        \notag\\
        &\leq
        \frac{4\cdot 2^{2{p}}\cdot 120^{2{p}-2}\gamma^{2{p}} \sqrt{L}^{2{p}}
        \exp(-\gamma\mu T) \sigma^{2{p}}
        A^{2{p}-2} (K+1)
        \exp(\gamma\mu {p} K)}{\sqrt{\Delta}^{2{p}-2}}
        \notag\\
        &\overset{\eqref{ch2:eq:clipped_SGD_PL_step_size_chapter_2}}{\leq}
        \frac{1}{5}\exp(-\gamma\mu T) \Delta.
        \label{ch2:eq:bound_7_SGD_PL_chapter_2}
\end{align}

 Now, we have the upper bounds for  $\circledOne, \circledTwo, \circledThree, \circledFour, \circledFive$. In particular, probability event $E_{T-1}$ implies
    \begin{gather*}
        \Delta_T \overset{\eqref{ch2:eq:SGD_PL_12345_bound_chapter_2}}{\leq} \exp(-\gamma\mu T) \Delta + \circledOne + \circledTwo + \circledThree + \circledFour + \circledFive,\\
        \circledTwo \overset{\eqref{ch2:eq:bound_4_SGD_PL_chapter_2}}{\leq} \frac{1}{5}\exp(-\gamma\mu T)\Delta,\quad \circledThree \overset{\eqref{ch2:eq:bound_5_SGD_PL_chapter_2}}{\leq} \frac{1}{5}\exp(-\gamma\mu T)\Delta,\quad \circledFive \overset{\eqref{ch2:eq:bound_7_SGD_PL_chapter_2}}{\leq} \frac{1}{5}\exp(-\gamma\mu T)\Delta,\\
        \sum\limits_{l=0}^{T-1}\sigma_l^2 \overset{\eqref{ch2:eq:bound_3_variances_SGD_PL_chapter_2}}{\leq} \frac{\exp(-2\gamma\mu T)\Delta^2}{150A},\quad \sum\limits_{l=0}^{T-1}\widehat\sigma_l^2 \overset{\eqref{ch2:eq:bound_6_variances_SGD_PL_chapter_2}}{\leq}  \frac{\exp(-2\gamma\mu T)\Delta^2}{150A}.
    \end{gather*}
    Moreover, we also have (see \eqref{ch2:eq:bound_1_SGD_PL_chapter_2}, \eqref{ch2:eq:bound_6_SGD_PL_chapter_2} and our induction assumption)
     \begin{gather*}
        \Prob\{E_{T-1}\} \geq 1 - \frac{(T-1)\beta}{K+1},\\
        \Prob\{E_{\circledOne}\} \geq 1 - \frac{\beta}{2(K+1)}, \quad \Prob\{E_{\circledFour}\} \geq 1 - \frac{\beta}{2(K+1)},
    \end{gather*}
    where
    \begin{align}
        E_{\circledOne}& = \left\{\text{either } \sum\limits_{l=0}^{T-1}\sigma_l^2 > \frac{\exp(-2\gamma\mu T) \Delta^2}{150A} \text{ or } |\circledOne| \leq \frac{1}{5}{\exp(-\gamma\mu T)} \Delta\right\},\notag\\
        E_{\circledFour}& = \left\{\text{either } \sum\limits_{l=0}^{T-1}\widehat\sigma_l^2 > \frac{\exp(-2\gamma\mu T) \Delta^2}{150A} \text{ or } |\circledFour| \leq \frac{1}{5}{\exp(-\gamma\mu T)} \Delta\right\}.\notag
    \end{align}
    Thus, probability event $E_{T-1} \cap E_{\circledOne} \cap E_{\circledFour}$ implies
    \begin{align*}
        \Delta_T & \overset{\eqref{ch2:eq:SGD_PL_12345_bound_chapter_2}}{\leq} \exp(-\gamma\mu T) \Delta + \circledOne + \circledTwo + \circledThree + \circledFour + \circledFive\\
        & \leq 2\exp(-\gamma\mu T) \Delta,
    \end{align*}
    which is equivalent to \eqref{ch2:eq:induction_inequality_SGD_PL_chapter_2} for $t = T$, and
    \begin{equation}
        \Prob\{E_T\} \geq \Prob\{E_{T-1} \cap E_{\circledOne} \cap E_{\circledFour}\} = 1 - \Prob\{\overline{E}_{T-1} \cup \overline{E}_{\circledOne} \cup \overline{E}_{\circledFour}\} \geq 1 - \frac{T\beta}{K+1}. \notag
    \end{equation}
    This finishes the inductive part of our proof, i.e., for all $k = 0,1,\ldots,K+1$ we have $\Prob\{E_k\} \geq 1 - k\beta/(K+1)$. In particular, for $k = K+1$ we have that with probability at least $1 - \beta$
    \begin{equation}
        f(x^{K+1}) - f^{\inf} \leq 2\exp(-\gamma\mu (K+1))\Delta. \notag
    \end{equation}
    Finally, if 
    \begin{align*}
        \gamma &= \min\left\{\frac{1}{500 L A}, \frac{\ln(B_K)}{\mu(K+1)}\right\},\\
        B_K &= \max\left\{2, \frac{(K+1)^{2(p-1)/p}\mu^2\Delta}{264600^{2/p}L\sigma^2A^{2(p-1)/p}\ln^2(B_K)} \right\}\\
        &= \cO\biggl(\max\biggl\{2, \frac{K^{2(p-1)/p}\mu^2 \Delta}{L\sigma^2\ln^{2(p-1)/p}(K/\beta)}\\
        &\qquad\qquad \cdot \frac{1}{\ln^2\left(\max\left\{2, \frac{K^{2(p-1)/p}\mu^2 \Delta}{L\sigma^2\ln^{2(p-1)/p}(K/\beta)} \right\}\right)} \biggr\}\biggr),
    \end{align*}
    then with probability at least $1-\beta$
    \begin{align*}
        f(x^{K+1}) - f^{\inf}
        &\leq 2\exp(-\gamma\mu (K+1))\Delta\\
        &= 2\Delta\max\left\{\exp\left(-\frac{\mu(K+1)}{500 L A}\right), \frac{1}{B_K}\right\}\\
        &= \cO\Biggl(\max\Biggl\{
        \Delta\exp\left(- \frac{\mu K}{L \ln(K/\beta)}\right),
        \frac{L\sigma^2\ln^{2(p-1)/p}(K/\beta)}{K^{2(p-1)/p}\mu^2}\\
        &\qquad\qquad\cdot\ln^2\left(\max\left\{2,
        \frac{K^{2(p-1)/p}\mu^2\Delta}{L\sigma^2\ln^{2(p-1)/p}(K/\beta)}
        \right\}\right)
        \Biggr\}\Biggr).
    \end{align*}
    To get $f(x^{K+1}) - f^{\inf} \leq \varepsilon$ with probability at least $1-\beta$ it is sufficient to choose $K$ such that both terms in the maximum above are $\cO(\varepsilon)$. This leads to
    \begin{align*}
        K &= \cO\biggl(\max\biggl\{\frac{L}{\mu}\ln\left(\frac{\Delta}{\varepsilon}\right)\ln\left(\frac{L}{\mu \beta}\ln\frac{\Delta}{\varepsilon}\right),\notag\\
        &\quad \left(\frac{L\sigma^2}{\mu^2\varepsilon}\right)^{p/(2(p-1))}\ln \left(\frac{1}{\beta} \left(\frac{L\sigma^2}{\mu^2\varepsilon}\right)^{p/(2(p-1))}\right)\ln^{p/(p-1)}\left(B_\varepsilon\right)\biggr\}\biggr),
    \end{align*}
    where
    \begin{equation*}
        B_\varepsilon = \max\left\{2, \frac{\Delta}{\varepsilon \ln \left(\frac{1}{\beta} \left(\frac{L\sigma^2}{\mu^2\varepsilon}\right)^{p/(2(p-1))}\right)}\right\}.
    \end{equation*}
    This concludes the proof.
\end{proof}

\subsection{Convex Functions}

Now, we focus on the case of convex functions. We start with the following lemma.
\begin{lemma}\label{ch2:lem:main_opt_lemma_clipped_SGD_convex_chapter_2}
    Let Assumptions~\ref{ch2:as:L_smoothness_chapter_2} and \ref{ch2:as:str_cvx_chapter_2} with $\mu = 0$ hold on $Q = B_{2R}(x^*)$, where $R \geq \|x^0 - x^*\|$, and let stepsize  $\gamma$ satisfy $\gamma \leq 1/(2L)$. If $x^{k} \in Q$ for all $k = 0,1,\ldots,K$, $K \ge 0$, then after $K+1$ iterations of \algname{clipped-SGD} we have
    \begin{align}
       \gamma\left(f(\overline{x}^K) -f(x^*)\right)
       &\leq \frac{\|x^0 - x^*\|^2 - \|x^{K+1} - x^*\|^2}{K+1}
       \notag\\
       &\quad - \frac{2\gamma}{K+1}\sum\limits_{k=0}^K
       \langle x^k - x^* - \gamma \nabla f(x^k), \theta_k \rangle
       \notag\\
       &\quad + \frac{\gamma^2}{K+1}\sum\limits_{k=0}^K \|\theta_k\|^2,
       \label{ch2:eq:main_opt_lemma_clipped_SGD_star_convex_chapter_2}\\
       \overline{x}^K
       &= \frac{1}{K+1}\sum\limits_{k=0}^K x^k.
       \label{ch2:eq:x_avg_clipped_SGD_chapter_2}
    \end{align}
    where $\theta_k$ is defined in \eqref{ch2:eq:theta_k_def_clipped_SGD_non_convex_chapter_2}.
\end{lemma}
\begin{proof}
    Using $x^{k+1} = x^k - \gamma \thesischtwotnabla f(x^{k};\xi^k)$, we derive for all $k = 0,1, \ldots, K$ that
    \begin{align*}
        \|x^{k+1} - x^*\|^2
        & = \|x^k - x^*\|^2 -2\gamma \langle x^k - x^*, \thesischtwotnabla f(x^{k};\xi^k)\rangle \\
        &\quad + \gamma^2\|\thesischtwotnabla f(x^k;\xi^k)\|^2\\
        &= \|x^k - x^*\|^2
        -2\gamma \langle x^k - x^*, \nabla f(x^{k})\rangle
        - 2\gamma \langle x^k - x^*, \theta_k\rangle\\
        &\quad + \gamma^2\|\nabla f(x^k) + \theta_k\|^2\\
        &\leq
        \|x^k - x^*\|^2
        -2\gamma\left(f(x^k) - f(x^*)\right)\\
        &\quad - 2\gamma\langle x^k - x^* - \gamma \nabla f(x^k), \theta_k \rangle
        + \gamma^2\|\nabla f(x^k)\|^2 + \gamma^2 \|\theta_k\|^2\\
        &\leq
        \|x^k - x^*\|^2
        -2\gamma\left(1 - \gamma L\right)\left(f(x^k) - f(x^*)\right)\\
        &\quad - 2\gamma\langle x^k - x^* - \gamma \nabla f(x^k), \theta_k \rangle
        + \gamma^2 \|\theta_k\|^2\\
        &\leq
        \|x^k - x^*\|^2 -\gamma \left(f(x^k) - f(x^*)\right)\\
        &\quad - 2\gamma\langle x^k - x^* - \gamma \nabla f(x^k), \theta_k \rangle
        + \gamma^2 \|\theta_k\|^2.
    \end{align*}
    where the second inequality uses convexity of $f$ (Assumption~\ref{ch2:as:str_cvx_chapter_2} with $\mu = 0$), the third one uses \eqref{ch2:eq:L_smoothness_cor_2_chapter_2} with $f^{\inf} = f(x^*)$, and the last one uses $\gamma \leq 1/(2L)$. Summing up the above inequalities for $k = 0,1,\ldots, K$ and rearranging the terms, we get
    \begin{align*}
        \MoveEqLeft[1] \frac{\gamma}{K+1}\sum\limits_{k=0}^K\left(f(x^k) - f(x^*)\right) \\
        &\leq \frac{1}{K+1}\sum\limits_{k=0}^K
        \left(\|x^k - x^*\|^2 - \|x^{k+1} - x^*\|^2\right) \\
        &\quad - \frac{2\gamma}{K+1}\sum\limits_{k=0}^K
        \langle x^k - x^* - \gamma \nabla f(x^k), \theta_k \rangle\\
        &\quad+ \frac{\gamma^2}{K+1}\sum\limits_{k=0}^K \|\theta_k\|^2\\
        &= \frac{\|x^0 - x^*\|^2 - \|x^{K+1} - x^*\|^2}{K+1} \\
        &\quad - \frac{2\gamma}{K+1}\sum\limits_{k=0}^K
        \langle x^k - x^* - \gamma \nabla f(x^k), \theta_k \rangle \\
        &\quad + \frac{\gamma^2}{K+1}\sum\limits_{k=0}^K \|\theta_k\|^2.  
    \end{align*}
    Finally, we use the definition of $\overline{x}^K$ and Jensen's inequality and get the result.
\end{proof}

Using this lemma we prove the main convergence result for \algname{clipped-SGD}.
\begin{theorem}[Case 3 from Theorem~\ref{ch2:thm:clipped_SGD_main_theorem_chapter_2}]\label{ch2:thm:clipped_SGD_convex_case_appendix_chapter_2}
    Let Assumptions~\ref{ch2:as:bounded_alpha_moment_chapter_2}, \ref{ch2:as:L_smoothness_chapter_2} and \ref{ch2:as:str_cvx_chapter_2} with $\mu = 0$ hold on $Q = B_{2R}(x^*)$, where $R \geq \|x^0 - x^*\|$. Let $K > 0$ and $\beta \in (0,1]$ be such that $A = \ln(4(K+1)/\beta) \geq 1$, and
    \begin{gather}
        \gamma \leq \min\left\{\frac{1}{120L A}, \; \frac{R}{108^{1/p}\cdot 20 \sigma  (K+1)^{1/p}A^{(p-1)/p}}\right\},\label{ch2:eq:clipped_SGD_step_size_cvx_chapter_2}\\
        \lambda_{k} \equiv \lambda = \frac{R}{40 \gamma A}, \label{ch2:eq:clipped_SGD_clipping_level_cvx_chapter_2}
    \end{gather}
    Then, after $K+1$ iterations of \algname{clipped-SGD} the iterates with probability at least $1 - \beta$ satisfy
    \begin{equation}
        f(\overline{x}^K) - f(x^*) \leq \frac{2R^2}{\gamma(K+1)} \quad \text{and} \quad  \{x^k\}_{k=0}^{K+1} \subseteq B_{\sqrt{2}R}(x^*). \label{ch2:eq:clipped_SGD_convex_case_appendix_chapter_2}
    \end{equation}
    In particular, when $\gamma$ equals the minimum from \eqref{ch2:eq:clipped_SGD_step_size_cvx_chapter_2}, then the iterates produced by \algname{clipped-SGD} after $K+1$ iterations with probability at least $1-\beta$ satisfy
    \begin{equation}
        f(\overline{x}^K) - f(x^*) = \cO\left(\max\left\{\frac{LR^2 \ln(K/\beta)}{K}, \frac{\sigma R\ln^{(p-1)/p}(K/\beta)}{K^{(p-1)/p}}\right\}\right), \label{ch2:eq:clipped_SGD_convex_case_2_appendix_chapter_2}
    \end{equation}
    meaning that \algname{clipped-SGD} requires the following number of iterations/oracle calls to achieve $f(\overline{x}^K) - f(x^*) \leq \varepsilon$ with probability at least $1 - \beta$:
    \begin{equation}
        K = \cO\left(\max\left\{\frac{LR^2}{\varepsilon}\ln\frac{LR^2}{\varepsilon\beta}, \left(\frac{\sigma R}{\varepsilon}\right)^{p/(p-1)} \ln \left(\frac{1}{\beta} \left(\frac{\sigma R}{\varepsilon}\right)^{p/(p-1)}\right) \right\}\right). \label{ch2:eq:clipped_SGD_convex_case_complexity_appendix_chapter_2}
    \end{equation}
\end{theorem}
\begin{proof}
    Let $R_k = \|x^k - x^*\|$ for all $k\geq 0$. Next, our goal is to show by induction that $R_{t} \leq \sqrt{2}R$ with high probability, which allows to apply the result of Lemma~\ref{ch2:lem:main_opt_lemma_clipped_SGD_convex_chapter_2} and then use Bernstein's inequality to estimate the stochastic part of the upper-bound. More precisely, for each $k = 0,\ldots, K+1$ we consider probability event $E_k$ defined as follows: inequalities
    \begin{align}
        \MoveEqLeft[1] - 2\gamma\sum\limits_{l=0}^{t-1} \langle x^l - x^* - \gamma \nabla f(x^l), \theta_l \rangle  + \gamma^2\sum\limits_{l=0}^{t-1} \|\theta_l\|^2 \notag \\
        & \leq  R^2, \label{ch2:eq:clipped_SGD_convex_induction_inequality_1_chapter_2}\\
        R_t & \leq \sqrt{2}R \label{ch2:eq:clipped_SGD_convex_induction_inequality_2_chapter_2}
    \end{align}
    hold for all $t = 0,1,\ldots, k$ simultaneously. We want to prove via induction that $\Prob\{E_k\} \geq 1 - k\beta/(K+1)$ for all $k = 0,1,\ldots, K+1$. For $k = 0$ the statement is trivial. Assume that the statement is true for some $k = T - 1 \leq K$: $\Prob\{E_{T-1}\} \geq 1 - (T-1)\beta/(K+1)$. One needs to prove that $\Prob\{E_{T}\} \geq 1 - T\beta/(K+1)$. First, we notice that probability event $E_{T-1}$ implies that $x^t \in B_{\sqrt{2}R}(x^*)$ for all $t = 0,1, \ldots, T-1$. Moreover, $E_{T-1}$ implies
    \begin{align*}
        \|x^T - x^*\|
        &= \|x^{T-1} - x^*
        - \gamma \thesischtwotnabla f(x^{T-1};\xi^{T-1})\|\\
        &\leq \|x^{T-1} - x^*\|
        + \gamma\|\thesischtwotnabla f(x^{T-1};\xi^{T-1})\|\\
        &\leq \sqrt{2}R + \gamma \lambda
        \overset{\eqref{ch2:eq:clipped_SGD_clipping_level_cvx_chapter_2}}{\leq} 2R,
    \end{align*}
    i.e., $x^0, x^1, \ldots, x^T \in B_{2R}(x^*)$. Therefore, $E_{T-1}$ implies $\{x^k\}_{k=0}^{T} \subseteq Q$, meaning that the assumptions of Lemma~\ref{ch2:lem:main_opt_lemma_clipped_SGD_convex_chapter_2} are satisfied and we have 
    \begin{align}
        \gamma \left(f(\overline{x}^{t-1}) -f(x^*)\right)
        &\leq \frac{\|x^0 - x^*\|^2 - \|x^{t} - x^*\|^2}{t}
        \notag\\
        &\quad - \frac{2\gamma}{t}\sum\limits_{l=0}^{t-1}
        \langle x^l - x^* - \gamma \nabla f(x^l), \theta_l \rangle
        \notag\\
        &\quad + \frac{\gamma^2}{t}\sum\limits_{l=0}^{t-1} \|\theta_l\|^2.
        \label{ch2:eq:clipped_SGD_convex_technical_1_chapter_2}
    \end{align}
    for all $t = 1, \ldots, T$ simultaneously, and for all $t = 1, \ldots, T-1$ this probability event also implies that
    \begin{align}
        f(\overline{x}^{t-1}) -f(x^*)
        & \leq \frac{1}{\gamma t}\Bigg(
        R^2 - 2\gamma\sum\limits_{l=0}^{t-1}
        \langle x^l - x^* \notag \\
        &\quad - \gamma \nabla f(x^l), \theta_l \rangle \notag\\
        &\hspace{7em}
        + \gamma^2\sum\limits_{l=0}^{t-1} \|\theta_l\|^2\Bigg)
        \notag\\
        &\overset{\eqref{ch2:eq:clipped_SGD_convex_induction_inequality_1_chapter_2}}{\leq}
        \frac{2R^2}{\gamma t}.
        \label{ch2:eq:clipped_SGD_convex_technical_1_1_chapter_2}
    \end{align}
    Taking into account that $f(\overline{x}^{T-1}) -f(x^*) \geq 0$, we also derive from \eqref{ch2:eq:clipped_SGD_convex_technical_1_chapter_2} that $E_{T-1}$ implies
    \begin{equation}
        R_T^2 \leq R^2 - 2\gamma\sum\limits_{l=0}^{T-1} \langle x^l - x^* - \gamma \nabla f(x^l), \theta_l \rangle  + \gamma^2\sum\limits_{l=0}^{T-1} \|\theta_l\|^2. \label{ch2:eq:clipped_SGD_convex_technical_2_chapter_2}
    \end{equation}
    Next, we define random vectors
    \begin{equation}
        u^{t} = \begin{cases} x^t - x^* - \gamma \nabla f(x^t),& \text{if } \|x^t - x^* - \gamma \nabla f(x^t)\| \leq 2R,\\ 0,&\text{otherwise}, \end{cases} \notag
    \end{equation}
    for all $t = 0,1, \ldots, T-1$. By definition, these random vectors are bounded with probability 1
    \begin{equation}
        \|u^{t}\| \leq 2R. \label{ch2:eq:clipped_SGD_convex_technical_6_chapter_2}
    \end{equation}
    Moreover, for $t = 0, \ldots, T-1$ event $E_{T-1}$ implies
    \begin{align}
        \|\nabla f(x^{t})\| & \overset{\eqref{ch2:eq:L_smoothness_chapter_2}}{\leq} L\|x^t - x^*\|  \overset{\eqref{ch2:eq:clipped_SGD_convex_induction_inequality_2_chapter_2}}{\leq}  \sqrt{2}LR \overset{\eqref{ch2:eq:clipped_SGD_step_size_cvx_chapter_2},\eqref{ch2:eq:clipped_SGD_clipping_level_cvx_chapter_2}}{\leq} \frac{\lambda}{2},\label{ch2:eq:clipped_SGD_convex_technical_4_chapter_2}\\
        \|x^t - x^* - \gamma \nabla f(x^t)\| & \leq \|x^t - x^*\| + \gamma \|\nabla f(x^t)\| \notag \\
        & \overset{\eqref{ch2:eq:clipped_SGD_convex_technical_4_chapter_2}}{\leq} \sqrt{2}R (1+L\gamma) \overset{\eqref{ch2:eq:clipped_SGD_step_size_cvx_chapter_2}}{\leq} 2R. \notag
    \end{align}
    Next, we define the unbiased part and the bias of $\theta_{t}$ as $\theta_{t}^u$ and $\theta_{t}^b$, respectively:
    \begin{align}
\theta_{t}^u & = \thesischtwotnabla f(x^{t};\xi^t) - \Exp_{\xi^t}\left[\thesischtwotnabla f(x^{t};\xi^t)\right],\quad \theta_{t}^b \notag \\
& = \Exp_{\xi^t}\left[\thesischtwotnabla f(x^{t};\xi^t)\right] - \nabla f(x^{t}). \label{ch2:eq:clipped_SGD_convex_theta_u_b_chapter_2}
    \end{align}
    We notice that $\theta_{t} = \theta_{t}^u + \theta_{t}^b$. Using new notation, $u^t = x^t - x^* - \gamma\nabla f(x^t)$ and $\|\theta_t\|^2 \leq 2\|\theta_t^u\|^2 + 2\|\theta_t^b\|^2$, we get that $E_{T-1}$ implies
    \begin{align}
        R_T^2 
        & \leq R^2 \underbrace{-2\gamma\sum\limits_{t=0}^{T-1}\langle \theta_t^u, u^{t}\rangle}_{\circledOne}  \underbrace{-2\gamma\sum\limits_{t=0}^{T-1}\langle \theta_t^b, u^{t}\rangle}_{\circledTwo} + \underbrace{2\gamma^2\sum\limits_{t=0}^{T-1}\left(\left\|\theta_{t}^u\right\|^2 - \Exp_{\xi^t}\left[\left\|\theta_{t}^u\right\|^2\right]\right)}_{\circledThree}\notag\\
        &\quad\quad + \underbrace{2\gamma^2\sum\limits_{t=0}^{T-1}\Exp_{\xi^t}\left[\left\|\theta_{t}^u\right\|^2\right]}_{\circledFour} + \underbrace{2\gamma^2\sum\limits_{t=0}^{T-1}\left\|\theta_{t}^b\right\|^2}_{\circledFive}. \label{ch2:eq:clipped_SGD_convex_technical_7_chapter_2}
    \end{align}

    It remains to derive good enough high-probability upper-bounds for the terms $\circledOne, \circledTwo, \circledThree, \circledFour, \circledFive$, i.e., to finish our inductive proof we need to show that $\circledOne + \circledTwo + \circledThree + \circledFour + \circledFive \leq R^2$ with high probability. In the subsequent parts of the proof, we will need to use many times the bounds for the norm and second moments of $\theta_{t}^u$ and $\theta_{t}^b$. First, by definition of clipping operator, we have with probability $1$ that
     \begin{equation}
        \|\theta_{t}^u\| \leq 2\lambda. \label{ch2:eq:clipped_SGD_convex_norm_theta_u_bound_chapter_2}
    \end{equation}
    Moreover, since $E_{T-1}$ implies that $\|\nabla f(x^{t})\| \leq \lambda/2$ for $t = 0,1, \ldots, T-1$ (see \eqref{ch2:eq:clipped_SGD_convex_technical_4_chapter_2}), then, in view of Lemma~\ref{ch2:lem:bias_and_variance_clip_chapter_2}, we have that $E_{T-1}$ implies
    \begin{align}
        \|\theta_{t}^b\| & \leq \frac{2^{p}\sigma^{p}}{\lambda^{{p}-1}}, \label{ch2:eq:clipped_SGD_convex_norm_theta_b_bound_chapter_2} \\
        \Exp_{\xi^t}\left[\|\theta_{t}^u\|^2\right] & \leq 18 \lambda^{2-{p}}\sigma^{p}. \label{ch2:eq:clipped_SGD_convex_second_moment_theta_u_bound_chapter_2}
    \end{align}

    \textbf{Upper bound for $\circledOne$.} By definition of $\theta_{t}^u$, we have $\Exp_{\xi^t}[\theta_{t}^u] = 0$ and
    \begin{equation}
        \Exp_{\xi^t}\left[-2\gamma\langle\theta_t^u, u^{t}\rangle\right] = 0. \notag
    \end{equation}
    Next, sum $\circledOne$ has bounded with probability $1$ terms:
    \begin{equation}
        |2\gamma\left\la \theta_{t}^u, u^{t}\right\ra| \leq 2\gamma \|\theta_{t}^u\| \cdot \|u^{t}\| \overset{\eqref{ch2:eq:clipped_SGD_convex_technical_6_chapter_2},\eqref{ch2:eq:clipped_SGD_convex_norm_theta_u_bound_chapter_2}}{\leq} 8\gamma \lambda R\overset{\eqref{ch2:eq:clipped_SGD_clipping_level_cvx_chapter_2}}{=} \frac{R^2}{5A} \eqdef c. \label{ch2:eq:clipped_SGD_convex_technical_8_chapter_2} 
    \end{equation}
    The summands also have bounded conditional variances $\sigma_t^2$ (the conditional expectations $\Exp_{\xi^t}$ of their squares):
    \begin{align}
\sigma_t^2 & \leq \Exp_{\xi^t}\left[4\gamma^2\|\theta_{t}^u\|^2\cdot \|u^{t}\|^2\right] \notag \\
& \overset{\eqref{ch2:eq:clipped_SGD_convex_technical_6_chapter_2}}{\leq} 16\gamma^2 R^2 \Exp_{\xi^t}\left[\|\theta_{t}^u\|^2\right]. \label{ch2:eq:clipped_SGD_convex_technical_9_chapter_2}
    \end{align}
    In other words, we showed that the summands of $\circledOne$ form a bounded martingale difference sequence with bounded conditional variances $\{\sigma_t^2\}_{t=0}^{T-1}$. Next, we apply Bernstein's inequality (Lemma~\ref{ch2:lem:Bernstein_ineq_chapter_2}) with $X_t = -2\gamma \left\la \theta_{t}^u, u^{t}\right\ra$, parameter $c$ as in \eqref{ch2:eq:clipped_SGD_convex_technical_8_chapter_2}, $b = \frac{R^2}{5}$, $H = \frac{R^4}{150A}$:
    \begin{align*}
\MoveEqLeft[1] \Prob\left\{|\circledOne| > \frac{R^2}{5}\quad \text{and}\quad \sum\limits_{t=0}^{T-1} \sigma_{t}^2 \leq \frac{R^4}{150A}\right\} \\
& \leq 2\exp\left(- \frac{b^2}{2H + 2cb/3}\right) = \frac{\beta}{2(K+1)}.
    \end{align*}
    Equivalently, we have $\Prob\{ E_{\circledOne} \} \geq 1 - \beta/(2(K+1))$ for
    \begin{equation}
        E_{\circledOne} = \left\{ \text{either }  \sum\limits_{t=0}^{T-1} \sigma_{t}^2 > \frac{R^4}{150A} \text{ or } |\circledOne| \leq \frac{R^2}{5}\right\}. \label{ch2:eq:clipped_SGD_convex_sum_1_upper_bound_chapter_2}
    \end{equation}
    In addition, $E_{T-1}$ implies that
    \begin{align}
        \sum\limits_{t=0}^{T-1} \sigma_{t}^2 & \overset{\eqref{ch2:eq:clipped_SGD_convex_technical_9_chapter_2}}{\leq} 16\gamma^2 R^2 \sum\limits_{t=0}^{T-1}  \Exp_{\xi^t}\left[\|\theta_{t}^u\|^2\right] \notag \\
        & \overset{\eqref{ch2:eq:clipped_SGD_convex_second_moment_theta_u_bound_chapter_2}}{\leq} 288\gamma^2 R^2 \sigma^{{p}}T \lambda^{2-{p}} \notag\\
        & \overset{\eqref{ch2:eq:clipped_SGD_clipping_level_cvx_chapter_2}}{=} \frac{9 \cdot 40^{{p}}R^{4-{p}}\sigma^{p} T\gamma^{{p}}}{50A^{2-{p}}} \overset{\eqref{ch2:eq:clipped_SGD_step_size_cvx_chapter_2}}{\leq} \frac{R^4}{150A}. \label{ch2:eq:clipped_SGD_convex_sum_1_variance_bound_chapter_2}
    \end{align}

    \textbf{Upper bound for $\circledTwo$.} From $E_{T-1}$ it follows that
\begin{align}
        \circledTwo & = -2\gamma\sum\limits_{t=0}^{T-1}\langle \theta_t^b, u^{t} \rangle \leq 2\gamma\sum\limits_{t=0}^{T-1}\|\theta_{t}^b\|\cdot \|u^{t}\| \notag \\
        & \overset{\eqref{ch2:eq:clipped_SGD_convex_technical_6_chapter_2},\eqref{ch2:eq:clipped_SGD_convex_norm_theta_b_bound_chapter_2}}{\leq}  \frac{4\cdot 2^{p} \gamma \sigma^{p} T R}{\lambda^{{p}-1}} \notag\\
        & \overset{\eqref{ch2:eq:clipped_SGD_clipping_level_cvx_chapter_2}}{=} \frac{80^{p}}{10} \cdot \frac{\sigma^{p} T R^{2-{p}}\gamma^{p}}{A^{1-{p}}}\overset{\eqref{ch2:eq:clipped_SGD_step_size_cvx_chapter_2}}{\leq} \frac{R^2}{5}. \label{ch2:eq:clipped_SGD_convex_sum_2_upper_bound_chapter_2}
\end{align}

\textbf{Upper bound for $\circledThree$.} First, we have
    \begin{equation}
        \Exp_{\xi^t}\left[2\gamma^2\left(\left\|\theta_{t}^u\right\|^2 - \Exp_{\xi^t}\left[\left\|\theta_{t}^u\right\|^2\right]\right)\right] = 0. \notag
    \end{equation}
    Next, sum $\circledThree$ has bounded with probability $1$ terms:
    \begin{align}
        \MoveEqLeft[1] \left|2\gamma^2\left(\left\|\theta_{t}^u\right\|^2 - \Exp_{\xi^t}\left[\left\|\theta_{t}^u\right\|^2\right]\right)\right| \notag \\
        & \leq 2\gamma^2\left( \|\theta_{t}^u\|^2 +   \Exp_{\xi^t}\left[\left\|\theta_{t}^u\right\|^2\right]\right)\notag\\
        & \overset{\eqref{ch2:eq:clipped_SGD_convex_norm_theta_u_bound_chapter_2}}{\leq} 16\gamma^2\lambda^2\overset{\eqref{ch2:eq:clipped_SGD_clipping_level_cvx_chapter_2}}{=} \frac{R^2}{100A^2} \le \frac{R^2}{5A}\eqdef c. \label{ch2:eq:clipped_SGD_convex_technical_10_chapter_2}
    \end{align}
    The summands also have bounded conditional variances $\widetilde\sigma_t^2$ (the conditional expectations $\Exp_{\xi^t}$ of their squares):
    \begin{align}
        \widetilde\sigma_t^2 & \overset{\eqref{ch2:eq:clipped_SGD_convex_technical_10_chapter_2}}{\leq} \frac{R^2}{5A} \Exp_{\xi^t}\left[2\gamma^2\left|\left\|\theta_{t}^u\right\|^2 - \Exp_{\xi^t}\left[\left\|\theta_{t}^u\right\|^2\right]\right|\right] \notag \\
        & \leq \frac{4\gamma^2 R^2}{5A} \Exp_{\xi^t}\left[\|\theta_{t}^u\|^2\right], \label{ch2:eq:clipped_SGD_convex_technical_11_chapter_2}
    \end{align}
    where in \eqref{ch2:eq:clipped_SGD_convex_technical_10_chapter_2} we used $A \geq 1$. In other words, we showed that the summands of $\circledThree$ form a bounded martingale difference sequence with bounded conditional variances $\{\widetilde\sigma_t^2\}_{t=0}^{T-1}$. Next, we apply Bernstein's inequality (Lemma~\ref{ch2:lem:Bernstein_ineq_chapter_2}) with $X_t$ equal to the summands of $\circledThree$, parameter $c$ as in \eqref{ch2:eq:clipped_SGD_convex_technical_10_chapter_2}, $b = \frac{R^2}{5}$, $H = \frac{R^4}{150A}$:
    \begin{align*}
\MoveEqLeft[1] \Prob\left\{|\circledThree| > \frac{R^2}{5}\quad \text{and}\quad \sum\limits_{t=0}^{T-1} \widetilde\sigma_{t}^2 \leq \frac{R^4}{150A}\right\} \\
& \leq 2\exp\left(- \frac{b^2}{2H + 2cb/3}\right) = \frac{\beta}{2(K+1)}.
    \end{align*}
    Equivalently, we have $\Prob\{ E_{\circledThree} \} \geq 1 - \beta/(2(K+1))$ for
    \begin{equation}
        E_{\circledThree} = \left\{ \text{either }  \sum\limits_{t=0}^{T-1} \widetilde\sigma_{t}^2 > \frac{R^4}{150A} \text{ or } |\circledThree| \leq \frac{R^2}{5}\right\}. \label{ch2:eq:clipped_SGD_convex_sum_3_upper_bound_chapter_2}
    \end{equation}
    In addition, $E_{T-1}$ implies that
    \begin{align}
        \sum\limits_{t=0}^{T-1} \widetilde\sigma_{t}^2 & \overset{\eqref{ch2:eq:clipped_SGD_convex_technical_11_chapter_2}}{\leq} \frac{4\gamma^2R^2}{5A} \sum\limits_{t=0}^{T-1}  \Exp_{\xi^t}\left[\|\theta_{t}^u\|^2\right] \notag \\
        & \overset{\eqref{ch2:eq:clipped_SGD_convex_second_moment_theta_u_bound_chapter_2}}{\leq}  
        \frac{72\gamma^2R^2\lambda^{2-{p}}\sigma^{{p}}T}{5A} \notag\\
        & \overset{\eqref{ch2:eq:clipped_SGD_clipping_level_cvx_chapter_2}}{=}
        \frac{9\cdot 40^{p}}{1000}\cdot\frac{\sigma^{p} T R^{4-{p}}\gamma^{{p}}}{A^{3-{p}}} \overset{\eqref{ch2:eq:clipped_SGD_step_size_cvx_chapter_2}}{\leq} \frac{R^4}{150A}. \label{ch2:eq:clipped_SGD_convex_sum_3_variance_bound_chapter_2}
    \end{align}

    \textbf{Upper bound for $\circledFour$.} From $E_{T-1}$ it follows that
\begin{align}
        \circledFour & = 2\gamma^2\sum\limits_{t=0}^{T-1}\Exp_{\xi^t}\left[\left\|\theta_{t}^u\right\|^2\right] \overset{\eqref{ch2:eq:clipped_SGD_convex_second_moment_theta_u_bound_chapter_2}}{\leq} 36\gamma^2\lambda^{2-{p}}\sigma^{{p}}T \notag \\
        & \overset{\eqref{ch2:eq:clipped_SGD_clipping_level_cvx_chapter_2}}{=}
        \frac{9\cdot 40^{{p}}}{400}\cdot\frac{\gamma^{{p}}\sigma^{{p}}TR^{2-{p}}}{A^{2-{p}}} \overset{\eqref{ch2:eq:clipped_SGD_step_size_cvx_chapter_2}}{\leq} \frac{R^2}{5}. \label{ch2:eq:clipped_SGD_convex_sum_4_upper_bound_chapter_2}
\end{align}

\textbf{Upper bound for $\circledFive$.} From $E_{T-1}$ it follows that
\begin{align}
        \circledFive & = 2\gamma^2\sum\limits_{t=0}^{T-1}\left\|\theta_{t}^b\right\|^2 \overset{\eqref{ch2:eq:clipped_SGD_convex_norm_theta_b_bound_chapter_2}}{\leq} \frac{2 \cdot 4^{p} \sigma^{2{p}}T\gamma^2}{\lambda^{2({p}-1)}} \notag \\
        & \overset{\eqref{ch2:eq:clipped_SGD_clipping_level_cvx_chapter_2}}{=}   \frac{6400^{p}}{800}\cdot \frac{\sigma^{2{p}}T\gamma^{2{p}}R^{2-2{p}}}{A^{2(1-{p})}} \overset{\eqref{ch2:eq:clipped_SGD_step_size_cvx_chapter_2}}{\leq} \frac{R^2}{5}. \label{ch2:eq:clipped_SGD_convex_sum_5_upper_bound_chapter_2}
\end{align}

Now, we have the upper bounds for  $\circledOne, \circledTwo, \circledThree, \circledFour, \circledFive$. In particular, probability event $E_{T-1}$ implies
\begin{gather*}
        R_T^2 \overset{\eqref{ch2:eq:clipped_SGD_convex_technical_7_chapter_2}}{\leq} R^2 + \circledOne + \circledTwo + \circledThree + \circledFour + \circledFive,\\
        \circledTwo \overset{\eqref{ch2:eq:clipped_SGD_convex_sum_2_upper_bound_chapter_2}}{\leq} \frac{R^2}{5},\quad \circledFour \overset{\eqref{ch2:eq:clipped_SGD_convex_sum_4_upper_bound_chapter_2}}{\leq} \frac{R^2}{5}, \quad \circledFive \overset{\eqref{ch2:eq:clipped_SGD_convex_sum_5_upper_bound_chapter_2}}{\leq} \frac{R^2}{5},\\
        \sum\limits_{t=0}^{T-1} \sigma_t^2 \overset{\eqref{ch2:eq:clipped_SGD_convex_sum_1_variance_bound_chapter_2}}{\leq} \frac{R^4}{150A},\quad \sum\limits_{t=0}^{T-1} \widetilde\sigma_t^2 \overset{\eqref{ch2:eq:clipped_SGD_convex_sum_3_variance_bound_chapter_2}}{\leq} \frac{R^4}{150A}.
\end{gather*}
    Moreover, we also have (see \eqref{ch2:eq:clipped_SGD_convex_sum_1_upper_bound_chapter_2}, \eqref{ch2:eq:clipped_SGD_convex_sum_3_upper_bound_chapter_2} and our induction assumption)
\begin{equation*}
        \Prob\{E_{T-1}\} \geq 1 - \frac{(T-1)\beta}{K+1},\quad \Prob\{E_{\circledOne}\} \geq 1 - \frac{\beta}{2(K+1)},\quad \Prob\{E_{\circledThree}\} \geq 1 - \frac{\beta}{2(K+1)},
\end{equation*}
    where 
\begin{align*}
        E_{\circledOne} & = \left\{ \text{either} \quad  \sum\limits_{t=0}^{T-1} \sigma_{t}^2 > \frac{R^4}{150A} \quad \text{or}\quad |\circledOne| \leq \frac{R^2}{5}\right\},\\
        E_{\circledThree} & = \left\{ \text{either} \quad  \sum\limits_{t=0}^{T-1} \widetilde\sigma_{t}^2 > \frac{R^4}{150A} \quad \text{or}\quad |\circledThree| \leq \frac{R^2}{5}\right\}.
\end{align*}
    Thus, probability event $E_{T-1} \cap E_{\circledOne} \cap E_{\circledThree}$ implies
\begin{align}
        R_T^2 & \leq R^2 + \frac{R^2}{5} + \frac{R^2}{5} + \frac{R^2}{5} + \frac{R^2}{5} + \frac{R^2}{5} = 2R^2, \notag
\end{align}
    which is equivalent to \eqref{ch2:eq:clipped_SGD_convex_induction_inequality_1_chapter_2} and \eqref{ch2:eq:clipped_SGD_convex_induction_inequality_2_chapter_2} for $t = T$, and 
\begin{align*}
        \Prob\{E_T\}
        &\geq \Prob\left\{E_{T-1} \cap E_{\circledOne} \cap E_{\circledThree}\right\} \\
        &= 1 - \Prob\left\{\overline{E}_{T-1} \cup \overline{E}_{\circledOne} \cup \overline{E}_{\circledThree}\right\} \\
        &\geq 1 - \Prob\{\overline{E}_{T-1}\} - \Prob\{\overline{E}_{\circledOne}\} - \Prob\{\overline{E}_{\circledThree}\} \\
        &\geq 1 - \frac{T\beta}{K+1}.
\end{align*}
    This finishes the inductive part of our proof, i.e., for all $k = 0,1,\ldots,K+1$ we have $\Prob\{E_k\} \geq 1 - k\beta/(K+1)$. In particular, for $k = K+1$ we have that with probability at least $1 - \beta$
\begin{equation*}
        f(\overline{x}^K) - f(x^*) \overset{\eqref{ch2:eq:clipped_SGD_convex_technical_1_1_chapter_2}}{\leq}\frac{2R^2}{\gamma(K+1)}
\end{equation*}
    and $\{x^k\}_{k=0}^{K+1} \subseteq B_{\sqrt{2}R}(x^*)$, which follows from \eqref{ch2:eq:clipped_SGD_convex_induction_inequality_2_chapter_2}.
    
Finally, if
\begin{equation*}
        \gamma = \min\left\{\frac{1}{120L A}, \; \frac{R}{108^{1/p}\cdot 20 \sigma  (K+1)^{1/p}A^{(p-1)/p}}\right\},
\end{equation*}
then with probability at least $1-\beta$
\begin{align*}
        f(\overline{x}^K) - f(x^*) & \leq \frac{2R^2}{\gamma(K+1)}\\
        & = \max\left\{\frac{240 LR^2A}{K+1}, \frac{40\cdot 108^{1/p}\sigma R(K+1)^{1/p}A^{(p-1)/p}}{K+1}\right\}\notag\\
        & = \cO\left(\max\left\{\frac{L R^2 \ln(K/\beta)}{K}, \frac{ \sigma R\ln^{(p-1)/p}(K/\beta)}{K^{(p-1)/p}}\right\}\right).
\end{align*}
To get $f(\overline{x}^K) - f(x^*) \leq \varepsilon$ with probability at least $1-\beta$ it is sufficient to choose $K$ such that both terms in the maximum above are $\cO(\varepsilon)$. This leads to
\begin{equation*}
         K = \cO\left(\max\left\{\frac{LR^2}{\varepsilon} \ln\frac{L R^2}{\varepsilon\beta}, \left(\frac{\sigma R}{\varepsilon}\right)^{p/(p-1)} \ln \left(\frac{1}{\beta} \left(\frac{\sigma R}{\varepsilon}\right)^{p/(p-1)}\right) \right\}\right),
\end{equation*}
which concludes the proof.

\end{proof}

\subsection{Quasi-Strongly Convex Functions}

Finally, we consider \algname{clipped-SGD} under smoothness and quasi-strong convexity assumptions. As the next lemma shows, the gradient of such function is quasi-strongly monotone and star-cocoercive operator.

\begin{lemma}
    Consider differentiable function $f: \R^d \to \R$. If $f$ satisfies Assumption~\ref{ch2:as:QSC_chapter_2} on some set $Q$ with parameter $\mu$, then operator $F(x) = \nabla f(x)$ satisfies Assumption~\ref{ch2:as:QSM_chapter_2} on $Q$ with parameter $\mu/2$. If $f$ satisfies Assumptions~\ref{ch2:as:L_smoothness_chapter_2} and \ref{ch2:as:QSC_chapter_2} with $\mu = 0$ on some set $Q$, then operator $F(x) = \nabla f(x)$ satisfies Assumption~\ref{ch2:as:star-cocoercivity_chapter_2} on $Q$ with $\ell = 2L$.
\end{lemma}
\begin{proof}
    We start with the first part. Assumption~\ref{ch2:as:QSC_chapter_2} on set $Q$ means that for any $x \in Q$
    \begin{equation*}
        f(x^*) \geq f(x) + \langle \nabla f(x), x^* - x \rangle + \frac{\mu}{2}\|x- x^*\|^2.
    \end{equation*}
    For $F(x) = \nabla f(x)$ it implies that for all $x \in Q$
    \begin{equation*}
        \langle F(x), x - x^* \rangle \geq f(x) - f(x^*) + \frac{\mu}{2}\|x - x^*\|^2 \geq \frac{\mu}{2}\|x - x^*\|^2,
    \end{equation*}
    i.e., Assumption~\ref{ch2:as:QSM_chapter_2} holds on $Q$ with parameter $\mu/2$ for operator $F(x)$.

    Next, we prove the second part. Assume that $f$ satisfies Assumptions~\ref{ch2:as:L_smoothness_chapter_2} and \ref{ch2:as:QSC_chapter_2} with $\mu = 0$ on some set $Q$. Our goal is to show that $F(x) = \nabla f(x)$ satisfies Assumption~\ref{ch2:as:star-cocoercivity_chapter_2} on $Q$. In view of \textcite[Lemma C.6]{gorbunov2022extragradient}, this is equivalent to showing that operator $\thesischtwoId - F/L$ is non-expansive around $x^*$, i.e., we need to show that $\|(\thesischtwoId - F/L)(x) - (\thesischtwoId - F/L)(x^*)\| \leq \|x-x^*\|$ for any $x \in Q$. We have
    \begin{align*}
        \MoveEqLeft \left\|\left(\thesischtwoId - \frac{1}{L}F\right)(x)
        - \left(\thesischtwoId - \frac{1}{L}F\right)(x^*)\right\|^2\\
        &= \left\|x-x^* - \frac{1}{L}F(x)\right\|^2\\
        &=\|x - x^*\|^2
        - \frac{2}{L}\langle x - x^*, F(x) \rangle
        + \frac{1}{L^2}\|F(x)\|^2\\
        &= \|x - x^*\|^2
        - \frac{2}{L}\langle x - x^*, \nabla f(x) \rangle
        + \frac{1}{L^2}\|\nabla f(x)\|^2\\
        &\leq
        \|x - x^*\|^2
        - \frac{2}{L}\left(f(x) - f(x^*)\right)\\
        &\quad + \frac{2}{L}\left(f(x) - f(x^*)\right)\\
        &= \|x-x^*\|^2.
    \end{align*}
    where the inequality uses \eqref{ch2:eq:QSC_chapter_2} with $\mu = 0$, $F(x^*) = 0$ and \eqref{ch2:eq:L_smoothness_cor_2_chapter_2} with $f^{\inf} = f(x^*)$. This finishes the proof.
\end{proof}

Therefore, using the result of Theorem~\ref{ch2:thm:main_result_str_mon_SGDA_chapter_2} with $\ell := 2L$ and $\mu := \mu/2$, we get the convergence result for \algname{clipped-SGD} under smoothness and quasi-strong convexity assumptions.

\begin{theorem}[Case 4 in Theorem~\ref{ch2:thm:clipped_SGD_main_theorem_chapter_2}]\label{ch2:thm:main_result_QSC_SGD_chapter_2}
    Let Assumptions~\ref{ch2:as:bounded_alpha_moment_chapter_2}, \ref{ch2:as:L_smoothness_chapter_2}, \ref{ch2:as:QSC_chapter_2} hold for $Q = B_{2R}(x^*) = \{x\in\R^d\mid \|x - x^*\| \leq 2R\}$, where $R \geq \|x^0 - x^*\|$. Let $K \geq 0$ and $\beta \in (0,1]$ be such that $A = \ln(4(K+1)/\beta) \geq 1$, and
    \begin{align}
        0< \gamma & \leq \min\left\{\frac{1}{800 LA}, \frac{2\ln(B_K)}{\mu(K+1)}\right\}, \label{ch2:eq:gamma_SGD_QSC_chapter_2}\\
        B_K & = \max\left\{2, \frac{(K+1)^{2(p-1)/p}\mu^2R^2}{4\cdot54000^{2/p}\sigma^2A^{2(p-1)/p}\ln^2(B_K)} \right\} \label{ch2:eq:B_K_SGD_QSC_1_chapter_2} \\
        & = \cO\biggl(\max\biggl\{2, \frac{K^{2(p-1)/p}\mu^2R^2}{\sigma^2\ln^{2(p-1)/p}(K/\beta)}\notag\\
        &\qquad\qquad \cdot \frac{1}{\ln^2\left(\max\left\{2, \frac{K^{2(p-1)/p}\mu^2R^2}{\sigma^2\ln^{2(p-1)/p}(K/\beta)} \right\}\right)} \biggr\}\biggr), \label{ch2:eq:B_K_SGD_QSC_2_chapter_2} \\
        \lambda_k & = \frac{\exp(-\gamma(\mu/2)(1 + k/2))R}{120\gamma A}, \label{ch2:eq:lambda_SGD_QSC_chapter_2}
    \end{align}
    where $K$ is large enough that the maximum in \eqref{ch2:eq:B_K_SGD_QSC_1_chapter_2} is attained at its second term. Then, after $K+1$ iterations the iterates produced by \algname{clipped-SGD} with probability at least $1 - \beta$ satisfy 
    \begin{equation}
        \|x^{K+1} - x^*\|^2 \leq 2\exp(-\gamma(\mu/2)(K+1))R^2. \label{ch2:eq:main_result_str_cvx_SGD_chapter_2}
    \end{equation}
    In particular, when $\gamma$ equals the minimum from \eqref{ch2:eq:gamma_SGD_QSC_chapter_2}, then the iterates produced by \algname{clipped-SGD} after $K+1$ iterations with probability at least $1-\beta$ satisfy
    \begin{align}
        \MoveEqLeft \|x^{K+1} - x^*\|^2\notag\\
        &= \cO\Biggl(\max\Biggl\{
        R^2\exp\left(- \frac{\mu K}{L \ln(K/\beta)}\right),
        \frac{\sigma^2\ln^{2(p-1)/p}(K/\beta)}{K^{2(p-1)/p}\mu^2}\notag\\
        &\qquad\qquad\cdot\ln^2\left(\max\left\{2,
        \frac{K^{2(p-1)/p}\mu^2R^2}
        {\sigma^2\ln^{2(p-1)/p}(K/\beta)}
        \right\}\right)
        \Biggr\}\Biggr),
        \label{ch2:eq:clipped_SGD_QSC_case_2_appendix_chapter_2}
    \end{align}
    meaning that \algname{clipped-SGD} requires the following number of iterations/oracle calls to achieve $\|x^{K+1} - x^*\|^2 \leq \varepsilon$ with probability at least $1 - \beta$:
    \begin{align}
        K &= \cO\biggl(\max\biggl\{\frac{L}{\mu}\ln\left(\frac{R^2}{\varepsilon}\right)\ln\left(\frac{L}{\mu \beta}\ln\frac{R^2}{\varepsilon}\right),\notag\\
        &\quad \left(\frac{\sigma^2}{\mu^2\varepsilon}\right)^{p/(2(p-1))}\ln \left(\frac{1}{\beta} \left(\frac{\sigma^2}{\mu^2\varepsilon}\right)^{p/(2(p-1))}\right)\ln^{p/(p-1)}\left(B_\varepsilon\right)\biggr\}\biggr), \label{ch2:eq:clipped_SGD_QSC_case_complexity_appendix_chapter_2}
    \end{align}
    where
    \begin{equation*}
        B_\varepsilon = \max\left\{2, \frac{R^2}{\varepsilon \ln \left(\frac{1}{\beta} \left(\frac{4\sigma^2}{\mu^2\varepsilon}\right)^{p/(2(p-1))}\right)}\right\}.
    \end{equation*}
\end{theorem}

\clearpage

\section{Additional Proofs for \algname{clipped-SSTM} and \algname{R-clipped-SSTM}}

In this section, we provide the complete formulation of the main results for \algname{clipped-SSTM} and \algname{R-clipped-SSTM} and the corresponding proofs. For brevity, we will use the following notation: $\thesischtwotnabla f(x^{k+1};\xi^k) = \clip\left(\nabla f(x^{k+1};\xi^k), \lambda_k\right)$.

\begin{algorithm}[h]
\caption{Clipped Stochastic Similar Triangles Method (\algname{clipped-SSTM}) \parencite{gorbunov2020stochastic}}
\label{ch2:alg:clipped-SSTM_chapter_2}   
\begin{algorithmic}[1]
\REQUIRE starting point $x^0$, number of iterations $K$, stepsize parameter $a>0$, clipping levels $\{\lambda_k\}_{k=0}^{K-1}$, smoothness constant $L$.
\STATE Set $A_0 = a_0 = 0$, $y^0 = z^0 = x^0$
\FOR{$k=0,\ldots, K-1$}
\STATE Set $a_{k+1} = \frac{k+2}{2aL}$, $A_{k+1} = A_k + a_{k+1}$
\STATE $x^{k+1} = \frac{A_k y^k + a_{k+1} z^k}{A_{k+1}}$
\STATE Compute $\thesischtwotnabla f(x^{k+1};\xi^k) = \clip\left(\nabla f(x^{k+1};\xi^k), \lambda_k\right)$ using a fresh sample $\xi^k \sim \thesischtwocD_k$
\STATE $z^{k+1} = z^k - a_{k+1}\thesischtwotnabla f(x^{k+1};\xi^k)$
\STATE $y^{k+1} = \frac{A_k y^k + a_{k+1} z^{k+1}}{A_{k+1}}$
\ENDFOR
\ENSURE $y^K$ 
\end{algorithmic}
\end{algorithm}

\subsection{Convex Functions}\label{ch2:appendix:clipped_SSTM_chapter_2}

We start with the following lemma, which is a special case of Lemma~6 from \textcite{gorbunov2021near}. This result can be seen as the ``optimization'' part of the analysis of \algname{clipped-SSTM}: the proof follows the same steps as the analysis of deterministic Similar Triangles Method \parencite{gasnikov2016universal, dvurechenskii2018decentralize} and separates stochasticity from the deterministic part of the method. 

\begin{lemma}[Special case of Lemma~4.1 from \textcite{gorbunov2021near}]\label{ch2:lem:main_opt_lemma_clipped_SSTM_chapter_2}
    Let Assumptions~\ref{ch2:as:L_smoothness_chapter_2} and \ref{ch2:as:str_cvx_chapter_2} with $\mu = 0$ hold on $Q = B_{3R}(x^*)$, where $R \ge \|x^0 - x^*\|$, and let stepsize parameter $a$ satisfy $a\ge 1$. If $x^{k}, y^k, z^k \in B_{3R}(x^*)$ for all $k = 0,1, \ldots, N$, $N \ge 0$, then after $N$ iterations of \algname{clipped-SSTM} for all $z\in B_{3R}(x^*)$ we have
    \begin{align}
        A_N\left(f(y^N) - f(z)\right) & \le \frac{1}{2}\|z^0 - z\|^2 - \frac{1}{2}\|z^{N} - z\|^2 \notag \\
        &\quad + \sum\limits_{k=0}^{N-1}a_{k+1}\left\la \theta_{k+1}, z - z^{k} + a_{k+1} \nabla f(x^{k+1})\right\ra \notag\\
        &\quad\quad + \sum\limits_{k=0}^{N-1}a_{k+1}^2\left\|\theta_{k+1}\right\|^2,\label{ch2:eq:main_opt_lemma_clipped_SSTM_chapter_2}\\
        \theta_{k+1} & \eqdef \thesischtwotnabla f(x^{k+1};\xi^k) - \nabla f(x^{k+1}).\label{ch2:eq:theta_k+1_def_clipped_SSTM_chapter_2}
    \end{align}
\end{lemma}
\begin{proof}
    For completeness, we provide the full proof. Using the update rule \eqref{ch2:eq:clipped_SSTM_z_update_chapter_2} for $z^{k+1}$, we get that for all $z\in B_{3R}(x^*)$ and $k = 0,1, \ldots, N-1$
    \begin{align}
        \MoveEqLeft[1] a_{k+1}\left\la\thesischtwotnabla f(x^{k+1};\xi^k), z^k - z \right\ra \notag \\
        &= a_{k+1}\left\la \thesischtwotnabla f(x^{k+1};\xi^k),
        z^k - z^{k+1}\right\ra
        \notag\\
        &\quad + a_{k+1}\left\la \thesischtwotnabla f(x^{k+1};\xi^k),
        z^{k+1} - z\right\ra
        \notag\\
        &= a_{k+1}\left\la \thesischtwotnabla f(x^{k+1};\xi^k),
        z^k - z^{k+1}\right\ra
        \notag\\
        &\quad + \left\la z^{k+1}-z^k, z - z^{k+1}\right\ra
        \notag\\
        &= a_{k+1}\left\la \thesischtwotnabla f(x^{k+1};\xi^k),
        z^k - z^{k+1}\right\ra
        - \frac{1}{2}\|z^k - z^{k+1}\|^2
        \notag\\
        &\quad + \frac{1}{2}\|z^k - z\|^2
        - \frac{1}{2}\|z^{k+1}-z\|^2.
        \label{ch2:eq:main_olc_clipped_SSTM_technical_2_chapter_2}
    \end{align}
    where in the last step we apply $2\langle a, b \rangle = \|a + b\|^2 - \|a\|^2 - \|b\|^2$ with $a = z^{k+1} - z^k$ and $b = z - z^{k+1}$. The update rules \eqref{ch2:eq:clipped_SSTM_y_update_chapter_2} and \eqref{ch2:eq:clipped_SSTM_x_update_chapter_2} give the following formula:
    \begin{align}
         y^{k+1}
         &= \frac{A_k y^k + {a}_{k+1}z^{k+1}}{A_{k+1}}
         \notag\\
         &= \frac{A_k y^k + {a}_{k+1}z^{k}}{A_{k+1}}
         + \frac{a_{k+1}}{A_{k+1}}\left(z^{k+1}-z^k\right)
         \notag\\
         &= x^{k+1}
         + \frac{a_{k+1}}{A_{k+1}}\left(z^{k+1}-z^k\right).
         \label{ch2:eq:main_olc_clipped_SSTM_technical_3_chapter_2}
    \end{align}
    It implies
    \begin{align*}
        \MoveEqLeft[1] a_{k+1}\left\la\thesischtwotnabla f(x^{k+1};\xi^k), z^k - z \right\ra \\
        &\overset{\eqref{ch2:eq:theta_k+1_def_clipped_SSTM_chapter_2},\eqref{ch2:eq:main_olc_clipped_SSTM_technical_2_chapter_2}}{\le}
        a_{k+1}\left\la \nabla f(x^{k+1}), z^{k} - z^{k+1}\right\ra
        \notag\\
        &\quad - \frac{1}{2}\|z^k - z^{k+1}\|^2
        \notag\\
        &\quad + a_{k+1}\left\la \theta_{k+1}, z^{k} - z^{k+1}\right\ra
        \notag\\
        &\quad + \frac{1}{2}\|z^k - z\|^2
        - \frac{1}{2}\|z^{k+1}-z\|^2
        \notag\\
        &\overset{\eqref{ch2:eq:main_olc_clipped_SSTM_technical_3_chapter_2}}{=}
        A_{k+1}\left\la \nabla f(x^{k+1}), x^{k+1} - y^{k+1}\right\ra
        \notag\\
        &\quad - \frac{1}{2}\|z^k - z^{k+1}\|^2
        \notag\\
        &\quad + a_{k+1}\left\la \theta_{k+1}, z^{k} - z^{k+1}\right\ra
        \notag\\
        &\quad + \frac{1}{2}\|z^k - z\|^2
        - \frac{1}{2}\|z^{k+1}-z\|^2
        \notag\\
        &\overset{\eqref{ch2:eq:L_smoothness_cor_1_chapter_2}}{\le}
        A_{k+1}\left(f(x^{k+1}) - f(y^{k+1})\right)
        \notag\\
        &\quad + \frac{A_{k+1}L}{2}\|x^{k+1}-y^{k+1}\|^2
        - \frac{1}{2}\|z^k - z^{k+1}\|^2
        \notag\\
        &\quad + a_{k+1}\left\la \theta_{k+1}, z^{k} - z^{k+1}\right\ra
        \notag\\
        &\quad + \frac{1}{2}\|z^k - z\|^2
        - \frac{1}{2}\|z^{k+1}-z\|^2
        \notag\\
        &\overset{\eqref{ch2:eq:main_olc_clipped_SSTM_technical_3_chapter_2}}{=}
        A_{k+1}\left(f(x^{k+1}) - f(y^{k+1})\right)
        \notag\\
        &\quad + \frac{1}{2}\left(\frac{a_{k+1}^2L}{A_{k+1}} - 1\right)
        \|z^k - z^{k+1}\|^2
        \notag\\
        &\quad + a_{k+1}\left\la \theta_{k+1}, z^{k} - z^{k+1}\right\ra
        \notag\\
        &\quad + \frac{1}{2}\|z^k - z\|^2
        - \frac{1}{2}\|z^{k+1}-z\|^2,
        \notag
    \end{align*}
    where in the third inequality we use $x^{k+1}, y^{k+1} \in B_{3R}(x^*)$. Since $A_{k+1} \ge aLa_{k+1}^2$ (Lemma~\ref{ch2:lem:alpha_k_A_K_lemma_chapter_2}) and $a\ge 1$ we can continue our derivation as follows:
    \begin{align}
        \MoveEqLeft[1] a_{k+1}\left\la\thesischtwotnabla f(x^{k+1};\xi^k), z^k - z \right\ra \notag \\
        & \le A_{k+1}\left(f(x^{k+1}) - f(y^{k+1})\right) + a_{k+1}\left\la \theta_{k+1}, z^{k} - z^{k+1}\right\ra\notag\\
        &\quad\quad + \frac{1}{2}\|z^k - z\|^2 - \frac{1}{2}\|z^{k+1}-z\|^2. \label{ch2:eq:main_olc_clipped_SSTM_technical_4_chapter_2}
    \end{align}
    Convexity of $f$ gives
    \begin{align}
        \MoveEqLeft[1] \left\la\thesischtwotnabla f(x^{k+1};\xi^k), y^k - x^{k+1}\right\ra \notag \\
        & \overset{\eqref{ch2:eq:theta_k+1_def_clipped_SSTM_chapter_2}}{=} \left\la\nabla f(x^{k+1}), y^k - x^{k+1}\right\ra + \left\la \theta_{k+1}, y^k - x^{k+1}\right\ra\notag\\
        & \le f(y^k) - f(x^{k+1}) + \left\la \theta_{k+1}, y^k - x^{k+1}\right\ra.\label{ch2:eq:main_olc_clipped_SSTM_technical_5_chapter_2}
    \end{align}
    The definition of $x^{k+1}$ \eqref{ch2:eq:clipped_SSTM_x_update_chapter_2} implies
    \begin{equation}
        a_{k+1}\left(x^{k+1} - z^k\right) = A_k\left(y^k - x^{k+1}\right)\label{ch2:eq:main_olc_clipped_SSTM_technical_6_chapter_2}
    \end{equation}
    since $A_{k+1} = A_k + a_{k+1}$. Putting all inequalities together, we derive that
    \begin{align*}
        \MoveEqLeft[1] a_{k+1}\left\la\thesischtwotnabla f(x^{k+1};\xi^k), x^{k+1} - z \right\ra \\
        &= a_{k+1}\left\la\thesischtwotnabla f(x^{k+1};\xi^k),
        x^{k+1} - z^k \right\ra \\
        &\quad + a_{k+1}\left\la\thesischtwotnabla f(x^{k+1};\xi^k),
        z^k - z \right\ra\\
        &\overset{\eqref{ch2:eq:main_olc_clipped_SSTM_technical_6_chapter_2}}{=}
        A_{k}\left\la\thesischtwotnabla f(x^{k+1};\xi^k),
        y^{k} - x^{k+1} \right\ra \\
        &\quad + a_{k+1}\left\la\thesischtwotnabla f(x^{k+1};\xi^k),
        z^k - z \right\ra\\
        &\overset{\eqref{ch2:eq:main_olc_clipped_SSTM_technical_5_chapter_2},\eqref{ch2:eq:main_olc_clipped_SSTM_technical_4_chapter_2}}{\le}
        A_k\left(f(y^k)-f(x^{k+1})\right) \\
        &\quad + A_k\left\la \theta_{k+1}, y^k - x^{k+1}\right\ra\\
        &\quad + A_{k+1}\left(f(x^{k+1}) - f(y^{k+1})\right)\\
        &\quad + a_{k+1}\left\la \theta_{k+1}, z^{k} - z^{k+1}\right\ra \\
        &\quad + \frac{1}{2}\|z^k - z\|^2
        - \frac{1}{2}\|z^{k+1}-z\|^2\\
        &\overset{\eqref{ch2:eq:main_olc_clipped_SSTM_technical_6_chapter_2}}{=}
        A_kf(y^k) - A_{k+1}f(y^{k+1}) \\
        &\quad + a_{k+1}\left\la \theta_{k+1}, x^{k+1}-z^k\right\ra\\
        &\quad + {a}_{k+1}f(x^{k+1})
        + a_{k+1}\left\la \theta_{k+1}, z^{k} - z^{k+1}\right\ra\\
        &\quad + \frac{1}{2}\|z^k - z\|^2
        - \frac{1}{2}\|z^{k+1}-z\|^2\\
        &= A_kf(y^k) - A_{k+1}f(y^{k+1}) + {a}_{k+1}f(x^{k+1})\\
        &\quad + a_{k+1}\left\la \theta_{k+1}, x^{k+1} - z^{k+1}\right\ra\\
        &\quad + \frac{1}{2}\|z^k - z\|^2
        - \frac{1}{2}\|z^{k+1}-z\|^2.
    \end{align*}
    Rearranging the terms, we get
    \begin{align*}
        \MoveEqLeft A_{k+1}f(y^{k+1}) - A_kf(y^k)\\
        &\le a_{k+1}\left(f(x^{k+1}) + \left\la\thesischtwotnabla f(x^{k+1};\xi^k), z-x^{k+1}\right\ra\right)\\
        &\quad + \frac{1}{2}\|z^k - z\|^2 - \frac{1}{2}\|z^{k+1} - z\|^2
        + a_{k+1}\left\la \theta_{k+1}, x^{k+1} - z^{k+1}\right\ra\\
        & \overset{\eqref{ch2:eq:theta_k+1_def_clipped_SSTM_chapter_2}}{=}
        a_{k+1}\left(f(x^{k+1}) + \left\la\nabla f(x^{k+1}), z-x^{k+1}\right\ra\right) \\
        &\quad + a_{k+1}\left\la\theta_{k+1}, z-x^{k+1}\right\ra\\
        &\quad + \frac{1}{2}\|z^k - z\|^2 - \frac{1}{2}\|z^{k+1} - z\|^2
        + a_{k+1}\left\la \theta_{k+1}, x^{k+1} - z^{k+1}\right\ra\\
        & \le {a}_{k+1}f(z) + \frac{1}{2}\|z^k - z\|^2 - \frac{1}{2}\|z^{k+1} - z\|^2 \\
        &\quad + a_{k+1}\left\la \theta_{k+1}, z - z^{k+1}\right\ra,
    \end{align*}
    where in the last inequality we use the convexity of $f$. Taking into account $A_0 = a_0 = 0$ and $A_{N} = \sum_{k=0}^{N-1}a_{k+1}$ we sum up these inequalities for $k = 0, \ldots, N-1$ and get
    \begin{align*}
        A_Nf(y^N)
        &\le A_N f(z) + \frac{1}{2}\|z^0 - z\|^2
        - \frac{1}{2}\|z^{N} - z\|^2 \\
        &\quad + \sum\limits_{k=0}^{N-1}a_{k+1}
        \left\la \theta_{k+1}, z - z^{k+1}\right\ra\\
        &= A_N f(z) + \frac{1}{2}\|z^0 - z\|^2
        - \frac{1}{2}\|z^{N} - z\|^2 \\
        &\quad + \sum\limits_{k=0}^{N-1}a_{k+1}
        \left\la \theta_{k+1}, z - z^{k}
        + a_{k+1} \thesischtwotnabla f(x^{k+1};\xi^k)\right\ra \\
        &\overset{\eqref{ch2:eq:theta_k+1_def_clipped_SSTM_chapter_2}}{=}
        A_N f(z) + \frac{1}{2}\|z^0 - z\|^2
        - \frac{1}{2}\|z^{N} - z\|^2 \\
        &\quad + \sum\limits_{k=0}^{N-1}a_{k+1}
        \left\la \theta_{k+1}, z - z^{k}
        + a_{k+1}\nabla f(x^{k+1})\right\ra\\
        &\quad + \sum\limits_{k=0}^{N-1}a_{k+1}^2
        \left\|\theta_{k+1}\right\|^2,
    \end{align*}
    which concludes the proof.
\end{proof}

Using this lemma we prove the main convergence result for \algname{clipped-SSTM}.
\begin{theorem}[Full version of Theorem~\ref{ch2:thm:clipped_SSTM_main_theorem_chapter_2}]\label{ch2:thm:clipped_SSTM_convex_case_appendix_chapter_2}
    Let Assumptions~\ref{ch2:as:bounded_alpha_moment_chapter_2}, \ref{ch2:as:L_smoothness_chapter_2} and \ref{ch2:as:str_cvx_chapter_2} with $\mu = 0$ hold on $Q = B_{3R}(x^*)$, where $R \geq \|x^0 - x^*\|$, and
    \begin{gather}
        a \geq \max\left\{48600\ln^2(4K/\beta), \frac{900\sigma(K+1)K^{1/p}\ln^{(p-1)/p}(4K/\beta)}{LR}\right\},\label{ch2:eq:clipped_SSTM_parameter_a_chapter_2}\\
        \lambda_{k} = \frac{R}{30a_{k+1}\ln(4K/\beta)}, \label{ch2:eq:clipped_SSTM_clipping_level_chapter_2}
    \end{gather}
    for some $K > 0$ and $\beta \in (0,1]$ such that $\ln(4K/\beta) \geq 1$. Then, after $K$ iterations of \algname{clipped-SSTM} the iterates with probability at least $1 - \beta$ satisfy
    \begin{align}
        & f(y^K) - f(x^*) \leq \frac{6aL R^2}{K(K+3)} \quad \text{and} \notag\\
        & \{x^k\}_{k=0}^{K+1}, \{z^k\}_{k=0}^{K}, \{y^k\}_{k=0}^K \subseteq B_{2R}(x^*). \label{ch2:eq:clipped_SSTM_convex_case_appendix_chapter_2}
    \end{align}
    In particular, when parameter $a$ equals the maximum from \eqref{ch2:eq:clipped_SSTM_parameter_a_chapter_2}, then the iterates produced by \algname{clipped-SSTM} after $K$ iterations with probability at least $1-\beta$ satisfy
    \begin{equation}
        f(y^K) - f(x^*) = \cO\left(\max\left\{\frac{LR^2\ln^2(K/\beta)}{K^2}, \frac{\sigma R \ln^{(p-1)/p}(K/\beta)}{K^{(p-1)/p}}\right\}\right), \label{ch2:eq:clipped_SSTM_convex_case_2_appendix_chapter_2}
    \end{equation}
    meaning that \algname{clipped-SSTM} requires the following number of iterations/oracle calls to achieve $f(y^K) - f(x^*) \leq \varepsilon$ with probability at least $1 - \beta$:
    \begin{equation}
        K = \cO\left(\max\left\{\sqrt{\frac{LR^2}{\varepsilon}}\ln\frac{LR^2}{\varepsilon\beta}, \left(\frac{\sigma R}{\varepsilon}\right)^{p/(p-1)}\ln\left( \frac{1}{\beta} \left(\frac{\sigma R}{\varepsilon}\right)^{p/(p-1)}\right)\right\}\right). \label{ch2:eq:clipped_SSTM_convex_case_complexity_appendix_chapter_2}
    \end{equation}
\end{theorem}
\begin{proof}
    The proof starts similarly to the proof of Theorem~4.1 from \textcite{gorbunov2021near}. Let $R_k = \|z^k - x^*\|$, $\widetilde{R}_0 = R_0$, $\widetilde{R}_{k+1} = \max\{\widetilde{R}_k, R_{k+1}\}$ for all $k\geq 0$. We first show by induction that for all $k\geq 0$ the iterates $x^{k+1}, z^k, y^k$ lie in $B_{\widetilde{R}_k}(x^*)$. The induction base is trivial since $y^0 = z^0$, $\widetilde{R}_0 = R_0$, and $x^1 = \tfrac{A_0 y^0 + a_1 z^0}{A_1} = z^0$. Next, assume that $x^{l}, z^{l-1}, y^{l-1} \in B_{\widetilde{R}_{l-1}}(x^*)$ for some $l\ge 1$. By definitions of $R_l$ and $\widetilde{R}_l$ we have that $z^l \in B_{R_{l}}(x^*)\subseteq B_{\widetilde{R}_{l}}(x^*)$. Since $y^l$ is a convex combination of $y^{l-1}\in B_{\widetilde{R}_{l-1}}(x^*)\subseteq B_{\widetilde{R}_{l}}(x^*)$, $z^l\in B_{\widetilde{R}_{l}}(x^*)$ and $B_{\widetilde{R}_{l}}(x^*)$ is a convex set we conclude that $y^l \in B_{\widetilde{R}_{l}}(x^*)$. Finally, since $x^{l+1}$ is a convex combination of $y^l$ and $z^l$ we have that $x^{l+1}$ lies in $B_{\widetilde{R}_{l}}(x^*)$ as well.
    
    Next, our goal is to show by induction that $\widetilde{R}_{l} \leq 3R$ with high probability, which allows us to apply the result of Lemma~\ref{ch2:lem:main_opt_lemma_clipped_SSTM_chapter_2} and then use Bernstein's inequality to estimate the stochastic part of the upper-bound. More precisely, for each $k = 0,\ldots, K$ we consider probability event $E_k$ defined as follows: inequalities
    \begin{gather}
        \sum\limits_{l=0}^{t-1}a_{l+1}\left\la \theta_{l+1}, x^* - z^{l} + a_{l+1}\nabla f(x^{l+1})\right\ra + \sum\limits_{l=0}^{t-1}a_{l+1}^2\left\|\theta_{l+1}\right\|^2 \leq R^2, \label{ch2:eq:clipped_SSTM_induction_inequality_1_chapter_2}\\
        R_t \leq 2R \label{ch2:eq:clipped_SSTM_induction_inequality_2_chapter_2}
    \end{gather}
    hold for all $t = 0,1,\ldots, k$ simultaneously. We want to prove via induction that $\Prob\{E_k\} \geq 1 - k\beta/K$ for all $k = 0,1,\ldots, K$. For $k = 0$ the statement is trivial: the left-hand side of \eqref{ch2:eq:clipped_SSTM_induction_inequality_1_chapter_2} equals zero and $R \geq R_0$ by definition. Assume that the statement is true for some $k = T - 1 \leq K-1$: $\Prob\{E_{T-1}\} \geq 1 - (T-1)\beta/K$. One needs to prove that $\Prob\{E_{T}\} \geq 1 - T\beta/K$. First, we notice that probability event $E_{T-1}$ implies that $\widetilde{R}_t \leq 2R$ for all $t = 0,1, \ldots, T-1$. Moreover, it implies that
    \begin{align*}
\|z^{T} - x^*\| & \overset{\eqref{ch2:eq:clipped_SSTM_z_update_chapter_2}}{\leq} \|z^{T-1} - x^*\| + a_{T}\|\thesischtwotnabla f(x^T;\xi^{T-1})\| \\
& \leq 2R + a_{T}\lambda_{T-1} \overset{\eqref{ch2:eq:clipped_SSTM_clipping_level_chapter_2}}{\leq} 3R.
    \end{align*}
    Therefore, $E_{T-1}$ implies $\{x^k\}_{k=0}^{T}, \{z^k\}_{k=0}^{T}, \{y^k\}_{k=0}^T \subseteq B_{3R}(x^*)$, meaning that the assumptions of Lemma~\ref{ch2:lem:main_opt_lemma_clipped_SSTM_chapter_2} are satisfied and we have
    \begin{align}
        A_t\left(f(y^t) - f(x^*)\right)
        &\le \frac{1}{2}R_0^2 - \frac{1}{2}R_t^2 \notag\\
        &\quad + \sum\limits_{l=0}^{t-1}a_{l+1}
        \left\la \theta_{l+1},
        x^* - z^{l} + a_{l+1} \nabla f(x^{l+1})\right\ra \notag\\
        &\quad + \sum\limits_{l=0}^{t-1}a_{l+1}^2
        \left\|\theta_{l+1}\right\|^2
        \label{ch2:eq:clipped_SSTM_technical_1_chapter_2}
    \end{align}
    for all $t = 0,1, \ldots, T$ simultaneously, and for all $t = 1, \ldots, T-1$ this probability event also implies that
    \begin{equation}
        f(y^t) - f(x^*) \overset{\eqref{ch2:eq:clipped_SSTM_induction_inequality_1_chapter_2}, \eqref{ch2:eq:clipped_SSTM_technical_1_chapter_2}}{\leq} \frac{\frac{1}{2}R_0^2 - \frac{1}{2}R_t^2 + R^2}{A_t} \leq \frac{3R^2}{2A_t} = \frac{6aLR^2}{t(t+3)}. \label{ch2:eq:clipped_SSTM_technical_1_1_chapter_2}
    \end{equation}
    Taking into account that $f(y^T) - f(x^*) \geq 0$, we also derive that $E_{T-1}$ implies
    \begin{align}
        R_T^2 & \leq R_0^2 + \underbrace{2\sum\limits_{t=0}^{T-1}a_{t+1}\left\la \theta_{t+1}, x^* - z^{t} + a_{t+1} \nabla f(x^{t+1})\right\ra + 2\sum\limits_{t=0}^{T-1}a_{t+1}^2\left\|\theta_{t+1}\right\|^2}_{2B_T}\notag\\
        & \leq R^2 + 2B_T. \label{ch2:eq:clipped_SSTM_technical_2_chapter_2}
    \end{align}
    Before we estimate $B_T$, we need to derive a few useful inequalities. We start with showing that $E_{T-1}$ implies $\|\nabla f(x^{t+1})\| \leq \lambda_t/2$ for all $t = 0,1, \ldots, T-1$. For $t = 0$ we have $x^1 = x^0$ and
    \begin{align}
\|\nabla f(x^{1})\| & = \|\nabla f(x^0)\| \overset{\eqref{ch2:eq:L_smoothness_chapter_2}}{\leq} L\|x^0 - x^*\| \leq \frac{R}{aa_1} \notag \\
& = \frac{\lambda_0}{2}\cdot \frac{60\ln(4K/\beta)}{a} \overset{\eqref{ch2:eq:clipped_SSTM_parameter_a_chapter_2}}{\leq} \frac{\lambda_0}{2}. \label{ch2:eq:clipped_SSTM_technical_3_chapter_2}
    \end{align}
    Next, for $t = 1, \ldots, T-1$ event $E_{T-1}$ implies
    \begin{align}
        \|\nabla f(x^{t+1})\| & \leq \|\nabla f(x^{t+1}) - \nabla f(y^t)\| + \|\nabla f(y^t)\| \notag\\
        & \overset{\eqref{ch2:eq:L_smoothness_chapter_2}, \eqref{ch2:eq:L_smoothness_cor_2_chapter_2}}{\leq} L\|x^{t+1} - y^t\| + \sqrt{2L\left(f(y^t) - f(x^*)\right)} \notag\\
        & \overset{\eqref{ch2:eq:main_olc_clipped_SSTM_technical_6_chapter_2},\eqref{ch2:eq:clipped_SSTM_technical_1_1_chapter_2}}{\leq} \frac{L{a}_{t+1}}{A_t}\|x^{t+1} - z^t\| + \sqrt{\frac{12aL^2R^2}{t(t+3)}}\notag\\
        & \leq \frac{4LR{a}_{t+1}}{A_t} + \sqrt{\frac{12aL^2R^2}{t(t+3)}} \notag\\
        & =  \frac{R}{60a_{t+1}\ln(4K/\beta)}\Biggl(\frac{240L{a}_{t+1}^2\ln(4K/\beta)}{A_t}\notag\\
        &\qquad\qquad + 60\sqrt{\frac{12aL^2a_{t+1}^2\ln^2(4K/\beta)}{t(t+3)}}\Biggr)\notag\\
        & \overset{\eqref{ch2:eq:A_k+1_explicit_chapter_2},\eqref{ch2:eq:clipped_SSTM_clipping_level_chapter_2}}{\leq} \frac{\lambda_t}{2}\Biggl(\frac{240L\left(\frac{t+2}{2aL}\right)^2\ln(4K/\beta)}{\frac{t(t+3)}{4aL}}\notag\\
        &\qquad\qquad + 60\sqrt{\frac{12aL^2\left(\frac{t+2}{2aL}\right)^2\ln^2(4K/\beta)}{t(t+3)}}\Biggr) \notag\\
        & = \frac{\lambda_t}{2} \left(\frac{240(t+2)^2\ln(4K/\beta)}{t(t+3)a} + 60\sqrt{\frac{3(t+2)^2\ln^2(4K/\beta)}{t(t+3)a}}\right) \notag\\
        & \leq \frac{\lambda_t}{2} \left(\frac{540\ln(4K/\beta)}{a} + \frac{90\sqrt{3}\ln(4K/\beta)}{\sqrt{a}}\right) \overset{\eqref{ch2:eq:clipped_SSTM_parameter_a_chapter_2}}{\leq} \frac{\lambda_t}{2}, \label{ch2:eq:clipped_SSTM_technical_4_chapter_2}
    \end{align}
    where in the last row we use $\frac{(t+2)^2}{t(t+3)} \leq \frac{9}{4}$ for all $t \geq 1$. Therefore, probability event $E_{T-1}$ implies that
    \begin{align}
        \MoveEqLeft[1] \|x^* - z^{t} + a_{t+1} \nabla f(x^{t+1})\| \notag \\
        &\leq \|x^* - z^{t}\| + a_{t+1}\|\nabla f(x^{t+1})\|
        \notag\\
        &\overset{\eqref{ch2:eq:clipped_SSTM_induction_inequality_2_chapter_2},
        \eqref{ch2:eq:clipped_SSTM_technical_3_chapter_2},
        \eqref{ch2:eq:clipped_SSTM_technical_4_chapter_2}}{\leq}
        2R + \frac{R}{60\ln(4K/\beta)}
        \leq 3R.
        \label{ch2:eq:clipped_SSTM_technical_5_chapter_2}
    \end{align}
    for all $t = 0,1, \ldots, T-1$. Next, we define random vectors
    \begin{equation}
        u^{t} = \begin{cases} x^* - z^{t} + a_{t+1} \nabla f(x^{t+1}),& \text{if } \|x^* - z^{t} + a_{t+1} \nabla f(x^{t+1})\| \leq 3R,\\ 0,&\text{otherwise}, \end{cases} \notag
    \end{equation}
    for all $t = 0,1, \ldots, T-1$. By definition these random vectors are bounded with probability 1
    \begin{equation}
        \|u^{t}\| \leq 3R \label{ch2:eq:clipped_SSTM_technical_6_chapter_2}
    \end{equation}
    and probability event $E_{T-1}$ implies that $u^{t} = x^* - z^{t} + a_{t+1} \nabla f(x^{t+1})$ for all $t = 0,1, \ldots, T-1$. Then, from $E_{T-1}$ it follows that
    \begin{align}
        B_T & = \sum\limits_{t=0}^{T-1}a_{t+1}\left\la \theta_{t+1}, u^{t}\right\ra + \sum\limits_{t=0}^{T-1}a_{t+1}^2\left\|\theta_{t+1}\right\|^2. \notag
    \end{align}
    Next, we define the unbiased part and the bias of $\theta_{t+1}$ as $\theta_{t+1}^u$ and $\theta_{t+1}^b$ respectively:
    \begin{equation}
        \begin{aligned}
            \theta_{t+1}^u &= \thesischtwotnabla f(x^{t+1};\xi^t) - \Exp_{\xi^t}\left[\thesischtwotnabla f(x^{t+1};\xi^t)\right],\\
            \theta_{t+1}^b &= \Exp_{\xi^t}\left[\thesischtwotnabla f(x^{t+1};\xi^t)\right] - \nabla f(x^{t+1}).
        \end{aligned}\label{ch2:eq:clipped_SSTM_theta_u_b_chapter_2}
    \end{equation}
    We notice that $\theta_{t+1} = \theta_{t+1}^u + \theta_{t+1}^b$. Using new notation and $\|\theta_{t+1}\|^2 \leq 2\|\theta_{t+1}^u\|^2 + 2\|\theta_{t+1}^b\|^2$, we get that $E_{T-1}$ implies
    \begin{align}
        B_T & = \sum\limits_{t=0}^{T-1}a_{t+1}\left\la \theta_{t+1}^u + \theta_{t+1}^b, u^{t}\right\ra + \sum\limits_{t=0}^{T-1}a_{t+1}^2\left\|\theta_{t+1}^u + \theta_{t+1}^b\right\|^2 \notag\\
        & \leq \underbrace{\sum\limits_{t=0}^{T-1}a_{t+1}\left\la \theta_{t+1}^u, u^{t}\right\ra}_{\circledOne} + \underbrace{\sum\limits_{t=0}^{T-1}a_{t+1}\left\la \theta_{t+1}^b, u^{t}\right\ra}_{\circledTwo}\notag\\
        &\quad + \underbrace{2\sum\limits_{t=0}^{T-1}a_{t+1}^2\left(\left\|\theta_{t+1}^u\right\|^2 - \Exp_{\xi^t}\left[\left\|\theta_{t+1}^u\right\|^2\right]\right)}_{\circledThree}\notag\\
        &\quad\quad + \underbrace{2\sum\limits_{t=0}^{T-1}a_{t+1}^2\Exp_{\xi^t}\left[\left\|\theta_{t+1}^u\right\|^2\right]}_{\circledFour} + \underbrace{2\sum\limits_{t=0}^{T-1}a_{t+1}^2\left\|\theta_{t+1}^b\right\|^2}_{\circledFive}. \label{ch2:eq:clipped_SSTM_technical_7_chapter_2}
    \end{align}
    It remains to derive good enough high-probability upper-bounds for the terms $\circledOne, \circledTwo, \circledThree, \circledFour, \circledFive$, i.e., to finish our inductive proof we need to show that $\circledOne + \circledTwo + \circledThree + \circledFour + \circledFive \leq R^2$ with high probability. In the subsequent parts of the proof, we will need to use many times the bounds for the norm and second moments of $\theta_{t+1}^u$ and $\theta_{t+1}^b$. First, by definition of clipping operator, we have with probability $1$ that
    \begin{equation}
        \|\theta_{t+1}^u\| \leq 2\lambda_{t}. \label{ch2:eq:clipped_SSTM_norm_theta_u_bound_chapter_2}
    \end{equation}
    Moreover, since $E_{T-1}$ implies that $\|\nabla f(x^{t+1})\| \leq \lambda_t/2$ for $t = 0,1, \ldots, T-1$ (see \eqref{ch2:eq:clipped_SSTM_technical_3_chapter_2} and \eqref{ch2:eq:clipped_SSTM_technical_4_chapter_2}), then, in view of Lemma~\ref{ch2:lem:bias_and_variance_clip_chapter_2}, we have that $E_{T-1}$ implies
    \begin{align}
        \|\theta_{t+1}^b\| & \leq \frac{2^{p}\sigma^{p}}{\lambda_{t}^{{p}-1}}, \label{ch2:eq:clipped_SSTM_norm_theta_b_bound_chapter_2} \\
        \Exp_{\xi^t}\left[\|\theta_{t+1}^u\|^2\right] & \leq 18 \lambda_t^{2-{p}}\sigma^{p}. \label{ch2:eq:clipped_SSTM_second_moment_theta_u_bound_chapter_2}
    \end{align}

    \textbf{Upper bound for $\circledOne$.} By definition of $\theta_{t+1}^u$, we have $\Exp_{\xi^t}[\theta_{t+1}^u] = 0$ and
    \begin{equation}
        \Exp_{\xi^t}\left[a_{t+1}\left\la \theta_{t+1}^u, u^{t}\right\ra\right] = 0. \notag
    \end{equation}
    Next, sum $\circledOne$ has bounded with probability $1$ terms:
    \begin{align}
|a_{t+1}\left\la \theta_{t+1}^u, u^{t}\right\ra| & \leq a_{t+1} \|\theta_{t+1}^u\| \cdot \|u^{t}\| \notag \\
& \overset{\eqref{ch2:eq:clipped_SSTM_technical_6_chapter_2},\eqref{ch2:eq:clipped_SSTM_norm_theta_u_bound_chapter_2}}{\leq} 6a_{t+1}\lambda_t R \overset{\eqref{ch2:eq:clipped_SSTM_clipping_level_chapter_2}}{=} \frac{R^2}{5\ln(4K/\beta)} \eqdef c. \label{ch2:eq:clipped_SSTM_technical_8_chapter_2}
    \end{align}
    The summands also have bounded conditional variances $\sigma_t^2$ (the conditional expectations $\Exp_{\xi^t}$ of their squares):
    \begin{align}
\sigma_t^2 & \leq \Exp_{\xi^t}\left[a_{t+1}^2\|\theta_{t+1}^u\|^2\cdot \|u^{t}\|^2\right] \notag \\
& \overset{\eqref{ch2:eq:clipped_SSTM_technical_6_chapter_2}}{\leq} 9a_{t+1}^2R^2 \Exp_{\xi^t}\left[\|\theta_{t+1}^u\|^2\right]. \label{ch2:eq:clipped_SSTM_technical_9_chapter_2}
    \end{align}
    In other words, we showed that the summands of $\circledOne$ form a bounded martingale difference sequence with bounded conditional variances $\{\sigma_t^2\}_{t=0}^{T-1}$. Next, we apply Bernstein's inequality (Lemma~\ref{ch2:lem:Bernstein_ineq_chapter_2}) with $X_t = a_{t+1}\left\la \theta_{t+1}^u, u^{t}\right\ra$, parameter $c$ as in \eqref{ch2:eq:clipped_SSTM_technical_8_chapter_2}, $b = \frac{R^2}{5}$, $H = \frac{R^4}{150\ln(4K/\beta)}$:
    \begin{align*}
\MoveEqLeft[1] \Prob\left\{|\circledOne| > \frac{R^2}{5}\quad \text{and}\quad \sum\limits_{t=0}^{T-1} \sigma_{t}^2 \leq \frac{R^4}{150\ln(4K/\beta)}\right\} \\
& \leq 2\exp\left(- \frac{b^2}{2H + 2cb/3}\right) = \frac{\beta}{2K}.
    \end{align*}
    Equivalently, we have
    \begin{align}
\Prob\left\{ E_{\circledOne} \right\} & \geq 1 - \frac{\beta}{2K},\quad \text{for}\quad E_{\circledOne} \notag \\
& = \left\{ \text{either} \quad  \sum\limits_{t=0}^{T-1} \sigma_{t}^2 > \frac{R^4}{150\ln(4K/\beta)} \quad \text{or}\quad |\circledOne| \leq \frac{R^2}{5}\right\}. \label{ch2:eq:clipped_SSTM_sum_1_upper_bound_chapter_2}
    \end{align}
    In addition, $E_{T-1}$ implies that
    \begin{align}
        \sum\limits_{t=0}^{T-1} \sigma_{t}^2 & \overset{\eqref{ch2:eq:clipped_SSTM_technical_9_chapter_2}}{\leq} 9R^2 \sum\limits_{t=0}^{T-1} a_{t+1}^2 \Exp_{\xi^t}\left[\|\theta_{t+1}^u\|^2\right] \notag \\
        & \overset{\eqref{ch2:eq:clipped_SSTM_second_moment_theta_u_bound_chapter_2}}{\leq} 162\sigma^{p} R^2 \sum\limits_{t=0}^{T-1} a_{t+1}^2\lambda_t^{2-{p}} \notag\\
        & \overset{\eqref{ch2:eq:clipped_SSTM_clipping_level_chapter_2}}{=} \frac{162\sigma^{p} R^{4-{p}}}{30^{2-{p}}\ln^{2-{p}}(4K/\beta)}\sum\limits_{t=0}^{T-1} a_{t+1}^{p} \notag \\
        & = \frac{162\sigma^{p} R^{4-{p}}}{30^{2-{p}}\cdot 2^{p} a^{p} L^{p}\ln^{2-{p}}(4K/\beta)}\sum\limits_{t=0}^{T-1} (t+2)^{p} \notag\\
        & \leq \frac{1}{a^{p}} \cdot \frac{162\sigma^{p} R^{4-{p}} T(T+1)^{p}}{30^{2-{p}}\cdot 2^{p} L^{p}\ln^{2-{p}}(4K/\beta)} \notag \\
        & \overset{\eqref{ch2:eq:clipped_SSTM_parameter_a_chapter_2}}{\leq} \frac{R^4}{150 \ln(4K/\beta)}. \label{ch2:eq:clipped_SSTM_sum_1_variance_bound_chapter_2}
    \end{align}

    \textbf{Upper bound for $\circledTwo$.} From $E_{T-1}$ it follows that
    \begin{align}
        \circledTwo
        &\leq \sum\limits_{t=0}^{T-1}a_{t+1}\|\theta_{t+1}^b\|\cdot \|u^{t}\|
        \notag\\
        &\overset{\eqref{ch2:eq:clipped_SSTM_technical_6_chapter_2},\eqref{ch2:eq:clipped_SSTM_norm_theta_b_bound_chapter_2}}{\leq}
        3R\cdot 2^{p}\sigma^{p} \sum\limits_{t=0}^{T-1}\frac{a_{t+1}}{\lambda_{t}^{{p}-1}}
        \notag\\
        &\overset{\eqref{ch2:eq:clipped_SSTM_clipping_level_chapter_2}}{\leq}
        12R\sigma^{p} \cdot
        \frac{30^{{p}-1}\ln^{{p}-1}(4K/\beta)}{R^{{p}-1}}
        \sum\limits_{t=0}^{T-1}a_{t+1}^{{p}}
        \notag\\
        &\leq
        \frac{360\sigma^{p} R^{2 - {p}}\ln^{{p}-1}(4K/\beta)}
        {2^{p} a^{{p}} L^{p}}
        \sum\limits_{t=0}^{T-1} (t+2)^{p}
        \notag\\
        &\leq \frac{1}{a^{p}} \cdot
        \frac{180\sigma^{p} R^{2-{p}} T(T+1)^{p}\ln^{{p}-1}(4K/\beta)}{L^{p}}
        \overset{\eqref{ch2:eq:clipped_SSTM_parameter_a_chapter_2}}{\leq}
        \frac{R^2}{5}.
        \label{ch2:eq:clipped_SSTM_sum_2_upper_bound_chapter_2}
    \end{align}

    \textbf{Upper bound for $\circledThree$.} First, we have
    \begin{equation}
        \Exp_{\xi^t}\left[2a_{t+1}^2\left(\left\|\theta_{t+1}^u\right\|^2 - \Exp_{\xi^t}\left[\left\|\theta_{t+1}^u\right\|^2\right]\right)\right] = 0. \notag
    \end{equation}
    Next, sum $\circledThree$ has bounded with probability $1$ terms:
    \begin{align}
        &\left|2a_{t+1}^2\left(\left\|\theta_{t+1}^u\right\|^2
        - \Exp_{\xi^t}\left[\left\|\theta_{t+1}^u\right\|^2\right]\right)\right|
        \notag\\
        &\quad \leq 2a_{t+1}^2\left(
        \|\theta_{t+1}^u\|^2
        + \Exp_{\xi^t}\left[\left\|\theta_{t+1}^u\right\|^2\right]\right)
        \notag\\
        &\quad \overset{\eqref{ch2:eq:clipped_SSTM_norm_theta_u_bound_chapter_2}}{\leq}
        16a_{t+1}^2\lambda_t^2
        \overset{\eqref{ch2:eq:clipped_SSTM_clipping_level_chapter_2}}{\leq}
        \frac{R^2}{5\ln(4K/\beta)}
        \eqdef c.
        \label{ch2:eq:clipped_SSTM_technical_10_chapter_2}
    \end{align}
    The summands also have bounded conditional variances
    \begin{equation*}
        \widetilde\sigma_t^2 \eqdef
        \Exp_{\xi^t}\left[
        4a_{t+1}^4\left(
        \left\|\theta_{t+1}^u\right\|^2
        - \Exp_{\xi^t}\left[\left\|\theta_{t+1}^u\right\|^2\right]
        \right)^2\right]:
    \end{equation*}
    \begin{align}
        \widetilde\sigma_t^2
        &\overset{\eqref{ch2:eq:clipped_SSTM_technical_10_chapter_2}}{\leq}
        \frac{R^2}{5 \ln(4K/\beta)}
        \Exp_{\xi^t}\left[
        2a_{t+1}^2\left|
        \left\|\theta_{t+1}^u\right\|^2
        - \Exp_{\xi^t}\left[\left\|\theta_{t+1}^u\right\|^2\right]
        \right|\right]
        \notag\\
        &\leq a_{t+1}^2R^2
        \Exp_{\xi^t}\left[\|\theta_{t+1}^u\|^2\right].
        \label{ch2:eq:clipped_SSTM_technical_11_chapter_2}
    \end{align}
    since $\ln(4K/\beta) \geq 1$. In other words, we showed that the summands of $\circledThree$ form a bounded martingale difference sequence with bounded conditional variances $\{\widetilde\sigma_t^2\}_{t=0}^{T-1}$. Next, we apply Bernstein's inequality (Lemma~\ref{ch2:lem:Bernstein_ineq_chapter_2}) with $X_t$ equal to the summands of $\circledThree$, parameter $c$ as in \eqref{ch2:eq:clipped_SSTM_technical_10_chapter_2}, $b = \frac{R^2}{5}$, $H = \frac{R^4}{150\ln(4K/\beta)}$:
    \begin{align*}
\MoveEqLeft[1] \Prob\left\{|\circledThree| > \frac{R^2}{5}\quad \text{and}\quad \sum\limits_{t=0}^{T-1} \widetilde\sigma_{t}^2 \leq \frac{R^4}{150\ln(4K/\beta)}\right\} \\
& \leq 2\exp\left(- \frac{b^2}{2H + 2cb/3}\right) = \frac{\beta}{2K}.
    \end{align*}
    Equivalently, we have
    \begin{align}
\Prob\left\{ E_{\circledThree} \right\} & \geq 1 - \frac{\beta}{2K},\quad \text{for}\quad E_{\circledThree} \notag \\
& = \left\{ \text{either} \quad  \sum\limits_{t=0}^{T-1} \widetilde\sigma_{t}^2 > \frac{R^4}{150\ln(4K/\beta)} \quad \text{or}\quad |\circledThree| \leq \frac{R^2}{5}\right\}. \label{ch2:eq:clipped_SSTM_sum_3_upper_bound_chapter_2}
    \end{align}
    In addition, $E_{T-1}$ implies that
    \begin{align}
        \sum\limits_{t=0}^{T-1} \widetilde\sigma_{t}^2
        &\overset{\eqref{ch2:eq:clipped_SSTM_technical_11_chapter_2}}{\leq}
        R^2 \sum\limits_{t=0}^{T-1} a_{t+1}^2
        \Exp_{\xi^t}\left[\|\theta_{t+1}^u\|^2\right]
        \notag\\
        &\leq 9R^2 \sum\limits_{t=0}^{T-1} a_{t+1}^2
        \Exp_{\xi^t}\left[\|\theta_{t+1}^u\|^2\right]
        \notag\\
        &\overset{\eqref{ch2:eq:clipped_SSTM_sum_1_variance_bound_chapter_2}}{\leq}
        \frac{R^4}{150 \ln(4K/\beta)}.
        \label{ch2:eq:clipped_SSTM_sum_3_variance_bound_chapter_2}
    \end{align}

    \textbf{Upper bound for $\circledFour$.} From $E_{T-1}$ it follows that
    \begin{align}
        \circledFour
        &= 2\sum\limits_{t=0}^{T-1}a_{t+1}^2
        \Exp_{\xi^t}\left[\left\|\theta_{t+1}^u\right\|^2\right]
        \notag\\
        &\leq \frac{1}{R^2}\cdot 9R^2
        \sum\limits_{t=0}^{T-1}a_{t+1}^2
        \Exp_{\xi^t}\left[\left\|\theta_{t+1}^u\right\|^2\right]
        \notag\\
        &\overset{\eqref{ch2:eq:clipped_SSTM_sum_1_variance_bound_chapter_2}}{\leq}
        \frac{R^2}{150\ln(4K/\beta)}
        \leq \frac{R^2}{5}.
        \label{ch2:eq:clipped_SSTM_sum_4_upper_bound_chapter_2}
    \end{align}

    \textbf{Upper bound for $\circledFive$.} From $E_{T-1}$ it follows that
    \begin{align}
        \circledFive
        &= 2\sum\limits_{t=0}^{T-1}a_{t+1}^2
        \left\|\theta_{t+1}^b\right\|^2
        \notag\\
        &\leq 2^{2{p} + 1}\sigma^{2{p}}
        \sum\limits_{t=0}^{T-1}
        \frac{a_{t+1}^2}{\lambda_{t}^{2{p} - 2}}
        \notag\\
        &\overset{\eqref{ch2:eq:clipped_SSTM_clipping_level_chapter_2}}{=}
        \frac{2^{2{p} + 1}\cdot 30^{2{p}-2}\sigma^{2{p}}
        \ln^{2{p}-2}(4K/\beta)}{R^{2{p}-2}}
        \sum\limits_{t=0}^{T-1}a_{t+1}^{2{p}}
        \notag\\
        &=
        \frac{2^{2{p} + 1}\cdot 30^{2{p}-2}\sigma^{2{p}}
        \ln^{2{p}-2}(4K/\beta)}
        {2^{2{p}}a^{2{p}} L^{2{p}} R^{2{p}-2}}
        \sum\limits_{t=0}^{T-1}(t+2)^{2{p}}
        \notag\\
        &\leq \frac{1}{a^{2{p}}} \cdot
        \frac{1800 \sigma^{2{p}}T(T+1)^{2{p}}
        \ln^{2{p}-2}(4K/\beta)}
        {L^{2{p}}R^{2{p}-2}}
        \notag\\
        &\overset{\eqref{ch2:eq:clipped_SSTM_parameter_a_chapter_2}}{\leq}
        \frac{R^2}{5}.
        \label{ch2:eq:clipped_SSTM_sum_5_upper_bound_chapter_2}
    \end{align}

    Now, we have the upper bounds for  $\circledOne, \circledTwo, \circledThree, \circledFour, \circledFive$. In particular, probability event $E_{T-1}$ implies
    \begin{gather*}
        B_T \overset{\eqref{ch2:eq:clipped_SSTM_technical_7_chapter_2}}{\leq} \circledOne + \circledTwo + \circledThree + \circledFour + \circledFive,\\
        \circledTwo \overset{\eqref{ch2:eq:clipped_SSTM_sum_2_upper_bound_chapter_2}}{\leq} \frac{R^2}{5},\quad \circledFour \overset{\eqref{ch2:eq:clipped_SSTM_sum_4_upper_bound_chapter_2}}{\leq} \frac{R^2}{5}, \quad \circledFive \overset{\eqref{ch2:eq:clipped_SSTM_sum_5_upper_bound_chapter_2}}{\leq} \frac{R^2}{5},\\
        \sum\limits_{t=0}^{T-1} \sigma_t^2 \overset{\eqref{ch2:eq:clipped_SSTM_sum_1_variance_bound_chapter_2}}{\leq} \frac{R^4}{150 \ln(4K/\beta)},\quad \sum\limits_{t=0}^{T-1} \widetilde\sigma_t^2 \overset{\eqref{ch2:eq:clipped_SSTM_sum_3_variance_bound_chapter_2}}{\leq} \frac{R^4}{150 \ln(4K/\beta)}.
    \end{gather*}
    Moreover, we also have (see \eqref{ch2:eq:clipped_SSTM_sum_1_upper_bound_chapter_2}, \eqref{ch2:eq:clipped_SSTM_sum_3_upper_bound_chapter_2} and our induction assumption)
    \begin{equation*}
        \Prob\{E_{T-1}\} \geq 1 - \frac{(T-1)\beta}{K},\quad \Prob\{E_{\circledOne}\} \geq 1 - \frac{\beta}{2K},\quad \Prob\{E_{\circledThree}\} \geq 1 - \frac{\beta}{2K},
    \end{equation*}
    where 
    \begin{align*}
        E_{\circledOne} & = \left\{ \text{either} \quad  \sum\limits_{t=0}^{T-1} \sigma_{t}^2 > \frac{R^4}{150\ln(4K/\beta)} \quad \text{or}\quad |\circledOne| \leq \frac{R^2}{5}\right\},\\
        E_{\circledThree} & = \left\{ \text{either} \quad  \sum\limits_{t=0}^{T-1} \widetilde\sigma_{t}^2 > \frac{R^4}{150\ln(4K/\beta)} \quad \text{or}\quad |\circledThree| \leq \frac{R^2}{5}\right\}.
    \end{align*}
    Thus, probability event $E_{T-1} \cap E_{\circledOne} \cap E_{\circledThree}$ implies
    \begin{align}
        B_T & \leq \frac{R^2}{5} + \frac{R^2}{5} + \frac{R^2}{5} + \frac{R^2}{5} + \frac{R^2}{5} = R^2, \notag\\
        R_T^2 & \overset{\eqref{ch2:eq:clipped_SSTM_technical_2_chapter_2}}{\leq} R^2 + 2R^2 \leq (2R)^2, \notag
    \end{align}
    which is equivalent to \eqref{ch2:eq:clipped_SSTM_induction_inequality_1_chapter_2} and \eqref{ch2:eq:clipped_SSTM_induction_inequality_2_chapter_2} for $t = T$, and 
    \begin{align*}
        \Prob\{E_T\}
        &\geq \Prob\left\{E_{T-1}
        \cap E_{\circledOne}\cap E_{\circledThree}\right\}\\
        &= 1 - \Prob\left\{\overline{E}_{T-1}
        \cup \overline{E}_{\circledOne}\cup \overline{E}_{\circledThree}\right\}\\
        &\geq 1 - \Prob\{\overline{E}_{T-1}\}
        - \Prob\{\overline{E}_{\circledOne}\}
        - \Prob\{\overline{E}_{\circledThree}\}\\
        &\geq 1 - \frac{T\beta}{K}.
    \end{align*}
    This finishes the inductive part of our proof, i.e., for all $k = 0,1,\ldots,K$ we have $\Prob\{E_k\} \geq 1 - k\beta/K$. In particular, for $k = K$ we have that with probability at least $1 - \beta$
    \begin{equation*}
        f(y^K) - f(x^*) \overset{\eqref{ch2:eq:clipped_SSTM_technical_1_1_chapter_2}}{\leq} \frac{6aLR^2}{K(K+3)}
    \end{equation*}
    and $\{x^k\}_{k=0}^{K+1}, \{z^k\}_{k=0}^{K}, \{y^k\}_{k=0}^K \subseteq B_{2R}(x^*)$, which follows from \eqref{ch2:eq:clipped_SSTM_induction_inequality_2_chapter_2}.
    
    Finally, if $a$ equals the maximum in \eqref{ch2:eq:clipped_SSTM_parameter_a_chapter_2},
    then with probability at least $1-\beta$
    \begin{align*}
        \MoveEqLeft f(y^K) - f(x^*) \leq \frac{6aLR^2}{K(K+3)} \\
        & = \max\left\{\frac{291600LR^2\ln^2(4K/\beta)}{K(K+3)}, \frac{5400\sigma R(K+1)K^{1/p}\ln^{(p-1)/p}(4K/\beta)}{K(K+3)}\right\}\\
        & = \cO\left(\max\left\{\frac{LR^2\ln^2(K/\beta)}{K^2}, \frac{\sigma R \ln^{(p-1)/p}(K/\beta)}{K^{(p-1)/p}}\right\}\right).
    \end{align*}
    To get $f(y^K) - f(x^*) \leq \varepsilon$ with probability at least $1-\beta$ it is sufficient to choose $K$ such that both terms in the maximum above are $\cO(\varepsilon)$. This leads to
    \begin{equation*}
         K = \cO\left(\max\left\{\sqrt{\frac{LR^2}{\varepsilon}}\ln\frac{LR^2}{\varepsilon\beta}, \left(\frac{\sigma R}{\varepsilon}\right)^{p/(p-1)}\ln\left( \frac{1}{\beta} \left(\frac{\sigma R}{\varepsilon}\right)^{p/(p-1)}\right)\right\}\right)
    \end{equation*}
    that concludes the proof.
\end{proof}

\subsection{Strongly Convex Functions}\label{ch2:appendix:R_clipped_SSTM_chapter_2}
In the strongly convex case, we consider the restarted version of \algname{clipped-SSTM} (\algname{R-clipped-SSTM}). The main result is summarized below.

\begin{algorithm}[h]
\caption{Restarted \algname{clipped-SSTM} (\algname{R-clipped-SSTM}) \parencite{gorbunov2020stochastic}}
\label{ch2:alg:R-clipped-SSTM_chapter_2}   
\begin{algorithmic}[1]
\REQUIRE starting point $x^0$, number of restarts $\tau$, number of steps of \algname{clipped-SSTM} between restarts $\{K_t\}_{t=1}^{\tau}$, stepsize parameters $\{a_t\}_{t=1}^\tau$, clipping levels $\{\lambda_{k}^1\}_{k=0}^{K_1-1}$, $\{\lambda_{k}^2\}_{k=0}^{K_2-1}$, \ldots, $\{\lambda_{k}^\tau\}_{k=0}^{K_\tau - 1}$, smoothness constant $L$.
\STATE $\hat{x}^0 = x^0$
\FOR{$t=1,\ldots, \tau$}
\STATE Run \algname{clipped-SSTM} (Algorithm~\ref{ch2:alg:clipped-SSTM_chapter_2}) for $K_t$ iterations with stepsize parameter $a_{t}$, clipping levels $\{\lambda_{k}^t\}_{k=0}^{K_t-1}$, and starting point $\hat{x}^{t-1}$. Define the output of \algname{clipped-SSTM} by $\hat{x}^{t}$.
\ENDFOR
\ENSURE $\hat{x}^\tau$ 
\end{algorithmic}
\end{algorithm}

\begin{theorem}[Full version of Theorem~\ref{ch2:thm:R_clipped_SSTM_main_theorem_chapter_2}]\label{ch2:thm:R_clipped_SSTM_main_theorem_appendix_chapter_2}
    Let Assumptions~\ref{ch2:as:bounded_alpha_moment_chapter_2}, \ref{ch2:as:L_smoothness_chapter_2}, \ref{ch2:as:str_cvx_chapter_2} with $\mu > 0$ hold for $Q = B_{3R}(x^*)$, where $R \geq \|x^0 - x^*\|$ and \algname{R-clipped-SSTM} runs \algname{clipped-SSTM} $\tau$ times. Let
    \begin{align}
        K_t
        &= \Biggl\lceil
        \max\Biggl\{1080\sqrt{\frac{LR_{t-1}^2}{\varepsilon_t}}
        \ln\frac{2160\sqrt{LR_{t-1}^2}\tau}{\sqrt{\varepsilon_t}\beta},\notag\\
        &\qquad 2\left(\frac{5400\sigma R_{t-1}}{\varepsilon_t}\right)^{p/(p-1)}
        \ln\left(\frac{4\tau}{\beta}
        \left(\frac{5400\sigma R_{t-1}}{\varepsilon_t}\right)^{p/(p-1)}\right)\Biggr\}
        \Biggr\rceil,
        \label{ch2:eq:R_clipped_SSTM_K_t_chapter_2}\\
        \varepsilon_t
        &= \frac{\mu R_{t-1}^2}{4},\quad
        R_{t-1} = \frac{R}{2^{(t-1)/2}},
        \quad \tau = \left\lceil \log_2 \frac{\mu R^2}{2\varepsilon} \right\rceil,
        \quad \ln\frac{4K_t\tau}{\beta} \geq 1,
        \label{ch2:eq:R_clipped_SSTM_epsilon_R_t_tau_chapter_2}\\
        a_t
        &= \max\left\{48600\ln^2\frac{4K_t\tau}{\beta},
        \frac{900\sigma(K_t+1)K_t^{1/p}
        \ln^{(p-1)/p}(4K_t\tau/\beta)}{L R_{t-1}}\right\},
        \label{ch2:eq:R_clipped_SSTM_parameter_a_chapter_2}\\
        \lambda_{k}^t
        &= \frac{R_{t-1}}{30a_{k+1}^t \ln(4K_t\tau/\beta)}
        \label{ch2:eq:R_clipped_SSTM_clipping_level_chapter_2}
    \end{align}
    for $t = 1, \ldots, \tau$, where $a_{k+1}^t = (k+2)/(2a_tL)$ are the step sizes of \algname{clipped-SSTM} at stage $t$. Then to guarantee $f(\hat x^\tau) - f(x^*) \leq \varepsilon$ with probability $\geq 1 - \beta$ \algname{R-clipped-SSTM} requires
    \begin{align}
        \cO\Biggl(
        \max\Biggl\{&\sqrt{\frac{L}{\mu}}
        \ln\left(\frac{\mu R^2}{\varepsilon}\right)
        \ln\left(\frac{\sqrt{L}}{\sqrt{\mu}\beta}
        \ln\left(\frac{\mu R^2}{\varepsilon}\right)\right),\notag\\
        &\left(\frac{\sigma^2}{\mu \varepsilon}\right)^{p/(2(p-1))}
        \ln\left(\frac{1}{\beta}
        \left(\frac{\sigma^2}{\mu \varepsilon}\right)^{p/(2(p-1))}
        \ln\left(\frac{\mu R^2}{\varepsilon}\right)\right)\Biggr\}
        \Biggr)
        \label{ch2:eq:R_clipped_SSTM_main_result_appendix_chapter_2}
    \end{align}
    iterations/oracle calls. Moreover, with probability $\geq 1-\beta$ the iterates of \algname{R-clipped-SSTM} at stage $t$ stay in the ball $B_{2R_{t-1}}(x^*)$.
\end{theorem}
\begin{proof}
    We show by induction that for any $t = 1,\ldots, \tau$ with probability at least $1 - t\beta/\tau$ inequalities
    \begin{equation}
        f(\hat{x}^l) - f(x^*) \leq \varepsilon_l, \quad \|\hat x^l - x^*\|^2 \leq  R_l^2 = \frac{R^2}{2^l}
    \end{equation}
    hold for $l = 1,\ldots, t$ simultaneously. First, we prove the base of the induction. Theorem~\ref{ch2:thm:clipped_SSTM_convex_case_appendix_chapter_2} (applied with $R_0 = R$ and $\beta/\tau$ in place of $\beta$) implies that with probability at least $1 - \beta/\tau$
    \begin{align*}
        \MoveEqLeft f(\hat x^1) - f(x^*)
        \leq \frac{6a_1LR^2}{K_1(K_1+3)}\\
        &= \max\Biggl\{\frac{291600LR^2\ln^2(4K_1\tau/\beta)}{K_1(K_1+3)},\\
        &\qquad\qquad \frac{5400\sigma R(K_1+1)K_1^{1/p}\ln^{(p-1)/p}(4K_1\tau/\beta)}{K_1(K_1+3)}\Biggr\}\\
        &\leq \max\left\{\frac{291600LR^2\ln^2(4K_1\tau/\beta)}{K_1^2},
        \frac{5400\sigma R\ln^{(p-1)/p}(4K_1\tau/\beta)}{K_1^{(p-1)/p}}\right\}\\
        &\leq \varepsilon_1 = \frac{\mu R^2}{4},
    \end{align*}
    where the equality follows from \eqref{ch2:eq:R_clipped_SSTM_parameter_a_chapter_2} and the last inequality follows from \eqref{ch2:eq:R_clipped_SSTM_K_t_chapter_2}. Due to the strong convexity,
    \begin{equation*}
        \|\hat x^1 - x^*\|^2 \leq \frac{2(f(\hat x^1) - f(x^*))}{\mu} \leq \frac{R^2}{2} = R_1^2.
    \end{equation*}
    The base of the induction is proven. Now, assume that the statement holds for some $t = T < \tau$, i.e., with probability at least $1 - T\beta/\tau$ inequalities
    \begin{equation}
        f(\hat{x}^l) - f(x^*) \leq \varepsilon_l, \quad \|\hat x^l - x^*\|^2 \leq  R_l^2 = \frac{R^2}{2^l}
    \end{equation}
    hold for $l = 1,\ldots, T$ simultaneously. In particular, with probability at least $1 - T\beta/\tau$ we have $\|\hat x^T - x^*\|^2 \leq R_T^2$. Applying Theorem~\ref{ch2:thm:clipped_SSTM_convex_case_appendix_chapter_2} (with $R_T$ in place of $R$ and $\beta/\tau$ in place of $\beta$) and using union bound for probability events, we get that with probability at least $1 - (T+1)\beta/\tau$
    \begin{align}
        \MoveEqLeft f(\hat x^{T+1}) - f(x^*)
        \leq \frac{6a_{T+1}LR_{T}^2}{K_{T+1}(K_{T+1}+3)}
        \notag\\
        &= \max\Biggl\{\frac{291600LR_{T}^2\ln^2(4K_{T+1}\tau/\beta)}
        {K_{T+1}(K_{T+1}+3)},\notag\\
        &\qquad\qquad \frac{5400\sigma R_T(K_{T+1}+1)K_{T+1}^{1/p}
        \ln^{(p-1)/p}(4K_{T+1}\tau/\beta)}
        {K_{T+1}(K_{T+1}+3)}\Biggr\}
        \notag\\
        &\leq \max\left\{\frac{291600LR_T^2\ln^2(4K_{T+1}\tau/\beta)}{K_{T+1}^2},
        \frac{5400\sigma R_T\ln^{(p-1)/p}(4K_{T+1}\tau/\beta)}
        {K_{T+1}^{(p-1)/p}}\right\}
        \notag\\
        &\leq \varepsilon_{T+1} = \frac{\mu R_T^2}{4},
        \notag
    \end{align}
    where we again use \eqref{ch2:eq:R_clipped_SSTM_parameter_a_chapter_2} and \eqref{ch2:eq:R_clipped_SSTM_K_t_chapter_2}. Due to the strong convexity,
    \begin{equation*}
        \|\hat x^{T+1} - x^*\|^2 \leq \frac{2(f(\hat x^{T+1}) - f(x^*))}{\mu} \leq \frac{R_T^2}{2} = R_{T+1}^2.
    \end{equation*}
    Thus, we finished the inductive part of the proof. In particular, with probability at least $1 - \beta$ inequalities
    \begin{equation}
        f(\hat{x}^l) - f(x^*) \leq \varepsilon_l, \quad \|\hat x^l - x^*\|^2 \leq  R_l^2 = \frac{R^2}{2^l}\notag
    \end{equation}
    hold for $l = 1,\ldots, \tau$ simultaneously, which gives for $l = \tau$ that with probability at least $1 - \beta$
    \begin{equation*}
        f(\hat{x}^\tau) - f(x^*) \leq \varepsilon_\tau = \frac{\mu R_{\tau-1}^2}{4} = \frac{\mu R^2}{2^{\tau+1}} \overset{\eqref{ch2:eq:R_clipped_SSTM_epsilon_R_t_tau_chapter_2}}{\leq} \varepsilon.
    \end{equation*}
    It remains to calculate the overall number of oracle calls during all runs of \algname{clipped-SSTM}. We have
    \begin{align*}
        \MoveEqLeft \sum\limits_{t=1}^\tau K_t\\
        &= \cO\Biggl(\sum\limits_{t=1}^\tau \max\Biggl\{
            \sqrt{\frac{LR_{t-1}^2}{\varepsilon_t}}
            \ln\left(\frac{\sqrt{LR_{t-1}^2}\tau}{\sqrt{\varepsilon_t}\beta}\right),\\
        &\qquad\qquad \left(\frac{\sigma R_{t-1}}{\varepsilon_t}\right)^{p/(p-1)}\ln\left(\frac{\tau}{\beta}
            \left(\frac{\sigma R_{t-1}}{\varepsilon_t}\right)^{p/(p-1)}\right)
            \Biggr\} \Biggr)\\
        &= \cO\Biggl(\sum\limits_{t=1}^\tau \max\Biggl\{
            \sqrt{\frac{L}{\mu}}\ln\left(\frac{\sqrt{L}\tau}{\sqrt{\mu}\beta}\right),\\
        &\qquad\qquad \left(\frac{\sigma}{\mu R_{t-1}}\right)^{p/(p-1)}\ln\left(\frac{\tau}{\beta}
            \left(\frac{\sigma }{\mu R_{t-1}}\right)^{p/(p-1)}\right)
            \Biggr\} \Biggr)\\
        &= \cO\Biggl(\max\Biggl\{
            \tau \sqrt{\frac{L}{\mu}}\ln\left(\frac{\sqrt{L}\tau}{\sqrt{\mu}\beta}\right),\\
        &\qquad\qquad \sum\limits_{t=1}^\tau             \left(\frac{\sigma \cdot 2^{t/2}}{\mu R}\right)^{p/(p-1)}\ln\left(\frac{\tau}{\beta}
            \left(\frac{\sigma \cdot 2^{t/2}}{\mu R}\right)^{p/(p-1)}\right)
            \Biggr\} \Biggr)\\
        &= \cO\Biggl(\max\Biggl\{
            \sqrt{\frac{L}{\mu}}\ln\left(\frac{\mu R^2}{\varepsilon}\right)
            \ln\left(\frac{\sqrt{L}}{\sqrt{\mu}\beta}
            \ln\left(\frac{\mu R^2}{\varepsilon}\right)\right),\\
        &\qquad\qquad \left(\frac{\sigma}{\mu R}\right)^{p/(p-1)}\ln\left(\frac{\tau}{\beta}
            \left(\frac{\sigma \cdot 2^{\tau/2}}{\mu R}\right)^{p/(p-1)}\right)
            \sum\limits_{t=1}^\tau 2^{pt/(2(p-1))}
            \Biggr\} \Biggr)\\
        &= \cO\Biggl(\max\Biggl\{
            \sqrt{\frac{L}{\mu}}\ln\left(\frac{\mu R^2}{\varepsilon}\right)
            \ln\left(\frac{\sqrt{L}}{\sqrt{\mu}\beta}
            \ln\left(\frac{\mu R^2}{\varepsilon}\right)\right),\\
        &\qquad\qquad \left(\frac{\sigma}{\mu R}\right)^{p/(p-1)}\ln\left(\frac{\tau}{\beta}
            \left(\frac{\sigma}{\mu R}\right)^{p/(p-1)}
            \cdot 2^{p\tau/(2(p-1))}\right)
            2^{p \tau/(2(p-1))}
            \Biggr\} \Biggr)\\
        &= \cO\Biggl(\max\Biggl\{
            \sqrt{\frac{L}{\mu}}\ln\left(\frac{\mu R^2}{\varepsilon}\right)
            \ln\left(\frac{\sqrt{L}}{\sqrt{\mu}\beta}
            \ln\left(\frac{\mu R^2}{\varepsilon}\right)\right),\\
        &\qquad\qquad \left(\frac{\sigma^2}{\mu \varepsilon}\right)^{p/(2(p-1))}
            \ln\left(\frac{1}{\beta}
            \left(\frac{\sigma^2}{\mu \varepsilon}\right)^{p/(2(p-1))}
            \ln\left(\frac{\mu R^2}{\varepsilon}\right)\right)
            \Biggr\} \Biggr),
    \end{align*}
    which concludes the proof.
\end{proof}

\clearpage

\section{Additional Proofs for \algname{clipped-SEG}}\label{ch2:appendix:SEG_chapter_2}
In this section, we provide the complete formulation of the main results for \algname{clipped-SEG} and the corresponding proofs. For brevity, we will use the following notation: $\thesischtwotF(x^{k};\xi_1^k) = \clip\left(F(x^{k};\xi_1^k), \lambda_k\right)$ and $\thesischtwotF(\thesischtwotx^{k};\xi_2^k) = \clip\left(F(\thesischtwotx^{k};\xi_2^k), \lambda_k\right)$.

\begin{algorithm}[h]
\caption{Clipped Stochastic Extragradient (\algname{clipped-SEG}) \parencite{gorbunov2022clipped}}
\label{ch2:alg:clipped-SEG_chapter_2}   
\begin{algorithmic}[1]
\REQUIRE starting point $x^0$, number of iterations $K$, stepsize $\gamma > 0$, clipping levels $\{\lambda_k\}_{k=0}^{K}$.
\FOR{$k=0,\ldots, K$}
\STATE Compute $\thesischtwotF(x^{k};\xi_1^k) = \clip\left(F(x^{k};\xi_1^k), \lambda_k\right)$ using a fresh sample $\xi_1^k \sim \thesischtwocD_k$
\STATE $\thesischtwotx^{k} = x^k - \gamma \thesischtwotF(x^{k};\xi_1^k)$
\STATE Compute $\thesischtwotF(\thesischtwotx^{k};\xi_2^k) = \clip\left(F(\thesischtwotx^{k};\xi_2^k), \lambda_k\right)$ using a fresh sample $\xi_2^k \sim \thesischtwocD_k$
\STATE $x^{k+1} = x^k - \gamma \thesischtwotF(\thesischtwotx^{k};\xi_2^k)$
\ENDFOR
\ENSURE $x^{K+1}$ or $\thesischtwotx_{\avg}^K = \frac{1}{K+1}\sum_{k=0}^K \thesischtwotx^k$
\end{algorithmic}
\end{algorithm}

\subsection{Monotone Problems}
We start with the following lemma derived by \textcite{gorbunov2022extragradient}. Since this lemma handles only deterministic part of the algorithm, the proof is the same as in the original work.
\begin{lemma}[Lemma C.1 from \textcite{gorbunov2022extragradient}]\label{ch2:lem:optimization_lemma_gap_SEG_chapter_2}
    Let Assumptions~\ref{ch2:as:L_Lip_chapter_2} and \ref{ch2:as:monotonicity_chapter_2} hold for $Q = B_{4R}(x^*)$, where $R \geq \|x^0 - x^*\|$ and $0 < \gamma \leq 1/(\sqrt{2}L)$. If $x^k$ and $\thesischtwotx^k$ lie in $B_{4R}(x^*)$ for all $k = 0,1,\ldots, K$ for some $K\geq 0$, then for all $u \in B_{4R}(x^*)$ the iterates produced by \algname{clipped-SEG} satisfy
    \begin{align}
        \langle F(u), \thesischtwotx^K_{\avg} - u\rangle & \leq \frac{\|x^0 - u\|^2 - \|x^{K+1} - u\|^2}{2\gamma(K+1)} \notag \\
        &\quad + \frac{\gamma}{2(K+1)}\sum\limits_{k=0}^K\left(\|\theta_k\|^2 + 2\|\omega_k\|^2\right) \notag\\
        &\quad\quad + \frac{1}{K+1}\sum\limits_{k=0}^K\langle x^k - u - \gamma F(\thesischtwotx^k), \theta_k\rangle, \label{ch2:eq:optimization_lemma_SEG_chapter_2}\\
        \thesischtwotx^K_{\avg} & \eqdef \frac{1}{K+1}\sum\limits_{k=0}^{K}\thesischtwotx^k, \label{ch2:eq:tx_avg_SEG_chapter_2}\\
        \theta_k & \eqdef F(\thesischtwotx^k) - \thesischtwotF(\thesischtwotx^k;\xi_2^k), \label{ch2:eq:theta_k_SEG_chapter_2}\\
        \omega_k & \eqdef F(x^k) - \thesischtwotF(x^k;\xi_1^k). \label{ch2:eq:omega_k_SEG_chapter_2}
    \end{align}
\end{lemma}

Using this lemma we prove the main convergence result for \algname{clipped-SEG} in the monotone case.
\begin{theorem}[Case 1 in Theorem~\ref{ch2:thm:clipped_SEG_main_theorem_chapter_2}]\label{ch2:thm:main_result_gap_SEG_chapter_2}
    Let Assumptions~\ref{ch2:as:bounded_alpha_moment_chapter_2}, \ref{ch2:as:L_Lip_chapter_2}, \ref{ch2:as:monotonicity_chapter_2} hold for $Q = B_{4R}(x^*)$, where $R \geq \|x^0 - x^*\|$. Let $K \geq 0$ and $\beta \in (0,1]$ be such that $A = \ln(6(K+1)/\beta) \geq 1$, and
    \begin{align}
        0< \gamma & \leq \min\left\{\frac{1}{160L A}, \frac{20^{(2-p)/p}R}{10800^{1/p}(K+1)^{1/p}\sigma A^{(p-1)/p}}\right\}, \label{ch2:eq:gamma_SEG_chapter_2}\\
        \lambda_{k} \equiv \lambda & = \frac{R}{20\gamma A}, \label{ch2:eq:lambda_SEG_chapter_2}
    \end{align}
    Then, after $K+1$ iterations the iterates produced by \algname{clipped-SEG} with probability at least $1 - \beta$ satisfy 
    \begin{align}
        & \gap_R(\thesischtwotx_{\avg}^K) \leq \frac{9R^2}{2\gamma(K+1)} \quad \text{and} \notag\\
        & \{x^k\}_{k=0}^{K+1} \subseteq B_{3R}(x^*), \{\thesischtwotx^k\}_{k=0}^{K+1} \subseteq B_{4R}(x^*), \label{ch2:eq:main_result_chapter_2}
    \end{align}
    where $\thesischtwotx_{\avg}^K$ is defined in \eqref{ch2:eq:tx_avg_SEG_chapter_2}. In particular, when $\gamma$ equals the minimum from \eqref{ch2:eq:gamma_SEG_chapter_2}, then the iterates produced by \algname{clipped-SEG} after $K+1$ iterations with probability at least $1-\beta$ satisfy
    \begin{equation}
        \gap_R(\thesischtwotx_{\avg}^K) = \cO\left(\max\left\{\frac{LR^2\ln(K/\beta)}{K}, \frac{\sigma R \ln^{(p-1)/p}(K/\beta)}{K^{(p-1)/p}}\right\}\right), \label{ch2:eq:clipped_SEG_monotone_case_2_appendix_chapter_2}
    \end{equation}
    meaning that \algname{clipped-SEG} requires the following number of iterations/oracle calls to achieve $\gap_R(\thesischtwotx_{\avg}^K) \leq \varepsilon$ with probability at least $1 - \beta$:
    \begin{equation}
        K = \cO\left(\max\left\{\frac{LR^2}{\varepsilon}\ln\frac{LR^2}{\varepsilon\beta}, \left(\frac{\sigma R}{\varepsilon}\right)^{p/(p-1)}\ln \left(\frac{1}{\beta}\left(\frac{\sigma R}{\varepsilon}\right)^{p/(p-1)}\right)\right\}\right). \label{ch2:eq:clipped_SEG_monotone_case_complexity_appendix_chapter_2}
    \end{equation}
\end{theorem}
\begin{proof}
    The proof follows similar steps as the proof of Theorem~C.1 from \textcite{gorbunov2022clipped}. The key difference is related to the application of Bernstein inequality and estimating biases and variances of stochastic terms.
    
    Let $R_k = \|x^k - x^*\|$ for all $k\geq 0$. As in the previous results, the proof is based on the induction argument and showing that the iterates do not leave some ball around the solution with high probability. More precisely, for each $k = 0,1,\ldots,K+1$ we consider probability event $E_k$ as follows: inequalities
    \begin{gather}
        \begin{aligned}
            A_t \eqdef \max\limits_{u \in B_{R}(x^*)}\biggl\{&\|x^0 - u\|^2 + 2\gamma \sum\limits_{l = 0}^{t-1} \langle x^l - u - \gamma F(\thesischtwotx^l), \theta_l \rangle\\
            &+ \gamma^2 \sum\limits_{l=0}^{t-1}\left(\|\theta_l\|^2 + 2\|\omega_l\|^2\right)\biggr\} \leq 9R^2,
        \end{aligned} \label{ch2:eq:induction_inequality_1_SEG_chapter_2}\\
        \left\|\gamma\sum\limits_{l=0}^{t-1}\theta_l\right\| \leq R \label{ch2:eq:induction_inequality_2_SEG_chapter_2}
    \end{gather}
    hold for $t = 0,1,\ldots,k$ simultaneously. We want to prove $\Prob\{E_k\} \geq  1 - k\beta/(K+1)$ for all $k = 0,1,\ldots,K+1$ by induction. The base of the induction is trivial: for $k=0$ we have $\|x^0 - u\|^2 \leq 2\|x^0 - x^*\|^2 + 2\|x^* - u\|^2 \leq 4R^2 \leq 9R^2$ and $\|\gamma\sum_{l=0}^{k-1}\theta_l\| = 0$ for any $u \in B_R(x^*)$. Next, assume that for $k = T-1 \leq K$ the statement holds: $\Prob\{E_{T-1}\} \geq  1 - (T-1)\beta/(K+1)$. Given this, we need to prove $\Prob\{E_{T}\} \geq  1 - T\beta/(K+1)$. We start with showing that $E_{T-1}$ implies $R_t \leq 3R$ for all $t = 0,1, \ldots, T$ (also by induction). For $t = 0$ this is already shown. Now, assume that $R_t \leq 3R$ for all $t = 0,1,\ldots,t'$ for some $t' < T$. Then for $t = 0,1,\ldots, t'$
    \begin{align}
        \|\thesischtwotx^t - x^*\| & = \|x^t - x^* - \gamma \thesischtwotF(x^t;\xi_1^t)\| \notag \\
        & \leq \|x^t - x^*\| + \gamma \|\thesischtwotF(x^t;\xi_1^t)\| \notag\\
        & \leq \|x^t - x^*\| + \gamma\lambda \overset{\eqref{ch2:eq:lambda_SEG_chapter_2}}{\leq} 3R + \frac{R}{20A}  \leq 4R. \label{ch2:eq:gap_thm_SEG_technical_1_chapter_2}
    \end{align}
    Therefore, the conditions of Lemma~\ref{ch2:lem:optimization_lemma_gap_SEG_chapter_2} are satisfied and we have that $E_{T-1}$ implies
    \begin{align*}
        \MoveEqLeft \max\limits_{u \in B_{R}(x^*)}\left\{2\gamma(t' + 1)\langle F(u), \thesischtwotx^{t'}_{\avg} - u\rangle + \|x^{t'+1} - u\|^2\right\}\\
        &\leq \max\limits_{u \in B_{R}(x^*)}\left\{ \|x^0 - u\|^2 + 2\gamma \sum\limits_{l = 0}^{t'} \langle x^l - u - \gamma F(\thesischtwotx^l), \theta_l \rangle\right\}\\
        &\quad + \gamma^2 \sum\limits_{l=0}^{t'}\left(\|\theta_l\|^2 + 2\|\omega_l\|^2\right)
        \overset{\eqref{ch2:eq:induction_inequality_1_SEG_chapter_2}}{\leq} 9R^2,
    \end{align*}
    meaning that (take $u = x^*$ and use $F(x^*) = 0$)
    \begin{align*}
\|x^{t' + 1} - x^*\|^2 & \leq \max\limits_{u \in B_{R}(x^*)}\left\{2\gamma(t' + 1)\langle F(u), \thesischtwotx^{t'}_{\avg} - u\rangle + \|x^{t'+1} - u\|^2\right\} \\
& \leq 9R^2,
    \end{align*}
    i.e., $R_{t'+1} \leq 3R$. In other words, we derived that probability event $E_{T-1}$ implies $R_{t} \leq 3R$ and
    \begin{equation}
        \max\limits_{u \in B_{R}(x^*)}\left\{2\gamma(t + 1)\langle F(u), \thesischtwotx^{t}_{\avg} - u\rangle + \|x^{t+1} - u\|^2\right\} \leq 9R^2 \label{ch2:eq:gap_thm_SEG_technical_1_5_chapter_2}
    \end{equation}
    for all $t = 0, 1, \ldots, T-1$. In addition, due to \eqref{ch2:eq:gap_thm_SEG_technical_1_chapter_2} $E_{T-1}$ also implies that $\|\thesischtwotx^t - x^*\| \leq 4R$ for all $t = 0, 1, \ldots, T-1$. Thus, $E_{T-1}$ implies
    \begin{align}
        \|x^t - x^* - \gamma F(\thesischtwotx^t)\| & \leq \|x^t - x^*\| + \gamma\|F(\thesischtwotx^t)\| \notag \\
        & \overset{\eqref{ch2:eq:L_Lip_chapter_2}}{\leq} 3R + \gamma L\|\thesischtwotx^t - x^*\| \notag\\
        & \overset{\eqref{ch2:eq:gap_thm_SEG_technical_1_chapter_2}}{\leq} 3R + 4R\gamma L \overset{\eqref{ch2:eq:gamma_SEG_chapter_2}}{\leq} 5R, \label{ch2:eq:gap_thm_SEG_technical_2_chapter_2}
    \end{align}
    for all $t = 0, 1, \ldots, T-1$. Next, we introduce random vectors
    \begin{equation*}
        u^{t} = \begin{cases}x^t - x^* - \gamma F(\thesischtwotx^t),& \text{if } \|x^t - x^* - \gamma F(\thesischtwotx^t)\| \leq 5R,\\ 0,& \text{otherwise,} \end{cases}
    \end{equation*}
    for all $t = 0, 1, \ldots, T-1$. These vectors are bounded almost surely:
     \begin{equation}
        \|u^{t}\| \leq 5R  \label{ch2:eq:gap_thm_SEG_technical_3_chapter_2}
    \end{equation}
    for all $t = 0, 1, \ldots, T-1$. Moreover, due to \eqref{ch2:eq:gap_thm_SEG_technical_2_chapter_2}, probability event $E_{T-1}$ implies $u^{t} = x^t - x^* - \gamma F(\thesischtwotx^t)$ for all $t = 0, 1, \ldots, T-1$ and
    \begin{align*}
        A_T
        &\leq \max\limits_{u \in B_{R}(x^*)}
        \left\{ \|x^0 - u\|^2
        + 2\gamma \sum\limits_{l = 0}^{T-1}
        \langle x^* - u, \theta_l \rangle \right\}\\
        &\quad + 2\gamma \sum\limits_{l = 0}^{T-1}
        \langle x^l - x^* - \gamma F(\thesischtwotx^l), \theta_l \rangle\\
        &\quad + \gamma^2 \sum\limits_{l=0}^{T-1}
        \left(\|\theta_l\|^2 + 2\|\omega_l\|^2\right)\\
        &\leq 4R^2
        + 2\gamma\max\limits_{u \in B_{R}(x^*)}
        \left\{\left\langle x^* - u,
        \sum\limits_{l = 0}^{T-1}\theta_l \right\rangle \right\}\\
        &\quad + 2\gamma \sum\limits_{l = 0}^{T-1}
        \langle u^{l}, \theta_l \rangle\\
        &\quad + \gamma^2 \sum\limits_{l=0}^{T-1}
        \left(\|\theta_l\|^2 + 2\|\omega_l\|^2\right)\\
        &= 4R^2
        + 2\gamma R \left\|\sum\limits_{l = 0}^{T-1}\theta_l\right\|\\
        &\quad + 2\gamma \sum\limits_{l = 0}^{T-1}
        \langle u^{l}, \theta_l \rangle\\
        &\quad + \gamma^2 \sum\limits_{l=0}^{T-1}
        \left(\|\theta_l\|^2 + 2\|\omega_l\|^2\right),
    \end{align*}
    where $A_T$ is defined in \eqref{ch2:eq:induction_inequality_1_SEG_chapter_2}.
    
    To handle the sums appeared in the right-hand side of the previous inequality we consider unbiased and biased parts of $\theta_l, \omega_l$:
     \begin{gather}
        \theta_l^u \eqdef \Exp_{\xi_2^l}\left[\thesischtwotF(\thesischtwotx^l;\xi_2^l)\right] - \thesischtwotF(\thesischtwotx^l;\xi_2^l),\quad \theta_l^b \eqdef F(\thesischtwotx^l) - \Exp_{\xi_2^l}\left[\thesischtwotF(\thesischtwotx^l;\xi_2^l)\right], \label{ch2:eq:gap_thm_SEG_technical_4_chapter_2}\\
        \omega_l^u \eqdef \Exp_{\xi_1^l}\left[\thesischtwotF(x^l;\xi_1^l)\right] - \thesischtwotF(x^l;\xi_1^l),\quad \omega_l^b \eqdef F(x^l) - \Exp_{\xi_1^l}\left[\thesischtwotF(x^l;\xi_1^l)\right], \label{ch2:eq:gap_thm_SEG_technical_5_chapter_2}
    \end{gather}
    for all $l = 0,\ldots, T-1$. By definition we have $\theta_l = \theta_l^u + \theta_l^b$, $\omega_l = \omega_l^u + \omega_l^b$ for all $l = 0,\ldots, T-1$. Therefore, $E_{T-1}$ implies
    \begin{align}
        A_T & \leq 4R^2 + 2\gamma R \left\|\sum\limits_{l = 0}^{T-1}\theta_l\right\| + \underbrace{2\gamma \sum\limits_{l = 0}^{T-1} \langle u^{l}, \theta_l^u \rangle}_{\circledOne} + \underbrace{2\gamma \sum\limits_{l = 0}^{T-1} \langle u^{l}, \theta_l^b \rangle}_{\circledTwo}\notag\\
        &\quad\quad + \underbrace{2\gamma^2 \sum\limits_{l=0}^{T-1}\left(\Exp_{\xi_2^l}\left[\|\theta_l^u\|^2\right] + 2\Exp_{\xi_1^l}\left[\|\omega_l^u\|^2\right]\right)}_{\circledThree} \notag\\
        &\quad\quad + \underbrace{2\gamma^2 \sum\limits_{l=0}^{T-1}\left(\|\theta_l^u\|^2 + 2\|\omega_l^u\|^2 - \Exp_{\xi_2^l}\left[\|\theta_l^u\|^2\right] - 2\Exp_{\xi_1^l}\left[\|\omega_l^u\|^2\right]\right)}_{\circledFour}\notag\\
        &\quad\quad + \underbrace{2\gamma^2 \sum\limits_{l=0}^{T-1}\left(\|\theta_l^b\|^2 + 2\|\omega_l^b\|^2\right)}_{\circledFive},\label{ch2:eq:gap_thm_SEG_technical_6_chapter_2}
    \end{align}
    where we also apply inequality $\|a+b\|^2 \leq 2\|a\|^2 + 2\|b\|^2$ holding for all $a,b \in \R^d$ to upper bound $\|\theta_l\|^2$ and $\|\omega_l\|^2$. It remains to derive good enough high-probability upper-bounds for the terms $2\gamma R \left\|\sum_{l = 0}^{T-1}\theta_l\right\|, \circledOne, \circledTwo, \circledThree, \circledFour, \circledFive$, i.e., to finish our inductive proof we need to show that $2\gamma R \left\|\sum_{l = 0}^{T-1}\theta_l\right\| + \circledOne + \circledTwo + \circledThree + \circledFour + \circledFive \leq 5R^2$ with high probability. In the subsequent parts of the proof, we will need to use many times the bounds for the norm and second moments of $\theta_{l}^u$, $\theta_{l}^b$, $\omega_l^u$ and $\omega_l^b$. First, by definition of clipping operator we have with probability $1$ that 
    \begin{equation}
        \|\theta_l^u\| \leq 2\lambda,\quad \|\omega_l^u\| \leq 2\lambda. \label{ch2:eq:theta_omega_magnitude_chapter_2}
    \end{equation}
    Moreover, since $E_{T-1}$ implies that
    \begin{gather*}
        \|F(x^l)\| \overset{\eqref{ch2:eq:L_Lip_chapter_2}}{\leq} L\|x^l - x^*\| \leq 3LR \overset{\eqref{ch2:eq:gamma_SEG_chapter_2}}{\leq} \frac{R}{40\gamma A} \overset{\eqref{ch2:eq:lambda_SEG_chapter_2}}{=} \frac{\lambda}{2},\\
        \|F(\thesischtwotx^l)\| \overset{\eqref{ch2:eq:L_Lip_chapter_2}}{\leq} L\|\thesischtwotx^l - x^*\| \overset{\eqref{ch2:eq:gap_thm_SEG_technical_1_chapter_2}}{\leq} 4LR \overset{\eqref{ch2:eq:gamma_SEG_chapter_2}}{\leq} \frac{R}{40\gamma A} \overset{\eqref{ch2:eq:lambda_SEG_chapter_2}}{=} \frac{\lambda}{2}
    \end{gather*}
    for $l = 0,1, \ldots, T-1$. Then, in view of Lemma~\ref{ch2:lem:bias_and_variance_clip_chapter_2}, we have that $E_{T-1}$ implies
    \begin{gather}
        \left\|\theta_l^b\right\| \leq \frac{2^{p}\sigma^{p}}{\lambda^{{p}-1}},\quad \left\|\omega_l^b\right\| \leq \frac{2^{p}\sigma^{p}}{\lambda^{{p}-1}}, \label{ch2:eq:bias_theta_omega_chapter_2}\\
        \Exp_{\xi_2^l}\left[\left\|\theta_l\right\|^2\right] \leq 18 \lambda^{2-{p}}\sigma^{p},\quad \Exp_{\xi_1^l}\left[\left\|\omega_l\right\|^2\right] \leq 18 \lambda^{2-{p}}\sigma^{p}, \label{ch2:eq:distortion_theta_omega_chapter_2}\\
        \Exp_{\xi_2^l}\left[\left\|\theta_l^u\right\|^2\right] \leq 18 \lambda^{2-{p}}\sigma^{p},\quad \Exp_{\xi_1^l}\left[\left\|\omega_l^u\right\|^2\right] \leq 18 \lambda^{2-{p}}\sigma^{p}, \label{ch2:eq:variance_theta_omega_chapter_2}
    \end{gather}
    for all $l = 0,1, \ldots, T-1$.

    \paragraph{Upper bound for $\circledOne$.} By definition of $\theta_{l}^u$, we have $\Exp_{\xi_2^l}[\theta_{l}^u] = 0$ and
    \begin{equation*}
        \Exp_{\xi_2^l}\left[2\gamma \langle u^{l}, \theta_l^u \rangle\right] = 0.
    \end{equation*}
    Next, sum $\circledOne$ has bounded with probability $1$ terms:
    \begin{equation}
        |2\gamma\langle u^{l}, \theta_l^u \rangle| \leq 2\gamma \|u^{l}\|\cdot \|\theta_l^u\| \overset{\eqref{ch2:eq:gap_thm_SEG_technical_3_chapter_2},\eqref{ch2:eq:theta_omega_magnitude_chapter_2}}{\leq} 20\gamma R\lambda \overset{\eqref{ch2:eq:lambda_SEG_chapter_2}}{=} \frac{R^2}{A} \eqdef c. \label{ch2:eq:gap_thm_SEG_technical_6_5_chapter_2}
    \end{equation}
    The summands also have bounded conditional variances $\sigma_l^2$ (the conditional expectations $\Exp_{\xi_2^l}$ of their squares):
    \begin{align}
\sigma_l^2 & \leq \Exp_{\xi_2^l}\left[4\gamma^2 \|u^{l}\|^2\cdot \|\theta_l^u\|^2\right] \notag \\
& \overset{\eqref{ch2:eq:gap_thm_SEG_technical_3_chapter_2}}{\leq} 100\gamma^2 R^2 \Exp_{\xi_2^l}\left[\|\theta_l^u\|^2\right]. \label{ch2:eq:gap_thm_SEG_technical_7_chapter_2}
    \end{align}
    In other words, we showed that the summands of $\circledOne$ form a bounded martingale difference sequence with bounded conditional variances $\{\sigma_l^2\}_{l = 0}^{T-1}$. Next, we apply Bernstein's inequality (Lemma~\ref{ch2:lem:Bernstein_ineq_chapter_2}) with $X_l = 2\gamma \langle u^{l}, \theta_l^u \rangle$, parameter $c$ as in \eqref{ch2:eq:gap_thm_SEG_technical_6_5_chapter_2}, $b = R^2$, $H = \tfrac{R^4}{6A}$:
    \begin{equation*}
        \Prob\left\{|\circledOne| > R^2 \text{ and } \sum\limits_{l=0}^{T-1}\sigma_l^2 \leq \frac{R^4}{6A}\right\} \leq 2\exp\left(- \frac{b^2}{2H + 2cb/3}\right) = \frac{\beta}{3(K+1)}.
    \end{equation*}
    Equivalently, we have
    \begin{align}
\Prob\{E_{\circledOne}\} & \geq 1 - \frac{\beta}{3(K+1)},\quad \text{for}\quad E_{\circledOne} \notag \\
& = \left\{\text{either} \quad \sum\limits_{l=0}^{T-1}\sigma_l^2 > \frac{R^4}{6A}\quad \text{or}\quad |\circledOne| \leq R^2\right\}. \label{ch2:eq:bound_1_gap_SEG_chapter_2}
    \end{align}
    In addition, $E_{T-1}$ implies that
    \begin{align}
        \sum\limits_{l=0}^{T-1}\sigma_l^2 & \overset{\eqref{ch2:eq:gap_thm_SEG_technical_7_chapter_2}}{\leq} 100\gamma^2 R^2 \sum\limits_{l=0}^{T-1} \Exp_{\xi_2^l}\left[\|\theta_l^u\|^2\right] \notag \\
        & \overset{\eqref{ch2:eq:variance_theta_omega_chapter_2}, T \leq K+1}{\leq} 1800(K+1)\gamma^2 R^2 \lambda^{2-{p}} \sigma^{p} \notag \\
        & \overset{\eqref{ch2:eq:lambda_SEG_chapter_2}}{=} \frac{1800 (K+1)\gamma^{{p}}\sigma^{p} R^{4-{p}}}{20^{2-{p}} A^{2-{p}}} \overset{\eqref{ch2:eq:gamma_SEG_chapter_2}}{\leq} \frac{R^4}{6A}. \label{ch2:eq:bound_1_variances_gap_SEG_chapter_2}
    \end{align}

    \textbf{Upper bound for $\circledTwo$.} From $E_{T-1}$ it follows that
    \begin{align}
        \circledTwo & \leq 2\gamma \sum\limits_{l=0}^{T-1}\|u^{l}\| \cdot \|\theta_l^b\| \overset{\eqref{ch2:eq:gap_thm_SEG_technical_3_chapter_2}, \eqref{ch2:eq:bias_theta_omega_chapter_2}, T \leq K+1}{\leq} \frac{10\cdot 2^{p}(K+1)\gamma R\sigma^{p}}{\lambda^{{p}-1}} \notag \\
        & \overset{\eqref{ch2:eq:lambda_SEG_chapter_2}}{=} \frac{10\cdot 2^{p}\cdot 20^{{p}-1} (K+1)\gamma^{p} \sigma^{p} A^{{p}-1}}{R^{{p}-2}} \overset{\eqref{ch2:eq:gamma_SEG_chapter_2}}{\leq} R^2. \label{ch2:eq:bound_2_variances_gap_SEG_chapter_2}
    \end{align}

    \textbf{Upper bound for $\circledThree$.} From $E_{T-1}$ it follows that
    \begin{align}
        2\gamma^2 \sum\limits_{l=0}^{T-1}
        \Exp_{\xi_2^l}[\|\theta_l^u\|^2]
        &\overset{\eqref{ch2:eq:distortion_theta_omega_chapter_2}, T \leq K+1}{\leq}
        36\gamma^2 (K+1)\lambda^{2-{p}}\sigma^{p}
        \notag\\
        &\overset{\eqref{ch2:eq:lambda_SEG_chapter_2}}{=}
        \frac{36\gamma^{p} (K+1)\sigma^{p}R^{2-{p}}}
        {20^{2-{p}} A^{2-{p}}}
        \overset{\eqref{ch2:eq:gamma_SEG_chapter_2}}{\leq}
        \frac{1}{12} R^2,
        \label{ch2:eq:sum_theta_squared_bound_gap_SEG_chapter_2}\\
        4\gamma^2 \sum\limits_{l=0}^{T-1}
        \Exp_{\xi_1^l}[\|\omega_l^u\|^2]
        &\overset{\eqref{ch2:eq:distortion_theta_omega_chapter_2}, T \leq K+1}{\leq}
        72\gamma^2 (K+1)\lambda^{2-{p}}\sigma^{p}
        \notag\\
        &\overset{\eqref{ch2:eq:lambda_SEG_chapter_2}}{=}
        \frac{72\gamma^{p} (K+1)\sigma^{p}R^{2-{p}}}
        {20^{2-{p}} A^{2-{p}}}
        \overset{\eqref{ch2:eq:gamma_SEG_chapter_2}}{\leq}
        \frac{1}{12} R^2,
        \label{ch2:eq:sum_omega_squared_bound_gap_SEG_chapter_2}\\
        \circledThree
        &\overset{\eqref{ch2:eq:sum_theta_squared_bound_gap_SEG_chapter_2},
        \eqref{ch2:eq:sum_omega_squared_bound_gap_SEG_chapter_2}}{\leq}
        \frac{1}{6}R^2.
        \label{ch2:eq:bound_3_variances_gap_SEG_chapter_2}
    \end{align}

    \paragraph{Upper bound for $\circledFour$.} By the construction we have
    \begin{equation*}
        2\gamma^2\Exp_{\xi_1^l,\xi_2^l}\left[\|\theta_l^u\|^2 + 2\|\omega_l^u\|^2 - \Exp_{\xi_2^l}\left[\|\theta_l^u\|^2\right] - 2\Exp_{\xi_1^l}\left[\|\omega_l^u\|^2\right]\right] = 0.
    \end{equation*}
    Next, sum $\circledFour$ has bounded with probability $1$ terms:
    \begin{align}
        &2\gamma^2\left|\|\theta_l^u\|^2 + 2\|\omega_l^u\|^2
        - \Exp_{\xi_2^l}\left[\|\theta_l^u\|^2\right]
        - 2\Exp_{\xi_1^l}\left[\|\omega_l^u\|^2\right] \right|
        \notag\\
        &\quad \leq 2\gamma^2 \|\theta_l^u\|^2
        + 2\gamma^2\Exp_{\xi_2^l}\left[\|\theta_l^u\|^2\right]
        \notag\\
        &\qquad + 4\gamma^2 \|\omega_l^u\|^2
        + 4\gamma^2\Exp_{\xi_1^l}\left[\|\omega_l^u\|^2\right]
        \notag\\
        &\quad \overset{\eqref{ch2:eq:theta_omega_magnitude_chapter_2}}{\leq}
        48\gamma^2 \lambda^2
        \notag\\
        &\quad \overset{\eqref{ch2:eq:lambda_SEG_chapter_2}}{\leq}
        \frac{R^2}{6A}
        \eqdef c.
        \label{ch2:eq:gap_thm_SEG_technical_6_5_1_chapter_2}
    \end{align}
    The summands also have bounded conditional variances
    \begin{equation*}
        \widetilde\sigma_l^2 \eqdef
        4\gamma^4\Exp_{\xi_1^l, \xi_2^l}\left[
        \left|
        \begin{aligned}
        &\|\theta_l^u\|^2 + 2\|\omega_l^u\|^2
        - \Exp_{\xi_2^l}\left[\|\theta_l^u\|^2\right]\\
        &\qquad
        - 2\Exp_{\xi_1^l}\left[\|\omega_l^u\|^2\right]
        \end{aligned}
        \right|^2\right]:
    \end{equation*}
    \begin{align}
        \widetilde\sigma_l^2 & \overset{\eqref{ch2:eq:gap_thm_SEG_technical_6_5_1_chapter_2}}{\leq} \frac{\gamma^2R^2}{3A}\Exp_{\xi_1^l,\xi_2^l}\left[\left|\|\theta_l^u\|^2 + 2\|\omega_l^u\|^2 - \Exp_{\xi_2^l}\left[\|\theta_l^u\|^2\right] - 2\Exp_{\xi_1^l}\left[\|\omega_l^u\|^2\right] \right|\right]\notag\\
        & \leq \frac{2\gamma^2R^2}{3A} \Exp_{\xi_1^l,\xi_2^l}\left[\|\theta_l^u\|^2 + 2\|\omega_l^u\|^2 \right]. \label{ch2:eq:gap_thm_SEG_technical_7_1_chapter_2}
    \end{align}
    In other words, we showed that the summands of $\circledFour$ form a bounded martingale difference sequence with bounded conditional variances $\{\widetilde\sigma_l^2\}_{l = 0}^{T-1}$. Next, we apply Bernstein's inequality (Lemma~\ref{ch2:lem:Bernstein_ineq_chapter_2}) with $X_l$ equal to the summands of $\circledFour$, parameter $c$ as in \eqref{ch2:eq:gap_thm_SEG_technical_6_5_1_chapter_2}, $b = \tfrac{1}{6}R^2$, $H = \tfrac{R^4}{216A}$:
    \begin{align*}
\MoveEqLeft[1] \Prob\left\{|\circledFour| > \frac{1}{6}R^2 \text{ and } \sum\limits_{l=0}^{T-1}\widetilde\sigma_l^2 \leq \frac{R^4}{216A}\right\} \\
& \leq 2\exp\left(- \frac{b^2}{2H + 2cb/3}\right) = \frac{\beta}{3(K+1)}.
    \end{align*}
    Equivalently, we have
    \begin{align}
\Prob\{E_{\circledFour}\} & \geq 1 - \frac{\beta}{3(K+1)},\quad \text{for}\quad E_{\circledFour} \notag \\
& = \left\{\text{either} \quad \sum\limits_{l=0}^{T-1}\widetilde\sigma_l^2 > \frac{R^4}{216A}\quad \text{or}\quad |\circledFour| \leq \frac{1}{6}R^2\right\}. \label{ch2:eq:bound_1_gap_SEG_1_chapter_2}
    \end{align}
    In addition, $E_{T-1}$ implies that
    \begin{align}
        \sum\limits_{l=0}^{T-1}\widetilde\sigma_l^2 & \overset{\eqref{ch2:eq:gap_thm_SEG_technical_7_1_chapter_2}}{\leq} \frac{2\gamma^2R^2}{3A} \sum\limits_{l=0}^{T-1} \Exp_{\xi_1^l,\xi_2^l}\left[\|\theta_l^u\|^2 + 2\|\omega_l^u\|^2 \right]\notag\\
        & \overset{\eqref{ch2:eq:variance_theta_omega_chapter_2}, T \leq K+1}{\leq} \frac{36(K+1)\gamma^2 R^2 \lambda^{2-{p}} \sigma^{p}}{A} \notag\\
        & \overset{\eqref{ch2:eq:lambda_SEG_chapter_2}}{\leq} \frac{36(K+1)\gamma^{p} R^{4-{p}} \sigma^{p}}{20^{2-{p}}A^{3-{p}}}  \overset{\eqref{ch2:eq:gamma_SEG_chapter_2}}{\leq} \frac{R^4}{216A}. \label{ch2:eq:bound_1_variances_gap_SEG_1_chapter_2}
    \end{align}

    \textbf{Upper bound for $\circledFive$.} From $E_{T-1}$ it follows that
    \begin{align}
        \circledFive & = 2\gamma^2 \sum\limits_{l=0}^{T-1}\left(\|\theta_l^b\|^2 + 2\|\omega_l^b\|^2\right) \notag \\
        & \overset{\eqref{ch2:eq:bias_theta_omega_chapter_2}, T\leq K+1}{\leq} \frac{6\cdot 2^{2{p}}\gamma^2\sigma^{2{p}} (K+1)}{\lambda^{2{p}-2}} \notag\\
        & \overset{\eqref{ch2:eq:lambda_SEG_chapter_2}}{=}  \frac{6\cdot 2^{2{p}}\cdot 20^{2{p} - 2}\gamma^{2{p}}\sigma^{2{p}} (K+1) A^{2{p}-2}}{R^{2{p}-2}} \overset{\eqref{ch2:eq:gamma_SEG_chapter_2}}{\leq} \frac{1}{6}R^2. \label{ch2:eq:nsjcbdjhcbfjdhfbj_chapter_2}
    \end{align}

    \paragraph{Upper bound for $2\gamma R \left\|\sum_{l=0}^{T-1} \theta_l\right\|$.} To upper-bound this sum, we introduce new random vectors:
    \begin{equation*}
        \zeta_l = \begin{cases} \gamma \sum\limits_{r=0}^{l-1}\theta_r,& \text{if } \left\|\gamma \sum\limits_{r=0}^{l-1}\theta_r\right\| \leq R,\\ 0, & \text{otherwise} \end{cases}
    \end{equation*}
    for $l = 0, 1, \ldots, T-1$. These vectors are bounded with probability $1$:
    \begin{equation}
        \|\zeta_l\| \leq R.  \label{ch2:eq:gap_thm_SEG_technical_8_chapter_2}
    \end{equation}
    Therefore, taking into account \eqref{ch2:eq:induction_inequality_2_SEG_chapter_2}, we derive that $E_{T-1}$ implies
    \begin{align}
        2\gamma R\left\|\sum\limits_{l = 0}^{T-1}\theta_l\right\| & = 2R\sqrt{\gamma^2\left\|\sum\limits_{l = 0}^{T-1}\theta_l\right\|^2}\notag\\
        & = 2R\sqrt{\gamma^2\sum\limits_{l=0}^{T-1}\|\theta_l\|^2 + 2\gamma\sum\limits_{l=0}^{T-1}\left\langle \gamma\sum\limits_{r=0}^{l-1} \theta_r, \theta_l \right\rangle} \notag\\
        & = 2R \sqrt{\gamma^2\sum\limits_{l=0}^{T-1}\|\theta_l\|^2 + 2\gamma \sum\limits_{l=0}^{T-1} \langle \zeta_l, \theta_l\rangle} \notag\\
        & \overset{\eqref{ch2:eq:gap_thm_SEG_technical_4_chapter_2}}{\leq} 2R \sqrt{\circledThree + \circledFour + \circledFive + \underbrace{2\gamma \sum\limits_{l=0}^{T-1} \langle \zeta_l, \theta_l^u\rangle}_{\circledSix} + \underbrace{2\gamma \sum\limits_{l=0}^{T-1} \langle \zeta_l, \theta_l^b\rangle}_{\circledSeven}}. \label{ch2:eq:norm_sum_theta_bound_gap_SEG_chapter_2}
    \end{align}
    Similarly to the previous parts of the proof, we bound $\circledSix$ and $\circledSeven$.

    \paragraph{Upper bound for $\circledSix$.} By definition of $\theta_{l}^u$, we have $\Exp_{\xi_2^l}[\theta_{l}^u] = 0$ and
    \begin{equation*}
        \Exp_{\xi_2^l}\left[2\gamma \langle \zeta_l, \theta_l^u \rangle\right] = 0.
    \end{equation*}
    Next, sum $\circledSix$ has bounded with probability $1$ terms:
    \begin{equation}
        |2\gamma\langle \zeta_l, \theta_l^u \rangle| \leq 2\gamma \|\zeta_l\|\cdot \|\theta_l^u\| \overset{\eqref{ch2:eq:gap_thm_SEG_technical_8_chapter_2},\eqref{ch2:eq:theta_omega_magnitude_chapter_2}}{\leq} 4\gamma R\lambda \overset{\eqref{ch2:eq:lambda_SEG_chapter_2}}{\leq} \frac{R^2}{4A} \eqdef c. \label{ch2:eq:gap_thm_SEG_technical_8_5_chapter_2}
    \end{equation}
    The summands also have bounded conditional variances $\hat\sigma_l^2$ (the conditional expectations $\Exp_{\xi_2^l}$ of their squares):
    \begin{align}
\hat\sigma_l^2 & \leq \Exp_{\xi_2^l}\left[4\gamma^2 \|\zeta_l\|^2\cdot \|\theta_l^u\|^2\right] \notag \\
& \overset{\eqref{ch2:eq:gap_thm_SEG_technical_8_chapter_2}}{\leq} 4\gamma^2 R^2 \Exp_{\xi_2^l}\left[\|\theta_l^u\|^2\right]. \label{ch2:eq:gap_thm_SEG_technical_9_chapter_2}
    \end{align}
    In other words, we showed that the summands of $\circledSix$ form a bounded martingale difference sequence with bounded conditional variances $\{\hat \sigma_l^2\}_{l = 0}^{T-1}$. Next, we apply Bernstein's inequality (Lemma~\ref{ch2:lem:Bernstein_ineq_chapter_2}) with $X_l = 2\gamma \langle \zeta_l, \theta_l^u \rangle$, parameter $c$ as in \eqref{ch2:eq:gap_thm_SEG_technical_8_5_chapter_2}, $b = \tfrac{R^2}{4}$, $H = \tfrac{R^4}{96A}$:
    \begin{align*}
\MoveEqLeft[1] \Prob\left\{|\circledSix| > \frac{1}{4}R^2 \text{ and } \sum\limits_{l=0}^{T-1}\hat\sigma_l^2 \leq \frac{R^4}{96A}\right\} \\
& \leq 2\exp\left(- \frac{b^2}{2H + 2cb/3}\right) = \frac{\beta}{3(K+1)}.
    \end{align*}
    Equivalently, we have $\Prob\{E_{\circledSix}\} \geq 1 - \beta/(3(K+1))$ for
    \begin{equation}
        E_{\circledSix} = \left\{ \text{either} \quad \sum\limits_{l=0}^{T-1}\hat\sigma_l^2 > \frac{R^4}{96A}\quad \text{or}\quad |\circledSix| \leq \frac{1}{4}R^2\right\}. \label{ch2:eq:bound_4_gap_SEG_chapter_2}
    \end{equation}
    In addition, $E_{T-1}$ implies that
    \begin{align}
        \sum\limits_{l=0}^{T-1}\hat\sigma_l^2 & \overset{\eqref{ch2:eq:gap_thm_SEG_technical_9_chapter_2}}{\leq} 4\gamma^2 R^2 \sum\limits_{l=0}^{T-1} \Exp_{\xi_2^l}\left[\|\theta_l^u\|^2\right] \notag \\
        & \overset{\eqref{ch2:eq:variance_theta_omega_chapter_2}, T \leq K+1}{\leq} 72(K+1)\gamma^2 R^2 \lambda^{2-{p}} \sigma^{p} \notag\\
        & \overset{\eqref{ch2:eq:lambda_SEG_chapter_2}}{=} \frac{72(K+1)\gamma^{p} R^{4-{p}} \sigma^{p}}{20^{2-{p}} A^{2-{p}}} \overset{\eqref{ch2:eq:gamma_SEG_chapter_2}}{\leq} \frac{R^4}{96A}. \label{ch2:eq:bound_4_variances_gap_SEG_chapter_2}
    \end{align}

    \textbf{Upper bound for $\circledSeven$.} From $E_{T-1}$ it follows that
    \begin{align}
        \circledSeven & \leq 2\gamma \sum\limits_{l=0}^{T-1}\|\zeta_l\| \cdot \|\theta_l^b\| \overset{\eqref{ch2:eq:gap_thm_SEG_technical_8_chapter_2}, \eqref{ch2:eq:bias_theta_omega_chapter_2}, T \leq K+1}{\leq} \frac{2^{{p}+1}(K+1)\gamma R\sigma^{p}}{\lambda^{{p}-1}} \notag \\
        & \overset{\eqref{ch2:eq:lambda_SEG_chapter_2}}{=} \frac{2^{{p}+1}\cdot 20^{{p}-1}(K+1)\gamma^{p} \sigma^{p} A^{{p}-1}}{R^{{p}-2}} \overset{\eqref{ch2:eq:gamma_SEG_chapter_2}}{\leq} \frac{1}{4}R^2. \label{ch2:eq:bound_5_variances_gap_SEG_chapter_2}
    \end{align}

    Now, we have the upper bounds for  $2\gamma R \left\|\sum_{l = 0}^{T-1}\theta_l\right\|, \circledOne, \circledTwo, \circledThree, \circledFour, \circledFive, \circledSix, \circledSeven$. In particular, probability event $E_{T-1}$ implies
    \begin{gather*}
        A_T \overset{\eqref{ch2:eq:gap_thm_SEG_technical_6_chapter_2}}{\leq} 4R^2 + 2\gamma R\left\|\sum\limits_{l=0}^{T-1} \theta_l\right\| + \circledOne + \circledTwo + \circledThree + \circledFour + \circledFive,\\
        2\gamma R\left\|\sum\limits_{l=0}^{T-1} \theta_l\right\| \overset{\eqref{ch2:eq:norm_sum_theta_bound_gap_SEG_chapter_2}}{\leq} 2R \sqrt{\circledThree + \circledFour + \circledFive + \circledSix + \circledSeven},\\
        \circledTwo \overset{\eqref{ch2:eq:bound_2_variances_gap_SEG_chapter_2}}{\leq} R^2,\quad \circledThree \overset{\eqref{ch2:eq:bound_3_variances_gap_SEG_chapter_2}}{\leq} \frac{1}{6}R^2,\quad \circledFive \overset{\eqref{ch2:eq:nsjcbdjhcbfjdhfbj_chapter_2}}{\leq} \frac{1}{6}R^2,\quad \circledSeven \overset{\eqref{ch2:eq:bound_5_variances_gap_SEG_chapter_2}}{\leq} \frac{1}{4}R^2,\\
        \sum\limits_{l=0}^{T-1}\sigma_l^2 \overset{\eqref{ch2:eq:bound_1_variances_gap_SEG_chapter_2}}{\leq}  \frac{R^4}{6A},\quad \sum\limits_{l=0}^{T-1}\widetilde\sigma_l^2 \overset{\eqref{ch2:eq:bound_1_variances_gap_SEG_1_chapter_2}}{\leq} \frac{R^4}{216A},\quad \sum\limits_{l=0}^{T-1}\hat\sigma_l^2 \overset{\eqref{ch2:eq:bound_4_variances_gap_SEG_chapter_2}}{\leq}  \frac{R^4}{96A}.
    \end{gather*}
    Moreover, we also have (see \eqref{ch2:eq:bound_1_gap_SEG_chapter_2}, \eqref{ch2:eq:bound_1_gap_SEG_1_chapter_2}, \eqref{ch2:eq:bound_4_gap_SEG_chapter_2} and our induction assumption)
    \begin{gather*}
        \Prob\{E_{T-1}\} \geq 1 - \frac{(T-1)\beta}{K+1},\\
        \Prob\{E_{\circledOne}\} \geq 1 - \frac{\beta}{3(K+1)}, \quad \Prob\{E_{\circledFour}\} \geq 1 - \frac{\beta}{3(K+1)}, \quad \Prob\{E_{\circledSix}\} \geq 1 - \frac{\beta}{3(K+1)} ,
    \end{gather*}
    where
    \begin{align}
        E_{\circledOne}& = \left\{\text{either} \quad \sum\limits_{l=0}^{T-1}\sigma_l^2 > \frac{R^4}{6A}\quad \text{or}\quad |\circledOne| \leq R^2\right\},\notag\\
        E_{\circledFour}& = \left\{\text{either} \quad \sum\limits_{l=0}^{T-1}\widetilde\sigma_l^2 > \frac{R^4}{216A}\quad \text{or}\quad |\circledFour| \leq \frac{1}{6}R^2\right\},\notag\\
        E_{\circledSix}& = \left\{\text{either} \quad \sum\limits_{l=0}^{T-1}\hat\sigma_l^2 > \frac{R^4}{96A}\quad \text{or}\quad |\circledSix| \leq \frac{1}{4}R^2\right\}.\notag
    \end{align}
    Thus, probability event $E_{T-1} \cap E_{\circledOne} \cap E_{\circledFour} \cap E_{\circledSix}$ implies
    \begin{align}
        \left\|\gamma\sum\limits_{l=0}^{T-1} \theta_l\right\| & \leq \sqrt{\frac{1}{6}R^2 + \frac{1}{6}R^2 + \frac{1}{6}R^2 + \frac{1}{4}R^2 + \frac{1}{4}R^2} = R, \label{ch2:eq:gap_thm_SEG_technical_10_chapter_2} \\
        A_T & \leq 4R^2 + 2R\sqrt{\frac{1}{6}R^2 + \frac{1}{6}R^2 + \frac{1}{6}R^2 + \frac{1}{4}R^2 + \frac{1}{4}R^2}\notag\\
        &\quad\quad + R^2 + R^2 + \frac{1}{6}R^2 + \frac{1}{6}R^2 + \frac{1}{6}R^2\notag\\
        & \leq 9R^2, \label{ch2:eq:gap_thm_SEG_technical_11_chapter_2}
    \end{align}
    which is equivalent to \eqref{ch2:eq:induction_inequality_1_SEG_chapter_2} and \eqref{ch2:eq:induction_inequality_2_SEG_chapter_2} for $t = T$, and 
    \begin{align*}
\Prob\{E_T\} & \geq \Prob\{E_{T-1} \cap E_{\circledOne} \cap E_{\circledFour} \cap E_{\circledSix}\} \\
& = 1 - \Prob\{\overline{E}_{T-1} \cup \overline{E}_{\circledOne} \cup \overline{E}_{\circledFour} \cup \overline{E}_{\circledSix}\} \geq 1 - \frac{T\beta}{K+1}.
    \end{align*}
    This finishes the inductive part of our proof, i.e., for all $k = 0,1,\ldots,K+1$ we have $\Prob\{E_k\} \geq 1 - k\beta/(K+1)$. In particular, for $k = K+1$ we have that with probability at least $1 - \beta$
    \begin{align*}
        \gap_{R}(\thesischtwotx^{K}_{\avg}) & = \max\limits_{u \in B_{R}(x^*)}\left\{\langle F(u), \thesischtwotx^{K}_{\avg} - u\rangle\right\}\\
        & \leq \frac{1}{2\gamma(K + 1)}\max\limits_{u \in B_{R}(x^*)}\left\{2\gamma(K + 1)\langle F(u), \thesischtwotx^{K}_{\avg} - u\rangle + \|x^{K+1} - u\|^2\right\}\\
        & \overset{\eqref{ch2:eq:gap_thm_SEG_technical_1_5_chapter_2}}{\leq} \frac{9R^2}{2\gamma(K+1)}.
    \end{align*}
    Finally, if $\gamma$ equals the minimum in \eqref{ch2:eq:gamma_SEG_chapter_2},
    then with probability at least $1-\beta$
    \begin{align*}
        \gap_{R}(\thesischtwotx^{K}_{\avg}) & \leq \frac{9R^2}{2\gamma(K+1)} \\
        & = \max\left\{\frac{720 LR^2A}{K+1}, \frac{9\cdot 10800^{1/p} \sigma R A^{(p-1)/p}}{2\cdot 20^{(2-p)/p} (K+1)^{(p-1)/p}}\right\}\\
        & = \cO\left(\max\left\{\frac{LR^2\ln(K/\beta)}{K}, \frac{\sigma R \ln^{(p-1)/p}(K/\beta)}{K^{(p-1)/p}}\right\}\right).
    \end{align*}
    To get $\gap_{R}(\thesischtwotx^{K}_{\avg}) \leq \varepsilon$ with probability at least $1-\beta$ it is sufficient to choose $K$ such that both terms in the maximum above are $\cO(\varepsilon)$. This leads to
    \begin{equation*}
         K = \cO\left(\max\left\{\frac{LR^2}{\varepsilon}\ln\frac{LR^2}{\varepsilon\beta}, \left(\frac{\sigma R}{\varepsilon}\right)^{p/(p-1)}\ln \left(\frac{1}{\beta}\left(\frac{\sigma R}{\varepsilon}\right)^{p/(p-1)}\right)\right\}\right)
    \end{equation*}
    that concludes the proof.
\end{proof}

\subsection{Quasi-Strongly Monotone Problems}
As in the monotone case, we use another lemma from \textcite{gorbunov2022clipped} that handles the deterministic part of \algname{clipped-SEG} in the quasi-strongly monotone case.
\begin{lemma}[Lemma~C.3 from \textcite{gorbunov2022clipped}]\label{ch2:lem:optimization_lemma_str_mon_SEG_chapter_2}
    Let Assumptions~\ref{ch2:as:L_Lip_chapter_2}, \ref{ch2:as:QSM_chapter_2} hold for $Q = B_{3R}(x^*) = \{x\in\R^d\mid \|x - x^*\| \leq 3R\}$, where $R \geq \|x^0 - x^*\|$, and $0 < \gamma \leq 1/(2(L+2\mu))$. If $x^k$ and $\thesischtwotx^k$ lie in $B_{3R}(x^*)$ for all $k = 0,1,\ldots, K$ for some $K\geq 0$, then the iterates produced by \algname{clipped-SEG} satisfy
    \begin{align}
        \|x^{K+1} - x^*\|^2 & \leq (1 - \gamma \mu)^{K+1}\|x^0 - x^*\|^2 \notag \\
        &\quad - 4\gamma^3 \mu \sum\limits_{k=0}^K (1-\gamma\mu)^{K-k}\langle F(x^k), \omega_k \rangle \notag\\
        &\quad\quad + 2\gamma \sum\limits_{k=0}^K (1-\gamma\mu)^{K-k} \langle x^k - x^* - \gamma F(\thesischtwotx^k), \theta_k \rangle\notag\\
        &\quad\quad + \gamma^2 \sum\limits_{k=0}^K (1-\gamma\mu)^{K-k} \left(\|\theta_k\|^2 + 4\|\omega_k\|^2\right), \label{ch2:eq:optimization_lemma_SEG_str_mon_chapter_2}
    \end{align}
    where $\theta_k, \omega_k$ are defined in \eqref{ch2:eq:theta_k_SEG_chapter_2}, \eqref{ch2:eq:omega_k_SEG_chapter_2}.
\end{lemma}

Using this lemma we prove the main convergence result for \algname{clipped-SEG} in the quasi-strongly monotone case.

\begin{theorem}[Case 2 in Theorem~\ref{ch2:thm:clipped_SEG_main_theorem_chapter_2}]\label{ch2:thm:main_result_str_mon_SEG_chapter_2}
    Let Assumptions~\ref{ch2:as:bounded_alpha_moment_chapter_2}, \ref{ch2:as:L_Lip_chapter_2}, \ref{ch2:as:QSM_chapter_2} hold for $Q = B_{3R}(x^*) = \{x\in\R^d\mid \|x - x^*\| \leq 3R\}$, where $R \geq \|x^0 - x^*\|$. Let $K \geq 0$ and $\beta \in (0,1]$ be such that $A = \ln(6(K+1)/\beta) \geq 1$, and
    \begin{align}
        0< \gamma & \leq \min\left\{\frac{1}{650 L A}, \frac{\ln(B_K)}{\mu(K+1)}\right\}, \label{ch2:eq:gamma_SEG_str_mon_chapter_2}\\
        B_K & = \max\left\{2, \frac{(K+1)^{2(p-1)/p}\mu^2R^2}{264600^{2/p}\sigma^2A^{2(p-1)/p}\ln^2(B_K)} \right\} \label{ch2:eq:B_K_SEG_str_mon_1_chapter_2} \\
        & = \cO\biggl(\max\biggl\{2, \frac{K^{2(p-1)/p}\mu^2R^2}{\sigma^2\ln^{2(p-1)/p}(K/\beta)}\notag\\
        &\qquad\qquad \cdot \frac{1}{\ln^2\left(\max\left\{2, \frac{K^{2(p-1)/p}\mu^2R^2}{\sigma^2\ln^{2(p-1)/p}(K/\beta)} \right\}\right)} \biggr\}\biggr), \label{ch2:eq:B_K_SEG_str_mon_2_chapter_2} \\
        \lambda_k & = \frac{\exp(-\gamma\mu(1 + k/2))R}{120\gamma A}, \label{ch2:eq:lambda_SEG_str_mon_chapter_2}
    \end{align}
    where $K$ is large enough that the maximum in \eqref{ch2:eq:B_K_SEG_str_mon_1_chapter_2} is attained at its second term. Then, after $K+1$ iterations the iterates produced by \algname{clipped-SEG} with probability at least $1 - \beta$ satisfy 
    \begin{equation}
        \|x^{K+1} - x^*\|^2 \leq 2\exp(-\gamma\mu(K+1))R^2. \label{ch2:eq:main_result_str_mon_chapter_2}
    \end{equation}
    In particular, when $\gamma$ equals the minimum from \eqref{ch2:eq:gamma_SEG_str_mon_chapter_2}, then the iterates produced by \algname{clipped-SEG} after $K+1$ iterations with probability at least $1-\beta$ satisfy
    \begin{align}
        \MoveEqLeft \|x^{K+1} - x^*\|^2\notag\\
        &= \cO\Biggl(\max\Biggl\{
        R^2\exp\left(- \frac{\mu K}{L \ln(K/\beta)}\right),
        \frac{\sigma^2\ln^{2(p-1)/p}(K/\beta)}{K^{2(p-1)/p}\mu^2}\notag\\
        &\qquad\qquad\cdot\ln^2\left(\max\left\{2,
        \frac{K^{2(p-1)/p}\mu^2R^2}{\sigma^2\ln^{2(p-1)/p}(K/\beta)}
        \right\}\right)
        \Biggr\}\Biggr),
        \label{ch2:eq:clipped_SEG_str_monotone_case_2_appendix_chapter_2}
    \end{align}
    meaning that \algname{clipped-SEG} requires the following number of iterations/oracle calls to achieve $\|x^{K+1} - x^*\|^2 \leq \varepsilon$ with probability at least $1 - \beta$:
    \begin{align}
        K &= \cO\biggl(\max\biggl\{\frac{L}{\mu}\ln\left(\frac{R^2}{\varepsilon}\right)\ln\left(\frac{L}{\mu \beta}\ln\frac{R^2}{\varepsilon}\right),\notag\\
        &\quad \left(\frac{\sigma^2}{\mu^2\varepsilon}\right)^{p/(2(p-1))}\ln \left(\frac{1}{\beta} \left(\frac{\sigma^2}{\mu^2\varepsilon}\right)^{p/(2(p-1))}\right)\ln^{p/(p-1)}\left(B_\varepsilon\right)\biggr\}\biggr), \label{ch2:eq:clipped_SEG_str_monotone_case_complexity_appendix_chapter_2}
    \end{align}
    where
    \begin{equation*}
        B_\varepsilon = \max\left\{2, \frac{R^2}{\varepsilon \ln \left(\frac{1}{\beta} \left(\frac{\sigma^2}{\mu^2\varepsilon}\right)^{p/(2(p-1))}\right)}\right\}.
    \end{equation*}
\end{theorem}
\begin{proof}
    Again, we will closely follow the proof of Theorem~C.3 from \textcite{gorbunov2022clipped} and the main difference will be reflected in the application of Bernstein inequality and estimating biases and variances of stochastic terms.
    
    Let $R_k = \|x^k - x^*\|$ for all $k\geq 0$. As in the previous results, the proof is based on the induction argument and showing that the iterates do not leave some ball around the solution with high probability. More precisely, for each $k = 0,1,\ldots,K+1$ we consider probability event $E_k$ as follows: inequalities
    \begin{equation}
        R_t^2 \leq 2 \exp(-\gamma\mu t) R^2 \label{ch2:eq:induction_inequality_str_mon_SEG_chapter_2}
    \end{equation}
    hold for $t = 0,1,\ldots,k$ simultaneously. We want to prove $\Prob\{E_k\} \geq  1 - k\beta/(K+1)$ for all $k = 0,1,\ldots,K+1$ by induction. The base of the induction is trivial: for $k=0$ we have $R_0^2 \leq R^2 < 2R^2$ by definition. Next, assume that for $k = T-1 \leq K$ the statement holds: $\Prob\{E_{T-1}\} \geq  1 - (T-1)\beta/(K+1)$. Given this, we need to prove $\Prob\{E_{T}\} \geq  1 - T\beta/(K+1)$. Since $R_t^2 \leq 2\exp(-\gamma\mu t) R^2 \leq 9R^2$, we have $x^t \in B_{3R}(x^*)$, where operator $F$ is $L$-Lipschitz. Thus, $E_{T-1}$ implies
    \begin{align}
        \|F(x^t)\| & \leq L\|x^t - x^*\| \overset{\eqref{ch2:eq:induction_inequality_str_mon_SEG_chapter_2}}{\leq} \sqrt{2}L\exp(- \gamma\mu t/2)R \overset{\eqref{ch2:eq:gamma_SEG_str_mon_chapter_2},\eqref{ch2:eq:lambda_SEG_str_mon_chapter_2}}{\leq} \frac{\lambda_t}{2} \label{ch2:eq:operator_bound_x_t_SEG_str_mon_chapter_2}
    \end{align}
    and
    \begin{align}
        \|\omega_t\|^2
        &\leq 2\|\thesischtwotF(x^t;\xi_1)\|^2 + 2\|F(x^t)\|^2
        \overset{\eqref{ch2:eq:operator_bound_x_t_SEG_str_mon_chapter_2}}{\leq}
        \frac{5}{2}\lambda_t^2 \notag\\
        &\overset{\eqref{ch2:eq:lambda_SEG_str_mon_chapter_2}}{\leq}
        \frac{\exp(-\gamma\mu t)R^2}{4\gamma^2}
        \label{ch2:eq:omega_bound_x_t_SEG_str_mon_chapter_2}
    \end{align}
    for all $t = 0, 1, \ldots, T-1$, where we use that $\|a+b\|^2 \leq 2\|a\|^2 + 2\|b\|^2$ holding for all $a,b \in \R^d$. 
    
    Next, we need to prove that $E_{T-1}$ implies $\|\thesischtwotx^t - x^*\| \leq 3R$ and show several useful inequalities related to $\theta_t$. By Lipschitzness of $F$, probability event $E_{T-1}$ implies
    \begin{align}
        \|\thesischtwotx^t - x^*\|^2 & = \|x^t - x^* - \gamma \thesischtwotF(x^t;\xi_1)\|^2 \notag \\
        & \leq  2\|x^t - x^*\|^2 + 2\gamma^2\|\thesischtwotF(x^t;\xi_1)\|^2 \notag \\
        & \leq 2R_t^2 + 4\gamma^2\|F(x^t)\|^2 + 4\gamma^2\|\omega_t\|^2 \notag \\
        & \overset{\eqref{ch2:eq:L_Lip_chapter_2}}{\leq} 2(1 + 2\gamma^2 L^2)R_t^2 + 4\gamma^2\|\omega_t\|^2 \notag \\
        & \overset{\eqref{ch2:eq:gamma_SEG_str_mon_chapter_2},\eqref{ch2:eq:omega_bound_x_t_SEG_str_mon_chapter_2}}{\leq} 7\exp(-\gamma\mu t) R^2 \leq 9R^2 \label{ch2:eq:tilde_x_distance_SEG_str_mon_chapter_2}
    \end{align}
    and
    \begin{align}
        \|F(\thesischtwotx^t)\| & \leq L\|\thesischtwotx^t - x^*\| \leq \sqrt{7}L\exp(- \gamma\mu t/2)R \overset{\eqref{ch2:eq:gamma_SEG_str_mon_chapter_2},\eqref{ch2:eq:lambda_SEG_str_mon_chapter_2}}{\leq} \frac{\lambda_t}{2} \label{ch2:eq:operator_bound_tx_t_SEG_str_mon_chapter_2}
    \end{align}
    for all $t = 0, 1, \ldots, T-1$. Therefore, $E_{T-1}$ implies that $x^t, \thesischtwotx^t \in B_{3R}(x^*)$ for all $t = 0, 1, \ldots, T-1$. Using Lemma~\ref{ch2:lem:optimization_lemma_str_mon_SEG_chapter_2} and $(1 - \gamma\mu)^T \leq \exp(-\gamma\mu T)$, we obtain that $E_{T-1}$ implies
    \begin{align}
        R_T^2 & \leq \exp(-\gamma\mu T)R^2 - 4\gamma^3 \mu \sum\limits_{l=0}^{T-1} (1-\gamma\mu)^{T-1-l}\langle F(x^l), \omega_l \rangle\notag\\
        &\quad\quad + 2\gamma \sum\limits_{l=0}^{T-1} (1-\gamma\mu)^{T-1-l} \langle x^l - x^* - \gamma F(\thesischtwotx^l), \theta_l \rangle\notag\\
        &\quad\quad + \gamma^2 \sum\limits_{l=0}^{T-1} (1-\gamma\mu)^{T-1-l} \left(\|\theta_l\|^2 + 4\|\omega_l\|^2\right). \notag
    \end{align}
    To handle the sums above, we introduce a new notation:
    \begin{gather}
        \zeta_t = \begin{cases} F(x^t),& \text{if } \|F(x^t)\| \leq \sqrt{2}L\exp(-\gamma\mu t/2)R,\\ 0,& \text{otherwise}, \end{cases} \label{ch2:eq:zeta_t_SEG_str_mon_chapter_2}\\
        u^{t} = \begin{cases} x^t - x^* - \gamma F(\thesischtwotx^t),& \text{if } \|x^t - x^* - \gamma F(\thesischtwotx^t)\|\\ & \qquad \leq \sqrt{7}(1 + \gamma L) \exp(- \gamma\mu t/2)R,\\ 0,& \text{otherwise}, \end{cases} \label{ch2:eq:eta_t_SEG_str_mon_chapter_2}
    \end{gather}
    for $t = 0, 1, \ldots, T-1$. These vectors are bounded almost surely:
    \begin{equation}
        \begin{aligned}
            \|\zeta_t\| &\leq \sqrt{2}L{\exp(-\gamma\mu t/2)}R,\\
            \|u^{t}\| &\leq \sqrt{7}(1 + \gamma L){\exp(-\gamma\mu t/2)}R
        \end{aligned}\label{ch2:eq:zeta_t_eta_t_bound_SEG_str_mon_chapter_2} 
    \end{equation}
    for all $t = 0, 1, \ldots, T-1$. We also notice that, due to \eqref{ch2:eq:operator_bound_x_t_SEG_str_mon_chapter_2}, event $E_{T-1}$ implies the bound $\|F(x^t)\| \leq \sqrt{2}L\exp(-\gamma\mu t/2)R$ and
    \begin{align*}
        \|x^t - x^* - \gamma F(\thesischtwotx^t)\| & \leq \|x^t - x^*\| + \gamma \|F(\thesischtwotx^t)\|\\
        & \overset{\eqref{ch2:eq:tilde_x_distance_SEG_str_mon_chapter_2},\eqref{ch2:eq:operator_bound_tx_t_SEG_str_mon_chapter_2}}{\leq} \sqrt{7}(1 + \gamma L)\exp(-\gamma\mu t/2)R
    \end{align*}
    for $t = 0, 1, \ldots, T-1$. In other words, $E_{T-1}$ implies $\zeta_t = F(x^t)$ and $u^{t} = x^t - x^* - \gamma F(\thesischtwotx^t)$ for all $t = 0,1, \ldots, T-1$, meaning that from $E_{T-1}$ it follows that
    \begin{align}
        R_T^2 & \leq \exp(-\gamma\mu T)R^2 - 4\gamma^3 \mu \sum\limits_{l=0}^{T-1} (1-\gamma\mu)^{T-1-l}\langle \zeta_l, \omega_l \rangle\notag\\
        & \quad\quad + 2\gamma \sum\limits_{l=0}^{T-1} (1-\gamma\mu)^{T-1-l} \langle u^{l}, \theta_l \rangle \notag \\
        &\quad + \gamma^2 \sum\limits_{l=0}^{T-1} (1-\gamma\mu)^{T-1-l} \left(\|\theta_l\|^2 + 4\|\omega_l\|^2\right). \notag
    \end{align}

    To handle the sums appeared on the right-hand side of the previous inequality we consider unbiased and biased parts of $\theta_l, \omega_l$:
    \begin{gather}
        \theta_l^u \eqdef \Exp_{\xi_2^l}\left[\thesischtwotF(\thesischtwotx^l;\xi_2^l)\right] - \thesischtwotF(\thesischtwotx^l;\xi_2^l), \notag\\
        \theta_l^b \eqdef F(\thesischtwotx^l) - \Exp_{\xi_2^l}\left[\thesischtwotF(\thesischtwotx^l;\xi_2^l)\right], \label{ch2:eq:theta_unbias_bias_SEG_str_mon_chapter_2}\\
        \omega_l^u \eqdef \Exp_{\xi_1^l}\left[\thesischtwotF(x^l;\xi_1^l)\right] - \thesischtwotF(x^l;\xi_1^l), \notag\\
        \omega_l^b \eqdef F(x^l) - \Exp_{\xi_1^l}\left[\thesischtwotF(x^l;\xi_1^l)\right], \label{ch2:eq:omega_unbias_bias_SEG_str_mon_chapter_2}
    \end{gather}
    for all $l = 0,\ldots, T-1$. By definition we have $\theta_l = \theta_l^u + \theta_l^b$, $\omega_l = \omega_l^u + \omega_l^b$ for all $l = 0,\ldots, T-1$. Therefore, $E_{T-1}$ implies
    \begin{align}
        R_T^2 & \leq \exp(-\gamma\mu T) R^2 \underbrace{- 4\gamma^3 \mu \sum\limits_{l=0}^{T-1} (1-\gamma\mu)^{T-1-l}\langle \zeta_l, \omega_l^u \rangle}_{\circledOne}\notag\\
        &\quad\quad \underbrace{- 4\gamma^3 \mu \sum\limits_{l=0}^{T-1} (1-\gamma\mu)^{T-1-l}\langle \zeta_l, \omega_l^b \rangle}_{\circledTwo}\notag\\
        &\quad\quad + \underbrace{2\gamma \sum\limits_{l=0}^{T-1} (1-\gamma\mu)^{T-1-l} \langle u^{l}, \theta_l^u \rangle}_{\circledThree} + \underbrace{2\gamma \sum\limits_{l=0}^{T-1} (1-\gamma\mu)^{T-1-l} \langle u^{l}, \theta_l^b \rangle}_{\circledFour} \notag\\
        &\quad\quad + \underbrace{2\gamma^2 \sum\limits_{l=0}^{T-1} (1-\gamma\mu)^{T-1-l} \left(\Exp_{\xi_2^l}\left[\|\theta_l^u\|^2\right] + 4\Exp_{\xi_1^l}\left[\|\omega_l^u\|^2\right]\right)}_{\circledFive} \notag\\
        &\quad\quad + \underbrace{2\gamma^2 \sum\limits_{l=0}^{T-1} (1-\gamma\mu)^{T-1-l} \left(\|\theta_l^u\|^2 + 4\|\omega_l^u\|^2 -\Exp_{\xi_2^l}\left[\|\theta_l^u\|^2\right] - 4\Exp_{\xi_1^l}\left[\|\omega_l^u\|^2\right]\right)}_{\circledSix}\notag\\
        &\quad\quad + \underbrace{2\gamma^2 \sum\limits_{l=0}^{T-1} (1-\gamma\mu)^{T-1-l} \left(\|\theta_l^b\|^2 + 4\|\omega_l^b\|^2\right)}_{\circledSeven}. \label{ch2:eq:SEG_str_mon_1234567_bound_chapter_2}
    \end{align}
    where we also apply inequality $\|a+b\|^2 \leq 2\|a\|^2 + 2\|b\|^2$ holding for all $a,b \in \R^d$ to upper bound $\|\theta_l\|^2$ and $\|\omega_l\|^2$. It remains to derive good enough high-probability upper-bounds for the terms $\circledOne, \circledTwo, \circledThree, \circledFour, \circledFive, \circledSix, \circledSeven$, i.e., to finish our inductive proof we need to show that $\circledOne + \circledTwo + \circledThree + \circledFour + \circledFive + \circledSix + \circledSeven \leq \exp(-\gamma\mu T) R^2$ with high probability. In the subsequent parts of the proof, we will need to use many times the bounds for the norm and second moments of $\theta_{l}^u$, $\theta_{l}^b$, $\omega_l^u$ and $\omega_l^b$. First, by definition of clipping operator, we have with probability $1$ that
    \begin{equation}
        \|\theta_l^u\| \leq 2\lambda_l,\quad \|\omega_l^u\| \leq 2\lambda_l.\label{ch2:eq:theta_omega_magnitude_str_mon_chapter_2}
    \end{equation}
    Moreover, since $E_{T-1}$ implies that $\|F(x^l)\| \leq \lambda_l/2$ and $\|F(\thesischtwotx^l)\| \leq \lambda_l/2$ for all $l = 0,1, \ldots, T-1$ (see \eqref{ch2:eq:operator_bound_x_t_SEG_str_mon_chapter_2} and \eqref{ch2:eq:operator_bound_tx_t_SEG_str_mon_chapter_2}), from Lemma~\ref{ch2:lem:bias_and_variance_clip_chapter_2} we also have that $E_{T-1}$ implies
    \begin{gather}
        \left\|\theta_l^b\right\| \leq \frac{2^{{p}}\sigma^{p}}{\lambda_l^{{p}-1}},\quad \left\|\omega_l^b\right\| \leq \frac{2^{{p}}\sigma^{p}}{\lambda_l^{{p}-1}}, \label{ch2:eq:bias_theta_omega_str_mon_chapter_2}\\
        \Exp_{\xi_2^l}\left[\left\|\theta_l\right\|^2\right] \leq 18 \lambda_l^{2-{p}}\sigma^{p},\quad \Exp_{\xi_1^l}\left[\left\|\omega_l\right\|^2\right] \leq 18 \lambda_l^{2-{p}}\sigma^{p}, \label{ch2:eq:distortion_theta_omega_str_mon_chapter_2}\\
        \Exp_{\xi_2^l}\left[\left\|\theta_l^u\right\|^2\right] \leq 18 \lambda_l^{2-{p}}\sigma^{p},\quad \Exp_{\xi_1^l}\left[\left\|\omega_l^u\right\|^2\right] \leq 18 \lambda_l^{2-{p}}\sigma^{p}, \label{ch2:eq:variance_theta_omega_str_mon_chapter_2}
    \end{gather}
    for all $l = 0,1, \ldots, T-1$.

    \paragraph{Upper bound for $\circledOne$.} By definition of $\omega_{l}^u$, we have $\Exp_{\xi_1^l}[\omega_{l}^u] = 0$ and
    \begin{equation*}
        \Exp_{\xi_1^l}\left[-4\gamma^3\mu (1-\gamma\mu)^{T-1-l} \langle \zeta_l, \omega_l^u \rangle\right] = 0.
    \end{equation*}
    Next, sum $\circledOne$ has bounded with probability $1$ terms:
    \begin{align}
        \MoveEqLeft[1] |-4\gamma^3\mu (1-\gamma\mu)^{T-1-l} \langle \zeta_l, \omega_l^u \rangle | \notag \\
        & \leq 4\gamma^3\mu \exp(-\gamma\mu (T - 1 - l)) \|\zeta_l\|\cdot \|\omega_l^u\|\notag\\
        & \overset{\eqref{ch2:eq:zeta_t_eta_t_bound_SEG_str_mon_chapter_2},\eqref{ch2:eq:theta_omega_magnitude_str_mon_chapter_2}}{\leq} 8\sqrt{2}\gamma^3\mu L \exp(-\gamma\mu (T - 1 - l/2)) R \lambda_l\notag\\
        & \overset{\eqref{ch2:eq:gamma_SEG_str_mon_chapter_2},\eqref{ch2:eq:lambda_SEG_str_mon_chapter_2}}{\leq} \frac{\exp(-\gamma\mu T)R^2}{7A}  \eqdef c. \label{ch2:eq:SEG_str_mon_technical_1_1_chapter_2}
    \end{align}
    The summands also have bounded conditional variances
    \begin{equation*}
        \sigma_l^2 \eqdef
        \Exp_{\xi_1^l}\left[
        16\gamma^6\mu^2 (1-\gamma\mu)^{2T-2-2l}
        \langle \zeta_l, \omega_l^u \rangle^2
        \right],
    \end{equation*}
    which satisfy
    \begin{align}
        \sigma_l^2 & \leq \Exp_{\xi_1^l}\left[16\gamma^6\mu^2 \exp(-\gamma\mu (2T - 2 - 2l)) \|\zeta_l\|^2\cdot \|\omega_l^u\|^2\right]\notag\\
        & \overset{\eqref{ch2:eq:zeta_t_eta_t_bound_SEG_str_mon_chapter_2}}{\leq} 36\gamma^6 \mu^2L^2 \exp(-\gamma\mu (2T - 2 - l)) R^2 \Exp_{\xi_1^l}\left[\|\omega_l^u\|^2\right]\notag\\
        & \overset{\eqref{ch2:eq:gamma_SEG_str_mon_chapter_2}}{\leq} \frac{4\gamma^2\exp(-\gamma\mu(2T - l))R^2}{2809A} \Exp_{\xi_1^l}\left[\|\omega_l^u\|^2\right]. \label{ch2:eq:SEG_str_mon_technical_1_2_chapter_2}
    \end{align}
    In other words, we showed that the summands of $\circledOne$ form a bounded martingale difference sequence with bounded conditional variances $\{\sigma_l^2\}_{l = 0}^{T-1}$. Next, we apply Bernstein's inequality (Lemma~\ref{ch2:lem:Bernstein_ineq_chapter_2}) with $X_l$ equal to the summands of $\circledOne$, parameter $c$ as in \eqref{ch2:eq:SEG_str_mon_technical_1_1_chapter_2}, $b = \tfrac{1}{7}\exp(-\gamma\mu T) R^2$, $H = \tfrac{\exp(-2 \gamma\mu T) R^4}{294A}$:
    \begin{align*}
        \MoveEqLeft \Prob\left\{|\circledOne| > \frac{1}{7}\exp(-\gamma\mu T) R^2 \text{ and } \sum\limits_{l=0}^{T-1}\sigma_l^2 \leq \frac{\exp(-2\gamma\mu T) R^4}{294A}\right\}\\
        &\leq 2\exp\left(- \frac{b^2}{2H + 2cb/3}\right) = \frac{\beta}{3(K+1)}.
    \end{align*}
    Equivalently, we have $\Prob\{E_{\circledOne}\} \geq 1 - \frac{\beta}{3(K+1)}$ for
    \begin{equation}
            E_{\circledOne} =
            \left\{\text{either } \sum\limits_{l=0}^{T-1}\sigma_l^2 > \frac{\exp(-2\gamma\mu T) R^4}{294A} \text{ or } |\circledOne| \leq \frac{1}{7}\exp(-\gamma\mu T) R^2\right\}.
        \label{ch2:eq:bound_1_SEG_str_mon_chapter_2}
    \end{equation}
    In addition, $E_{T-1}$ implies that
    \begin{align}
        \sum\limits_{l=0}^{T-1}\sigma_l^2 & \overset{\eqref{ch2:eq:SEG_str_mon_technical_1_2_chapter_2}}{\leq} \frac{4\gamma^2\exp(-2\gamma\mu T)R^2}{2809A} \sum\limits_{l=0}^{T-1} \frac{\Exp_{\xi_1^l}\left[\|\omega_l^u\|^2\right]}{\exp(-\gamma\mu l)}\notag\\
        & \overset{\eqref{ch2:eq:variance_theta_omega_str_mon_chapter_2}, T \leq K+1}{\leq} \frac{72\gamma^2\exp(-2\gamma\mu T) R^2 \sigma^{p}}{2809A} \sum\limits_{l=0}^{K} \frac{\lambda_{l}^{2-{p}}}{\exp(-\gamma\mu l)}\notag\\
        & \overset{\eqref{ch2:eq:lambda_SEG_str_mon_chapter_2}}{\leq} \frac{72\gamma^{p}\exp(-2\gamma\mu T) R^{4-{p}} \sigma^{p}}{2809\cdot 120^{2-{p}} A^{3-{p}}} \sum\limits_{l=0}^{K} \frac{1}{\exp(-\gamma\mu l)} \cdot \left(\exp(-\gamma\mu(1 + l/2))\right)^{2-{p}} \notag\\
        & \leq \frac{72\gamma^{p}\exp(-2\gamma\mu T) R^{4-{p}} \sigma^{p}}{2809\cdot 120^{2-{p}} A^{3-{p}}} \sum\limits_{l=0}^{K} \exp(\gamma\mu({p}-2)) \cdot \exp\left(\frac{\gamma\mu{p} l}{2}\right) \notag\\
        & \leq \frac{72\gamma^{p}\exp(-2\gamma\mu T) R^{4-{p}} \sigma^{p} (K+1) \exp\left(\frac{\gamma\mu{p} K}{2}\right)}{2809\cdot 120^{2-{p}} A^{3-{p}}}  \notag\\
        & \overset{\eqref{ch2:eq:gamma_SEG_str_mon_chapter_2}}{\leq} \frac{\exp(-2\gamma\mu T)R^4}{294A}, \label{ch2:eq:bound_1_variances_SEG_str_mon_chapter_2}
    \end{align}
    where we also show that $E_{T-1}$ implies
    \begin{equation}
        \gamma^2 R^2\sum\limits_{l=0}^{K} \frac{\lambda_{l}^{2-{p}}}{\exp(-\gamma\mu l)} \leq \frac{\gamma^{p} R^{4-{p}}(K+1)\exp(\frac{\gamma\mu{p} K}{2})}{120^{2-{p}}A^{2-{p}}}. \label{ch2:eq:useful_inequality_on_lambda_SEG_str_mon_chapter_2}
    \end{equation}

    \textbf{Upper bound for $\circledTwo$.} From $E_{T-1}$ it follows that
    \begin{align}
        \circledTwo
        &\leq 4\gamma^3 \mu \sum\limits_{l=0}^{T-1}
        \exp(-\gamma\mu (T-1-l)) \|\zeta_l\| \cdot \|\omega_l^b\|
        \notag\\
        &\overset{\eqref{ch2:eq:zeta_t_eta_t_bound_SEG_str_mon_chapter_2},
        \eqref{ch2:eq:bias_theta_omega_str_mon_chapter_2}}{\leq}
        \begin{multlined}[t]
        2^{2+{p}}\cdot\sqrt{2} \exp(-\gamma\mu (T-1))
        \gamma^3 \mu L R\\
        {}\sum\limits_{l=0}^{T-1}
        \frac{\sigma^{p}}{\lambda_l^{{p}-1} \exp(-\gamma\mu l/2)}
        \end{multlined}
        \notag\\
        &\overset{\eqref{ch2:eq:lambda_SEG_str_mon_chapter_2}}{=}
        \begin{multlined}[t]
        \frac{2^{2+{p}}\cdot 120^{{p}-1}\sqrt{2}
        \exp(-\gamma\mu (T-1)) \gamma^{2+{p}} \mu L \sigma^{p}
        A^{{p}-1}}{R^{{p}-2}}\\
        {}\sum\limits_{l=0}^{T-1}
        \frac{1}{\exp\left(-\gamma\mu(1 + l/2)\right)^{{p}-1}
        \cdot\exp(-\gamma\mu l/2)}
        \end{multlined}
        \notag\\
        &\overset{T \leq K+1}{\leq}
        \begin{multlined}[t]
        \frac{2^{3+{p}}\cdot 120^{{p}-1}\sqrt{2}
        \exp(-\gamma\mu (T-1)) \gamma^{2+{p}} \mu L \sigma^{p}
        A^{{p}-1}}{R^{{p}-2}}\\
        {}\sum\limits_{l=0}^{K}
        \exp\left(\frac{\gamma\mu {p} l}{2}\right)
        \end{multlined}
        \notag\\
        &\leq
        \frac{2^{3+{p}}\cdot 120^{{p}-1}\sqrt{2}
            \exp(-\gamma\mu (T-1)) \gamma^{2+{p}} \mu L \sigma^{p}}{R^{{p}-2}}
        \notag\\
        &\qquad \cdot A^{{p}-1}
            (K+1)\exp\left(\frac{\gamma\mu {p} K}{2}\right)
        \notag\\
        &\overset{\eqref{ch2:eq:gamma_SEG_str_mon_chapter_2}}{\leq}
        \frac{1}{7}\exp(-\gamma\mu T) R^2,
        \label{ch2:eq:bound_2_SEG_str_mon_chapter_2}
    \end{align}
    where we also show that $E_{T-1}$ implies
    \begin{equation}
        \gamma R\sum\limits_{l=0}^{T-1} \frac{1}{\lambda_l^{{p}-1}\exp(-\gamma\mu l/2)} \leq \frac{120^{{p}-1}\gamma^{p} (K+1)\exp(\frac{\gamma\mu {p} K}{2})A^{{p}-1}}{R^{{p}-2}}. \label{ch2:eq:useful_inequality_on_lambda_SEG_str_mon_1_chapter_2}
    \end{equation}

    \paragraph{Upper bound for $\circledThree$.} By definition of $\theta_l^u$, we have $\Exp_{\xi_2^l}[\theta_{l}^u] = 0$ and
    \begin{equation*}
        \Exp_{\xi_2^l}\left[2\gamma (1-\gamma\mu)^{T-1-l} \langle u^{l}, \theta_l^u \rangle\right] = 0.
    \end{equation*}
    Next, sum $\circledThree$ has bounded with probability $1$ terms:
    \begin{align}
        \MoveEqLeft[1] |2\gamma (1-\gamma\mu)^{T-1-l} \langle u^{l}, \theta_l^u \rangle | \notag \\
        & \leq 2\gamma\exp(-\gamma\mu (T - 1 - l)) \|u^{l}\|\cdot \|\theta_l^u\|\notag\\
        & \overset{\eqref{ch2:eq:zeta_t_eta_t_bound_SEG_str_mon_chapter_2},\eqref{ch2:eq:theta_omega_magnitude_str_mon_chapter_2}}{\leq} 4\sqrt{7}\gamma (1 + \gamma L) \exp(-\gamma\mu (T - 1 - l/2)) R \lambda_l\notag\\
        & \overset{\eqref{ch2:eq:gamma_SEG_str_mon_chapter_2},\eqref{ch2:eq:lambda_SEG_str_mon_chapter_2}}{\leq} \frac{\exp(-\gamma\mu T)R^2}{7A} \eqdef c. \label{ch2:eq:SEG_str_mon_technical_3_1_chapter_2}
    \end{align}
    The summands also have bounded conditional variances $\widetilde\sigma_l^2$ (the conditional expectations $\Exp_{\xi_2^l}$ of their squares):
    \begin{align}
        \widetilde\sigma_l^2 & \leq \Exp_{\xi_2^l}\left[4\gamma^2\exp(-\gamma\mu (2T - 2 - 2l)) \|u^{l}\|^2\cdot \|\theta_l^u\|^2\right]\notag\\
        & \overset{\eqref{ch2:eq:zeta_t_eta_t_bound_SEG_str_mon_chapter_2}}{\leq} 49\gamma^2 (1 + \gamma L)^2 \exp(-\gamma\mu (2T - 2 - l)) R^2 \Exp_{\xi_2^l}\left[\|\theta_l^u\|^2\right]\notag\\
        & \overset{\eqref{ch2:eq:gamma_SEG_str_mon_chapter_2}}{\leq} 50\gamma^2\exp(-\gamma\mu (2T - l))R^2 \Exp_{\xi_2^l}\left[\|\theta_l^u\|^2\right]. \label{ch2:eq:SEG_str_mon_technical_3_2_chapter_2}
    \end{align}
    In other words, we showed that the summands of $\circledThree$ form a bounded martingale difference sequence with bounded conditional variances $\{\widetilde\sigma_l^2\}_{l = 0}^{T-1}$. Next, we apply Bernstein's inequality (Lemma~\ref{ch2:lem:Bernstein_ineq_chapter_2}) with $X_l$ equal to the summands of $\circledThree$, parameter $c$ as in \eqref{ch2:eq:SEG_str_mon_technical_3_1_chapter_2}, $b = \tfrac{1}{7}\exp(-\gamma\mu T) R^2$, $H = \tfrac{\exp(-2 \gamma\mu T) R^4}{294A}$:
    \begin{align*}
        &\Prob\left\{
        |\circledThree| > \frac{1}{7}\exp(-\gamma\mu T) R^2
        \ \text{and}\
        \sum\limits_{l=0}^{T-1}\widetilde\sigma_l^2
        \leq \frac{\exp(- 2\gamma\mu T) R^4}{294A}
        \right\}\\
        &\qquad \leq
        2\exp\left(- \frac{b^2}{2H + 2cb/3}\right)
        = \frac{\beta}{3(K+1)}.
    \end{align*}
    Equivalently, we have $\Prob\{E_{\circledThree}\} \geq 1 - \frac{\beta}{3(K+1)}$ for
    \begin{equation}
        E_{\circledThree} = \left\{\text{either } \sum\limits_{l=0}^{T-1}\widetilde\sigma_l^2 > \frac{\exp(- 2\gamma\mu T) R^4}{294A} \text{ or } |\circledThree| \leq \frac{1}{7}{\exp(-\gamma\mu T)} R^2\right\}. \label{ch2:eq:bound_3_SEG_str_mon_chapter_2}
    \end{equation}
    In addition, $E_{T-1}$ implies that
    \begin{align}
        \sum\limits_{l=0}^{T-1}\widetilde\sigma_l^2 & \overset{\eqref{ch2:eq:SEG_str_mon_technical_3_2_chapter_2}}{\leq} 50\gamma^2\exp(- 2\gamma\mu T)R^2\sum\limits_{l=0}^{T-1} \frac{\Exp_{\xi_2^l}\left[\|\theta_l^u\|^2\right]}{\exp(-\gamma\mu l)}\notag\\ 
        & \overset{\eqref{ch2:eq:variance_theta_omega_str_mon_chapter_2}, T \leq K+1}{\leq} 900\gamma^2\exp(-2\gamma\mu T) R^2 \sigma^{p} \sum\limits_{l=0}^{K} \frac{\lambda_l^{2-{p}}}{\exp(-\gamma\mu l)}\notag\\
        & \overset{\eqref{ch2:eq:useful_inequality_on_lambda_SEG_str_mon_chapter_2}}{\leq} \frac{900\gamma^{p}\exp(-2\gamma\mu T) R^{4-{p}} \sigma^{p} (K+1)\exp(\frac{\gamma\mu{p} K}{2})}{120^{2-{p}} A^{2-{p}}}\notag\\
        & \overset{\eqref{ch2:eq:gamma_SEG_str_mon_chapter_2}}{\leq} \frac{\exp(-2\gamma\mu T)R^4}{294A}. \label{ch2:eq:bound_3_variances_SEG_str_mon_chapter_2}
    \end{align}

    \paragraph{Upper bound for $\circledFour$.} From $E_{T-1}$ it follows that
    \begin{align}
        \circledFour & \leq 2\gamma \exp(-\gamma\mu (T-1)) \sum\limits_{l=0}^{T-1} \frac{\|u^{l}\|\cdot \|\theta_l^b\|}{\exp(-\gamma\mu l)}\notag\\
        & \overset{\eqref{ch2:eq:zeta_t_eta_t_bound_SEG_str_mon_chapter_2}, \eqref{ch2:eq:bias_theta_omega_str_mon_chapter_2}}{\leq} 2^{1+{p}}\sqrt{7} \gamma (1+\gamma L) \exp(-\gamma\mu (T-1)) R \sigma^{p} \sum\limits_{l=0}^{T-1} \frac{1}{\lambda_l^{{p}-1} \exp(-\gamma\mu l/2)}\notag\\
        & \overset{\eqref{ch2:eq:useful_inequality_on_lambda_SEG_str_mon_1_chapter_2}}{\leq} \frac{2^{3+{p}}\cdot 120^{{p}-1}\sqrt{7} \gamma^{p}(1+\gamma L) \exp(-\gamma\mu T) (K+1) \exp\left(\frac{\gamma\mu{p} K}{2}\right) A^{{p}-1}}{R^{{p}-2}} \notag \\
        & \overset{\eqref{ch2:eq:gamma_SEG_str_mon_chapter_2}}{\leq} \frac{1}{7}\exp(-\gamma\mu T) R^2. \label{ch2:eq:bound_4_SEG_str_mon_chapter_2}
    \end{align}

    \paragraph{Upper bound for $\circledFive$.} From $E_{T-1}$ it follows that
    \begin{align}
        \circledFive & = 2\gamma^2 \exp(-\gamma\mu (T-1)) \sum\limits_{l=0}^{T-1} \frac{\Exp_{\xi_2^l}\left[\|\theta_l^u\|^2\right] + 4\Exp_{\xi_1^l}\left[\|\omega_l^u\|^2\right]}{\exp(-\gamma\mu l)} \notag\\
        & \overset{\eqref{ch2:eq:variance_theta_omega_str_mon_chapter_2}}{\leq} 180\gamma^2\exp(-\gamma\mu (T-1)) \sigma^{p}\sum\limits_{l=0}^{T-1} \frac{\lambda_l^{2-{p}}}{\exp(-\gamma\mu l)} \notag\\
        & \overset{\eqref{ch2:eq:useful_inequality_on_lambda_SEG_str_mon_chapter_2}}{\leq}  \frac{180\gamma^{p} R^{2-{p}}\exp(-\gamma\mu (T-1)) \sigma^{p} (K+1)\exp(\frac{\gamma\mu{p} K}{2})}{120^{2-{p}} A^{2-{p}}} \notag\\
        & \overset{\eqref{ch2:eq:gamma_SEG_str_mon_chapter_2}}{\leq} \frac{1}{7} \exp(-\gamma\mu T) R^2. \label{ch2:eq:bound_5_SEG_str_mon_chapter_2}
    \end{align}

    \paragraph{Upper bound for $\circledSix$.} First, we have
    \begin{equation*}
        2\gamma^2 (1-\gamma\mu)^{T-1-l}\Exp_{\xi_1^l, \xi_2^l}\left[\|\theta_l^u\|^2 + 4\|\omega_l^u\|^2 -\Exp_{\xi_2^l}\left[\|\theta_l^u\|^2\right] - 4\Exp_{\xi_1^l}\left[\|\omega_l^u\|^2\right]\right] = 0.
    \end{equation*}
    Next, sum $\circledSix$ has bounded with probability $1$ terms:
    \begin{align}
        \MoveEqLeft[1] 2\gamma^2 (1-\gamma\mu)^{T-1-l}\left| \|\theta_l^u\|^2 + 4\|\omega_l^u\|^2 -\Exp_{\xi_2^l}\left[\|\theta_l^u\|^2\right] - 4\Exp_{\xi_1^l}\left[\|\omega_l^u\|^2\right] \right| \notag \\
        & \overset{\eqref{ch2:eq:theta_omega_magnitude_str_mon_chapter_2}}{\leq} \frac{80\gamma^2 \exp(-\gamma\mu T) \lambda_l^2}{\exp(-\gamma\mu (1+l))}\notag\\
        & \overset{\eqref{ch2:eq:lambda_SEG_str_mon_chapter_2}}{\leq} \frac{\exp(-\gamma\mu T)R^2}{7A}\notag\\
        & \eqdef c. \label{ch2:eq:SEG_str_mon_technical_6_1_chapter_2}
    \end{align}
    The summands also have conditional variances
    \begin{equation*}
        \widehat\sigma_l^2 \eqdef \Exp_{\xi_1^l,\xi_2^l}\left[4\gamma^4 (1-\gamma\mu)^{2T-2-2l} \left| \|\theta_l^u\|^2 + 4\|\omega_l^u\|^2 -\Exp_{\xi_2^l}\left[\|\theta_l^u\|^2\right] - 4\Exp_{\xi_1^l}\left[\|\omega_l^u\|^2\right] \right|^2\right]
    \end{equation*}
    that are bounded
    \begin{align}
        \widehat\sigma_l^2
        &\overset{\eqref{ch2:eq:SEG_str_mon_technical_6_1_chapter_2}}{\leq}
        \frac{2\gamma^2\exp(-2\gamma\mu T)R^2}{7\exp(-\gamma\mu (1+l))A}
        \notag\\[-0.25em]
        &\qquad {}\times
        \Exp_{\xi_1^l,\xi_2^l}\left[
        \left|\|\theta_l^u\|^2 + 4\|\omega_l^u\|^2
        -\Exp_{\xi_2^l}\left[\|\theta_l^u\|^2\right] - 4\Exp_{\xi_1^l}\left[\|\omega_l^u\|^2\right]\right|\right]
        \notag\\
        &\leq
        \frac{4\gamma^2\exp(-2\gamma\mu T)R^2}{7\exp(-\gamma\mu (1+l))A}
        \Exp_{\xi_1^l,\xi_2^l}\left[
        \|\theta_l^u\|^2 + 4\|\omega_l^u\|^2\right].
        \label{ch2:eq:SEG_str_mon_technical_6_2_chapter_2}
    \end{align}
    In other words, the summands of $\circledSix$ form a bounded martingale difference sequence with bounded conditional variances $\{\widehat\sigma_l^2\}_{l = 0}^{T-1}$. Next, we apply Bernstein's inequality (Lemma~\ref{ch2:lem:Bernstein_ineq_chapter_2}) with $X_l$ equal to the summands of $\circledSix$, $b = \frac{1}{7}\exp(-\gamma\mu T)R^2$, $H = \frac{\exp(-2\gamma\mu T)R^4}{294A}$,
    and parameter $c$ as in
    \eqref{ch2:eq:SEG_str_mon_technical_6_1_chapter_2}:
    \begin{align*}
        &\Prob\left\{
        |\circledSix| > \frac{1}{7}\exp(-\gamma\mu T) R^2
        \ \text{and}\
        \sum\limits_{l=0}^{T-1}\widehat\sigma_l^2
        \leq \frac{\exp(-2\gamma\mu T) R^4}{294A}
        \right\}\\
        &\qquad \leq
        2\exp\left(- \frac{b^2}{2H + 2cb/3}\right)
        = \frac{\beta}{3(K+1)}.
    \end{align*}
    Equivalently, we have $\Prob\{E_{\circledSix}\} \geq 1 - \frac{\beta}{3(K+1)}$ for
    \begin{equation}
        E_{\circledSix} = \left\{\text{either } \sum\limits_{l=0}^{T-1}\widehat\sigma_l^2 > \frac{\exp(-2\gamma\mu T) R^4}{294A} \text{ or } |\circledSix| \leq \frac{1}{7}{\exp(-\gamma\mu T)} R^2\right\}. \label{ch2:eq:bound_6_SEG_str_mon_chapter_2}
    \end{equation}
    In addition, $E_{T-1}$ implies that
    \begin{align}
        \sum\limits_{l=0}^{T-1}\widehat\sigma_l^2 & \overset{\eqref{ch2:eq:SEG_str_mon_technical_6_2_chapter_2}}{\leq} \frac{4\gamma^2\exp(-\gamma\mu (2T-1))R^2}{7A} \sum\limits_{l=0}^{T-1} \frac{\Exp_{\xi_1^l,\xi_2^l}\left[\|\theta_l^u\|^2 + 4\|\omega_l^u\|^2\right]}{\exp(-\gamma\mu l)}\notag\\
        & \overset{\eqref{ch2:eq:variance_theta_omega_str_mon_chapter_2}, T \leq K+1}{\leq} \frac{360\gamma^2\exp(-\gamma\mu (2T-1)) R^2 \sigma^{p}}{7A} \sum\limits_{l=0}^{K} \frac{\lambda_l^{2-{p}}}{\exp(-\gamma\mu l)}\notag\\
        & \overset{\eqref{ch2:eq:useful_inequality_on_lambda_SEG_str_mon_chapter_2}}{\leq} \frac{360\gamma^{p}\exp(-\gamma\mu (2T-1)) R^{4-{p}} \sigma^{p} (K+1)\exp(\frac{\gamma\mu{p} K}{2})}{7\cdot 120^{2-{p}}A^{3-{p}}} \notag\\
        & \overset{\eqref{ch2:eq:gamma_SEG_str_mon_chapter_2}}{\leq} \frac{\exp(-2\gamma\mu T)R^4}{294A}. \label{ch2:eq:bound_6_variances_SEG_str_mon_chapter_2}
    \end{align}

    \paragraph{Upper bound for $\circledSeven$.} From $E_{T-1}$ it follows that
    \begin{align}
        \circledSeven
        &= 2\gamma^2 \sum\limits_{l=0}^{T-1}
        \exp(-\gamma\mu (T-1-l))
        \left(\|\theta_l^b\|^2 + 4\|\omega_l^b\|^2\right)
        \notag\\
        &\overset{\eqref{ch2:eq:bias_theta_omega_str_mon_chapter_2}}{\leq}
        10\cdot 2^{2{p}}\gamma^2 \exp(-\gamma\mu (T-1)) \sigma^{2{p}}
        \sum\limits_{l=0}^{T-1}
        \frac{1}{\lambda_l^{2{p}-2} \exp(-\gamma\mu l)}
        \notag\\
        &\overset{\eqref{ch2:eq:lambda_SEG_str_mon_chapter_2}, T \leq K+1}{\leq}
        \frac{20\cdot 2^{2{p}}\cdot 120^{2{p}-2}\gamma^{2{p}}
        \exp(-\gamma\mu T) \sigma^{2{p}}}
        {R^{2{p}-2}}
        \notag\\
        &\qquad {}\times A^{2{p}-2}
        \sum\limits_{l=0}^{K}
        \exp\left(\gamma\mu(2{p}-2)\left(1 + \frac{l}{2}\right)\right)
        \exp(\gamma\mu l)
        \notag\\
        &\leq
        \frac{40\cdot 2^{2{p}}\cdot 120^{2{p}-2}\gamma^{2{p}}
        \exp(-\gamma\mu T) \sigma^{2{p}}}
        {R^{2{p}-2}}
        A^{2{p}-2}
        \sum\limits_{l=0}^{K} \exp(\gamma\mu {p} l)
        \notag\\
        &\leq
        \frac{40\cdot 2^{2{p}}\cdot 120^{2{p}-2}\gamma^{2{p}}
        \exp(-\gamma\mu T) \sigma^{2{p}}}
        {R^{2{p}-2}}
        \notag\\
        &\qquad {}\times A^{2{p}-2}
        (K+1) \exp(\gamma\mu {p} K)
        \notag\\
        &\overset{\eqref{ch2:eq:gamma_SEG_str_mon_chapter_2}}{\leq}
        \frac{1}{7}\exp(-\gamma\mu T) R^2.
        \label{ch2:eq:bound_7_SEG_str_mon_chapter_2}
    \end{align}

    Now, we have the upper bounds for  $\circledOne, \circledTwo, \circledThree, \circledFour, \circledFive, \circledSix, \circledSeven$. In particular, probability event $E_{T-1}$ implies
    \begin{gather*}
        R_T^2 \overset{\eqref{ch2:eq:SEG_str_mon_1234567_bound_chapter_2}}{\leq} \exp(-\gamma\mu T) R^2 + \circledOne + \circledTwo + \circledThree + \circledFour + \circledFive + \circledSix + \circledSeven,\\
        \circledTwo \overset{\eqref{ch2:eq:bound_2_SEG_str_mon_chapter_2}}{\leq} \frac{1}{7}\exp(-\gamma\mu T)R^2,\\
        \circledFour \overset{\eqref{ch2:eq:bound_4_SEG_str_mon_chapter_2}}{\leq} \frac{1}{7}\exp(-\gamma\mu T)R^2,\\
        \circledFive \overset{\eqref{ch2:eq:bound_5_SEG_str_mon_chapter_2}}{\leq} \frac{1}{7}\exp(-\gamma\mu T)R^2,\\
        \circledSeven \overset{\eqref{ch2:eq:bound_7_SEG_str_mon_chapter_2}}{\leq} \frac{1}{7}\exp(-\gamma\mu T)R^2,\\
        \sum\limits_{l=0}^{T-1}\sigma_l^2 \overset{\eqref{ch2:eq:bound_1_variances_SEG_str_mon_chapter_2}}{\leq}  \frac{\exp(-2\gamma\mu T)R^4}{294A},\\
        \sum\limits_{l=0}^{T-1}\widetilde\sigma_l^2 \overset{\eqref{ch2:eq:bound_3_variances_SEG_str_mon_chapter_2}}{\leq} \frac{\exp(-2\gamma\mu T)R^4}{294A},\\
        \sum\limits_{l=0}^{T-1}\widehat\sigma_l^2 \overset{\eqref{ch2:eq:bound_6_variances_SEG_str_mon_chapter_2}}{\leq}  \frac{\exp(-2\gamma\mu T)R^4}{294A}.
    \end{gather*}
     Moreover, we also have (see \eqref{ch2:eq:bound_1_SEG_str_mon_chapter_2}, \eqref{ch2:eq:bound_3_SEG_str_mon_chapter_2}, \eqref{ch2:eq:bound_6_SEG_str_mon_chapter_2} and our induction assumption)
     \begin{gather*}
        \Prob\{E_{T-1}\} \geq 1 - \frac{(T-1)\beta}{K+1},\\
        \Prob\{E_{\circledOne}\} \geq 1 - \frac{\beta}{3(K+1)}, \quad \Prob\{E_{\circledThree}\} \geq 1 - \frac{\beta}{3(K+1)}, \quad \Prob\{E_{\circledSix}\} \geq 1 - \frac{\beta}{3(K+1)} ,
    \end{gather*}
    where
    \begin{align}
        E_{\circledOne}& = \left\{\text{either } \sum\limits_{l=0}^{T-1}\sigma_l^2 > \frac{\exp(-2\gamma\mu T) R^4}{294A} \text{ or } |\circledOne| \leq \frac{1}{7}{\exp(-\gamma\mu T)} R^2\right\},\notag\\
        E_{\circledThree}& = \left\{\text{either } \sum\limits_{l=0}^{T-1}\widetilde\sigma_l^2 > \frac{\exp(-2\gamma\mu T) R^4}{294A} \text{ or } |\circledThree| \leq \frac{1}{7}{\exp(-\gamma\mu T)} R^2\right\},\notag\\
        E_{\circledSix}& = \left\{\text{either } \sum\limits_{l=0}^{T-1}\widehat\sigma_l^2 > \frac{\exp(-2\gamma\mu T) R^4}{294A} \text{ or } |\circledSix| \leq \frac{1}{7}{\exp(-\gamma\mu T)} R^2\right\}.\notag
    \end{align}
    Thus, probability event $E_{T-1} \cap E_{\circledOne} \cap E_{\circledThree} \cap E_{\circledSix}$ implies
    \begin{align*}
        R_T^2 & \overset{\eqref{ch2:eq:SEG_str_mon_1234567_bound_chapter_2}}{\leq} \exp(-\gamma\mu T) R^2 + \circledOne + \circledTwo + \circledThree + \circledFour + \circledFive + \circledSix + \circledSeven\\
        & \leq 2\exp(-\gamma\mu T) R^2,
    \end{align*}
    which is equivalent to \eqref{ch2:eq:induction_inequality_str_mon_SEG_chapter_2} for $t = T$, and
    \begin{align}
\Prob\{E_T\} & \geq \Prob\{E_{T-1} \cap E_{\circledOne} \cap E_{\circledThree} \cap E_{\circledSix}\} \notag \\
& = 1 - \Prob\{\overline{E}_{T-1} \cup \overline{E}_{\circledOne} \cup \overline{E}_{\circledThree} \cup \overline{E}_{\circledSix}\} \geq 1 - \frac{T\beta}{K+1}. \notag
    \end{align}
    This finishes the inductive part of our proof, i.e., for all $k = 0,1,\ldots,K+1$ we have $\Prob\{E_k\} \geq 1 - k\beta/(K+1)$. In particular, for $k = K+1$ we have that with probability at least $1 - \beta$
    \begin{equation}
        \|x^{K+1} - x^*\|^2 \leq 2\exp(-\gamma\mu (K+1))R^2. \notag
    \end{equation}
    Finally, if
    \begin{align*}
        \gamma & = \min\left\{\frac{1}{650 L A}, \frac{\ln(B_K)}{\mu(K+1)}\right\}, \notag\\
        B_K & = \max\left\{2, \frac{(K+1)^{2(p-1)/p}\mu^2R^2}{264600^{2/p}\sigma^2A^{2(p-1)/p}\ln^2(B_K)} \right\}  \\
        & = \cO\biggl(\max\biggl\{2, \frac{K^{2(p-1)/p}\mu^2R^2}{\sigma^2\ln^{2(p-1)/p}(K/\beta)}\notag\\
        &\qquad\qquad \cdot \frac{1}{\ln^2\left(\max\left\{2, \frac{K^{2(p-1)/p}\mu^2R^2}{\sigma^2\ln^{2(p-1)/p}(K/\beta)} \right\}\right)} \biggr\}\biggr),
    \end{align*}
    then with probability at least $1-\beta$
    \begin{align*}
        \MoveEqLeft \|x^{K+1} - x^*\|^2
        \leq 2\exp(-\gamma\mu (K+1))R^2\\
        &= 2R^2\max\left\{
        \exp\left(-\frac{\mu(K+1)}{650 L A}\right),
        \frac{1}{B_K} \right\}\\
        &= \cO\Biggl(\max\Biggl\{
        R^2\exp\left(- \frac{\mu K}{L \ln(K/\beta)}\right),
        \frac{\sigma^2\ln^{2(p-1)/p}(K/\beta)}{K^{2(p-1)/p}\mu^2}\\
        &\qquad\qquad\cdot\ln^2\left(\max\left\{2,
        \frac{K^{2(p-1)/p}\mu^2R^2}{\sigma^2\ln^{2(p-1)/p}(K/\beta)}
        \right\}\right)
        \Biggr\}\Biggr).
    \end{align*}
    To get $\|x^{K+1} - x^*\|^2 \leq \varepsilon$ with probability at least $1-\beta$ it is sufficient to choose $K$ such that both terms in the maximum above are $\cO(\varepsilon)$. This leads to
    \begin{align*}
        K &= \cO\biggl(\max\biggl\{\frac{L}{\mu}\ln\left(\frac{R^2}{\varepsilon}\right)\ln\left(\frac{L}{\mu \beta}\ln\frac{R^2}{\varepsilon}\right),\notag\\
        &\quad \left(\frac{\sigma^2}{\mu^2\varepsilon}\right)^{p/(2(p-1))}\ln \left(\frac{1}{\beta} \left(\frac{\sigma^2}{\mu^2\varepsilon}\right)^{p/(2(p-1))}\right)\ln^{p/(p-1)}\left(B_\varepsilon\right)\biggr\}\biggr),
    \end{align*}
    where
    \begin{equation*}
        B_\varepsilon = \max\left\{2, \frac{R^2}{\varepsilon \ln \left(\frac{1}{\beta} \left(\frac{\sigma^2}{\mu^2\varepsilon}\right)^{p/(2(p-1))}\right)}\right\}.
    \end{equation*}
    This concludes the proof.
\end{proof}

\clearpage

\section{Additional Proofs for \algname{clipped-SGDA}}\label{ch2:appendix:SGDA_chapter_2}

In this section, we provide the complete formulation of the main results for \algname{clipped-SGDA} and the corresponding proofs. For brevity, we will use the following notation: $\thesischtwotF(x^{k};\xi^k) = \clip\left(F(x^{k};\xi^k), \lambda_k\right)$.

\begin{algorithm}[h]
\caption{Clipped Stochastic Gradient Descent Ascent (\algname{clipped-SGDA}) \parencite{gorbunov2022clipped}}
\label{ch2:alg:clipped-SGDA_chapter_2}   
\begin{algorithmic}[1]
\REQUIRE starting point $x^0$, number of iterations $K$, stepsize $\gamma > 0$, clipping levels $\{\lambda_k\}_{k=0}^{K}$.
\FOR{$k=0,\ldots, K$}
\STATE Compute $\thesischtwotF(x^{k};\xi^k) = \clip\left(F(x^{k};\xi^k), \lambda_k\right)$ using a fresh sample $\xi^k \sim \thesischtwocD_k$
\STATE $x^{k+1} = x^k - \gamma \thesischtwotF(x^{k};\xi^k)$
\ENDFOR
\ENSURE $x^{K+1}$ or $x_{\avg}^K = \frac{1}{K+1}\sum_{k=0}^K x^k$
\end{algorithmic}
\end{algorithm}

\subsection{Monotone Star-Cocoercive Problems}

We start with the following lemma derived by \textcite{gorbunov2022extragradient}. Since this lemma handles only deterministic part of the algorithm, the proof is the same as in the original work.
\begin{lemma}[Lemma D.1 from \textcite{gorbunov2022extragradient}]\label{ch2:lem:optimization_lemma_gap_SGDA_chapter_2}
    Let Assumptions~\ref{ch2:as:monotonicity_chapter_2} and \ref{ch2:as:star-cocoercivity_chapter_2} hold for $Q = B_{3R}(x^*)$, where $R \geq \|x^0 - x^*\|$ and $0 < \gamma \leq 2/\ell$. If $x^k$ lies in $B_{3R}(x^*)$ for all $k = 0,1,\ldots, K$ for some $K\geq 0$, then for all $u \in B_{3R}(x^*)$ the iterates produced by \algname{clipped-SGDA} satisfy
    \begin{align}
        \langle F(u), x^K_{\avg} - u\rangle & \leq \frac{\|x^0 - u\|^2 - \|x^{K+1} - u\|^2}{2\gamma(K+1)} \notag \\
        &\quad + \frac{\gamma}{2(K+1)}\sum\limits_{k=0}^K\left(\|F(x^k)\|^2 + \|\omega_k\|^2\right) \notag\\
        &\quad\quad + \frac{1}{K+1}\sum\limits_{k=0}^K\langle x^k - u - \gamma F(x^k), \omega_k\rangle, \label{ch2:eq:optimization_lemma_SGDA_chapter_2}\\
        x^K_{\avg} & \eqdef \frac{1}{K+1}\sum\limits_{k=0}^{K}x^k, \label{ch2:eq:x_avg_SGDA_chapter_2}\\
        \omega_k & \eqdef F(x^k) - \thesischtwotF(x^k;\xi^k). \label{ch2:eq:omega_k_SGDA_chapter_2}
    \end{align}
\end{lemma}

Also we need to use the following lemma to estimate the term $\sum\limits_{k=0}^K\|F(x^k)\|^2$ from the right hand side of \eqref{ch2:eq:optimization_lemma_SGDA_chapter_2} in the proof of the main theorem.

\begin{lemma}[Lemma D.2 from \textcite{gorbunov2022extragradient}]\label{ch2:lem:optimization_lemma_SGDA_chapter_2}
    Let Assumption~\ref{ch2:as:star-cocoercivity_chapter_2} hold for $Q = B_{3R}(x^*)$, where $R \geq R_0 \eqdef \|x^0 - x^*\|$ and $0 < \gamma \leq 2/\ell$. If $x^k$ lies in $B_{3R}(x^*)$ for all $k = 0,1,\ldots, K$ for some $K\geq 0$, then the iterates produced by \algname{clipped-SGDA} satisfy
    \begin{align}
        \MoveEqLeft[1] \frac{\gamma}{K+1}\left(\frac{2}{\ell} - \gamma\right)
        \sum\limits_{k=0}^K\|F(x^k)\|^2 \notag \\
        &\leq
        \frac{\|x^0 - x^*\|^2 - \|x^{K+1} - x^*\|^2}{K+1}
        \notag\\
        &\quad + \frac{2\gamma}{K+1}\sum\limits_{k=0}^K
        \langle x^k - x^* - \gamma F(x^k), \omega_k\rangle
        \notag\\
        &\quad + \frac{\gamma^2}{K+1}
        \sum\limits_{k=0}^K\|\omega_k\|^2.
        \label{ch2:eq:optimization_lemma_norm_SGDA_chapter_2}
    \end{align}
    where $\omega_k$ is defined in \eqref{ch2:eq:omega_k_SGDA_chapter_2}.
\end{lemma}

Using those lemmas, we prove the main convergence result for \algname{clipped-SGDA} in the monotone star-cocoercive case.

\begin{theorem}[Case 1 in Theorem~\ref{ch2:thm:clipped_SGDA_main_theorem_chapter_2}]\label{ch2:thm:main_result_gap_SGDA_chapter_2}
    Let Assumptions~\ref{ch2:as:bounded_alpha_moment_chapter_2}, \ref{ch2:as:monotonicity_chapter_2}, \ref{ch2:as:star-cocoercivity_chapter_2} hold for $Q = B_{3R}(x^*)$, where $R \geq \|x^0 - x^*\|$. Let $K \geq 0$ and $\beta \in (0,1]$ be such that $A = \ln(6(K+1)/\beta) \geq 1$, and
    \begin{align}
        0< \gamma & \leq \min\left\{\frac{1}{170\ell A}, \frac{R}{97200^{1/p}(K+1)^{1/p}\sigma A^{(p-1)/p}}\right\}, \label{ch2:eq:gamma_SGDA_chapter_2}\\
        \lambda_{k} \equiv \lambda & = \frac{R}{60\gamma A}, \label{ch2:eq:lambda_SGDA_chapter_2}
    \end{align}
    Then, after $K+1$ iterations the iterates produced by \algname{clipped-SGDA} with probability at least $1 - \beta$ satisfy 
    \begin{equation}
        \gap_R(x_{\avg}^K) \leq \frac{5R^2}{\gamma(K+1)} \quad \text{and}\quad \{x^k\}_{k=0}^{K+1} \subseteq B_{3R}(x^*), \label{ch2:eq:main_result_SGDA_chapter_2}
    \end{equation}
    where $x_{\avg}^K$ is defined in \eqref{ch2:eq:x_avg_SGDA_chapter_2}. In particular, when $\gamma$ equals the minimum from \eqref{ch2:eq:gamma_SGDA_chapter_2}, then the iterates produced by \algname{clipped-SGDA} after $K+1$ iterations with probability at least $1-\beta$ satisfy
    \begin{equation}
        \gap_R(x_{\avg}^K) = \cO\left(\max\left\{\frac{\ell R^2\ln(K/\beta)}{K}, \frac{\sigma R \ln^{(p-1)/p}(K/\beta)}{K^{(p-1)/p}}\right\}\right), \label{ch2:eq:clipped_SGDA_monotone_case_2_appendix_chapter_2}
    \end{equation}
    meaning that \algname{clipped-SGDA} requires the following number of iterations/oracle calls to achieve $\gap_R(x_{\avg}^K) \leq \varepsilon$ with probability at least $1 - \beta$:
    \begin{equation}
        K = \cO\left(\max\left\{\frac{\ell R^2}{\varepsilon}\ln\frac{\ell R^2}{\varepsilon\beta}, \left(\frac{\sigma R}{\varepsilon}\right)^{p/(p-1)}\ln\left(\frac{1}{\beta}\left(\frac{\sigma R}{\varepsilon}\right)^{p/(p-1)}\right)\right\}\right). \label{ch2:eq:clipped_SGDA_monotone_case_complexity_appendix_chapter_2}
    \end{equation}
\end{theorem}

\begin{proof}
    The proof follows similar steps as the proof of Theorem~D.1 from \textcite{gorbunov2022clipped}. The key difference is related to the application of Bernstein inequality and estimating biases and variances of stochastic terms.
    
    Let $R_k = \|x^k - x^*\|$ for all $k\geq 0$. As in the previous results, the proof is based on the induction argument and showing that the iterates do not leave some ball around the solution with high probability. More precisely, for each $k = 0,1,\ldots,K+1$ we consider probability event $E_k$ as follows: inequalities
    \begin{gather}
        \|x^t - x^*\|^2 \leq 2R^2 \quad \text{and}\quad \gamma\left\|\sum\limits_{l=0}^{t-1} \omega_l\right\| \leq R \label{ch2:eq:induction_inequality_SGDA_gap_chapter_2}
    \end{gather}
    hold for $t = 0,1,\ldots,k$ simultaneously. We want to prove that $\Prob\{E_k\} \geq  1 - k\beta/(K+1)$ for all $k = 0,1,\ldots,K+1$ by induction. The base of the induction is trivial: for $k = 0$ we have $R_0^2 \leq 2R^2$ by definition and $\sum_{l=0}^{-1} \omega_l = 0$. Next, assume that the statement holds for $k = T \leq K$, i.e., we have $\Prob\{E_{T}\} \geq 1 - T\beta/(K+1)$. Given this, we need to prove that $\Prob\{E_{T+1}\} \geq 1 - (T+1)\beta/(K+1)$. Since probability event $E_{T}$ implies $R^2_t \leq 2R^2$, we have $x^t \in B_{2R}(x^*)$ for all $t = 0, 1, \ldots, T$. According to this, the assumptions of Lemma~\ref{ch2:lem:optimization_lemma_SGDA_chapter_2} hold and  $E_{T}$ implies ($\gamma < 1/\ell$)
    \begin{align}
        \MoveEqLeft[1] \frac{\gamma}{\ell(T+1)}\sum\limits_{t=0}^T\|F(x^t)\|^2 \notag \\
        & \leq \frac{\|x^0 - x^*\|^2 - \|x^{T+1} - x^*\|^2}{T+1}\notag\\
        &\quad\quad + \frac{2\gamma}{T+1}\sum\limits_{t=0}^T\langle x^t - x^* - \gamma F(x^t), \omega_t\rangle + \frac{\gamma^2}{T+1}\sum\limits_{t=0}^T\|\omega_t\|^2
        \label{ch2:eq:thm_SGDA_technical_gap_0_chapter_2}
    \end{align}
    and by $\ell$-star-cocoersivity we have 
    \begin{align}
        \|F(x^t)\| & \leq \ell\|x^t - x^*\| \overset{\eqref{ch2:eq:induction_inequality_SGDA_gap_chapter_2}}{\leq} \sqrt{2}\ell R \overset{\eqref{ch2:eq:gamma_SGDA_chapter_2},\eqref{ch2:eq:lambda_SGDA_chapter_2}}{\leq} \frac{\lambda}{2} \label{ch2:eq:operator_bound_x_t_SGDA_gap_chapter_2}
    \end{align}
    for all $t = 0, 1, \ldots, T$.
    Using \eqref{ch2:eq:thm_SGDA_technical_gap_0_chapter_2}, we obtain  
    \begin{equation*}
       R_{T+1}^2 \leq R_0^2 + 2\gamma \sum\limits_{t=0}^T\langle x^t - x^* - \gamma F(x^t), \omega_t\rangle + \gamma^2\sum\limits_{t=0}^T\|\omega_t\|^2.
    \end{equation*}
    Due to \eqref{ch2:eq:operator_bound_x_t_SGDA_gap_chapter_2}, we have 
    \begin{align}
        \|x^t - x^* - \gamma F(x^t)\| & \leq \|x^t - x^*\| + \gamma\|F(x^t)\| \notag \\
        & \overset{\eqref{ch2:eq:star-cocoercivity_chapter_2},\eqref{ch2:eq:induction_inequality_SGDA_gap_chapter_2}}{\leq} 2R + \gamma \ell\|x^t - x^*\| \notag\\
        & \overset{\eqref{ch2:eq:induction_inequality_SGDA_gap_chapter_2}}{\leq} 2R + 2R\gamma \ell \overset{\eqref{ch2:eq:gamma_SGDA_chapter_2}}{\leq} 3R, \label{ch2:eq:thm_SGDA_technical_gap_2_chapter_2}
    \end{align}
    for all $t = 0, 1, \ldots, T$. To handle the sum above, we introduce a new vector
    \begin{equation*}
        u^{t} = \begin{cases}x^t - x^* - \gamma F(x^t),& \text{if } \|x^t - x^* - \gamma F(x^t)\| \leq 3R,\\ 0,& \text{otherwise,} \end{cases}
    \end{equation*}
    for all $t = 0, 1, \ldots, T$. This vector $u^{t}$ is bounded with probability $1$:
    \begin{equation}
        \|u^{t}\| \leq 3R  \label{ch2:eq:thm_SGDA_technical_gap_3_chapter_2}
    \end{equation}
    for all $t = 0, 1, \ldots, T$. We also notice that probability event $E_{T}$ implies $u^{t} = x^t - x^* - \gamma F(x^t)$ for all $t = 0, 1, \ldots, T$. Thus, thanks to \eqref{ch2:eq:thm_SGDA_technical_gap_2_chapter_2}, $E_{T}$ implies
    \begin{equation*}
       R_{T+1}^2 \leq R^2 + 2\gamma\sum\limits_{t=0}^T\langle u^{t}, \omega_t\rangle + \gamma^2\sum\limits_{t=0}^T\|\omega_t\|^2.
    \end{equation*}
    To handle the sums appeared on the right-hand side of the previous inequality we consider unbiased and biased parts of $\omega_t$:
    \begin{gather}
        \omega_t^u \eqdef \Exp_{\xi^t}\left[\thesischtwotF(x^t;\xi^t)\right] - \thesischtwotF(x^t;\xi^t),\quad \omega_t^b \eqdef F(x^t) - \Exp_{\xi^t}\left[\thesischtwotF(x^t;\xi^t)\right] \label{ch2:eq:thm_SGDA_technical_gap_4_chapter_2}
    \end{gather}
    for all $t = 0, \ldots, T$. Also, by definition we have $\omega_t = \omega_t^u + \omega_t^b$ for all $t = 0, \ldots, T$. Therefore, $E_{T}$ implies
    \begin{align}
        R_{T+1}^2 & \leq R^2 + \underbrace{2\gamma \sum\limits_{t = 0}^{T} \langle u^{t}, \omega_t^u \rangle}_{\circledOne} + \underbrace{2\gamma \sum\limits_{t = 0}^{T} \langle u^{t}, \omega_t^b \rangle}_{\circledTwo} + \underbrace{2\gamma^2 \sum\limits_{t=0}^{T}\left(\Exp_{\xi^t}\left[\|\omega_t^u\|^2\right] \right)}_{\circledThree} \notag\\
        &\quad\quad + \underbrace{2\gamma^2 \sum\limits_{t=0}^{T}\left(\|\omega_t^u\|^2 - \Exp_{\xi^t}\left[\|\omega_t^u\|^2\right]\right)}_{\circledFour} + \underbrace{2\gamma^2 \sum\limits_{t=0}^{T}\left(\|\omega_t^b\|^2\right)}_{\circledFive}.\label{ch2:eq:thm_SGDA_technical_gap_5_chapter_2}
    \end{align}
    We notice that the above inequality does not rely on monotonicity of $F$.
    
    According to the induction assumption, from probability event $E_T$ we have $x^t \in B_{2R}(x^*)$ for all $t = 0, 1, \ldots, T$. Thus, the assumptions of Lemma~\ref{ch2:lem:optimization_lemma_gap_SGDA_chapter_2} hold and  probability event $E_{T}$ implies
    \begin{align}
        2\gamma(T+1)\gap_R(x^T_{\avg}) & \leq \max\limits_{u\in B_R(x^*)}\left\{\|x^0 - u\|^2 + 2\gamma \sum\limits_{t=0}^T\langle x^t - u - \gamma F(x^t), \omega_t\rangle\right\}\notag\\
        &\quad\quad + \gamma^2\sum\limits_{t=0}^T\left(\|F(x^t)\|^2 + \|\omega_t\|^2\right),  \notag\\
        & = \max\limits_{u\in B_R(x^*)}\left\{\|x^0 - u\|^2 + 2\gamma \sum\limits_{t=0}^T\langle x^* - u, \omega_t\rangle\right\}\notag\\
        &\quad\quad + 2\gamma \sum\limits_{t=0}^T\langle x^t - x^* - \gamma F(x^t), \omega_t\rangle\notag\\
        &\quad\quad + \gamma^2\sum\limits_{t=0}^T\left(\|F(x^t)\|^2 + \|\omega_t\|^2\right).  \notag
    \end{align}
    As we mentioned before, $E_T$ implies $u^{t} = x^t - x^* - \gamma F(x^t)$ for all $t = 0, 1, \ldots, T$ as well as \eqref{ch2:eq:thm_SGDA_technical_gap_0_chapter_2} and $\gamma < 1/\ell$. Due to that, probability event $E_T$ implies
    \begin{align}
        2\gamma(T+1)\gap_R(x^T_{\avg}) & \leq \max\limits_{u\in B_R(x^*)}\left\{\|x^0 - u\|^2\right\} \notag \\
        &\quad + 2\gamma \max\limits_{u\in B_R(x^*)}\left\{\sum\limits_{t=0}^T\langle x^* - u, \omega_t\rangle\right\} \notag\\
        &\quad\quad + 2\gamma \sum\limits_{t=0}^T\langle u^{t}, \omega_t\rangle + \frac{\gamma}{\ell}\sum\limits_{t=0}^T\|F(x^t)\|^2 + \gamma^2\sum\limits_{t=0}^T\|\omega_t\|^2 \notag\\
        & \leq 4R^2 + 2\gamma \max\limits_{u\in B_R(x^*)}\left\{\left\langle x^* - u, \sum\limits_{t=0}^T\omega_t\right\rangle\right\} \notag\\
        &\quad\quad + R^2 + 4\gamma \sum\limits_{t=0}^T\langle u^{t}, \omega_t\rangle + 2\gamma^2\sum\limits_{t=0}^T\|\omega_t\|^2\notag\\
        & \leq 5R^2 + 2\gamma R\left\|\sum\limits_{t=0}^T\omega_t\right\| \notag\\
        &\quad + 2\cdot\left(\circledOne + \circledTwo + \circledThree + \circledFour + \circledFive\right),\label{ch2:eq:thm_SGDA_technical_gap_5_1_chapter_2}
    \end{align}
    where we also apply inequality $\|a+b\|^2 \leq 2\|a\|^2 + 2\|b\|^2$ holding for all $a,b \in \R^d$ to upper bound $\|\omega_t\|^2$.
    
    It remains to derive good enough high-probability upper-bounds for the terms $\circledOne, \circledTwo, \circledThree, \circledFour, \circledFive$  and $2\gamma R\left\|\sum_{t=0}^{T}\omega_t\right\|$, i.e., to finish our inductive proof we need to show that $\circledOne + \circledTwo + \circledThree + \circledFour + \circledFive \leq R^2$ and $2\gamma R\left\|\sum_{t=0}^{T}\omega_t\right\| \leq 2R^2$ with high probability. In the subsequent parts of the proof, we will need to use many times the bounds for the norm and second moments of $\omega_t^u, \omega_t^b$. First, by Lemma~\ref{ch2:lem:bias_variance_chapter_2}, we have with probability $1$ that
    \begin{equation}
        \|\omega_t^u\| \leq 2\lambda \label{ch2:eq:omega_magnitude_SGDA_gap_chapter_2}
    \end{equation}
    for all $t = 0,1, \ldots, T$. Moreover, due to Lemma~\ref{ch2:lem:bias_variance_chapter_2}, we also have that $E_{T}$ implies 
    \begin{gather}
         \left\|\omega_t^b\right\| \leq \frac{2^{{p}}\sigma^{{p}}}{\lambda^{{p}-1}}, \label{ch2:eq:bias_omega_SGDA_gap_chapter_2}\\
         \Exp_{\xi^t}\left[\left\|\omega^b_t\right\|^2\right] \leq 18\lambda^{2-{p}}\sigma^{{p}}, \label{ch2:eq:distortion_omega_SGDA_gap_chapter_2}\\
         \Exp_{\xi^t}\left[\left\|\omega_t^u\right\|^2\right] \leq 18\lambda^{2-{p}}\sigma^{{p}} \label{ch2:eq:variance_omega_SGDA_gap_chapter_2}
    \end{gather}
    for all $k = 0,1, \ldots, K$.
    
    \paragraph{Upper bound for $\circledOne$.} By definition of $\omega^u_t$, we have $\Exp_{\xi^t}[\omega_t^u] = 0$ and 
    \begin{equation*}
        \Exp_{\xi^t}\left[2\gamma\langle u^{t}, \omega_t^u \rangle\right] = 0.
    \end{equation*}
    Next, the sum $\circledOne$ has bounded with probability $1$ term:
    \begin{equation}
        |2\gamma\langle u^{t}, \omega_t^u \rangle | \leq 2\gamma  \|u^{t}\|\cdot \|\omega_t^u\| \overset{\eqref{ch2:eq:thm_SGDA_technical_gap_3_chapter_2},\eqref{ch2:eq:omega_magnitude_SGDA_gap_chapter_2}}{\leq} 12 \gamma R \lambda \overset{\eqref{ch2:eq:lambda_SGDA_chapter_2}}{\leq} \frac{R^2}{5A} \eqdef c. \label{ch2:eq:SGDA_neg_mon_technical_1_1_chapter_2}
    \end{equation}
    Moreover, these summands also have bounded conditional variances $\sigma_t^2$ (the conditional expectations $\Exp_{\xi^t}$ of their squares):
    \begin{align}
\sigma_t^2 & \leq \Exp_{\xi^t}\left[4\gamma^2 \|u^{t}\|^2\cdot \|\omega_t^u\|^2\right] \notag \\
& \overset{\eqref{ch2:eq:thm_SGDA_technical_gap_3_chapter_2}}{\leq} 36\gamma^2 R^2 \Exp_{\xi^t}\left[\|\omega_t^u\|^2\right]. \label{ch2:eq:SGDA_neg_mon_technical_1_2_chapter_2}
    \end{align}
    In other words, we showed that the summands of $\circledOne$ form a bounded martingale difference sequence with bounded conditional variances $\{\sigma_t^2\}_{t \geq 0}$. Next, we apply Bernstein's inequality (Lemma~\ref{ch2:lem:Bernstein_ineq_chapter_2}) with $X_t = 2\gamma \langle u^{t}, \omega_t^u \rangle$,  parameter $c$ as in \eqref{ch2:eq:SGDA_neg_mon_technical_1_1_chapter_2}, $b = \frac{R^2}{5}$, $H = \tfrac{R^4}{150A}$:
    \begin{equation*}
        \Prob\left\{|\circledOne| > \frac{R^2}{5} \text{ and } \sum\limits_{t=0}^{T}\sigma_t^2 \leq \frac{R^4}{150A}\right\} \leq 2\exp\left(- \frac{b^2}{2H + 2cb/3}\right) = \frac{\beta}{3(K+1)}.
    \end{equation*}
    Equivalently, we have $\Prob\{E_{\circledOne}\} \geq 1 - \tfrac{\beta}{3(K+1)}$ for
    \begin{equation}
        E_{\circledOne} = \left\{\text{either } \sum\limits_{t=0}^{T}\sigma_t^2 > \frac{R^4}{150A} \text{ or } |\circledOne| \leq \frac{R^2}{5}\right\}. \label{ch2:eq:bound_1_SGDA_neg_mon_chapter_2}
    \end{equation}
    In addition, $E_{T}$ implies that
    \begin{align}
        \sum\limits_{t=0}^{T}\sigma_t^2 & \overset{\eqref{ch2:eq:SGDA_neg_mon_technical_1_2_chapter_2}}{\leq} 36\gamma^2R^2\sum\limits_{t=0}^{T} \Exp_{\xi^t}\left[\|\omega_t^u\|^2\right] \notag\\
        & \overset{\eqref{ch2:eq:variance_omega_SGDA_gap_chapter_2}, T \leq K+1}{\leq} 648\gamma^2 R^2 \sigma^{{p}} (K+1) \lambda^{2-{p}} \notag\\
        & \overset{\eqref{ch2:eq:lambda_SGDA_chapter_2}}{\leq} 648 \gamma^{{p}}R^{4-{p}}\sigma^{{p}} (K+1)A^{{p}-2} \notag\\
        & \overset{\eqref{ch2:eq:gamma_SGDA_chapter_2}}{\leq} \frac{R^4}{150A}. \label{ch2:eq:bound_1_variances_SGDA_neg_mon_chapter_2}
    \end{align}
    
    \paragraph{Upper bound for $\circledTwo$.} From $E_{T}$ it follows that
    \begin{align}
        \circledTwo & \leq 2\gamma \sum\limits_{t=0}^{T}\|u^{t}\| \cdot \|\omega_t^b\| \overset{\eqref{ch2:eq:thm_SGDA_technical_gap_3_chapter_2},\eqref{ch2:eq:bias_omega_SGDA_gap_chapter_2}, T \leq K+1}{\leq} 6\cdot 2^{{p}}\gamma R(K+1) \frac{\sigma^{{p}}}{\lambda^{{p}-1}}\notag\\
        & \overset{\eqref{ch2:eq:lambda_SGDA_chapter_2}}{=} 12 \cdot 120^{{p}-1}\gamma^{{p}}\sigma^{{p}}R^{2-{p}} (K+1)A^{{p}-1} \overset{\eqref{ch2:eq:gamma_SGDA_chapter_2}}{\leq} \frac{R^2}{5}. \label{ch2:eq:bound_2_SGDA_neg_mon_chapter_2}
    \end{align}

    \paragraph{Upper bound for $\circledThree$.}  From $E_{T}$ it follows that
    \begin{align}
        \circledThree & =  2\gamma^2 \sum\limits_{t = 0}^{T}\Exp_{\xi^t}\left[\|\omega_t^u\|^2\right] \overset{\eqref{ch2:eq:variance_omega_SGDA_gap_chapter_2}, T \leq K+1}{\leq} 36\gamma^2 \lambda^{2-{p}} \sigma^{{p}} (K+1) \notag\\
        & \overset{\eqref{ch2:eq:lambda_SGDA_chapter_2}}{\leq} 36 \gamma^{{p}} R^{2-{p}} \sigma^{{p}} (K+1)  A^{{p}-2}
        \overset{\eqref{ch2:eq:gamma_SGDA_chapter_2}}{\leq} \frac{R^2}{5}. \label{ch2:eq:bound_3_SGDA_neg_mon_chapter_2}
    \end{align}
    
    \paragraph{Upper bound for $\circledFour$.} First, we have
    \begin{equation*}
        2\gamma^2\Exp_{\xi^t}\left[\|\omega_t^u\|^2 - \Exp_{\xi^t}\left[\|\omega_t^u\|^2\right]\right] = 0.
    \end{equation*}
    Next, the sum $\circledFour$ has bounded with probability $1$ terms:
    \begin{align}
        \MoveEqLeft[1] 2\gamma^2 \left|\|\omega_t^u\|^2 - \Exp_{\xi^t}\left[\|\omega_t^u\|^2\right] \right| \notag \\
        & \leq 2\gamma^2\left( \|\omega_t^u\|^2 + \Exp_{\xi^t}\left[\|\omega_t^u\|^2\right] \right) \overset{\eqref{ch2:eq:omega_magnitude_SGDA_gap_chapter_2}}{\leq} 16\gamma^2 \lambda^2 \notag\\
        & \overset{\eqref{ch2:eq:lambda_SGDA_chapter_2}}{\leq} \frac{R^2}{225A} \leq \frac{R^2}{5A} \eqdef c. \label{ch2:eq:SGDA_neg_mon_technical_4_1_chapter_2}
    \end{align}
    The summands also have conditional variances $\widetilde\sigma_t^2$ (the conditional expectations $\Exp_{\xi^t}$ of their squares) that are bounded
    \begin{align}
\widetilde\sigma_t^2 & \overset{\eqref{ch2:eq:SGDA_neg_mon_technical_4_1_chapter_2}}{\leq} \frac{2\gamma^2 R^2}{225A} \Exp_{\xi^t}\left[\left| \|\omega_t^u\|^2 - \Exp_{\xi^t}\left[\|\omega_t^u\|^2\right] \right|\right] \notag \\
& \leq \frac{4\gamma^2 R^2}{225A} \Exp_{\xi^t}\left[\|\omega_t^u\|^2\right]. \label{ch2:eq:SGDA_neg_mon_technical_4_2_chapter_2}
    \end{align}
    In other words, we showed that the summands of $\circledFour$ form a bounded martingale difference sequence with bounded conditional variances $\{\widetilde\sigma_t^2\}_{t \geq 0}$. Next, we apply Bernstein's inequality (Lemma~\ref{ch2:lem:Bernstein_ineq_chapter_2}) with $X_t$ equal to the summands of $\circledFour$, parameter $c$ as in \eqref{ch2:eq:SGDA_neg_mon_technical_4_1_chapter_2}, $b = \frac{R^2}{5}$, $H = \tfrac{R^4}{150A}$:
    \begin{align*}
\MoveEqLeft[1] \Prob\left\{|\circledFour| > \frac{R^2}{5} \text{ and } \sum\limits_{t=0}^{T}\widetilde\sigma_t^2 \leq \frac{R^4}{150A}\right\} \\
& \leq 2\exp\left(- \frac{b^2}{2H + 2cb/3}\right) = \frac{\beta}{3(K+1)}.
    \end{align*}
    Equivalently, we have $\Prob\{E_{\circledFour}\} \geq 1 - \tfrac{\beta}{3(K+1)}$ for
    \begin{equation}
        E_{\circledFour} = \left\{\text{either } \sum\limits_{t=0}^{T}\widetilde\sigma_t^2 > \frac{R^4}{150A} \text{ or } |\circledFour| \leq \frac{R^2}{5}\right\}. \label{ch2:eq:bound_4_SGDA_neg_mon_chapter_2}
    \end{equation}
    In addition, $E_{T}$ implies that
    \begin{align}
        \sum\limits_{t=0}^{T}\widetilde\sigma_t^2 & \overset{\eqref{ch2:eq:SGDA_neg_mon_technical_4_2_chapter_2}}{\leq} \frac{4\gamma^2 R^2}{225A}\sum\limits_{t=0}^{T} \Exp_{\xi^t}\left[\|\omega_t^u\|^2\right] \notag \\
        & \overset{\eqref{ch2:eq:variance_omega_SGDA_gap_chapter_2}, T \leq K+1}{\leq} \frac{8\gamma^2 R^2 (K+1)}{25A} \lambda^{2-{p}} \sigma^{{p}} \notag\\
        & \overset{\eqref{ch2:eq:lambda_SGDA_chapter_2}}{\leq} \frac{8}{25}\gamma^{{p}}
         R^{4-{p}} (K+1) \sigma^{{p}}A^{{p}-3}\notag\\
        & \overset{\eqref{ch2:eq:gamma_SGDA_chapter_2}}{\leq}\frac{R^4}{150A}. \label{ch2:eq:bound_4_variances_SGDA_neg_mon_chapter_2}
    \end{align}

    \paragraph{Upper bound for $\circledFive$.} From $E_{T}$ it follows that
    \begin{align}
        \circledFive & =  2\gamma^2 \sum\limits_{t = 0}^{T}\|\omega_t^b\|^2  \overset{\eqref{ch2:eq:bias_omega_SGDA_gap_chapter_2}, T \leq K+1}{\leq} 2^{2{p}+1}\gamma^2 (K+1) \frac{\sigma^{2{p}}}{\lambda^{2{p}-2}}\notag\\
        & \overset{\eqref{ch2:eq:lambda_SGDA_chapter_2}}{=} 2^{2{p}+1}\cdot 60^{2{p}-2} \gamma^{2{p}} (K+1) \frac{\sigma^{2{p}}}{R^{2{p}-2}}A^{2{p}-2} \notag\\
        & \overset{\eqref{ch2:eq:gamma_SGDA_chapter_2}}{\leq} \frac{R^2}{5}. \label{ch2:eq:bound_5_SGDA_neg_mon_chapter_2}
    \end{align}

    \paragraph{Upper bound for $\gamma \left\|\sum_{t=0}^{T}\omega_t\right\|$.} To estimate this term from above, we consider a new vector:
    \begin{equation*}
        \zeta_l = \begin{cases} \gamma \sum\limits_{r=0}^{l-1}\omega_r,& \text{if } \left\|\gamma \sum\limits_{r=0}^{l-1}\omega_r\right\| \leq R,\\ 0, & \text{otherwise} \end{cases}
    \end{equation*}
    for $l = 0, 1, \ldots, T$. This vector is bounded almost surely:
    \begin{equation}
        \|\zeta_l\| \leq R.  \label{ch2:eq:gap_thm_SGDA_technical_8_chapter_2}
    \end{equation}
    Thus, by \eqref{ch2:eq:induction_inequality_SGDA_gap_chapter_2}, probability event $E_{T}$ implies
    \begin{align}
        \gamma\left\|\sum\limits_{l = 0}^{T}\omega_l\right\| & = \sqrt{\gamma^2\left\|\sum\limits_{l = 0}^{T}\omega_l\right\|^2}\notag\\
        & = \sqrt{\gamma^2\sum\limits_{l=0}^{T}\|\omega_l\|^2 + 2\gamma\sum\limits_{l=0}^{T}\left\langle \gamma\sum\limits_{r=0}^{l-1} \omega_r, \omega_l \right\rangle} \notag\\
        & = \sqrt{\gamma^2\sum\limits_{l=0}^{T}\|\omega_l\|^2 + 2\gamma \sum\limits_{l=0}^{T} \langle \zeta_l, \omega_l\rangle} \notag\\
        & \overset{\eqref{ch2:eq:thm_SGDA_technical_gap_5_chapter_2}}{\leq} \sqrt{\circledThree + \circledFour + \circledFive + \underbrace{2\gamma \sum\limits_{l=0}^{T} \langle \zeta_l, \omega_l^u\rangle}_{\circledSix} + \underbrace{2\gamma \sum\limits_{l=0}^{T} \langle \zeta_l, \omega_l^b\rangle}_{\circledSeven}}. \label{ch2:eq:norm_sum_omega_bound_gap_SGDA_chapter_2}
    \end{align}
    Following similar steps as before, we bound $\circledSix$ and $\circledSeven$.

    \paragraph{Upper bound for $\circledSix$.} By definition of $\omega_t^u$, we have  $\Exp_{\xi^t}[\omega_t^u] = 0$ and 
    \begin{equation*}
        \Exp_{\xi^t}\left[2\gamma\langle \zeta_t, \omega_t^u \rangle\right] = 0.
    \end{equation*}
    Next, sum $\circledSix$ has bounded with probability $1$ terms:
    \begin{equation}
        |2\gamma\langle \zeta_t, \omega_t^u \rangle | \leq 2\gamma  \|\zeta_t\|\cdot \|\omega_t^u\| 
        \overset{\eqref{ch2:eq:gap_thm_SGDA_technical_8_chapter_2},\eqref{ch2:eq:omega_magnitude_SGDA_gap_chapter_2}}{\leq} 4 \gamma R \lambda \overset{\eqref{ch2:eq:lambda_SGDA_chapter_2}}{\leq} \frac{R^2}{5A} \eqdef c. \label{ch2:eq:SGDA_neg_mon_technical_6_1_chapter_2}
    \end{equation}
    The summands also have bounded conditional variances $\widehat\sigma_t^2$ (the conditional expectations $\Exp_{\xi^t}$ of their squares):
    \begin{align}
\widehat\sigma_t^2 & \leq \Exp_{\xi^t}\left[4\gamma^2 \|\zeta_t\|^2\cdot \|\omega_t^u\|^2\right] \notag \\
& \overset{\eqref{ch2:eq:gap_thm_SGDA_technical_8_chapter_2}}{\leq} 4\gamma^2 R^2 \Exp_{\xi^t}\left[\|\omega_t^u\|^2\right]. \label{ch2:eq:SGDA_neg_mon_technical_6_2_chapter_2}
    \end{align}
    In other words, we showed that the summands of $\circledSix$ form a bounded martingale difference sequence with bounded conditional variances $\{\widehat\sigma_t^2\}_{t \geq 0}$. Applying Bernstein's inequality (Lemma~\ref{ch2:lem:Bernstein_ineq_chapter_2}) with $X_t$ equal to the summands of $\circledSix$, parameter $c$ as in \eqref{ch2:eq:SGDA_neg_mon_technical_6_1_chapter_2}, $b = \frac{R^2}{5}$, $H = \tfrac{R^4}{150A}$:
    \begin{align*}
\MoveEqLeft[1] \Prob\left\{|\circledSix| > \frac{R^2}{5} \text{ and } \sum\limits_{t=0}^{T}\widehat\sigma_t^2 \leq \frac{R^4}{150A}\right\} \\
& \leq 2\exp\left(- \frac{b^2}{2H + 2cb/3}\right) = \frac{\beta}{3(K+1)}.
    \end{align*}
    Equivalently, we have $\Prob\{E_{\circledSix}\} \geq 1 - \tfrac{\beta}{3(K+1)}$ for
    \begin{equation}
        E_{\circledSix} = \left\{\text{either } \sum\limits_{t=0}^{T}\widehat\sigma_t^2 > \frac{R^4}{150A} \text{ or } |\circledSix| \leq \frac{R^2}{5}\right\}. \label{ch2:eq:bound_6_SGDA_neg_mon_chapter_2}
    \end{equation}
    In addition, $E_{T}$ implies that
    \begin{align}
        \sum\limits_{t=0}^{T}\widehat\sigma_t^2 & \overset{\eqref{ch2:eq:SGDA_neg_mon_technical_6_2_chapter_2}}{\leq} 4\gamma^2R^2\sum\limits_{t=0}^{T} \Exp_{\xi^t}\left[\|\omega_t^u\|^2\right] \notag\\
        & \overset{\eqref{ch2:eq:variance_omega_SGDA_gap_chapter_2}, T \leq K+1}{\leq} 72\gamma^2 R^2 \sigma^{{p}} (K+1) \lambda^{2-{p}} \notag\\
        & \overset{\eqref{ch2:eq:lambda_SGDA_chapter_2}}{\leq} 72 \gamma^{{p}}R^{4-{p}}\sigma^{{p}} (K+1)A^{{p}-2} \notag\\
        & \overset{\eqref{ch2:eq:gamma_SGDA_chapter_2}}{\leq} \frac{R^4}{150A} . \label{ch2:eq:bound_6_variances_SGDA_neg_mon_chapter_2}
    \end{align}
    
    \paragraph{Upper bound for $\circledSeven$.} From $E_{T}$ it follows that
    \begin{align}
        \circledSeven & \leq 2\gamma \sum\limits_{t=0}^{T}\|\zeta_t\| \cdot \|\omega_t^b\| \overset{\eqref{ch2:eq:gap_thm_SGDA_technical_8_chapter_2},\eqref{ch2:eq:bias_omega_SGDA_gap_chapter_2}, T \leq K+1}{\leq} 8\cdot 2^{{p}}\gamma R(K+1) \frac{\sigma^{{p}}}{\lambda^{{p}-1}}\notag\\
        & \overset{\eqref{ch2:eq:lambda_SGDA_chapter_2}}{=} 16 \cdot 120^{{p}-1}\gamma^{{p}}\sigma^{{p}}R^{2-{p}}(K+1)A^{{p}-1} \overset{\eqref{ch2:eq:gamma_SGDA_chapter_2}}{\leq} \frac{R^2}{5}. \label{ch2:eq:bound_7_SGDA_neg_mon_chapter_2}
    \end{align}

    Now, we have the upper bounds for  $\circledOne, \circledTwo, \circledThree, \circledFour, \circledFive, \circledSix, \circledSeven$. In particular, probability event $E_{T}$ implies
    \begin{gather*}
        R_{T+1}^2 \overset{\eqref{ch2:eq:thm_SGDA_technical_gap_5_chapter_2}}{\leq} R^2 + \circledOne + \circledTwo + \circledThree + \circledFour + \circledFive ,\\
        2\gamma(T+1)\gap_R(x^T_{\avg}) \overset{\eqref{ch2:eq:thm_SGDA_technical_gap_5_1_chapter_2}}{\leq} 5R^2 + 2\gamma R\left\|\sum\limits_{t=0}^T\omega_t\right\| + 2\cdot\left(\circledOne + \circledTwo + \circledThree + \circledFour + \circledFive\right), \\
        \gamma\left\|\sum\limits_{l = 0}^{T}\omega_l\right\| \overset{\eqref{ch2:eq:norm_sum_omega_bound_gap_SGDA_chapter_2}}{\leq} \sqrt{\circledThree + \circledFour + \circledFive + \circledSix + \circledSeven},\\
        \circledTwo \overset{\eqref{ch2:eq:bound_2_SGDA_neg_mon_chapter_2}}{\leq} \frac{R^2}{5},\quad \circledThree \overset{\eqref{ch2:eq:bound_3_SGDA_neg_mon_chapter_2}}{\leq} \frac{R^2}{5},\quad \circledFive \overset{\eqref{ch2:eq:bound_5_SGDA_neg_mon_chapter_2}}{\leq} \frac{R^2}{5},\quad \circledSeven \overset{\eqref{ch2:eq:bound_7_SGDA_neg_mon_chapter_2}}{\leq} \frac{R^2}{5},\\
        \sum\limits_{t=0}^{T}\sigma_t^2 \overset{\eqref{ch2:eq:bound_1_variances_SGDA_neg_mon_chapter_2}}{\leq}  \frac{R^4}{150A},\quad \sum\limits_{t=0}^{T}\widetilde\sigma_t^2 \overset{\eqref{ch2:eq:bound_4_variances_SGDA_neg_mon_chapter_2}}{\leq} \frac{R^4}{150A},\quad \sum\limits_{t=0}^{T}\widehat\sigma_t^2 \overset{\eqref{ch2:eq:bound_6_variances_SGDA_neg_mon_chapter_2}}{\leq} \frac{R^4}{150A}.
    \end{gather*}
    Moreover, we also have (see \eqref{ch2:eq:bound_1_SGDA_neg_mon_chapter_2}, \eqref{ch2:eq:bound_4_SGDA_neg_mon_chapter_2}, \eqref{ch2:eq:bound_6_SGDA_neg_mon_chapter_2} and our induction assumption)
    \begin{gather*}
        \Prob\{E_{T}\} \geq 1 - \frac{T\beta}{K+1},\\
        \Prob\{E_{\circledOne}\} \geq 1 - \frac{\beta}{3(K+1)}, \quad \Prob\{E_{\circledFour}\} \geq 1 - \frac{\beta}{3(K+1)}, \quad \Prob\{E_{\circledSix}\} \geq 1 - \frac{\beta}{3(K+1)},
    \end{gather*}
    where
    \begin{align}
        E_{\circledOne}& = \left\{\text{either } \sum\limits_{t=0}^{T}\sigma_t^2 > \frac{R^4}{150A} \text{ or } |\circledOne| \leq \frac{R^2}{5}\right\},\notag\\
        E_{\circledFour}& = \left\{\text{either } \sum\limits_{t=0}^{T}\widetilde\sigma_t^2 > \frac{R^4}{150A} \text{ or } |\circledFour| \leq \frac{R^2}{5}\right\},\notag\\
        E_{\circledSix}& = \left\{\text{either } \sum\limits_{t=0}^{T}\widehat\sigma_t^2 > \frac{R^4}{150A} \text{ or } |\circledSix| \leq \frac{R^2}{5}\right\}.\notag
    \end{align}
    Thus, probability event $E_{T} \cap E_{\circledOne} \cap E_{\circledFour} \cap E_{\circledSix}$ implies
    \begin{align*}
        R_{T+1}^2 & \leq R^2 + \circledOne + \circledTwo + \circledThree + \circledFour + \circledFive \leq 2R^2,\\
        \gamma\left\|\sum\limits_{l = 0}^{T}\omega_l\right\| & \leq \sqrt{\circledThree + \circledFour + \circledFive + \circledSix + \circledSeven} \leq R,\\
        2\gamma(T+1)\gap_R(x^T_{\avg}) & \leq 5R^2 + 2\gamma R\left\|\sum\limits_{t=0}^T\omega_t\right\| + 2\cdot\left(\circledOne + \circledTwo + \circledThree + \circledFour + \circledFive\right)\\
        & \leq 10R^2,
    \end{align*}
    which gives \eqref{ch2:eq:induction_inequality_SGDA_gap_chapter_2} for $t = T+1$, and 
    \begin{align}
\Prob\{E_{T+1}\} & \geq \Prob\{E_{T} \cap E_{\circledOne} \cap E_{\circledFour} \cap E_{\circledSix}\} = 1 - \Prob\{\overline{E}_{T} \cup \overline{E}_{\circledOne} \cup \overline{E}_{\circledFour} \cup \overline{E}_{\circledSix}\} \notag \\
& \geq 1 - \frac{(T+1)\beta}{K+1}. \notag
    \end{align}
    This finishes the inductive part of our proof, i.e., for all $k = 0,1,\ldots,K+1$ we have $\Prob\{E_k\} \geq 1 - k\beta/(K+1)$.  In particular, for $k = K+1$ we have that with probability at least $1 - \beta$
    \begin{equation*}
        \gap_R(x_{\avg}^K) \leq \frac{5R^2}{\gamma(K+1)}.
    \end{equation*}
    Finally, if $\gamma$ equals the minimum in \eqref{ch2:eq:gamma_SGDA_chapter_2},
    then with probability at least $1-\beta$
    \begin{align*}
        \gap_{R}(x^{K}_{\avg}) & \leq \frac{5R^2}{\gamma(K+1)} \\
        & = \max\left\{\frac{850 \ell R^2A}{K+1}, \frac{5 \cdot97200^{1/p}\sigma R A^{(p-1)/p}}{(K+1)^{(p-1)/p}}\right\}\\
        & = \cO\left(\max\left\{\frac{\ell R^2\ln(K/\beta)}{K}, \frac{\sigma R \ln^{(p-1)/p}(K/\beta)}{K^{(p-1)/p}}\right\}\right).
    \end{align*}
    To get $\gap_{R}(x^{K}_{\avg}) \leq \varepsilon$ with probability at least $1-\beta$ it is sufficient to choose $K$ such that both terms in the maximum above are $\cO(\varepsilon)$. This leads to
    \begin{equation*}
         K = \cO\left(\max\left\{\frac{\ell R^2}{\varepsilon}\ln\frac{\ell R^2}{\varepsilon\beta}, \left(\frac{\sigma R}{\varepsilon}\right)^{p/(p-1)}\ln\left(\frac{1}{\beta}\left(\frac{\sigma R}{\varepsilon}\right)^{p/(p-1)}\right)\right\}\right)
    \end{equation*}
    that concludes the proof.
\end{proof}

\subsection{Star-Cocoercive Problems}

\begin{theorem}[Case 2 in Theorem~\ref{ch2:thm:clipped_SGDA_main_theorem_chapter_2}]\label{ch2:thm:main_result_non_mon_SGDA_chapter_2}
    Let Assumptions~\ref{ch2:as:bounded_alpha_moment_chapter_2}, \ref{ch2:as:star-cocoercivity_chapter_2} hold for $Q = B_{2R}(x^*)$, where $R \geq \|x^0 - x^*\|$. Let $K \geq 0$ and $\beta \in (0,1]$ be such that $A = \ln(4(K+1)/\beta) \geq 1$, and
    \begin{align}
        0< \gamma & \leq \min\left\{\frac{1}{170\ell A}, \frac{R}{97200^{1/p}(K+1)^{1/p}\sigma A^{(p-1)/p}}\right\}, \label{ch2:eq:gamma_SGDA_non_mon_chapter_2}\\
        \lambda_{k} \equiv \lambda & = \frac{R}{60\gamma A}, \label{ch2:eq:lambda_SGDA_non_mon_chapter_2}
    \end{align}
    Then, after $K+1$ iterations the iterates produced by \algname{clipped-SGDA} with probability at least $1 - \beta$ satisfy 
    \begin{equation}
        \frac{1}{K+1}\sum^K_{k=0}\|F(x^k)\|^2 \leq \frac{2\ell R^2}{\gamma(K+1)} . \label{ch2:eq:main_result_non_mon_chapter_2}
    \end{equation}
     In particular, when $\gamma$ equals the minimum from \eqref{ch2:eq:gamma_SGDA_non_mon_chapter_2}, then the iterates produced by \algname{clipped-SGDA} after $K+1$ iterations with probability at least $1-\beta$ satisfy
    \begin{align}
         \MoveEqLeft \frac{1}{K+1}\sum^K_{k=0}\|F(x^k)\|^2 \notag\\
         &= \cO\left(\max\left\{\frac{\ell^2 R^2\ln(K/\beta)}{K+1}, \frac{\ell \sigma R \ln^{(p-1)/p}(K/\beta)}{K^{(p-1)/p}}\right\}\right), \label{ch2:eq:clipped_SGDA_non_monotone_case_2_appendix_chapter_2}
    \end{align}
    meaning that \algname{clipped-SGDA} requires the following number of iterations/oracle calls to achieve $\frac{1}{K+1}\sum_{k=0}^K\|F(x^k)\|^2 \leq \varepsilon$ with probability at least $1 - \beta$:
    \begin{align}
        K = \cO\Biggl(\max\Biggl\{&\frac{\ell^2 R^2}{\varepsilon}\ln\frac{\ell^2 R^2}{\varepsilon\beta}, \notag\\
        &\left(\frac{\ell \sigma R}{\varepsilon}\right)^{p/(p-1)}\ln\left( \frac{1}{\beta} \left(\frac{\ell \sigma R}{\varepsilon}\right)^{p/(p-1)}\right)\Biggr\}\Biggr). \label{ch2:eq:clipped_SGDA_non_monotone_case_complexity_appendix_chapter_2}
    \end{align}
\end{theorem}
\begin{proof}
    Again, we will closely follow the proof of Theorem~D.2 from \textcite{gorbunov2022clipped} and the main difference will be reflected in the application of Bernstein inequality and estimating biases and variances of stochastic terms.

    Let $R_k = \|x^k - x^*\|$ for all $k\geq 0$. As the previous result, the proof is based on the induction argument and showing that the iterates do not leave some ball around the solution with high probability. More precisely, for each $k = 0,\ldots, K+1$ we define probability event $E_k$ as follows: inequalities
    \begin{gather}
        \|x^t - x^*\|^2 \leq 2R^2, \label{ch2:eq:induction_inequality_SGDA_chapter_2}
    \end{gather}
    hold for $t = 0,1,\ldots,k$ simultaneously. We want to prove that $\Prob\{E_k\} \geq  1 - k\beta/(K+1)$ for all $k = 0,1,\ldots,K+1$ by induction. One of the important things is that inequalities \eqref{ch2:eq:thm_SGDA_technical_gap_0_chapter_2} and \eqref{ch2:eq:thm_SGDA_technical_gap_5_chapter_2} are obtained without assuming monotonicity of $F$. Thus, if we do exactly the same steps as in the proof of Theorem~\ref{ch2:thm:main_result_gap_SGDA_chapter_2} (up to the replacement of $\ln \frac{6(K+1)}{\beta}$ by $A$), we gain that
    \begin{gather*}
        R_{T+1}^2 \overset{\eqref{ch2:eq:thm_SGDA_technical_gap_5_chapter_2}}{\leq} R^2 + \circledOne + \circledTwo + \circledThree + \circledFour + \circledFive ,\\
        \circledTwo \overset{\eqref{ch2:eq:bound_2_SGDA_neg_mon_chapter_2}}{\leq} \frac{R^2}{5},\quad \circledThree \overset{\eqref{ch2:eq:bound_3_SGDA_neg_mon_chapter_2}}{\leq} \frac{R^2}{5},\quad \circledFive \overset{\eqref{ch2:eq:bound_5_SGDA_neg_mon_chapter_2}}{\leq} \frac{R^2}{5},\\
        \sum\limits_{t=0}^{T}\sigma_t^2 \overset{\eqref{ch2:eq:bound_1_variances_SGDA_neg_mon_chapter_2}}{\leq}  \frac{R^4}{150A},\quad \sum\limits_{t=0}^{T}\widetilde\sigma_t^2 \overset{\eqref{ch2:eq:bound_4_variances_SGDA_neg_mon_chapter_2}}{\leq} \frac{R^4}{150A}.
    \end{gather*}
    Moreover, we also have (see \eqref{ch2:eq:bound_1_SGDA_neg_mon_chapter_2}, \eqref{ch2:eq:bound_4_SGDA_neg_mon_chapter_2} and our induction assumption)
    \begin{gather*}
        \Prob\{E_{T}\} \geq 1 - \frac{T\beta}{K+1},\\
        \Prob\{E_{\circledOne}\} \geq 1 - \frac{\beta}{2(K+1)}, \quad \Prob\{E_{\circledFour}\} \geq 1 - \frac{\beta}{2(K+1)},
    \end{gather*}
    where 
    \begin{align}
        E_{\circledOne}& = \left\{\text{either } \sum\limits_{t=0}^{T}\sigma_t^2 > \frac{R^4}{150A} \text{ or } |\circledOne| \leq \frac{R^2}{5}\right\},\notag\\
        E_{\circledFour}& = \left\{\text{either } \sum\limits_{t=0}^{T}\widetilde\sigma_t^2 > \frac{R^4}{150A} \text{ or } |\circledFour| \leq \frac{R^2}{5}\right\}.\notag
    \end{align}
    Thus, probability event $E_{T} \cap E_{\circledOne} \cap E_{\circledFour}$ implies
    \begin{equation*}
        R_{T+1}^2 \leq R^2 + \circledOne + \circledTwo + \circledThree + \circledFour + \circledFive \leq 2R^2,
    \end{equation*}
    and 
    \begin{equation}
        \Prob\{E_{T+1}\} \geq \Prob\{E_{T} \cap E_{\circledOne} \cap E_{\circledFour}\} = 1 - \Prob\{\overline{E}_{T} \cup \overline{E}_{\circledOne} \cup \overline{E}_{\circledFour}\} \geq 1 - \frac{(T+1)\beta}{K+1}.
    \end{equation}
    This finishes the inductive part of our proof, i.e., for all $k = 0,1,\ldots,K+1$ we have $\Prob\{E_k\} \geq 1 - k\beta/(K+1)$. In particular, for $k = K+1$ we have that with probability at least $1 - \beta$
    \begin{align*}
        \frac{1}{K+1}\sum\limits_{k=0}^{K}\|F(x^k)\|^2 & \overset{\eqref{ch2:eq:thm_SGDA_technical_gap_0_chapter_2}}{\leq} \frac{\ell(R^2 - R_{K+1}^2)}{\gamma(K+1)} \\
        &\quad + \frac{\ell(\circledOne + \circledTwo + \circledThree + \circledFour + \circledFive)}{\gamma(K+1)}\\
        & \leq \frac{2\ell R^2}{\gamma (K+1)}.
    \end{align*}
    Finally, if $\gamma$ equals the minimum in \eqref{ch2:eq:gamma_SGDA_non_mon_chapter_2},
    then with probability at least $1-\beta$
    \begin{align*}
        \frac{1}{K+1}\sum\limits_{k=0}^{K}\|F(x^k)\|^2 & \leq \frac{2\ell R^2}{\gamma(K+1)} \\
        & = \max\left\{\frac{340 \ell^2 R^2A}{K+1}, \frac{2 \cdot97200^{1/p}\ell\sigma R A^{(p-1)/p}}{(K+1)^{(p-1)/p}}\right\}\\
        & = \cO\left(\max\left\{\frac{\ell^2 R^2\ln(K/\beta)}{K}, \frac{\ell \sigma R \ln^{(p-1)/p}(K/\beta)}{K^{(p-1)/p}}\right\}\right).
    \end{align*}
    To get $\frac{1}{K+1}\sum\limits_{k=0}^{K}\|F(x^k)\|^2 \leq \varepsilon$ with probability at least $1-\beta$ it is sufficient to choose $K$ such that both terms in the maximum above are $\cO(\varepsilon)$. This leads to
    \begin{equation*}
         K = \cO\left(\max\left\{\frac{\ell^2 R^2}{\varepsilon}\ln\frac{\ell^2 R^2}{\varepsilon\beta}, \left(\frac{\ell \sigma R}{\varepsilon}\right)^{p/(p-1)}\ln\left( \frac{1}{\beta} \left(\frac{\ell \sigma R}{\varepsilon}\right)^{p/(p-1)}\right)\right\}\right)
    \end{equation*}
    that concludes the proof.

\end{proof}

\subsection{Quasi-Strongly Monotone Star-Cocoercive Problems}

As in the monotone case, we use another lemma from \textcite{gorbunov2022clipped} that handles the deterministic part of \algname{clipped-SGDA} in the quasi-strongly monotone case.
\begin{lemma}[Lemma~D.3 from \textcite{gorbunov2022clipped}]\label{ch2:lem:optimization_lemma_str_mon_SGDA_chapter_2}
    Let Assumptions~\ref{ch2:as:QSM_chapter_2}, \ref{ch2:as:star-cocoercivity_chapter_2}  hold for $Q = B_{2R}(x^*) = \{x\in\R^d\mid \|x - x^*\| \leq 2R\}$, where $R \geq \|x^0 - x^*\|$, and $0 < \gamma \leq 1/\ell$. If $x^k$ lie in $B_{2R}(x^*)$ for all $k = 0,1,\ldots, K$ for some $K\geq 0$, then the iterates produced by \algname{clipped-SGDA} satisfy
    \begin{align}
        \|x^{K+1} - x^*\|^2 & \leq (1 - \gamma \mu)^{K+1}\|x^0 - x^*\|^2 \notag \\
        &\quad + 2\gamma \sum\limits_{k=0}^K (1-\gamma\mu)^{K-k} \langle x^k - x^* - \gamma F(x^k), \omega_k \rangle \notag\\
        &\quad\quad + \gamma^2 \sum\limits_{k=0}^K (1-\gamma\mu)^{K-k}\|\omega_k\|^2, \label{ch2:eq:optimization_lemma_SGDA_str_mon_chapter_2}
    \end{align}
    where $\omega_k$ are defined in \eqref{ch2:eq:omega_k_SGDA_chapter_2}.
\end{lemma}

Using this lemma we prove the main convergence result for \algname{clipped-SGDA} in the quasi-strongly monotone case.

\begin{theorem}[Case 2 in Theorem~\ref{ch2:thm:clipped_SGDA_main_theorem_chapter_2}]\label{ch2:thm:main_result_str_mon_SGDA_chapter_2}
    Let Assumptions~\ref{ch2:as:bounded_alpha_moment_chapter_2}, \ref{ch2:as:QSM_chapter_2}, \ref{ch2:as:star-cocoercivity_chapter_2} hold for $Q = B_{2R}(x^*) = \{x\in\R^d\mid \|x - x^*\| \leq 2R\}$, where $R \geq \|x^0 - x^*\|$. Let $K \geq 0$ and $\beta \in (0,1]$ be such that $A = \ln(4(K+1)/\beta) \geq 1$, and
    \begin{align}
        0< \gamma & \leq \min\left\{\frac{1}{400 \ell A}, \frac{\ln(B_K)}{\mu(K+1)}\right\}, \label{ch2:eq:gamma_SGDA_str_mon_chapter_2}\\
        B_K & = \max\left\{2, \frac{(K+1)^{2(p-1)/p}\mu^2R^2}{54000^{2/p}\sigma^2A^{2(p-1)/p}\ln^2(B_K)} \right\} \label{ch2:eq:B_K_SGDA_str_mon_1_chapter_2} \\
        \mathcal D_K
        & \eqdef
        \max\left\{
            2,
            \frac{K^{2(p-1)/p}\mu^2R^2}{\sigma^2\ln^{2(p-1)/p}(K/\beta)}
        \right\},\notag\\
        B_K
        & =
        \cO\left(
            \max\left\{
                2,
                \frac{K^{2(p-1)/p}\mu^2R^2}{\sigma^2\ln^{2(p-1)/p}(K/\beta)
                \ln^2\left(\mathcal D_K\right)}
            \right\}
        \right), \label{ch2:eq:B_K_SGDA_str_mon_2_chapter_2} \\
        \lambda_k & = \frac{\exp(-\gamma\mu(1 + k/2))R}{120\gamma A}, \label{ch2:eq:lambda_SGDA_str_mon_chapter_2}
    \end{align}
    where $K$ is large enough that the maximum in \eqref{ch2:eq:B_K_SGDA_str_mon_1_chapter_2} is attained at its second term. Then, after $K+1$ iterations the iterates produced by \algname{clipped-SGDA} with probability at least $1 - \beta$ satisfy 
    \begin{equation}
        \|x^{K+1} - x^*\|^2 \leq 2\exp(-\gamma\mu(K+1))R^2. \label{ch2:eq:main_result_str_mon_SGDA_chapter_2}
    \end{equation}
    In particular, when $\gamma$ equals the minimum from \eqref{ch2:eq:gamma_SGDA_str_mon_chapter_2}, then the iterates produced by \algname{clipped-SGDA} after $K+1$ iterations with probability at least $1-\beta$ satisfy
    \begin{align}
        \MoveEqLeft \|x^{K+1} - x^*\|^2\notag\\
        &= \cO\Biggl(\max\Biggl\{R^2\exp\left(- \frac{\mu K}{\ell \ln(K/\beta)}\right),\notag\\
        &\qquad\qquad \frac{\sigma^2\ln^{2(p-1)/p}(K/\beta)\ln^2\left(\mathcal D_K\right)}{K^{2(p-1)/p}\mu^2}\Biggr\}\Biggr).
        \label{ch2:eq:clipped_SGDA_str_monotone_case_2_appendix_chapter_2}
    \end{align}
    meaning that \algname{clipped-SGDA} requires the following number of iterations/oracle calls to achieve $\|x^{K+1} - x^*\|^2 \leq \varepsilon$ with probability at least $1 - \beta$:
    \begin{align}
        K &= \cO\biggl(\max\biggl\{\frac{\ell}{\mu}\ln\left(\frac{R^2}{\varepsilon}\right)\ln\left(\frac{\ell}{\mu \beta}\ln\frac{R^2}{\varepsilon}\right),\notag\\
        &\quad \left(\frac{\sigma^2}{\mu^2\varepsilon}\right)^{p/(2(p-1))}\ln \left(\frac{1}{\beta} \left(\frac{\sigma^2}{\mu^2\varepsilon}\right)^{p/(2(p-1))}\right)\ln^{p/(p-1)}\left(B_\varepsilon\right)\biggr\}\biggr), \label{ch2:eq:clipped_SGDA_str_monotone_case_complexity_appendix_chapter_2}
    \end{align}
    where
    \begin{equation*}
        B_\varepsilon = \max\left\{2, \frac{R^2}{\varepsilon \ln \left(\frac{1}{\beta} \left(\frac{\sigma^2}{\mu^2\varepsilon}\right)^{p/(2(p-1))}\right)}\right\}.
    \end{equation*}
\end{theorem}

\begin{proof}
    Again, we will closely follow the proof of Theorem~D.3 from \textcite{gorbunov2022clipped} and the main difference will be reflected in the application of Bernstein inequality and estimating biases and variances of stochastic terms.
    
    Let $R_k = \|x^k - x^*\|$ for all $k\geq 0$. As in the previous results, the proof is based on the induction argument and showing that the iterates do not leave some ball around the solution with high probability. More precisely, for each $k = 0,1,\ldots,K+1$ we consider probability event $E_k$ as follows: inequalities
    \begin{equation}
        R_t^2 \leq 2 \exp(-\gamma\mu t) R^2 \label{ch2:eq:induction_inequality_str_mon_SGDA_chapter_2}
    \end{equation}
    hold for $t = 0,1,\ldots,k$ simultaneously. We want to prove $\Prob\{E_k\} \geq  1 - k\beta/(K+1)$ for all $k = 0,1,\ldots,K+1$ by induction. The base of the induction is trivial: for $k=0$ we have $R_0^2 \leq R^2 < 2R^2$ by definition. Next, assume that for $k = T-1 \leq K$ the statement holds: $\Prob\{E_{T-1}\} \geq  1 - (T-1)\beta/(K+1)$. Given this, we need to prove $\Prob\{E_{T}\} \geq  1 - T\beta/(K+1)$. Since $R_t^2 \leq 2\exp(-\gamma\mu t) R^2 \leq 2R^2$, we have $x^t \in B_{2R}(x^*)$, where operator $F$ is $\ell$-star-cocoersive. Thus, $E_{T-1}$ implies
    \begin{align}
        \|F(x^t)\| & \leq \ell\|x^t - x^*\| \overset{\eqref{ch2:eq:induction_inequality_str_mon_SGDA_chapter_2}}{\leq} \sqrt{2}\ell\exp(- \gamma\mu t/2)R \overset{\eqref{ch2:eq:gamma_SGDA_str_mon_chapter_2},\eqref{ch2:eq:lambda_SGDA_str_mon_chapter_2}}{\leq} \frac{\lambda_t}{2} \label{ch2:eq:operator_bound_x_t_SGDA_str_mon_chapter_2}
    \end{align}
    and
    \begin{align}
        \|\omega_t\|^2 & \leq 2\|\thesischtwotF(x^t;\xi)\|^2 + 2\|F(x^t)\|^2 \overset{\eqref{ch2:eq:operator_bound_x_t_SGDA_str_mon_chapter_2}}{\leq} \frac{5}{2}\lambda_t^2 \overset{\eqref{ch2:eq:lambda_SGDA_str_mon_chapter_2}}{\leq} \frac{\exp(-\gamma\mu t)R^2}{4\gamma^2} \label{ch2:eq:omega_bound_x_t_SGDA_str_mon_chapter_2}
    \end{align}
    for all $t = 0, 1, \ldots, T-1$, where we use that $\|a+b\|^2 \leq 2\|a\|^2 + 2\|b\|^2$ holding for all $a,b \in \R^d$. 

    Using Lemma~\ref{ch2:lem:optimization_lemma_str_mon_SGDA_chapter_2} and $(1 - \gamma\mu)^T \leq \exp(-\gamma\mu T)$, we obtain that $E_{T-1}$ implies
    \begin{align}
        R_T^2 & \leq \exp(-\gamma\mu T)R^2 + 2\gamma \sum\limits_{t=0}^{T-1} (1-\gamma\mu)^{T-1-t} \langle x^t - x^* \notag \\
        &\quad - \gamma F(x^t), \omega_t \rangle \notag\\
        &\quad\quad + \gamma^2 \sum\limits_{t=0}^{T-1} (1-\gamma\mu)^{T-1-t}\|\omega_t\|^2. \notag
    \end{align}
    To handle the sums above, we introduce a new notation:
    \begin{gather}
        u^{t} = \begin{cases} x^t - x^* - \gamma F(x^t),& \text{if } \|x^t - x^* - \gamma F(x^t)\|\\ & \quad \leq \sqrt{2}(1 + \gamma \ell) \exp(- \gamma\mu t/2)R,\\ 0,& \text{otherwise}, \end{cases} \label{ch2:eq:eta_t_SGDA_str_mon_chapter_2}
    \end{gather}
    for $t = 0, 1, \ldots, T-1$. This vector is bounded almost surely:
    \begin{equation}
         \|u^{t}\| \leq \sqrt{2}(1 + \gamma \ell)\exp(-\gamma\mu t/2)R \label{ch2:eq:eta_t_bound_SGDA_str_mon_chapter_2} 
    \end{equation}
    for all $t = 0, 1, \ldots, T-1$. We also notice that, due to \eqref{ch2:eq:operator_bound_x_t_SGDA_str_mon_chapter_2}, event $E_{T-1}$ implies the bound $\|F(x^t)\| \leq \sqrt{2}\ell\exp(-\gamma\mu t/2)R$ and
    \begin{align*}
        \|x^t - x^* - \gamma F(x^t)\| & \leq \|x^t - x^*\| + \gamma \|F(x^t)\|\\
        & \overset{\eqref{ch2:eq:operator_bound_x_t_SGDA_str_mon_chapter_2}}{\leq} \sqrt{2}(1 + \gamma \ell)\exp(-\gamma\mu t/2)R
    \end{align*}
    for $t = 0, 1, \ldots, T-1$. In other words, $E_{T-1}$ implies $u^{t} = x^t - x^* - \gamma F(x^t)$ for all $t = 0,1, \ldots, T-1$, meaning that from $E_{T-1}$ it follows that
    \begin{align}
        R_T^2 & \leq \exp(-\gamma\mu T)R^2 + 2\gamma \sum\limits_{t=0}^{T-1} (1-\gamma\mu)^{T-1-t} \langle u^{t}, \omega_t \rangle \notag \\
        &\quad + \gamma^2 \sum\limits_{t=0}^{T-1} (1-\gamma\mu)^{T-1-t} \|\omega_t\|^2. \notag
    \end{align}

    To handle the sums appeared on the right-hand side of the previous inequality we consider unbiased and biased parts of $ \omega_t$:
    \begin{gather}
        \omega_t^u \eqdef \Exp_{\xi^t}\left[\thesischtwotF(x^t;\xi^t)\right] - \thesischtwotF(x^t;\xi^t),\quad \omega_t^b \eqdef F(x^t) - \Exp_{\xi^t}\left[\thesischtwotF(x^t;\xi^t)\right], \label{ch2:eq:omega_unbias_bias_SGDA_str_mon_chapter_2}
    \end{gather}
    for all $t = 0, \ldots, T-1$. By definition we have $\omega_t = \omega_t^u + \omega_t^b$ for all $t = 0, \ldots, T-1$. Therefore, $E_{T-1}$ implies
    \begin{align}
        R_T^2
        &\leq \exp(-\gamma\mu T) R^2
        + \underbrace{
            2\gamma \sum\limits_{t=0}^{T-1}
            (1-\gamma\mu)^{T-1-t}
            \langle u^{t}, \omega_t^u \rangle
        }_{\circledOne}
        \notag\\
        &\quad
        + \underbrace{
            2\gamma \sum\limits_{t=0}^{T-1}
            (1-\gamma\mu)^{T-1-t}
            \langle u^{t}, \omega_t^b \rangle
        }_{\circledTwo}
        \notag\\
        &\quad
        + \underbrace{
            2\gamma^2 \sum\limits_{t=0}^{T-1}
            (1-\gamma\mu)^{T-1-t}
            \Exp_{\xi^t}\left[\|\omega^u_t\|^2\right]
        }_{\circledThree}
        \notag\\
        &\quad
        + \underbrace{
            2\gamma^2 \sum\limits_{t=0}^{T-1}
            (1-\gamma\mu)^{T-1-t}
            \left(
                \|\omega^u_t\|^2
                - \Exp_{\xi^t}\left[\|\omega^u_t\|^2\right]
            \right)
        }_{\circledFour}
        \notag\\
        &\quad
        + \underbrace{
            2\gamma^2 \sum\limits_{t=0}^{T-1}
            (1-\gamma\mu)^{T-1-t}
            \|\omega^b_t\|^2
        }_{\circledFive}.
        \label{ch2:eq:SGDA_str_mon_12345_bound_chapter_2}
    \end{align}
    where we also apply inequality $\|a+b\|^2 \leq 2\|a\|^2 + 2\|b\|^2$ holding for all $a,b \in \R^d$ to upper bound $\|\omega_t\|^2$. It remains to derive good enough high-probability upper-bounds for the terms $\circledOne, \circledTwo, \circledThree, \circledFour, \circledFive$, i.e., to finish our inductive proof we need to show that $\circledOne + \circledTwo + \circledThree + \circledFour + \circledFive  \leq \exp(-\gamma\mu T) R^2$ with high probability. In the subsequent parts of the proof, we will need to use many times the bounds for the norm and second moments of $\omega_{t}^u$ and $\omega_{t}^b$. First, by Lemma \ref{ch2:lem:bias_and_variance_clip_chapter_2}, we have with probability $1$ that
    \begin{equation}
        \|\omega_t^u\| \leq 2\lambda_t.\label{ch2:eq:omega_magnitude_str_mon_chapter_2}
    \end{equation}
    Moreover, since $E_{T-1}$ implies that $\|F(x^t)\| \leq \lambda_t/2$ for all $t = 0,1, \ldots, T-1$ (see \eqref{ch2:eq:operator_bound_x_t_SGDA_str_mon_chapter_2}), from Lemma~\ref{ch2:lem:bias_and_variance_clip_chapter_2} we also have that $E_{T-1}$ implies
    \begin{gather}
        \left\|\omega_t^b\right\| \leq \frac{2^{{p}}\sigma^{p}}{\lambda_t^{{p}-1}}, \label{ch2:eq:bias_omega_str_mon_chapter_2}\\
        \Exp_{\xi^t}\left[\left\|\omega_t\right\|^2\right] \leq 18 \lambda_t^{2-{p}}\sigma^{p}, \label{ch2:eq:distortion_omega_str_mon_chapter_2}\\
        \Exp_{\xi^t}\left[\left\|\omega_t^u\right\|^2\right] \leq 18 \lambda_t^{2-{p}}\sigma^{p}, \label{ch2:eq:variance_omega_str_mon_chapter_2}
    \end{gather}
    for all $t = 0,1, \ldots, T-1$.

\paragraph{Upper bound for $\circledOne$.} By definition of $\omega_t^u$, we have $\Exp_{\xi^t}[\omega_{t}^u] = 0$ and
    \begin{equation*}
        \Exp_{\xi^t}\left[2\gamma (1-\gamma\mu)^{T-1-t} \langle u^{t}, \omega_t^u \rangle\right] = 0.
    \end{equation*}
    Next, sum $\circledOne$ has bounded with probability $1$ terms:
    \begin{align}
        \MoveEqLeft[1] |2\gamma (1-\gamma\mu)^{T-1-t} \langle u^{t}, \omega_t^u \rangle | \notag \\
        & \leq 2\gamma\exp(-\gamma\mu (T - 1 - t)) \|u^{t}\|\cdot \|\omega_t^u\|\notag\\
        & \overset{\eqref{ch2:eq:eta_t_bound_SGDA_str_mon_chapter_2},\eqref{ch2:eq:omega_magnitude_str_mon_chapter_2}}{\leq} 4\sqrt{2}\gamma (1 + \gamma \ell) \exp(-\gamma\mu (T - 1 - t/2)) R \lambda_t\notag\\
        & \overset{\eqref{ch2:eq:gamma_SGDA_str_mon_chapter_2},\eqref{ch2:eq:lambda_SGDA_str_mon_chapter_2}}{\leq} \frac{\exp(-\gamma\mu T)R^2}{5A} \eqdef c. \label{ch2:eq:SGDA_str_mon_technical_1_1_chapter_2}
    \end{align}
    The summands also have bounded conditional variances $\sigma_t^2$ (the conditional expectations $\Exp_{\xi^t}$ of their squares):
    \begin{align}
        \sigma_t^2 & \leq \Exp_{\xi^t}\left[4\gamma^2\exp(-\gamma\mu (2T - 2 - 2t)) \|u^{t}\|^2\cdot \|\omega_t^u\|^2\right]\notag\\
        & \overset{\eqref{ch2:eq:eta_t_bound_SGDA_str_mon_chapter_2}}{\leq} 8\gamma^2 (1 + \gamma \ell)^2 \exp(-\gamma\mu (2T - 2 - t)) R^2 \Exp_{\xi^t}\left[\|\omega_t^u\|^2\right]\notag\\
        & \overset{\eqref{ch2:eq:gamma_SGDA_str_mon_chapter_2}}{\leq} 10\gamma^2\exp(-\gamma\mu (2T - t))R^2 \Exp_{\xi^t}\left[\|\omega_t^u\|^2\right]. \label{ch2:eq:SGDA_str_mon_technical_1_2_chapter_2}
    \end{align}

    In other words, we showed that the summands of $\circledOne$ form a bounded martingale difference sequence with bounded conditional variances $\{\sigma_t^2\}_{t = 0}^{T-1}$. Next, we apply Bernstein's inequality (Lemma~\ref{ch2:lem:Bernstein_ineq_chapter_2}) with $X_t$ equal to the summands of $\circledOne$, parameter $c$ as in \eqref{ch2:eq:SGDA_str_mon_technical_1_1_chapter_2}, and
    \begin{equation*}
        b = \frac{1}{5}\exp(-\gamma\mu T) R^2,
        \qquad
        F = \frac{\exp(-2 \gamma\mu T) R^4}{300A}:
    \end{equation*}
    \begin{equation*}
        \Prob\left\{|\circledOne| > b \quad\text{and}\quad \sum\limits_{t=0}^{T-1}\sigma_t^2 \leq F\right\}
        \leq 2\exp\left(- \frac{b^2}{2F + 2cb/3}\right)
        \leq \frac{\beta}{2(K+1)}.
    \end{equation*}
    Equivalently, we have $\Prob\{E_{\circledOne}\} \geq 1 - \frac{\beta}{2(K+1)}$ for
    \begin{equation}
            E_{\circledOne} =
            \left\{\text{either } \sum\limits_{t=0}^{T-1}\sigma_t^2 > \frac{\exp(-2\gamma\mu T) R^4}{300A} \text{ or } |\circledOne| \leq \frac{1}{5}\exp(-\gamma\mu T) R^2\right\}.
        \label{ch2:eq:bound_1_SGDA_str_mon_chapter_2}
    \end{equation}
    In addition, $E_{T-1}$ implies that
    \begin{align}
        \sum\limits_{t=0}^{T-1}\sigma_t^2 & \overset{\eqref{ch2:eq:SGDA_str_mon_technical_1_2_chapter_2}}{\leq} 10\gamma^2\exp(- 2\gamma\mu T)R^2\sum\limits_{t=0}^{T-1} \frac{\Exp_{\xi^t}\left[\|\omega_t^u\|^2\right]}{\exp(-\gamma\mu t)}\notag\\ 
        & \overset{\eqref{ch2:eq:variance_omega_str_mon_chapter_2}, T \leq K+1}{\leq} 180\gamma^2\exp(-2\gamma\mu T) R^2 \sigma^{p} \sum\limits_{t=0}^{K} \frac{\lambda_t^{2-{p}}}{\exp(-\gamma\mu t)}\notag\\
        & \overset{\eqref{ch2:eq:lambda_SGDA_str_mon_chapter_2}}{\leq} \frac{180\gamma^{p}\exp(-2\gamma\mu T) R^{4-{p}} \sigma^{p} (K+1)\exp(\frac{\gamma\mu{p} K}{2})}{120^{2-{p}} A^{2-{p}}}\notag\\
        & \overset{\eqref{ch2:eq:gamma_SGDA_str_mon_chapter_2}}{\leq} \frac{\exp(-2\gamma\mu T)R^4}{300A}. \label{ch2:eq:bound_1_variances_SGDA_str_mon_chapter_2}
    \end{align}

    \paragraph{Upper bound for $\circledTwo$.} From $E_{T-1}$ it follows that
    \begin{align}
        \circledTwo & \leq 2\gamma \exp(-\gamma\mu (T-1)) \sum\limits_{t=0}^{T-1} \frac{\|u^{t}\|\cdot \|\omega_t^b\|}{\exp(-\gamma\mu t)}\notag\\
        & \overset{\eqref{ch2:eq:eta_t_bound_SGDA_str_mon_chapter_2}, \eqref{ch2:eq:bias_omega_str_mon_chapter_2}}{\leq} 2^{1+{p}}\sqrt{2} \gamma (1+\gamma \ell) \exp(-\gamma\mu (T-1)) R \sigma^{p} \sum\limits_{t=0}^{T-1} \frac{1}{\lambda_t^{{p}-1} \exp(-\gamma\mu t/2)}\notag\\
        & \overset{\eqref{ch2:eq:lambda_SGDA_str_mon_chapter_2}}{\leq} \frac{2^{3+{p}}\cdot 120^{{p}-1}\sqrt{2} \gamma^{p}(1+\gamma \ell) \exp(-\gamma\mu T) \sigma^{{p}} (K+1) \exp\left(\frac{\gamma\mu{p} K}{2}\right) A^{{p}-1}}{R^{{p}-2}} \notag \\
        & \overset{\eqref{ch2:eq:gamma_SGDA_str_mon_chapter_2}}{\leq} \frac{1}{5}\exp(-\gamma\mu T) R^2. \label{ch2:eq:bound_2_SGDA_str_mon_chapter_2}
    \end{align}

    \paragraph{Upper bound for $\circledThree$.} From $E_{T-1}$ it follows that
    \begin{align}
        \circledThree & = 2\gamma^2 \exp(-\gamma\mu (T-1)) \sum\limits_{t=0}^{T-1} \frac{\Exp_{\xi^t}\left[\|\omega_t^u\|^2\right]}{\exp(-\gamma\mu t)} \notag\\
        & \overset{\eqref{ch2:eq:variance_omega_str_mon_chapter_2}}{\leq} 144\gamma^2\exp(-\gamma\mu (T-1)) \sigma^{p}\sum\limits_{t=0}^{T-1} \frac{\lambda_t^{2-{p}}}{\exp(-\gamma\mu t)} \notag\\
        & \overset{\eqref{ch2:eq:lambda_SGDA_str_mon_chapter_2}}{\leq}  \frac{144\gamma^{p} R^{2-{p}}\exp(-\gamma\mu (T-1)) \sigma^{p} (K+1)\exp(\frac{\gamma\mu{p} K}{2})}{120^{2-{p}} A^{2-{p}}} \notag\\
        & \overset{\eqref{ch2:eq:gamma_SGDA_str_mon_chapter_2}}{\leq} \frac{1}{5} \exp(-\gamma\mu T) R^2. \label{ch2:eq:bound_3_SGDA_str_mon_chapter_2}
    \end{align}

    \paragraph{Upper bound for $\circledFour$.} First, we have
    \begin{equation*}
        2\gamma^2 (1-\gamma\mu)^{T-1-t}\Exp_{\xi^t}\left[\|\omega_t^u\|^2  -\Exp_{\xi^t}\left[\|\omega_t^u\|^2\right] \right] = 0.
    \end{equation*}
    Next, sum $\circledFour$ has bounded with probability $1$ terms:
    \begin{align}
        \MoveEqLeft[1] 2\gamma^2 (1-\gamma\mu)^{T-1-t}\left| \|\omega_t^u\|^2  -\Exp_{\xi^t}\left[\|\omega_t^u\|^2\right] \right| \notag \\
        & \overset{\eqref{ch2:eq:omega_magnitude_str_mon_chapter_2}}{\leq} \frac{16\gamma^2 \exp(-\gamma\mu T) \lambda_t^2}{\exp(-\gamma\mu (1+t))}\notag\\
        & \overset{\eqref{ch2:eq:lambda_SGDA_str_mon_chapter_2}}{\leq} \frac{\exp(-\gamma\mu T)R^2}{5A}\notag\\
        & \eqdef c. \label{ch2:eq:SGDA_str_mon_technical_4_1_chapter_2}
    \end{align}
    The summands also have conditional variances
    \begin{equation*}
        \widetilde\sigma_t^2 \eqdef \Exp_{\xi^t}\left[4\gamma^4 (1-\gamma\mu)^{2T-2-2t} \left|\|\omega_t^u\|^2  -\Exp_{\xi^t}\left[\|\omega_t^u\|^2\right] \right|^2 \right]
    \end{equation*}
    that are bounded
    \begin{align}
        \widetilde\sigma_t^2 & \overset{\eqref{ch2:eq:SGDA_str_mon_technical_4_1_chapter_2}}{\leq} \frac{2\gamma^2\exp(-2\gamma\mu T)R^2}{5\exp(-\gamma\mu (1+t))A} \Exp_{\xi^t}\left[ \left|\|\omega_t^u\|^2  -\Exp_{\xi^t}\left[\|\omega_t^u\|^2\right] \right|\right]\notag\\
        & \leq \frac{4\gamma^2\exp(-2\gamma\mu T)R^2}{5\exp(-\gamma\mu (1+t))A} \Exp_{\xi^t}\left[\|\omega_t^u\|^2\right]. \label{ch2:eq:SGDA_str_mon_technical_4_2_chapter_2}
    \end{align}
    In other words, we showed that the summands of $\circledFour$ form a bounded martingale difference sequence with bounded conditional variances $\{\widetilde\sigma_t^2\}_{t = 0}^{T-1}$. Next, we apply Bernstein's inequality (Lemma~\ref{ch2:lem:Bernstein_ineq_chapter_2}) with $X_t$ equal to the summands of $\circledFour$, parameter $c$ as in \eqref{ch2:eq:SGDA_str_mon_technical_4_1_chapter_2}, $b = \tfrac{1}{5}\exp(-\gamma\mu T) R^2$, $H = \tfrac{\exp(-2 \gamma\mu T) R^4}{300A}$:
    \begin{align*}
        &\Prob\left\{|\circledFour| > \frac{1}{5}\exp(-\gamma\mu T) R^2
        \text{ and } \sum\limits_{t=0}^{T-1}\widetilde\sigma_t^2
        \leq \frac{\exp(-2\gamma\mu T) R^4}{300A}\right\}\\
        &\qquad \leq
        2\exp\left(- \frac{b^2}{2H + 2cb/3}\right)
        \leq \frac{\beta}{2(K+1)}.
    \end{align*}
    Equivalently, we have $\Prob\{E_{\circledFour}\} \geq 1 - \frac{\beta}{2(K+1)}$ for
    \begin{equation}
        E_{\circledFour} = \left\{\text{either } \sum\limits_{t=0}^{T-1}\widetilde\sigma_t^2 > \frac{\exp(-2\gamma\mu T) R^4}{300A} \text{ or } |\circledFour| \leq \frac{1}{5}\exp(-\gamma\mu T) R^2\right\}. \label{ch2:eq:bound_4_SGDA_str_mon_chapter_2}
    \end{equation}
    In addition, $E_{T-1}$ implies that
    \begin{align}
        \sum\limits_{t=0}^{T-1}\widetilde\sigma_t^2 & \overset{\eqref{ch2:eq:SGDA_str_mon_technical_4_2_chapter_2}}{\leq} \frac{4\gamma^2\exp(-\gamma\mu (2T-1))R^2}{5A} \sum\limits_{t=0}^{T-1} \frac{\Exp_{\xi^t}\left[\|\omega_t^u\|^2\right]}{\exp(-\gamma\mu t)}\notag\\
        & \overset{\eqref{ch2:eq:variance_omega_str_mon_chapter_2}, T \leq K+1}{\leq} \frac{72\gamma^2\exp(-\gamma\mu (2T-1)) R^2 \sigma^{p}}{5A} \sum\limits_{t=0}^{K} \frac{\lambda_t^{2-{p}}}{\exp(-\gamma\mu t)}\notag\\
        & \overset{\eqref{ch2:eq:lambda_SGDA_str_mon_chapter_2}}{\leq} \frac{72\gamma^{p}\exp(-\gamma\mu (2T-1)) R^{4-{p}} \sigma^{p} (K+1)\exp(\frac{\gamma\mu{p} K}{2})}{5\cdot 120^{2-{p}}A^{3-{p}}} \notag\\
        & \overset{\eqref{ch2:eq:gamma_SGDA_str_mon_chapter_2}}{\leq} \frac{\exp(-2\gamma\mu T)R^4}{300A}. \label{ch2:eq:bound_4_variances_SGDA_str_mon_chapter_2}
    \end{align}

    \paragraph{Upper bound for $\circledFive$.} From $E_{T-1}$ it follows that
    \begin{align}
        \circledFive
        &= 2\gamma^2 \sum\limits_{t=0}^{T-1}
        \exp(-\gamma\mu (T-1-t)) \left(\|\omega_t^b\|^2\right)
        \notag\\
        &\overset{\eqref{ch2:eq:bias_omega_str_mon_chapter_2}}{\leq}
        2\cdot 2^{2{p}}\gamma^2 \exp(-\gamma\mu (T-1)) \sigma^{2{p}}
        \sum\limits_{t=0}^{T-1}
        \frac{1}{\lambda_t^{2{p}-2} \exp(-\gamma\mu t)}
        \notag\\
        &\overset{\eqref{ch2:eq:lambda_SGDA_str_mon_chapter_2},
        T \leq K+1}{\leq}
        \frac{2\cdot 2^{2{p}}\cdot 120^{2{p}-2}\gamma^{2{p}}
        \exp(-\gamma\mu (T-1)) \sigma^{2{p}}
        A^{2{p}-2}}{R^{2{p}-2}}
        \notag\\
        &\quad \times
        \sum\limits_{t=0}^{K}
        \exp\left(\gamma\mu(2{p}-2)\left(1 + \frac{t}{2}\right)\right)
        \exp(\gamma\mu t)
        \notag\\
        &\leq
        \frac{4\cdot 2^{2{p}}\cdot 120^{2{p}-2}\gamma^{2{p}}
        \exp(-\gamma\mu (T-3)) \sigma^{2{p}}
        A^{2{p}-2}}{R^{2{p}-2}}
        \sum\limits_{t=0}^{K} \exp(\gamma\mu {p} t)
        \notag\\
        &\leq
        \frac{4\cdot 2^{2{p}}\cdot 120^{2{p}-2}\gamma^{2{p}}
        \exp(-\gamma\mu (T-3)) \sigma^{2{p}}
        A^{2{p}-2}
        (K+1) \exp(\gamma\mu {p} K)}{R^{2{p}-2}}
        \notag\\
        &\overset{\eqref{ch2:eq:gamma_SGDA_str_mon_chapter_2}}{\leq}
        \frac{1}{5}\exp(-\gamma\mu T) R^2.
        \label{ch2:eq:bound_5_SGDA_str_mon_chapter_2}
    \end{align}

    Now, we have the upper bounds for  $\circledOne, \circledTwo, \circledThree, \circledFour, \circledFive$. In particular, probability event $E_{T-1}$ implies
    \begin{gather*}
        R_T^2 \overset{\eqref{ch2:eq:SGDA_str_mon_12345_bound_chapter_2}}{\leq} \exp(-\gamma\mu T) R^2 + \circledOne + \circledTwo + \circledThree + \circledFour + \circledFive ,\\
        \circledTwo \overset{\eqref{ch2:eq:bound_2_SGDA_str_mon_chapter_2}}{\leq} \frac{1}{5}\exp(-\gamma\mu T)R^2,\quad \circledThree \overset{\eqref{ch2:eq:bound_3_SGDA_str_mon_chapter_2}}{\leq} \frac{1}{5}\exp(-\gamma\mu T)R^2,\\ \circledFive \overset{\eqref{ch2:eq:bound_5_SGDA_str_mon_chapter_2}}{\leq} \frac{1}{5}\exp(-\gamma\mu T)R^2,\\
        \sum\limits_{t=0}^{T-1}\sigma_t^2 \overset{\eqref{ch2:eq:bound_1_variances_SGDA_str_mon_chapter_2}}{\leq}  \frac{\exp(-2\gamma\mu T)R^4}{300A},\quad \sum\limits_{t=0}^{T-1}\widetilde\sigma_t^2 \overset{\eqref{ch2:eq:bound_4_variances_SGDA_str_mon_chapter_2}}{\leq} \frac{\exp(-2\gamma\mu T)R^4}{300A}.
    \end{gather*}
     Moreover, we also have (see \eqref{ch2:eq:bound_1_SGDA_str_mon_chapter_2}, \eqref{ch2:eq:bound_4_SGDA_str_mon_chapter_2} and our induction assumption)
     \begin{gather*}
        \Prob\{E_{T-1}\} \geq 1 - \frac{(T-1)\beta}{K+1},\\
        \Prob\{E_{\circledOne}\} \geq 1 - \frac{\beta}{2(K+1)}, \quad \Prob\{E_{\circledFour}\} \geq 1 - \frac{\beta}{2(K+1)},
    \end{gather*}
    where
    \begin{align}
        E_{\circledOne}& =  \left\{\text{either } \sum\limits_{t=0}^{T-1}\sigma_t^2 > \frac{\exp(- 2\gamma\mu T) R^4}{300A} \text{ or } |\circledOne| \leq \frac{1}{5}\exp(-\gamma\mu T) R^2\right\},\notag\\
        E_{\circledFour}& = \left\{\text{either } \sum\limits_{t=0}^{T-1}\widetilde\sigma_t^2 > \frac{\exp(-2\gamma\mu T) R^4}{300A} \text{ or } |\circledFour| \leq \frac{1}{5}\exp(-\gamma\mu T) R^2\right\}.\notag
    \end{align}
     Thus, probability event $E_{T-1} \cap E_{\circledOne} \cap E_{\circledFour} $ implies
    \begin{align*}
        R_T^2 & \overset{\eqref{ch2:eq:SGDA_str_mon_12345_bound_chapter_2}}{\leq} \exp(-\gamma\mu T) R^2 + \circledOne + \circledTwo + \circledThree + \circledFour + \circledFive\\
        & \leq 2\exp(-\gamma\mu T) R^2,
    \end{align*}
    which is equivalent to \eqref{ch2:eq:induction_inequality_str_mon_SGDA_chapter_2} for $t = T$, and
    \begin{equation}
        \Prob\{E_T\} \geq \Prob\{E_{T-1} \cap E_{\circledOne} \cap E_{\circledFour} \} = 1 - \Prob\{\overline{E}_{T-1} \cup \overline{E}_{\circledOne} \cup \overline{E}_{\circledFour} \} \geq 1 - \frac{T\beta}{K+1}. \notag
    \end{equation}
    This finishes the inductive part of our proof, i.e., for all $k = 0,1,\ldots,K+1$ we have $\Prob\{E_k\} \geq 1 - k\beta/(K+1)$. In particular, for $k = K+1$ we have that with probability at least $1 - \beta$
    \begin{equation}
        \|x^{K+1} - x^*\|^2 \leq 2\exp(-\gamma\mu (K+1))R^2. \notag
    \end{equation}
    Finally, if
    \begin{align*}
        \gamma
        &= \min\left\{\frac{1}{400 \ell A},
        \frac{\ln(B_K)}{\mu(K+1)}\right\},\\
        B_K & = \max\left\{2,
        \frac{(K+1)^{2(p-1)/p}\mu^2R^2}{54000^{2/p}\sigma^2A^{2(p-1)/p}\ln^2(B_K)} \right\} \\
        & = \cO\left(\max\left\{2,
        \frac{K^{2(p-1)/p}\mu^2R^2}{\sigma^2\ln^{2(p-1)/p}(K/\beta)
        \ln^2\left(\mathcal D_K\right)}\right\}\right),
    \end{align*}
    then with probability at least $1-\beta$
    \begin{align*}
        \MoveEqLeft \|x^{K+1} - x^*\|^2
        \leq 2\exp(-\gamma\mu (K+1))R^2\\
        &= 2R^2\max\left\{
        \exp\left(-\frac{\mu(K+1)}{400 \ell A}\right),
        \frac{1}{B_K} \right\}\\
        &= \cO\Biggl(\max\Biggl\{R^2\exp\left(- \frac{\mu K}{\ell \ln(K/\beta)}\right),
        \frac{\sigma^2\ln^{2(p-1)/p}(K/\beta)\ln^2\left(\mathcal D_K\right)}{K^{2(p-1)/p}\mu^2}\Biggr\}\Biggr).
    \end{align*}
    To get $\|x^{K+1} - x^*\|^2 \leq \varepsilon$ with probability at least $1-\beta$ it is sufficient to choose $K$ such that both terms in the maximum above are $\cO(\varepsilon)$. This leads to
    \begin{align*}
         K &= \cO\biggl(\max\biggl\{\frac{\ell}{\mu}\ln\left(\frac{R^2}{\varepsilon}\right)
         \ln\left(\frac{\ell}{\mu \beta}\ln\frac{R^2}{\varepsilon}\right),\\
         &\quad \left(\frac{\sigma^2}{\mu^2\varepsilon}\right)^{p/(2(p-1))}
         \ln \left(\frac{1}{\beta}\left(\frac{\sigma^2}{\mu^2\varepsilon}\right)^{p/(2(p-1))}\right)
         \ln^{p/(p-1)}\left(B_\varepsilon\right)\biggr\}\biggr),
    \end{align*}
    where
    \begin{equation*}
        B_\varepsilon = \max\left\{2, \frac{R^2}{\varepsilon \ln \left(\frac{1}{\beta} \left(\frac{\sigma^2}{\mu^2\varepsilon}\right)^{p/(2(p-1))}\right)}\right\}.
    \end{equation*}
    This concludes the proof.
\end{proof}
  \endgroup

\chapter{Supplementary Material and Proofs for Chapter 3}
  \begingroup
  \thesischaptersetup{3}%
  \stepEightAAppendixColors{3}%

\section{Extra Related Work}\label{ch3:appendix:extra_related_chapter_3}

\paragraph{Non-convex case.} \textcite{li2020high} analyze the high-probability convergence rate of \algname{SGD} for finding first-order stationary points for smooth non-convex unconstrained problems. The first high-probability result under Assumption~\ref{ch3:as:bounded_alpha_moment_chapter_3} for the same class of functions is derived by \textcite{cutkosky2021high}. However, the result of \textcite{cutkosky2021high} relies on the additional assumption that the gradients are bounded. \textcite{sadiev2023high} remove the bounded gradient assumption but derive a slightly worse rate. \textcite{nguyen2023high} improve the result and achieve the same rate as that in \textcite{cutkosky2021high} without assuming boundedness of the gradients. It is worth mentioning that \textcite{cutkosky2021high, sadiev2023high, nguyen2023high} derive their main results for the methods that use gradient clipping.

\paragraph{Gradient clipping} is a very useful algorithmic tool in the training of deep neural networks \parencite{pascanu2013difficulty, goodfellow2016deep}. Gradient clipping also has some good theoretical properties, e.g., it can be useful for minimization of $(L_0, L_1)$-smooth functions \parencite{zhang2020gradient}, in differential privacy \parencite{abadi2016deep}, Byzantine-robustness \parencite{karimireddy2021learning}. Moreover, as we already mentioned, almost all existing high-probability results that do not rely on the light-tailed noise assumption are derived for the methods with clipping. Recently, \textcite{sadiev2023high} theoretically showed that \algname{SGD} has worse high-probability convergence than \algname{clipped-SGD} even when the noise in the gradient has bounded variance.

\clearpage

\section{Auxiliary and Technical Results}\label{ch3:appendix:aux_results_chapter_3}

\paragraph{Bernstein inequality.} In the final stages of our proofs, we need to estimate certain sums of random variables. The main tool that we use to handle such sums is {\it Bernstein inequality for martingale differences} \parencite{bennett1962probability,dzhaparidze2001bernstein,freedman1975tail}.
\begin{lemma}\label{ch3:lem:Bernstein_ineq_chapter_3}
    Let the sequence of random variables $\{X_i\}_{i\ge 1}$ form a martingale difference sequence, i.e.\ $\Exp\left[X_i\mid X_{i-1},\ldots, X_1\right] = 0$ for all $i \ge 1$. Assume that conditional variances $\sigma_i^2\eqdef\Exp\left[X_i^2\mid X_{i-1},\ldots, X_1\right]$ exist and are bounded and also assume that there exists deterministic constant $c>0$ such that $|X_i| \le c$ almost surely for all $i\ge 1$. Then for all $b > 0$, $H > 0$ and $n\ge 1$
    \begin{equation}
        \Prob\left\{\Big|\sum\limits_{i=1}^nX_i\Big| > b \text{ and } \sum\limits_{i=1}^n\sigma_i^2 \le H\right\} \le 2\exp\left(-\frac{b^2}{2H + 2cb/3}\right).
    \end{equation}
\end{lemma}

\paragraph{Impact of clipping on the bias and variance.} The following lemma also helps to handle the aforementioned sums of random variables.

\begin{lemma}[Lemma~5.1 from \textcite{sadiev2023high}]\label{ch3:lem:bias_and_variance_clip_chapter_3}
    Let $X$ be a random vector in $\R^d$ and $\thesischthreetX = \clip(X,\lambda)$. Then, $\|\thesischthreetX - \Exp[\thesischthreetX]\| \leq 2\lambda$. 
    Moreover, if for some $\sigma \geq 0$ and ${p} \in (1,2]$ we have $\Exp[X] = x\in\R^d$, $\Exp[\|X - x\|^{p}] \leq \sigma^{p}$, 
    and $\|x\| \leq \lambda/2$, then
    \begin{align}
        \left\|\Exp[\thesischthreetX] - x\right\| & \leq \frac{2^{p}\sigma^{p}}{\lambda^{{p}-1}}, \label{ch3:eq:bias_X_main_chapter_3}\\
        \Exp\left[\left\|\thesischthreetX - \Exp[\thesischthreetX]\right\|^2\right] & \leq 18 \lambda^{2-{p}} \sigma^{p}. \label{ch3:eq:variance_X_main_chapter_3}
    \end{align}
\end{lemma}

\paragraph{Intuition behind the choice of clipping level in the distributed case.} To better illustrate why we increase clipping level $n$ times, we prove the following lemma.

\begin{lemma}\label{ch3:lem:bias_and_variance_clip_distributed_chapter_3}
    Let $X_1, X_2, \ldots, X_n$ be independent random vectors in $\R^d$ and $\thesischthreetX_i = \clip(X_i,\lambda)$ for all $i \in [n]$. Then, for $\thesischthreetX = (1/n)\sum_{i=1}^n \thesischthreetX_i$ we have  $\|\thesischthreetX - \Exp[\thesischthreetX]\| \leq 2\lambda$. Moreover, if for some $\sigma \geq 0$ and ${p} \in (1,2]$ we have $\Exp[X_i] = x_i\in\R^d$, $\Exp[\|X_i - x_i\|^{p}] \leq \sigma^{p}$, and $\|x_i\| \leq \lambda/2$ for all $i \in [n]$, then for $x = (1/n)\sum_{i=1}^n x_i$ the following inequalities hold
    \begin{align}
        \left\|\Exp[\thesischthreetX] - x\right\| & \leq \frac{2^{p}\sigma^{p}}{\lambda^{{p}-1}}, \label{ch3:eq:bias_X_main_distr_chapter_3}\\
        \Exp\left[\left\|\thesischthreetX - \Exp[\thesischthreetX]\right\|^2\right] & \leq \frac{18 \lambda^{2-{p}} \sigma^{p}}{n}. \label{ch3:eq:variance_X_main_distr_chapter_3}
    \end{align}
\end{lemma}
\begin{proof}
    From Lemma~\ref{ch3:lem:bias_and_variance_clip_chapter_3} we have for all $i \in [n]$ that $\|\thesischthreetX_i - \Exp[\thesischthreetX_i]\| \leq 2\lambda$ and
    \begin{align}
        \left\|\Exp[\thesischthreetX_i] - x_i\right\| & \leq \frac{2^{p}\sigma^{p}}{\lambda^{{p}-1}}, \label{ch3:eq:bias_X_main_1111111_chapter_3}\\
        \Exp\left[\left\|\thesischthreetX_i - \Exp[\thesischthreetX_i]\right\|^2\right] & \leq 18 \lambda^{2-{p}} \sigma^{p}. \label{ch3:eq:variance_X_main_1111111_chapter_3}
    \end{align}
    Jensen's inequality implies
    \begin{align*}
        \left\|\thesischthreetX - \Exp[\thesischthreetX]\right\| & = \left\| \frac{1}{n}\sum\limits_{i=1}^n \left( \thesischthreetX_i - \Exp[\thesischthreetX_i]\right) \right\| \\
        & \leq \frac{1}{n}\sum\limits_{i=1}^n \left\|\thesischthreetX_i - \Exp[\thesischthreetX_i]\right\| \leq 2\lambda,\\
        \left\|\Exp[\thesischthreetX] - x\right\| & = \left\| \frac{1}{n}\sum\limits_{i=1}^n \left(\Exp[\thesischthreetX_i] - x_i\right) \right\| \\
        & \leq \frac{1}{n}\sum\limits_{i=1}^n \left\|\Exp[\thesischthreetX_i] - x_i\right\| \overset{\eqref{ch3:eq:bias_X_main_1111111_chapter_3}}{\leq} \frac{2^{p}\sigma^{p}}{\lambda^{{p}-1}}.
    \end{align*}
    Finally, using the independence of $\thesischthreetX_1,\ldots, \thesischthreetX_n$, we derive
    \begin{align*}
        \Exp\left[\left\|\thesischthreetX - \Exp[\thesischthreetX]\right\|^2\right] & = \Exp\left[\left\| \frac{1}{n}\sum\limits_{i=1}^n \left( \thesischthreetX_i - \Exp[\thesischthreetX_i]\right) \right\|^2\right] \\
        & = \frac{1}{n^2}\sum\limits_{i=1}^n\Exp\left[\left\|\thesischthreetX_i - \Exp[\thesischthreetX_i]\right\|^2\right]\\
        & \overset{\eqref{ch3:eq:variance_X_main_1111111_chapter_3}}{\leq} \frac{18 \lambda^{2-{p}} \sigma^{p}}{n},
    \end{align*}
    which concludes the proof.
\end{proof}

From \eqref{ch3:eq:bias_X_main_distr_chapter_3}-\eqref{ch3:eq:variance_X_main_distr_chapter_3}, we see that number of workers $n$ appears differently in the bound on bias and variance. However, if we replace $\lambda$ with $n\lambda$, then both bounds will transform to \eqref{ch3:eq:bias_X_main_chapter_3}-\eqref{ch3:eq:variance_X_main_chapter_3} respectively with $\sigma^{p} = \sigma^{p}/n^{{p}-1}$ (in other words, bias and variance will have the same dependence on $n$). These observations hint that the complexity bounds for distributed methods should be similar to the ones proven for non-distributed methods (in the unconstrained case) by \textcite{sadiev2023high} up to the replacement of $\sigma^{p}$ with $\sigma^{p}/n^{{p}-1}$. Nevertheless, our analysis of the distributed case does not rely on Lemma~\ref{ch3:lem:bias_and_variance_clip_distributed_chapter_3} and has some important differences with the single-node case (even when $\Psi \equiv 0$).

\paragraph{Useful inequality related to prox-operator.} In the analysis of \algname{DProx-clipped-SGDA-shift}, we use the following standard result.

\begin{lemma}[Theorem 6.39 (iii) from \textcite{beck2017first}]\label{ch3:lem:prox_lemma_chapter_3}
    Let $\Psi$ be a proper lower semicontinuous convex function and $x^+ = \thesischthreeproxkPsi(x)$. Then for all $y \in \R^d$ the following inequality holds:
    \begin{equation*}
        \langle x^+ - x, y - x^+ \rangle \geq \gamma \left(\Psi(x^+) - \Psi(y)\right).
    \end{equation*}
\end{lemma}

\clearpage

\section{Composite Distributed Variational Inequalities}\label{ch3:appendix:VIs_chapter_3}

In this section, we provide further details on the results we obtained for variational inequalities.

\subsection{Assumptions}

In addition to Assumption~\ref{ch3:as:bounded_alpha_moment_VI_chapter_3}, we consider the following assumptions in the case of variational inequalities.

\paragraph{Assumptions on $F_i$.} We make standard assumptions on $\{F_i\}_{i\in [n]}$. The first one is Lipschitzness. 

\begin{assumption}\label{ch3:as:L_Lip_chapter_3}
    There exist some set $Q \subseteq \R^d$ such that operators $F_i$ are $L$-Lipschitz:
    \begin{equation}
        \|F_i(x) - F_i(y)\| \leq L\|x - y\| \quad \forall \, x, y \in Q, \, i\in [n]. \label{ch3:eq:L_Lip_chapter_3}
    \end{equation}
\end{assumption}

Next, for each particular result, we make one or two of the following assumptions.

\begin{assumption}\label{ch3:as:monotonicity_chapter_3}
    There exist some set $Q \subseteq \R^d$ such that $F$ is monotone on $Q$:
    \begin{equation}
        \langle F(x) - F(y), x - y \rangle \geq  0 \quad \forall \, x, y \in Q.\label{ch3:eq:monotonicity_chapter_3}
    \end{equation}
\end{assumption}

\begin{assumption}\label{ch3:as:str_monotonicity_chapter_3}
    There exist some set $Q \subseteq \R^d$ such that $F$ is $(\mu, x^*)$-quasi strongly monotone on $Q$ for some $\mu \geq 0$ and any solution $x^*$ of \eqref{ch3:eq:main_problem_chapter_3}:
    \begin{gather}
        \langle F(x) - F(x^*) , x - x^* \rangle \geq \mu\|x - x^*\|^2, \quad \forall\, x\in Q.
        \label{ch3:eq:str_monotonicity_chapter_3}
    \end{gather}
\end{assumption}

\begin{assumption}\label{ch3:as:star_cocoercivity_chapter_3}
    There exist some set $Q \subseteq \R^d$ such that $\{F_i\}_{i\in[n]}$ are $(\ell, x^*)$-star-cocoercive on $Q$ for some $\ell > 0$ and any solution $x^*$ of \eqref{ch3:eq:main_problem_chapter_3}:
    \begin{gather}
        \|F_i(x) - F_i(x^*)\|^2 \leq \ell\langle F_i(x) - F_i(x^*), x - x^* \rangle, \quad \forall\, x \in Q,\, i\in [n]. \label{ch3:eq:star_cocoercivity_chapter_3}
    \end{gather}
\end{assumption}

\begin{assumption}\label{ch3:as:cocoercivity_chapter_3}
    There exist some set $Q \subseteq \R^d$ such that $F$ is $\ell$-cocoercive on $Q$ for some $\ell > 0$:
    \begin{gather}
        \|F(x) - F(y)\|^2 \leq \ell\langle F(x) - F(y), x - y \rangle, \quad \forall\, x,y \in Q. \label{ch3:eq:cocoercivity_chapter_3}
    \end{gather}
\end{assumption}

Assumption~\ref{ch3:as:monotonicity_chapter_3} is a standard assumption for the literature on VIPs. Quasi-strong monotonicity \parencite{mertikopoulos2019learning, song2020optimistic, loizou2021stochastic} is weaker than standard strong monotonicity\footnote{Operator $F$ is called $\mu$-strongly monotone on $Q$ if $\langle F(x) - F(y), x - y \rangle \geq  \mu\|x - y\|^2$.} and star-cocoercivity is weaker than standard cocoercivity (Assumption~\ref{ch3:as:cocoercivity_chapter_3}), which implies monotonicity and Lipschitzness but not vice versa. Both conditions \eqref{ch3:eq:str_monotonicity_chapter_3} and \eqref{ch3:eq:star_cocoercivity_chapter_3} imply neither monotonicity nor Lipschitzness \parencite{loizou2021stochastic}.

\subsection{\algname{DProx-clipped-SGDA-shift}}

For composite variational inequalities, we start with Distributed \algname{Prox-clipped-SGDA-shift} (\algname{DProx-clipped-SGDA-shift}) that is defined in \eqref{ch3:eq:Dprox_clipped_SGD_shift_chapter_3}-\eqref{ch3:eq:Dprox_clipped_SGD_shift_shifts_chapter_3} with the following change: $\hat\Delta_i^k = \clip\left(F(x^k;\xi_i^k) - h_i^k, \lambda_k\right)$, 
where $\xi_1^k,\ldots, \xi_n^k$ are sampled independently from each other and previous steps. For the proposed method, we derive the following result.
\begin{theorem}[Convergence of \algname{DProx-clipped-SGDA-shift}]\label{ch3:thm:D_prox_clipped_SGDA_main_chapter_3}
    Let $K \geq 1$, $\beta \in (0,1)$, and
    \begin{equation*}
        A = \ln \frac{48n(K+1)}{\beta}, \quad V \geq \|x^0 - x^*\|^2 + \frac{25600\gamma^2A^2}{n^2}\sum_{i=1}^n \|F_i(x^*)\|^2.
    \end{equation*}
    \textbf{Case 1.} Let Assumptions~\ref{ch3:as:bounded_alpha_moment_chapter_3}, \ref{ch3:as:str_monotonicity_chapter_3} with $\mu > 0$, and \ref{ch3:as:star_cocoercivity_chapter_3} hold for $Q = B_{3\sqrt{V}}(x^*)$. Assume that
    \begin{align*}
        0 < \nu & = \cO\left(\frac{1}{\sqrt{n}A}\right),\\
        0 < \gamma & = \cO\left(\min\left\{\frac{1}{\sqrt{n}A\mu}, \frac{1}{\ell A}, \frac{\ln(B_K)}{\mu(K+1)}\right\}\right),\\
        B_K & = \Theta\left(\max\left\{2, \frac{(K+1)^{2({p}-1)/{p}} \mu^2 n^{2({p}-1)/{p}} V}{\sigma^2 A^{2({p}-1)/{p}}\ln^2(B_K)}\right\}\right),\\
        \lambda_k & = \Theta\left(\frac{n\exp(-\gamma\mu(1+k/2))\sqrt{V}}{\gamma A}\right).
    \end{align*}
    \textbf{Case 2.} Let Assumptions~\ref{ch3:as:bounded_alpha_moment_chapter_3}, \ref{ch3:as:monotonicity_chapter_3},  and \ref{ch3:as:star_cocoercivity_chapter_3} hold for $Q = B_{3\sqrt{V}}(x^*)$. Assume that
    \begin{align*}
        \nu & = 0,\\
        0 < \gamma & = \cO\left(\min\left\{\frac{1}{\ell A}, \frac{n^{({p}-1)/{p}}\sqrt{V}}{\sigma K^{1/{p}}A^{({p}-1)/{p}}}\right\}\right),\\
        \lambda_k = \lambda & = \Theta\left(\frac{n\sqrt{V}}{\gamma A}\right).
    \end{align*}
    Then to guarantee $\|x^K - x^*\|^2 \leq \varepsilon$ in \textbf{Case 1} and
    \begin{equation*}
        \gap_{\sqrt{V}}(x_{\text{avg}}^K)
        = \max_{y \in B_{\sqrt{V}}(x^*)}
        \left\{\langle F(y), x_{\text{avg}}^K - y \rangle
        + \Psi(x_{\text{avg}}^K) - \Psi(y)\right\}
        \leq \varepsilon
    \end{equation*}
    in \textbf{Case 2}, with
    $x_{\text{avg}}^K = \frac{1}{K+1}\sum_{k=0}^K x^k$,
    with probability $\geq 1 - \beta$,
    \algname{DProx-clipped-SGDA-shift} requires the following number of iterations/oracle calls per worker:
    \begin{align}
        \textbf{Case 1:}&\quad \widetilde\cO\left(\max\left\{\frac{\ell}{\mu}, \frac{1}{n}\left(\frac{\sigma^2}{\mu^2\varepsilon}\right)^{{p}/(2({p}-1))}\right\}\right), \label{ch3:eq:d_prox_clipped_SGDA_shift_main_result_QSM_chapter_3}\\
        \textbf{Case 2:}&\quad \widetilde\cO\left(\max\left\{\frac{\ell V}{\varepsilon}, \frac{1}{n}\left(\frac{\sigma \sqrt{V}}{\varepsilon}\right)^{{p}/({p}-1)}\right\}\right). \label{ch3:eq:d_prox_clipped_SGDA_shift_main_result_mon_chapter_3}
    \end{align}
\end{theorem}

As in the case of minimization, in the single-node case, the derived results coincide with ones known for \algname{clipped-SGD} in the unconstrained case \parencite{sadiev2023high} up to the difference between $V$ and $\|x^0 - x^*\|^2$. In the distributed case, we also observe the benefits of parallelization.

\subsection{\algname{DProx-clipped-SEG-shift}}

Finally, we propose a distributed version of \algname{clipped-SEG} for composite VIPs (\algname{DProx-clipped-SEG-shift}):
\begin{align}
    \thesischthreetx^k &= \thesischthreeprox_{\gamma\Psi}\left(x^k - \gamma \thesischthreetg^k\right),\;\; \thesischthreetg^k = \frac{1}{n}\sum\limits_{i=1}^n \thesischthreetg_i^k, \notag\\
    \thesischthreetg_i^k &= \thesischthreethh_{i}^k + \tilde\Delta_i^k,\;\; \thesischthreethh_i^{k+1} = \thesischthreethh_i^k + \nu \tilde \Delta_i^k\label{ch3:eq:clipped_SEG_extr_step_chapter_3}\\
    x^{k+1} &= \thesischthreeprox_{\gamma\Psi}\left(x^k - \gamma \thesischthreehg^k\right),\;\; \thesischthreehg^k = \frac{1}{n}\sum\limits_{i=1}^n \thesischthreehg_i^k, \notag\\
    \thesischthreehg_i^k &= \thesischthreehh_{i}^k + \hat\Delta_i^k,\;\; \thesischthreehh_i^{k+1} = \thesischthreehh_i^k + \nu \hat \Delta_i^k \label{ch3:eq:clipped_SEG_upd_step_chapter_3}
\vspace{-0.3cm}
\end{align}
where $\tilde\Delta_i^k = \clip(F_i(x^k;\xi_{1,i}^k) - \thesischthreethh_i^k, \lambda_k)$, $\hat\Delta_i^k = \clip(F_i(\thesischthreetx^k;\xi_{2,i}^k) - \thesischthreehh_i^k, \lambda_k)$ and $\xi_{1,1}^k,\ldots, \xi_{1,n}^k, \xi_{2,1}^k,\ldots, \xi_{2,n}^k$ are sampled independently from each other and previous steps. For the proposed method, we derive the following result.
\begin{theorem}[Convergence of \algname{DProx-clipped-SEG-shift}]\label{ch3:thm:D_prox_clipped_SEG_main_chapter_3}
    Let $K \geq 1$, $\beta \in (0,1)$, and
    \begin{equation*}
        A = \ln \frac{48n(K+1)}{\beta}, \quad V \geq \|x^0 - x^*\|^2 + \frac{409600\gamma^2 A^2}{n^2}\sum_{i=1}^n \|F_i(x^*)\|^2.
    \end{equation*}
    \textbf{Case 1.} Let Assumptions~\ref{ch3:as:bounded_alpha_moment_chapter_3}, \ref{ch3:as:L_Lip_chapter_3}, and \ref{ch3:as:str_monotonicity_chapter_3} with $\mu > 0$ hold for $Q = B_{3\sqrt{V}}(x^*)$. Assume that
    \begin{align*}
        \nu & = \gamma\mu,\\
        0 < \gamma & = \cO\left(\min\left\{\frac{1}{\mu A^2}, \frac{1}{L}, \frac{\sqrt{n}}{L A}, \frac{\ln(B_K)}{\mu(K+1)}\right\}\right),\\
        B_K & = \Theta\left(\max\left\{2, \frac{(K+1)^{2({p}-1)/{p}} \mu^2 n^{2({p}-1)/{p}} V}{\sigma^2 A^{2({p}-1)/{p}}\ln^2(B_K)}\right\}\right),\\
        \lambda_k & = \Theta\left(\frac{n\exp(-\gamma\mu(1+k/4))\sqrt{V}}{\gamma A}\right).
    \end{align*}
    \textbf{Case 2.} Let Assumptions~\ref{ch3:as:bounded_alpha_moment_chapter_3}, \ref{ch3:as:L_Lip_chapter_3}, and \ref{ch3:as:monotonicity_chapter_3} hold for $Q = B_{4n\sqrt{V}}(x^*)$. Assume that
    \begin{align*}
        \nu & = 0,\\
        0 < \gamma & = \cO\left(\min\left\{\frac{1}{L A}, \frac{n^{({p}-1)/{p}}\sqrt{V}}{\sigma K^{1/{p}}A^{({p}-1)/{p}}}\right\}\right),\\
        \lambda_k = \lambda & = \Theta\left(\frac{n\sqrt{V}}{\gamma A}\right).
    \end{align*}
    Then to guarantee $\|x^K - x^*\|^2 \leq \varepsilon$ in \textbf{Case 1} and
    \begin{equation*}
        \gap_{\sqrt{V}}(\thesischthreetx_{\text{avg}}^K)
        = \max_{y \in B_{\sqrt{V}}(x^*)}
        \left\{\langle F(y), \thesischthreetx_{\text{avg}}^K - y \rangle
        + \Psi(\thesischthreetx_{\text{avg}}^K) - \Psi(y)\right\}
        \leq \varepsilon
    \end{equation*}
    in \textbf{Case 2}, with
    $\thesischthreetx_{\text{avg}}^K
    = \frac{1}{K+1}\sum_{k=0}^K \thesischthreetx^k$,
    with probability $\geq 1 - \beta$,
    \algname{DProx-clipped-SEG-shift} requires the following number of iterations/oracle calls per worker:
    \begin{align}
        \textbf{Case 1:}&\quad \widetilde\cO\left(\max\left\{\frac{L}{\mu}, \frac{1}{n}\left(\frac{\sigma^2}{\mu^2\varepsilon}\right)^{{p}/(2({p}-1))}\right\}\right), \label{ch3:eq:d_prox_clipped_SEG_shift_main_result_QSM_chapter_3}\\
        \textbf{Case 2:}&\quad \widetilde\cO\left(\max\left\{\frac{L V}{\varepsilon}, \frac{1}{n}\left(\frac{\sigma \sqrt{V}}{\varepsilon}\right)^{{p}/({p}-1)}\right\}\right). \label{ch3:eq:d_prox_clipped_SEG_shift_main_result_mon_chapter_3}
    \end{align}
\end{theorem}

The main properties of the above result are similar to the ones of the result for \algname{DProx-clipped-SGDA-shift}. The only difference is that the methods (\algname{DProx-clipped-SGDA/SEG-shift}) are analyzed for different classes of problems and, thus, complement each other. According to the known lower bounds, our upper bound \eqref{ch3:eq:d_prox_clipped_SEG_shift_main_result_QSM_chapter_3} has optimal dependence on $\varepsilon$ up to logarithmic factors.

\clearpage

\section{Additional Proofs for \algname{Prox-clipped-SGD-star}}\label{ch3:appendix:prox_clipped_SGD_star_chapter_3}

This section provides the complete formulations of our results for \algname{Prox-clipped-SGD-star} and rigorous proofs. We start with the following result -- a generalization of Lemma E.7 from \textcite{sadiev2023high} to the composite distributed problems.
\begin{lemma}\label{ch3:lem:QSC_to_QSM_chapter_3}
    Consider differentiable function $f: \R^d \to \R$ having a finite-sum structure \eqref{ch3:eq:distributed_problem_chapter_3}. If $f$ satisfies Assumption~\ref{ch3:as:QSC_chapter_3} on some set $Q$ with parameter $\mu$ and $D_f(x,x^*) \geq 0$ for all\footnote{For example $D_f(x,x^*) \geq 0$ when $f$ is convex or when $\Psi(x) = 0$. We notice that Assumption~\ref{ch3:as:L_smoothness_chapter_3} implies $D_f(x,x^*)\geq 0$ since the right-hand side of \eqref{ch3:eq:L_smoothness_cor_2_chapter_3} equals $D_f(x,x^*)$ after averaging.} $x\in Q$, then operator $F(x) = \nabla f(x)$ satisfies Assumption~\ref{ch3:as:str_monotonicity_chapter_3} on $Q$ with parameter $\mu/2$. If $f_1,\ldots,f_n$ satisfy Assumption~\ref{ch3:as:L_smoothness_chapter_3} and \ref{ch3:as:QSC_chapter_3} with $\mu = 0$ on some set $Q$, then operator $F(x) = \nabla f(x)$ satisfies Assumption~\ref{ch3:as:star_cocoercivity_chapter_3} on $Q$ with $\ell = 2L$.
\end{lemma}
\begin{proof}
    Let Assumption~\ref{ch3:as:QSC_chapter_3} hold on some set $Q$ and $D_f(x,x^*) \geq 0$ for all $x\in Q$. Then, averaging inequalities \eqref{ch3:eq:QSC_chapter_3}, we get that for all $x \in Q$
    \begin{equation*}
        f(x^*) \geq f(x) + \langle \nabla f(x), x^* - x \rangle+ \frac{\mu}{2}\|x - x^*\|^2, 
    \end{equation*}
    implying for $F(x) = \nabla f(x)$ that
    \begin{align*}
        \langle F(x) - F(x^*), x  - x^* \rangle & \geq D_f(x,x^*) + \frac{\mu}{2}\|x - x^*\|^2 \geq \frac{\mu}{2}\|x - x^*\|^2,
    \end{align*}
    meaning that Assumption~\ref{ch3:as:str_monotonicity_chapter_3} is satisfied with parameter $\mu/2$.

    It remains to show the second part of the lemma. Let Assumptions~\ref{ch3:as:L_smoothness_chapter_3} and \ref{ch3:as:QSC_chapter_3} with $\mu = 0$ hold on some set $Q$. We need to show that operators $F_i(x) = \nabla f_i(x)$, $i=1,\ldots,n$ satisfy Assumption~\ref{ch3:as:star_cocoercivity_chapter_3} on $Q$ with $\ell = 2L$. Guided by \textcite[Lemma C.6]{gorbunov2022extragradient} and \textcite[Lemma E.7]{sadiev2023high}, we derive
    \begin{align}
        \left\|x - x^* -\frac{1}{L}(F_i(x) - F_i(x^*))\right\|^2 & = \|x - x^*\|^2 - \frac{2}{L}\langle x - x^*, F_i(x) - F_i(x^*) \rangle\notag\\
        &\quad\quad + \frac{1}{L^2}\|F_i(x) - F_i(x^*)\|^2 \label{ch3:eq:cjsidcncnvdf_chapter_3}\\
        & = \|x - x^*\|^2 - \frac{2}{L}\langle x - x^*, \nabla f_i(x) - \nabla f_i(x^*) \rangle\notag\\
        &\quad\quad + \frac{1}{L^2}\|\nabla f_i(x) - \nabla f_i(x^*)\|^2 \notag\\
        & \overset{\eqref{ch3:eq:L_smoothness_cor_2_chapter_3}}{\leq}  \|x - x^*\|^2 - \frac{2}{L} \langle x - x^*, \nabla f_i(x) \rangle \notag\\
        &\quad\quad + \frac{2}{L}\left(f_i(x) - f_i(x^*)\right) \notag\\
        & \overset{\eqref{ch3:eq:QSC_chapter_3}}{\leq} \|x - x^*\|^2. \label{ch3:eq:vhjdfbhjvdfbjv_chapter_3}
    \end{align}
    From \eqref{ch3:eq:cjsidcncnvdf_chapter_3} and \eqref{ch3:eq:vhjdfbhjvdfbjv_chapter_3} we get
    \begin{equation*}
        \|x - x^*\|^2 - \frac{2}{L}\langle x - x^*, F_i(x) - F_i(x^*) \rangle + \frac{1}{L^2}\|F_i(x) - F_i(x^*)\|^2 \leq \|x - x^*\|^2
    \end{equation*}
    that is equivalent to \eqref{ch3:eq:star_cocoercivity_chapter_3} with $\ell = 2L$.
\end{proof}

Therefore, for smooth quasi-strongly convex $f$ such that $D_f(x,x^*)\geq 0$ ($n =1$) we can consider operator $F(x) = \nabla f(x)$ and VI formulation instead. In this case, the method is equivalent to \algname{Prox-clipped-SGDA-star}:
\begin{gather*}
    x^{k+1} = \thesischthreeproxkPsi \left(x^k - \gamma \thesischthreetg^k\right),\;\; \thesischthreetg^k = F(x^*) + \clip\left(F(x^k;\xi^k) - F(x^*), \lambda_k\right)\\
    \thesischthreehg^k = \clip\left(F(x^k;\xi^k) - F(x^*), \lambda_k\right).
\end{gather*}

The following lemma is the main ``optimization'' part of the analysis of \algname{Prox-clipped-SGDA-star}.

\begin{lemma}\label{ch3:lem:optimization_lemma_str_mon_SGDA_chapter_3}
    Let $n=1$, Assumptions~\ref{ch3:as:str_monotonicity_chapter_3}, \ref{ch3:as:star_cocoercivity_chapter_3} hold for $Q = B_{2R}(x^*)$, where $R \geq R_0 \eqdef \|x^0 - x^*\|$, and $0 < \gamma \leq 1/\ell$. If $x^k$ lies in $B_{2R}(x^*)$ for all $k = 0,1,\ldots, K$ for some $K\geq 0$, then the iterates produced by \algname{Prox-clipped-SGDA-star} satisfy
    \begin{align}
        \|x^{K+1} - x^*\|^2
        &\leq (1 - \gamma \mu)^{K+1}\|x^0 - x^*\|^2
        + \gamma^2 \sum\limits_{k=0}^K
        (1-\gamma\mu)^{K-k} \|\omega_k\|^2 \notag\\
        &\quad + 2\gamma \sum\limits_{k=0}^K
        (1-\gamma\mu)^{K-k} \notag\\
        &\qquad \langle x^k - x^* - \gamma (F(x^k)-F(x^*)),
        \omega_k \rangle,
        \label{ch3:eq:optimization_lemma_prox_SGDA_str_mon_chapter_3}\\
        \omega_k
        &\eqdef F(x^k) - F(x^*)-\thesischthreehg^k.
        \label{ch3:eq:omega_k_prox_SGDA_chapter_3}
    \end{align}
\end{lemma}
\begin{proof}
     Using the update rule of \algname{Prox-clipped-SGDA-star}, we obtain
     \begin{align*}
        \MoveEqLeft \|x^{k+1} - x^*\|^2\\
        &= \|\thesischthreeproxkPsi\left(x^k - \gamma\thesischthreetg^k\right) - \thesischthreeproxkPsi\left(x^* - \gamma F(x^*)\right)\|^2\\
        &\le \|x^k - x^* - \gamma(\thesischthreetg^k - F(x^*))\|^2\\
        &= \|x^k - x^*\|^2 -2\gamma \langle x^k - x^*, \thesischthreehg^k \rangle + \gamma^2\|\thesischthreehg^k\|^2\\
        &\overset{\eqref{ch3:eq:omega_k_prox_SGDA_chapter_3}}{=} \|x^k - x^*\|^2 - 2\gamma \langle x^k - x^*, F(x^k)-F(x^*) \rangle - 2\gamma^2\langle F(x^k) - F(x^*), \omega_k\rangle\\
        &\quad  + 2\gamma \langle x^k - x^*, \omega_k \rangle  + \gamma^2\|F(x^k) - F(x^*)\|^2 + \gamma^2\|\omega_k\|^2\\
        &\overset{\eqref{ch3:eq:star_cocoercivity_chapter_3}}{\leq}  \|x^k - x^*\|^2 + 2\gamma \langle x^k - x^*, \omega_k \rangle - 2 \gamma^2 \langle F(x^k) - F(x^*), \omega_k\rangle\\
        &\quad  - 2\gamma \left(1 - \frac{\gamma \ell}{2}\right) \langle x^k - x^*, F(x^k)  -F(x^*)\rangle + \gamma^2\|\omega_k\|^2\\
        &\overset{\eqref{ch3:eq:str_monotonicity_chapter_3}}{\leq}\|x^k - x^*\|^2 + 2\gamma \langle x^k - x^* - \gamma (F(x^k) - F(x^*)), \omega_k \rangle\\
        &\quad  - 2\gamma \mu \left(1 - \frac{\gamma \ell}{2}\right) \|x^k - x^*\|^2 + \gamma^2\|\omega_k\|^2\\
        &\leq(1-\gamma \mu)\|x^k - x^*\|^2 + 2\gamma \langle x^k - x^* - \gamma (F(x^k) - F(x^*)), \omega_k \rangle + \gamma^2\|\omega_k\|^2.
    \end{align*}
    where in the last two inequalities we also use $\gamma \le 1/\ell$. Unrolling the recurrence, we obtain \eqref{ch3:eq:optimization_lemma_prox_SGDA_str_mon_chapter_3}.
\end{proof}

\begin{theorem}\label{ch3:thm:main_result_prox_clipped_SGDA_star_chapter_3}
Let $n =1$, Assumptions~\ref{ch3:as:str_monotonicity_chapter_3}, \ref{ch3:as:star_cocoercivity_chapter_3}, hold for $Q = B_{2R}(x^*) = \{x\in\R^d\mid \|x - x^*\| \leq 2R\}$ for any $x \in B_{2R}(x^*)$, where $R \geq \|x^0 - x^*\|$, and
    \begin{align}
        0< \gamma & \leq \min\left\{\frac{1}{400 \ell\ln(4(K+1)/\beta)}, \frac{\ln(B_K)}{\mu(K+1)}\right\}, \label{ch3:eq:gamma_SGDA_prox_clipped_chapter_3}\\
        B_K & = \max\left\{2, \frac{120^{2(2-p)/p}(K+1)^{2(p-1)/p}\mu^2R^2}{27000^{2/p}\sigma^2\ln^{2(p-1)/p}(4(K+1)/\beta)\ln^2(B_K)} \right\} \label{ch3:eq:B_K_SGDA_prox_clipped_1_chapter_3} \\
        & = \cO\left(\max\left\{2, \frac{Q_K}{\ln^2(\max\{2, Q_K\})} \right\}\right), \label{ch3:eq:B_K_SGDA_prox_clipped_2_chapter_3} \\
        Q_K & \eqdef \frac{K^{2(p-1)/p}\mu^2R^2}{\sigma^2\ln^{2(p-1)/p}(K/\beta)}, \notag\\
        \lambda_k & = \frac{\exp(-\gamma\mu(1 + k/2))R}{120\gamma \ln(4(K+1)/\beta)}, \label{ch3:eq:lambda_SGDA_prox_clipped_chapter_3}
    \end{align}
    for some $K \geq 0$ and $\beta \in (0,1]$ such that $\ln(4(K+1)/\beta) \geq 1$. Then, after $K$ iterations the iterates produced by \algname{Prox-clipped-SGDA-star} with probability at least $1 - \beta$ satisfy 
    \begin{equation}
        \|x^{K+1} - x^*\|^2 \leq 2\exp(-\gamma\mu(K+1))R^2. \label{ch3:eq:main_result_str_mon_SGDA_chapter_3}
    \end{equation}
    In particular, when $\gamma$ equals the minimum from \eqref{ch3:eq:gamma_SGDA_prox_clipped_chapter_3}, then the iterates produced by \algname{Prox-clipped-SGDA-star} after $K$ iterations with probability at least $1-\beta$ satisfy
    \begin{align}
       R_K^2 & = \cO\biggl(\max\biggl\{R^2\exp\left(- \frac{\mu K}{\ell \ln(K/\beta)}\right), \notag\\
        &\qquad \frac{\sigma^2\ln^{2(p-1)/p}(K/\beta)\ln^2(\max\{2, Q_K\})}{K^{2(p-1)/p}\mu^2}\biggr\}\biggr),  \label{ch3:eq:clipped_SGDA_prox_clipped_case_1_chapter_3}
    \end{align}
    meaning that to achieve $R_K^2 = \|x^{K} - x^*\|^2 \leq \varepsilon$ with probability at least $1 - \beta$ \algname{Prox-clipped-SGDA-star} requires
    \begin{align}
        K & = \cO\biggl(\max\biggl\{\frac{\ell}{\mu}\ln\left(\frac{R^2}{\varepsilon}\right)\ln\left(\frac{\ell}{\mu \beta}\ln\frac{R^2}{\varepsilon}\right), \notag\\
        &\qquad \left(\frac{\sigma^2}{\mu^2\varepsilon}\right)^{p/(2(p-1))}\ln \left(\frac{1}{\beta} \left(\frac{\sigma^2}{\mu^2\varepsilon}\right)^{p/(2(p-1))}\right)\ln^{p/(p-1)}\left(B_\varepsilon\right)\biggr\}\biggr)  \label{ch3:eq:clipped_SGDA_prox_clipped_case_complexity_appendix_chapter_3}
    \end{align}
    iterations/oracle calls, where
    \begin{equation*}
        B_\varepsilon = \max\left\{2, \frac{R^2}{\varepsilon \ln \left(\frac{1}{\beta} \left(\frac{\sigma^2}{\mu^2\varepsilon}\right)^{p/(2(p-1))}\right)}\right\}.
    \end{equation*}
\end{theorem}
\begin{proof}
    Let $R_k = \|x^k - x^*\|$ for all $k\geq 0$. Our proof is induction-based: by induction, we show that the iterates of the method stay in some ball around the solution with high probability. To formulate the statement rigorously, we introduce probability event $E_k$ for each $k = 0,1,\ldots, K+1$ as follows: inequalities
    \begin{equation}
        R_t^2 \leq 2 \exp(-\gamma\mu t) R^2 \label{ch3:eq:induction_inequality_prox_clipped_SGDA_star_chapter_3}
    \end{equation}
    hold for $t = 0,1,\ldots,k$ simultaneously. We will prove by induction that $\Prob\{E_k\} \geq  1 - k\beta/(K+1)$ for all $k = 0,1,\ldots,K+1$. The base of the induction follows immediately by the definition of $R$. Next, assume that for $k = T-1 \leq K$ the statement holds: $\Prob\{E_{T-1}\} \geq  1 - (T-1)\beta/(K+1)$. Given this, we need to prove $\Prob\{E_{T}\} \geq  1 - T\beta/(K+1)$. Since $R_t^2 \leq 2\exp(-\gamma\mu t) R^2 \leq 2R^2$, we have $x^t \in B_{2R}(x^*)$ for $t = 0, 1, \ldots, T-1$, where operator $F$ is $\ell$-star-cocoercive. Thus, $E_{T-1}$ implies
    \begin{align}
        \|F(x^t) - F(x^*)\| & \leq \ell\|x^t - x^*\| \overset{\eqref{ch3:eq:induction_inequality_prox_clipped_SGDA_star_chapter_3}}{\leq} \sqrt{2}\ell\exp(- \gamma\mu t/2)R \overset{\eqref{ch3:eq:gamma_SGDA_prox_clipped_chapter_3},\eqref{ch3:eq:lambda_SGDA_prox_clipped_chapter_3}}{\leq} \frac{\lambda_t}{2} \label{ch3:eq:operator_bound_x_t_SGDA_prox_clipped_chapter_3}
    \end{align}
    and
    \begin{align}
        \|\omega_t\|^2 & \leq 2\|F(x^t)-F(x^*)\|^2 + 2\|\thesischthreehg^t\|^2\overset{\eqref{ch3:eq:operator_bound_x_t_SGDA_prox_clipped_chapter_3}}{\leq} \frac{5}{2}\lambda_t^2 \overset{\eqref{ch3:eq:lambda_SGDA_prox_clipped_chapter_3}}{\leq} \frac{\exp(-\gamma\mu t)R^2}{4\gamma^2} \label{ch3:eq:omega_bound_x_t_SGDA_prox_clipped_chapter_3}
    \end{align}
    for all $t = 0, 1, \ldots, T-1$, where we use that $\|a+b\|^2 \leq 2\|a\|^2 + 2\|b\|^2$ holding for all $a,b \in \R^d$. 
    This means that we can apply Lemma~\ref{ch3:lem:optimization_lemma_str_mon_SGDA_chapter_3} and $(1 - \gamma\mu)^T \leq \exp(-\gamma\mu T)$: $E_{T-1}$ implies
    \begin{align}
        R_T^2 & \leq \exp(-\gamma\mu T)R^2 + 2\gamma \sum\limits_{t=0}^{T-1} (1-\gamma\mu)^{T-1-t} \langle x^t - x^* - \gamma (F(x^t)-F(x^*)), \omega_t \rangle\notag\\
        &\quad\quad + \gamma^2 \sum\limits_{t=0}^{T-1} (1-\gamma\mu)^{T-1-t}\|\omega_t\|^2. \notag
    \end{align}
    Before we proceed, we introduce a new notation:
    \begin{gather}
        u^{t} = \begin{cases} \underbrace{x^t - x^* - \gamma (F(x^t)-F(x^*))}_{\hat{\eta}_t},& \text{if } \|\hat{\eta}_t\| \leq \sqrt{2}(1 + \gamma \ell) \exp(- \gamma\mu t/2)R,\\ 0,& \text{otherwise}, \end{cases} \label{ch3:eq:eta_t_SGDA_str_mon_chapter_3}
    \end{gather}
    for $t = 0, 1, \ldots, T-1$. Random vectors $\{u^{t}\}_{t=0}^{T-1}$ are bounded almost surely:
    \begin{equation}
         \|u^{t}\| \leq \sqrt{2}(1 + \gamma \ell)\exp(-\gamma\mu t/2)R \label{ch3:eq:eta_t_bound_SGDA_prox_clipped_chapter_3} 
    \end{equation}
    for all $t = 0, 1, \ldots, T-1$. We also notice that $E_{T-1}$ implies $\|F(x^t)-F(x^*)\| \leq \sqrt{2}\ell\exp(-\gamma\mu t/2)R$ (due to \eqref{ch3:eq:operator_bound_x_t_SGDA_prox_clipped_chapter_3}) and
    \begin{align*}
        \|x^t - x^* - \gamma (F(x^t)-F(x^*))\| & \leq \|x^t - x^*\| + \gamma \|F(x^t)-F(x^*)\|\\
        & \overset{\eqref{ch3:eq:operator_bound_x_t_SGDA_prox_clipped_chapter_3}}{\leq} \sqrt{2}(1 + \gamma \ell)\exp(-\gamma\mu t/2)R
    \end{align*}
    for $t = 0, 1, \ldots, T-1$. Therefore, $E_{T-1}$ implies  $u^{t} = x^t - x^* - \gamma (F(x^t)-F(x^*))$ for all $t = 0, 1, \ldots, T-1$ and from $E_{T-1}$ it follows that
    \begin{align}
        R_T^2 & \leq \exp(-\gamma\mu T)R^2 + 2\gamma \sum\limits_{t=0}^{T-1} (1-\gamma\mu)^{T-1-t} \langle u^{t}, \omega_t \rangle\notag\\
        &\quad\quad + \gamma^2 \sum\limits_{t=0}^{T-1} (1-\gamma\mu)^{T-1-t} \|\omega_t\|^2. \notag
    \end{align}
    For convenience, we define unbiased and biased parts of $ \omega_t$:
    \begin{gather}
        \omega_t^u \eqdef \Exp_{\xi^t}\left[\thesischthreehg^t\right] - \thesischthreehg^t,\quad \omega_t^b \eqdef F(x^t) - F(x^*) - \Exp_{\xi^t}\left[\thesischthreehg^t\right], \label{ch3:eq:omega_unbias_bias_SGDA_str_mon_chapter_3}
    \end{gather}
    for all $t = 0, 1, \ldots, T-1$. By definition we have $\omega_t = \omega_t^u + \omega_t^b$ for all $t = 0, \ldots, T-1$. Therefore, $E_{T-1}$ implies
    \begin{align}
        R_T^2 & \leq \exp(-\gamma\mu T) R^2 + \underbrace{2\gamma \sum\limits_{t=0}^{T-1} (1-\gamma\mu)^{T-1-t} \langle u^{t}, \omega_t^u \rangle}_{\circledOne} \notag\\ &\quad\quad + \underbrace{2\gamma \sum\limits_{t=0}^{T-1} (1-\gamma\mu)^{T-1-t} \langle u^{t}, \omega_t^b \rangle}_{\circledTwo} 
        + \underbrace{2\gamma^2 \sum\limits_{t=0}^{T-1} (1-\gamma\mu)^{T-1-t} \Exp_{\xi^t}\left[\|\omega^u_t\|^2\right]}_{\circledThree}\notag\\
        &\quad\quad + \underbrace{2\gamma^2 \sum\limits_{t=0}^{T-1} (1-\gamma\mu)^{T-1-t}\left( \|\omega^u_t\|^2 -  \Exp_{\xi^t}\left[\|\omega^u_t\|^2\right]\right)}_{\circledFour}\notag\\
        &\quad\quad + \underbrace{2\gamma^2 \sum\limits_{t=0}^{T-1} (1-\gamma\mu)^{T-1-t} \|\omega^b_t\|^2}_{\circledFive}. \label{ch3:eq:SGDA_prox_clipped_12345_bound_chapter_3}
    \end{align}
    where we also use inequality $\|a+b\|^2 \leq 2\|a\|^2 + 2\|b\|^2$ holding for all $a,b \in \R^d$ to upper bound $\|\omega_t\|^2$. To derive high-probability bounds for $\circledOne, \circledTwo, \circledThree, \circledFour, \circledFive$ we need to establish several useful inequalities related to $\omega_{i,t}^u, \omega_{i,t}^b$. First, by definition of clipping
    \begin{equation}
        \|\omega_t^u\| \leq 2\lambda_t.\label{ch3:eq:omega_magnitude_prox_clipped_chapter_3}
    \end{equation}
    Next, $E_{T-1}$ implies that $\|F(x^t)-F(x^*)\| \leq \lambda_t/2$ for all $t = 0, 1, \ldots, T-1$ (see \eqref{ch3:eq:operator_bound_x_t_SGDA_prox_clipped_chapter_3}). Therefore, from Lemma~\ref{ch3:lem:bias_and_variance_clip_chapter_3} we also have that $E_{T-1}$ implies
    \begin{gather}
        \left\|\omega_t^b\right\| \leq \frac{2^{{p}}\sigma^{{p}}}{\lambda_t^{{p}-1}}, \label{ch3:eq:bias_omega_prox_clipped_chapter_3}\\
        \Exp_{\xi^t}\left[\left\|\omega_t^b\right\|^2\right] \leq 18 \lambda_t^{2-{p}} \sigma^{{p}}, \label{ch3:eq:distortion_omega_prox_clipped_chapter_3}\\
        \Exp_{\xi^t}\left[\left\|\omega_t^u\right\|^2\right] \leq 18 \lambda_t^{2-{p}} \sigma^{{p}}, \label{ch3:eq:variance_omega_prox_clipped_chapter_3}
    \end{gather}
    for all $t = 0, 1, \ldots, T-1$.

In the bounds below we repeatedly use $T \le K+1$ to extend sums over $t \le T-1$ to $t \le K$.

\paragraph{Upper bound for $\circledOne$.} To estimate this sum, we will use Bernstein's inequality. The summands have conditional expectations equal to zero:
    \begin{equation*}
        \Exp_{\xi^t}\left[2\gamma (1-\gamma\mu)^{T-1-t} \langle u^{t}, \omega_t^u \rangle\right] = 0.
    \end{equation*}
    Next, the summands are bounded:
    \begin{align}
        |2\gamma (1-\gamma\mu)^{T-1-t} \langle u^{t}, \omega_t^u \rangle | & \leq 2\gamma\exp(-\gamma\mu (T - 1 - t)) \|u^{t}\|\cdot \|\omega_t^u\|\notag\\
        & \overset{\eqref{ch3:eq:eta_t_bound_SGDA_prox_clipped_chapter_3},\eqref{ch3:eq:omega_magnitude_prox_clipped_chapter_3}}{\leq} 4\sqrt{2}\gamma (1 + \gamma \ell) \exp(-\gamma\mu (T - 1 - t/2)) R \lambda_t\notag\\
        & \overset{\eqref{ch3:eq:gamma_SGDA_prox_clipped_chapter_3},\eqref{ch3:eq:lambda_SGDA_prox_clipped_chapter_3}}{\leq} \frac{\exp(-\gamma\mu T)R^2}{5\ln(4(K+1)/\beta)} \eqdef c. \label{ch3:eq:SGDA_str_mon_technical_1_1_chapter_3}
    \end{align}
Finally, conditional variances $\sigma_t^2 \eqdef \Exp_{\xi^t}\left[4\gamma^2 (1-\gamma\mu)^{2T-2-2t} \langle u^{t}, \omega_t^u \rangle^2\right]$ of the summands are bounded:
    \begin{align}
        \sigma_t^2 & \leq \Exp_{\xi^t}\left[4\gamma^2\exp(-\gamma\mu (2T - 2 - 2t)) \|u^{t}\|^2\cdot \|\omega_t^u\|^2\right]\notag\\
        & \overset{\eqref{ch3:eq:eta_t_bound_SGDA_prox_clipped_chapter_3}}{\leq} 8\gamma^2 (1 + \gamma \ell)^2 \exp(-\gamma\mu (2T - 2 - t)) R^2 \Exp_{\xi^t}\left[\|\omega_t^u\|^2\right]\notag\\
        & \overset{\eqref{ch3:eq:gamma_SGDA_prox_clipped_chapter_3}}{\leq} 10\gamma^2\exp(-\gamma\mu (2T - t))R^2 \Exp_{\xi^t}\left[\|\omega_t^u\|^2\right]. \label{ch3:eq:SGDA_str_mon_technical_1_2_chapter_3}
    \end{align}
Applying Bernstein's inequality (Lemma~\ref{ch3:lem:Bernstein_ineq_chapter_3}) with the constant $c$ defined in \eqref{ch3:eq:SGDA_str_mon_technical_1_1_chapter_3}, $b = \exp(-\gamma\mu T) R^2/5$,
    \begin{equation*}
        X_t = 2\gamma (1-\gamma\mu)^{T-1-t} \langle u^{t}, \omega_t^u \rangle, \quad H = \frac{\exp(-2 \gamma\mu T) R^4}{150\ln(4(K+1)/\beta)},
    \end{equation*}
    we get
    \begin{align*}
        \MoveEqLeft[1] \Prob\left\{|\circledOne| > \frac{1}{5}\exp(-\gamma\mu T) R^2 \text{ and } \sum\limits_{t=0}^{T-1}\sigma_t^2 \leq \frac{\exp(- 2\gamma\mu T) R^4}{150\ln(4(K+1)/\beta)}\right\} \\
        & \leq 2\exp\left(- \frac{b^2}{2H + 2cb/3}\right) = \frac{\beta}{2(K+1)}.
    \end{align*}
The above is equivalent to $\Prob\{E_{\circledOne}\} \geq 1 - \frac{\beta}{2(K+1)}$ for
    \begin{equation}
        E_{\circledOne} = \left\{\text{either} \quad \sum\limits_{t=0}^{T-1}\sigma_t^2 > \frac{\exp(- 2\gamma\mu T) R^4}{150\ln(4(K+1)/\beta)}\quad \text{or}\quad |\circledOne| \leq \frac{1}{5}\exp(-\gamma\mu T) R^2\right\}. \label{ch3:eq:bound_1_SGDA_prox_clipped_chapter_3}
    \end{equation}
Moreover, $E_{T-1}$ implies
    \begin{align}
        \sum\limits_{t=0}^{T-1}\sigma_t^2 & \overset{\eqref{ch3:eq:SGDA_str_mon_technical_1_2_chapter_3}}{\leq} 10\gamma^2\exp(- 2\gamma\mu T)R^2\sum\limits_{t=0}^{T-1} \frac{\Exp_{\xi^t}\left[\|\omega_t^u\|^2\right]}{\exp(-\gamma\mu t)}\notag\\ 
        & \overset{\eqref{ch3:eq:variance_omega_prox_clipped_chapter_3}}{\leq} 180\gamma^2\exp(-2\gamma\mu T) R^2 \sigma^{{p}} \sum\limits_{t=0}^{K} \frac{\lambda_t^{2-{p}}}{\exp(-\gamma\mu t)}\notag\\
        & \overset{\eqref{ch3:eq:lambda_SGDA_prox_clipped_chapter_3}}{\leq}\frac{ 180\gamma^{{p}}\exp(-2\gamma\mu T) R^{4-{p}} \sigma^{{p}} (K+1)\exp(\frac{\gamma \mu {p} K}{2})}{120^{2-{p}}\ln^{2-{p}}\frac{4(K+1)}{\beta}}\notag\\
        & \overset{\eqref{ch3:eq:gamma_SGDA_prox_clipped_chapter_3}}{\leq} \frac{\exp(-2\gamma\mu T)R^4}{150\ln(4(K+1)/\beta)}. \label{ch3:eq:bound_1_variances_SGDA_prox_clipped_chapter_3}
    \end{align}

 \paragraph{Upper bound for $\circledTwo$.} Probability event $E_{T-1}$ implies
    \begin{align}
        \circledTwo & \leq 2\gamma \exp(-\gamma\mu (T-1)) \sum\limits_{t=0}^{T-1} \frac{\|u^{t}\|\cdot \|\omega_t^b\|}{\exp(-\gamma\mu t)}\notag\\
        & \overset{\eqref{ch3:eq:eta_t_bound_SGDA_prox_clipped_chapter_3}, \eqref{ch3:eq:bias_omega_prox_clipped_chapter_3}}{\leq} 2^{1+{p}}\sqrt{2} \gamma (1+\gamma \ell) \exp(-\gamma\mu (T-1)) R \sigma^{{p}} \sum\limits_{t=0}^{T-1} \frac{1}{\lambda_t^{{p}-1} \exp(-\gamma\mu t/2)}\notag\\
        & \overset{\eqref{ch3:eq:lambda_SGDA_prox_clipped_chapter_3}}{\leq} \frac{2^{1+{p}}120^{{p}-1}\sqrt{2} \gamma^{{p}}\sigma^{{p}} R^{2-{p}} (1+\gamma \ell)\exp(-\gamma\mu (T-1)) (K+1) \exp\left(\frac{\gamma\mu {p} K}{2} \right)}{\ln^{1-{p}}\tfrac{4(K+1)}{\beta}} \notag \\
        & \overset{\eqref{ch3:eq:gamma_SGDA_prox_clipped_chapter_3}}{\leq} \frac{1}{5}\exp(-\gamma\mu T) R^2. \label{ch3:eq:bound_2_SGDA_prox_clipped_chapter_3}
    \end{align}

\paragraph{Upper bound for $\circledThree$.} Probability event $E_{T-1}$ implies
    \begin{align}
        \circledThree & = 2\gamma^2 \exp(-\gamma\mu (T-1)) \sum\limits_{t=0}^{T-1} \frac{\Exp_{\xi^t}\left[\|\omega_t^u\|^2\right]}{\exp(-\gamma\mu t)} \notag\\
        & \overset{\eqref{ch3:eq:variance_omega_prox_clipped_chapter_3}}{\leq} 36\gamma^2\exp(-\gamma\mu (T-1)) \sigma^{{p}}\sum\limits_{t=0}^{T-1} \frac{\lambda^{2-{p}}_t}{\exp(-\gamma\mu t)} \notag\\
        & \overset{\eqref{ch3:eq:lambda_SGDA_prox_clipped_chapter_3}}{\leq}  \frac{36\gamma^{{p}} R^{2-{p}} \exp(-\gamma\mu (T-1)) \sigma^{{p}} (K+1)\exp(\frac{\gamma\mu{p} K}{2})}{120^{2-{p}}\ln^{2-{p}}\frac{4(K+1)}{\beta}} \notag\\
        & \overset{\eqref{ch3:eq:gamma_SGDA_prox_clipped_chapter_3}}{\leq} \frac{1}{5} \exp(-\gamma\mu T) R^2. \label{ch3:eq:bound_3_SGDA_prox_clipped_chapter_3}
    \end{align}

\paragraph{Upper bound for $\circledFour$.} To estimate this sum, we will use Bernstein's inequality. The summands have conditional expectations equal to zero:
    \begin{equation*}
        2\gamma^2 (1-\gamma\mu)^{T-1-t}\Exp_{\xi^t}\left[\|\omega_t^u\|^2  -\Exp_{\xi^t}\left[\|\omega_t^u\|^2\right] \right] = 0.
    \end{equation*}
    Next, the summands are bounded:
    \begin{align}
        2\gamma^2 (1-\gamma\mu)^{T-1-t}\left| \|\omega_t^u\|^2  -\Exp_{\xi^t}\left[\|\omega_t^u\|^2\right] \right| 
        & \overset{\eqref{ch3:eq:omega_magnitude_prox_clipped_chapter_3}}{\leq} \frac{16\gamma^2 \exp(-\gamma\mu T) \lambda_t^2}{\exp(-\gamma\mu (t+1))}\notag\\
        & \overset{\eqref{ch3:eq:lambda_SGDA_prox_clipped_chapter_3}}{\leq} \frac{\exp(-\gamma\mu T)R^2}{5\ln(4(K+1)/\beta)}\notag\\
        & \eqdef c. \label{ch3:eq:SGDA_prox_clipped_technical_4_1_chapter_3}
    \end{align}
    Finally, conditional variances
    \begin{equation*}
        \widetilde\sigma_t^2 \eqdef \Exp_{\xi^t}\left[4\gamma^4 (1-\gamma\mu)^{2T-2-2t} \left|\|\omega_t^u\|^2  -\Exp_{\xi^t}\left[\|\omega_t^u\|^2\right] \right|^2 \right]
    \end{equation*}
    of the summands are bounded:
    \begin{align}
        \widetilde\sigma_t^2 & \overset{\eqref{ch3:eq:SGDA_prox_clipped_technical_4_1_chapter_3}}{\leq} \frac{2\gamma^2\exp(-2\gamma\mu T)R^2}{5\exp(-\gamma\mu (1+t))\ln(4(K+1)/\beta)} \notag\\
        &\qquad \cdot \Exp_{\xi^t}\left[ \left|\|\omega_t^u\|^2  -\Exp_{\xi^t}\left[\|\omega_t^u\|^2\right] \right|\right]\notag\\
        & \leq \frac{4\gamma^2\exp(-2\gamma\mu T)R^2}{5\exp(-\gamma\mu (1+t))\ln(4(K+1)/\beta)} \notag\\
        &\qquad \cdot \Exp_{\xi^t}\left[\|\omega_t^u\|^2\right]. \label{ch3:eq:SGDA_str_mon_technical_4_2_chapter_3}
    \end{align}
    Applying Bernstein's inequality (Lemma~\ref{ch3:lem:Bernstein_ineq_chapter_3}) with the constant $c$ defined in \eqref{ch3:eq:SGDA_prox_clipped_technical_4_1_chapter_3}, $b = \exp(-\gamma\mu T) R^2/5$,
    \begin{equation*}
        X_t = 2\gamma^2 (1-\gamma\mu)^{T-1-t}\left( \|\omega_t^u\|^2 -\Exp_{\xi^t}\left[\|\omega_t^u\|^2\right]\right), \quad H = \frac{\exp(-2 \gamma\mu T) R^4}{150\ln(4(K+1)/\beta)},
    \end{equation*}
    we get
    \begin{align*}
        \MoveEqLeft[1] \Prob\left\{|\circledFour| > \frac{1}{5}\exp(-\gamma\mu T) R^2 \text{ and } \sum\limits_{t=0}^{T-1}\widetilde\sigma_t^2 \leq \frac{\exp(-2\gamma\mu T) R^4}{150\ln(4(K+1)/\beta)}\right\} \\
        & \leq 2\exp\left(- \frac{b^2}{2H + 2cb/3}\right) = \frac{\beta}{2(K+1)}.
    \end{align*}
    The above is equivalent to $\Prob\{E_{\circledFour}\} \geq 1 - \frac{\beta}{2(K+1)}$ for
    \begin{equation}
        E_{\circledFour} = \left\{\text{either} \quad \sum\limits_{t=0}^{T-1}\widetilde\sigma_t^2 > \frac{\exp(-2\gamma\mu T) R^4}{150\ln(4(K+1)/\beta)}\quad \text{or}\quad |\circledFour| \leq \frac{1}{5}\exp(-\gamma\mu T) R^2\right\}. \label{ch3:eq:bound_4_SGDA_prox_clipped_chapter_3}
    \end{equation}
    Moreover, $E_{T-1}$ implies
    \begin{align}
        \sum\limits_{t=0}^{T-1}\widetilde\sigma_t^2 & \overset{\eqref{ch3:eq:SGDA_str_mon_technical_4_2_chapter_3}}{\leq} \frac{4\gamma^2\exp(-\gamma\mu (2T-1))R^2}{5\ln(4(K+1)/\beta)} \sum\limits_{t=0}^{T-1} \frac{\Exp_{\xi^t}\left[\|\omega_t^u\|^2\right]}{\exp(-\gamma\mu t)}\notag\\ & \overset{\eqref{ch3:eq:variance_omega_prox_clipped_chapter_3}}{\leq} \frac{72\gamma^2\exp(-\gamma\mu (2T-1)) R^2 \sigma^{{p}}}{5\ln(4(K+1)/\beta)} \sum\limits_{t=0}^{K} \frac{\lambda_t^{2-{p}}}{\exp(-\gamma\mu t)}\notag\\
        & \overset{\eqref{ch3:eq:lambda_SGDA_prox_clipped_chapter_3}}{\leq} \frac{72\gamma^{{p}}\exp(-\gamma\mu (2T-1)) R^{4-{p}} \sigma^{{p}} (K+1)\exp(\frac{\gamma\mu{p} K}{2})}{5\cdot120^{2-{p}}\ln^{3-{p}}\tfrac{4(K+1)}{\beta}} \notag\\
        & \overset{\eqref{ch3:eq:gamma_SGDA_prox_clipped_chapter_3}}{\leq} \frac{\exp(-2\gamma\mu T)R^4}{150\ln(4(K+1)/\beta)}. \label{ch3:eq:bound_4_variances_SGDA_prox_clipped_chapter_3}
    \end{align}

\paragraph{Upper bound for $\circledFive$.} Probability event $E_{T-1}$ implies
    \begin{align}
        \circledFive & =  2\gamma^2 \sum\limits_{t=0}^{T-1} \exp(-\gamma\mu (T-1-t)) \|\omega_t^b\|^2\notag\\
        & \overset{\eqref{ch3:eq:bias_omega_prox_clipped_chapter_3}}{\leq} 2\cdot 2^{2 {p}} \gamma^2 \sigma^{2 {p}}  \exp(-\gamma\mu (T-1)) \sum\limits_{t=0}^{T-1} \frac{1}{\lambda_t^{2{p}-2} \exp(-\gamma\mu t)} \notag\\
        & \overset{\eqref{ch3:eq:lambda_SGDA_prox_clipped_chapter_3}}{\leq} \frac{2\cdot 2^{2{p}}120^{2{p}-2}\gamma^{2{p}} \sigma^{2{p}} \exp(-\gamma\mu (T-3))  \ln^{2{p}-2}\tfrac{4(K+1)}{\beta}}{R^{2{p}-2}} \sum\limits_{t=0}^{K} \exp\left(\gamma\mu{p} t\right)\notag\\
        & \leq \frac{2\cdot 2^{2{p}}120^{2{p}-2}\gamma^{2{p}} \sigma^{2{p}} \exp(-\gamma\mu (T-3))  \ln^{2{p}-2}\tfrac{4(K+1)}{\beta}(K+1)\exp(\gamma\mu{p} K)}{R^{2{p}-2}}\notag\\
        & \overset{\eqref{ch3:eq:gamma_SGDA_prox_clipped_chapter_3}}{\leq} \frac{1}{5}\exp(-\gamma\mu T) R^2. \label{ch3:eq:bound_5_SGDA_prox_clipped_chapter_3}
    \end{align}

    That is, we derive the upper bounds for  $\circledOne, \circledTwo, \circledThree, \circledFour, \circledFive$. More precisely, $E_{T-1}$ implies
    \begin{gather*}
        R_T^2 \overset{\eqref{ch3:eq:SGDA_prox_clipped_12345_bound_chapter_3}}{\leq} \exp(-\gamma\mu T) R^2 + \circledOne + \circledTwo + \circledThree + \circledFour + \circledFive ,\\
        \circledTwo \overset{\eqref{ch3:eq:bound_2_SGDA_prox_clipped_chapter_3}}{\leq} \frac{1}{5}\exp(-\gamma\mu T)R^2,\quad \circledThree \overset{\eqref{ch3:eq:bound_3_SGDA_prox_clipped_chapter_3}}{\leq} \frac{1}{5}\exp(-\gamma\mu T)R^2,\quad \circledFive \overset{\eqref{ch3:eq:bound_5_SGDA_prox_clipped_chapter_3}}{\leq} \frac{1}{5}\exp(-\gamma\mu T)R^2,\\
        \sum\limits_{t=0}^{T-1}\sigma_t^2 \overset{\eqref{ch3:eq:bound_1_variances_SGDA_prox_clipped_chapter_3}}{\leq}  \frac{\exp(-2\gamma\mu T)R^4}{150\ln(4(K+1)/\beta)},\quad \sum\limits_{t=0}^{T-1}\widetilde\sigma_t^2 \overset{\eqref{ch3:eq:bound_4_variances_SGDA_prox_clipped_chapter_3}}{\leq} \frac{\exp(-2\gamma\mu T)R^4}{150\ln(4(K+1)/\beta)}.
    \end{gather*}
     In addition, we also establish (see \eqref{ch3:eq:bound_1_SGDA_prox_clipped_chapter_3}, \eqref{ch3:eq:bound_4_SGDA_prox_clipped_chapter_3} and our induction assumption)
     \begin{gather*}
        \Prob\{E_{T-1}\} \geq 1 - \frac{(T-1)\beta}{K+1},\\
        \Prob\{E_{\circledOne}\} \geq 1 - \frac{\beta}{2(K+1)}, \quad \Prob\{E_{\circledFour}\} \geq 1 - \frac{\beta}{2(K+1)}.
    \end{gather*}
    where
    \begin{align}
        E_{\circledOne}& =  \left\{\text{either} \quad \sum\limits_{t=0}^{T-1}\sigma_t^2 > \frac{\exp(- 2\gamma\mu T) R^4}{150\ln(4(K+1)/\beta)}\quad \text{or}\quad |\circledOne| \leq \frac{1}{5}\exp(-\gamma\mu T) R^2\right\},\notag\\
        E_{\circledFour}& = \left\{\text{either} \quad \sum\limits_{t=0}^{T-1}\widetilde\sigma_t^2 > \frac{\exp(-2\gamma\mu T) R^4}{150\ln(4(K+1)/\beta)}\quad \text{or}\quad |\circledFour| \leq \frac{1}{5}\exp(-\gamma\mu T) R^2\right\}.\notag
    \end{align}
    Therefore, probability event $E_{T-1} \cap E_{\circledOne} \cap E_{\circledFour} $ implies
    \begin{align*}
        R_T^2 & \overset{\eqref{ch3:eq:SGDA_prox_clipped_12345_bound_chapter_3}}{\leq} \exp(-\gamma\mu T) R^2 + \circledOne + \circledTwo + \circledThree + \circledFour + \circledFive\\
        & \leq 2\exp(-\gamma\mu T) R^2,
    \end{align*}
    which is equivalent to \eqref{ch3:eq:induction_inequality_prox_clipped_SGDA_star_chapter_3} for $t = T$. Moreover,
    \begin{equation}
        \Prob\{E_T\} \geq \Prob\{E_{T-1} \cap E_{\circledOne} \cap E_{\circledFour} \} = 1 - \Prob\{\overline{E}_{T-1} \cup \overline{E}_{\circledOne} \cup \overline{E}_{\circledFour} \} \geq 1 - \frac{T\beta}{K+1}. \notag
    \end{equation}
    In other words, we showed that $\Prob\{E_k\} \geq 1 - k\beta/(K+1)$ for all $k = 0,1,\ldots,K+1$. For $k = K+1$ we have that with probability at least $1 - \beta$
    \begin{equation}
        \|x^{K+1} - x^*\|^2 \leq 2\exp(-\gamma\mu (K+1))R^2. \notag
    \end{equation}
    
    Finally, if 
    \begin{align*}
        \gamma & = \min\left\{\frac{1}{400 \ell \ln(4(K+1)/\beta)}, \frac{\ln(B_K)}{\mu(K+1)}\right\}, \\
        B_K & = \max\left\{2, \frac{120^{2(2-p)/p}(K+1)^{2(p-1)/p}\mu^2R^2}{27000^{2/p}\sigma^2\ln^{2(p-1)/p}(4(K+1)/\beta)\ln^2(B_K)} \right\} \\
        & = \cO\left(\max\left\{2, \frac{Q_K}{\ln^2(\max\{2, Q_K\})} \right\}\right)
    \end{align*}
    then with probability at least $1-\beta$
    \begin{align*}
        \|x^{K+1} - x^*\|^2 & \leq 2\exp(-\gamma\mu (K+1))R^2\\
        & = 2R^2\max\left\{\exp\left(-\frac{\mu(K+1)}{400 \ell \ln(4(K+1)/\beta)}\right), \frac{1}{B_K} \right\}\\
        & = \cO\biggl(\max\biggl\{R^2\exp\left(- \frac{\mu K}{\ell \ln(K/\beta)}\right), \\
        &\qquad \frac{\sigma^2\ln^{2(p-1)/p}(K/\beta)\ln^2(\max\{2, Q_K\})}{K^{2(p-1)/p}\mu^2}\biggr\}\biggr).
    \end{align*}
    To get $\|x^{K+1} - x^*\|^2 \leq \varepsilon$ with probability at least $1-\beta$, $K$ should be
    \begin{align*}
         K & = \cO\biggl(\max\biggl\{\frac{\ell}{\mu}\ln\left(\frac{R^2}{\varepsilon}\right)\ln\left(\frac{\ell}{\mu \beta}\ln\frac{R^2}{\varepsilon}\right), \\
        &\qquad \left(\frac{\sigma^2}{\mu^2\varepsilon}\right)^{p/(2(p-1))}\ln \left(\frac{1}{\beta} \left(\frac{\sigma^2}{\mu^2\varepsilon}\right)^{p/(2(p-1))}\right)\ln^{p/(p-1)}\left(B_\varepsilon\right)\biggr\}\biggr),
    \end{align*}
    where
    \begin{equation*}
        B_\varepsilon = \max\left\{2, \frac{R^2}{\varepsilon \ln \left(\frac{1}{\beta} \left(\frac{\sigma^2}{\mu^2\varepsilon}\right)^{p/(2(p-1))}\right)}\right\}.
    \end{equation*}
\end{proof}

\clearpage

\section{Additional Proofs for \algname{DProx-clipped-SGD-shift}}\label{ch3:appendix:dprox_clipped_SGD_shift_chapter_3}

In this section, we give the complete formulations of our results for \algname{DProx-clipped-SGD-shift} and rigorous proofs. For the readers' convenience, the method's update rule is repeated below:
\begin{gather*}
    x^{k+1} = \thesischthreeproxkPsi \left(x^k - \gamma \thesischthreetg^k\right),\;\; \text{where}\;\; \thesischthreetg^k = \frac{1}{n}\sum\limits_{i=1}^n \thesischthreetg_i^k,\;\; \thesischthreetg_i^k = h_i^k + \hat\Delta_i^k, \\
    h_i^{k+1} = h_i^k + \nu \hat \Delta_i^k,\quad \hat\Delta_i^k = \clip\left(\nabla f_i(x^k;\xi_i^k) - h_i^k, \lambda_k\right). 
\end{gather*}

\begin{lemma}\label{ch3:lem:main_opt_lemma_prox_clipped_SGD_convex_chapter_3}
    Let Assumptions~\ref{ch3:as:L_smoothness_chapter_3} and \ref{ch3:as:str_cvx_chapter_3} with $\mu = 0$ hold on $Q = B_{3n\sqrt{V}}(x^*)$, where
    \begin{equation*}
        A = \ln\frac{48n(K+1)}{\beta}, \quad V \geq \|x^0 - x^*\|^2 + \frac{36864\gamma^2A^2}{n^2}\sum_{i=1}^n\|\nabla f_i(x^*)\|^2,
    \end{equation*}
    and let stepsize $\gamma$ satisfy $\gamma \leq 1/L$. If $x^{k} \in Q$ for all $k = 0,1,\ldots,K+1$, $K \ge 0$, then after $K$ iterations of \algname{DProx-clipped-SGD-shift} we have
	    \begin{align}
	       2\gamma\left(\Phi(\overline{x}^{K+1}) -\Phi(x^*)\right)
	       &\leq \frac{\|x^0 - x^*\|^2 - \|x^{K+1} - x^*\|^2}{K+1}\notag\\
	       &\quad + \frac{2\gamma}{K+1}\sum\limits_{k=0}^{K}
	       \langle \omega_k, \Hat{x}^k - x^* \rangle\notag\\
	       &\quad + \frac{2\gamma^2}{K+1}\sum\limits_{k=0}^{K}
	       \|\omega_k\|^2 ,\label{ch3:eq:main_opt_lemma_clipped_SGD_star_convex_chapter_3}\\
	       \overline{x}^{K+1}
	       &\eqdef \frac{1}{K+1}\sum\limits_{k=0}^{K} x^{k+1},
	       \label{ch3:eq:x_avg_prox_clipped_SGD_chapter_3}\\
	       \Hat{x}^k
	       &\eqdef \thesischthreeproxkPsi \left(x^k - \gamma \nabla f(x^k)\right),
	       \label{ch3:eq:x_hat_prox_clipped_SGD_chapter_3}\\
	       \omega_k
	       &\eqdef \nabla f(x^k) - \thesischthreetg^k.
	       \label{ch3:eq:omega_k_prox_clipped_SGD_chapter_3}
	    \end{align}
\end{lemma}
\begin{proof}
    Using Lemma~C.2 from \textcite{khaled2020unified} with $p=x^{k+1}, \ y = x^k - \gamma \thesischthreetg^k, \ x=x^k$, we derive for all $k = 0,1, \ldots, K$ that
    \begin{align*}
        2\gamma\left(\Phi(x^{k+1}) - \Phi(x^*)\right) & \leq \|x^{k} - x^*\|^2 - \|x^{k+1} - x^*\|^2\\
        &\quad - 2\gamma \langle \thesischthreetg^k - \nabla f(x^k), x^{k+1}-x^*\rangle.
    \end{align*}
    Next, we obtain the following inequality
    \begin{equation*}
    \begin{aligned}
        &-2\gamma\langle \thesischthreetg^k - \nabla f(x^k),
        x^{k+1} - x^*\rangle\\
        &= -2\gamma \langle \thesischthreetg^k - \nabla f(x^k),
        \Hat{x}^k - x^*\rangle\\
        &\quad +2\gamma \langle \thesischthreetg^k - \nabla f(x^k),
        \Hat{x}^k - x^{k+1}\rangle\\
        &\overset{\eqref{ch3:eq:omega_k_prox_clipped_SGD_chapter_3}}{\leq}
        2\gamma \langle \omega_k, \Hat{x}^k - x^*\rangle\\
        &\quad + 2\gamma \|\thesischthreetg^k - \nabla f(x^k)\|
        \cdot\|\Hat{x}^k - x^{k+1}\|\\
        &\overset{\eqref{ch3:eq:x_hat_prox_clipped_SGD_chapter_3}}{=}
        2\gamma \langle \omega_k, \Hat{x}^k - x^*\rangle\\
        &\quad + 2\gamma \|\thesischthreetg^k - \nabla f(x^k)\|\\
        &\qquad \cdot\|\thesischthreeproxkPsi
        \left(x^k - \gamma \nabla f(x^k)\right)
        - \thesischthreeproxkPsi
        \left(x^k - \gamma \thesischthreetg^k\right)\|\\
        &\overset{\eqref{ch3:eq:omega_k_prox_clipped_SGD_chapter_3}}{\leq}
        2\gamma \langle \omega_k, \Hat{x}^k - x^*\rangle
        + 2 \gamma^2 \|\omega_k\|^2.
    \end{aligned}
    \end{equation*}
    Putting all together we get 
    \begin{align*}
        2\gamma\left(\Phi(x^{k+1}) - \Phi(x^*)\right) & \leq \|x^{k} - x^*\|^2 - \|x^{k+1} - x^*\|^2\\
        &\quad + 2\gamma \langle \omega_k, \Hat{x}^k - x^*\rangle + 2 \gamma^2 \|\omega_k\|^2.
    \end{align*}
    Summing up the above inequalities for $k = 0,1,\ldots, K$, we get
    \begin{align*}
        \MoveEqLeft[1] \frac{2\gamma}{K+1}\sum\limits_{k=0}^{K}\left(\Phi(x^{k+1}) - \Phi(x^*)\right)\\
        & \leq \frac{1}{K+1}\sum\limits_{k=0}^{K}\left(\|x^k - x^*\|^2 - \|x^{k+1} - x^*\|^2\right)\\
        &\quad + \frac{2\gamma}{K+1}\sum\limits_{k=0}^{K} \langle \omega_k, \Hat{x}^k - x^* \rangle + \frac{2\gamma^2}{K+1}\sum\limits_{k=0}^{K} \|\omega_k\|^2\\
        & = \frac{\|x^0 - x^*\|^2 - \|x^{K+1} - x^*\|^2}{K+1}\\
        &\quad + \frac{2\gamma}{K+1}\sum\limits_{k=0}^{K} \langle \omega_k, \Hat{x}^k - x^* \rangle + \frac{2\gamma^2}{K+1}\sum\limits_{k=0}^{K} \|\omega_k\|^2.
    \end{align*}
    Finally, we use the definition of $\overline{x}^{K+1}$ and Jensen's inequality and get the result.
\end{proof}

\begin{theorem}\label{ch3:thm:prox_clipped_SGD_convex_case_chapter_3}
    Let Assumptions~\ref{ch3:as:L_smoothness_chapter_3} and \ref{ch3:as:str_cvx_chapter_3} with $\mu = 0$ hold on $Q = B_{3n\sqrt{V}}(x^*)$, where
    \begin{equation*}
        A = \ln\frac{48n(K+1)}{\beta}, \quad V \geq \|x^0 - x^*\|^2 + \frac{36864\gamma^2A^2}{n^2}\sum_{i=1}^n\|\nabla f_i(x^*)\|^2,
    \end{equation*}
    and $\nu = 0$, $h_1^0 = \ldots = h_n^0 = 0$,
    \begin{align}
        \gamma & \leq \min\left\{\frac{1}{360LA}, \frac{R\sqrt{n}}{192A \zeta_*}, \frac{\sqrt{V}n^{(p-1)/p}}{27^{1/p}\cdot 48 \sigma  (K+1)^{1/p}A^{(p-1)/p}}\right\},\label{ch3:eq:prox_clipped_SGD_step_size_cvx_chapter_3}\\
        \lambda_{k} & = \lambda = \frac{n\sqrt{V}}{48 \gamma A}, \label{ch3:eq:prox_clipped_SGD_clipping_level_cvx_chapter_3}
    \end{align}
    where $R \geq \|x^0 - x^*\|$, $\zeta_* \geq \sqrt{(1/n)\sum_{i=1}^n \|\nabla f_i(x^*)\|^2}$, $K \geq 0$ and $\beta \in (0,1]$. Then, after $K+1$ iterations of \algname{DProx-clipped-SGD-shift} the iterates with probability at least $1 - \beta$ satisfy
    \begin{equation}
        \Phi(\overline{x}^{K+1}) - \Phi(x^*) \leq \frac{V}{\gamma(K+1)} \quad \text{and} \quad  \{x^k\}_{k=0}^{K+1} \subseteq B_{3n\sqrt{V}}(x^*). \label{ch3:eq:prox_clipped_SGD_convex_case_chapter_3}
    \end{equation}
    In particular, the second entry of the minimum in \eqref{ch3:eq:prox_clipped_SGD_step_size_cvx_chapter_3} makes $V = 2R^2$ admissible, and when $\gamma$ equals the minimum from \eqref{ch3:eq:prox_clipped_SGD_step_size_cvx_chapter_3}, then the iterates produced by \algname{DProx-clipped-SGD-shift} after $K+1$ iterations with probability at least $1-\beta$ satisfy
    \begin{align}
        \Phi(\overline{x}^{K+1}) - \Phi(x^*) = \cO\biggl(\max\biggl\{ & \frac{LR^2 \ln(nK/\beta)}{K}, \frac{R\zeta_* \ln(nK/\beta)}{\sqrt{n}K}, \notag\\
        & \frac{\sigma R\ln^{(p-1)/p}(nK/\beta)}{n^{(p-1)/p}K^{(p-1)/p}}\biggr\}\biggr), \label{ch3:eq:prox_clipped_SGD_convex_case_2_chapter_3}
    \end{align}
    meaning that to achieve $\Phi(\overline{x}^{K+1}) - \Phi(x^*) \leq \varepsilon$ with probability at least $1 - \beta$ \algname{DProx-clipped-SGD-shift} requires
    \begin{align}
        K = \cO\biggl(\max\biggl\{ & \frac{LR^2}{\varepsilon}\ln\frac{nLR^2}{\beta\varepsilon}, \frac{R\zeta_*}{\sqrt{n}\varepsilon}\ln\frac{\sqrt{n}R\zeta_*}{\beta\varepsilon}, \notag\\
        & \left(\frac{\sigma R}{\varepsilon n^{(p-1)/p}}\right)^{p/(p-1)} \ln \left(\frac{1}{\beta} \left(\frac{\sigma R}{\varepsilon}\right)^{p/(p-1)}\right) \biggr\}\biggr) \label{ch3:eq:prox_clipped_SGD_convex_case_complexity_chapter_3}
    \end{align}
    iterations/oracle calls.
\end{theorem}
\begin{proof}
 The key idea behind the proof is similar to the one used by \textcite{gorbunov2022clipped, sadiev2023high}: we prove by induction that the iterates do not leave some ball and the sums decrease as $1/(K+1)$. To formulate the statement rigorously, we introduce probability event $E_k$ for each $k = 0,1,\ldots, K+1$ as follows: inequalities
    \begin{equation}
        \underbrace{\|x^0 - x^*\|^2 + 2\gamma\sum\limits_{l=0}^{t-1} \langle \omega_l, \Hat{x}^l - x^* \rangle  + 2\gamma^2\sum\limits_{l=0}^{t-1} \|\omega_l\|^2}_{A_t} \leq 2V, \label{ch3:eq:prox_clipped_SGD_convex_induction_inequality_1_chapter_3}
    \end{equation}
    \begin{gather}
        \left\|\frac{\gamma}{n}\sum\limits_{i=1}^{r-1}\omega_{i,t-1}^u\right\| \leq \frac{\sqrt{V}}{2}
    \end{gather}
hold for $t = 0,1,\ldots,k$ and $r = 1, 2, \ldots, n$ simultaneously, where
    \begin{gather}
        \omega_l = \omega_l^u + \omega_l^b, \label{ch3:eq:prox_clipped_SGD_technical_4_full_theta_omega_chapter_3}\\
        \omega_l^u \eqdef \frac{1}{n}\sum\limits_{i=1}^n \omega_{i,l}^u,\quad \omega_l^b \eqdef \frac{1}{n}\sum\limits_{i=1}^n \omega_{i,l}^b, \\
        \omega_{i,l}^u \eqdef \Exp_{\xi_{i}^l}\left[\thesischthreetg_{i}^l\right] - \thesischthreetg_{i}^l,\quad \omega_{i,l}^b \eqdef \nabla f_i(x^l) - \Exp_{\xi_{i}^l}\left[\thesischthreetg_{i}^l\right] \quad \forall\; i\in [n]. \label{ch3:eq:prox_clipped_SGD_convex_theta_u_b_chapter_3}
    \end{gather}
    We will prove by induction that $\Prob\{E_k\} \geq  1 - k\beta/(K+1)$ for all $k = 0,1,\ldots,K+1$. The base of induction follows immediately: $\|x^0 - x^*\|^2 \leq V < 2V$ and for $k = 0$ we have  $\|\frac{\gamma}{n}\sum_{i=1}^{r-1} \omega_{i,k-1}^u\| = 0$ since $\omega_{i,-1}^u = 0$. Next, we assume that the statement holds for $k = T-1 \leq K$, i.e., $\Prob\{E_{T-1}\} \geq  1 - (T-1)\beta/(K+1)$. Let us show that it also holds for $k = T$, i.e., $\Prob\{E_{T}\} \geq  1 - T\beta/(K+1)$.
    
    To proceed, we need to show that $E_{T-1}$ implies $\|x^t - x^*\| \leq 3n\sqrt{V}$ for all $t = 0,1, \ldots, T$. First, for $t = 0,1, \ldots, T-1$ probability event $E_{T-1}$ implies (in view of, $\Phi(\overline{x}^t) - \Phi(x^*) \geq 0$)
    \begin{equation}
        \|x^t - x^*\|^2 \overset{\eqref{ch3:eq:main_opt_lemma_clipped_SGD_star_convex_chapter_3}}{\leq} A_t \overset{\eqref{ch3:eq:prox_clipped_SGD_convex_induction_inequality_1_chapter_3}}{\leq} 2V. \label{ch3:eq:jsknfdbvufbvudfu_chapter_3}
    \end{equation}
    Next, by definition of $V$ we have
    \begin{align}
        \|\nabla f(x^*)\| & = \sqrt{\|\nabla f(x^*)\|^2} \leq \sqrt{\sum\limits_{i=1}^n \|\nabla f_i(x^*)\|^2} \notag \\
        & \leq \frac{n\sqrt{V}}{192\gamma A}. \label{ch3:eq:msodfnvidvidbsi_chapter_3}
    \end{align}
    Then, for $t = T$ we have that $E_{T-1}$ implies
    \begin{align*}
        \|x^T - x^*\| & = \|\thesischthreeproxkPsi (x^{T-1} - \gamma \thesischthreetg^{T-1}) - \thesischthreeproxkPsi (x^* - \gamma \nabla f(x^*))\|\\
        & \leq \|x^{T-1} - \gamma \thesischthreetg^{T-1} - x^* + \gamma \nabla f(x^*)\| \\
        & \leq \|x^{T-1} - x^*\| + \gamma\|\thesischthreetg^{T-1}\| + \gamma\|\nabla f(x^*)\|\\
        & \overset{\eqref{ch3:eq:jsknfdbvufbvudfu_chapter_3},\eqref{ch3:eq:msodfnvidvidbsi_chapter_3}}{\leq} \left(\sqrt{2} + \frac{n}{192A}\right)\sqrt{V} + \gamma\lambda \overset{\eqref{ch3:eq:prox_clipped_SGD_clipping_level_cvx_chapter_3}}{\leq} 3n\sqrt{V}.
    \end{align*}

    This means that $E_{T-1}$ implies $x^t \in B_{3n\sqrt{V}}(x^*)$ for $t = 0,1, \ldots, T$ and we can apply Lemma~\ref{ch3:lem:main_opt_lemma_prox_clipped_SGD_convex_chapter_3}: $E_{T-1}$ implies
     \begin{align}
        2\gamma \left(\Phi(\overline{x}^{T}) -\Phi(x^*)\right) & \leq \frac{\|x^0 - x^*\|^2 - \|x^T - x^*\|^2}{T}\notag\\
       &\quad\quad + \frac{2\gamma}{T}\sum\limits_{l=0}^{T-1} \langle \omega_l, \Hat{x}^l - x^* \rangle  + \frac{2\gamma^2}{T}\sum\limits_{l=0}^{T-1} \|\omega_l\|^2 \notag\\
       & \leq \frac{A_T}{T}.\label{ch3:eq:prox_clipped_SGD_convex_technical_1_chapter_3}
    \end{align}
    Before we proceed, we introduce a new notation:
    \begin{equation}
        u^{t} = \begin{cases} \Hat{x}^t - x^*,& \text{if } \|\Hat{x}^t - x^*\| \leq 2\sqrt{V},\\ 0,&\text{otherwise}, \end{cases} \notag
    \end{equation}
    for all $t = 0, 1, \ldots, T-1$. Random vectors $\{u^{t}\}_{t=0}^{T-1}$ are bounded almost surely:
    \begin{equation}
        \|u^{t}\| \leq 2\sqrt{V}. \label{ch3:eq:prox_clipped_SGD_convex_technical_6_chapter_3}
    \end{equation}
    for all $t = 0, 1, \ldots, T-1$. In addition, $E_{T-1}$ implies for all $t = 0, 1, \ldots, T-1$ that
    \begin{align*}
        \|\Hat{x}^{t} - x^*\| & = \|\thesischthreeproxkPsi \left(x^{t} - \gamma \nabla f(x^{t})\right) - \thesischthreeproxkPsi \left(x^* - \gamma \nabla f(x^*)\right)\|\\
        & \leq \|x^{t} - x^* - \gamma(\nabla f(x^{t}) - \nabla f(x^*))\| \\
        & \leq \|x^{t} - x^*\| + \gamma \|\nabla f(x^{t}) - \nabla f(x^*)\|\\
        & \overset{\eqref{ch3:eq:L_smoothness_chapter_3}}{\leq} (1+L\gamma)\|x^{t} - x^*\| \overset{\eqref{ch3:eq:prox_clipped_SGD_step_size_cvx_chapter_3}}{\leq} \frac{361}{360}\|x^{t} - x^*\| \overset{\eqref{ch3:eq:jsknfdbvufbvudfu_chapter_3}}{\leq} 2\sqrt{V}.
    \end{align*}
    meaning that $u^{t} = \Hat{x}^t - x^*$ follows from $E_{T-1}$ for all $t = 0, 1, \ldots, T-1$. Thus, $E_{T-1}$ implies
    \begin{align}
        A_T & \overset{\eqref{ch3:eq:prox_clipped_SGD_convex_induction_inequality_1_chapter_3}}{=} \|x^0 - x^*\|^2 + 2\gamma\sum\limits_{l=0}^{T-1} \langle \omega_l, \Hat{x}^l - x^* \rangle  + 2\gamma^2\sum\limits_{l=0}^{T-1} \|\omega_l\|^2\notag\\
        & \leq V + 2\gamma\sum\limits_{l=0}^{T-1} \langle \omega_l, u^{l} \rangle  + 2\gamma^2\sum\limits_{l=0}^{T-1} \|\omega_l\|^2.\label{ch3:eq:prox_clipped_SGD_convex_technical_2_chapter_3}
    \end{align}
    Using the notation from \eqref{ch3:eq:prox_clipped_SGD_technical_4_full_theta_omega_chapter_3}-\eqref{ch3:eq:prox_clipped_SGD_convex_theta_u_b_chapter_3}, we can rewrite $\|\omega_l\|^2$ as
    \begin{align}
        \|\omega_l\|^2 & \leq  2\|\omega_l^u\|^2 + 2\|\omega_l^b\|^2 \notag \\
        & = \frac{2}{n^2}\left\|\sum\limits_{i=1}^n \omega_{i,l}^u\right\|^2 + 2\|\omega_l^b\|^2 \notag \\
        & = \frac{2}{n^2}\sum\limits_{i=1}^n\|\omega_{i,l}^u\|^2 + \frac{4}{n^2}\sum\limits_{j=2}^n\left\langle \sum\limits_{i=1}^{j-1} \omega_{i,l}^u, \omega_{j,l}^u \right\rangle + 2\|\omega_l^b\|^2 .\label{ch3:eq:omega_bound_x_t_prox_clipped_SGD_convex_chapter_3}
    \end{align}
    Putting all together, we obtain that $E_{T-1}$ implies
	    \begin{align}
	        A_T
	        &\leq V
	        + \underbrace{\frac{2\gamma}{n}\sum\limits_{l=0}^{T-1}\sum\limits_{i=1}^n
	        \langle \omega_{i,l}^u, u^{l}\rangle}_{\circledOne}
	        + \underbrace{2\gamma\sum\limits_{l=0}^{T-1}
	        \langle \omega_l^b, u^{l}\rangle}_{\circledTwo}\notag\\
	        &\quad + \underbrace{\frac{4\gamma^2}{n^2}
	        \sum\limits_{l=0}^{T-1}\sum\limits_{i=1}^n
	        \left(\left\|\omega_{i,l}^u\right\|^2
	        - \Exp_{\xi^l_i}\left[\left\|\omega_{i,l}^u\right\|^2\right]\right)}_{\circledThree}\notag\\
	        &\quad + \underbrace{\frac{4\gamma^2}{n^2}
	        \sum\limits_{l=0}^{T-1}\sum\limits_{i=1}^n
	        \Exp_{\xi^l_i}\left[\left\|\omega_{i,l}^u\right\|^2\right]}_{\circledFour}
	        + \underbrace{4\gamma^2\sum\limits_{l=0}^{T-1}
	        \|\omega_{l}^b\|^2}_{\circledFive}\notag\\
	        &\quad + \underbrace{\frac{8\gamma^2}{n^2}
	        \sum\limits_{l=0}^{T-1}\sum\limits_{j=2}^n
	        \left\langle \sum\limits_{i=1}^{j-1}\omega_{i,l}^u,
	        \omega_{j,l}^u\right\rangle}_{\circledSix}.
	        \label{ch3:eq:prox_clipped_SGD_convex_technical_7_chapter_3}
	    \end{align}
    To finish the proof, it remains to estimate $\circledOne, \circledTwo, \circledThree, \circledFour, \circledFive, \circledSix$ with high probability. More precisely, the goal is to prove that $\circledOne + \circledTwo + \circledThree + \circledFour + \circledFive + \circledSix \leq V$ with high probability. Before we proceed, we need to derive several useful inequalities related to $ \omega_{i,l}^u, \omega_{l}^b$. First of all, we have
    \begin{equation}
       \|\omega_{i,l}^u\| \leq 2\lambda \label{ch3:eq:prox_clipped_SGD_convex_norm_theta_u_bound_chapter_3}
    \end{equation}
    by definition of the clipping operator. Next, probability event $E_{T-1}$ implies 
    \begin{align}
        \|\nabla f_i(x^{t})\| & \leq \|\nabla f_i(x^t) - \nabla f_i(x^*)\| + \|\nabla f_i(x^*)\| \notag \\
        & \overset{\eqref{ch3:eq:L_smoothness_chapter_3}}{\leq} L\|x^t - x^*\| + \sqrt{\sum\limits_{i=1}^n\|\nabla f_i(x^*)\|^2} \notag\\ 
        & \leq \sqrt{2}L\sqrt{V} + \frac{n\sqrt{V}}{192\gamma A} \notag \\
        & \leq \frac{n\sqrt{V}}{96\gamma A} \leq \frac{\lambda}{2}. \label{ch3:eq:prox_clipped_SGD_convex_technical_4_chapter_3}
    \end{align}for $t = 0, 1, \ldots, T-1$ and $i\in [n]$. Therefore, Lemma~\ref{ch3:lem:bias_and_variance_clip_chapter_3} and $E_{T-1}$ imply
    \begin{gather}
        \left\|\omega_l^b\right\| \leq \frac{1}{n}\sum\limits_{i=1}^n \|\omega_{i,l}^b\| \leq \frac{2^{p}\sigma^{p}}{\lambda^{{p}-1}}, \label{ch3:eq:prox_clipped_SGD_convex_norm_theta_b_bound_chapter_3}\\
        \Exp_{\xi_{i}^l}\left[\left\|\omega_{i,l}^u\right\|^2\right] \leq 18 \lambda^{2-{p}}\sigma^{p}, \label{ch3:eq:prox_clipped_SGD_convex_second_moment_theta_u_bound_chapter_3}
    \end{gather}
    for all $l = 0,1, \ldots, T-1$ and $i\in [n]$.

\textbf{Upper bound for $\circledOne$.} 
To estimate this sum, we will use Bernstein's inequality. The summands have conditional expectations equal to zero:
    \begin{equation}
        \Exp_{\xi^l_i}\left[\frac{2\gamma}{n}\langle\omega_{i,l}^u, u^{l}\rangle\right] = \frac{2\gamma}{n} \left\langle u^{l}, \Exp_{\xi_{i}^l}[\omega_{i,l}^u] \right\rangle =  0. \notag
    \end{equation}
    Moreover, for all $l = 0,\ldots, T-1$ random vectors $\{\omega_{i,l}^u\}_{i=1}^n$ are independent. Thus, sequence $\left\{\frac{2\gamma}{n} \langle u^{l}, \omega_{i,l}^u \rangle\right\}_{l,i = 0,1}^{T-1,n}$ is a martingale difference sequence. Next, the summands are bounded:
    \begin{align}
\left|\frac{2\gamma}{n}\left\la \omega_{i,l}^u, u^{l}\right\ra\right| & \leq \frac{2\gamma}{n} \|\omega_{i,l}^u\| \cdot \|u^{l}\| \overset{\eqref{ch3:eq:prox_clipped_SGD_convex_technical_6_chapter_3},\eqref{ch3:eq:prox_clipped_SGD_convex_norm_theta_u_bound_chapter_3}}{\leq} \frac{8\gamma \lambda \sqrt{V}}{n} \notag \\
& \overset{\eqref{ch3:eq:prox_clipped_SGD_clipping_level_cvx_chapter_3}}{=} \frac{V}{6A} \eqdef c. \label{ch3:eq:prox_clipped_SGD_convex_technical_8_chapter_3}
    \end{align}
    Finally, conditional variances $\sigma_{i,l}^2 \eqdef \Exp_{\xi^l_i}[\frac{4\gamma^2}{n^2}\langle\omega_{i,l}^u, u^{l}\rangle^2]$ of the summands are bounded:
    \begin{align}
\sigma_{i,t}^2 & \leq \Exp_{\xi^t_i}\left[\frac{4\gamma^2}{n^2}\|\omega_{i,t}^u\|^2\cdot \|u^{t}\|^2\right] \notag \\
& \overset{\eqref{ch3:eq:prox_clipped_SGD_convex_technical_6_chapter_3}}{\leq} \frac{16\gamma^2 V}{n^2} \Exp_{\xi^t_i}\left[\|\omega_{i,t}^u\|^2\right]. \label{ch3:eq:prox_clipped_SGD_convex_technical_9_chapter_3}
    \end{align}
    Applying Bernstein's inequality (Lemma~\ref{ch3:lem:Bernstein_ineq_chapter_3}) with $X_{i,l} = \frac{2\gamma}{n} \langle u^{l}, \omega_{i,l}^u \rangle$, constant $c$ defined in \eqref{ch3:eq:prox_clipped_SGD_convex_technical_8_chapter_3}, $b = \frac{V}{6}$, $H = \frac{V^2}{216A}$, we get
    \begin{align*}
\MoveEqLeft[1] \Prob\left\{|\circledOne| > \frac{V}{6}\quad \text{and}\quad \sum\limits_{l=0}^{T-1}\sum\limits_{i=1}^n \sigma_{i,l}^2 \leq \frac{V^2}{216A}\right\} \\
& \leq 2\exp\left(- \frac{b^2}{2H + 2cb/3}\right) = \frac{\beta}{24n(K+1)}.
    \end{align*}
	    The above is equivalent to 
	    \begin{equation}
	        \begin{aligned}
	            \Prob\left\{ E_{\circledOne} \right\}
	            &\geq 1 - \frac{\beta}{24n(K+1)},\quad \text{for}\quad\\
	            E_{\circledOne}
	            &= \left\{\text{either } \sum\limits_{l=0}^{T-1}\sum\limits_{i=1}^n \sigma_{i,l}^2 > \frac{V^2}{216A} \text{ or } |\circledOne| \leq \frac{V}{6}\right\}.
	        \end{aligned}
	        \label{ch3:eq:prox_clipped_SGD_convex_sum_1_upper_bound_chapter_3}
	    \end{equation}
    Moreover, $E_{T-1}$ implies
    \begin{align}
        \sum\limits_{l=0}^{T-1}\sum\limits_{i=1}^n \sigma_{i,l}^2 & \overset{\eqref{ch3:eq:prox_clipped_SGD_convex_technical_9_chapter_3}}{\leq} \frac{16\gamma^2 V}{n^2} \sum\limits_{l=0}^{T-1} \sum\limits_{i=1}^n \Exp_{\xi^l_i}\left[\|\omega_{i,l}^u\|^2\right] \notag \\
        & \overset{\eqref{ch3:eq:prox_clipped_SGD_convex_second_moment_theta_u_bound_chapter_3}}{\leq} \frac{288\gamma^2 V \sigma^{{p}}T \lambda^{2-{p}}}{n} \notag\\
        & \overset{\eqref{ch3:eq:prox_clipped_SGD_clipping_level_cvx_chapter_3}}{=} \frac{ 48^{{p}}\sqrt{V}^{4-{p}}\sigma^{p} T\gamma^{{p}}}{8 n^{{p}-1}A^{2-{p}}} \notag \\
        & \overset{\eqref{ch3:eq:prox_clipped_SGD_step_size_cvx_chapter_3}}{\leq} \frac{V^2}{216A}. \label{ch3:eq:prox_clipped_SGD_convex_sum_1_variance_bound_chapter_3}
    \end{align}
    
 \textbf{Upper bound for $\circledTwo$.} Probability event $E_{T-1}$ implies
    \begin{align}
        \circledTwo & = 2\gamma\sum\limits_{l=0}^{T-1}\langle \omega_l^b, u^{l} \rangle \leq 2\gamma\sum\limits_{l=0}^{T-1}\|\omega_{l}^b\|\cdot \|u^{l}\| \notag \\
        & \overset{\eqref{ch3:eq:prox_clipped_SGD_convex_technical_6_chapter_3},\eqref{ch3:eq:prox_clipped_SGD_convex_norm_theta_b_bound_chapter_3}}{\leq}  \frac{4\cdot 2^{p} \gamma \sigma^{p} T \sqrt{V}}{\lambda^{{p}-1}} \notag \\
        & \overset{\eqref{ch3:eq:prox_clipped_SGD_clipping_level_cvx_chapter_3}}{=} \frac{96^{p}}{12} \cdot \frac{\sigma^{p} T \sqrt{V}^{2-{p}}\gamma^{p} A^{{p}-1}}{n^{{p}-1}} \overset{\eqref{ch3:eq:prox_clipped_SGD_step_size_cvx_chapter_3}}{\leq} \frac{V}{6}. \label{ch3:eq:prox_clipped_SGD_convex_sum_2_upper_bound_chapter_3}
\end{align}

\textbf{Upper bound for $\circledThree$.} To estimate this sum, we will use Bernstein's inequality. The summands have conditional expectations equal to zero:
    \begin{equation}
        \Exp_{\xi^l_i}\left[\frac{4\gamma^2}{n^2}\left(\left\|\omega_{i,l}^u\right\|^2 - \Exp_{\xi^l_i}\left[\left\|\omega_{i,l}^u\right\|^2\right]\right)\right] = 0. \notag
    \end{equation}
    Moreover, for all $l = 0,\ldots, T-1$ random vectors $\{\omega_{i,l}^u\}_{i=1}^n$ are independent. Thus, sequence $\left\{\frac{4\gamma^2}{n^2}\left(\|\omega_{i,l}^u\|^2 - \Exp_{\xi_{i}^l}\left[\|\omega_{i,l}^u\|^2\right]\right)\right\}_{l,i = 0,1}^{T-1,n}$ is a martingale difference sequence. Next, the summands are bounded:
    \begin{align}
        \left|\frac{4\gamma^2}{n^2}\left(\left\|\omega_{i,l}^u\right\|^2 - \Exp_{\xi^l_i}\left[\left\|\omega_{i,l}^u\right\|^2\right]\right)\right| & \leq \frac{4\gamma^2}{n^2}\left( \|\omega_{i,l}^u\|^2 +   \Exp_{\xi^l_i}\left[\left\|\omega_{i,l}^u\right\|^2\right]\right)\notag\\
        & \overset{\eqref{ch3:eq:prox_clipped_SGD_convex_norm_theta_u_bound_chapter_3}}{\leq} \frac{32\gamma^2\lambda^2}{n^2}\overset{\eqref{ch3:eq:prox_clipped_SGD_clipping_level_cvx_chapter_3}}{=} \frac{V}{72A^2} \notag\\
        & \le \frac{V}{12A}\eqdef c. \label{ch3:eq:prox_clipped_SGD_convex_technical_10_chapter_3}
    \end{align}
    Finally, conditional variances  $$\widetilde\sigma_{i,t}^2 \eqdef \Exp_{\xi^t_i}\left[\frac{16\gamma^4}{n^4}\left(\left\|\omega_{i,t}^u\right\|^2 - \Exp_{\xi^t_i}\left[\left\|\omega_{i,t}^u\right\|^2\right]\right)^2\right]$$
    of the summands are bounded:
    \begin{align}
        \widetilde\sigma_{i,t}^2 & \overset{\eqref{ch3:eq:prox_clipped_SGD_convex_technical_10_chapter_3}}{\leq} \frac{V}{12A} \Exp_{\xi^t_i}\left[\frac{4\gamma^2}{n^2}\left|\left\|\omega_{i,t}^u\right\|^2 - \Exp_{\xi^t_i}\left[\left\|\omega_{i,t}^u\right\|^2\right]\right|\right] \notag\\
        & \leq \frac{2\gamma^2 V}{3n^2A} \Exp_{\xi^t_i}\left[\|\omega_{i,t}^u\|^2\right]. \label{ch3:eq:prox_clipped_SGD_convex_technical_11_chapter_3}
    \end{align}
    Applying Bernstein's inequality (Lemma~\ref{ch3:lem:Bernstein_ineq_chapter_3}) with $X_{i,l} = (4\gamma^2/n^2)(\|\omega_{i,l}^u\|^2 - \Exp_{\xi_{i}^l}[\|\omega_{i,l}^u\|^2])$, constant $c$ defined in \eqref{ch3:eq:prox_clipped_SGD_convex_technical_10_chapter_3}, $b = \frac{V}{12}$, $H = \frac{V^2}{864A}$, we get
    \begin{align*}
\MoveEqLeft[1] \Prob\left\{|\circledThree| > \frac{V}{12}\quad \text{and}\quad \sum\limits_{l=0}^{T-1}\sum\limits_{i=1}^n \widetilde\sigma_{i,l}^2 \leq \frac{V^2}{864A}\right\} \\
& \leq 2\exp\left(- \frac{b^2}{2H + 2cb/3}\right) = \frac{\beta}{24n(K+1)}.
    \end{align*}
    The above is equivalent to 
    \begin{equation}
        \begin{aligned}
        \Prob\left\{ E_{\circledThree} \right\}
        &\geq 1 - \frac{\beta}{24n(K+1)},\quad \text{for}\quad\\
        E_{\circledThree}
        &= \left\{\text{either } \sum\limits_{l=0}^{T-1}\sum\limits_{i=1}^n \widetilde\sigma_{i,l}^2 > \frac{V^2}{864A} \text{ or } |\circledThree| \leq \frac{V}{12}\right\}.
        \end{aligned}
        \label{ch3:eq:prox_clipped_SGD_convex_sum_3_upper_bound_chapter_3}
    \end{equation}
    Moreover, $E_{T-1}$ implies
    \begin{align}
        \sum\limits_{l=0}^{T-1}\sum\limits_{i=1}^n \widetilde\sigma_{i,l}^2 & \overset{\eqref{ch3:eq:prox_clipped_SGD_convex_technical_11_chapter_3}}{\leq} \frac{2\gamma^2V}{3n^2A} \sum\limits_{l=0}^{T-1}\sum\limits_{i=1}^n  \Exp_{\xi^l_i}\left[\|\omega_{i,l}^u\|^2\right] \notag \\
        & \overset{\eqref{ch3:eq:prox_clipped_SGD_convex_second_moment_theta_u_bound_chapter_3}}{\leq}  
        \frac{12\gamma^2V\lambda^{2-{p}}\sigma^{{p}}T}{nA} \notag\\
        & \overset{\eqref{ch3:eq:prox_clipped_SGD_clipping_level_cvx_chapter_3}}{=}
        \frac{ 48^{p}\cdot 12}{48^2}\cdot\frac{\sigma^{p} T \sqrt{V}^{4-{p}}\gamma^{{p}}}{n^{{p}-1}A^{3-{p}}} \notag \\
        & \overset{\eqref{ch3:eq:prox_clipped_SGD_step_size_cvx_chapter_3}}{\leq} \frac{V^2}{864A}. \label{ch3:eq:prox_clipped_SGD_convex_sum_3_variance_bound_chapter_3}
    \end{align}

\textbf{Upper bound for $\circledFour$.} Probability event $E_{T-1}$ implies
\begin{align}
\circledFour & = \frac{4\gamma^2}{n^2}\sum\limits_{l=0}^{T-1}\sum\limits_{i=1}^n\Exp_{\xi^l_i}\left[\left\|\omega_{i,l}^u\right\|^2\right] \notag \\
& \overset{\eqref{ch3:eq:prox_clipped_SGD_convex_second_moment_theta_u_bound_chapter_3}}{\leq} \frac{72\gamma^2\lambda^{2-{p}}\sigma^{{p}}T}{n} \notag \\
& \overset{\eqref{ch3:eq:prox_clipped_SGD_clipping_level_cvx_chapter_3}}{=}
        \frac{48^{{p}}\gamma^{{p}}\sigma^{{p}}T\sqrt{V}^{2-{p}}}{32n^{{p}-1}A^{2-{p}}} \overset{\eqref{ch3:eq:prox_clipped_SGD_step_size_cvx_chapter_3}}{\leq} \frac{V}{12}. \label{ch3:eq:prox_clipped_SGD_convex_sum_4_upper_bound_chapter_3}
    \end{align}

\textbf{Upper bound for $\circledFive$.} Probability event $E_{T-1}$ implies
\begin{align}
\circledFive & = 4\gamma^2\sum\limits_{l=0}^{T-1}\left\|\omega_{l}^b\right\|^2 \overset{\eqref{ch3:eq:prox_clipped_SGD_convex_norm_theta_b_bound_chapter_3}}{\leq} \frac{4^{({p}+1)} \sigma^{2{p}}T\gamma^2}{\lambda^{2({p}-1)}} \notag \\
& \overset{\eqref{ch3:eq:prox_clipped_SGD_clipping_level_cvx_chapter_3}}{=}   \frac{9216^{p}}{576}\cdot \frac{\sigma^{2{p}}T\gamma^{2{p}}A^{2({p}-1)}}{n^{2({p}-1)}V^{{p}-1}} \notag \\
& \overset{\eqref{ch3:eq:prox_clipped_SGD_step_size_cvx_chapter_3}}{\leq} \frac{V}{6}. \label{ch3:eq:prox_clipped_SGD_convex_sum_5_upper_bound_chapter_3}
    \end{align}

\paragraph{Upper bounds for $\circledSix$.} This sum requires a more refined analysis. We introduce new vectors:
    \begin{equation}
        \delta_{j}^l = \begin{cases}
            \frac{\gamma}{n}\sum\limits_{i=1}^{j-1} \omega_{i,l}^u,& \text{if } \left\|\frac{\gamma}{n}\sum\limits_{i=1}^{j-1} \omega_{i,l}^u\right\| \leq \frac{\sqrt{V}}{2},\\
            0,& \text{otherwise,}
        \end{cases}\label{ch3:eq:SGD_zeta_delta_def_chapter_3}
    \end{equation}
    for all $j \in [n]$ and $l = 0,\ldots, T-1$. Then, by definition
    \begin{equation}
        \|\delta_j^l\| \leq \frac{\sqrt{V}}{2} \label{ch3:eq:SGD_zeta_delta_bound_chapter_3}
    \end{equation}
    and
    \begin{align}
        \circledSix & = \underbrace{\frac{8\gamma}{n} \sum\limits_{l=0}^{T-1} \sum\limits_{j=2}^n \left\langle \delta_j^l, \omega_{j,l}^u \right\rangle}_{\circledSix'} + \frac{8\gamma}{n} \sum\limits_{l=0}^{T-1} \sum\limits_{j=2}^n \left\langle \frac{\gamma}{n}\sum\limits_{i=1}^{j-1} \omega_{i,l}^u - \delta_j^l, \omega_{j,l}^u \right\rangle. \label{ch3:eq:SGD_extra_sums_seven_distributed_chapter_3}
    \end{align}
	    We also note here that $E_{T-1}$ implies
	    \begin{multline}
	        \frac{8\gamma}{n} \sum\limits_{l=0}^{T-1} \sum\limits_{j=2}^n
	        \left\langle \frac{\gamma}{n}\sum\limits_{i=1}^{j-1} \omega_{i,l}^u
	        - \delta_j^l, \omega_{j,l}^u \right\rangle\\
	        = \frac{8\gamma}{n}  \sum\limits_{j=2}^n
	        \left\langle \frac{\gamma}{n}\sum\limits_{i=1}^{j-1} \omega_{i,T-1}^u
	        - \delta_j^{T-1}, \omega_{j,T-1}^u \right\rangle.
	        \label{ch3:eq:SGD_extra_sums_seven_distributed_1_chapter_3}
	    \end{multline}
    
\paragraph{Upper bound for $\circledSix'$.} To estimate this sum, we will use Bernstein's inequality. The summands have conditional expectations equal to zero:
    \begin{equation*}
    \Exp_{\xi_{j}^l}\left[\frac{8\gamma}{n}\left\langle \delta_j^l, \omega_{j,l}^u \right\rangle\right] = \frac{8\gamma}{n}\left\langle \delta_j^l, \Exp_{\xi_{j}^l}[\omega_{j,l}^u] \right\rangle = 0.
    \end{equation*}
	    Moreover, for all $l = 0,\ldots, T-1$ random vectors $\{\omega_{i,l}^u\}_{i=1}^n$ are independent. Thus, sequence $\left\{\frac{8\gamma}{n}\left\langle \delta_j^l, \omega_{j,l}^u \right\rangle\right\}_{l,j = 0,2}^{T-1,n}$ is a martingale difference sequence. Next, the summands are bounded:
	    \begin{align}
	        \left|\frac{8\gamma}{n}
	        \left\langle \delta_j^l, \omega_{j,l}^u \right\rangle \right|
	        &\leq \frac{8\gamma}{n} \left\|\delta_j^l\right\|
	        \cdot \|\omega_{j,l}^u\|\notag\\
	        &\overset{\eqref{ch3:eq:SGD_zeta_delta_bound_chapter_3},
	        \eqref{ch3:eq:prox_clipped_SGD_convex_norm_theta_u_bound_chapter_3}}{\leq}
	        \frac{8\gamma}{n}\cdot \frac{\sqrt{V}}{2}\cdot 2\lambda\notag\\
	        &= \frac{V}{6A} \eqdef c.
	        \label{ch3:eq:gap_thm_SGD_technical_6_5_1234_chapter_3}
	    \end{align}
    Finally, conditional variances $(\sigma_{j,l}')^2 \eqdef \Exp_{\xi_{j}^l}\left[\frac{64\gamma^2}{n^2} \langle \delta_j^l, \omega_{j,l}^u \rangle^2\right]$ of the summands are bounded:
    \begin{align}
(\sigma_{j,l}')^2 & \leq \Exp_{\xi_{j}^l}\left[\frac{64\gamma^2 }{n^2}\|\delta_j^l\|^2\cdot \|\omega_{j,l}^u\|^2\right] \notag \\
& \overset{\eqref{ch3:eq:SGD_zeta_delta_bound_chapter_3}}{\leq} \frac{16\gamma^2 V}{n^2} \Exp_{\xi_{j}^l}\left[\|\omega_{j,l}^u\|^2\right]. \label{ch3:eq:gap_thm_SGD_technical_7_1234_chapter_3}
    \end{align}
    Applying Bernstein's inequality (Lemma~\ref{ch3:lem:Bernstein_ineq_chapter_3}) with $X_{j,l} = \frac{8\gamma}{n}\left\langle \delta_j^l, \omega_{j,l}^u \right\rangle$, constant $c$ defined in \eqref{ch3:eq:gap_thm_SGD_technical_6_5_1234_chapter_3}, $b = \frac{V}{6}$, $H = \frac{V^2}{216A}$, we get
    \begin{align*}
\MoveEqLeft[1] \Prob\left\{|\circledSix'| > \frac{V}{6} \text{ and } \sum\limits_{l=0}^{T-1}\sum\limits_{j=2}^n(\sigma_{j,l}')^2 \leq \frac{V^2}{216A}\right\} \\
& \leq 2\exp\left(- \frac{b^2}{2H + 2cb/3}\right) = \frac{\beta}{24n(K+1)}.
    \end{align*}
    The above is equivalent to 
    \begin{equation}
        \begin{aligned}
        \Prob\left\{ E_{\circledSix'} \right\}
        &\geq 1 - \frac{\beta}{24n(K+1)},\quad \text{for}\quad\\
        E_{\circledSix'}
        &= \left\{\text{either } \sum\limits_{l=0}^{T-1}\sum\limits_{j=2}^n(\sigma_{j,l}')^2 > \frac{V^2}{216A} \text{ or } |\circledSix'| \leq \frac{V}{6}\right\}.
        \end{aligned}
        \label{ch3:eq:bound_7_gap_SGD_chapter_3}
    \end{equation}
    Moreover, $E_{T-1}$ implies
    \begin{align}
        \sum\limits_{l=0}^{T-1}\sum\limits_{j=2}^n(\sigma_{j,l}')^2 & \overset{\eqref{ch3:eq:gap_thm_SGD_technical_7_1234_chapter_3}}{\leq} \frac{16\gamma^2 V}{n^2} \sum\limits_{l=0}^{T-1} \sum\limits_{j=2}^n \Exp_{\xi_{j}^l}\left[\|\omega_{j,l}^u\|^2\right] \notag \\
        & \overset{\eqref{ch3:eq:prox_clipped_SGD_convex_second_moment_theta_u_bound_chapter_3}}{\leq} \frac{288(K+1)\gamma^2 V \lambda^{2-{p}} \sigma^{p}}{n} \notag \\
        & \overset{\eqref{ch3:eq:prox_clipped_SGD_clipping_level_cvx_chapter_3}}{\leq} \frac{288 (K+1)\gamma^{{p}}\sigma^{p} V^{2-{p}/2}}{48^{2-{p}} n^{{p}-1} A^{2-{p}}} \notag \\
        & \overset{\eqref{ch3:eq:prox_clipped_SGD_step_size_cvx_chapter_3}}{\leq} \frac{V^2}{216A}. \label{ch3:eq:bound_7_variances_gap_SGD_chapter_3}
    \end{align}
    where in the second step we also use $T \leq K+1$.

That is, we derive the upper bounds for  $\circledOne, \circledTwo, \circledThree, \circledFour, \circledFive, \circledSix$. More precisely, $E_{T-1}$ implies
\begin{gather*}
        A_T \overset{\eqref{ch3:eq:prox_clipped_SGD_convex_technical_7_chapter_3}}{\leq} V + \circledOne + \circledTwo + \circledThree + \circledFour + \circledFive + \circledSix,\\
        \circledSix \overset{\eqref{ch3:eq:SGD_extra_sums_seven_distributed_chapter_3}}{=} \circledSix' + \frac{8\gamma}{n} \sum\limits_{j=2}^n \left\langle \frac{\gamma}{n}\sum\limits_{i=1}^{j-1} \omega_{i,T-1}^u - \delta_j^{T-1}, \omega_{j,T-1}^u \right\rangle,\\
        \circledTwo \overset{\eqref{ch3:eq:prox_clipped_SGD_convex_sum_2_upper_bound_chapter_3}}{\leq} \frac{V}{6},\quad \circledFour \overset{\eqref{ch3:eq:prox_clipped_SGD_convex_sum_4_upper_bound_chapter_3}}{\leq} \frac{V}{12}, \quad \circledFive \overset{\eqref{ch3:eq:prox_clipped_SGD_convex_sum_5_upper_bound_chapter_3}}{\leq} \frac{V}{6},\\
        \sum\limits_{l=0}^{T-1}\sum\limits_{i=1}^n \sigma_{i,l}^2 \overset{\eqref{ch3:eq:prox_clipped_SGD_convex_sum_1_variance_bound_chapter_3}}{\leq} \frac{V^2}{216A},\quad \sum\limits_{l=0}^{T-1}\sum\limits_{i=1}^n \widetilde\sigma_{i,l}^2 \overset{\eqref{ch3:eq:prox_clipped_SGD_convex_sum_3_variance_bound_chapter_3}}{\leq} \frac{V^2}{864A},\\
        \sum\limits_{l=0}^{T-1}\sum\limits_{j=2}^n(\sigma_{j,l}')^2 \overset{\eqref{ch3:eq:bound_7_variances_gap_SGD_chapter_3}}{\leq} \frac{V^2}{216A}.
\end{gather*}
    In addition, we also establish (see \eqref{ch3:eq:prox_clipped_SGD_convex_sum_1_upper_bound_chapter_3}, \eqref{ch3:eq:prox_clipped_SGD_convex_sum_3_upper_bound_chapter_3}, \eqref{ch3:eq:bound_7_gap_SGD_chapter_3} and our induction assumption)
\begin{gather*}
        \Prob\{E_{T-1}\} \geq 1 - \frac{(T-1)\beta}{K+1},\\
        \Prob\{E_{\circledOne}\} \geq 1 - \frac{\beta}{24n(K+1)}, \quad \Prob\{E_{\circledThree}\} \geq 1 - \frac{\beta}{24n(K+1)},\\
        \Prob\{E_{\circledSix'}\} \geq 1 - \frac{\beta}{24n(K+1)},
    \end{gather*}
    where
\begin{align*}
        E_{\circledOne} & = \left\{ \text{either } \sum\limits_{l=0}^{T-1}\sum\limits_{i=1}^n \sigma_{i,l}^2 > \frac{V^2}{216A} \text{ or } |\circledOne| \leq \frac{V}{6}\right\},\\
        E_{\circledThree} & = \left\{ \text{either } \sum\limits_{l=0}^{T-1}\sum\limits_{i=1}^n \widetilde\sigma_{i,l}^2 > \frac{V^2}{864A} \text{ or } |\circledThree| \leq \frac{V}{12}\right\}, \\
        E_{\circledSix'} & = \left\{ \text{either } \sum\limits_{l=0}^{T-1} \sum\limits^n_{j=2}(\sigma_{j,l}')^2  > \frac{V^2}{216A} \text{ or } |\circledSix'| \leq \frac{V}{6}\right\}.
\end{align*}
     Therefore, probability event $E_{T-1} \cap E_{\circledOne} \cap E_{\circledThree} \cap E_{\circledSix'}$ implies

     \begin{align}
        A_T & \leq V + \frac{V}{6} + \frac{V}{6} + \frac{V}{12} + \frac{V}{12} + \frac{V}{6} + \frac{V}{6}\notag\\
        &\quad\quad + \frac{8\gamma}{n} \sum\limits_{j=2}^n \left\langle \frac{\gamma}{n}\sum\limits_{i=1}^{j-1} \omega_{i,T-1}^u - \delta_j^{T-1}, \omega_{j,T-1}^u \right\rangle \notag \\
        & \leq 2V + \frac{8\gamma}{n} \sum\limits_{j=2}^n \left\langle \frac{\gamma}{n}\sum\limits_{i=1}^{j-1} \omega_{i,T-1}^u - \delta_j^{T-1}, \omega_{j,T-1}^u \right\rangle \label{ch3:eq:gap_thm_SGD_technical_11_chapter_3}
    \end{align}
    for $t = T$. 

    In the final part of the proof, we will show that $\frac{\gamma}{n}\sum\limits_{i=1}^{j-1} \omega_{i,T-1}^u = \delta_j^{T-1}$ with high probability. In particular, we consider probability event $\widetilde{E}_{T-1,j}$ defined as follows: inequalities
    \begin{equation}
         \left\|\frac{\gamma}{n}\sum\limits_{i=1}^{r-1} \omega_{i,T-1}^u\right\| \leq \frac{\sqrt{V}}{2} \notag
    \end{equation}
    hold for $r = 2, \ldots, j$ simultaneously. We want to show that $\Prob\{E_{T-1} \cap \widetilde{E}_{T-1,j}\} \geq 1 - \frac{(T-1)\beta}{K+1} - \frac{j\beta}{8n(K+1)}$ for all $j = 2, \ldots, n$. For $j = 2$ the statement is trivial since
    \begin{equation*}
         \left\|\frac{\gamma}{n} \omega_{1,T-1}^u\right\| \overset{\eqref{ch3:eq:prox_clipped_SGD_convex_norm_theta_u_bound_chapter_3}}{\leq} \frac{2\gamma\lambda}{n} \leq \frac{\sqrt{V}}{2}.
    \end{equation*}
    Next, we assume that the statement holds for some $j = m-1 < n$, i.e.,  $\Prob\{E_{T-1}\cap\widetilde{E}_{T-1,m-1}\} \geq 1 - \frac{(T-1)\beta}{K+1} - \frac{(m-1)\beta}{8n(K+1)}$. Our goal is to prove that $\Prob\{E_{T-1}\cap\widetilde{E}_{T-1,m}\} \geq 1 - \frac{(T-1)\beta}{K+1} - \frac{m\beta}{8n(K+1)}$. We have
    \begin{align}
        \left\|\frac{\gamma}{n}\sum\limits_{i=1}^{m-1} \omega_{i,T-1}^u\right\| & = \sqrt{\frac{\gamma^2}{n^2}\left\|\sum\limits_{i=1}^{m-1} \omega_{i,T-1}^u\right\|^2} \notag\\
        & = \sqrt{\frac{\gamma^2}{n^2}\sum\limits_{i=1}^{m-1}\|\omega_{i,T-1}^u\|^2 + \frac{2\gamma}{n}\sum\limits_{i=1}^{m-1}\left\langle \frac{\gamma}{n}\sum\limits_{r=1}^{i-1}\omega_{r,T-1}^u , \omega_{i,T-1}^u \right\rangle} \notag\\
        & \leq \sqrt{\frac{\gamma^2}{n^2}\sum\limits_{l=0}^{T-1}\sum\limits_{i=1}^{m-1}\|\omega_{i,l}^u\|^2 + \frac{2\gamma}{n}\sum\limits_{i=1}^{m-1}\left\langle \frac{\gamma}{n}\sum\limits_{r=1}^{i-1}\omega_{r,T-1}^u , \omega_{i,T-1}^u \right\rangle}. \notag
    \end{align}

    Next, we introduce a new notation:
    \begin{gather*}
        \rho_{i,T-1}' = \begin{cases}
            \frac{\gamma}{n}\sum\limits_{r=1}^{i-1}\omega_{r,T-1}^u,& \text{if } \left\|\frac{\gamma}{n}\sum\limits_{r=1}^{i-1}\omega_{r,T-1}^u\right\| \leq \frac{\sqrt{V}}{2},\\
            0,& \text{otherwise}
        \end{cases}
    \end{gather*}
    for $i = 1,\ldots,m-1$. By definition, we have
    \begin{equation}
        \|\rho_{i,T-1}'\| \leq \frac{\sqrt{V}}{2} \label{ch3:eq:SGD_bound_rho_chapter_3}
    \end{equation}
    for $i = 1,\ldots,m-1$. Moreover, $\widetilde{E}_{T-1,m-1}$ implies $\rho_{i,T-1}' = \frac{\gamma}{n}\sum\limits_{r=1}^{i-1}\omega_{r,T-1}^u$ for $i = 1,\ldots,m-1$ and
    \begin{align}
        \left\|\frac{\gamma}{n}\sum\limits_{i=1}^{m-1} \omega_{i,l}^u\right\| & \leq \sqrt{\circledThree + \circledFour + \circledSeven}, \notag
    \end{align}
    where
    \begin{gather*}
         \circledSeven =  \frac{2\gamma}{n}\sum\limits_{i=1}^{m-1}\left\langle \rho_{i,T-1}' , \omega_{i,T-1}^u \right\rangle.
    \end{gather*}
    It remains to estimate $\circledSeven$.

     \textbf{Upper bound for $\circledSeven$.} To estimate this sum, we will use Bernstein's inequality. The summands have conditional expectations equal to zero:
    \begin{align*}
\MoveEqLeft[1] \Exp_{\xi_{i}^{T-1}}\left[\frac{2\gamma}{n} \langle \rho_{i,T-1}', \omega_{i,T-1}^u \rangle\right] \\
& = \frac{2\gamma}{n} \left\langle \rho_{i,T-1}', \Exp_{\xi_{i}^{T-1}}[\omega_{i,T-1}^u] \right\rangle = 0,
    \end{align*}
    since random vectors $\{\omega_{i,T-1}^u\}_{i=1}^n$ are independent. Thus, for $i = 1, \ldots, m-1$ the summands $(2\gamma/n) \langle \rho_{i,T-1}', \omega_{i,T-1}^u \rangle$ form a martingale difference sequence. Next, the summands are bounded:
    \begin{align}
\left|\frac{2\gamma}{n}\langle \rho_{i,T-1}', \omega_{i,T-1}^u \rangle \right| & \leq \frac{2\gamma}{n} \|\rho_{i,T-1}'\|\cdot \|\omega_{i,T-1}^u\| \notag \\
& \overset{\eqref{ch3:eq:SGD_bound_rho_chapter_3},\eqref{ch3:eq:prox_clipped_SGD_convex_norm_theta_u_bound_chapter_3}}{\leq} \frac{2\gamma}{n} \sqrt{V}\lambda \overset{\eqref{ch3:eq:prox_clipped_SGD_clipping_level_cvx_chapter_3}}{=} \frac{V}{24A} \notag \\
& \eqdef c. \label{ch3:eq:gap_thm_SGD_technical_6_5_10'_chapter_3}
    \end{align}
    Finally, conditional variances $(\widetilde\sigma_{i,T-1}')^2 \eqdef \Exp_{\xi_{i}^{T-1}}\left[\frac{4\gamma^2}{n^2} \langle \rho_{i,T-1}', \omega_{i,T-1}^u \rangle^2\right]$ of the summands are bounded:
    \begin{align}
(\widetilde\sigma_{i,T-1}')^2 & \leq \Exp_{\xi_{i}^{T-1}}\left[\frac{4\gamma^2 }{n^2}\|\rho_{i,T-1}'\|^2\cdot \|\omega_{i,T-1}^u\|^2\right] \notag \\
& \overset{\eqref{ch3:eq:SGD_bound_rho_chapter_3}}{\leq} \frac{\gamma^2 V}{n^2} \Exp_{\xi_{i}^{T-1}}\left[\|\omega_{i,T-1}^u\|^2\right]. \label{ch3:eq:gap_thm_SGD_technical_7_10'_chapter_3}
    \end{align}
    Applying Bernstein's inequality (Lemma~\ref{ch3:lem:Bernstein_ineq_chapter_3}) with $X_{i} = \frac{2\gamma}{n} \langle \rho_{i,T-1}', \omega_{i,T-1}^u \rangle$, constant $c$ defined in \eqref{ch3:eq:gap_thm_SGD_technical_6_5_10'_chapter_3}, $b = \frac{V}{24}$, $H = \frac{V^2}{3456A}$, we get
    \begin{align*}
\MoveEqLeft[1] \Prob\left\{|\circledSeven| > \frac{V}{24} \text{ and } \sum\limits_{i=1}^{m-1}(\widetilde\sigma_{i,T-1}')^2 \leq \frac{V^2}{3456A}\right\} \\
& \leq 2\exp\left(- \frac{b^2}{2H + 2cb/3}\right) = \frac{\beta}{24n(K+1)}.
    \end{align*}
    The above is equivalent to 
    \begin{equation}
        \begin{aligned}
        \Prob\left\{ E_{\circledSeven} \right\}
        &\geq 1 - \frac{\beta}{24n(K+1)},\quad \text{for}\quad\\
        E_{\circledSeven}
        &= \left\{\text{either } \sum\limits_{i=1}^{m-1}(\widetilde\sigma_{i,T-1}')^2 > \frac{V^2}{3456A} \text{ or } |\circledSeven| \leq \frac{V}{24}\right\}.
        \end{aligned}
        \label{ch3:eq:bound_10_gap_SGD_chapter_3}
    \end{equation}
    Moreover, $E_{T-1}$ implies
    \begin{align}
        \sum\limits_{i=1}^{m-1}(\widetilde\sigma_{i,T-1}')^2 & \overset{\eqref{ch3:eq:gap_thm_SGD_technical_7_10'_chapter_3}}{\leq} \frac{\gamma^2 V}{n^2} \sum\limits_{i=1}^n \Exp_{\xi_{i}^{T-1}}\left[\|\omega_{i,T-1}^u\|^2\right] \notag \\
        & \overset{\eqref{ch3:eq:prox_clipped_SGD_convex_second_moment_theta_u_bound_chapter_3}}{\leq} \frac{18\gamma^2 V \lambda^{2-{p}} \sigma^{p}}{n} \notag \\
        & \overset{\eqref{ch3:eq:prox_clipped_SGD_clipping_level_cvx_chapter_3}}{\leq} \frac{18 \gamma^{{p}}\sigma^{p} V^{2-{p}/2}}{48^{2-{p}} n^{{p}-1} A^{2-{p}}} \notag \\
        & \overset{\eqref{ch3:eq:prox_clipped_SGD_step_size_cvx_chapter_3}}{\leq} \frac{V^2}{3456A}. \label{ch3:eq:bound_10_variances_gap_SGD_chapter_3}
    \end{align}

    Putting all together we get that $E_{T-1}\cap \widetilde{E}_{T-1,m-1}$ implies
    \begin{gather*}
         \left\|\frac{\gamma}{n}\sum\limits_{i=1}^{m-1} \omega_{i,T-1}^u\right\| \leq \sqrt{\circledThree + \circledFour + \circledSeven},\quad \circledFour \overset{\eqref{ch3:eq:prox_clipped_SGD_convex_sum_4_upper_bound_chapter_3}}{\leq} \frac{V}{12},\\
        \sum\limits_{l=0}^{T-1}\sum\limits_{i=1}^n\widetilde\sigma_{i,l}^2 \overset{\eqref{ch3:eq:prox_clipped_SGD_convex_sum_3_variance_bound_chapter_3}}{\leq} \frac{V^2}{864A},\quad 
        \sum\limits_{i=1}^{m-1}(\widetilde\sigma_{i,T-1}')^2 \leq \frac{V^2}{3456A}.
    \end{gather*}
    In addition, we also establish (see \eqref{ch3:eq:prox_clipped_SGD_convex_sum_3_upper_bound_chapter_3}, \eqref{ch3:eq:bound_10_gap_SGD_chapter_3} and our induction assumption)
    \begin{gather*}
        \Prob\{E_{T-1}\cap \widetilde{E}_{T-1,m-1}\} \geq 1 - \frac{(T-1)\beta}{K+1} - \frac{(m-1)\beta}{8n(K+1)},\\
        \Prob\{E_{\circledThree}\} \geq 1 - \frac{\beta}{24n(K+1)}, \quad \Prob\{E_{\circledSeven}\} \geq 1 - \frac{\beta}{24n(K+1)}
    \end{gather*}
    where
    \begin{align}
        E_{\circledThree} & = \left\{\text{either } \sum\limits_{l=0}^{T-1}\sum\limits_{i=1}^n\widetilde\sigma_{i,l}^2 > \frac{V^2}{864A} \text{ or } |\circledThree| \leq \frac{V}{12}\right\}, \notag\\
        E_{\circledSeven} & = \left\{\text{either } \sum\limits_{i=1}^{m-1}(\widetilde\sigma_{i,T-1}')^2 > \frac{V^2}{3456A} \text{ or } |\circledSeven| \leq \frac{V}{24}\right\} \notag
    \end{align}
    Therefore, probability event $E_{T-1} \cap \widetilde{E}_{T-1,m-1} \cap E_{\circledThree} \cap E_{\circledSeven}$ implies
    \begin{gather*}
        \left\|\frac{\gamma}{n}\sum\limits_{i=1}^{m-1} \omega_{i,T-1}^u\right\| \leq \sqrt{\frac{V}{12} + \frac{V}{12} + \frac{V}{24}} \leq \frac{\sqrt{V}}{2}.
    \end{gather*}
     This implies $\widetilde E_{T-1,m}$ and
    \begin{align*}
        \Prob\{E_{T-1} \cap \widetilde{E}_{T-1,m}\} & \geq \Prob\{E_{T-1} \cap \widetilde{E}_{T-1,m-1} \cap E_{\circledThree} \cap E_{\circledSeven} \} \\
        & = 1 - \Prob\left\{\overline{E_{T-1} \cap \widetilde{E}_{T-1,m-1}} \cup \overline{E}_{\circledThree} \cup \overline{E}_{\circledSeven} \right\}\\
        & \geq 1 - \frac{(T-1)\beta}{K+1} - \frac{m\beta}{8n(K+1)}.
    \end{align*}
    Therefore, for all $m = 2,\ldots,n$ the statement holds and, in particular, $\Prob\{E_{T-1} \cap \widetilde{E}_{T-1, n}\} \geq 1 - \frac{(T-1)\beta}{K+1} - \frac{\beta}{8(K+1)}$. Taking into account \eqref{ch3:eq:gap_thm_SGD_technical_11_chapter_3}, we conclude that $E_{T-1} \cap \widetilde{E}_{T-1, n} \cap E_{\circledOne} \cap E_{\circledThree} \cap E_{\circledSix'}$ implies
    \begin{equation*}
        A_T \leq 2V
    \end{equation*}
    that is equivalent to \eqref{ch3:eq:prox_clipped_SGD_convex_induction_inequality_1_chapter_3} for $t = T$. Moreover,
    \begin{align*}
        \Prob\left\{E_T\right\} & \geq \Prob\left\{E_{T-1} \cap \widetilde{E}_{T-1,n} \cap E_{\circledOne} \cap E_{\circledThree} \cap E_{\circledSix'}\right\} \\
        & = 1 - \Prob\left\{\overline{E_{T-1} \cap \widetilde{E}_{n}} \cup \overline{E}_{\circledOne} \cup \overline{E}_{\circledThree} \cup \overline{E}_{\circledSix'}\right\}\\
        & \geq 1 - \frac{(T-1)\beta}{K+1} - \frac{\beta}{8(K+1)} - 3\cdot \frac{\beta}{24n(K+1)} \geq 1 - \frac{T\beta}{K+1}.
    \end{align*}
     In other words, we showed that $\Prob\{E_k\} \geq 1 - k\beta/(K+1)$ for all $k = 0,1,\ldots,K+1$. For $k = K+1$ we have that with probability at least $1 - \beta$
\begin{equation*}
        \Phi(\overline{x}^{K+1}) - \Phi(x^*) \overset{\eqref{ch3:eq:prox_clipped_SGD_convex_technical_1_chapter_3},\eqref{ch3:eq:prox_clipped_SGD_convex_induction_inequality_1_chapter_3}}{\leq}\frac{V}{\gamma(K+1)}.
\end{equation*}
Finally, if $V = 2R^2$ and $\gamma$ equals the minimum in \eqref{ch3:eq:prox_clipped_SGD_step_size_cvx_chapter_3}, then with probability at least $1-\beta$
\begin{align*}
        \Phi(\overline{x}^{K+1}) - \Phi(x^*) & \leq \frac{V}{\gamma(K+1)}\\
        & = \max\biggl\{\frac{720 LR^2A}{K+1}, \frac{384A R \zeta_*}{\sqrt{n}(K+1)},\\
        &\qquad\quad \frac{48\sqrt{2}\cdot 27^{1/{p}}\sigma R A^{({p}-1)/{p}}}{n^{({p}-1)/{p}}(K+1)^{({p}-1)/{p}}}\biggr\}\\
        & = \cO\biggl(\max\biggl\{\frac{L R^2 \ln(nK/\beta)}{K}, \frac{R\zeta_* \ln(nK/\beta)}{\sqrt{n}K},\\
        &\qquad\quad \frac{ \sigma R\ln^{({p}-1)/{p}}(nK/\beta)}{n^{({p}-1)/{p}}K^{({p}-1)/{p}}}\biggr\}\biggr).
\end{align*}
To get $\Phi(\overline{x}^{K+1}) - \Phi(x^*) \leq \varepsilon$ with probability at least $1-\beta$ it is sufficient to choose $K$ such that all terms in the maximum above are $\cO(\varepsilon)$. This leads to
    \begin{align*}
        K = \cO\biggl(\max\biggl\{ & \frac{LR^2}{\varepsilon}\ln\frac{nLR^2}{\beta\varepsilon}, \frac{R\zeta_*}{\sqrt{n}\varepsilon}\ln\frac{\sqrt{n}R\zeta_*}{\beta\varepsilon},\\
        & \left(\frac{\sigma R}{\varepsilon n^{(p-1)/p}}\right)^{p/(p-1)} \ln \left(\frac{1}{\beta} \left(\frac{\sigma R}{\varepsilon}\right)^{p/(p-1)}\right) \biggr\}\biggr),
    \end{align*}
which concludes the proof.
\end{proof}

In view of Lemma~\ref{ch3:lem:QSC_to_QSM_chapter_3}, the result in the quasi-strongly convex case for \algname{DProx-clipped-SGD-shift} follows from our result for \algname{DProx-clipped-SGDA-shift}.

\clearpage

\section{Additional Proofs for \algname{DProx-clipped-SSTM-shift}}\label{ch3:appendix:dprox_clipped_SSTM_shift_chapter_3}

In this section, we provide the complete formulations of our results for \algname{DProx-clipped-SSTM-shift} and proofs. For the readers' convenience, the method's update rule is repeated below: $x^0 = y^0 = z^0$, $A_0 = a_0 = 0$, $a_{k+1} = (k+2)/(2aL)$, $A_{k+1} = A_k + a_{k+1}$ and
\begin{gather*}
    x^{k+1} = \frac{A_k y^k + {a}_{k+1}z^k}{A_{k+1}},\quad z^{k+1} = \thesischthreeprox_{a_{k+1}\Psi}\left(z^k - a_{k+1}\thesischthreetg(x^{k+1})\right), \\
    \thesischthreetg(x^{k+1}) = \frac{1}{n}\sum\limits_{i=1}^n \thesischthreetg_i(x^{k+1}),\;\; \thesischthreetg_i(x^{k+1}) = h_i^k + \hat\Delta_i^k, \\
    h_i^{k+1} = h_i^k + \nu_k \hat \Delta_i^k,\quad \hat\Delta_i^k = \clip\left(\nabla f_i(x^{k+1};\xi_i^k) - h_i^k, \lambda_k\right), \\
    y^{k+1} = \frac{A_k y^k + {a}_{k+1}z^{k+1}}{A_{k+1}} 
\end{gather*}
where $\xi_1^k,\ldots, \xi_n^k$ are sampled independently from each other and previous steps.

\subsection{Convex case}
The following lemma is the main ``optimization'' part of the analysis of \algname{DProx-clipped-SSTM-shift}.
\begin{lemma}
\label{ch3:lem:SSTM_optimization_lemma_convex_chapter_3}
    Let Assumptions~\ref{ch3:as:bounded_alpha_moment_chapter_3}, \ref{ch3:as:L_smoothness_chapter_3} and \ref{ch3:as:str_cvx_chapter_3}($\mu = 0$) hold on $Q = B_{5n\sqrt{M}}(x^*)$, where
    \begin{equation*}
        M \geq \|x^0-x^*\|^2+C^2{a}^2_{K_0+1}\frac{1}{n}\sum^n_{i=1}\|\nabla f_i(x^*)\|^2
    \end{equation*}
    for some $C > 0$ with $C^2 n \geq 1$. Let $x^{k},~y^k,~z^k$ lie in $B_{5n\sqrt{M}}(x^*)$ for all $k = 0,1,\dots, K$ for some $K\geq 0$. Additionally, let parameters of \algname{DProx-clipped-SSTM-shift} satisfy
    \begin{equation}
        \label{ch3:eq:a_K0_D_prox_clipped_SSTM_shift_lemma_chapter_3}
        a \geq \max\left\{2, 2C^2 + \frac{2}{n}\right\},\quad K_0 = \left\lceil 2C^2n\right\rceil;
    \end{equation}
    \begin{equation}
         \label{ch3:eq:nu_k_D_prox_clipped_SSTM_shift_lemma_chapter_3}
         \nu_k = \begin{cases}
                    \frac{(k+2)^2}{C^2(K_0+2)^2n} ,& \text{if $k < K_0$};\\
                    \frac{2k+5}{(k+3)^2} ,& \text{if $k \geq K_0$};
    \end{cases},     
    \end{equation}
    
    then the iterates produced by \algname{DProx-clipped-SSTM-shift} satisfy 
    \begin{align}
        A_K(\Phi(y^K)-\Phi(x^*)) &\leq \frac{1}{2}M_0-\frac{1}{2}M_K + \sum^{K-1}_{k=0}a_{k+1}\la\omega_{k+1}, x^*-z^k\ra \notag\\
        &\quad + \sum^{K-1}_{k=0}{a}^2_{k+1}\|\omega_{k+1}\|^2 + \sum^{K-1}_{k=0}\sum^{n}_{i=1}\frac{{a}^2_{k+1}}{n^2}\|\omega_{i,k+1}\|^2,
    \end{align}
    where Lyapunov function $M_k$ is defined as follows 
    \begin{equation}
        \label{ch3:eq:Lyapunov_function_D_prox_clipped_SSTM_shift_chapter_3}
        M_k = \|z^k-x^*\|^2 +C^2\widetilde{a}^2_{k+1}\frac{1}{n}\sum^n_{i=1}\|h^k_i-h^*_i\|^2,
    \end{equation}
    where 
    \begin{equation}
        \widetilde{a}_{k+1} = \begin{cases}
                    {a}_{K_0+1}& \text{if $k < K_0$};\\
                    a_{k+1} & \text{if $k \geq K_0$};
    \end{cases},     
    \end{equation}
     and $\omega_{k+1}$ is defined as follows 
    \begin{align}
        \omega_{i,k+1} &\eqdef \thesischthreetg_i(x^{k+1}) - \nabla f_i(x^{k+1}), \notag\\
        \omega_{k+1} &\eqdef \frac{1}{n}\sum^n_{i=1}\omega_{i,k+1}. \label{ch3:eq:omega_k+1_def_dprox_clipped_SSTM_shift_chapter_3}
    \end{align}
\end{lemma}
\begin{proof}
  By optimality condition for the problem \eqref{ch3:eq:D_prox_clipped_SSTM_shift_x_z_chapter_3}, we have for any $z$
\begin{align*}
    \MoveEqLeft[1] a_{k+1}\left\la \thesischthreetg(x^{k+1}), z^k-z\right\ra \\
    & \leq a_{k+1}(\Psi(z)-\Psi(z^{k+1})) +a_{k+1}\left\la \thesischthreetg(x^{k+1}), z^k-z^{k+1}\right\ra \\
    &\quad+ \frac{1}{2}\|z^k-z\|^2 - \frac{1}{2}\|z^{k+1}-z\|^2 - \frac{1}{2}\|z^{k+1}-z^k\|^2  \\
    & \overset{\eqref{ch3:eq:omega_k+1_def_dprox_clipped_SSTM_shift_chapter_3}}{\leq} a_{k+1}(\Psi(z)-\Psi(z^{k+1})) +a_{k+1}\left\la \omega_{k+1}, z^k-z^{k+1}\right\ra \\
        &\quad +a_{k+1}\left\la \nabla f(x^{k+1}), z^k-z^{k+1}\right\ra\\
    &\quad + \frac{1}{2}\|z^k-z\|^2 - \frac{1}{2}\|z^{k+1}-z\|^2 - \frac{1}{2}\|z^{k+1}-z^k\|^2 
\end{align*}
Using $A_{k+1}(y^{k+1}-x^{k+1})= a_{k+1}(z^{k+1}-z^k)$, we get 
\begin{align*}
    \MoveEqLeft[1] a_{k+1}\left\la \thesischthreetg(x^{k+1}), z^k-z\right\ra \\
    & \leq a_{k+1}(\Psi(z)-\Psi(z^{k+1})) +a_{k+1}\left\la \omega_{k+1}, z^k-z^{k+1}\right\ra \\
        &\quad +A_{k+1}\left\la \nabla f(x^{k+1}), x^{k+1}-y^{k+1}\right\ra\\
    &\quad + \frac{1}{2}\|z^k-z\|^2 - \frac{1}{2}\|z^{k+1}-z\|^2 - \frac{1}{2}\|z^{k+1}-z^k\|^2 \\
    & \leq a_{k+1}(\Psi(z)-\Psi(z^{k+1})) +a_{k+1}\left\la \omega_{k+1}, z^k-z^{k+1}\right\ra \\
        &\quad +A_{k+1}\left(f(x^{k+1})-f(y^{k+1}) + \frac{L}{2}\|y^{k+1}-x^{k+1}\|^2\right)\\
    &\quad + \frac{1}{2}\|z^k-z\|^2 - \frac{1}{2}\|z^{k+1}-z\|^2 - \frac{1}{2}\|z^{k+1}-z^k\|^2 \\
    & = a_{k+1}(\Psi(z)-\Psi(z^{k+1})) +a_{k+1}\left\la \omega_{k+1}, z^k-z^{k+1}\right\ra \\
        &\quad +A_{k+1}\left(f(x^{k+1})-f(y^{k+1})\right) + \frac{a_{k+1}^2L}{2A_{k+1}}\|z^{k+1}-z^{k}\|^2\\
    &\quad + \frac{1}{2}\|z^k-z\|^2 - \frac{1}{2}\|z^{k+1}-z\|^2 - \frac{1}{2}\|z^{k+1}-z^k\|^2 \\
    & = a_{k+1}(\Psi(z)-\Psi(z^{k+1})) +a_{k+1}\left\la \omega_{k+1}, z^k-z^{k+1}\right\ra \\
        &\quad +A_{k+1}\left(f(x^{k+1})-f(y^{k+1})\right) + \frac{1}{2}\|z^k-z\|^2 - \frac{1}{2}\|z^{k+1}-z\|^2\\
    &\quad  - \frac{1}{2}\left(1- \frac{a_{k+1}^2L}{A_{k+1}}\right)\|z^{k+1}-z^k\|^2 
\end{align*}
where in the second inequality we use $L$-smoothness of $f$. 
Using Young's inequality with a parameter $D > 0$, we have 
\begin{align*}
    \MoveEqLeft[1] a_{k+1}\left\la \thesischthreetg(x^{k+1}), z^k-z\right\ra \\
    & \leq a_{k+1}(\Psi(z)-\Psi(z^{k+1})) +a_{k+1}\frac{D}{2}\|\omega_{k+1}\|^2+ \frac{a_{k+1}}{2D} \|z^k-z^{k+1}\|^2 \\ 
    &\quad +A_{k+1}\left(f(x^{k+1})-f(y^{k+1})\right) + \frac{1}{2}\|z^k-z\|^2 - \frac{1}{2}\|z^{k+1}-z\|^2\\
    &\quad  - \frac{1}{2}\left(1- \frac{a_{k+1}^2L}{A_{k+1}}\right)\|z^{k+1}-z^k\|^2 \\
    & =  a_{k+1}(\Psi(z)-\Psi(z^{k+1})) +{a}^2_{k+1}\|\omega_{k+1}\|^2+ \frac{1}{4} \|z^k-z^{k+1}\|^2 \\ 
    &\quad +A_{k+1}\left(f(x^{k+1})-f(y^{k+1})\right) + \frac{1}{2}\|z^k-z\|^2 - \frac{1}{2}\|z^{k+1}-z\|^2\\
    &\quad  - \frac{1}{2}\left(1- \frac{a_{k+1}^2L}{A_{k+1}}\right)\|z^{k+1}-z^k\|^2 
\end{align*}
where in the equality we take $D = 2a_{k+1}$. Now, $a_{k+1}^2L/A_{k+1} = (k+2)^2/(a(k+1)(k+4)) \leq 1/a$, so $a \geq 2$ gives $1/2- a_{k+1}^2L/A_{k+1} \geq 0$ and 
\begin{align}
\label{ch3:eq:hdfjdofodpf_chapter_3}
    \MoveEqLeft[1] a_{k+1}\left\la \thesischthreetg(x^{k+1}), z^k-z\right\ra \notag \\
    & \leq a_{k+1}(\Psi(z)-\Psi(z^{k+1})) +{a}^2_{k+1}\|\omega_{k+1}\|^2 \notag\\ 
    &\quad +A_{k+1}\left(f(x^{k+1})-f(y^{k+1})\right) + \frac{1}{2}\|z^k-z\|^2 - \frac{1}{2}\|z^{k+1}-z\|^2\\
    &\quad  - \frac{1}{2}\left(\frac{1}{2}- \frac{a_{k+1}^2L}{A_{k+1}}\right)\|z^{k+1}-z^k\|^2 \notag\\
    & \leq a_{k+1}(\Psi(z)-\Psi(z^{k+1})) +{a}^2_{k+1}\|\omega_{k+1}\|^2 +A_{k+1}\left(f(x^{k+1})-f(y^{k+1})\right) \notag\\
    &\quad+ \frac{1}{2}\|z^k-z\|^2 - \frac{1}{2}\|z^{k+1}-z\|^2.
\end{align}

To continue the proof, we have to mention that 
\begin{align}
\label{ch3:eq:dfjfimo_chapter_3}
    \left\la \thesischthreetg(x^{k+1}), y^k-x^{k+1}\right\ra & \overset{\eqref{ch3:eq:omega_k+1_def_dprox_clipped_SSTM_shift_chapter_3}}{=}
    \left\la \nabla f(x^{k+1}), y^k-x^{k+1}\right\ra + \left\la \omega_{k+1}, y^k-x^{k+1}\right\ra \notag\\
    & \leq  f(y^k) - f(x^{k+1}) + \left\la \omega_{k+1}, y^k-x^{k+1}\right\ra,
\end{align}
where in the last inequality we used convexity of $f$. Also, by convexity of $\Psi$ and definition of $y^{k+1}$, we have 
 \begin{equation}
\begin{aligned}
    \Psi(y^{k+1})
    &= \Psi\left(\frac{A_k}{A_{k+1}}y^{k}
    + \frac{a_{k+1}}{A_{k+1}} z^{k+1}\right)\\
    &\leq \frac{A_k}{A_{k+1}}\Psi(y^{k})
    + \frac{a_{k+1}}{A_{k+1}} \Psi(z^{k+1});\\
    -a_{k+1}\Psi(z^{k+1})
    &\leq - A_{k+1}\Psi(y^{k+1}) + A_k \Psi(y^k).
\end{aligned}
\label{ch3:eq:hfjfk_chapter_3}
\end{equation}

Thus, we acquire 
\begin{align*}
    \MoveEqLeft[1] a_{k+1} \left\la \thesischthreetg(x^{k+1}), x^{k+1}- z\right\ra \\
    & = a_{k+1} \left\la \thesischthreetg(x^{k+1}), x^{k+1}- z^k\right\ra + a_{k+1} \left\la \thesischthreetg(x^{k+1}), z^k-z\right\ra \\
    & = A_k \left\la \thesischthreetg(x^{k+1}), y^k - x^{k+1}\right\ra + a_{k+1} \left\la \thesischthreetg(x^{k+1}), z^k-z\right\ra
\end{align*}
where the last equation is true due to that $a_{k+1}(x^{k+1}-z^k) = A_k(y^k-x^{k+1})$. By \eqref{ch3:eq:hdfjdofodpf_chapter_3}, \eqref{ch3:eq:dfjfimo_chapter_3}, we get 
\begin{align*}
    a_{k+1} \left\la \thesischthreetg(x^{k+1}), x^{k+1}- z\right\ra
    &\leq A_k( f(y^k) - f(x^{k+1}))
    + A_k\left\la \omega_{k+1}, y^k-x^{k+1}\right\ra\\
    &\quad + a_{k+1}(\Psi(z)-\Psi(z^{k+1}))
    +A_{k+1}\left(f(x^{k+1})-f(y^{k+1})\right)\\
    &\quad +{a}^2_{k+1}\|\omega_{k+1}\|^2
    + \frac{1}{2}\|z^k-z\|^2 - \frac{1}{2}\|z^{k+1}-z\|^2\\
    &\overset{\eqref{ch3:eq:hfjfk_chapter_3}}{\leq}
    A_k( f(y^k) - f(x^{k+1}))
    + A_k\left\la \omega_{k+1}, y^k-x^{k+1}\right\ra\\
    &\quad + a_{k+1}\Psi(z)- A_{k+1}\Psi(y^{k+1})
    + A_k \Psi(y^k)\\
    &\quad +A_{k+1}\left(f(x^{k+1})-f(y^{k+1})\right)
    +{a}^2_{k+1}\|\omega_{k+1}\|^2\\
    &\quad + \frac{1}{2}\|z^k-z\|^2
    - \frac{1}{2}\|z^{k+1}-z\|^2.
\end{align*}
By definition of function $\Phi(\cdot)$ \eqref{ch3:eq:min_problem_chapter_3}, we have 
\begin{align*}
    \MoveEqLeft[1] a_{k+1} \left\la \thesischthreetg(x^{k+1}), x^{k+1}- z\right\ra \\
    & \leq A_k \Phi(y^k) - A_{k+1} \Phi(y^{k+1})  + A_k\left\la \omega_{k+1}, y^k-x^{k+1}\right\ra\\
    &\quad+ a_{k+1}\Psi(z) + (A_{k+1} - A_k)f(x^{k+1}) \\
    &\quad +{a}^2_{k+1}\|\omega_{k+1}\|^2 + \frac{1}{2}\|z^k-z\|^2 - \frac{1}{2}\|z^{k+1}-z\|^2\\
    & = A_k \Phi(y^k) - A_{k+1} \Phi(y^{k+1}) + a_{k+1}\left\la \omega_{k+1},x^{k+1} - z^k\right\ra\\
    &\quad+ a_{k+1}\Psi(z) +  {a}_{k+1}f(x^{k+1})\\
    &\quad +{a}^2_{k+1}\|\omega_{k+1}\|^2 + \frac{1}{2}\|z^k-z\|^2 - \frac{1}{2}\|z^{k+1}-z\|^2
\end{align*}
where in the equality we used $a_{k+1}(x^{k+1}-z^k) = A_k(y^k-x^{k+1})$ and $A_{k+1} = A_k + a_{k+1}$. Making a small rearrangement, we derive
\begin{align}
    A_{k+1} \Phi(y^{k+1}) - A_k \Phi(y^k)
    &\leq \frac{1}{2}\|z^k-z\|^2 - \frac{1}{2}\|z^{k+1}-z\|^2
    + a_{k+1}\Psi(z) \notag\\
    &\quad + {a}_{k+1}f(x^{k+1})
    + a_{k+1} \left\la \thesischthreetg(x^{k+1}),z - x^{k+1}\right\ra \notag\\
    &\quad + a_{k+1}\left\la \omega_{k+1},x^{k+1} - z^k\right\ra
    + {a}^2_{k+1}\|\omega_{k+1}\|^2 \notag\\
    &\overset{\eqref{ch3:eq:omega_k+1_def_dprox_clipped_SSTM_shift_chapter_3}}{=}
    \frac{1}{2}\|z^k-z\|^2 - \frac{1}{2}\|z^{k+1}-z\|^2
    + a_{k+1}\Psi(z) \notag\\
    &\quad + a_{k+1}\frac{1}{n}\sum^n_{i=1}f_i(x^{k+1}) \notag\\
    &\quad + a_{k+1} \left\la \frac{1}{n}\sum^n_{i=1}\nabla f_i(x^{k+1}),
    z - x^{k+1}\right\ra \notag\\
    &\quad + a_{k+1} \left\la \omega_{k+1} ,z - x^{k+1}\right\ra
    + a_{k+1}\left\la \omega_{k+1},x^{k+1} - z^k\right\ra \notag\\
    &\quad + {a}^2_{k+1}\|\omega_{k+1}\|^2 \notag\\
    &\leq \frac{1}{2}\|z^k-z\|^2 - \frac{1}{2}\|z^{k+1}-z\|^2
    + a_{k+1}\Psi(z) + a_{k+1} f(z) \notag\\
    &\quad + a_{k+1} \left\la \omega_{k+1} ,z - z^k\right\ra
    + {a}^2_{k+1}\|\omega_{k+1}\|^2 \notag\\
    &\quad -\frac{a_{k+1}}{2Ln}\sum^n_{i=1}
    \|\nabla f_i(x^{k+1}) - \nabla f_i(z)\|^2,
    \label{ch3:eq:uyerueio_chapter_3}
\end{align}
where in the last inequality we used $L$-smoothness and convexity of each $f_i$. Now we consider the sequences of $h^k_i$, produced by the method, for any $i \in [n]$. Denoting $h^*_i = \nabla f_i(x^*)$, we have 
\begin{align}
    \|h^{k+1}_i-h^*_i\|^2
    &\overset{\eqref{ch3:eq:D_prox_clipped_SSTM_shift_estimator_2_chapter_3}}{=}
    \|h^{k}_i-h^*_i\|^2 \notag\\
    &\quad + 2\nu_k \left\la \hat\Delta^k_i, h^k_i-h^*_i\right\ra
    + \nu^2_k\|\hat\Delta^k_i\|^2 \notag\\
    &= \|h^{k}_i-h^*_i\|^2 \notag\\
    &\quad + 2\nu_k \left\la \thesischthreetg_i(x^{k+1}) -h^k_i,
    h^k_i-h^*_i\right\ra
    + \nu^2_k\|\thesischthreetg_i(x^{k+1}) -h^k_i\|^2 \notag\\
    &\leq \|h^{k}_i-h^*_i\|^2 \notag\\
    &\quad + 2\nu_k \left\la \thesischthreetg_i(x^{k+1})-h^k_i,
    h^k_i-h^*_i\right\ra
    + \nu_k\|\thesischthreetg_i(x^{k+1})-h^k_i\|^2 \notag\\
    &= \|h^{k}_i-h^*_i\|^2 \notag\\
    &\quad + \nu_k \left\la \thesischthreetg_i(x^{k+1})-h^k_i,
    \thesischthreetg_i(x^{k+1}) + h^k_i- 2h^*_i\right\ra \notag\\
    &\leq (1-\nu_k)\|h^{k}_i-h^*_i\|^2
    + \nu_k\|\thesischthreetg_i(x^{k+1})-h^*_i\|^2 \notag\\
    &\leq (1-\nu_k)\|h^k_i-h^*_i\|^2 \notag\\
    &\quad + 2\nu_k\|\thesischthreetg_i(x^{k+1})-\nabla f_i(x^{k+1})\|^2
    + 2\nu_k\|\nabla f_i(x^{k+1})-h^*_i\|^2 \notag\\
    &\overset{\eqref{ch3:eq:omega_k+1_def_dprox_clipped_SSTM_shift_chapter_3}}{=}
    (1-\nu_k)\|h^{k}_i-h^*_i\|^2
    + 2\nu_k\|\omega_{i,k+1}\|^2 \notag\\
    &\quad + 2\nu_k\|\nabla f_i(x^{k+1})-\nabla f_i(x^{\star})\|^2.
    \label{ch3:eq:djkfdjksd_chapter_3}
\end{align}
where in the first inequality we use $\nu_k \leq 1$, then the identity $\la u - v, u + v - 2w\ra = \|u - w\|^2 - \|v - w\|^2$ and $\|u + v\|^2 \leq 2\|u\|^2 + 2\|v\|^2$. Summing up \eqref{ch3:eq:djkfdjksd_chapter_3} by $i$ from $1$ to $n$, we obtain 
\begin{align}
\label{ch3:eq:jdjldsdslocxnxcmn_chapter_3}
    \frac{1}{n}\sum^n_{i=1}\|h^{k+1}_i-h^*_i\|^2 & \leq (1-\nu_k)\frac{1}{n}\sum^n_{i=1}\|h^{k}_i-h^*_i\|^2 + \frac{2\nu_k}{n}\sum^n_{i=1}\|\omega_{i,k+1}\|^2\notag\\
    &\quad + \frac{2\nu_k}{n}\sum^n_{i=1}\|\nabla f_i(x^{k+1})-\nabla f_i(x^{\star})\|^2.
\end{align}

Combining inequality \eqref{ch3:eq:uyerueio_chapter_3}, where we take $z = x^*$, and inequality \eqref{ch3:eq:jdjldsdslocxnxcmn_chapter_3} multiplied by $\frac{1}{2}C^2\widetilde{a}^2_{k+2}$, we get 
\begin{align*}
    \MoveEqLeft A_{k+1}\left( \Phi(y^{k+1})-\Phi(x^*)\right)\\
    &\leq A_{k}\left( \Phi(y^{k})-\Phi(x^*)\right) + \frac{1}{2}\|z^k-x^*\|^2\\
    &\quad + \frac{1}{2}C^2\widetilde{a}^2_{k+1}\frac{1}{n}\sum^n_{i=1} \|h^k_i-h^*_i\|^2 \\
    &\quad - \frac{1}{2}\|z^{k+1}-x^*\|^2 - \frac{1}{2}C^2\widetilde{a}^2_{k+2}\frac{1}{n}\sum^n_{i=1} \|h^{k+1}_i-h^*_i\|^2\\
    &\quad + \frac{1}{2}(1-\nu_k)C^2\widetilde{a}^2_{k+2}\frac{1}{n}\sum^n_{i=1} \|h^k_i-h^*_i\|^2\\
    &\quad - \frac{1}{2}C^2\widetilde{a}^2_{k+1}\frac{1}{n}\sum^n_{i=1} \|h^k_i-h^*_i\|^2 + a_{k+1} \left\la \omega_{k+1} , x^*- z^k\right\ra\\
    &\quad +{a}^2_{k+1}\|\omega_{k+1}\|^2 + \frac{1}{2}C^2\widetilde{a}^2_{k+2}\frac{2\nu_k}{n}\sum^n_{i=1}\|\omega_{i,k+1}\|^2\\
    &\quad - \left(\frac{a_{k+1}}{2Ln} - \frac{1}{n}\nu_kC^2\widetilde{a}^2_{k+2}\right)\sum^n_{i=1}\|\nabla f_i(x^{k+1}) - \nabla f_i(x^*)\|^2.
\end{align*}
By the selection of parameters \eqref{ch3:eq:a_K0_D_prox_clipped_SSTM_shift_lemma_chapter_3}, \eqref{ch3:eq:nu_k_D_prox_clipped_SSTM_shift_lemma_chapter_3} and definition of Lyapunov function $M_k$ \eqref{ch3:eq:Lyapunov_function_D_prox_clipped_SSTM_shift_chapter_3}, we have 
\begin{align*}
    A_{k+1}\left( \Phi(y^{k+1})-\Phi(x^*)\right) & \leq A_{k}\left( \Phi(y^{k})-\Phi(x^*)\right) + \frac{1}{2}M_k - \frac{1}{2}M_{k+1}\\
    &\quad + a_{k+1} \left\la \omega_{k+1} , x^*- z^k\right\ra +{a}^2_{k+1}\|\omega_{k+1}\|^2\\
    &\quad + \frac{{a}^2_{k+1}}{n^2}\sum^n_{i=1}\|\omega_{i,k+1}\|^2.
\end{align*}
Here we use three consequences of \eqref{ch3:eq:a_K0_D_prox_clipped_SSTM_shift_lemma_chapter_3} and \eqref{ch3:eq:nu_k_D_prox_clipped_SSTM_shift_lemma_chapter_3}:
\begin{gather*}
    (1-\nu_k)\widetilde{a}_{k+2}^2 \leq \widetilde{a}_{k+1}^2,\quad C^2\widetilde{a}_{k+2}^2\nu_k \leq \frac{a_{k+1}^2}{n},\quad C^2\widetilde{a}_{k+2}^2\nu_k \leq \frac{a_{k+1}}{2L}.
\end{gather*}
For $k < K_0$ we have $\widetilde{a}_{k+1} = \widetilde{a}_{k+2} = a_{K_0+1}$, which gives the first inequality; moreover, $C^2\widetilde{a}_{k+2}^2\nu_k = (k+2)^2/(4a^2L^2n) = a_{k+1}^2/n$, and this quantity equals $a_{k+1}(k+2)/(2aLn) \leq a_{k+1}/(2L)$, since $k+2 \leq K_0+1 \leq 2C^2n + 2 \leq an$. For $k \geq K_0$ we have $\widetilde{a}_{k+1} = a_{k+1}$, $\widetilde{a}_{k+2} = a_{k+2}$ and $(1-\nu_k)(k+3)^2 = (k+2)^2$, which gives the first inequality with equality; moreover, $C^2\widetilde{a}_{k+2}^2\nu_k = C^2(2k+5)/(4a^2L^2)$, which is at most $a_{k+1}^2/n$ because $(k+2)^2 \geq C^2n(2k+5)$ for $k \geq 2C^2n$, and at most $a_{k+1}/(2L)$ because $a \geq 2C^2 + 2/n \geq C^2(2k+5)/(k+2)$. Finally, $C^2n \geq 1$ gives $\nu_k \leq 1$, which was used above. Summing up the previous inequality by $k$ from $0$ to $K-1$, we finish the proof. 

\end{proof}

\begin{theorem}
\label{ch3:thm:D_prox_clipped_SSTM_shift_full_chapter_3}
    Let Assumptions~\ref{ch3:as:bounded_alpha_moment_chapter_3}, \ref{ch3:as:L_smoothness_chapter_3} and \ref{ch3:as:str_cvx_chapter_3}($\mu = 0$) hold on $Q = B_{5n\sqrt{M}}(x^*)$, where
    \begin{equation*}
        M \geq \|x^0-x^*\|^2+C^2{a}^2_{K_0+1}\frac{1}{n}\sum^n_{i=1}\|\nabla f_i(x^*)\|^2,
    \end{equation*}
    $C = (864/\sqrt{n})\ln(10nK/\beta)$, $h_1^0 = \ldots = h_n^0 = 0$, $\nu_k$ is given by \eqref{ch3:eq:nu_k_D_prox_clipped_SSTM_shift_lemma_chapter_3}, and
    \begin{align}
        a & \geq \max\biggl\{2, \frac{8\cdot 3^5\cdot 72^4}{n} \ln^4(10nK/\beta), \notag\\
        &\qquad\qquad \frac{ 18\cdot 6^5\sigma K^{1/p}(K+1)\ln^{(p-1)/p}(10nK/\beta)}{\sqrt{M}Ln^{(p-1)/p}} \biggr\},\label{ch3:eq:Dprox_clipped_SSTM_shift_parameter_a_chapter_3}\\
        \lambda_{k} & = \frac{n\sqrt{M}}{72\widetilde{a}_{k+1}\ln(10nK/\beta)}, \label{ch3:eq:Dprox_clipped_SSTM_shift_clipping_level_chapter_3}
    \end{align}
    for some $K \geq K_0 = \lceil 2C^2n\rceil$ and $\beta \in (0,1]$ such that $\ln(10nK/\beta) \geq 1$. Then, after $K$ iterations of~\algname{DProx-clipped-SSTM-shift} the following inequality holds with probability at least $1 - \beta$ 
    \begin{align}
        & \Phi(y^K) - \Phi(x^*) \leq \frac{6aL M}{K(K+3)} \quad \text{and} \notag\\
        & \{x^k\}_{k=0}^{K+1}, \{z^k\}_{k=0}^{K}, \{y^k\}_{k=0}^K \subseteq B_{2\sqrt{M}}(x^*). \label{ch3:eq:Dprox_clipped_SSTM_shift_chapter_3}
    \end{align}
    In particular, when parameter $a$ equals the maximum from \eqref{ch3:eq:Dprox_clipped_SSTM_shift_parameter_a_chapter_3}, then after $K$ iterations of \algname{DProx-clipped-SSTM-shift}, we have with probability at least $1-\beta$
    \begin{align}
        \Phi(y^K) - \Phi(x^*) = \cO\biggl(\max\biggl\{ & \frac{LM}{K^2}, \frac{LM\ln^4(nK/\beta)}{nK^2}, \notag\\
        & \frac{\sigma \sqrt{M} \ln^{(p-1)/p}(nK/\beta)}{(nK)^{(p-1)/p}}\biggr\}\biggr), \label{ch3:eq:Dprox_clipped_SSTM_shift_2_chapter_3}
    \end{align}
    i.e. achieve $\Phi(y^K) - \Phi(x^*) \leq \varepsilon$ with probability at least $1 - \beta$ \algname{DProx-clipped-SSTM-shift} requires
    \begin{align}
        K = \cO\biggl(\max\biggl\{ & \sqrt{\frac{LM}{\varepsilon}}, \sqrt{\frac{LM}{\varepsilon n}}\ln^2\frac{nLM}{\varepsilon\beta},\notag\\
        & \frac{1}{n}\left(\frac{\sigma \sqrt{M}}{\varepsilon}\right)^{p/(p-1)}\ln\frac{\sigma \sqrt{M}}{\varepsilon\beta}\biggr\}\biggr) \label{ch3:eq:Dprox_clipped_SSTM_shift_complexity_chapter_3}
    \end{align}
    iterations/oracle calls per worker.
\end{theorem}
\begin{proof}
    The key idea behind the proof is similar to the one used by \textcite{gorbunov2021near,sadiev2023high}. We prove by induction that the iterates do not leave some ball and $\Phi(y^K)-\Phi(x^*)$ decreases as $1/(K(K+3))$.
    
    Firstly, we denote $R_k = \|z^k - x^*\|$, $\widetilde{R}_0 = R_0$, $\widetilde{R}_{k+1} = \max\{\widetilde{R}_k, R_{k+1}\}$ for all $k\geq 0$, and now we show by induction that for all $k\geq 0$ the iterates $x^{k+1}, z^k, y^k$ lie in $B_{\widetilde{R}_k}(x^*)$. The induction base is trivial since $y^0 = z^0$, $\widetilde{R}_0 = R_0$, and $x^1 = (A_0 y^0 + a_1 z^0)/A_1 = z^0$. Next, we assume this statement is true for some  $l \geq 1$: $x^{l}, z^{l-1}, y^{l-1} \in B_{\widetilde{R}_{l-1}}(x^*)$. According to definitions of $R_l$ and $\widetilde{R}_l$, we obtain $z^l \in B_{R_{l}}(x^*)\subseteq B_{\widetilde{R}_{l}}(x^*)$. Due to that $y^l$ is a convex combination of $y^{l-1}\in B_{\widetilde{R}_{l-1}}(x^*)\subseteq B_{\widetilde{R}_{l}}(x^*)$, $z^l\in B_{\widetilde{R}_{l}}(x^*)$ and $B_{\widetilde{R}_{l}}(x^*)$ is a convex set, we have that $y^l \in B_{\widetilde{R}_{l}}(x^*)$. Finally, since $x^{l+1}$ is a convex combination of $y^l$ and $z^l$, we conclude $x^{l+1}$ lies in $B_{\widetilde{R}_{l}}(x^*)$ as well.

    Now to formulate the statement rigorously, we introduce probability event $E_k$ for each $k = 0,\ldots, K$ as follows: inequalities
    \begin{gather}
         \underbrace{2\sum\limits^{t-1}_{l=0}a_{l+1}\la\omega_{l+1}, x^*-z^l\ra  + 2\sum\limits^{t-1}_{l=0}{a}^2_{l+1}\|\omega_{l+1}\|^2 + 2\sum\limits^{t-1}_{l=0}\sum^{n}_{i=1}\frac{{a}^2_{l+1}}{n^2}\|\omega_{i,l+1}\|^2}_{B_t} \leq M, \label{ch3:eq:SSTM_induction_inequality_1_chapter_3}\\
        R_t\leq  \sqrt{M_t} \leq 2\sqrt{M} \label{ch3:eq:SSTM_induction_inequality_2_chapter_3},\\
        \left\|\frac{a_{t}}{n}\sum^{r-1}_{i=1}\omega^u_{i,t}\right\| \leq \frac{\sqrt{M}}{2}\label{ch3:eq:SSTM_induction_inequality_3_chapter_3}
    \end{gather}
    hold for  $t = 0,1,\ldots, k$ and $r=1,2\dots,n$ simultaneously, where 
     \begin{gather}
        \omega_{l+1} = \omega_{l+1}^u +\omega^b_{l+1}, \label{ch3:eq:SSTM_omega_full_chapter_3}\\
        \omega_{l+1}^u \eqdef \frac{1}{n}\sum_{i=1}^n \omega_{i,l+1}^u, \quad \omega_{l+1}^b \eqdef \frac{1}{n}\sum_{i=1}^n \omega_{i,l+1}^b, \label{ch3:eq:SSTM_omega_u_omega_b_chapter_3}\\
        \omega_{i,l+1}^u \eqdef \thesischthreetg_i(x^{l+1})  - \Exp_{\xi^l_i}\left[\thesischthreetg_i(x^{l+1}) \right], \notag\\
        \omega_{i,l+1}^b \eqdef \Exp_{\xi^l_i}\left[\thesischthreetg_i(x^{l+1}) \right] - \nabla f_i(x^{l+1}) \quad \forall i \in [n]. \label{ch3:eq:SSTM_omega_i_chapter_3}
    \end{gather}

     We want to show via induction $\widetilde{R}_{l} \leq 5n\sqrt{M}$ with high probability, which allows us to apply the result of Lemma~\ref{ch3:lem:SSTM_optimization_lemma_convex_chapter_3} and  Bernstein's inequality to estimate the stochastic part of the upper-bound. After that, we will prove by induction that $\Prob\{E_k\} \geq 1 - k\beta/K$ for all $k = 0,1,\ldots, K$. The base induction follows immediately: the left-hand side of \eqref{ch3:eq:SSTM_induction_inequality_1_chapter_3} equals zero and $M \geq M_0$ by definition, and for $k=0$ we have  $\left\|\frac{a_{0}}{n}\sum\limits^r_{i=1}\omega^u_{i,0}\right\| = 0$, since $a_0 = 0$. Next we assume that the statement holds for some $k = T - 1 \leq K-1$: $\Prob\{E_{T-1}\} \geq 1 - (T-1)\beta/K$. Let us show that $\Prob\{E_{T}\} \geq 1 - T\beta/K$. 
     
     To proceed, we note that $E_{T-1}$ implies $R_t \leq 2\sqrt{M}$, and hence $\widetilde{R}_t \leq 2\sqrt{M}$ and $x^{t+1}, y^t, z^t \in B_{2\sqrt{M}}(x^*)$, for all $t = 0,1,\ldots, T-1$. For $t = T$, the event $E_{T-1}$ implies
    \begin{align*}
        \|z^{t} - x^{*}\|& = \|\thesischthreeprox_{a_{t}\Psi}\left(z^{t-1} - a_{t}\thesischthreetg (x^{t})\right) - \thesischthreeprox_{a_{t}\Psi}\left(x^{*} - a_{t}\nabla f(x^{*})\right)\| \\
        & \leq\|z^{t-1} - x^* -a_{t}\left(\thesischthreetg(x^{t})- \nabla f(x^{*})\right)\|\\
        & \leq\|z^{t-1} - x^*\| + a_{t}\|\thesischthreetg(x^t)-h^{t-1}\| + a_{t}\|h^{t-1}- h^*\| \\
        & \leq \left(1 + \frac{1}{C}\right)\sqrt{M_{t-1}} + a_{t}\lambda_{t-1}\\
        & \overset{\eqref{ch3:eq:Dprox_clipped_SSTM_shift_clipping_level_chapter_3},\eqref{ch3:eq:SSTM_induction_inequality_2_chapter_3}}{\leq} 2\left(1 + \frac{\sqrt{n}}{864}\right)\sqrt{M} + \frac{n\sqrt{M}}{72} \leq 5n\sqrt{M},
    \end{align*}
    where $h^{t-1} = \frac{1}{n}\sum_{i=1}^n h_i^{t-1}$ and $h^* = \nabla f(x^*)$; in the fourth step we use $\|h^{t-1} - h^*\|^2 \leq \frac{1}{n}\sum_{i=1}^n\|h_i^{t-1} - h_i^*\|^2$, $a_t \leq \widetilde{a}_t$ and $\|\thesischthreetg_i(x^t) - h_i^{t-1}\| = \|\hat\Delta_i^{t-1}\| \leq \lambda_{t-1}$, and in the last step $1/C \leq \sqrt{n}/864$.

    This means that $x^t, z^t, y^t \in B_{5n\sqrt{M}}(x^*)$ for $t=0,1,\dots, T$, and we can apply Lemma~\ref{ch3:lem:SSTM_optimization_lemma_convex_chapter_3} with $K = t \leq T$: $E_{T-1}$ implies
    \begin{align}
        A_{t}\left(\Phi(y^{t}) - \Phi(x^*)\right) & \leq \frac{1}{2}M_0 - \frac{1}{2}M_{t} + \sum^{t-1}_{l=0}a_{l+1}\la\omega_{l+1}, x^*-z^l\ra  + \sum^{t-1}_{l=0}{a}^2_{l+1}\|\omega_{l+1}\|^2\notag\\
        &\quad + \sum^{t-1}_{l=0}\sum^{n}_{i=1}\frac{{a}^2_{l+1}}{n^2}\|\omega_{i,l+1}\|^2\notag\\
        & = \frac{1}{2}M_0 -\frac{1}{2}M_{t} +\frac{1}{2}B_{t}. \label{ch3:eq:D_prox_clipped_SSTM_shift_technical_1_chapter_3}
    \end{align}
    For $t \leq T-1$, inequality \eqref{ch3:eq:SSTM_induction_inequality_1_chapter_3} gives $B_t \leq M$, hence $M_{t} \leq M_0 + M \leq 2M$ and
    \begin{equation}
        \Phi(y^{t}) - \Phi(x^*) \leq \frac{M}{A_{t}} \leq \frac{3M}{2A_{t}} = \frac{6aLM}{t(t+3)} \label{ch3:eq:D_prox_clipped_SSTM_shift_technical_1_1_chapter_3}
    \end{equation}
    for all $t = 1, \ldots, T-1$, where we use $A_t = t(t+3)/(4aL)$.  Before we proceed, we introduce a new notation:
    \begin{equation*}
        u^{t} = \begin{cases} x^* - z^t,& \text{if } \|x^* - z^t\| \leq 2\sqrt{M},\\ 0,& \text{otherwise,} \end{cases}
    \end{equation*}
    for all $t = 0, 1, \ldots, T-1$. Random vectors $\{u^{t}\}_{t=0}^{T-1}$ are bounded almost surely:
     \begin{equation}
        \|u^{t}\| \leq 2\sqrt{M}  \label{ch3:eq:gap_thm_SSTM_technical_3_chapter_3}
    \end{equation}
    for all $t = 0, 1, \ldots, T-1$. In addition, $u^{t} = x^* - z^t$ follows from $E_{T-1}$ for all $t = 0,1, \ldots, T-1$ and, thus, $E_{T-1}$ implies 
    \begin{align}
        B_T & = 2 \sum^{T-1}_{k=0}a_{k+1}\la\omega_{k+1}, x^*-z^k\ra  + 2\sum^{T-1}_{k=0}{a}^2_{k+1}\|\omega_{k+1}\|^2 \notag\\
        &\quad + 2\sum^{T-1}_{k=0}\sum^{n}_{i=1}\frac{{a}^2_{k+1}}{n^2}\|\omega_{i,k+1}\|^2 \notag\\
        & = 2 \sum^{T-1}_{k=0}a_{k+1}\la\omega_{k+1}, u^{k}\ra  + 2\sum^{T-1}_{k=0}{a}^2_{k+1}\|\omega_{k+1}\|^2 \notag\\
        &\quad + 2\sum^{T-1}_{k=0}\sum^{n}_{i=1}\frac{{a}^2_{k+1}}{n^2}\|\omega_{i,k+1}\|^2 . \label{ch3:eq:D_prox_clipped_SSTM_shift_technical_2_chapter_3}
    \end{align}
    Using the notation from \eqref{ch3:eq:SSTM_omega_full_chapter_3}-\eqref{ch3:eq:SSTM_omega_i_chapter_3}, we can rewrite $\|\omega_{k+1}\|^2$ and $\|\omega_{i,k+1}\|^2$ as
    \begin{align}
        \|\omega_{k+1}\|^2 & \leq \frac{2}{n^2}\sum\limits_{i=1}^n\|\omega_{i,k+1}^u\|^2 + \frac{4}{n^2}\sum\limits_{j=2}^n\left\langle \sum\limits_{i=1}^{j-1} \omega_{i,k+1}^u, \omega_{j,k+1}^u \right\rangle + 2\|\omega_{k+1}^b\|^2. \label{ch3:eq:SSTM_sum_omega_representation_chapter_3}
    \end{align}
    Putting all together, we obtain that $E_{T-1}$  implies
    \begin{align}
        B_T 
        & \leq \underbrace{2\sum\limits_{k=0}^{T-1}\sum^n_{i=1}\frac{a_{k+1}}{n}\left\la \omega_{i,k+1}^u, u^{k}\right\ra}_{\circledOne} + \underbrace{2\sum\limits_{k=0}^{T-1}\sum^n_{i=1}\frac{a_{k+1}}{n}\left\la \omega_{i,k+1}^b, u^{k}\right\ra}_{\circledTwo} \notag\\
        &\quad\quad+ \underbrace{8\sum\limits_{k=0}^{T-1}\sum^n_{i=1}\frac{a_{k+1}^2}{n^2}\left(\left\|\omega_{i,k+1}^u\right\|^2 - \Exp_{\xi^k_i}\left[\left\|\omega_{i,k+1}^u\right\|^2\right]\right)}_{\circledThree}\notag\\
        &\quad\quad + \underbrace{8\sum\limits_{k=0}^{T-1}\sum^n_{i=1}\frac{a_{k+1}^2}{n^2}\Exp_{\xi^k}\left[\left\|\omega_{i,k+1}^u\right\|^2\right]}_{\circledFour} + \underbrace{8\sum\limits_{k=0}^{T-1}\sum^n_{i=1}\frac{a_{k+1}^2}{n}\left\|\omega_{i,k+1}^b\right\|^2}_{\circledFive}\notag\\
        &\quad\quad +\underbrace{8\sum\limits_{k=0}^{T-1}\sum^n_{j=2}\frac{a_{k+1}^2}{n^2}\left\la\sum^{j-1}_{i=1}\omega^u_{i,k+1}, \omega^u_{j, k+1}\right\ra}_{\circledSix}. \label{ch3:eq:D_prox_clipped_SSTM_shift_technical_7_chapter_3}
    \end{align}

    To finish the proof, it remains to estimate  $\circledOne, \circledTwo, \circledThree, \circledFour, \circledFive, \circledSix$ with high probability. More precisely, the goal to prove  that $\circledOne + \circledTwo + \circledThree + \circledFour + \circledFive + \circledSix  \leq M$ with high probability.
    Before we proceed, we need to derive several useful inequalities related to $\omega^u_{i,k+1},~\omega^b_{i,k+1}$. First of all, we have 
    \begin{equation}
        \|\omega_{i,k+1}^u\| \leq 2\lambda_{k}. \label{ch3:eq:D_prox_clipped_SSTM_shift_norm_omega_u_bound_chapter_3}
    \end{equation}
    by definition of the clipping operator. Next, probability event $E_{T-1}$ implies that for $t = 0,1,\ldots,\min\{K_0, T-1\}$
    \begin{align}
        \|\nabla f_i(x^{t+1}) - h^t_i\| &\leq \|\nabla f_i(x^{t+1})- \nabla f_i(x^{*})\| + \|h^t_i -h^*_i\| \notag\\
        &\leq L\|x^{t+1} - x^*\| + \frac{\sqrt{n}}{C\widetilde{a}_{t+1}}\sqrt{C^2\widetilde{a}^2_{t+1}\frac{1}{n}\sum^n_{i=1}\|h^t_i-h^*_i\|^2}\notag\\
        &\leq 2L\sqrt{M} + \frac{2\sqrt{nM}}{C\widetilde{a}_{t+1}} \notag\\
        &\overset{\eqref{ch3:eq:Dprox_clipped_SSTM_shift_clipping_level_chapter_3}}{=} \frac{\lambda_t}{2}\left(\frac{144 (K_0+2)\ln(10nK/\beta)}{an} + \frac{288\ln(10nK/\beta)}{C\sqrt{n}}\right) \notag\\
        &\overset{\eqref{ch3:eq:Dprox_clipped_SSTM_shift_parameter_a_chapter_3}}{\leq} \frac{\lambda_t}{2}, \label{ch3:eq:D_prox_clipped_SSTM_shift_technical_3_chapter_3}
    \end{align}
    where in the third step we use $x^{t+1} \in B_{2\sqrt{M}}(x^*)$ and $M_t \leq 4M$, and in the last two steps $\widetilde{a}_{t+1} = a_{K_0+1} = (K_0+2)/(2aL)$ for $t \leq K_0$, $C\sqrt{n} = 864\ln(10nK/\beta)$ and $K_0 + 2 \leq 2C^2n + 3 \leq 1.5\cdot 10^6\ln^2(10nK/\beta)$. Next, for $t = K_0+1, \ldots, T-1$ event $E_{T-1}$ implies
    \begin{align}
        \MoveEqLeft[1] \|\nabla f_i(x^{t+1})-h^t_i\| \notag\\
        &\leq \|\nabla f_i(x^{t+1}) - \nabla f_i(y^t)\| + \|\nabla f_i(y^t) - \nabla f_i(x^{*})\| + \|h^t_i - h^*_i\| \notag\\
        &\leq L\|x^{t+1} - y^t\| + \sqrt{2L\left(f_i(y^t) - f_i(x^*) - \la \nabla f_i(x^{*}), y^t-x^{*}\ra\right)} \notag\\
        &\quad + \|h^t_i - h^*_i\| \notag\\
        &\leq L\|x^{t+1} - y^t\| + \sqrt{2nL\left(\Phi(y^t) - \Phi(x^*)\right)} + \sqrt{\sum^n_{i=1}\|h^t_i - h^*_i\|^2}\\
        &\overset{\eqref{ch3:eq:D_prox_clipped_SSTM_shift_technical_1_1_chapter_3}}{\leq} \frac{L{a}_{t+1}}{A_t}\|x^{t+1} - z^t\| + \sqrt{\frac{12anL^2M}{t(t+3)}} + \frac{2\sqrt{nM}}{C\widetilde{a}_{t+1}}\notag\\
        &\leq \frac{4L\sqrt{M}a_{t+1}}{A_t} + \sqrt{\frac{12anL^2M}{t(t+3)}} + \frac{2\sqrt{nM}}{C\widetilde{a}_{t+1}} \notag\\
        &\overset{\eqref{ch3:eq:Dprox_clipped_SSTM_shift_clipping_level_chapter_3}}{=} \frac{\lambda_t}{2}\biggl(\frac{576(t+2)^2\ln(10nK/\beta)}{an\, t(t+3)} + \frac{72\sqrt{12}(t+2)\ln(10nK/\beta)}{\sqrt{an\, t(t+3)}}\notag\\
        &\qquad\quad + \frac{288\ln(10nK/\beta)}{C\sqrt{n}}\biggr)\notag\\
        &\leq \frac{\lambda_t}{2}\left(\frac{1296\ln(10nK/\beta)}{an} + \frac{108\sqrt{12}\ln(10nK/\beta)}{\sqrt{an}} + \frac{1}{3}\right) \notag\\
        &\overset{\eqref{ch3:eq:Dprox_clipped_SSTM_shift_parameter_a_chapter_3}}{\leq} \frac{\lambda_t}{2}, \label{ch3:eq:D_prox_clipped_SSTM_shift_technical_4_chapter_3}
    \end{align}
    where in the third step we use $- \la \nabla f(x^{*}), y^t-x^{*}\ra \leq \Psi(y^t) - \Psi(x^{*})$, in the fourth step $\|x^{t+1} - y^t\| = (a_{t+1}/A_t)\|x^{t+1} - z^t\|$ and $M_t \leq 4M$, in the fifth step $x^{t+1}, z^t \in B_{2\sqrt{M}}(x^*)$, and in the last steps $\widetilde{a}_{t+1} = a_{t+1} = (t+2)/(2aL)$ for $t > K_0$, $A_t = t(t+3)/(4aL)$, $(t+2)^2/(t(t+3)) \leq 9/4$ for $t \geq 1$ and $C\sqrt{n} = 864\ln(10nK/\beta)$.

    Therefore, Lemma~\ref{ch3:lem:bias_and_variance_clip_chapter_3} and  $E_{T-1}$ imply
    \begin{align}
        \|\omega_{i,k+1}^b\| & \leq \frac{2^{p}\sigma^{p}}{\lambda_{k}^{{p}-1}}, \label{ch3:eq:D_prox_clipped_SSTM_shift_norm_omega_b_bound_chapter_3} \\
        \Exp_{\xi^k_i}\left[\|\omega_{i,k+1}^u\|^2\right] & \leq 18 \lambda_k^{2-{p}}\sigma^{p}. \label{ch3:eq:D_prox_clipped_SSTM_shift_second_moment_omega_u_bound_chapter_3}
    \end{align}

    \paragraph{Upper bound for $\circledOne$.} To estimate this sum, we will use Bernstein's inequality. The summands have conditional equal to zero, since  $\Exp_{\xi^k_i}[\omega_{i,k+1}^u] = 0$:
    \begin{equation}
        \Exp_{\xi^k_i}\left[\frac{a_{k+1}}{n}\left\la \omega_{i,k+1}^u, u^{k}\right\ra\right] = 0. \notag
    \end{equation}
    Moreover, for all $k = 0, \ldots, T-1$ random vectors $\{\omega_{i,k+1}^u\}_{i=1}^n$ are independent. Thus, sequence $\{2(a_{k+1}/n)\la \omega_{i,k+1}^u, u^{k}\ra\}_{k,i=0,1}^{T-1,n}$ is a martingale difference sequence. Next, the summands are bounded:
    \begin{align}
        \left|2\frac{a_{k+1}}{n}\left\la \omega_{i,k+1}^u, u^{k}\right\ra\right| & \leq 2\frac{a_{k+1}}{n} \|\omega_{i,k+1}^u\| \cdot \|u^{k}\| \notag \\
        & \overset{\eqref{ch3:eq:gap_thm_SSTM_technical_3_chapter_3},\eqref{ch3:eq:D_prox_clipped_SSTM_shift_norm_omega_u_bound_chapter_3}}{\leq} 8\frac{a_{k+1}}{n}\lambda_k \sqrt{M} \notag\\
        & \overset{\eqref{ch3:eq:Dprox_clipped_SSTM_shift_clipping_level_chapter_3}}{=} \frac{8{a}_{k+1}M}{72\widetilde{a}_{k+1}\ln(10nK/\beta)}\leq \frac{M}{9\ln(10nK/\beta)} \eqdef c. \label{ch3:eq:D_prox_clipped_SSTM_shift_technical_8_chapter_3} 
    \end{align}
    
    Finally, conditional variances $\sigma_{i,k}^2 \eqdef \Exp_{\xi^k_i}\left[4\frac{a_{k+1}^2}{n^2}\left\la \omega_{i,k+1}^u, u^{k}\right\ra^2\right]$ of the summands are bounded:
    \begin{align}
\sigma_{i,k}^2 & \leq \Exp_{\xi^k_i}\left[4\frac{a_{k+1}^2}{n^2}\|\omega_{i,k+1}^u\|^2\cdot \|u^{k}\|^2\right] \notag \\
& \leq  16\frac{a_{k+1}^2}{n^2}M \Exp_{\xi^k_i}\left[\|\omega_{i,k+1}^u\|^2\right]. \label{ch3:eq:D_prox_clipped_SSTM_shift_technical_9_chapter_3}
    \end{align}
    Applying Bernstein's inequality (Lemma~\ref{ch3:lem:Bernstein_ineq_chapter_3}) with $X_{i,k} = 2\frac{a_{k+1}}{n}\left\la \omega_{i,k+1}^u, u^{k}\right\ra$, parameter $c$ as in \eqref{ch3:eq:D_prox_clipped_SSTM_shift_technical_8_chapter_3}, $b = \frac{M}{6}$, $H = \frac{M^2}{6^3\ln(10nK/\beta)}$:
    \begin{align*}
\MoveEqLeft[1] \Prob\left\{|\circledOne| > \frac{M}{6}\quad \text{and}\quad \sum\limits_{k=0}^{T-1}\sum^n_{i=1} \sigma_{i,k}^2 \leq \frac{M^2}{6^3\ln(10nK/\beta)}\right\} \\
& \leq 2\exp\left(- \frac{b^2}{2H + 2cb/3}\right) \leq \frac{\beta}{5nK}.
    \end{align*}
    The above is equivalent to 
    \begin{align}
\Prob\left\{ E_{\circledOne} \right\} & \geq 1 - \frac{\beta}{5nK},\quad \text{for}\quad E_{\circledOne} \notag \\
& = \left\{ \text{either} \quad  \sum\limits_{k=0}^{T-1}\sum^n_{i=1} \sigma_{i,k}^2 > \frac{M^2}{6^3\ln(10nK/\beta)} \quad \text{or}\quad |\circledOne| \leq \frac{M}{6}\right\}. \label{ch3:eq:D_prox_clipped_SSTM_shift_sum_1_upper_bound_chapter_3}
    \end{align}
    Moreover, $E_{T-1}$ implies that
    \begin{align}
        \sum\limits_{k=0}^{T-1}\sum^n_{i=1} \sigma_{i, k}^2 & \overset{\eqref{ch3:eq:D_prox_clipped_SSTM_shift_technical_9_chapter_3}}{\leq} 16 M \sum\limits_{k=0}^{T-1}\sum^n_{i=1} \frac{a_{k+1}^2}{n^2} \Exp_{\xi^k_i}\left[\|\omega_{i,k+1}^u\|^2\right] \notag \\
        & \overset{\eqref{ch3:eq:D_prox_clipped_SSTM_shift_second_moment_omega_u_bound_chapter_3}}{\leq} 288\sigma^{p} M \sum\limits_{k=0}^{T-1}\sum^n_{i=1} \frac{a_{k+1}^2}{n^2}\lambda_k^{2-{p}} \notag\\
        & \overset{\eqref{ch3:eq:Dprox_clipped_SSTM_shift_clipping_level_chapter_3}}{\leq} \frac{288\sigma^{p} M^{2-{p}/2}}{72^{2-{p}}\ln^{2-{p}}(10nK/\beta)}\sum\limits_{k=0}^{T-1} \frac{a_{k+1}^{2}}{n^{{p}-1}\widetilde{a}_{k+1}^{2-{p}}} \notag \\
        & \leq \frac{288\sigma^{p} M^{2-{p}/2}}{72^{2-{p}}\ln^{2-{p}}(10nK/\beta)}\sum\limits_{k=0}^{T-1} \frac{a_{k+1}^{{p}}}{n^{{p}-1}} \notag\\
        & \leq \frac{288\sigma^{p} M^{2-{p}/2}}{n^{{p}-1}72^{2-{p}}\cdot 2^{p} a^{p} L^{p}\ln^{2-{p}}(10nK/\beta)}\sum\limits_{k=0}^{T-1} (k+2)^{p} \notag\\
        & \leq \frac{1}{a^{p}} \cdot \frac{144\sigma^{p} M^{2-{p}/2} T(T+1)^{p}}{n^{{p}-1} L^{p}\ln^{2-{p}}(10nK/\beta)} \notag \\
        & \overset{\eqref{ch3:eq:Dprox_clipped_SSTM_shift_parameter_a_chapter_3}}{\leq} \frac{M^2}{6^3 \ln(10nK/\beta)}. \label{ch3:eq:prox_clipped_SSTM_shift_sum_1_variance_bound_chapter_3}
    \end{align}

     \paragraph{Upper bound for $\circledTwo$.} Probability event $E_{T-1}$ implies
    \begin{align}
        \circledTwo & \leq 2\sum\limits_{k=0}^{T-1}\sum^n_{i=1}\frac{a_{k+1}}{n}\|\omega_{i,k+1}^b\|\cdot \|u^{k}\| \notag\\
        & \overset{\eqref{ch3:eq:gap_thm_SSTM_technical_3_chapter_3},\eqref{ch3:eq:D_prox_clipped_SSTM_shift_norm_omega_b_bound_chapter_3}}{\leq} 4\sqrt{M}\cdot 2^{p}\sigma^{p} \sum\limits_{k=0}^{T-1}\frac{a_{k+1}}{\lambda_{k}^{{p}-1}}\notag\\
        & \overset{\eqref{ch3:eq:Dprox_clipped_SSTM_shift_clipping_level_chapter_3}}{\leq} \frac{16\cdot 72^{{p}-1}M^{1-{p}/2}\sigma^{p}\ln^{{p}-1}(10nK/\beta)}{n^{{p}-1}} \sum\limits_{k=0}^{T-1}a_{k+1}\widetilde{a}_{k+1}^{{p}-1}\notag\\
        & = \frac{16\cdot 72^{{p}-1}\sigma^{p} M^{1 - {p}/2}\ln^{{p}-1}(10nK/\beta)}{2^{p} a^{{p}} L^{p}n^{{p}-1}}\sum\limits_{k=0}^{T-1} (k+2)\max\{K_0+2, k+2\}^{{p}-1} \notag\\
        & \leq \frac{1}{a^{p}} \cdot \frac{16\cdot 72^{{p}-1}\sigma^{p} M^{1 - {p}/2} K(K+1)^{p}\ln^{{p}-1}(10nK/\beta)}{n^{{p}-1}L^{p}} \notag\\
        & \overset{\eqref{ch3:eq:Dprox_clipped_SSTM_shift_parameter_a_chapter_3}}{\leq} \frac{M}{6}. \label{ch3:eq:D_prox_clipped_SSTM_shift_sum_2_upper_bound_chapter_3}
    \end{align}
    where in the last but one step we use $T \leq K$, $K_0 \leq K$ and $(K+2)^{{p}-1} \leq 2^{{p}-1}(K+1)^{{p}-1}$.

    \paragraph{Upper bound for $\circledThree$.} To estimate this sum, we will use Bernstein's inequality. The summands have conditional expectations equal to zero: 
    \begin{equation}
        \frac{8a_{k+1}^2}{n^2}\Exp_{\xi^k_i}\left[\left\|\omega_{i,k+1}^u\right\|^2 - \Exp_{\xi^k_i}\left[\left\|\omega_{i,k+1}^u\right\|^2\right]\right] = 0. \notag
    \end{equation}
    Moreover, for all $k = 0, \ldots, T-1$ random vectors $\left\{\omega^u_{i,k+1}\right\}^n_{i=1}$ are independent. Thus, sequence $\left\{ \frac{8a_{k+1}^2}{n^2}\left(\left\|\omega_{i,k+1}^u\right\|^2 - \Exp_{\xi^k_i}\left[\left\|\omega_{i,k+1}^u\right\|^2\right]\right) \right\}_{k,i = 0,1}^{T-1, n}$ is a martingale difference sequence. Next, the summands are bounded: 
    \begin{align}
        \MoveEqLeft \left|\frac{8a_{k+1}^2}{n^2}\left(\left\|\omega_{i,k+1}^u\right\|^2 - \Exp_{\xi^k_i}\left[\left\|\omega_{i,k+1}^u\right\|^2\right]\right)\right| \notag\\
        & \leq \frac{8a_{k+1}^2}{n^2}\left( \|\omega_{i,k+1}^u\|^2 +   \Exp_{\xi^k_i}\left[\left\|\omega_{i,k+1}^u\right\|^2\right]\right)\notag\\
        & \overset{\eqref{ch3:eq:D_prox_clipped_SSTM_shift_norm_omega_u_bound_chapter_3}}{\leq} \frac{64a_{k+1}^2\lambda_k^2}{n^2} \overset{\eqref{ch3:eq:Dprox_clipped_SSTM_shift_clipping_level_chapter_3}}{\leq} \frac{M}{9\ln(10nK/\beta)} \eqdef c. \label{ch3:eq:D_prox_clipped_SSTM_shift_technical_10_chapter_3}
    \end{align}

   Finally, conditional variances $$\widetilde\sigma_{i,k}^2 \eqdef \Exp_{\xi^k_i}\left[\frac{64 a_{k+1}^4}{n^4}\left(\left\|\omega_{i,k+1}^u\right\|^2 - \Exp_{\xi^k_i}\left[\left\|\omega_{i,k+1}^u\right\|^2\right]\right)^2\right]$$
   of the summands are bounded: 
    \begin{align}
        \widetilde\sigma_{i,k}^2 & \overset{\eqref{ch3:eq:D_prox_clipped_SSTM_shift_technical_10_chapter_3}}{\leq} \frac{8 a_{k+1}^2 M}{9 n^2 \ln(10nK/\beta)} \Exp_{\xi^k_i}\left[\left|\left\|\omega_{i,k+1}^u\right\|^2 - \Exp_{\xi^k_i}\left[\left\|\omega_{i,k+1}^u\right\|^2\right]\right|\right]\notag\\ 
        & \leq \frac{16a_{k+1}^2M}{9n^2} \Exp_{\xi^k_i}\left[\|\omega_{i,k+1}^u\|^2\right]. \label{ch3:eq:D_prox_clipped_SSTM_shift_technical_11_chapter_3}
    \end{align}
    
    Applying Bernstein's inequality (Lemma~\ref{ch3:lem:Bernstein_ineq_chapter_3}) with
    \begin{equation*}
        \widetilde X_{i,k} = \frac{8a_{k+1}^2}{n^2}\left(\left\|\omega_{i,k+1}^u\right\|^2 - \Exp_{\xi^k_i}\left[\left\|\omega_{i,k+1}^u\right\|^2\right]\right),
    \end{equation*}
    parameter $c$ defined in \eqref{ch3:eq:D_prox_clipped_SSTM_shift_technical_10_chapter_3}, $b = \frac{M}{9}$, $H = \frac{M^2}{6\cdot9^2\ln(10nK/\beta)}$:
    \begin{multline*}
        \Prob\left\{|\circledThree| > \frac{M}{9}\quad \text{and}\quad \sum\limits_{k=0}^{T-1}\sum^n_{i=1} \widetilde\sigma_{i,k}^2 \leq \frac{M^2}{6\cdot 9^2\ln(10nK/\beta)}\right\}\\
        \leq 2\exp\left(- \frac{b^2}{2H + 2cb/3}\right) = \frac{\beta}{5nK}.
    \end{multline*}
    The above is equivalent to
    \begin{align}
\Prob\left\{ E_{\circledThree} \right\} & \geq 1 - \frac{\beta}{5nK},\quad \text{for}\quad E_{\circledThree} \notag \\
& = \left\{ \text{either} \quad  \sum\limits_{k=0}^{T-1}\sum^n_{i=1} \widetilde\sigma_{i,k}^2 > \frac{M^2}{6\cdot 9^2\ln(10nK/\beta)} \quad \text{or}\quad |\circledThree| \leq \frac{M}{9}\right\}. \label{ch3:eq:D_prox_clipped_SSTM_shift_sum_3_upper_bound_chapter_3}
    \end{align}
    Moreover, $E_{T-1}$ implies 
    \begin{align}
        \sum\limits_{k=0}^{T-1}\sum^n_{i=1} \widetilde\sigma_{i,k}^2 & \overset{\eqref{ch3:eq:D_prox_clipped_SSTM_shift_technical_11_chapter_3}}{\leq} \frac{16}{9}M \sum\limits_{k=0}^{T-1} \sum^n_{i=1}\frac{a_{k+1}^2}{n^2} \Exp_{\xi^k_i}\left[\|\omega_{i,k+1}^u\|^2\right] \notag \\
        & \overset{\eqref{ch3:eq:prox_clipped_SSTM_shift_sum_1_variance_bound_chapter_3}}{\leq} \frac{M^2}{6 \cdot 9^2\ln(10nK/\beta)}. \label{ch3:eq:D_prox_clipped_SSTM_shift_sum_3_variance_bound_chapter_3}
    \end{align}

    \paragraph{Upper bound for $\circledFour$.} Probability event $E_{T-1}$ implies
    \begin{align}
        \circledFour & = 8\sum\limits_{k=0}^{T-1}\sum^n_{i=1}\frac{a_{k+1}^2}{n^2}\Exp_{\xi^k_i}\left[\left\|\omega_{i,k+1}^u\right\|^2\right] \notag \\
        & \leq \frac{1}{M}\cdot 8M\sum\limits_{k=0}^{T-1}\sum^n_{i=1}\frac{a_{k+1}^2}{n^2}\Exp_{\xi^k_i}\left[\left\|\omega_{i,k+1}^u\right\|^2\right] \notag\\
        & \overset{\eqref{ch3:eq:prox_clipped_SSTM_shift_sum_1_variance_bound_chapter_3}}{\leq} \frac{M}{6^3\ln(10nK/\beta)} \leq \frac{M}{9}. \label{ch3:eq:D_prox_clipped_SSTM_shift_sum_4_upper_bound_chapter_3}
    \end{align}

    \paragraph{Upper bound for $\circledFive$.} Probability event $E_{T-1}$ implies
    \begin{align}
        \circledFive & = 8\sum\limits_{k=0}^{T-1}\sum^n_{i=1}\frac{a_{k+1}^2}{n}\left\|\omega_{i,k+1}^b\right\|^2 \notag \\
        & \leq 2^{2{p} + 3}\sigma^{2{p}} \sum\limits_{k=0}^{T-1}\frac{a_{k+1}^2}{\lambda_{k}^{2{p} - 2}} \notag\\
        & \overset{\eqref{ch3:eq:Dprox_clipped_SSTM_shift_clipping_level_chapter_3}}{=} \frac{2^{2{p} + 3}\cdot 72^{2{p}-2}\sigma^{2{p}} \ln^{2{p}-2}(10nK/\beta)}{n^{2{p}-2}M^{{p}-1}} \sum\limits_{k=0}^{T-1}a_{k+1}^{2}\widetilde{a}_{k+1}^{2{p}-2} \notag\\
        & =   \frac{8\cdot 72^{2{p}-2}\sigma^{2{p}} \ln^{2{p}-2}(10nK/\beta)}{a^{2{p}} n^{2{p}-2} L^{2{p}} M^{{p}-1}} \sum\limits_{k=0}^{T-1}(k+2)^{2}\max\{K_0+2, k+2\}^{2{p}-2} \notag\\
        & \leq \frac{1}{a^{2{p}}} \cdot \frac{32\cdot72^{2{p}-2} \sigma^{2{p}}K(K+1)^{2{p}}\ln^{2{p}-2}(10nK/\beta)}{n^{2{p}-2}L^{2{p}}M^{{p}-1}} \notag \\
        & \overset{\eqref{ch3:eq:Dprox_clipped_SSTM_shift_parameter_a_chapter_3}}{\leq} \frac{M}{6}. \label{ch3:eq:D_prox_clipped_SSTM_shift_sum_5_upper_bound_chapter_3}
    \end{align}

    \paragraph{Upper bound for $\circledSix$.} This sum requires more refined analysis. We introduce a new vector: 
    \begin{equation}
        \chi_{j}^k = \begin{cases}
            \frac{a_{k+1}}{n}\sum\limits_{i=1}^{j-1} \omega_{i,k+1}^u,& \text{if } \left\|\frac{a_{k+1}}{n} \sum\limits_{i=1}^{j-1} \omega_{i,k+1}^u\right\| \leq \frac{\sqrt{M}}{2},\\
            0,& \text{otherwise,}
        \end{cases}\label{ch3:eq:prox_clipped_SSTM_chi_def_chapter_3}
    \end{equation}
    Then, by definition
    \begin{equation}
         \|\chi_j^k\| \leq \frac{\sqrt{M}}{2} \label{ch3:eq:prox_clipped_SGDA_chi_bound_chapter_3}
    \end{equation}
    and
    \begin{align}
        \circledSix & = \underbrace{8 \sum\limits_{k=0}^{T-1} \sum\limits_{j=2}^n \frac{a_{k+1}}{n}\left\langle \chi_j^k, \omega_{j,k+1}^u \right\rangle}_{\circledSix'} \notag\\
        &\quad + 8\sum\limits_{k=0}^{T-1} \sum\limits_{j=2}^n \frac{a_{k+1}}{n}\left\langle \frac{a_{k+1}}{n}\sum\limits_{i=1}^{j-1} \omega_{i,k+1}^u - \chi_j^k, \omega_{j,k+1}^u \right\rangle. \label{ch3:eq:prox_clipped_SSTM_extra_sums_six_distributed_chapter_3}
    \end{align}

    We also note here that $E_{T-1}$ implies
    \begin{align}
        \MoveEqLeft[1] 8\sum\limits_{k=0}^{T-1} \sum\limits_{j=2}^n \frac{a_{k+1}}{n}\left\langle \frac{a_{k+1}}{n}\sum\limits_{i=1}^{j-1} \omega_{i,k+1}^u - \chi_j^k, \omega_{j,k+1}^u \right\rangle \notag\\
        & = 8\sum\limits_{j=2}^n \frac{a_{T}}{n}\left\langle \frac{a_{T}}{n}\sum\limits_{i=1}^{j-1} \omega_{i,T}^u - \chi_j^{T-1}, \omega_{j,T}^u \right\rangle. \label{ch3:eq:prox_SSTM_extra_sum_six_bound_chapter_3}
    \end{align}

    \paragraph{Upper bound for $\circledSix'$.} To estimate this sum, we will use Bernstein's inequality. The summands have conditional expectations equal to zero:
    \begin{align}
\MoveEqLeft[1] \Exp_{\xi^k_j}\left[\frac{8a_{k+1}}{n}\langle\chi_j^k, \omega^u_{j,k+1} \rangle \right] \notag \\
& = \frac{8a_{k+1}}{n}\la \chi_j^k, \Exp_{\xi^k_j}\left[\omega^u_{j,k+1}\right]\ra = 0. \notag
    \end{align}

    Moreover, for all $k = 0, \ldots, T-1$ random vectors $\{\omega_{i,k+1}^u\}_{i=1}^n$ are independent. Thus, sequence $\{(8a_{k+1}/n)\langle\chi_j^k, \omega^u_{j,k+1} \rangle\}_{k,j = 0,2}^{T-1,n}$ is a martingale difference sequence. Next, the summands are bounded:
    \begin{align}
\MoveEqLeft[1] \left|\frac{8a_{k+1}}{n}\la \chi^k_j, \omega^u_{j,k+1}\ra\right| \notag \\
& \leq \frac{8a_{k+1}}{n} \|\chi^k_j\|\|\omega^u_{j,k+1}\| \overset{\eqref{ch3:eq:D_prox_clipped_SSTM_shift_norm_omega_u_bound_chapter_3},\eqref{ch3:eq:prox_clipped_SGDA_chi_bound_chapter_3}}{\leq} \frac{8a_{k+1}}{n} \cdot \frac{\sqrt{M}}{2}\cdot 2\lambda_{k} \notag \\
& \leq \frac{M}{6 \ln(10nK/\beta)} \eqdef c. \label{ch3:eq:prox_clipped_SSTM_shift_technical_6_1_chapter_3}
    \end{align}

    Finally, conditional variances $$\hat{\sigma}_{j,k}^2 \eqdef \Exp_{\xi^k_j}\left[\frac{64{a}^2_{k+1}}{n^2} \la \chi_j^k, \omega_{j,k+1}^u\ra^2\right]$$ of summands are bounded: 
    \begin{align}
\hat{\sigma}_{j,k}^2 & \leq \frac{64{a}^2_{k+1}}{n^2}\Exp_{\xi^k_j}\left[\|\chi^k_j\|^2\|\omega^u_{j,k+1}\|^2\right] \notag \\
& \leq \frac{16a_{k+1}^2 M}{n^2}\Exp_{\xi^k_{j}}\left[\|\omega_{j,k+1}^u\|^2\right]. \label{ch3:eq:prox_clipped_SSTM_shift_technical_6_2_chapter_3}
    \end{align}

    Applying Bernstein's inequality (Lemma~\ref{ch3:lem:Bernstein_ineq_chapter_3}) with $X_{j,k} =\frac{8a_{k+1}}{n}\left\la \chi^k_j, \omega^u_{j,k+1}\right\ra$, constant $c$ defined in \eqref{ch3:eq:prox_clipped_SSTM_shift_technical_6_1_chapter_3}, $b = \frac{M}{6}$, $H = \frac{M^2}{6^3\ln(10nK/\beta)}$, we get
    \begin{align*}
\MoveEqLeft[1] \Prob\left\{|\circledSix'| > \frac{M}{6} \text{ and } \sum\limits_{k=0}^{T-1}\sum\limits_{j=2}^n\hat{\sigma}_{j,k}^2 \leq \frac{M^2}{6^3\ln(10nK/\beta)}\right\} \\
& \leq 2\exp\left(- \frac{b^2}{2H + 2cb/3}\right) = \frac{\beta}{5nK}.
    \end{align*}
    The above is equivalent to 
    \begin{align}
\Prob\{E_{\circledSix'}\} & \geq 1 - \frac{\beta}{5nK},\; \text{ for }\; E_{\circledSix'} \notag \\
& = \left\{\text{either} ~~ \sum\limits_{k=0}^{T-1}\sum\limits_{j=2}^n\hat{\sigma}_{j,k}^2 > \frac{M^2}{6^3\ln(10nK/\beta)} \quad\text{ or } ~~|\circledSix'| \leq \frac{M}{6}\right\}. \label{ch3:eq:prox_clipped_SSTM_shift_technical_6_3_chapter_3}
    \end{align}
    
    Moreover, $E_{T-1}$ implies 
    \begin{align}
\sum^{T-1}_{k=0}\sum^{n}_{j=2}\hat{\sigma}^2_{j,k} & \overset{\eqref{ch3:eq:prox_clipped_SSTM_shift_technical_6_2_chapter_3}}{\leq} 16 M\sum^{T-1}_{k=0}\sum^{n}_{j=1}\frac{{a}^2_{k+1}}{n^2} \Exp_{\xi^k_{j}}\left[\|\omega_{j,k+1}^u\|^2\right] \notag \\
& \overset{\eqref{ch3:eq:prox_clipped_SSTM_shift_sum_1_variance_bound_chapter_3}}{\leq} \frac{M^2}{6^3\ln(10nK/\beta)} \label{ch3:eq:prox_clipped_SSTM_shift_sum_6_prime_variance_bound_chapter_3}
    \end{align}

    That is, we derive the upper bounds for  $ \circledOne, \circledTwo, \circledThree, \circledFour, \circledFive, \circledSix$. More precisely, $E_{T-1}$ implies

    \begin{gather*}
        B_T \overset{\eqref{ch3:eq:D_prox_clipped_SSTM_shift_technical_7_chapter_3}}{\leq} \circledOne + \circledTwo + \circledThree + \circledFour + \circledFive + \circledSix,\\
        \circledSix \overset{\eqref{ch3:eq:prox_clipped_SSTM_extra_sums_six_distributed_chapter_3}, \eqref{ch3:eq:prox_SSTM_extra_sum_six_bound_chapter_3}}{=} \circledSix' + 8 \sum\limits_{j=2}^n \frac{a_{T}}{n}\left\langle \frac{a_{T}}{n}\sum\limits_{i=1}^{j-1} \omega_{i,T}^u - \chi_j^{T-1}, \omega_{j,T}^u \right\rangle,\\
        \circledTwo \overset{\eqref{ch3:eq:D_prox_clipped_SSTM_shift_sum_2_upper_bound_chapter_3}}{\leq} \frac{M}{6},\quad \circledFour \overset{\eqref{ch3:eq:D_prox_clipped_SSTM_shift_sum_4_upper_bound_chapter_3}}{\leq} \frac{M}{9},\quad \circledFive \overset{\eqref{ch3:eq:D_prox_clipped_SSTM_shift_sum_5_upper_bound_chapter_3}}{\leq} \frac{M}{6},\\
        \sum\limits_{k=0}^{T-1}\sum^n_{i=1} \sigma_{i, k}^2 \overset{\eqref{ch3:eq:prox_clipped_SSTM_shift_sum_1_variance_bound_chapter_3}}{\leq} \frac{M^2}{6^3\ln(10nK/\beta)},\quad  \sum\limits_{k=0}^{T-1}\sum^n_{i=1} \widetilde{\sigma}_{i, k}^2 \overset{\eqref{ch3:eq:D_prox_clipped_SSTM_shift_sum_3_variance_bound_chapter_3}}{\leq} \frac{M^2}{6\cdot 9^2\ln(10nK/\beta)},\\
        \sum\limits_{k=0}^{T-1}\sum\limits_{j=2}^n\hat{\sigma}_{j,k}^2 \overset{\eqref{ch3:eq:prox_clipped_SSTM_shift_sum_6_prime_variance_bound_chapter_3}}{\leq} \frac{M^2}{6^3\ln(10nK/\beta)}.
    \end{gather*}

    In addition, we also establish (see \eqref{ch3:eq:D_prox_clipped_SSTM_shift_sum_1_upper_bound_chapter_3}, \eqref{ch3:eq:D_prox_clipped_SSTM_shift_sum_3_upper_bound_chapter_3}, \eqref{ch3:eq:prox_clipped_SSTM_shift_technical_6_3_chapter_3} and our induction assumption):
    \begin{gather*}
        \Prob\{E_{T-1}\} \geq 1 - \frac{(T-1)\beta}{K},\\
        \Prob\{E_{\circledOne}\} \geq 1 - \frac{\beta}{5nK}, \quad \Prob\{E_{\circledThree}\} \geq 1 - \frac{\beta}{5nK},\quad  \Prob\{E_{\circledSix'}\} \geq 1 - \frac{\beta}{5nK},
    \end{gather*}
    where 
    \begin{align}
        E_{\circledOne} & = \left\{ \text{either} \quad  \sum\limits_{k=0}^{T-1}\sum^n_{i=1} \sigma_{i,k}^2 > \frac{M^2}{6^3\ln(10nK/\beta)} \quad \text{or}\quad |\circledOne| \leq \frac{M}{6}\right\}, \notag\\
        E_{\circledThree} & = \left\{ \text{either} \quad  \sum\limits_{k=0}^{T-1}\sum^n_{i=1} \widetilde\sigma_{i,k}^2 > \frac{M^2}{6\cdot 9^2\ln(10nK/\beta)} \quad \text{or}\quad |\circledThree| \leq \frac{M}{9}\right\}, \notag\\
        E_{\circledSix'} & = \left\{\text{either} \quad\sum\limits_{k=0}^{T-1}\sum\limits_{j=2}^n\hat{\sigma}_{j,k}^2 > \frac{M^2}{6^3\ln(10nK/\beta)}\quad \text{ or } \quad|\circledSix'| \leq \frac{M}{6}\right\}. \notag 
    \end{align}

    Therefore, probability event $E_{T-1} \cap E_{\circledOne} \cap E_{\circledThree} \cap E_{\circledSix'}$ implies 
    \begin{align}
        B_T & \leq \frac{M}{6} + \frac{M}{6} + \frac{M}{9} + \frac{M}{9}+ \frac{M}{6}+ \frac{M}{6}  \notag\\
        &\quad+ 8\sum\limits_{k=0}^{T-1} \sum\limits_{j=2}^n \frac{a_{k+1}}{n}\left\langle \frac{a_{k+1}}{n}\sum\limits_{i=1}^{j-1} \omega_{i,k+1}^u - \chi_j^k, \omega_{j,k+1}^u \right\rangle \notag\\
        & \leq M +  8 \sum\limits_{j=2}^n \frac{a_{T}}{n}\left\langle \frac{a_{T}}{n}\sum\limits_{i=1}^{j-1} \omega_{i,T}^u - \chi_j^{T-1}, \omega_{j,T}^u \right\rangle. \label{ch3:eq:prox_clipped_SSTM_technical_8_chapter_3}
    \end{align}

    In the final part of the proof, we will show that $\frac{a_{k+1}}{n}\sum\limits_{i=1}^{j-1} \omega_{i,k+1}^u = \chi_j^k$ with high probability. In particular, we consider probability event $\widetilde{E}_{T-1,j}$ defined as follows: inequalities
    \begin{equation}
       \left\|\frac{a_{T}}{n}\sum\limits_{i=1}^{r-1} \omega_{i,T}^u\right\| \leq \frac{\sqrt{M}}{2}
    \end{equation}
    hold for $r = 2, \ldots, j$ simultaneously. We want to show that $\Prob\{E_{T-1} \cap \widetilde{E}_{T-1,j}\} \geq 1 - \frac{(T-1)\beta}{K} - \frac{2j\beta}{5nK}$ for all $j = 2, \ldots, n$. For $j = 2$ the statement is trivial since
    \begin{equation*}
         \left\|\frac{a_{T}}{n} \omega_{1,T}^u\right\| \overset{\eqref{ch3:eq:D_prox_clipped_SSTM_shift_norm_omega_u_bound_chapter_3}}{\leq} \frac{2a_{T}\lambda_{T-1}}{n} \leq \frac{\sqrt{M}}{2}.
    \end{equation*}
    
    Next, we assume that the statement holds for some $j = m-1 < n$, i.e.,  $\Prob\{E_{T-1}\cap\widetilde{E}_{T-1,m-1}\} \geq 1 - \frac{(T-1)\beta}{K} - \frac{2(m-1)\beta}{5nK}$. Our goal is to prove that $\Prob\{E_{T-1}\cap\widetilde{E}_{T-1,m}\} \geq 1 - \frac{(T-1)\beta}{K} - \frac{2m\beta}{5nK}$. First, we consider $\left\|\frac{a_{T}}{n}\sum\limits_{i=1}^{m-1} \omega_{i,T}^u\right\|$: 
    \begin{align}
        \left\|\frac{a_{T}}{n}\sum\limits_{i=1}^{m-1} \omega_{i,T}^u\right\| & = \sqrt{\frac{{a}^2_{T}}{n^2}\left\|\sum\limits_{i=1}^{m-1} \omega_{i,T}^u\right\|^2} \notag\\
        & = \sqrt{\frac{{a}^2_{T}}{n^2}\sum\limits_{i=1}^{m-1}\|\omega_{i,T}^u\|^2 + \frac{2a_{T}}{n}\sum\limits_{i=1}^{m-1}\left\langle \frac{a_{T}}{n}\sum\limits_{r=1}^{i-1}\omega_{r,T}^u , \omega_{i,T}^u \right\rangle} \notag\\
        & \leq \sqrt{\sum\limits_{k=0}^{T-1}\sum\limits_{i=1}^{m-1}\frac{a_{k+1}^2}{n^2}\|\omega_{i,k+1}^u\|^2 + \frac{2a_{T}}{n}\sum\limits_{i=1}^{m-1}\left\langle \frac{a_{T}}{n}\sum\limits_{r=1}^{i-1}\omega_{r,T}^u , \omega_{i,T}^u \right\rangle}. \notag
    \end{align}
    Next, we introduce a new notation:
    \begin{equation*}
        \rho_{i,T-1} = \begin{cases}
            \frac{a_{T}}{n}\sum\limits_{r=1}^{i-1}\omega_{r,T}^u,& \text{if } \left\|\frac{a_{T}}{n}\sum\limits_{r=1}^{i-1}\omega_{r,T}^u\right\| \leq \frac{\sqrt{M}}{2},\\
            0,& \text{otherwise}
        \end{cases}
    \end{equation*}
    for $i = 1,\ldots,m-1$. By definition, we have
    \begin{equation}
        \|\rho_{i,T-1}\| \leq \frac{\sqrt{M}}{2} \label{ch3:eq:prox_clipped_SSTM_bound_rho_chapter_3}
    \end{equation}
    for $i = 1,\ldots,m-1$. Moreover, $\widetilde{E}_{T-1,m-1}$ implies $\rho_{i,T-1} = \frac{a_{T}}{n}\sum\limits_{r=1}^{i-1}\omega_{r,T}^u$ for $i = 1,\ldots,m-1$ and
    \begin{align}
        \left\|\frac{a_{T}}{n}\sum\limits_{i=1}^{m-1} \omega_{i,T}^u\right\| & \leq \sqrt{\circledThree + \circledFour + \circledSeven}, \notag
    \end{align}
    where
    \begin{gather*}
        \circledSeven =  \frac{2a_{T}}{n}\sum\limits_{i=1}^{m-1}\left\langle \rho_{i,T-1} , \omega_{i,T}^u \right\rangle.
    \end{gather*}
    It remains to estimate $\circledSeven$.

    \paragraph{Upper bound for $\circledSeven$.} To estimate this sum, we will use Bernstein's inequality. The summands have conditional expectations equal to zero:
    \begin{align*}
\MoveEqLeft[1] \Exp_{\xi_{i}^{T-1}}\left[\frac{2a_{T}}{n} \langle \rho_{i,T-1}, \omega_{i,T}^u \rangle\right] \\
& = \frac{2a_{T}}{n} \left\langle \rho_{i,T-1}, \Exp_{\xi_{i}^{T-1}}[\omega_{i,T}^u] \right\rangle = 0,
    \end{align*}
    since random vectors $\{\omega_{i,T}^u\}_{i=1}^n$ are independent. Thus, for $i = 1, \ldots, m-1$ the summands $(2a_T/n) \langle \rho_{i,T-1}, \omega_{i,T}^u \rangle$ form a martingale difference sequence. Next, the summands are bounded:
    \begin{align}
\left|\frac{2a_{T}}{n}\langle \rho_{i,T-1}, \omega_{i,T}^u \rangle \right| & \leq \frac{2a_T}{n} \|\rho_{i,T-1}\|\cdot \|\omega_{i,T}^u\| \notag \\
& \overset{\eqref{ch3:eq:D_prox_clipped_SSTM_shift_norm_omega_u_bound_chapter_3}, \eqref{ch3:eq:prox_clipped_SSTM_bound_rho_chapter_3}}{\leq} \frac{2a_T}{n} \sqrt{M}\lambda_{T-1} \overset{\eqref{ch3:eq:Dprox_clipped_SSTM_shift_clipping_level_chapter_3}}{=} \frac{M}{36\ln(10nK/\beta)} \notag \\
& \eqdef c. \label{ch3:eq:prox_SSTM_technical_7_1_chapter_3}
    \end{align}
    Finally, conditional variances $\Bar{\sigma}_{i,T-1}^2 \eqdef \Exp_{\xi_{i}^{T-1}}\left[\frac{4a_{T}^2}{n^2} \langle \rho_{i,T-1}, \omega_{i,T}^u \rangle^2\right]$ of the summands are bounded:
    \begin{align}
\Bar{\sigma}_{i,T-1}^2 & \leq \Exp_{\xi_{i}^{T-1}}\left[\frac{4{a}^2_T }{n^2}\|\rho_{i,T-1}\|^2\cdot \|\omega_{i,T}^u\|^2\right] \notag \\
& \overset{\eqref{ch3:eq:prox_clipped_SSTM_bound_rho_chapter_3}}{\leq} \frac{{a}^2_T M}{n^2} \Exp_{\xi_{i}^{T-1}}\left[\|\omega_{i,T}^u\|^2\right]. \label{ch3:eq:prox_SSTM_technical_7_2_chapter_3}
    \end{align}

    Applying Bernstein's inequality (Lemma~\ref{ch3:lem:Bernstein_ineq_chapter_3}) with $X_{i} = \frac{2a_T}{n} \langle \rho_{i,T-1}, \omega_{i,T}^u \rangle$, constant $c$ defined in \eqref{ch3:eq:prox_SSTM_technical_7_1_chapter_3}, $b = \frac{M}{36}$, $H = \frac{M^2}{6^5\ln(10nK/\beta)}$, we get
    \begin{align*}
\MoveEqLeft[1] \Prob\left\{|\circledSeven| > \frac{M}{36} \text{ and } \sum\limits_{i=1}^{m-1}\bar{\sigma}_{i,T-1}^2 \leq \frac{M^2}{6^5\ln(10nK/\beta)}\right\} \\
& \leq 2\exp\left(- \frac{b^2}{2H + 2cb/3}\right) = \frac{\beta}{5nK}.
    \end{align*}
    The above is equivalent to 
    \begin{align}
\Prob\{E_{\circledSeven}\} & \geq 1 - \frac{\beta}{5nK},\; \text{ for }\; E_{\circledSeven} \notag \\
& = \left\{\text{either } \sum\limits_{i=1}^{m-1}\bar{\sigma}_{i,T-1}^2 > \frac{M^2}{6^5\ln(10nK/\beta)} \text{ or } |\circledSeven| \leq \frac{M}{36}\right\}. \label{ch3:eq:prox_SSTM_technical_7_3_chapter_3}
    \end{align}

    Moreover, $E_{T-1}$ implies 
    \begin{align}
        \sum^{m-1}_{i=1}\bar{\sigma}_{i,T-1}^2 & \overset{\eqref{ch3:eq:prox_SSTM_technical_7_2_chapter_3}}{\leq} \sum^{m-1}_{i=1} \frac{{a}^2_T M}{n^2} \Exp_{\xi_{i}^{T-1}}\left[\|\omega_{i,T}^u\|^2\right] \notag \\
        & \overset{\eqref{ch3:eq:D_prox_clipped_SSTM_shift_second_moment_omega_u_bound_chapter_3}}{\leq} 18 (m-1)\frac{{a}^2_T}{n^2} M \sigma^{{p}}\lambda_{T-1}^{2-{p}} \notag\\
        & \overset{\eqref{ch3:eq:Dprox_clipped_SSTM_shift_clipping_level_chapter_3}}{\leq} 18 (m-1) \frac{{a}^2_T}{n^2} \cdot \frac{n^{2-{p}}\sigma^{p} M^{2-{p}/2}}{72^{2-{p}}\widetilde{a}_T^{2-{p}}} \ln^{{p}-2}(10nK/\beta) \notag\\
        & \leq 18\cdot 6^5 \frac{a_K^{p}}{n^{{p}-1}} \left(\frac{\sigma}{\sqrt{M}}\right)^{{p}}\ln^{{p}-1}(10nK/\beta) \cdot \frac{M^2}{6^5\ln(10nK/\beta)}\notag\\
        & \leq \frac{1}{a^{{p}}} \frac{18\cdot 6^5}{2^{{p}}} \frac{1}{n^{{p}-1}}\left(\frac{\sigma}{L\sqrt{M}}\right)^{{p}} (K+1)^{{p}} \ln^{{p}-1}(10nK/\beta) \cdot \frac{M^2}{6^5\ln(10nK/\beta)}\notag\\
        & \overset{\eqref{ch3:eq:Dprox_clipped_SSTM_shift_parameter_a_chapter_3}}{\leq} \frac{M^2}{6^5\ln(10nK/\beta)}. \label{ch3:eq:prox_clipped_SSTM_shift_sum_7_variance_bound_chapter_3}
    \end{align}

    Putting all together we get that $E_{T-1}\cap \widetilde{E}_{T-1,m-1}$ implies
    \begin{gather*}
        \left\|\frac{a_{T}}{n}\sum\limits_{i=1}^{m-1} \omega_{i,T}^u\right\| \leq \sqrt{\circledThree + \circledFour + \circledSeven},\quad \circledFour \overset{\eqref{ch3:eq:D_prox_clipped_SSTM_shift_sum_4_upper_bound_chapter_3}}{\leq} \frac{M}{9}, \\
        \sum\limits_{k=0}^{T-1}\sum^n_{i=1} \widetilde{\sigma}_{i, k}^2 \overset{\eqref{ch3:eq:D_prox_clipped_SSTM_shift_sum_3_variance_bound_chapter_3}}{\leq} \frac{M^2}{6\cdot 9^2\ln(10nK/\beta)},\quad \sum^{m-1}_{i=1}\bar{\sigma}_{i,T-1}^2 \overset{\eqref{ch3:eq:prox_clipped_SSTM_shift_sum_7_variance_bound_chapter_3}}{\leq} \frac{M^2}{6^5\ln(10nK/\beta)}
    \end{gather*}
    In addition, we also establish (see \eqref{ch3:eq:D_prox_clipped_SSTM_shift_sum_3_upper_bound_chapter_3}, \eqref{ch3:eq:prox_SSTM_technical_7_3_chapter_3} and our induction assumption): 
    \begin{gather*}
        \Prob\{E_{T-1}\cap \widetilde{E}_{T-1,m-1}\} \geq 1 - \frac{(T-1)\beta}{K} - \frac{2(m-1)\beta}{5nK},\\
        \Prob\{E_{\circledThree}\} \geq 1 - \frac{\beta}{5nK}, \quad \Prob\{E_{\circledSeven}\} \geq 1 - \frac{\beta}{5nK},
    \end{gather*}
    where
    \begin{align}
        E_{\circledThree} & = \left\{\text{either } \sum\limits_{k=0}^{T-1}\sum\limits_{i=1}^n\widetilde\sigma_{i,k}^2 > \frac{M^2}{6\cdot 9^2\ln(10nK/\beta)} \text{ or } |\circledThree| \leq \frac{M}{9}\right\},\notag\\
        E_{\circledSeven} & = \left\{\text{either } \sum\limits_{i=1}^{m-1}\bar{\sigma}_{i,T-1}^2 > \frac{M^2}{6^5\ln(10nK/\beta)} \text{ or } |\circledSeven| \leq \frac{M}{36}\right\} \notag
    \end{align}
    Therefore, probability event $E_{T-1}\cap \widetilde E_{T-1, m-1} \cap E_{\circledThree} \cap E_{\circledSeven}$ implies
    \begin{equation*}
        \left\|\frac{a_{T}}{n}\sum\limits_{i=1}^{m-1} \omega_{i,T}^u\right\| \leq \sqrt{\frac{M}{9} + \frac{M}{9} + \frac{M}{36}} \leq \frac{\sqrt{M}}{2}.
    \end{equation*}
    This implies $\widetilde E_{T-1,m}$ and
    \begin{align*}
        \Prob\{E_{T-1} \cap \widetilde{E}_{T-1,m}\} & \geq \Prob\{E_{T-1}\cap \widetilde E_{T-1, m-1} \cap E_{\circledThree} \cap E_{\circledSeven}\} \\
        & = 1 - \Prob\left\{\overline{E_{T-1} \cap \widetilde{E}_{T-1,m-1}} \cup \overline{E}_{\circledThree} \cup \overline{E}_{\circledSeven}\right\}\\
        & \geq 1 - \frac{(T-1)\beta}{K} - \frac{2m\beta}{5nK}.
    \end{align*}
    Therefore, for all $m = 2,\ldots,n$ the statement holds and, in particular, $\Prob\{E_{T-1} \cap \widetilde{E}_{T-1, n}\} \geq 1 - \frac{(T-1)\beta}{K} - \frac{2\beta}{5K}$. Taking into account \eqref{ch3:eq:prox_clipped_SSTM_technical_8_chapter_3}, we conclude that $E_{T-1} \cap \widetilde{E}_{T-1, n} \cap E_{\circledOne} \cap E_{\circledThree} \cap E_{\circledSix'} \cap E_{\circledSeven}$ implies
    \begin{equation*}
        B_T \leq M
    \end{equation*}
    that is equivalent to \eqref{ch3:eq:SSTM_induction_inequality_1_chapter_3} for $t = T$. Together with \eqref{ch3:eq:D_prox_clipped_SSTM_shift_technical_1_chapter_3} for $t = T$, this gives $M_T \leq M_0 + B_T \leq 2M$, i.e., \eqref{ch3:eq:SSTM_induction_inequality_2_chapter_3} for $t = T$, and $\widetilde{E}_{T-1,n}$ gives \eqref{ch3:eq:SSTM_induction_inequality_3_chapter_3} for $t = T$. Moreover,
    \begin{align*}
        \Prob\left\{E_T\right\} & \geq \Prob\left\{E_{T-1} \cap \widetilde{E}_{T-1, n} \cap E_{\circledOne} \cap E_{\circledThree} \cap E_{\circledSix'} \cap E_{\circledSeven} \right\} \\
        & = 1 - \Prob\left\{\overline{E_{T-1} \cap \widetilde{E}_{n}} \cup \overline{E}_{\circledOne} \cup \overline{E}_{\circledThree} \cup \overline{E}_{\circledSix'} \cup \overline{E}_{\circledSeven}\right\}\\
        & \geq 1 - \frac{(T-1)\beta}{K} - \frac{2\beta}{5K} - 3\cdot \frac{\beta}{5nK} \geq 1 - \frac{T\beta}{K}.
    \end{align*}

    Finally, if $a$ equals the maximum in \eqref{ch3:eq:Dprox_clipped_SSTM_shift_parameter_a_chapter_3},
    then with probability at least $1-\beta$
    \begin{align*}
        \Phi(y^K) - \Phi(x^*) & \leq \frac{6aLM}{K(K+3)} \\
        & = \max \biggl\{\frac{12LM}{K(K+3)}, \frac{162\cdot 72^5LM\ln^4(10nK/\beta)}{nK(K+3)},\\
        &\quad \qquad\frac{3\cdot 6^7 \sigma\sqrt{M} (K+1)\ln^{(p-1)/p}(10nK/\beta)}{(K+3)(nK)^{(p-1)/p}}\biggr\}\\
        & = \cO\left(\max\left\{\frac{LM}{K^2}, \frac{LM\ln^4(nK/\beta)}{nK^2},\frac{\sigma \sqrt{M} \ln^{(p-1)/p}(nK/\beta)}{(nK)^{(p-1)/p}}\right\}\right).
    \end{align*}
    To get $\Phi(y^K) - \Phi(x^*) \leq \varepsilon$ with probability $1-\beta$, $K$ should be 
    \begin{align*}
        \cO\biggl(\max\biggl\{ & \sqrt{\frac{LM}{\varepsilon}}, \sqrt{\frac{LM}{\varepsilon n}}\ln^2\frac{nLM}{\varepsilon\beta},\\\
        & \frac{1}{n}\left(\frac{\sigma \sqrt{M}}{\varepsilon}\right)^{p/(p-1)}\ln\frac{\sigma \sqrt{M}}{\varepsilon\beta}\biggr\}\biggr)
    \end{align*}
    that concludes the proof.
\end{proof}

\subsection{Strongly Convex case}\label{ch3:appendix:SSTM_str_cvx_chapter_3}

In this section, we provide the complete formulation of our result for \algname{R-DProx-clipped-SSTM-shift} (a restarted version for \algname{DProx-clipped-SSTM-shift}) and proofs. We should mention that the results for \algname{DProx-clipped-SSTM-shift}, Theorem~\ref{ch3:thm:D_prox_clipped_SSTM_shift_full_chapter_3} and Lemma~\ref{ch3:lem:SSTM_optimization_lemma_convex_chapter_3}, can be proven in the same way for arbitrary initial shifts $h_i^0$ if
\begin{equation*}
    M \geq \|x^0-x^*\|^2 + C^2{a}^2_{K_0+1}\frac{1}{n}\sum^n_{i=1}\|h_i^0  - \nabla f_i(x^*)\|^2.
\end{equation*}

For the readers' convenience, the method’s update rule is repeated below:

\begin{algorithm}[h]
\caption{Restarted \algname{DProx-clipped-SSTM-shift} (\algname{R-DProx-clipped-SSTM-shift})}
\label{ch3:alg:R-DProx-clipped-SSTM-shift_chapter_3}   
\begin{algorithmic}[1]
\REQUIRE starting point $x^0$, number of restarts $\tau$, number of steps of \algname{DProx-clipped-SSTM-shift} between restarts $\{K_t\}_{t=1}^{\tau}$, stepsize parameters $\{a_t\}_{t=1}^\tau$, clipping levels $\{\lambda_{k}^1\}_{k=0}^{K_1-1}$, $\{\lambda_{k}^2\}_{k=0}^{K_2-1}$, \ldots, $\{\lambda_{k}^\tau\}_{k=0}^{K_\tau - 1}$, smoothness constant $L$, constants $\{N_t\}^{\tau}_{t=1}$, initial shifts $h_1^0, \ldots, h_n^0$.
\STATE $\hat{x}^0 = x^0$
\FOR{$t=1,\ldots, \tau$}
\STATE Run \algname{DProx-clipped-SSTM-shift} for $K_t$ iterations with stepsize parameter $a_{t}$, clipping levels $\{\lambda_{k}^t\}_{k=0}^{K_t-1}$, starting point $\hat{x}^{t-1}$, and initial shifts $h_i^0 = \nabla f_i(\hat{x}^{t-1})$ if $t \geq 2$. Define the output of \algname{DProx-clipped-SSTM-shift} by $\hat{x}^{t}$.
\ENDFOR
\ENSURE $\hat{x}^\tau$ 
\end{algorithmic}
\end{algorithm}

\begin{theorem}\label{ch3:thm:R_DProx_clipped_SSTM_shift_main_theorem_appendix_chapter_3}
    Let Assumptions~\ref{ch3:as:bounded_alpha_moment_chapter_3}, \ref{ch3:as:L_smoothness_chapter_3}, \ref{ch3:as:str_cvx_chapter_3} with $\mu > 0$ hold for $Q = B_{5n\sqrt{M}}(x^*)$, where
    \begin{equation*}
        M \geq \|x^0-x^*\|^2 + C^2_1(a^1_{N_1+1})^2\frac{1}{n}\sum^n_{i=1}\|h_i^0  - \nabla f_i(x^*)\|^2,
    \end{equation*}
    and \algname{R-DProx-clipped-SSTM-shift} runs \algname{DProx-clipped-SSTM-shift} $\tau$ times. For $t = 1, \ldots, \tau$, let $A_t = \ln(10nK_t\tau/\beta)$, $C_t = (864/\sqrt{n})A_t$, $N_t = \lceil 2C^2_tn\rceil$, let $a^t_{k+1} = (k+2)/(2a_tL)$ and $\widetilde{a}^t_{k+1}$ be the stepsizes of stage $t$ (defined as in Section~\ref{ch3:appendix:dprox_clipped_SSTM_shift_chapter_3} with $a = a_t$ and $K_0 = N_t$), let $\nu^t_k$ be given by \eqref{ch3:eq:nu_k_D_prox_clipped_SSTM_shift_lemma_chapter_3} with $C = C_t$ and $K_0 = N_t$, and let
    \begin{align}
        K_t & \geq \max\biggl\{N_t, \sqrt{\frac{96L}{\mu}}, 216\sqrt{2}\cdot 72^2 A_t^2\sqrt{\frac{L}{n\mu}}, \notag\\
        &\qquad\qquad \sqrt{96\cdot 864^3}\,A_t^{3/2}\sqrt{\frac{L}{\sqrt{n}\,\mu}}, \frac{A_t}{n}\left(\frac{2\cdot 6^9\sigma}{\mu\sqrt{M_{t-1}}}\right)^{p/(p-1)}\biggr\}, \label{ch3:eq:R_clipped_SSTM_K_t_chapter_3} \\
        \varepsilon_t & = \frac{\mu M_{t-1}}{4},\quad M_{t-1} = \frac{M}{2^{t-1}}, \quad \tau = \left\lceil \log_2 \frac{\mu M}{2\varepsilon} \right\rceil,\quad A_t \geq 1, \label{ch3:eq:R_clipped_SSTM_epsilon_R_t_tau_chapter_3} \\
        a_t & = \max\biggl\{2, \frac{8\cdot 3^5\cdot 72^4}{n} A_t^4, \frac{C_t(N_t+2)}{2}, \notag\\
        &\qquad\qquad \frac{ 18\cdot 6^5\sigma K_t^{1/p}(K_t+1)A_t^{(p-1)/p}}{\sqrt{M_{t-1}}Ln^{(p-1)/p}} \biggr\},\label{ch3:eq:R_Dprox_clipped_SSTM_shift_parameter_a_chapter_3}\\
        \lambda^t_{k} & = \frac{n\sqrt{M_{t-1}}}{72\widetilde{a}^t_{k+1}A_t}. \label{ch3:eq:R_Dprox_clipped_SSTM_shift_clipping_level_chapter_3}
    \end{align}
    If each $K_t$ is the smallest integer satisfying \eqref{ch3:eq:R_clipped_SSTM_K_t_chapter_3}, then to achieve $\Phi(\hat x^\tau) - \Phi(x^*) \leq \varepsilon$ with probability at least $ 1 - \beta$ \algname{R-DProx-clipped-SSTM-shift} requires
    \begin{align}
        \sum_{t=1}^\tau K_t = \cO\biggl(&\max\left\{\sqrt{\frac{L}{\mu}}, A^2\sqrt{\frac{L}{n\mu}}, A^{3/2}\sqrt{\frac{L}{\sqrt{n}\,\mu}}, A^2\right\}\ln\frac{\mu M}{\varepsilon} \notag\\
        & + \frac{A}{n}\left(\frac{\sigma^2}{\mu \varepsilon}\right)^{p/(2(p-1))}\biggr) \label{ch3:eq:R_clipped_SSTM_main_result_appendix_chapter_3}
    \end{align}
    iterations/oracle calls per worker, where $A = \max_t A_t$. Moreover, with probability $\geq 1-\beta$ the iterates of \algname{R-DProx-clipped-SSTM-shift} at stage $t$ stay in the ball $B_{2\sqrt{M_{t-1}}}(x^*)$.
\end{theorem}
\begin{proof}
    The key idea behind the proof is similar to the one used by \textcite{gorbunov2021near,sadiev2023high}. We prove by induction that for any $t = 1,\ldots, \tau$ with probability at least $1 - t\beta/\tau$ inequalities
    \begin{equation}
        \Phi(\hat{x}^l) - \Phi(x^*) \leq \varepsilon_l, \quad \hat{M}_l \leq  M_l = \frac{M}{2^l}
    \end{equation}
    hold for $l = 1,\ldots, t$ simultaneously (for $l = \tau$ only the first inequality is claimed), where, for $l = 0, 1, \ldots, \tau-1$,
    \begin{equation*}
        \hat{M}_l = \|\hat{x}^l-x^*\|^2 +C^2_{l+1}(a^{l+1}_{N_{l+1}+1})^2\frac{1}{n}\sum^n_{i=1}\|h_i^{0,l+1}-\nabla f_i(x^*)\|^2
    \end{equation*}
    is the initial value of the Lyapunov function \eqref{ch3:eq:Lyapunov_function_D_prox_clipped_SSTM_shift_chapter_3} at stage $l+1$, where $h_i^{0,l+1}$ are the initial shifts of stage $l+1$; thus $h_i^{0,1} = h_i^0$ and $h_i^{0,l+1} = \nabla f_i(\hat x^l)$ for $l \geq 1$. By assumption, $\hat{M}_0 \leq M = M_0$, and for $l \geq 1$
    \begin{equation*}
        \hat{M}_l \overset{\eqref{ch3:eq:L_smoothness_chapter_3}}{\leq} \underbrace{\left(1+C^2_{l+1}(a^{l+1}_{N_{l+1}+1})^2L^2\right)}_{ \eqdef H_{l+1}}\|\hat{x}^l-x^*\|^2.
    \end{equation*}
    Since $a^t_{N_t+1} = (N_t+2)/(2a_tL)$ and $a_t \geq C_t(N_t+2)/2$ by \eqref{ch3:eq:R_Dprox_clipped_SSTM_shift_parameter_a_chapter_3}, we have
    \begin{equation}
        \label{ch3:eq:G_t_chapter_3}
        H_t = 1 + \frac{C_t^2(N_t+2)^2}{4a_t^2} \leq 2 \quad \text{for all } t = 1, \ldots, \tau.
    \end{equation}
    Moreover, $K_t^2 \geq 48a_tL/\mu$ for all $t$. Indeed, the second and the third entries of the maximum in \eqref{ch3:eq:R_clipped_SSTM_K_t_chapter_3} give $K_t^2 \geq 48eL/\mu$ for the first two entries $e$ of the maximum in \eqref{ch3:eq:R_Dprox_clipped_SSTM_shift_parameter_a_chapter_3}; the fourth entry gives the same for $e = C_t(N_t+2)/2 \leq 2C_t^3n = 2\cdot 864^3A_t^3/\sqrt{n}$; and the last entry gives $K_t^{(p-1)/p} \geq 96\cdot 18\cdot 6^5\sigma A_t^{(p-1)/p}/(\mu\sqrt{M_{t-1}}n^{(p-1)/p})$, which together with $K_t + 1 \leq 2K_t$ gives $K_t^2 \geq 48eL/\mu$ for the last entry $e$ of \eqref{ch3:eq:R_Dprox_clipped_SSTM_shift_parameter_a_chapter_3}.

    Let $t \in \{1,\ldots,\tau\}$ and assume that $\hat{M}_{t-1} \leq M_{t-1}$ holds with probability at least $1 - (t-1)\beta/\tau$ together with the statement for $l \leq t-1$ (for $t = 1$ this is true with probability $1$). Stage $t$ is \algname{DProx-clipped-SSTM-shift} started at $\hat{x}^{t-1}$, its initial Lyapunov function equals $\hat{M}_{t-1} \leq M_{t-1}$, and its parameters satisfy the assumptions of Theorem~\ref{ch3:thm:D_prox_clipped_SSTM_shift_full_chapter_3} with $M_{t-1}$, $K_t \geq N_t$ and $\beta/\tau$ in place of $M$, $K \geq K_0$ and $\beta$. Therefore, Theorem~\ref{ch3:thm:D_prox_clipped_SSTM_shift_full_chapter_3} and the union bound imply that with probability at least $1 - t\beta/\tau$ the iterates of stage $t$ stay in $B_{2\sqrt{M_{t-1}}}(x^*)$ and
    \begin{equation}
        H_{t+1}\left(\Phi(\hat x^t) - \Phi(x^*)\right) \overset{\eqref{ch3:eq:G_t_chapter_3}}{\leq} \frac{12a_tLM_{t-1}}{K_t(K_t+3)} \leq \frac{12a_tLM_{t-1}}{K_t^2} \leq \frac{\mu M_{t-1}}{4} = \varepsilon_t, \label{ch3:eq:R_clipped_SSTM_stage_bound_chapter_3}
    \end{equation}
    where for $t = \tau$ we set $H_{\tau+1} = 1$, and in the last inequality we use $K_t^2 \geq 48a_tL/\mu$. Since $H_{t+1} \geq 1$, this gives $\Phi(\hat x^t) - \Phi(x^*) \leq \varepsilon_t$, and, for $t < \tau$, the strong convexity gives
    \begin{equation*}
        \hat{M}_{t}\leq H_{t+1}\|\hat x^{t} - x^*\|^2 \leq \frac{2H_{t+1}(\Phi(\hat x^{t}) - \Phi(x^*))}{\mu} \leq \frac{2\varepsilon_t}{\mu} = \frac{M_{t-1}}{2} = M_{t}.
    \end{equation*}
    This finishes the induction. In particular, for $t = \tau$ we get that with probability at least $1 - \beta$
    \begin{equation*}
        \Phi(\hat{x}^\tau) - \Phi(x^*) \leq \varepsilon_\tau =  \frac{\mu M_{\tau-1}}{4} = \frac{\mu M}{2^{\tau+1}}\overset{\eqref{ch3:eq:R_clipped_SSTM_epsilon_R_t_tau_chapter_3}}{\leq} \varepsilon.
    \end{equation*}
    It remains to calculate the overall number of oracle calls during all runs of \algname{DProx-clipped-SSTM-shift}. Since each $K_t$ is the smallest integer satisfying \eqref{ch3:eq:R_clipped_SSTM_K_t_chapter_3}, $N_t \leq 2\cdot 864^2A^2 + 1$, $M_{t-1} = M/2^{t-1}$ and $2^\tau \leq \mu M/\varepsilon$, we have
    \begin{align*}
        \sum\limits_{t=1}^\tau K_t
        &\leq \tau\biggl(2 + 2\cdot 864^2A^2 + \sqrt{\frac{96L}{\mu}} + 216\sqrt{2}\cdot 72^2 A^2\sqrt{\frac{L}{n\mu}}\\
        &\qquad + \sqrt{96\cdot 864^3}\,A^{3/2}\sqrt{\frac{L}{\sqrt{n}\,\mu}}\biggr)\\
        &\quad + \frac{A}{n}\left(\frac{2\cdot 6^9\sigma}{\mu\sqrt{M}}\right)^{p/(p-1)}\sum\limits_{t=1}^\tau 2^{(t-1)p/(2(p-1))}\\
        &\leq \tau\biggl(2 + 2\cdot 864^2A^2 + \sqrt{\frac{96L}{\mu}} + 216\sqrt{2}\cdot 72^2 A^2\sqrt{\frac{L}{n\mu}}\\
        &\qquad + \sqrt{96\cdot 864^3}\,A^{3/2}\sqrt{\frac{L}{\sqrt{n}\,\mu}}\biggr)\\
        &\quad + \frac{A}{n}\left(\frac{2\cdot 6^9\sigma}{\mu\sqrt{M}}\right)^{p/(p-1)}\left(\frac{\mu M}{\varepsilon}\right)^{p/(2(p-1))},
    \end{align*}
    where in the last step we use $\sum_{t=1}^\tau q^{t-1} \leq q^\tau$ for $q = 2^{p/(2(p-1))} \geq 2$. Since $\tau \leq \log_2(\mu M/\varepsilon)$ and $(\sigma/(\mu\sqrt{M}))^{p/(p-1)}(\mu M/\varepsilon)^{p/(2(p-1))} = (\sigma^2/(\mu\varepsilon))^{p/(2(p-1))}$, this gives \eqref{ch3:eq:R_clipped_SSTM_main_result_appendix_chapter_3}, which concludes the proof.
\end{proof}

\clearpage

\section{Additional Proofs for \algname{DProx-clipped-SGDA-shift}}\label{ch3:appendix:dprox_clipped_SGDA_shift_chapter_3}

\subsection{Cocoercive case}
In this section, we give the complete formulations of our results for \algname{DProx-clipped-SGDA-shift} and rigorous proofs. For the readers' convenience, the method's update rule is repeated below:
\begin{gather*}
    x^{k+1} = \thesischthreeproxkPsi \left(x^k - \gamma \thesischthreetg^k\right),\;\; \text{where}\;\; \thesischthreetg^k = \frac{1}{n}\sum\limits_{i=1}^n \thesischthreetg_i^k,\;\; \thesischthreetg_i^k = h_i^k + \hat\Delta_i^k, \\
    h_i^{k+1} = h_i^k + \nu \hat \Delta_i^k,\quad \hat\Delta_i^k = \clip\left(F_i(x^k;\xi_i^k) - h_i^k, \lambda_k\right). 
\end{gather*}

\begin{lemma}\label{ch3:lem:optimization_lemma_gap_prox_clipped_SGDA_chapter_3}
    Let Assumptions~\ref{ch3:as:monotonicity_chapter_3}, \ref{ch3:as:star_cocoercivity_chapter_3} and \ref{ch3:as:cocoercivity_chapter_3} hold for $Q = B_{3\sqrt{V}}(x^*)$, where
    \begin{gather*}
        A = \ln\frac{48n(K+1)}{\beta}, \quad V \geq \|x^0 - x^*\|^2 + \frac{25600\gamma^2A^2}{n^2}\sum_{i=1}^n\|F_i(x^*)\|^2,
    \end{gather*}
    and $0 < \gamma \leq 1/\ell$. If $x^k$ lies in $B_{3\sqrt{V}}(x^*)$ for all $k = 0,1,\ldots, K-1$ for some $K\geq 0$, then for all $u \in B_{3\sqrt{V}}(x^*)$ the iterates produced by \algname{DProx-clipped-SGDA-shift} satisfy
    \begin{align}
        \MoveEqLeft[1] \langle F(u), x^K_{\avg} - u\rangle + \Psi(x^K_{\avg}) - \Psi(u) \notag \\
        & \leq \frac{\|x^0 - u\|^2 - \|x^{K} - u\|^2}{2\gamma K} + \frac{\gamma}{K}\sum\limits_{k=0}^{K-1} \|\omega_k\|^2\notag\\
        &\quad\quad + \frac{1}{K}\sum\limits_{k=0}^{K-1}\langle x^k - u , \omega_k\rangle, \label{ch3:eq:optimization_lemma_prox_clipped_SGDA_chapter_3}\\
        x^K_{\avg} & \eqdef \frac{1}{K}\sum\limits_{k=0}^{K-1}x^{k+1}, \label{ch3:eq:x_avg_prox_clipped_SGDA_chapter_3}\\
        \omega_k & \eqdef F(x^k) - \thesischthreetg^k. \label{ch3:eq:omega_k_rpox_clipped_SGDA_chapter_3}
    \end{align}
\end{lemma}
\begin{proof}
    The proof of this lemma follows the proof of Theorem D.3 from \textcite{beznosikov2022stochastic}. For completeness, we provide here the full proof. We start with the application of Lemma~\ref{ch3:lem:prox_lemma_chapter_3} with $x^+= x^{k+1}$, $x = x^k - \gamma \thesischthreetg^k$, and $y = u$ for arbitrary $u \in B_{3\sqrt{V}}(x^*)$:
    \begin{equation*}
        \langle x^{k+1} - x^k + \gamma \thesischthreetg^k, u - x^{k+1} \rangle \geq \gamma\left(\Psi(x^{k+1}) - \Psi(u)\right).
    \end{equation*}
    Rearranging the terms, we get
    \begin{align*}
        2\gamma\left(\Psi(x^{k+1}) - \Psi(u)\right) & \leq 2\gamma\langle \thesischthreetg^k, u - x^k \rangle + 2\langle x^{k+1} - x^k, u - x^k \rangle \\
        &\quad\quad + 2\langle x^{k+1} - x^k + \gamma \thesischthreetg^k, x^k - x^{k+1} \rangle 
    \end{align*}
    implying
    \begin{align*}
        \MoveEqLeft[1] 2\gamma\left(\langle F(x^k), x^k - u \rangle + \Psi(x^{k+1}) - \Psi(u)\right) \\
        & \leq 2\langle x^{k+1} - x^k, u - x^k \rangle + 2\gamma\langle F(x^k) - \thesischthreetg^k, x^k - u\rangle \\
        &\quad\quad + 2\langle x^{k+1} - x^k + \gamma \thesischthreetg^k, x^k - x^{k+1} \rangle\\
        & = \|x^{k+1} - x^k\|^2 + \|x^k - u\|^2 - \|x^{k+1} - u\|^2\\
        &\quad\quad + 2\gamma\langle F(x^k) - \thesischthreetg^k, x^k - u\rangle \\
        &\quad\quad - 2\|x^{k+1} - x^k\|^2 + 2\gamma\langle \thesischthreetg^k, x^k - x^{k+1} \rangle\\
        & = \|x^k - u\|^2 - \|x^{k+1} - u\|^2 - \|x^{k+1} - x^k\|^2\\
        &\quad\quad + 2\gamma\langle F(x^k) - \thesischthreetg^k, x^k - u\rangle\\
        &\quad\quad + 2\gamma\langle F(u), x^k - x^{k+1} \rangle\\
        &\quad\quad + 2\gamma\langle \thesischthreetg^k - F(u), x^k - x^{k+1} \rangle\\
        & \leq \|x^k - u\|^2 - \|x^{k+1} - u\|^2\\
        &\quad\quad + 2\gamma\langle F(x^k) - \thesischthreetg^k, x^k - u\rangle\\
        &\quad\quad + 2\gamma\langle F(u), x^k - x^{k+1} \rangle + \gamma^2\|\thesischthreetg^k - F(u)\|^2,
    \end{align*}
    where in the last step we apply $2\gamma\langle \thesischthreetg^k - F(u), x^k - x^{k+1} \rangle \leq \gamma^2\|\thesischthreetg^k - F(u)\|^2 + \|x^k - x^{k+1}\|^2$. Adding $2\gamma\langle F(u), x^{k+1} - u \rangle - 2\gamma\langle F(x^k), x^k - u \rangle$ to the both sides, we derive
    \begin{equation*}
    \begin{aligned}
        &2\gamma\left(\langle F(u), x^{k+1} - u \rangle
        + \Psi(x^{k+1}) - \Psi(u)\right)\\
        &\leq \|x^k - u\|^2 - \|x^{k+1} - u\|^2\\
        &\quad + 2\gamma\langle F(u) - \thesischthreetg^k, x^k - u\rangle\\
        &\quad + \gamma^2\|\thesischthreetg^k - F(u)\|^2\\
        &= \|x^k - u\|^2 - \|x^{k+1} - u\|^2\\
        &\quad - 2\gamma\langle F(x^k) - F(u), x^k - u\rangle\\
        &\quad + \gamma^2\|\thesischthreetg^k - F(u)\|^2
        + 2\gamma\langle F(x^k) - \thesischthreetg^k, x^k - u\rangle\\
        &\overset{\eqref{ch3:eq:cocoercivity_chapter_3}}{\leq}
        \|x^k - u\|^2 - \|x^{k+1} - u\|^2\\
        &\quad - \frac{2\gamma}{\ell}\|F(x^k) - F(u)\|^2
        + 2\gamma^2\|F(x^k) - F(u)\|^2\\
        &\quad + 2\gamma\langle F(x^k) - \thesischthreetg^k, x^k - u\rangle\\
        &\quad + 2\gamma^2\|F(x^k) - \thesischthreetg^k\|^2\\
        &\leq \|x^k - u\|^2 - \|x^{k+1} - u\|^2\\
        &\quad+ 2\gamma\langle \omega_k, x^k - u\rangle
        + 2\gamma^2\|\omega_k\|^2.
    \end{aligned}
    \end{equation*}
    Next, we sum up the above inequality for $k=0,1,\ldots,K-1$ and divide both sides by $2\gamma K$:
    \begin{align*}
         \MoveEqLeft[1] \frac{1}{K}\sum\limits_{k=0}^{K-1}\left(\langle F(u), x^{k+1} - u\rangle + \Psi(x^{k+1}) - \Psi(u)\right) \\
         & \leq \frac{\|x^0 - u\|^2 - \|x^{K} - u\|^2}{2\gamma K} + \frac{\gamma}{K}\sum\limits_{k=0}^{K-1} \|\omega_k\|^2\notag\\
        &\quad\quad + \frac{1}{K}\sum\limits_{k=0}^{K-1}\langle x^k - u , \omega_k\rangle,
    \end{align*}
    To finish the proof, we need to use Jensen's inequality $\Psi\left(\frac{1}{K}\sum\limits_{k=0}^{K-1}x^{k+1}\right) \leq \frac{1}{K}\sum\limits_{k=0}^{K-1}\Psi(x^{k+1})$:
    \begin{align*}
         \MoveEqLeft[1] \langle F(u), x^K_{\avg} - u\rangle + \Psi(x^K_{\avg}) - \Psi(u) \\
         & \leq \frac{\|x^0 - u\|^2 - \|x^{K} - u\|^2}{2\gamma K} + \frac{\gamma}{K}\sum\limits_{k=0}^{K-1} \|\omega_k\|^2\notag\\
        &\quad\quad + \frac{1}{K}\sum\limits_{k=0}^{K-1}\langle x^k - u , \omega_k\rangle,
    \end{align*}
    where 
    $x^K_{\avg} = \frac{1}{K}\sum\limits_{k=0}^{K-1}x^{k+1}.$
\end{proof}

\begin{theorem}\label{ch3:thm:main_result_gap_prox_clipped_SGDA_chapter_3}
    Let Assumptions~\ref{ch3:as:monotonicity_chapter_3}, \ref{ch3:as:star_cocoercivity_chapter_3}, and \ref{ch3:as:cocoercivity_chapter_3} hold for $Q = B_{3\sqrt{V}}(x^*)$, where
    \begin{gather*}
        A = \ln\frac{48n(K+1)}{\beta}, \quad V \geq \|x^0 - x^*\|^2 + \frac{25600\gamma^2A^2}{n^2}\sum_{i=1}^n\|F_i(x^*)\|^2,
    \end{gather*}
    $\nu = 0$, $h_1^0 = \ldots = h_n^0 = 0$, and
    \begin{align}
        0< \gamma & \leq \min\left\{\frac{1}{480\ell A}, \frac{\sqrt{V} n^{(p-1)/p}}{86400^{1/p}(K+1)^{1/p}\sigma A^{(p-1)/p}}\right\}, \label{ch3:eq:gamma_prox_clipped_SGDA_chapter_3}\\
        \lambda_{k} \equiv \lambda & = \frac{n\sqrt{V}}{40\gamma A}, \label{ch3:eq:lambda_prox_clipped_SGDA_chapter_3}
    \end{align}
    for some $K \geq 0$ and $\beta \in (0,1]$. Then, after $K$ iterations the iterates produced by \algname{DProx-clipped-SGDA-shift} with probability at least $1 - \beta$ satisfy 
    \begin{equation}
        \gap_{\sqrt{V}}(x_{\avg}^{K+1}) \leq \frac{4V}{\gamma(K+1)} \quad \text{and}\quad \{x^k\}_{k=0}^{K+1} \subseteq B_{3\sqrt{V}}(x^*), \label{ch3:eq:main_result_prox_clipped_SGDA_chapter_3}
    \end{equation}
    where $x_{\avg}^{K+1}$ is defined in \eqref{ch3:eq:x_avg_prox_clipped_SGDA_chapter_3}. In particular, when $\gamma$ equals the minimum from \eqref{ch3:eq:gamma_prox_clipped_SGDA_chapter_3}, then the iterates produced by \algname{DProx-clipped-SGDA-shift} after $K$ iterations with probability at least $1-\beta$ satisfy
    \begin{equation}
        \gap_{\sqrt{V}}(x_{\avg}^{K+1}) = \cO\left(\max\left\{\frac{\ell V\ln(nK/\beta)}{K}, \frac{\sigma \sqrt{V} \ln^{(p-1)/p}(nK/\beta)}{n^{(p-1)/p}K^{(p-1)/p}}\right\}\right), \label{ch3:eq:prox_clipped_SGDA_monotone_case_2_chapter_3}
    \end{equation}
    meaning that to achieve $\gap_{\sqrt{V}}(x_{\avg}^{K+1}) \leq \varepsilon$ with probability at least $1 - \beta$ \algname{DProx-clipped-SGDA-shift} requires
    \begin{align}
        K = \cO\biggl(\max\biggl\{ & \frac{\ell V}{\varepsilon}\ln\frac{n\ell V}{\varepsilon\beta}, \notag\\
        & \frac{1}{n}\left(\frac{\sigma \sqrt{V}}{\varepsilon}\right)^{p/(p-1)}\ln\left(\frac{1}{\beta}\left(\frac{\sigma \sqrt{V}}{\varepsilon}\right)^{p/(p-1)}\right)\biggr\}\biggr) \label{ch3:eq:prox_clipped_SGDA_monotone_case_complexity_chapter_3}
    \end{align}
    iterations/oracle calls.
\end{theorem}

\begin{proof}
    The key idea behind the proof is similar to the one used by \textcite{gorbunov2022clipped, sadiev2023high}: we prove by induction that the iterates do not leave some ball and the sums decrease as $1/(K+1)$. To formulate the statement rigorously, we introduce probability event $E_k$ for each $k = 0,1,\ldots, K+1$ as follows: inequalities
    \begin{gather}
        \underbrace{\max\limits_{u \in B_{\sqrt{V}}(x^*)}\left\{ \|x^0 - u\|^2 + 2\gamma \sum\limits_{l = 0}^{t-1} \langle x^l - u, \omega_l \rangle + 2\gamma^2 \sum\limits_{l=0}^{t-1}\|\omega_l\|^2\right\}}_{A_t} \leq 8V, \label{ch3:eq:induction_inequality_1_prox_clipped_SGDA_chapter_3}\\
        \left\|\gamma\sum\limits_{l=0}^{t-1}\omega_l\right\| \leq \sqrt{V} \label{ch3:eq:induction_inequality_prox_clipped_SGDA_gap_chapter_3},\\
        \left\|\frac{\gamma}{n}\sum\limits_{i=1}^{r-1}\omega_{i,t-1}^u\right\| \leq \frac{\sqrt{V}}{2} \label{ch3:eq:induction_inequality_prox_clipped_SGDA_gap_2_chapter_3}
    \end{gather}
    hold for $t = 0,1,\ldots,k$ and $r = 1,2,\ldots, n$ simultaneously, where
    \begin{gather}
        \omega_l = \omega_l^u + \omega_l^b, \label{ch3:eq:gap_thm_prox_clipped_SGDA_technical_4_full_theta_omega_chapter_3}\\
         \omega_l^u \eqdef \frac{1}{n}\sum\limits_{i=1}^n \omega_{i,l}^u,\quad \omega_l^b \eqdef \frac{1}{n}\sum\limits_{i=1}^n \omega_{i,l}^b, \\
        \omega_{i,l}^u \eqdef \Exp_{\xi_{i}^l}\left[\thesischthreetg_{i}^l\right] - \thesischthreetg_{i}^l,\quad \omega_{i,l}^b \eqdef F_i(x^l) - \Exp_{\xi_{i}^l}\left[\thesischthreetg_{i}^l\right] \quad \forall\; i\in [n]. \label{ch3:eq:thm_prox_clipped_SGDA_technical_gap_4_chapter_3}
    \end{gather}
    We will prove by induction that $\Prob\{E_k\} \geq  1 - k\beta/(K+1)$ for all $k = 0,1,\ldots,K+1$. The base of induction follows immediately: for all $u \in B_{\sqrt{V}}(x^*)$ we have $\|x^0 - u\|^2 \leq 2\|x^0 - x^*\|^2 + 2\|x^* - u\|^2 \leq 4V < 8V$ and for $k = 0$ we have $\|\gamma\sum_{l=0}^{k-1} \omega_l\| = 0$. Next, we assume that the statement holds for $k = T-1 \leq K$, i.e., $\Prob\{E_{T-1}\} \geq  1 - (T-1)\beta/(K+1)$. Let us show that it also holds for $k = T$, i.e., $\Prob\{E_{T}\} \geq  1 - T\beta/(K+1)$. 

    To proceed, we need to show that $E_{T-1}$ implies $\|x^t - x^*\| \leq 3\sqrt{V}$ for all $t = 0,1, \ldots, T-1$. We will use the induction argument as well. The base is already proven. Next, we assume that $\|x^t - x^*\| \leq 3\sqrt{V}$ for all $t = 0,1,\ldots,t'$ for some $t' < T-1$. This means that $x^t \in B_{3\sqrt{V}}(x^*)$ for $t = 0,1,\ldots, t'$ and we can apply Lemma~\ref{ch3:lem:optimization_lemma_gap_prox_clipped_SGDA_chapter_3}: $E_{T-1}$ implies
    \begin{align*}
        & \max_{u \in B_{\sqrt{V}}(x^*)}\Big\{ 2\gamma(t'+1)\Big(\langle F(u), x^{t'+1}_{\avg} - u\rangle + \Psi(x^{t'+1}_{\avg}) - \Psi(u)\Big)\\
        &\qquad + \|x^{t'+1} - u\|^2  \Big\}\\
        &\leq \max_{u \in B_{\sqrt{V}}(x^*)}\left\{\|x^0 - u\|^2 + 2\gamma \sum\limits_{l = 0}^{t'} \langle x^l - u, \omega_l \rangle \right\}\\
        &\quad + 2\gamma^2 \sum\limits_{l=0}^{t'}\|\omega_l\|^2\\
        &\overset{\eqref{ch3:eq:induction_inequality_1_prox_clipped_SGDA_chapter_3}}{\leq} 8V.
    \end{align*}
    that gives
    \begin{align*}
        \|x^{t'+1} - x^*\|^2 & \leq \max_{u\in B_{\sqrt{V}}(x^*)}\Big\{ 2\gamma(t'+1)\Big(\langle F(u), x^{t'+1}_{\avg} - u\rangle + \Psi(x^{t'+1}_{\avg}) - \Psi(u)\Big)\\
        &\qquad + \|x^{t'+1} - u\|^2  \Big\}\\
        & \leq 8V,
    \end{align*}
    where we take $u = x^*$ and use $\langle F(x^*), x - x^*\rangle + \Psi(x) - \Psi(x^*) \geq 0$ for all $x$, since $x^*$ is a solution. Hence $\|x^{t'+1} - x^*\| \leq 2\sqrt{2V} \leq 3\sqrt{V}$. That is, we showed that $E_{T-1}$ implies $\|x^t - x^*\| \leq 3\sqrt{V}$ and
    \begin{equation}
        \max_{u \in B_{\sqrt{V}}(x^*)}\left\{ 2\gamma t\left(\langle F(u), \thesischthreetx^{t}_{\avg} - u\rangle + \Psi(\thesischthreetx^{t}_{\avg}) - \Psi(u)\right) + \|x^{t} - u\|^2  \right\} \leq 8V \label{ch3:eq:gap_thm_prox_clipped_SGDA_technical_1_5_chapter_3}
    \end{equation}
    for all $t = 1, \ldots, T-1$. Before we proceed, we introduce a new notation:
    \begin{equation*}
        u^{t} = \begin{cases}x^t - x^* ,& \text{if } \|x^t - x^*\| \leq 3\sqrt{V},\\ 0,& \text{otherwise,} \end{cases}
    \end{equation*}
    for all $t = 0, 1, \ldots, T-1$. Random vectors $\{u^{t}\}_{t=0}^{T-1}$ are bounded almost surely:
     \begin{equation}
        \|u^{t}\| \leq 3\sqrt{V}  \label{ch3:eq:thm_prox_clipped_SGDA_technical_gap_3_chapter_3}
    \end{equation}
    for all $t = 0, 1, \ldots, T-1$. In addition, $u^{t} = x^t - x^*$ follows from $E_{T-1}$ for all $t = 0, 1, \ldots, T-1$ and, thus, $E_{T-1}$ implies
   \begin{align}
        A_T & \overset{\eqref{ch3:eq:induction_inequality_1_prox_clipped_SGDA_chapter_3}}{=} \max\limits_{u \in B_{\sqrt{V}}(x^*)}\left\{ \|x^0 - u\|^2 + 2\gamma \sum\limits_{l = 0}^{T-1} \langle x^* - u, \omega_l \rangle \right\} \notag\\
        &\quad + 2\gamma \sum\limits_{l = 0}^{T-1} \langle x^l - x^*, \omega_l \rangle + 2\gamma^2 \sum\limits_{l=0}^{T-1}\|\omega_l\|^2\notag\\
        & \leq 4V + 2\gamma\max\limits_{u \in B_{\sqrt{V}}(x^*)}\left\{\left\langle x^* - u, \sum\limits_{l = 0}^{T-1}\omega_l \right\rangle \right\} \notag\\
        &\quad + 2\gamma \sum\limits_{l = 0}^{T-1} \langle u^{l}, \omega_l \rangle + 2\gamma^2 \sum\limits_{l=0}^{T-1}\|\omega_l\|^2\notag\\
        & = 4V + 2\gamma \sqrt{V} \left\|\sum\limits_{l = 0}^{T-1}\omega_l\right\| + 2\gamma \sum\limits_{l = 0}^{T-1} \langle u^{l}, \omega_l \rangle \notag\\
        &\quad + 2\gamma^2 \sum\limits_{l=0}^{T-1}\|\omega_l\|^2. \notag
    \end{align}
    Using the notation from \eqref{ch3:eq:gap_thm_prox_clipped_SGDA_technical_4_full_theta_omega_chapter_3}-\eqref{ch3:eq:thm_prox_clipped_SGDA_technical_gap_4_chapter_3}, we can rewrite $\|\omega_l\|^2$ as
    \begin{align}
        \|\omega_l\|^2 & \leq  2\|\omega_l^u\|^2 + 2\|\omega_l^b\|^2 \notag \\
        & = \frac{2}{n^2}\left\|\sum\limits_{i=1}^n \omega_{i,l}^u\right\|^2 + 2\|\omega_l^b\|^2 \notag \\
        & = \frac{2}{n^2}\sum\limits_{i=1}^n\|\omega_{i,l}^u\|^2 + \frac{4}{n^2}\sum\limits_{j=2}^n\left\langle \sum\limits_{i=1}^{j-1} \omega_{i,l}^u, \omega_{j,l}^u \right\rangle + 2\|\omega_l^b\|^2. \label{ch3:eq:omega_bound_x_t_SGDA_3_prox_clipped_chapter_3}
    \end{align}
    Putting all together, we obtain that $E_{T-1}$ implies
    \begin{align}
        A_T & \leq4V + 2\gamma \sqrt{V} \left\|\sum\limits_{l = 0}^{T-1}\omega_l\right\| + \underbrace{\frac{2\gamma}{n} \sum\limits_{l = 0}^{T-1}\sum\limits_{i=1}^n \langle u^{l}, \omega_{i,l}^u \rangle}_{\circledOne} + \underbrace{2\gamma \sum\limits_{l = 0}^{T-1} \langle u^{l}, \omega_l^b \rangle}_{\circledTwo}\notag\\
        &\quad\quad + \underbrace{\frac{4\gamma^2}{n^2} \sum\limits_{l=0}^{T-1} \sum\limits_{i=1}^n\Exp_{\xi_{i}^l}\left[\|\omega_{i,l}^u\|^2\right]}_{\circledThree} + \underbrace{\frac{4\gamma^2}{n^2} \sum\limits_{l=0}^{T-1} \sum\limits_{i=1}^n\left(\|\omega_{i,l}^u\|^2 - \Exp_{\xi_{i}^l}\left[\|\omega_{i,l}^u\|^2\right]\right)}_{\circledFour}\notag\\
        &\quad\quad + \underbrace{4\gamma^2 \sum\limits_{l=0}^{T-1}\|\omega_l^b\|^2}_{\circledFive} + \underbrace{\frac{8\gamma^2}{n^2} \sum\limits_{l=0}^{T-1} \sum\limits_{j=2}^n \left\langle \sum\limits_{i=1}^{j-1} \omega_{i,l}^u, \omega_{j,l}^u \right\rangle}_{\circledSix}.\label{ch3:eq:thm_prox_clipped_SGDA_technical_gap_5_chapter_3}
    \end{align}
    To finish the proof, it remains to estimate $2\gamma \sqrt{V} \left\|\sum_{l = 0}^{T-1}\omega_l\right\|, \circledOne, \circledTwo, \circledThree, \circledFour, \circledFive, \circledSix$ with high probability. More precisely, the goal is to prove that $2\gamma \sqrt{V} \left\|\sum_{l = 0}^{T-1}\omega_l\right\| + \circledOne + \circledTwo + \circledThree + \circledFour + \circledFive + \circledSix \leq 4V$ with high probability. Before we proceed, we need to derive several useful inequalities related to $ \omega_{i,l}^u, \omega_{l}^b$. First of all, we have
    \begin{equation}
        \|\omega_{i,l}^u\| \leq 2\lambda \label{ch3:eq:omega_magnitude_prox_clipped_SGDA_gap_chapter_3}
    \end{equation}
    by definition of the clipping operator. Next, probability event $E_{T-1}$ implies 
    \begin{align}
        \|F_i(x^l)\| & \leq \|F_i(x^l) - F_i(x^*)\| + \|F_i(x^*)\| \notag \\
        & \leq \ell\|x^l - x^*\| + \sqrt{\sum\limits_{i=1}^n\|F_i(x^*)\|^2} \notag\\
        & \leq 3\ell \sqrt{V} + \frac{n\sqrt{V}}{160\gamma A} \overset{\eqref{ch3:eq:gamma_prox_clipped_SGDA_chapter_3}}{\leq} \frac{n\sqrt{V}}{80\gamma A} \overset{\eqref{ch3:eq:lambda_prox_clipped_SGDA_chapter_3}}{=} \frac{\lambda}{2} \label{ch3:eq:operator_bound_x_t_prox_clipped_SGDA_gap_chapter_3}
    \end{align}
    for $l = 0,1,\ldots,T-1$ and $i\in [n]$, where in the second inequality we use \eqref{ch3:eq:star_cocoercivity_chapter_3} and the Cauchy--Schwarz inequality.
    Therefore, Lemma~\ref{ch3:lem:bias_and_variance_clip_chapter_3} and $E_{T-1}$ imply
    \begin{gather}
        \left\|\omega_l^b\right\| \leq \frac{1}{n}\sum\limits_{i=1}^n \|\omega_{i,l}^b\| \leq \frac{2^{p}\sigma^{p}}{\lambda^{{p}-1}}, \label{ch3:eq:bias_omega_prox_clipped_SGDA_gap_chapter_3}\\
        \Exp_{\xi_{i}^l}\left[\left\|\omega_{i,l}^u\right\|^2\right] \leq 18 \lambda^{2-{p}}\sigma^{p}, \label{ch3:eq:variance_omega_prox_clipped_SGDA_gap_chapter_3}
    \end{gather}
    for all $l = 0,1, \ldots, T-1$ and $i\in [n]$.    
    
 \paragraph{Upper bound for $\circledOne$.} To estimate this sum, we will use Bernstein's inequality. The summands have conditional expectations equal to zero:
    \begin{equation*}
        \Exp_{\xi^l_i}\left[\frac{2\gamma}{n}\langle u^{l}, \omega_{i,l}^u \rangle\right] = \frac{2\gamma}{n} \left\langle u^{l}, \Exp_{\xi_{i}^l}[\omega_{i,l}^u] \right\rangle = 0.
    \end{equation*}
    Moreover, for all $l = 0,\ldots, T-1$ random vectors $\{\omega_{i,l}^u\}_{i=1}^n$ are independent. Thus, sequence $\left\{\frac{2\gamma}{n} \langle u^{l}, \omega_{i,l}^u \rangle\right\}_{l,i = 0,1}^{T-1,n}$ is a martingale difference sequence. Next, the summands are bounded:
    \begin{align}
\left|\frac{2\gamma}{n}\langle u^{l}, \omega_{i,l}^u \rangle \right| & \leq \frac{2\gamma}{n} \|u^{l}\|\cdot \|\omega_{i,l}^u\| \overset{\eqref{ch3:eq:thm_prox_clipped_SGDA_technical_gap_3_chapter_3},\eqref{ch3:eq:omega_magnitude_prox_clipped_SGDA_gap_chapter_3}}{\leq} \frac{12 \gamma \sqrt{V} \lambda}{n} \notag \\
& \overset{\eqref{ch3:eq:lambda_prox_clipped_SGDA_chapter_3}}{\leq} \frac{3V}{10A} \eqdef c. \label{ch3:eq:SGDA_neg_mon_technical_1_1_chapter_3}
    \end{align}
    Finally, conditional variances $\sigma_{i,l}^2 \eqdef \Exp_{\xi^l_i}\left[\frac{4\gamma^2}{n^2} \langle u^{l}, \omega_{i,l}^u \rangle^2\right]$ of the summands are bounded:
    \begin{align}
\sigma_{i,l}^2 & \leq \Exp_{\xi^l_i}\left[\frac{4\gamma^2}{n^2} \|u^{l}\|^2\cdot \|\omega_{i,l}^u\|^2\right] \notag \\
& \overset{\eqref{ch3:eq:thm_prox_clipped_SGDA_technical_gap_3_chapter_3}}{\leq} \frac{36\gamma^2 V}{n^2} \Exp_{\xi^l_i}\left[\|\omega_{i,l}^u\|^2\right]. \label{ch3:eq:prox_clipped_SGDA_neg_mon_technical_1_2_chapter_3}
    \end{align}
    Applying Bernstein's inequality (Lemma~\ref{ch3:lem:Bernstein_ineq_chapter_3}) with $X_{i,l} = \frac{2\gamma}{n} \langle u^{l}, \omega_{i,l}^u \rangle$, constant $c$ defined in \eqref{ch3:eq:SGDA_neg_mon_technical_1_1_chapter_3}, $b = \frac{3V}{10}$, $H = \tfrac{3V^2}{200A}$, we get
    \begin{align*}
\MoveEqLeft[1] \Prob\left\{|\circledOne| > \frac{3V}{10} \text{ and } \sum\limits_{l=0}^{T-1}\sum\limits_{i=1}^n\sigma_{i,l}^2 \leq \frac{3V^2}{200A}\right\} \\
& \leq 2\exp\left(- \frac{b^2}{2H + 2cb/3}\right) = \frac{\beta}{24n(K+1)}.
    \end{align*}
    The above is equivalent to 
    \begin{align}
\Prob\{E_{\circledOne}\} & \geq 1 - \tfrac{\beta}{24n(K+1)}, \text{ for }E_{\circledOne} \notag \\
& = \left\{\text{either} \quad \sum\limits_{l=0}^{T-1}\sum\limits_{i=1}^n\sigma_{i,l}^2 > \frac{3V^2}{200A}\quad \text{or}\quad |\circledOne| \leq \frac{3V}{10}\right\}. \label{ch3:eq:bound_1_prox_clipped_SGDA_neg_mon_chapter_3}
    \end{align}
    Moreover, $E_{T-1}$ implies
    \begin{align}
        \sum\limits_{l=0}^{T-1}\sum\limits_{i=1}^n\sigma_{i,l}^2 & \overset{\eqref{ch3:eq:prox_clipped_SGDA_neg_mon_technical_1_2_chapter_3}}{\leq} \frac{36\gamma^2V}{n^2}\sum\limits_{l=0}^{T-1} \sum\limits_{i=1}^n\Exp_{\xi^l_i}\left[\|\omega_{i,l}^u\|^2\right] \notag\\
        & \overset{\eqref{ch3:eq:variance_omega_prox_clipped_SGDA_gap_chapter_3}}{\leq} \frac{648\gamma^2 V \sigma^{{p}} (K+1) \lambda^{2-{p}}}{n} \notag\\
        & \overset{\eqref{ch3:eq:lambda_prox_clipped_SGDA_chapter_3}}{\leq} \frac{648 \gamma^{{p}}\sqrt{V}^{4-{p}}\sigma^{{p}} (K+1)A^{{p}-2}}{40^{2-{p}}n^{{p}-1}} \notag\\
        & \overset{\eqref{ch3:eq:gamma_prox_clipped_SGDA_chapter_3}}{\leq} \frac{3V^2}{200A}. \label{ch3:eq:bound_1_variances_prox_clipped_SGDA_neg_mon_chapter_3}
    \end{align}
    
    \paragraph{Upper bound for $\circledTwo$.} Probability event $E_{T-1}$ implies
    \begin{align}
        \circledTwo & \leq 2\gamma \sum\limits_{l=0}^{T-1}\|u^{l}\| \cdot \|\omega_l^b\| \overset{\eqref{ch3:eq:thm_prox_clipped_SGDA_technical_gap_3_chapter_3},\eqref{ch3:eq:bias_omega_prox_clipped_SGDA_gap_chapter_3}}{\leq} 6\cdot 2^{{p}}\gamma \sqrt{V}(K+1) \frac{\sigma^{{p}}}{\lambda^{{p}-1}}\notag\\
        & \overset{\eqref{ch3:eq:lambda_prox_clipped_SGDA_chapter_3}}{=} \frac{6\cdot 40^{{p}-1}\cdot2^{{p}}}{n^{{p}-1}}\gamma^{{p}}\sigma^{{p}}\sqrt{V}^{2-{p}} (K+1)A^{{p}-1} \overset{\eqref{ch3:eq:gamma_prox_clipped_SGDA_chapter_3}}{\leq} \frac{3V}{100}. \label{ch3:eq:bound_2_prox_clipped_SGDA_neg_mon_chapter_3}
    \end{align}
    
    In the bounds for $\circledOne$ and $\circledTwo$ above and for $\circledThree$--$\circledFive$, $\circledSix'$ and $\circledSeven$--$\circledEight$ below we also use $T \leq K+1$.

    \paragraph{Upper bound for $\circledThree$.}  Probability event $E_{T-1}$ implies
    \begin{align}
        \circledThree & =  \frac{4\gamma^2}{n^2} \sum\limits_{l = 0}^{T-1}\sum\limits_{i=1}^n\Exp_{\xi^l_i}\left[\|\omega_{i,l}^u\|^2\right] \notag \\
        & \overset{\eqref{ch3:eq:variance_omega_prox_clipped_SGDA_gap_chapter_3}}{\leq} \frac{72\gamma^2 \lambda^{2-{p}} \sigma^{{p}} (K+1)}{n} \notag\\
        & \overset{\eqref{ch3:eq:lambda_prox_clipped_SGDA_chapter_3}}{\leq} \frac{72}{40^{2-{p}} n^{{p}-1}} \gamma^{{p}} \sqrt{V}^{2-{p}} \sigma^{{p}} (K+1)  A^{{p}-2}
        \overset{\eqref{ch3:eq:gamma_prox_clipped_SGDA_chapter_3}}{\leq} \frac{3V}{100}. \label{ch3:eq:bound_3_prox_clipped_SGDA_neg_mon_chapter_3}
    \end{align}
    
    \paragraph{Upper bound for $\circledFour$.} To estimate this sum, we will use Bernstein's inequality. The summands have conditional expectations equal to zero:
    \begin{equation*}
        \frac{4\gamma^2}{n^2}\Exp_{\xi^l_i}\left[\|\omega_{i,l}^u\|^2 - \Exp_{\xi^l_i}\left[\|\omega_{i,l}^u\|^2\right]\right] = 0.
    \end{equation*}
    Moreover, for all $l = 0,\ldots, T-1$ random vectors  $\{\omega_{i,l}^u\}_{i=1}^n$ are independent. Thus, sequence $\left\{\frac{4\gamma^2}{n^2}\left(\|\omega_{i,l}^u\|^2 -\Exp_{\xi_{i}^l}\left[\|\omega_{i,l}^u\|^2\right]\right)\right\}_{l,i = 0,1}^{T-1,n}$ is a martingale difference sequence. Next, the summands are bounded:
    \begin{align}
        \MoveEqLeft[1] \frac{4\gamma^2}{n^2} \left|\|\omega_{i,l}^u\|^2 - \Exp_{\xi^l_i}\left[\|\omega_{i,l}^u\|^2\right] \right| \notag \\
        & \leq \frac{4\gamma^2}{n^2}\left( \|\omega_{i,l}^u\|^2 + \Exp_{\xi^l_i}\left[\|\omega_{i,l}^u\|^2\right] \right) \overset{\eqref{ch3:eq:omega_magnitude_prox_clipped_SGDA_gap_chapter_3}}{\leq} \frac{32\gamma^2 \lambda^2}{n^2} \notag\\
        & \overset{\eqref{ch3:eq:lambda_prox_clipped_SGDA_chapter_3}}{\leq} \frac{V}{20A^2} \leq \frac{V}{10A} \eqdef c. \label{ch3:eq:prox_clipped_SGDA_neg_mon_technical_4_1_chapter_3}
    \end{align}
    Finally, conditional variances $$\widetilde\sigma_{i,l}^2 \eqdef \frac{16\gamma^4}{n^4}\Exp_{\xi^l_i}\left[\left(\|\omega_{i,l}^u\|^2 - \Exp_{\xi^l_i}\left[\|\omega_{i,l}^u\|^2\right]\right)^2\right]$$ of the summands are bounded:
    \begin{align}
\widetilde\sigma_{i,l}^2 & \overset{\eqref{ch3:eq:prox_clipped_SGDA_neg_mon_technical_4_1_chapter_3}}{\leq} \frac{\gamma^2 V}{5n^2A^2} \Exp_{\xi^l_i}\left[\left| \|\omega_{i,l}^u\|^2 - \Exp_{\xi^l_i}\left[\|\omega_{i,l}^u\|^2\right] \right|\right] \notag \\
& \leq \frac{2\gamma^2 V}{5n^2A^2} \Exp_{\xi^l_i}\left[\|\omega_{i,l}^u\|^2\right]. \label{ch3:eq:prox_clipped_SGDA_neg_mon_technical_4_2_chapter_3}
    \end{align}
    Applying Bernstein's inequality (Lemma~\ref{ch3:lem:Bernstein_ineq_chapter_3}) with  $X_{i,l} = (4\gamma^2/n^2)(\|\omega_{i,l}^u\|^2 - \Exp_{\xi^l_i}[\|\omega_{i,l}^u\|^2])$, constant $c$ defined in \eqref{ch3:eq:prox_clipped_SGDA_neg_mon_technical_4_1_chapter_3}, $b = \frac{V}{10}$, $H = \tfrac{V^2}{600A}$, we get
    \begin{align*}
\MoveEqLeft[1] \Prob\left\{|\circledFour| > \frac{V}{10} \text{ and } \sum\limits_{l=0}^{T-1}\sum\limits_{i=1}^n\widetilde\sigma_{i,l}^2 \leq \frac{V^2}{600A}\right\} \\
& \leq 2\exp\left(- \frac{b^2}{2H + 2cb/3}\right) = \frac{\beta}{24n(K+1)}.
    \end{align*}
    The above is equivalent to 
    \begin{align}
\Prob\{E_{\circledFour}\} & \geq 1 - \tfrac{\beta}{24n(K+1)}, \text{ for }E_{\circledFour} \notag \\
& = \left\{\text{either} \quad \sum\limits_{l=0}^{T-1}\sum\limits_{i=1}^n\widetilde\sigma_{i,l}^2 > \frac{V^2}{600A}\quad \text{or}\quad |\circledFour| \leq \frac{V}{10}\right\}. \label{ch3:eq:bound_4_prox_clipped_SGDA_neg_mon_chapter_3}
    \end{align}
    Moreover, $E_{T-1}$ implies
    \begin{align}
        \sum\limits_{l=0}^{T-1}\sum\limits_{i=1}^n\widetilde\sigma_{i,l}^2 & \overset{\eqref{ch3:eq:prox_clipped_SGDA_neg_mon_technical_4_2_chapter_3}}{\leq} \frac{2\gamma^2 V}{5n^2A^2}\sum\limits_{l=0}^{T-1}\sum\limits_{i=1}^n \Exp_{\xi^l_i}\left[\|\omega_{i,l}^u\|^2\right] \notag \\
        & \overset{\eqref{ch3:eq:variance_omega_prox_clipped_SGDA_gap_chapter_3}}{\leq} \frac{36\gamma^2 V (K+1)}{5 n A} \lambda^{2-{p}} \sigma^{{p}} \notag\\
        & \overset{\eqref{ch3:eq:lambda_prox_clipped_SGDA_chapter_3}}{\leq} \frac{9\cdot 40^{{p}}\sqrt{2}^{{p}}}{2000 n^{{p}-1}}\gamma^{{p}}
         \sqrt{V}^{4-{p}} (K+1) \sigma^{{p}}A^{{p}-4}\notag\\
        & \overset{\eqref{ch3:eq:gamma_prox_clipped_SGDA_chapter_3}}{\leq}\frac{V^2}{600A}. \label{ch3:eq:bound_4_variances_prox_clipped_SGDA_neg_mon_chapter_3}
    \end{align}
    
    \paragraph{Upper bound for $\circledFive$.} Probability event $E_{T-1}$ implies
    \begin{align}
        \circledFive & =  4\gamma^2 \sum\limits_{l = 0}^{T-1}\|\omega_l^b\|^2  \overset{\eqref{ch3:eq:bias_omega_prox_clipped_SGDA_gap_chapter_3}}{\leq} 2^{2{p}+2}\gamma^2 (K+1) \frac{\sigma^{2{p}}}{\lambda^{2{p}-2}}\notag\\
        & \overset{\eqref{ch3:eq:lambda_prox_clipped_SGDA_chapter_3}}{\leq} \frac{12800^{{p}}}{800}\gamma^{2{p}} (K+1) \frac{\sigma^{2{p}}}{n^{2{p}-2}\sqrt{V}^{2{p}-2}}A^{2{p}-2} \notag\\
        & \overset{\eqref{ch3:eq:gamma_prox_clipped_SGDA_chapter_3}}{\leq} \frac{V}{10}. \label{ch3:eq:bound_5_prox_clipped_SGDA_neg_mon_chapter_3}
    \end{align}

     \paragraph{Upper bound for $\circledSix$.} This sum requires more refined analysis. We introduce new vectors:
     \begin{equation}
        \delta_{j}^l = \begin{cases}
            \frac{\gamma}{n}\sum\limits_{i=1}^{j-1} \omega_{i,l}^u,& \text{if } \left\|\frac{\gamma}{n}\sum\limits_{i=1}^{j-1} \omega_{i,l}^u\right\| \leq \frac{\sqrt{V}}{2},\\
            0,& \text{otherwise,}
        \end{cases}\label{ch3:eq:prox_clipped_SGDA_delta_def_chapter_3}
    \end{equation}
    for all $j \in [n]$ and $l = 0,\ldots, T-1$. Then, by definition
    \begin{equation}
         \|\delta_j^l\| \leq \frac{\sqrt{V}}{2} \label{ch3:eq:prox_clipped_SGDA_delta_bound_chapter_3}
    \end{equation}
    and

\begin{align}
        \circledSix & = \underbrace{\frac{8\gamma}{n} \sum\limits_{l=0}^{T-1} \sum\limits_{j=2}^n \left\langle \delta_j^l, \omega_{j,l}^u \right\rangle}_{\circledSix'} + \frac{8\gamma}{n} \sum\limits_{l=0}^{T-1} \sum\limits_{j=2}^n \left\langle \frac{\gamma}{n}\sum\limits_{i=1}^{j-1} \omega_{i,l}^u - \delta_j^l, \omega_{j,l}^u \right\rangle. \label{ch3:eq:prox_clipped_SGDA_extra_sums_six_distributed_chapter_3}
    \end{align}

    We also note here that $E_{T-1}$ implies
    \begin{align}
        \MoveEqLeft \frac{8\gamma}{n} \sum\limits_{l=0}^{T-1} \sum\limits_{j=2}^n \left\langle \frac{\gamma}{n}\sum\limits_{i=1}^{j-1} \omega_{i,l}^u - \delta_j^l, \omega_{j,l}^u \right\rangle \notag\\
        & = \frac{8\gamma}{n}  \sum\limits_{j=2}^n \left\langle \frac{\gamma}{n}\sum\limits_{i=1}^{j-1} \omega_{i,T-1}^u - \delta_j^{T-1}, \omega_{j,T-1}^u \right\rangle. \label{ch3:eq:SGDA_extra_sums_six_distributed_1_chapter_3}
    \end{align}

    \paragraph{Upper bound for $\circledSix'$.} To estimate this sum, we will use Bernstein's inequality. The summands have conditional expectations equal to zero:

    \begin{equation*}
        \Exp_{\xi_{j}^l}\left[\frac{8\gamma}{n}\left\langle \delta_j^l, \omega_{j,l}^u \right\rangle\right] = \frac{8\gamma}{n}\left\langle \delta_j^l, \Exp_{\xi_{j}^l}[\omega_{j,l}^u] \right\rangle = 0.
    \end{equation*}

    Moreover, for all $l = 0,\ldots, T-1$ random vectors $\{\omega_{i,l}^u\}_{i=1}^n$ are independent. Thus, sequence $\left\{\frac{8\gamma}{n}\left\langle \delta_j^l, \omega_{j,l}^u \right\rangle\right\}_{l,j = 0,2}^{T-1,n}$ is a martingale difference sequence. Next, the summands are bounded:
    \begin{align}
\left|\frac{8\gamma}{n}\left\langle \delta_j^l, \omega_{j,l}^u \right\rangle \right| & \leq \frac{8\gamma}{n} \left\|\delta_j^l\right\| \cdot \|\omega_{j,l}^u\| \notag \\
& \overset{\eqref{ch3:eq:prox_clipped_SGDA_delta_bound_chapter_3},\eqref{ch3:eq:omega_magnitude_prox_clipped_SGDA_gap_chapter_3}}{\leq} \frac{8\gamma}{n}\cdot \frac{\sqrt{V}}{2}\cdot 2\lambda \leq \frac{V}{5A} \eqdef c. \label{ch3:eq:gap_thm_prox_clipped_SGDA_technical_6_5_1234_chapter_3}
    \end{align}
    Finally, conditional variances $(\sigma_{j,l}')^2 \eqdef \Exp_{\xi_{j}^l}\left[\frac{64\gamma^2}{n^2} \langle \delta_j^l, \omega_{j,l}^u \rangle^2\right]$ of the summands are bounded:
    \begin{align}
(\sigma_{j,l}')^2 & \leq \Exp_{\xi_{j}^l}\left[\frac{64\gamma^2 }{n^2}\|\delta_j^l\|^2\cdot \|\omega_{j,l}^u\|^2\right] \notag \\
& \overset{\eqref{ch3:eq:prox_clipped_SGDA_delta_bound_chapter_3}}{\leq} \frac{16\gamma^2 V}{n^2} \Exp_{\xi_{j}^l}\left[\|\omega_{j,l}^u\|^2\right]. \label{ch3:eq:gap_thm_prox_clipped_SGDA_technical_7_1234_chapter_3}
    \end{align}
    Applying Bernstein's inequality (Lemma~\ref{ch3:lem:Bernstein_ineq_chapter_3}) with $X_{j,l} = \frac{8\gamma}{n}\left\langle \delta_j^l, \omega_{j,l}^u \right\rangle$, constant $c$ defined in \eqref{ch3:eq:gap_thm_prox_clipped_SGDA_technical_6_5_1234_chapter_3}, $b = \frac{V}{5}$, $H = \frac{V^2}{150A}$, we get
    \begin{align*}
\MoveEqLeft[1] \Prob\left\{|\circledSix'| > \frac{V}{5} \text{ and } \sum\limits_{l=0}^{T-1}\sum\limits_{j=2}^n(\sigma_{j,l}')^2 \leq \frac{V^2}{150A}\right\} \\
& \leq 2\exp\left(- \frac{b^2}{2H + 2cb/3}\right) = \frac{\beta}{24n(K+1)}.
    \end{align*}
    The above is equivalent to 
    \begin{align}
\Prob\{E_{\circledSix'}\} & \geq 1 - \frac{\beta}{24n(K+1)},\; \text{ for }\; E_{\circledSix'} \notag \\
& = \left\{\text{either } \sum\limits_{l=0}^{T-1}\sum\limits_{j=2}^n(\sigma_{j,l}')^2 > \frac{V^2}{150A} \text{ or } |\circledSix'| \leq \frac{V}{5}\right\}. \label{ch3:eq:bound_7_gap_prox_clipped_SGDA_chapter_3}
    \end{align}
    Moreover, $E_{T-1}$ implies
    \begin{align}
        \sum\limits_{l=0}^{T-1}\sum\limits_{j=2}^n(\sigma_{j,l}')^2 & \overset{\eqref{ch3:eq:gap_thm_prox_clipped_SGDA_technical_7_1234_chapter_3}}{\leq} \frac{16\gamma^2 V}{n^2} \sum\limits_{l=0}^{T-1} \sum\limits_{j=2}^n \Exp_{\xi_{j}^l}\left[\|\omega_{j,l}^u\|^2\right] \notag \\
        & \overset{\eqref{ch3:eq:variance_omega_prox_clipped_SGDA_gap_chapter_3}}{\leq} \frac{288(K+1)\gamma^2 V \lambda^{2-{p}} \sigma^{p}}{n} \notag \\
        & \overset{\eqref{ch3:eq:lambda_prox_clipped_SGDA_chapter_3}}{\leq} \frac{288 (K+1)\gamma^{{p}}\sigma^{p} V^{2-\frac{{p}}{2}}}{40^{2-{p}} n^{{p}-1} A^{2-{p}}} \overset{\eqref{ch3:eq:gamma_prox_clipped_SGDA_chapter_3}}{\leq} \frac{V^2}{150A}. \label{ch3:eq:bound_6_variances_gap_prox_clipped_SGDA_chapter_3}
    \end{align}

    \paragraph{Upper bound for $2\gamma \sqrt{V} \left\|\sum_{l=0}^{T-1}\omega_t\right\|$.} We introduce new random vectors:
    \begin{equation*}
        \zeta_l = \begin{cases} \gamma \sum\limits_{r=0}^{l-1}\omega_r,& \text{if } \left\|\gamma \sum\limits_{r=0}^{l-1}\omega_r\right\| \leq \sqrt{V},\\ 0, & \text{otherwise} \end{cases}
    \end{equation*}
    for $l = 1, 2, \ldots, T-1$. With probability $1$ we have
    \begin{equation}
        \|\zeta_l\| \leq \sqrt{V}.  \label{ch3:eq:gap_thm_prox_clipped_SGDA_technical_8_chapter_3}
    \end{equation}
    Using this and \eqref{ch3:eq:induction_inequality_prox_clipped_SGDA_gap_chapter_3}, we obtain that $E_{T-1}$ implies
    \begin{align}
        \MoveEqLeft \gamma\left\|\sum\limits_{l = 0}^{T-1}\omega_l\right\| \notag\\
        & = \sqrt{\gamma^2\left\|\sum\limits_{l = 0}^{T-1}\omega_l\right\|^2}\notag\\
        & = \sqrt{\gamma^2\sum\limits_{l=0}^{T-1}\|\omega_l\|^2 + 2\gamma\sum\limits_{l=0}^{T-1}\left\langle \gamma\sum\limits_{r=0}^{l-1} \omega_r, \omega_l \right\rangle} \notag\\
        & = \sqrt{\gamma^2\sum\limits_{l=0}^{T-1}\|\omega_l\|^2 + 2\gamma \sum\limits_{l=0}^{T-1} \langle \zeta_l, \omega_l\rangle} \notag\\
        & \overset{\eqref{ch3:eq:thm_prox_clipped_SGDA_technical_gap_5_chapter_3}}{\leq} \sqrt{\frac{1}{2}\left(\circledThree + \circledFour + \circledFive + \circledSix\right) + \underbrace{\frac{2\gamma}{n} \sum\limits_{l=0}^{T-1}\sum\limits_{i=1}^n \langle \zeta_l, \omega_{i,l}^u\rangle}_{\circledSeven} + \underbrace{2\gamma \sum\limits_{l=0}^{T-1} \langle \zeta_l, \omega_l^b\rangle}_{\circledEight}}. \label{ch3:eq:norm_sum_omega_bound_gap_prox_clipped_SGDA_chapter_3}
    \end{align}
    
    \paragraph{Upper bound for $\circledSeven$.} 
    To estimate this sum, we will use Bernstein's inequality. The summands have conditional expectations equal to zero:
    
    \begin{equation*}
        \Exp_{\xi_i^l}\left[\frac{2\gamma}{n}\langle \zeta_l, \omega_{i,l}^u \rangle\right] = \frac{2\gamma}{n} \left\langle \zeta_l, \Exp_{\xi_{i}^l}[\omega_{i,l}^u] \right\rangle = 0.
    \end{equation*}

    Moreover, for all $l = 0,\ldots, T-1$ random vectors $\{\omega_{i,l}^u\}_{i=1}^n$ are independent. Thus, sequence $\left\{\frac{2\gamma}{n} \langle \zeta_l, \omega_{i,l}^u \rangle\right\}_{l,i = 0,1}^{T-1,n}$ is a martingale difference sequence. Next, the summands are bounded:
    \begin{equation}
        \left|\frac{2\gamma}{n}\langle \zeta_l, \omega_{i,l}^u \rangle \right| \leq \frac{2\gamma}{n}  \|\zeta_l\|\cdot \|\omega_{i,l}^u\| 
        \overset{\eqref{ch3:eq:gap_thm_prox_clipped_SGDA_technical_8_chapter_3},\eqref{ch3:eq:omega_magnitude_prox_clipped_SGDA_gap_chapter_3}}{\leq} \frac{4 \gamma}{n} \sqrt{V} \lambda \overset{\eqref{ch3:eq:lambda_prox_clipped_SGDA_chapter_3}}{\leq} \frac{V}{5A} \eqdef c. \label{ch3:eq:prox_clipped_SGDA_neg_mon_technical_6_1_chapter_3}
    \end{equation}
    Finally, conditional variances $\widehat\sigma_{i,l}^2 \eqdef \Exp_{\xi_i^l}\left[\frac{4\gamma^2}{n^2} \langle \zeta_l, \omega_{i,l}^u \rangle^2\right]$ of the summands are bounded:
    
    \begin{align}
\widehat\sigma_{i,l}^2 & \leq \Exp_{\xi^l_i}\left[\frac{4\gamma^2}{n^2} \|\zeta_l\|^2\cdot \|\omega_{i,l}^u\|^2\right] \notag \\
& \overset{\eqref{ch3:eq:gap_thm_prox_clipped_SGDA_technical_8_chapter_3}}{\leq} \frac{4\gamma^2}{n^2} V \Exp_{\xi^l_i}\left[\|\omega_{i,l}^u\|^2\right]. \label{ch3:eq:prox_clipped_SGDA_neg_mon_technical_6_2_chapter_3}
    \end{align}
    Applying Bernstein's inequality (Lemma~\ref{ch3:lem:Bernstein_ineq_chapter_3}) with $X_{i,l} = \frac{2\gamma}{n} \langle \zeta_l, \omega_{i,l}^u \rangle$, constant $c$ defined in \eqref{ch3:eq:prox_clipped_SGDA_neg_mon_technical_6_1_chapter_3}, $b = \frac{V}{5}$, $H = \tfrac{V^2}{150A}$, we get
    \begin{align*}
\MoveEqLeft[1] \Prob\left\{|\circledSeven| > \frac{V}{5} \text{ and } \sum\limits_{l=0}^{T-1}\sum\limits_{i=1}^n\widehat\sigma_{i,l}^2 \leq \frac{V^2}{150A}\right\} \\
& \leq 2\exp\left(- \frac{b^2}{2H + 2cb/3}\right) = \frac{\beta}{24n(K+1)}.
    \end{align*}
    The above is equivalent to 
    \begin{align}
\Prob\{E_{\circledSeven}\} & \geq 1 - \tfrac{\beta}{24n(K+1)} \text{ for }E_{\circledSeven} \notag \\
& = \left\{\text{either} \quad \sum\limits_{l=0}^{T-1}\sum\limits_{i=1}^n\widehat\sigma_{i,l}^2 > \frac{V^2}{150A}\quad \text{or}\quad |\circledSeven| \leq \frac{V}{5}\right\}. \label{ch3:eq:bound_6_prox_clipped_SGDA_neg_mon_chapter_3}
    \end{align}
    Moreover, $E_{T-1}$ implies
    \begin{align}
        \sum\limits_{l=0}^{T-1}\sum\limits_{i=1}^n\widehat\sigma_{i,l}^2 & \overset{\eqref{ch3:eq:prox_clipped_SGDA_neg_mon_technical_6_2_chapter_3}}{\leq} \frac{4\gamma^2}{n^2}V\sum\limits_{l=0}^{T-1}\sum\limits_{i=1}^n \Exp_{\xi^l_i}\left[\|\omega_{i,l}^u\|^2\right] \notag\\
        & \overset{\eqref{ch3:eq:variance_omega_prox_clipped_SGDA_gap_chapter_3}}{\leq} \frac{72\gamma^2 V \sigma^{{p}} (K+1) \lambda^{2-{p}}}{n} \notag\\
        & \overset{\eqref{ch3:eq:lambda_prox_clipped_SGDA_chapter_3}}{\leq} \frac{72 \gamma^{{p}}\sigma^{{p}} V^{2-{p}/2} (K+1)A^{{p}-2}}{40^{2-{p}}n^{{p}-1}} \notag\\
        & \overset{\eqref{ch3:eq:gamma_prox_clipped_SGDA_chapter_3}}{\leq} \frac{V^2}{150A} . \label{ch3:eq:bound_6_variances_prox_clipped_SGDA_neg_mon_chapter_3}
    \end{align}
    
    \paragraph{Upper bound for $\circledEight$.} Probability event $E_{T-1}$ implies
    \begin{align}
        \circledEight & \leq 2\gamma \sum\limits_{l=0}^{T-1}\|\zeta_l\| \cdot \|\omega_l^b\| \overset{\eqref{ch3:eq:gap_thm_prox_clipped_SGDA_technical_8_chapter_3},\eqref{ch3:eq:bias_omega_prox_clipped_SGDA_gap_chapter_3}}{\leq} 2\cdot 2^{{p}}\gamma \sqrt{V}(K+1) \frac{\sigma^{{p}}}{\lambda^{{p}-1}}\notag\\
        & \overset{\eqref{ch3:eq:lambda_prox_clipped_SGDA_chapter_3}}{=} \frac{2^{{p}+1}\cdot 40^{{p}-1}}{n^{{p}-1}}\gamma^{{p}}\sigma^{{p}}V^{1-{p}/2}(K+1)A^{{p}-1} \overset{\eqref{ch3:eq:gamma_prox_clipped_SGDA_chapter_3}}{\leq} \frac{V}{5}. \label{ch3:eq:bound_7_prox_clipped_SGDA_neg_mon_chapter_3}
    \end{align}

That is, we derive the upper bounds for  $2\gamma \sqrt{V}\left\|\sum_{l=0}^{T-1}\omega_l\right\|, \circledOne, \circledTwo, \circledThree, \circledFour, \circledFive, \circledSix$. More precisely, $E_{T-1}$ implies

    \begin{gather*}
        A_T \overset{\eqref{ch3:eq:thm_prox_clipped_SGDA_technical_gap_5_chapter_3}}{\leq} 4V + 2\gamma \sqrt{V}\left\|\sum\limits_{l=0}^{T-1}\omega_l\right\| + \circledOne + \circledTwo + \circledThree + \circledFour + \circledFive + \circledSix, \\
        \circledSix \overset{\eqref{ch3:eq:prox_clipped_SGDA_extra_sums_six_distributed_chapter_3}}{=} \circledSix' + \frac{8\gamma}{n} \sum\limits_{j=2}^n \left\langle \frac{\gamma}{n}\sum\limits_{i=1}^{j-1} \omega_{i,T-1}^u - \delta_j^{T-1}, \omega_{j,T-1}^u \right\rangle, \\
        2\gamma\sqrt{V}\left\|\sum\limits_{l = 0}^{T-1}\omega_l\right\| \overset{\eqref{ch3:eq:norm_sum_omega_bound_gap_prox_clipped_SGDA_chapter_3}}{\leq} 2\sqrt{V}\sqrt{\frac{1}{2}\left(\circledThree + \circledFour + \circledFive + \circledSix\right) + \circledSeven + \circledEight},\\
        \circledTwo \overset{\eqref{ch3:eq:bound_2_prox_clipped_SGDA_neg_mon_chapter_3}}{\leq} \frac{3V}{100},\quad \circledThree \overset{\eqref{ch3:eq:bound_3_prox_clipped_SGDA_neg_mon_chapter_3}}{\leq} \frac{3V}{100},\quad \circledFive \overset{\eqref{ch3:eq:bound_5_prox_clipped_SGDA_neg_mon_chapter_3}}{\leq} \frac{V}{10},\quad \circledEight \overset{\eqref{ch3:eq:bound_7_prox_clipped_SGDA_neg_mon_chapter_3}}{\leq} \frac{V}{5},\\
        \sum\limits_{l=0}^{T-1}\sum\limits_{i=1}^n\sigma_{i,l}^2 \overset{\eqref{ch3:eq:bound_1_variances_prox_clipped_SGDA_neg_mon_chapter_3}}{\leq}  \frac{3V^2}{200A},\quad \sum\limits_{l=0}^{T-1}\sum\limits_{i=1}^n\widetilde\sigma_{i,l}^2 \overset{\eqref{ch3:eq:bound_4_variances_prox_clipped_SGDA_neg_mon_chapter_3}}{\leq} \frac{V^2}{600A}, \\
        \sum\limits_{l=0}^{T-1}\sum\limits_{i=1}^n\widehat\sigma_{i,l}^2 \overset{\eqref{ch3:eq:bound_6_variances_prox_clipped_SGDA_neg_mon_chapter_3}}{\leq} \frac{V^2}{150A}, \quad \sum\limits_{l=0}^{T-1}\sum\limits_{j=2}^n(\sigma_{j,l}')^2  \overset{\eqref{ch3:eq:bound_6_variances_gap_prox_clipped_SGDA_chapter_3}}{\leq} \frac{V^2}{150A}.
    \end{gather*}
    In addition, we also establish (see \eqref{ch3:eq:bound_1_prox_clipped_SGDA_neg_mon_chapter_3}, \eqref{ch3:eq:bound_4_prox_clipped_SGDA_neg_mon_chapter_3}, \eqref{ch3:eq:bound_7_prox_clipped_SGDA_neg_mon_chapter_3}, \eqref{ch3:eq:bound_7_gap_prox_clipped_SGDA_chapter_3} and our induction assumption)
    \begin{gather*}
        \Prob\{E_{T-1}\} \geq 1 - \frac{(T-1)\beta}{K+1},\\
        \Prob\{E_{\circledOne}\} \geq 1 - \frac{\beta}{24n(K+1)}, \quad \Prob\{E_{\circledFour}\} \geq 1 - \frac{\beta}{24n(K+1)}, \\ \Prob\{E_{\circledSix'}\} \geq 1 - \frac{\beta}{24n(K+1)}, \quad \Prob\{E_{\circledSeven}\} \geq 1 - \frac{\beta}{24n(K+1)},
    \end{gather*}
    where
    \begin{align}
        E_{\circledOne}& = \left\{\text{either} \quad \sum\limits_{l=0}^{T-1}\sum\limits_{i=1}^n\sigma_{i,l}^2 > \frac{3V^2}{200A}\quad \text{or}\quad |\circledOne| \leq \frac{3V}{10}\right\},\notag\\
        E_{\circledFour}& = \left\{\text{either} \quad \sum\limits_{l=0}^{T-1}\sum\limits_{i=1}^n\widetilde\sigma_{i,l}^2 > \frac{V^2}{600A}\quad \text{or}\quad |\circledFour| \leq \frac{V}{10}\right\},\notag\\
        E_{\circledSix'} & = \left\{\text{either } \sum\limits_{l=0}^{T-1}\sum\limits_{j=2}^n(\sigma_{j,l}')^2 > \frac{V^2}{150A} \text{ or } |\circledSix'| \leq \frac{V}{5}\right\},\notag\\
        E_{\circledSeven}& = \left\{\text{either} \quad \sum\limits_{l=0}^{T-1}\sum\limits_{i=1}^n\widehat\sigma_{i,l}^2 > \frac{V^2}{150A}\quad \text{or}\quad |\circledSeven| \leq \frac{V}{5}\right\}.\notag
    \end{align}
    Therefore, denoting
    \begin{equation*}
        S \eqdef \frac{\gamma}{n} \sum\limits_{j=2}^n \left\langle \frac{\gamma}{n}\sum\limits_{i=1}^{j-1} \omega_{i,T-1}^u - \delta_j^{T-1}, \omega_{j,T-1}^u \right\rangle,
    \end{equation*}
    we get that probability event $E_{T-1} \cap E_{\circledOne} \cap E_{\circledFour} \cap E_{\circledSix'} \cap E_{\circledSeven}$ implies

    \begin{align}
        \left\|\gamma\sum\limits_{l = 0}^{T-1}\omega_l\right\| & \leq \sqrt{\frac{1}{2}\left(\frac{3V}{100} + \frac{V}{10} + \frac{V}{10} + \frac{V}{5}\right) + \frac{V}{5} + \frac{V}{5} + 4S} \notag\\
        & \leq \sqrt{V} + \sqrt{4S},\label{ch3:eq:gap_thm_prox_clipped_SGDA_technical_10_chapter_3}\\
        A_T & \leq 4V + 2V + 2\sqrt{V}\sqrt{4S} \notag\\
        &\quad + \frac{3V}{10} + \frac{3V}{100} + \frac{3V}{100} + \frac{V}{10} + \frac{V}{5} + \frac{V}{5} + 8S \notag\\
        & \leq 8V + 2\sqrt{V}\sqrt{4S} + 8S.\label{ch3:eq:gap_thm_prox_clipped_SGDA_technical_11_chapter_3}
    \end{align}
    
    In the final part of the proof, we will show that $\frac{\gamma}{n}\sum\limits_{i=1}^{j-1} \omega_{i,T-1}^u = \delta_j^{T-1}$ with high probability. In particular, we consider probability event $\widetilde{E}_{T-1,j}$ defined as follows: inequalities
    \begin{equation}
       \left\|\frac{\gamma}{n}\sum\limits_{i=1}^{r-1} \omega_{i,T-1}^u\right\| \leq \frac{\sqrt{V}}{2}
    \end{equation}
    hold for $r = 2, \ldots, j$ simultaneously. We want to show that $\Prob\{E_{T-1} \cap \widetilde{E}_{T-1,j}\} \geq 1 - \frac{(T-1)\beta}{K+1} - \frac{j\beta}{8n(K+1)}$ for all $j = 2, \ldots, n$. For $j = 2$ the statement is trivial since
    \begin{equation*}
         \left\|\frac{\gamma}{n} \omega_{1,T-1}^u\right\| \overset{\eqref{ch3:eq:omega_magnitude_prox_clipped_SGDA_gap_chapter_3}}{\leq} \frac{2\gamma\lambda}{n} \leq \frac{\sqrt{V}}{2}.
    \end{equation*}
    Next, we assume that the statement holds for some $j = m-1 < n$, i.e.,  $\Prob\{E_{T-1}\cap\widetilde{E}_{T-1,m-1}\} \geq 1 - \frac{(T-1)\beta}{K+1} - \frac{(m-1)\beta}{8n(K+1)}$. Our goal is to prove that $\Prob\{E_{T-1}\cap\widetilde{E}_{T-1,m}\} \geq 1 - \frac{(T-1)\beta}{K+1} - \frac{m\beta}{8n(K+1)}$. First, we consider $\left\|\frac{\gamma}{n}\sum\limits_{i=1}^{m-1} \omega_{i,T-1}^u\right\|$:
    \begin{align}
        \left\|\frac{\gamma}{n}\sum\limits_{i=1}^{m-1} \omega_{i,T-1}^u\right\| & = \sqrt{\frac{\gamma^2}{n^2}\left\|\sum\limits_{i=1}^{m-1} \omega_{i,T-1}^u\right\|^2} \notag\\
        & = \sqrt{\frac{\gamma^2}{n^2}\sum\limits_{i=1}^{m-1}\|\omega_{i,T-1}^u\|^2 + \frac{2\gamma}{n}\sum\limits_{i=1}^{m-1}\left\langle \frac{\gamma}{n}\sum\limits_{r=1}^{i-1}\omega_{r,T-1}^u , \omega_{i,T-1}^u \right\rangle} \notag\\
        & \leq \sqrt{\frac{\gamma^2}{n^2}\sum\limits_{l=0}^{T-1}\sum\limits_{i=1}^{m-1}\|\omega_{i,l}^u\|^2 + \frac{2\gamma}{n}\sum\limits_{i=1}^{m-1}\left\langle \frac{\gamma}{n}\sum\limits_{r=1}^{i-1}\omega_{r,T-1}^u , \omega_{i,T-1}^u \right\rangle}. \notag
    \end{align}
    Next, we introduce a new notation:
    \begin{equation*}
        \rho_{i,T-1} = \begin{cases}
            \frac{\gamma}{n}\sum\limits_{r=1}^{i-1}\omega_{r,T-1}^u,& \text{if } \left\|\frac{\gamma}{n}\sum\limits_{r=1}^{i-1}\omega_{r,T-1}^u\right\| \leq \frac{\sqrt{V}}{2},\\
            0,& \text{otherwise}
        \end{cases}
    \end{equation*}
    for $i = 1,\ldots,m-1$. By definition, we have
    \begin{equation}
        \|\rho_{i,T-1}\| \leq \frac{\sqrt{V}}{2} \label{ch3:eq:prox_clipped_SGDA_bound_rho_chapter_3}
    \end{equation}
    for $i = 1,\ldots,m-1$. Moreover, $\widetilde{E}_{T-1,m-1}$ implies $\rho_{i,T-1} = \frac{\gamma}{n}\sum\limits_{r=1}^{i-1}\omega_{r,T-1}^u$ for $i = 1,\ldots,m-1$ and
    \begin{align}
        \left\|\frac{\gamma}{n}\sum\limits_{i=1}^{m-1} \omega_{i,T-1}^u\right\| & \leq \sqrt{\circledThree + \circledFour + \circledNine}, \notag
    \end{align}
    where
    \begin{gather*}
        \circledNine =  \frac{2\gamma}{n}\sum\limits_{i=1}^{m-1}\left\langle \rho_{i,T-1} , \omega_{i,T-1}^u \right\rangle.
    \end{gather*}
    It remains to estimate $\circledNine$.
    
    \textbf{Upper bound for $\circledNine$.} To estimate this sum, we will use Bernstein's inequality. The summands have conditional expectations equal to zero:
    \begin{align*}
\MoveEqLeft[1] \Exp_{\xi_{i}^{T-1}}\left[\frac{2\gamma}{n} \langle \rho_{i,T-1}, \omega_{i,T-1}^u \rangle\right] \\
& = \frac{2\gamma}{n} \left\langle \rho_{i,T-1}, \Exp_{\xi_{i}^{T-1}}[\omega_{i,T-1}^u] \right\rangle = 0
    \end{align*}
    since random vectors $\{\omega_{i,T-1}^u\}_{i=1}^n$ are independent. Thus, for $i = 1, \ldots, m-1$ the summands $(2\gamma/n) \langle \rho_{i,T-1}, \omega_{i,T-1}^u \rangle$ form a martingale difference sequence. Next, the summands are bounded:
    \begin{align}
\left|\frac{2\gamma}{n}\langle \rho_{i,T-1}, \omega_{i,T-1}^u \rangle \right| & \leq \frac{2\gamma}{n} \|\rho_{i,T-1}\|\cdot \|\omega_{i,T-1}^u\| \notag \\
& \overset{\eqref{ch3:eq:prox_clipped_SGDA_bound_rho_chapter_3},\eqref{ch3:eq:omega_magnitude_prox_clipped_SGDA_gap_chapter_3}}{\leq} \frac{2\gamma}{n} \sqrt{V}\lambda \overset{\eqref{ch3:eq:lambda_prox_clipped_SGDA_chapter_3}}{=} \frac{V}{20A} \eqdef c. \label{ch3:eq:gap_thm_prox_clipped_SGDA_technical_6_5_10_chapter_3}
    \end{align}
    Finally, conditional variances $(\widehat\sigma_{i,T-1}')^2 \eqdef \Exp_{\xi_{i}^{T-1}}\left[\frac{4\gamma^2}{n^2} \langle \rho_{i,T-1}, \omega_{i,T-1}^u \rangle^2\right]$ of the summands are bounded:
    \begin{align}
(\widehat\sigma_{i,T-1}')^2 & \leq \Exp_{\xi_{i}^{T-1}}\left[\frac{4\gamma^2 }{n^2}\|\rho_{i,T-1}\|^2\cdot \|\omega_{i,T-1}^u\|^2\right] \notag \\
& \overset{\eqref{ch3:eq:prox_clipped_SGDA_bound_rho_chapter_3}}{\leq} \frac{\gamma^2 V}{n^2} \Exp_{\xi_{i}^{T-1}}\left[\|\omega_{i,T-1}^u\|^2\right]. \label{ch3:eq:gap_thm_prox_clipped_SGDA_technical_7_10_chapter_3}
    \end{align}
    Applying Bernstein's inequality (Lemma~\ref{ch3:lem:Bernstein_ineq_chapter_3}) with $X_{i} = \frac{2\gamma}{n} \langle \rho_{i,T-1}, \omega_{i,T-1}^u \rangle$, constant $c$ defined in \eqref{ch3:eq:gap_thm_prox_clipped_SGDA_technical_6_5_10_chapter_3}, $b = \frac{V}{20}$, $H = \frac{V^2}{2400A}$, we get
    \begin{align*}
\MoveEqLeft[1] \Prob\left\{|\circledNine| > \frac{V}{20} \text{ and } \sum\limits_{i=1}^{m-1}(\widehat\sigma_{i,T-1}')^2 \leq \frac{V^2}{2400A}\right\} \\
& \leq 2\exp\left(- \frac{b^2}{2H + 2cb/3}\right) = \frac{\beta}{24n(K+1)}.
    \end{align*}
    The above is equivalent to 
    \begin{align}
\Prob\{E_{\circledNine}\} & \geq 1 - \frac{\beta}{24n(K+1)},\; \text{ for }\; E_{\circledNine} \notag \\
& = \left\{\text{either } \sum\limits_{i=1}^{m-1}(\widehat\sigma_{i,T-1}')^2 > \frac{V^2}{2400A} \text{ or } |\circledNine| \leq \frac{V}{20}\right\}. \label{ch3:eq:bound_10_gap_prox_clipped_SGDA_chapter_3}
    \end{align}
    Moreover, $E_{T-1}$ implies
    \begin{align}
        \sum\limits_{i=1}^{m-1}(\widehat\sigma_{i,T-1}')^2 & \overset{\eqref{ch3:eq:gap_thm_prox_clipped_SGDA_technical_7_10_chapter_3}}{\leq} \frac{\gamma^2 V}{n^2} \sum\limits_{i=1}^n \Exp_{\xi_{i}^{T-1}}\left[\|\omega_{i,T-1}^u\|^2\right] \notag \\
        & \overset{\eqref{ch3:eq:variance_omega_prox_clipped_SGDA_gap_chapter_3}}{\leq} \frac{18\gamma^2 V \lambda^{2-{p}} \sigma^{p}}{n} \notag \\
        & \overset{\eqref{ch3:eq:lambda_prox_clipped_SGDA_chapter_3}}{\leq} \frac{18 \gamma^{{p}}\sigma^{p} V^{2-\frac{{p}}{2}}}{40^{2-{p}} n^{{p}-1} A^{2-{p}}} \overset{\eqref{ch3:eq:gamma_prox_clipped_SGDA_chapter_3}}{\leq} \frac{V^2}{2400A}. \label{ch3:eq:bound_10_variances_gap_prox_clipped_SGDA_chapter_3}
    \end{align}
    Putting all together we get that $E_{T-1}\cap \widetilde{E}_{T-1, m-1}$ implies
    \begin{gather*}
        \left\|\frac{\gamma}{n}\sum\limits_{i=1}^{m-1} \omega_{i,T-1}^u\right\| \leq \sqrt{\circledThree + \circledFour + \circledNine},\quad \circledThree \overset{\eqref{ch3:eq:bound_3_prox_clipped_SGDA_neg_mon_chapter_3}}{\leq} \frac{V}{10},\\
        \sum\limits_{l=0}^{T-1}\sum\limits_{i=1}^n\widetilde\sigma_{i,l}^2 \overset{\eqref{ch3:eq:bound_4_variances_prox_clipped_SGDA_neg_mon_chapter_3}}{\leq} \frac{V^2}{600A},\quad \sum\limits_{i=1}^{m-1}(\widehat\sigma_{i,T-1}')^2 \leq \frac{V^2}{2400A}
    \end{gather*}
    In addition, we also establish (see \eqref{ch3:eq:bound_4_prox_clipped_SGDA_neg_mon_chapter_3}, \eqref{ch3:eq:bound_10_gap_prox_clipped_SGDA_chapter_3} and our induction assumption)
    \begin{gather*}
        \Prob\{E_{T-1}\cap \widetilde{E}_{T-1,m-1}\} \geq 1 - \frac{(T-1)\beta}{K+1} - \frac{(m-1)\beta}{8n(K+1)},\\
        \Prob\{E_{\circledFour}\} \geq 1 - \frac{\beta}{24n(K+1)}, \quad \Prob\{E_{\circledNine}\} \geq 1 - \frac{\beta}{24n(K+1)},
    \end{gather*}
    where
    \begin{align}
        E_{\circledFour} & = \left\{\text{either } \sum\limits_{l=0}^{T-1}\sum\limits_{i=1}^n\widetilde\sigma_{i,l}^2 > \frac{V^2}{600A} \text{ or } |\circledFour| \leq \frac{V}{10}\right\},\notag\\
        E_{\circledNine} & = \left\{\text{either } \sum\limits_{i=1}^{m-1}(\widehat\sigma_{i,T-1}')^2 > \frac{V^2}{2400A} \text{ or } |\circledNine| \leq \frac{V}{20}\right\}. \notag
    \end{align}
    Therefore, probability event $E_{T-1} \cap \widetilde{E}_{T-1, m-1} \cap E_{\circledFour} \cap E_{\circledNine}$ implies
    \begin{equation*}
        \left\|\frac{\gamma}{n}\sum\limits_{i=1}^{m-1} \omega_{i,T-1}^u\right\| \leq \sqrt{\frac{V}{10} + \frac{V}{10} + \frac{V}{20}} = \frac{\sqrt{V}}{2}.
    \end{equation*} 
    This implies $\widetilde E_{T-1,m}$ and
    \begin{align*}
        \Prob\{E_{T-1} \cap \widetilde{E}_{T-1,m}\} & \geq \Prob\{E_{T-1} \cap \widetilde{E}_{T-1,m-1} \cap E_{\circledFour} \cap E_{\circledNine} \} \\
        & = 1 - \Prob\left\{\overline{E_{T-1} \cap \widetilde{E}_{T-1,m-1}} \cup \overline{E}_{\circledFour} \cup \overline{E}_{\circledNine} \right\}\\
        & \geq 1 - \frac{(T-1)\beta}{K+1} - \frac{m\beta}{8n(K+1)}.
    \end{align*}
    Therefore, for all $m = 2,\ldots,n$ the statement holds and, in particular, $\Prob\{E_{T-1} \cap \widetilde{E}_{T-1,n}\} \geq 1 - \frac{(T-1)\beta}{K+1} - \frac{\beta}{8(K+1)}$. Taking into account  \eqref{ch3:eq:gap_thm_prox_clipped_SGDA_technical_10_chapter_3} and \eqref{ch3:eq:gap_thm_prox_clipped_SGDA_technical_11_chapter_3}, we conclude that $E_{T-1} \cap \widetilde{E}_{T-1,n}$ implies (since $S = 0$ on $\widetilde{E}_{T-1,n}$)  
    \begin{equation*}
       \left\|\gamma\sum\limits_{l = 0}^{T-1}\omega_l\right\| \leq \sqrt{V},\quad  A_T \leq 8V,
    \end{equation*}
    which is equivalent to \eqref{ch3:eq:induction_inequality_1_prox_clipped_SGDA_chapter_3} and \eqref{ch3:eq:induction_inequality_prox_clipped_SGDA_gap_chapter_3} for $t = T$; moreover, $\widetilde{E}_{T-1,n}$ gives \eqref{ch3:eq:induction_inequality_prox_clipped_SGDA_gap_2_chapter_3} for $t = T$. Finally,
    \begin{align*}
        \Prob\left\{E_T\right\} & \geq \Prob\left\{E_{T-1} \cap \widetilde{E}_{T-1,n} \cap E_{\circledOne} \cap E_{\circledFour} \cap E_{\circledSix'} \cap E_{\circledSeven} \cap E_{\circledNine}\right\} \\
        & = 1 - \Prob\left\{\overline{E_{T-1} \cap \widetilde{E}_{T-1,n}} \cup \overline{E}_{\circledOne} \cup \overline{E}_{\circledFour} \cup \overline{E}_{\circledSix'} \cup \overline{E}_{\circledSeven} \cup \overline{E}_{\circledNine} \right\}\\
        & \geq 1 - \frac{(T-1)\beta}{K+1} - \frac{\beta}{8(K+1)} - 4\cdot \frac{\beta}{24n(K+1)} \geq 1 - \frac{T\beta}{K+1}.
    \end{align*}

    In other words, we showed that $\Prob\{E_k\} \geq 1 - k\beta/(K+1)$ for all $k = 0,1,\ldots,K+1$. For $k = K+1$ we have that with probability at least $1 - \beta$

    \begin{align*}
        \gap_{\sqrt{V}}(x_{\avg}^{K+1})
        & \overset{\eqref{ch3:eq:gap_thm_prox_clipped_SGDA_technical_1_5_chapter_3}}{\leq} \frac{4V}{\gamma(K+1)}.
    \end{align*}

    Finally, if $\gamma$ equals the minimum in \eqref{ch3:eq:gamma_prox_clipped_SGDA_chapter_3}, then with probability at least $1-\beta$
    \begin{align*}
        \gap_{\sqrt{V}}(x^{K+1}_{\avg}) & \leq \frac{4V}{\gamma(K+1)} \\
        & = \max\left\{\frac{1920 \ell V A}{K+1}, \frac{4\cdot 86400^{1/p} \sigma \sqrt{V} A^{(p-1)/p}}{n^{(p-1)/p}(K+1)^{(p-1)/p}}\right\}\\
        &= \cO\left(\max\left\{\frac{\ell V\ln(nK/\beta)}{K}, \frac{\sigma \sqrt{V} \ln^{(p-1)/p}(nK/\beta)}{n^{(p-1)/p}K^{(p-1)/p}}\right\}\right).
    \end{align*}
    To get $\gap_{\sqrt{V}}(x^{K+1}_{\avg}) \leq \varepsilon$ with probability at least $1-\beta$ it is sufficient to choose $K$ such that both terms in the maximum above are $\cO(\varepsilon)$. This leads to
    \begin{align*}
         K = \cO\biggl(\max\biggl\{ & \frac{\ell V}{\varepsilon}\ln\frac{n\ell V}{\varepsilon\beta}, \\
         & \frac{1}{n}\left(\frac{\sigma \sqrt{V}}{\varepsilon}\right)^{p/(p-1)}\ln\left(\frac{1}{\beta}\left(\frac{\sigma \sqrt{V}}{\varepsilon}\right)^{p/(p-1)}\right)\biggr\}\biggr),
    \end{align*}
    that concludes the proof.
\end{proof}

\subsection{Quasi-Strongly Monotone case}

\begin{lemma}\label{ch3:lem:optimization_lemma_str_mon_SGDA_2_chapter_3}
    Let Assumptions~\ref{ch3:as:str_monotonicity_chapter_3}, \ref{ch3:as:star_cocoercivity_chapter_3} hold for $Q = B_{\sqrt{2V}}(x^*)$, where
    \begin{gather*}
        A = \ln\frac{48n(K+1)}{\beta},\quad V \geq \|x^0 - x^*\|^2 + \frac{9\cdot 10^6\gamma^2A^2}{n^2}\sum_{i=1}^n\|F_i(x^*)\|^2,
    \end{gather*}
    $0 < \gamma \leq 1/(\ell + 18\cdot 10^6 \nu A^2\ell/n)$ and $\gamma\mu \leq \nu \leq 1/(18\cdot 10^6A^2)$. If $x^k$ lies in $B_{\sqrt{2V}}(x^*)$ for all $k = 0,1,\ldots, K$ for some $K\geq 0$, then the iterates produced by \algname{DProx-clipped-SGDA-shift} satisfy
    \begin{align}
        V_{K+1} & \leq (1 - \gamma \mu)^{K+1}V_0 + \frac{2\gamma}{n} \sum\limits_{k=0}^K\sum\limits_{i=1}^n (1-\gamma\mu)^{K-k}\langle x^k - x^* - \gamma (F(x^k)-h^*), \omega_{i,k} \rangle\notag\\
        &\quad\quad  + \frac{ \gamma^2  }{n^2}\sum\limits_{k=0}^K\sum\limits_{i=1}^n (1-\gamma \mu)^{K-k} \|\omega_{i,k}\|^2  + \gamma^2 \sum\limits_{k=0}^K (1-\gamma\mu)^{K-k} \|\omega_k\|^2, \label{ch3:eq:optimization_lemma_prox_SGDA_str_mon_2_chapter_3}
    \end{align}
    where $V_k = \|x^k - x^*\|^2 + (9\cdot 10^6\gamma^2A^2/n^2)\sum_{i=1}^n\|h^k_i - h^*_i\|^2$, $h_i^* = F_i(x^*)$, and $\omega_k, \omega_k^u, \omega_k^b, \omega_{k,i}^u, \omega_{k,i}^b$ are defined in \eqref{ch3:eq:gap_thm_prox_clipped_SGDA_technical_4_full_theta_omega_chapter_3}-\eqref{ch3:eq:thm_prox_clipped_SGDA_technical_gap_4_chapter_3}.
\end{lemma}

\begin{proof}
     Using the update rule of \algname{DProx-clipped-SGDA-shift} and
    $\omega_k = F(x^k) - \thesischthreetg^k$ we obtain
     \begin{align*}
        \|x^{k+1} - x^*\|^2 & = \|\thesischthreeproxkPsi\left(x^k - \gamma\thesischthreetg^k\right) - \thesischthreeproxkPsi\left(x^* - \gamma h^*\right)\|^2\\
        & \le \|x^k - x^* - \gamma(\thesischthreetg^k - h^*)\|^2\\
        & = \|x^k - x^*\|^2 -2\gamma \langle x^k - x^*, \thesischthreetg^k - h^* \rangle + \gamma^2\|\thesischthreetg^k - h^*\|^2\\
        & = \|x^k - x^*\|^2 - 2\gamma \langle x^k - x^*, F(x^k)-h^* \rangle - 2\gamma^2\langle F(x^k) - h^*, \omega_k\rangle\\
        &\quad\quad  + 2\gamma \langle x^k - x^*, \omega_k \rangle  + \gamma^2\|F(x^k) - h^*\|^2 + \gamma^2\|\omega_k\|^2,
    \end{align*}
    Next, let us recall that
    \begin{gather*}
    h_i^{k+1} = h_i^k + \nu \hat \Delta_i^k,\quad \hat\Delta_i^k = \clip\left(F_i(x^k;\xi_i^k) - h_i^k, \lambda_k\right),\\
    \thesischthreetg_i^k = h_i^k + \hat\Delta_i^k,\quad \omega_{i,k} = F_i(x^k) - \thesischthreetg_i^k.
\end{gather*}
    Then, $\forall i \in [n]$ we have
    \begin{align*}
        \|h^{k+1}_i - h^*_i\|^2 & = \|h^k_i - h^*_i + \nu \hat\Delta^k_i\|^2 \\
        & = \|h^k_i - h^*_i\|^2 + 2\nu \langle h^k_i - h^*_i, \hat\Delta^k_i\rangle + \nu^2\|\hat\Delta^k_i\|^2\\
        & = \|h^k_i - h^*_i\|^2 + 2\nu \langle h^k_i - h^*_i, \thesischthreetg^k_i - h_i^k\rangle + \nu^2\|\thesischthreetg^k_i - h_i^k\|^2\\
        & \leq \|h^k_i - h^*_i\|^2 + 2\nu \langle h^k_i - h^*_i, \thesischthreetg^k_i - h_i^k\rangle + \nu\|\thesischthreetg^k_i - h_i^k\|^2\\
        & = \|h^k_i - h^*_i\|^2 + \nu \langle \thesischthreetg^k_i - h_i^k, \thesischthreetg^k_i + h^k_i - 2h^*_i \rangle \\
        & = (1-\nu)\|h^k_i - h^*_i\|^2 + \nu\|\thesischthreetg^k_i - h_i^*\|^2\\
        & \le  (1-\nu)\|h^k_i - h^*_i\|^2 + 2\nu\|\thesischthreetg^k_i - F_i(x^k)\|^2 + 2\nu\|F_i(x^k) - h_i^*\|^2\\
        & = (1-\nu)\|h^k_i - h^*_i\|^2 + 2\nu\|\omega_{i,k}\|^2 + 2\nu \|F_i(x^k) - h^*_i\|^2,
    \end{align*}
    where in the inequality we use $\nu \leq 1$.
    Let us consider the following stepsize condition
    \begin{equation}
        0 < \gamma \leq \frac{1}{\ell + \frac{18\cdot 10^6 \nu A^2\ell}{n}}.\label{ch3:eq:gamma_lemma_prox_SGDA_str_mon_2_chapter_3}
    \end{equation}
    Lyapunov function
    \begin{equation*}
        V_k = \|x^k - x^*\|^2 + \frac{9\cdot 10^6\gamma^2A^2}{n^2}\sum\limits_{i=1}^n\|h^k_i - h^*_i\|^2.
    \end{equation*}
    \begin{align*}
        V_{k+1}
        &\le \|x^k - x^*\|^2
        - 2\gamma \langle x^k - x^*, F(x^k)-h^* \rangle\\
        &\quad - 2\gamma^2\langle F(x^k) - h^*, \omega_k\rangle\\
        &\quad + 2\gamma \langle x^k - x^*, \omega_k \rangle
        + \gamma^2\|F(x^k) - h^*\|^2
        + \gamma^2\|\omega_k\|^2\\
        &\quad + \frac{9\cdot 10^6\gamma^2
        A^2}{n^2}
        \sum\limits_{i=1}^n \Bigr[
        (1-\nu)\|h^k_i - h^*_i\|^2\\
        &\qquad + 2\nu \|\omega_{i,k}\|^2
        + 2\nu \|F_i(x^k) - h^*_i\|^2\Bigl]\\
        &\overset{\eqref{ch3:eq:star_cocoercivity_chapter_3}}{\le}
        \|x^k-x^*\|^2\\
        &\quad + (1-\nu)
        \frac{9\cdot 10^6\gamma^2
        A^2}{n^2}
        \sum\limits_{i=1}^n\|h^k_i - h^*_i\|^2\\
        &\quad -2 \gamma \left(
        1 - \frac{\gamma \ell}{2}
        - \frac{\nu n}{\gamma}\cdot
        \frac{9\cdot 10^6\gamma^2
        A^2}{n^2}
        \ell\right)
        \langle x^k - x^*, F(x^k)-h^* \rangle\\
        &\quad + \frac{2\gamma}{n}\sum\limits_{i=1}^n
        \langle x^k - x^* - \gamma(F(x^k) - h^*),
        \omega_{i,k} \rangle\\
        &\quad + \gamma^2\|\omega_k\|^2
        + \frac{\gamma^2}{n^2}\sum\limits_{i=1}^n\|\omega_{i,k}\|^2\\
        \displaybreak[3]
        &\overset{ \eqref{ch3:eq:gamma_lemma_prox_SGDA_str_mon_2_chapter_3}}{\le}
        \|x^k-x^*\|^2\\
        &\quad + (1-\nu)
        \frac{9\cdot 10^6\gamma^2
        A^2}{n^2}
        \sum\limits_{i=1}^n\|h^k_i - h^*_i\|^2\\
        &\quad - \gamma \langle x^k - x^*, F(x^k)-h^* \rangle\\
        &\quad + \frac{2\gamma}{n}\sum\limits_{i=1}^n
        \langle x^k - x^* - \gamma(F(x^k) - h^*),
        \omega_{i,k} \rangle\\
        &\quad + \gamma^2\|\omega_k\|^2
        + \frac{\gamma^2}{n^2}\sum\limits_{i=1}^n\|\omega_{i,k}\|^2\\
        &\overset{\eqref{ch3:eq:str_monotonicity_chapter_3}}{\le}
        (1-\gamma \mu) \|x^k - x^*\|^2\\
        &\quad + (1-\nu)
        \frac{9\cdot 10^6\gamma^2
        A^2}{n^2}
        \sum\limits_{i=1}^n\|h^k_i - h^*_i\|^2\\
        &\quad + \frac{2\gamma}{n}\sum\limits_{i=1}^n
        \langle x^k - x^* - \gamma(F(x^k) - h^*),
        \omega_{i,k} \rangle\\
        &\quad + \gamma^2\|\omega_k\|^2
        + \frac{\gamma^2}{n^2}\sum\limits_{i=1}^n\|\omega_{i,k}\|^2\\
        &\le
        (1 - \gamma \mu) V_k\\
        &\quad + \frac{2\gamma}{n}\sum\limits_{i=1}^n
        \langle x^k - x^* - \gamma(F(x^k) - h^*),
        \omega_{i,k} \rangle\\
        &\quad + \gamma^2\|\omega_k\|^2
        + \frac{\gamma^2}{n^2}\sum\limits_{i=1}^n\|\omega_{i,k}\|^2.
    \end{align*}
    Here the first inequality uses \eqref{ch3:eq:star_cocoercivity_chapter_3}, which gives $\sum_{i=1}^n\|F_i(x^k) - h_i^*\|^2 \leq n\ell\langle x^k - x^*, F(x^k) - h^*\rangle$ and $\|F(x^k) - h^*\|^2 \leq \ell\langle x^k - x^*, F(x^k) - h^*\rangle$, together with $2\nu\cdot 9\cdot 10^6A^2 \leq 1$; the second one uses \eqref{ch3:eq:gamma_lemma_prox_SGDA_str_mon_2_chapter_3}, which gives $1 - \gamma\ell/2 - 9\cdot 10^6\nu\gamma A^2\ell/n \geq 1/2$; and the last one uses $\gamma\mu \leq \nu$. Unrolling the recurrence, we obtain \eqref{ch3:eq:optimization_lemma_prox_SGDA_str_mon_2_chapter_3}.
\end{proof}

\begin{theorem}\label{ch3:thm:main_result_prox_clipped_SGDA_2_chapter_3}
Let Assumptions~\ref{ch3:as:str_monotonicity_chapter_3}, \ref{ch3:as:star_cocoercivity_chapter_3} hold for $Q = B_{\sqrt{2V}}(x^*)$, where
    \begin{gather*}
        A = \ln\frac{48n(K+1)}{\beta},\quad V \geq \|x^0 - x^*\|^2 + \frac{9\cdot 10^6\gamma^2A^2}{n^2}\sum_{i=1}^n\|F_i(x^*)\|^2,
    \end{gather*}
    let $R \geq \|x^0 - x^*\|$ and $\zeta_* \geq \sqrt{(1/n)\sum_{i=1}^n\|F_i(x^*)\|^2}$, let the initial shifts be $h_1^0 = \ldots = h_n^0 = 0$ and the shift parameter be $\nu = 1/{(18\cdot 10^6A^2)}$, and let
    \begin{align}
        0< \gamma &\leq \min\left\{\frac{1}{4096 \ell A}, \frac{\sqrt{n}R}{3000 \zeta_* A}, \frac{1}{18\cdot 10^6\mu A^2}, \frac{\ln(B_K)}{\mu(K+1)}\right\}, \label{ch3:eq:gamma_SGDA_2_prox_clipped_chapter_3}\\
        Q_K &= \frac{(K+1)^{2(p-1)/p}\mu^2V n^{2(p-1)/p}}{\sigma^2A^{2(p-1)/p}}, \notag\\
        B_K\ln^2(B_K) &= \frac{Q_K}{(768\cdot 576)^{2/p}}, \label{ch3:eq:B_K_SGDA_2_prox_clipped_1_chapter_3} \\
        B_K &= \Theta\left(\frac{Q_K}{\ln^2 Q_K}\right), \label{ch3:eq:B_K_SGDA_2_prox_clipped_2_chapter_3}\\
        \lambda_k &= \frac{n\exp(-\gamma\mu(1 + k/2))\sqrt{V}}{256\sqrt{2}\gamma A}, \label{ch3:eq:lambda_SGDA_2_prox_clipped_chapter_3}
    \end{align}
    for some $K \geq 0$ and $\beta \in (0,1]$ such that $B_K \geq 2$. Then, after $K$ iterations the iterates produced by \algname{DProx-clipped-SGDA-shift} with probability at least $1 - \beta$ satisfy 
    \begin{equation}
        V_{K+1} \leq 2\exp(-\gamma\mu(K+1))V, \label{ch3:eq:main_result_str_mon_SGDA_2_chapter_3}
    \end{equation}
    where $V_k = \|x^k - x^*\|^2 + (9\cdot 10^6\gamma^2A^2/n^2)\sum_{i=1}^n\|h^k_i - h^*_i\|^2$ and $h_i^* = F_i(x^*)$. 
    In particular, the second entry of the minimum in \eqref{ch3:eq:gamma_SGDA_2_prox_clipped_chapter_3} makes $V = 2R^2$ admissible, and when $\gamma$ equals the minimum from \eqref{ch3:eq:gamma_SGDA_2_prox_clipped_chapter_3}, then the iterates produced by \algname{DProx-clipped-SGDA-shift} after $K$ iterations with probability at least $1-\beta$ satisfy
    \begin{align}
       V_{K+1} = \cO\biggl(\max\biggl\{ & R^2\exp\left(- \frac{\mu K}{\ell \ln(nK/\beta)}\right), R^2\exp\left(- \frac{\mu \sqrt{n} RK}{\zeta_* \ln(nK/\beta)}\right), \notag\\
       & R^2\exp\left(- \frac{K}{\ln^2(nK/\beta)}\right), \notag\\
       & \frac{\sigma^2\ln^{2(p-1)/p}(nK/\beta)\ln^2B_K}{K^{2(p-1)/p}\mu^2 n^{2(p-1)/p}}\biggr\}\biggr), \label{ch3:eq:clipped_SGDA_2_prox_clipped_case_1_chapter_3}
    \end{align}
    meaning that to achieve $V_{K+1} \leq \varepsilon$ with probability at least $1 - \beta$ \algname{DProx-clipped-SGDA-shift} requires
    \begin{align}
        K = \cO\biggl(\max\biggl\{ & \frac{\ell}{\mu}\ln\frac{R^2}{\varepsilon}\ln\left(\frac{n\ell}{\mu \beta}\ln\frac{R^2}{\varepsilon}\right), \notag\\
        & \frac{\zeta_*}{\sqrt{n}R\mu}\ln\frac{R^2}{\varepsilon}\ln\left(\frac{\sqrt{n}\zeta_*}{R\mu \beta}\ln\frac{R^2}{\varepsilon}\right), \ln\frac{R^2}{\varepsilon}\ln^2\left(\frac{n}{\beta}\ln\frac{R^2}{\varepsilon}\right), \notag\\
        & \frac{1}{n}\left(\frac{\sigma^2}{\mu^2\varepsilon}\right)^{p/(2(p-1))}\ln \left(\frac{1}{\beta} \left(\frac{\sigma^2}{\mu^2\varepsilon }\right)^{p/(2(p-1))}\right) \notag\\
        & \qquad\cdot\ln^{p/(p-1)}\left(B_\varepsilon\right)\biggr\}\biggr) \label{ch3:eq:clipped_SGDA_2_prox_clipped_case_complexity_appendix_chapter_3}
    \end{align}
    iterations/oracle calls, where
    \begin{equation*}
         B_\varepsilon = \max\left\{2, \frac{2R^2}{\varepsilon \ln \left(\frac{1}{\beta} \left(\sigma^2/(\mu^2\varepsilon)\right)^{p/(2(p-1))}\right)}\right\}.
    \end{equation*}
\end{theorem}
\begin{proof}
The Lyapunov function has the following form
    \begin{equation*}
        V_k = \|x^k - x^*\|^2 + \frac{9\cdot 10^6\gamma^2A^2}{n^2}\sum\limits_{i=1}^n\|h^k_i - h^*_i\|^2.
    \end{equation*} Similar to previous results, our proof is induction-based.  To formulate the statement rigorously, we introduce probability event $E_k$ for each $k = 0,1,\ldots, K+1$ as follows: inequalities
    \begin{equation}
        V_t \leq 2 \exp(-\gamma\mu t) V \label{ch3:eq:induction_inequality_prox_clipped_SGDA_2_chapter_3}
    \end{equation}
    \begin{gather}
         \left\|\frac{\gamma}{n}\sum\limits_{i=1}^{r-1}\omega_{i,t-1}^u\right\| \leq \exp\left(-\frac{\gamma\mu (t-1)}{2}\right)\frac{\sqrt{V}}{2} \label{ch3:eq:induction_inequality_str_mon_SGDA_2_chapter_3}
    \end{gather}
    hold for $t = 0,1,\ldots,k$ and $r = 1, 2, \ldots, n$ simultaneously. We will prove by induction that $\Prob\{E_k\} \geq  1 - k\beta/(K+1)$ for all $k = 0,1,\ldots,K+1$. The base of induction follows immediately by the definition of $V$. Next, we assume that the statement holds for $k = T-1 \leq K$, i.e., $\Prob\{E_{T-1}\} \geq  1 - (T-1)\beta/(K+1)$. Let us show that it also holds for $k = T$, i.e., $\Prob\{E_{T}\} \geq  1 - T\beta/(K+1)$.

    Inequality \eqref{ch3:eq:induction_inequality_prox_clipped_SGDA_2_chapter_3} gives $\|x^t - x^*\|^2 \leq V_t \leq 2V$, i.e., $E_{T-1}$ implies $x^t \in B_{\sqrt{2V}}(x^*)$ for $t = 0, \ldots, T-1$. Moreover, the choice of $\nu$ and \eqref{ch3:eq:gamma_SGDA_2_prox_clipped_chapter_3} give $\gamma\mu \leq \nu$ and $\gamma \leq 1/(2\ell) \leq 1/(\ell + 18\cdot 10^6\nu A^2\ell/n)$. This means that we can apply Lemma~\ref{ch3:lem:optimization_lemma_str_mon_SGDA_2_chapter_3}: $E_{T-1}$ implies

    \begin{align*}
        V_{T} & \leq (1 - \gamma \mu)^{T}V + \frac{2\gamma}{n} \sum\limits_{t=0}^{T-1}\sum\limits_{i=1}^n (1-\gamma\mu)^{T-1-t}\langle x^t - x^* - \gamma (F(x^t)-h^*), \omega_{i,t} \rangle\notag\\
        &\quad\quad  + \frac{\gamma^2}{n^2}\sum\limits_{t=0}^{T-1}\sum\limits_{i=1}^n (1-\gamma \mu)^{T-1-t} \|\omega_{i,t}\|^2 + \gamma^2 \sum\limits_{t=0}^{T-1} (1-\gamma\mu)^{T-1-t} \|\omega_t\|^2.
    \end{align*}
Before we proceed, we introduce a new notation:
\begin{gather}
        r_t = \begin{cases} y_t,& \text{if } \|y_t\| \leq 2\sqrt{2} \exp(- \gamma\mu t/2)\sqrt{V},\\ 0,& \text{otherwise}, \end{cases} \label{ch3:eq:xi_t_SGDA_2_str_mon_chapter_3}
\end{gather}
where $y_t = x^t - x^* - \gamma (F(x^t)-h^*)$,
for $t = 0, 1, \ldots, T-1$. Random vectors $\{r_t\}_{t=0}^{T-1}$ are bounded almost surely:
\begin{equation}
         \|r_t\| \leq 2\sqrt{2}\exp(-\gamma\mu t/2)\sqrt{V} \label{ch3:eq:xi_t_bound_SGDA_2_prox_clipped_chapter_3} 
\end{equation}
for all $t = 0, 1, \ldots, T-1$. 
    In addition, $r_t = x^t - x^* - \gamma (F(x^t)-h^*)$ follows from $E_{T-1}$ for all $t = 0, 1, \ldots, T-1$, since $\|x^t - x^* - \gamma(F(x^t) - h^*)\| \leq (1 + \gamma\ell)\|x^t - x^*\|$ by \eqref{ch3:eq:star_cocoercivity_chapter_3}, and, thus, using $(1-\gamma\mu)^T \leq \exp(-\gamma\mu T)$ and $V_0 \leq V$, we get that $E_{T-1}$ implies
\begin{align}
        V_{T} & \leq \exp(-\gamma\mu T)V + \underbrace{\frac{2\gamma}{n} \sum\limits_{t=0}^{T-1}\sum\limits_{i=1}^n (1-\gamma\mu)^{T-1-t}\langle r_t, \omega^u_{i,t} \rangle}_{\circledOne} \notag\\
        &\quad\quad + \underbrace{\frac{2\gamma}{n} \sum\limits_{t=0}^{T-1}\sum\limits_{i=1}^n (1-\gamma\mu)^{T-1-t}\langle r_t, \omega^b_{i,t} \rangle}_{\circledTwo}\notag\\
        &\quad\quad  + \underbrace{\frac{4\gamma^2}{n^2}  \sum\limits_{t=0}^{T-1}\sum\limits_{i=1}^n (1-\gamma \mu)^{T-1-t} \left[\|\omega^u_{i,t}\|^2 - \Exp_{\xi^t}[\|\omega^u_{i,t}\|^2] \right]}_{\circledThree}\notag\\
        &\quad\quad  + \underbrace{\frac{4\gamma^2}{n^2} \sum\limits_{t=0}^{T-1}\sum\limits_{i=1}^n (1-\gamma \mu)^{T-1-t} \Exp_{\xi^t}[\|\omega^u_{i,t}\|^2]}_{\circledFour} \notag\\
        &\quad\quad  + \underbrace{\frac{4\gamma^2}{n} \sum\limits_{t=0}^{T-1}\sum\limits_{i=1}^n (1-\gamma \mu)^{T-1-t} \|\omega^b_{i,t}\|^2}_{\circledFive}\notag\\
        &\quad\quad + \underbrace{\frac{4\gamma^2}{n^2} \sum\limits_{t=0}^{T-1}\sum\limits_{j=2}^n (1-\gamma\mu)^{T-1-t} \left\langle \sum\limits_{i=1}^{j-1} \omega_{i,t}^u, \omega_{j,t}^u\right\rangle}_{\circledSix}. \label{ch3:eq:SGDA_2_prox_clipped_12345_bound_chapter_3}
\end{align}
To derive high-probability bounds for $\circledOne, \circledTwo, \circledThree, \circledFour, \circledFive, \circledSix$ we need to establish several useful inequalities related to $\omega_{i,t}^u, \omega_{i,t}^b$. First, by definition of clipping
    \begin{equation}
         \|\omega_{i,t}^u\| \leq 2\lambda_t. \label{ch3:eq:omega_magnitude_prox_clipped_2_chapter_3}
    \end{equation}
Next, we notice that $E_{T-1}$ implies
\begin{align}
        \|F_i(x^t) - h_i^t\|
        &\le \|F_i(x^t) - h^*_i\| + \|h^t_i- h^*_i\| \notag\\
        &\overset{\eqref{ch3:eq:star_cocoercivity_chapter_3}}{\le}
        \ell\|x^t - x^*\|
        + \sqrt{\sum\limits_{i=1}^n\|h^t_i -  h^*_i\|^2}\notag\\
        &\le \left( \ell
        + \frac{n}{3000\gamma A}\right)
        \sqrt{V_t}\notag\\
        &\overset{\eqref{ch3:eq:induction_inequality_prox_clipped_SGDA_2_chapter_3}}{\leq}
        \sqrt{2}\left( \ell
        + \frac{n}{3000\gamma A}\right)
        \exp(- \gamma\mu t/2)\sqrt{V}\notag\\
        &\overset{\eqref{ch3:eq:gamma_SGDA_2_prox_clipped_chapter_3},
        \eqref{ch3:eq:lambda_SGDA_2_prox_clipped_chapter_3}}{\leq}
        \frac{\lambda_t}{2}.
        \label{ch3:eq:operator_bound_x_t_SGDA_2_prox_clipped_chapter_3}
    \end{align}
for $t = 0,1, \ldots, T-1$ and $i\in [n]$. Therefore, one can apply Lemma~\ref{ch3:lem:bias_and_variance_clip_chapter_3} and get \begin{gather}
         \left\|\omega_t^b\right\| \leq \frac{1}{n}\sum\limits_{i=1}^n \|\omega_{i,t}^b\| \leq \frac{2^{p}\sigma^{p}}{\lambda_t^{{p}-1}}, \label{ch3:eq:bias_omega_prox_clipped_2_chapter_3}\\
         \Exp_{\xi_{i}^t}\left[\left\|\omega_{i,t}^u\right\|^2\right] \leq 18 \lambda_t^{2-{p}}\sigma^{p}, \label{ch3:eq:variance_omega_prox_clipped_2_chapter_3}
    \end{gather}
    for all $t = 0,1, \ldots, T-1$ and $i\in [n]$. Finally, in the bounds below we use $(1-\gamma\mu)^{s} \leq \exp(-\gamma\mu s)$, $T \leq K+1$, and
    \begin{equation}
        \gamma\mu \leq \frac{1}{4096},\quad \exp\left(\frac{\gamma\mu{p}K}{2}\right) \leq B_K^{{p}/2},\quad \gamma \leq \frac{\ln(B_K)}{\mu(K+1)},\label{ch3:eq:nu_theorem_prox_SGDA_str_mon_2_chapter_3}
    \end{equation}
    which follow from \eqref{ch3:eq:gamma_SGDA_2_prox_clipped_chapter_3} and $\mu \leq \ell$; together with \eqref{ch3:eq:B_K_SGDA_2_prox_clipped_1_chapter_3}, they give each of the inequalities marked by \eqref{ch3:eq:gamma_SGDA_2_prox_clipped_chapter_3} below.

\paragraph{Upper bound for $\circledOne$.}
To estimate this sum, we will use Bernstein's inequality. The summands have conditional expectations equal to zero:
\begin{align*}
\MoveEqLeft[1] \Exp_{\xi_{i}^l}\left[\frac{2\gamma}{n} (1-\gamma\mu)^{T-1-l} \langle r_l, \omega_{i,l}^u \rangle\right] \\
& = \frac{2\gamma}{n} (1-\gamma\mu)^{T-1-l}\left\langle r_l, \Exp_{\xi_{i}^l}[\omega_{i,l}^u] \right\rangle = 0.
    \end{align*}
Moreover, for all $l = 0,\ldots, T-1$ random vectors $\{\omega_{i,l}^u\}_{i=1}^n$ are independent. Thus, sequence $\{(2\gamma/n)(1-\gamma\mu)^{T-1-l} \langle r_l, \omega_{i,l}^u \rangle\}_{l,i = 0,1}^{T-1,n}$ is a martingale difference sequence. Next, the summands are bounded:
    \begin{align}
        \left|\frac{2\gamma}{n} (1-\gamma\mu)^{T - 1 - l} \langle r_l, \omega_{i,l}^u \rangle \right| & \leq \frac{2\gamma}{n}\exp(-\gamma\mu (T - 1 - l)) \|r_l\|\cdot \|\omega_{i,l}^u\|\notag\\
        & \overset{\eqref{ch3:eq:xi_t_bound_SGDA_2_prox_clipped_chapter_3},\eqref{ch3:eq:omega_magnitude_prox_clipped_2_chapter_3}}{\leq} \frac{8\sqrt{2}}{n}\gamma  \exp(-\gamma\mu (T - 1 - l/2)) \sqrt{V} \lambda_l\notag\\
        & \overset{\eqref{ch3:eq:lambda_SGDA_2_prox_clipped_chapter_3}}{\leq} \frac{\exp(-\gamma\mu T)V}{8A} \eqdef c. \label{ch3:eq:SGDA_2_prox_clipped_technical_1_1_xi_chapter_3}
    \end{align}
    Finally, conditional variances $\sigma_{i,l}^2 \eqdef \Exp_{\xi_{i}^l}\left[\frac{4\gamma^2}{n^2}(1-\gamma\mu)^{2(T-1-l)} \langle r_l, \omega_{i,l}^u \rangle^2\right]$ of the summands are bounded:
    \begin{align}
        \sigma_{i,l}^2 & \leq \Exp_{\xi^l_i}\left[\frac{4\gamma^2}{n^2}\exp(-2\gamma\mu (T - 1 - l)) \|r_l\|^2\cdot \|\omega_{i,l}^u\|^2\right]\notag\\
        & \overset{\eqref{ch3:eq:xi_t_bound_SGDA_2_prox_clipped_chapter_3}}{\leq} \frac{32\gamma^2}{n^2}  \exp(-\gamma\mu (2T - 2 - l)) V \Exp_{\xi^l_i}\left[\|\omega_{i,l}^u\|^2\right]. \label{ch3:eq:SGDA_2_prox_clipped_technical_1_2_xi_chapter_3}
    \end{align}
Applying Bernstein's inequality (Lemma~\ref{ch3:lem:Bernstein_ineq_chapter_3}) with $X_{i,l}$ equal to the summands of $\circledOne$, constant $c$ defined in \eqref{ch3:eq:SGDA_2_prox_clipped_technical_1_1_xi_chapter_3}, $b = \frac{\exp\left(-\gamma\mu T\right)V}{8}$, $H = \frac{\exp\left(-2\gamma\mu T\right)V^2}{384A}$, we get

    \begin{align*}
        &\Prob\left\{|\circledOne| > \frac{\exp(-\gamma\mu T) V}{8} \text{ and } \sum\limits_{l=0}^{T-1}\sum\limits_{i=1}^{n}\sigma_{i,l}^2 \leq \frac{\exp(- 2\gamma\mu T) V^2} {384A}\right\}\\
        &\leq 2\exp\left(- \frac{b^2}{2H + 2cb/3}\right)\\
        &= \frac{\beta}{24n(K+1)}.
    \end{align*}
    The above is equivalent to $\Prob\{E_{\circledOne}\} \geq 1 - \frac{\beta}{24n(K+1)}$ for
    \begin{align}
         E_{\circledOne} = \Biggl\{&\text{either} \quad \sum\limits_{l=0}^{T-1}\sum\limits_{i=1}^{n}\sigma_{i,l}^2 > \frac{\exp(- 2\gamma\mu T) V^2}{384A} \notag\\
         &\text{or}\quad |\circledOne| \leq \frac{\exp(-\gamma\mu T) V}{8}\Biggr\}. \label{ch3:eq:bound_1_SGDA_2_prox_clipped_chapter_3}
    \end{align}
    Moreover, $E_{T-1}$ implies
    \begin{align}
        \sum\limits_{l=0}^{T-1}\sum\limits_{i=1}^n\sigma_{i,l}^2 & \overset{\eqref{ch3:eq:SGDA_2_prox_clipped_technical_1_2_xi_chapter_3}}{\leq} \frac{32\gamma^2}{n^2}\exp(- 2\gamma\mu (T-1))V\sum\limits_{l=0}^{T-1}\sum\limits_{i=1}^n \frac{\Exp_{\xi^l_i}\left[\|\omega_{i,l}^u\|^2\right]}{\exp(-\gamma\mu l)}\notag\\ 
        & \overset{\eqref{ch3:eq:variance_omega_prox_clipped_2_chapter_3}}{\leq} \frac{576\gamma^2}{n}\exp(-2\gamma\mu (T-1)) V \sigma^{p} \sum\limits_{l=0}^{T-1} \frac{\lambda_l^{2-{p}}}{\exp(-\gamma\mu l)}\notag\\
         & \overset{\eqref{ch3:eq:lambda_SGDA_2_prox_clipped_chapter_3}}{\leq} \frac{576\gamma^{p}\exp(-2\gamma\mu (T-1)) V^{2-{p}/2} \sigma^{p} (K+1)\exp(\gamma\mu{p} K/2)}{(256\sqrt{2})^{2-{p}} n^{{p}-1}A^{2-{p}}}\notag\\
        & \overset{\eqref{ch3:eq:gamma_SGDA_2_prox_clipped_chapter_3}}{\leq} \frac{\exp(-2\gamma\mu T)V^2}{384A}. \label{ch3:eq:bound_1_variances_SGDA_2_prox_clipped_chapter_3}
    \end{align}
    
\paragraph{Upper bound for $\circledTwo$.} Probability event $E_{T-1}$ implies
    \begin{align}
        \circledTwo
        &\leq \frac{2\gamma}{n} \exp(-\gamma\mu (T-1))
        \sum\limits_{l=0}^{T-1}\sum\limits_{i=1}^n
        \frac{\|r_l\|\cdot \|\omega_{i,l}^b\|}{\exp(-\gamma\mu l)}
        \notag\\
        &\overset{\eqref{ch3:eq:xi_t_bound_SGDA_2_prox_clipped_chapter_3},
        \eqref{ch3:eq:bias_omega_prox_clipped_2_chapter_3}}{\leq}
        2^{2+{p}}\sqrt{2} \gamma \exp(-\gamma\mu (T-1))
        \sqrt{V} \sigma^{p}
        \sum\limits_{l=0}^{T-1}
        \frac{1}{\lambda_l^{{p}-1}
        \exp(-\gamma\mu l/2)}
        \notag\\
        &\overset{\eqref{ch3:eq:lambda_SGDA_2_prox_clipped_chapter_3}}{\leq}
        \frac{(512\sqrt{2})^{{p}}}{64}\cdot
        \frac{\gamma^{p} \sigma^{{p}} \exp(-\gamma\mu (T-1)) (K+1)\exp(\gamma\mu{p}K/2)\exp(\gamma\mu{p})}
        {n^{{p}-1}V^{{p}/2-1}A^{1-{p}}}
        \notag\\
        &\overset{\eqref{ch3:eq:gamma_SGDA_2_prox_clipped_chapter_3}}{\leq}
        \frac{\exp(-\gamma\mu T) V}{8}.
        \label{ch3:eq:bound_2_SGDA_2_prox_clipped_chapter_3}
    \end{align}

\paragraph{Upper bound for $\circledThree$.} To estimate this sum, we will use Bernstein's inequality. The summands have conditional expectations equal to zero:
    \begin{equation*}
        \Exp_{\xi^l_i}\left[\frac{4\gamma^2}{n^2}(1-\gamma\mu)^{T-1-l}\left[\|\omega_{i,l}^u\|^2  -\Exp_{\xi^l_i}\left[\|\omega_{i,l}^u\|^2\right]\right] \right] = 0.
    \end{equation*}
    Moreover, for all $l = 0,\ldots, T-1$ random vectors $\{\omega_{i,l}^u\}_{i=1}^n$ are independent. Thus, sequence $\{(4\gamma^2/n^2)(1-\gamma\mu)^{T-1-l}(\|\omega_{i,l}^u\|^2 - \Exp_{\xi_{i}^l}[\|\omega_{i,l}^u\|^2])\}_{l,i = 0,1}^{T-1,n}$ is a martingale difference sequence. Next, the summands are bounded:
    \begin{align}
        \frac{4\gamma^2}{n^2} (1-\gamma\mu)^{T-1-l}\left| \|\omega_{i,l}^u\|^2  -\Exp_{\xi^l_i}\left[\|\omega_{i,l}^u\|^2\right] \right| & \overset{\eqref{ch3:eq:omega_magnitude_prox_clipped_2_chapter_3}}{\leq}
        \frac{32\gamma^2\lambda_l^2}{n^2} \frac{\exp(-\gamma\mu T) }{\exp(-\gamma\mu (1+l))}\notag\\
        & \overset{\eqref{ch3:eq:lambda_SGDA_2_prox_clipped_chapter_3}}{\leq} \frac{\exp(-\gamma\mu (T+1))V}{256A^2}\notag\\
        & \leq \frac{\exp(-\gamma\mu T)V}{8A}
        \eqdef c. \label{ch3:eq:SGDA_2_prox_clipped_technical_4_1_chapter_3}
    \end{align}

    Finally, conditional variances 
    \begin{equation*}
        \widetilde\sigma_{i,l}^2 \eqdef \Exp_{\xi^l_i}\left[\frac{16\gamma^4}{n^4}  (1-\gamma\mu)^{2T-2-2l} \left|\|\omega_{i,l}^u\|^2  -\Exp_{\xi^l_i}\left[\|\omega_{i,l}^u\|^2\right] \right|^2 \right]
    \end{equation*}
    of the summands are bounded:
    \begin{align}
        \widetilde\sigma_{i,l}^2 & \overset{\eqref{ch3:eq:SGDA_2_prox_clipped_technical_4_1_chapter_3}}{\leq} \frac{4\gamma^2\exp(-2\gamma\mu T)V}{8n^2\exp(-\gamma\mu (1+l))A} \Exp_{\xi^l_i}\left[ \left|\|\omega_{i,l}^u\|^2  -\Exp_{\xi^l_i}\left[\|\omega_{i,l}^u\|^2\right] \right|\right]\notag\\
        & \leq \frac{\gamma^2\exp(-2\gamma\mu T)V}{n^2\exp(-\gamma\mu (1+l))A} \Exp_{\xi^l_i}\left[\|\omega_{i,l}^u\|^2\right]. \label{ch3:eq:SGDA_2_prox_clipped_technical_4_2_chapter_3}
    \end{align}
    Applying Bernstein's inequality (Lemma~\ref{ch3:lem:Bernstein_ineq_chapter_3}) with
    \begin{equation*}
        X_{i,l} = \frac{4\gamma^2}{n^2} (1-\gamma\mu)^{T-1-l}
        \left[\|\omega_{i,l}^u\|^2
        -\Exp_{\xi^l_i}\left[\|\omega_{i,l}^u\|^2\right]\right],
    \end{equation*}
    constant $c$ defined in \eqref{ch3:eq:SGDA_2_prox_clipped_technical_4_1_chapter_3}, $b = \frac{\exp\left(-\gamma\mu T\right)V}{8}$, $H = \frac{\exp\left(-2\gamma\mu T\right)V^2}{384A}$, we get
    \begin{align*}
        &\Prob\left\{|\circledThree| > \frac{\exp(-\gamma\mu T) V}{8} \text{ and } \sum\limits_{l=0}^{T-1}\sum\limits_{i=1}^n\widetilde\sigma_{i,l}^2 \leq \frac{\exp(-2\gamma\mu T) V^2} {384A}\right\}\\
        &\leq 2\exp\left(- \frac{b^2}{2H + 2cb/3}\right)\\
        &= \frac{\beta}{24n(K+1)}.
    \end{align*}
    The above is equivalent to $\Prob\{E_{\circledThree}\} \geq 1 - \frac{\beta}{24n(K+1)}$ for
    \begin{align}
        E_{\circledThree} = \Biggl\{&\text{either} \quad \sum\limits_{l=0}^{T-1}\sum\limits_{i=1}^n\widetilde\sigma_{i,l}^2 > \frac{\exp(-2\gamma\mu T) V^2}{384A} \notag\\
        &\text{or}\quad |\circledThree| \leq \frac{\exp(-\gamma\mu T) V}{8}\Biggr\}. \label{ch3:eq:bound_5_SGDA_2_prox_clipped_chapter_3}
    \end{align}
    Moreover, $E_{T-1}$ implies
    \begin{align}
        \sum\limits_{l=0}^{T-1}\sum\limits_{i=1}^n\widetilde\sigma_{i,l}^2 & \overset{\eqref{ch3:eq:SGDA_2_prox_clipped_technical_4_2_chapter_3}}{\leq} \frac{\gamma^2\exp(-\gamma\mu (2T-1))V}{n^2A} \sum\limits_{l=0}^{T-1}\sum\limits_{i=1}^n \frac{\Exp_{\xi^l_i}\left[\|\omega_{i,l}^u\|^2\right]}{\exp(-\gamma\mu l)}\notag\\
        & \overset{\eqref{ch3:eq:variance_omega_prox_clipped_2_chapter_3}}{\leq} \frac{18\gamma^2\exp(-\gamma\mu (2T-1)) V \sigma^{p}}{nA} \sum\limits_{l=0}^{T-1} \frac{\lambda_l^{2-{p}}}{\exp(-\gamma\mu l)}\notag\\
        & \overset{\eqref{ch3:eq:lambda_SGDA_2_prox_clipped_chapter_3}}{\leq} 
        \frac{18\gamma^{p}\exp(-\gamma\mu (2T-1)) V^{2-{p}/2} \sigma^{p} (K+1)\exp(\gamma\mu{p} K/2)}{(256\sqrt{2})^{2-{p}}n^{{p}-1}A^{3-{p}}} \notag\\
        & \overset{\eqref{ch3:eq:gamma_SGDA_2_prox_clipped_chapter_3}}{\leq} \frac{\exp(-2\gamma\mu T)V^2}{384A}. \label{ch3:eq:bound_5_variances_SGDA_2_prox_clipped_chapter_3}
    \end{align}

\paragraph{Upper bound for $\circledFour$.} Probability event $E_{T-1}$ implies
    \begin{align}
        \circledFour & = \frac{4\gamma^2}{n^2} \sum\limits_{l=0}^{T-1}\sum\limits_{i=1}^n (1-\gamma \mu)^{T-1-l} \Exp_{\xi^l_i}[\|\omega^u_{i,l}\|^2]\notag\\
        & \overset{\eqref{ch3:eq:variance_omega_prox_clipped_2_chapter_3}}{\leq} \frac{72\gamma^2\exp(-\gamma\mu (T-1)) \sigma^{p}}{n}\sum\limits_{l=0}^{T-1} \frac{\lambda_l^{2-{p}}}{\exp(-\gamma\mu l)} \notag\\
        & \overset{\eqref{ch3:eq:lambda_SGDA_2_prox_clipped_chapter_3}}{\leq} \frac{72\gamma^{p} V^{1-{p}/2}\exp(-\gamma\mu (T-1)) \sigma^{p} (K+1)\exp(\gamma\mu{p} K/2)}{(256\sqrt{2})^{2-{p}}n^{{p}-1} A^{2-{p}}} \notag\\
        & \overset{\eqref{ch3:eq:gamma_SGDA_2_prox_clipped_chapter_3}}{\leq} \frac{\exp(-\gamma\mu T) V}{8} . \label{ch3:eq:bound_6_SGDA_2_prox_clipped_chapter_3}
    \end{align}

\paragraph{Upper bound for $\circledFive$.} Probability event $E_{T-1}$ implies
    \begin{align}
        \circledFive & =  \frac{4\gamma^2}{n} \sum\limits_{l=0}^{T-1}\sum\limits_{i=1}^n (1-\gamma \mu)^{T-1-l} \|\omega^b_{i,l}\|^2\notag\\
        & \overset{\eqref{ch3:eq:bias_omega_prox_clipped_2_chapter_3}}{\leq} 4\cdot 2^{2{p}}\gamma^2 \exp(-\gamma\mu (T-1)) \sigma^{2{p}} \sum\limits_{l=0}^{T-1} \frac{1}{\lambda_l^{2{p}-2} \exp(-\gamma\mu l)} \notag\\
        & \overset{\eqref{ch3:eq:lambda_SGDA_2_prox_clipped_chapter_3}}{\leq}  4^{{p}+1}(256\sqrt{2})^{2{p}-2}\frac{\gamma^{2{p}} \sigma^{2{p}} A^{2({p}-1)} (K+1)}{n^{2({p}-1)}V^{{p}-1}}\notag\\
        &\qquad \cdot\exp\left(-\gamma\mu (T-1) + \gamma\mu {p} K + \gamma\mu(2{p}-2)\right)\notag\\
        & \overset{\eqref{ch3:eq:gamma_SGDA_2_prox_clipped_chapter_3}}{\leq} \frac{\exp(-\gamma\mu T) V}{8}. \label{ch3:eq:bound_7_SGDA_2_prox_clipped_chapter_3}
    \end{align}

\paragraph{Upper bounds for $\circledSix$.} This sum requires more refined analysis. We introduce new vectors:
\begin{align}
        \delta_{j}^l = \begin{cases}
            \frac{\gamma}{n}\sum\limits_{i=1}^{j-1} \omega_{i,l}^u,& \text{if } \left\|\frac{\gamma}{n}\sum\limits_{i=1}^{j-1} \omega_{i,l}^u\right\| \leq \exp\left(-\frac{\gamma\mu l}{2}\right)\frac{\sqrt{V}}{2},\\
            0,& \text{otherwise,}
        \end{cases}\label{ch3:eq:SGDA_zeta_delta_def_str_mon_2_chapter_3}
\end{align}
    for all $j \in [n]$ and $l = 0,\ldots, T-1$. Then, by definition
    \begin{equation}
        \|\delta_j^l\| \leq \exp\left(-\frac{\gamma\mu l}{2}\right)\frac{\sqrt{V}}{2} \label{ch3:eq:SGDA_zeta_delta_bound_str_mon_chapter_3}
    \end{equation}
    and
    \begin{align}
        \circledSix & = \underbrace{\frac{4\gamma}{n} \sum\limits_{l=0}^{T-1} \sum\limits_{j=2}^n (1-\gamma\mu)^{T-1-l}\left\langle \delta_j^l, \omega_{j,l}^u \right\rangle}_{\circledSix'}\notag\\
        &\quad\quad + \frac{4\gamma}{n} \sum\limits_{l=0}^{T-1} \sum\limits_{j=2}^n (1-\gamma\mu)^{T-1-l} \left\langle \frac{\gamma}{n}\sum\limits_{i=1}^{j-1} \omega_{i,l}^u - \delta_j^l, \omega_{j,l}^u \right\rangle. \label{ch3:eq:SGDA_extra_sums_seven_distributed_str_mon_chapter_3}
    \end{align}
We also note here that $E_{T-1}$ implies
    \begin{align}
        \MoveEqLeft[1] \frac{4\gamma}{n} \sum\limits_{l=0}^{T-1} \sum\limits_{j=2}^n (1-\gamma\mu)^{T-1-l}\left\langle \frac{\gamma}{n}\sum\limits_{i=1}^{j-1} \omega_{i,l}^u - \delta_j^l, \omega_{j,l}^u \right\rangle \notag\\
        &= \frac{4\gamma}{n}  \sum\limits_{j=2}^n \left\langle \frac{\gamma}{n}\sum\limits_{i=1}^{j-1} \omega_{i,T-1}^u - \delta_j^{T-1}, \omega_{j,T-1}^u \right\rangle. \label{ch3:eq:SGDA_extra_sums_seven_distributed_1_str_mon_chapter_3}
    \end{align}
\paragraph{Upper bound for $\circledSix'$.} To estimate this sum, we will use Bernstein's inequality. The summands have conditional expectations equal to zero:
    \begin{align*}
\MoveEqLeft[1] \Exp_{\xi_{j}^l}\left[\frac{4\gamma}{n}(1-\gamma\mu)^{T-1-l} \langle \delta_j^l, \omega_{j,l}^u \rangle\right] \\
& = \frac{4\gamma}{n} (1-\gamma\mu)^{T-1-l}\left\langle \delta_j^l, \Exp_{\xi_{j}^l}[\omega_{j,l}^u] \right\rangle = 0.
    \end{align*}
    Moreover, for all $l = 0,\ldots, T-1$ random vectors $\{\omega_{j,l}^u\}_{j=1}^n$ are independent. Thus, sequence $\{(4\gamma/n)(1-\gamma\mu)^{T-1-l} \langle \delta_j^l, \omega_{j,l}^u \rangle\}_{l,j = 0,2}^{T-1,n}$ is a martingale difference sequence. Next, the summands are bounded:
    \begin{align}
        \left|\frac{4\gamma}{n} (1-\gamma\mu)^{T-1-l}\langle \delta_j^l, \omega_{j,l}^u \rangle \right| & \leq \frac{4\gamma}{n}\exp\left(-\gamma\mu(T-1-l)\right) \|\delta_j^l\|\cdot \|\omega_{j,l}^u\|\notag\\
        & \overset{\eqref{ch3:eq:SGDA_zeta_delta_bound_str_mon_chapter_3},\eqref{ch3:eq:omega_magnitude_prox_clipped_2_chapter_3}}{\leq} \frac{4\sqrt{V}\gamma \exp\left(-\gamma\mu(T-1)\right)}{n} \exp\left(\frac{\gamma\mu l}{2}\right)\lambda_l\notag\\
        & \overset{\eqref{ch3:eq:lambda_SGDA_2_prox_clipped_chapter_3}}{=} \frac{\exp\left(-\gamma\mu T\right)V}{16\sqrt{2}A}\notag\\
        & \leq \frac{\exp\left(-\gamma\mu T\right)V}{8A} \eqdef c. \label{ch3:eq:str_mon_thm_SGDA_technical_6_5_seven_chapter_3}
    \end{align}
    Finally, conditional variances $(\sigma_{j,l}')^2 \eqdef \Exp_{\xi_{j}^l}\left[\frac{16\gamma^2}{n^2}(1-\gamma\mu)^{2(T-1-l)} \langle \delta_j^l, \omega_{j,l}^u \rangle^2\right]$ of the summands are bounded:
    \begin{align}
        (\sigma_{j,l}')^2 & \leq \Exp_{\xi_{j}^l}\left[\frac{16\gamma^2 }{n^2}\exp\left(-\gamma\mu(2T-2-2l)\right)\|\delta_j^l\|^2\cdot \|\omega_{j,l}^u\|^2\right] \notag\\
        & \overset{\eqref{ch3:eq:SGDA_zeta_delta_bound_str_mon_chapter_3}}{\leq} \frac{4\gamma^2 V\exp\left(-\gamma\mu\left(2T-2-l\right)\right)}{n^2} \Exp_{\xi_{j}^l}\left[\|\omega_{j,l}^u\|^2\right]. \label{ch3:eq:str_mon_thm_SGDA_technical_7_seven_chapter_3}
    \end{align}
    Applying Bernstein's inequality (Lemma~\ref{ch3:lem:Bernstein_ineq_chapter_3}) with $X_{j,l}$ equal to the summands of $\circledSix'$, constant $c$ defined in \eqref{ch3:eq:str_mon_thm_SGDA_technical_6_5_seven_chapter_3}, $b = \frac{\exp\left(-\gamma\mu T\right)V}{8}$, $H = \frac{\exp\left(-2\gamma\mu T\right)V^2}{384A}$, we get
    \begin{align*}
        &\Prob\left\{|\circledSix'| > \frac{\exp\left(-\gamma\mu T\right)V}{8} \text{ and } \sum\limits_{l=0}^{T-1}\sum\limits_{j=2}^n (\sigma_{j,l}')^2 \leq \frac{\exp\left(-2\gamma\mu T\right)V^2}{384A}\right\}\\
        &\leq 2\exp\left(- \frac{b^2}{2H + 2cb/3}\right) = \frac{\beta}{24n(K+1)}.
    \end{align*}
    The above is equivalent to $\Prob\{E_{\circledSix'}\} \geq 1 - \frac{\beta}{24n(K+1)}$ for
    \begin{align}
        E_{\circledSix'} = \Biggl\{&\text{either } \sum\limits_{l=0}^{T-1}\sum\limits_{j=2}^n (\sigma_{j,l}')^2 > \frac{\exp\left(-2\gamma\mu T\right)V^2}{384A} \notag\\
        &\text{or } |\circledSix'| \leq \frac{\exp\left(-\gamma\mu T\right)V}{8}\Biggr\}. \label{ch3:eq:bound_6_str_mon_SGDA_seven_chapter_3}
    \end{align}
    Moreover, $E_{T-1}$ implies
    \begin{align}
        \sum\limits_{l=0}^{T-1}\sum\limits_{j=2}^n (\sigma_{j,l}')^2 & \overset{\eqref{ch3:eq:str_mon_thm_SGDA_technical_7_seven_chapter_3}}{\leq} \frac{4\gamma^2 V\exp\left(-\gamma\mu\left(2T-2\right)\right)}{n^2} \sum\limits_{l=0}^{T-1}\exp\left(\gamma\mu l\right) \sum\limits_{i=1}^n \Exp_{\xi_{i}^l}\left[\|\omega_{i,l}^u\|^2\right]\notag\\
        & \overset{\eqref{ch3:eq:variance_omega_prox_clipped_2_chapter_3}}{\leq}  \frac{72\gamma^2 V\exp\left(-\gamma\mu\left(2T-2\right)\right)\sigma^{p}}{n} \sum\limits_{l=0}^{T-1}\exp\left(\gamma\mu l\right) \lambda_l^{2-{p}}\notag\\
        & \overset{\eqref{ch3:eq:lambda_SGDA_2_prox_clipped_chapter_3}}{\leq} \frac{72\gamma^{p} V^{2 - \frac{{p}}{2}}\exp\left(-2\gamma\mu (T-1)\right)\sigma^{p}}{(256\sqrt{2})^{2-{p}} n^{{p}-1} A^{2-{p}}} \sum\limits_{l=0}^{T-1}\exp\left(\frac{\gamma\mu l {p}}{2}\right)\notag\\
        & \leq \frac{72\gamma^{p} V^{2 - \frac{{p}}{2}}\exp\left(-2\gamma\mu (T-1)\right)\sigma^{p} (K+1) \exp\left(\frac{\gamma\mu K {p}}{2}\right)}{(256\sqrt{2})^{2-{p}} n^{{p}-1} A^{2-{p}}}\notag\\
        & \overset{\eqref{ch3:eq:gamma_SGDA_2_prox_clipped_chapter_3}}{\leq} \frac{\exp\left(-2\gamma\mu T\right)V^2}{384A}.\label{ch3:eq:bound_6_variances_str_mon_SGDA_seven_chapter_3}
    \end{align}

That is, we derive the upper bounds for  $\circledOne, \circledTwo, \circledThree, \circledFour, \circledFive, \circledSix$. More precisely, $E_{T-1}$ implies
    \begin{gather*}
        V_T \overset{\eqref{ch3:eq:SGDA_2_prox_clipped_12345_bound_chapter_3}}{\leq} \exp\left(-\gamma\mu T\right)V + \circledOne + \circledTwo + \circledThree + \circledFour + \circledFive + \circledSix,\\
         \circledSix \overset{\eqref{ch3:eq:SGDA_extra_sums_seven_distributed_str_mon_chapter_3}}{=} \circledSix' + \frac{4\gamma}{n}  \sum\limits_{j=2}^n \left\langle \frac{\gamma}{n}\sum\limits_{i=1}^{j-1} \omega_{i,T-1}^u - \delta_j^{T-1}, \omega_{j,T-1}^u \right\rangle,\\
         \circledTwo \overset{\eqref{ch3:eq:bound_2_SGDA_2_prox_clipped_chapter_3}}{\leq} \frac{\exp\left(-\gamma\mu T\right) V}{8},\quad
         \circledFour \overset{\eqref{ch3:eq:bound_6_SGDA_2_prox_clipped_chapter_3}}{\leq} \frac{\exp\left(-\gamma\mu T\right) V}{8},\\ \circledFive \overset{\eqref{ch3:eq:bound_7_SGDA_2_prox_clipped_chapter_3}}{\leq} \frac{\exp\left(-\gamma\mu T\right) V}{8}, \\
         \sum\limits_{l=0}^{T-1}\sum\limits_{i=1}^n\sigma_{i,l}^2 \overset{\eqref{ch3:eq:bound_1_variances_SGDA_2_prox_clipped_chapter_3}}{\leq} \frac{\exp\left(-2\gamma\mu T\right)V^2}{384A},\quad 
         \sum\limits_{l=0}^{T-1}\sum\limits_{i=1}^n \widetilde \sigma_{i,l}^2 \overset{\eqref{ch3:eq:bound_5_variances_SGDA_2_prox_clipped_chapter_3}}{\leq} \frac{\exp\left(-2\gamma\mu T\right)V^2}{384A},\\ \sum\limits_{l=0}^{T-1}\sum\limits_{j=2}^n (\sigma_{j,l}')^2 \overset{\eqref{ch3:eq:bound_6_variances_str_mon_SGDA_seven_chapter_3}}{\leq} \frac{\exp\left(-2\gamma\mu T\right)V^2}{384A}.
    \end{gather*}
    In addition, we also establish (see \eqref{ch3:eq:bound_1_SGDA_2_prox_clipped_chapter_3}, \eqref{ch3:eq:bound_5_SGDA_2_prox_clipped_chapter_3}, \eqref{ch3:eq:bound_6_str_mon_SGDA_seven_chapter_3}, and our induction assumption)
    \begin{gather*}
        \Prob\{E_{T-1}\} \geq 1 - \frac{(T-1)\beta}{K+1},\quad \Prob\{E_{\circledOne}\} \geq 1 - \frac{\beta}{24n(K+1)}, \\ \Prob\{E_{\circledThree}\} \geq 1 - \frac{\beta}{24n(K+1)},\quad \Prob\{E_{\circledSix'}\} \geq 1 - \frac{\beta}{24n(K+1)},
    \end{gather*}
    where
    \begin{align*}
        E_{\circledOne} & = \left\{\text{either } \sum\limits_{l=0}^{T-1}\sum\limits_{i=1}^n\sigma_{i,l}^2 > \frac{\exp\left(-2\gamma\mu T\right)V^2}{384A} \text{ or } |\circledOne| \leq \frac{\exp\left(-\gamma\mu T\right)V}{8}\right\},\\
        E_{\circledThree} & = \left\{\text{either } \sum\limits_{l=0}^{T-1}\sum\limits_{i=1}^n\widetilde \sigma_{i,l}^2 > \frac{\exp\left(-2\gamma\mu T\right)V^2}{384A} \text{ or } |\circledThree| \leq \frac{\exp\left(-\gamma\mu T\right)V}{8}\right\},\\
        E_{\circledSix'} & = \left\{\text{either } \sum\limits_{l=0}^{T-1}\sum\limits_{j=2}^n (\sigma_{j,l}')^2 > \frac{\exp\left(-2\gamma\mu T\right)V^2}{384A} \text{ or } |\circledSix'| \leq \frac{\exp\left(-\gamma\mu T\right)V}{8}\right\}.
    \end{align*}
	    Therefore, probability event $E_{T-1} \cap E_{\circledOne}  \cap E_{\circledThree} \cap E_{\circledSix'}$ implies
	    \begin{align}
	        V_T
	        &\leq \exp\left(-\gamma\mu T\right)V
	        \underbrace{\left(1 + \frac{1}{8} + \frac{1}{8} + \frac{1}{8}
	        + \frac{1}{8} + \frac{1}{8} + \frac{1}{8}\right)}_{< 2}\notag\\
	        &\quad + \frac{4\gamma}{n}  \sum\limits_{j=2}^n
	        \left\langle \frac{\gamma}{n}\sum\limits_{i=1}^{j-1} \omega_{i,T-1}^u - \delta_j^{T-1}, \omega_{j,T-1}^u \right\rangle.
	        \label{ch3:eq:SGDA_new_chapter_3}
	    \end{align}

To finish the proof, we need to show that  $\frac{\gamma}{n}\sum\limits_{i=1}^{j-1} \omega_{i,T-1}^u = \delta_j^{T-1}$ with high probability. In particular, we consider probability event $\widetilde{E}_{T-1,j}$ defined as follows: inequalities
    \begin{equation}
         \left\|\frac{\gamma}{n}\sum\limits_{i=1}^{r-1} \omega_{i,T-1}^u\right\| \leq \exp\left(-\frac{\gamma\mu (T-1)}{2}\right)\frac{\sqrt{V}}{2} \notag
    \end{equation}
    hold for $r = 2, \ldots, j$ simultaneously. We want to show that $\Prob\{E_{T-1} \cap \widetilde{E}_{T-1,j}\} \geq 1 - \frac{(T-1)\beta}{K+1} - \frac{j\beta}{8n(K+1)}$ for all $j = 2, \ldots, n$. For $j = 2$ the statement is trivial since
    \begin{align*}
        \left\|\frac{\gamma}{n} \omega_{1,T-1}^u\right\| & \overset{\eqref{ch3:eq:omega_magnitude_prox_clipped_2_chapter_3}}{\leq} \frac{2\gamma\lambda_{T-1}}{n} \leq \exp\left(-\frac{\gamma\mu (T-1)}{2}\right)\frac{\sqrt{V}}{2}.
    \end{align*}
    Next, we assume that the statement holds for some $j = m-1 < n$, i.e.,  $\Prob\{E_{T-1}\cap\widetilde{E}_{T-1,m-1}\} \geq 1 - \frac{(T-1)\beta}{K+1} - \frac{(m-1)\beta}{8n(K+1)}$. Our goal is to prove that $\Prob\{E_{T-1}\cap\widetilde{E}_{T-1,m}\} \geq 1 - \frac{(T-1)\beta}{K+1} - \frac{m\beta}{8n(K+1)}$. First, we consider $\left\|\frac{\gamma}{n}\sum\limits_{i=1}^{m-1} \omega_{i,T-1}^u\right\|$:
    \begin{align*}
        \left\|\frac{\gamma}{n}\sum\limits_{i=1}^{m-1} \omega_{i,T-1}^u\right\|
        &= \sqrt{\frac{\gamma^2}{n^2}
        \left\|\sum\limits_{i=1}^{m-1} \omega_{i,T-1}^u\right\|^2}\\
        &= \sqrt{
        \begin{aligned}
        &\frac{\gamma^2}{n^2}\sum\limits_{i=1}^{m-1}
        \|\omega_{i,T-1}^u\|^2\\
        &+ \frac{2\gamma}{n}\sum\limits_{i=1}^{m-1}
        \left\langle \frac{\gamma}{n}\sum\limits_{r=1}^{i-1}
        \omega_{r,T-1}^u , \omega_{i,T-1}^u \right\rangle
        \end{aligned}}\\
        &\leq \sqrt{
        \begin{aligned}
        &\frac{\gamma^2}{n^2}\sum\limits_{t=0}^{T-1}
        (1-\gamma\mu)^{T-1-t}
        \sum\limits_{i=1}^{m-1}\|\omega_{i,t}^u\|^2\\
        &+ \frac{2\gamma}{n}\sum\limits_{i=1}^{m-1}
        \left\langle \frac{\gamma}{n}\sum\limits_{r=1}^{i-1}
        \omega_{r,T-1}^u , \omega_{i,T-1}^u \right\rangle
        \end{aligned}}.
    \end{align*}

    Next, we introduce a new notation:
    \begin{gather*}
        \rho_{i,T-1}' = \begin{cases}
            \frac{\gamma}{n}\sum\limits_{r=1}^{i-1}\omega_{r,T-1}^u,& \text{if } \left\|\frac{\gamma}{n}\sum\limits_{r=1}^{i-1}\omega_{r,T-1}^u\right\| \leq \exp\left(-\frac{\gamma\mu (T-1)}{2}\right)\frac{\sqrt{V}}{2},\\
            0,& \text{otherwise}
        \end{cases}
    \end{gather*}
    for $i = 1,\ldots,m-1$. By definition, we have
    \begin{equation}
         \|\rho_{i,T-1}'\| \leq \exp\left(-\frac{\gamma\mu (T-1)}{2}\right)\frac{\sqrt{V}}{2} \label{ch3:eq:SGDA_bound_rho_str_mon_chapter_3}
    \end{equation}
    for $i = 1,\ldots,m-1$. Moreover, $\widetilde{E}_{T-1,m-1}$ implies $\rho_{i,T-1}' = \frac{\gamma}{n}\sum\limits_{r=1}^{i-1}\omega_{r,T-1}^u$ for $i = 1,\ldots,m-1$ and
    \begin{align}
        \left\|\frac{\gamma}{n}\sum\limits_{i=1}^{m-1} \omega_{i,T-1}^u\right\| & \leq \sqrt{\frac{1}{4}\left(\circledThree + \circledFour\right) + \circledSeven}, \notag
    \end{align}
    where
    \begin{gather*}
        \circledSeven =  \frac{2\gamma}{n}\sum\limits_{i=1}^{m-1}\left\langle \rho_{i,T-1}' , \omega_{i,T-1}^u \right\rangle.
    \end{gather*}
    It remains to estimate $\circledSeven$.

\paragraph{Upper bound for $\circledSeven$.} To estimate this sum, we will use Bernstein's inequality. The summands have conditional expectations equal to zero:
    \begin{align*}
\MoveEqLeft[1] \Exp_{\xi_{i}^{T-1}}\left[\frac{2\gamma}{n} \langle \rho_{i,T-1}', \omega_{i,T-1}^u \rangle\right] \\
& = \frac{2\gamma}{n} \left\langle \rho_{i,T-1}', \Exp_{\xi_{i}^{T-1}}[\omega_{i,T-1}^u] \right\rangle = 0.
    \end{align*}
    Thus, for $i = 1,\ldots,m-1$ the summands $(2\gamma/n)\langle \rho_{i,T-1}', \omega_{i,T-1}^u \rangle$ form a martingale difference sequence. Next, the summands are bounded:
    \begin{align}
        \left|\frac{2\gamma}{n} \langle \rho_{i,T-1}', \omega_{i,T-1}^u \rangle \right| & \leq \frac{2\gamma}{n}\|\rho_{i,T-1}'\|\cdot \|\omega_{i,T-1}^u\|\notag\\
        & \overset{\eqref{ch3:eq:SGDA_bound_rho_str_mon_chapter_3},\eqref{ch3:eq:omega_magnitude_prox_clipped_2_chapter_3}}{\leq} \frac{2\sqrt{V}\gamma \exp\left(-\frac{\gamma\mu(T-1)}{2}\right)}{n} \lambda_{T-1}\notag\\
        & \overset{\eqref{ch3:eq:lambda_SGDA_2_prox_clipped_chapter_3}}{=} \frac{\exp\left(-\gamma\mu T\right)V}{32\sqrt{2}A} \notag\\
        & \leq \frac{\exp\left(-\gamma\mu T\right)V}{8A} \eqdef c. \label{ch3:eq:str_mon_thm_SGDA_technical_6_5_eight'_chapter_3}
    \end{align}
    Finally, conditional variances $(\widetilde\sigma_{i,T-1}')^2 \eqdef \Exp_{\xi_{i}^{T-1}}\left[\frac{4\gamma^2}{n^2} \langle \rho_{i,T-1}', \omega_{i,T-1}^u \rangle^2\right]$ of the summands are bounded:
    \begin{align}
        (\widetilde\sigma_{i,T-1}')^2 & \leq \Exp_{\xi_{i}^{T-1}}\left[\frac{4\gamma^2 }{n^2}\|\rho_{i,T-1}'\|^2\cdot \|\omega_{i,T-1}^u\|^2\right] \notag\\
        & \overset{\eqref{ch3:eq:SGDA_bound_rho_str_mon_chapter_3}}{\leq} \frac{\gamma^2 V\exp\left(-\gamma\mu(T-1)\right)}{n^2} \Exp_{\xi_{i}^{T-1}}\left[\|\omega_{i,T-1}^u\|^2\right]. \label{ch3:eq:str_mon_thm_SGDA_technical_7_eight'_chapter_3}
    \end{align}
    Applying Bernstein's inequality (Lemma~\ref{ch3:lem:Bernstein_ineq_chapter_3}) with $X_{i} = \frac{2\gamma}{n}\langle \rho_{i,T-1}', \omega_{i,T-1}^u \rangle$, constant $c$ defined in \eqref{ch3:eq:str_mon_thm_SGDA_technical_6_5_eight'_chapter_3}, $b = \frac{\exp\left(-\gamma\mu T\right)V}{8}$, $H = \frac{\exp\left(-2\gamma\mu T\right)V^2}{384A}$, we get
    \begin{align*}
        &\Prob\left\{|\circledSeven| > \frac{\exp\left(-\gamma\mu T\right)V}{8} \text{ and } \sum\limits_{i=1}^{m-1}(\widetilde\sigma_{i,T-1}')^2 \leq \frac{\exp\left(-2\gamma\mu T\right)V^2}{384A}\right\}\\
        &\leq 2\exp\left(- \frac{b^2}{2H + 2cb/3}\right) = \frac{\beta}{24n(K+1)}.
    \end{align*}
    The above is equivalent to $\Prob\{E_{\circledSeven}\} \geq 1 - \frac{\beta}{24n(K+1)}$ for
    \begin{equation}
        E_{\circledSeven} = \left\{\text{either } \sum\limits_{i=1}^{m-1}(\widetilde\sigma_{i,T-1}')^2 > \frac{\exp\left(-2\gamma\mu T\right)V^2}{384A} \text{ or } |\circledSeven| \leq \frac{\exp\left(-\gamma\mu T\right)V}{8}\right\}. \label{ch3:eq:bound_1_str_mon_SGDA_eight'_chapter_3}
    \end{equation}
    Moreover, $E_{T-1}$ implies
    \begin{align}
        \sum\limits_{i=1}^{m-1}(\widetilde\sigma_{i,T-1}')^2 & \overset{\eqref{ch3:eq:str_mon_thm_SGDA_technical_7_eight'_chapter_3}}{\leq} \frac{\gamma^2 V\exp\left(-\gamma\mu\left(T-1\right)\right)}{n^2} \sum\limits_{i=1}^{m-1} \Exp_{\xi_{i}^{T-1}}\left[\|\omega_{i,T-1}^u\|^2\right]\notag\\
        & \overset{\eqref{ch3:eq:variance_omega_prox_clipped_2_chapter_3}}{\leq}  \frac{18\gamma^2 V\exp\left(-\gamma\mu\left(T-1\right)\right)\sigma^{p}}{n} \lambda_{T-1}^{2-{p}}\notag\\
        & \overset{\eqref{ch3:eq:lambda_SGDA_2_prox_clipped_chapter_3}}{\leq} \frac{18\gamma^{p} V^{2 - \frac{{p}}{2}}\exp\left(-\gamma\mu (T-1)(4-{p})/2\right)\sigma^{p}}{(256\sqrt{2})^{2-{p}} n^{{p}-1} A^{2-{p}}}\notag\\
        & \leq \frac{18\gamma^{p} V^{2 - \frac{{p}}{2}}\exp\left(-2\gamma\mu T\right)\sigma^{p} \exp\left(\gamma\mu K {p}/2 + 2\gamma\mu\right)}{(256\sqrt{2})^{2-{p}} n^{{p}-1} A^{2-{p}}}\notag\\
        & \overset{\eqref{ch3:eq:gamma_SGDA_2_prox_clipped_chapter_3}}{\leq} \frac{\exp\left(-2\gamma\mu T\right)V^2}{384A}.\label{ch3:eq:bound_1_variances_str_mon_SGDA_eight'_chapter_3}
    \end{align}
Putting all together we get that $E_{T-1}\cap \widetilde{E}_{T-1,m-1}$ implies
    \begin{gather*}
        \left\|\frac{\gamma}{n}\sum\limits_{i=1}^{m-1} \omega_{i,T-1}^u\right\| \leq \sqrt{\frac{1}{4}\left(\circledThree + \circledFour\right) + \circledSeven},\\ \circledFour \overset{\eqref{ch3:eq:bound_6_SGDA_2_prox_clipped_chapter_3}}{\leq} \frac{\exp\left(-\gamma\mu T\right)V}{8},\quad
        \sum\limits_{l=0}^{T-1}\sum\limits_{i=1}^n\widetilde\sigma_{i,l}^2 \overset{\eqref{ch3:eq:bound_5_variances_SGDA_2_prox_clipped_chapter_3}}{\leq} \frac{\exp\left(-2\gamma\mu T\right)V^2}{384A},\\
        \sum\limits_{i=1}^{m-1}(\widetilde\sigma_{i,T-1}')^2 \leq \frac{\exp\left(-2\gamma\mu T\right)V^2}{384A}.
    \end{gather*}
In addition, we also establish (see \eqref{ch3:eq:bound_5_SGDA_2_prox_clipped_chapter_3}, \eqref{ch3:eq:bound_1_str_mon_SGDA_eight'_chapter_3} and our induction assumption)
    \begin{gather*}
        \Prob\{E_{T-1}\cap \widetilde{E}_{T-1,m-1}\} \geq 1 - \frac{(T-1)\beta}{K+1} - \frac{(m-1)\beta}{8n(K+1)},\\
        \Prob\{E_{\circledThree}\} \geq 1 - \frac{\beta}{24n(K+1)}, \quad \Prob\{E_{\circledSeven}\} \geq 1 - \frac{\beta}{24n(K+1)}
    \end{gather*}
    where
    \begin{align}
        E_{\circledThree} & = \left\{\text{either } \sum\limits_{l=0}^{T-1}\sum\limits_{i=1}^n\widetilde\sigma_{i,l}^2 > \frac{\exp\left(-2\gamma\mu T\right)V^2}{384A} \text{ or } |\circledThree| \leq \frac{\exp\left(-\gamma\mu T\right)V}{8}\right\}, \notag\\
        E_{\circledSeven} & = \left\{\text{either } \sum\limits_{i=1}^{m-1}(\widetilde\sigma_{i,T-1}')^2 > \frac{\exp\left(-2\gamma\mu T\right)V^2}{384A} \text{ or } |\circledSeven| \leq \frac{\exp\left(-\gamma\mu T\right)V}{8}\right\}. \notag
    \end{align}
    Therefore, probability event $E_{T-1} \cap \widetilde{E}_{T-1,m-1} \cap E_{\circledThree} \cap E_{\circledSeven}$ implies
    \begin{align*}
        \left\|\frac{\gamma}{n}\sum\limits_{i=1}^{m-1} \omega_{i,T-1}^u\right\| & \leq \exp\left(-\frac{\gamma\mu(T-1)}{2}\right)\sqrt{V}\sqrt{\frac{1}{32} + \frac{1}{32} + \frac{1}{8}}\\
        & \leq \frac{1}{2}\exp\left(-\frac{\gamma\mu(T-1)}{2}\right)\sqrt{V}.
    \end{align*}
This implies $\widetilde E_{T-1,m}$ and
    \begin{align*}
        \Prob\{E_{T-1} \cap \widetilde{E}_{T-1,m}\} & \geq \Prob\{E_{T-1} \cap \widetilde{E}_{T-1,m-1} \cap E_{\circledThree} \cap E_{\circledSeven}\} \\
        & = 1 - \Prob\left\{\overline{E_{T-1} \cap \widetilde{E}_{T-1,m-1}}  \cup \overline{E}_{\circledThree} \cup \overline{E}_{\circledSeven}\right\}\\
        & \geq 1 - \frac{(T-1)\beta}{K+1} - \frac{m\beta}{8n(K+1)}.
    \end{align*}
    Therefore, for all $m = 2,\ldots,n$ the statement holds and, in particular, $\Prob\{E_{T-1} \cap \widetilde{E}_{T-1, n}\} \geq 1 - \frac{(T-1)\beta}{K+1} - \frac{\beta}{8(K+1)}$, i.e., \eqref{ch3:eq:induction_inequality_str_mon_SGDA_2_chapter_3} holds. Taking into account \eqref{ch3:eq:SGDA_new_chapter_3}, we conclude that $E_{T-1} \cap \widetilde{E}_{T-1, n} \cap E_{\circledOne} \cap E_{\circledThree} \cap E_{\circledSix'}$ implies
    \begin{equation*}
        V_T \leq 2\exp\left(-\gamma\mu T\right)V
    \end{equation*}
    that is equivalent to \eqref{ch3:eq:induction_inequality_prox_clipped_SGDA_2_chapter_3} for $t = T$. Moreover,
    \begin{align*}
        \Prob\left\{E_T\right\} & \geq \Prob\left\{E_{T-1} \cap \widetilde{E}_{T-1,n} \cap E_{\circledOne} \cap E_{\circledThree} \cap E_{\circledSix'} \right\} \\
        & = 1 - \Prob\left\{\overline{E_{T-1} \cap \widetilde{E}_{n}} \cup \overline{E}_{\circledOne} \cup \overline{E}_{\circledThree} \cup \overline{E}_{\circledSix'}\right\}\\
        & \geq 1 - \frac{(T-1)\beta}{K+1} - \frac{\beta}{8(K+1)} - 3\cdot \frac{\beta}{24n(K+1)} \geq 1 - \frac{T\beta}{K+1}.
    \end{align*}

    In other words, we showed that $\Prob\{E_k\} \geq 1 - k\beta/(K+1)$ for all $k = 0,1,\ldots,K+1$. For $k = K+1$ we have that with probability at least $1 - \beta$
    \begin{equation}
        \|x^{K+1} - x^*\|^2 \leq 2\exp(-\gamma\mu (K+1))V. \notag
    \end{equation}
    Finally, if $V = 2R^2$ and $\gamma$ equals the minimum in \eqref{ch3:eq:gamma_SGDA_2_prox_clipped_chapter_3}, then with probability at least $1-\beta$
    \begin{align*}
        \|x^{K+1} - x^*\|^2 &\leq 2\exp(-\gamma\mu (K+1))V\\
        &= 2V\max\biggl\{\exp\left(-\frac{\mu(K+1)}{4096 \ell A}\right), \exp\left(- \frac{\mu \sqrt{n} R(K+1)}{3000\zeta_* A}\right),\\
        &\qquad\qquad \exp\left(-\frac{K+1}{18\cdot 10^6A^2}\right), \frac{1}{B_K} \biggr\}\\
        &= \cO\biggl(\max\biggl\{R^2\exp\left(- \frac{\mu K}{\ell \ln(nK/\beta)}\right), R^2\exp\left(- \frac{\mu \sqrt{n} RK}{\zeta_* \ln(nK/\beta)}\right),\\
        &\qquad\qquad R^2\exp\left(- \frac{K}{\ln^2(nK/\beta)}\right), \frac{\sigma^2\ln^{2(p-1)/p}(nK/\beta)\ln^2B_K}{K^{2(p-1)/p}\mu^2 n^{2(p-1)/p}}\biggr\}\biggr),
    \end{align*}
    where in the last step we use \eqref{ch3:eq:B_K_SGDA_2_prox_clipped_1_chapter_3}. To get $\|x^{K+1} - x^*\|^2 \leq \varepsilon$ with probability $\geq 1 - \beta$, $K$ should be
    \begin{align*}
        K = \cO\biggl(\max\biggl\{ & \frac{\ell}{\mu}\ln\frac{R^2}{\varepsilon}\ln\left(\frac{n\ell}{\mu \beta}\ln\frac{R^2}{\varepsilon}\right),\\
        & \frac{\zeta_*}{\sqrt{n}R\mu}\ln\frac{R^2}{\varepsilon}\ln\left(\frac{\sqrt{n}\zeta_*}{R\mu \beta}\ln\frac{R^2}{\varepsilon}\right), \ln\frac{R^2}{\varepsilon}\ln^2\left(\frac{n}{\beta}\ln\frac{R^2}{\varepsilon}\right),\\
        & \frac{1}{n}\left(\frac{\sigma^2}{\mu^2\varepsilon}\right)^{p/(2(p-1))}\ln \left(\frac{1}{\beta} \left(\frac{\sigma^2}{\mu^2\varepsilon }\right)^{p/(2(p-1))}\right)\ln^{p/(p-1)}\left(B_\varepsilon\right)\biggr\}\biggr),
    \end{align*}
    where
    \begin{equation*}
         B_\varepsilon = \max\left\{2, \frac{2R^2}{\varepsilon \ln \left(\frac{1}{\beta} \left(\sigma^2/(\mu^2\varepsilon)\right)^{p/(2(p-1))}\right)}\right\}.
    \end{equation*}
    This concludes the proof.
\end{proof}

\clearpage

\section{Additional Proofs for \algname{DProx-clipped-SEG-shift}}\label{ch3:appendix:dprox_clipped_SEG_shift_chapter_3}

In this section, we give the complete formulations of our results for \algname{DProx-clipped-SEG-shift} and rigorous proofs. For the readers' convenience, the method's update rule is repeated below:
\begin{gather*}
    \thesischthreetx^k = \thesischthreeprox_{\gamma\Psi}\left(x^k - \gamma \thesischthreetg^k\right),\;\; \thesischthreetg^k = \frac{1}{n}\sum\limits_{i=1}^n \thesischthreetg_i^k,\;\; \thesischthreetg_i^k = \thesischthreethh_{i}^k + \tilde\Delta_i^k,\;\; \thesischthreethh_i^{k+1} = \thesischthreethh_i^k + \nu \tilde \Delta_i^k,\\
    x^{k+1} = \thesischthreeprox_{\gamma\Psi}\left(x^k - \gamma \thesischthreehg^k\right),\;\; \thesischthreehg^k = \frac{1}{n}\sum\limits_{i=1}^n \thesischthreehg_i^k,\;\; \thesischthreehg_i^k = \thesischthreehh_{i}^k + \hat\Delta_i^k,\;\; \thesischthreehh_i^{k+1} = \thesischthreehh_i^k + \nu \hat \Delta_i^k, \\
    \tilde\Delta_i^k = \clip(F_i(x^k;\xi_{1,i}^k) - \thesischthreethh_i^k, \lambda_k),\;\; \hat\Delta_i^k = \clip(F_i(\thesischthreetx^k;\xi_{2,i}^k) - \thesischthreehh_i^k, \lambda_k).
\end{gather*}

\subsection{Monotone Case}

The following lemma is the main ``optimization'' part of the analysis of \algname{DProx-clipped-SEG-shift}.

\begin{lemma}\label{ch3:lem:SEG_optimization_lemma_monotone_chapter_3}
    Let Assumptions \ref{ch3:as:L_Lip_chapter_3} and \ref{ch3:as:monotonicity_chapter_3} hold for $Q = B_{4n\sqrt{V}}(x^*)$, where
    \begin{gather*}
        A = \ln\frac{48n(K+1)}{\beta},\quad V \geq \|x^0 - x^*\|^2 + \frac{409600\gamma^2A^2}{n^2}\sum_{i=1}^n\|F_i(x^*)\|^2,
    \end{gather*}
    and $0 < \gamma \leq 1/(\sqrt{12}L)$. If $x^k$ and $\thesischthreetx^k$ lie in $B_{4n\sqrt{V}}(x^*)$ for all $k = 0,1,\ldots,K$ for some $K \geq 0$, then for all $u \in B_{4n\sqrt{V}}(x^*)$ the iterates produced by \algname{DProx-clipped-SEG-shift} satisfy
    \begin{align}
        \langle F(u), \thesischthreetx^K_{\avg} - u\rangle + \Psi(\thesischthreetx^K_{\avg}) - \Psi(u) & \leq \frac{\|x^0 - u\|^2 - \|x^{K+1}-u\|^2}{2\gamma(K+1)} \notag\\
        &\quad\quad + \frac{\gamma}{K+1}\sum\limits_{k=0}^K (3\|\omega_k\|^2 + 4\|\theta_k\|^2) \notag\\
        &\quad\quad + \frac{1}{K+1}\sum\limits_{k=0}^K \langle \theta_k, x^k - u \rangle, \label{ch3:eq:optimization_lemma_SEG_chapter_3}\\
        \thesischthreetx^K_{\avg} & \eqdef \frac{1}{K+1}\sum\limits_{k=0}^{K}\thesischthreetx^k, \label{ch3:eq:tx_avg_SEG_chapter_3}\\
        \theta_k & \eqdef F(\thesischthreetx^k) - \thesischthreehg^k, \label{ch3:eq:theta_k_SEG_chapter_3}\\
        \omega_k & \eqdef F(x^k) - \thesischthreetg^k. \label{ch3:eq:omega_k_SEG_chapter_3}
    \end{align}
\end{lemma}
\begin{proof}
    Since $\thesischthreetx^k = \thesischthreeprox_{\gamma\Psi}\left(x^k - \gamma \thesischthreetg^k\right)$ and $x^{k+1} = \thesischthreeprox_{\gamma\Psi}\left(x^k - \gamma \thesischthreehg^k\right)$, we have $x^k - \gamma \thesischthreetg^k - \thesischthreetx^k \in \gamma \partial \Psi(\thesischthreetx^k)$ and $x^k - \gamma \thesischthreehg^k - x^{k+1} \in \gamma \partial\Psi(x^{k+1})$. By definition of the subgradient, we have $\forall u \in \R^d$
    \begin{align}
        \gamma\left(\Psi(\thesischthreetx^{k}) - \Psi(x^{k+1})\right) & \leq \langle \thesischthreetx^k - x^k + \gamma \thesischthreetg^k, x^{k+1} - \thesischthreetx^k \rangle, \notag\\
        \gamma\left(\Psi(x^{k+1}) - \Psi(u)\right) & \leq \langle x^{k+1} - x^k + \gamma \thesischthreehg^k, u - x^{k+1} \rangle. \notag
    \end{align}
    Summing up the above inequalities, we get
    \begin{align}
        \MoveEqLeft[1] \gamma\left(\Psi(\thesischthreetx^k) - \Psi(u)\right) \notag \\
        & \leq \langle \thesischthreetx^k - x^k, x^{k+1} - \thesischthreetx^k \rangle + \langle x^{k+1} - x^k, u - x^{k+1} \rangle \notag\\
        &\quad\quad + \gamma \langle \thesischthreetg^k - \thesischthreehg^k, x^{k+1} - \thesischthreetx^k \rangle + \gamma \langle \thesischthreehg^k, u - \thesischthreetx^k \rangle. \label{ch3:eq:ddmvklfkdvnfkvdn_chapter_3}
    \end{align}
    Since
    \begin{align}
        \langle \thesischthreetx^k - x^k, x^{k+1} - \thesischthreetx^k \rangle & = \frac{1}{2}\|x^{k+1} - x^k\|^2 - \frac{1}{2}\|\thesischthreetx^k - x^k\|^2 - \frac{1}{2}\|x^{k+1} - \thesischthreetx^k\|^2, \notag\\
        \langle x^{k+1} - x^k, u - x^{k+1} \rangle & = \frac{1}{2}\|x^k - u\|^2 - \frac{1}{2}\|x^{k+1} - x^k\|^2 - \frac{1}{2}\|x^{k+1} - u\|^2, \notag
    \end{align}
    we can rewrite \eqref{ch3:eq:ddmvklfkdvnfkvdn_chapter_3} as follows
    \begin{align}
        \MoveEqLeft[1] \gamma\left(\langle F(\thesischthreetx^k), \thesischthreetx^k - u \rangle +  \Psi(\thesischthreetx^k) - \Psi(u)\right) \notag \\
        & \leq \frac{1}{2}\|x^k - u\|^2 - \frac{1}{2}\|x^{k+1} - u\|^2  - \frac{1}{2}\|\thesischthreetx^k - x^k\|^2 \notag\\
        &\quad\quad - \frac{1}{2}\|x^{k+1} - \thesischthreetx^k\|^2 + \gamma \langle \thesischthreetg^k - \thesischthreehg^k, x^{k+1} - \thesischthreetx^k \rangle\notag\\
        &\quad\quad + \gamma \langle \theta_k, \thesischthreetx^k - u \rangle. \label{ch3:eq:shvghaxsacshdbjchjd_chapter_3}
    \end{align}
    Next, we upper-bound $\gamma \langle \thesischthreetg^k - \thesischthreehg^k, x^{k+1} - \thesischthreetx^k \rangle$ using Young's inequality, stating that $\langle a, b \rangle \leq \frac{1}{2\eta}\|a\|^2 + \frac{\eta}{2}\|b\|^2$ for all $a,b\in\R^d$ and $\eta > 0$, and Jensen's inequality for the squared norm:
    \begin{align}
        \MoveEqLeft[1] \gamma \langle \thesischthreetg^k - \thesischthreehg^k, x^{k+1} - \thesischthreetx^k \rangle \notag \\
        & \leq \gamma^2 \|\thesischthreetg^k - \thesischthreehg^k\|^2 + \frac{1}{4}\|x^{k+1} - \thesischthreetx^k\|^2 \notag\\
        & = \gamma^2\|F(x^k) - F(\thesischthreetx^k) - \omega_k + \theta_k\|^2\notag\\
        &\quad + \frac{1}{4}\|x^{k+1} - \thesischthreetx^k\|^2\notag\\
        & \leq 3\gamma^2 \|F(x^k) - F(\thesischthreetx^k)\|^2 + 3\gamma^2 \|\omega_k\|^2 + 3\gamma^2\|\theta_k\|^2\notag\\
        &\quad + \frac{1}{4}\|x^{k+1} - \thesischthreetx^k\|^2\notag\\
        & \overset{\eqref{ch3:eq:L_Lip_chapter_3}}{\leq} 3\gamma^2 L^2 \|x^k - \thesischthreetx^k\|^2 + 3\gamma^2 \|\omega_k\|^2 + 3\gamma^2\|\theta_k\|^2\notag\\
        &\quad + \frac{1}{4}\|x^{k+1} - \thesischthreetx^k\|^2 .\label{ch3:eq:njksdnjnsdkjncikds_chapter_3}
    \end{align}
    Plugging \eqref{ch3:eq:njksdnjnsdkjncikds_chapter_3} in \eqref{ch3:eq:shvghaxsacshdbjchjd_chapter_3}, we derive for all $u \in \R^d$
    \begin{align}
        \MoveEqLeft[1] \gamma\left(\langle F(\thesischthreetx^k), \thesischthreetx^k - u \rangle +  \Psi(\thesischthreetx^k) - \Psi(u)\right) \notag \\
        & \leq \frac{1}{2}\|x^k - u\|^2 - \frac{1}{2}\|x^{k+1} - u\|^2  \notag\\
        &\quad\quad - \frac{1}{2}\left(1 - 6\gamma^2 L^2\right)\|\thesischthreetx^k - x^k\|^2 - \frac{1}{4}\|x^{k+1} - \thesischthreetx^k\|^2 \notag\\
        &\quad\quad + 3\gamma^2 \|\omega_k\|^2 + 3\gamma^2\|\theta_k\|^2  + \gamma \langle \theta_k, \thesischthreetx^k - u \rangle. \label{ch3:eq:important_optimization_ineq_SEG_chapter_3}
    \end{align}
    We notice that the above inequality does not rely on monotonicity. Next, we apply monotonicity and get that for all $u \in B_{4n\sqrt{V}}(x^*)$:
    \begin{align}
        \MoveEqLeft[1] \gamma \left(\langle F(u), \thesischthreetx^k - u \rangle +  \Psi(\thesischthreetx^k) - \Psi(u)\right) \notag \\
        & \leq \frac{1}{2}\|x^k - u\|^2 - \frac{1}{2}\|x^{k+1} - u\|^2  \notag\\
        &\quad\quad - \frac{1}{2}\left(1 - 6\gamma^2 L^2\right)\|\thesischthreetx^k - x^k\|^2 + \gamma\langle \theta_k, \thesischthreetx^k - x^k \rangle\notag\\
        &\quad\quad + 3\gamma^2 \|\omega_k\|^2 + 3\gamma^2\|\theta_k\|^2  + \gamma \langle \theta_k, x^k - u \rangle \notag\\
        & \leq \frac{1}{2}\|x^k - u\|^2 - \frac{1}{2}\|x^{k+1} - u\|^2  \notag\\
        &\quad\quad - \frac{1}{2}\left(\frac{1}{2} - 6\gamma^2 L^2\right)\|\thesischthreetx^k - x^k\|^2 \notag\\
        &\quad\quad + 3\gamma^2 \|\omega_k\|^2 + 4\gamma^2\|\theta_k\|^2  + \gamma \langle \theta_k, x^k - u \rangle, \notag
    \end{align}
    where in the last step we apply $\gamma\langle \theta_k, \thesischthreetx^k - x^k \rangle \leq \gamma^2\|\theta_k\|^2 + \frac{1}{4}\|\thesischthreetx^k - x^k\|^2$. Since $\gamma \leq 1/(\sqrt{12}L)$, we have
    \begin{align}
        \MoveEqLeft[1] \gamma \left(\langle F(u), \thesischthreetx^k - u \rangle +  \Psi(\thesischthreetx^k) - \Psi(u)\right) \notag \\
        & \leq \frac{1}{2}\|x^k - u\|^2 - \frac{1}{2}\|x^{k+1} - u\|^2  \notag\\
        &\quad\quad + 3\gamma^2 \|\omega_k\|^2 + 4\gamma^2\|\theta_k\|^2  + \gamma \langle \theta_k, x^k - u \rangle, \notag
    \end{align}
    Summing up the above inequalities for $k = 0,1,\ldots, K$ and dividing both sides by $\gamma(K+1)$, we obtain
    \begin{align}
        \frac{1}{K+1}\sum\limits_{k=0}^{K}\left(\langle F(u), \thesischthreetx^k - u \rangle +  \Psi(\thesischthreetx^k) - \Psi(u)\right) & \leq \frac{\|x^0 - u\|^2 - \|x^{K+1}-u\|^2}{2\gamma(K+1)} \notag\\
        &\quad\quad + \frac{\gamma}{K+1}\sum\limits_{k=0}^{K} (3\|\omega_k\|^2 + 4\|\theta_k\|^2) \notag\\
        &\quad\quad + \frac{1}{K+1}\sum\limits_{k=0}^{K} \langle \theta_k, x^k - u \rangle. \notag
    \end{align}
    It remains to use that $\frac{1}{K+1}\sum_{k=0}^K \langle F(u), \thesischthreetx^k - u \rangle$ equals $\langle F(u), \thesischthreetx^K_{\avg} - u \rangle$ and that, by Jensen's inequality, $\Psi(\thesischthreetx^K_{\avg}) \leq \frac{1}{K+1}\sum_{k=0}^K\Psi(\thesischthreetx^k)$.
\end{proof}

Next, we proceed with the full statement of our main result for \algname{DProx-clipped-SEG-shift} in the monotone case.

\begin{theorem}[Case 2 from Theorem~\ref{ch3:thm:D_prox_clipped_SEG_main_chapter_3}]\label{ch3:thm:D_prox_clipped_SEG_monotone_appendix_chapter_3}
    Let Assumptions \ref{ch3:as:bounded_alpha_moment_chapter_3}, \ref{ch3:as:L_Lip_chapter_3} and \ref{ch3:as:monotonicity_chapter_3} hold for $Q = B_{4n\sqrt{V}}(x^*)$, where
    \begin{gather*}
        A = \ln\frac{48n(K+1)}{\beta},\quad V \geq \|x^0 - x^*\|^2 + \frac{409600\gamma^2A^2}{n^2}\sum_{i=1}^n\|F_i(x^*)\|^2,
    \end{gather*}
    and
    \begin{align}
        \gamma & \leq \min\left\{\frac{1}{1920 L A}, \frac{60^{(2-p)/p}\sqrt{V}n^{(p-1)/p}}{97200^{1/p}(K+1)^{1/p}\sigma A^{(p-1)/p}}\right\}, \label{ch3:eq:gamma_SEG_monotone_chapter_3}\\
        \lambda_k \equiv \lambda & = \frac{n\sqrt{V}}{60\gamma A}, \label{ch3:eq:lambda_SEG_monotone_chapter_3}\\
        \nu & = 0 \label{ch3:eq:nu_SEG_monotone_chapter_3}
    \end{align}
    for some $K \geq 1$ and $\beta \in (0,1]$. Then, after $K$ iterations of \algname{DProx-clipped-SEG-shift}, the following inequality holds with probability at least $1 - \beta$:
    \begin{align}
        & \gap_{\sqrt{V}}(\thesischthreetx_{\avg}^K) \leq \frac{9V}{2\gamma(K+1)} \quad \text{and} \notag\\
        & \{x^k\}_{k=0}^{K+1} \subseteq B_{3\sqrt{V}}(x^*), \{\thesischthreetx^k\}_{k=0}^{K+1} \subseteq B_{4n\sqrt{V}}(x^*), \label{ch3:eq:main_result_SEG_monotone_chapter_3}
    \end{align}
    where $\thesischthreetx_{\avg}^K$ is defined in \eqref{ch3:eq:tx_avg_SEG_chapter_3}. In particular, when $\gamma$ equals the minimum from \eqref{ch3:eq:gamma_SEG_monotone_chapter_3}, then after $K$ iterations of \algname{DProx-clipped-SEG-shift}, we have with probability at least $1-\beta$
    \begin{equation}
        \gap_{\sqrt{V}}(\thesischthreetx_{\avg}^K) = \cO\left(\max\left\{\frac{LV\ln(nK/\beta)}{K}, \frac{\sigma \sqrt{V} \ln^{(p-1)/p}(nK/\beta)}{n^{(p-1)/p}K^{(p-1)/p}}\right\}\right), \label{ch3:eq:clipped_SEG_monotone_rate_appendix_chapter_3}
    \end{equation}
    i.e., to achieve $\gap_{\sqrt{V}}(\thesischthreetx_{\avg}^K) \leq \varepsilon$ with probability at least $1 - \beta$ \algname{DProx-clipped-SEG-shift} needs
    \begin{align}
        K = \cO\biggl(\max\biggl\{ & \frac{LV}{\varepsilon}\ln\frac{nLV}{\varepsilon\beta}, \frac{1}{n}\left(\frac{\sigma \sqrt{V}}{\varepsilon}\right)^{p/(p-1)}\ln \frac{\sigma \sqrt{V}}{\varepsilon\beta}\biggr\}\biggr) \label{ch3:eq:clipped_SEG_monotone_complexity_appendix_chapter_3}
    \end{align}
    iterations/oracle calls per worker.
\end{theorem}
\begin{proof}
    The key idea behind the proof is similar to the one used by \textcite{gorbunov2022clipped, sadiev2023high}: we prove by induction that the iterates do not leave some ball and the sums decrease as $1/(K+1)$. To formulate the statement rigorously, we introduce probability event $E_k$ for each $k = 0,1,\ldots, K+1$ as follows: inequalities
    \begin{gather}
        \underbrace{\begin{aligned}
            &\max\limits_{u \in B_{\sqrt{V}}(x^*)}\biggl\{ \|x^0 - u\|^2 + 2\gamma \sum\limits_{l = 0}^{t-1} \langle x^l - u, \theta_l \rangle \\
            &\qquad\qquad + \gamma^2 \sum\limits_{l=0}^{t-1}\left(8\|\theta_l\|^2 + 6\|\omega_l\|^2\right)\biggr\}
        \end{aligned}}_{A_t} \leq 9V, \label{ch3:eq:induction_inequality_1_SEG_chapter_3}\\
        \left\|\gamma\sum\limits_{l=0}^{t-1}\theta_l\right\| \leq \sqrt{V}, \label{ch3:eq:induction_inequality_2_SEG_chapter_3}\\
        \left\|\frac{\gamma}{n}\sum\limits_{i=1}^{r-1}\theta_{i,t-1}^u\right\| \leq \frac{\sqrt{V}}{2}, \quad \left\|\frac{\gamma}{n}\sum\limits_{i=1}^{r-1}\omega_{i,t-1}^u\right\| \leq \frac{\sqrt{V}}{2}
    \end{gather}
    hold for $t = 0,1,\ldots,k$ and $r = 1, 2, \ldots, n$ simultaneously, where
    \begin{gather}
        \theta_l = \theta_l^u + \theta_l^b,\quad  \omega_l = \omega_l^u + \omega_l^b, \label{ch3:eq:gap_thm_SEG_technical_4_full_theta_omega_chapter_3}\\
        \theta_l^u \eqdef \frac{1}{n}\sum\limits_{i=1}^n \theta_{i,l}^u,\quad \theta_l^b \eqdef \frac{1}{n}\sum\limits_{i=1}^n \theta_{i,l}^b,\quad \omega_l^u \eqdef \frac{1}{n}\sum\limits_{i=1}^n \omega_{i,l}^u,\quad \omega_l^b \eqdef \frac{1}{n}\sum\limits_{i=1}^n \omega_{i,l}^b, \\
        \theta_{i,l}^u \eqdef \Exp_{\xi_{2,i}^l}\left[\thesischthreehg_{i}^l\right] - \thesischthreehg_{i}^l,\quad \theta_{i,l}^b \eqdef F_i(\thesischthreetx^l) - \Exp_{\xi_{2,i}^l}\left[\thesischthreehg_{i}^l\right] \quad \forall\; i\in [n], \label{ch3:eq:gap_thm_SEG_technical_4_chapter_3}\\
        \omega_{i,l}^u \eqdef \Exp_{\xi_{1,i}^l}\left[\thesischthreetg_{i}^l\right] - \thesischthreetg_{i}^l,\quad \omega_{i,l}^b \eqdef F_i(x^l) - \Exp_{\xi_{1,i}^l}\left[\thesischthreetg_{i}^l\right] \quad \forall\; i\in [n]. \label{ch3:eq:gap_thm_SEG_technical_5_chapter_3}
    \end{gather}

    We will prove by induction that $\Prob\{E_k\} \geq  1 - k\beta/(K+1)$ for all $k = 0,1,\ldots,K+1$. The base of induction follows immediately: for all $u \in B_{\sqrt{V}}(x^*)$ we have $\|x^0 - u\|^2 \leq 2\|x^0 - x^*\|^2 + 2\|x^* - u\|^2 \leq 4V < 9V$ and for $k = 0$ we have $\|\gamma\sum_{l=0}^{k-1} \theta_l\| = 0$, $\|\frac{\gamma}{n}\sum_{i=1}^{r-1} \theta_{i,k-1}^u\| = \|\frac{\gamma}{n}\sum_{i=1}^{r-1} \omega_{i,k-1}^u\| = 0$ since $\theta_{i,-1}^u = \omega_{i,-1}^u = 0$. Next, we assume that the statement holds for $k = T-1 \leq K$, i.e., $\Prob\{E_{T-1}\} \geq  1 - (T-1)\beta/(K+1)$. Let us show that it also holds for $k = T$, i.e., $\Prob\{E_{T}\} \geq  1 - T\beta/(K+1)$. 
    
    To proceed, we need to show that $E_{T-1}$ implies $\|x^t - x^*\| \leq 3\sqrt{V}$ for all $t = 0,1, \ldots, T-1$. We will use the induction argument as well. The base is already proven. Next, we assume that $\|x^t - x^*\| \leq 3\sqrt{V}$ for all $t = 0,1,\ldots,t'$ for some $t' < T-1$. Then 
    \begin{align}
        \|F(x^*)\| & = \sqrt{\|F(x^*)\|^2} \leq \sqrt{\sum\limits_{i=1}^n \|F_i(x^*)\|^2} \notag \\
        & \leq \frac{n\sqrt{V}}{160\gamma A} < \lambda \label{ch3:eq:njscnjdnsdcibv_chapter_3}
    \end{align}
    and for $t = 0,1,\ldots, t'$
    \begin{align}
        \|\thesischthreetx^t - x^*\| & = \|\thesischthreeprox_{\gamma \Psi}(x^t - \gamma \thesischthreetg^t) - \thesischthreeprox_{\gamma \Psi}(x^* - \gamma F(x^*))\| \notag\\
        & \leq \|x^t - x^* - \gamma(\thesischthreetg^t - F(x^*))\| \leq \|x^t - x^*\| + \gamma \|\thesischthreetg^t - F(x^*)\|\notag\\
        & \leq \|x^t - x^*\| + \gamma(\|\thesischthreetg^t\| + \|F(x^*)\|) \overset{\eqref{ch3:eq:lambda_SEG_monotone_chapter_3}, \eqref{ch3:eq:njscnjdnsdcibv_chapter_3}}{\leq} 3\sqrt{V} + 2\gamma\lambda \notag \\
        & \leq 3\sqrt{V} + \frac{n\sqrt{V}}{30A} \notag\\
        & \leq 4n\sqrt{V}. \label{ch3:eq:gap_thm_SEG_technical_1_chapter_3}
    \end{align}
    This means that $x^t \in B_{4n\sqrt{V}}(x^*)$ and $\thesischthreetx^t \in B_{4n\sqrt{V}}(x^*)$ for $t = 0,1,\ldots, t'$, and we can apply Lemma~\ref{ch3:lem:SEG_optimization_lemma_monotone_chapter_3}: $E_{T-1}$ implies
    \begin{align*}
        & \max_{B_{\sqrt{V}}(x^*)}\Big\{ 2\gamma(t'+1)\Big(\langle F(u), \thesischthreetx^{t'}_{\avg} - u\rangle + \Psi(\thesischthreetx^{t'}_{\avg}) - \Psi(u)\Big)\\
        &\qquad + \|x^{t'+1} - u\|^2  \Big\}\\
        &\leq \max_{B_{\sqrt{V}}(x^*)}\left\{\|x^0 - u\|^2 + 2\gamma \sum\limits_{l = 0}^{t'} \langle x^l - u, \theta_l \rangle \right\}\\
        &\quad + \gamma^2 \sum\limits_{l=0}^{t'}\left(8\|\theta_l\|^2 + 6\|\omega_l\|^2\right)\\
        &\overset{\eqref{ch3:eq:induction_inequality_1_SEG_chapter_3}}{\leq} 9V
    \end{align*}
    that gives
    \begin{align*}
        \|x^{t'+1} - x^*\|^2 & \leq \max_{B_{\sqrt{V}}(x^*)}\Big\{ 2\gamma(t'+1)\Big(\langle F(u), \thesischthreetx^{t'}_{\avg} - u\rangle + \Psi(\thesischthreetx^{t'}_{\avg}) - \Psi(u)\Big)\\
        &\qquad + \|x^{t'+1} - u\|^2  \Big\}\\
        & \leq 9V,
    \end{align*}
    where we take $u = x^*$ and use $\langle F(x^*), x - x^*\rangle + \Psi(x) - \Psi(x^*) \geq 0$ for all $x$. That is, we showed that $E_{T-1}$ implies $\|x^t - x^*\| \leq 3\sqrt{V}$, $\|\thesischthreetx^t - x^*\| \leq 4n\sqrt{V}$ and
    \begin{equation}
        \max_{B_{\sqrt{V}}(x^*)}\Big\{ 2\gamma(t+1)\Big(\langle F(u), \thesischthreetx^{t}_{\avg} - u\rangle + \Psi(\thesischthreetx^{t}_{\avg}) - \Psi(u)\Big) + \|x^{t+1} - u\|^2  \Big\} \leq 9V \label{ch3:eq:gap_thm_SEG_technical_1_5_chapter_3}
    \end{equation}
    for all $t = 0, 1, \ldots, T-2$. Before we proceed, we introduce a new notation:
    \begin{equation*}
        u^{t} = \begin{cases}x^t - x^*,& \text{if } \|x^t - x^*\| \leq 3\sqrt{V},\\ 0,& \text{otherwise,} \end{cases}
    \end{equation*}
    for all $t = 0, 1, \ldots, T-1$. Random vectors $\{u^{t}\}_{t=0}^{T-1}$ are bounded almost surely:
     \begin{equation}
        \|u^{t}\| \leq 3\sqrt{V}  \label{ch3:eq:gap_thm_SEG_technical_3_chapter_3}
    \end{equation}
    for all $t = 0, 1, \ldots, T-1$. In addition, $u^{t} = x^t - x^*$ follows from $E_{T-1}$ for all $t = 0, 1, \ldots, T-1$ and, thus, $E_{T-1}$ implies
    \begin{align}
        A_T & \overset{\eqref{ch3:eq:induction_inequality_1_SEG_chapter_3}}{=} \max\limits_{u \in B_{\sqrt{V}}(x^*)}\left\{ \|x^0 - u\|^2 + 2\gamma \sum\limits_{l = 0}^{T-1} \langle x^* - u, \theta_l \rangle \right\} \notag\\
        &\quad + 2\gamma \sum\limits_{l = 0}^{T-1} \langle x^l - x^*, \theta_l \rangle + \gamma^2 \sum\limits_{l=0}^{T-1}\left(8\|\theta_l\|^2 + 6\|\omega_l\|^2\right)\notag\\
        & \leq 4V + 2\gamma\max\limits_{u \in B_{\sqrt{V}}(x^*)}\left\{\left\langle x^* - u, \sum\limits_{l = 0}^{T-1}\theta_l \right\rangle \right\} \notag\\
        &\quad + 2\gamma \sum\limits_{l = 0}^{T-1} \langle u^{l}, \theta_l \rangle + \gamma^2 \sum\limits_{l=0}^{T-1}\left(8\|\theta_l\|^2 + 6\|\omega_l\|^2\right)\notag\\
        & = 4V + 2\gamma \sqrt{V} \left\|\sum\limits_{l = 0}^{T-1}\theta_l\right\| + 2\gamma \sum\limits_{l = 0}^{T-1} \langle u^{l}, \theta_l \rangle \notag\\
        &\quad + \gamma^2 \sum\limits_{l=0}^{T-1}\left(8\|\theta_l\|^2 + 6\|\omega_l\|^2\right). \notag
    \end{align}

    Using the notation from \eqref{ch3:eq:gap_thm_SEG_technical_4_full_theta_omega_chapter_3}-\eqref{ch3:eq:gap_thm_SEG_technical_5_chapter_3}, we can rewrite $\|\theta_l\|^2$ as
    \begin{align}
        \|\theta_l\|^2 & \leq  2\|\theta_l^u\|^2 + 2\|\theta_l^b\|^2 \notag \\
        & = \frac{2}{n^2}\left\|\sum\limits_{i=1}^n \theta_{i,l}^u\right\|^2 + 2\|\theta_l^b\|^2 \notag \\
        & = \frac{2}{n^2}\sum\limits_{i=1}^n\|\theta_{i,l}^u\|^2 + \frac{4}{n^2}\sum\limits_{j=2}^n\left\langle \sum\limits_{i=1}^{j-1} \theta_{i,l}^u, \theta_{j,l}^u \right\rangle + 2\|\theta_l^b\|^2 \label{ch3:eq:SEG_sum_theta_representation_chapter_3}
    \end{align}
    and, similarly,
    it holds for $\|\omega_l\|^2$.
    Putting all together, we obtain that $E_{T-1}$ implies
    \begin{align}
        A_T & \leq4V + 2\gamma \sqrt{V} \left\|\sum\limits_{l = 0}^{T-1}\theta_l\right\| + \underbrace{\frac{2\gamma}{n} \sum\limits_{l = 0}^{T-1}\sum\limits_{i=1}^n \langle u^{l}, \theta_{i,l}^u \rangle}_{\circledOne} + \underbrace{2\gamma \sum\limits_{l = 0}^{T-1} \langle u^{l}, \theta_l^b \rangle}_{\circledTwo}\notag\\
        &\quad\quad + \underbrace{\frac{2\gamma^2}{n^2} \sum\limits_{l=0}^{T-1} \sum\limits_{i=1}^n\left(8\Exp_{\xi_{2,i}^l}\left[\|\theta_{i,l}^u\|^2\right] + 6\Exp_{\xi_{1,i}^l}\left[\|\omega_{i,l}^u\|^2\right]\right)}_{\circledThree} \notag\\
        &\quad\quad + \underbrace{\frac{2\gamma^2}{n^2} \sum\limits_{l=0}^{T-1} \sum\limits_{i=1}^n\left(8\|\theta_{i,l}^u\|^2 + 6\|\omega_{i,l}^u\|^2 - 8\Exp_{\xi_{2,i}^l}\left[\|\theta_{i,l}^u\|^2\right] - 6\Exp_{\xi_{1,i}^l}\left[\|\omega_{i,l}^u\|^2\right]\right)}_{\circledFour}\notag\\
        &\quad\quad + \underbrace{2\gamma^2 \sum\limits_{l=0}^{T-1}\left(8\|\theta_l^b\|^2 + 6\|\omega_l^b\|^2\right)}_{\circledFive} + \underbrace{\frac{32\gamma^2}{n^2} \sum\limits_{l=0}^{T-1} \sum\limits_{j=2}^n \left\langle \sum\limits_{i=1}^{j-1} \theta_{i,l}^u, \theta_{j,l}^u \right\rangle}_{\circledSix}\notag\\
        &\quad\quad + \underbrace{\frac{24\gamma^2}{n^2} \sum\limits_{l=0}^{T-1} \sum\limits_{j=2}^n \left\langle \sum\limits_{i=1}^{j-1} \omega_{i,l}^u, \omega_{j,l}^u \right\rangle}_{\circledSeven}.\label{ch3:eq:gap_thm_SEG_technical_6_chapter_3}
    \end{align}
    To finish the proof, it remains to estimate $2\gamma \sqrt{V} \left\|\sum_{l = 0}^{T-1}\theta_l\right\|, \circledOne, \circledTwo, \circledThree, \circledFour, \circledFive, \circledSix, \circledSeven$ with high probability. More precisely, the goal is to prove that $2\gamma \sqrt{V} \left\|\sum_{l = 0}^{T-1}\theta_l\right\| + \circledOne + \circledTwo + \circledThree + \circledFour + \circledFive + \circledSix + \circledSeven \leq 5V$ with high probability. Before we proceed, we need to derive several useful inequalities related to $\theta_{i,l}^u, \omega_{i,l}^u, \theta_{l}^b, \omega_{l}^b$. First of all, we have
    \begin{equation}
        \|\theta_{i,l}^u\| \leq 2\lambda,\quad \|\omega_{i,l}^u\| \leq 2\lambda \label{ch3:eq:theta_omega_magnitude_chapter_3}
    \end{equation}
    by definition of the clipping operator. Next, probability event $E_{T-1}$ implies
    \begin{align*}
        \|F_i(x^l)\| & \leq \|F_i(x^l) - F_i(x^*)\| + \|F_i(x^*)\| \\
        & \overset{\eqref{ch3:eq:L_Lip_chapter_3}}{\leq} L\|x^l - x^*\| + \sqrt{\sum\limits_{i=1}^n \|F_i(x^*)\|^2} \\
        & \leq 3L\sqrt{V} + \frac{n\sqrt{V}}{160\gamma A} \overset{\eqref{ch3:eq:gamma_SEG_monotone_chapter_3}}{\leq} \frac{n\sqrt{V}}{120\gamma A} \overset{\eqref{ch3:eq:lambda_SEG_monotone_chapter_3}}{=} \frac{\lambda}{2},\\
        \|F_i(\thesischthreetx^l)\| & \leq \|F_i(\thesischthreetx^l) - F_i(x^*)\| + \|F_i(x^*)\| \\
        & \overset{\eqref{ch3:eq:L_Lip_chapter_3}}{\leq} L\|\thesischthreetx^l - x^*\| + \sqrt{\sum\limits_{i=1}^n \|F_i(x^*)\|^2} \\
        & \leq 4Ln\sqrt{V} + \frac{n\sqrt{V}}{160\gamma A} \overset{\eqref{ch3:eq:gamma_SEG_monotone_chapter_3}}{\leq} \frac{n\sqrt{V}}{120\gamma A} \overset{\eqref{ch3:eq:lambda_SEG_monotone_chapter_3}}{=} \frac{\lambda}{2}
    \end{align*}
    for $l = 0,1,\ldots,T-1$ and $i\in [n]$. Therefore, Lemma~\ref{ch3:lem:bias_and_variance_clip_chapter_3} and $E_{T-1}$ imply
    \begin{gather}
        \left\|\theta_l^b\right\| \leq \frac{1}{n}\sum\limits_{i=1}^n \|\theta_{i,l}^b\| \leq \frac{2^{p}\sigma^{p}}{\lambda^{{p}-1}},\quad \left\|\omega_l^b\right\| \leq \frac{1}{n}\sum\limits_{i=1}^n \|\omega_{i,l}^b\| \leq \frac{2^{p}\sigma^{p}}{\lambda^{{p}-1}}, \label{ch3:eq:bias_theta_omega_chapter_3}\\
        \Exp_{\xi_{2,i}^l}\left[\left\|\theta_{i,l}^u\right\|^2\right] \leq 18 \lambda^{2-{p}}\sigma^{p},\quad \Exp_{\xi_{1,i}^l}\left[\left\|\omega_{i,l}^u\right\|^2\right] \leq 18 \lambda^{2-{p}}\sigma^{p}, \label{ch3:eq:variance_theta_omega_chapter_3}
    \end{gather}
    for all $l = 0,1, \ldots, T-1$ and $i\in [n]$.

    \paragraph{Upper bound for $\circledOne$.} To estimate this sum, we will use Bernstein's inequality. The summands have conditional expectations equal to zero:
    \begin{equation*}
        \Exp_{\xi_{2,i}^l}\left[\frac{2\gamma}{n} \langle u^{l}, \theta_{i,l}^u \rangle\right] = \frac{2\gamma}{n} \left\langle u^{l}, \Exp_{\xi_{2,i}^l}[\theta_{i,l}^u] \right\rangle = 0.
    \end{equation*}
	    Moreover, for all $l = 0,\ldots, T-1$ random vectors $\{\theta_{i,l}^u\}_{i=1}^n$ are independent. Thus, sequence $\left\{\frac{2\gamma}{n} \langle u^{l}, \theta_{i,l}^u \rangle\right\}_{l,i = 0,1}^{T-1,n}$ is a martingale difference sequence. Next, the summands are bounded:
	    \begin{align}
	        \left|\frac{2\gamma}{n}
	        \langle u^{l}, \theta_{i,l}^u \rangle \right|
	        &\leq \frac{2\gamma}{n} \|u^{l}\|
	        \cdot \|\theta_{i,l}^u\|\notag\\
	        &\overset{\eqref{ch3:eq:gap_thm_SEG_technical_3_chapter_3},
	        \eqref{ch3:eq:theta_omega_magnitude_chapter_3}}{\leq}
	        \frac{12\gamma}{n} \sqrt{V}\lambda\notag\\
	        &\overset{\eqref{ch3:eq:lambda_SEG_monotone_chapter_3}}{\leq}
	        \frac{3V}{10A} \eqdef c.
	        \label{ch3:eq:gap_thm_SEG_technical_6_5_chapter_3}
	    \end{align}
    Finally, conditional variances $\sigma_{i,l}^2 \eqdef \Exp_{\xi_{2,i}^l}\left[\frac{4\gamma^2}{n^2} \langle u^{l}, \theta_{i,l}^u \rangle^2\right]$ of the summands are bounded:
    \begin{align}
\sigma_{i,l}^2 & \leq \Exp_{\xi_{2,i}^l}\left[\frac{4\gamma^2 }{n^2}\|u^{l}\|^2\cdot \|\theta_{i,l}^u\|^2\right] \notag \\
& \overset{\eqref{ch3:eq:gap_thm_SEG_technical_3_chapter_3}}{\leq} \frac{36\gamma^2 V}{n^2} \Exp_{\xi_{2,i}^l}\left[\|\theta_{i,l}^u\|^2\right]. \label{ch3:eq:gap_thm_SEG_technical_7_chapter_3}
    \end{align}
    Applying Bernstein's inequality (Lemma~\ref{ch3:lem:Bernstein_ineq_chapter_3}) with $X_{i,l} = \frac{2\gamma}{n} \langle u^{l}, \theta_{i,l}^u \rangle$, constant $c$ defined in \eqref{ch3:eq:gap_thm_SEG_technical_6_5_chapter_3}, $b = \frac{3V}{10}$, $H = \frac{3V^2}{200A}$, we get
    \begin{align*}
\MoveEqLeft[1] \Prob\left\{|\circledOne| > \frac{3V}{10} \text{ and } \sum\limits_{l=0}^{T-1}\sum\limits_{i=1}^n\sigma_{i,l}^2 \leq \frac{3V^2}{200A}\right\} \\
& \leq 2\exp\left(- \frac{b^2}{2H + 2cb/3}\right) = \frac{\beta}{24n(K+1)}.
    \end{align*}
    The above is equivalent to 
    \begin{align}
\Prob\{E_{\circledOne}\} & \geq 1 - \frac{\beta}{24n(K+1)},\; \text{ for }\; E_{\circledOne} \notag \\
& = \left\{\text{either } \sum\limits_{l=0}^{T-1}\sum\limits_{i=1}^n\sigma_{i,l}^2 > \frac{3V^2}{200A} \text{ or } |\circledOne| \leq \frac{3V}{10}\right\}. \label{ch3:eq:bound_1_gap_SEG_chapter_3}
    \end{align}
    Moreover, $E_{T-1}$ implies
    \begin{align}
        \sum\limits_{l=0}^{T-1}\sum\limits_{i=1}^n\sigma_{i,l}^2 & \overset{\eqref{ch3:eq:gap_thm_SEG_technical_7_chapter_3}}{\leq} \frac{36\gamma^2 V}{n^2} \sum\limits_{l=0}^{T-1} \sum\limits_{i=1}^n \Exp_{\xi_{2,i}^l}\left[\|\theta_{i,l}^u\|^2\right] \notag \\
        & \overset{\eqref{ch3:eq:variance_theta_omega_chapter_3}}{\leq} \frac{648(K+1)\gamma^2 V \lambda^{2-{p}} \sigma^{p}}{n} \notag \\
        & \overset{\eqref{ch3:eq:lambda_SEG_monotone_chapter_3}}{\leq} \frac{648 (K+1)\gamma^{{p}}\sigma^{p} V^{2-\frac{{p}}{2}}}{60^{2-{p}} n^{{p}-1} A^{2-{p}}} \overset{\eqref{ch3:eq:gamma_SEG_monotone_chapter_3}}{\leq} \frac{3V^2}{200A}. \label{ch3:eq:bound_1_variances_gap_SEG_chapter_3}
    \end{align}

    \paragraph{Upper bound for $\circledTwo$.} Probability event $E_{T-1}$ implies
    \begin{align}
        \circledTwo & \leq 2\gamma \sum\limits_{l=0}^{T-1}\|u^{l}\| \cdot \|\theta_l^b\| \overset{\eqref{ch3:eq:gap_thm_SEG_technical_3_chapter_3}, \eqref{ch3:eq:bias_theta_omega_chapter_3}}{\leq} \frac{6\cdot 2^{p}(K+1)\gamma \sqrt{V}\sigma^{p}}{\lambda^{{p}-1}} \notag \\
        & \overset{\eqref{ch3:eq:lambda_SEG_monotone_chapter_3}}{=} \frac{6\cdot 2^{p}\cdot 60^{{p}-1} (K+1)\gamma^{p} \sigma^{p} A^{{p}-1}}{n^{{p}-1}V^{\frac{{p}}{2}-1}} \overset{\eqref{ch3:eq:gamma_SEG_monotone_chapter_3}}{\leq} \frac{3V}{100}. \label{ch3:eq:bound_2_variances_gap_SEG_chapter_3}
    \end{align}

    In these bounds and in the bounds below we also use $T \leq K+1$.

    \paragraph{Upper bound for $\circledThree$.} Probability event $E_{T-1}$ implies
    \begin{align}
        \MoveEqLeft[1] \frac{16\gamma^2}{n^2} \sum\limits_{l=0}^{T-1}\sum\limits_{i=1}^n\Exp_{\xi_{2,i}^l}[\|\theta_{i,l}^u\|^2] \notag \\
        & \overset{\eqref{ch3:eq:variance_theta_omega_chapter_3}}{\leq} \frac{288\gamma^2 (K+1)\lambda^{2-{p}}\sigma^{p}}{n} \notag \\
        & \overset{\eqref{ch3:eq:lambda_SEG_monotone_chapter_3}}{=} \frac{288\gamma^{p} (K+1)\sigma^{p} V^{1-\frac{{p}}{2}}}{60^{2-{p}} n^{{p}-1} A^{2-{p}}} \notag\\
        & \overset{\eqref{ch3:eq:gamma_SEG_monotone_chapter_3}}{\leq} \frac{3}{100} V, \label{ch3:eq:sum_theta_squared_bound_gap_SEG_chapter_3}\\
        \MoveEqLeft[1] \frac{12\gamma^2}{n^2} \sum\limits_{l=0}^{T-1}\sum\limits_{i=1}^n\Exp_{\xi_{1,i}^l}[\|\omega_{i,l}^u\|^2] \notag \\
        & \overset{\eqref{ch3:eq:variance_theta_omega_chapter_3}}{\leq} \frac{216\gamma^2 (K+1)\lambda^{2-{p}}\sigma^{p}}{n} \notag \\
        & \overset{\eqref{ch3:eq:lambda_SEG_monotone_chapter_3}}{=} \frac{216\gamma^{p} (K+1)\sigma^{p} V^{1-\frac{{p}}{2}}}{60^{2-{p}} n^{{p}-1} A^{2-{p}}} \notag\\
        & \overset{\eqref{ch3:eq:gamma_SEG_monotone_chapter_3}}{\leq} \frac{1}{50} V, \label{ch3:eq:sum_omega_squared_bound_gap_SEG_chapter_3}\\
        \circledThree & \overset{\eqref{ch3:eq:sum_theta_squared_bound_gap_SEG_chapter_3}, \eqref{ch3:eq:sum_omega_squared_bound_gap_SEG_chapter_3}}{\leq} \frac{1}{20}V. \label{ch3:eq:bound_3_variances_gap_SEG_chapter_3}
    \end{align}

    \paragraph{Upper bound for $\circledFour$.} 
    To estimate this sum, we will use Bernstein's inequality. The summands
    \begin{align*}
        Z_{i,l} \eqdef \frac{2\gamma^2}{n^2}\Bigl( & 8\|\theta_{i,l}^u\|^2 + 6\|\omega_{i,l}^u\|^2 \\
        & - 8\Exp_{\xi_{2,i}^l}\left[\|\theta_{i,l}^u\|^2\right] - 6\Exp_{\xi_{1,i}^l}\left[\|\omega_{i,l}^u\|^2\right]\Bigr)
    \end{align*}
    have conditional expectations equal to zero: $\Exp_{\xi_{1,i}^l,\xi_{2,i}^l}[Z_{i,l}] = 0$. Moreover, for all $l = 0,\ldots, T-1$ random vectors $\{\theta_{i,l}^u\}_{i=1}^n$, $\{\omega_{i,l}^u\}_{i=1}^n$ are independent. Thus, sequence $\{Z_{i,l}\}_{l,i = 0,1}^{T-1,n}$ is a martingale difference sequence. Next, the summands are bounded:
    \begin{align}
        |Z_{i,l}| &\leq \frac{16\gamma^2}{n^2} \left(\|\theta_{i,l}^u\|^2 + \Exp_{\xi_{2,i}^l}\left[\|\theta_{i,l}^u\|^2\right]\right) \notag\\
        &\quad + \frac{12\gamma^2}{n^2} \left(\|\omega_{i,l}^u\|^2 + \Exp_{\xi_{1,i}^l}\left[\|\omega_{i,l}^u\|^2\right]\right) \notag\\
        &\overset{\eqref{ch3:eq:theta_omega_magnitude_chapter_3}}{\leq} \frac{224\gamma^2 \lambda^2}{n^2} \overset{\eqref{ch3:eq:lambda_SEG_monotone_chapter_3}}{\leq} \frac{V}{6A} \eqdef c.\label{ch3:eq:gap_thm_SEG_technical_6_5_1_chapter_3}
    \end{align}
    Finally, conditional variances $\widetilde\sigma_{i,l}^2 \eqdef \Exp_{\xi_{1,i}^l, \xi_{2,i}^l}[Z_{i,l}^2]$ of the summands are bounded:
    \begin{align}
        \widetilde\sigma_{i,l}^2
        &\overset{\eqref{ch3:eq:gap_thm_SEG_technical_6_5_1_chapter_3}}{\leq}
        \frac{V}{6A}\Exp_{\xi_{1,i}^l, \xi_{2,i}^l}\left[|Z_{i,l}|\right] \notag\\
        &\leq
        \frac{4\gamma^2V}{3n^2A}
        \Exp_{\xi_{1,i}^l, \xi_{2,i}^l}\left[
        4\|\theta_{i,l}^u\|^2 + 3\|\omega_{i,l}^u\|^2 \right].
        \label{ch3:eq:gap_thm_SEG_technical_7_1_chapter_3}
    \end{align}
    Applying Bernstein's inequality (Lemma~\ref{ch3:lem:Bernstein_ineq_chapter_3}) with $X_{i,l} = Z_{i,l}$, $b = V/6$, $H = V^2/(216A)$ and constant $c$ defined in \eqref{ch3:eq:gap_thm_SEG_technical_6_5_1_chapter_3}, we get
    \begin{align*}
        &\Prob\left\{|\circledFour| > \frac{V}{6} \text{ and } \sum\limits_{l=0}^{T-1}\sum\limits_{i=1}^n \widetilde\sigma_{i,l}^2 \leq \frac{V^2}{216A}\right\}\\
        &\qquad \leq
        2\exp\left(- \frac{b^2}{2H + 2cb/3}\right)
        = \frac{\beta}{24n(K+1)}.
    \end{align*}
    The above is equivalent to 
    \begin{equation}
        \begin{aligned}
        \Prob\{E_{\circledFour}\}
        &\geq 1 - \frac{\beta}{24n(K+1)},\quad \text{for}\\
        E_{\circledFour}
        &= \left\{\text{either } \sum\limits_{l=0}^{T-1}\sum\limits_{i=1}^n \widetilde\sigma_{i,l}^2 > \frac{V^2}{216A} \text{ or } |\circledFour| \leq \frac{V}{6}\right\}.
        \end{aligned}
        \label{ch3:eq:bound_1_gap_SEG_1_chapter_3}
    \end{equation}
	    Moreover, $E_{T-1}$ implies
	    \begin{align}
	        \sum\limits_{l=0}^{T-1}\sum\limits_{i=1}^n\widetilde\sigma_{i,l}^2
	        &\overset{\eqref{ch3:eq:gap_thm_SEG_technical_7_1_chapter_3}}{\leq}
	        \frac{4\gamma^2V}{3n^2A}
	        \sum\limits_{l=0}^{T-1}\sum\limits_{i=1}^n
	        \Exp_{\xi_{1,i}^l, \xi_{2,i}^l}
	        \left[4\|\theta_{i,l}^u\|^2 + 3\|\omega_{i,l}^u\|^2 \right]\notag\\
	        &\overset{\eqref{ch3:eq:variance_theta_omega_chapter_3}}{\leq}
	        \frac{168(K+1)\gamma^2 V \lambda^{2-{p}} \sigma^{p}}
	        {nA}\notag\\
	        &\overset{\eqref{ch3:eq:lambda_SEG_monotone_chapter_3}}{\leq}
	        \frac{168(K+1)\gamma^{p} V^{2-\frac{{p}}{2}} \sigma^{p}}
	        {60^{2-{p}}n^{{p}-1}A^{3-{p}}}\notag\\
	        &\overset{\eqref{ch3:eq:gamma_SEG_monotone_chapter_3}}{\leq}
	        \frac{V^2}{216A}.
	        \label{ch3:eq:bound_1_variances_gap_SEG_1_chapter_3}
	    \end{align}

    \paragraph{Upper bound for $\circledFive$.} Probability event $E_{T-1}$ implies
    \begin{align}
        \circledFive & = 2\gamma^2 \sum\limits_{l=0}^{T-1}\left(8\|\theta_l^b\|^2 + 6\|\omega_l^b\|^2\right) \notag \\
        & \overset{\eqref{ch3:eq:bias_theta_omega_chapter_3}}{\leq} \frac{28\cdot 2^{2{p}}\gamma^2\sigma^{2{p}} (K+1)}{\lambda^{2{p}-2}} \notag\\
        & \overset{\eqref{ch3:eq:lambda_SEG_monotone_chapter_3}}{=}  \frac{28\cdot 2^{2{p}}\cdot 60^{2{p} - 2}\gamma^{2{p}}\sigma^{2{p}} (K+1) A^{2{p}-2}}{n^{2{p}-2}V^{{p}-1}} \overset{\eqref{ch3:eq:gamma_SEG_monotone_chapter_3}}{\leq} \frac{V}{6}. \label{ch3:eq:nsjcbdjhcbfjdhfbj_chapter_3}
    \end{align}

    \paragraph{Upper bounds for $\circledSix$ and $\circledSeven$.} These sums require more refined analysis. We introduce new vectors:
    \begin{equation}
    \begin{aligned}
        \zeta_{j}^l &= \begin{cases}
            \frac{\gamma}{n}\sum\limits_{i=1}^{j-1} \theta_{i,l}^u,& \text{if } \left\|\frac{\gamma}{n}\sum\limits_{i=1}^{j-1} \theta_{i,l}^u\right\| \leq \frac{\sqrt{V}}{2},\\
            0,& \text{otherwise,}
        \end{cases}, \\
        \delta_{j}^l &= \begin{cases}
            \frac{\gamma}{n}\sum\limits_{i=1}^{j-1} \omega_{i,l}^u,& \text{if } \left\|\frac{\gamma}{n}\sum\limits_{i=1}^{j-1} \omega_{i,l}^u\right\| \leq \frac{\sqrt{V}}{2},\\
            0,& \text{otherwise,}
        \end{cases}
    \end{aligned}\label{ch3:eq:SEG_zeta_delta_def_chapter_3}
    \end{equation}
    for all $j \in [n]$ and $l = 0,\ldots, T-1$. Then, by definition
    \begin{equation}
        \|\zeta_j^l\| \leq \frac{\sqrt{V}}{2},\quad \|\delta_j^l\| \leq \frac{\sqrt{V}}{2} \label{ch3:eq:SEG_zeta_delta_bound_chapter_3}
    \end{equation}
    and
    \begin{align}
        \circledSix & = \underbrace{\frac{32\gamma}{n} \sum\limits_{l=0}^{T-1} \sum\limits_{j=2}^n \left\langle \zeta_j^l, \theta_{j,l}^u \right\rangle}_{\circledSix'} + \frac{32\gamma}{n} \sum\limits_{l=0}^{T-1} \sum\limits_{j=2}^n \left\langle \frac{\gamma}{n}\sum\limits_{i=1}^{j-1} \theta_{i,l}^u - \zeta_j^l, \theta_{j,l}^u \right\rangle, \label{ch3:eq:SEG_extra_sums_six_distributed_chapter_3} \\
        \circledSeven & = \underbrace{\frac{24\gamma}{n} \sum\limits_{l=0}^{T-1} \sum\limits_{j=2}^n \left\langle \delta_j^l, \omega_{j,l}^u \right\rangle}_{\circledSeven'} + \frac{24\gamma}{n} \sum\limits_{l=0}^{T-1} \sum\limits_{j=2}^n \left\langle \frac{\gamma}{n}\sum\limits_{i=1}^{j-1} \omega_{i,l}^u - \delta_j^l, \omega_{j,l}^u \right\rangle. \label{ch3:eq:SEG_extra_sums_seven_distributed_chapter_3}
    \end{align}
    We also note here that $E_{T-1}$ implies
    \begin{align}
        \MoveEqLeft \frac{32\gamma}{n} \sum\limits_{l=0}^{T-1} \sum\limits_{j=2}^n \left\langle \frac{\gamma}{n}\sum\limits_{i=1}^{j-1} \theta_{i,l}^u - \zeta_j^l, \theta_{j,l}^u \right\rangle \notag\\
        & = \frac{32\gamma}{n}  \sum\limits_{j=2}^n \left\langle \frac{\gamma}{n}\sum\limits_{i=1}^{j-1} \theta_{i,T-1}^u - \zeta_j^{T-1}, \theta_{j,T-1}^u \right\rangle, \label{ch3:eq:SEG_extra_sums_six_distributed_1_chapter_3}\\
        \MoveEqLeft \frac{24\gamma}{n} \sum\limits_{l=0}^{T-1} \sum\limits_{j=2}^n \left\langle \frac{\gamma}{n}\sum\limits_{i=1}^{j-1} \omega_{i,l}^u - \delta_j^l, \omega_{j,l}^u \right\rangle \notag\\
        & = \frac{24\gamma}{n}  \sum\limits_{j=2}^n \left\langle \frac{\gamma}{n}\sum\limits_{i=1}^{j-1} \omega_{i,T-1}^u - \delta_j^{T-1}, \omega_{j,T-1}^u \right\rangle. \label{ch3:eq:SEG_extra_sums_seven_distributed_1_chapter_3}
    \end{align}
    
    \paragraph{Upper bound for $\circledSix'$.} To estimate this sum, we will use Bernstein's inequality. The summands have conditional expectations equal to zero:
    \begin{equation*}
    \Exp_{\xi_{2,j}^l}\left[\frac{32\gamma}{n}\left\langle \zeta_j^l, \theta_{j,l}^u \right\rangle\right] = \frac{32\gamma}{n}\left\langle \zeta_j^l, \Exp_{\xi_{2,j}^l}[\theta_{j,l}^u] \right\rangle = 0.
    \end{equation*}
    Moreover, for all $l = 0,\ldots, T-1$ random vectors $\{\theta_{i,l}^u\}_{i=1}^n$ are independent. Thus, sequence $\left\{\frac{32\gamma}{n}\left\langle \zeta_j^l, \theta_{j,l}^u \right\rangle\right\}_{l,j = 0,2}^{T-1,n}$ is a martingale difference sequence. Next, the summands are bounded:
    \begin{align}
\left|\frac{32\gamma}{n}\left\langle \zeta_j^l, \theta_{j,l}^u \right\rangle \right| & \leq \frac{32\gamma}{n} \left\|\zeta_j^l\right\| \cdot \|\theta_{j,l}^u\| \notag \\
& \overset{\eqref{ch3:eq:SEG_zeta_delta_bound_chapter_3},\eqref{ch3:eq:theta_omega_magnitude_chapter_3}}{\leq} \frac{32\gamma}{n}\cdot \frac{\sqrt{V}}{2}\cdot 2\lambda \overset{\eqref{ch3:eq:lambda_SEG_monotone_chapter_3}}{\leq} \frac{4V}{5A} \eqdef c. \label{ch3:eq:gap_thm_SEG_technical_6_5_123_chapter_3}
    \end{align}
    Finally, conditional variances $\widehat\sigma_{j,l}^2 \eqdef \Exp_{\xi_{2,j}^l}\left[\frac{1024\gamma^2}{n^2} \langle \zeta_j^l, \theta_{j,l}^u \rangle^2\right]$ of the summands are bounded:
    \begin{align}
\widehat\sigma_{j,l}^2 & \leq \Exp_{\xi_{2,j}^l}\left[\frac{1024\gamma^2 }{n^2}\|\zeta_j^l\|^2\cdot \|\theta_{j,l}^u\|^2\right] \notag \\
& \overset{\eqref{ch3:eq:SEG_zeta_delta_bound_chapter_3}}{\leq} \frac{256\gamma^2 V}{n^2} \Exp_{\xi_{2,j}^l}\left[\|\theta_{j,l}^u\|^2\right]. \label{ch3:eq:gap_thm_SEG_technical_7_123_chapter_3}
    \end{align}
    Applying Bernstein's inequality (Lemma~\ref{ch3:lem:Bernstein_ineq_chapter_3}) with $X_{j,l} = \frac{32\gamma}{n}\left\langle \zeta_j^l, \theta_{j,l}^u \right\rangle$, constant $c$ defined in \eqref{ch3:eq:gap_thm_SEG_technical_6_5_123_chapter_3}, $b = \frac{4V}{5}$, $H = \frac{8V^2}{75A}$, we get
    \begin{align*}
\MoveEqLeft[1] \Prob\left\{|\circledSix'| > \frac{4V}{5} \text{ and } \sum\limits_{l=0}^{T-1}\sum\limits_{j=2}^n\widehat\sigma_{j,l}^2 \leq \frac{8V^2}{75A}\right\} \\
& \leq 2\exp\left(- \frac{b^2}{2H + 2cb/3}\right) = \frac{\beta}{24n(K+1)}.
    \end{align*}
    The above is equivalent to 
    \begin{align}
\Prob\{E_{\circledSix'}\} & \geq 1 - \frac{\beta}{24n(K+1)},\; \text{ for }\; E_{\circledSix'} \notag \\
& = \left\{\text{either } \sum\limits_{l=0}^{T-1}\sum\limits_{j=2}^n\widehat\sigma_{j,l}^2 > \frac{8V^2}{75A} \text{ or } |\circledSix'| \leq \frac{4V}{5}\right\}. \label{ch3:eq:bound_6_gap_SEG_chapter_3}
    \end{align}
    Moreover, $E_{T-1}$ implies
    \begin{align}
        \sum\limits_{l=0}^{T-1}\sum\limits_{j=2}^n\widehat\sigma_{j,l}^2 & \overset{\eqref{ch3:eq:gap_thm_SEG_technical_7_123_chapter_3}}{\leq} \frac{256\gamma^2 V}{n^2} \sum\limits_{l=0}^{T-1} \sum\limits_{j=2}^n \Exp_{\xi_{2,j}^l}\left[\|\theta_{j,l}^u\|^2\right] \notag \\
        & \overset{\eqref{ch3:eq:variance_theta_omega_chapter_3}}{\leq} \frac{4608(K+1)\gamma^2 V \lambda^{2-{p}} \sigma^{p}}{n} \notag \\
        & \overset{\eqref{ch3:eq:lambda_SEG_monotone_chapter_3}}{\leq} \frac{4608 (K+1)\gamma^{{p}}\sigma^{p} V^{2-\frac{{p}}{2}}}{60^{2-{p}} n^{{p}-1} A^{2-{p}}} \overset{\eqref{ch3:eq:gamma_SEG_monotone_chapter_3}}{\leq} \frac{8V^2}{75A}. \label{ch3:eq:bound_6_variances_gap_SEG_chapter_3}
    \end{align}

    \paragraph{Upper bound for $\circledSeven'$.} To estimate this sum, we will use Bernstein's inequality. The summands have conditional expectations equal to zero:
    \begin{equation*}
    \Exp_{\xi_{1,j}^l}\left[\frac{24\gamma}{n}\left\langle \delta_j^l, \omega_{j,l}^u \right\rangle\right] = \frac{24\gamma}{n}\left\langle \delta_j^l, \Exp_{\xi_{1,j}^l}[\omega_{j,l}^u] \right\rangle = 0.
    \end{equation*}
    Moreover, for all $l = 0,\ldots, T-1$ random vectors $\{\omega_{i,l}^u\}_{i=1}^n$ are independent. Thus, sequence $\left\{\frac{24\gamma}{n}\left\langle \delta_j^l, \omega_{j,l}^u \right\rangle\right\}_{l,j = 0,2}^{T-1,n}$ is a martingale difference sequence. Next, the summands are bounded:
    \begin{align}
\left|\frac{24\gamma}{n}\left\langle \delta_j^l, \omega_{j,l}^u \right\rangle \right| & \leq \frac{24\gamma}{n} \left\|\delta_j^l\right\| \cdot \|\omega_{j,l}^u\| \notag \\
& \overset{\eqref{ch3:eq:SEG_zeta_delta_bound_chapter_3},\eqref{ch3:eq:theta_omega_magnitude_chapter_3}}{\leq} \frac{24\gamma}{n}\cdot \frac{\sqrt{V}}{2}\cdot 2\lambda \overset{\eqref{ch3:eq:lambda_SEG_monotone_chapter_3}}{\leq} \frac{3V}{5A} \eqdef c. \label{ch3:eq:gap_thm_SEG_technical_6_5_1234_chapter_3}
    \end{align}
    Finally, conditional variances $(\sigma_{j,l}')^2 \eqdef \Exp_{\xi_{1,j}^l}\left[\frac{576\gamma^2}{n^2} \langle \delta_j^l, \omega_{j,l}^u \rangle^2\right]$ of the summands are bounded:
    \begin{align}
(\sigma_{j,l}')^2 & \leq \Exp_{\xi_{1,j}^l}\left[\frac{576\gamma^2 }{n^2}\|\delta_j^l\|^2\cdot \|\omega_{j,l}^u\|^2\right] \notag \\
& \overset{\eqref{ch3:eq:SEG_zeta_delta_bound_chapter_3}}{\leq} \frac{144\gamma^2 V}{n^2} \Exp_{\xi_{1,j}^l}\left[\|\omega_{j,l}^u\|^2\right]. \label{ch3:eq:gap_thm_SEG_technical_7_1234_chapter_3}
    \end{align}
    Applying Bernstein's inequality (Lemma~\ref{ch3:lem:Bernstein_ineq_chapter_3}) with $X_{j,l} = \frac{24\gamma}{n}\left\langle \delta_j^l, \omega_{j,l}^u \right\rangle$, constant $c$ defined in \eqref{ch3:eq:gap_thm_SEG_technical_6_5_1234_chapter_3}, $b = \frac{3V}{5}$, $H = \frac{3V^2}{50A}$, we get
    \begin{align*}
\MoveEqLeft[1] \Prob\left\{|\circledSeven'| > \frac{3V}{5} \text{ and } \sum\limits_{l=0}^{T-1}\sum\limits_{j=2}^n(\sigma_{j,l}')^2 \leq \frac{3V^2}{50A}\right\} \\
& \leq 2\exp\left(- \frac{b^2}{2H + 2cb/3}\right) = \frac{\beta}{24n(K+1)}.
    \end{align*}
    The above is equivalent to 
    \begin{align}
\Prob\{E_{\circledSeven'}\} & \geq 1 - \frac{\beta}{24n(K+1)},\; \text{ for }\; E_{\circledSeven'} \notag \\
& = \left\{\text{either } \sum\limits_{l=0}^{T-1}\sum\limits_{j=2}^n(\sigma_{j,l}')^2 > \frac{3V^2}{50A} \text{ or } |\circledSeven'| \leq \frac{3V}{5}\right\}. \label{ch3:eq:bound_7_gap_SEG_chapter_3}
    \end{align}
    Moreover, $E_{T-1}$ implies
    \begin{align}
        \sum\limits_{l=0}^{T-1}\sum\limits_{j=2}^n(\sigma_{j,l}')^2 & \overset{\eqref{ch3:eq:gap_thm_SEG_technical_7_1234_chapter_3}}{\leq} \frac{144\gamma^2 V}{n^2} \sum\limits_{l=0}^{T-1} \sum\limits_{j=2}^n \Exp_{\xi_{1,j}^l}\left[\|\omega_{j,l}^u\|^2\right] \notag \\
        & \overset{\eqref{ch3:eq:variance_theta_omega_chapter_3}}{\leq} \frac{2592(K+1)\gamma^2 V \lambda^{2-{p}} \sigma^{p}}{n} \notag \\
        & \overset{\eqref{ch3:eq:lambda_SEG_monotone_chapter_3}}{\leq} \frac{2592 (K+1)\gamma^{{p}}\sigma^{p} V^{2-\frac{{p}}{2}}}{60^{2-{p}} n^{{p}-1} A^{2-{p}}} \overset{\eqref{ch3:eq:gamma_SEG_monotone_chapter_3}}{\leq} \frac{3V^2}{50A}. \label{ch3:eq:bound_7_variances_gap_SEG_chapter_3}
    \end{align}

    \paragraph{Upper bound for $2\gamma \sqrt{V} \left\|\sum_{l=0}^{T-1} \theta_l\right\|$.} We introduce new random vectors:
    \begin{equation*}
        \zeta_l' = \begin{cases} \gamma \sum\limits_{r=0}^{l-1}\theta_r,& \text{if } \left\|\gamma \sum\limits_{r=0}^{l-1}\theta_r\right\| \leq \sqrt{V},\\ 0, & \text{otherwise} \end{cases}
    \end{equation*}
    for $l = 1, 2, \ldots, T-1$. With probability $1$ we have
    \begin{equation}
        \|\zeta_l'\| \leq \sqrt{V}.  \label{ch3:eq:gap_thm_SEG_technical_8_chapter_3}
    \end{equation}
   Using this and \eqref{ch3:eq:induction_inequality_2_SEG_chapter_3}, we obtain that $E_{T-1}$ implies
    \begin{align}
        2\gamma \sqrt{V}\left\|\sum\limits_{l = 0}^{T-1}\theta_l\right\|
        &= 2\sqrt{V}\sqrt{\gamma^2
        \left\|\sum\limits_{l = 0}^{T-1}\theta_l\right\|^2}\notag\\
        &= 2\sqrt{V}\sqrt{\gamma^2\sum\limits_{l=0}^{T-1}\|\theta_l\|^2
        + 2\gamma\sum\limits_{l=0}^{T-1}
        \left\langle \gamma\sum\limits_{r=0}^{l-1} \theta_r,
        \theta_l \right\rangle} \notag\\
        &= 2\sqrt{V} \sqrt{\gamma^2\sum\limits_{l=0}^{T-1}\|\theta_l\|^2
        + 2\gamma \sum\limits_{l=0}^{T-1}
        \langle \zeta_l', \theta_l\rangle} \notag\\
        &\overset{\eqref{ch3:eq:gap_thm_SEG_technical_6_chapter_3}}{\leq}
        2\sqrt{V} \sqrt{
        \begin{aligned}
        &\frac{\circledThree + \circledFour + \circledFive + \circledSix + \circledSeven}{8}\\
        &+ \underbrace{\frac{2\gamma}{n}
        \sum\limits_{l=0}^{T-1}\sum\limits_{i=1}^n
        \langle \zeta_l', \theta_{i,l}^u\rangle}_{\circledEight}\\
        &+ \underbrace{2\gamma \sum\limits_{l=0}^{T-1}
        \langle \zeta_l', \theta_l^b\rangle}_{\circledNine}
        \end{aligned}}.
        \label{ch3:eq:norm_sum_theta_bound_gap_SEG_chapter_3}
    \end{align}

    \paragraph{Upper bound for $\circledEight$.} To estimate this sum, we will use Bernstein's inequality. The summands have conditional expectations equal to zero:
    \begin{equation*}
        \Exp_{\xi_{2,i}^l}\left[\frac{2\gamma}{n} \langle \zeta_l', \theta_{i,l}^u \rangle\right] = \frac{2\gamma}{n} \left\langle \zeta_l', \Exp_{\xi_{2,i}^l}[\theta_{i,l}^u] \right\rangle = 0.
    \end{equation*}
    Moreover, for all $l = 0,\ldots, T-1$ random vectors $\{\theta_{i,l}^u\}_{i=1}^n$ are independent. Thus, sequence $\left\{\frac{2\gamma}{n} \langle \zeta_l', \theta_{i,l}^u \rangle\right\}_{l,i = 0,1}^{T-1,n}$ is a martingale difference sequence. Next, the summands are bounded:
    \begin{align}
        \left|\frac{2\gamma}{n}\langle \zeta_l', \theta_{i,l}^u \rangle \right|
        &\leq \frac{2\gamma}{n} \|\zeta_l'\|\cdot \|\theta_{i,l}^u\| \notag\\
        &\overset{\eqref{ch3:eq:gap_thm_SEG_technical_8_chapter_3},\eqref{ch3:eq:theta_omega_magnitude_chapter_3}}{\leq}
        \frac{4\gamma}{n} \sqrt{V}\lambda \notag\\
        &\overset{\eqref{ch3:eq:lambda_SEG_monotone_chapter_3}}{\leq}
        \frac{V}{10A} \eqdef c.
        \label{ch3:eq:gap_thm_SEG_technical_6_5_8_chapter_3}
    \end{align}
    Finally, conditional variances $(\widetilde\sigma_{i,l}')^2 \eqdef \Exp_{\xi_{2,i}^l}\left[\frac{4\gamma^2}{n^2} \langle \zeta_l', \theta_{i,l}^u \rangle^2\right]$ of the summands are bounded:
    \begin{align}
(\widetilde\sigma_{i,l}')^2 & \leq \Exp_{\xi_{2,i}^l}\left[\frac{4\gamma^2 }{n^2}\|\zeta_l'\|^2\cdot \|\theta_{i,l}^u\|^2\right] \notag \\
& \overset{\eqref{ch3:eq:gap_thm_SEG_technical_8_chapter_3}}{\leq} \frac{4\gamma^2 V}{n^2} \Exp_{\xi_{2,i}^l}\left[\|\theta_{i,l}^u\|^2\right]. \label{ch3:eq:gap_thm_SEG_technical_7_8_chapter_3}
    \end{align}
    Applying Bernstein's inequality (Lemma~\ref{ch3:lem:Bernstein_ineq_chapter_3}) with $X_{i,l} = \frac{2\gamma}{n} \langle \zeta_l', \theta_{i,l}^u \rangle$, constant $c$ defined in \eqref{ch3:eq:gap_thm_SEG_technical_6_5_8_chapter_3}, $b = \frac{V}{10}$, $H = \frac{V^2}{600A}$, we get
    \begin{align*}
\MoveEqLeft[1] \Prob\left\{|\circledEight| > \frac{V}{10} \text{ and } \sum\limits_{l=0}^{T-1}\sum\limits_{i=1}^n(\widetilde\sigma_{i,l}')^2 \leq \frac{V^2}{600A}\right\} \\
& \leq 2\exp\left(- \frac{b^2}{2H + 2cb/3}\right) = \frac{\beta}{24n(K+1)}.
    \end{align*}
    The above is equivalent to 
    \begin{align}
\Prob\{E_{\circledEight}\} & \geq 1 - \frac{\beta}{24n(K+1)},\; \text{ for }\; E_{\circledEight} \notag \\
& = \left\{\text{either } \sum\limits_{l=0}^{T-1}\sum\limits_{i=1}^n(\widetilde\sigma_{i,l}')^2 > \frac{V^2}{600A} \text{ or } |\circledEight| \leq \frac{V}{10}\right\}. \label{ch3:eq:bound_8_gap_SEG_chapter_3}
    \end{align}
    Moreover, $E_{T-1}$ implies
    \begin{align}
        \sum\limits_{l=0}^{T-1}\sum\limits_{i=1}^n(\widetilde\sigma_{i,l}')^2 & \overset{\eqref{ch3:eq:gap_thm_SEG_technical_7_8_chapter_3}}{\leq} \frac{4\gamma^2 V}{n^2} \sum\limits_{l=0}^{T-1} \sum\limits_{i=1}^n \Exp_{\xi_{2,i}^l}\left[\|\theta_{i,l}^u\|^2\right] \notag \\
        & \overset{\eqref{ch3:eq:variance_theta_omega_chapter_3}}{\leq} \frac{72(K+1)\gamma^2 V \lambda^{2-{p}} \sigma^{p}}{n} \notag \\
        & \overset{\eqref{ch3:eq:lambda_SEG_monotone_chapter_3}}{\leq} \frac{72 (K+1)\gamma^{{p}}\sigma^{p} V^{2-\frac{{p}}{2}}}{60^{2-{p}} n^{{p}-1} A^{2-{p}}} \overset{\eqref{ch3:eq:gamma_SEG_monotone_chapter_3}}{\leq} \frac{V^2}{600A}. \label{ch3:eq:bound_8_variances_gap_SEG_chapter_3}
    \end{align}

    \paragraph{Upper bound for $\circledNine$.} Probability event $E_{T-1}$ implies
    \begin{align}
        \circledNine & \leq 2\gamma \sum\limits_{l=0}^{T-1}\|\zeta_l'\| \cdot \|\theta_l^b\| \overset{\eqref{ch3:eq:gap_thm_SEG_technical_8_chapter_3}, \eqref{ch3:eq:bias_theta_omega_chapter_3}}{\leq} \frac{2\cdot 2^{p}(K+1)\gamma \sqrt{V}\sigma^{p}}{\lambda^{{p}-1}} \notag \\
        & \overset{\eqref{ch3:eq:lambda_SEG_monotone_chapter_3}}{=} \frac{2\cdot 2^{p}\cdot 60^{{p}-1} (K+1)\gamma^{p} \sigma^{p} A^{{p}-1}}{n^{{p}-1}V^{\frac{{p}}{2}-1}} \overset{\eqref{ch3:eq:gamma_SEG_monotone_chapter_3}}{\leq} \frac{V}{100}. \label{ch3:eq:bound_9_variances_gap_SEG_chapter_3}
    \end{align}

    That is, we derive the upper bounds for  $2\gamma \sqrt{V} \left\|\sum_{l = 0}^{T-1}\theta_l\right\|, \circledOne, \circledTwo, \circledThree, \circledFour, \circledFive, \circledSix, \circledSeven$. More precisely, $E_{T-1}$ implies
    \begin{gather*}
        A_T \overset{\eqref{ch3:eq:gap_thm_SEG_technical_6_chapter_3}}{\leq} 4V + 2\gamma \sqrt{V}\left\|\sum\limits_{l=0}^{T-1} \theta_l\right\| + \circledOne + \circledTwo + \circledThree + \circledFour + \circledFive + \circledSix + \circledSeven,\\
        \circledSix \overset{\eqref{ch3:eq:SEG_extra_sums_six_distributed_chapter_3}}{=} \circledSix' + \frac{32\gamma}{n} \sum\limits_{j=2}^n \left\langle \frac{\gamma}{n}\sum\limits_{i=1}^{j-1} \theta_{i,T-1}^u - \zeta_j^{T-1}, \theta_{j,T-1}^u \right\rangle,\\
        \circledSeven \overset{\eqref{ch3:eq:SEG_extra_sums_seven_distributed_chapter_3}}{=} \circledSeven' + \frac{24\gamma}{n} \sum\limits_{j=2}^n \left\langle \frac{\gamma}{n}\sum\limits_{i=1}^{j-1} \omega_{i,T-1}^u - \delta_j^{T-1}, \omega_{j,T-1}^u \right\rangle,\\
        2\gamma \sqrt{V}\left\|\sum\limits_{l=0}^{T-1} \theta_l\right\| \overset{\eqref{ch3:eq:norm_sum_theta_bound_gap_SEG_chapter_3}}{\leq} 2\sqrt{V} \sqrt{\frac{\circledThree + \circledFour + \circledFive + \circledSix + \circledSeven}{8} + \circledEight + \circledNine},\\
        \circledTwo \overset{\eqref{ch3:eq:bound_2_variances_gap_SEG_chapter_3}}{\leq} \frac{3V}{100},\quad \circledThree \overset{\eqref{ch3:eq:bound_3_variances_gap_SEG_chapter_3}}{\leq} \frac{V}{20},\quad \circledFive \overset{\eqref{ch3:eq:nsjcbdjhcbfjdhfbj_chapter_3}}{\leq} \frac{V}{6},\quad \circledNine \overset{\eqref{ch3:eq:bound_9_variances_gap_SEG_chapter_3}}{\leq} \frac{V}{100},\\
        \sum\limits_{l=0}^{T-1}\sum\limits_{i=1}^n\sigma_{i,l}^2 \overset{\eqref{ch3:eq:bound_1_variances_gap_SEG_chapter_3}}{\leq}  \frac{3V^2}{200A},\quad \sum\limits_{l=0}^{T-1}\sum\limits_{i=1}^n\widetilde\sigma_{i,l}^2 \overset{\eqref{ch3:eq:bound_1_variances_gap_SEG_1_chapter_3}}{\leq} \frac{V^2}{216A},\quad\\
        \sum\limits_{l=0}^{T-1}\sum\limits_{j=2}^n\widehat\sigma_{j,l}^2 \overset{\eqref{ch3:eq:bound_6_variances_gap_SEG_chapter_3}}{\leq}  \frac{8V^2}{75A},\\
        \sum\limits_{l=0}^{T-1}\sum\limits_{j=2}^n(\sigma_{j,l}')^2 \overset{\eqref{ch3:eq:bound_7_variances_gap_SEG_chapter_3}}{\leq} \frac{3V^2}{50A},\quad \sum\limits_{l=0}^{T-1}\sum\limits_{i=1}^n(\widetilde\sigma_{i,l}')^2 \leq \frac{V^2}{600A}.
    \end{gather*}
    In addition, we also establish (see \eqref{ch3:eq:bound_1_gap_SEG_chapter_3}, \eqref{ch3:eq:bound_1_gap_SEG_1_chapter_3}, \eqref{ch3:eq:bound_6_gap_SEG_chapter_3}, \eqref{ch3:eq:bound_7_gap_SEG_chapter_3}, \eqref{ch3:eq:bound_8_gap_SEG_chapter_3} and our induction assumption)
    \begin{gather*}
        \Prob\{E_{T-1}\} \geq 1 - \frac{(T-1)\beta}{K+1},\quad \Prob\{E_{\circledOne}\} \geq 1 - \frac{\beta}{24n(K+1)},\\
        \Prob\{E_{\circledFour}\} \geq 1 - \frac{\beta}{24n(K+1)}, \quad \Prob\{E_{\circledSix'}\} \geq 1 - \frac{\beta}{24n(K+1)},\\
        \Prob\{E_{\circledSeven'}\} \geq 1 - \frac{\beta}{24n(K+1)},\quad \Prob\{E_{\circledEight}\} \geq 1 - \frac{\beta}{24n(K+1)},
    \end{gather*}
    where
    \begin{align}
        E_{\circledOne} & = \left\{\text{either } \sum\limits_{l=0}^{T-1}\sum\limits_{i=1}^n\sigma_{i,l}^2 > \frac{3V^2}{200A} \text{ or } |\circledOne| \leq \frac{3V}{10}\right\},\notag\\
        E_{\circledFour} & = \left\{\text{either } \sum\limits_{l=0}^{T-1}\sum\limits_{i=1}^n\widetilde\sigma_{i,l}^2 > \frac{V^2}{216A} \text{ or } |\circledFour| \leq \frac{V}{6}\right\},\notag\\
        E_{\circledSix'} & = \left\{\text{either } \sum\limits_{l=0}^{T-1}\sum\limits_{j=2}^n\widehat\sigma_{j,l}^2 > \frac{8V^2}{75A} \text{ or } |\circledSix'| \leq \frac{4V}{5}\right\},\notag\\
        E_{\circledSeven'} & = \left\{\text{either } \sum\limits_{l=0}^{T-1}\sum\limits_{j=2}^n(\sigma_{j,l}')^2 > \frac{3V^2}{50A} \text{ or } |\circledSeven'| \leq \frac{3V}{5}\right\}\notag \\
        E_{\circledEight} & = \left\{\text{either } \sum\limits_{l=0}^{T-1}\sum\limits_{i=1}^n(\widetilde\sigma_{i,l}')^2 > \frac{V^2}{600A} \text{ or } |\circledEight| \leq \frac{V}{10}\right\} \notag
    \end{align}
    Therefore, denoting
    \begin{gather*}
        S_\theta \eqdef \frac{\gamma}{n} \sum\limits_{j=2}^n \left\langle \frac{\gamma}{n}\sum\limits_{i=1}^{j-1} \theta_{i,T-1}^u - \zeta_j^{T-1}, \theta_{j,T-1}^u \right\rangle,\\
        S_\omega \eqdef \frac{\gamma}{n} \sum\limits_{j=2}^n \left\langle \frac{\gamma}{n}\sum\limits_{i=1}^{j-1} \omega_{i,T-1}^u - \delta_j^{T-1}, \omega_{j,T-1}^u \right\rangle,
    \end{gather*}
    we get that probability event $E_{T-1} \cap E_{\circledOne} \cap E_{\circledFour} \cap E_{\circledSix'} \cap E_{\circledSeven'} \cap E_{\circledEight}$ implies
    \begin{align}
        \left\|\gamma\sum\limits_{l=0}^{T-1} \theta_l\right\| & \leq \sqrt{\frac{1}{8}\left(\frac{V}{20} + \frac{V}{6} + \frac{V}{6} + \frac{4V}{5} + \frac{3V}{5}\right) + \frac{V}{10} + \frac{V}{100}} + \sqrt{4S_\theta + 3S_\omega} \notag\\
        & \leq \sqrt{V} + \sqrt{4S_\theta + 3S_\omega}, \label{ch3:eq:gap_thm_SEG_technical_10_chapter_3}
    \end{align}
    and
    \begin{align}
        A_T & \leq 4V + 2V + 2\sqrt{V}\sqrt{4S_\theta + 3S_\omega} \notag\\
        &\quad + \frac{3V}{10} + \frac{3V}{100} + \frac{V}{20} + \frac{V}{6} + \frac{V}{6} + \frac{4V}{5} + \frac{3V}{5} + 32S_\theta + 24S_\omega \notag\\
        & \leq 9V + 2\sqrt{V}\sqrt{4S_\theta + 3S_\omega} + 32S_\theta + 24S_\omega, \label{ch3:eq:gap_thm_SEG_technical_11_chapter_3}
    \end{align}

    In the final part of the proof, we will show that $\frac{\gamma}{n}\sum\limits_{i=1}^{j-1} \theta_{i,T-1}^u = \zeta_j^{T-1}$ and $\frac{\gamma}{n}\sum\limits_{i=1}^{j-1} \omega_{i,T-1}^u = \delta_j^{T-1}$ with high probability. In particular, we consider probability event $\widetilde{E}_{T-1,j}$ defined as follows: inequalities
    \begin{align}
\left\|\frac{\gamma}{n}\sum\limits_{i=1}^{r-1} \theta_{i,T-1}^u\right\| & \leq \frac{\sqrt{V}}{2}, \quad \left\|\frac{\gamma}{n}\sum\limits_{i=1}^{r-1} \omega_{i,T-1}^u\right\| \notag \\
& \leq \frac{\sqrt{V}}{2} \notag
    \end{align}
    hold for $r = 2, \ldots, j$ simultaneously. We want to show that $\Prob\{E_{T-1} \cap \widetilde{E}_{T-1,j}\} \geq 1 - \frac{(T-1)\beta}{K+1} - \frac{j\beta}{8n(K+1)}$ for all $j = 2, \ldots, n$. For $j = 2$ the statement is trivial since
    \begin{equation*}
        \left\|\frac{\gamma}{n} \theta_{1,T-1}^u\right\| \overset{\eqref{ch3:eq:theta_omega_magnitude_chapter_3}}{\leq} \frac{2\gamma\lambda}{n} \leq \frac{\sqrt{V}}{2},\quad \left\|\frac{\gamma}{n} \omega_{1,T-1}^u\right\| \overset{\eqref{ch3:eq:theta_omega_magnitude_chapter_3}}{\leq} \frac{2\gamma\lambda}{n} \leq \frac{\sqrt{V}}{2}.
    \end{equation*}
    Next, we assume that the statement holds for some $j = m-1 < n$, i.e.,  $\Prob\{E_{T-1}\cap\widetilde{E}_{T-1,m-1}\} \geq 1 - \frac{(T-1)\beta}{K+1} - \frac{(m-1)\beta}{8n(K+1)}$. Our goal is to prove that $\Prob\{E_{T-1}\cap\widetilde{E}_{T-1,m}\} \geq 1 - \frac{(T-1)\beta}{K+1} - \frac{m\beta}{8n(K+1)}$. First, we consider $\left\|\frac{\gamma}{n}\sum\limits_{i=1}^{m-1} \theta_{i,T-1}^u\right\|$:
    \begin{align}
        \left\|\frac{\gamma}{n}\sum\limits_{i=1}^{m-1} \theta_{i,T-1}^u\right\| & = \sqrt{\frac{\gamma^2}{n^2}\left\|\sum\limits_{i=1}^{m-1} \theta_{i,T-1}^u\right\|^2} \notag\\
        & = \sqrt{\frac{\gamma^2}{n^2}\sum\limits_{i=1}^{m-1}\|\theta_{i,T-1}^u\|^2 + \frac{2\gamma}{n}\sum\limits_{i=1}^{m-1}\left\langle \frac{\gamma}{n}\sum\limits_{r=1}^{i-1}\theta_{r,T-1}^u , \theta_{i,T-1}^u \right\rangle} \notag\\
        & \leq \sqrt{\frac{\gamma^2}{n^2}\sum\limits_{l=0}^{T-1}\sum\limits_{i=1}^{m-1}\|\theta_{i,l}^u\|^2 + \frac{2\gamma}{n}\sum\limits_{i=1}^{m-1}\left\langle \frac{\gamma}{n}\sum\limits_{r=1}^{i-1}\theta_{r,T-1}^u , \theta_{i,T-1}^u \right\rangle}. \notag
    \end{align}
    Similarly, we have
    \begin{align}
        \left\|\frac{\gamma}{n}\sum\limits_{i=1}^{m-1} \omega_{i,T-1}^u\right\| & = \sqrt{\frac{\gamma^2}{n^2}\left\|\sum\limits_{i=1}^{m-1} \omega_{i,T-1}^u\right\|^2} \notag\\
        & = \sqrt{\frac{\gamma^2}{n^2}\sum\limits_{i=1}^{m-1}\|\omega_{i,T-1}^u\|^2 + \frac{2\gamma}{n}\sum\limits_{i=1}^{m-1}\left\langle \frac{\gamma}{n}\sum\limits_{r=1}^{i-1}\omega_{r,T-1}^u , \omega_{i,T-1}^u \right\rangle} \notag\\
        & \leq \sqrt{\frac{\gamma^2}{n^2}\sum\limits_{l=0}^{T-1}\sum\limits_{i=1}^{m-1}\|\omega_{i,l}^u\|^2 + \frac{2\gamma}{n}\sum\limits_{i=1}^{m-1}\left\langle \frac{\gamma}{n}\sum\limits_{r=1}^{i-1}\omega_{r,T-1}^u , \omega_{i,T-1}^u \right\rangle}. \notag
    \end{align}

    Next, we introduce a new notation:
    \begin{gather*}
        \rho_{i,T-1} = \begin{cases}
            \frac{\gamma}{n}\sum\limits_{r=1}^{i-1}\theta_{r,T-1}^u,& \text{if } \left\|\frac{\gamma}{n}\sum\limits_{r=1}^{i-1}\theta_{r,T-1}^u\right\| \leq \frac{\sqrt{V}}{2},\\
            0,& \text{otherwise}
        \end{cases},\\
        \rho_{i,T-1}' = \begin{cases}
            \frac{\gamma}{n}\sum\limits_{r=1}^{i-1}\omega_{r,T-1}^u,& \text{if } \left\|\frac{\gamma}{n}\sum\limits_{r=1}^{i-1}\omega_{r,T-1}^u\right\| \leq \frac{\sqrt{V}}{2},\\
            0,& \text{otherwise}
        \end{cases}
    \end{gather*}
    for $i = 1,\ldots,m-1$. By definition, we have
    \begin{equation}
        \|\rho_{i,T-1}\| \leq \frac{\sqrt{V}}{2},\quad \|\rho_{i,T-1}'\| \leq \frac{\sqrt{V}}{2} \label{ch3:eq:SEG_bound_rho_chapter_3}
    \end{equation}
    for $i = 1,\ldots,m-1$. Moreover, $\widetilde{E}_{T-1,m-1}$ implies $\rho_{i,T-1} = \frac{\gamma}{n}\sum\limits_{r=1}^{i-1}\theta_{r,T-1}^u$, $\rho_{i,T-1}' = \frac{\gamma}{n}\sum\limits_{r=1}^{i-1}\omega_{r,T-1}^u$ for $i = 1,\ldots,m-1$ and
    \begin{align}
        \left\|\frac{\gamma}{n}\sum\limits_{i=1}^{m-1} \theta_{i,T-1}^u\right\| & \leq \sqrt{\circledThree + \circledFour + \circledTen}, \notag\\
        \left\|\frac{\gamma}{n}\sum\limits_{i=1}^{m-1} \omega_{i,T-1}^u\right\| & \leq \sqrt{\circledThree + \circledFour + \circledTen'}, \notag
    \end{align}
    where
    \begin{gather*}
        \circledTen =  \frac{2\gamma}{n}\sum\limits_{i=1}^{m-1}\left\langle \rho_{i,T-1} , \theta_{i,T-1}^u \right\rangle,\quad \circledTen' =  \frac{2\gamma}{n}\sum\limits_{i=1}^{m-1}\left\langle \rho_{i,T-1}' , \omega_{i,T-1}^u \right\rangle.
    \end{gather*}
    It remains to estimate $\circledTen$ and $\circledTen'$.

    \textbf{Upper bound for $\circledTen$.} To estimate this sum, we will use Bernstein's inequality. The summands have conditional expectations equal to zero:
    \begin{align*}
\MoveEqLeft[1] \Exp_{\xi_{2,i}^{T-1}}\left[\frac{2\gamma}{n} \langle \rho_{i,T-1}, \theta_{i,T-1}^u \rangle\right] \\
& = \frac{2\gamma}{n} \left\langle \rho_{i,T-1}, \Exp_{\xi_{2,i}^{T-1}}[\theta_{i,T-1}^u] \right\rangle = 0,
    \end{align*}
    since random vectors $\{\theta_{i,T-1}^u\}_{i=1}^n$ are independent. Thus, for $i = 1, \ldots, m-1$ the summands $(2\gamma/n) \langle \rho_{i,T-1}, \theta_{i,T-1}^u \rangle$ form a martingale difference sequence. Next, the summands are bounded:
    \begin{align}
        \left|\frac{2\gamma}{n}\langle \rho_{i,T-1},
        \theta_{i,T-1}^u \rangle \right|
        &\leq \frac{2\gamma}{n} \|\rho_{i,T-1}\|
        \cdot \|\theta_{i,T-1}^u\| \notag\\
        &\overset{\eqref{ch3:eq:SEG_bound_rho_chapter_3},\eqref{ch3:eq:theta_omega_magnitude_chapter_3}}{\leq}
        \frac{2\gamma}{n} \sqrt{V}\lambda \notag\\
        &\overset{\eqref{ch3:eq:lambda_SEG_monotone_chapter_3}}{=}
        \frac{V}{30A} \eqdef c.
        \label{ch3:eq:gap_thm_SEG_technical_6_5_10_chapter_3}
    \end{align}
    Finally, conditional variances $(\widehat\sigma_{i,T-1}')^2 \eqdef \Exp_{\xi_{2,i}^{T-1}}\left[\frac{4\gamma^2}{n^2} \langle \rho_{i,T-1}, \theta_{i,T-1}^u \rangle^2\right]$ of the summands are bounded:
    \begin{align}
(\widehat\sigma_{i,T-1}')^2 & \leq \Exp_{\xi_{2,i}^{T-1}}\left[\frac{4\gamma^2 }{n^2}\|\rho_{i,T-1}\|^2\cdot \|\theta_{i,T-1}^u\|^2\right] \notag \\
& \overset{\eqref{ch3:eq:SEG_bound_rho_chapter_3}}{\leq} \frac{\gamma^2 V}{n^2} \Exp_{\xi_{2,i}^{T-1}}\left[\|\theta_{i,T-1}^u\|^2\right]. \label{ch3:eq:gap_thm_SEG_technical_7_10_chapter_3}
    \end{align}
    Applying Bernstein's inequality (Lemma~\ref{ch3:lem:Bernstein_ineq_chapter_3}) with $X_{i} = \frac{2\gamma}{n} \langle \rho_{i,T-1}, \theta_{i,T-1}^u \rangle$, constant $c$ defined in \eqref{ch3:eq:gap_thm_SEG_technical_6_5_10_chapter_3}, $b = \frac{V}{30}$, $H = \frac{V^2}{5400A}$, we get
    \begin{align*}
\MoveEqLeft[1] \Prob\left\{|\circledTen| > \frac{V}{30} \text{ and } \sum\limits_{i=1}^{m-1}(\widehat\sigma_{i,T-1}')^2 \leq \frac{V^2}{5400A}\right\} \\
& \leq 2\exp\left(- \frac{b^2}{2H + 2cb/3}\right) = \frac{\beta}{24n(K+1)}.
    \end{align*}
    The above is equivalent to 
    \begin{align}
\Prob\{E_{\circledTen}\} & \geq 1 - \frac{\beta}{24n(K+1)},\; \text{ for }\; E_{\circledTen} \notag \\
& = \left\{\text{either } \sum\limits_{i=1}^{m-1}(\widehat\sigma_{i,T-1}')^2 > \frac{V^2}{5400A} \text{ or } |\circledTen| \leq \frac{V}{30}\right\}. \label{ch3:eq:bound_10_gap_SEG_chapter_3}
    \end{align}
    Moreover, $E_{T-1}$ implies
    \begin{align}
        \sum\limits_{i=1}^{m-1}(\widehat\sigma_{i,T-1}')^2
        &\overset{\eqref{ch3:eq:gap_thm_SEG_technical_7_10_chapter_3}}{\leq}
        \frac{\gamma^2 V}{n^2} \sum\limits_{i=1}^n
        \Exp_{\xi_{2,i}^{T-1}}\left[\|\theta_{i,T-1}^u\|^2\right]\notag\\
        &\overset{\eqref{ch3:eq:variance_theta_omega_chapter_3}}{\leq}
        \frac{18\gamma^2 V \lambda^{2-{p}} \sigma^{p}}{n}\notag\\
        &\overset{\eqref{ch3:eq:lambda_SEG_monotone_chapter_3}}{\leq}
        \frac{18 \gamma^{{p}}\sigma^{p} V^{2-\frac{{p}}{2}}}
        {60^{2-{p}} n^{{p}-1} A^{2-{p}}}\notag\\
        &\overset{\eqref{ch3:eq:gamma_SEG_monotone_chapter_3}}{\leq}
        \frac{V^2}{5400A}.
        \label{ch3:eq:bound_10_variances_gap_SEG_chapter_3}
    \end{align}

    \textbf{Upper bound for $\circledTen'$.} To estimate this sum, we will use Bernstein's inequality. The summands have conditional expectations equal to zero:
    \begin{align*}
\MoveEqLeft[1] \Exp_{\xi_{1,i}^{T-1}}\left[\frac{2\gamma}{n} \langle \rho_{i,T-1}', \omega_{i,T-1}^u \rangle\right] \\
& = \frac{2\gamma}{n} \left\langle \rho_{i,T-1}', \Exp_{\xi_{1,i}^{T-1}}[\omega_{i,T-1}^u] \right\rangle = 0,
    \end{align*}
    since random vectors $\{\omega_{i,T-1}^u\}_{i=1}^n$ are independent. Thus, for $i = 1, \ldots, m-1$ the summands $(2\gamma/n) \langle \rho_{i,T-1}', \omega_{i,T-1}^u \rangle$ form a martingale difference sequence. Next, the summands are bounded:
    \begin{align}
        \left|\frac{2\gamma}{n}\langle \rho_{i,T-1}',
        \omega_{i,T-1}^u \rangle \right|
        &\leq \frac{2\gamma}{n} \|\rho_{i,T-1}'\|
        \cdot \|\omega_{i,T-1}^u\| \notag\\
        &\overset{\eqref{ch3:eq:SEG_bound_rho_chapter_3},\eqref{ch3:eq:theta_omega_magnitude_chapter_3}}{\leq}
        \frac{2\gamma}{n} \sqrt{V}\lambda \notag\\
        &\overset{\eqref{ch3:eq:lambda_SEG_monotone_chapter_3}}{=}
        \frac{V}{30A} \eqdef c.
        \label{ch3:eq:gap_thm_SEG_technical_6_5_10'_chapter_3}
    \end{align}
    Finally, conditional variances $(\widetilde\sigma_{i,T-1}')^2 \eqdef \Exp_{\xi_{1,i}^{T-1}}\left[\frac{4\gamma^2}{n^2} \langle \rho_{i,T-1}', \omega_{i,T-1}^u \rangle^2\right]$ of the summands are bounded:
    \begin{align}
(\widetilde\sigma_{i,T-1}')^2 & \leq \Exp_{\xi_{1,i}^{T-1}}\left[\frac{4\gamma^2 }{n^2}\|\rho_{i,T-1}'\|^2\cdot \|\omega_{i,T-1}^u\|^2\right] \notag \\
& \overset{\eqref{ch3:eq:SEG_bound_rho_chapter_3}}{\leq} \frac{\gamma^2 V}{n^2} \Exp_{\xi_{1,i}^{T-1}}\left[\|\omega_{i,T-1}^u\|^2\right]. \label{ch3:eq:gap_thm_SEG_technical_7_10'_chapter_3}
    \end{align}
    Applying Bernstein's inequality (Lemma~\ref{ch3:lem:Bernstein_ineq_chapter_3}) with $X_{i} = \frac{2\gamma}{n} \langle \rho_{i,T-1}', \omega_{i,T-1}^u \rangle$, constant $c$ defined in \eqref{ch3:eq:gap_thm_SEG_technical_6_5_10'_chapter_3}, $b = \frac{V}{30}$, $H = \frac{V^2}{5400A}$, we get
    \begin{align*}
\MoveEqLeft[1] \Prob\left\{|\circledTen'| > \frac{V}{30} \text{ and } \sum\limits_{i=1}^{m-1}(\widetilde\sigma_{i,T-1}')^2 \leq \frac{V^2}{5400A}\right\} \\
& \leq 2\exp\left(- \frac{b^2}{2H + 2cb/3}\right) = \frac{\beta}{24n(K+1)}.
    \end{align*}
    The above is equivalent to 
    \begin{equation}
        \begin{aligned}
        \Prob\{E_{\circledTen'}\}
        &\geq 1 - \frac{\beta}{24n(K+1)},\; \text{ for }\;\\
        E_{\circledTen'}
        &= \left\{\text{either } \sum\limits_{i=1}^{m-1}(\widetilde\sigma_{i,T-1}')^2 > \frac{V^2}{5400A} \text{ or } |\circledTen'| \leq \frac{V}{30}\right\}.
        \end{aligned}
        \label{ch3:eq:bound_10_gap_SEG'_chapter_3}
    \end{equation}
    Moreover, $E_{T-1}$ implies
    \begin{align}
        \sum\limits_{i=1}^{m-1}(\widetilde\sigma_{i,T-1}')^2 & \overset{\eqref{ch3:eq:gap_thm_SEG_technical_7_10'_chapter_3}}{\leq} \frac{\gamma^2 V}{n^2} \sum\limits_{i=1}^n \Exp_{\xi_{1,i}^{T-1}}\left[\|\omega_{i,T-1}^u\|^2\right] \notag \\
        & \overset{\eqref{ch3:eq:variance_theta_omega_chapter_3}}{\leq} \frac{18\gamma^2 V \lambda^{2-{p}} \sigma^{p}}{n} \notag \\
        & \overset{\eqref{ch3:eq:lambda_SEG_monotone_chapter_3}}{\leq} \frac{18 \gamma^{{p}}\sigma^{p} V^{2-\frac{{p}}{2}}}{60^{2-{p}} n^{{p}-1} A^{2-{p}}} \overset{\eqref{ch3:eq:gamma_SEG_monotone_chapter_3}}{\leq} \frac{V^2}{5400A}. \label{ch3:eq:bound_10_variances_gap_SEG'_chapter_3}
    \end{align}

    Putting all together we get that $E_{T-1}\cap \widetilde{E}_{T-1, m-1}$ implies

    \begin{gather*}
        \left\|\frac{\gamma}{n}\sum\limits_{i=1}^{m-1} \theta_{i,T-1}^u\right\| \leq \sqrt{\circledThree + \circledFour + \circledTen},\\
        \left\|\frac{\gamma}{n}\sum\limits_{i=1}^{m-1} \omega_{i,T-1}^u\right\| \leq \sqrt{\circledThree + \circledFour + \circledTen'},\quad \circledThree \overset{\eqref{ch3:eq:bound_3_variances_gap_SEG_chapter_3}}{\leq} \frac{V}{20},\\
        \sum\limits_{l=0}^{T-1}\sum\limits_{i=1}^n\widetilde\sigma_{i,l}^2 \overset{\eqref{ch3:eq:bound_1_variances_gap_SEG_1_chapter_3}}{\leq} \frac{V^2}{216A},\\
        \sum\limits_{i=1}^{m-1}(\widehat\sigma_{i,T-1}')^2 \leq \frac{V^2}{5400A},\quad \sum\limits_{i=1}^{m-1}(\widetilde\sigma_{i,T-1}')^2 \leq \frac{V^2}{5400A}.
    \end{gather*}
    In addition, we also establish (see \eqref{ch3:eq:bound_1_gap_SEG_1_chapter_3}, \eqref{ch3:eq:bound_10_gap_SEG_chapter_3}, \eqref{ch3:eq:bound_10_gap_SEG'_chapter_3} and our induction assumption)
    \begin{gather*}
        \Prob\{E_{T-1}\cap \widetilde{E}_{T-1,m-1}\} \geq 1 - \frac{(T-1)\beta}{K+1} - \frac{(m-1)\beta}{8n(K+1)},\\
        \Prob\{E_{\circledFour}\} \geq 1 - \frac{\beta}{24n(K+1)}, \quad \Prob\{E_{\circledTen}\} \geq 1 - \frac{\beta}{24n(K+1)},\\
        \Prob\{E_{\circledTen'}\} \geq 1 - \frac{\beta}{24n(K+1)}
    \end{gather*}
    where
    \begin{align}
        E_{\circledFour} & = \left\{\text{either } \sum\limits_{l=0}^{T-1}\sum\limits_{i=1}^n\widetilde\sigma_{i,l}^2 > \frac{V^2}{216A} \text{ or } |\circledFour| \leq \frac{V}{6}\right\},\notag\\
        E_{\circledTen} & = \left\{\text{either } \sum\limits_{i=1}^{m-1}(\widehat\sigma_{i,T-1}')^2 > \frac{V^2}{5400A} \text{ or } |\circledTen| \leq \frac{V}{30}\right\}, \notag\\
        E_{\circledTen'} & = \left\{\text{either } \sum\limits_{i=1}^{m-1}(\widetilde\sigma_{i,T-1}')^2 > \frac{V^2}{5400A} \text{ or } |\circledTen'| \leq \frac{V}{30}\right\} \notag
    \end{align}
    Therefore, probability event $E_{T-1} \cap \widetilde{E}_{T-1,m-1} \cap E_{\circledFour} \cap E_{\circledTen} \cap E_{\circledTen'}$ implies
    \begin{gather*}
        \left\|\frac{\gamma}{n}\sum\limits_{i=1}^{m-1} \theta_{i,T-1}^u\right\| \leq \sqrt{\frac{V}{20} + \frac{V}{6} + \frac{V}{30}} = \frac{\sqrt{V}}{2},\\
        \left\|\frac{\gamma}{n}\sum\limits_{i=1}^{m-1} \omega_{i,T-1}^u\right\| \leq \sqrt{\frac{V}{20} + \frac{V}{6} + \frac{V}{30}} = \frac{\sqrt{V}}{2}.
    \end{gather*}
    This implies $\widetilde E_{T-1,m}$ and
    \begin{align*}
        \Prob\{E_{T-1} \cap \widetilde{E}_{T-1,m}\} & \geq \Prob\{E_{T-1} \cap \widetilde{E}_{T-1,m-1} \cap E_{\circledFour} \cap E_{\circledTen} \cap E_{\circledTen'}\} \\
        & = 1 - \Prob\left\{\overline{E_{T-1} \cap \widetilde{E}_{T-1,m-1}} \cup \overline{E}_{\circledFour} \cup \overline{E}_{\circledTen} \cup \overline{E}_{\circledTen'}\right\}\\
        & \geq 1 - \frac{(T-1)\beta}{K+1} - \frac{m\beta}{8n(K+1)}.
    \end{align*}
    Therefore, for all $m = 2,\ldots,n$ the statement holds and, in particular, $\Prob\{E_{T-1} \cap \widetilde{E}_{T-1, n}\} \geq 1 - \frac{(T-1)\beta}{K+1} - \frac{\beta}{8(K+1)}$. Taking into account \eqref{ch3:eq:gap_thm_SEG_technical_11_chapter_3}, we conclude that $E_{T-1} \cap \widetilde{E}_{T-1, n} \cap E_{\circledOne} \cap E_{\circledFour} \cap E_{\circledSix'} \cap E_{\circledSeven'} \cap E_{\circledEight}$ implies
    \begin{equation*}
       \left\|\gamma\sum\limits_{l=0}^{T-1} \theta_l\right\| \leq \sqrt{V},\quad A_T \leq 9V
    \end{equation*}
    (since $S_\theta = S_\omega = 0$ on $\widetilde{E}_{T-1,n}$), which is equivalent to \eqref{ch3:eq:induction_inequality_1_SEG_chapter_3} and \eqref{ch3:eq:induction_inequality_2_SEG_chapter_3} for $t = T$; moreover, $\widetilde{E}_{T-1,n}$ gives the last two inequalities in the definition of $E_T$. Finally,
    \begin{align*}
        \Prob\left\{E_T\right\} & \geq \Prob\left\{E_{T-1} \cap \widetilde{E}_{T-1,n} \cap E_{\circledOne} \cap E_{\circledFour} \cap E_{\circledSix'} \cap E_{\circledSeven'} \cap E_{\circledEight} \right\} \\
        & = 1 - \Prob\left\{\overline{E_{T-1} \cap \widetilde{E}_{n}} \cup \overline{E}_{\circledOne} \cup \overline{E}_{\circledFour} \cup \overline{E}_{\circledSix'} \cup \overline{E}_{\circledSeven'} \cup \overline{E}_{\circledEight}\right\}\\
        & \geq 1 - \frac{(T-1)\beta}{K+1} - \frac{\beta}{8(K+1)} - 5\cdot \frac{\beta}{24n(K+1)} \geq 1 - \frac{T\beta}{K+1}.
    \end{align*}

    In other words, we showed that $\Prob\{E_k\} \geq 1 - k\beta/(K+1)$ for all $k = 0,1,\ldots,K+1$. For $k = K+1$ we have that with probability at least $1 - \beta$
    \begin{align*}
        \gap_{\sqrt{V}}(\thesischthreetx^{K}_{\avg})
        &= \max_{B_{\sqrt{V}}(x^*)}\left\{\langle F(u), \thesischthreetx^{K}_{\avg} - u\rangle + \Psi(\thesischthreetx^{K}_{\avg}) - \Psi(u)\right\}\\
        &\leq \frac{1}{2\gamma(K + 1)}\max_{B_{\sqrt{V}}(x^*)}\Big\{2\gamma(K+1)\Big(\langle F(u), \thesischthreetx^{K}_{\avg} - u\rangle\\
        &\qquad\qquad + \Psi(\thesischthreetx^{K}_{\avg}) - \Psi(u)\Big) + \|x^{K+1} - u\|^2\Big\}\\
        &\overset{\eqref{ch3:eq:gap_thm_SEG_technical_1_5_chapter_3}}{\leq}
        \frac{9V}{2\gamma(K+1)}.
    \end{align*}
    Finally, if $\gamma$ equals the minimum in \eqref{ch3:eq:gamma_SEG_monotone_chapter_3}, then with probability at least $1-\beta$
    \begin{align*}
        \gap_{\sqrt{V}}(\thesischthreetx^{K}_{\avg}) & \leq \frac{9V}{2\gamma(K+1)} \\
        & = \max\left\{\frac{8640 LV A}{K+1}, \frac{9\cdot 97200^{1/p}\sigma \sqrt{V} A^{(p-1)/p}}{2\cdot 60^{(2-p)/p} n^{(p-1)/p}(K+1)^{(p-1)/p}}\right\}\\
        &= \cO\left(\max\left\{\frac{LV\ln(nK/\beta)}{K}, \frac{\sigma \sqrt{V} \ln^{(p-1)/p}(nK/\beta)}{n^{(p-1)/p}K^{(p-1)/p}}\right\}\right).
    \end{align*}
    To get $\gap_{\sqrt{V}}(\thesischthreetx^{K}_{\avg}) \leq \varepsilon$ with probability $\geq 1 - \beta$, $K$ should be
    \begin{align*}
         K = \cO\biggl(\max\biggl\{ & \frac{LV}{\varepsilon}\ln\frac{nLV}{\varepsilon\beta}, \frac{1}{n}\left(\frac{\sigma \sqrt{V}}{\varepsilon}\right)^{p/(p-1)}\ln \frac{\sigma \sqrt{V}}{\varepsilon\beta}\biggr\}\biggr),
    \end{align*}
    that concludes the proof.
\end{proof}

\subsection{Quasi-Strongly Monotone Case}

We start with the following lemma.

\begin{lemma}\label{ch3:lem:SEG_optimization_lemma_QSM_chapter_3}
    Let Assumptions \ref{ch3:as:L_Lip_chapter_3} and \ref{ch3:as:str_monotonicity_chapter_3} hold for $Q = B_{4n\sqrt{V}}(x^*)$, where
    \begin{gather*}
        A = \ln\frac{48n(K+1)}{\beta},\quad V \geq \|x^0 - x^*\|^2 + \frac{72\cdot 10^6\gamma^2A^2}{n^2}\sum_{i=1}^n\|F_i(x^*)\|^2,
    \end{gather*}
    the initial shifts are $\thesischthreethh_i^0 = \thesischthreehh_i^0 = 0$, $\nu = \gamma\mu$, and
    \begin{equation*}
        0 < \gamma \leq \min\left\{\frac{1}{6L}, \frac{\sqrt{n}}{15000 L A}, \frac{1}{72\cdot 10^6 \mu A^2}\right\}.
    \end{equation*}
    If $x^k$ and $\thesischthreetx^k$ lie in $B_{4n\sqrt{V}}(x^*)$ for all $k = 0,1,\ldots,K$ for some $K \geq 0$, then the iterates produced by \algname{DProx-clipped-SEG-shift} satisfy
    \begin{align}
       V_{K+1} & \leq \exp\left(-\frac{\gamma\mu}{2}(K+1)\right) V + 2\gamma \sum\limits_{k=0}^{K}\exp\left(-\frac{\gamma\mu}{2}(K-k)\right)\langle \theta_k, x^k - x^* \rangle\notag\\
       &\quad\quad + \frac{\gamma^2}{n^2}\sum\limits_{k=0}^{K}\sum\limits_{i=1}^n\exp\left(-\frac{\gamma\mu}{2}(K-k)\right)\left(18\|\theta_{i,k}^u\|^2 + 14\|\omega_{i,k}^u\|^2\right)\notag\\
        &\quad\quad  + \frac{\gamma^2}{n}\sum\limits_{k=0}^{K}\sum\limits_{i=1}^n\exp\left(-\frac{\gamma\mu}{2}(K-k)\right)\left(18\|\theta_{i,k}^b\|^2 + 14\|\omega_{i,k}^b\|^2\right)\notag\\
        &\quad\quad + \frac{32\gamma^2}{n^2}\sum\limits_{k=0}^{K}\sum\limits_{j=2}^n\exp\left(-\frac{\gamma\mu}{2}(K-k)\right)\left\langle \sum\limits_{i=1}^{j-1} \theta_{i,k}^u, \theta_{j,k}^u \right\rangle \notag\\
        &\quad\quad + \frac{24\gamma^2}{n^2}\sum\limits_{k=0}^{K}\sum\limits_{j=2}^n\exp\left(-\frac{\gamma\mu}{2}(K-k)\right)\left\langle \sum\limits_{i=1}^{j-1} \omega_{i,k}^u, \omega_{j,k}^u \right\rangle, \label{ch3:eq:optimization_lemma_SEG_QSM_chapter_3}
    \end{align}
    where $V_k = \|x^k - x^*\|^2 + \frac{36\cdot 10^6\gamma^2A^2}{n^2}\sum_{i=1}^n\left(\|\thesischthreethh_{i}^k - F_i(x^*)\|^2 + \|\thesischthreehh_{i}^k - F_i(x^*)\|^2\right)$ and $\theta_k$, $\theta_{i,k}^u$, $\theta_{i,k}^b$, $\omega_k$, $\omega_{i,k}^u$, $\omega_{i,k}^b$ are defined in \eqref{ch3:eq:theta_k_SEG_chapter_3}, \eqref{ch3:eq:omega_k_SEG_chapter_3}, and \eqref{ch3:eq:gap_thm_SEG_technical_4_chapter_3}-\eqref{ch3:eq:gap_thm_SEG_technical_5_chapter_3}.
\end{lemma}
\begin{proof}
    From \eqref{ch3:eq:important_optimization_ineq_SEG_chapter_3} with $u = x^*$ we have
    \begin{align*}
        \MoveEqLeft[1] \gamma\left(\langle F(\thesischthreetx^k), \thesischthreetx^k - x^* \rangle +  \Psi(\thesischthreetx^k) - \Psi(x^*)\right) \\
        & \leq \frac{1}{2}\|x^k - x^*\|^2 - \frac{1}{2}\|x^{k+1} - x^*\|^2  \notag\\
        &\quad\quad - \frac{1}{2}\left(1 - 6\gamma^2 L^2\right)\|\thesischthreetx^k - x^k\|^2 - \frac{1}{4}\|x^{k+1} - \thesischthreetx^k\|^2 \notag\\
        &\quad\quad + 3\gamma^2 \|\omega_k\|^2 + 3\gamma^2\|\theta_k\|^2  + \gamma \langle \theta_k, \thesischthreetx^k - x^* \rangle.
    \end{align*}
    Using quasi-strong monotonicity of $F$, the fact that $-F(x^*) \in \partial \Psi(x^*)$, and convexity of $\Psi(x^*)$, we derive
    \begin{align*}
        2\gamma\mu\|\thesischthreetx^k - x^*\|^2 & \leq 2\gamma \langle F(\thesischthreetx^k) - F(x^*), \thesischthreetx^k - x^* \rangle \\
        & \leq 2\gamma\left(\langle F(\thesischthreetx^k), \thesischthreetx^k - x^* \rangle +  \Psi(\thesischthreetx^k) - \Psi(x^*)\right)\\
        & \leq \|x^k - x^*\|^2 - \|x^{k+1} - x^*\|^2   - \left(1 - 6\gamma^2 L^2\right)\|\thesischthreetx^k - x^k\|^2\\
        & \quad\quad - \frac{1}{2}\|x^{k+1} - \thesischthreetx^k\|^2 + 6\gamma^2 \|\omega_k\|^2 + 6\gamma^2\|\theta_k\|^2 \\
        &\quad + 2\gamma \langle \theta_k, \thesischthreetx^k - x^* \rangle.
    \end{align*}
    Next, we apply $\|\thesischthreetx^k - x^*\|^2 \geq \frac{1}{2}\|x^k - x^*\|^2 - \|\thesischthreetx^k - x^k\|^2$ and rearrange the terms:
    \begin{align}
        \|x^{k+1} - x^*\|^2 & \leq \left(1 - \frac{\gamma\mu}{2}\right)\|x^k - x^*\|^2 - \gamma\mu\|\thesischthreetx^k - x^*\|^2 \notag \\
        &\quad - \left(1-\gamma\mu - 6\gamma^2 L^2\right)\|\thesischthreetx^k - x^k\|^2 \notag\\
        &\quad\quad  + 6\gamma^2 \|\omega_k\|^2 + 6\gamma^2\|\theta_k\|^2  + 2\gamma \langle \theta_k, \thesischthreetx^k - x^* \rangle\notag.
    \end{align}
    Since $2\gamma\langle \theta_k, \thesischthreetx^k - x^* \rangle = 2\gamma\langle \theta_k, x^k - x^* \rangle +  2\gamma\langle \theta_k, \thesischthreetx^k - x^k \rangle$ and $2\gamma\langle \theta_k, \thesischthreetx^k - x^k \rangle \leq 2\gamma^2\|\theta_k\|^2 + \frac{1}{2}\|\thesischthreetx^k - x^k\|^2$, we have
    \begin{align}
        \|x^{k+1} - x^*\|^2 & \leq \left(1 - \frac{\gamma\mu}{2}\right)\|x^k - x^*\|^2 - \gamma\mu\|\thesischthreetx^k - x^*\|^2 \notag \\
        &\quad - \left(\frac{1}{2}-\gamma\mu - 6\gamma^2 L^2\right)\|\thesischthreetx^k - x^k\|^2 \notag\\
        &\quad\quad  + 6\gamma^2 \|\omega_k\|^2 + 8\gamma^2\|\theta_k\|^2  + 2\gamma \langle \theta_k, x^k - x^* \rangle. \label{ch3:eq:nshjcbjdbhcjsd_chapter_3}
    \end{align}
    Now, we move on to the shifts: for all $i\in [n]$ (for convenience, we use the new notation: $h_i^* = F_i(x^*)$ for all $i\in [n]$)
    \begin{align}
        \|\thesischthreethh_i^{k+1}-h^*_i\|^2 & = \|\thesischthreethh^{k}_i-h^*_i\|^2 \notag \\
        &\quad + 2\nu \left\la \tilde \Delta^k_i, \thesischthreethh^k_i-h^*_i\right\ra + \nu^2\|\tilde \Delta^k_i\|^2 \notag\\
    & = \|\thesischthreethh^{k}_i-h^*_i\|^2 + 2\nu \left\la \thesischthreetg_i^k -\thesischthreethh^k_i, \thesischthreethh^k_i-h^*_i\right\ra + \nu^2\|\thesischthreetg_i^k \notag \\
    &\quad -\thesischthreethh^k_i\|^2 \notag\\
    & \leq \|\thesischthreethh^{k}_i-h^*_i\|^2 + 2\nu \left\la \thesischthreetg_i^k-\thesischthreethh^k_i, \thesischthreethh^k_i-h^*_i\right\ra + \nu\|\thesischthreetg_i^k-\thesischthreethh^k_i\|^2\notag\\
    & =\|\thesischthreethh^{k}_i-h^*_i\|^2 + \nu \left\la \thesischthreetg_i^k-\thesischthreethh^k_i,  \thesischthreetg_i^k + \thesischthreethh^k_i- 2h^*_i\right\ra\notag\\
    & \leq (1-\nu)\|\thesischthreethh^{k}_i-h^*_i\|^2 + \nu\|\thesischthreetg_i^k-h^*_i\|^2 \notag\\
    & \leq (1-\nu)\|\thesischthreethh^k_i-h^*_i\|^2 + 2\nu\|\thesischthreetg_i^k- F_i(x^k)\|^2 \notag \\
    &\quad + 2\nu\|F_i(x^{k})-h^*_i\|^2 \notag\\
    & \leq (1-\nu)\|\thesischthreethh^{k}_i-h^*_i\|^2 + 2\nu\|\omega_{i,k}\|^2 \notag \\
    &\quad + 2\nu\|F_i(x^{k})-h_i^*\|^2 \notag\\
    & = (1-\gamma\mu)\|\thesischthreethh^{k}_i-h^*_i\|^2 + 2\gamma\mu\|\omega_{i,k}\|^2 \notag \\
    &\quad + 2\gamma\mu\|F_i(x^{k})-h_i^*\|^2 \notag\\
    & \overset{\eqref{ch3:eq:L_Lip_chapter_3}}{\leq}  (1-\gamma\mu)\|\thesischthreethh^{k}_i-h^*_i\|^2 + 2\gamma\mu\|\omega_{i,k}\|^2 + 2\gamma\mu L^2\|x^{k} \notag \\
    &\quad - x^*\|^2 \notag\\
    & \leq (1-\gamma\mu)\|\thesischthreethh^{k}_i-h^*_i\|^2 + 2\gamma\mu\|\omega_{i,k}\|^2 + 4\gamma\mu L^2\|\thesischthreetx^{k} \notag \\
    &\quad - x^*\|^2 \notag\\
    &\quad\quad +  4\gamma\mu L^2\|\thesischthreetx^{k} - x^k\|^2\label{ch3:eq:nvbhjsdbhdbfhjvbdjhvdjbvcs_chapter_3}
    \end{align}
    where in the first inequality we use $\nu \leq 1$, and, similarly,
    \begin{align}
        \|\thesischthreehh_i^{k+1}-h^*_i\|^2 & \leq (1-\gamma\mu)\|\thesischthreehh^{k}_i-h^*_i\|^2 \notag \\
        &\quad + 2\gamma\mu\|\theta_{i,k}\|^2 + 2\gamma\mu\|F_i(\thesischthreetx^{k})-h_i^*\|^2 \notag\\
        & \overset{\eqref{ch3:eq:L_Lip_chapter_3}}{\leq} (1-\gamma\mu)\|\thesischthreehh^{k}_i-h^*_i\|^2 + 2\gamma\mu\|\theta_{i,k}\|^2 + 2\gamma\mu L^2\|\thesischthreetx^{k} \notag \\
        &\quad -x^*\|^2. \label{ch3:eq:nvbhjsdbhdbfhjvbdjhvdjbvcs_2_chapter_3}
    \end{align}

    Summing up \eqref{ch3:eq:nshjcbjdbhcjsd_chapter_3}, \eqref{ch3:eq:nvbhjsdbhdbfhjvbdjhvdjbvcs_chapter_3}, and \eqref{ch3:eq:nvbhjsdbhdbfhjvbdjhvdjbvcs_2_chapter_3}, we derive
    \begin{align}
        V_{k+1} & \leq \left(1 - \frac{\gamma\mu}{2}\right)\|x^k - x^*\|^2\notag\\
        & \quad\quad + (1 - \gamma\mu) \frac{36\cdot 10^6\gamma^2A^2}{n^2}\notag \\
        &\qquad\quad \cdot\sum_{i=1}^n\left(\|\thesischthreethh_i^k - h_i^*\|^2 + \|\thesischthreehh_i^k - h_i^*\|^2\right) \notag\\
        &\quad\quad - \left(\gamma\mu - \frac{216\cdot 10^6 \gamma^3\mu L^2A^2}{n}\right)\|\thesischthreetx^k - x^*\|^2\notag\\
        &\quad\quad - \left(\frac{1}{2}-\gamma\mu - 6\gamma^2 L^2 - \frac{144\cdot 10^6\gamma^3\mu L^2A^2}{n}\right)\|\thesischthreetx^k - x^k\|^2 \notag\\
        &\quad\quad  + 6\gamma^2 \|\omega_k\|^2 + 8\gamma^2\|\theta_k\|^2  + 2\gamma \langle \theta_k, x^k - x^* \rangle\notag\\
        &\quad\quad + \frac{72\cdot10^6\gamma^3\mu A^2}{n^2}\sum_{i=1}^n\left(\|\theta_{i,k}\|^2 + \|\omega_{i,k}\|^2\right) \notag\\
        & \leq \left(1 - \frac{\gamma\mu}{2}\right) V_k + 6\gamma^2 \|\omega_k\|^2 + 8\gamma^2\|\theta_k\|^2  + 2\gamma \langle \theta_k, x^k - x^* \rangle\notag\\
        &\quad\quad + \frac{\gamma^2}{n^2}\sum\limits_{i=1}^n\left(\|\theta_{i,k}\|^2 + \|\omega_{i,k}\|^2\right) \notag\\
        & \overset{\eqref{ch3:eq:SEG_sum_theta_representation_chapter_3},\eqref{ch3:eq:gap_thm_SEG_technical_4_full_theta_omega_chapter_3}}{\leq} \exp\left(-\frac{\gamma\mu}{2}\right) V_k + 2\gamma \langle \theta_k, x^k - x^* \rangle \notag \\
        &\quad + \frac{\gamma^2}{n^2}\sum\limits_{i=1}^n\left(18\|\theta_{i,k}^u\|^2 + 14\|\omega_{i,k}^u\|^2\right) \notag\\
        & \quad\quad  + \frac{\gamma^2}{n}\sum\limits_{i=1}^n\left(18\|\theta_{i,k}^b\|^2 + 14\|\omega_{i,k}^b\|^2\right) \notag \\
        &\quad + \frac{32\gamma^2}{n^2}\sum\limits_{j=2}^n\left\langle \sum\limits_{i=1}^{j-1} \theta_{i,k}^u, \theta_{j,k}^u \right\rangle \notag\\
        &\quad\quad + \frac{24\gamma^2}{n^2}\sum\limits_{j=2}^n\left\langle \sum\limits_{i=1}^{j-1} \omega_{i,k}^u, \omega_{j,k}^u \right\rangle.\notag
    \end{align}
    Here the second inequality uses $\gamma \leq \sqrt{n}/(15000LA)$, $\gamma \leq 1/(6L)$ and $72\cdot 10^6\gamma\mu A^2 \leq 1$, which make the coefficients of $\|\thesischthreetx^k - x^*\|^2$ and $\|\thesischthreetx^k - x^k\|^2$ nonpositive and give $72\cdot 10^6\gamma^3\mu A^2/n^2 \leq \gamma^2/n^2$, and the last one uses $1 - \gamma\mu/2 \leq \exp(-\gamma\mu/2)$ and $16/n + 2/n^2 \leq 18/n$, $12/n + 2/n^2 \leq 14/n$. Unrolling the recurrence and using $V_0 \leq V$ (the initial shifts are zero), we get the result.
\end{proof}

Next, we proceed with the full statement of our main result for \algname{DProx-clipped-SEG-shift} in the quasi-strongly monotone case.

\begin{theorem}[Case 1 from Theorem~\ref{ch3:thm:D_prox_clipped_SEG_main_chapter_3}]\label{ch3:thm:D_prox_clipped_SEG_QSM_appendix_chapter_3}
    Let Assumptions \ref{ch3:as:bounded_alpha_moment_chapter_3}, \ref{ch3:as:L_Lip_chapter_3} and \ref{ch3:as:str_monotonicity_chapter_3} hold for $Q = B_{4n\sqrt{V}}(x^*)$, where
    \begin{gather*}
        A = \ln\frac{48n(K+1)}{\beta},\quad V \geq \|x^0 - x^*\|^2 + \frac{72\cdot 10^6\gamma^2A^2}{n^2}\sum_{i=1}^n\|F_i(x^*)\|^2,
    \end{gather*}
    let the initial shifts be $\thesischthreethh_i^0 = \thesischthreehh_i^0 = 0$ for all $i \in [n]$, and let
    \begin{align}
        0< \gamma &\leq \min\left\{\frac{1}{72\cdot 10^6 \mu A^2}, \frac{1}{6L}, \frac{\sqrt{n}}{15000 L A}, \frac{2\ln(B_K)}{\mu(K+1)}\right\}, \label{ch3:eq:gamma_SEG_str_mon_chapter_3}\\
        Q_K &= \frac{n^{2(p-1)/p}(K+1)^{2(p-1)/p}\mu^2V}{\sigma^2A^{2(p-1)/p}}, \notag\\
        B_K\ln^2(B_K) &= \frac{Q_K}{(3\cdot 10^7)^{2/p}}, \label{ch3:eq:B_K_SEG_str_mon_1_chapter_3} \\
        B_K &= \Theta\left(\frac{Q_K}{\ln^2 Q_K}\right), \label{ch3:eq:B_K_SEG_str_mon_2_chapter_3} \\
        \lambda_k &= \frac{n\exp(-\gamma\mu(1 + k/4))\sqrt{V}}{300\gamma A}, \label{ch3:eq:lambda_SEG_str_mon_chapter_3}\\
        \nu &= \gamma\mu \label{ch3:eq:nu_SEG_str_mon_chapter_3}
    \end{align}
    for some $K \geq 1$ and $\beta \in (0,1]$ such that $B_K \geq 2$. Then, after $K$ iterations of \algname{DProx-clipped-SEG-shift}, the following inequality holds with probability at least $1 - \beta$:
    \begin{equation}
        \|x^{K+1} - x^*\|^2 \leq 2\exp\left(-\frac{\gamma\mu(K+1)}{2}\right)V. \label{ch3:eq:main_result_str_mon_chapter_3}
    \end{equation}
    In particular, when $\gamma$ equals the minimum from \eqref{ch3:eq:gamma_SEG_str_mon_chapter_3}, then after $K$ iterations of \algname{DProx-clipped-SEG-shift}, we have with probability at least $1-\beta$ that 
    \begin{eqnarray*}
        \|x^{K+1} - x^*\|^2 &=& \cO\biggl(\max\biggl\{V\exp\left(- \frac{K}{\ln^2(nK/\beta)}\right), V\exp\left(- \frac{\mu  K}{L}\right),\notag \\
        && V\exp\left(- \frac{\mu \sqrt{n} K}{L \ln(nK/\beta)}\right), \frac{\sigma^2\ln^{2(p-1)/p}(nK/\beta)\ln^2B_K}{n^{2(p-1)/p}K^{2(p-1)/p}\mu^2}\biggr\}\biggr),\label{ch3:eq:clipped_SEG_str_monotone_case_2_appendix_chapter_3}
    \end{eqnarray*}
    i.e., to achieve $\|x^{K} - x^*\|^2 \leq \varepsilon$ with probability at least $1 - \beta$ \algname{DProx-clipped-SEG-shift} needs
    \begin{align}
        K = \cO\biggl(\max\biggl\{ & \left(\frac{L}{\sqrt{n}\mu}+\ln\left(\frac{nL}{\mu \beta}\ln\frac{V}{\varepsilon}\right)\right)\ln\frac{V}{\varepsilon}\ln\left(\frac{nL}{\mu \beta}\ln\frac{V}{\varepsilon}\right), \frac{L}{\mu}\ln\frac{V}{\varepsilon},\notag\\
        & \frac{1}{n}\left(\frac{\sigma^2}{\mu^2\varepsilon}\right)^{p/(2(p-1))}\ln \left(\frac{n}{\beta} \left(\frac{\sigma^2}{\mu^2\varepsilon}\right)^{p/(2(p-1))}\right) \notag\\
        & \qquad\cdot\ln^{p/(p-1)}\left(B_\varepsilon\right)\biggr\}\biggr) \label{ch3:eq:clipped_SEG_str_monotone_case_complexity_appendix_chapter_3}
    \end{align}
    iterations/oracle calls per worker, where
    \begin{equation*}
        B_\varepsilon = \max\left\{2, \frac{V}{\varepsilon \ln \left(\frac{1}{\beta} \left(\sigma^2/(\mu^2\varepsilon)\right)^{p/(2(p-1))}\right)}\right\}.
    \end{equation*}
\end{theorem}
\begin{proof}
    Similar to previous results, our proof is induction-based.  To formulate the statement rigorously, we introduce probability event $E_k$ for each $k = 0,1,\ldots, K+1$ as follows: inequalities
    \begin{gather}
        V_t \leq 2 \exp\left(-\frac{\gamma\mu t}{2}\right) V \label{ch3:eq:induction_inequality_str_mon_SEG_chapter_3}\\
         \left\|\frac{\gamma}{n}\sum\limits_{i=1}^{r-1}\theta_{i,t-1}^u\right\| \leq \exp\left(-\frac{\gamma\mu (t-1)}{4}\right)\frac{\sqrt{V}}{2}, \label{ch3:eq:induction_inequality_str_mon_SEG_1_chapter_3}\\
         \left\|\frac{\gamma}{n}\sum\limits_{i=1}^{r-1}\omega_{i,t-1}^u\right\| \leq \exp\left(-\frac{\gamma\mu (t-1)}{4}\right)\frac{\sqrt{V}}{2} \label{ch3:eq:induction_inequality_str_mon_SEG_2_chapter_3}
    \end{gather}
    hold for $t = 0,1,\ldots,k$ and $r = 1, 2, \ldots, n$ simultaneously. We will prove by induction that $\Prob\{E_k\} \geq  1 - k\beta/(K+1)$ for all $k = 0,1,\ldots,K+1$. The base of induction follows immediately by the definition of $V$. Next, we assume that the statement holds for $k = T-1 \leq K$, i.e., $\Prob\{E_{T-1}\} \geq  1 - (T-1)\beta/(K+1)$. Let us show that it also holds for $k = T$, i.e., $\Prob\{E_{T}\} \geq  1 - T\beta/(K+1)$.

    Similarly to the monotone case, one can show that due to our choice of the clipping level, we have that $E_{T-1}$ implies $x^t, \thesischthreetx^t \in B_{4n\sqrt{V}}(x^*)$ for $t = 0, \ldots, T-1$. Indeed, for $t = 0,1, \ldots, T-1$ inequality \eqref{ch3:eq:induction_inequality_str_mon_SEG_chapter_3} gives $x^t \in B_{2\sqrt{V}}(x^*)$. Next, for $\thesischthreetx^t$, $t = 0, \ldots, T-1$, event $E_{T-1}$ implies
    \begin{align*}
        \|\thesischthreetx^t - x^*\| & = \|\thesischthreeprox_{\gamma\Psi}(x^t - \gamma \thesischthreetg^t) - \thesischthreeprox_{\gamma \Psi}(x^* - \gamma F(x^*))\| \\
        & \leq \|x^t - x^* - \gamma (\thesischthreetg^t - F(x^*))\| \\
        & \leq \|x^t - x^*\| + \gamma \left\|\frac{1}{n}\sum\limits_{i=1}^n (\thesischthreethh_i^t - F_i(x^*)) + \frac{1}{n}\sum\limits_{i=1}^n \tilde\Delta_i^t \right\| \\
        & \leq 2\sqrt{V} + \gamma \sqrt{\frac{1}{n}\sum\limits_{i=1}^n \|\thesischthreethh_i^t - F_i(x^*)\|^2} + \frac{\gamma}{n}\sum\limits_{i=1}^n \|\tilde \Delta_i^t\|\\
        & \leq 2\sqrt{V} + \frac{\sqrt{n}}{6000A}\sqrt{V_t} + \gamma \lambda_t\\
        & \overset{\eqref{ch3:eq:induction_inequality_str_mon_SEG_chapter_3}}{\leq} \left(2 + \frac{20n + \sqrt{2n}}{6000A}\right)\sqrt{V} \leq 4n\sqrt{V}.
    \end{align*}

    This means that we can apply Lemma~\ref{ch3:lem:SEG_optimization_lemma_QSM_chapter_3}: $E_{T-1}$ implies
    \begin{align}
       V_{T} & \leq \exp\left(-\frac{\gamma\mu}{2}T\right) V + 2\gamma \sum\limits_{t=0}^{T-1}\exp\left(-\frac{\gamma\mu}{2}(T-1-t)\right)\langle \theta_t, x^t - x^* \rangle\notag\\
       &\quad\quad + \frac{\gamma^2}{n^2}\sum\limits_{t=0}^{T-1}\sum\limits_{i=1}^n\exp\left(-\frac{\gamma\mu}{2}(T-1-t)\right)\left(18\|\theta_{i,t}^u\|^2 + 14\|\omega_{i,t}^u\|^2\right)\notag\\
        &\quad\quad  + \gamma^2\sum\limits_{t=0}^{T-1}\sum\limits_{i=1}^n\exp\left(-\frac{\gamma\mu}{2}(T-1-t)\right)\left(18\|\theta_{i,t}^b\|^2 + 14\|\omega_{i,t}^b\|^2\right)\notag\\
        &\quad\quad + \frac{32\gamma^2}{n^2}\sum\limits_{t=0}^{T-1}\sum\limits_{j=2}^n\exp\left(-\frac{\gamma\mu}{2}(T-1-t)\right)\left\langle \sum\limits_{i=1}^{j-1} \theta_{i,t}^u, \theta_{j,t}^u \right\rangle \notag\\
        &\quad\quad + \frac{24\gamma^2}{n^2}\sum\limits_{t=0}^{T-1}\sum\limits_{j=2}^n\exp\left(-\frac{\gamma\mu}{2}(T-1-t)\right)\left\langle \sum\limits_{i=1}^{j-1} \omega_{i,t}^u, \omega_{j,t}^u \right\rangle. \notag
    \end{align}
    Before we proceed, we introduce a new notation:
    \begin{equation*}
        u^{t} =
        \begin{cases}
        x^t - x^*,&
        \begin{aligned}
        &\text{if } \|x^t - x^*\|\\[-0.2em]
        &\quad \leq \exp\left(-\frac{\gamma\mu t}{4}\right)\sqrt{2V},
        \end{aligned}\\
        0,& \text{otherwise,}
        \end{cases}
    \end{equation*}
    for all $t = 0, 1, \ldots, T-1$. Random vectors $\{u^{t}\}_{t=0}^{T-1}$ are bounded almost surely:
     \begin{equation}
        \|u^{t}\| \leq \exp\left(-\frac{\gamma\mu t}{4}\right)\sqrt{2V}  \label{ch3:eq:eta_t_bound_SEG_str_mon_chapter_3}
    \end{equation}
    for all $t = 0, 1, \ldots, T-1$. In addition, $u^{t} = x^t - x^*$ follows from $E_{T-1}$ for all $t = 0, 1, \ldots, T-1$ (since $\|x^t - x^*\|^2 \leq V_t$) and, thus, $E_{T-1}$ implies
    \begin{align}
       V_{T}
       &\leq \exp\left(-\frac{\gamma\mu}{2}T\right)V
       +\circledOne+\circledTwo+\circledThree+\circledFour
       \notag\\
       &\quad+\circledFive+\circledSix+\circledSeven.
       \label{ch3:eq:ndvjinuibvbsvdhbvdhbfussdncdni_chapter_3}\\
       \circledOne
       &\eqdef
       \frac{2\gamma}{n}
       \sum\limits_{l=0}^{T-1}\sum\limits_{i=1}^n
       \exp\left(-\frac{\gamma\mu}{2}(T-1-l)\right)
       \langle \theta_{i,l}^u, u^{l} \rangle,
       \notag\\
       \circledTwo
       &\eqdef
       2\gamma\sum\limits_{l=0}^{T-1}
       \exp\left(-\frac{\gamma\mu}{2}(T-1-l)\right)
       \langle \theta_{l}^b, u^{l} \rangle,
       \notag\\
       \circledThree
       &\eqdef
       \frac{\gamma^2}{n^2}
       \sum\limits_{l=0}^{T-1}\sum\limits_{i=1}^n
       \exp\left(-\frac{\gamma\mu}{2}(T-1-l)\right)
       \notag\\[-0.25em]
       &\qquad {}\times
       \left(
       18\Exp_{\xi_{2,i}^l}\left[\|\theta_{i,l}^u\|^2\right]
       +14\Exp_{\xi_{1,i}^l}\left[\|\omega_{i,l}^u\|^2\right]
       \right),
       \notag\\
       \circledFour
       &\eqdef
       \frac{\gamma^2}{n^2}
       \sum\limits_{l=0}^{T-1}\sum\limits_{i=1}^n
       \exp\left(-\frac{\gamma\mu}{2}(T-1-l)\right)
       \notag\\[-0.25em]
       & \qquad {}\times
       \Bigl(
       18\|\theta_{i,l}^u\|^2+14\|\omega_{i,l}^u\|^2 \notag \\
       &\quad -18\Exp_{\xi_{2,i}^l}\left[\|\theta_{i,l}^u\|^2\right] \notag\\[-0.25em]
       &\hspace{15em}
       -14\Exp_{\xi_{1,i}^l}\left[\|\omega_{i,l}^u\|^2\right]
       \Bigr),
       \notag\\
       \circledFive
       &\eqdef
       \frac{\gamma^2}{n}
       \sum\limits_{l=0}^{T-1}\sum\limits_{i=1}^n
       \exp\left(-\frac{\gamma\mu}{2}(T-1-l)\right)
       \left(18\|\theta_{i,l}^b\|^2+14\|\omega_{i,l}^b\|^2\right),
       \notag\\
       \circledSix
       &\eqdef
       \frac{32\gamma^2}{n^2}
       \sum\limits_{l=0}^{T-1}\sum\limits_{j=2}^n
       \exp\left(-\frac{\gamma\mu}{2}(T-1-l)\right)
       \left\langle
       \sum\limits_{i=1}^{j-1}\theta_{i,l}^u,\theta_{j,l}^u
       \right\rangle,
       \notag\\
       \circledSeven
       &\eqdef
       \frac{24\gamma^2}{n^2}
       \sum\limits_{l=0}^{T-1}\sum\limits_{j=2}^n
       \exp\left(-\frac{\gamma\mu}{2}(T-1-l)\right)
       \left\langle
       \sum\limits_{i=1}^{j-1}\omega_{i,l}^u,\omega_{j,l}^u
       \right\rangle.
       \notag
    \end{align}

    To derive high-probability bounds for $\circledOne, \circledTwo, \circledThree, \circledFour, \circledFive, \circledSix, \circledSeven$ we need to establish several useful inequalities related to $\theta_{i,l}^u, \theta_{i,l}^b, \omega_{i,l}^u, \omega_{i,l}^b$. First, by definition of clipping
    \begin{equation}
        \|\theta_{i,l}^u\| \leq 2\lambda_l,\quad \|\omega_{i,l}^u\| \leq 2\lambda_l. \label{ch3:eq:theta_omega_magnitude_str_mon_chapter_3}
    \end{equation}
    Next, we notice that $E_{T-1}$ implies
    \begin{align*}
        \|F_i(x^l) - \thesischthreethh_i^l\| & \leq \|F_i(x^l) - F_i(x^*)\| + \|\thesischthreethh_i^l - F_i(x^*)\|\\
        & \overset{\eqref{ch3:eq:L_Lip_chapter_3}}{\leq} L\|x^l - x^*\| + \sqrt{\sum\limits_{j=1}^n \|\thesischthreethh_i^l - F_i(x^*)\|^2}\\
        & \leq L\sqrt{V_l} + \frac{n\sqrt{V_l}}{6000\gamma A}\\
        & \overset{\eqref{ch3:eq:induction_inequality_str_mon_SEG_chapter_3}}{\leq} \sqrt{2}\left(L + \frac{n}{6000\gamma A}\right)\exp\left(-\frac{\gamma\mu l}{4}\right)\sqrt{V} \overset{\eqref{ch3:eq:gamma_SEG_str_mon_chapter_3}, \eqref{ch3:eq:lambda_SEG_str_mon_chapter_3}}{\leq} \frac{\lambda_l}{2}
    \end{align*}
    and
    \begin{align*}
        \|F_i(\thesischthreetx^l) - \thesischthreehh_i^l\|
        & \leq \|F_i(\thesischthreetx^l) - F_i(x^*)\| + \|\thesischthreehh_i^l \\
        &\quad - F_i(x^*)\|\\
        &\overset{\eqref{ch3:eq:L_Lip_chapter_3}}{\leq}
        L\|\thesischthreetx^l - x^*\|
        + \sqrt{\sum\limits_{j=1}^n
        \|\thesischthreehh_i^l - F_i(x^*)\|^2}\\
        & \leq
        L\|\thesischthreeprox_{\gamma\Psi}(x^l - \gamma \thesischthreetg^l) - \thesischthreeprox_{\gamma \Psi}(x^* - \gamma F(x^*))\| \\
        &\quad + \frac{n\sqrt{V_l}}
        {6000\gamma A}\\
        &\leq L\|x^l - x^* - \gamma (\thesischthreetg^l - F(x^*))\|
        + \frac{n\sqrt{V_l}}
        {6000\gamma A}\\
        & \leq L\|x^l - x^*\| \\
        &\quad + L\gamma \left\|
        \frac{1}{n}\sum\limits_{i=1}^n
        (\thesischthreethh_i^l - F_i(x^*))
        + \frac{1}{n}\sum\limits_{i=1}^n \tilde\Delta_i^l
        \right\| \\
        &\quad + \frac{n\sqrt{V_l}}
        {6000\gamma A}\\
        &\leq
        \left(L + \frac{n}
        {6000\gamma A}\right)\sqrt{V_l}\\
        &\quad + L\gamma
        \sqrt{\frac{1}{n}\sum\limits_{i=1}^n
        \|\thesischthreethh_i^l - F_i(x^*)\|^2}
        + \frac{L\gamma}{n}\sum\limits_{i=1}^n
        \|\tilde \Delta_i^l\|\\
        &\leq
        \left(L + \frac{n + L\gamma\sqrt{n}}
        {6000\gamma A}\right)\sqrt{V_l}
        + L\gamma \lambda_t\\
        &\overset{\eqref{ch3:eq:induction_inequality_str_mon_SEG_chapter_3}}{\leq}
        \sqrt{2}\left(L + \frac{n + L\gamma\sqrt{n}}
        {6000\gamma A}\right)
        \exp\left(-\frac{\gamma\mu l}{4}\right)\sqrt{V}
        + L\gamma \lambda_l\\
        &\overset{\eqref{ch3:eq:gamma_SEG_str_mon_chapter_3},
        \eqref{ch3:eq:lambda_SEG_str_mon_chapter_3}}{\leq}
        \frac{\lambda_l}{2}
    \end{align*}
    for $l = 0,1,\ldots,T-1$ and $i\in [n]$. Therefore, one can apply Lemma~\ref{ch3:lem:bias_and_variance_clip_chapter_3} and get
    \begin{gather}
        \left\|\theta_l^b\right\| \leq \frac{1}{n}\sum\limits_{i=1}^n \|\theta_{i,l}^b\| \leq \frac{2^{p}\sigma^{p}}{\lambda_l^{{p}-1}},\quad \left\|\omega_l^b\right\| \leq \frac{1}{n}\sum\limits_{i=1}^n \|\omega_{i,l}^b\| \leq \frac{2^{p}\sigma^{p}}{\lambda_l^{{p}-1}}, \label{ch3:eq:bias_theta_omega_str_mon_chapter_3}\\
        \Exp_{\xi_{2,i}^l}\left[\left\|\theta_{i,l}^u\right\|^2\right] \leq 18 \lambda_l^{2-{p}}\sigma^{p},\quad \Exp_{\xi_{1,i}^l}\left[\left\|\omega_{i,l}^u\right\|^2\right] \leq 18 \lambda_l^{2-{p}}\sigma^{p}, \label{ch3:eq:variance_theta_omega_str_mon_chapter_3}
    \end{gather}
    for all $l = 0,1, \ldots, T-1$ and $i\in [n]$.

    \paragraph{Upper bound for $\circledOne$.} To estimate this sum, we will use Bernstein's inequality. The summands have conditional expectations equal to zero:
    \begin{align*}
\MoveEqLeft[1] \Exp_{\xi_{2,i}^l}\left[\frac{2\gamma}{n}\exp\left(-\frac{\gamma\mu}{2}(T-1-l)\right) \langle u^{l}, \theta_{i,l}^u \rangle\right] \\
& = \frac{2\gamma}{n} \exp\left(-\frac{\gamma\mu}{2}(T-1-l)\right)\left\langle u^{l}, \Exp_{\xi_{2,i}^l}[\theta_{i,l}^u] \right\rangle = 0.
    \end{align*}
    Moreover, for all $l = 0,\ldots, T-1$ random vectors $\{\theta_{i,l}^u\}_{i=1}^n$ are independent. Thus, the summands $\frac{2\gamma}{n}\exp\left(-\frac{\gamma\mu}{2}(T-1-l)\right) \langle u^{l}, \theta_{i,l}^u \rangle$ form a martingale difference sequence. Next, the summands are bounded:
    \begin{align}
        \MoveEqLeft[1] \left|\frac{2\gamma}{n} \exp\left(-\frac{\gamma\mu}{2}(T-1-l)\right)\langle u^{l}, \theta_{i,l}^u \rangle \right| \notag \\
        & \leq \frac{2\gamma}{n}\exp\left(-\frac{\gamma\mu}{2}(T-1-l)\right) \|u^{l}\|\cdot \|\theta_{i,l}^u\|\notag\\
        & \overset{\eqref{ch3:eq:eta_t_bound_SEG_str_mon_chapter_3},\eqref{ch3:eq:theta_omega_magnitude_str_mon_chapter_3}}{\leq} \frac{4\sqrt{2V}\gamma \exp\left(-\frac{\gamma\mu(T-1)}{2}\right)}{n} \exp\left(\frac{\gamma\mu l}{4}\right)\lambda_l\notag\\
        & \overset{\eqref{ch3:eq:lambda_SEG_str_mon_chapter_3}}{\leq} \frac{\exp\left(-\frac{\gamma\mu T}{2}\right)V}{50A} \eqdef c. \label{ch3:eq:str_mon_thm_SEG_technical_6_5_chapter_3}
    \end{align}
    Finally, conditional variances $\sigma_{i,l}^2 \eqdef \Exp_{\xi_{2,i}^l}\left[\frac{4\gamma^2}{n^2}\exp\left(-\gamma\mu(T-1-l)\right) \langle u^{l}, \theta_{i,l}^u \rangle^2\right]$ of the summands are bounded:
    \begin{align}
        \sigma_{i,l}^2 & \leq \Exp_{\xi_{2,i}^l}\left[\frac{4\gamma^2 }{n^2}\exp\left(-\gamma\mu(T-1-l)\right)\|u^{l}\|^2\cdot \|\theta_{i,l}^u\|^2\right] \notag\\
        & \overset{\eqref{ch3:eq:eta_t_bound_SEG_str_mon_chapter_3}}{\leq} \frac{8\gamma^2 V\exp\left(-\gamma\mu\left(T-1-\frac{l}{2}\right)\right)}{n^2} \Exp_{\xi_{2,i}^l}\left[\|\theta_{i,l}^u\|^2\right]. \label{ch3:eq:str_mon_thm_SEG_technical_7_chapter_3}
    \end{align}
    Applying Bernstein's inequality (Lemma~\ref{ch3:lem:Bernstein_ineq_chapter_3}) with $X_{i,l}$ equal to the summands above, constant $c$ defined in \eqref{ch3:eq:str_mon_thm_SEG_technical_6_5_chapter_3}, $b = \frac{\exp\left(-\frac{\gamma\mu T}{2}\right)V}{50}$, $H = \frac{\exp\left(-\gamma\mu T\right)V^2}{15000A}$, we get
    \begin{align*}
        \MoveEqLeft[1] \Prob\left\{|\circledOne| > \frac{\exp\left(-\frac{\gamma\mu T}{2}\right)V}{50} \text{ and } \sum\limits_{l=0}^{T-1}\sum\limits_{i=1}^n\sigma_{i,l}^2 \leq \frac{\exp\left(-\gamma\mu T\right)V^2}{15000A}\right\}\\
        & \leq 2\exp\left(- \frac{b^2}{2H + 2cb/3}\right)\\
        & = \frac{\beta}{24n(K+1)}.
    \end{align*}
    The above is equivalent to $\Prob\{E_{\circledOne}\} \geq 1 - \frac{\beta}{24n(K+1)}$ for
    \begin{equation}
        E_{\circledOne} = \left\{\text{either } \sum\limits_{l=0}^{T-1}\sum\limits_{i=1}^n\sigma_{i,l}^2 > \frac{\exp\left(-\gamma\mu T\right)V^2}{15000A} \text{ or } |\circledOne| \leq \frac{\exp\left(-\frac{\gamma\mu T}{2}\right)V}{50}\right\}. \label{ch3:eq:bound_1_str_mon_SEG_chapter_3}
    \end{equation}
    Moreover, $E_{T-1}$ implies
    \begin{align}
        \sum\limits_{l=0}^{T-1}\sum\limits_{i=1}^n\sigma_{i,l}^2 & \overset{\eqref{ch3:eq:str_mon_thm_SEG_technical_7_chapter_3}}{\leq} \frac{8\gamma^2 V\exp\left(-\gamma\mu\left(T-1\right)\right)}{n^2} \sum\limits_{l=0}^{T-1}\exp\left(\frac{\gamma\mu l}{2}\right) \sum\limits_{i=1}^n \Exp_{\xi_{2,i}^l}\left[\|\theta_{i,l}^u\|^2\right]\notag\\
        & \overset{\eqref{ch3:eq:variance_theta_omega_str_mon_chapter_3}}{\leq}  \frac{144\gamma^2 V\exp\left(-\gamma\mu\left(T-1\right)\right)\sigma^{p}}{n} \sum\limits_{l=0}^{T-1}\exp\left(\frac{\gamma\mu l}{2}\right) \lambda_l^{2-{p}}\notag\\
        & \overset{\eqref{ch3:eq:lambda_SEG_str_mon_chapter_3}}{\leq} \frac{144\gamma^{p} V^{2 - \frac{{p}}{2}}\exp\left(-\gamma\mu (T-1)\right)\sigma^{p}}{300^{2-{p}} n^{{p}-1} A^{2-{p}}} \sum\limits_{l=0}^{T-1}\exp\left(\frac{\gamma\mu l {p}}{4}\right)\notag\\
        & \leq \frac{144\gamma^{p} V^{2 - \frac{{p}}{2}}\exp\left(-\gamma\mu (T-1)\right)\sigma^{p} (K+1) \exp\left(\frac{\gamma\mu K {p}}{4}\right)}{300^{2-{p}} n^{{p}-1} A^{2-{p}}}\notag\\
        & \overset{\eqref{ch3:eq:gamma_SEG_str_mon_chapter_3}}{\leq} \frac{\exp\left(-\gamma\mu T\right)V^2}{15000A}.\label{ch3:eq:bound_1_variances_str_mon_SEG_chapter_3}
    \end{align}

    \paragraph{Upper bound for $\circledTwo$.} Probability event $E_{T-1}$ implies
    \begin{align}
        \circledTwo & \leq 2\gamma \sum\limits_{l=0}^{T-1}\exp\left(-\frac{\gamma\mu(T-1-l)}{2}\right)\|u^{l}\| \cdot \|\theta_l^b\| \\
        & \overset{\eqref{ch3:eq:eta_t_bound_SEG_str_mon_chapter_3}, \eqref{ch3:eq:bias_theta_omega_str_mon_chapter_3}}{\leq} 2^{{p}+1}\gamma\sigma^{p}\sqrt{2V}\exp\left(-\frac{\gamma\mu(T-1)}{2}\right) \sum\limits_{l=0}^{T-1}\exp\left(\frac{\gamma\mu l}{4}\right)\frac{1}{\lambda_l^{{p}-1}}\notag\\
        & \overset{\eqref{ch3:eq:lambda_SEG_str_mon_chapter_3}}{\leq} \frac{2\sqrt{2}\cdot 2^{{p}}\cdot 300^{{p}-1} \exp\left(-\frac{\gamma\mu (T-1)}{2}\right) e^{\gamma\mu({p}-1)} (K+1) \exp\left(\frac{\gamma\mu K {p}}{4}\right)\gamma^{p} \sigma^{p} A^{{p}-1}}{n^{{p}-1}V^{\frac{{p}}{2}-1}} \notag\\
        & \overset{\eqref{ch3:eq:gamma_SEG_str_mon_chapter_3}}{\leq} \frac{3\exp\left(-\frac{\gamma\mu T}{2}\right) V}{100}. \label{ch3:eq:bound_2_variances_str_mon_SEG_chapter_3}
    \end{align}

    \paragraph{Upper bound for $\circledThree$.} Probability event $E_{T-1}$ implies
    \begin{align}
        &\frac{18\gamma^2}{n^2}
        \sum\limits_{l=0}^{T-1}\sum\limits_{i=1}^n
        \exp\left(-\frac{\gamma\mu(T-1-l)}{2}\right)
        \Exp_{\xi_{2,i}^l}\left[\|\theta_{i,l}^u\|^2\right]
        \notag\\
        &\quad
        \overset{\eqref{ch3:eq:bound_1_variances_str_mon_SEG_chapter_3}}{\leq}
        \frac{\exp\left(-\frac{\gamma\mu T}{2}\right)V}{100}
        \notag
    \end{align}
    and, similarly,
    \begin{align}
        &\frac{14\gamma^2}{n^2}
        \sum\limits_{l=0}^{T-1}\sum\limits_{i=1}^n
        \exp\left(-\frac{\gamma\mu(T-1-l)}{2}\right)
        \Exp_{\xi_{1,i}^l}\left[\|\omega_{i,l}^u\|^2\right]
        \notag\\
        &\quad
        \overset{\eqref{ch3:eq:bound_1_variances_str_mon_SEG_chapter_3}}{\leq}
        \frac{\exp\left(-\frac{\gamma\mu T}{2}\right)V}{100}
        \notag
    \end{align}
    that give
    \begin{equation}
        \circledThree \leq \frac{\exp\left(-\frac{\gamma\mu T}{2}\right)V}{50}. \label{ch3:eq:bound_3_variances_str_mon_SEG_chapter_3}
    \end{equation}

    \paragraph{Upper bound for $\circledFour$.} 
    To estimate this sum, we will use Bernstein's inequality. Let
    \begin{align*}
        Z_{i,l} \eqdef{} & 18\|\theta_{i,l}^u\|^2 + 14\|\omega_{i,l}^u\|^2\\
        & - 18\Exp_{\xi_{2,i}^l}\left[\|\theta_{i,l}^u\|^2\right] - 14\Exp_{\xi_{1,i}^l}\left[\|\omega_{i,l}^u\|^2\right],
    \end{align*}
    so that the summands of $\circledFour$ are $(\gamma^2/n^2)\exp(-\gamma\mu(T-1-l)/2)Z_{i,l}$. They have conditional expectations equal to zero, since $\Exp_{\xi_{1,i}^l,\xi_{2,i}^l}[Z_{i,l}] = 0$. Moreover, for all $l = 0,\ldots, T-1$ random vectors $\{\theta_{i,l}^u\}_{i=1}^n$, $\{\omega_{i,l}^u\}_{i=1}^n$ are independent. Thus, these summands form a martingale difference sequence. Next, they are bounded:
    \begin{align}
        \frac{\gamma^2}{n^2}\exp\left(-\frac{\gamma\mu}{2}(T-1-l)\right)|Z_{i,l}| &\overset{\eqref{ch3:eq:theta_omega_magnitude_str_mon_chapter_3}}{\leq} \frac{256\gamma^2 \lambda_l^2}{n^2}\exp\left(-\frac{\gamma\mu}{2}(T-1-l)\right) \notag\\
        &\overset{\eqref{ch3:eq:lambda_SEG_str_mon_chapter_3}}{\leq} \frac{\exp\left(-\frac{\gamma\mu T}{2}\right)V}{6A} \eqdef c.\label{ch3:eq:str_mon_thm_SEG_technical_6_5_1_chapter_3}
    \end{align}
    Finally, conditional variances $\widetilde\sigma_{i,l}^2 \eqdef (\gamma^4/n^4)\exp(-\gamma\mu(T-1-l))\Exp_{\xi_{1,i}^l,\xi_{2,i}^l}[Z_{i,l}^2]$ of the summands are bounded:
    \begin{align}
        \widetilde\sigma_{i,l}^2 & \overset{\eqref{ch3:eq:str_mon_thm_SEG_technical_6_5_1_chapter_3}}{\leq} \frac{\gamma^2\exp\left(-\frac{\gamma\mu }{2}(2T-1-l)\right)V}{6n^2A}\Exp_{\xi_{1,i}^l, \xi_{2,i}^l}\left[|Z_{i,l}|\right]\notag\\
        & \leq \frac{2\gamma^2\exp\left(-\frac{\gamma\mu }{2}(2T-1-l)\right)V}{3n^2A} \Exp_{\xi_{1,i}^l, \xi_{2,i}^l}\left[9\|\theta_{i,l}^u\|^2 + 7\|\omega_{i,l}^u\|^2 \right]. \label{ch3:eq:str_mon_thm_SEG_technical_7_1_chapter_3}
    \end{align}
    Applying Bernstein's inequality (Lemma~\ref{ch3:lem:Bernstein_ineq_chapter_3}) with $X_{i,l}$ equal to the summands of $\circledFour$, constant $c$ defined in \eqref{ch3:eq:str_mon_thm_SEG_technical_6_5_1_chapter_3}, $b = \frac{\exp\left(-\frac{\gamma\mu T}{2}\right)V}{5}$, $H = \frac{\exp\left(-\gamma\mu T\right)V^2}{150A}$, we get
    \begin{align*}
        \MoveEqLeft[1] \Prob\left\{|\circledFour| > \frac{\exp\left(-\frac{\gamma\mu T}{2}\right)V}{5} \text{ and } \sum\limits_{l=0}^{T-1}\sum\limits_{i=1}^n\widetilde\sigma_{i,l}^2 \leq \frac{\exp\left(-\gamma\mu T\right)V^2}{150A}\right\}\\
        & \leq 2\exp\left(- \frac{b^2}{2H + 2cb/3}\right)\\
        & = \frac{\beta}{24n(K+1)}.
    \end{align*}
    The above is equivalent to $\Prob\{E_{\circledFour}\} \geq 1 - \frac{\beta}{24n(K+1)}$ for
    \begin{equation}
        E_{\circledFour} = \left\{\text{either } \sum\limits_{l=0}^{T-1}\sum\limits_{i=1}^n\widetilde\sigma_{i,l}^2 > \frac{\exp\left(-\gamma\mu T\right)V^2}{150A} \text{ or } |\circledFour| \leq \frac{\exp\left(-\frac{\gamma\mu T}{2}\right)V}{5}\right\}. \label{ch3:eq:bound_1_str_mon_SEG_1_chapter_3}
    \end{equation}
    Moreover, $E_{T-1}$ implies
    \begin{align}
        \sum\limits_{l=0}^{T-1}\sum\limits_{i=1}^n\widetilde\sigma_{i,l}^2 & \overset{\eqref{ch3:eq:str_mon_thm_SEG_technical_7_1_chapter_3}}{\leq} \frac{2\exp\left(-\gamma\mu (T-\frac{1}{2})\right)\gamma^2V}{3n^2A} \notag\\
        &\qquad \cdot \sum\limits_{l=0}^{T-1}\exp\left(\frac{\gamma\mu l}{2}\right)\sum\limits_{i=1}^n \Exp_{\xi_{1,i}^l, \xi_{2,i}^l}\left[9\|\theta_{i,l}^u\|^2 + 7\|\omega_{i,l}^u\|^2 \right]\notag\\
        & = \frac{\exp\left(-\frac{\gamma\mu T}{2}\right)V}{3A}\cdot\circledThree \overset{\eqref{ch3:eq:bound_3_variances_str_mon_SEG_chapter_3}}{\leq} \frac{\exp\left(-\gamma\mu T\right)V^2}{150A}. \label{ch3:eq:bound_1_variances_str_mon_SEG_1_chapter_3}
    \end{align}

    \paragraph{Upper bound for $\circledFive$.} Probability event $E_{T-1}$ implies
    \begin{align}
        \circledFive & \overset{\eqref{ch3:eq:bias_theta_omega_str_mon_chapter_3}}{\leq} 32\gamma^2\sum\limits_{l=0}^{T-1}\exp\left(-\frac{\gamma\mu}{2}(T-1-l)\right) \frac{2^{2{p}}\sigma^{2{p}}}{\lambda_{l}^{2{p}-2}} \notag \\
        & \overset{\eqref{ch3:eq:lambda_SEG_str_mon_chapter_3}, \eqref{ch3:eq:gamma_SEG_str_mon_chapter_3}}{\leq}  \frac{\exp\left(-\frac{\gamma\mu T}{2}\right)V}{6}. \label{ch3:eq:bound_5_variances_str_mon_SEG_chapter_3}
    \end{align}

    \paragraph{Upper bounds for $\circledSix$ and $\circledSeven$.} These sums require more refined analysis. We introduce new vectors:
    \begin{align}
        \zeta_{j}^l = \begin{cases}
            \frac{\gamma}{n}\sum\limits_{i=1}^{j-1} \theta_{i,l}^u,& \text{if } \left\|\frac{\gamma}{n}\sum\limits_{i=1}^{j-1} \theta_{i,l}^u\right\| \leq \exp\left(-\frac{\gamma\mu l}{4}\right)\frac{\sqrt{V}}{2},\\
            0,& \text{otherwise,}
        \end{cases}, \label{ch3:eq:SEG_zeta_delta_def_str_mon_1_chapter_3} \\
        \delta_{j}^l = \begin{cases}
            \frac{\gamma}{n}\sum\limits_{i=1}^{j-1} \omega_{i,l}^u,& \text{if } \left\|\frac{\gamma}{n}\sum\limits_{i=1}^{j-1} \omega_{i,l}^u\right\| \leq \exp\left(-\frac{\gamma\mu l}{4}\right)\frac{\sqrt{V}}{2},\\
            0,& \text{otherwise,}
        \end{cases}\label{ch3:eq:SEG_zeta_delta_def_str_mon_2_chapter_3}
    \end{align}
    for all $j \in [n]$ and $l = 0,\ldots, T-1$. Then, by definition
    \begin{equation}
        \|\zeta_j^l\| \leq \exp\left(-\frac{\gamma\mu l}{4}\right)\frac{\sqrt{V}}{2},\quad \|\delta_j^l\| \leq \exp\left(-\frac{\gamma\mu l}{4}\right)\frac{\sqrt{V}}{2} \label{ch3:eq:SEG_zeta_delta_bound_str_mon_chapter_3}
    \end{equation}
    and
    \begin{align}
        \circledSix & = \underbrace{\frac{32\gamma}{n} \sum\limits_{l=0}^{T-1} \sum\limits_{j=2}^n \exp\left(-\frac{\gamma\mu}{2}(T-1-l)\right)\left\langle \zeta_j^l, \theta_{j,l}^u \right\rangle}_{\circledSix'}\notag\\
        &\quad\quad + \frac{32\gamma}{n} \sum\limits_{l=0}^{T-1} \sum\limits_{j=2}^n\exp\left(-\frac{\gamma\mu}{2}(T-1-l)\right) \left\langle \frac{\gamma}{n}\sum\limits_{i=1}^{j-1} \theta_{i,l}^u - \zeta_j^l, \theta_{j,l}^u \right\rangle, \label{ch3:eq:SEG_extra_sums_six_distributed_str_mon_chapter_3} \\
        \circledSeven & = \underbrace{\frac{24\gamma}{n} \sum\limits_{l=0}^{T-1} \sum\limits_{j=2}^n \exp\left(-\frac{\gamma\mu}{2}(T-1-l)\right)\left\langle \delta_j^l, \omega_{j,l}^u \right\rangle}_{\circledSeven'}\notag\\
        &\quad\quad + \frac{24\gamma}{n} \sum\limits_{l=0}^{T-1} \sum\limits_{j=2}^n \exp\left(-\frac{\gamma\mu}{2}(T-1-l)\right) \left\langle \frac{\gamma}{n}\sum\limits_{i=1}^{j-1} \omega_{i,l}^u - \delta_j^l, \omega_{j,l}^u \right\rangle. \label{ch3:eq:SEG_extra_sums_seven_distributed_str_mon_chapter_3}
    \end{align}
    We also note here that $E_{T-1}$ implies
    \begin{align}
        \MoveEqLeft[1] \frac{32\gamma}{n} \sum\limits_{l=0}^{T-1} \sum\limits_{j=2}^n\exp\left(-\frac{\gamma\mu}{2}(T-1-l)\right) \notag \\
        & \left\langle \frac{\gamma}{n}\sum\limits_{i=1}^{j-1} \theta_{i,l}^u - \zeta_j^l, \theta_{j,l}^u \right\rangle \notag \\
        &= \frac{32\gamma}{n}  \sum\limits_{j=2}^n \left\langle \frac{\gamma}{n}\sum\limits_{i=1}^{j-1} \theta_{i,T-1}^u - \zeta_j^{T-1}, \theta_{j,T-1}^u \right\rangle, \label{ch3:eq:SEG_extra_sums_six_distributed_1_str_mon_chapter_3}\\
        \MoveEqLeft[1] \frac{24\gamma}{n} \sum\limits_{l=0}^{T-1} \sum\limits_{j=2}^n \exp\left(-\frac{\gamma\mu}{2}(T-1-l)\right) \notag \\
        &\left\langle \frac{\gamma}{n}\sum\limits_{i=1}^{j-1} \omega_{i,l}^u - \delta_j^l, \omega_{j,l}^u \right\rangle \notag\\
        &= \frac{24\gamma}{n}  \sum\limits_{j=2}^n \left\langle \frac{\gamma}{n}\sum\limits_{i=1}^{j-1} \omega_{i,T-1}^u - \delta_j^{T-1}, \omega_{j,T-1}^u \right\rangle. \label{ch3:eq:SEG_extra_sums_seven_distributed_1_str_mon_chapter_3}
    \end{align}

    \paragraph{Upper bound for $\circledSix'$.} To estimate this sum, we will use Bernstein's inequality. The summands have conditional expectations equal to zero:
    \begin{align*}
\MoveEqLeft[1] \Exp_{\xi_{2,j}^l}\left[\frac{32\gamma}{n}\exp\left(-\frac{\gamma\mu}{2}(T-1-l)\right) \langle \zeta_j^l, \theta_{j,l}^u \rangle\right] \\
& = \frac{32\gamma}{n} \exp\left(-\frac{\gamma\mu}{2}(T-1-l)\right)\left\langle \zeta_j^l, \Exp_{\xi_{2,j}^l}[\theta_{j,l}^u] \right\rangle = 0.
    \end{align*}
    Moreover, for all $l = 0,\ldots, T-1$ random vectors $\{\theta_{j,l}^u\}_{j=1}^n$ are independent. Thus, the summands $\frac{32\gamma}{n}\exp\left(-\frac{\gamma\mu}{2}(T-1-l)\right) \langle \zeta_j^l, \theta_{j,l}^u \rangle$ form a martingale difference sequence. Next, the summands are bounded:
    \begin{align}
        \MoveEqLeft[1] \left|\frac{32\gamma}{n} \exp\left(-\frac{\gamma\mu}{2}(T-1-l)\right)\langle \zeta_j^l, \theta_{j,l}^u \rangle \right| \notag \\
        & \leq \frac{32\gamma}{n}\exp\left(-\frac{\gamma\mu}{2}(T-1-l)\right) \|\zeta_j^l\|\cdot \|\theta_{j,l}^u\|\notag\\
        & \overset{\eqref{ch3:eq:SEG_zeta_delta_bound_str_mon_chapter_3},\eqref{ch3:eq:theta_omega_magnitude_str_mon_chapter_3}}{\leq} \frac{32\sqrt{V}\gamma \exp\left(-\frac{\gamma\mu(T-1)}{2}\right)}{n} \exp\left(\frac{\gamma\mu l}{4}\right)\lambda_l\notag\\
        & \overset{\eqref{ch3:eq:lambda_SEG_str_mon_chapter_3}}{\leq} \frac{\exp\left(-\frac{\gamma\mu T}{2}\right)V}{9A} \eqdef c. \label{ch3:eq:str_mon_thm_SEG_technical_6_5_six_chapter_3}
    \end{align}
    Finally, conditional variances $\widehat \sigma_{j,l}^2 \eqdef \Exp_{\xi_{2,j}^l}\left[\frac{1024\gamma^2}{n^2}\exp\left(-\gamma\mu(T-1-l)\right) \langle \zeta_j^l, \theta_{j,l}^u \rangle^2\right]$ of the summands are bounded:
    \begin{align}
        \widehat \sigma_{j,l}^2 & \leq \Exp_{\xi_{2,j}^l}\left[\frac{1024\gamma^2 }{n^2}\exp\left(-\gamma\mu(T-1-l)\right)\|\zeta_j^l\|^2\cdot \|\theta_{j,l}^u\|^2\right] \notag\\
        & \overset{\eqref{ch3:eq:SEG_zeta_delta_bound_str_mon_chapter_3}}{\leq} \frac{256\gamma^2 V\exp\left(-\gamma\mu\left(T-1-\frac{l}{2}\right)\right)}{n^2} \Exp_{\xi_{2,j}^l}\left[\|\theta_{j,l}^u\|^2\right]. \label{ch3:eq:str_mon_thm_SEG_technical_7_six_chapter_3}
    \end{align}
    Applying Bernstein's inequality (Lemma~\ref{ch3:lem:Bernstein_ineq_chapter_3}) with $X_{j,l}$ equal to the summands above, constant $c$ defined in \eqref{ch3:eq:str_mon_thm_SEG_technical_6_5_six_chapter_3}, $b = \frac{\exp\left(-\frac{\gamma\mu T}{2}\right)V}{9}$, $H = \frac{\exp\left(-\gamma\mu T\right)V^2}{600A}$, we get
    \begin{align*}
        \MoveEqLeft[1] \Prob\left\{|\circledSix^{\prime}| > \frac{\exp\left(-\frac{\gamma\mu T}{2}\right)V}{9} \text{ and } \sum\limits_{l=0}^{T-1}\sum\limits_{j=2}^n\widehat \sigma_{j,l}^2 \leq \frac{\exp\left(-\gamma\mu T\right)V^2}{600A}\right\}\\
        & \leq 2\exp\left(- \frac{b^2}{2H + 2cb/3}\right)\\
        & = \frac{\beta}{24n(K+1)}.
    \end{align*}
    The above is equivalent to $\Prob\{E_{\circledSix'}\} \geq 1 - \frac{\beta}{24n(K+1)}$ for
    \begin{align}
        E_{\circledSix'} = \Biggl\{&\text{either } \sum\limits_{l=0}^{T-1}\sum\limits_{j=2}^n\widehat \sigma_{j,l}^2 > \frac{\exp\left(-\gamma\mu T\right)V^2}{600A} \notag\\
        &\text{or } |\circledSix^{\prime}| \leq \frac{\exp\left(-\frac{\gamma\mu T}{2}\right)V}{9}\Biggr\}. \label{ch3:eq:bound_6_str_mon_SEG_six_chapter_3}
    \end{align}
    Moreover, $E_{T-1}$ implies
    \begin{align}
        \sum\limits_{l=0}^{T-1}\sum\limits_{j=2}^n \widehat \sigma_{j,l}^2 & \overset{\eqref{ch3:eq:str_mon_thm_SEG_technical_7_six_chapter_3}}{\leq} \frac{256\gamma^2 V\exp\left(-\gamma\mu\left(T-1\right)\right)}{n^2} \sum\limits_{l=0}^{T-1}\exp\left(\frac{\gamma\mu l}{2}\right) \sum\limits_{i=1}^n \Exp_{\xi_{2,i}^l}\left[\|\theta_{i,l}^u\|^2\right]\notag\\
        & \overset{\eqref{ch3:eq:variance_theta_omega_str_mon_chapter_3}}{\leq}  \frac{4608\gamma^2 V\exp\left(-\gamma\mu\left(T-1\right)\right)\sigma^{p}}{n} \sum\limits_{l=0}^{T-1}\exp\left(\frac{\gamma\mu l}{2}\right) \lambda_l^{2-{p}}\notag\\
        & \overset{\eqref{ch3:eq:lambda_SEG_str_mon_chapter_3}}{\leq} \frac{4608\gamma^{p} V^{2 - \frac{{p}}{2}}\exp\left(-\gamma\mu (T-1)\right)\sigma^{p}}{300^{2-{p}} n^{{p}-1} A^{2-{p}}} \sum\limits_{l=0}^{T-1}\exp\left(\frac{\gamma\mu l {p}}{4}\right)\notag\\
        & \leq \frac{4608\gamma^{p} V^{2 - \frac{{p}}{2}}\exp\left(-\gamma\mu (T-1)\right)\sigma^{p} (K+1) \exp\left(\frac{\gamma\mu K {p}}{4}\right)}{300^{2-{p}} n^{{p}-1} A^{2-{p}}}\notag\\
        & \overset{\eqref{ch3:eq:gamma_SEG_str_mon_chapter_3}}{\leq} \frac{\exp\left(-\gamma\mu T\right)V^2}{600A}.\label{ch3:eq:bound_6_variances_str_mon_SEG_six_chapter_3}
    \end{align}

    \paragraph{Upper bound for $\circledSeven'$.} To estimate this sum, we will use Bernstein's inequality. The summands have conditional expectations equal to zero:
    \begin{multline*}
        \Exp_{\xi_{1,j}^l}\left[\frac{24\gamma}{n}\exp\left(-\frac{\gamma\mu}{2}(T-1-l)\right) \langle \delta_j^l, \omega_{j,l}^u \rangle\right]\\
        = \frac{24\gamma}{n} \exp\left(-\frac{\gamma\mu}{2}(T-1-l)\right)\left\langle \delta_j^l, \Exp_{\xi_{1,j}^l}[\omega_{j,l}^u] \right\rangle = 0.
    \end{multline*}
    Moreover, for all $l = 0,\ldots, T-1$ random vectors $\{\omega_{j,l}^u\}_{j=1}^n$ are independent. Thus, the summands $\frac{24\gamma}{n}\exp\left(-\frac{\gamma\mu}{2}(T-1-l)\right) \langle \delta_j^l, \omega_{j,l}^u \rangle$ form a martingale difference sequence. Next, the summands are bounded:
    \begin{align}
        \MoveEqLeft[1] \left|\frac{24\gamma}{n} \exp\left(-\frac{\gamma\mu}{2}(T-1-l)\right)\langle \delta_j^l, \omega_{j,l}^u \rangle \right| \notag \\
        & \leq \frac{24\gamma}{n}\exp\left(-\frac{\gamma\mu}{2}(T-1-l)\right) \|\delta_j^l\|\cdot \|\omega_{j,l}^u\|\notag\\
        & \overset{\eqref{ch3:eq:SEG_zeta_delta_bound_str_mon_chapter_3},\eqref{ch3:eq:theta_omega_magnitude_str_mon_chapter_3}}{\leq} \frac{24\sqrt{V}\gamma \exp\left(-\frac{\gamma\mu(T-1)}{2}\right)}{n} \exp\left(\frac{\gamma\mu l}{4}\right)\lambda_l\notag\\
        & \overset{\eqref{ch3:eq:lambda_SEG_str_mon_chapter_3}}{\leq} \frac{\exp\left(-\frac{\gamma\mu T}{2}\right)V}{10A} \eqdef c. \label{ch3:eq:str_mon_thm_SEG_technical_6_5_seven_chapter_3}
    \end{align}
    Finally, conditional variances $(\sigma_{j,l}')^2 \eqdef \Exp_{\xi_{1,j}^l}\left[\frac{576\gamma^2}{n^2}\exp\left(-\gamma\mu(T-1-l)\right) \langle \delta_j^l, \omega_{j,l}^u \rangle^2\right]$ of the summands are bounded:
    \begin{align}
        (\sigma_{j,l}')^2 & \leq \Exp_{\xi_{1,j}^l}\left[\frac{576\gamma^2 }{n^2}\exp\left(-\gamma\mu(T-1-l)\right)\|\delta_j^l\|^2\cdot \|\omega_{j,l}^u\|^2\right] \notag\\
        & \overset{\eqref{ch3:eq:SEG_zeta_delta_bound_str_mon_chapter_3}}{\leq} \frac{288\gamma^2 V\exp\left(-\gamma\mu\left(T-1-\frac{l}{2}\right)\right)}{n^2} \Exp_{\xi_{1,j}^l}\left[\|\omega_{j,l}^u\|^2\right]. \label{ch3:eq:str_mon_thm_SEG_technical_7_seven_chapter_3}
    \end{align}
    Applying Bernstein's inequality (Lemma~\ref{ch3:lem:Bernstein_ineq_chapter_3}) with $X_{j,l}$ equal to the summands above, constant $c$ defined in \eqref{ch3:eq:str_mon_thm_SEG_technical_6_5_seven_chapter_3}, $b = \frac{\exp\left(-\frac{\gamma\mu T}{2}\right)V}{10}$, $H = \frac{\exp\left(-\gamma\mu T\right)V^2}{600A}$, we get
    \begin{align*}
        \MoveEqLeft[1] \Prob\left\{|\circledSeven^{\prime}| > \frac{\exp\left(-\frac{\gamma\mu T}{2}\right)V}{10} \text{ and } \sum\limits_{l=0}^{T-1}\sum\limits_{j=2}^n (\sigma_{j,l}^{\prime})^2 \leq \frac{\exp\left(-\gamma\mu T\right)V^2}{600A}\right\}\\
        & \leq 2\exp\left(- \frac{b^2}{2H + 2cb/3}\right)\\
        & = \frac{\beta}{24n(K+1)}.
    \end{align*}
    The above is equivalent to $\Prob\{E_{\circledSeven'}\} \geq 1 - \frac{\beta}{24n(K+1)}$ for
    \begin{align}
        E_{\circledSeven'} = \Biggl\{&\text{either } \sum\limits_{l=0}^{T-1}\sum\limits_{j=2}^n (\sigma_{j,l}^{\prime})^2 > \frac{\exp\left(-\gamma\mu T\right)V^2}{600A} \notag\\
        &\text{or } |\circledSeven^{\prime}| \leq \frac{\exp\left(-\frac{\gamma\mu T}{2}\right)V}{10}\Biggr\}. \label{ch3:eq:bound_6_str_mon_SEG_seven_chapter_3}
    \end{align}
    Moreover, $E_{T-1}$ implies
    \begin{align}
        \MoveEqLeft[1] \sum\limits_{l=0}^{T-1}\sum\limits_{j=2}^n (\sigma_{j,l}')^2 \notag \\
        & \overset{\eqref{ch3:eq:str_mon_thm_SEG_technical_7_seven_chapter_3}}{\leq} \frac{288\gamma^2 V\exp\left(-\gamma\mu\left(T-1\right)\right)}{n^2} \sum\limits_{l=0}^{T-1}\exp\left(\frac{\gamma\mu l}{2}\right) \sum\limits_{i=1}^n \Exp_{\xi_{1,i}^l}\left[\|\omega_{i,l}^u\|^2\right]\notag\\
        & \overset{\eqref{ch3:eq:variance_theta_omega_str_mon_chapter_3}}{\leq}  \frac{5184\gamma^2 V\exp\left(-\gamma\mu\left(T-1\right)\right)\sigma^{p}}{n} \sum\limits_{l=0}^{T-1}\exp\left(\frac{\gamma\mu l}{2}\right) \lambda_l^{2-{p}}\notag\\
        & \overset{\eqref{ch3:eq:lambda_SEG_str_mon_chapter_3}}{\leq} \frac{5184\gamma^{p} V^{2 - \frac{{p}}{2}}\exp\left(-\gamma\mu (T-1)\right)\sigma^{p}}{300^{2-{p}} n^{{p}-1} A^{2-{p}}} \sum\limits_{l=0}^{T-1}\exp\left(\frac{\gamma\mu l {p}}{4}\right)\notag\\
        & \leq \frac{5184\gamma^{p} V^{2 - \frac{{p}}{2}}\exp\left(-\gamma\mu (T-1)\right)\sigma^{p} (K+1) \exp\left(\frac{\gamma\mu K {p}}{4}\right)}{300^{2-{p}} n^{{p}-1} A^{2-{p}}}\notag\\
        & \overset{\eqref{ch3:eq:gamma_SEG_str_mon_chapter_3}}{\leq} \frac{\exp\left(-\gamma\mu T\right)V^2}{600A}.\label{ch3:eq:bound_6_variances_str_mon_SEG_seven_chapter_3}
    \end{align}

    That is, we derive the upper bounds for  $\circledOne, \circledTwo, \circledThree, \circledFour, \circledFive, \circledSix, \circledSeven$. More precisely, $E_{T-1}$ implies
    \begin{gather*}
        V_T \overset{\eqref{ch3:eq:ndvjinuibvbsvdhbvdhbfussdncdni_chapter_3}}{\leq} \exp\left(-\frac{\gamma\mu}{2}T\right)V + \circledOne + \circledTwo + \circledThree + \circledFour + \circledFive + \circledSix + \circledSeven,\\
         \circledSix \overset{\eqref{ch3:eq:SEG_extra_sums_six_distributed_str_mon_chapter_3}}{=}  \circledSix' + \frac{32\gamma}{n}  \sum\limits_{j=2}^n \left\langle \frac{\gamma}{n}\sum\limits_{i=1}^{j-1} \theta_{i,T-1}^u - \zeta_j^{T-1}, \theta_{j,T-1}^u \right\rangle,\\
         \circledSeven \overset{\eqref{ch3:eq:SEG_extra_sums_seven_distributed_str_mon_chapter_3}}{=} \circledSeven' + \frac{24\gamma}{n}  \sum\limits_{j=2}^n \left\langle \frac{\gamma}{n}\sum\limits_{i=1}^{j-1} \omega_{i,T-1}^u - \delta_j^{T-1}, \omega_{j,T-1}^u \right\rangle,\\
         \circledTwo \overset{\eqref{ch3:eq:bound_2_variances_str_mon_SEG_chapter_3}}{\leq} \frac{3\exp\left(-\frac{\gamma\mu T}{2}\right) V}{100},\quad \circledThree \overset{\eqref{ch3:eq:bound_3_variances_str_mon_SEG_chapter_3}}{\leq} \frac{\exp\left(-\frac{\gamma\mu T}{2}\right) V}{50},\quad \circledFive \overset{\eqref{ch3:eq:bound_5_variances_str_mon_SEG_chapter_3}}{\leq} \frac{\exp\left(-\frac{\gamma\mu T}{2}\right) V}{6},\\
         \sum\limits_{l=0}^{T-1}\sum\limits_{i=1}^n\sigma_{i,l}^2 \overset{\eqref{ch3:eq:bound_1_variances_str_mon_SEG_chapter_3}}{\leq} \frac{\exp\left(-\gamma\mu T\right)V^2}{15000A},\quad \sum\limits_{l=0}^{T-1}\sum\limits_{i=1}^n\widetilde\sigma_{i,l}^2 \overset{\eqref{ch3:eq:bound_1_variances_str_mon_SEG_1_chapter_3}}{\leq} \frac{\exp\left(-\gamma\mu T\right)V^2}{150A},\\
         \sum\limits_{l=0}^{T-1}\sum\limits_{j=2}^n \widehat \sigma_{j,l}^2 \overset{\eqref{ch3:eq:bound_6_variances_str_mon_SEG_six_chapter_3}}{\leq} \frac{\exp\left(-\gamma\mu T\right)V^2}{600A},\quad \sum\limits_{l=0}^{T-1}\sum\limits_{j=2}^n (\sigma_{j,l}')^2 \overset{\eqref{ch3:eq:bound_6_variances_str_mon_SEG_seven_chapter_3}}{\leq} \frac{\exp\left(-\gamma\mu T\right)V^2}{600A}.
    \end{gather*}
    In addition, we also establish (see \eqref{ch3:eq:bound_1_str_mon_SEG_chapter_3}, \eqref{ch3:eq:bound_1_str_mon_SEG_chapter_3}, \eqref{ch3:eq:bound_6_str_mon_SEG_six_chapter_3}, \eqref{ch3:eq:bound_6_str_mon_SEG_six_chapter_3}, and our induction assumption)
    \begin{gather*}
        \Prob\{E_{T-1}\} \geq 1 - \frac{(T-1)\beta}{K+1},\quad \Prob\{E_{\circledOne}\} \geq 1 - \frac{\beta}{24n(K+1)}, \\
        \Prob\{E_{\circledFour}\} \geq 1 - \frac{\beta}{24n(K+1)}, \quad \Prob\{E_{\circledSix'}\} \geq 1 - \frac{\beta}{24n(K+1)},\\
        \Prob\{E_{\circledSeven'}\} \geq 1 - \frac{\beta}{24n(K+1)},
    \end{gather*}
    where
    \begin{align*}
        E_{\circledOne} & = \left\{\text{either } \sum\limits_{l=0}^{T-1}\sum\limits_{i=1}^n\sigma_{i,l}^2 > \frac{\exp\left(-\gamma\mu T\right)V^2}{15000A} \text{ or } |\circledOne| \leq \frac{\exp\left(-\frac{\gamma\mu T}{2}\right)V}{50}\right\},\\
        E_{\circledFour} & = \left\{\text{either } \sum\limits_{l=0}^{T-1}\sum\limits_{i=1}^n\widetilde\sigma_{i,l}^2 > \frac{\exp\left(-\gamma\mu T\right)V^2}{150A} \text{ or } |\circledFour| \leq \frac{\exp\left(-\frac{\gamma\mu T}{2}\right)V}{5}\right\},\\
        E_{\circledSix'} & = \left\{\text{either } \sum\limits_{l=0}^{T-1}\sum\limits_{j=2}^n\widehat \sigma_{j,l}^2 > \frac{\exp\left(-\gamma\mu T\right)V^2}{600A} \text{ or } |\circledSix^{\prime}| \leq \frac{\exp\left(-\frac{\gamma\mu T}{2}\right)V}{9}\right\},\\
        E_{\circledSeven'} & = \left\{\text{either } \sum\limits_{l=0}^{T-1}\sum\limits_{j=2}^n (\sigma_{j,l}^{\prime})^2 > \frac{\exp\left(-\gamma\mu T\right)V^2}{600A} \text{ or } |\circledSeven^{\prime}| \leq \frac{\exp\left(-\frac{\gamma\mu T}{2}\right)V}{10}\right\}.
    \end{align*}
    Therefore, denoting
    \begin{gather*}
        S_\theta \eqdef \frac{\gamma}{n} \sum\limits_{j=2}^n \left\langle \frac{\gamma}{n}\sum\limits_{i=1}^{j-1} \theta_{i,T-1}^u - \zeta_j^{T-1}, \theta_{j,T-1}^u \right\rangle,\\
        S_\omega \eqdef \frac{\gamma}{n} \sum\limits_{j=2}^n \left\langle \frac{\gamma}{n}\sum\limits_{i=1}^{j-1} \omega_{i,T-1}^u - \delta_j^{T-1}, \omega_{j,T-1}^u \right\rangle,
    \end{gather*}
    we get that probability event $E_{T-1} \cap E_{\circledOne} \cap E_{\circledFour} \cap E_{\circledSix'} \cap E_{\circledSeven'}$ implies
    \begin{align}
        V_T & \leq \exp\left(-\frac{\gamma\mu}{2}T\right)V \underbrace{\left(1 + \frac{1}{50} + \frac{3}{100} + \frac{1}{50} + \frac{1}{5} + \frac{1}{6} + \frac{1}{9} + \frac{1}{10}\right)}_{\leq 2} \notag\\
        &\quad\quad + 32S_\theta + 24S_\omega. \label{ch3:eq:nvsdibhjbsjdbvshbcjs_chapter_3}
    \end{align}

    To finish the proof, we need to show that $\frac{\gamma}{n}\sum\limits_{i=1}^{j-1} \theta_{i,T-1}^u = \zeta_j^{T-1}$ and $\frac{\gamma}{n}\sum\limits_{i=1}^{j-1} \omega_{i,T-1}^u = \delta_j^{T-1}$ with high probability. In particular, we consider probability event $\widetilde{E}_{T-1,j}$ defined as follows: inequalities
    \begin{align}
\MoveEqLeft[1] \left\|\frac{\gamma}{n}\sum\limits_{i=1}^{r-1} \theta_{i,T-1}^u\right\| \notag \\
& \leq \exp\left(-\frac{\gamma\mu (T-1)}{4}\right)\frac{\sqrt{V}}{2}, \quad \left\|\frac{\gamma}{n}\sum\limits_{i=1}^{r-1} \omega_{i,T-1}^u\right\| \notag \\
& \leq \exp\left(-\frac{\gamma\mu (T-1)}{4}\right)\frac{\sqrt{V}}{2} \notag
    \end{align}
    hold for $r = 2, \ldots, j$ simultaneously. We want to show that $\Prob\{E_{T-1} \cap \widetilde{E}_{T-1,j}\} \geq 1 - \frac{(T-1)\beta}{K+1} - \frac{j\beta}{8n(K+1)}$ for all $j = 2, \ldots, n$. For $j = 2$ the statement is trivial since
    \begin{align*}
        \left\|\frac{\gamma}{n} \theta_{1,T-1}^u\right\| & \overset{\eqref{ch3:eq:theta_omega_magnitude_str_mon_chapter_3}}{\leq} \frac{2\gamma\lambda_{T-1}}{n} \leq \exp\left(-\frac{\gamma\mu (T-1)}{4}\right)\frac{\sqrt{V}}{2},\\
        \left\|\frac{\gamma}{n} \omega_{1,T-1}^u\right\| & \overset{\eqref{ch3:eq:theta_omega_magnitude_chapter_3}}{\leq} \frac{2\gamma\lambda_{T-1}}{n} \leq \exp\left(-\frac{\gamma\mu (T-1)}{4}\right)\frac{\sqrt{V}}{2}.
    \end{align*}
    Next, we assume that the statement holds for some $j = m-1 < n$, i.e.,  $\Prob\{E_{T-1}\cap\widetilde{E}_{T-1,m-1}\} \geq 1 - \frac{(T-1)\beta}{K+1} - \frac{(m-1)\beta}{8n(K+1)}$. Our goal is to prove that $\Prob\{E_{T-1}\cap\widetilde{E}_{T-1,m}\} \geq 1 - \frac{(T-1)\beta}{K+1} - \frac{m\beta}{8n(K+1)}$. First, we consider $\left\|\frac{\gamma}{n}\sum\limits_{i=1}^{m-1} \theta_{i,T-1}^u\right\|$:
    \begin{align*}
        \left\|\frac{\gamma}{n}\sum\limits_{i=1}^{m-1} \theta_{i,T-1}^u\right\|
        &= \sqrt{\frac{\gamma^2}{n^2}\left\|\sum\limits_{i=1}^{m-1} \theta_{i,T-1}^u\right\|^2}\\
        &= \sqrt{\begin{aligned}
            &\frac{\gamma^2}{n^2}\sum\limits_{i=1}^{m-1}\|\theta_{i,T-1}^u\|^2\\
            &+ \frac{2\gamma}{n}\sum\limits_{i=1}^{m-1}\left\langle \frac{\gamma}{n}\sum\limits_{r=1}^{i-1}\theta_{r,T-1}^u , \theta_{i,T-1}^u \right\rangle
        \end{aligned}}\\
        &\leq \sqrt{\begin{aligned}
            &\frac{\gamma^2}{n^2}\sum\limits_{l=0}^{T-1}\exp\left(-\frac{\gamma\mu(T-1-l)}{2}\right)\sum\limits_{i=1}^{m-1}\|\theta_{i,l}^u\|^2\\
            &+ \frac{2\gamma}{n}\sum\limits_{i=1}^{m-1}\left\langle \frac{\gamma}{n}\sum\limits_{r=1}^{i-1}\theta_{r,T-1}^u , \theta_{i,T-1}^u \right\rangle
        \end{aligned}}.
    \end{align*}
    Similarly, we have
    \begin{align*}
        \left\|\frac{\gamma}{n}\sum\limits_{i=1}^{m-1} \omega_{i,T-1}^u\right\|
        &= \sqrt{\frac{\gamma^2}{n^2}\left\|\sum\limits_{i=1}^{m-1} \omega_{i,T-1}^u\right\|^2}\\
        &= \sqrt{\begin{aligned}
            &\frac{\gamma^2}{n^2}\sum\limits_{i=1}^{m-1}\|\omega_{i,T-1}^u\|^2\\
            &+ \frac{2\gamma}{n}\sum\limits_{i=1}^{m-1}\left\langle \frac{\gamma}{n}\sum\limits_{r=1}^{i-1}\omega_{r,T-1}^u , \omega_{i,T-1}^u \right\rangle
        \end{aligned}}\\
        &\leq \sqrt{\begin{aligned}
            &\frac{\gamma^2}{n^2}\sum\limits_{l=0}^{T-1}\exp\left(-\frac{\gamma\mu(T-1-l)}{2}\right)\sum\limits_{i=1}^{m-1}\|\omega_{i,l}^u\|^2\\
            &+ \frac{2\gamma}{n}\sum\limits_{i=1}^{m-1}\left\langle \frac{\gamma}{n}\sum\limits_{r=1}^{i-1}\omega_{r,T-1}^u , \omega_{i,T-1}^u \right\rangle
        \end{aligned}}.
    \end{align*}

    Next, we introduce a new notation:
    \begin{gather*}
        \rho_{i,T-1} = \begin{cases}
            \frac{\gamma}{n}\sum\limits_{r=1}^{i-1}\theta_{r,T-1}^u,& \text{if } \left\|\frac{\gamma}{n}\sum\limits_{r=1}^{i-1}\theta_{r,T-1}^u\right\| \leq \exp\left(-\frac{\gamma\mu (T-1)}{4}\right)\frac{\sqrt{V}}{2},\\
            0,& \text{otherwise}
        \end{cases},\\
        \rho_{i,T-1}' = \begin{cases}
            \frac{\gamma}{n}\sum\limits_{r=1}^{i-1}\omega_{r,T-1}^u,& \text{if } \left\|\frac{\gamma}{n}\sum\limits_{r=1}^{i-1}\omega_{r,T-1}^u\right\| \leq \exp\left(-\frac{\gamma\mu (T-1)}{4}\right)\frac{\sqrt{V}}{2},\\
            0,& \text{otherwise}
        \end{cases}
    \end{gather*}
    for $i = 1,\ldots,m-1$. By definition, we have
    \begin{align}
\|\rho_{i,T-1}\| & \leq \exp\left(-\frac{\gamma\mu (T-1)}{4}\right)\frac{\sqrt{V}}{2},\quad \|\rho_{i,T-1}'\| \notag \\
& \leq \exp\left(-\frac{\gamma\mu (T-1)}{4}\right)\frac{\sqrt{V}}{2} \label{ch3:eq:SEG_bound_rho_str_mon_chapter_3}
    \end{align}
    for $i = 1,\ldots,m-1$. Moreover, $\widetilde{E}_{T-1,m-1}$ implies $\rho_{i,T-1} = \frac{\gamma}{n}\sum\limits_{r=1}^{i-1}\theta_{r,T-1}^u$, $\rho_{i,T-1}' = \frac{\gamma}{n}\sum\limits_{r=1}^{i-1}\omega_{r,T-1}^u$ for $i = 1,\ldots,m-1$ and
    \begin{align}
        \left\|\frac{\gamma}{n}\sum\limits_{i=1}^{m-1} \theta_{i,T-1}^u\right\| & \leq \sqrt{\circledThree + \circledFour + \circledEight}, \notag\\
        \left\|\frac{\gamma}{n}\sum\limits_{i=1}^{m-1} \omega_{i,T-1}^u\right\| & \leq \sqrt{\circledThree + \circledFour + \circledEight'}, \notag
    \end{align}
    where
    \begin{gather*}
        \circledEight =  \frac{2\gamma}{n}\sum\limits_{i=1}^{m-1}\left\langle \rho_{i,T-1} , \theta_{i,T-1}^u \right\rangle,\quad \circledEight' =  \frac{2\gamma}{n}\sum\limits_{i=1}^{m-1}\left\langle \rho_{i,T-1}' , \omega_{i,T-1}^u \right\rangle.
    \end{gather*}
    It remains to estimate $\circledEight$ and $\circledEight'$.

    \paragraph{Upper bound for $\circledEight$.} To estimate this sum, we will use Bernstein's inequality. The summands have conditional expectations equal to zero:
    \begin{align*}
\MoveEqLeft[1] \Exp_{\xi_{2,i}^{T-1}}\left[\frac{2\gamma}{n} \langle \rho_{i,T-1}, \theta_{i,T-1}^u \rangle\right] \\
& = \frac{2\gamma}{n} \left\langle \rho_{i,T-1}, \Exp_{\xi_{2,i}^{T-1}}[\theta_{i,T-1}^u] \right\rangle = 0.
    \end{align*}
    Thus, for $i = 1, \ldots, m-1$ the summands $(2\gamma/n)\langle \rho_{i,T-1}, \theta_{i,T-1}^u \rangle$ form a martingale difference sequence. Next, the summands are bounded:
    \begin{align}
        \left|\frac{2\gamma}{n} \langle \rho_{i,T-1}, \theta_{i,T-1}^u \rangle \right| & \leq \frac{2\gamma}{n}\|\rho_{i,T-1}\|\cdot \|\theta_{i,T-1}^u\|\notag\\
        & \overset{\eqref{ch3:eq:SEG_bound_rho_str_mon_chapter_3},\eqref{ch3:eq:theta_omega_magnitude_str_mon_chapter_3}}{\leq} \frac{2\sqrt{V}\gamma \exp\left(-\frac{\gamma\mu(T-1)}{4}\right)}{n} \lambda_{T-1}\notag\\
        & \overset{\eqref{ch3:eq:lambda_SEG_str_mon_chapter_3}}{\leq} \frac{\exp\left(-\frac{\gamma\mu (T-1)}{2}\right)V}{80A} \eqdef c. \label{ch3:eq:str_mon_thm_SEG_technical_6_5_eight_chapter_3}
    \end{align}
    Finally, conditional variances $(\widehat\sigma_{i,T-1}')^2 \eqdef \Exp_{\xi_{2,i}^{T-1}}\left[\frac{4\gamma^2}{n^2} \langle \rho_{i,T-1}, \theta_{i,T-1}^u \rangle^2\right]$ of the summands are bounded:
    \begin{align}
        (\widehat\sigma_{i,T-1}')^2 & \leq \Exp_{\xi_{2,i}^{T-1}}\left[\frac{4\gamma^2 }{n^2}\|\rho_{i,T-1}\|^2\cdot \|\theta_{i,T-1}^u\|^2\right] \notag\\
        & \overset{\eqref{ch3:eq:SEG_bound_rho_str_mon_chapter_3}}{\leq} \frac{\gamma^2 V\exp\left(-\frac{\gamma\mu(T-1)}{2}\right)}{n^2} \Exp_{\xi_{2,i}^{T-1}}\left[\|\theta_{i,T-1}^u\|^2\right]. \label{ch3:eq:str_mon_thm_SEG_technical_7_eight_chapter_3}
    \end{align}
    Applying Bernstein's inequality (Lemma~\ref{ch3:lem:Bernstein_ineq_chapter_3}) with $X_{i} = \frac{2\gamma}{n}\langle \rho_{i,T-1}, \theta_{i,T-1}^u \rangle$, constant $c$ defined in \eqref{ch3:eq:str_mon_thm_SEG_technical_6_5_eight_chapter_3}, $b = \frac{\exp\left(-\frac{\gamma\mu (T-1)}{2}\right)V}{80}$, $H = \frac{\exp\left(-\gamma\mu (T-1)\right)V^2}{38400A}$, we get
    \begin{align*}
        &\Prob\left\{|\circledEight| > \frac{\exp\left(-\frac{\gamma\mu (T-1)}{2}\right)V}{80} \text{ and } \sum\limits_{i=1}^{m-1}(\widehat\sigma_{i,T-1}')^2 \leq \frac{\exp\left(-\gamma\mu (T-1)\right)V^2} {38400A}\right\}\\
        &\leq 2\exp\left(- \frac{b^2}{2H + 2cb/3}\right)\\
        &= \frac{\beta}{24n(K+1)}.
    \end{align*}
    The above is equivalent to $\Prob\{E_{\circledEight}\} \geq 1 - \frac{\beta}{24n(K+1)}$ for
    \begin{equation}
        \begin{aligned}
        E_{\circledEight} = \biggl\{&\text{either } \sum\limits_{i=1}^{m-1}(\widehat\sigma_{i,T-1}^{\prime})^2 > \frac{\exp\left(-\gamma\mu (T-1)\right)V^2}{38400A}\\
        &\text{or } |\circledEight| \leq \frac{\exp\left(-\frac{\gamma\mu (T-1)}{2}\right)V}{80}\biggr\}.
        \end{aligned}\label{ch3:eq:bound_1_str_mon_SEG_eight_chapter_3}
    \end{equation}
    Moreover, $E_{T-1}$ implies
    \begin{align}
        \sum\limits_{i=1}^{m-1}(\widehat\sigma_{i,T-1}')^2 & \overset{\eqref{ch3:eq:str_mon_thm_SEG_technical_7_eight_chapter_3}}{\leq} \frac{\gamma^2 V\exp\left(-\frac{\gamma\mu(T-1)}{2}\right)}{n^2} \sum\limits_{i=1}^{m-1} \Exp_{\xi_{2,i}^{T-1}}\left[\|\theta_{i,T-1}^u\|^2\right]\notag\\
        & \overset{\eqref{ch3:eq:variance_theta_omega_str_mon_chapter_3}}{\leq}  \frac{18\gamma^2 V\exp\left(-\frac{\gamma\mu(T-1)}{2}\right)\sigma^{p}}{n} \lambda_{T-1}^{2-{p}}\notag\\
        & \overset{\eqref{ch3:eq:lambda_SEG_str_mon_chapter_3}}{\leq} \frac{18\gamma^{p} V^{2 - \frac{{p}}{2}}\exp\left(-\gamma\mu (T-1)\right)\sigma^{p}}{300^{2-{p}} n^{{p}-1} A^{2-{p}}} \exp\left(\frac{\gamma\mu (T-1) {p}}{4}\right)\notag\\
        & \overset{T-1 \leq K}{\leq} \frac{18\gamma^{p} V^{2 - \frac{{p}}{2}}\exp\left(-\gamma\mu (T-1)\right)\sigma^{p} \exp\left(\frac{\gamma\mu K {p}}{4}\right)}{300^{2-{p}} n^{{p}-1} A^{2-{p}}}\notag\\
        & \overset{\eqref{ch3:eq:gamma_SEG_str_mon_chapter_3}}{\leq} \frac{\exp\left(-\gamma\mu (T-1)\right)V^2}{38400A}.\label{ch3:eq:bound_1_variances_str_mon_SEG_eight_chapter_3}
    \end{align}

    \paragraph{Upper bound for $\circledEight'$.} To estimate this sum, we will use Bernstein's inequality. The summands have conditional expectations equal to zero:
    \begin{align*}
\MoveEqLeft[1] \Exp_{\xi_{1,i}^{T-1}}\left[\frac{2\gamma}{n} \langle \rho_{i,T-1}', \omega_{i,T-1}^u \rangle\right] \\
& = \frac{2\gamma}{n} \left\langle \rho_{i,T-1}', \Exp_{\xi_{1,i}^{T-1}}[\omega_{i,T-1}^u] \right\rangle = 0.
    \end{align*}
    Thus, for $i = 1, \ldots, m-1$ the summands $(2\gamma/n)\langle \rho_{i,T-1}', \omega_{i,T-1}^u \rangle$ form a martingale difference sequence. Next, the summands are bounded:
    \begin{align}
        \MoveEqLeft[1] \left|\frac{2\gamma}{n} \langle \rho_{i,T-1}', \omega_{i,T-1}^u \rangle \right| \notag \\
        & \leq \frac{2\gamma}{n}\|\rho_{i,T-1}'\|\cdot \|\omega_{i,T-1}^u\|\notag\\
        & \overset{\eqref{ch3:eq:SEG_bound_rho_str_mon_chapter_3},\eqref{ch3:eq:theta_omega_magnitude_str_mon_chapter_3}}{\leq} \frac{2\sqrt{V}\gamma \exp\left(-\frac{\gamma\mu(T-1)}{4}\right)}{n} \lambda_{T-1}\notag\\
        & \overset{\eqref{ch3:eq:lambda_SEG_str_mon_chapter_3}}{=} \frac{\exp\left(-\frac{\gamma\mu (T-1)}{2}\right)V}{80A} \eqdef c. \label{ch3:eq:str_mon_thm_SEG_technical_6_5_eight'_chapter_3}
    \end{align}
    Finally, conditional variances $(\widetilde\sigma_{i,T-1}')^2 \eqdef \Exp_{\xi_{1,i}^{T-1}}\left[\frac{4\gamma^2}{n^2} \langle \rho_{i,T-1}', \omega_{i,T-1}^u \rangle^2\right]$ of the summands are bounded:
    \begin{align}
        (\widetilde\sigma_{i,T-1}')^2 & \leq \Exp_{\xi_{1,i}^{T-1}}\left[\frac{4\gamma^2 }{n^2}\|\rho_{i,T-1}'\|^2\cdot \|\omega_{i,T-1}^u\|^2\right] \notag\\
        & \overset{\eqref{ch3:eq:SEG_bound_rho_str_mon_chapter_3}}{\leq} \frac{\gamma^2 V\exp\left(-\frac{\gamma\mu(T-1)}{2}\right)}{n^2} \Exp_{\xi_{1,i}^{T-1}}\left[\|\omega_{i,T-1}^u\|^2\right]. \label{ch3:eq:str_mon_thm_SEG_technical_7_eight'_chapter_3}
    \end{align}
    Applying Bernstein's inequality (Lemma~\ref{ch3:lem:Bernstein_ineq_chapter_3}) with $X_{i} = \frac{2\gamma}{n}\langle \rho_{i,T-1}', \omega_{i,T-1}^u \rangle$, constant $c$ defined in \eqref{ch3:eq:str_mon_thm_SEG_technical_6_5_eight'_chapter_3}, $b = \frac{\exp\left(-\frac{\gamma\mu (T-1)}{2}\right)V}{80}$, $H = \frac{\exp\left(-\gamma\mu (T-1)\right)V^2}{38400A}$, we get
    \begin{align*}
        &\Prob\left\{|\circledEight'| > \frac{\exp\left(-\frac{\gamma\mu (T-1)}{2}\right)V}{80} \text{ and } \sum\limits_{i=1}^{m-1}(\widetilde\sigma_{i,T-1}')^2 \leq \frac{\exp\left(-\gamma\mu (T-1)\right)V^2} {38400A}\right\}\\
        &\leq 2\exp\left(- \frac{b^2}{2H + 2cb/3}\right)\\
        &= \frac{\beta}{24n(K+1)}.
    \end{align*}
    The above is equivalent to $\Prob\{E_{\circledEight'}\} \geq 1 - \frac{\beta}{24n(K+1)}$ for
    \begin{equation}
        \begin{aligned}
        E_{\circledEight'} = \biggl\{&\text{either } \sum\limits_{i=1}^{m-1}(\widetilde\sigma_{i,T-1}^{\prime})^2 > \frac{\exp\left(-\gamma\mu (T-1)\right)V^2}{38400A}\\
        &\text{or } |\circledEight^{\prime}| \leq \frac{\exp\left(-\frac{\gamma\mu (T-1)}{2}\right)V}{80}\biggr\}.
        \end{aligned}\label{ch3:eq:bound_1_str_mon_SEG_eight'_chapter_3}
    \end{equation}
    Moreover, $E_{T-1}$ implies
    \begin{align}
        \sum\limits_{i=1}^{m-1}(\widetilde\sigma_{i,T-1}')^2 & \overset{\eqref{ch3:eq:str_mon_thm_SEG_technical_7_eight'_chapter_3}}{\leq} \frac{\gamma^2 V\exp\left(-\frac{\gamma\mu(T-1)}{2}\right)}{n^2} \sum\limits_{i=1}^{m-1} \Exp_{\xi_{1,i}^{T-1}}\left[\|\omega_{i,T-1}^u\|^2\right]\notag\\
        & \overset{\eqref{ch3:eq:variance_theta_omega_str_mon_chapter_3}}{\leq}  \frac{18\gamma^2 V\exp\left(-\frac{\gamma\mu(T-1)}{2}\right)\sigma^{p}}{n} \lambda_{T-1}^{2-{p}}\notag\\
        & \overset{\eqref{ch3:eq:lambda_SEG_str_mon_chapter_3}}{\leq} \frac{18\gamma^{p} V^{2 - \frac{{p}}{2}}\exp\left(-\gamma\mu (T-1)\right)\sigma^{p}}{300^{2-{p}} n^{{p}-1} A^{2-{p}}} \exp\left(\frac{\gamma\mu (T-1) {p}}{4}\right)\notag\\
        & \overset{T-1 \leq K}{\leq} \frac{18\gamma^{p} V^{2 - \frac{{p}}{2}}\exp\left(-\gamma\mu (T-1)\right)\sigma^{p} \exp\left(\frac{\gamma\mu K {p}}{4}\right)}{300^{2-{p}} n^{{p}-1} A^{2-{p}}}\notag\\
        & \overset{\eqref{ch3:eq:gamma_SEG_str_mon_chapter_3}}{\leq} \frac{\exp\left(-\gamma\mu (T-1)\right)V^2}{38400A}.\label{ch3:eq:bound_1_variances_str_mon_SEG_eight'_chapter_3}
    \end{align}

    Putting all together we get that $E_{T-1}\cap \widetilde{E}_{T-1,m-1}$ implies
    \begin{gather*}
        \left\|\frac{\gamma}{n}\sum\limits_{i=1}^{m-1} \theta_{i,T-1}^u\right\| \leq \sqrt{\circledThree + \circledFour + \circledEight},\\
        \left\|\frac{\gamma}{n}\sum\limits_{i=1}^{m-1} \omega_{i,T-1}^u\right\| \leq \sqrt{\circledThree + \circledFour + \circledEight'},\\
        \circledThree \overset{\eqref{ch3:eq:bound_3_variances_str_mon_SEG_chapter_3}}{\leq} \frac{\exp\left(-\frac{\gamma\mu T}{2}\right)V}{50},\\
        \sum\limits_{l=0}^{T-1}\sum\limits_{i=1}^n\widetilde\sigma_{i,l}^2 \overset{\eqref{ch3:eq:bound_1_variances_str_mon_SEG_1_chapter_3}}{\leq} \frac{\exp\left(-\gamma\mu T\right)V^2}{150A},\\
        \sum\limits_{i=1}^{m-1}(\widehat\sigma_{i,T-1}')^2 \leq \frac{\exp\left(-\gamma\mu (T-1)\right)V^2}{38400A},\\
        \sum\limits_{i=1}^{m-1}(\widetilde\sigma_{i,T-1}')^2 \leq \frac{\exp\left(-\gamma\mu (T-1)\right)V^2}{38400A}.
    \end{gather*}
    In addition, we also establish (see \eqref{ch3:eq:bound_1_str_mon_SEG_1_chapter_3}, \eqref{ch3:eq:bound_1_str_mon_SEG_eight_chapter_3}, \eqref{ch3:eq:bound_1_str_mon_SEG_eight'_chapter_3} and our induction assumption)
    \begin{gather*}
        \Prob\{E_{T-1}\cap \widetilde{E}_{T-1,m-1}\} \geq 1 - \frac{(T-1)\beta}{K+1} - \frac{(m-1)\beta}{8n(K+1)},\\
        \Prob\{E_{\circledFour}\} \geq 1 - \frac{\beta}{24n(K+1)}, \quad \Prob\{E_{\circledEight}\} \geq 1 - \frac{\beta}{24n(K+1)},\\
        \Prob\{E_{\circledEight'}\} \geq 1 - \frac{\beta}{24n(K+1)}
    \end{gather*}
    where
    \begin{align}
        E_{\circledFour} & = \biggl\{\text{either } \sum\limits_{l=0}^{T-1}\sum\limits_{i=1}^n\widetilde\sigma_{i,l}^2 > \frac{\exp\left(-\gamma\mu T\right)V^2}{150A}\notag\\
        &\hphantom{{}= \biggl\{}\text{or } |\circledFour| \leq \frac{\exp\left(-\frac{\gamma\mu T}{2}\right)V}{5}\biggr\},\notag\\
        E_{\circledEight} & = \biggl\{\text{either } \sum\limits_{i=1}^{m-1}(\widehat\sigma_{i,T-1}^{\prime})^2 > \frac{\exp\left(-\gamma\mu (T-1)\right)V^2}{38400A}\notag\\
        &\hphantom{{}= \biggl\{}\text{or } |\circledEight| \leq \frac{\exp\left(-\frac{\gamma\mu (T-1)}{2}\right)V}{80}\biggr\}, \notag\\
        E_{\circledEight'} & = \biggl\{\text{either } \sum\limits_{i=1}^{m-1}(\widetilde\sigma_{i,T-1}^{\prime})^2 > \frac{\exp\left(-\gamma\mu (T-1)\right)V^2}{38400A}\notag\\
        &\hphantom{{}= \biggl\{}\text{or } |\circledEight^{\prime}| \leq \frac{\exp\left(-\frac{\gamma\mu (T-1)}{2}\right)V}{80}\biggr\}. \notag
    \end{align}
    Therefore, probability event $E_{T-1} \cap \widetilde{E}_{T-1,m-1} \cap E_{\circledFour} \cap E_{\circledEight} \cap E_{\circledEight'}$ implies
    \begin{align*}
        \left\|\frac{\gamma}{n}\sum\limits_{i=1}^{m-1} \theta_{i,T-1}^u\right\| & \leq \sqrt{\frac{1}{50} + \frac{1}{5} + \frac{1}{80}}\exp\left(-\frac{\gamma\mu(T-1)}{4}\right)\sqrt{V}\\
        & \leq \frac{1}{2}\exp\left(-\frac{\gamma\mu(T-1)}{4}\right)\sqrt{V},\\
        \left\|\frac{\gamma}{n}\sum\limits_{i=1}^{m-1} \omega_{i,T-1}^u\right\| & \leq \sqrt{\frac{1}{50} + \frac{1}{5} + \frac{1}{80}}\exp\left(-\frac{\gamma\mu(T-1)}{4}\right)\sqrt{V}\\
        & \leq \frac{1}{2}\exp\left(-\frac{\gamma\mu(T-1)}{4}\right)\sqrt{V}.
    \end{align*}
    This implies $\widetilde E_{T-1,m}$ and
    \begin{align*}
        \Prob\{E_{T-1} \cap \widetilde{E}_{T-1,m}\} & \geq \Prob\{E_{T-1} \cap \widetilde{E}_{T-1,m-1} \cap E_{\circledFour} \cap E_{\circledEight} \cap E_{\circledEight'}\} \\
        & = 1 - \Prob\left\{\overline{E_{T-1} \cap \widetilde{E}_{T-1,m-1}} \cup \overline{E}_{\circledFour} \cup \overline{E}_{\circledEight} \cup \overline{E}_{\circledEight'}\right\}\\
        & \geq 1 - \frac{(T-1)\beta}{K+1} - \frac{m\beta}{8n(K+1)}.
    \end{align*}
    Therefore, for all $m = 2,\ldots,n$ the statement holds and, in particular, $\Prob\{E_{T-1} \cap \widetilde{E}_{T-1, n}\} \geq 1 - \frac{(T-1)\beta}{K+1} - \frac{\beta}{8(K+1)}$, i.e., \eqref{ch3:eq:induction_inequality_str_mon_SEG_1_chapter_3} and \eqref{ch3:eq:induction_inequality_str_mon_SEG_2_chapter_3} hold. Taking into account \eqref{ch3:eq:nvsdibhjbsjdbvshbcjs_chapter_3} and $S_\theta = S_\omega = 0$ on $\widetilde{E}_{T-1,n}$, we conclude that $E_{T-1} \cap \widetilde{E}_{T-1, n} \cap E_{\circledOne} \cap E_{\circledFour} \cap E_{\circledSix'} \cap E_{\circledSeven'}$ implies
    \begin{equation*}
        V_T \leq 2\exp\left(-\frac{\gamma\mu}{2}T\right)V
    \end{equation*}
    that is equivalent to \eqref{ch3:eq:induction_inequality_str_mon_SEG_chapter_3} for $t = T$. Moreover,
    \begin{align*}
        \Prob\left\{E_T\right\} & \geq \Prob\left\{E_{T-1} \cap \widetilde{E}_{T-1,n} \cap E_{\circledOne} \cap E_{\circledFour} \cap E_{\circledSix'} \cap E_{\circledSeven'} \right\} \\
        & = 1 - \Prob\left\{\overline{E_{T-1} \cap \widetilde{E}_{n}} \cup \overline{E}_{\circledOne} \cup \overline{E}_{\circledFour} \cup \overline{E}_{\circledSix'} \cup \overline{E}_{\circledSeven'}\right\}\\
        & \geq 1 - \frac{(T-1)\beta}{K+1} - \frac{\beta}{8(K+1)} - 4\cdot \frac{\beta}{24n(K+1)} \geq 1 - \frac{T\beta}{K+1}.
    \end{align*}

    In other words, we showed that $\Prob\{E_k\} \geq 1 - k\beta/(K+1)$ for all $k = 0,1,\ldots,K+1$. For $k = K+1$ we have that with probability at least $1 - \beta$
    \begin{equation*}
       \|x^{K+1} - x^*\|^2 \leq V_{K+1} \leq 2\exp\left(-\frac{\gamma\mu(K+1)}{2}\right)V.
    \end{equation*}

    Finally, if $\gamma$ equals the minimum in \eqref{ch3:eq:gamma_SEG_str_mon_chapter_3},
    then with probability at least $1-\beta$
    \begin{align*}
        \|x^{K+1} - x^*\|^2 & \leq 2\exp\left(-\frac{\gamma\mu(K+1)}{2}\right)V\\
        & = 2V\max\Bigg\{\exp\left(-\frac{K+1}{144\cdot 10^6A^2}\right), \exp\left(-\frac{\mu(K+1)}{12 L}\right),\\
        &\qquad\qquad \exp\left(-\frac{\mu\sqrt{n}(K+1)}{30000 L A}\right), \frac{1}{B_K} \Bigg\}\\
        & = \cO\Bigg(\max\Bigg\{V\exp\left(- \frac{K}{\ln^2(nK/\beta)}\right), V\exp\left(- \frac{\mu  K}{L}\right),\\
        &\qquad\qquad V\exp\left(- \frac{\mu \sqrt{n} K}{L \ln(nK/\beta)}\right), \frac{\sigma^2\ln^{2(p-1)/p}(nK/\beta)\ln^2B_K}{n^{2(p-1)/p}K^{2(p-1)/p}\mu^2}\Bigg\}\Bigg).
    \end{align*}
    To get $\|x^{K+1} - x^*\|^2 \leq \varepsilon$ with probability $\geq 1 - \beta$, $K$ should be
    \begin{align}
        K = \cO\biggl(\max\biggl\{ & \left(\frac{L}{\sqrt{n}\mu}+\ln\left(\frac{nL}{\mu \beta}\ln\frac{V}{\varepsilon}\right)\right)\ln\frac{V}{\varepsilon}\ln\left(\frac{nL}{\mu \beta}\ln\frac{V}{\varepsilon}\right), \frac{L}{\mu}\ln\frac{V}{\varepsilon},\notag\\
        & \frac{1}{n}\left(\frac{\sigma^2}{\mu^2\varepsilon}\right)^{p/(2(p-1))}\ln \left(\frac{n}{\beta} \left(\frac{\sigma^2}{\mu^2\varepsilon}\right)^{p/(2(p-1))}\right)\ln^{p/(p-1)}\left(B_\varepsilon\right)\biggr\}\biggr), \notag
    \end{align}
    where
    \begin{equation*}
        B_\varepsilon = \max\left\{2, \frac{V}{\varepsilon \ln \left(\frac{1}{\beta} \left(\frac{\sigma^2}{\mu^2\varepsilon}\right)^{p/(2(p-1))}\right)}\right\}.
    \end{equation*}
\end{proof}

\newpage
\section{Numerical Experiments}
In this section we provide numerical experiments for the following simple problem:
\begin{equation}
    \label{ch3:eq:experiment_problem_chapter_3}
    \min_{x \in B_r(\hat{x})} f(x),
\end{equation}
where a radius  $r=1$, a central point $\hat{x} = (3,3,\dots,3)^{\top} \in \mathbb{R}^{10}$, and $f(x) = \frac{1}{2}\|x\|^2$, $f(x;\xi) = \frac{1}{2}\|x\|^2 +\langle \xi, x\rangle$, where $\xi$ comes from the symmetric Levy ${p}$-stable distribution ${p} = \frac{3}{2}$. We use the following parameters: $\gamma = 0.001$, $x^0 = \hat{x} + r \frac{e}{\|e\|}$, where $e = (1,1,\dots,1)^{\top}$. We tried three values of $\lambda$: $0.1$, $0.01$ and $0.001$. 

\begin{figure}[h]
    \centering
    \includegraphics[width=0.32\textwidth]{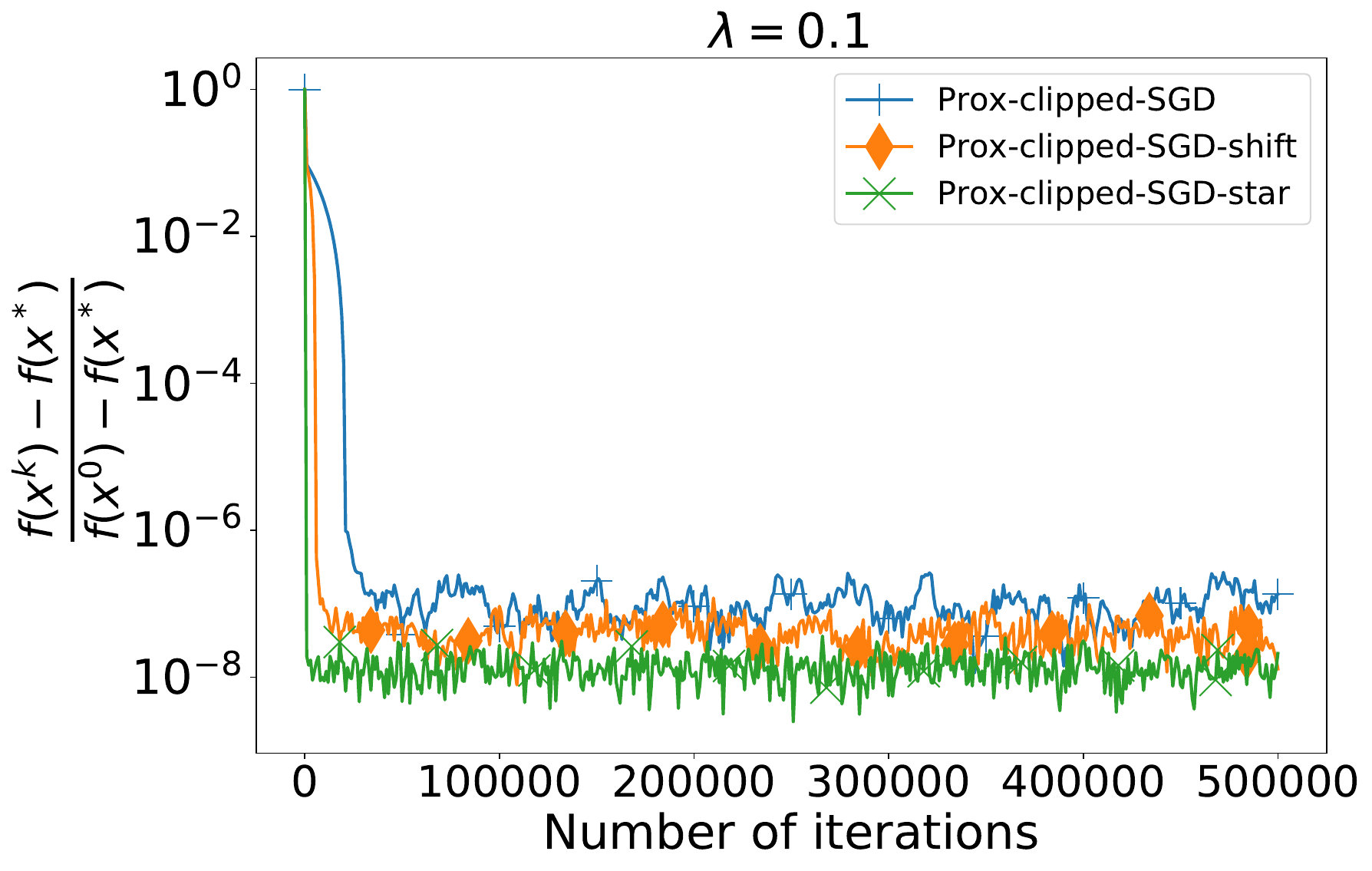}
    \includegraphics[width=0.32\textwidth]{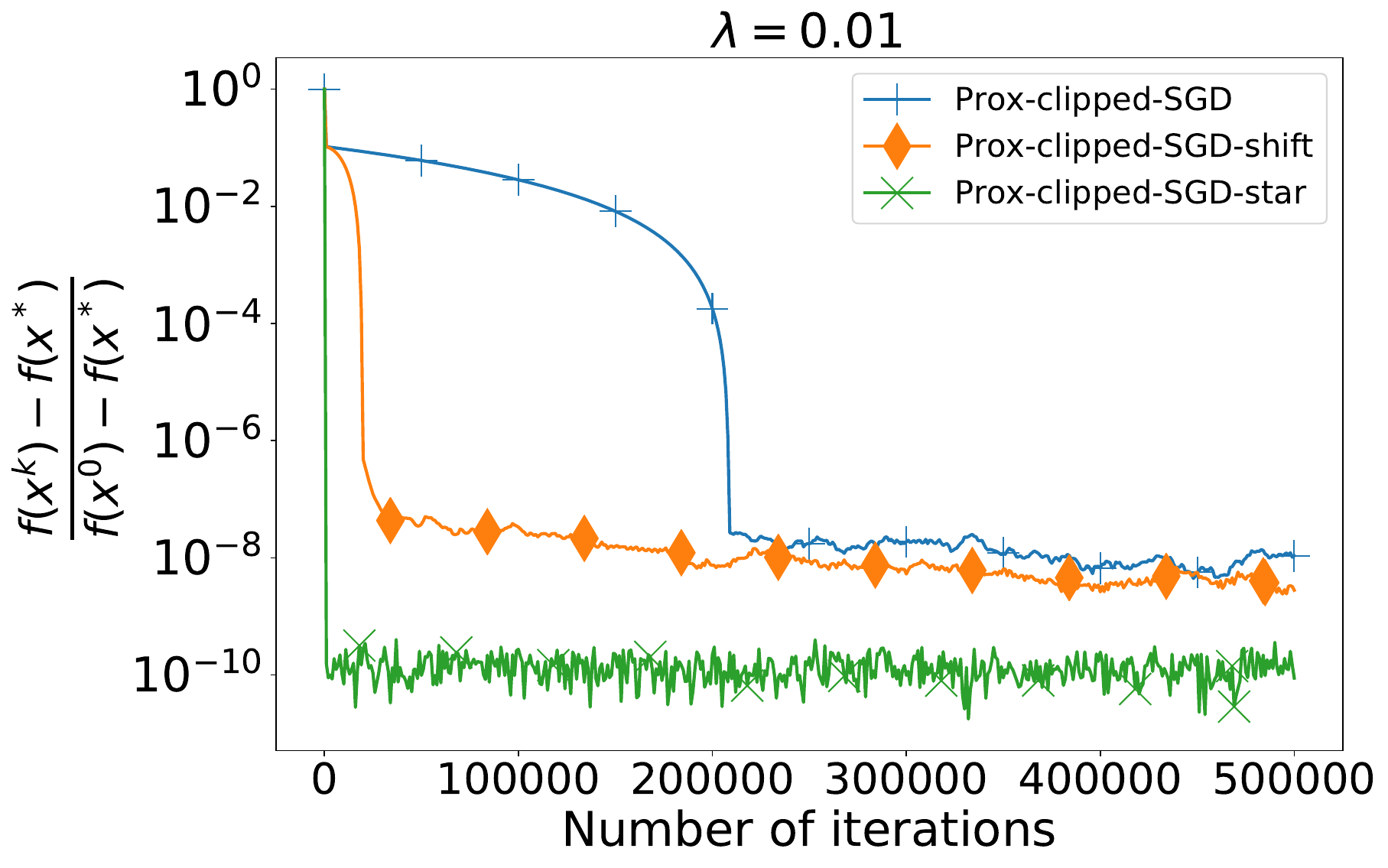}
    \includegraphics[width=0.32\textwidth]{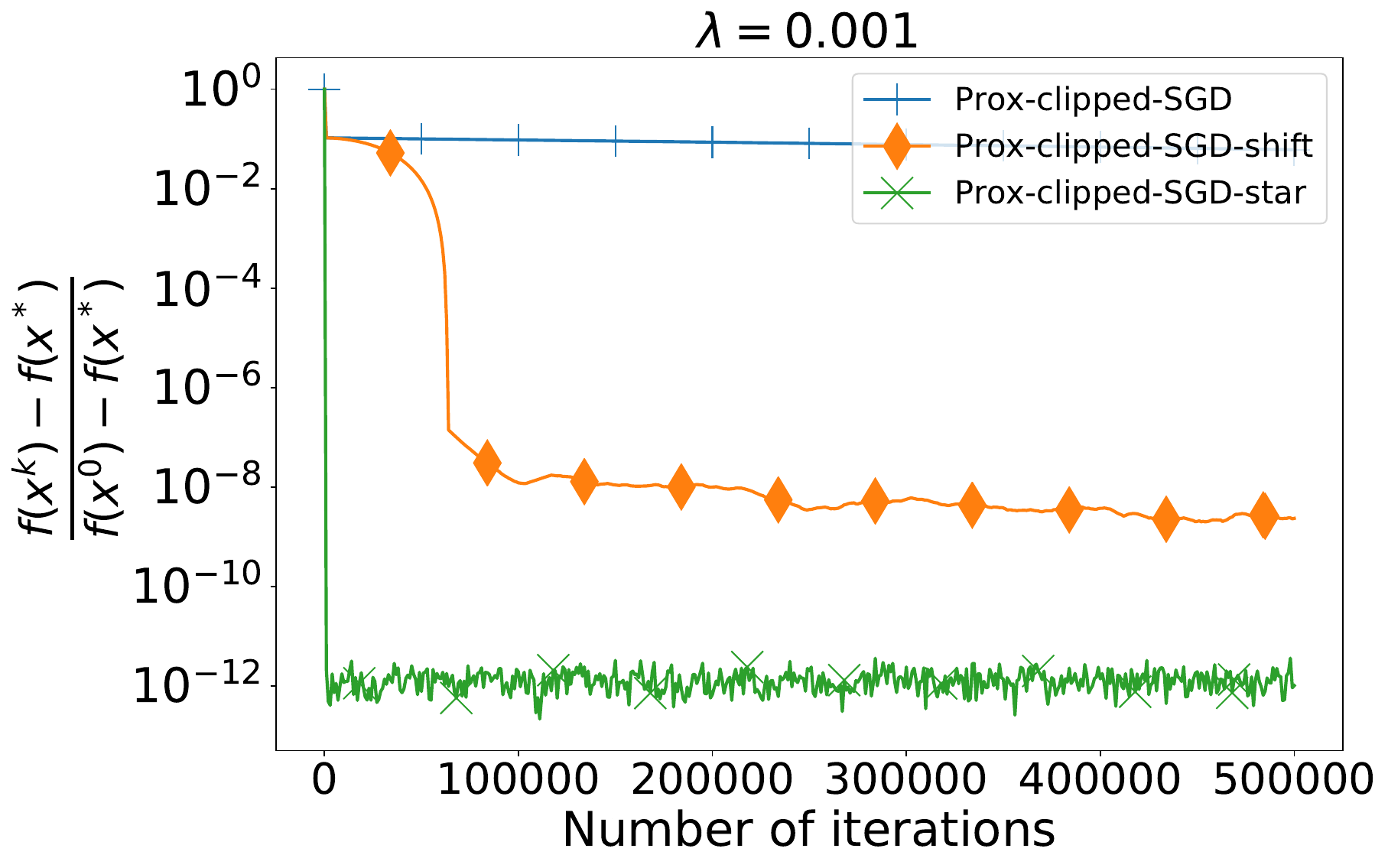}
    \caption{Comparison between performances of \algname{Prox-clipped-SGD}, \algname{Prox-clipped-SGD-star}, \algname{Prox-clipped-SGD-shift} in solving problem \eqref{ch3:eq:experiment_problem_chapter_3} with fixed clipping level for each of them $\lambda \in \{0.1, 0.01, 0.001\}$.}
    \label{ch3:fig:exp_1_chapter_3}
\end{figure}

In our numerical experiments (see Figure~\ref{ch3:fig:exp_1_chapter_3}), we observe that the na\"ive \algname{Prox-clipped-SGD} converges slower than \algname{Prox-clipped-SGD-star} and \algname{Prox-clipped-SGD-shift}. Moreover, when the clipping level is small \algname{Prox-clipped-SGD} converges extremely slow, while \algname{Prox-clipped-SGD-shifts} takes some time to learn the shift and then converges to much better accuracy. We also see that the smaller clipping level is, the better accuracy \algname{Prox-clipped-SGD-star} achieves. For \algname{Prox-clipped-SGD-shift} we observe the same phenomenon when we reduce $\lambda$  from $0.1$ to $0.01$ . We expect the improvement in the accuracy even further if we decrease the stepsizes $\gamma$ and $\nu$ .
  \endgroup

\chapter{Supplementary Material and Proofs for Chapter 4}
  \begingroup
  \thesischaptersetup{4}%
  \stepEightAAppendixColors{4}%
  
\newpage
\section{Additional Related Work}
\paragraph{Gradient clipping and normalization.} Gradient clipping has been applied successfully also in zero-order optimization \parencite{kornilov2024accelerated}, bandit and RL literature \parencite{bubeck2013bandits,cayci2024provably}, online learning \parencite{zhang2022parameter} and differential privacy \parencite{abadi2016deep,sha2024clip}. Some other works, which use normalization and gradient clipping under heavy-tailed noise include \parencite{NonlinearHP2023Armacki,puchkin2024breaking,liu2024nonconvex}.

\section{About Normalized SGD and Different Variants of Momentum}
In this section, we present a general framework for momentum-based methods. The study of momentum has a long history, particularly in convex optimization, where \textcite{nesterov1} demonstrated that an appropriate use of momentum in gradient descent can accelerate convergence. In contrast, our focus is on momentum in the setting of stochastic non-convex optimization. The general update rule is given by
\begin{equation*}
x^{k+1} = x^k - {\gamma} \frac{g^k}{\norm{g^k}},
\end{equation*}
where $g^k$ denotes the momentum term. Different choices of $g^k$ give rise to different algorithms, several of which we now describe.

\textbf{\algname{NSGD}} is the standard stochastic gradient method. Setting $g^k = \nabla f(x^k;\xi^{k})$ in the general framework immediately recovers \algname{NSGD}. It is well known that \algname{NSGD} converges only to a neighborhood of the solution, with radius proportional to $\sigma$ \parencite{yang2023two_sides,Hubler2024clip_to_norm}. One possible remedy is to employ mini-batching, which reduces the variance. However, this strategy cannot be applied if only a single sample per iteration is allowed.

\textbf{\algname{NSGDM}} is the solution of the problem described above and was proposed and analyzed by \textcite{MomentumImprovesNSGD_Cutkosky_2020}. The momentum term is defined as
\begin{equation*}
g^k = (1-\eta)g^{k-1} + \eta \nabla f(x^k;\xi^{k}).
\end{equation*}
In their work, the authors proved that \algname{NSGDM} converges to an $\varepsilon$-stationary point and achieves a sample complexity of $\cO(\varepsilon^{-4})$, which is optimal according to \textcite{arjevani2023lower}. Subsequently, \textcite{Hubler2024clip_to_norm} extended this result to the setting where the stochastic gradients have bounded $p$-th moments, deriving a sample complexity of $\cO(\varepsilon^{-(3p-2)/(p-1)})$.

\textbf{\algname{NSGDM} with clipping} was proposed for the general setting where the stochastic gradients have bounded $p$-th moments. The momentum term is modified as
\begin{equation*}
g^k = (1-\eta)g^{k-1} + \eta\,\texttt{clip}\left(\nabla f(x^k;\xi^{k}), \lambda\right).
\end{equation*}
\textcite{cutkosky2021high} established high-probability convergence guarantees and showed that this method achieves a sample complexity of $\cO(\varepsilon^{-(3p-2)/(p-1)})$, which is optimal and matches the lower bounds established by \textcite{zhang2020adaptive}. It is worth noting that the clipping operator is not strictly necessary: normalization of the momentum term is sufficient to guarantee convergence under Assumption~\ref{ch4:ass:p_BCM_chapter_4} \parencite{Hubler2024clip_to_norm,liu2024nonconvex}.

\textbf{\algname{NIGT}} is a variant of \algname{NSGD} that incorporates momentum and implicit gradient transport, proposed by \textcite{MomentumImprovesNSGD_Cutkosky_2020}. In this method, the momentum term is modified by changing the point at which the stochastic gradient is evaluated:
\begin{align*}
g^k & = (1-\eta)g^{k-1} + \eta \nabla f(y^k;\xi^{k}), \quad \text{where} ~~ y^k \\
& = x^k + \frac{1-\eta}{\eta}(x^k - x^{k-1}).
    \end{align*}
Thanks to this modification, and under the assumption that the objective function has a Lipschitz continuous Hessian, \textcite{MomentumImprovesNSGD_Cutkosky_2020} proved that \algname{NIGT} achieves a sample complexity of $\cO\left(\varepsilon^{-7/2}\right)$, which improves upon that of \algname{NSGDM}.

\textbf{\algname{NIGT} with clipping} was introduced by \textcite{cutkosky2021high} to enable high-probability analysis under the assumption of bounded $p$-th moments of the stochastic gradients. The only modification compared to standard \algname{NIGT} is the inclusion of a clipping operator, which is required from a theoretical standpoint. As shown by \textcite{cutkosky2021high}, \algname{NIGT} with clipping achieves a sample complexity of $\cO\left(\varepsilon^{-(5p-3)/(2p-2)}\right)$, which is optimal according to Theorem~\ref{ch4:thm:lower_bounds_chapter_4}.

\textbf{\algname{NSGD} with MVR} is a variant of \algname{NSGD} that incorporates momentum variance reduction (MVR). The \algname{STORM}/MVR method was proposed by \textcite{cutkosky-orabona19} to improve the sample complexity compared to standard SGD. More recently, \textcite{he2025complexity} established that \algname{NSGD} with MVR achieves a sample complexity of $\cO\left(\varepsilon^{-(2p-1)/(p-1)}\right)$, which matches ours in terms of dependence on $\varepsilon$ (Theorem~\ref{ch4:thm:nsgdm_hess_chapter_4}).

The methods discussed above rely solely on first-order information. A natural question is whether incorporating second-order information can further improve the sample complexity.

\textbf{\algname{NSGDHess}} is a variant of Normalized SGD that incorporates Hessian-corrected momentum, proposed by \textcite{BetterSGDUsingSOM_Tran_2021}. In this method, a stochastic Hessian-vector product is introduced into the momentum update:
\begin{equation*}
g^k = (1-\eta)\left(g^{k-1} + \nabla^2 f(x^k;\xi^{k})\cdot (x^k-x^{k-1}) \right) + \eta \nabla f(y^k;\xi^{k}).
\end{equation*}

\textbf{Refined \algname{NSGDHess}.} In our refined version, the momentum term is modified by adjusting the point at which the stochastic Hessian is evaluated and by using an independent fresh sample $\hat{\xi}^{k} \sim {\cD}$ for the Hessian-vector product computation:
\begin{align*}
    \hat{x}^k & = {q^k} x^k + (1-{q^k})x^{k-1},\quad \text{where  }~~ {q^k} \sim {\cU}([0,1]), \\
    g^k & = (1-\eta)\left(g^{k-1} + {\nabla^2 f(\hat{x}^k;\hat{\xi}^{k})}\cdot (x^k-x^{k-1}) \right) + \eta \nabla f(x^k;\xi^{k}).
\end{align*}

\newpage 
\section{Technical Lemmas}
We list the following technical lemmas, which are useful for our analysis in the subsequent sections.
\begin{lemma}[Lemma 10 from \textcite{Hubler2024clip_to_norm}]
\label{ch4:lem:von_Bahr_inequality_chapter_4}
    Let $p \in [1,2]$ and let $X_1,\dots, X_n \in \R^d$ be a martingale difference sequence, i.e.,
    \begin{equation*}
        \Exp\left[X_j\mid X_{j-1},\dots, X_1\right] = 0 \ \text{a.s. for all } j=1,\dots, n,
    \end{equation*}
    satisfying 
    \begin{equation*}
        \Exp\left[\norm{X_j}^p\right] < \infty \quad \quad \text{for all } {j}=1,\dots, n.
    \end{equation*}
    Define $S_n \eqdef \sum^n_{j=1}X_j$, then 
    \begin{equation*}
        \Exp\left[\norm{S_n}^p\right] \leq 2\sum^n_{j=1}\Exp\left[\norm{X_j}^p\right].
    \end{equation*}
\end{lemma}

\begin{lemma}[Lemma 10 from \textcite{cutkosky2021high}]
\label{ch4:lem:inequlity_from_cutkosky_chapter_4}
    Let  $X_1,\dots, X_n \in \R^d$  be a sequence of random vectors. Define the sequence of real numbers $w_1,\dots, w_n$ recursively 
    \begin{enumerate}
        \item $w_0 = 0$.
        \item If $\sum^{j-1}_{i=1}X_i \neq 0$, then we set: 
        \begin{equation*}
            w_j = {\sign}\left(\sum^{j-1}_{i=1}w_i\right)\frac{\left\la \sum^{j-1}_{i=1}X_i, X_j\right\ra}{\norm{\sum^{j-1}_{i=1}X_i}}.
        \end{equation*}
        \item If $\sum^{j-1}_{i=1}X_i = 0$, set $w_j = 0$.
    \end{enumerate}
    Then  $|w_j| \leq \norm{X_j}$ for all $j=1,\dots, n$, and, with $S_n \eqdef \sum^n_{j=1}X_j$, 
    \begin{equation}
        \norm{S_n}\leq \left|\sum^n_{j=1} w_j\right| + \sqrt{\max_{j \in [n]}\norm{X_j}^2 + \sum^n_{j=1}\norm{X_j}^2}.
    \end{equation}
\end{lemma}

\begin{lemma}[Bernstein inequality]
\label{ch4:lem:Bernstain_inequality_chapter_4}
    Let  $X_1,\dots, X_n \in \R$  be a martingale difference sequence, i.e., $\Exp\left[X_j\mid X_{j-1},\dots, X_1\right] = 0$ a.s. for all $j=1,\dots, n$. Assume that conditional variances $\sigma^2_{j} \eqdef \Exp\left[X^2_j| X_{j-1},\dots, X_1\right]$ exist and are bounded, and assume that there exists a deterministic constant $c > 0$ such that $|X_j| \leq c$ almost surely for all  $j=1,\dots, n$. Define $S_n \eqdef {\sum^n_{j=1}X_j}$; then for all $b > 0$ and $G> 0$,
    \begin{equation}
        \Prob\left\{\left|S_n\right| > b \text{ and } \sum^n_{j=1}\sigma^2_{j} \leq G\right\} \leq 2\exp\left( - \frac{b^2}{2G + 2cb/3}\right).
    \end{equation}
\end{lemma}

\begin{lemma}[Lemma 5.1 from \textcite{sadiev2023high}]
\label{ch4:lem:bound_variance_magnitude_chapter_4}
    Let $\lambda > 0$ and $X \in \R^d$ be a random vector and $\widetilde{X} \eqdef \texttt{clip}\left(X, \lambda\right)$. Then, $\norm{\widetilde{X} - \Exp{\widetilde{X}}} \leq 2\lambda$. Moreover, if for some $\sigma \geq 0$ and $p \in {(1,2]}$ we have $\Exp{X} = x \in \R^d$, $\Exp{\norm{X-x}^p}\leq \sigma^p$, and $\norm{x} \leq \lambda/2$, then 
    {
    \begin{align*}
        \norm{\Exp{\widetilde{X}} -x}
        &\leq \frac{2^p\sigma^p}{\lambda^{p-1}},\\
        \Exp{\norm{\widetilde{X} - x}^2}
        &\leq 18\lambda^{2-p}\sigma^p,
        \qquad
        \Exp{\norm{\widetilde{X} - \Exp{\widetilde{X}}}^2}
        \leq 18\lambda^{2-p}\sigma^p.
    \end{align*}
    }
\end{lemma}

\newpage
\section{Additional Proofs for Section~\ref{ch4:sec:lower_bounds_chapter_4}}
\label{ch4:app:missing_proofs_lower_bounds_chapter_4}
In this section we provide lower bounds for {a particular} class of functions and the oracle class. Our proof is inspired by the work of \textcite{arjevani2020second}.

\textbf{Auxiliary Lemmas.} Now we state two helpful lemmas to prove lower bound. The first lemma is about the number of iteration to reveal all coordinates, when zero-respecting algorithms use the information received from the oracle forming \textit{probability-$\rho$ zero-chain}. 
\begin{definition}
    A collection of derivative estimators ${\nabla^1 f}(x,\xi)$, ${\nabla^2 f}(x,\xi)$, $\dots$, ${\nabla^q f}(x,\xi)$ for a function $f$ forms a probability-$\rho$ zero-chain if 
    \begin{equation*}
        \Prob\left\{\exists~ x ~|~ \text{prog}\left({\nabla^1 f}(x,\xi), \dots, {\nabla^q f}(x,\xi)\right) = \text{prog}_{1/4}(x) +1 \right\} \leq \rho 
    \end{equation*}
    and 
    \begin{equation*}
        \Prob\left\{\exists~ x ~|~ \text{prog}\left({\nabla^1 f}(x,\xi), \dots, {\nabla^q f}(x,\xi)\right) = \text{prog}_{1/4}(x) + i \right\} = 0,~ i > 1.
    \end{equation*}
\end{definition}
The second lemma is about the main properties of the \textit{worst} function from the class ${\cF}(\Delta, L_{1:q})$. 

\begin{lemma}[\textcite{arjevani2020second}]
\label{ch4:lem:lower_bound_1_chapter_4}
    Let ${\{\nabla^r f(x;\xi)\}_{r\in[q]}}$ be a collection of probability-$\rho$ zero-chain derivative estimators for $f: \R^{T}\rightarrow \R$, and let
    {
    \[
        \text{O}^q_f
        =\left({\nabla^r f}(x,\xi)\right)_{r\in[q]}.
    \]
    }
    Let ${\{x^{(t)}_{A[\text{O}_f^q]}\}}$ be a sequence of queries produced by a zero-respecting algorithm $A$ interacting with $\text{O}^q_f$. Then, with probability at least $1-\beta$,
    \begin{equation*}
        \text{prog}\left(x^{(t)}_{A\left[O^q_f\right]}\right) < T, \qquad \text{ for all } \quad t \leq \frac{T-\log\frac{1}{\beta}}{2\rho} .
    \end{equation*}
\end{lemma}

\begin{lemma}[\textcite{carmon2020lower}]
\label{ch4:lem:lower_bound_2_chapter_4}
Let $h: \R^T\rightarrow \R$ be the following function:
\begin{equation*}
    h(x) = -\Psi(1)\Phi(x_1) + \sum^T_{i=2}\left(\Psi(-x_{i-1})\Phi(-x_i) - \Psi(x_{i-1})\Phi(x_i)\right),
\end{equation*}
where 
\begin{equation*}
    \Psi(x) = \begin{cases}
        0,& x \leq \frac{1}{2};\\
        \exp\left(1-\frac{1}{(2x-1)^2}\right),& x > \frac{1}{2},
    \end{cases}
    \quad\text{and}\quad \Phi(x) = \sqrt{e}\int\limits^x_{-\infty} e^{-t^2/2}\, dt. 
\end{equation*}
Then the function $h$ satisfies the following properties: 
    \begin{itemize}
        \item $h(0) - \inf\limits_x h(x) \leq \Delta_0 T$, where $\Delta_0 = 12$. 
        \item For $q \geq 1$, the $q$-th order derivatives of $h$ are $\ell_q$-Lipschitz continuous, where $\ell_q \leq \exp\left(\frac{5}{2}q\log q + cq\right)$ for a numerical constant $c < \infty$. 
        \item For all $x \in \R^{T}$, $q \in \mathbb{N}
        $ and $i \in [T]$, we have $\|\nabla^q_i h(x)\|_{\text{op}} \leq \ell_{q-1}$, where $\ell_0 = 23$. 
        \item For all $x \in \R^T$ and $q \in \mathbb{N} $, $\text{prog}\left(\nabla^q h(x)\right) \leq \text{prog}_{1/2} (x) +1 $.
        \item For all $x \in \R^T$, if $\text{prog}_1(x) < T$, then $\|\nabla h(x)\| \geq |\nabla_{\text{prog}_1(x)+1} h(x)| > 1$.
    \end{itemize}
\end{lemma}

\begin{lemma}
\label{ch4:lem:lower_bound_3_chapter_4}
{
    For each $r \in [q]$, define the derivative estimator
    \begin{equation*}
        \left[{\nabla^r h}(x,\xi)\right]_i \eqdef \left(1 + \mathbb{I}\left\{ i > \text{prog}_{1/4}(x)\right\}\left(\frac{\xi}{\rho}-1\right)\right) \cdot \nabla^r_i h(x),
    \end{equation*}
    where $\xi \sim \operatorname{Bernoulli}(\rho)$. Then,
    $\{\nabla^r h(x,\xi)\}_{r \in [q]}$ forms a probability-$\rho$
    zero-chain. Moreover, for any $p \in [1,2]$ and every $r \in [q]$,
    {
    \begin{align*}
        \Exp\left[{\nabla^r h}(x,\xi)\right]
        &= {\nabla^r h}(x),\\
        \MoveEqLeft[1] \Exp\left[
        \left\|{\nabla^r h}(x,\xi)-{\nabla^r h}(x)\right\|_{\mathrm{op}}^p
        \right] \\
        &\leq \frac{2\ell_{r-1}^p}{\rho^{p-1}},
        \qquad x\in\R^T.
    \end{align*}
    }
}
\end{lemma}
\begin{proof}
{
    Since $\Exp[\xi/\rho-1]=0$, the estimator is unbiased:
    \begin{equation*}
        \Exp\left[{\nabla^r h}(x,\xi)\right] = (1+0)\cdot \nabla^r h(x) = \nabla^r h(x).
    \end{equation*}
    Let $m=\operatorname{prog}_{1/4}(x)$. By
    Lemma~\ref{ch4:lem:lower_bound_2_chapter_4},
    $\operatorname{prog}(\nabla^r h(x))\leq m+1$, and hence
    \begin{equation*}
        \left[{\nabla^r h}(x,\xi)\right]_i
        = \nabla^r_i h(x)=0,
        \qquad i>m+1.
    \end{equation*}
    Thus the random multiplier leaves every already active subtensor
    $i\leq m$ unchanged and can modify only the next potentially active
    subtensor $i=m+1$. This subtensor is revealed only when $\xi=1$.
    Consequently,
    \begin{equation*}
        \Prob\left\{\exists x~\middle|~
        \operatorname{prog}\left({\nabla^1 h}(x,\xi),\dots,
        {\nabla^q h}(x,\xi)\right)
        =\operatorname{prog}_{1/4}(x)+1\right\}\leq\rho,
    \end{equation*}
    while the probability of revealing any subtensor beyond $m+1$ is
    zero. Therefore the collection is a probability-$\rho$ zero-chain.

    It remains to bound the central moment. Only the $(m+1)$st
    subtensor can differ from its deterministic counterpart, so for
    any $p\in[1,2]$ and $r\in[q]$,
    \begin{align*}
        \MoveEqLeft[1] \Exp\left[
            \left\|{\nabla^r h}(x,\xi)-\nabla^r h(x)\right\|_{\mathrm{op}}^p
        \right] \\
        &\leq
        \Exp\left[\left|\frac{\xi}{\rho}-1\right|^p\right]
        \left\|\nabla^r_{m+1}h(x)\right\|_{\mathrm{op}}^p\\
        &\leq
        \ell_{r-1}^p
        \frac{\rho(1-\rho)^p+(1-\rho)\rho^p}{\rho^p}\\
        &\leq \frac{2\ell_{r-1}^p}{\rho^{p-1}}.
    \end{align*}
}
\end{proof}

\textbf{Main Theorem.} Now we are ready to state and prove the main theorem in this section, i.e. Theorem~\ref{ch4:thm:lower_bounds_chapter_4}. 

\begin{proof}
First of all, we fix all parameters: $q\in \mathbb{N}, \Delta, L_{1:q}, \sigma_{1:q} > 0 $  and $\varepsilon > 0$. Next, we rescale our function $h$ as follows: $h^{\star}(x) = \nu h(\chi x)$, where $\nu> 0$ and $\chi >0$.
According to Lemma~\ref{ch4:lem:lower_bound_2_chapter_4}, selecting $T = \left\lfloor \frac{\Delta}{\nu\Delta_0}\right\rfloor$, we have  for any $r \in [q-1]$ 
\begin{align*}
    h^{\star}(0) - \inf\limits_x h^{\star}(x) & = \nu\left(h(0) - \inf\limits_x h( \chi x)\right) \leq \nu \Delta_0 T \leq \Delta;\\
    \|\nabla^{r+1} h^{\star}(x)\|_{\text{op}} & = \nu \chi^{r+1}\|\nabla^{r+1} h(\chi x)\| \leq \nu \chi^{r+1}\ell_r;\\
    \|\nabla h^{\star}(x)\| & = \nu\chi\|\nabla h(\chi x)\| \geq \nu \chi, ~\forall~x:~\text{prog}_1(\chi x) < T.  
\end{align*}
Since the scaled estimators ${\nabla^r h^{\star}}(x,\xi) = \nu\chi^r{\nabla^r h}(\chi x,\xi)$, $r\in [q]$, form a probability-$\rho$ zero-chain, Lemma~\ref{ch4:lem:lower_bound_1_chapter_4} with $\beta = 1/2$ gives $\text{prog}_1(\chi x^{t}_{A[O_h^q]}) < T$ with probability at least $1/2$ for all $t\leq (T-1)/(2\rho)$ (note that $\log 2 < 1$). Hence  
\begin{align*}
\Exp\left[\left\|\nabla h^{\star}(x^{t}_{A[O_h^q]})\right\|\right] & = \nu \chi \Exp\left[\left\|\nabla  h(\chi x^{t}_{A[O_h^q]})\right\|\right] \\
& \geq \frac{\nu \chi}{2},\quad \forall~t \leq \frac{T-1}{2\rho}.
    \end{align*}
The $p$-th  central moment of the scaled derivatives estimators is bounded as 
\begin{align*}
    \MoveEqLeft[1] \Exp\left[\left\|{\nabla^r  h^{\star}}(x,\xi) - \nabla^r  h^{\star} (x)\right\|_{\text{op}}^{p}\right] \\
    & \leq \nu^{p}\chi^{rp} \Exp\left[\left\|{\nabla^r h}(\chi x,\xi) - \nabla^r h (\chi x)\right\|_{\text{op}}^{p}\right]\\
    & \leq \frac{2\nu^{p}\chi^{rp}\ell_{r-1}^{p}}{\rho^{p-1}},
\end{align*}
where in the last inequality we applied Lemma~\ref{ch4:lem:lower_bound_3_chapter_4}.
Thus, we have the set of constraints for parameters $\nu$ and $\chi$ for all $r \in [q]$:
\begin{equation*}
    \begin{cases}
        \nu \Delta_0 T \leq \Delta; \\
        \nu \chi^{r+1}\ell_r \leq L_r;\\
        \frac{\nu\chi}{2} \geq \varepsilon; \\
        \nu^{p}\chi^{p r}\ell^{p}_{r-1}\frac{2}{\rho^{p-1}} \leq \sigma^{p}_r.
    \end{cases}
\end{equation*}
We will resolve the system of inequalities step by step. The first step is to set $\nu = \frac{2\varepsilon}{\chi}$.  For $r = 1$, we have 
\begin{align*}
&\nu^{p}\chi^{p}\ell_0^{p} \cdot \frac{2}{\rho^{p-1}} \leq \sigma^{p}_1\\
&\quad \Rightarrow\quad \rho = \min\left\{\left(\frac{\varepsilon}{\sigma_1}\right)^{p/(p-1)}\cdot 2^{(p+1)/(p-1)}\ell_0^{p/(p-1)}, 1\right\} .
    \end{align*} 
To set  $\chi$, we need to find its specific value, which can satisfy the following constraints: for any $r \in \{2, \dots, q\} $ and  any $r' \in [q]$ 
\begin{gather*}
\chi^{r+1} \leq \frac{L_r}{\nu \ell_r} = \frac{\chi L_r}{2\varepsilon\ell_{r}},\\
\nu^{p}\chi^{p r}\ell^{p}_{r-1}\frac{2}{\rho^{p-1}} \overset{(\ast)}{\leq} \chi^{p(r-1)}\left(\frac{\sigma_1 \ell_{r-1}}{\ell_0}\right)^p \leq \sigma^{p}_r,
\end{gather*}
where in the inequality $(\ast)$ we used  $\frac{2}{\rho^{p-1}}\leq \frac{\sigma^p_1}{\nu^p\chi^p\ell_0^p}$. Therefore, we obtain 
\begin{equation*}
    \chi = \min_{{r' \in [q];}~{r \in \{2,\dots, q\}}}\min\left\{ \left(\frac{\ell_0\sigma_r}{\ell_{r-1}\sigma_1}\right)^{1/(r-1)}, \left(\frac{L_{r'}}{2\varepsilon\ell_{r'}}\right)^{1/r'}\right\}.
\end{equation*}
Then, assuming $T \geq 3$ and $\left(\frac{\varepsilon}{\sigma_1}\right)^{p/(p-1)}\cdot 2^{(p+1)/(p-1)}\ell_0^{p/(p-1)} \leq 1$, we get
\begin{align*}
    \frac{T-1}{2\rho} & =\frac{1}{2\rho}\left(\left\lfloor\frac{\Delta\chi}{2\Delta_0\varepsilon}\right\rfloor -1 \right) \geq  \frac{\Delta \chi}{8\rho\Delta_0\varepsilon} \\
    & = \left(\frac{\sigma_1}{\varepsilon\ell_0}\right)^{p/(p-1)} \cdot \frac{\Delta}{2^{(4p-2)/(p-1)}\Delta_0\varepsilon}\\
    &\quad\cdot \min_{{r' \in [q];}~{r \in \{2,\dots, q\}}}\min\left\{ \left(\frac{\ell_0\sigma_r}{\ell_{r-1}\sigma_1}\right)^{1/(r-1)}, \left(\frac{L_{r'}}{2\varepsilon\ell_{r'}}\right)^{1/r'}\right\}\\
    & \geq\Omega(1)\cdot\frac{\Delta}{\varepsilon}\left(\frac{\sigma_1}{\varepsilon}\right)^{p/(p-1)}  \min\left\{ \min_{r \in \{2,\dots, q\}}\left(\frac{\sigma_r}{\sigma_1}\right)^{1/(r-1)}, \min_{r' \in [q]}\left(\frac{L_{r'}}{\varepsilon}\right)^{1/r'}\right\}.
\end{align*}
Since $h$ and $h^{\star}$ are functions of $T$ arguments, the dimension of the problem is equal to $T$. Here $\Omega(1)$ hides factors that depend only on $p$ and $q$. This concludes the proof. 

\end{proof}

\newpage
\section{Additional Proofs for Section~\ref{ch4:sec:NSGD_SOM_conv_chapter_4}}
\label{ch4:app:missing_proofs_NSGD_SOM_chapter_4}
In this section we study Normalized SGD with Hessian correction (\algname{NSGDHess}) (see Algorithm~\ref{ch4:alg:NSGD_SOM_chapter_4}). 

\textbf{Auxiliary Lemmas.} To show the convergence guarantees for Algorithm~\ref{ch4:alg:NSGD_SOM_chapter_4}, we state and prove auxiliary lemmas. The first one is well-known \textit{Descent Lemma} for Normalized SGD. A proof of this lemma is given by \textcite{MomentumImprovesNSGD_Cutkosky_2020,Hubler2024clip_to_norm}, but for completeness we restate and reprove it. The second one is devoted to the bounds on \textit{error} term from \textit{Descent Lemma}. 

\begin{lemma}
\label{ch4:lem:descent_lemma_chapter_4}
    Let Assumptions~\ref{ch4:assum:lower_bounded_chapter_4} and~\ref{ch4:assum:l_smooth_chapter_4} hold. Then for any selection of stepsize $\gamma >0$ the iterates $\{x^{k}\}_{k=0}^{K}$ generated by Algorithm~\ref{ch4:alg:NSGD_SOM_chapter_4} satisfy
    \begin{equation}
        \gamma \sum^{K-1}_{k=0}\|\nabla f(x^{k})\| + \Delta_{K} \leq \Delta_0 +2\gamma \sum^{K-1}_{k=0}\|\hat{e}_k\| + \frac{\gamma^2 L K}{2},
    \end{equation}
    where functional gap is defined as $\Delta_k \eqdef f(x^{k}) -f^{\inf}$, error term is defined as $\hat{e}_k \eqdef g^{k} - \nabla f(x^{k})$.
\end{lemma}
\begin{proof}
    According to the update rule for $x^{k}$ and Assumption~\ref{ch4:assum:l_smooth_chapter_4}, we have 
    \begin{align*}
    f(x^{k+1}) & \leq f(x^{k}) +\la \nabla f(x^{k}), x^{k+1} -x^{k} \ra + \frac{L}{2}\|x^{k+1} - x^{k}\|^2 \\
    & = f(x^{k}) - \gamma \left\la \nabla f(x^{k}), \frac{g^{k}}{\|g^{k}\|} \right\ra + \frac{\gamma^2 L}{2}\\
    & = f(x^{k}) - \gamma \|g^{k}\| - \gamma \left\la \nabla f(x^{k}) - g^{k}, \frac{g^{k}}{\|g^{k}\|} \right\ra + \frac{\gamma^2 L}{2}\\
    & \leq f(x^{k}) - \gamma \|g^{k}\| +\gamma \|\nabla f(x^{k}) - g^{k}\| + \frac{\gamma^2 L}{2}\\
    & \leq f(x^{k}) - \gamma \|\nabla f(x^{k})\| +2\gamma \|\nabla f(x^{k}) - g^{k}\| + \frac{\gamma^2 L}{2}.
\end{align*}
Using notation for the functional gap and the error term, we get
\begin{equation*}
    \gamma \|\nabla f(x^{k})\| + \Delta_{k+1} \leq \Delta_k +2\gamma \|\hat{e}_k\| + \frac{\gamma^2 L}{2}.
\end{equation*}
Summing over $k$ from $0$ to $K-1$, we obtain 
\begin{equation}
    \gamma \sum^{K-1}_{k=0}\|\nabla f(x^{k})\| + \Delta_{K} \leq \Delta_0 +2\gamma \sum^{K-1}_{k=0}\|\hat{e}_k\| + \frac{\gamma^2 L K}{2}.
\end{equation}
\end{proof}

Next, we derive the lemma for error control for the Hessian corrected momentum estimator (lines $4-7$ in Algorithm~\ref{ch4:alg:NSGD_SOM_chapter_4}). Similar lemma for the special case $p=2$ appeared previously in \parencite{salehkaleybar-et-al22,Fatkhullin_SPGM_FND_2023}. For the case $p<2$, similar recursion was derived by \textcite{Hubler2024clip_to_norm} for first-order momentum. Now we extend this idea to the case of second-order momentum. 
\begin{lemma}\label{ch4:lem:dev_bound_chapter_4}
{
    Let Assumptions~\ref{ch4:assum:l_smooth_chapter_4},~\ref{ch4:ass:p_BCM_chapter_4}, and~\ref{ch4:ass:p_BCM_hess_chapter_4} hold. Then, for all $k \geq 0$,
\begin{equation}\label{ch4:eq:dev_bound_chapter_4}
    \Exp\left[\norm{\hat e_k}\right]
    \leq
    (1-\eta)^k\Exp\left[\norm{\hat e_0}\right]
    +2\eta^{(p-1)/p}\sigma
    +12\gamma(L+\sigma_h)\eta^{-1/p},
\end{equation}
    where $\hat e_k\eqdef g^k-\nabla f(x^k)$.
}
\end{lemma}
\begin{proof}
{
For $k\geq1$, let $\mathcal F_{k-1}$ contain all randomness used
to construct $x^k$ and $g^{k-1}$. Conditional on
$\mathcal F_{k-1}$, the variables ${q^k}$, $\xi^k$, and
$\hat\xi^k$ are fresh and mutually independent. Define
\[
    e_k\eqdef\nabla f(x^k;\xi^k)-\nabla f(x^k)
\]
and
\[
    \hat S_k\eqdef
    \nabla^2 f(\hat x^k;\hat\xi^k)(x^k-x^{k-1})
    -\bigl(\nabla f(x^k)-\nabla f(x^{k-1})\bigr).
\]
Using the uniform interpolation identity for $\hat x^k$, we have
\[
    \Exp[e_k\mid\mathcal F_{k-1}]=0,
    \qquad
    \Exp[\hat S_k\mid\mathcal F_{k-1}]=0,
    \qquad
    \Exp[\norm{e_k}^p\mid\mathcal F_{k-1}]\leq\sigma^p.
\]
Moreover, $\norm{x^k-x^{k-1}}\leq\gamma$, and therefore
\begin{align*}
    \Exp[\norm{\hat S_k}^p\mid\mathcal F_{k-1}]
    &\leq
    3\norm{\nabla f(x^k)-\nabla f(x^{k-1})}^p\\
    &\quad+
    3\Exp\left[
        \norm{\nabla^2f(\hat x^k)(x^k-x^{k-1})}^p
        \,\middle|\,\mathcal F_{k-1}
    \right]\\
    &\quad+
    3\Exp\left[
        \norm{
            (\nabla^2f(\hat x^k;\hat\xi^k)-\nabla^2f(\hat x^k))
            (x^k-x^{k-1})
        }^p
        \,\middle|\,\mathcal F_{k-1}
    \right]\\
    &\leq 6L^p\gamma^p+3\sigma_h^p\gamma^p.
\end{align*}

The estimator update gives, for $k\geq1$,
\[
    \hat e_k
    =(1-\eta)\hat e_{k-1}+\eta e_k+(1-\eta)\hat S_k.
\]
Unrolling this recursion yields
\begin{equation}\label{ch4:eq:hat_et_unrolled_chapter_4}
    \hat e_k
    =(1-\eta)^k\hat e_0
    +\eta\sum_{j=0}^{k-1}(1-\eta)^{k-j-1}e_{j+1}
    +\sum_{j=0}^{k-1}(1-\eta)^{k-j}\hat S_{j+1}.
\end{equation}
Consequently,
\begin{align}
\Exp[\norm{\hat e_k}]
&\leq
(1-\eta)^k\Exp[\norm{\hat e_0}]
+\Exp\left[
    \norm{\eta\sum_{j=0}^{k-1}(1-\eta)^{k-j-1}e_{j+1}}
\right]\notag\\
&\quad+
\Exp\left[
    \norm{\sum_{j=0}^{k-1}(1-\eta)^{k-j}\hat S_{j+1}}
\right].
\label{ch4:eq:dhsbhdsdsiwp_chapter_4}
\end{align}

The weighted sequences in both sums are martingale differences
with respect to the filtration above. Jensen's inequality and
Lemma~\ref{ch4:lem:von_Bahr_inequality_chapter_4} therefore give
\begin{align}
&\Exp\left[
    \norm{\eta\sum_{j=0}^{k-1}(1-\eta)^{k-j-1}e_{j+1}}
\right]\notag\\
&\quad\leq
\left(
    2\eta^p\sum_{j=0}^{k-1}
    (1-\eta)^{p(k-j-1)}\Exp[\norm{e_{j+1}}^p]
\right)^{1/p}\notag\\
&\quad\leq 2\eta^{(p-1)/p}\sigma,
\label{ch4:eq:bjhdvbfuioiapa_chapter_4}
\end{align}
where we used
$\sum_{j=0}^{k-1}(1-\eta)^{p(k-j-1)}\leq\eta^{-1}$.
Similarly—and with no extra factor of $\eta$ inside the norm—
\begin{align}
&\Exp\left[
    \norm{\sum_{j=0}^{k-1}(1-\eta)^{k-j}\hat S_{j+1}}
\right]\notag\\
&\quad\leq
\left(
    2\sum_{j=0}^{k-1}
    (1-\eta)^{p(k-j)}
    \Exp[\norm{\hat S_{j+1}}^p]
\right)^{1/p}\notag\\
&\quad\leq
12\gamma(L+\sigma_h)\eta^{-1/p}.
\label{ch4:eq:hdsbhfjbshbhf_chapter_4}
\end{align}
Substituting the last two estimates into
\eqref{ch4:eq:dhsbhdsdsiwp_chapter_4} proves
\eqref{ch4:eq:dev_bound_chapter_4}.
}
\end{proof}

\textbf{Main Theorem.} Now we are ready to {state and prove} the main convergence theorem for Algorithm~\ref{ch4:alg:NSGD_SOM_chapter_4}. 
\begin{theorem}
\label{ch4:thm:nsgdm_hess_app_chapter_4}
{
Let Assumptions~\ref{ch4:assum:lower_bounded_chapter_4},~\ref{ch4:assum:l_smooth_chapter_4},~\ref{ch4:ass:p_BCM_chapter_4}, and~\ref{ch4:ass:p_BCM_hess_chapter_4} hold, and define
\[
    \Delta_0\eqdef f(x^0)-f^{\inf},
    \qquad
    \mathcal E_0\geq
    \Exp\left[\norm{g^0-\nabla f(x^0)}\right].
\]
For $K\geq1$, choose
\[
    \eta=\min\{1,\eta_{\mathrm{eff}}\},
    \qquad
    \gamma=
    \sqrt{\frac{\Delta_0\eta^{1/p}}{(L+\sigma_h)K}},
\]
where
\[
    \eta_{\mathrm{eff}}
    \eqdef
    \max\left\{
        \left(\frac{\mathcal E_0}{K\sigma}\right)^{p/(2p-1)},
        \left(\frac{\Delta_0(L+\sigma_h)}{K\sigma^2}\right)^{p/(2p-1)}
    \right\}.
\]
Then the iterates of Algorithm~\ref{ch4:alg:NSGD_SOM_chapter_4}
satisfy
\begin{align*}
    \frac1K\sum_{k=0}^{K-1}
    \Exp\left[\norm{\nabla f(x^k)}\right]
    =\cO\Bigg(
        &\sqrt{\frac{\Delta_0(L+\sigma_h)}{K}}
        +\sigma
        \left(
            \frac{\Delta_0(L+\sigma_h)}{K\sigma^2}
        \right)^{(p-1)/(2p-1)}\\
        &+\frac{\mathcal E_0}{K}
        +\sigma
        \left(
            \frac{\mathcal E_0}{K\sigma}
        \right)^{(p-1)/(2p-1)}
    \Bigg).
\end{align*}
}
\end{theorem}
\begin{proof}
{
Applying Lemma~\ref{ch4:lem:dev_bound_chapter_4} and summing the
geometric term gives
\[
    \sum_{k=0}^{K-1}\Exp[\norm{\hat e_k}]
    \leq
    \frac{\Exp[\norm{\hat e_0}]}{\eta}
    +2K\eta^{(p-1)/p}\sigma
    +12K\gamma(L+\sigma_h)\eta^{-1/p}.
\]
Lemma~\ref{ch4:lem:descent_lemma_chapter_4} therefore implies
\begin{align*}
    \frac1K\sum_{k=0}^{K-1}
    \Exp[\norm{\nabla f(x^k)}]
    &\leq
    \frac{\Delta_0}{\gamma K}
    +\frac{2\Exp[\norm{\hat e_0}]}{\eta K}
    +4\eta^{(p-1)/p}\sigma\\
    &\quad+
    25\gamma(L+\sigma_h)\eta^{-1/p}.
\end{align*}
Substituting the stated value of $\gamma$ and using
$\Exp[\norm{\hat e_0}]\leq\mathcal E_0$ yields
\[
    \frac1K\sum_{k=0}^{K-1}
    \Exp[\norm{\nabla f(x^k)}]
    =
    \cO\left(
        \sqrt{\frac{\Delta_0(L+\sigma_h)}{\eta^{1/p}K}}
        +\frac{\mathcal E_0}{\eta K}
        +\eta^{(p-1)/p}\sigma
    \right).
\]
Finally, substituting
$\eta=\min\{1,\eta_{\mathrm{eff}}\}$ and considering separately
the cases $\eta_{\mathrm{eff}}\leq1$ and
$\eta_{\mathrm{eff}}>1$ gives the claimed bound.
}
\end{proof}

Now we investigate how different choices of initial estimator $g^{0}$ affect the total sample complexity bound.
\begin{corollary}
{
Let the assumptions and parameter choices of
Theorem~\ref{ch4:thm:nsgdm_hess_app_chapter_4} hold.
\begin{enumerate}
    \item If $g^0=\nabla f(x^0)$, then $\mathcal E_0=0$, and
    the sample complexity is
    \[
        \cO\left(
            \frac{\Delta_0(L+\sigma_h)}{\varepsilon^2}
            +
            \frac{\Delta_0(L+\sigma_h)}{\varepsilon^2}
            \left(\frac{\sigma}{\varepsilon}\right)^{1/(p-1)}
        \right).
    \]
    \item If $g^0=0$, then
    $\mathcal E_0\leq\sqrt{2L\Delta_0}$, and the sample complexity is
    \[
        \cO\left(
            \frac{\Delta_0(L+\sigma_h)}{\varepsilon^2}
            +
            \frac{\Delta_0(L+\sigma_h)}{\varepsilon^2}
            \left(\frac{\sigma}{\varepsilon}\right)^{1/(p-1)}
            +
            \frac{\sqrt{\Delta_0L}\,\sigma}{\varepsilon^2}
            \left(\frac{\sigma}{\varepsilon}\right)^{1/(p-1)}
        \right).
    \]
    \item If
    \[
        g^0=\frac1{B_{\mathrm{init}}}
        \sum_{j=1}^{B_{\mathrm{init}}}\nabla f(x^0;\xi_j^0),
        \qquad
        B_{\mathrm{init}}
        =
        \max\left\{
            1,
            \left(\frac{\sigma}{\varepsilon}\right)^{p/(p-1)},
            \left(\frac{\sigma}{\varepsilon}\right)^{p/(2p-1)}
        \right\},
    \]
    then
    $\mathcal E_0\leq
    2\sigma/B_{\mathrm{init}}^{(p-1)/p}$, and the total sample
    complexity, including the initial batch, is
    \[
        \cO\left(
            \frac{\Delta_0(L+\sigma_h)}{\varepsilon^2}
            +
            \frac{\Delta_0(L+\sigma_h)}{\varepsilon^2}
            \left(\frac{\sigma}{\varepsilon}\right)^{1/(p-1)}
            +
            \frac{\sigma}{\varepsilon}
            \left(\frac{\sigma}{\varepsilon}\right)^{1/(p-1)}
            +
            \left(\frac{\sigma}{\varepsilon}\right)^{p/(2p-1)}
        \right).
    \]
\end{enumerate}
}
\end{corollary}
\begin{proof}
{
The first two claims follow from
Theorem~\ref{ch4:thm:nsgdm_hess_app_chapter_4} after substituting
$\mathcal E_0=0$ and
$\mathcal E_0\leq\sqrt{2L\Delta_0}$, respectively.

For the mini-batch initialization, Lemma~\ref{ch4:lem:von_Bahr_inequality_chapter_4}
and the $p$-moment assumption give
\begin{align*}
    \Exp[\norm{\hat e_0}]
    &\leq
    \left(\Exp[\norm{\hat e_0}^p]\right)^{1/p}\\
    &\leq
    \frac{2}{B_{\mathrm{init}}}
    \left(
        \sum_{j=1}^{B_{\mathrm{init}}}
        \Exp\left[
            \norm{
                \nabla f(x^0;\xi_j^0)-\nabla f(x^0)
            }^p
        \right]
    \right)^{1/p}\\
    &\leq
    \frac{2\sigma}{B_{\mathrm{init}}^{(p-1)/p}}.
\end{align*}
Substitution into Theorem~\ref{ch4:thm:nsgdm_hess_app_chapter_4},
followed by solving each term for the requirement that it be at
most $\varepsilon$, gives the displayed complexities. The last
case also includes the $B_{\mathrm{init}}$ samples used to form
$g^0$.
}
\end{proof}

\newpage
\section{Additional Proofs for Section~\ref{ch4:sec:high_probability_analysis_chapter_4}}
\label{ch4:app:missing_proofs_hp_chapter_4}

Now we start the high-probability convergence analysis for Algorithm~\ref{ch4:alg:NSGD_SOM_clipped_chapter_4}. We divide our analysis into two parts: Optimization part, where we prove descent lemma, and High-Probability part, where we use concentration inequality to bound several terms from Descent Lemma. Finally, combining results from both parts, we prove the main results of this section via induction. The idea is based on work of \textcite{sadiev2023high} and \textcite{UnboundedClippedNSGDM2023Liu}.

\subsection{Analysis: Optimization Part}
 We start with \textit{Descent Lemma}. 

\begin{lemma}
\label{ch4:lem:descent_lemma_clipped_sgd_chapter_4}
    Let Assumptions~\ref{ch4:assum:lower_bounded_chapter_4},~\ref{ch4:assum:l_smooth_chapter_4} hold. Then Algorithm~\ref{ch4:alg:NSGD_SOM_clipped_chapter_4} with stepsize $\gamma > 0$ and momentum parameter $\eta \in (0,1)$  generates iterates $\left\{x^{k}\right\}_{k=0}^{K}$ satisfying  the following inequality
    {
    \begin{align*}
        \gamma \sum^{K-1}_{k=0}\|\nabla f(x^{k})\| + \Delta_{K}
        \leq{}& \Delta_0 + \frac{\gamma^2 L K}{2}
        + \frac{3\gamma}{\eta} \sqrt{L\Delta_0}\\
        &+2\gamma\eta \sum^{K-1}_{k=1}
        \left\|\sum^k_{j=1}(1-\eta)^{k-j}\theta_j\right\|\\
        &+2\gamma\sum^{K-1}_{k=1}
        \left\|\sum^k_{j=1}(1-\eta)^{k-j+1}\omega_j\right\|.
    \end{align*}
    }
    {where, for every $j\in[K-1]$,}
{
\begin{align}
    \theta_j
    &\eqdef \texttt{clip}\left(\nabla f(x^{j};\xi^j), \lambda\right)
    - \nabla f(x^{j}),
    \label{ch4:eq:definiton_theta_chapter_4}\\
    \omega_j
    &\eqdef \texttt{clip}\left(
        \nabla^2 f(\hat{x}^j;\hat{\xi}^{j})(x^{j} - x^{j-1}),
        \lambda_h
    \right)\notag\\
    &\qquad-\left(\nabla f(x^{j})-\nabla f(x^{j-1})\right).
    \label{ch4:eq:definition_omega_chapter_4}
\end{align}
}
\end{lemma}

 \begin{proof}
First, we notice that Lemma~\ref{ch4:lem:descent_lemma_chapter_4} holds for Algorithm~\ref{ch4:alg:NSGD_SOM_clipped_chapter_4} and can be proven in the same way as for Algorithm~\ref{ch4:alg:NSGD_SOM_chapter_4}, since the gradient updates in {both methods} are identical except for $x^{1}$, and the update rule for the {momentum estimate} does not play any role in the proof of Lemma~\ref{ch4:lem:descent_lemma_chapter_4}. Therefore, the iterates $\left\{x^{k}\right\}_{k=0}^{K}$ of Algorithm~\ref{ch4:alg:NSGD_SOM_clipped_chapter_4} satisfy
\begin{equation}
\label{ch4:eq:descent_lemma_2_chapter_4}
    \gamma \sum^{K-1}_{k=0}\|\nabla f(x^{k})\| + \Delta_{K} \leq \Delta_0 +2\gamma \sum^{K-1}_{k=0}\|\hat{e}_k\| + \frac{\gamma^2 L K}{2},
\end{equation}
where $x^{0}=x^{1}$, and hence $f(x^{1})-f^{\inf}=\Delta_0$. Next we bound $\|\hat{e}_k\|$ in almost the same way as we did in Lemma~\ref{ch4:lem:dev_bound_chapter_4}. By update rule for momentum parameter, we have 
\begin{align*}
    \hat e_k & = g^{k} - \nabla f(x^{k})\\
    & =  (1-\eta)\left( g^{k-1} + {\texttt{clip}\left(\nabla^2 f(\hat{x}^k;\hat{\xi}^{k})(x^{k} - x^{k-1}), \lambda_h \right)}\right)\\
    &\quad\quad  + \eta \texttt{clip}\left(\nabla f(x^{k};\xi^{k}), \lambda\right) - \nabla f(x^{k})\\
    & \overset{\eqref{ch4:eq:definiton_theta_chapter_4}, \eqref{ch4:eq:definition_omega_chapter_4}}{=} (1-\eta)\hat e_{k-1} + \eta \theta_k +(1-\eta)\omega_k\\
    & = (1-\eta)^k \hat e_0 + \eta\sum^{k}_{j=1}(1-\eta)^{k-j}\theta_j+ \sum^k_{j=1}(1-\eta)^{k-j+1}\omega_j.
\end{align*}
Taking the norm, we obtain the following bound 
{
\begin{align}
    \|\hat{e}_k\|
    \leq{}&(1-\eta)^k\|\hat{e}_0\|
    + \eta \left\|\sum^k_{j=1}(1-\eta)^{k-j}\theta_j\right\|\notag\\
    &+\left\|\sum^k_{j=1}(1-\eta)^{k-j+1}\omega_j\right\|\notag\\
    \leq{}&(1-\eta)^k \sqrt{2L\Delta_0}
    + \eta \left\|\sum^k_{j=1}(1-\eta)^{k-j}\theta_j\right\|\notag\\
    &+\left\|\sum^k_{j=1}(1-\eta)^{k-j+1}\omega_j\right\|.
    \label{ch4:eq:error_norm_bound_chapter_4}
\end{align}
}
where in the last inequality we used the following chain of inequalities 
{
\begin{align*}
    \|\hat{e}_0\|
    &=\|g^{0}-\nabla f(x^{0})\|
      =\|\nabla f(x^{0})\|
      =\|\nabla f(x^{1})\|\\
    &\leq \sqrt{2L\bigl(f(x^{1})-f^{\inf}\bigr)}
      \leq \sqrt{2L\Delta_0}.
\end{align*}
}
    Plugging \eqref{ch4:eq:error_norm_bound_chapter_4} into \eqref{ch4:eq:descent_lemma_2_chapter_4}, we acquire
    {
    \begin{align*}
        \MoveEqLeft[1] \gamma \sum_{k=0}^{K-1}\|\nabla f(x^{k})\|+\Delta_K
        \leq{} \\
        &\Delta_0
        +2\gamma\sum_{k=0}^{K-1}\|\hat e_k\|
        +\frac{\gamma^2LK}{2}\\
        \leq{}&\Delta_0+\frac{\gamma^2LK}{2}
        +2\gamma\sum_{k=0}^{K-1}(1-\eta)^k\sqrt{2L\Delta_0}\\
        &+2\gamma\eta\sum_{k=1}^{K-1}
        \left\|\sum_{j=1}^k(1-\eta)^{k-j}\theta_j\right\|\\
        &+2\gamma\sum_{k=1}^{K-1}
        \left\|\sum_{j=1}^k(1-\eta)^{k-j+1}\omega_j\right\|\\
        \leq{}&\Delta_0+\frac{\gamma^2LK}{2}
        +\frac{3\gamma}{\eta}\sqrt{L\Delta_0}\\
        &+2\gamma\eta\sum_{k=1}^{K-1}
        \left\|\sum_{j=1}^k(1-\eta)^{k-j}\theta_j\right\|\\
        &+2\gamma\sum_{k=1}^{K-1}
        \left\|\sum_{j=1}^k(1-\eta)^{k-j+1}\omega_j\right\|.
    \end{align*}
    }
where in the last inequality we used 
\begin{align*}
2\gamma \sum^{K-1}_{k=0}(1-\eta)^k \sqrt{2L\Delta_0} & \leq  2\gamma \sum^{\infty}_{k=0}(1-\eta)^k \sqrt{2L\Delta_0} \\
& \leq \frac{2\sqrt{2}\gamma}{1-(1-\eta)}\sqrt{L\Delta_0} \\
& \leq \frac{3\gamma}{\eta}\sqrt{L\Delta_0}.
    \end{align*}
\end{proof}

\subsection{Analysis: Statistical Part}

According to Lemma~\ref{ch4:lem:descent_lemma_clipped_sgd_chapter_4}, we have two new terms
\begin{equation*}
    \left\|\sum^k_{j=1}(1-\eta)^{k-j}\theta_j\right\| \quad \text{and}\quad  \left\|\sum^k_{j=1}(1-\eta)^{k-j+1}\omega_j\right\|.
\end{equation*}
{
To bound both terms, let $\mathcal F_{k-1}$ contain all randomness
used to construct $x^k$ and $g^{k-1}$. Conditional on
$\mathcal F_{k-1}$, the variables ${q^k}$, $\xi^k$, and
$\hat\xi^k$ are fresh and mutually independent. Following prior work
\parencite{gorbunov2020stochastic,sadiev2023high,UnboundedClippedNSGDM2023Liu},
decompose $\theta_k$ and $\omega_k$ into conditionally centered and
bias terms:
\begin{align}
    \theta_k & = \theta_k^u+\theta_k^b,\quad \text{where} \notag\\
    \theta^u_k & \eqdef \texttt{clip}\left(\nabla f(x^{k};\xi^{k}), \lambda\right) \notag \\
    &\quad - \Exp\left[\texttt{clip}\left(\nabla f(x^{k};\xi^{k}), \lambda\right)\mid\mathcal F_{k-1}\right], \label{ch4:eq:theta_unbiased_chapter_4}\\
    \theta^b_k & \eqdef \Exp\left[\texttt{clip}\left(\nabla f(x^{k};\xi^{k}), \lambda\right)\mid\mathcal F_{k-1}\right] - \nabla f(x^{k}); \label{ch4:eq:theta_biased_chapter_4}\\
    \omega_k & = \omega^u_k +\omega_k^b,\quad \text{where} \notag\\
    \omega^u_k & \eqdef \texttt{clip}\left(\nabla^2 f(\hat{x}^k;\hat{\xi}^{k})(x^{k} - x^{k-1}), \lambda_h \right) \notag\\
    &\quad\quad- \Exp\left[\texttt{clip}\left(\nabla^2 f(\hat{x}^k;\hat{\xi}^{k})(x^{k} - x^{k-1}), \lambda_h \right)\mid\mathcal F_{k-1}\right], \label{ch4:eq:omega_unbiased_chapter_4}\\
    \omega^b_k & \eqdef \Exp\left[\texttt{clip}\left(\nabla^2 f(\hat{x}^k;\hat{\xi}^{k})(x^{k} - x^{k-1}), \lambda_h \right)\mid\mathcal F_{k-1}\right] \notag\\
    &\quad\quad - \left(\nabla f(x^{k}) - \nabla f(x^{k-1})\right). \label{ch4:eq:omega_biased_chapter_4}
\end{align}
Under Assumption~\ref{ch4:ass:p_BCM_chapter_4}, if
$\lambda\geq2\norm{\nabla f(x^k)}$,
Lemma~\ref{ch4:lem:bound_variance_magnitude_chapter_4} implies
\begin{equation}
\label{ch4:eq:theta_u_theta_b_variance_and_magnitude_chapter_4}
    \Exp\left[\norm{\theta^u_k}^2\mid\mathcal F_{k-1}\right]
    \leq18\lambda^{2-p}\sigma^p,
    \qquad
    \norm{\theta^b_k}\leq2^p\lambda^{1-p}\sigma^p.
\end{equation}
Similarly, define
\[
    \hat S_k\eqdef
    \nabla^2f(\hat x^k;\hat\xi^k)(x^k-x^{k-1})
    -\bigl(\nabla f(x^k)-\nabla f(x^{k-1})\bigr).
\]
The uniform interpolation identity gives
{
\begin{align*}
    &\Exp\left[
    \nabla^2 f(\hat{x}^k;\hat{\xi}^{k})(x^{k}-x^{k-1})
    \,\middle|\,\mathcal F_{k-1}
    \right]\\
    &\qquad=
    \int_0^1\nabla^2 f\bigl(qx^k+(1-q)x^{k-1}\bigr)
    (x^k-x^{k-1})\,dq\\
    &\qquad=
    \nabla f(x^k)-\nabla f(x^{k-1}).
\end{align*}
}
Thus, $\Exp{\hat S_k\mid\mathcal F_{k-1}}=0$ and
\[
    \Exp{\norm{\hat S_k}^p\mid\mathcal F_{k-1}}
    \leq(6L^p+3\sigma_h^p)\gamma^p.
\]
If $\lambda_h\geq2L\gamma$, then
Lemma~\ref{ch4:lem:bound_variance_magnitude_chapter_4} yields
\begin{equation}
\label{ch4:eq:omega_u_omega_b_variance_and_magnitude_chapter_4}
    \begin{aligned}
    \Exp\left[\norm{\omega^u_k}^2\mid\mathcal F_{k-1}\right]
    &\leq18\lambda_h^{2-p}\gamma^p(6L^p+3\sigma_h^p),\\
    \norm{\omega^b_k}
    &\leq2^p\lambda_h^{1-p}\gamma^p(6L^p+3\sigma_h^p).
    \end{aligned}
\end{equation}
}

\vspace{0.5cm}

\begin{lemma}
\label{ch4:lem:technical_lemma_hp_1_chapter_4}
{
Let Assumption~\ref{ch4:ass:p_BCM_chapter_4} hold. For
$\beta'\in(0,\frac12]$, $k\in[K]$, and
$\lambda\geq\sigma/\eta^{1/p}$, define the localization event
\[
    \mathcal H_k
    \eqdef
    \left\{
    \max_{1\leq j\leq k}\norm{\nabla f(x^j)}
    \leq\frac{\lambda}{2}
    \right\}.
\]
Then the events defined in
\eqref{ch4:eq:def_event_one_chapter_4} and
\eqref{ch4:eq:def_event_two_chapter_4} satisfy
\[
    \Prob\left(
    \mathcal H_k\cap\overline{E_{\circledOne,k}}
    \right)\leq\beta',
    \qquad
    \Prob\left(
    \mathcal H_k\cap\overline{E_{\circledTwo,k}}
    \right)\leq\beta'.
\]
Consequently,
\[
    \Prob\left(
    \mathcal H_k
    \cap
    \left\{
    \left\|\sum_{j=1}^k(1-\eta)^{k-j}\theta_j\right\|
    >
    22\lambda\log\frac{2}{\beta'}
    \right\}
    \right)
    \leq2\beta'.
\]
}
\end{lemma}
\begin{proof} 
    We start with bounding $\left\|\sum^k_{j=1}(1-\eta)^{k-j}\theta_j\right\|$. 
    
    \paragraph{Upper bound  for $\left\|\sum^k_{j=1}(1-\eta)^{k-j}\theta_j\right\|$ } By \eqref{ch4:eq:theta_unbiased_chapter_4} and \eqref{ch4:eq:theta_biased_chapter_4}, we have 
    \begin{align*}
        \left\|\sum\limits^k_{j=1}(1-\eta)^{k-j}\theta_j\right\| & \leq \left\|\sum\limits^k_{j=1}(1-\eta)^{k-j}\theta_j^u\right\| + \underbrace{\left\|\sum\limits^k_{j=1}(1-\eta)^{k-j}\theta_j^b\right\|}_{\circledFour}.
    \end{align*}

    {For fixed $k$, define
    \[
        Y_j^k\eqdef
        \norm{(1-\eta)^{k-j}\theta_j^u}^2
        -
        \Exp\left[
            \norm{(1-\eta)^{k-j}\theta_j^u}^2
            \,\middle|\,\mathcal F_{j-1}
        \right],
        \qquad j\in[k].
    \]
    {For compactness, set
    \[
        S_{j-1}^k\eqdef
        \sum_{i=1}^{j-1}(1-\eta)^{k-i}\theta_i^u.
    \]
    Define the scalar sequence $V_j^k$ recursively by}}
{
\begin{equation*}
    V_j^k \eqdef
    \begin{cases}
        0,
        & j=0,\\[1mm]
        \displaystyle
        \sign\left(\sum_{i=1}^{j-1}V_i^k\right)
        \frac{\left\langle S_{j-1}^k,
        (1-\eta)^{k-j}\theta_j^u\right\rangle}{\norm{S_{j-1}^k}},
        & j\neq0\ \text{and}\ S_{j-1}^k\neq0,\\[3mm]
        0,
        & j\neq0\ \text{and}\ S_{j-1}^k=0.
    \end{cases}
\end{equation*}
}
{
For every $j\in[k]$, define the prefix event
\[
    \mathcal H_j
    \eqdef
    \left\{
    \max_{1\leq i\leq j}\norm{\nabla f(x^i)}
    \leq\frac{\lambda}{2}
    \right\},
\]
and let
\[
    I_j\eqdef\mathbf{1}_{\mathcal H_j},
    \qquad
    \overline V_j^k\eqdef I_jV_j^k,
    \qquad
    \overline Y_j^k\eqdef I_jY_j^k.
\]
Since $x^j$ is $\mathcal F_{j-1}$-measurable,
$I_j$ is also $\mathcal F_{j-1}$-measurable. Hence
$(\overline V_j^k)_{j=1}^k$ and
$(\overline Y_j^k)_{j=1}^k$ are martingale-difference
sequences. Moreover, on $\mathcal H_k$ the events
$\mathcal H_j$ hold for all $j\leq k$, and therefore
\[
    \overline V_j^k=V_j^k,
    \qquad
    \overline Y_j^k=Y_j^k,
    \qquad j\in[k].
\]
}
    Then to bound $\left\|\sum\limits^k_{j=1}(1-\eta)^{k-j}\theta_j^u\right\|$ we use Lemma~\ref{ch4:lem:inequlity_from_cutkosky_chapter_4} and obtain 
    {
    \begin{align*}
        \left\|\sum_{j=1}^k(1-\eta)^{k-j}\theta_j^u\right\|
        \leq{}&
        \left|\sum_{j=1}^kV_j^k\right|\\
        &+\left(
        \max_{j\in[k]}\left\|(1-\eta)^{k-j}\theta_j^u\right\|^2
        +\sum_{j=1}^k\left\|(1-\eta)^{k-j}\theta_j^u\right\|^2
        \right)^{1/2}\\
        \leq{}&
        \left|\sum_{j=1}^kV_j^k\right|
        +\left(
        2\sum_{j=1}^k\left\|(1-\eta)^{k-j}\theta_j^u\right\|^2
        \right)^{1/2}\\
        \leq{}&
        \left|\underbrace{\sum_{j=1}^kV_j^k}_{\circledOne}\right|\\
        &+\left(
        2\underbrace{\sum_{j=1}^kY_j^k}_{\circledTwo}
        +2\underbrace{
        \sum_{j=1}^k
        \Exp\left[
        \left\|(1-\eta)^{k-j}\theta_j^u\right\|^2
        \,\middle|\,\mathcal F_{j-1}
        \right]}_{\circledThree}
        \right)^{1/2}.
    \end{align*}
    }

\textbf{Upper bound for \circledOne.}
{The sequence $V_1^k,\dots,V_k^k$ is a martingale
difference sequence with respect to
$(\mathcal F_j)_{j\geq0}$ because the coefficient multiplying
$\theta_j^u$ is $\mathcal F_{j-1}$-measurable and
$\Exp[\theta_j^u\mid\mathcal F_{j-1}]=0$. Moreover,
Lemma~\ref{ch4:lem:inequlity_from_cutkosky_chapter_4} gives}
\begin{equation}
    {|V_j^k| \leq \norm{(1-\eta)^{k-j}\theta_j^u} \leq \norm{\theta_j^u} \leq 2\lambda, \quad \text{where } c_1 \eqdef 2\lambda.}
\end{equation}
{
Let
\[
    \sigma_j^2
    \eqdef
    \Exp\left[(V_j^k)^2\mid\mathcal F_{j-1}\right],
    \qquad
    \overline\sigma_j^2
    \eqdef
    \Exp\left[(\overline V_j^k)^2\mid\mathcal F_{j-1}\right]
    =I_j\sigma_j^2.
\]
Since $I_j=1$ only on $\mathcal H_j$, the clipping bound
\eqref{ch4:eq:theta_u_theta_b_variance_and_magnitude_chapter_4}
gives
\begin{align}
    \sum_{j=1}^k\overline\sigma_j^2
    &\leq
    \sum_{j=1}^k
    I_j(1-\eta)^{2(k-j)}
    \Exp\left[
    \norm{\theta_j^u}^2
    \,\middle|\,\mathcal F_{j-1}
    \right]\notag\\
    &\leq
    18\lambda^{2-p}\sigma^p
    \sum_{j=1}^k(1-\eta)^{2(k-j)}
    \leq
    \frac{18\lambda^{2-p}\sigma^p}{\eta}
    \eqdef G_1.
    \label{ch4:eq:qwertyuaisja_chapter_4}
\end{align}
Moreover, $|\overline V_j^k|\leq c_1\eqdef2\lambda$. Thus,
Lemma~\ref{ch4:lem:Bernstain_inequality_chapter_4} with
\[
    b_1
    \eqdef
    \left(
    \frac{c_1}{3}
    +
    \sqrt{\frac{c_1^2}{9}+2G_1}
    \right)
    \log\frac{2}{\beta'}
\]
implies
\begin{equation}
\label{ch4:eq:localized_event_one_chapter_4}
    \Prob\left\{
    \left|\sum_{j=1}^k\overline V_j^k\right|>b_1
    \ \text{and}\
    \sum_{j=1}^k\overline\sigma_j^2
    \leq G_1\log\frac{2}{\beta'}
    \right\}
    \leq\beta'.
\end{equation}
}

\textbf{Upper bound for \circledTwo.}
{The sequence $Y_1^k,\dots,Y_k^k$ is a martingale
difference sequence with respect to $(\mathcal F_j)_{j\geq0}$.
Moreover, Lemma~\ref{ch4:lem:bound_variance_magnitude_chapter_4}
gives}
\begin{align}
|Y_j^k| & \leq \norm{(1-\eta)^{k-j}\theta_j^u}^2 + {\Exp_{\xi_j}\left[\norm{(1-\eta)^{k-j}\theta_j^u}^2\right]} \notag \\
& \leq 4\lambda^2 + 4 \lambda^2 = 8\lambda^2. \label{ch4:eq:bound_on_Y_chapter_4}
    \end{align}
{
Set $c_2\eqdef8\lambda^2$ and define
\[
    \widetilde\sigma_j^2
    \eqdef
    \Exp\left[(Y_j^k)^2\mid\mathcal F_{j-1}\right],
    \qquad
    \overline{\widetilde\sigma}_j^2
    \eqdef
    \Exp\left[(\overline Y_j^k)^2\mid\mathcal F_{j-1}\right]
    =I_j\widetilde\sigma_j^2.
\]
On $\mathcal H_j$, \eqref{ch4:eq:bound_on_Y_chapter_4} and
\eqref{ch4:eq:theta_u_theta_b_variance_and_magnitude_chapter_4}
give
\[
    \widetilde\sigma_j^2
    \leq
    16\lambda^2(1-\eta)^{2(k-j)}
    \Exp\left[
    \norm{\theta_j^u}^2
    \,\middle|\,\mathcal F_{j-1}
    \right].
\]
It follows that
\begin{align*}
    \sum_{j=1}^k\overline{\widetilde\sigma}_j^2
    &\leq
    16\lambda^2
    \sum_{j=1}^k
    I_j(1-\eta)^{2(k-j)}
    \Exp\left[
    \norm{\theta_j^u}^2
    \,\middle|\,\mathcal F_{j-1}
    \right]\\
    &\leq
    \frac{16\cdot18\lambda^{4-p}\sigma^p}{\eta}
    \eqdef G_2.
\end{align*}
Since $|\overline Y_j^k|\leq c_2$, Lemma~
\ref{ch4:lem:Bernstain_inequality_chapter_4} with
\[
    b_2
    \eqdef
    \left(
    \frac{c_2}{3}
    +
    \sqrt{\frac{c_2^2}{9}+2G_2}
    \right)
    \log\frac{2}{\beta'}
\]
implies
\begin{equation}
\label{ch4:eq:localized_event_two_chapter_4}
    \Prob\left\{
    \left|\sum_{j=1}^k\overline Y_j^k\right|>b_2
    \ \text{and}\
    \sum_{j=1}^k\overline{\widetilde\sigma}_j^2
    \leq G_2\log\frac{2}{\beta'}
    \right\}
    \leq\beta'.
\end{equation}
}

\textbf{Upper bound for \circledThree.} {On $\mathcal H_k$, the stopped and original conditional variances agree. Hence, by \eqref{ch4:eq:qwertyuaisja_chapter_4},}
{
\begin{align*}
    \circledThree
    &=\sum_{j=1}^k
    \Exp\left[
    \left\|(1-\eta)^{k-j}\theta_j^u\right\|^2
    \,\middle|\,\mathcal F_{j-1}
    \right]\\
    &\overset{\eqref{ch4:eq:qwertyuaisja_chapter_4}}{\leq}
    \frac{18\lambda^{2-p}\sigma^p}{\eta}.
\end{align*}
}

\textbf{Upper bound on \circledFour.} {On $\mathcal H_k$, the clipping bias bound holds for every $j\in[k]$, and therefore}
{
\begin{align*}
    \circledFour
    &=\left\|\sum_{j=1}^k(1-\eta)^{k-j}\theta_j^b\right\|\\
    &\leq\sum_{j=1}^k(1-\eta)^{k-j}\norm{\theta_j^b}\\
    &\overset{\eqref{ch4:eq:theta_u_theta_b_variance_and_magnitude_chapter_4}}{\leq}
    4\lambda^{1-p}\sigma^p
    \sum_{j=1}^k(1-\eta)^{k-j}
    \leq\frac{4\lambda^{1-p}\sigma^p}{\eta}.
\end{align*}
}

To sum up, we introduce event $E_{\circledOne, k}$ as follows 
\begin{equation}
\label{ch4:eq:def_event_one_chapter_4}
    E_{\circledOne, k} \eqdef \left\{\left|\circledOne\right| \leq b_1 \text{ or } \sum^k_{j=1}\sigma^2_j > G_1\log\frac{2}{\beta'}\right\},
\end{equation}
where $c_1 = 2\lambda$, $G_1 = \frac{18\lambda^{2-p}\sigma^p}{\eta}$, $b_1 = \left(\frac{1}{3}c_1 + \sqrt{\frac{1}{9}c_1^2 + 2G_1}\right)\log\frac{2}{\beta'}$.
{On $\mathcal H_k$, we have
$\overline V_j^k=V_j^k$ and
$\overline\sigma_j^2=\sigma_j^2$ for every $j\leq k$.
Thus, \eqref{ch4:eq:localized_event_one_chapter_4} yields
\[
    \Prob\left(
    \mathcal H_k\cap\overline{E_{\circledOne,k}}
    \right)\leq\beta'.
\]}
The bound on $b_1$ 
\begin{align*}
    b_1 & = \left(\frac{1}{3}c_1 + \sqrt{\frac{1}{9}c_1^2 + 2G_1}\right)\log\frac{2}{\beta'} \leq \left(\frac{2}{3}c_1 + \sqrt{2G_1}\right)\log\frac{2}{\beta'}\\
    & = \left(\frac{4}{3}\lambda + \sqrt{36\frac{\lambda^{2-p}\sigma^p}{\eta}}\right)\log\frac{2}{\beta'} =  \lambda\left(\frac{4}{3}+6\sqrt{\frac{1}{\eta}\left(\frac{\sigma}{\lambda}\right)^p}\right)\log\frac{2}{\beta'}.
\end{align*}

Also we define event $E_{\circledTwo, k}$ as follows 
\begin{equation}
\label{ch4:eq:def_event_two_chapter_4}
    E_{\circledTwo, k} = \left\{\left|\circledTwo\right| \leq b_2 \text{ or } \sum^k_{j=1}\widetilde{\sigma}^2_j > G_2\log\frac{2}{\beta'}\right\},
\end{equation}
where $c_2 = 8\lambda^2$, $G_2 = \frac{16\cdot 18 \lambda^{4-p}\sigma^p}{\eta}$, $b_2 = \left(\frac{1}{3}c_2 + \sqrt{\frac{1}{9}c^2_2 + 2G_2}\right)\log\frac{2}{\beta'}$.
{On $\mathcal H_k$, we also have
$\overline Y_j^k=Y_j^k$ and
$\overline{\widetilde\sigma}_j^2=\widetilde\sigma_j^2$.
Therefore, \eqref{ch4:eq:localized_event_two_chapter_4} yields
\[
    \Prob\left(
    \mathcal H_k\cap\overline{E_{\circledTwo,k}}
    \right)\leq\beta'.
\]}
We adjust the bound on $b_2$:
\begin{align*}
    b_2 & = \left(\frac{1}{3}c_2 + \sqrt{\frac{1}{9}c_2^2 + 2G_2}\right)\log\frac{2}{\beta'} \leq  \left(\frac{2}{3}c_2 + \sqrt{2G_2}\right)\log\frac{2}{\beta'}\\
    & = \left(\frac{16}{3}\lambda^2 + \sqrt{\frac{16 \cdot 36\lambda^{4-p}\sigma^p}{\eta}}\right)\log\frac{2}{\beta'} = \lambda^2\left(\frac{16}{3}+24\sqrt{\frac{1}{\eta}\left(\frac{\sigma}{\lambda}\right)^p}\right)\log\frac{2}{\beta'}.
\end{align*}

Thus, {on $\mathcal H_k$,} the event $E_{\circledOne, k} \cap E_{\circledTwo, k}$ implies 
\begin{align*}
    \MoveEqLeft[1] \left\|\sum\limits^k_{j=1}(1-\eta)^{k-j}\theta_j\right\| \\
    & \leq \left\|\sum\limits^k_{j=1}(1-\eta)^{k-j}\theta_j^u\right\| + \left\|\sum\limits^k_{j=1}(1-\eta)^{k-j}\theta_j^b\right\|\\
    & \leq \left|\circledOne\right| + \sqrt{2\cdot \circledTwo + 2\cdot \circledThree} + \circledFour\\
    & \leq \lambda\left(\frac{4}{3}+6\sqrt{\frac{1}{\eta}\left(\frac{\sigma}{\lambda}\right)^p}\right)\log\frac{2}{\beta'} + 4\lambda \cdot \frac{1}{\eta}\left(\frac{\sigma}{\lambda}\right)^p\\
    &\quad + \sqrt{2\lambda^2\left(\frac{16}{3}+24\sqrt{\frac{1}{\eta}\left(\frac{\sigma}{\lambda}\right)^p}\right)\log\frac{2}{\beta'} + \frac{36 \lambda^{2-p}\sigma^p}{\eta}}.
\end{align*}
Taking $\lambda \geq \frac{\sigma}{\eta^{1/p}}$, we have {on $\mathcal H_k$} the event $E_{\circledOne, k} \cap E_{\circledTwo, k}$ implies 
\begin{align*}
    \MoveEqLeft[1] \left\|\sum\limits^k_{j=1}(1-\eta)^{k-j}\theta_j\right\| \\
    & \leq \lambda\left(\frac{4}{3}+ 6 + 4 + \sqrt{
    \frac{32}{3} + 84}\right)\log\frac{2}{\beta'} \leq 22\lambda\log\frac{2}{\beta'}.
\end{align*}
{
Combining the two localized failure bounds with a union bound gives
\[
    \Prob\left(
    \mathcal H_k
    \cap
    \left\{
    \left\|\sum_{j=1}^k(1-\eta)^{k-j}\theta_j\right\|
    >
    22\lambda\log\frac{2}{\beta'}
    \right\}
    \right)
    \leq2\beta'.
\]
}
This concludes the proof.
\end{proof}

\begin{lemma}
\label{ch4:lem:technical_lemma_hp_2_chapter_4}
    {Let Assumptions~\ref{ch4:assum:l_smooth_chapter_4}
    and~\ref{ch4:ass:p_BCM_hess_chapter_4} hold. For any
    $\beta''\in(0,\frac12]$ and $k\in[K]$, suppose
    \[
        \lambda_h\geq
        \max\left\{
            2\gamma L,
            \frac{\gamma(L+\sigma_h)}{\eta^{1/p}}
        \right\}.
    \]
    Then, with probability at least $1-2\beta''$,}
    \begin{equation*}
        \left\|\sum^k_{j=1}(1-\eta)^{k-j+1}\omega_j\right\| \leq  59\lambda_h\log\frac{2}{\beta''}.
\end{equation*}
\end{lemma}

\begin{proof}
 By \eqref{ch4:eq:omega_unbiased_chapter_4} and \eqref{ch4:eq:omega_biased_chapter_4}, we have 
    \begin{align}
        \left\|\sum\limits^k_{j=1}(1-\eta)^{k-j+1}\omega_j\right\| & \leq \left\|\sum\limits^k_{j=1}(1-\eta)^{k-j+1}\omega_j^u\right\|\notag\\
        &\quad + \underbrace{\left\|\sum\limits^k_{j=1}(1-\eta)^{k-j+1}\omega_j^b\right\|}_{\circledEight}. \label{ch4:eq:eight_def_chapter_4}
    \end{align}
    {For fixed $k$, define
    \[
        Z_j^k\eqdef
        \norm{(1-\eta)^{k-j+1}\omega_j^u}^2
        -
        \Exp\left[
            \norm{(1-\eta)^{k-j+1}\omega_j^u}^2
            \,\middle|\,\mathcal F_{j-1}
        \right],
        \qquad j\in[k].
    \]
    Also define $W_j^k$ recursively by}
\begin{equation*}
    W^k_j \eqdef \begin{cases}
        0,& \text{if } j = 0;\\
        {\sign}\left(\sum\limits^{j-1}_{i=1} W^k_i\right) A_j, & \text{if } j\neq 0 \text{ and } \sum\limits^{j-1}_{i=1}(1-\eta)^{k-i+1}\omega_i^u \neq 0;\\
        0, &  \text{if } j\neq 0 \text{ and } \sum\limits^{j-1}_{i=1}(1-\eta)^{k-i+1}\omega_i^u = 0.
    \end{cases}
\end{equation*}
where $A_j$ is defined as 
\begin{equation*}
    A_j \eqdef \frac{\left\la \sum\limits^{j-1}_{i=1}(1-\eta)^{k-i+1}\omega_i^u, (1-\eta)^{k-j+1}\omega_j^u \right\ra}{\left\|\sum\limits^{j-1}_{i=1}(1-\eta)^{k-i+1}\omega_i^u\right\|}.
\end{equation*}
    Then to bound $\left\|\sum\limits^k_{j=1}(1-\eta)^{k-j+1}\omega_j^u\right\|$ we use Lemma~\ref{ch4:lem:inequlity_from_cutkosky_chapter_4} and obtain 
    {
    \begin{align*}
        \left\|\sum_{j=1}^k(1-\eta)^{k-j+1}\omega_j^u\right\|
        \leq{}&
        \left|\sum_{j=1}^k W_j^k\right|\\
        &+\biggl(
        \max_{j\in[k]}\left\|(1-\eta)^{k-j+1}\omega_j^u\right\|^2\\
        &\qquad+\sum_{j=1}^k\left\|(1-\eta)^{k-j+1}\omega_j^u\right\|^2
        \biggr)^{1/2}\\
        \leq{}&
        \left|\sum_{j=1}^k W_j^k\right|
        +\left(
        2\sum_{j=1}^k\left\|(1-\eta)^{k-j+1}\omega_j^u\right\|^2
        \right)^{1/2}\\
        \leq{}&
        \left|\underbrace{\sum_{j=1}^k W_j^k}_{\circledFive}\right|\\
        &+\left(
        2\underbrace{\sum_{j=1}^k Z_j^k}_{\circledSix}
        +2\underbrace{
        \sum_{j=1}^k
        \Exp\left[
        \left\|(1-\eta)^{k-j+1}\omega_j^u\right\|^2
        \,\middle|\,\mathcal F_{j-1}
        \right]}_{\circledSeven}
        \right)^{1/2}.
    \end{align*}
    }

\textbf{Upper bound for \circledFive.}
{The sequence $W_1^k,\dots,W_k^k$ is a martingale
difference sequence with respect to $(\mathcal F_j)_{j\geq0}$.
Moreover, Lemma~\ref{ch4:lem:inequlity_from_cutkosky_chapter_4}
gives}
\begin{equation}
\label{ch4:eq:mdskcmsdl_chapter_4}
    {|W_j^k| \leq \norm{(1-\eta)^{k-j+1}\omega_j^u} \leq \norm{\omega_j^u} \leq 2\lambda_h, \quad \text{where } c_3 \eqdef 2\lambda_h.}
\end{equation}
{Let
$\bar\sigma_j^2\eqdef\Exp[(W_j^k)^2\mid\mathcal F_{j-1}]$.
Lemma~\ref{ch4:lem:Bernstain_inequality_chapter_4} implies}
\begin{equation*}
    {\Prob\left\{\left|\circledFive\right| > b_3 \text{ and } \sum^k_{j=1}\bar\sigma^2_j \leq G_3 \log \frac{2}{\beta''}\right\} \leq 2\exp\left(- \frac{b^2_3}{2G_3\log\frac{2}{\beta''} + \frac{2c_3b_3}{3} }\right) = \beta'',}
\end{equation*}
{where the last identity holds for
$b_3 = \left(\frac{1}{3}c_3+\sqrt{\frac{1}{9}c_3^2 + 2G_3}\right)\log \frac{2}{\beta''}$.
To define $G_3$, we bound $\sum_{j=1}^k\bar\sigma_j^2$:}
\begin{align}
    {\sum^k_{j=1}\bar{\sigma}^2_j} & = {\sum^k_{j=1} \Exp\left[(W^k_j)^2\mid\mathcal F_{j-1}\right]
    \overset{\eqref{ch4:eq:mdskcmsdl_chapter_4}}{\leq}  \sum^k_{j=1} (1-\eta)^{2(k-j+1)}\Exp\left[\|\omega_j^u\|^2\mid\mathcal F_{j-1}\right]} \notag \\
    & \overset{\eqref{ch4:eq:omega_u_omega_b_variance_and_magnitude_chapter_4}}{\leq}  \sum^k_{j=1} (1-\eta)^{2(k-j+1)} \cdot 18\lambda^{2-p}_h\left(6L^p+3\sigma^p_h\right)\gamma^p \notag \\
    & \leq \frac{18 \lambda^{2-p}_h(6L^p+3\sigma^p_h)\gamma^p}{1-(1-\eta)^2} \notag\\
    & = \frac{18 \lambda^{2-p}_h\left(6L^p+3\sigma^p_h\right)\gamma^p}{\eta(2-\eta)} \notag \\
    & \leq \frac{18 \lambda^{2-p}_h\left(6L^p+3\sigma^p_h\right)\gamma^p}{\eta} = G_3. \label{ch4:eq:mjkcniwecuidhhjs_chapter_4}
\end{align}

\textbf{Upper bound for \circledSix.}
{The sequence $Z_1^k,\dots,Z_k^k$ is a martingale
difference sequence with respect to $(\mathcal F_j)_{j\geq0}$.
Moreover, Lemma~\ref{ch4:lem:bound_variance_magnitude_chapter_4}
gives}
\begin{align}
   |Z_j^k| & \leq \norm{(1-\eta)^{k-j+1}\omega_j^u}^2 + \Exp\left[\norm{(1-\eta)^{k-j+1}\omega_j^u}^2\mid\mathcal F_{j-1}\right] \notag\\
   & \leq 4\lambda^2_h + 4 \lambda^2_h = 8\lambda^2_h,\label{ch4:eq:bound_on_Z_chapter_4}
\end{align}
{where $c_4 \eqdef 8\lambda_h^2$. Define
$\widehat{\sigma}_j^2\eqdef
\Exp[(Z_j^k)^2\mid\mathcal F_{j-1}]$. Then}
\begin{align}
    {\widehat{\sigma}^2_j}  & \overset{\eqref{ch4:eq:bound_on_Z_chapter_4}}{\leq} {8\lambda^2_h \cdot \Exp\left[\left|\left\|(1-\eta)^{k-j+1}\omega_j^u\right\|^2 - \Exp\left[\left\|(1-\eta)^{k-j+1}\omega_j^u\right\|^2\mid\mathcal F_{j-1}\right]\right|\middle|\mathcal F_{j-1}\right]} \notag\\
    & \leq {16\lambda^2_h \Exp\left[\left\|(1-\eta)^{k-j+1}\omega_j^u\right\|^2\mid\mathcal F_{j-1}\right].} \label{ch4:eq:lsmclksnciooeiw_chapter_4}
\end{align}
Then, Lemma~\ref{ch4:lem:Bernstain_inequality_chapter_4} implies 
\begin{equation*}
    {\Prob\left\{\left|\circledSix\right| > b_4 \text{ and } \sum^k_{j=1} \widehat{\sigma}^2_j \leq G_4 \log \frac{2}{\beta''}\right\} \leq 2\exp\left(-\frac{b^2_4}{2G_4\log\frac{2}{\beta''} +\frac{2c_4b_4}{3}}\right) = \beta'',}
\end{equation*}
{where the last identity holds for
$b_4 = \left(\frac{1}{3}c_4+\sqrt{\frac{1}{9}c_4^2 + 2G_4}\right)\log \frac{2}{\beta''}$.
To define $G_4$, we bound $\sum_{j=1}^k\widehat{\sigma}_j^2$:}
\begin{align*}
    {\sum^k_{j=1} \widehat{\sigma}^2_j} & \overset{\eqref{ch4:eq:lsmclksnciooeiw_chapter_4}}{\leq} {16\lambda^2_h \sum^k_{j=1}(1-\eta)^{2(k-j+1)}\Exp\left[\left\|\omega_j^u\right\|^2\mid\mathcal F_{j-1}\right]}\\
    & \leq 16\lambda^2_h\cdot \frac{18\lambda_h^{2-p}\left(6L^p+3\sigma^p_h\right)\gamma^p}{\eta} \\
    & = \frac{16\cdot 18 \lambda^{4-p}_h(6L^p+3\sigma^p_h)\gamma^p}{\eta} = G_4.
\end{align*}

\textbf{Upper bound for \circledSeven.}
{By \eqref{ch4:eq:mjkcniwecuidhhjs_chapter_4}, with probability one,}
\begin{equation*}
    {\circledSeven
    = \sum^k_{j=1} \Exp\left[\left\|(1-\eta)^{k-j+1}\omega_j^u\right\|^2\mid\mathcal F_{j-1}\right]
    \leq \frac{18 \lambda_h^{2-p}(6L^p+3\sigma_h^p)\gamma^p}{\eta}.}
\end{equation*}

\textbf{Upper bound on \circledEight.}  By definition of $\omega^b_j$, with probability $1$ we have 
\begin{align*}
    \circledEight & \overset{\eqref{ch4:eq:eight_def_chapter_4}}{=} \left\|\sum\limits^k_{j=1} (1-\eta)^{k-j+1}\omega^b_j\right\| \leq \sum^k_{j=1}(1-\eta)^{k-j+1}\|\omega^b_j\|\\
    & \overset{\eqref{ch4:eq:omega_u_omega_b_variance_and_magnitude_chapter_4}}{\leq} 4\lambda_h^{1-p}\left(6L^p +3\sigma^p_h\right)\gamma^p\sum^k_{j=1}(1-\eta)^{k-j+1} \\
    & \leq \frac{4\lambda_h^{1-p}\left(6L^p+3\sigma^p_h\right)\gamma^p}{\eta}.
\end{align*}

To sum up, we introduce event $E_{\circledFive, k}$ as follows 
\begin{equation}
\label{ch4:eq:def_event_5_chapter_4}
    {E_{\circledFive, k} = \left\{\left|\circledFive\right| \leq b_3 \text{ or } \sum^k_{j=1}\bar{\sigma}^2_j > G_3\log\frac{2}{\beta''}\right\},}
\end{equation}
where $c_3 = 2\lambda_h$, $G_3= \frac{18\lambda^{2-p}_h\left(6L^p+3\sigma^p_h\right)\gamma^p}{\eta}$, $b_3 = \left(\frac{1}{3}c_3 + \sqrt{\frac{1}{9}c_3^2 + 2G_3}\right)\log\frac{2}{\beta''}$.
We have shown that $\Prob\left\{E_{\circledFive, k}\right\} \geq 1-\beta''$. We adjust the bound on $b_3$:
\begin{align*}
    {b_3} & = {\left(\frac{1}{3}c_3 + \sqrt{\frac{1}{9}c_3^2 + 2G_3}\right)\log\frac{2}{\beta''} \leq \left(\frac{2}{3}c_3 + \sqrt{2G_3}\right)\log\frac{2}{\beta''}}\\
    & = \left(\frac{4}{3}\lambda_h + \sqrt{36\frac{\lambda^{2-p}_h\left(6L^p + 3\sigma^p_h\right)\gamma^p}{\eta}}\right)\log\frac{2}{\beta''}\\
    & =  {\lambda_h\left(\frac{4}{3}+6\sqrt{\frac{(6L^p+3\sigma^p_h)\gamma^p}{\eta\lambda_h^p}}\right)\log\frac{2}{\beta''}.}
\end{align*}

Also we introduce event $E_{\circledSix, k}$ as follows 
\begin{equation}
\label{ch4:eq:def_event_6_chapter_4}
    E_{\circledSix, k} = \left\{\left|\circledSix\right| \leq b_4 \text{ or } \sum^k_{j=1}\widehat{\sigma}^2_j > G_4\log\frac{2}{\beta''}\right\},
\end{equation}
where $c_4 = 8\lambda^2_h$, $G_4 = \frac{16\cdot 18 \lambda^{4-p}_h\left(6L^p+3\sigma^p_h\right)\gamma^p}{\eta}$, $b_4 = \left(\frac{1}{3}c_4 + \sqrt{\frac{1}{9}c^2_4 + 2G_4}\right)\log\frac{2}{\beta''}$. We have shown that $\Prob\left\{E_{\circledSix, k}\right\} \geq 1-\beta''$. We adjust the bound on $b_4$:
\begin{align*}
    {b_4} & = {\left(\frac{1}{3}c_4 + \sqrt{\frac{1}{9}c_4^2 + 2G_4}\right)\log\frac{2}{\beta''} \leq  \left(\frac{2}{3}c_4 + \sqrt{2G_4}\right)\log\frac{2}{\beta''}}\\
    & = \left(\frac{16}{3}\lambda^2_h + \sqrt{\frac{16 \cdot 36\lambda_h^{4-p}\left(6L^p+3\sigma^p_h\right)\gamma^p}{\eta}}\right)\log\frac{2}{\beta''} \\
    & = \lambda^2_h\left(\frac{16}{3}+24\sqrt{\frac{\left(6L^p+3\sigma^p_h\right)\gamma^p}{\eta\lambda^p_h}}\right)\log\frac{2}{\beta''}.
\end{align*}

Thus, the event $E_{\circledFive, k} \cap E_{\circledSix, k}$ implies 
\begin{align*}
    \MoveEqLeft[1] \left\|\sum\limits^k_{j=1}(1-\eta)^{k-j+1}\omega_j\right\| \\
    & \leq \left\|\sum\limits^k_{j=1}(1-\eta)^{k-j+1}\omega_j^u\right\| + \left\|\sum\limits^k_{j=1}(1-\eta)^{k-j+1}\omega_j^b\right\|\\
    & \leq \left|\circledFive\right| + \sqrt{2\cdot \circledSix + 2\cdot \circledSeven} + \circledEight\\
    & \leq \lambda_h\left(\frac{4}{3}+6\sqrt{\frac{\left(6L^p+3\sigma^p_h\right)\gamma^p}{\eta\lambda^p_h}}\right)\log\frac{2}{\beta''} \\
    &\quad+ 4\lambda_h \cdot \frac{\left(6L^p+3\sigma^p_h\right)\gamma^p}{\eta\lambda^p_h}\\
    &\quad + \Bigg(2\lambda^2_h\left(\frac{16}{3}+24\sqrt{\frac{\left(6L^p+3\sigma^p_h\right)\gamma^p}{\eta\lambda^p_h}}\right)\log\frac{2}{\beta''} \\
    &\quad \qquad\qquad+36 \lambda^{2}_h\frac{(6L^p+3\sigma^p_h)\gamma^p}{\eta\lambda^p_h}\Bigg)^{1/2}.
\end{align*}
Assuming that $\lambda_h \geq \frac{\gamma(L+\sigma_h)}{\eta^{1/p}}$, we have the event $E_{\circledFive, k} \cap E_{\circledSix, k}$ implies 
\begin{align*}
    \MoveEqLeft[1] \left\|\sum\limits^k_{j=1}(1-\eta)^{k-j+1}\omega_j\right\| \\
    & \leq \lambda_h\left(\frac{4}{3}+ 6\sqrt{6} + 24 + \sqrt{
    \frac{32}{3} + 48\sqrt{6}+216}\right)\log\frac{2}{\beta''}\\
    & \leq 59\lambda_h\log\frac{2}{\beta''}.
\end{align*}
By the union bound, $\Prob\left\{E_{\circledFive, k} \cap E_{\circledSix, k}\right\} \geq 1-2\beta''$, which proves the lemma.

\end{proof}

\newpage
\paragraph{Main Theorem} Now we are ready to state and prove the main results of this section: high-probability convergence guarantees for Algorithm~\ref{ch4:alg:NSGD_SOM_clipped_chapter_4}. 

\begin{theorem}
\label{ch4:thm:main_theorem_hp_app_chapter_4}
Let Assumptions~\ref{ch4:assum:lower_bounded_chapter_4},~\ref{ch4:assum:l_smooth_chapter_4},~\ref{ch4:ass:p_BCM_chapter_4} and ~\ref{ch4:ass:p_BCM_hess_chapter_4} hold, and let $\Delta_0 {\eqdef} f(x^{0}) - f^{\inf}$. If the parameters of Algorithm~\ref{ch4:alg:NSGD_SOM_clipped_chapter_4} are chosen with stepsize
\begin{equation}
\label{ch4:eq:stepsize_hp_part_chapter_4}
{
\begin{aligned}
\gamma {=} \min\Bigg\{&
\frac{1}{2}\sqrt{\frac{\Delta_0}{LK}},\,
\frac{\eta}{12}\sqrt{\frac{\Delta_0}{L}},\,
\frac{1}{1408\eta K\log\frac{8K}{\beta}}\sqrt{\frac{\Delta_0}{L}},\\
&
\frac{\Delta_0}{352\sigma\eta^{(p-1)/p}K\log\frac{8K}{\beta}},\,
\sqrt{\frac{\Delta_0\eta^{1/p}}{968(L+\sigma_h)K\log\frac{8K}{\beta}}}
\Bigg\},
\end{aligned}
}
\end{equation}
{and with momentum parameter and algorithmic clipping inputs}
\begin{equation}
\label{ch4:eq:momentum_clipping_level_chapter_4}
{
\begin{aligned}
    \eta
    &= \frac{1}{K^{p/(2p-1)}},\\
    \lambda
    &= \max\left\{4\sqrt{L\Delta_0},
       \frac{\sigma}{\eta^{1/p}}\right\},\\
    \bar\lambda_h
    &= \frac{2(L+\sigma_h)}{\eta^{1/p}},
    \qquad
    \lambda_h \eqdef \gamma\bar\lambda_h
     = \frac{2\gamma(L+\sigma_h)}{\eta^{1/p}},
\end{aligned}
}
\end{equation}
for some $K \geq 1$ and $\beta \in (0,1]$ such that $\log\frac{8K}{\beta} \geq 1$. {Here, \(\bar\lambda_h\) is the clipping level supplied to Algorithm~\ref{ch4:alg:NSGD_SOM_clipped_chapter_4}, whereas \(\lambda_h=\gamma\bar\lambda_h\) is the equivalent clipping level for the unscaled Hessian-vector product used in the analysis.} Then, after $K$ iterations of Algorithm~\ref{ch4:alg:NSGD_SOM_clipped_chapter_4}, the iterates {satisfy, with probability at least $1-\beta$,}
\begin{equation}
    \frac{1}{K}\sum^{K-1}_{k=0}\norm{\nabla f(x^{k})} \leq \frac{2\Delta_0}{\gamma K}.
\end{equation}
{Consequently, with probability at least $1-\beta$,}
\begin{equation}
    {\frac{1}{K}\sum^{K-1}_{k=0}\norm{\nabla f(x^{k})} = \cO\left(\frac{\sqrt{\Delta_0(L+\sigma_h)}+\sigma}{K^{(p-1)/(2p-1)}}\log \frac{K}{\beta}\right),}
\end{equation}
implying that to achieve $\frac{1}{K}\sum^{K-1}_{k=0}\norm{\nabla f(x^{k})} \leq \varepsilon$ with probability at least $1-\beta$ Algorithm~\ref{ch4:alg:NSGD_SOM_clipped_chapter_4} requires 
\begin{equation}
{
\begin{aligned}
K = \cO\Bigg(&
\left(\frac{\sqrt{\Delta_0(L+\sigma_h)}+\sigma}{\varepsilon}\right)^{(2p-1)/(p-1)}\\
&{}\times
\log^{(2p-1)/(p-1)}\left(
\frac{1}{\beta}
\left(\frac{\sqrt{\Delta_0(L+\sigma_h)}+\sigma}{\varepsilon}\right)^{(2p-1)/(p-1)}
\right)
\Bigg).
\end{aligned}
}
\end{equation}
{iterations, and therefore the same order of stochastic-oracle calls.}
\end{theorem}
\begin{proof} {Recall that the proofs of Lemmas~\ref{ch4:lem:technical_lemma_hp_1_chapter_4} and~\ref{ch4:lem:technical_lemma_hp_2_chapter_4} establish}
\begin{align*}
    E_{\circledOne, k} & \overset{\eqref{ch4:eq:def_event_one_chapter_4}}{=} \left\{\left|\circledOne\right| \leq b_1 \text{ or } \sum^k_{j=1}\sigma^2_j > G_1\log\frac{2}{\beta'}\right\},\\
    E_{\circledTwo, k} & \overset{\eqref{ch4:eq:def_event_two_chapter_4}}{=} \left\{\left|\circledTwo\right| \leq b_2 \text{ or } \sum^k_{j=1}\widetilde{\sigma}^2_j > G_2\log\frac{2}{\beta'}\right\},\\
    E_{\circledFive, k} & \overset{\eqref{ch4:eq:def_event_5_chapter_4}}{=} \left\{\left|\circledFive\right| \leq b_3 \text{ or } \sum^k_{j=1}\bar{\sigma}^2_j > G_3\log\frac{2}{\beta''}\right\},\\
     E_{\circledSix, k} & \overset{\eqref{ch4:eq:def_event_6_chapter_4}}{=} \left\{\left|\circledSix\right| \leq b_4 \text{ or } \sum^k_{j=1}\widehat{\sigma}^2_j > G_4\log\frac{2}{\beta''}\right\},
\end{align*}
where
{
\begin{align*}
c_1&=2\lambda,
&G_1&=\frac{18\lambda^{2-p}\sigma^p}{\eta},\\
b_1&\leq
\lambda\left(\frac{4}{3}
+6\sqrt{\frac{1}{\eta}\left(\frac{\sigma}{\lambda}\right)^p}\right)
\log\frac{2}{\beta'},\\[1mm]
c_2&=8\lambda^2,
&G_2&=\frac{16\cdot18\lambda^{4-p}\sigma^p}{\eta},\\
b_2&\leq
\lambda^2\left(\frac{16}{3}
+24\sqrt{\frac{1}{\eta}\left(\frac{\sigma}{\lambda}\right)^p}\right)
\log\frac{2}{\beta'},\\[1mm]
c_3&=2\lambda_h,
&G_3&=\frac{18\lambda_h^{2-p}(6L^p+3\sigma_h^p)\gamma^p}{\eta},\\
b_3&\leq
\lambda_h\left(\frac{4}{3}
+6\sqrt{\frac{(6L^p+3\sigma_h^p)\gamma^p}{\eta\lambda_h^p}}\right)
\log\frac{2}{\beta''},\\[1mm]
c_4&=8\lambda_h^2,
&G_4&=\frac{16\cdot18\lambda_h^{4-p}(6L^p+3\sigma_h^p)\gamma^p}{\eta},\\
b_4&\leq
\lambda_h^2\left(\frac{16}{3}
+24\sqrt{\frac{(6L^p+3\sigma_h^p)\gamma^p}{\eta\lambda_h^p}}\right)
\log\frac{2}{\beta''}.
\end{align*}
}
Moreover, according to Lemma~\ref{ch4:lem:technical_lemma_hp_1_chapter_4} and Lemma~\ref{ch4:lem:technical_lemma_hp_2_chapter_4}, we have  
{
\begin{align*}
\Prob\left(
\mathcal H_k\cap\overline{E_{\circledOne,k}}
\right)&\leq\beta',
&
\Prob\left(
\mathcal H_k\cap\overline{E_{\circledTwo,k}}
\right)&\leq\beta',\\
\Prob\left\{E_{\circledFive,k}\right\}&\geq1-\beta'',
&\Prob\left\{E_{\circledSix,k}\right\}&\geq1-\beta''.
\end{align*}
}

The idea of the proof is based on {a technique from prior works} \parencite{gorbunov2020stochastic, sadiev2023high, UnboundedClippedNSGDM2023Liu}: via mathematical induction we plan to show that for any $\tau \in \{0,1,\dots,K-1\}$ the event 
    \begin{equation*}
        G_{\tau} = E_{\tau}\cap E_{1,\tau} \cap E_{2,\tau} \cap E_{3,\tau} \cap E_{4,\tau}
    \end{equation*}
    holds with probability at least $1-\frac{\tau\beta}{K}$, where event $E_{\tau}$ is defined as 
    \begin{equation*}
        {E_{\tau} = \left\{\gamma\sum^k_{s=0}\|\nabla f(x^s)\| + \Delta_{k+1} \leq 2\Delta_0,\quad \forall~k\leq \tau\right\},}
    \end{equation*}
    and 
    \begin{align*}
E_{1,\tau} & = \bigcap^{\tau}_{k=1} E_{\circledOne, k},\quad E_{2,\tau} = \bigcap^{\tau}_{k=1} E_{\circledTwo, k},\quad E_{3,\tau} \\
& = \bigcap^{\tau}_{k=1} E_{\circledFive, k},\quad E_{4,\tau} = \bigcap^{\tau}_{k=1} E_{\circledSix, k}.
    \end{align*}

\textbf{Basis of induction.}
{Since $x^1=x^0$, we have $\Delta_1=\Delta_0$.
Furthermore, $L$-smoothness and lower boundedness give
$\norm{\nabla f(x^0)}\leq\sqrt{2L\Delta_0}$. The first
stepsize restriction in \eqref{ch4:eq:stepsize_hp_part_chapter_4}
therefore implies
\[
    \gamma\norm{\nabla f(x^0)}+\Delta_1
    \leq
    \frac{\sqrt2}{2\sqrt K}\Delta_0+\Delta_0
    \leq2\Delta_0.
\]
Thus $G_0=E_0$ holds with probability one.}

\textbf{Step of induction.} Assume that the induction hypothesis is true for $\tau-1$:  $\Prob\left\{G_{\tau-1}\right\} \geq 1 - ((\tau-1)\beta)/K$. One needs to prove $\Prob\left\{G_{\tau}\right\} \geq 1 - \tau\beta/K$.
{
The event $E_{\tau-1}$ implies, for every $j\leq\tau$,
\[
    \Delta_j
    \leq
    \gamma\sum_{s=0}^{j-1}\norm{\nabla f(x^s)}
    +\Delta_j
    \leq2\Delta_0.
\]
Therefore, by $L$-smoothness and the choice
$\lambda\geq4\sqrt{L\Delta_0}$,
\[
    \norm{\nabla f(x^j)}
    \leq\sqrt{2L\Delta_j}
    \leq2\sqrt{L\Delta_0}
    \leq\frac{\lambda}{2},
    \qquad j\leq\tau.
\]
Consequently,
\[
    G_{\tau-1}\subseteq E_{\tau-1}\subseteq\mathcal H_\tau.
\]
}
We want to mention that 
{
\begin{align*}
    &E_{\tau-1}\cap E_{1,\tau}\cap E_{2,\tau}
    \cap E_{3,\tau}\cap E_{4,\tau}\\
    &\qquad=
    G_{\tau-1}\cap E_{\circledOne,\tau}
    \cap E_{\circledTwo,\tau}
    \cap E_{\circledFive,\tau}
    \cap E_{\circledSix,\tau}.
\end{align*}
}
Then, we have 
{
\begin{align*}
    &\Prob\left\{
    E_{\tau-1}\cap E_{1,\tau}\cap E_{2,\tau}
    \cap E_{3,\tau}\cap E_{4,\tau}
    \right\}\\
    &\quad=
    \Prob\left\{
    G_{\tau-1}\cap E_{\circledOne,\tau}
    \cap E_{\circledTwo,\tau}
    \cap E_{\circledFive,\tau}
    \cap E_{\circledSix,\tau}
    \right\}\\
    &\quad\geq
    \Prob\left\{G_{\tau-1}\right\}
    -\Prob\left\{
    G_{\tau-1}\cap\overline{E_{\circledOne,\tau}}
    \right\}\\
    &\qquad
    -\Prob\left\{
    G_{\tau-1}\cap\overline{E_{\circledTwo,\tau}}
    \right\}
    -\Prob\left\{
    G_{\tau-1}\cap\overline{E_{\circledFive,\tau}}
    \right\}\\
    &\qquad
    -\Prob\left\{
    G_{\tau-1}\cap\overline{E_{\circledSix,\tau}}
    \right\}\\
    &\quad\geq
    1-\frac{(\tau-1)\beta}{K}-2\beta'-2\beta''\\
    &\quad=
    1-\frac{\tau\beta}{K}
    +\left(\frac{\beta}{K}-2\beta'-2\beta''\right)
    =1-\frac{\tau\beta}{K}.
\end{align*}
}
where we have selected $\beta' = \beta'' = \frac{\beta}{4K}$. {The second inequality uses
$G_{\tau-1}\subseteq\mathcal H_\tau$, the two localized bounds
from Lemma~\ref{ch4:lem:technical_lemma_hp_1_chapter_4}, and the
unconditional bounds from
Lemma~\ref{ch4:lem:technical_lemma_hp_2_chapter_4}.}
Assuming $E_{\tau-1}\cap E_{1,\tau} \cap E_{2,\tau} \cap E_{3,\tau} \cap E_{4,\tau}$ happens, Lemma~\ref{ch4:lem:descent_lemma_clipped_sgd_chapter_4} implies
\begin{align*}
    \gamma\sum^{\tau}_{k=0}\|\nabla f(x^{k})\| + \Delta_{\tau+1} & \leq \Delta_0 + \frac{3\gamma}{\eta}\sqrt{L\Delta_0} + \frac{\gamma^2 L (\tau+1)}{2}\\
    &\quad + 2\gamma\eta\sum^{\tau}_{k=1}\left\|\sum^k_{j=1}(1-\eta)^{k-j} \theta_j\right\| \\
    &\quad+ 2\gamma \sum^{\tau}_{k=1}\left\|\sum^k_{j=1}(1-\eta)^{k-j+1}\omega_j\right\|\\
    & \leq \Delta_0 + \frac{3\gamma}{\eta}\sqrt{L\Delta_0} + \frac{\gamma^2 L (\tau+1)}{2}\\
    &\quad + 44\gamma\eta(\tau +1)\lambda\log\frac{8K}{\beta} + 118\gamma(\tau+1) \lambda_h\log\frac{8K}{\beta}.
\end{align*}

By setting clipping levels \eqref{ch4:eq:momentum_clipping_level_chapter_4} and stepsize \eqref{ch4:eq:stepsize_hp_part_chapter_4}, we have 
\begin{equation*}
    \gamma\sum^{\tau}_{k=0}\|\nabla f(x^{k})\| + \Delta_{\tau+1} \leq \Delta_0 +\frac{\Delta_0}{4} +\frac{\Delta_0}{4} +\frac{\Delta_0}{4} +\frac{\Delta_0}{4} = 2\Delta_0.
\end{equation*}
Therefore, we obtain the following. The event
$E_{\tau-1}\cap E_{1,\tau}\cap E_{2,\tau}\cap E_{3,\tau}\cap E_{4,\tau}$
implies $E_\tau$, and hence it is a subset of $G_\tau$.
\begin{align*}
    \Prob\left\{G_{\tau}\right\} & = \Prob\left\{E_{\tau}\cap E_{1,\tau}\cap E_{2,\tau} \cap E_{3,\tau} \cap E_{4, \tau}\right\}\\
    & \geq \Prob\left\{E_{\tau-1}\cap E_{1,\tau}\cap E_{2,\tau} \cap E_{3,\tau} \cap E_{4, \tau}\right\} \\
    & \geq 1-\frac{\tau\beta}{K}.
\end{align*}
Thus, we have 
\begin{equation*}
    {\Prob\left\{E_{K-1}\right\}
    \geq \Prob\left\{G_{K-1}\right\}
    \geq 1-\frac{(K-1)\beta}{K}
    \geq1-\beta,}
\end{equation*}
or equivalently with probability at least $1-\beta$ 
\begin{align*}
\MoveEqLeft[1] \gamma\sum^{K-1}_{k=0}\|\nabla f(x^{k})\| + \Delta_{K} \\
& \leq 2\Delta_0~~\Rightarrow~~ \frac{1}{K}\sum^{K-1}_{k=0}\|\nabla f(x^{k})\| \\
& \leq \frac{2\Delta_0}{\gamma K}.
    \end{align*}
{Substituting the chosen stepsize from
\eqref{ch4:eq:stepsize_hp_part_chapter_4} into the preceding
inequality gives}
\begin{align*}
    \frac{1}{K}\sum^{K-1}_{k=0}\|\nabla f(x^{k})\| & \leq \frac{2\Delta_0}{\gamma K}\\
    & = \mathcal{O}\Bigg(\max\Bigg\{\sqrt{\frac{L\Delta_0}{K}}, \frac{\sqrt{L\Delta_0}}{\eta K}, \eta \sqrt{L\Delta_0}\log \frac{K}{\beta}, \\
    &\quad\qquad\qquad\quad\sigma\eta^{(p-1)/p}\log \frac{K}{\beta}, \sqrt{\frac{\Delta_0(L+\sigma_h)\log\frac{K}{\beta}}{K\eta^{1/p}}}\Bigg\}\Bigg).
\end{align*}

{Selecting
$\eta = K^{-p/(2p-1)}$
as in \eqref{ch4:eq:momentum_clipping_level_chapter_4}, we obtain}
\begin{align*}
    \MoveEqLeft[1] {\frac{1}{K}\sum^{K-1}_{k=0}\|\nabla f(x^{k})\|} \\
    & =\mathcal{O}\Bigg(\max\Bigg\{\sqrt{\frac{L\Delta_0}{K}}, \frac{\sqrt{L\Delta_0}}{K^{(p-1)/(2p-1)}},  \frac{\sqrt{L\Delta_0}}{K^{(p-1)/(2p-1)}}\log \frac{K}{\beta},\\
    &\quad\qquad\qquad\quad \frac{\sigma}{K^{(p-1)/(2p-1)}}\log\frac{K}{\beta}, \sqrt{\frac{\Delta_0(L+\sigma_h)\log\frac{K}{\beta}}{K^{(2p-2)/(2p-1)}}}\Bigg\}\Bigg)\\
    & = {\mathcal{O}\Bigg(\frac{\sqrt{\Delta_0(L+\sigma_h)}+\sigma}{K^{(p-1)/(2p-1)}}\log \frac{K}{\beta} \Bigg).}
\end{align*}

\textbf{Discussion.} According to Theorem~\ref{ch4:thm:main_theorem_hp_chapter_4}, the sample complexity of Algorithm~\ref{ch4:alg:NSGD_SOM_clipped_chapter_4} is 
\begin{equation*}
    {\widetilde\cO\left(\left(\frac{\sqrt{\Delta_0(L+\sigma_h)}+\sigma}{\varepsilon}\right)^{(2p-1)/(p-1)}\right),}
\end{equation*}
which has the same dependence on $\varepsilon$ as lower bounds do (see Theorem~\ref{ch4:thm:lower_bounds_chapter_4}). However, the obtained upper bound has {worse} dependence on parameters $L,\sigma,\sigma_h, \Delta_0$ than lower bounds have. {We believe that this is caused} by the analysis technique, since Algorithm~\ref{ch4:alg:NSGD_SOM_chapter_4} has better convergence rate in terms of parameter dependence. This observation raises an interesting question: Is it possible to conduct high-probability analysis for Algorithm~\ref{ch4:alg:NSGD_SOM_chapter_4}, or equivalently, Algorithm~\ref{ch4:alg:NSGD_SOM_clipped_chapter_4} without clipping?

\end{proof}

\newpage
\section{Experimental Results}

We solve a simple quadratic problem, $\frac{1}{2}\|x\|^2$, in dimension $d=10$, with synthetic noise sampled from a two-sided Pareto distribution. Recall that {a} two-sided Pareto distribution with tail index $\bar p$ satisfies Assumption~\ref{ch4:ass:p_BCM_chapter_4} with $p < \bar p$. We consider different values of the tail index $\bar p$ ranging from $1.1$ to $2.0$. We start each algorithm from the same point sampled uniformly at random from the standard normal distribution. Our empirical analysis consists of three experimental settings:

\begin{figure}[ht]
    \begin{minipage}[t]{0.48\linewidth}
        \centering
        \includegraphics[width=\linewidth]{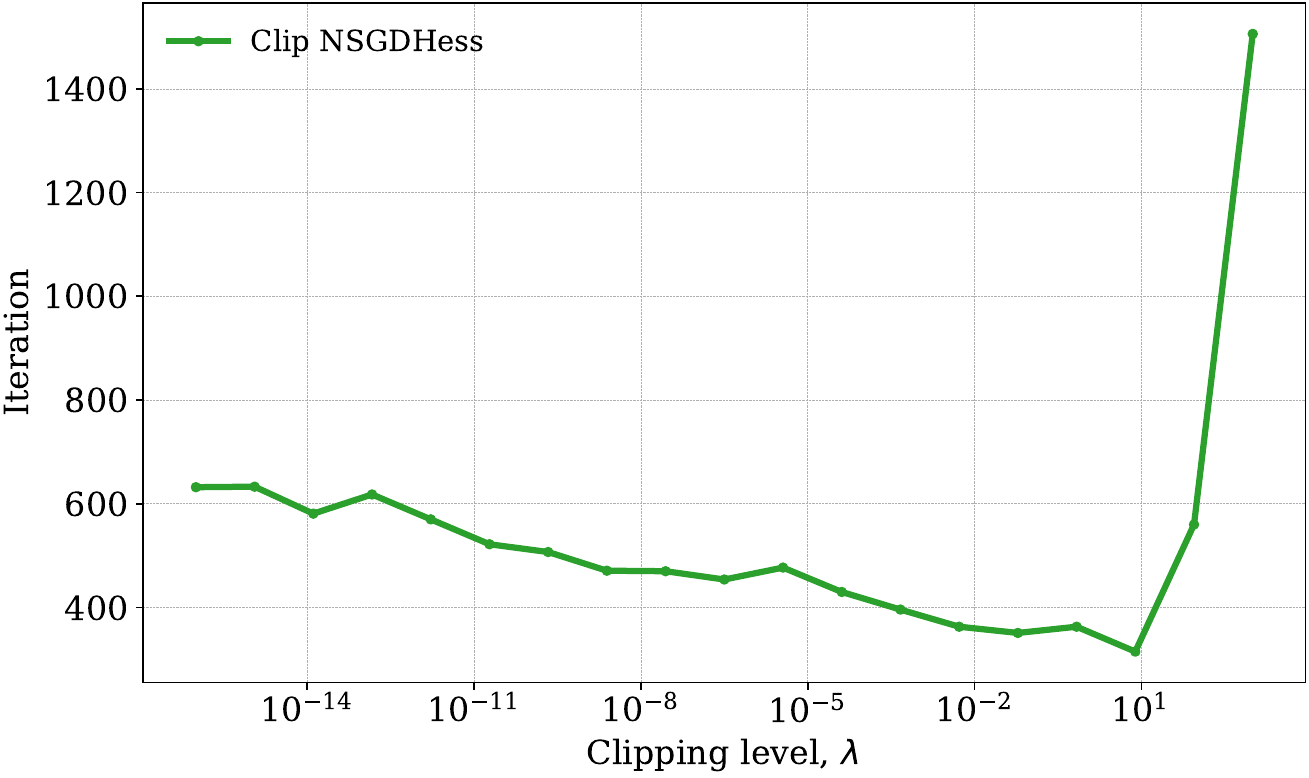}
        \caption{Effect of the {algorithmic HVP clipping level $\bar\lambda_h = \lambda$} on the iteration complexity. The plot shows the number of iterations required for {\algname{Clip NSGDHess}} to find a point with $\norm{\nabla f(x)} \leq 3/2$. For extremely small and large values of $\lambda,$ more iterations are needed. The recommended value for this task is ${\bar\lambda_h = 10}$.
        }
        \label{ch4:fig:sensitivity_clipping_chapter_4}
    \end{minipage}%
    \hfill
    \begin{minipage}[t]{0.48\linewidth}
        \centering
        \includegraphics[width=\linewidth]{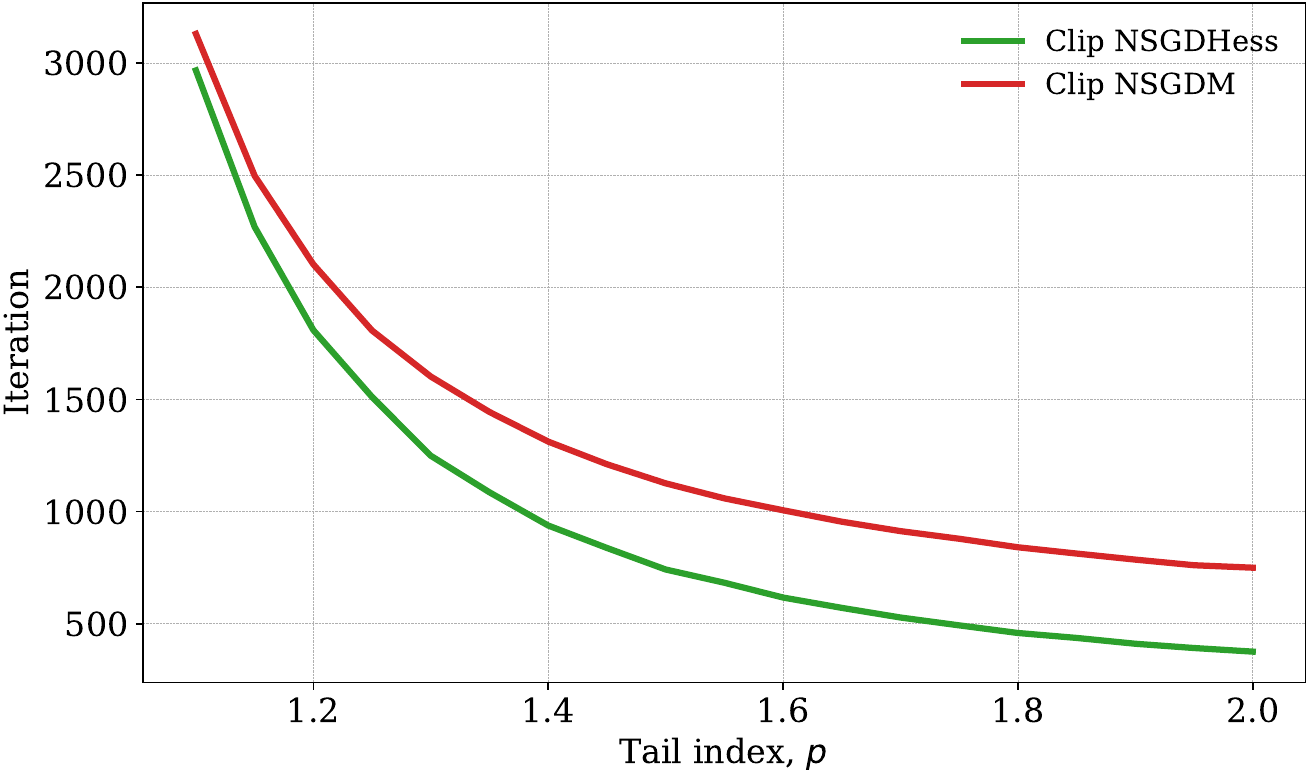}
        \caption{Number of iterations needed for \algname{Clip NSGDHess} and \algname{Clip NSGDM} {for varying tail indices} to find a point with $\norm{\nabla f(x)} \leq 3/2$ starting with the same initial point. The performance of both algorithms {degrades gradually as the tail index decreases}. The iteration complexity of the second-order algorithm, \algname{Clip NSGDHess}, is uniformly better for all values of $p \in [1.1, 2].$}
        \label{ch4:fig:complexity_comparison_chapter_4}
    \end{minipage}
\end{figure}

\paragraph{1. Effect of Clipping Level on Iteration Complexity.} 
We run Algorithm~\ref{ch4:alg:NSGD_SOM_clipped_chapter_4} (\algname{Clip NSGDHess}) with ${\bar\lambda_h = \lambda}$ until it {reaches} the target accuracy $\|\nabla f(x)\| \leq 3/2$. The stepsize and momentum parameters are set to $\gamma = 0.01$ and ${\eta = 0.2}$. {The equivalent clipping level for the unscaled Hessian-vector product is \(\lambda_h=\gamma\bar\lambda_h\).} Clipping levels are varied from $10^{-16}$ to $10^{3}$ in multiplicative steps of $10$. In Figure~\ref{ch4:fig:sensitivity_clipping_chapter_4}, we observe that intermediate clipping ($\lambda=10$) yields the lowest iteration complexity; see Table~\ref{ch4:tab:iter_comple_lambda_h_chapter_4} for precise numerical values. Extremely small or large values lead to slower convergence. The proposed algorithm tolerates very small values of $\lambda$, but very large values (e.g., $\lambda > 10^3$) result in very slow convergence. This empirically confirms the need for Hessian clipping to stabilize and accelerate convergence.

\begin{table}[h]
\centering
\caption{Iteration complexity depending on the clipping level ${\lambda = \bar\lambda_h}$ for Algorithm~\ref{ch4:alg:NSGD_SOM_clipped_chapter_4} (\algname{Clip NSGDHess}), cf. Figure~\ref{ch4:fig:sensitivity_clipping_chapter_4}.
}
\label{ch4:tab:iter_comple_lambda_h_chapter_4}

\begin{tabular}{c|ccccc}
 Clipping Level, $\lambda$ & $10^{-16}$ & $10^{-8}$ & $10^{1}$ & $10^{2}$ & $10^{3}$ \\
\hline
Iterations & 632& 471& 315& 560& 1506 \\
\end{tabular}
\vspace{0.2cm}

\end{table}

\paragraph{2. Comparison with \algname{Clip NSGDM} {for Varying Tail Indices}.}
We compare \algname{Clip NSGDHess} with \algname{Clip NSGDM}, which was proposed by \textcite{cutkosky2021high}. For \algname{Clip NSGDM}, we used the theoretical parameter choices: stepsize ${\gamma = K^{-(2p-1)/(3p-2)}}$ and momentum ${\eta = K^{-p/(3p-2)}}$. For both methods we fixed ${K=4000}$, $\lambda=0.5$, and ${\bar\lambda_h}=0.05$. These experiments were conducted to complement Figure~\ref{ch4:fig:complexities_chapter_4}, which presents the theoretical iteration complexities. As shown in Figure~\ref{ch4:fig:complexity_comparison_chapter_4}, the empirical results align with the theory: \algname{Clip NSGDHess} consistently outperforms \algname{Clip NSGDM}, and its relative advantage (the ratio of the iteration counts) grows as $p$ increases. This provides strong empirical support for our theoretical findings.

\paragraph{3. Sensitivity to Gradient Clipping for Fixed Hessian Clipping.}
In this set of experiments, we test \algname{Clip NSGDHess} (Algorithm~\ref{ch4:alg:NSGD_SOM_clipped_chapter_4}). Our goal is to study how iteration complexity depends on the choice of the gradient clipping level $\lambda$ when the Hessian clipping level is fixed. Specifically, we test three different values of the Hessian clipping level, and vary the gradient clipping threshold $\lambda$:
\begin{itemize}
\item For ${\bar\lambda_h = 0.01}$: $\lambda$ ranged from $10^{-4}$ to $10^{2}$ (step $10$).
\item For ${\bar\lambda_h = 1.0}$: $\lambda$ ranged from $10^{-2}$ to $10^{4}$ (step $10$).
\item For ${\bar\lambda_h = 100}$: $\lambda$ ranged from $10^{0}$ to $10^{6}$ (step $10$).
\end{itemize}
The results are summarized in Figure~\ref{ch4:fig:dmdskkcs_chapter_4}.

\begin{figure}[ht]
    \begin{minipage}[t]{0.48\linewidth}
        \centering
        \includegraphics[width=\linewidth]{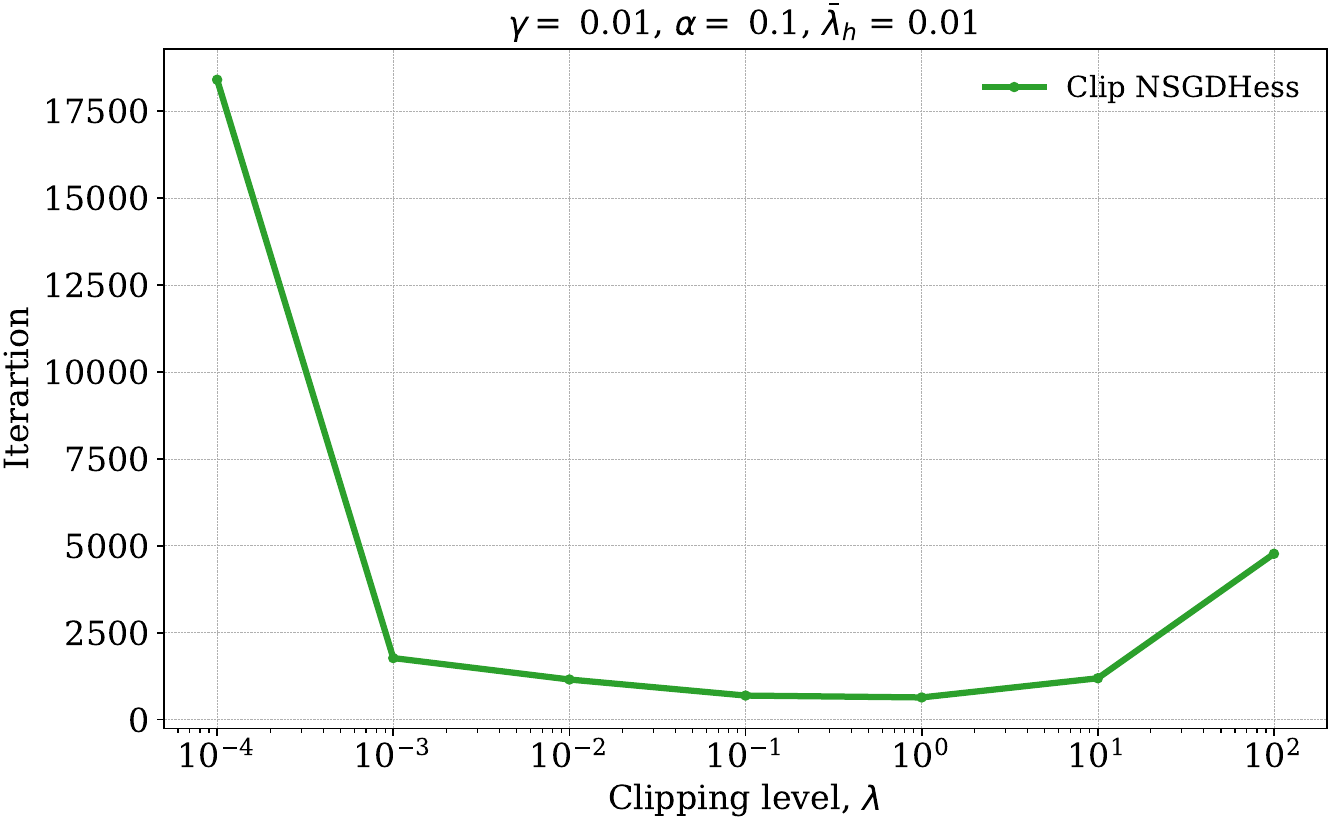}
    \end{minipage}%
    \hfill
    \begin{minipage}[t]{0.48\linewidth}
        \centering
        \includegraphics[width=\linewidth]{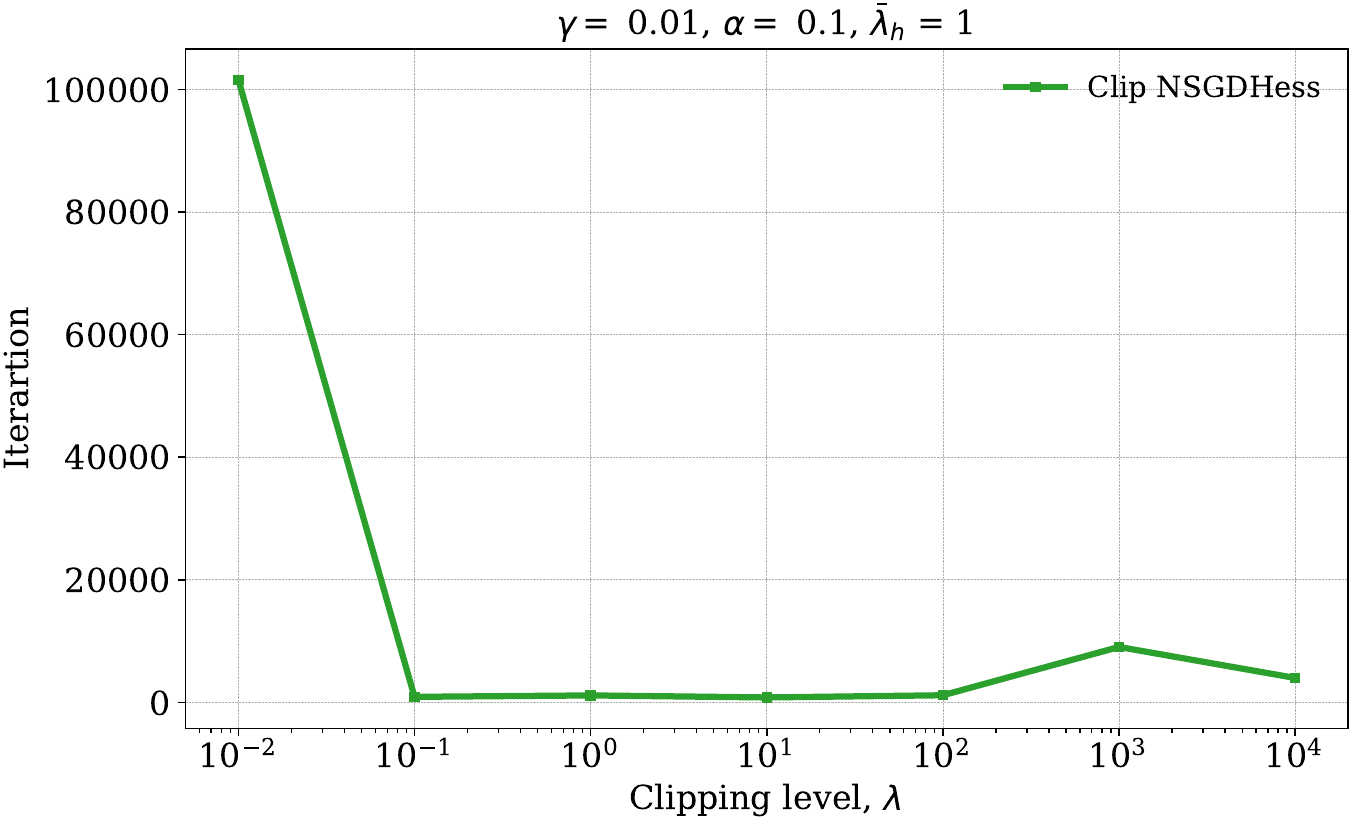}
    \end{minipage}
    \hfill
    \begin{center}
    \begin{minipage}[h]{0.48\linewidth}
        \centering
        \includegraphics[width=\linewidth]{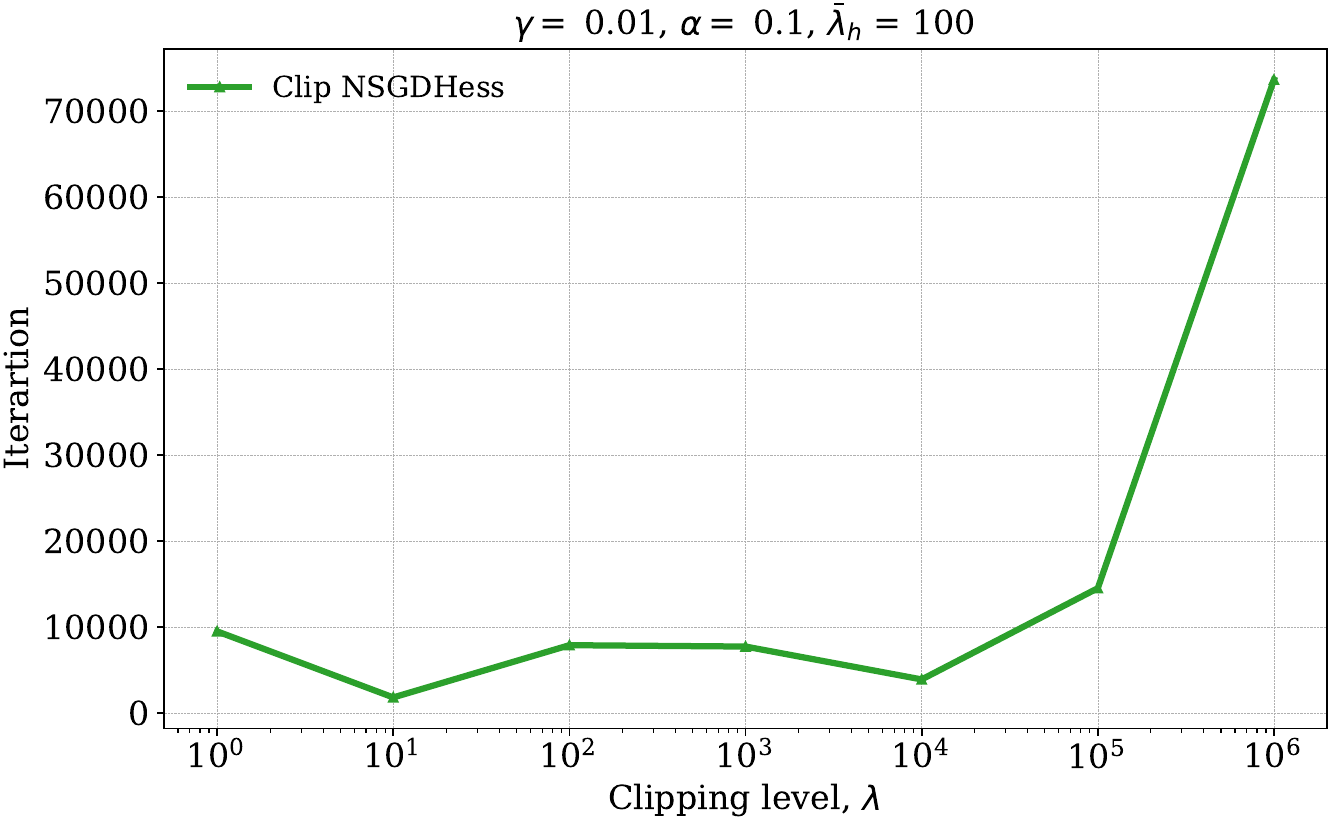}
    \end{minipage}
    \end{center}
    \caption{Iteration complexity of \algname{Clip NSGDHess} (Algorithm~\ref{ch4:alg:NSGD_SOM_clipped_chapter_4}) depending on gradient clipping for the three different fixed values of Hessian clipping ${\bar\lambda_h \in \{0.01, 1, 100\}}$.}
    \label{ch4:fig:dmdskkcs_chapter_4}
\end{figure}

In all cases, the algorithm is tested until $\|\nabla f(x)\| \leq 1/2$. We observe that intermediate values of $\lambda$ yield the best performance across all three regimes (small, medium, and large Hessian clipping levels). As shown in Figure~\ref{ch4:fig:dmdskkcs_chapter_4}, both extremely small and extremely large values of $\lambda$ lead to higher iteration complexity, while moderate choices minimize it. The optimal values are:
{
\begin{align*}
\lambda &= 1.0
&&\text{for } {\bar\lambda_h} = 0.01,\\
\lambda &= 0.1
&&\text{for } {\bar\lambda_h} = 1.0,\\
\lambda &= 10.0
&&\text{for } {\bar\lambda_h} = 100.0.
\end{align*}
}

These findings confirm that moderate clipping levels are crucial for achieving the best performance aligning with our theoretical findings in Theorem~\ref{ch4:thm:main_theorem_hp_chapter_4}.
  \endgroup

\chapter{Supplementary Material and Proofs for Chapter 5}
  \begingroup
  \thesischaptersetup{5}%
  \stepEightAAppendixColors{5}%
  
\clearpage
\begingroup
\setlength{\emergencystretch}{3em}

\section{Notation}\label{ch5:appdx-sec:notation_chapter_5}

{\color{navy}Following~\textcite{pmlr-v125-arjevani20a}, an $r^{\textnormal{th}}$-order tensor ($r\ge0$) is an element of $\R^{\otimes^r d}$. We use $r$ for tensor and oracle order so that the moment exponents $p$ and $q$ retain their thesis-wide meanings. By convention, orders $0$, $1$, and $2$ correspond to scalars, vectors, and matrices. If $r\ge1$, write $T=(T_1,\ldots,T_d)$, where $T_i$ is obtained by fixing the first index to $i$:
\begin{thesischfivefitdisplay}
    [T_i]_{j_1,\ldots,j_{r-1}}=T_{i,j_1,\ldots,j_{r-1}}.
\end{thesischfivefitdisplay}
When $r=2$, $[T]_{i,\cdot}$ and $[T]_{\cdot,j}$ denote the $i$th row and $j$th column. Only first- and second-order tensors are used in the algorithms of this chapter.}

For any integer $n > 0$, we define the index set $[n] \eqdef \thesischfiveens{1, 2, \ldots, n}$. Let $d \ge 1$ denote the ambient dimension. We use $\thesischfiveps{\cdot}{\cdot}$ to represent the standard Euclidean inner product on $\R^d$, i.e., $\thesischfiveps{x}{y} = \sum_{i=1}^d x_i y_i$, and $\norm{\cdot}$ to denote the associated $\ell^2$–norm, $\norm{x} = \sqrt{\thesischfiveps{x}{x}}$. For a matrix $A \in \R^{d \times d}$, we denote by $\normop{A}$ its operator (spectral) norm, defined as $\normop{A} = \sup_{\norm{x}=1} \norm{A x}$. We write {\color{navy}$\nabla f(\cdot;\cdot)$ and $\nabla^2 f(\cdot;\cdot)$} for the stochastic gradient and Hessian oracles, respectively. For any two real numbers $a,b$, we use $a \wedge b \eqdef \min\thesischfiveens{a,b}$ and $a \vee b \eqdef \max\thesischfiveens{a,b}$. We adopt the standard asymptotic notation $\cO(\cdot)$ and $\Omega(\cdot)$ to denote upper and lower bounds on growth rates. To prevent ambiguity, constants related to first-order (gradient) quantities are indexed by the subscript $1$ (e.g., $L_1$, $\sigma_1$), whereas those related to second-order (Hessian) quantities carry the subscript $2$ (e.g., $L_2$, $\sigma_2$).

We denote by $\thesischfiveProba{E}$ the probability of an event $E$ and by $\thesischfiveE{X}$ the expectation of a random variable $X$. Conditional probability and conditional expectation are written, respectively, as $\thesischfiveprobac{E}{\mathcal{F}}$ and $\thesischfiveExpCond{X}{\mathcal{F}}$, where $\mathcal{F}$ denotes a $\sigma$–algebra or a conditioning event. The notation $\thesischfiveExpSub{\xi}{\cdot}$ (or equivalently {\color{navy}$\thesischfiveExpSub{\xi \sim \mathcal P}{\cdot}$}) indicates that the expectation is taken with respect to the randomness of $\xi$, with law {\color{navy}$\mathcal P$}.

\section{Additional Definitions}\label{ch5:appdx-sec:lower-bound-def_chapter_5}
We present in this section the formal setup (function, oracle and algorithm/protocol classes) we considered and under which our (dimension-free) lower and upper bounds are derived.

\subsection{The Setup}

\subsubsection{Function Class}
The lower bounds developed in this work apply to algorithms that find an  $\thesischfiveeps$--stationary point of (nonconvex) functions. All functions considered here are defined from $\R^d$ to $\R$ where $d \ge 1$ is an integer. Depending on the assumptions (Assumption~\ref{ch5:ass:mean-squared-smoothness_chapter_5}, Assumption~\ref{ch5:ass:mean-squared-smoothness-2_chapter_5}, and Assumption~\ref{ch5:ass:L-lipschitz-hessians_chapter_5}), the class of functions may vary. Formally, we define
    \begin{thesischfivefitdisplay} \mathcal{F}(\Delta) \eqdef \thesischfiveenstq{f \in \mathcal{C}^1(\R^d, \R)}{f(0) - \inf_{x \in \R^d} f(x) \le \Delta}, \thesischfivenumberthis\label{ch5:106ec3da-9b9d-4af4-acb0-1a325cdaabae_chapter_5} \end{thesischfivefitdisplay}
which is the class of continuously differentiable and $\Delta$--bounded functions. If the considered function further satisfies some (standard) regularity conditions (see~Assumption~\ref{ch5:ass:L-lipschitz-gradients_chapter_5} and Assumption~\ref{ch5:ass:L-lipschitz-hessians_chapter_5}) we let
    \begin{thesischfivefitalignat}{2}
        \MoveEqLeft[1] \mathcal{F}(\Delta, L_1) \eqdef \Bigl\{ f \in \mathcal{C}^1(\R^d, \R) \, : \, f(0) - \inf_{x \in \R^d} f(x) \le \Delta, \\
        &\quad \norm{\nabla f(x) - \nabla f(y)} \le L_1 \norm{x - y} \text{ for all $x, y \in \R^d$} \Bigr\},
    \end{thesischfivefitalignat}
and
    \begin{thesischfivefitalignat}{2}
        \MoveEqLeft[1] \mathcal{F}(\Delta, L_1, L_2) \eqdef \Bigl\{ f \in \mathcal{C}^2(\R^d, \R) \, : \, f(0) - \inf_{x \in \R^d} f(x) \le \Delta, \\
        &\quad \norm{\nabla f(x) - \nabla f(y)} \le L_1 \norm{x - y}, \\
        &\quad \normop{\nabla^2 f(x) - \nabla^2 f(y)} \le L_2 \norm{x - y} \text{ for all $x, y \in \R^d$} \Bigr\},
    \end{thesischfivefitalignat}
when access to second-order information is possible. The {\color{navy}Lipschitz} constants $L_1$ and $L_2$ appearing above are measured with respect to the canonical (Euclidean) $\ell^2$--norm over $\R^d$ for the gradients and points while, for the Hessian we use the operator norm $\normop{\cdot}$ induced by $\norm{\cdot}$. We recall in~Definition~\ref{ch5:appdx-def:operator-norm_chapter_5} the definition of $\normop{\cdot}$.

The dimension $d \ge 1$ appearing in the above definitions will be made explicit in the proofs (see~Section~\ref{ch5:appdx-sec:proofs-lower-bound_chapter_5}) and it may {\color{navy}depend} on problem-specific parameters, e.g., $\thesischfiveeps$, $\Delta$ and the smoothness constants.

Note that the function class~\eqref{ch5:106ec3da-9b9d-4af4-acb0-1a325cdaabae_chapter_5} as stated is too broad and lacks regularity, e.g., {\color{navy}Lipschitz-continuous} gradients. Nonetheless in some cases, appropriate assumptions on the stochastic gradients (for instance, \textit{mean-squared smoothness}, see~Assumption~\ref{ch5:ass:mean-squared-smoothness_chapter_5}) are enough to ensure the underlying function itself is smooth, i.e., belongs to $\mathcal{F}(\Delta, L)$ for some constant $L \ge 0$ (for the exact derivation of this fact, see~\eqref{ch5:94828cd2-cf02-492b-a6ec-40b3012f66c1_chapter_5}).

\subsubsection{Oracle Class}

\begin{definition}[{Stochastic $r^{\textnormal{th}}$-order Oracles~\textcite{pmlr-v125-arjevani20a}}]\label{ch5:appdx-def:stochastic-pth-order-oracle_chapter_5}
    {\color{navy}Given integers $d,r\ge1$ and $f\in\mathcal C^r(\R^d,\R)$, let $\mathcal O_r(f)$ denote the class of stochastic $r$th-order oracles. An element is a pair $(\mathcal P_\xi,\texttt O_f^r)$, where $\mathcal P_\xi$ is a probability law on a measurable space $\mathcal Z$ and
    \begin{thesischfivefitdisplay}
        \texttt O_f^r(x,\xi)
        =\bigl(f(x),\nabla f(x;\xi),\ldots,\nabla^r f(x;\xi)\bigr).
    \end{thesischfivefitdisplay}
    The oracle is unbiased: for every $s\in[r]$,
    \begin{thesischfivefitdisplay}
        \thesischfiveExpSub{\xi\sim\mathcal P_\xi}{\nabla^s f(x;\xi)}
        =\nabla^s f(x).
    \end{thesischfivefitdisplay}}
\end{definition}
Furthermore, if some of the derivative estimator (here, only the gradient estimator $\nabla f(x;\xi)$ or the Hessian estimator $\nabla^2 f(x;\xi)$ are used) satisfies additional properties or assumption depending on some parameters $\sigma_1$ (Assumption~\ref{ch5:ass:p-bounded-central-moment-gradient_chapter_5}), $\sigma_2$ (Assumption~\ref{ch5:ass:q-bounded-central-moment-hessian_chapter_5}), $\bar{L}$ (Assumption~\ref{ch5:ass:mean-squared-smoothness_chapter_5}), $\delta$ (Assumption~\ref{ch5:ass:mean-squared-smoothness-2_chapter_5}) and so on, the oracle class will be denoted by
    \begin{thesischfivefitdisplay} {\color{navy}\mathcal{O}_r(f, \sigma_1, \sigma_2\ldots)}, \thesischfivenumberthis\label{ch5:440058b9-2a5c-4a0d-a158-4aaff45f095b_chapter_5} \end{thesischfivefitdisplay}
with all parameters listed inside the parenthesis, in an arbitrary order. As each assumption defines unambiguously its own parameters, we keep the notation~\eqref{ch5:440058b9-2a5c-4a0d-a158-4aaff45f095b_chapter_5} for simplicity. {\color{navy}Quantities such as $\delta^q$ and $\sigma_1^p$ are the actual powers appearing in the moment bounds; the oracle order is carried separately by the subscript $r$.}

\subsubsection{Optimization Protocol and Algorithm Class}

First, let us introduce some important definitions.
\begin{definition}[Support of a Vector/Tensor]\label{ch5:appdx-def:support_chapter_5}
    Let $d \ge 1$ be an integer, the support of a vector $x = (x_1, \ldots, x_d) \in \R^d$ is the set
        \begin{thesischfivefitdisplay} \thesischfivesupp(x) \eqdef \thesischfiveenstq{i \in [d]}{x_i \neq 0}, \end{thesischfivefitdisplay}
    i.e., the set of all indices $i \in [d]$ for which $x$ has a nonzero $i^{\textnormal{th}}$ coordinate.

    {\color{navy}For a given $r$th-order tensor $T=(T_1,\ldots,T_d)\in\R^{\otimes^r d}$, its support is defined as}
        \begin{thesischfivefitdisplay} \thesischfivesupp(T) \eqdef \thesischfiveenstq{i \in [d]}{T_i \neq 0}, \end{thesischfivefitdisplay}
    {\color{navy}where $T_i$ is the $i$th subtensor (an $(r-1)$st-order tensor), e.g., the $i$th row of the matrix $T$ when $r=2$.}

    {\color{navy}For a finite tuple or batch $(Z_1,\ldots,Z_m)$ of vectors or tensors, define
    \begin{thesischfivefitdisplay}
        \thesischfivesupp(Z_1,\ldots,Z_m)
        \eqdef\bigcup_{\ell=1}^m\thesischfivesupp(Z_\ell),
        \qquad
        \thesischfiveprog_0(Z_1,\ldots,Z_m)
        \eqdef\max_{\ell\in[m]}\thesischfiveprog_0(Z_\ell).
    \end{thesischfivefitdisplay}}
\end{definition}

\begin{definition}[Progress, \protect\say{$\thesischfiveprog$}]
    Let $d \ge 1$ be an integer, for any $x = (x_1, \ldots, x_d) \in \R^d$ and any $\alpha \in \thesischfiveintfo{0}{+\infty}$ we define
        \begin{thesischfivefitdisplay} \thesischfiveprog_{\alpha}(x) \eqdef \max\thesischfiveenstq{i \in [d]}{\abs{x_i} > \alpha}, \end{thesischfivefitdisplay}
    and $\thesischfiveprog_{\alpha}(x) = 0$ if $\abs{x_i} \le \alpha$ for all $i \in [d]$.

    {\color{navy}If $\alpha=0$, for a given $r$th-order tensor $T=(T_1,\ldots,T_d)\in\R^{\otimes^r d}$, its \say{$\thesischfiveprog$} is defined as}
        \begin{thesischfivefitdisplay} \thesischfiveprog_0(T) \eqdef \max\thesischfiveenstq{i \in [d]}{T_i \neq 0}, \end{thesischfivefitdisplay}
    {\color{navy}and is set to $0$ if there is no index $i\in[d]$ for which $T_i\neq0$.}
\end{definition}
Notably, if $\alpha = 0$ then $\thesischfiveprog_0(x)$ is the largest index at which $x \in \R^d$ has a nonzero coordinate. $\thesischfiveprog_0$ will be used to capture the rate at which new coordinates are \say{discovered}. Initially, all coordinates are set to $0$ and, as we progressively acquire information from the queries to the oracle, the union of the support of the oracle responses grows. {\color{navy}The growth rate is quantified using $\thesischfiveprog$ and controlled through the notion of a \emph{zero-chain}, recalled formally in Definitions~\ref{ch5:appdx-def:deterministic-zero-chain_chapter_5} and~\ref{ch5:appdx-def:probability-zero-chain_chapter_5}.}

We recall below some technical notions from the paper of~\textcite{arjevani2022lower}.
\paragraph{Optimization Protocol:} {\color{navy}Let $b\ge1$ denote the number of simultaneous queries per oracle round. At round $i\ge1$, the algorithm selects
$\mathbf x^{(i)}=(x_1^{(i)},\ldots,x_b^{(i)})$. The oracle then draws one fresh seed $\xi^{(i)}\sim\mathcal P_\xi$, independently across rounds, and reuses that seed throughout the batch:
\begin{thesischfivefitdisplay}
    \texttt O_f^{r,b}(\mathbf x^{(i)},\xi^{(i)})
    \eqdef\bigl(\texttt O_f^r(x_1^{(i)},\xi^{(i)}),\ldots,
    \texttt O_f^r(x_b^{(i)},\xi^{(i)})\bigr).
\end{thesischfivefitdisplay}
Thus the draws are independent across rounds, not across points within one batch. Algorithm~\ref{ch5:algo:nsgd-mvr_chapter_5} has $b=2$.}

\paragraph{Optimization Algorithms:} a (randomized) algorithm $\texttt{A}$ consists of a probability law $\mathcal P_{\rm alg}$ on a measurable set $\mathcal R$, a seed $r_{\rm alg}\sim\mathcal P_{\rm alg}$ drawn at the beginning, and measurable mappings $\{\texttt A^{(i)}\}_{i\ge1}$. {\color{navy}The resulting query sequence is denoted by $\{\mathbf x_{\texttt A[\texttt O_f^{r,b}]}^{(i)}\}_{i\ge1}$ and each $\mathbf x^{(i)}$ is measurable with respect to the algorithmic seed and all responses through round $i-1$.}

We define {\color{navy}$\mathcal{A}_{\texttt{rand}}(b)$} as the class of all algorithms that follow the aforementioned protocol with a batch size of {\color{navy}$b$} queries per round.

\begin{definition}[{Zero-Respecting Algorithm~\textcite[Definition~1]{arjevani2022lower,pmlr-v125-arjevani20a}}]\label{ch5:appdx-def:zero-especting-algo_chapter_5}
    {\color{navy}A stochastic $r$th-order algorithm $\texttt A$ is zero-respecting if, for every $f\in\mathcal C^r(\R^d,\R)$, every oracle $\texttt O_f^{r,b}$, and every round $i\ge1$,
    \begin{thesischfivefitdisplay}
        \thesischfivesupp\!\left(\mathbf x_{\texttt A[\texttt O_f^{r,b}]}^{(i)}\right)
        \subseteq
        \bigcup_{j=1}^{i-1}
        \thesischfivesupp\!\left(
        \texttt O_f^{r,b}\!\left(
        \mathbf x_{\texttt A[\texttt O_f^{r,b}]}^{(j)},\xi^{(j)}
        \right)\right)
        \quad\text{almost surely}.
    \end{thesischfivefitdisplay}
    We denote this class by
    $\mathcal A_{\texttt{zr}}(b)\subseteq\mathcal A_{\texttt{rand}}(b)$.}
\end{definition}
In other words, a zero-respecting algorithm cannot modify coordinates where no information is known; that is, its queries at each round have support in the union of the supports of all previous oracle responses.

{\color{navy}When the batch arity is fixed, we suppress $b$ from $\mathcal A_{\texttt{zr}}(b)$ to lighten notation.}

\paragraph{Proof Strategy: Zero-Chains}
The main motivation for building the hard instance~\eqref{ch5:f9d4cc59-a955-4e12-997e-7870c5314204_chapter_5} is the notion of zero-chain, which in some sense, are functions for which it is \say{hard} to uncover new, nonzero, coordinates.

\begin{definition}[{Deterministic $r^{\textnormal{th}}$-order Zero-Chain~\textcite[Definition~3]{carmon2020lower}}]\label{ch5:appdx-def:deterministic-zero-chain_chapter_5}
    {\color{navy}Given $r\ge1$ and $f\in\mathcal C^r(\R^d,\R)$, we say that $f$ is an $r$th-order zero-chain if, for every $x\in\R^d$ and $i\in[d]$,
    \begin{thesischfivefitdisplay}
        \thesischfivesupp(x)\subseteq\{1,\ldots,i-1\}
        \quad\Longrightarrow\quad
        \bigcup_{s=1}^r\thesischfivesupp(\nabla^s f(x))
        \subseteq\{1,\ldots,i\}.
        \thesischfivenumberthis\label{ch5:15294226-4f58-484e-b80c-22c08bec83d5-0_chapter_5}
    \end{thesischfivefitdisplay}}

    {\color{navy}We say that $f$ is a zero-chain if it is an $r$th-order zero-chain for every integer $r\ge1$.}
\end{definition}
In other word,~Definition~\ref{ch5:appdx-def:deterministic-zero-chain_chapter_5} tells us that given $x \in \R^d$, we can \say{discover} at most one new coordinate when accessing the gradient of $f$ at $x$ (or any high-order derivatives of $f$). When dealing with stochastic estimators instead, we can extend the previous~Definition~\ref{ch5:appdx-def:deterministic-zero-chain_chapter_5} as follows:
\begin{definition}[{Probability-$\theta$ $r^{\textnormal{th}}$-order Zero-Chain~\textcite[Definition~2]{pmlr-v125-arjevani20a}}]\label{ch5:appdx-def:probability-zero-chain_chapter_5}
    {\color{navy}Given $r\ge1$, $\theta\in(0,1]$, $f\in\mathcal C^r(\R^d,\R)$, and unbiased derivative estimators through order $r$, we call them a probability-$\theta$ $r$th-order zero-chain (at the hard-instance threshold $1/4$) if}
        \begin{thesischfivefitdisplay}{\color{navy}\thesischfiveProba{\exists x \mid \thesischfiveprog_0\left( \nabla f(x;\xi), \ldots, \nabla^r f(x;\xi) \right) > \thesischfiveprog_{\frac14}(x) + 1} = 0}, \thesischfivenumberthis\label{ch5:15294226-4f58-484e-b80c-22c08bec83d5-1_chapter_5} \end{thesischfivefitdisplay}
    and
        \begin{thesischfivefitdisplay} {\color{navy}\thesischfiveProba{\exists x \mid \thesischfiveprog_0\left( \nabla f(x;\xi), \ldots, \nabla^r f(x;\xi) \right) = \thesischfiveprog_{\frac14}(x) + 1} \le \theta}, \thesischfivenumberthis\label{ch5:15294226-4f58-484e-b80c-22c08bec83d5-2_chapter_5} \end{thesischfivefitdisplay}

    {\color{navy}We say probability-$\theta$ zero-chain when this holds for every order $r\ge1$.}
\end{definition}
In the above~Definition~\ref{ch5:appdx-def:probability-zero-chain_chapter_5}, condition~\eqref{ch5:15294226-4f58-484e-b80c-22c08bec83d5-1_chapter_5} is the analogue of condition~\eqref{ch5:15294226-4f58-484e-b80c-22c08bec83d5-0_chapter_5} when we have only access to noisy derivatives estimators of $f$. In addition, condition~\eqref{ch5:15294226-4f58-484e-b80c-22c08bec83d5-2_chapter_5} tells us that we have a \say{small} chance to discover a new coordinate upon querying the oracle, i.e., the added noise behaves \say{adversarially} and can slow down the discovery of new coordinates (in expectation).

An important properties of probability-$\theta$ zero-chain is the following lemma.
\begin{lemma}[{\textcite[Lemma~16]{pmlr-v125-arjevani20a}}]\label{ch5:appdx-lem:zero-chain-progress_chapter_5}
    {\color{navy}Let $r\ge1$, $\theta\in(0,1]$, and let the unbiased derivative estimators through order $r$ form a probability-$\theta$ zero-chain. Let $\texttt O_f^r$ denote the associated oracle.}
    
    {\color{navy}Let $\{\mathbf x_{\texttt A[\texttt O_f^r]}^{(i)}\}_{i\ge1}$ be the queries of any zero-respecting algorithm. For $\beta\in(0,1)$ satisfying $N>\log(1/\beta)$, with probability at least $1-\beta$,
    \begin{thesischfivefitdisplay}
        \thesischfiveprog_0\!\left(\mathbf x_{\texttt A[\texttt O_f^r]}^{(t)}\right)<N
        \quad\text{for all}\quad
        t\le\frac{N-\log(1/\beta)}{2\theta}.
    \end{thesischfivefitdisplay}
    If $N\le\log(1/\beta)$, the displayed time range is empty and the assertion is vacuous.}
\end{lemma}
The proof can be found in~\textcite[Lemma~16]{pmlr-v125-arjevani20a} and~\textcite[Lemma~1]{arjevani2020tight}.

\subsubsection{Complexity Measures}

As in~\textcite{arjevani2022lower,pmlr-v125-arjevani20a}, we develop lower bounds on the \emph{distributional complexity} for finding an $\thesischfiveeps$--stationary point, {\color{navy}which, in turn, implies lower bounds on the \emph{minimax complexity}} for finding such stationary points. Formally, following~\textcite{arjevani2022lower}
    \begin{thesischfivefitalignat}{2}
        \MoveEqLeft[1] {\color{navy}\mathfrak{m}^{\texttt{zr}}_{\thesischfiveeps}(b, \Delta,\ldots)} \\
        & {\color{navy}\eqdef \sup_{f \in \mathcal{F}(\Delta,\ldots)} \sup_{(\mathcal P_{\xi}, \texttt{O}_f) \in \mathcal{O}(f,\ldots)} \inf_{\texttt{A} \in \mathcal{A}_{\texttt{zr}}(b)} \inf \thesischfiveenstq{t \ge 1}{\thesischfiveExpSub{\texttt{O}_f, \texttt{A}}{\norm{\nabla f\left( \mathbf{x}_{{\normalfont\texttt{A}}[{\normalfont\texttt{O}_f}]}^{(t)} \right)}} \le \thesischfiveeps}},
    \end{thesischfivefitalignat}
where the expectation is taken over the randomness in the oracle $\texttt{O}_f$ and in the algorithm $\texttt{A} \in \mathcal{A}_{\texttt{zr}}$, if any.

As in the definition of the oracle class, complexity measures may depend on various parameters depending on the assumptions considered. The convention is to list all involved parameters (in an arbitrary order), and there is no ambiguity in doing so.

\subsection{The Worst-Case Function}\label{ch5:appdx-subsec:hard-instance_chapter_5}

In this section, we recall the \say{worst-case} function introduced in~\textcite{carmon2020lower,arjevani2022lower} and used to prove our lower bounds. {\color{navy}This function, or a variation adapted to the target function class, underlies many lower bounds in stochastic nonconvex optimization}~\parencite{arjevani2022lower,NEURIPS2024_dd850be1,10.5555/3737916.3739635,islamov2025safeef,sun2025improved,sadiev2025second}. Formally, given an integer $N \ge 1$ which denotes the dimension in which the hard instance $f_N \colon \R^N \to \R$ lies, we define
    \begin{thesischfivefitdisplay} f_N \colon x \mapsto -\Psi(1) \Phi(x_1) + \sum_{i = 2}^N \left[ \Psi(-x_{i - 1}) \Phi(-x_i) - \Psi(x_{i - 1}) \Phi(x_i) \right], \thesischfivenumberthis\label{ch5:f9d4cc59-a955-4e12-997e-7870c5314204_chapter_5} \end{thesischfivefitdisplay}
where the inner components $\Psi$ and $\Phi$ are defined, for any $t \in \R$, as
    \begin{thesischfivefitdisplay} \Psi(t) \eqdef \begin{cases}
        0, & \text{if $t \le \frac{1}{2}$;} \\
        \exp\left( 1 - \frac{1}{(2 t - 1)^2} \right), & \text{if $t > \frac{1}{2}$;}
    \end{cases} \quad \text{ and } \quad \Phi(t) \eqdef \sqrt{e} \int_{-\infty}^t e^{-\frac{s^2}{2}} \odif{s}. \end{thesischfivefitdisplay}
The particular design~\eqref{ch5:f9d4cc59-a955-4e12-997e-7870c5314204_chapter_5} of $f_N$ is such that it is a (deterministic) zero-chain (see~Definition~\ref{ch5:appdx-def:deterministic-zero-chain_chapter_5}) and the gradients of $f_N$ are large unless all coordinates are large (i.e., $\thesischfiveprog_1(x) \ge N$); these two keys properties are recalled below ($5.$ and $6.$). Several other properties of the function $f_N$ are stated in the next lemma.
\begin{lemma}[Properties of the Hard Instance~\textcite{carmon2020lower,arjevani2022lower}]\label{ch5:appdx-lem:properties-hard-instance_chapter_5}
    For any integer $N \ge 1$, the function $f_N$ satisfies:
    \begin{enumerate}
        \item $f_N(0) - \inf_{x \in \R^N} f_N(x) \le \Delta_0 N$ where $\Delta_0 \eqdef 12$,

        \item The gradient of $f_N$ is $\ell_1$--Lipschitz continuous with $\ell_1 \eqdef 152$,

        \item {\color{navy}For all $x \in \R^N$, we have $\norminf{\nabla f_N(x)} \le \gamma_{\infty}$, where $\gamma_{\infty} \eqdef 23$. In particular, $\norm{\nabla f_N(x)} \le 23\sqrt N$.}

        \item {\color{navy}There exists a universal constant $0<c<+\infty$} such that for every integer $p \ge 1$, the $p$--th order derivatives of $f_N$ are $\ell_p$--Lipschitz continuous for some $\ell_p \le \exp\left( \frac{5}{2} p \log p + c p \right)$,

        \item For all $x \in \R^N$ we have $\thesischfiveprog_0(\nabla f_N(x)) \le \thesischfiveprog_{\frac{1}{2}}(x) + 1$,

        \item {\color{navy}For all $x \in \R^N$, if $\thesischfiveprog_1(x) < N$, then $\norm{\nabla f_N(x)} > 1$.}
    \end{enumerate}
\end{lemma}
In the proofs of our lower bounds, we mostly focus on properties $1.$, $2.$, $3.$ and $4.$ which, after a specific rescaling, are enough to ensure the hard instance $f_N$ belongs to the right class of functions. {\color{navy}Properties $5.$ and $6.$ are not used explicitly in this chapter, but they preserve the mechanism of the hard instance from~\textcite{carmon2020lower,arjevani2022lower}. Consequently, zero-respecting algorithms behave as required and the original lower-bound strategy remains applicable.}

\subsection{The Stochastic Gradient and Hessian Estimator}

We recall below from~\textcite{arjevani2022lower,pmlr-v125-arjevani20a} the gradient and Hessian estimators used:
\begin{definition}[Gradient Estimator]\label{ch5:appdx-def:gradient-estimator_chapter_5}
    {\color{navy}Given an integer $N \ge 1$ and $\theta \in \thesischfiveintof{0}{1}$, the gradient estimator $g_N \colon \R^N \times \mathcal{Z} \to \R^N$ of the hard instance $f_N$ in~\eqref{ch5:f9d4cc59-a955-4e12-997e-7870c5314204_chapter_5} is defined coordinate-wise by
        \begin{thesischfivefitdisplay} \left[ g_N(x;\xi) \right]_i \eqdef \left[ \nabla f_N(x) \right]_i \left( 1 + \mathbb{I}\!\thesischfiveens{i > \thesischfiveprog_{\frac{1}{4}}(x)} \left( \frac{\xi}{\theta} - 1 \right) \right). \end{thesischfivefitdisplay}
    Here $i \in [N]$, $x \in \R^N$, and $\xi \sim \thesischfiveBer(\theta)$, so the oracle law is $\mathcal P_\xi=\thesischfiveBer(\theta)$ on $\mathcal{Z} = \thesischfiveens{0, 1}$.}
\end{definition}

\begin{definition}[Hessian Estimator]\label{ch5:appdx-def:hessian-estimator_chapter_5}
    {\color{navy}Given an integer $N \ge 1$ and $\theta \in \thesischfiveintof{0}{1}$, the Hessian estimator $\nabla g_N \colon \R^N \times \mathcal{Z} \to \R^{N\times N}$ of the hard instance $f_N$ in~\eqref{ch5:f9d4cc59-a955-4e12-997e-7870c5314204_chapter_5} is defined row-wise by
        \begin{thesischfivefitdisplay} \left[ \nabla g_N(x;\xi) \right]_{i, \cdot} \eqdef \left( 1 + \mathbb{I}\!\thesischfiveens{i > \thesischfiveprog_{\frac{1}{4}}(x)} \left( \frac{\xi}{\theta} - 1 \right) \right) \left[ \nabla^2 f_N(x) \right]_{i, \cdot}. \end{thesischfivefitdisplay}
    Here $i \in [N]$, $x \in \R^N$, and $\xi \sim \thesischfiveBer(\theta)$.}
\end{definition}
Note that, for all $x \in \R^N$ and all $i > \thesischfiveprog_{\frac{1}{4}}(x) + 1$ we have $\left[ \nabla f_N(x) \right]_i = 0$ by~Lemma~\ref{ch5:appdx-lem:properties-hard-instance_chapter_5} (property $5.$) so, only the specific coordinate at $i = \thesischfiveprog_{\frac{1}{4}}(x) + 1$ is noisy.

In the case of~Assumption~\ref{ch5:ass:mean-squared-smoothness_chapter_5} and Assumption~\ref{ch5:ass:mean-squared-smoothness-2_chapter_5}, i.e., additional regularity conditions on the gradient estimator, we define, as in~\textcite{arjevani2022lower} the \say{smoothed} indicator $\Theta_i$ of $\mathbb{I}\!\thesischfiveens{i > \thesischfiveprog_{1/4}(\cdot)}$ as:
    \begin{thesischfivefitdisplay} {\color{navy}\Theta_i \colon x \mapsto \Gamma\left( 1 - \left( \sum_{k = i}^N \Gamma\left(\abs{x_k}\right)^2 \right)^{1/2} \right) = \Gamma\left( 1 - \norm{\left( \Gamma(\abs{x_i}), \ldots, \Gamma(\abs{x_N}) \right)} \right)}, \end{thesischfivefitdisplay}
where $\Gamma \colon \R \to \R$ is any $\mathcal{C}^{\infty}$, non-decreasing Lipschitz continuous function with $\Gamma(t) = 0$ for all $t \le 1/4$ and $\Gamma(t) = 1$ for all $t \ge 1/2$. For instance, we can take, as in~\textcite{arjevani2022lower}
    \begin{align*}
        \Gamma(t) & \eqdef \frac{\int_{1/4}^t \Lambda(s) \odif{s}}{\int_{1/4}^{1/2} \Lambda(s) \odif{s}}, \quad \text{where} \\
        \Lambda(t) & \eqdef \begin{cases}
        0, & \text{if $t \le 1/4$ or $t \ge 1/2$;} \\
        \exp\left( -\dfrac{1}{100 \left( t - 1/4 \right) \left( 1/2 - t \right)} \right), & \text{if $1/4 < t < 1/2$,}
    \end{cases}
    \end{align*}
for any $t \in \R$. Notably, the smoothed indicator $\Theta_i$ satisfies for all $x \in \R^N$
    \begin{thesischfivefitdisplay} \mathbb{I}\!\thesischfiveens{i > \thesischfiveprog_{1/4}(x)} \le \Theta_i(x) \le \mathbb{I}\!\thesischfiveens{i > \thesischfiveprog_{1/2}(x)}, \end{thesischfivefitdisplay}
hence
    \begin{thesischfivefitdisplay} \Theta_i(x) = \begin{cases}
        1, & \text{for all $i > \thesischfiveprog_{\frac{1}{4}}(x)$;} \\
        0, & \text{for all $i \le \thesischfiveprog_{\frac{1}{2}}(x)$.}
    \end{cases} \end{thesischfivefitdisplay}

Overall, this gives the following \say{smoothed} gradient estimator:
\begin{definition}[\protect\say{Smoothed} Gradient Estimator]\label{ch5:appdx-def:smoothed-gradient-estimator_chapter_5}
    {\color{navy}Given an integer $N \ge 1$ and $\theta \in \thesischfiveintof{0}{1}$, the smoothed gradient estimator $g_N \colon \R^N \times \mathcal{Z} \to \R^N$ of the hard instance $f_N$ in~\eqref{ch5:f9d4cc59-a955-4e12-997e-7870c5314204_chapter_5} is defined coordinate-wise by
        \begin{thesischfivefitdisplay} \left[ g_N(x;\xi) \right]_i \eqdef \left[ \nabla f_N(x) \right]_i \left( 1 + \Theta_i(x) \left( \frac{\xi}{\theta} - 1 \right) \right). \end{thesischfivefitdisplay}
    Here $i \in [N]$, $x \in \R^N$, and $\xi \sim \thesischfiveBer(\theta)$.}
\end{definition}

The gradient and Hessian estimators defined above are all probability-$\theta$ zero-chain.

\newpage
\section{Proofs of the Lower Bounds}\label{ch5:appdx-sec:proofs-lower-bound_chapter_5}

\subsection{A Lower Bound in the \textit{Global Stochastic Model} under Heavy-Tail Noise}

In this part, we extend~\textcite[Lemma~11]{arjevani2022lower} in the setting where we only assume the variance of the gradient estimator to have bounded $p$--th moment for some $p \in \thesischfiveintof{1}{2}$. This result holds in any dimension $d \ge 1$.

\begin{lemma}\label{ch5:appdx-lem:lower-bound-global-stichastic-model-mean-squared-smoothness-p-BCM_chapter_5}
    Under~Assumption~\ref{ch5:ass:lower-boundedness_chapter_5}, Assumption~\ref{ch5:ass:p-bounded-central-moment-gradient_chapter_5}, and Assumption~\ref{ch5:ass:mean-squared-smoothness_chapter_5} and as long as $0 < \thesischfiveeps \le \sqrt{\bar{L} \Delta}/8$, the number of samples required to obtain an $\thesischfiveeps$--stationary point in the global stochastic model defined above is
        \begin{thesischfivefitdisplay} \Omega(1) \cdot \left( \frac{\sigma_1}{\thesischfiveeps} \right)^{p/(p - 1)}, \end{thesischfivefitdisplay}
    where $p \in \thesischfiveintof{1}{2}$.
\end{lemma}

In particular, this lower bound does not depend on the exponent $q$ from~Assumption~\ref{ch5:ass:mean-squared-smoothness_chapter_5}.

\begin{remark}\label{ch5:appdx-rem:improve-lower-bound-global-stichastic-model-mean-squared-smoothness-p-BCM_chapter_5}
    More precisely,~Lemma~\ref{ch5:appdx-lem:lower-bound-global-stichastic-model-mean-squared-smoothness-p-BCM_chapter_5} shows that
        \begin{thesischfivefitdisplay} \mathfrak{m}^{\normalfont\texttt{zr}}_{\thesischfiveeps} \left( {\color{navy}b}, \Delta, \bar{L}, \sigma_1^p \right) \ge \Omega(1) \cdot \max\thesischfiveens{1, \left( \frac{\sigma_1}{\thesischfiveeps} \right)^{p/(p - 1)}}. \thesischfivenumberthis\label{ch5:df3dc971-577b-4af1-818f-ae694fee0100_chapter_5} \end{thesischfivefitdisplay}
\end{remark}

\begin{remark}\label{ch5:appdx-rem:improve-lower-bound-global-stichastic-model-mean-squared-smoothness-2-p-BCM_chapter_5}
    It is also worth mentioning that the same lower bound holds under~Assumption~\ref{ch5:ass:lower-boundedness_chapter_5}, Assumption~\ref{ch5:ass:L-lipschitz-gradients_chapter_5}, Assumption~\ref{ch5:ass:p-bounded-central-moment-gradient_chapter_5}, and Assumption~\ref{ch5:ass:mean-squared-smoothness-2_chapter_5} since, by our choice of function~\eqref{ch5:af23ded8-f6c2-4361-8fc4-1e72f4395d6d_chapter_5} for any $s \in \thesischfiveens{-1, 1}$ and any $x, y \in \R^d$ we have
        \begin{thesischfivefitalignat}{2}
            \MoveEqLeft[1] \thesischfiveExpSub{\xi \sim \Prob^s}{\norm{\left[ \nabla f_d^{\star}(x, \xi) - \nabla f_d^{\star}(y, \xi) \right] - \left[ \nabla f_{d, s}^{\star}(x) - \nabla f_{d, s}^{\star}(y) \right]}^q} \\
            \thesischfiveoversetlab{\eqref{ch5:70110a40-c2b6-4b8f-b790-e09561a73316_chapter_5}}&{=} {\color{navy}\chi}^q \thesischfiveExpSub{\xi \sim \Prob^s}{\norm{\left[ \nabla f_d({\color{navy}\chi} x, \xi) - \nabla f_d({\color{navy}\chi} y, \xi) \right] - \left[ \nabla f({\color{navy}\chi} x;d, s) - \nabla f({\color{navy}\chi} y;d, s) \right]}^q} \\
            \thesischfiveoversetlab{\eqref{ch5:2130180a-a3a9-4040-8a04-214e99e99948_chapter_5}}&{=} {\color{navy}\chi}^q \thesischfiveExpSub{\xi \sim \Prob^s}{\norm{\left[ L ({\color{navy}\chi} x - \boldsymbol{\xi}) - L({\color{navy}\chi} y - \boldsymbol{\xi}) \right] - \left[ L ({\color{navy}\chi} x - \theta_s) - L ({\color{navy}\chi} y - \theta_s) \right]}^q} \\
            & = 0,
        \end{thesischfivefitalignat}
    hence~Assumption~\ref{ch5:ass:mean-squared-smoothness-2_chapter_5} is satisfied. Moreover, taking $L = 2 \Delta/r^2$ is enough to ensure that $f_{d, s}^{\star}$ has $\Delta$--bounded sub-optimality (see~\eqref{ch5:f38cccea-2eb2-49b1-910d-1eb1671cb5ff_chapter_5}). On the other hand, for any $x, y \in \R^d$ we have
        \begin{thesischfivefitalignatnum}{2}
            \norm{\nabla f_{d, s}^{\star}(x) - \nabla f_{d, s}^{\star}(y)} \thesischfiveoversetlab{\eqref{ch5:70110a40-c2b6-4b8f-b790-e09561a73316_chapter_5}}&{=} {\color{navy}\chi} \norm{\nabla f({\color{navy}\chi} x;d, s) - \nabla f({\color{navy}\chi} y;d, s)} \\
            & = {\color{navy}\chi}^2 L \norm{x - y}, \label{ch5:1446c17c-7865-4a08-be55-9a576a51828f_chapter_5}
        \end{thesischfivefitalignatnum}
    {\color{navy}and it suffices to take $0<\chi\le\sqrt{L_1/L}$} to ensure the function $f_{d, s}^{\star}$ has $L_1$--Lipschitz gradients. The rest of the proof (\textbf{steps 3}, \textbf{4} and \textbf{5}) is the same, using $L_1$ instead of $\bar{L}$.
\end{remark}

\begin{remark}\label{ch5:appdx-rem:improve-lower-bound-global-stichastic-model-mean-squared-smoothness-2-hessian_chapter_5}
    {\color{navy}As a final remark before proving~Lemma~\ref{ch5:appdx-lem:lower-bound-global-stichastic-model-mean-squared-smoothness-p-BCM_chapter_5}, the construction used for~\eqref{ch5:df3dc971-577b-4af1-818f-ae694fee0100_chapter_5} remains valid even when the oracle supplies higher-order information, such as the Hessian. Indeed, $f_d$ is quadratic, so for any $s \in \thesischfiveens{-1,1}$ and any $x,y\in\R^d$ we have}
        \begin{thesischfivefitalignat}{2}
            \MoveEqLeft[1] \thesischfiveExpSub{\xi \sim \Prob^s}{\normop{\nabla^2 f^{\star}_d(x, \xi) - \nabla^2 f^{\star}_{d, s}(x)}^q} \\
            & = {\color{navy}\chi}^{2 q} \thesischfiveExpSub{\xi \sim \Prob^s}{\normop{\nabla^2 f_d({\color{navy}\chi} x, \xi) - \nabla^2 f({\color{navy}\chi} x;d, s)}^q} \\
            & = {\color{navy}\chi}^{2 q} \thesischfiveExpSub{\xi \sim \Prob^s}{\normop{L I_d - L I_d}^q} \\
            & = 0,
        \end{thesischfivefitalignat}
    where $I_d \in \R^{d \times d}$ is the identity matrix. Hence,~Assumption~\ref{ch5:ass:q-bounded-central-moment-hessian_chapter_5} is satisfied. Additionally, since $\nabla^2 f^{\star}_{d, s}(x) = L {\color{navy}\chi}^2 I_d$ for any $x \in \R^d$ then the Hessian of $f^{\star}_{d, s}$ is $L_2$--Lipschitz continuous for any constant $L_2 \ge 0$. Combining this observation with the arguments in the previous remark (see~\eqref{ch5:1446c17c-7865-4a08-be55-9a576a51828f_chapter_5}), we obtain that $f^{\star}_{d, s} \in \mathcal{F}(\Delta, L_1, L_2)$ for any $d \ge 1$ and any $s \in \thesischfiveens{-1, 1}$, as long as we choose
        \begin{thesischfivefitdisplay} L = \frac{2 \Delta}{r^2} \quad \text{and} \quad 0 < {\color{navy}\chi} \le \sqrt{\frac{L_1}{L}}. \end{thesischfivefitdisplay}
    The rest of the proof (\textbf{steps 3}, \textbf{4} and \textbf{5}) is the same.
\end{remark}

\begin{proof}
    Let the accuracy parameter $\thesischfiveeps > 0$, initial sub-optimality $\Delta \ge 0$, the mean-squared smoothness parameter $\bar{L}$, and the variance parameter $\sigma_1 \ge 0$ be fixed, and $L > 0$ to be specified later. Our proof follows the same steps as in~\textcite{arjevani2020tight} except that, taking inspiration from~\textcite{hübler2025gradientclippingnormalizationheavy} and the gradient oracle from~\textcite{arjevani2022lower} as recalled in~Definition~\ref{ch5:appdx-def:gradient-estimator_chapter_5}, we take two Bernoulli distributions instead of normal distributions.

    Let us fix $d \ge 1$ an integer and consider the following family of functions $f_d \colon \R^d \times \thesischfiveens{-1, 0, 1} \to \R$ defined as
        \begin{thesischfivefitdisplay} f_d(x, \xi) \eqdef \frac{L}{2} \left( \sqnorm{x} - 2 x_1 \xi + r^2 \right), \thesischfivenumberthis\label{ch5:af23ded8-f6c2-4361-8fc4-1e72f4395d6d_chapter_5} \end{thesischfivefitdisplay}
    where $r \in \thesischfiveintof{0}{1}$ is a fixed parameter (to be specified later), $(x, \xi) \in \R^d \times \thesischfiveens{-1, 0, 1}$ and $x = (x_1, \ldots, x_d)$ are its coordinates. Then, for ${\color{navy}\chi} > 0$, we rescale the function $f_d$ as
        \begin{thesischfivefitdisplay} f_d^{\star} (x, \xi) \colon x \mapsto f_d({\color{navy}\chi} x, \xi). \thesischfivenumberthis\label{ch5:70110a40-c2b6-4b8f-b790-e09561a73316_chapter_5} \end{thesischfivefitdisplay}

    Next, we define the two probability distributions $\Prob^1 = \thesischfiveBer(r)$ and $\Prob^{-1} = -\thesischfiveBer(r)$ (whose support is $\thesischfiveens{-1, 0}$) and we let $\theta_s \eqdef (r s, 0, \ldots, 0) \in \R^d$ for $s \in \thesischfiveens{-1, 1}$. Additionally, for any $s \in \thesischfiveens{-1, 1}$, we define the function $f_{d, s}^{\star} \colon \R^d \to \R$ for all $x \in \R^d$ as
        \begin{align*}
            f_{d, s}^{\star}(x) & \eqdef \thesischfiveExpSub{\xi \sim \Prob^s}{f_d^{\star}\left(x, \xi \right)} = \frac{L}{2} \left( \sqnorm{{\color{navy}\chi} x} - 2 {\color{navy}\chi} x_1 r s + r^2 \right) \\
            & = \frac{L}{2} \sqnorm{{\color{navy}\chi} x - \theta_s} \thesischfivereqdef f({\color{navy}\chi} x;d, s),
        \end{align*}
    where the first equality follows from the fact that, when $s = 1$ then $\xi \sim \thesischfiveBer(r)$ and $\thesischfiveE{\xi} = r = r s$, while when $s = -1$ then $\xi \sim -\thesischfiveBer(r)$ and $\thesischfiveE{\xi} = -r = r s$.

    \begin{itemize}
        \item \textbf{Step 1:} \textit{Ensuring $f_{d, s}^{\star}$ has $\Delta$--bounded initial sub-optimality (Assumption~\ref{ch5:ass:lower-boundedness_chapter_5})}.

        To ensure $f_{d, s}^{\star}$ satisfies~Assumption~\ref{ch5:ass:lower-boundedness_chapter_5}, let us compute the initial sub-optimality. Assuming the starting point is $x^0 = 0$, we have
            \begin{thesischfivefitdisplay} f^{\star}_{d, s}(0) - \inf_{x \in \R^d} f^{\star}_{d, s}(x) \thesischfiveoversetlab{\eqref{ch5:70110a40-c2b6-4b8f-b790-e09561a73316_chapter_5}}{=} f(0;d, s) = \frac{L r^2}{2}, \thesischfivenumberthis\label{ch5:f38cccea-2eb2-49b1-910d-1eb1671cb5ff_chapter_5} \end{thesischfivefitdisplay}
        {\color{navy}Thus, it suffices to take $L=2\Delta/r^2<+\infty$ (since $r>0$)} to ensure the function $f_{d, s}^{\star}$ has $\Delta$--bounded initial sub-optimality. In the next step, we will show that $f_{d, s}^{\star}$ also has $\bar{L}$--Lipschitz gradients, as a consequence of the mean-squared smoothness property.

        \item \textbf{Step 2:} \textit{The Oracle Class and~Assumption~\ref{ch5:ass:p-bounded-central-moment-gradient_chapter_5} and Assumption~\ref{ch5:ass:mean-squared-smoothness_chapter_5}}.

        Now, let us compute the $p$--th central moment, for any $x \in \R^d$ and any $s \in \thesischfiveens{-1, 1}$ we have
            \begin{thesischfivefitalignatnum}{1}
                &\thesischfiveExpSub{\xi \sim \Prob^s}
                {\norm{\nabla f_d^{\star}(x,\xi)-\nabla f_{d,s}^{\star}(x)}^p}
                \notag\\[-0.2em]
                &\quad\thesischfiveoversetlab{\eqref{ch5:70110a40-c2b6-4b8f-b790-e09561a73316_chapter_5}}{=}
                {\color{navy}\chi}^p\,
                \thesischfiveExpSub{\xi\sim\Prob^s}
                {\norm{\nabla f_d({\color{navy}\chi}x,\xi)
                -\nabla f({\color{navy}\chi}x;d,s)}^p}
                \notag\\
                &\quad=
                {\color{navy}\chi}^p\,
                \thesischfiveExpSub{\xi\sim\Prob^s}
                {\norm{L({\color{navy}\chi}x-\boldsymbol{\xi})
                -L({\color{navy}\chi}x-\theta_s)}^p}
                \notag\\
                &\quad=
                \left(L{\color{navy}\chi}\right)^p
                \thesischfiveExpSub{\xi\sim\Prob^s}{\abs{\xi-rs}^p}
                \notag\\
                &\quad=
                \left(L{\color{navy}\chi}r\right)^p
                \thesischfiveExpSub{\xi\sim\Prob^s}
                {\abs{\frac{\xi}{r}-s}^p}
                \notag\\
                &\quad=
                \left(L{\color{navy}\chi}r\right)^p
                \left((1-r)\abs{s}^p+r\abs{\frac{s}{r}-s}^p\right)
                \notag\\
                &\quad=
                \left(L{\color{navy}\chi}r\right)^p
                \left((1-r)+r\abs{\frac1r-1}^p\right)
                \notag\\
                &\quad=
                \left(L{\color{navy}\chi}r\right)^p(1-r)
                \left(1+\left(\frac{1-r}{r}\right)^{p-1}\right)
                \notag\\
                &\quad=
                \left(L{\color{navy}\chi}\right)^p r(1-r)
                \left(r^{p-1}+(1-r)^{p-1}\right)
                \notag\\
                &\quad\le
                2\left(L{\color{navy}\chi}\right)^p r(1-r),
                \label{ch5:279a6448-b491-4172-8fb0-c0db4abd8e3f_chapter_5}
            \end{thesischfivefitalignatnum}
        where in the second equality we let $\boldsymbol{\xi} = (\xi, 0, \ldots, 0) \in \R^d$ and use the gradients of $f_d$ and $f_{d, s}$, given for any $x \in \R^d$ by
            \begin{thesischfivefitdisplay} \nabla f_d(x, \xi) = L (x - \boldsymbol{\xi}) \quad \text{and} \quad \nabla f(x;d, s) = L (x - \theta_s), \thesischfivenumberthis\label{ch5:2130180a-a3a9-4040-8a04-214e99e99948_chapter_5} \end{thesischfivefitdisplay}
        notably, $\thesischfiveExpSub{\xi \sim \Prob^s}{\nabla f_d^{\star}(x, \xi)} = \nabla f_{d, s}^{\star} (x)$. In the fifth equality we use that $\xi = s$ with probability $r$ and $\xi = 0$ otherwise, in the sixth we use $\abs{s} = 1$, and in the last inequality we use both $p > 1$ and $r \in \thesischfiveintof{0}{1}$ to bound $r^{p - 1} + \left( 1 - r \right)^{p - 1}$ by $2$. Hence, following~\eqref{ch5:279a6448-b491-4172-8fb0-c0db4abd8e3f_chapter_5} we need to guarantee
            \begin{thesischfivefitdisplay} 2 \left( L {\color{navy}\chi} \right)^p r ( 1 - r) \le \sigma_1^p, \quad \text{that is,} \quad 2^{1/p} L {\color{navy}\chi} \left( r ( 1 - r) \right)^{1/p} \le \sigma_1. \thesischfivenumberthis\label{ch5:d98b1956-8d19-4272-af9a-9a2d85b182c3_chapter_5} \end{thesischfivefitdisplay}

        Now, computing the mean-squared smoothness constant, we have, for $q \in \thesischfiveintof{1}{2}$, $s \in \thesischfiveens{-1, 1}$ and any $x, y \in \R^d$
            \begin{thesischfivefitalignat}{2}
                \MoveEqLeft[1] \thesischfiveExpSub{\xi \sim \Prob^s}{\norm{\nabla f_d^{\star}(x, \xi) - \nabla f_d^{\star}(y, \xi)}^q} \\
                \thesischfiveoversetlab{\eqref{ch5:70110a40-c2b6-4b8f-b790-e09561a73316_chapter_5}}&{=} {\color{navy}\chi}^q \, \thesischfiveExpSub{\xi \sim \Prob^s}{\norm{\nabla f_d({\color{navy}\chi} x, \xi) - \nabla f_d({\color{navy}\chi} y, \xi)}^q} \\
                \thesischfiveoversetlab{\eqref{ch5:2130180a-a3a9-4040-8a04-214e99e99948_chapter_5}}&{=} \left( L {\color{navy}\chi} \right)^q \, \thesischfiveExpSub{\xi \sim \Prob^s}{\norm{{\color{navy}\chi} x - {\color{navy}\chi} y}^q} \\
                & = \left( L {\color{navy}\chi}^2 \right)^q \norm{x - y}^q,
            \end{thesischfivefitalignat}
        and it suffices to ensure
            \begin{thesischfivefitdisplay} \left( L {\color{navy}\chi}^2 \right)^q \le \bar{L}^q, \quad \text{that is,} \quad L {\color{navy}\chi}^2 \le \bar{L}, \thesischfivenumberthis\label{ch5:63deeb98-3d9a-4bca-9cc0-9acaba35425e_chapter_5} \end{thesischfivefitdisplay}
        so as to satisfy~Assumption~\ref{ch5:ass:mean-squared-smoothness_chapter_5}.

        \item \textbf{Step 3:} \textit{Choice of ${\color{navy}\chi}$, $r$}.

        It remains to choose ${\color{navy}\chi} > 0$ and a universal constant $r \in \thesischfiveintof{0}{1}$ satisfying the two inequalities~\eqref{ch5:d98b1956-8d19-4272-af9a-9a2d85b182c3_chapter_5} and~\eqref{ch5:63deeb98-3d9a-4bca-9cc0-9acaba35425e_chapter_5}. Additionally, {\color{navy}let $c'>0$ be a universal constant, to be fixed later, such that $0<\thesischfiveeps\le c'\sqrt{\bar L\Delta}$.} From inequality~\eqref{ch5:63deeb98-3d9a-4bca-9cc0-9acaba35425e_chapter_5} we have
            \begin{thesischfivefitdisplay} {\color{navy}\chi} \le \sqrt{\frac{\bar{L}}{L}} = \frac{1}{L} \sqrt{\bar{L} L} = \frac{\sqrt{2 \bar{L} \Delta}}{L r}, \end{thesischfivefitdisplay}
        hence for any constant $c' > 0$, taking ${\color{navy}\chi} = \sqrt{2}\,\thesischfiveeps/(c' L r) > 0$ gives
            \begin{thesischfivefitdisplay} {\color{navy}\chi} = \frac{\thesischfiveeps \sqrt{2}}{c' L r} \le \frac{c' \sqrt{\bar{L} \Delta} \sqrt{2}}{c' L r} = \frac{\sqrt{2 \bar{L} \Delta}}{L r}, \end{thesischfivefitdisplay}
        as desired. Plugging back this value in the inequality~\eqref{ch5:d98b1956-8d19-4272-af9a-9a2d85b182c3_chapter_5} we need to have
            \begin{thesischfivefitdisplay} \frac{2^{1/p} \sqrt{2}}{c'} \thesischfiveeps r^{\frac{1}{p} - 1} (1 - r)^{1/p} \le \sigma_1, \quad \text{i.e.,} \quad \frac{2^{1/p} \sqrt{2}}{c'} \left( \frac{\thesischfiveeps}{\sigma_1} \right) (1 - r)^{1/p} \le r^{(p - 1)/p}, \end{thesischfivefitdisplay}
        which is equivalent to
            \begin{thesischfivefitdisplay} C^{p/(p - 1)} \left( \frac{\thesischfiveeps}{\sigma_1} \right)^{p/(p - 1)} (1 - r)^{1/(p - 1)} \le r, \end{thesischfivefitdisplay}
        where we set $C \eqdef 2^{1/p} \sqrt{2}/c'$. Thus, as $1 - r \le 1$, it is enough to take
            \begin{thesischfivefitdisplay} r = \min\thesischfiveens{1, \left( \frac{C \thesischfiveeps}{\sigma_1} \right)^{p/(p - 1)}} > 0, \end{thesischfivefitdisplay}
        and now all the inequalities are satisfied so do~Assumption~\ref{ch5:ass:lower-boundedness_chapter_5}, Assumption~\ref{ch5:ass:p-bounded-central-moment-gradient_chapter_5}, and Assumption~\ref{ch5:ass:mean-squared-smoothness_chapter_5} hold.

        \item \textbf{Step 4:} \textit{Transforming the Optimization Problem into a Function Identification Problem}.

        {\color{navy}We now follow the two-point reduction of~\textcite{arjevani2022lower}. Draw $S$ uniformly from $\{-1,1\}$ and let an arbitrary algorithm observe $N$ i.i.d.\ samples $\xi^{1},\ldots,\xi^{N}\sim\mathcal P^{S}$ before returning $\widehat x\in\R^d$. Define
        \begin{thesischfivefitdisplay}
            \widehat S\in\operatorname*{arg\,min}_{s\in\{-1,1\}}
            \norm{\nabla f^\star_{d,s}(\widehat x)},
        \end{thesischfivefitdisplay}
        breaking ties measurably. On the event $\{\widehat S\ne S\}$, the definition of $\widehat S$ and the triangle inequality give
        \begin{align*}
            2\norm{\nabla f^\star_{d,S}(\widehat x)}
            &\ge
            \norm{\nabla f^\star_{d,S}(\widehat x)}
            +\norm{\nabla f^\star_{d,\widehat S}(\widehat x)}\\
            &\ge L\chi
            \norm{\theta_S-\theta_{\widehat S}}
            =2L\chi r.
        \end{align*}
        Consequently,
        \begin{thesischfivefitdisplay}
            \thesischfiveE{\norm{\nabla f^\star_{d,S}(\widehat x)}}
            \ge L\chi r\,\thesischfiveProba{\widehat S\ne S}.
            \thesischfivenumberthis\label{ch5:174ddb60-9406-45e1-b077-a317743c46a2_chapter_5}
        \end{thesischfivefitdisplay}}

        \item \textbf{Step 5:} \textit{Lower Bounding the Misidentification Probability $\thesischfiveProba{\widehat{S} \neq S}$}.

        {\color{navy}For $s\in\{-1,1\}$, let $\mathcal P_N^s$ be the law of $(\xi^1,\ldots,\xi^N)$ conditional on $S=s$; equivalently, it is the product law of $N$ independent copies of $sB$, where $B\sim\operatorname{Ber}(r)$. Le Cam's two-point inequality yields, for every estimator $\widehat S$,
            \begin{thesischfivefitdisplay}
                \thesischfiveProba{\widehat S\ne S}
                \ge \frac12\left(1-\operatorname{TV}
                (\mathcal P_N^1,\mathcal P_N^{-1})\right).
            \end{thesischfivefitdisplay}
        To compute this total variation distance, observe that the supports of $\mathcal P_N^1$ and $\mathcal P_N^{-1}$ are $\{0,1\}^N$ and $\{-1,0\}^N$, respectively, and intersect only at $\mathbf0=(0,\ldots,0)$. Therefore}
            \begin{thesischfivefitalignat}{2}
                {\color{navy}\operatorname{TV}(\mathcal P_N^1,\mathcal P_N^{-1})}
                &{\color{navy}= \frac{1}{2} \sum_{x \in \thesischfiveens{-1, 0, 1}^N} \abs{\mathcal P_N^1(x) - \mathcal P_N^{-1} (x)}} \\
                &{\color{navy}= \frac{1}{2} \sum_{x \in \thesischfiveens{-1, 0}^N \setminus \thesischfiveens{\mathbf{0}}} \mathcal P_N^{-1}(x) + \frac{1}{2} \sum_{x \in \thesischfiveens{0, 1}^N \setminus \thesischfiveens{\mathbf{0}}} \mathcal P_N^1(x)} \\
                & = \frac{1}{2} \left( 1 - (1 - r)^N \right) + \frac{1}{2} \left( 1 - (1 - r)^N \right) \\
                & = 1 - (1 - r)^N,
            \end{thesischfivefitalignat}
        Hence,
            \begin{thesischfivefitdisplay} \thesischfiveProba{\widehat{S} \neq S} \ge \frac{1}{2} ( 1 - r)^N, \end{thesischfivefitdisplay}
        and plugging this bound into~\eqref{ch5:174ddb60-9406-45e1-b077-a317743c46a2_chapter_5} gives
            \begin{thesischfivefitdisplay} \thesischfiveE{\norm{\nabla f^{\star}_{d, S} \left( \widehat{x} \right)}} \ge \frac{L {\color{navy}\chi} r}{2} (1 - r)^N. \end{thesischfivefitdisplay}
        Also, note that
            \begin{thesischfivefitalignat}{2}
                \MoveEqLeft[1] \max\thesischfiveens{\thesischfiveE{\norm{\nabla f^{\star}_{d, 1} \left( \widehat{x} \right)}}, \thesischfiveE{\norm{\nabla f^{\star}_{d, -1} \left( \widehat{x} \right)}}} \\
                & \ge \frac{1}{2} \left( \thesischfiveE{\norm{\nabla f^{\star}_{d, 1} \left( \widehat{x} \right)}} + \thesischfiveE{\norm{\nabla f^{\star}_{d, -1} \left( \widehat{x} \right)}} \right) = \thesischfiveE{\norm{\nabla f^{\star}_{d, S} \left( \widehat{x} \right)}}.
            \end{thesischfivefitalignat}

        It remains to lower bound adequately the quantity $\frac{L {\color{navy}\chi} r}{2} (1 - r)^N$. For this, we consider two cases:
            \begin{itemize}
                \item if $\sigma_1 > C \thesischfiveeps$ then $C \thesischfiveeps/\sigma_1 < 1$, so $r = (C \thesischfiveeps/\sigma_1)^{p/(p - 1)} < 1$ and
                    \begin{thesischfivefitalignat}{2}
                        \frac{L {\color{navy}\chi} r}{2} (1 - r)^N & = \frac{\thesischfiveeps \sqrt{2}}{2 c'} (1 - r)^N,
                    \end{thesischfivefitalignat}
                and,
                \begin{itemize}
                    \item assume we have $N \le 1/(2 r)$ then
                        \begin{thesischfivefitdisplay} {\color{navy}(1 - r)^N \ge 1-rN \ge \frac{1}{2}}, \end{thesischfivefitdisplay}
                    {\color{navy}where the first inequality follows from Bernoulli's inequality (Lemma~\ref{ch5:lem:bernoulli-inequality_chapter_5}). Thus,}
                        \begin{thesischfivefitdisplay} \max\thesischfiveens{\thesischfiveE{\norm{\nabla f^{\star}_{d, 1} \left( \widehat{x} \right)}}, \thesischfiveE{\norm{\nabla f^{\star}_{d, -1} \left( \widehat{x} \right)}}} \ge \frac{\thesischfiveeps \sqrt{2}}{4 c'}, \end{thesischfivefitdisplay}
                    {\color{navy}and it suffices to take $c'=\frac18$ to ensure}
                        \begin{thesischfivefitdisplay} \max\thesischfiveens{\thesischfiveE{\norm{\nabla f^{\star}_{d, 1} \left( \widehat{x} \right)}}, \thesischfiveE{\norm{\nabla f^{\star}_{d, -1} \left( \widehat{x} \right)}}} \ge 2 \thesischfiveeps \sqrt{2} > 2\thesischfiveeps. \end{thesischfivefitdisplay}
                    So, when $N \le 1/(2 r)$ it is not possible to reach an $\thesischfiveeps$--stationary point on both $f^{\star}_{d, 1}$ and on $f^{\star}_{d, -1}$.
                    
                    \item Hence, we deduce that we must have 
                        \begin{thesischfivefitdisplay} N > \frac{1}{2 r} = \frac{1}{2} \left( \frac{\sigma_1}{C \thesischfiveeps} \right)^{p/(p - 1)} = \Omega(1) \cdot \left( \frac{\sigma_1}{\thesischfiveeps} \right)^{p/(p - 1)} \ge \Omega(1), \end{thesischfivefitdisplay}
                    as desired,
                \end{itemize}

            \item if $\sigma_1 \le C \thesischfiveeps$ then $r = 1$ and
                \begin{thesischfivefitdisplay} \left( \frac{\sigma_1}{C \thesischfiveeps} \right)^{p/(p - 1)} \le 1, \end{thesischfivefitdisplay}
            and since for any $s \in \thesischfiveens{-1, 1}$ we have 
                \begin{thesischfivefitdisplay} \norm{\nabla f^{\star}_{d, s} (0)} \thesischfiveoversetlab{\eqref{ch5:2130180a-a3a9-4040-8a04-214e99e99948_chapter_5}}{=} L {\color{navy}\chi} r = \frac{\thesischfiveeps \sqrt{2}}{c'} = 8 \thesischfiveeps \sqrt{2} > 2 \thesischfiveeps, \end{thesischfivefitdisplay}
            so at least one query is required to reach an $\thesischfiveeps$--stationary point, i.e., $N \ge 1$ from where 
                \begin{thesischfivefitdisplay} N \ge 1 \ge \left( \frac{\sigma_1}{C \thesischfiveeps} \right)^{p/(p - 1)} = \Omega(1) \cdot \left( \frac{\sigma_1}{\thesischfiveeps} \right)^{p/(p - 1)}, \end{thesischfivefitdisplay}
            {\color{navy}which completes the proof of the lemma.}
            \end{itemize}
    \end{itemize}
\end{proof}

\begin{remark}\label{ch5:appdx-rem:lower-bound-global-stochastic-model-motivations_chapter_5}
    Let us expand on the proof strategy of the above lemma. First, we start by defining a stochastic function $f \colon (x, \xi) \mapsto \frac{L}{2} \left( \sqnorm{x} - 2 x_1 \xi + r^2 \right)$ and we then provide two probability distributions $\Prob^-$ and $\Prob^+$ for the random variable $\xi$ which allows to define two deterministic functions
        \begin{thesischfivefitdisplay} f_- \colon x \mapsto \thesischfiveExpSub{\xi \sim \Prob_-}{f(x;\xi)} \quad \text{and} \quad f_+ \colon x \mapsto \thesischfiveExpSub{\xi \sim \Prob_+}{f(x;\xi)}. \end{thesischfivefitdisplay}
    {\color{navy}The strategy is to randomize the initial function by choosing $f_-$ or $f_+$ with probability $\frac12$ and then lower-bound the gradient-norm error through the probability of misidentification. Equivalently, we ask how many samples $\xi^{(1)},\ldots,\xi^{(N)}$ are needed, in expectation, to distinguish $\Prob^-$ from $\Prob^+$.}

    Contrary to~\textcite{arjevani2022lower}, instead of using normal distributions we use two Bernoulli distributions $\thesischfiveBer(r)$ and $-\thesischfiveBer(r)$ where, after tuning the parameters so as to satisfy the different assumptions, we obtain
        \begin{thesischfivefitdisplay} r = \Theta\left( \min\thesischfiveens{1, \left( \frac{\thesischfiveeps}{\sigma_1} \right)^{p/(p - 1)}} \right). \end{thesischfivefitdisplay}
    The choice of Bernoulli distributions is actually inspired from the definition of the gradient estimator in the work of~\textcite{arjevani2022lower} where a Bernoulli distribution is used to probabilistically hide new coordinates (see~Definition~\ref{ch5:appdx-def:gradient-estimator_chapter_5}).
    
    {\color{navy}Under these Bernoulli distributions, the oracle returns zeros until it first returns $-1$ or $+1$. A zero does not distinguish $\thesischfiveBer(r)$ from $-\thesischfiveBer(r)$, whereas a nonzero observation identifies the distribution immediately. Since a nonzero observation occurs with probability $r$, distinguishing $\Prob^-$ from $\Prob^+$ requires, in expectation,
    $1/r=\Omega\!\left(1+\left(\sigma_1/\thesischfiveeps\right)^{p/(p-1)}\right)$
    oracle calls, which gives the claimed lower bound.}
\end{remark}

\begin{remark}
    Combining~\textcite[Lemma~11]{arjevani2022lower} and the above lemma, we deduce that under the general $p$--bounded central moment where $p > 1$ is any real number, we have
        \begin{thesischfivefitdisplay} \mathfrak{m}^{\normalfont\texttt{zr}}_{\thesischfiveeps} \left( {\color{navy}b}, \Delta, \bar{L}, \sigma_1^p \right) \ge \Omega(1) \cdot \left( \frac{\sigma_1}{\thesischfiveeps} \right)^{\max\{p/(p - 1),\, 2\}}, \thesischfivenumberthis\label{ch5:4f3b2ba5-4a06-44f5-b589-94a0605dc4ab_chapter_5} \end{thesischfivefitdisplay}
    {\color{navy}where the hidden constant in $\Omega(1)$ may depend on $p$. Hence,}
        \begin{itemize}
            \item {\color{navy}a bounded first central moment ($p=1$) is insufficient to guarantee an $\thesischfiveeps$--stationary point after any finite number $N$ of oracle queries},
            
            \item assuming high-order bounded moments ($p > 2$) does not bring additional information than in the bounded variance case ($p = 2$).
        \end{itemize}
\end{remark}

\begin{proof}[Proof of~\eqref{ch5:4f3b2ba5-4a06-44f5-b589-94a0605dc4ab_chapter_5}]
    The proof here applies not only to the setting of~Lemma~\ref{ch5:appdx-lem:lower-bound-global-stichastic-model-mean-squared-smoothness-p-BCM_chapter_5}, i.e.,~Assumption~\ref{ch5:ass:lower-boundedness_chapter_5}, Assumption~\ref{ch5:ass:p-bounded-central-moment-gradient_chapter_5}, and Assumption~\ref{ch5:ass:mean-squared-smoothness_chapter_5} but also to the settings discussed in~Remark~\ref{ch5:appdx-rem:improve-lower-bound-global-stichastic-model-mean-squared-smoothness-2-p-BCM_chapter_5} and Remark~\ref{ch5:appdx-rem:improve-lower-bound-global-stichastic-model-mean-squared-smoothness-2-hessian_chapter_5}.

    First, the case of an exponent $p \in \thesischfiveintof{1}{2}$ is already covered by~Lemma~\ref{ch5:appdx-lem:lower-bound-global-stichastic-model-mean-squared-smoothness-p-BCM_chapter_5}. For the case when $p > 2$, we reuse the exact same proof as in~\textcite[Lemma~11]{arjevani2022lower}. More precisely, instead of considering Bernoulli distributions for $\Prob^{-1}$ and $\Prob^1$, we fall back to normal distributions, that is, $\Prob^s = \mathcal{N}(r s, \sigma_1^2/c_p^2)$ for $s \in \thesischfiveens{-1, 1}$ and some well-chosen constant $c_p > 0$ (depending on the exponent $p$ and $\bar{L}$) that we fix later. Taking ${\color{navy}\chi} = 1$ and $L = \bar{L}$, as in~\textcite{arjevani2022lower}, it remains to bound the $p^{\textnormal{th}}$ central moment of the stochastic gradient, we have for any $x \in \R^d$ and any $s \in \thesischfiveens{-1, 1}$, by~\eqref{ch5:279a6448-b491-4172-8fb0-c0db4abd8e3f_chapter_5}
            \begin{thesischfivefitalignat}{2}
                \thesischfiveExpSub{\xi \sim \Prob^s}{\norm{\nabla f_d^{\star}(x, \xi) - \nabla f_{d, s}^{\star}(x)}^p} \thesischfiveoversetlab{\eqref{ch5:279a6448-b491-4172-8fb0-c0db4abd8e3f_chapter_5}}&{=} \bar{L}^p \thesischfiveExpSub{\xi \sim \Prob^s}{\abs{\xi - r s}^p} \\
                & = \bar{L}^p \left( \frac{\sigma_1 \sqrt{2}}{c_p} \right)^p \frac{\Gamma\left( \frac{p + 1}{2} \right)}{\sqrt{\pi}} \\
                & \le \sigma_1^p,
            \end{thesischfivefitalignat}
        where the second equality follows from Lemma~\ref{ch5:appdx-lem:absolute-central-moment-gaussian_chapter_5}, and in the last inequality we took
            \begin{thesischfivefitdisplay} c_p = \bar{L} \sqrt{2} \left( \frac{\Gamma\left( \frac{p + 1}{2} \right)}{\sqrt{\pi}} \right)^{1/p} = C_p \bar{L} > 0, \end{thesischfivefitdisplay}
        where $C_p = \sqrt{2}\, \bigl( \Gamma((p + 1)/2)/\sqrt{\pi} \bigr)^{1/p}$. The rest of the proof is similar to~\textcite[Lemma~10]{arjevani2022lower}.        
\end{proof}

\subsection{\texorpdfstring{{\color{navy}Some Bounds on the Gradient Estimator $\bar{g}_N$}}{Some Bounds on the Gradient Estimator g-bar-N}}

\begin{lemma}[Properties of the Gradient Estimator $\bar{g}_N$]\label{ch5:appdx-lem:properties-gradient-estimator-1_chapter_5}
    The stochastic gradient estimator $\bar{g}_N$ is a probability--$\theta$ zero-chain, is unbiased with respect to $\nabla f_N$ and satisfies
        \begin{thesischfivefitalignatnum}{2}
            \thesischfiveE{\norm{\bar{g}_N(x, \xi) - \nabla f_N(x)}^p} & \le \frac{2 \gamma_{\infty}^p (1 - \theta)}{\theta^{p - 1}}, \quad \text{and} \\
            \thesischfiveE{\norm{\bar{g}_N(x, \xi) - \bar{g}_N(y, \xi)}^q} & \le \frac{\bar{\ell}_1^q}{\theta^{q - 1}} \norm{x - y}^q, \label{ch5:9ff1cee4-5c5c-4c5b-854f-445bd4e1539d_chapter_5}
        \end{thesischfivefitalignatnum}
    and
        \begin{thesischfivefitdisplay} {\color{navy}\thesischfiveE{\norm{\left[ \bar{g}_N(x, \xi) - \bar{g}_N(y, \xi) \right] - \left[ \nabla f_N(x) - \nabla f_N(y) \right]}^q} \le \frac{\bar{\delta_1}^q (1 - \theta)}{\theta^{q - 1}} \norm{x - y}^q}, \thesischfivenumberthis\label{ch5:5ebfcdd7-34dd-4433-bb48-69a29f3c7dab_chapter_5} \end{thesischfivefitdisplay}
    for all $x, y \in \R^N$, where $p, q \in \thesischfiveintof{1}{2}$, $\gamma_{\infty}$ is defined in~Lemma~\ref{ch5:appdx-lem:properties-hard-instance_chapter_5}, $\bar{\ell}_1 \eqdef 2 \left( 2^q \left( 6^{2 q} \gamma_{\infty}^q + \ell_1^q \right) + \ell_1^q  \right)^{1/q} \ge \ell_1$ and $\bar{\delta}_1 \eqdef 4 \left( 6^{2 q} \gamma_{\infty}^q + \ell_1^q \right)^{1/q}$.
\end{lemma}

The proof of~Lemma~\ref{ch5:appdx-lem:properties-gradient-estimator-1_chapter_5} is largely inspired by~\textcite[Lemma~4]{arjevani2022lower}. Notably, we rely on~Lemma~\ref{ch5:lem:norm-power-alpha-inequality_chapter_5} which states the following  inequality
    \begin{thesischfivefitdisplay} \norm{a + b}^{\alpha} \le 2^{\alpha - 1} \left( \norm{a}^{\alpha} + \norm{b}^{\alpha} \right), \thesischfivenumberthis\label{ch5:12e46e5c-0bf2-4de8-ac73-920ea999660b_chapter_5} \end{thesischfivefitdisplay}
holds, for any vectors $a, b \in \R^d$ and any exponent $\alpha \ge 1$. This generalizes the well-known inequality $\sqnorm{a + b} \le 2( \sqnorm{a} + \sqnorm{b})$. We use inequality~\eqref{ch5:12e46e5c-0bf2-4de8-ac73-920ea999660b_chapter_5} in~\eqref{ch5:af4fdc92-39df-4ba0-b39e-528947561b51_chapter_5} as a \say{substitute} for the squared norm expansion (as used in~\textcite{arjevani2022lower}) which does not hold anymore with exponent $q \in \thesischfiveintof{1}{2}$ instead of $2$.

\begin{proof}
    Following the proof of~\textcite[Lemma~4]{arjevani2022lower}, for any $\xi$, the vector $\delta(x, \xi) \eqdef \bar{g}_N(x, \xi) - \nabla f_N(x)$ has at most one nonzero entry at coordinate $i_x = \thesischfiveprog_{\frac{1}{4}}(x) + 1$. Moreover, for any $i \in \thesischfiveInt{1}{N}$, the $i^{\textnormal{th}}$ entry $\delta_i(x, \xi)$ of $\delta(x, \xi)$ reads
        \begin{thesischfivefitdisplay} \delta_i(x, \xi) = \left[ \nabla f_N(x) \right]_i \Theta_i(x) \left( \frac{\xi}{\theta} - 1 \right), \end{thesischfivefitdisplay}
    where the function $\Theta_i$ is defined in~Definition~\ref{ch5:appdx-def:smoothed-gradient-estimator_chapter_5}. Hence, we have
        \begin{thesischfivefitalignatnum}{2}
            \thesischfiveE{\norm{\bar{g}_N(x, \xi) - \nabla f_N(x)}^p} & = \thesischfiveE{\abs{\delta_{i_x}(x, \xi)}^p} \\
            & = \abs{ \left[ \nabla f_N(x) \right]_{i_x}}^p \abs{\Theta_{i_x}(x)}^p \thesischfiveE{\abs{\frac{\xi}{\theta} - 1}^p} \\
            \thesischfiveoversetref{Lem.}{\ref{ch5:appdx-lem:properties-hard-instance_chapter_5}}&{\le} \gamma_{\infty}^p \abs{\Theta_{i_x}(x)}^p \thesischfiveE{\abs{\frac{\xi}{\theta} - 1}^p} \\
            & \le \gamma_{\infty}^p \, \thesischfiveE{\abs{\frac{\xi}{\theta} - 1}^p} \\
            \thesischfiveoversetrel{rel:91d78f4b-c07c-4b5e-974c-cea7451e140d}&{=} \gamma_{\infty}^p \left( (1 - \theta) + \theta \left( \frac{1}{\theta} - 1 \right)^p \right) \\
            & = \gamma_{\infty}^p (1 - \theta) \left( 1 + \left( \frac{1 - \theta}{\theta} \right)^{p - 1} \right) \\
            \thesischfiveoversetrel{rel:51b68383-0e41-40f2-a4b7-c07c7d81d151}&{\le} \frac{2 \gamma_{\infty}^p (1 - \theta)}{\theta^{p - 1}}, \label{ch5:717d22a1-a3cf-4c72-9b37-af698eef534b_chapter_5}
        \end{thesischfivefitalignatnum}
    where in~\thesischfiverelref{rel:91d78f4b-c07c-4b5e-974c-cea7451e140d} we use the fact that $\xi \sim \thesischfiveBer(\theta)$. In~\thesischfiverelref{rel:51b68383-0e41-40f2-a4b7-c07c7d81d151} we use $\theta \in \thesischfiveintff{0}{1}$ and $p \ge 1$ to bound $(1 - \theta)^{p - 1} \le 1$ and $1 \le \theta^{- (p - 1)}$. This establishes the first inequality from~\eqref{ch5:9ff1cee4-5c5c-4c5b-854f-445bd4e1539d_chapter_5}.

    Now, for the second inequality in~\eqref{ch5:9ff1cee4-5c5c-4c5b-854f-445bd4e1539d_chapter_5}, we use the fact that $\delta(\cdot, \xi) \eqdef \bar{g}_N(\cdot, \xi) - \nabla f_N(\cdot)$ has at most one nonzero coordinate, then if we let $i_x = \thesischfiveprog_{\frac{1}{4}}(x) + 1$ and $i_y = \thesischfiveprog_{\frac{1}{4}}(y) + 1$ we obtain
        \begin{thesischfivefitalignatnum}{2}
            \MoveEqLeft[1] \thesischfiveE{\norm{\bar{g}_N(x, \xi) - \bar{g}_N(y, \xi)}^q}\\
            & = \thesischfiveE{\norm{\left[ \delta(x, \xi) - \delta(y, \xi) \right] + \left[ \nabla f_N(x) - \nabla f_N(y) \right]}^q} \\
            & \le 2^{q - 1} \left( \thesischfiveE{\norm{\delta(x, \xi) - \delta(y, \xi)}^q} + \norm{\nabla f_N(x) - \nabla f_N(y)}^q \right) \\
            & \le 2^{q - 1} \Bigl( \thesischfiveE{\abs{\delta_{i_x}(x, \xi) - \delta_{i_x}(y, \xi)}^q} + \thesischfiveE{\abs{\delta_{i_y}(x, \xi) - \delta_{i_y}(y, \xi)}^q}\\
            &\qquad + \norm{\nabla f_N(x) - \nabla f_N(y)}^q \Bigr), \label{ch5:af4fdc92-39df-4ba0-b39e-528947561b51_chapter_5}
        \end{thesischfivefitalignatnum}
    where the first inequality follows from Lemma~\ref{ch5:lem:norm-power-alpha-inequality_chapter_5}. Moreover, for any integer $i \in \thesischfiveInt{1}{N}$ we have
        \begin{alignat}{2}
            \MoveEqLeft[1] \thesischfiveE{\abs{\delta_i(x, \xi) - \delta_i(y, \xi)}^q} \notag\\
            & = \abs{\left[ \nabla f_N(x) \right]_i \Theta_i(x) - \left[ \nabla f_N(y) \right]_i \Theta_i(y)}^q \thesischfiveE{\abs{\frac{\xi}{\theta} - 1}^q} \notag\\
            & \le \abs{\left[ \nabla f_N(x) \right]_i \Theta_i(x) - \left[ \nabla f_N(y) \right]_i \Theta_i(y)}^q \left( \frac{2}{\theta^{q - 1}} \right) \notag\\
            & = \abs{\left[ \nabla f_N(x) \right]_i \left( \Theta_i(x) - \Theta_i(y) \right) + \left( \left[ \nabla f_N(x) \right]_i - \left[ \nabla f_N(y) \right]_i \right) \Theta_i(y)}^q \left( \frac{2}{\theta^{q - 1}} \right) \notag\\
            & \le 2^{q - 1} \Bigl( \abs{\left[ \nabla f_N(x) \right]_i}^q \abs{\Theta_i(x) - \Theta_i(y)}^q \notag\\
            &\qquad + \abs{\left[ \nabla f_N(x) \right]_i - \left[ \nabla f_N(y) \right]_i}^q \abs{\Theta_i(y)}^q \Bigr) \left( \frac{2}{\theta^{q - 1}} \right) \notag\\
            & \le 2^{q - 1} \left( 6^{2 q} \abs{\left[ \nabla f_N(x) \right]_i}^q \norm{x - y}^q + \abs{\left[ \nabla f_N(x) \right]_i - \left[ \nabla f_N(y) \right]_i}^q \right) \left( \frac{2}{\theta^{q - 1}} \right) \notag\\
            & \le 2^{q - 1} \left( 6^{2 q} \gamma_{\infty}^q \norm{x - y}^q + \norm{\nabla f_N(x) - \nabla f_N(y)}^q \right) \left( \frac{2}{\theta^{q - 1}} \right), \label{ch5:78a8d3c0-b886-4e27-965b-695544cbd732_chapter_5}
        \end{alignat}
    where the first inequality uses our previous bound from~\eqref{ch5:717d22a1-a3cf-4c72-9b37-af698eef534b_chapter_5} (which we derived with exponent $p$ instead of $q$), the second follows from Lemma~\ref{ch5:lem:norm-power-alpha-inequality_chapter_5}, the third uses the fact that $\Theta_i$ is $6^2$--Lipschitz \parencite[see][]{arjevani2022lower} and $\abs{\Theta_i(\cdot)} \le 1$, and the last follows from Lemma~\ref{ch5:appdx-lem:properties-hard-instance_chapter_5}. Now, plugging back the bound~\eqref{ch5:78a8d3c0-b886-4e27-965b-695544cbd732_chapter_5} in~\eqref{ch5:af4fdc92-39df-4ba0-b39e-528947561b51_chapter_5} we obtain
        \begin{thesischfivefitalignat}{2}
            \MoveEqLeft[1] \thesischfiveE{\norm{\bar{g}_N(x, \xi) - \bar{g}_N(y, \xi)}^q}\\
            \thesischfiveoversetlab{\eqref{ch5:78a8d3c0-b886-4e27-965b-695544cbd732_chapter_5}+\eqref{ch5:af4fdc92-39df-4ba0-b39e-528947561b51_chapter_5}}&{\le} 2^{q - 1} \Bigl( 2^q \left( 6^{2 q} \gamma_{\infty}^q \norm{x - y}^q + \norm{\nabla f_N(x) - \nabla f_N(y)}^q \right) \left( \frac{2}{\theta^{q - 1}} \right)\\
            &\qquad + \norm{\nabla f_N(x) - \nabla f_N(y) }^q \Bigr) \\
            & \le 2^{q - 1} \left( 2^q \left( 6^{2 q} \gamma_{\infty}^q + \ell_1^q \right) \left( \frac{2}{\theta^{q - 1}} \right) + \ell_1^q \right) \norm{x - y}^q \\
            & \le 2^q \left( 2^q \left( 6^{2 q} \gamma_{\infty}^q + \ell_1^q \right) + \ell_1^q \right) \frac{\norm{x - y}^q}{\theta^{q - 1}},
        \end{thesischfivefitalignat}
    where the second inequality follows from Lemma~\ref{ch5:appdx-lem:properties-hard-instance_chapter_5}, and in the last one we use $1 \le 2/\theta^{q - 1}$ to factor it out. {\color{navy}Define}
        \begin{thesischfivefitdisplay}
            {\color{navy}\bar{\ell}_1 \eqdef
            2 \left(2^q\left(6^{2q}\gamma_{\infty}^q+\ell_1^q\right)
            +\ell_1^q\right)^{1/q}.}
        \end{thesischfivefitdisplay}
    {\color{navy}Then}
        \begin{thesischfivefitdisplay} \thesischfiveE{\norm{\bar{g}_N(x, \xi) - \bar{g}_N(y, \xi)}^q} \le \frac{\bar{\ell}_1^q}{\theta^{q - 1}} \norm{x - y}^q, \end{thesischfivefitdisplay}
    as desired.

    It remains to establish the third inequality~\eqref{ch5:5ebfcdd7-34dd-4433-bb48-69a29f3c7dab_chapter_5}. Combining the bounds~\eqref{ch5:af4fdc92-39df-4ba0-b39e-528947561b51_chapter_5} and~\eqref{ch5:78a8d3c0-b886-4e27-965b-695544cbd732_chapter_5} we have
        \begin{thesischfivefitalignat}{2}
            \MoveEqLeft[1] \thesischfiveE{\norm{\left[ \bar{g}_N(x, \xi) - \bar{g}_N(y, \xi) \right] - \left[ \nabla f_N(x) - \nabla f_N(y) \right]}^q}\\
            & = \thesischfiveE{\norm{\delta(x, \xi) - \delta(y, \xi)}^q} \\
            & \le 2^{q - 1} \left( \thesischfiveE{\abs{\delta_{i_x}(x, \xi) - \delta_{i_x}(y, \xi)}^q} + \thesischfiveE{\abs{\delta_{i_y}(x, \xi) - \delta_{i_y}(y, \xi)}^q} \right) \\
            \thesischfiveoversetlab{\eqref{ch5:717d22a1-a3cf-4c72-9b37-af698eef534b_chapter_5}+\eqref{ch5:78a8d3c0-b886-4e27-965b-695544cbd732_chapter_5}}&{\le} 4^q \left( 6^{2 q} \gamma_{\infty}^q + \ell_1^q \right) \frac{1 - \theta}{\theta^{q - 1}} \norm{x - y}^q \\
            & {\color{navy}= \frac{\bar{\delta_1}^q (1 - \theta)}{\theta^{q - 1}} \norm{x - y}^q},
        \end{thesischfivefitalignat}
    where the first inequality follows from Lemma~\ref{ch5:lem:norm-power-alpha-inequality_chapter_5}, and we define
        \begin{thesischfivefitdisplay}
            \bar{\delta}_1 \eqdef 4 \left( 6^{2 q} \gamma_{\infty}^q + \ell_1^q \right)^{1/q}.
        \end{thesischfivefitdisplay}
    {\color{navy}This completes the proof.}
\end{proof}

\subsection{Proof of~Theorem~\ref{ch5:thm:lower-bound-sofo-mean-squared-smoothness_chapter_5}}\label{ch5:appdx-subsec:proof-lower-bound-sofo-mean-squared-smoothness_chapter_5}

\begin{theorem}[Restatement of Theorem~\ref{ch5:thm:lower-bound-sofo-mean-squared-smoothness_chapter_5}]
     Given $\Delta, \bar{L} > 0$, $\sigma_1 \ge 0$ and $0 < \thesischfiveeps \le c_1 \sqrt{\bar{L} \Delta}$ for some universal constant $c_1 > 0$. Then, for any algorithm $A \in \mathcal{A}_{\texttt{zr}}$, there exists a function $f \in \mathcal{F}\left( \Delta \right)$, an oracle and a distribution $(O, \thesischfivecD) \in \mathcal{O}\left( f, \bar{L}^q, \sigma_1^p \right)$ satisfying~Assumption~\ref{ch5:ass:lower-boundedness_chapter_5}, Assumption~\ref{ch5:ass:mean-squared-smoothness_chapter_5}, and Assumption~\ref{ch5:ass:p-bounded-central-moment-gradient_chapter_5} such that

    \begin{thesischfivefitalignat}{2}
        \mathfrak{m}^{\normalfont\texttt{zr}}_{\thesischfiveeps}\left( {\color{navy}b}, \Delta, \bar{L}^q, \sigma^p_1 \right) \ge \Omega(1) \cdot \left( \left( \frac{\sigma_1}{\thesischfiveeps} \right)^{p/(p - 1)} + \frac{\bar{L} \Delta}{\thesischfiveeps^2} + \frac{\bar{L} \Delta}{\thesischfiveeps^2} \left( \frac{\sigma_1}{\thesischfiveeps} \right)^{p/(q (p - 1))} \right).
    \end{thesischfivefitalignat}
\end{theorem}

\begin{proof}
    Let $\Delta_0$, $\ell_1$, $\gamma_{\infty}$ and $\bar{\ell}_1$ be the numerical constants in~Lemma~\ref{ch5:appdx-lem:properties-hard-instance_chapter_5} and Lemma~\ref{ch5:appdx-lem:properties-gradient-estimator-1_chapter_5} respectively. Additionally, we let the accuracy parameter $\thesischfiveeps > 0$, initial sub-optimality $\Delta \ge 0$, the $q$--weak average smoothness parameter $\bar{L}$, and the variance parameter $\sigma_1 \ge 0$ be fixed, and $0 < L \le \bar{L}$ to be specified later. Then, for $\alpha, {\color{navy}\chi} > 0$ two positive real numbers, following~\textcite{arjevani2022lower}, we rescale the function $f_N$ as
        \begin{thesischfivefitdisplay} f^{\star}_N \colon x \mapsto \alpha f_N({\color{navy}\chi} x). \thesischfivenumberthis\label{ch5:6f55ce68-cc03-4d22-a262-4aa3367f948d_chapter_5} \end{thesischfivefitdisplay}

    \begin{itemize}
        \item \textbf{Step 1:} \textit{Ensuring $f^{\star}_N \in \mathcal{F} \left(\Delta, L \right)$}.

        To guarantee the rescaled function $f^{\star}_N$ belongs to the function class $\mathcal{F}(\Delta, L)$, let us compute the initial sub-optimality $\Delta$ and the smoothness constant $L$. Assuming the algorithm \texttt{A} starts at $x^0 = 0$ we have
            \begin{thesischfivefitdisplay} f^{\star}_N(0) - \inf_{x \in \R^N} f^{\star}_N(x) \thesischfiveoversetlab{\eqref{ch5:6f55ce68-cc03-4d22-a262-4aa3367f948d_chapter_5}}{=} \alpha \left( f_N(0) - \inf_{x \in \R^N} f_N(x) \right) \le \alpha \Delta_0 N, \end{thesischfivefitdisplay}
        where the inequality follows from Lemma~\ref{ch5:appdx-lem:properties-hard-instance_chapter_5}; {\color{navy}thus, it suffices to take $N = \thesischfiveFloor{\Delta/(\alpha \Delta_0)}$ to ensure $f^{\star}_N(0) - \inf_{x \in \R^N} f^{\star}_N(x) \le \Delta$.} Moreover, for any $x, y \in \R^N$,
            \begin{thesischfivefitalignat}{2}
                \norm{\nabla f^{\star}_N(x) - \nabla f^{\star}_N(y)} & = \alpha {\color{navy}\chi} \norm{\nabla f_N({\color{navy}\chi} x) - \nabla f_N({\color{navy}\chi} y)} \\
                & \le \alpha {\color{navy}\chi} \ell_1 \norm{{\color{navy}\chi} x - {\color{navy}\chi} y} \\
                & = \alpha {\color{navy}\chi}^2 \ell_1 \norm{x - y},
            \end{thesischfivefitalignat}
        where the inequality follows from Lemma~\ref{ch5:appdx-lem:properties-hard-instance_chapter_5}, and it suffices to take $\alpha = L/({\color{navy}\chi}^2 \ell_1) > 0$ to ensure the function $f^{\star}_N$ has $L$--Lipschitz gradients. Consequently, we have $f^{\star}_N \in \mathcal{F}(\Delta, L)$, as desired.

        \item \textbf{Step 2:} \textit{Analysis of the Protocol and Choice for ${\color{navy}\chi}$}.

        Following the proof of~\textcite[Theorem~1]{arjevani2022lower}, according to~Lemma~\ref{ch5:appdx-lem:properties-hard-instance_chapter_5}, for all points $x \in \R^N$ such that $\thesischfiveprog_0(x) < N$ we have $\thesischfiveprog_0({\color{navy}\chi} x) = \thesischfiveprog_0(x) < N$ so
            \begin{thesischfivefitdisplay} \norm{\nabla f^{\star}_N(x)} \thesischfiveoversetlab{\eqref{ch5:6f55ce68-cc03-4d22-a262-4aa3367f948d_chapter_5}}{=} \frac{L}{\ell_1 {\color{navy}\chi}} \norm{\nabla f_N({\color{navy}\chi} x)} > \frac{L}{\ell_1 {\color{navy}\chi}}, \thesischfivenumberthis\label{ch5:191cd1ea-8e5f-4b85-8e63-f31fd25f3b17_chapter_5} \end{thesischfivefitdisplay}
        where the inequality follows from Lemma~\ref{ch5:appdx-lem:properties-hard-instance_chapter_5}, and we need to guarantee that
            \begin{thesischfivefitdisplay} \norm{\nabla f^{\star}_N(x)} > 2\thesischfiveeps, \end{thesischfivefitdisplay}
        for all $x \in \R^N$ with $\thesischfiveprog_0(x) < N$ which, given~\eqref{ch5:191cd1ea-8e5f-4b85-8e63-f31fd25f3b17_chapter_5}, can be done if we set ${\color{navy}\chi} = L/(2 \ell_1 \thesischfiveeps)$.

        \item \textbf{Step 3:} \textit{The Oracle Class and~Assumption~\ref{ch5:ass:p-bounded-central-moment-gradient_chapter_5} and Assumption~\ref{ch5:ass:mean-squared-smoothness_chapter_5}}.

        It remains to choose the parameter $\theta \in \thesischfiveintof{0}{1}$ and constant $L$ such that the gradient estimator $\bar{g}^{\star}_N$ of $\nabla f^{\star}_N$ satisfies~Assumption~\ref{ch5:ass:p-bounded-central-moment-gradient_chapter_5} and Assumption~\ref{ch5:ass:mean-squared-smoothness_chapter_5}. Computing the $p$--th central moment of $\bar{g}^{\star}_N$ gives, for all $x \in \R^N$
            \begin{thesischfivefitalignatnum}{2}
                \thesischfiveE{\norm{\bar{g}^{\star}_N(x, \xi) - \nabla f^{\star}_N(x)}^p} \thesischfiveoversetlab{\eqref{ch5:6f55ce68-cc03-4d22-a262-4aa3367f948d_chapter_5}}&{=} \left( \alpha {\color{navy}\chi} \right)^p \thesischfiveE{\norm{\bar{g}_N({\color{navy}\chi} x, \xi) - \nabla f_N({\color{navy}\chi} x)}^p} \\
                & \le \frac{2 \left( \gamma_{\infty} \alpha {\color{navy}\chi} \right)^p (1 - \theta)}{\theta^{p - 1}} \\
                & \le \frac{\left( 2 \gamma_{\infty} \alpha {\color{navy}\chi} \right)^p}{\theta^{p - 1}}, \label{ch5:29bb8bb5-e6e5-482d-ab46-fab6357955ec_chapter_5}
            \end{thesischfivefitalignatnum}
        where the first inequality follows from Lemma~\ref{ch5:appdx-lem:properties-gradient-estimator-1_chapter_5}, and in the second we use the fact that $0 < \theta \le 1$ and $p \ge 1$ so that $2 \le 2^p$. {\color{navy}By~\eqref{ch5:29bb8bb5-e6e5-482d-ab46-fab6357955ec_chapter_5}, Assumption~\ref{ch5:ass:p-bounded-central-moment-gradient_chapter_5} is satisfied by taking $\theta = \min\thesischfiveens{1, \bar{\theta}}$, where}
        
            \begin{thesischfivefitdisplay} \bar{\theta}^{p - 1} \ge \left( \frac{2 \gamma_{\infty} \alpha {\color{navy}\chi}}{\sigma_1} \right)^p = \left( \frac{4 \gamma_{\infty} \thesischfiveeps}{\sigma_1} \right)^p, \,\, \text{ so } \,\, \bar{\theta} \ge \left( \frac{4 \gamma_{\infty} \thesischfiveeps}{\sigma_1} \right)^{p/(p - 1)} \end{thesischfivefitdisplay}
        where the equality uses the values of $\alpha$ and ${\color{navy}\chi}$ fixed earlier. Hence
            \begin{thesischfivefitdisplay} \theta = \min\thesischfiveens{1, \left( \frac{4 \gamma_{\infty} \thesischfiveeps}{\sigma_1} \right)^{p/(p - 1)}}. \thesischfivenumberthis\label{ch5:b3242270-31d2-4d6a-a567-27370b26f605-2_chapter_5} \end{thesischfivefitdisplay}
        {\color{navy}If this minimum equals $1$, the Bernoulli oracle is deterministic and the central moment is exactly zero. If $\theta<1$, the preceding moment bound applies directly. Thus the choice is valid in both branches.}

        Next, concerning the mean-squared smoothness assumption, we have
            \begin{thesischfivefitalignatnum}{2}
                \thesischfiveE{\norm{\bar{g}^{\star}_N(x, \xi) - \bar{g}^{\star}_N(y, \xi)}^q} \thesischfiveoversetlab{\eqref{ch5:6f55ce68-cc03-4d22-a262-4aa3367f948d_chapter_5}}&{=} \left( \alpha {\color{navy}\chi} \right)^q \thesischfiveE{\norm{\bar{g}_N({\color{navy}\chi} x, \xi) - \bar{g}_N({\color{navy}\chi} y, \xi)}^q} \\
                & \le \frac{\left( \alpha {\color{navy}\chi} \bar{\ell}_1 \right)^q}{\theta^{q - 1}} \norm{{\color{navy}\chi} x - {\color{navy}\chi} y}^q \\
                & = \left( \frac{\alpha {\color{navy}\chi}^2 \bar{\ell}_1}{\theta^{(q - 1)/q}} \right)^q \norm{x - y}^q \\
                & = \left( \frac{L \bar{\ell}_1}{\ell_1 \theta^{(q - 1)/q}} \right)^q \norm{x - y}^q, \label{ch5:b4ce77e4-45fb-40a3-9459-1426ca4e5db4-2_chapter_5}
            \end{thesischfivefitalignatnum}
        where the inequality follows from Lemma~\ref{ch5:appdx-lem:properties-gradient-estimator-1_chapter_5}, and the last equality uses $\alpha = L/({\color{navy}\chi}^2 \ell_1)$. Hence, from the upper bound~\eqref{ch5:b4ce77e4-45fb-40a3-9459-1426ca4e5db4-2_chapter_5} it suffices to take
            \begin{thesischfivefitdisplay} L = \frac{\ell_1 \bar{L}}{\bar{\ell}_1} \theta^{(q - 1)/q} \le \bar{L} \min\thesischfiveens{1, \left( \frac{4 \gamma_{\infty} \thesischfiveeps}{\sigma_1} \right)^{p (q - 1)/(q (p - 1))}} \le \bar{L}, \end{thesischfivefitdisplay}
        since $\ell_1 = 152 \le \bar{\ell}_1$ (see~Lemma~\ref{ch5:appdx-lem:properties-gradient-estimator-1_chapter_5}). This proves~Assumption~\ref{ch5:ass:mean-squared-smoothness_chapter_5} is satisfied by the gradient estimator $\bar{g}^{\star}_N$.

        \item \textbf{Step 4:} \textit{Lower Bounding $\mathfrak{m}^{\normalfont\texttt{zr}}_{\thesischfiveeps}\left( {\color{navy}b}, \Delta, \bar{L}^q, \sigma^p_1 \right)$}.

        {\color{navy}Continuing Step~2, Lemma~\ref{ch5:appdx-lem:zero-chain-progress_chapter_5} shows that, with probability at least $1/2$, for every integer $0 \le t \le (N - 1)/(2 \theta)$ and every $k \in [b]$,}
            \begin{thesischfivefitdisplay} \norm{\nabla f^{\star}_N\left( x^{(t, k)}_{\texttt{A}[\texttt{O}_f]} \right)} > 2\thesischfiveeps, \,\, \text{ hence } \,\, \thesischfiveE{\norm{\nabla f^{\star}_N\left( x^{(t, k)}_{\texttt{A}[\texttt{O}_f]} \right)}} > \thesischfiveeps, \end{thesischfivefitdisplay}
        from where it follows that
            \begin{thesischfivefitalignatnum}{1}
                &\mathfrak{m}^{\normalfont\texttt{zr}}_{\thesischfiveeps}
                \left({\color{navy}b},\Delta,\bar{L}^q,\sigma_1^p\right)
                >\frac{N-1}{2\theta}
                \notag\\
                &\quad=
                \frac{1}{2\theta}
                \left(\thesischfiveFloor{\frac{\Delta}{\alpha\Delta_0}}-1\right)
                \notag\\
                &\quad=
                \frac{1}{2\theta}
                \left(\thesischfiveFloor{
                \frac{{\color{navy}\chi}^2\ell_1\Delta}{L\Delta_0}}-1\right)
                \notag\\
                &\quad=
                \frac{1}{2\theta}
                \left(\thesischfiveFloor{
                \frac{L\Delta}{4\Delta_0\ell_1\thesischfiveeps^2}}-1\right)
                \notag\\
                &\quad=
                \frac{1}{2\theta}
                \left(\thesischfiveFloor{
                \frac{\bar{L}\Delta\theta^{(q-1)/q}}
                {4\Delta_0\bar{\ell}_1\thesischfiveeps^2}}-1\right),
                \label{ch5:b3242270-31d2-4d6a-a567-27370b26f605_chapter_5}
            \end{thesischfivefitalignatnum}
        we then distinguish two cases:
            \begin{itemize}
                \item if $\bar{L} \Delta \theta^{(q - 1)/q} / (4 \Delta_0 \bar{\ell}_1 \thesischfiveeps^2) \ge 3$ then, using the inequality $\thesischfiveFloor{x} - 1 \ge x/2$, valid for all real number $x \ge 3$ we obtain
                    \begin{thesischfivefitalignat}{2}
                        \mathfrak{m}^{\normalfont\texttt{zr}}_{\thesischfiveeps}\left( {\color{navy}b}, \Delta, \bar{L}^q, \sigma^p_1 \right) \thesischfiveoversetlab{\eqref{ch5:b3242270-31d2-4d6a-a567-27370b26f605_chapter_5}}&{>} \frac{1}{4 \theta} \cdot \frac{\bar{L} \Delta \theta^{(q - 1)/q}}{4 \Delta_0 \bar{\ell}_1 \thesischfiveeps^2} \\
                        & = \frac{1}{16 \Delta_0 \bar{\ell}_1} \cdot \frac{\bar{L} \Delta}{\thesischfiveeps^2 \theta^{1/q}} \\
                        \thesischfiveoversetlab{\eqref{ch5:b3242270-31d2-4d6a-a567-27370b26f605-2_chapter_5}}&{\ge} \frac{1}{32 \Delta_0 \bar{\ell}_1} \cdot \frac{\bar{L} \Delta}{\thesischfiveeps^2} \left[ \left( \frac{\sigma_1}{4 \gamma_{\infty} \thesischfiveeps} \right)^{p/(q (p - 1))} + 1 \right] \\
                        & = \Omega(1) \cdot \left( \frac{\bar{L} \Delta}{\thesischfiveeps^2} + \frac{\bar{L} \Delta}{\thesischfiveeps^2} \left( \frac{\sigma_1}{\thesischfiveeps} \right)^{p/(q (p - 1))} \right),
                    \end{thesischfivefitalignat}
                and combining the lower bound above with~Lemma~\ref{ch5:appdx-lem:lower-bound-global-stichastic-model-mean-squared-smoothness-p-BCM_chapter_5}, with the same choice for $c'$ as provided below, gives the desired results.

                \item otherwise, if $\bar{L} \Delta \theta^{(q - 1)/q} / (4 \Delta_0 \bar{\ell}_1 \thesischfiveeps^2) < 3$, choosing the universal constant $0 < c' = \left( 12 \Delta_0 \bar{\ell}_1 \right)^{-1/2} \le \left( 12 \Delta_0 \ell_1 \right)^{-1/2} \approx 0.0067592 < 1/8$ then the assumption 
                    \begin{thesischfivefitdisplay} 0 < \thesischfiveeps \le c' \sqrt{\bar{L} \Delta} = \sqrt{\frac{\bar{L} \Delta}{12 \Delta_0 \bar{\ell}_1}}, \end{thesischfivefitdisplay}
                precludes the possibility that $\theta = 1$ so $\sigma_1 \ge 4 \gamma_{\infty} \thesischfiveeps$ hence $\sigma_1/\thesischfiveeps \gtrsim 1$. Moreover, we have
                    \begin{thesischfivefitalignat}{2}
                        3 > \frac{\bar{L} \Delta}{4 \Delta_0 \bar{\ell}_1 \thesischfiveeps^2} \left( \frac{4 \gamma_{\infty} \thesischfiveeps}{\sigma_1} \right)^{p (q - 1)/(q (p - 1))},
                    \end{thesischfivefitalignat}
                and multiplying both sides by $(\sigma_1/\thesischfiveeps)^{p/(p - 1)} > 0$ leads to
                    \begin{thesischfivefitalignat}{2}
                        \left( \frac{\sigma_1}{\thesischfiveeps} \right)^{p/(p - 1)} & > \frac{\bar{L} \Delta \left( 4 \gamma_{\infty} \right)^{p (q - 1)/(q (p - 1))}}{12 \Delta_0 \bar{\ell}_1 \thesischfiveeps^2} \left( \frac{\sigma_1}{\thesischfiveeps} \right)^{p/(p - 1) - p (q - 1)/(q (p - 1))} \\
                        & = \frac{\bar{L} \Delta \left( 4 \gamma_{\infty} \right)^{p (q - 1)/(q (p - 1))}}{12 \Delta_0 \bar{\ell}_1 \thesischfiveeps^2} \left( \frac{\sigma_1}{\thesischfiveeps} \right)^{p/(q (p - 1))},
                    \end{thesischfivefitalignat}
                and using~Lemma~\ref{ch5:appdx-lem:lower-bound-global-stichastic-model-mean-squared-smoothness-p-BCM_chapter_5}, {\color{navy}since $0 < \thesischfiveeps \le \sqrt{\bar{L} \Delta}/8$, there exists a universal constant $C_p > 0$ (depending only on $p$) such that}
                    \begin{thesischfivefitdisplay} \mathfrak{m}^{\normalfont\texttt{zr}}_{\thesischfiveeps}\left( {\color{navy}b}, \Delta, \bar{L}^q, \sigma^p_1 \right) \ge C_p \left( \frac{\sigma_1}{\thesischfiveeps} \right)^{p/(p - 1)} \ge \Omega(1) \cdot \frac{\bar{L} \Delta}{\thesischfiveeps^2} \left( \frac{\sigma_1}{\thesischfiveeps} \right)^{p/(q (p - 1))} \gtrsim \frac{\bar{L} \Delta}{\thesischfiveeps^2}, \thesischfivenumberthis\label{ch5:e278b503-7290-42c8-a9cd-e396da8844d8_chapter_5} \end{thesischfivefitdisplay}
                hence, 
                    \begin{thesischfivefitdisplay} \mathfrak{m}^{\normalfont\texttt{zr}}_{\thesischfiveeps}\left( {\color{navy}b}, \Delta, \bar{L}^q, \sigma^p_1 \right) \ge \Omega(1) \cdot \left( \left( \frac{\sigma_1}{\thesischfiveeps} \right)^{p/(p - 1)} + \frac{\bar{L} \Delta}{\thesischfiveeps^2} + \frac{\bar{L} \Delta}{\thesischfiveeps^2} \left( \frac{\sigma_1}{\thesischfiveeps} \right)^{p/(q (p - 1))} \right), \end{thesischfivefitdisplay}
                which holds in both cases and concludes the proof of the theorem.
            \end{itemize}
    \end{itemize}
\end{proof}

\subsection{Proof of~Theorem~\ref{ch5:thm:lower-bound-sofo-mean-squared-smoothness-2_chapter_5}}\label{ch5:appdx-subsec:proof-lower-bound-sofo-mean-squared-smoothness-2_chapter_5}

\begin{theorem}[Restatement of Theorem~\ref{ch5:thm:lower-bound-sofo-mean-squared-smoothness-2_chapter_5}]
    Given $\Delta, L_1, \delta > 0$, $\sigma_1 \ge 0$ and $0 < \thesischfiveeps \le c_1 \sqrt{L_1 \Delta}$ for some universal constant $c_1 > 0$. Then, for any algorithm $A \in \mathcal{A}_{\texttt{zr}}$, there exists a function $f \in \mathcal{F}\left( \Delta, L_1 \right)$, an oracle and a distribution $(O, \thesischfivecD) \in \mathcal{O}\left( f, \delta^q, \sigma_1^p \right)$ satisfying~Assumption~\ref{ch5:ass:lower-boundedness_chapter_5}, Assumption~\ref{ch5:ass:L-lipschitz-gradients_chapter_5}, Assumption~\ref{ch5:ass:p-bounded-central-moment-gradient_chapter_5}, and Assumption~\ref{ch5:ass:mean-squared-smoothness-2_chapter_5} such that

    \begin{thesischfivefitalignat}{2}
        \MoveEqLeft[1] \mathfrak{m}^{\normalfont\texttt{zr}}_{\thesischfiveeps}\left( {\color{navy}b}, \Delta, L_1, \delta^q, \sigma_1^p \right)\\
        & \ge \Omega(1) \cdot \min\biggl\{ \frac{L_1 \Delta}{\thesischfiveeps^2} + \frac{L_1 \Delta}{\thesischfiveeps^2} \left( \frac{\sigma_1}{\thesischfiveeps} \right)^{p/(p - 1)},\\
        &\qquad \left( \frac{\sigma_1}{\thesischfiveeps} \right)^{p/(p - 1)} + \frac{(L_1 + \delta) \Delta}{\thesischfiveeps^2} + \frac{\delta \Delta}{\thesischfiveeps^2} \left( \frac{\sigma_1}{\thesischfiveeps} \right)^{p/(q (p - 1))} \biggr\}.
    \end{thesischfivefitalignat}
\end{theorem}

\begin{proof}
    The proof follows the same lines as the proof of~Theorem~\ref{ch5:thm:lower-bound-sofo-mean-squared-smoothness_chapter_5}.  Let $\Delta_0$, $\ell_1$, $\gamma_{\infty}$ and $\bar{\delta}_1$ be the numerical constants in~Lemma~\ref{ch5:appdx-lem:properties-hard-instance_chapter_5} and Lemma~\ref{ch5:appdx-lem:properties-gradient-estimator-1_chapter_5} respectively. Additionally, we let the accuracy parameter $\thesischfiveeps > 0$, initial sub-optimality $\Delta \ge 0$, the smoothness constant $L_1 \ge 0$, the $q$--weak average smoothness parameter $\delta$ (Assumption~\ref{ch5:ass:mean-squared-smoothness-2_chapter_5}), and the variance parameter $\sigma_1 \ge 0$ be fixed, and $0 < L \le L_1$ to be specified later. Then, for $\alpha, {\color{navy}\chi} > 0$ two positive real numbers, following~\textcite{arjevani2022lower}, we rescale the function $f_N$ as
        \begin{thesischfivefitdisplay} f^{\star}_N \colon x \mapsto \alpha f_N({\color{navy}\chi} x). \thesischfivenumberthis\label{ch5:6f55ce68-cc03-4d22-a262-4aa3367f948d-bis_chapter_5} \end{thesischfivefitdisplay}

    \begin{itemize}
        \item \textbf{Step 1:} \textit{Ensuring $f^{\star}_N \in \mathcal{F} \left(\Delta, L \right)$}.

        To guarantee the rescaled function $f^{\star}_N$ belongs to the function class $\mathcal{F}(\Delta, L)$, let us compute the initial sub-optimality $\Delta$ and the smoothness constant $L$. Assuming the algorithm \texttt{A} starts at $x^0 = 0$ we have
            \begin{thesischfivefitdisplay} f^{\star}_N(0) - \inf_{x \in \R^N} f^{\star}_N(x) \thesischfiveoversetlab{\eqref{ch5:6f55ce68-cc03-4d22-a262-4aa3367f948d-bis_chapter_5}}{=} \alpha \left( f_N(0) - \inf_{x \in \R^N} f_N(x) \right) \le \alpha \Delta_0 N, \end{thesischfivefitdisplay}
        where the inequality follows from Lemma~\ref{ch5:appdx-lem:properties-hard-instance_chapter_5}; {\color{navy}thus, it suffices to take $N = \thesischfiveFloor{\Delta/(\alpha \Delta_0)}$ to ensure $f^{\star}_N(0) - \inf_{x \in \R^N} f^{\star}_N(x) \le \Delta$.} Moreover, as done previously, for any $x, y \in \R^N$,
            \begin{thesischfivefitalignat}{2}
                \norm{\nabla f^{\star}_N(x) - \nabla f^{\star}_N(y)} & = \alpha {\color{navy}\chi} \norm{\nabla f_N({\color{navy}\chi} x) - \nabla f_N({\color{navy}\chi} y)} \\
                & \le \alpha {\color{navy}\chi} \ell_1 \norm{{\color{navy}\chi} x - {\color{navy}\chi} y} \\
                & = \alpha {\color{navy}\chi}^2 \ell_1 \norm{x - y},
            \end{thesischfivefitalignat}
        where the inequality follows from Lemma~\ref{ch5:appdx-lem:properties-hard-instance_chapter_5}, and it suffices to take $\alpha = L/({\color{navy}\chi}^2 \ell_1) > 0$ to ensure the function $f^{\star}_N$ has $L$--Lipschitz gradients. Consequently, we have $f^{\star}_N \in \mathcal{F}(\Delta, L)$, as desired.

        \item \textbf{Step 2:} \textit{Analysis of the Protocol and Choice for ${\color{navy}\chi}$}.

        Following the proof of~\textcite[Theorem~1]{arjevani2022lower}, according to~Lemma~\ref{ch5:appdx-lem:properties-hard-instance_chapter_5}, for all points $x \in \R^N$ such that $\thesischfiveprog_0(x) < N$ we have $\thesischfiveprog_0({\color{navy}\chi} x) = \thesischfiveprog_0(x) < N$ so
            \begin{thesischfivefitdisplay} \norm{\nabla f^{\star}_N(x)} \thesischfiveoversetlab{\eqref{ch5:6f55ce68-cc03-4d22-a262-4aa3367f948d-bis_chapter_5}}{=} \frac{L}{\ell_1 {\color{navy}\chi}} \norm{\nabla f_N({\color{navy}\chi} x)} > \frac{L}{\ell_1 {\color{navy}\chi}}, \thesischfivenumberthis\label{ch5:191cd1ea-8e5f-4b85-8e63-f31fd25f3b17-bis_chapter_5} \end{thesischfivefitdisplay}
        where the inequality follows from Lemma~\ref{ch5:appdx-lem:properties-hard-instance_chapter_5}, and we need to guarantee that
            \begin{thesischfivefitdisplay} \norm{\nabla f^{\star}_N(x)} > 2\thesischfiveeps, \thesischfivenumberthis\label{ch5:e8299abb-d1db-4ee4-9919-3b49ed05df54_chapter_5} \end{thesischfivefitdisplay}
        for all $x \in \R^N$ with $\thesischfiveprog_0(x) < N$ which, given~\eqref{ch5:191cd1ea-8e5f-4b85-8e63-f31fd25f3b17_chapter_5}, can be done if we set ${\color{navy}\chi} = L/(2 \ell_1 \thesischfiveeps)$.

        \item \textbf{Step 3:} \textit{The Oracle Class and~Assumption~\ref{ch5:ass:p-bounded-central-moment-gradient_chapter_5} and Assumption~\ref{ch5:ass:mean-squared-smoothness_chapter_5}}.

        It remains to choose the parameter $\theta \in \thesischfiveintof{0}{1}$ and constant $L$ such that the gradient estimator $\bar{g}^{\star}_N$ of $\nabla f^{\star}_N$ satisfies~Assumption~\ref{ch5:ass:p-bounded-central-moment-gradient_chapter_5} and Assumption~\ref{ch5:ass:mean-squared-smoothness-2_chapter_5}. Computing the $p$--th central moment of $\bar{g}^{\star}_N$ gives, for all $x \in \R^N$
            \begin{thesischfivefitalignatnum}{2}
                \thesischfiveE{\norm{\bar{g}^{\star}_N(x, \xi) - \nabla f^{\star}_N(x)}^p} \thesischfiveoversetlab{\eqref{ch5:6f55ce68-cc03-4d22-a262-4aa3367f948d-bis_chapter_5}}&{=} \left( \alpha {\color{navy}\chi} \right)^p \thesischfiveE{\norm{\bar{g}_N({\color{navy}\chi} x, \xi) - \nabla f_N({\color{navy}\chi} x)}^p} \\
                & \le \frac{2 \left( \gamma_{\infty} \alpha {\color{navy}\chi} \right)^p (1 - \theta)}{\theta^{p - 1}} \\
                & \le \frac{\left( 2 \gamma_{\infty} \alpha {\color{navy}\chi} \right)^p}{\theta^{p - 1}}, \label{ch5:29bb8bb5-e6e5-482d-ab46-fab6357955ec-bis_chapter_5}
            \end{thesischfivefitalignatnum}
        where the first inequality follows from Lemma~\ref{ch5:appdx-lem:properties-gradient-estimator-1_chapter_5}, and in the second we use the fact that $0 < \theta \le 1$ and $p \ge 1$ so that $2 \le 2^p$. {\color{navy}By~\eqref{ch5:29bb8bb5-e6e5-482d-ab46-fab6357955ec_chapter_5}, Assumption~\ref{ch5:ass:p-bounded-central-moment-gradient_chapter_5} is satisfied by taking $\theta = \min\thesischfiveens{1, \bar{\theta}}$, where}
            \begin{thesischfivefitdisplay} \bar{\theta}^{p - 1} \ge \left( \frac{2 \gamma_{\infty} \alpha {\color{navy}\chi}}{\sigma_1} \right)^p = \left( \frac{4 \gamma_{\infty} \thesischfiveeps}{\sigma_1} \right)^p, \,\, \text{ so } \,\, \bar{\theta} \ge \left( \frac{4 \gamma_{\infty} \thesischfiveeps}{\sigma_1} \right)^{p/(p - 1)} \end{thesischfivefitdisplay}
        where the equality uses the values of $\alpha$ and ${\color{navy}\chi}$ fixed earlier. Hence
            \begin{thesischfivefitdisplay} \theta = \min\thesischfiveens{1, \left( \frac{4 \gamma_{\infty} \thesischfiveeps}{\sigma_1} \right)^{p/(p - 1)}}. \thesischfivenumberthis\label{ch5:b3242270-31d2-4d6a-a567-27370b26f605-2-bis_chapter_5} \end{thesischfivefitdisplay}
        {\color{navy}Again, when $\theta=1$ the estimator is deterministic and both centered moment terms vanish; when $\theta<1$, the displayed upper bounds establish the required conditions.}

        Next, concerning the mean-squared smoothness assumption (Assumption~\ref{ch5:ass:mean-squared-smoothness-2_chapter_5}), we have
            \begin{alignat}{2}
                \MoveEqLeft[1] \thesischfiveE{\norm{\left[ \bar{g}^{\star}_N(x, \xi) - \bar{g}^{\star}_N(y, \xi) \right] - \left[ \nabla f_N^{\star}(x) - \nabla f_N^{\star}(y) \right]}^q} \notag\\
                \thesischfiveoversetlab{\eqref{ch5:6f55ce68-cc03-4d22-a262-4aa3367f948d-bis_chapter_5}}&{=} \left( \alpha {\color{navy}\chi} \right)^q \thesischfiveE{\norm{\left[ \bar{g}_N({\color{navy}\chi} x, \xi) - \bar{g}_N({\color{navy}\chi} y, \xi) \right] - \left[ \nabla f_N({\color{navy}\chi} x) - \nabla f_N({\color{navy}\chi} y) \right]}^q} \notag\\
                & \le \frac{\left( \alpha {\color{navy}\chi} \bar{\delta}_1 \right)^q (1 - \theta)}{\theta^{q - 1}} \norm{{\color{navy}\chi} x - {\color{navy}\chi} y}^q \notag\\
                & = \left( \frac{\alpha {\color{navy}\chi}^2 \bar{\delta}_1}{\theta^{(q - 1)/q}} \right)^q (1 - \theta) \norm{x - y}^q \notag\\
                & = \left( \frac{L \bar{\delta}_1}{\ell_1 \theta^{(q - 1)/q}} \right)^q (1 - \theta) \norm{x - y}^q, \label{ch5:b4ce77e4-45fb-40a3-9459-1426ca4e5db4-2-bis_chapter_5}
            \end{alignat}
        where the inequality follows from Lemma~\ref{ch5:appdx-lem:properties-gradient-estimator-1_chapter_5}, and the last equality uses $\alpha = L/({\color{navy}\chi}^2 \ell_1)$. Hence, from the upper bound~\eqref{ch5:b4ce77e4-45fb-40a3-9459-1426ca4e5db4-2_chapter_5} it suffices to take
            \begin{thesischfivefitdisplay} L \le \frac{\ell_1 \delta}{\bar{\delta}_1} \theta^{(q - 1)/q}, \end{thesischfivefitdisplay}
        and since we must have $L \le L_1$, we set
            \begin{thesischfivefitdisplay} L = \min\thesischfiveens{L_1, \frac{\ell_1 \delta}{\bar{\delta}_1} \theta^{(q - 1)/q}}. \thesischfivenumberthis\label{ch5:bf6dc5ba-5913-4acb-b17c-aa57b62b93b3_chapter_5} \end{thesischfivefitdisplay}
        This proves~Assumption~\ref{ch5:ass:mean-squared-smoothness-2_chapter_5} is satisfied by the gradient estimator $\bar{g}^{\star}_N$.

        \item \textbf{Step 4:} \textit{Lower Bounding $\mathfrak{m}^{\normalfont\texttt{zr}}_{\thesischfiveeps}\left( {\color{navy}b}, \Delta, L_1, \delta^q, \sigma_1^p \right)$}.

        {\color{navy}Continuing Step~2, Lemma~\ref{ch5:appdx-lem:zero-chain-progress_chapter_5} shows that, with probability at least $1/2$, for every integer $0 \le t \le (N - 1)/(2 \theta)$ and every $k \in [b]$,}
            \begin{thesischfivefitdisplay} \norm{\nabla f^{\star}_N\left( x^{(t, k)}_{\texttt{A}[\texttt{O}_f]} \right)} > 2\thesischfiveeps, \,\, \text{ hence } \,\, \thesischfiveE{\norm{\nabla f^{\star}_N\left( x^{(t, k)}_{\texttt{A}[\texttt{O}_f]} \right)}} > \thesischfiveeps, \end{thesischfivefitdisplay}
        from where it follows that
            \begin{thesischfivefitalignatnum}{1}
                &\mathfrak{m}^{\normalfont\texttt{zr}}_{\thesischfiveeps}\left( {\color{navy}b}, \Delta, L_1, \delta^q, \sigma_1^p \right) > \frac{N - 1}{2 \theta}\\
                &\quad= \frac{1}{2 \theta} \left( \thesischfiveFloor{\frac{\Delta}{\alpha \Delta_0}} - 1 \right)\\
                &\quad= \frac{1}{2 \theta} \left( \thesischfiveFloor{\frac{{\color{navy}\chi}^2 \ell_1 \Delta}{L \Delta_0}} - 1 \right)\\
                &\quad= \frac{1}{2 \theta} \left( \thesischfiveFloor{\frac{L \Delta}{4 \Delta_0 \ell_1 \thesischfiveeps^2}} - 1 \right), \label{ch5:b3242270-31d2-4d6a-a567-27370b26f605-bis_chapter_5}
            \end{thesischfivefitalignatnum}
        we then distinguish two cases:
            \begin{itemize}
                \item if $L \Delta / (4 \Delta_0 \ell_1 \thesischfiveeps^2) \ge 3$ then, using the inequality $\thesischfiveFloor{x} - 1 \ge x/2$, valid for all real number $x \ge 3$ we obtain
                    \begin{thesischfivefitalignatnum}{2}
                        \MoveEqLeft[1] \mathfrak{m}^{\normalfont\texttt{zr}}_{\thesischfiveeps}\left( {\color{navy}b}, \Delta, L_1, \delta^q, \sigma_1^p \right) \\
                        \thesischfiveoversetlab{\eqref{ch5:b3242270-31d2-4d6a-a567-27370b26f605-bis_chapter_5}}&{>} \frac{1}{4 \theta} \cdot \frac{L \Delta}{4 \Delta_0 \ell_1 \thesischfiveeps^2} \\
                        \thesischfiveoversetlab{\eqref{ch5:bf6dc5ba-5913-4acb-b17c-aa57b62b93b3_chapter_5}}&{=} \frac{1}{16 \Delta_0 \ell_1} \cdot \frac{\Delta}{\theta \thesischfiveeps^2} \min\thesischfiveens{L_1, \frac{\ell_1 \delta}{\bar{\delta}_1} \theta^{(q - 1)/q}} \\
                        & = \frac{1}{16 \Delta_0 \ell_1} \cdot \min\thesischfiveens{\frac{L_1 \Delta}{\theta \thesischfiveeps^2}, \frac{\ell_1 \delta \Delta}{\bar{\delta}_1 \thesischfiveeps^2} \theta^{-1/q}} \\
                        \thesischfiveoversetlab{\eqref{ch5:b3242270-31d2-4d6a-a567-27370b26f605-2-bis_chapter_5}}&{\ge} \frac{1}{32 \Delta_0 \ell_1} \cdot \min\biggl\{\frac{L_1 \Delta}{\thesischfiveeps^2} \left( 1 + \left( \frac{\sigma_1}{4 \gamma_{\infty} \thesischfiveeps} \right)^{p/(p - 1)} \right),\\
                        &\qquad \frac{\ell_1 \delta \Delta}{\bar{\delta}_1 \thesischfiveeps^2} \left( 1 + \left( \frac{\sigma_1}{4 \gamma_{\infty} \thesischfiveeps} \right)^{p/(q (p - 1))} \right)\biggr\} \\
                        & = \Omega(1) \cdot \min\biggl\{\frac{L_1 \Delta}{\thesischfiveeps^2} + \frac{L_1 \Delta}{\thesischfiveeps^2} \left( \frac{\sigma_1}{\thesischfiveeps} \right)^{p/(p - 1)},\\
                        &\qquad \frac{\delta \Delta}{\thesischfiveeps^2} + \frac{\delta \Delta}{\thesischfiveeps^2} \left( \frac{\sigma_1}{\thesischfiveeps} \right)^{p/(q (p - 1))}\biggr\}, \label{ch5:3b256059-e30a-4307-b1b7-a8e09886c0f1_chapter_5}
                    \end{thesischfivefitalignatnum}
                and combining the lower bound~\eqref{ch5:3b256059-e30a-4307-b1b7-a8e09886c0f1_chapter_5} with~Lemma~\ref{ch5:appdx-lem:lower-bound-global-stichastic-model-mean-squared-smoothness-p-BCM_chapter_5}, using the choice for $c'$ provided below, gives:
                    \begin{thesischfivefitalignatnum}{2}
                        \MoveEqLeft[1] \mathfrak{m}^{\normalfont\texttt{zr}}_{\thesischfiveeps}\left( {\color{navy}b}, \Delta, L_1, \delta^q, \sigma_1^p \right) \\
                        & \ge \Omega(1) \cdot \min\biggl\{\frac{L_1 \Delta}{\thesischfiveeps^2} + \frac{L_1 \Delta}{\thesischfiveeps^2} \left( \frac{\sigma_1}{\thesischfiveeps} \right)^{p/(p - 1)},\\
                        &\qquad \left( \frac{\sigma_1}{\thesischfiveeps} \right)^{p/(p - 1)} + \frac{\delta \Delta}{\thesischfiveeps^2} + \frac{\delta \Delta}{\thesischfiveeps^2} \left( \frac{\sigma_1}{\thesischfiveeps} \right)^{p/(q (p - 1))}\biggr\}, \label{ch5:5aea445e-d78f-404f-ac66-dbff30f0a00c_chapter_5}
                    \end{thesischfivefitalignatnum}
                since from the assumption $0 < \thesischfiveeps \le c' \sqrt{L_1 \Delta}$ we have
                    \begin{thesischfivefitdisplay} \frac{L_1 \Delta}{\thesischfiveeps^2} \left( \frac{\sigma_1}{\thesischfiveeps} \right)^{p/(p - 1)} \gtrsim \left( \frac{\sigma_1}{\thesischfiveeps} \right)^{p/(p - 1)}. \end{thesischfivefitdisplay}

                \item otherwise, if $L \Delta / (4 \Delta_0 \ell_1 \thesischfiveeps^2) < 3$, choosing the universal constant $0 < c' = \left( 12 \Delta_0 \ell_1 \right)^{-1/2} \approx 0.0067592 < 1/8$ then the assumption 
                    \begin{thesischfivefitdisplay} 0 < \thesischfiveeps \le c' \sqrt{L_1 \Delta} = \sqrt{\frac{L_1 \Delta}{12 \Delta_0 \ell_1}}, \end{thesischfivefitdisplay}
                precludes the possibility that $L = L_1$. {\color{navy}Although the argument used in the proof of~Theorem~\ref{ch5:thm:lower-bound-sofo-mean-squared-smoothness_chapter_5} remains applicable, here we use the alternative strategy described in~Remark~\ref{ch5:appdx-rem:improve-lower-bound-global-stichastic-model-mean-squared-smoothness-p-BCM_chapter_5}. Dividing both sides of $L \Delta / (4 \Delta_0 \ell_1 \thesischfiveeps^2) < 3$ by our choice of $\theta > 0$ gives}
                    \begin{thesischfivefitalignatnum}{2}
                        \frac{1}{12 \Delta_0 \ell_1} \cdot \frac{L \Delta}{\thesischfiveeps^2} \theta^{-1} & \le \theta^{-1} = \max\thesischfiveens{1, \left( \frac{\sigma_1}{4 \gamma_{\infty} \thesischfiveeps} \right)^{p/(p - 1)}} \\
                        & \le \max\thesischfiveens{1, \left( \frac{\sigma_1}{\thesischfiveeps} \right)^{p/(p - 1)}}, \label{ch5:2b14b2b5-0770-4ca0-b177-61bcd5125db3_chapter_5}
                    \end{thesischfivefitalignatnum}
                where the last inequality follows from the fact that $p/(p - 1) > 0$ and $4 \gamma_{\infty} = 92 \ge 1$. Using the definition of $L$ we obtain
                    \begin{thesischfivefitalignat}{2}
                        \MoveEqLeft[1] \max\thesischfiveens{1, \left( \frac{\sigma_1}{\thesischfiveeps} \right)^{p/(p - 1)}} \\
                        \thesischfiveoversetlab{\eqref{ch5:2b14b2b5-0770-4ca0-b177-61bcd5125db3_chapter_5}+\eqref{ch5:bf6dc5ba-5913-4acb-b17c-aa57b62b93b3_chapter_5}}&{\ge} \frac{1}{12 \Delta_0 \ell_1} \cdot \min\thesischfiveens{\frac{L_1 \Delta}{\theta \thesischfiveeps^2}, \frac{\ell_1 \delta \Delta}{\bar{\delta}_1 \thesischfiveeps^2} \theta^{-1/q}} \\
                        \thesischfiveoversetlab{\eqref{ch5:b3242270-31d2-4d6a-a567-27370b26f605-2-bis_chapter_5}}&{\ge} \frac{1}{24 \Delta_0 \ell_1} \cdot \min\biggl\{\frac{L_1 \Delta}{\thesischfiveeps^2} + \frac{L_1 \Delta}{\thesischfiveeps^2} \left( \frac{\sigma_1}{4 \gamma_{\infty} \thesischfiveeps} \right)^{p/(p - 1)},\\
                        &\qquad \frac{\ell_1 \delta \Delta}{\bar{\delta}_1 \thesischfiveeps^2} + \frac{\ell_1 \delta \Delta}{\bar{\delta}_1 \thesischfiveeps^2} \left( \frac{\sigma_1}{4 \gamma_{\infty} \thesischfiveeps} \right)^{p/(q (p - 1))}\biggr\} \\
                        & \ge C_{p, q} \min\biggl\{\frac{L_1 \Delta}{\thesischfiveeps^2} + \frac{L_1 \Delta}{\thesischfiveeps^2} \left( \frac{\sigma_1}{\thesischfiveeps} \right)^{p/(p - 1)},\\
                        &\qquad \frac{\delta \Delta}{\thesischfiveeps^2} + \frac{\delta \Delta}{\thesischfiveeps^2} \left( \frac{\sigma_1}{\thesischfiveeps} \right)^{p/(q (p - 1))}\biggr\},
                    \end{thesischfivefitalignat}
                {\color{navy}where $C_{p,q}>0$ is a universal constant depending only on $p$ and $q$. Moreover, by~Lemma~\ref{ch5:appdx-lem:lower-bound-global-stichastic-model-mean-squared-smoothness-p-BCM_chapter_5}, and more precisely Remark~\ref{ch5:appdx-rem:improve-lower-bound-global-stichastic-model-mean-squared-smoothness-2-p-BCM_chapter_5}, the condition $0 < \thesischfiveeps \le \sqrt{L_1 \Delta}/8$ implies that there exists a universal constant $C_p > 0$ (depending only on $p$) such that}
                    \begin{thesischfivefitdisplay} \mathfrak{m}^{\normalfont\texttt{zr}}_{\thesischfiveeps}\left( {\color{navy}b}, \Delta, L_1, \delta^q, \sigma_1^p \right) \ge C_p \max\thesischfiveens{1, \left( \frac{\sigma_1}{\thesischfiveeps} \right)^{p/(p - 1)}}, \end{thesischfivefitdisplay}
                hence
                    \begin{thesischfivefitalignatnum}{2}
                        \MoveEqLeft[1] \mathfrak{m}^{\normalfont\texttt{zr}}_{\thesischfiveeps}\left( {\color{navy}b}, \Delta, L_1, \delta^q, \sigma_1^p \right) \\
                        & \ge \Omega(1) \cdot \biggl( \left( \frac{\sigma_1}{\thesischfiveeps} \right)^{p/(p - 1)} + \min\biggl\{\frac{L_1 \Delta}{\thesischfiveeps^2} + \frac{L_1 \Delta}{\thesischfiveeps^2} \left( \frac{\sigma_1}{\thesischfiveeps} \right)^{p/(p - 1)},\\
                        &\qquad \frac{\delta \Delta}{\thesischfiveeps^2} + \frac{\delta \Delta}{\thesischfiveeps^2} \left( \frac{\sigma_1}{\thesischfiveeps} \right)^{p/(q (p - 1))}\biggr\} \biggr) \\
                        & \ge \Omega(1) \cdot \min\biggl\{\frac{L_1 \Delta}{\thesischfiveeps^2} + \frac{L_1 \Delta}{\thesischfiveeps^2} \left( \frac{\sigma_1}{\thesischfiveeps} \right)^{p/(p - 1)},\\
                        &\qquad \left( \frac{\sigma_1}{\thesischfiveeps} \right)^{p/(p - 1)} + \frac{\delta \Delta}{\thesischfiveeps^2} + \frac{\delta \Delta}{\thesischfiveeps^2} \left( \frac{\sigma_1}{\thesischfiveeps} \right)^{p/(q (p - 1))}\biggr\}, \label{ch5:5aea445e-d78f-404f-ac66-dbff30f0a00c-2_chapter_5}
                    \end{thesischfivefitalignatnum}
                as $0 < \thesischfiveeps \le c' \sqrt{L_1 \Delta}$ implies
                    \begin{thesischfivefitdisplay} \frac{L_1 \Delta}{\thesischfiveeps^2} \left( \frac{\sigma_1}{\thesischfiveeps} \right)^{p/(p - 1)} \gtrsim \left( \frac{\sigma_1}{\thesischfiveeps} \right)^{p/(p - 1)}, \end{thesischfivefitdisplay}
                and we can forget $(\sigma_1/\thesischfiveeps)^{p/(p - 1)}$ in the first term of the $\min$.
            \end{itemize}

            \item \textbf{Step 5:} \textit{A Last Bound: the Case $\theta = 1$}.

            Observe that, if instead of taking $\theta$ as in~\eqref{ch5:b3242270-31d2-4d6a-a567-27370b26f605-2-bis_chapter_5}, we choose directly $\theta = 1$ then, thanks to~Lemma~\ref{ch5:appdx-lem:properties-gradient-estimator-1_chapter_5} we immediately have
                \begin{thesischfivefitdisplay} \thesischfiveE{\norm{\bar{g}^{\star}_N(x, \xi) - \nabla f^{\star}_N(x)}^p} = 0, \end{thesischfivefitdisplay}
            and
                \begin{thesischfivefitdisplay} \thesischfiveE{\norm{\left[ \bar{g}^{\star}_N(x, \xi) - \bar{g}^{\star}_N(y, \xi) \right] - \left[ \nabla f_N^{\star}(x) - \nabla f_N^{\star}(y) \right]}^q} = 0, \end{thesischfivefitdisplay}
            so~Assumption~\ref{ch5:ass:p-bounded-central-moment-gradient_chapter_5} and~Assumption~\ref{ch5:ass:mean-squared-smoothness-2_chapter_5} are satisfied. Hence, if we set 
                \begin{thesischfivefitdisplay} N = \thesischfiveFloor{\frac{\Delta}{\alpha \Delta_0}}, \quad \alpha = \frac{L_1}{{\color{navy}\chi}^2 \ell_1}, \quad \text{ and } \quad {\color{navy}\chi} = \frac{L_1}{2 \ell_1 \thesischfiveeps}, \end{thesischfivefitdisplay}
            then $f_N^{\star} \in \mathcal{F}(\Delta, L_1)$ and the inequality~\eqref{ch5:e8299abb-d1db-4ee4-9919-3b49ed05df54_chapter_5} is satisfied. {\color{navy}Hence, with probability at least $1/2$, for every integer $0 \le t \le (N - 1)/(2 \theta)$ and every $k \in [b]$,}
            \begin{thesischfivefitdisplay} \norm{\nabla f^{\star}_N\left( x^{(t, k)}_{\texttt{A}[\texttt{O}_f]} \right)} > 2\thesischfiveeps, \,\, \text{ hence } \,\, \thesischfiveE{\norm{\nabla f^{\star}_N\left( x^{(t, k)}_{\texttt{A}[\texttt{O}_f]} \right)}} > \thesischfiveeps, \end{thesischfivefitdisplay}
        from where we obtain also
            \begin{thesischfivefitalignatnum}{1}
                &\mathfrak{m}^{\normalfont\texttt{zr}}_{\thesischfiveeps}\left( {\color{navy}b}, \Delta, L_1, \delta^q, \sigma_1^p \right) > \frac{N - 1}{2 \theta}\\
                &\quad= \frac{1}{2 \theta} \left( \thesischfiveFloor{\frac{\Delta}{\alpha \Delta_0}} - 1 \right)\\
                &\quad= \frac{1}{2 \theta} \left( \thesischfiveFloor{\frac{{\color{navy}\chi}^2 \ell_1 \Delta}{L_1 \Delta_0}} - 1 \right)\\
                &\quad= \frac{1}{2 \theta} \left( \thesischfiveFloor{\frac{L_1 \Delta}{4 \Delta_0 \ell_1 \thesischfiveeps^2}} - 1 \right), \label{ch5:b3242270-31d2-4d6a-a567-27370b26f605-bis-bis_chapter_5}
            \end{thesischfivefitalignatnum}
        and, since we assume $0 < \thesischfiveeps < c' \sqrt{L_1 \Delta} = \sqrt{L_1 \Delta / (12 \Delta_0 \ell_1)}$ this implies $L_1 \Delta / (4 \Delta_0 \ell_1 \thesischfiveeps^2) \ge 3$ thus
            \begin{thesischfivefitalignat}{2}
                \mathfrak{m}^{\normalfont\texttt{zr}}_{\thesischfiveeps}\left( {\color{navy}b}, \Delta, L_1, \delta^q, \sigma_1^p \right) \thesischfiveoversetlab{\eqref{ch5:b3242270-31d2-4d6a-a567-27370b26f605-bis-bis_chapter_5}}&{>} \frac{1}{4 \theta} \cdot \frac{L_1 \Delta}{4 \Delta_0 \ell_1 \thesischfiveeps^2} \\
                & = \Omega(1) \cdot \frac{L_1 \Delta}{\thesischfiveeps^2},
            \end{thesischfivefitalignat}
        and, combining this bound with~\eqref{ch5:5aea445e-d78f-404f-ac66-dbff30f0a00c_chapter_5} and~\eqref{ch5:5aea445e-d78f-404f-ac66-dbff30f0a00c-2_chapter_5} respectively leads to the desired result, i.e.,
            \begin{thesischfivefitalignat}{2}
                \MoveEqLeft[1] \mathfrak{m}^{\normalfont\texttt{zr}}_{\thesischfiveeps}\left( {\color{navy}b}, \Delta, L_1, \delta^q, \sigma_1^p \right) \\
                & \ge \Omega(1) \cdot \min\biggl\{ \frac{L_1 \Delta}{\thesischfiveeps^2} + \frac{L_1 \Delta}{\thesischfiveeps^2} \left( \frac{\sigma_1}{\thesischfiveeps} \right)^{p/(p - 1)},\\
                &\qquad \left( \frac{\sigma_1}{\thesischfiveeps} \right)^{p/(p - 1)} + \frac{(L_1 + \delta) \Delta}{\thesischfiveeps^2} + \frac{\delta \Delta}{\thesischfiveeps^2} \left( \frac{\sigma_1}{\thesischfiveeps} \right)^{p/(q (p - 1))} \biggr\}.
            \end{thesischfivefitalignat}
    \end{itemize}
\end{proof}

\subsection{Proof of~Theorem~\ref{ch5:thm:lower-bound-sofo-bounded-central-moments_chapter_5}}\label{ch5:appdx-subsec:proof-lower-bound-sofo-bounded-central-moments_chapter_5}

\begin{lemma}[Properties of the Gradient and Hessian Estimators $g_K$ and $\nabla g_K$]\label{ch5:appdx-lem:properties-gradient-estimator-2_chapter_5}
    The stochastic gradient estimator $g_K$ is a probability--$\theta$ zero-chain, is unbiased with respect to $\nabla f_N$ and satisfies
        \begin{align*}
            \thesischfiveE{\norm{\nabla f_N (x, \xi) - \nabla f_N(x)}^p} &\le \frac{2 \gamma_{\infty}^p (1 - \theta)}{\theta^{p - 1}}, \quad \text{and} \\
            \thesischfiveE{\normop{\nabla^2 f_N (x, \xi) - \nabla^2 f_N (x)}^q} &\le \frac{2 \ell_1^q (1 - \theta)}{\theta^{q - 1}},
        \end{align*}
    for all $x \in \R^N$, where $p \in \thesischfiveintof{1}{2}$, $q \in \thesischfiveintff{1}{2}$, $\gamma_{\infty}$ and $\ell_1$ are defined in~Lemma~\ref{ch5:appdx-lem:properties-hard-instance_chapter_5}.
\end{lemma}

The proof of~Lemma~\ref{ch5:appdx-lem:properties-gradient-estimator-2_chapter_5} is very similar to~Lemma~\ref{ch5:appdx-lem:properties-gradient-estimator-1_chapter_5} (and is simpler since we can directly bound $\mathbb{I}\!\{i > \thesischfiveprog_{1/4}(x)\}$ by $1$). For the bound on the difference $\nabla^2 f_N (x, \xi) - \nabla^2 f_N (x)$ we use~Lemma~\ref{ch5:appdx-lem:properties-hard-instance_chapter_5}, that is, $f_N$ has $\ell_1$--Lipschitz gradients along with the fact that $f_N$ is twice continuously differentiable which allows to bound the Hessian appropriately (see~Lemma~\ref{ch5:appdx-lem:lipschitz-gradients-implies-bounded-hessian_chapter_5}). Moreover, by definition of the operator norm $\normop{\cdot}$ (Definition~\ref{ch5:appdx-def:operator-norm_chapter_5}) we have
    \begin{thesischfivefitalignat}{2}
        \MoveEqLeft[1] \normop{\left( \left[ \nabla^2 f_N(x) \right]_{i, \cdot} \, \mathbb{I}\!\thesischfiveens{i > \thesischfiveprog_{1/4}(x)} \right)_{i \in [N]}} \\
        & = \sup_{y \in \R^d, \norm{y} = 1} \left( \sum_{i = 1}^N \mathbb{I}\!\thesischfiveens{i > \thesischfiveprog_{1/4}(x)} \abs{\thesischfiveT{\left[ \nabla^2 f_N(x) \right]_{i, \cdot}} y}^2 \right)^{1/2} \\
        & \le \sup_{y \in \R^d, \norm{y} = 1} \left( \sum_{i = 1}^N \abs{\thesischfiveT{\left[ \nabla^2 f_N(x) \right]_{i, \cdot}} y}^2 \right)^{1/2} \\
        & = \normop{\nabla^2 f_N(x)} \\
        & \le \ell_1,
    \end{thesischfivefitalignat}
where $\left[ \nabla^2 f_N(x) \right]_{i, \cdot}$ denotes the $i$--th row of the Hessian of $f_N$. The equality is Definition~\ref{ch5:appdx-def:operator-norm_chapter_5}; in the first inequality we use $\mathbb{I}\!\thesischfiveens{i > \thesischfiveprog_{1/4}(x)} \le 1$, and the last one follows from Lemma~\ref{ch5:appdx-lem:lipschitz-gradients-implies-bounded-hessian_chapter_5}.

\begin{theorem}[Restatement of Theorem~\ref{ch5:thm:lower-bound-sofo-bounded-central-moments_chapter_5}]
    Given $\Delta, L_1, L_2 > 0$, $\sigma_1, \sigma_2 \ge 0$ and $0 < \thesischfiveeps \le c_1 \min\{\sqrt{L_1 \Delta}, L_2^{1/3} \Delta^{2/3} \}$ for some universal constant $c_1 > 0$. Then, for any algorithm $A \in \mathcal{A}_{\texttt{zr}}$, there exists a function $f \in \mathcal{F}\left( \Delta, L_1, L_2 \right)$, an oracle and a distribution $(O, \thesischfivecD) \in \mathcal{O}\left( f, \sigma_1^p, \sigma_2^q \right)$ satisfying~Assumption~\ref{ch5:ass:lower-boundedness_chapter_5}, Assumption~\ref{ch5:ass:L-lipschitz-gradients_chapter_5}, Assumption~\ref{ch5:ass:p-bounded-central-moment-gradient_chapter_5}, Assumption~\ref{ch5:ass:L-lipschitz-hessians_chapter_5}, and Assumption~\ref{ch5:ass:q-bounded-central-moment-hessian_chapter_5} such that
    \begin{thesischfivefitalignat}{2}
        \MoveEqLeft[1] \mathfrak{m}^{\normalfont\texttt{zr}}_{\thesischfiveeps}\left( {\color{navy}b}, \Delta, L_1, L_2, \sigma^p_1, \sigma_2^q \right) \\
        & \ge \Omega(1) \cdot \min\biggl\{\frac{L_1 \Delta}{\thesischfiveeps^2} + \frac{L_1 \Delta}{\thesischfiveeps^2} \left( \frac{\sigma_1}{\thesischfiveeps} \right)^{p/(p - 1)},\ 
        \frac{L_2^{1/2} \Delta}{\thesischfiveeps^{3/2}} + \frac{L_2^{1/2} \Delta}{\thesischfiveeps^{3/2}} \left( \frac{\sigma_1}{\thesischfiveeps} \right)^{p/(p - 1)}, \\
        &\qquad \min\thesischfiveens{\frac{L_1 \Delta}{\thesischfiveeps^2}, \frac{L_2^{1/2} \Delta}{\thesischfiveeps^{3/2}}} + \frac{\Delta \sigma_2}{\thesischfiveeps^2} + \frac{\Delta \sigma_2}{\thesischfiveeps^2} \left( \frac{\sigma_1}{\thesischfiveeps} \right)^{p/(q (p - 1))} + \left( \frac{\sigma_1}{\thesischfiveeps} \right)^{p/(p - 1)}\biggr\}.
    \end{thesischfivefitalignat}
\end{theorem}

\begin{proof}
    Let $\Delta_0$, $\ell_1$ and $\gamma_{\infty}$ be the numerical constants in~Lemma~\ref{ch5:appdx-lem:properties-hard-instance_chapter_5} respectively. Additionally, we let the accuracy parameter $\thesischfiveeps > 0$, initial sub-optimality $\Delta \ge 0$, the Lipschitz constants $L_1, L_2 \ge 0$ of the gradients and Hessians of $f$ respectively, and the variance parameters $\sigma_1, \sigma_2 \ge 0$ of the stochastic gradients and Hessians be fixed. Then, for $\alpha, {\color{navy}\chi} > 0$ two positive real numbers, as in~\textcite{arjevani2022lower}, we rescale the function $f_N$ as
        \begin{thesischfivefitdisplay} f^{\star}_N \colon x \mapsto \alpha f_N({\color{navy}\chi} x). \thesischfivenumberthis\label{ch5:6f55ce68-cc03-4d22-a262-4aa3367f948d-2_chapter_5} \end{thesischfivefitdisplay}

    \begin{itemize}
        \item \textbf{Step 1:} \textit{Ensuring $f^{\star}_N \in \mathcal{F} \left(\Delta, L_1, L_2 \right)$}.

        To guarantee the rescaled function $f^{\star}_N$ belongs to the function class $\mathcal{F}(\Delta, L_1, L_2)$, let us compute the initial sub-optimality $\Delta$ and check if $f_N$ has $L_1$--Lipschitz gradients and $L_2$--Lipschitz Hessians. Assuming the algorithm \texttt{A} starts at $x^0 = 0$ we have
            \begin{thesischfivefitdisplay} f^{\star}_N(0) - \inf_{x \in \R^N} f^{\star}_N(x) \thesischfiveoversetlab{\eqref{ch5:6f55ce68-cc03-4d22-a262-4aa3367f948d-2_chapter_5}}{=} \alpha \left( f_N(0) - \inf_{x \in \R^N} f_N(x) \right) \le \alpha \Delta_0 N, \end{thesischfivefitdisplay}
        where the inequality follows from Lemma~\ref{ch5:appdx-lem:properties-hard-instance_chapter_5}; {\color{navy}thus, it suffices to take $N = \thesischfiveFloor{\Delta/(\alpha \Delta_0)}$ to ensure $f^{\star}_N(0) - \inf_{x \in \R^N} f^{\star}_N(x) \le \Delta$.} Moreover, for any $x, y \in \R^N$,
            \begin{thesischfivefitalignat}{2}
                \norm{\nabla f^{\star}_N(x) - \nabla f^{\star}_N(y)} & = \alpha {\color{navy}\chi} \norm{\nabla f_N({\color{navy}\chi} x) - \nabla f_N({\color{navy}\chi} y)} \\
                & \le \alpha {\color{navy}\chi} \ell_1 \norm{{\color{navy}\chi} x - {\color{navy}\chi} y} \\
                & = \alpha {\color{navy}\chi}^2 \ell_1 \norm{x - y},
            \end{thesischfivefitalignat}
        and
            \begin{thesischfivefitalignat}{2}
                \normop{\nabla^2 f^{\star}_N(x) - \nabla^2 f^{\star}_N(y)} & = \alpha {\color{navy}\chi}^2 \normop{\nabla^2 f_N({\color{navy}\chi} x) - \nabla^2 f_N({\color{navy}\chi} y)} \\
                & \le \alpha {\color{navy}\chi}^2 \ell_2 \norm{{\color{navy}\chi} x - {\color{navy}\chi} y} \\
                & = \alpha {\color{navy}\chi}^3 \ell_2 \norm{x - y},
            \end{thesischfivefitalignat}
        where both inequalities follow from Lemma~\ref{ch5:appdx-lem:properties-hard-instance_chapter_5}. So it suffices to take $0 < \alpha {\color{navy}\chi}^2 \le L_1/\ell_1$ and $0 < \alpha {\color{navy}\chi}^3 \le L_2/\ell_2$ to ensure the function $f^{\star}_N$ has $L_1$--Lipschitz gradients and $L_2$--Lipschitz Hessians. Consequently, we have $f^{\star}_N \in \mathcal{F}(\Delta, L_1, L_2)$, as desired.

        \item \textbf{Step 2:} \textit{Analysis of the Protocol and Choice for ${\color{navy}\chi}$}.

        Following the proof of~\textcite[Theorem~1]{arjevani2022lower}, according to~Lemma~\ref{ch5:appdx-lem:properties-hard-instance_chapter_5}, for all points $x \in \R^N$ such that $\thesischfiveprog_0(x) < N$ we have $\thesischfiveprog_0({\color{navy}\chi} x) = \thesischfiveprog_0(x) < N$ so
            \begin{thesischfivefitdisplay} \norm{\nabla f^{\star}_N(x)} \thesischfiveoversetlab{\eqref{ch5:6f55ce68-cc03-4d22-a262-4aa3367f948d-2_chapter_5}}{=} \alpha {\color{navy}\chi} \norm{\nabla f_N({\color{navy}\chi} x)} > \alpha {\color{navy}\chi}, \thesischfivenumberthis\label{ch5:191cd1ea-8e5f-4b85-8e63-f31fd25f3b17-2_chapter_5} \end{thesischfivefitdisplay}
        where the inequality follows from Lemma~\ref{ch5:appdx-lem:properties-hard-instance_chapter_5}, and we need to guarantee that
            \begin{thesischfivefitdisplay} \norm{\nabla f^{\star}_N(x)} > 2\thesischfiveeps, \thesischfivenumberthis\label{ch5:f01cde9c-e4f8-4c8b-b9b4-7e86b18da9c1_chapter_5} \end{thesischfivefitdisplay}
        for all $x \in \R^N$ with $\thesischfiveprog_0(x) < N$ which, given~\eqref{ch5:191cd1ea-8e5f-4b85-8e63-f31fd25f3b17-2_chapter_5}, can be done if we set $\alpha = 2 \thesischfiveeps/{\color{navy}\chi}$.

        \item \textbf{Step 3:} \textit{The Oracle Class and~Assumption~\ref{ch5:ass:p-bounded-central-moment-gradient_chapter_5} and Assumption~\ref{ch5:ass:q-bounded-central-moment-hessian_chapter_5}}.

        It remains to choose the parameter $\theta \in \thesischfiveintof{0}{1}$ and ${\color{navy}\chi} > 0$ such that $0 < \alpha {\color{navy}\chi}^2 = 2 \thesischfiveeps {\color{navy}\chi} \le L_1/\ell_1$, $0 < \alpha {\color{navy}\chi}^3 = 2 \thesischfiveeps {\color{navy}\chi}^2 \le L_2/\ell_2$ and the gradient estimator $g^{\star}_N$ of $\nabla f^{\star}_N$ and Hessian estimator $\nabla g^{\star}_N$ of $\nabla^2 f^{\star}_N$ satisfies~Assumption~\ref{ch5:ass:p-bounded-central-moment-gradient_chapter_5} and Assumption~\ref{ch5:ass:q-bounded-central-moment-hessian_chapter_5}. Computing the $p$--th central moment of $g^{\star}_N$ gives, for all $x \in \R^N$
            \begin{thesischfivefitalignatnum}{2}
                \thesischfiveE{\norm{g^{\star}_N(x, \xi) - \nabla f^{\star}_N(x)}^p} \thesischfiveoversetlab{\eqref{ch5:6f55ce68-cc03-4d22-a262-4aa3367f948d-2_chapter_5}}&{=} \left( \alpha {\color{navy}\chi} \right)^p \thesischfiveE{\norm{g_K({\color{navy}\chi} x, \xi) - \nabla f_N({\color{navy}\chi} x)}^p} \\
                & \le \frac{2 \left( \gamma_{\infty} \alpha {\color{navy}\chi} \right)^p (1 - \theta)}{\theta^{p - 1}} \\
                & \le \frac{\left( 2 \gamma_{\infty} \alpha {\color{navy}\chi} \right)^p}{\theta^{p - 1}} \\
                & = \frac{\left( 4 \gamma_{\infty} \thesischfiveeps \right)^p}{\theta^{p - 1}}, \label{ch5:29bb8bb5-e6e5-482d-ab46-fab6357955ec-2_chapter_5}
            \end{thesischfivefitalignatnum}
        where the first inequality follows from Lemma~\ref{ch5:appdx-lem:properties-gradient-estimator-2_chapter_5}, in the second we use the fact that $0 < \theta \le 1$ and $p \ge 1$ so that $2 \le 2^p$, and in the last equality we use $\alpha = 2 \thesischfiveeps/{\color{navy}\chi}$. Moreover, following the same lines as in~\eqref{ch5:29bb8bb5-e6e5-482d-ab46-fab6357955ec-2_chapter_5} we have
            \begin{thesischfivefitalignatnum}{2}
                \thesischfiveE{\norm{\nabla g^{\star}_N(x, \xi) - \nabla^2 f^{\star}_N(x)}^q} \thesischfiveoversetlab{\eqref{ch5:6f55ce68-cc03-4d22-a262-4aa3367f948d-2_chapter_5}}&{=} \left( \alpha {\color{navy}\chi}^2 \right)^q \thesischfiveE{\norm{\nabla g_K({\color{navy}\chi} x, \xi) - \nabla^2 f_N({\color{navy}\chi} x)}^q} \\
                & \le \frac{2 \left( \ell_1 \alpha {\color{navy}\chi}^2 \right)^q (1 - \theta)}{\theta^{q - 1}} \\
                & \le \frac{\left( 2 \ell_1 \alpha {\color{navy}\chi}^2 \right)^q}{\theta^{q - 1}}. \label{ch5:29bb8bb5-e6e5-482d-ab46-fab6357955ec-3_chapter_5}
            \end{thesischfivefitalignatnum}
        
        {\color{navy}By~\eqref{ch5:29bb8bb5-e6e5-482d-ab46-fab6357955ec-2_chapter_5}, Assumptions~\ref{ch5:ass:p-bounded-central-moment-gradient_chapter_5} and~\ref{ch5:ass:q-bounded-central-moment-hessian_chapter_5} are satisfied by taking $\theta = \min\thesischfiveens{1, \theta_1}$ such that}
            \begin{thesischfivefitdisplay} \theta_1^{p - 1} \ge \left( \frac{4 \gamma_{\infty} \thesischfiveeps}{\sigma_1} \right)^p, \quad \text{ so } \quad \theta_1 \ge \left( \frac{4 \gamma_{\infty} \thesischfiveeps}{\sigma_1} \right)^{p/(p - 1)}, \thesischfivenumberthis\label{ch5:46a820b5-b4a2-4c73-b77f-41e02836f2e0-2_chapter_5} \end{thesischfivefitdisplay}
        hence
            \begin{thesischfivefitdisplay} \theta = \min\thesischfiveens{1, \left( \frac{4 \gamma_{\infty} \thesischfiveeps}{\sigma_1} \right)^{p/(p - 1)}}, \thesischfivenumberthis\label{ch5:b3242270-31d2-4d6a-a567-27370b26f605-4_chapter_5} \end{thesischfivefitdisplay}
        {\color{navy}If $\theta=1$, both the gradient and Hessian estimators are deterministic, so their centered moments are zero. If $\theta<1$, the two preceding inequalities give the required gradient and Hessian moment bounds.}
        For~\eqref{ch5:29bb8bb5-e6e5-482d-ab46-fab6357955ec-3_chapter_5} we need $0 < \alpha {\color{navy}\chi}^2 \le \sigma_2 \theta^{(q - 1)/q}/(2 \ell_1)$, so we fix ${\color{navy}\chi}$ such that
            \begin{align*}
                \alpha {\color{navy}\chi}^2 = 2 \thesischfiveeps {\color{navy}\chi} &\le \min\thesischfiveens{\frac{L_1}{\ell_1}, \frac{\sigma_2 \theta^{(q - 1)/q}}{2 \ell_1}} \quad \text{and}\\
                \alpha {\color{navy}\chi}^3 = 2 \thesischfiveeps {\color{navy}\chi}^2 &\le \frac{L_2}{\ell_2},
            \end{align*}
        that is to say
            \begin{thesischfivefitdisplay} {\color{navy}\chi} = \min\thesischfiveens{\frac{L_1}{2 \ell_1 \thesischfiveeps}, \sqrt{\frac{L_2}{2 \ell_2 \thesischfiveeps}}, \frac{\sigma_2 \theta^{(q - 1)/q}}{4 \ell_1 \thesischfiveeps}}. \thesischfivenumberthis\label{ch5:979cbf5c-3220-43c2-aaca-043c521969e6_chapter_5} \end{thesischfivefitdisplay}

        \item \textbf{Step 4:} \textit{Lower Bounding $\mathfrak{m}^{\normalfont\texttt{zr}}_{\thesischfiveeps}\left( {\color{navy}b}, \Delta, L_1, L_2, \sigma^p_1, \sigma_2^q \right)$}.

        {\color{navy}Continuing Step~2, Lemma~\ref{ch5:appdx-lem:zero-chain-progress_chapter_5} shows that, with probability at least $1/2$, for every integer $0 \le t \le (N - 1)/(2 \theta)$ and every $k \in [b]$,}
            \begin{thesischfivefitdisplay} \norm{\nabla f^{\star}_N\left( x^{(t, k)}_{\texttt{A}[\texttt{O}_f]} \right)} > 2\thesischfiveeps, \quad \text{ hence } \quad \thesischfiveE{\norm{\nabla f^{\star}_N\left( x^{(t, k)}_{\texttt{A}[\texttt{O}_f]} \right)}} > \thesischfiveeps, \end{thesischfivefitdisplay}
        from where it follows that
            \begin{thesischfivefitalignatnum}{1}
                &\mathfrak{m}^{\normalfont\texttt{zr}}_{\thesischfiveeps}\left( {\color{navy}b}, \Delta, L_1, L_2, \sigma^p_1, \sigma_2^q \right) > \frac{N - 1}{2 \theta}\\
                &\quad= \frac{1}{2 \theta} \left( \thesischfiveFloor{\frac{\Delta}{\alpha \Delta_0}} - 1 \right)\\
                &\quad= \frac{1}{2 \theta} \left( \thesischfiveFloor{\frac{\Delta {\color{navy}\chi}}{2 \Delta_0 \thesischfiveeps}} - 1 \right), \label{ch5:b3242270-31d2-4d6a-a567-27370b26f605-3_chapter_5}
            \end{thesischfivefitalignatnum}
        we then distinguish two cases:
            \begin{itemize}
                \item if $\Delta {\color{navy}\chi} / (2 \Delta_0 \thesischfiveeps) \ge 3$ then, using the inequality $\thesischfiveFloor{x} - 1 \ge x/2$, valid for all real number $x \ge 3$, we obtain
                    \begin{alignat}{2}
                        \MoveEqLeft[1] \mathfrak{m}^{\normalfont\texttt{zr}}_{\thesischfiveeps}\left( {\color{navy}b}, \Delta, L_1, L_2, \sigma^p_1, \sigma_2^q \right) \notag\\
                        \thesischfiveoversetlab{\eqref{ch5:b3242270-31d2-4d6a-a567-27370b26f605-3_chapter_5}}&{>} \frac{1}{4 \theta} \cdot \frac{\Delta {\color{navy}\chi}}{2 \Delta_0 \thesischfiveeps} \notag\\
                        \thesischfiveoversetlab{\eqref{ch5:979cbf5c-3220-43c2-aaca-043c521969e6_chapter_5}}&{=} \frac{1}{8 \Delta_0} \cdot \frac{\Delta}{\thesischfiveeps} \min\thesischfiveens{\frac{L_1}{2 \ell_1 \thesischfiveeps \theta}, \frac{1}{\theta} \sqrt{\frac{L_2}{2 \ell_2 \thesischfiveeps}}, \frac{\sigma_2 \theta^{-1/q}}{4 \ell_1 \thesischfiveeps}} \notag\\
                        & \ge \frac{1}{16 \Delta_0 \ell_1} \cdot \frac{\Delta}{\thesischfiveeps} \min\thesischfiveens{\frac{L_1}{\thesischfiveeps \theta}, \frac{1}{\theta} \sqrt{\frac{2 L_2}{\ell_2 \thesischfiveeps}}, \frac{\sigma_2 \theta^{-1/q}}{2 \thesischfiveeps}} \notag\\
                        \thesischfiveoversetlab{\eqref{ch5:b3242270-31d2-4d6a-a567-27370b26f605-4_chapter_5}}&{\ge} \frac{1}{32 \Delta_0 \ell_1} \cdot \frac{\Delta}{\thesischfiveeps} \min\biggl\{\frac{L_1}{\thesischfiveeps} \left( 1 + \left( \frac{\sigma_1}{4 \gamma_{\infty} \thesischfiveeps} \right)^{p/(p - 1)} \right), \notag\\
                        &\qquad \sqrt{\frac{2 L_2}{\ell_2 \thesischfiveeps}} \left( 1 + \left( \frac{\sigma_1}{4 \gamma_{\infty} \thesischfiveeps} \right)^{p/(p - 1)} \right),\ \frac{\sigma_2}{2 \thesischfiveeps} \left( 1 + \left( \frac{\sigma_1}{4 \gamma_{\infty} \thesischfiveeps} \right)^{p/(q (p - 1))} \right)\biggr\} \notag\\
                        & = \Omega(1) \cdot \min\biggl\{\frac{L_1 \Delta}{\thesischfiveeps^2} + \frac{L_1 \Delta}{\thesischfiveeps^2} \left( \frac{\sigma_1}{\thesischfiveeps} \right)^{p/(p - 1)}, \notag\\
                        &\qquad \frac{L_2^{1/2} \Delta}{\thesischfiveeps^{3/2}} + \frac{L_2^{1/2} \Delta}{\thesischfiveeps^{3/2}} \left( \frac{\sigma_1}{\thesischfiveeps} \right)^{p/(p - 1)}, \notag\\
                        &\qquad \frac{\Delta \sigma_2}{\thesischfiveeps^2} + \frac{\Delta \sigma_2}{\thesischfiveeps^2} \left( \frac{\sigma_1}{\thesischfiveeps} \right)^{p/(q (p - 1))}\biggr\}, \label{ch5:3b256059-e30a-4307-b1b7-a8e09886c0f1-2_chapter_5}
                    \end{alignat}
                where in the second inequality we use $\ell_1 \ge 1$ for the second entry of the minimum. Combining the lower bound~\eqref{ch5:3b256059-e30a-4307-b1b7-a8e09886c0f1-2_chapter_5} with~Lemma~\ref{ch5:appdx-lem:lower-bound-global-stichastic-model-mean-squared-smoothness-p-BCM_chapter_5}, using the choice for $c'$ provided below, gives:
                    \begin{thesischfivefitalignatnum}{2}
                        \MoveEqLeft[1] \mathfrak{m}^{\normalfont\texttt{zr}}_{\thesischfiveeps}\left( {\color{navy}b}, \Delta, L_1, L_2, \sigma_1^p, \sigma_2^q \right) \\
                        & \ge \Omega(1) \cdot \min\biggl\{\frac{L_1 \Delta}{\thesischfiveeps^2} + \frac{L_1 \Delta}{\thesischfiveeps^2} \left( \frac{\sigma_1}{\thesischfiveeps} \right)^{p/(p - 1)},\\
                        &\qquad \frac{L_2^{1/2} \Delta}{\thesischfiveeps^{3/2}} + \frac{L_2^{1/2} \Delta}{\thesischfiveeps^{3/2}} \left( \frac{\sigma_1}{\thesischfiveeps} \right)^{p/(p - 1)},\\
                        &\qquad \left( \frac{\sigma_1}{\thesischfiveeps} \right)^{p/(p - 1)} + \frac{\Delta \sigma_2}{\thesischfiveeps^2} + \frac{\Delta \sigma_2}{\thesischfiveeps^2} \left( \frac{\sigma_1}{\thesischfiveeps} \right)^{p/(q (p - 1))}\biggr\}, \label{ch5:5aea445e-d78f-404f-ac66-dbff30f0a00c-3_chapter_5}
                    \end{thesischfivefitalignatnum}
                since from the assumption $0 < \thesischfiveeps \le c' \min\thesischfiveens{\sqrt{L_1 \Delta}, L_2^{1/3} \Delta^{2/3}}$ we have
                    \begin{align*}
                        \frac{L_1 \Delta}{\thesischfiveeps^2} \left( \frac{\sigma_1}{\thesischfiveeps} \right)^{p/(p - 1)} &\gtrsim \left( \frac{\sigma_1}{\thesischfiveeps} \right)^{p/(p - 1)} \quad \text{and} \\
                        \frac{L_2^{1/2} \Delta}{\thesischfiveeps^{3/2}} \left( \frac{\sigma_1}{\thesischfiveeps} \right)^{p/(p - 1)} &\gtrsim \left( \frac{\sigma_1}{\thesischfiveeps} \right)^{p/(p - 1)},
                    \end{align*}
                {\color{navy}so the factor $(\sigma_1/\thesischfiveeps)^{p/(p - 1)}$ may be omitted from the first two entries of the minimum.}

                \item otherwise, if $\Delta {\color{navy}\chi} / (2 \Delta_0 \thesischfiveeps) < 3$, we choose the universal constant $c'$ as
                    \begin{thesischfivefitdisplay} 0 < c' = \min\thesischfiveens{\left( 12 \Delta_0 \ell_1 \right)^{-1/2}, \left( 6 \sqrt{2} \Delta_0 \ell_2^{1/2} \right)^{-2/3}} < \frac{1}{8}, \thesischfivenumberthis\label{ch5:e4896e63-043b-4b68-b23e-bb9b172efed4_chapter_5} \end{thesischfivefitdisplay}
                 then, using~Remark~\ref{ch5:appdx-rem:improve-lower-bound-global-stichastic-model-mean-squared-smoothness-p-BCM_chapter_5}, dividing both sides of $\Delta {\color{navy}\chi} / (2 \Delta_0 \thesischfiveeps) < 3$ by our choice of $\theta > 0$ yields
                    \begin{thesischfivefitalignatnum}{2}
                        \frac{1}{6 \Delta_0} \cdot \frac{\Delta {\color{navy}\chi}}{\thesischfiveeps} \theta^{-1} & \le \theta^{-1} = \max\thesischfiveens{1, \left( \frac{\sigma_1}{4 \gamma_{\infty} \thesischfiveeps} \right)^{p/(p - 1)}} \\
                        & \le \max\thesischfiveens{1, \left( \frac{\sigma_1}{\thesischfiveeps} \right)^{p/(p - 1)}}, \label{ch5:2b14b2b5-0770-4ca0-b177-61bcd5125db3-2_chapter_5}
                    \end{thesischfivefitalignatnum}
                where the last inequality follows from the fact that $p/(p - 1) > 0$ and $4 \gamma_{\infty} = 92 \ge 1$. Using the definition of ${\color{navy}\chi}$ we obtain
                    \begin{thesischfivefitalignat}{2}
                        \MoveEqLeft[1] \max\thesischfiveens{1, \left( \frac{\sigma_1}{\thesischfiveeps} \right)^{p/(p - 1)}} \\
                        \thesischfiveoversetlab{\eqref{ch5:2b14b2b5-0770-4ca0-b177-61bcd5125db3-2_chapter_5}+\eqref{ch5:979cbf5c-3220-43c2-aaca-043c521969e6_chapter_5}}&{\ge} \frac{1}{6 \Delta_0} \cdot \min\thesischfiveens{\frac{L_1 \Delta}{2 \ell_1 \theta \thesischfiveeps^2}, \frac{\Delta}{\theta \thesischfiveeps}\sqrt{\frac{L_2}{2 \ell_2 \thesischfiveeps}}, \frac{\Delta \sigma_2 \theta^{-1/q}}{4 \ell_1 \thesischfiveeps^2}} \\
                        \thesischfiveoversetlab{\eqref{ch5:b3242270-31d2-4d6a-a567-27370b26f605-4_chapter_5}}&{\ge} c_0 \min\biggl\{\frac{L_1 \Delta}{\thesischfiveeps^2} + \frac{L_1 \Delta}{\thesischfiveeps^2} \left( \frac{\sigma_1}{4 \gamma_{\infty} \thesischfiveeps} \right)^{p/(p - 1)},\\
                        &\qquad \frac{L_2^{1/2} \Delta}{\thesischfiveeps^{3/2}} + \frac{L_2^{1/2} \Delta}{\thesischfiveeps^{3/2}} \left( \frac{\sigma_1}{4 \gamma_{\infty} \thesischfiveeps} \right)^{p/(p - 1)},\\
                        &\qquad \frac{\Delta \sigma_2}{\thesischfiveeps^2} + \frac{\Delta \sigma_2}{\thesischfiveeps^2} \left( \frac{\sigma_1}{4 \gamma_{\infty} \thesischfiveeps} \right)^{p/(q (p - 1))}\biggr\} \\
                        & \ge C_{p, q} \min\biggl\{\frac{L_1 \Delta}{\thesischfiveeps^2} + \frac{L_1 \Delta}{\thesischfiveeps^2} \left( \frac{\sigma_1}{\thesischfiveeps} \right)^{p/(p - 1)},\\
                        &\qquad \frac{L_2^{1/2} \Delta}{\thesischfiveeps^{3/2}} + \frac{L_2^{1/2} \Delta}{\thesischfiveeps^{3/2}} \left( \frac{\sigma_1}{\thesischfiveeps} \right)^{p/(p - 1)},\\
                        &\qquad \frac{\Delta \sigma_2}{\thesischfiveeps^2} + \frac{\Delta \sigma_2}{\thesischfiveeps^2} \left( \frac{\sigma_1}{\thesischfiveeps} \right)^{p/(q (p - 1))}\biggr\},
                    \end{thesischfivefitalignat}
                where $c_0 \eqdef 1/(48 \Delta_0 \max\{\ell_1, \ell_2^{1/2}\})$ {\color{navy}and $C_{p,q}>0$ is a universal constant depending only on $p$ and $q$. Moreover, by~Lemma~\ref{ch5:appdx-lem:lower-bound-global-stichastic-model-mean-squared-smoothness-p-BCM_chapter_5}, and more precisely Remark~\ref{ch5:appdx-rem:improve-lower-bound-global-stichastic-model-mean-squared-smoothness-2-hessian_chapter_5}, the condition $0 < \thesischfiveeps \le \sqrt{L_1 \Delta}/8$ implies that there exists a universal constant $C_p > 0$ (depending only on $p$) such that}
                    \begin{thesischfivefitdisplay} \mathfrak{m}^{\normalfont\texttt{zr}}_{\thesischfiveeps}\left( {\color{navy}b}, \Delta, L_1, L_2, \sigma_1^p, \sigma_2^q \right) \ge C_p \max\thesischfiveens{1, \left( \frac{\sigma_1}{\thesischfiveeps} \right)^{p/(p - 1)}}, \end{thesischfivefitdisplay}
                hence
                    \begin{thesischfivefitalignatnum}{2}
                        \MoveEqLeft[1] \mathfrak{m}^{\normalfont\texttt{zr}}_{\thesischfiveeps}\left( {\color{navy}b}, \Delta, L_1, L_2, \sigma_1^p, \sigma_2^q \right) \notag\\
                        & \ge \Omega(1) \cdot \biggl( \left( \frac{\sigma_1}{\thesischfiveeps} \right)^{p/(p - 1)} + \min\biggl\{\frac{L_1 \Delta}{\thesischfiveeps^2} + \frac{L_1 \Delta}{\thesischfiveeps^2} \left( \frac{\sigma_1}{\thesischfiveeps} \right)^{p/(p - 1)}, \notag\\
                        &\qquad \frac{L_2^{1/2} \Delta}{\thesischfiveeps^{3/2}} + \frac{L_2^{1/2} \Delta}{\thesischfiveeps^{3/2}} \left( \frac{\sigma_1}{\thesischfiveeps} \right)^{p/(p - 1)}, \notag\\
                        &\qquad \frac{\Delta \sigma_2}{\thesischfiveeps^2} + \frac{\Delta \sigma_2}{\thesischfiveeps^2} \left( \frac{\sigma_1}{\thesischfiveeps} \right)^{p/(q (p - 1))}\biggr\} \biggr) \notag\\
                        & \ge \Omega(1) \cdot \min\biggl\{\frac{L_1 \Delta}{\thesischfiveeps^2} + \frac{L_1 \Delta}{\thesischfiveeps^2} \left( \frac{\sigma_1}{\thesischfiveeps} \right)^{p/(p - 1)}, \notag\\
                        &\qquad \frac{L_2^{1/2} \Delta}{\thesischfiveeps^{3/2}} + \frac{L_2^{1/2} \Delta}{\thesischfiveeps^{3/2}} \left( \frac{\sigma_1}{\thesischfiveeps} \right)^{p/(p - 1)}, \notag\\
                        &\qquad \left( \frac{\sigma_1}{\thesischfiveeps} \right)^{p/(p - 1)} + \frac{\Delta \sigma_2}{\thesischfiveeps^2} + \frac{\Delta \sigma_2}{\thesischfiveeps^2} \left( \frac{\sigma_1}{\thesischfiveeps} \right)^{p/(q (p - 1))}\biggr\}, \label{ch5:5aea445e-d78f-404f-ac66-dbff30f0a00c-4_chapter_5}
                    \end{thesischfivefitalignatnum}
                as $0 < \thesischfiveeps \le c' \min\thesischfiveens{\sqrt{L_1 \Delta}, L_2^{1/3} \Delta^{2/3}}$ implies
                    \begin{align*}
                        \frac{L_1 \Delta}{\thesischfiveeps^2} \left( \frac{\sigma_1}{\thesischfiveeps} \right)^{p/(p - 1)} &\gtrsim \left( \frac{\sigma_1}{\thesischfiveeps} \right)^{p/(p - 1)} \quad \text{and} \\
                        \frac{L_2^{1/2} \Delta}{\thesischfiveeps^{3/2}} \left( \frac{\sigma_1}{\thesischfiveeps} \right)^{p/(p - 1)} &\gtrsim \left( \frac{\sigma_1}{\thesischfiveeps} \right)^{p/(p - 1)},
                    \end{align*}
                {\color{navy}so the factor $(\sigma_1/\thesischfiveeps)^{p/(p - 1)}$ may be omitted from the first two entries of the minimum.}
                    
            \end{itemize}

            \item \textbf{Step 5:} \textit{A Last Bound: the Case $\theta = 1$}.

            Observe that, if instead of taking $\theta$ as in~\eqref{ch5:b3242270-31d2-4d6a-a567-27370b26f605-4_chapter_5}, we choose directly $\theta = 1$ then, thanks to~Lemma~\ref{ch5:appdx-lem:properties-gradient-estimator-2_chapter_5} we immediately have
                \begin{thesischfivefitdisplay} \thesischfiveE{\norm{g^{\star}_N(x, \xi) - \nabla f^{\star}_N(x)}^p} = 0, \end{thesischfivefitdisplay}
            and
                \begin{thesischfivefitdisplay} \thesischfiveE{\norm{\nabla g^{\star}_N(x, \xi) - \nabla^2 f^{\star}_N(x)}^q} = 0, \end{thesischfivefitdisplay}
            so~Assumption~\ref{ch5:ass:p-bounded-central-moment-gradient_chapter_5} and~Assumption~\ref{ch5:ass:q-bounded-central-moment-hessian_chapter_5} are satisfied. Hence, if we set 
                \begin{thesischfivefitdisplay} N = \thesischfiveFloor{\frac{\Delta}{\alpha \Delta_0}}, \quad \alpha = \frac{2 \thesischfiveeps}{{\color{navy}\chi}}, \quad \text{ and } \quad {\color{navy}\chi} = \min\thesischfiveens{\frac{L_1}{2 \ell_1 \thesischfiveeps}, \sqrt{\frac{L_2}{2 \ell_2 \thesischfiveeps}}}, \end{thesischfivefitdisplay}
            then $f_N^{\star} \in \mathcal{F}(\Delta, L_1, L_2)$ and the inequality~\eqref{ch5:f01cde9c-e4f8-4c8b-b9b4-7e86b18da9c1_chapter_5} is satisfied. {\color{navy}Hence, with probability at least $1/2$, for every integer $0 \le t \le (N - 1)/(2 \theta)$ and every $k \in [b]$,}
            \begin{thesischfivefitdisplay} \norm{\nabla f^{\star}_N\left( x^{(t, k)}_{\texttt{A}[\texttt{O}_f]} \right)} > 2\thesischfiveeps, \quad \text{ hence } \quad \thesischfiveE{\norm{\nabla f^{\star}_N\left( x^{(t, k)}_{\texttt{A}[\texttt{O}_f]} \right)}} > \thesischfiveeps, \end{thesischfivefitdisplay}
        from where we obtain also
            \begin{thesischfivefitalignatnum}{1}
                &\mathfrak{m}^{\normalfont\texttt{zr}}_{\thesischfiveeps}\left( {\color{navy}b}, \Delta, L_1, L_2, \sigma^p_1, \sigma_2^q \right) > \frac{N - 1}{2 \theta}\\
                &\quad= \frac{1}{2 \theta} \left( \thesischfiveFloor{\frac{\Delta}{\alpha \Delta_0}} - 1 \right)\\
                &\quad= \frac{1}{2 \theta} \left( \thesischfiveFloor{\frac{\Delta {\color{navy}\chi}}{2 \Delta_0 \thesischfiveeps}} - 1 \right), \label{ch5:b3242270-31d2-4d6a-a567-27370b26f605-bis-bis-bis_chapter_5}
            \end{thesischfivefitalignatnum}
        and by our assumption on $\thesischfiveeps$ with the choice of $c'$ in~\eqref{ch5:e4896e63-043b-4b68-b23e-bb9b172efed4_chapter_5}, we have
            \begin{align*}
                0 < \thesischfiveeps &< c' \min\thesischfiveens{\sqrt{L_1 \Delta}, L_2^{1/3} \Delta^{2/3}} \\
                &\le \min\thesischfiveens{\sqrt{\frac{L_1 \Delta}{12 \Delta_0 \ell_1}}, \left( \frac{L_2^{1/2} \Delta}{6 \sqrt{2} \Delta_0 \ell_2^{1/2}} \right)^{2/3}},
            \end{align*}
        which is enough to imply the inequality
            \begin{thesischfivefitdisplay} \min\thesischfiveens{\frac{L_1 \Delta}{4 \Delta_0 \ell_1 \thesischfiveeps^2}, \frac{L_2^{1/2} \Delta}{2 \sqrt{2} \Delta_0 \ell_2^{1/2} \thesischfiveeps^{3/2}}} = \frac{\Delta {\color{navy}\chi}}{2 \Delta_0 \thesischfiveeps} \ge 3. \end{thesischfivefitdisplay}
        Thus, since $\theta = 1$, we have
            \begin{thesischfivefitalignat}{2}
                \MoveEqLeft[1] \mathfrak{m}^{\normalfont\texttt{zr}}_{\thesischfiveeps}\left( {\color{navy}b}, \Delta, L_1, L_2, \sigma_1^p, \sigma_2^q \right) \\
                \thesischfiveoversetlab{\eqref{ch5:b3242270-31d2-4d6a-a567-27370b26f605-bis-bis-bis_chapter_5}}&{>} \frac{1}{4 \theta} \cdot \frac{\Delta {\color{navy}\chi}}{2 \Delta_0 \thesischfiveeps} \\
                & = \frac{1}{8 \Delta_0} \cdot \frac{\Delta}{\thesischfiveeps}\min\thesischfiveens{\frac{L_1}{2 \ell_1 \thesischfiveeps}, \sqrt{\frac{L_2}{2 \ell_2 \thesischfiveeps}}} \\
                & = \Omega(1) \cdot \min\thesischfiveens{\frac{L_1 \Delta}{\thesischfiveeps^2}, \frac{L_2^{1/2} \Delta}{\thesischfiveeps^{3/2}}},
            \end{thesischfivefitalignat}
        and, combining this bound with~\eqref{ch5:5aea445e-d78f-404f-ac66-dbff30f0a00c-3_chapter_5} and~\eqref{ch5:5aea445e-d78f-404f-ac66-dbff30f0a00c-4_chapter_5} respectively leads to the desired result, i.e.,
            \begin{thesischfivefitalignat}{2}
                \MoveEqLeft[1] \mathfrak{m}^{\normalfont\texttt{zr}}_{\thesischfiveeps}\left( {\color{navy}b}, \Delta, L_1, L_2, \sigma^p_1, \sigma_2^q \right) \\
                & \ge \Omega(1) \cdot \min\biggl\{\frac{L_1 \Delta}{\thesischfiveeps^2} + \frac{L_1 \Delta}{\thesischfiveeps^2} \left( \frac{\sigma_1}{\thesischfiveeps} \right)^{p/(p - 1)},\ \frac{L_2^{1/2} \Delta}{\thesischfiveeps^{3/2}} + \frac{L_2^{1/2} \Delta}{\thesischfiveeps^{3/2}} \left( \frac{\sigma_1}{\thesischfiveeps} \right)^{p/(p - 1)}, \\
                &\qquad \min\thesischfiveens{\frac{L_1 \Delta}{\thesischfiveeps^2}, \frac{L_2^{1/2} \Delta}{\thesischfiveeps^{3/2}}} + \frac{\Delta \sigma_2}{\thesischfiveeps^2} + \frac{\Delta \sigma_2}{\thesischfiveeps^2} \left( \frac{\sigma_1}{\thesischfiveeps} \right)^{p/(q (p - 1))} + \left( \frac{\sigma_1}{\thesischfiveeps} \right)^{p/(p - 1)}\biggr\}.
            \end{thesischfivefitalignat}
    \end{itemize}
\end{proof}

\newpage
\section{Additional Proofs for Section~\ref{ch5:sec:optimal-method-nsgd-mvr_chapter_5}}\label{ch5:appdx:proof-nsgd-mvr_chapter_5}

\subsection{Auxiliary Lemmas}

\begin{lemma}[A Descent Lemma]\label{ch5:appdx-lem:descent-lemma-nsgd-mvr_chapter_5}
    Under~Assumption~\ref{ch5:ass:mean-squared-smoothness_chapter_5} and Assumption~\ref{ch5:ass:p-bounded-central-moment-gradient_chapter_5}, for any choice of stepsize $\gamma > 0$, and for any $k \in \thesischfiveInt{0}{K - 1}$ we have
        \begin{thesischfivefitdisplay} f(x^{k+1}) \le f(x^{k}) + 2 \gamma \norm{\hat{e}_k} - \gamma \norm{\nabla f(x^{k})} + \frac{\gamma^2 \bar{L}}{2}, \thesischfivenumberthis\label{ch5:2331e9e8-3c3e-4a5b-b3b6-50d71097a688_chapter_5} \end{thesischfivefitdisplay}
    where $\hat{e}_k \eqdef g^{k} - \nabla f(x^{k})$ is the error term.
\end{lemma}

\begin{proof}
    According to~Assumption~\ref{ch5:ass:mean-squared-smoothness_chapter_5} we know that the function $f$ is $\bar{L}$--smooth~\parencite{nesterov2018lectures} since by Jensen's inequality (Lemma~\ref{ch5:lem:jensen-inequality_chapter_5}) applied on the convex function $x \mapsto \norm{x}^q$ (since $q \ge 1$)
        \begin{thesischfivefitalignatnum}{2}
            \norm{\nabla f(x) - \nabla f(y)}^q & = \norm{\thesischfiveExpSub{\xi \sim \thesischfivecD}{\nabla f(x;\xi) - \nabla f(y;\xi)}}^q \\
            & \le \thesischfiveExpSub{\xi \sim \thesischfivecD}{\norm{\nabla f(x;\xi) - \nabla f(y;\xi)}^q} \\
            & \le \bar{L}^q \norm{x - y}^q, \label{ch5:94828cd2-cf02-492b-a6ec-40b3012f66c1_chapter_5}
        \end{thesischfivefitalignatnum}
    where the equality follows from Assumption~\ref{ch5:ass:p-bounded-central-moment-gradient_chapter_5} and the last inequality from Assumption~\ref{ch5:ass:mean-squared-smoothness_chapter_5}. Thus it holds
        \begin{thesischfivefitalignat}{2}
            f(x^{k+1}) & \le f(x^{k}) + \thesischfiveps{\nabla f(x^{k})}{x^{k+1} - x^{k}} + \frac{\bar{L}}{2}\sqnorm{x^{k+1} - x^{k}} \\
            & = f(x^{k}) - \gamma \thesischfiveps{\nabla f(x^{k})}{\unitdir{g^{k}}} + \frac{\gamma^2 \bar{L}}{2} \\
            & = f(x^{k}) - \gamma \thesischfiveps{\nabla f(x^{k}) - g^{k}}{\unitdir{g^{k}}} - \gamma \norm{g^{k}} + \frac{\gamma^2 \bar{L}}{2} \\
            & \le f(x^{k}) + \gamma \norm{\nabla f(x^{k}) - g^{k}} - \gamma \norm{g^{k}} + \frac{\gamma^2 \bar{L}}{2} \\
            & \le f(x^{k}) + 2 \gamma \norm{\nabla f(x^{k}) - g^{k}} - \gamma \norm{\nabla f(x^{k})} + \frac{\gamma^2 \bar{L}}{2} \\
            & = f(x^{k}) + 2 \gamma \norm{\hat{e}_k} - \gamma \norm{\nabla f(x^{k})} + \frac{\gamma^2 \bar{L}}{2},
        \end{thesischfivefitalignat}
    where the first inequality uses the $\bar{L}$--smoothness of $f$, in the first equality we use the update rule {\color{darkred}$x^{k+1}=x^k-\gamma g^k/\norm{g^k}$}, with the zero-numerator convention stated in Chapter~5, the second inequality is the Cauchy--Schwarz inequality (Lemma~\ref{ch5:lem:cauchy-schwarz_chapter_5}), and in the third we use the triangle inequality.
\end{proof}

\begin{lemma}[Another Descent Lemma]\label{ch5:appdx-lem:descent-lemma-nsgd-mvr-hess_chapter_5}
    Under~Assumption~\ref{ch5:ass:L-lipschitz-gradients_chapter_5} and Assumption~\ref{ch5:ass:p-bounded-central-moment-gradient_chapter_5}, for any choice of stepsize $\gamma > 0$, and for any $k \in \thesischfiveInt{0}{K - 1}$ we have
        \begin{thesischfivefitdisplay} f(x^{k+1}) \le f(x^{k}) + 2 \gamma \norm{\hat{e}_k} - \gamma \norm{\nabla f(x^{k})} + \frac{\gamma^2 L_1}{2}, \thesischfivenumberthis\label{ch5:2331e9e8-3c3e-4a5b-b3b6-50d71097a688-2_chapter_5} \end{thesischfivefitdisplay}
    where $\hat{e}_k \eqdef g^{k} - \nabla f(x^{k})$ is the error term.
\end{lemma}

\begin{proof}
    The proof is the same as in the previous descent lemma (Lemma~\ref{ch5:appdx-lem:descent-lemma-nsgd-mvr_chapter_5}) where now $f$ has $L_1$--Lipschitz continuous gradients (instead of $\bar{L}$).
\end{proof}

\begin{lemma}[Unrolling the Descent Lemma]\label{ch5:appdx-lem:unrolling-descent-lemma-nsgd-mvr_chapter_5}
    Under~Assumption~\ref{ch5:ass:lower-boundedness_chapter_5}, Assumption~\ref{ch5:ass:p-bounded-central-moment-gradient_chapter_5}, and Assumption~\ref{ch5:ass:mean-squared-smoothness_chapter_5}, for any choice of stepsize $\gamma > 0$ the iterates $\{x^{k}\}_{k \in \thesischfiveInt{0}{K}}$ produced by~\ref{ch5:algo:nsgd-mvr_chapter_5} satisfy
        \begin{thesischfivefitdisplay} \frac{1}{K} \sum_{k=0}^{K - 1} \norm{\nabla f(x^{k})} \le \frac{\Delta}{\gamma K} + \frac{2}{K} \sum_{k=0}^{K - 1} \norm{\hat{e}_k} + \frac{\gamma \bar{L}}{2}. \end{thesischfivefitdisplay}
\end{lemma}

\begin{proof}
    From the previous descent lemma (Lemma~\ref{ch5:appdx-lem:descent-lemma-nsgd-mvr_chapter_5}), summing inequality~\eqref{ch5:2331e9e8-3c3e-4a5b-b3b6-50d71097a688_chapter_5} over $k \in \thesischfiveInt{0}{K - 1}$ gives
        \begin{thesischfivefitdisplay} \gamma \sum_{k=0}^{K - 1} \norm{\nabla f(x^{k})} \le f(x^0) - f(x^K) + 2 \gamma \sum_{k=0}^{K - 1} \norm{\hat{e}_k} + \frac{\gamma^2 \bar{L} K}{2}, \thesischfivenumberthis\label{ch5:eaec4407-2c6f-461b-805b-b59a34c6eee9_chapter_5} \end{thesischfivefitdisplay}
    where we telescope the terms $f(x^{k}) - f(x^{k+1})$. Multiplying both sides of~\eqref{ch5:eaec4407-2c6f-461b-805b-b59a34c6eee9_chapter_5} by $1/(\gamma K)$ leads to
        \begin{thesischfivefitdisplay} \frac{1}{K} \sum_{k=0}^{K - 1} \norm{\nabla f(x^{k})} \le \frac{1}{\gamma K} \left( f(x^0) - f(x^K)  \right) + \frac{2}{K} \sum_{k=0}^{K - 1} \norm{\hat{e}_k} + \frac{\gamma \bar{L}}{2}, \end{thesischfivefitdisplay}
    and using~Assumption~\ref{ch5:ass:lower-boundedness_chapter_5} we obtain
        \begin{thesischfivefitdisplay} \frac{1}{K} \sum_{k=0}^{K - 1} \norm{\nabla f(x^{k})} \le \frac{\Delta}{\gamma K} + \frac{2}{K} \sum_{k=0}^{K - 1} \norm{\hat{e}_k} + \frac{\gamma \bar{L}}{2}, \end{thesischfivefitdisplay}
    as desired.
\end{proof}

If we assume the function $f$ has $L_1$--Lipschitz continuous gradients (Assumption~\ref{ch5:ass:L-lipschitz-gradients_chapter_5}) then~Lemma~\ref{ch5:appdx-lem:descent-lemma-nsgd-mvr-hess_chapter_5} holds and we can unroll it in the same way as we did above. For that reason, we only state the result and we omit the proof.
\begin{lemma}[Unrolling the Descent Lemma]\label{ch5:appdx-lem:unrolling-descent-lemma-nsgd-mvr-hess_chapter_5}
    Under~Assumption~\ref{ch5:ass:lower-boundedness_chapter_5}, Assumption~\ref{ch5:ass:L-lipschitz-gradients_chapter_5}, and Assumption~\ref{ch5:ass:p-bounded-central-moment-gradient_chapter_5}, for any choice of stepsize $\gamma > 0$ the iterates $\{x^{k}\}_{k \in \thesischfiveInt{0}{K}}$ produced by~\ref{ch5:algo:nsgd-mvr_chapter_5} satisfy
        \begin{thesischfivefitdisplay} \frac{1}{K} \sum_{k=0}^{K - 1} \norm{\nabla f(x^{k})} \le \frac{\Delta}{\gamma K} + \frac{2}{K} \sum_{k=0}^{K - 1} \norm{\hat{e}_k} + \frac{\gamma L_1}{2}. \end{thesischfivefitdisplay}
\end{lemma}

\begin{lemma}[Bounding the Error Term]\label{ch5:appdx-lem:bounding-error-term-nsgd-mvr_chapter_5}
    Under~Assumption~\ref{ch5:ass:p-bounded-central-moment-gradient_chapter_5} and Assumption~\ref{ch5:ass:mean-squared-smoothness_chapter_5}, for all $k \in \thesischfiveInt{0}{K - 1}$ we have
        \begin{thesischfivefitdisplay} \thesischfiveE{\norm{\hat{e}_k}} \le (1 - \eta)^k \thesischfiveE{\norm{\hat{e}_0}} + 2 \sigma_1 \eta^{(p - 1)/p} + 4 \gamma \bar{L} \eta^{-1/q}, \end{thesischfivefitdisplay}
    where $\hat{e}_k \eqdef g^{k} - \nabla f(x^{k})$.
\end{lemma}

\begin{proof}
    {\color{navy}By the update rule of Algorithm~\ref{ch5:algo:nsgd-mvr_chapter_5}, for every $k\in\thesischfiveInt{1}{K-1}$, the estimator $g^k$ and the error term $\hat e_k = g^k - \nabla f(x^k)$ satisfy
        \begin{align}
            g^k &= (1-\eta) \left( g^{k-1} + \nabla f(x^k;\xi^k) - \nabla f(x^{k-1};\xi^k) \right) + \eta \nabla f(x^k;\xi^k), \notag\\
            \hat e_k &= (1-\eta) \hat e_{k-1} + \eta e_k - (1-\eta) \hat S_k,
            \label{ch5:2f56c833-5913-40e8-a00c-20c7b1a4b32b_chapter_5}
        \end{align}
        where
        \begin{align*}
            e_k &\eqdef \nabla f(x^k;\xi^k) - \nabla f(x^k), \\
            \hat S_k &\eqdef [\nabla f(x^k) - \nabla f(x^{k-1})] - [\nabla f(x^k;\xi^k) - \nabla f(x^{k-1};\xi^k)].
        \end{align*}
        Set $\mathcal F_0\eqdef\sigma(g^0)$ and
        $\mathcal F_k\eqdef\sigma(g^0,\xi^1,\ldots,\xi^k)$ for $k\ge1$.
        Since $x^{k-1}$ and $x^k$ are $\mathcal F_{k-1}$-measurable and
        $\xi^k$ is fresh, both $e_k$ and $\hat S_k$ are
        $\mathcal F_k$-measurable and satisfy
        $\Exp{e_k\mid\mathcal F_{k-1}}=0$ and
        $\Exp{\hat S_k\mid\mathcal F_{k-1}}=0$.}

        Then, unrolling the recursion~\eqref{ch5:2f56c833-5913-40e8-a00c-20c7b1a4b32b_chapter_5} gives
            \begin{thesischfivefitdisplay} \hat{e}_k = (1 - \eta)^k \hat{e}_0 + \eta \sum_{j = 0}^{k - 1} (1 - \eta)^{k - j - 1} e_{j + 1} - \sum_{j = 0}^{k - 1} (1 - \eta)^{k  - j} \hat{S}_{j + 1}, \end{thesischfivefitdisplay}
        and taking the norm followed by the total expectation yields
            \begin{thesischfivefitalignat}{2}
                \thesischfiveE{\norm{\hat{e}_k}} & \le (1 - \eta)^k \thesischfiveE{\norm{\hat{e}_0}} + \eta \thesischfiveE{\norm{ \sum_{j = 0}^{k - 1} (1 - \eta)^{k - j - 1} e_{j + 1}}} \\
                &\quad + (1 - \eta) \thesischfiveE{\norm{\sum_{j = 0}^{k - 1} (1 - \eta)^{k  - j - 1} \hat{S}_{j + 1}}}.
            \end{thesischfivefitalignat}
        We now need to upper bound the last two terms of the previous inequality. For the first term we have
            \begin{thesischfivefitalignatnum}{2}
                \thesischfiveE{\norm{ \sum_{j = 0}^{k - 1} (1 - \eta)^{k - j - 1} e_{j + 1}}} & \le \left( \thesischfiveE{\norm{ \sum_{j = 0}^{k - 1} (1 - \eta)^{k - j - 1} e_{j + 1}}^p} \right)^{1/p} \\
                & \le \left( 2 \sum_{j = 0}^{k - 1} (1 - \eta)^{p (k - j -1)} \thesischfiveE{\norm{e_{j + 1}}^p} \right)^{1/p} \\
                & \le 2 \left( \sum_{j = 0}^{k - 1} (1 - \eta)^{p (k - j -1)} \sigma_1^p \right)^{1/p} \\
                & \le 2 \sigma_1 \left( \sum_{j = 0}^{k - 1} (1 - \eta)^{k - j - 1} \right)^{1/p} \\
                & \le 2 \sigma_1 \left( \sum_{j \ge 0} (1 - \eta)^j \right)^{1/p} \\
                & = 2 \sigma_1 \eta^{-1/p}, \label{ch5:cea5881d-8688-41f3-8083-df58a5981198_chapter_5}
            \end{thesischfivefitalignatnum}
        where the first inequality is Jensen's inequality (Lemma~\ref{ch5:lem:jensen-inequality_chapter_5}), the second follows from Lemma~\ref{ch5:appdx-technical-lem:app-von-Bahr-and-Essen_chapter_5} since the $e_{j+1}$ are martingale differences, the third from Assumption~\ref{ch5:ass:p-bounded-central-moment-gradient_chapter_5} and $2^{1/p} \le 2$, and the fourth from $(1 - \eta)^p \le 1 - \eta$. For the last term we have
            \begin{thesischfivefitalignatnum}{2}
                \thesischfiveE{\norm{\sum_{j = 0}^{k - 1} (1 - \eta)^{k  - j - 1} \hat{S}_{j + 1}}} & \le \left( \thesischfiveE{\norm{\sum_{j = 0}^{k - 1} (1 - \eta)^{k  - j - 1} \hat{S}_{j + 1}}^q} \right)^{1/q} \\
                & \le \left( 2 \sum_{j = 0}^{k - 1} (1 - \eta)^{q (k - j - 1)} \thesischfiveE{\norm{\hat{S}_{j + 1}}^q} \right)^{1/q}, \label{ch5:d98da1f4-086e-4004-a809-2eb471da42a6_chapter_5}
            \end{thesischfivefitalignatnum}
        by Jensen's inequality and Lemma~\ref{ch5:appdx-technical-lem:app-von-Bahr-and-Essen_chapter_5} again. Moreover,
            \begin{thesischfivefitalignatnum}{2}
                \MoveEqLeft[1] \thesischfiveE{\norm{\hat{S}_{j + 1}}^q} \\
                & = \thesischfiveE{\norm{\left[ \nabla f(x^{j+1}) - \nabla f(x^{j}) \right] - \left[ \nabla f(x^{j+1};\xi^{j+1}) - \nabla f(x^{j};\xi^{j+1}) \right]}^q} \\
                & \le 2^{q - 1} \Exp\Bigl[\norm{\nabla f(x^{j+1}) - \nabla f(x^{j})}^q \\
                &\qquad + \norm{\nabla f(x^{j+1};\xi^{j+1}) - \nabla f(x^{j};\xi^{j+1})}^q\Bigr] \\
                & \le 2^q \bar{L}^q \thesischfiveE{\norm{x^{j+1} - x^{j}}^q} \\
                & \le 2^q \gamma^q \bar{L}^q, \label{ch5:f79e21c1-04c0-497e-81e9-c10afd7f04f9_chapter_5}
            \end{thesischfivefitalignatnum}
        where the first inequality follows from Lemma~\ref{ch5:lem:norm-power-alpha-inequality_chapter_5}, the second from~\eqref{ch5:94828cd2-cf02-492b-a6ec-40b3012f66c1_chapter_5} and Assumption~\ref{ch5:ass:mean-squared-smoothness_chapter_5}, and the last from $\norm{x^{j+1} - x^{j}} \le \gamma$. Thus, using~\eqref{ch5:f79e21c1-04c0-497e-81e9-c10afd7f04f9_chapter_5} and $2^{1/q} \le 2$ we obtain
             \begin{thesischfivefitalignatnum}{2}
                \thesischfiveE{\norm{\sum_{j = 0}^{k - 1} (1 - \eta)^{k  - j - 1} \hat{S}_{j + 1}}} \thesischfiveoversetlab{\eqref{ch5:d98da1f4-086e-4004-a809-2eb471da42a6_chapter_5}}&{\le} 4 \gamma \bar{L} \left( \sum_{j = 0}^{k - 1} (1 - \eta)^{q (k - j - 1)} \right)^{1/q} \\
                & \le 4 \gamma \bar{L} \eta^{-1/q}. \label{ch5:91e3e64b-54a6-46f3-96a5-a9b4f24dc6ee_chapter_5}
            \end{thesischfivefitalignatnum}

        Then, combining the bounds~\eqref{ch5:cea5881d-8688-41f3-8083-df58a5981198_chapter_5} and~\eqref{ch5:91e3e64b-54a6-46f3-96a5-a9b4f24dc6ee_chapter_5} we have
            \begin{thesischfivefitdisplay} \thesischfiveE{\norm{\hat{e}_k}} \le (1 - \eta)^k \thesischfiveE{\norm{\hat{e}_0}} + 2 \sigma_1 \eta^{(p - 1)/p} + 4 \gamma \bar{L} \eta^{-1/q}, \end{thesischfivefitdisplay}
        {\color{navy}which completes the proof of the lemma.}
\end{proof}

\begin{remark}
    In~\eqref{ch5:f79e21c1-04c0-497e-81e9-c10afd7f04f9_chapter_5} we only use the $q$--weak average smoothness assumption (Assumption~\ref{ch5:ass:mean-squared-smoothness_chapter_5}) to achieve the bound $2^q \gamma^q \bar{L}^q$. {\color{navy}This assumption may instead be replaced by the combination of $\bar{L}$--Lipschitz continuity of $\nabla f$ and}
        \begin{thesischfivefitdisplay} \thesischfiveExpSub{\xi \sim \thesischfivecD}{\norm{\left[ \nabla f(x;\xi) - \nabla f(y;\xi) \right] - \left[ \nabla f(x) - \nabla f(y) \right]}^q} \le \delta^q \norm{x - y}^q, \end{thesischfivefitdisplay}
    for all $x, y \in \R^d$, where {\color{navy}$\delta \ge 0$ can be much smaller than $\bar L$} (see~Assumption~\ref{ch5:ass:mean-squared-smoothness-2_chapter_5}).
\end{remark}

\begin{corollary}[Bounding the Error Term: a Refined Version]\label{ch5:appdx-lem:bounding-error-term-nsgd-mvr-refined_chapter_5}
    Under~Assumption~\ref{ch5:ass:p-bounded-central-moment-gradient_chapter_5} and Assumption~\ref{ch5:ass:mean-squared-smoothness-2_chapter_5}, for all $k \in \thesischfiveInt{0}{K - 1}$ we have
        \begin{thesischfivefitdisplay} \thesischfiveE{\norm{\hat{e}_k}} \le (1 - \eta)^k \thesischfiveE{\norm{\hat{e}_0}} + 2 \sigma_1 \eta^{(p - 1)/p} + 2 \gamma \delta \eta^{-1/q}, \thesischfivenumberthis\label{ch5:b9474e8b-ac2e-49ce-9636-37a177622ba3_chapter_5} \end{thesischfivefitdisplay}
    where $\hat{e}_k \eqdef g^{k} - \nabla f(x^{k})$.
\end{corollary}

\begin{proof}
    {\color{navy}The first two terms in~\eqref{ch5:b9474e8b-ac2e-49ce-9636-37a177622ba3_chapter_5} are obtained exactly as in Lemma~\ref{ch5:appdx-lem:bounding-error-term-nsgd-mvr_chapter_5}. For the last term, we begin as in~\eqref{ch5:d98da1f4-086e-4004-a809-2eb471da42a6_chapter_5} and apply Assumption~\ref{ch5:ass:mean-squared-smoothness-2_chapter_5}:}
        \begin{align*}
            \color{navy}\thesischfiveE{\norm{\hat{S}_{j + 1}}^q} &\color{navy}= \thesischfiveE{\norm{\left[ \nabla f(x^{j+1}) - \nabla f(x^{j}) \right] - \left[ \nabla f(x^{j+1};\xi^{j+1}) - \nabla f(x^{j};\xi^{j+1}) \right]}^q} \\
            &\color{navy}\le \delta^q \thesischfiveE{\norm{x^{j+1} - x^{j}}^q} \le \gamma^q \delta^q,
        \end{align*}
    and plugging this new bound into~\eqref{ch5:d98da1f4-086e-4004-a809-2eb471da42a6_chapter_5} and following the steps that lead to~\eqref{ch5:91e3e64b-54a6-46f3-96a5-a9b4f24dc6ee_chapter_5} gives the inequality~\eqref{ch5:b9474e8b-ac2e-49ce-9636-37a177622ba3_chapter_5}, as claimed.
\end{proof}

\subsection{Proof of~Theorem~\ref{ch5:thm:nsgd-mvr-convergence-analysis_chapter_5}}

{\color{navy}Combining Lemmas~\ref{ch5:appdx-lem:unrolling-descent-lemma-nsgd-mvr_chapter_5} and~\ref{ch5:appdx-lem:bounding-error-term-nsgd-mvr_chapter_5} yields the in-expectation convergence analysis of Algorithm~\ref{ch5:algo:nsgd-mvr_chapter_5}.}

\begin{theorem}[Restatement of Theorem~\ref{ch5:thm:nsgd-mvr-convergence-analysis_chapter_5}]
	     Under~Assumption~\ref{ch5:ass:lower-boundedness_chapter_5}, Assumption~\ref{ch5:ass:p-bounded-central-moment-gradient_chapter_5}, and Assumption~\ref{ch5:ass:mean-squared-smoothness_chapter_5}, {\color{navy}choose, for a sufficiently large universal constant $C>0$,}
        \begin{thesischfivefitdisplay}
            {\color{navy}K=\left\lceil C\left(1+\left(\frac{\sigma_1}{\thesischfiveeps}\right)^{p/(p-1)}
            +\frac{\bar L\Delta}{\thesischfiveeps^2}
            +\frac{\bar L\Delta}{\thesischfiveeps^2}
            \left(\frac{\sigma_1}{\thesischfiveeps}\right)^{p/(q(p-1))}\right)\right\rceil .}
        \end{thesischfivefitdisplay}
        {\color{navy}Let the initial gradient estimate $g^0$ be given by}
        \begin{thesischfivefitdisplay} g^0 = \frac{1}{B_{\textnormal{init}}} \sum_{j = 0}^{B_{\textnormal{init}} - 1} \nabla f\left( x^0; \xi_{j}^{0} \right), \end{thesischfivefitdisplay}
    where the initial batch size and the stepsize are
        \begin{align*}
            \color{navy}B_{\textnormal{init}} &\color{navy}= \left\lceil\max\thesischfiveens{1, \left( \frac{2\sigma_1}{\thesischfiveeps} \right)^{p/(p - 1)}}\right\rceil, \\
            \gamma &= \sqrt{\frac{\Delta \eta^{1/q}}{\bar{L} K}},
        \end{align*}
    and the momentum parameter is $\eta = \min\thesischfiveens{1, \eta_{\textnormal{eff}}}$ with
        \begin{thesischfivefitdisplay} \eta_{\textnormal{eff}} = \max\thesischfiveens{\left( \frac{\thesischfiveeps}{\sigma_1 K} \right)^{p/(2p - 1)}, \left( \frac{\bar{L} \Delta}{\sigma_1^2 K} \right)^{pq/(p (2 q + 1) - 2 q)}}. \thesischfivenumberthis\label{ch5:2c0dbb02-323f-4955-a4f2-8f3b73a0f89d-3_chapter_5} \end{thesischfivefitdisplay}
    {\color{navy}When $\sigma_1=0$, we use the deterministic convention $\eta=1$ and $B_{\textnormal{init}}=1$.}
    {\color{navy}If $\bar L=0$, the assumptions imply $\nabla f\equiv0$ and the conclusion is immediate; the displayed parameter choice is therefore used for $\bar L>0$.}
    {\color{navy}For $x_{\mathrm{out}}\sim\mathcal U\{x^0,\ldots,x^{K-1}\}$, Algorithm~\ref{ch5:algo:nsgd-mvr_chapter_5} guarantees $\Exp{\norm{\nabla f(x_{\mathrm{out}})}}\le\thesischfiveeps$ with total sample complexity}
        \begin{thesischfivefitdisplay} \cO\left( {\color{navy}1+} \left( \frac{\sigma_1}{\thesischfiveeps} \right)^{p/(p - 1)} + \frac{\bar{L} \Delta}{\thesischfiveeps^2} + \frac{\bar{L} \Delta}{\thesischfiveeps^2} \left( \frac{\sigma_1}{\thesischfiveeps} \right)^{p/(q (p - 1))} \right). \end{thesischfivefitdisplay}
\end{theorem}

\begin{proof}
    According to~Lemma~\ref{ch5:appdx-lem:unrolling-descent-lemma-nsgd-mvr_chapter_5} we have
        \begin{thesischfivefitdisplay} \frac{1}{K} \sum_{k=0}^{K - 1} \norm{\nabla f(x^{k})} \le \frac{\Delta}{\gamma K} + \frac{2}{K} \sum_{k=0}^{K - 1} \norm{\hat{e}_k} + \frac{\gamma \bar{L}}{2}, \end{thesischfivefitdisplay}
    and using~Lemma~\ref{ch5:appdx-lem:bounding-error-term-nsgd-mvr_chapter_5} this yields
        \begin{align*}
            \frac{1}{K} \sum_{k=0}^{K - 1} \thesischfiveE{\norm{\hat{e}_k}} & \le \frac{1}{K} \sum_{k=0}^{K - 1} (1 - \eta)^k \thesischfiveE{\norm{\hat{e}_0}} + 2 \sigma_1 \eta^{(p - 1)/p} + 4 \gamma \bar{L} \eta^{-1/q} \\
            & \le \frac{\thesischfiveE{\norm{\hat{e}_0}}}{\eta K} + 2 \sigma_1 \eta^{(p - 1)/p} + 4 \gamma \bar{L} \eta^{-1/q},
        \end{align*}
    hence
        \begin{thesischfivefitalignatnum}{2}
            \MoveEqLeft[1] \frac{1}{K} \sum_{k=0}^{K - 1} \thesischfiveE{\norm{\nabla f(x^{k})}} \\
            & \le \frac{\Delta}{\gamma K} + \frac{2 \thesischfiveE{\norm{\hat{e}_0}}}{\eta K} + 4 \sigma_1 \eta^{(p - 1)/p} + 8 \gamma \bar{L} \eta^{-1/q} + \frac{\gamma \bar{L}}{2} \\
            & \le \frac{\Delta}{\gamma K} + \frac{2 \thesischfiveE{\norm{\hat{e}_0}}}{\eta K} + 4 \sigma_1 \eta^{(p - 1)/p} + 9 \gamma \bar{L} \eta^{-1/q}, \label{ch5:d53dc0d7-cefc-4bc1-88bf-08c8201be495_chapter_5}
        \end{thesischfivefitalignatnum}
    since $0 < \eta \le 1$. Now, using $\gamma = \sqrt{\Delta \eta^{1/q}/(\bar{L} K)}$ we have
        \begin{thesischfivefitdisplay} \frac{1}{K} \sum_{k=0}^{K - 1} \thesischfiveE{\norm{\nabla f(x^{k})}} = \cO\left( \eta^{-1/(2 q)} \sqrt{\frac{\bar{L} \Delta}{K}} + \frac{\thesischfiveE{\norm{\hat{e}_0}}}{\eta K} + \sigma_1 \eta^{(p - 1)/p} \right). \thesischfivenumberthis\label{ch5:36a6b109-a2ff-445f-937c-7040717cb89d_chapter_5} \end{thesischfivefitdisplay}
    Now, by our choice of $g^0$ we have
        \begin{thesischfivefitalignatnum}{2}
            \thesischfiveE{\norm{\hat{e}_0}} & \le \left( \thesischfiveE{\norm{\hat{e}_0}^p} \right)^{1/p} \\
            & = \left( \thesischfiveE{\norm{g^0 - \nabla f(x^0)}^p} \right)^{1/p} \\
            & \le \frac{2}{B_{\textnormal{init}}} \left( \sum_{j = 0}^{B_{\textnormal{init}} - 1} \thesischfiveE{\norm{{\color{navy}\nabla f(x^0;\xi_{j}^{0})} - \nabla f(x^0)}^p} \right)^{1/p} \\
            & \le \frac{2}{B_{\textnormal{init}}} \left( \sum_{j = 0}^{B_{\textnormal{init}} - 1} \sigma_1^p \right)^{1/p} \\
            & = \frac{2 \sigma_1}{B_{\textnormal{init}}^{(p - 1)/p}}, \label{ch5:6a9db682-cb1e-43bb-878f-e4b60ebd3183_chapter_5}
        \end{thesischfivefitalignatnum}
    where the first inequality is Jensen's inequality (Lemma~\ref{ch5:lem:jensen-inequality_chapter_5}), the second follows from Lemma~\ref{ch5:appdx-technical-lem:app-von-Bahr-and-Essen_chapter_5} and $2^{1/p} \le 2$, and the third from Assumption~\ref{ch5:ass:p-bounded-central-moment-gradient_chapter_5}. Since {\color{navy}$B_{\textnormal{init}} \ge (2\sigma_1/\thesischfiveeps)^{p/(p - 1)}$}, we have {\color{navy}$\thesischfiveE{\norm{\hat{e}_0}}\le\thesischfiveeps$}. Moreover, $\eta = \min\{1, \eta_{\textnormal{eff}}\}$ gives $\eta^{-1/(2q)} \le 1 + \eta_{\textnormal{eff}}^{-1/(2q)}$, $\eta^{-1} \le 1 + \eta_{\textnormal{eff}}^{-1}$ and $\eta \le \eta_{\textnormal{eff}}$. This gives
        \begin{thesischfivefitalignatnum}{2}
            \MoveEqLeft[1] \frac{1}{K} \sum_{k=0}^{K - 1} \thesischfiveE{\norm{\nabla f(x^{k})}} \\
            \thesischfiveoversetlab{\eqref{ch5:36a6b109-a2ff-445f-937c-7040717cb89d_chapter_5}}&{=} \cO\left( \eta^{-1/(2 q)} \sqrt{\frac{\bar{L} \Delta}{K}} + \frac{\thesischfiveE{\norm{\hat{e}_0}}}{\eta K} + \sigma_1 \eta^{(p - 1)/p} \right) \\
            & = \cO\left( \eta^{-1/(2 q)} \sqrt{\frac{\bar{L} \Delta}{K}} + \frac{\thesischfiveeps}{\eta K} + \sigma_1 \eta^{(p - 1)/p} \right) \\
            & = \cO\left( \sqrt{\frac{\bar{L} \Delta}{K}} + \eta_{\textnormal{eff}}^{-1/(2 q)} \sqrt{\frac{\bar{L} \Delta}{K}} + \frac{\thesischfiveeps}{K} + \frac{\thesischfiveeps}{K} \eta_{\textnormal{eff}}^{-1} + \sigma_1 \eta_{\textnormal{eff}}^{(p - 1)/p} \right) \\
            & = \cO\biggl( \frac{\thesischfiveeps}{K} + \sqrt{\frac{\bar{L} \Delta}{K}} \\
            &\qquad + \sigma_1 \left( \frac{\bar{L} \Delta}{\sigma_1^2 K} \right)^{q (p - 1)/(p (2 q + 1) - 2 q)} + \sigma_1 \left( \frac{\thesischfiveeps}{\sigma_1 K}  \right)^{(p - 1)/(2 p - 1)} \biggr), \label{ch5:0f0a2e3d-d04e-4e64-948f-d38d2619b1fb_chapter_5}
        \end{thesischfivefitalignatnum}
    where in the last equality we use the choice of $\eta_{\textnormal{eff}}$ from~\eqref{ch5:2c0dbb02-323f-4955-a4f2-8f3b73a0f89d-3_chapter_5}, which gives
        \begin{align*}
            \frac{\thesischfiveeps}{K} \eta_{\textnormal{eff}}^{-1} & \le \frac{\thesischfiveeps}{K} \left( \frac{\sigma_1 K}{\thesischfiveeps} \right)^{p/(2p - 1)} = \sigma_1 \left( \frac{\thesischfiveeps}{\sigma_1 K} \right)^{(p - 1)/(2 p - 1)}, \\
            \eta_{\textnormal{eff}}^{-1/(2 q)} \sqrt{\frac{\bar{L} \Delta}{K}} & \le \sigma_1 \left( \frac{\bar{L} \Delta}{\sigma_1^2 K} \right)^{q (p - 1)/(p (2 q + 1) - 2 q)},
        \end{align*}
    and
        \begin{thesischfivefitdisplay} \sigma_1 \eta_{\textnormal{eff}}^{(p - 1)/p} \le \sigma_1 \left( \frac{\thesischfiveeps}{\sigma_1 K} \right)^{(p - 1)/(2p - 1)} + \sigma_1 \left( \frac{\bar{L} \Delta}{\sigma_1^2 K} \right)^{q (p - 1)/(p (2 q + 1) - 2 q)}. \end{thesischfivefitdisplay}

    Finally, since $\thesischfiveeps/K \le \thesischfiveeps$, from the bound~\eqref{ch5:0f0a2e3d-d04e-4e64-948f-d38d2619b1fb_chapter_5} we deduce that the sample complexity of~\ref{ch5:algo:nsgd-mvr_chapter_5} is exactly
        \begin{thesischfivefitdisplay} \cO\left( {\color{navy}1+} \left( \frac{\sigma_1}{\thesischfiveeps} \right)^{p/(p - 1)} + \frac{\bar{L} \Delta}{\thesischfiveeps^2} + \frac{\bar{L} \Delta}{\thesischfiveeps^2} \left( \frac{\sigma_1}{\thesischfiveeps} \right)^{p/(q (p - 1))} \right). \end{thesischfivefitdisplay}
    {\color{navy}This proves the claim and matches the lower bound in Theorem~\ref{ch5:thm:lower-bound-sofo-mean-squared-smoothness_chapter_5}.}
\end{proof}

\subsection{Proof of~Theorem~\ref{ch5:thm:nsgd-mvr-convergence-analysis-2_chapter_5}}

\begin{theorem}[Restatement of Theorem~\ref{ch5:thm:nsgd-mvr-convergence-analysis-2_chapter_5}]
    Under~Assumption~\ref{ch5:ass:lower-boundedness_chapter_5}, Assumption~\ref{ch5:ass:L-lipschitz-gradients_chapter_5}, Assumption~\ref{ch5:ass:p-bounded-central-moment-gradient_chapter_5}, and Assumption~\ref{ch5:ass:mean-squared-smoothness-2_chapter_5}, {\color{navy}choose, for a sufficiently large universal constant $C>0$,}
        \begin{thesischfivefitdisplay}
            {\color{navy}K=\left\lceil C\left(1+\left(\frac{\sigma_1}{\thesischfiveeps}\right)^{p/(p-1)}
            +\frac{(L_1+\delta)\Delta}{\thesischfiveeps^2}
            +\frac{\delta\Delta}{\thesischfiveeps^2}
            \left(\frac{\sigma_1}{\thesischfiveeps}\right)^{p/(q(p-1))}\right)\right\rceil .}
        \end{thesischfivefitdisplay}
        {\color{navy}Let the initial gradient estimate $g^0$ be given by}
        \begin{thesischfivefitdisplay} g^0 = \frac{1}{B_{\textnormal{init}}} \sum_{j = 0}^{B_{\textnormal{init}} - 1} \nabla f\left( x^0; \xi_{j}^{0} \right), \end{thesischfivefitdisplay}
    where the initial batch size and the stepsize are
        \begin{align*}
            \color{navy}B_{\textnormal{init}} &\color{navy}= \left\lceil\max\thesischfiveens{1, \left( \frac{2\sigma_1}{\thesischfiveeps} \right)^{p/(p - 1)}}\right\rceil, \\
            \gamma &= \min\thesischfiveens{\sqrt{\frac{\Delta}{L_1 K}}, \sqrt{\frac{\Delta \eta^{1/q}}{\delta K}}},
        \end{align*}
    and the momentum parameter is $\eta = \min\thesischfiveens{1, \eta_{\textnormal{eff}}}$ with
        \begin{thesischfivefitdisplay} \eta_{\textnormal{eff}} = \max\thesischfiveens{\left( \frac{\thesischfiveeps}{\sigma_1 K} \right)^{p/(2p - 1)}, \left( \frac{\delta \Delta}{\sigma_1^2 K} \right)^{pq/(p (2 q + 1) - 2 q)}}. \thesischfivenumberthis\label{ch5:2c0dbb02-323f-4955-a4f2-8f3b73a0f89d-2-bis_chapter_5} \end{thesischfivefitdisplay}
    {\color{navy}When $\sigma_1=0$, we use the deterministic convention $\eta=1$ and $B_{\textnormal{init}}=1$.}
    {\color{navy}A stepsize term with zero denominator is interpreted as $+\infty$ and omitted from the minimum. If $L_1=\delta=0$, then $\nabla f\equiv0$ and the conclusion is immediate.}
    {\color{navy}For $x_{\mathrm{out}}\sim\mathcal U\{x^0,\ldots,x^{K-1}\}$, Algorithm~\ref{ch5:algo:nsgd-mvr_chapter_5} guarantees $\Exp{\norm{\nabla f(x_{\mathrm{out}})}}\le\thesischfiveeps$ with total sample complexity}
        \begin{thesischfivefitdisplay} \cO\left( {\color{navy}1+} \left( \frac{\sigma_1}{\thesischfiveeps} \right)^{p/(p - 1)} + \frac{(L_1 + \delta) \Delta}{\thesischfiveeps^2} + \frac{\delta \Delta}{\thesischfiveeps^2} \left( \frac{\sigma_1}{\thesischfiveeps} \right)^{p/(q (p - 1))} \right). \end{thesischfivefitdisplay}
\end{theorem}

\begin{proof}
    As the function $f$ has $L_1$--Lipschitz gradients by~Assumption~\ref{ch5:ass:L-lipschitz-gradients_chapter_5} then applying~Lemma~\ref{ch5:appdx-lem:unrolling-descent-lemma-nsgd-mvr-hess_chapter_5} gives
        \begin{thesischfivefitdisplay} \frac{1}{K} \sum_{k=0}^{K - 1} \norm{\nabla f(x^{k})} \le \frac{\Delta}{\gamma K} + \frac{2}{K} \sum_{k=0}^{K - 1} \norm{\hat{e}_k} + \frac{\gamma L_1}{2}, \end{thesischfivefitdisplay}
    and using~Lemma~\ref{ch5:appdx-lem:bounding-error-term-nsgd-mvr-refined_chapter_5} yields
        \begin{align*}
            \frac{1}{K} \sum_{k=0}^{K - 1} \thesischfiveE{\norm{\hat{e}_k}} & \le \frac{1}{K} \sum_{k=0}^{K - 1} (1 - \eta)^k \thesischfiveE{\norm{\hat{e}_0}} + 2 \sigma_1 \eta^{(p - 1)/p} + 2 \gamma \delta \eta^{-1/q} \\
            & \le \frac{\thesischfiveE{\norm{\hat{e}_0}}}{\eta K} + 2 \sigma_1 \eta^{(p - 1)/p} + 2 \gamma \delta \eta^{-1/q},
        \end{align*}
    hence
        \begin{thesischfivefitdisplay} \frac{1}{K} \sum_{k=0}^{K - 1} \thesischfiveE{\norm{\nabla f(x^{k})}} \le \frac{\Delta}{\gamma K} + \frac{2 \thesischfiveE{\norm{\hat{e}_0}}}{\eta K} + 4 \sigma_1 \eta^{(p - 1)/p} + 4 \gamma \delta \eta^{-1/q} + \frac{\gamma L_1}{2}. \thesischfivenumberthis\label{ch5:c6f0fec1-2964-442a-b2c3-7c25dfd2e653_chapter_5} \end{thesischfivefitdisplay}
    Now, using our choice of stepsize $\gamma$ we have
        \begin{thesischfivefitdisplay} \frac{1}{K} \sum_{k=0}^{K - 1} \thesischfiveE{\norm{\nabla f(x^{k})}} = \cO\left( \sqrt{\frac{L_1 \Delta}{K}} + \eta^{-1/(2 q)} \sqrt{\frac{\delta \Delta}{K}} + \frac{\thesischfiveE{\norm{\hat{e}_0}}}{\eta K} + \sigma_1 \eta^{(p - 1)/p} \right). \thesischfivenumberthis\label{ch5:36a6b109-a2ff-445f-937c-7040717cb89d-bis_chapter_5} \end{thesischfivefitdisplay}
    {\color{navy}Next, by our choice of $g^0$, as in the preceding proof, we have}
        \begin{thesischfivefitdisplay} \thesischfiveE{\norm{\hat{e}_0}} \thesischfiveoversetlab{\eqref{ch5:6a9db682-cb1e-43bb-878f-e4b60ebd3183_chapter_5}}{\le} \frac{2 \sigma_1}{B_{\textnormal{init}}^{(p - 1)/p}}, \end{thesischfivefitdisplay}
    and since {\color{navy}$B_{\textnormal{init}} \ge (2\sigma_1/\thesischfiveeps)^{p/(p - 1)}$}, we have {\color{navy}$\thesischfiveE{\norm{\hat{e}_0}}\le\thesischfiveeps$}. Using also $\eta = \min\{1, \eta_{\textnormal{eff}}\}$ as in the preceding proof, this gives
        \begin{thesischfivefitalignatnum}{2}
            \MoveEqLeft[1] \frac{1}{K} \sum_{k=0}^{K - 1} \thesischfiveE{\norm{\nabla f(x^{k})}} \\
            \thesischfiveoversetlab{\eqref{ch5:36a6b109-a2ff-445f-937c-7040717cb89d-bis_chapter_5}}&{=} \cO\left( \sqrt{\frac{L_1 \Delta}{K}} + \eta^{-1/(2 q)} \sqrt{\frac{\delta \Delta}{K}} + \frac{\thesischfiveE{\norm{\hat{e}_0}}}{\eta K} + \sigma_1 \eta^{(p - 1)/p} \right) \\
            & = \cO\left( \sqrt{\frac{L_1 \Delta}{K}} + \eta^{-1/(2 q)} \sqrt{\frac{\delta \Delta}{K}} + \frac{\thesischfiveeps}{\eta K} + \sigma_1 \eta^{(p - 1)/p} \right) \\
            & = \cO\left( \sqrt{\frac{L_1 \Delta}{K}} + \sqrt{\frac{\delta \Delta}{K}} + \eta_{\textnormal{eff}}^{-1/(2 q)} \sqrt{\frac{\delta \Delta}{K}} + \frac{\thesischfiveeps}{K} + \frac{\thesischfiveeps}{K} \eta_{\textnormal{eff}}^{-1} + \sigma_1 \eta_{\textnormal{eff}}^{(p - 1)/p} \right) \\
            & = \cO\biggl( \frac{\thesischfiveeps}{K} + \sqrt{\frac{L_1 \Delta}{K}} + \sqrt{\frac{\delta \Delta}{K}} \\
            &\qquad + \sigma_1 \left( \frac{\delta \Delta}{\sigma_1^2 K} \right)^{q (p - 1)/(p (2 q + 1) - 2 q)} + \sigma_1 \left( \frac{\thesischfiveeps}{\sigma_1 K}  \right)^{(p - 1)/(2 p - 1)} \biggr), \label{ch5:0f0a2e3d-d04e-4e64-948f-d38d2619b1fbc_chapter_5}
        \end{thesischfivefitalignatnum}
    where in the last equality we use the choice of $\eta_{\textnormal{eff}}$ from~\eqref{ch5:2c0dbb02-323f-4955-a4f2-8f3b73a0f89d-2-bis_chapter_5}, which gives
        \begin{align*}
            \frac{\thesischfiveeps}{K} \eta_{\textnormal{eff}}^{-1} & \le \frac{\thesischfiveeps}{K} \left( \frac{\sigma_1 K}{\thesischfiveeps} \right)^{p/(2p - 1)} = \sigma_1 \left( \frac{\thesischfiveeps}{\sigma_1 K} \right)^{(p - 1)/(2 p - 1)}, \\
            \eta_{\textnormal{eff}}^{-1/(2 q)} \sqrt{\frac{\delta \Delta}{K}} & \le \sigma_1 \left( \frac{\delta \Delta}{\sigma_1^2 K} \right)^{q (p - 1)/(p (2 q + 1) - 2 q)},
        \end{align*}
    and
        \begin{thesischfivefitdisplay} \sigma_1 \eta_{\textnormal{eff}}^{(p - 1)/p} \le \sigma_1 \left( \frac{\thesischfiveeps}{\sigma_1 K} \right)^{(p - 1)/(2p - 1)} + \sigma_1 \left( \frac{\delta \Delta}{\sigma_1^2 K} \right)^{q (p - 1)/(p (2 q + 1) - 2 q)}. \end{thesischfivefitdisplay}

    Finally, since $\thesischfiveeps/K \le \thesischfiveeps$, from the bound~\eqref{ch5:0f0a2e3d-d04e-4e64-948f-d38d2619b1fbc_chapter_5} we deduce that the sample complexity of~\ref{ch5:algo:nsgd-mvr_chapter_5} is exactly
        \begin{thesischfivefitdisplay} \cO\left( {\color{navy}1+} \left( \frac{\sigma_1}{\thesischfiveeps} \right)^{p/(p - 1)} + \frac{(L_1 + \delta) \Delta}{\thesischfiveeps^2} + \frac{\delta \Delta}{\thesischfiveeps^2} \left( \frac{\sigma_1}{\thesischfiveeps} \right)^{p/(q (p - 1))} \right), \end{thesischfivefitdisplay}
    as claimed, and, combining~\ref{ch5:algo:nsgd-mvr_chapter_5} with \algname{NSGD-Mom} is enough to match our lower bound from~Theorem~\ref{ch5:thm:lower-bound-sofo-mean-squared-smoothness-2_chapter_5}.
\end{proof}

\subsection{Proof of~Theorem~\ref{ch5:thm:nsgd-mvr-convergence-analysis-unknown-p-q_chapter_5}}

\begin{theorem}[Restatement of Theorem~\ref{ch5:thm:nsgd-mvr-convergence-analysis-unknown-p-q_chapter_5}]
    Under~Assumption~\ref{ch5:ass:lower-boundedness_chapter_5}, Assumption~\ref{ch5:ass:p-bounded-central-moment-gradient_chapter_5}, and Assumption~\ref{ch5:ass:mean-squared-smoothness_chapter_5}, let the initial gradient estimate $g^0 = \nabla f(x^0;\xi^{0})$, let the stepsize $\gamma = \sqrt{\Delta \eta/(\bar{L} K)}$, and let the momentum parameter $\eta = K^{-1/2} \in \thesischfiveintof{0}{1}$.
    {\color{navy}If $\bar L=0$, the assumptions imply $\nabla f\equiv0$ and the conclusion is immediate; hence the displayed stepsize is used for $\bar L>0$. Then Algorithm~\ref{ch5:algo:nsgd-mvr_chapter_5} guarantees}
        \begin{thesischfivefitdisplay} \frac{1}{K} \sum_{k=0}^{K - 1} \thesischfiveE{\norm{\nabla f(x^{k})}} = \cO\left( \frac{\sigma_1}{K^{(p - 1)/(2 p)}} + \frac{\sqrt{\bar{L} \Delta}}{K^{1/4}} \right). \end{thesischfivefitdisplay}
\end{theorem}

\begin{proof}
    Notice that the bound~\eqref{ch5:d53dc0d7-cefc-4bc1-88bf-08c8201be495_chapter_5} holds for any choice of the parameters $\gamma > 0$, $\eta \in \thesischfiveintof{0}{1}$ and $g^0 \in \R^d$. Hence, for our particular choice $g^0 = \nabla f(x^0;\xi^{0})$, $\gamma = \sqrt{\Delta \eta/(\bar{L} K)}$ and $\eta = K^{-1/2}$ we obtain
        \begin{thesischfivefitalignat}{2}
            \MoveEqLeft[1] \frac{1}{K} \sum_{k=0}^{K - 1} \thesischfiveE{\norm{\nabla f(x^{k})}} \\
            \thesischfiveoversetlab{\eqref{ch5:d53dc0d7-cefc-4bc1-88bf-08c8201be495_chapter_5}}&{=} \cO\left( \frac{\Delta}{\gamma K} + \frac{\thesischfiveE{\norm{\hat{e}_0}}}{\eta K} + \sigma_1 \eta^{(p - 1)/p} + \gamma \bar{L} \eta^{-1/q} \right) \\
            & = \cO\left( \eta^{-1/2} \sqrt{\frac{\bar{L} \Delta}{K}} + \frac{\thesischfiveE{\norm{\hat{e}_0}}}{K^{1/2}} + \frac{\sigma_1}{K^{(p - 1)/(2 p)}} + \eta^{-(2 - q)/(2 q)} \sqrt{\frac{\bar{L} \Delta}{K}} \right) \\
            & = \cO\left( \frac{\sqrt{\bar{L} \Delta}}{K^{1/4}} + \frac{\thesischfiveE{\norm{\hat{e}_0}}}{K^{1/2}} + \frac{\sigma_1}{K^{(p - 1)/(2 p)}} \right) \\
            & = \cO\left( \frac{\sigma_1}{K^{(p - 1)/(2 p)}} + \frac{\sigma_1}{K^{1/2}} + \frac{\sqrt{\bar{L} \Delta}}{K^{1/4}} \right) \\
            & = \cO\left( \frac{\sigma_1}{K^{(p - 1)/(2 p)}} + \frac{\sqrt{\bar{L} \Delta}}{K^{1/4}} \right),
        \end{thesischfivefitalignat}
    where in the third equality we use $\eta^{-1/2} \ge \eta^{-(2 - q)/(2 q)}$, since $\eta \in \thesischfiveintof{0}{1}$ and $(2 - q)/(2 q) \le 1/2$, and in the fourth we use the definition of $\hat{e}_0$ and $g^0$, i.e.,
        \begin{align*}
            \thesischfiveE{\norm{\hat{e}_0}} & = \thesischfiveE{\norm{\nabla f(x^0;\xi^{0}) - \nabla f(x^0)}} \\
            & \le \left( \thesischfiveE{\norm{\nabla f(x^0;\xi^{0}) - \nabla f(x^0)}^p} \right)^{1/p} \le \sigma_1,
        \end{align*}
    by Jensen's inequality (Lemma~\ref{ch5:lem:jensen-inequality_chapter_5}) and Assumption~\ref{ch5:ass:p-bounded-central-moment-gradient_chapter_5}. The last equality uses $(p - 1)/(2p) \le 1/2$.
    {\color{navy}This completes the proof of the theorem.}
\end{proof}

\subsection{Proof of~Theorem~\ref{ch5:thm:nsgd-mvr-convergence-analysis-unknown-p-q-2_chapter_5}}

\begin{theorem}[Restatement of Theorem~\ref{ch5:thm:nsgd-mvr-convergence-analysis-unknown-p-q-2_chapter_5}]
    Under~Assumption~\ref{ch5:ass:lower-boundedness_chapter_5}, Assumption~\ref{ch5:ass:L-lipschitz-gradients_chapter_5}, Assumption~\ref{ch5:ass:p-bounded-central-moment-gradient_chapter_5}, and Assumption~\ref{ch5:ass:mean-squared-smoothness-2_chapter_5}, let the initial gradient estimate $g^0 = \nabla f(x^0;\xi^{0})$, let the stepsize $\gamma = \min\{\sqrt{\Delta/(L_1 K)}, \sqrt{\Delta \eta/(\delta K)}\}$, and let the momentum parameter $\eta = K^{-1/2} \in \thesischfiveintof{0}{1}$.
    {\color{navy}A stepsize term with zero denominator is interpreted as $+\infty$ and omitted from the minimum. If $L_1=\delta=0$, then $\nabla f\equiv0$ and the conclusion is immediate. Otherwise, Algorithm~\ref{ch5:algo:nsgd-mvr_chapter_5} guarantees}
        \begin{thesischfivefitdisplay} \frac{1}{K} \sum_{k=0}^{K - 1} \thesischfiveE{\norm{\nabla f(x^{k})}} = \cO\left( \frac{\sigma_1}{K^{(p - 1)/(2 p)}} + \frac{\sqrt{\delta \Delta}}{K^{1/4}} + \frac{\sqrt{L_1 \Delta}}{K^{1/2}} \right). \thesischfivenumberthis\label{ch5:004535e0-7ad4-4957-8cc1-223ee18990af_chapter_5} \end{thesischfivefitdisplay}
\end{theorem}

\begin{proof}
    As in the previous theorem, note that the bound~\eqref{ch5:c6f0fec1-2964-442a-b2c3-7c25dfd2e653_chapter_5} holds for any choice of the parameters $\gamma > 0$, $\eta \in \thesischfiveintof{0}{1}$ and $g^0 \in \R^d$. Hence, for our particular choice $g^0 = \nabla f(x^0;\xi^{0})$, $\gamma = \min\{\sqrt{\Delta/(L_1 K)}, \sqrt{\Delta \eta/(\delta K)}\}$ and $\eta = K^{-1/2}$ we obtain, as in the previous proof,
        \begin{thesischfivefitalignat}{2}
            \MoveEqLeft[1] \frac{1}{K} \sum_{k=0}^{K - 1} \thesischfiveE{\norm{\nabla f(x^{k})}} \\
            \thesischfiveoversetlab{\eqref{ch5:c6f0fec1-2964-442a-b2c3-7c25dfd2e653_chapter_5}}&{=} \cO\left( \frac{\Delta}{\gamma K} + \frac{\thesischfiveE{\norm{\hat{e}_0}}}{\eta K} + \sigma_1 \eta^{(p - 1)/p} + \gamma \delta \eta^{-1/q} + \gamma L_1 \right) \\
            & = \cO\biggl( \sqrt{\frac{L_1 \Delta}{K}} + \eta^{-1/2} \sqrt{\frac{\delta \Delta}{K}} + \frac{\thesischfiveE{\norm{\hat{e}_0}}}{K^{1/2}} \\
            &\qquad + \frac{\sigma_1}{K^{(p - 1)/(2 p)}} + \eta^{-(2 - q)/(2 q)} \sqrt{\frac{\delta \Delta}{K}} \biggr) \\
            & = \cO\left( \sqrt{\frac{L_1 \Delta}{K}} + \frac{\sqrt{\delta \Delta}}{K^{1/4}} + \frac{\thesischfiveE{\norm{\hat{e}_0}}}{K^{1/2}} + \frac{\sigma_1}{K^{(p - 1)/(2 p)}} \right) \\
            & = \cO\left( \frac{\sigma_1}{K^{(p - 1)/(2 p)}} + \frac{\sigma_1}{K^{1/2}} + \frac{\sqrt{\delta \Delta}}{K^{1/4}} + \sqrt{\frac{L_1 \Delta}{K}}\right) \\
            & = \cO\left( \frac{\sigma_1}{K^{(p - 1)/(2 p)}} + \frac{\sqrt{\delta \Delta}}{K^{1/4}} + \sqrt{\frac{L_1 \Delta}{K}} \right),
        \end{thesischfivefitalignat}
    where in the third equality we use $\eta^{-1/2} \ge \eta^{-(2 - q)/(2 q)}$, since $\eta \in \thesischfiveintof{0}{1}$ and $(2 - q)/(2 q) \le 1/2$. {\color{navy}In the fourth, Jensen's inequality and Assumption~\ref{ch5:ass:p-bounded-central-moment-gradient_chapter_5} give
        \begin{thesischfivefitdisplay} \thesischfiveE{\norm{\hat{e}_0}} = \thesischfiveE{\norm{\nabla f(x^0;\xi^0) - \nabla f(x^0)}} \le \sigma_1, \end{thesischfivefitdisplay}
    exactly as in the preceding proof.} The last equality uses $(p - 1)/(2p) \le 1/2$.
    
    {\color{navy}This completes the proof of the theorem.}
\end{proof}

\newpage%
\section{Additional Proofs for Section~\ref{ch5:sec:near-optimal-method-nsgd-mvr-hess_chapter_5}}\label{ch5:appdx:proof-nsgd-mvr-hess_chapter_5}

\subsection{Auxiliary Lemmas}

Some of the auxiliary lemmas needed in this section (Lemma~\ref{ch5:appdx-lem:descent-lemma-nsgd-mvr-hess_chapter_5} and Lemma~\ref{ch5:appdx-lem:unrolling-descent-lemma-nsgd-mvr-hess_chapter_5}) have already been established in the previous section.
\begin{lemma}[Bounding the Error Term]\label{ch5:appdx-lem:bounding-error-term-nsgd-mvr-hess_chapter_5}
    Under~Assumption~\ref{ch5:ass:p-bounded-central-moment-gradient_chapter_5}, Assumption~\ref{ch5:ass:q-bounded-central-moment-hessian_chapter_5}, Assumption~\ref{ch5:ass:L-lipschitz-gradients_chapter_5}, and Assumption~\ref{ch5:ass:L-lipschitz-hessians_chapter_5}, for all $k \in \thesischfiveInt{0}{K - 1}$ we have
        \begin{thesischfivefitdisplay} \thesischfiveE{\norm{\hat{e}_k}} \le (1 - \eta)^k \thesischfiveE{\norm{\hat{e}_0}} + 2 \sigma_1 \eta^{(p - 1)/p} + 2 \gamma \sigma_2 \eta^{-1/q} + 4 \gamma \min\thesischfiveens{L_1, \gamma L_2} \eta^{-1/2}, \end{thesischfivefitdisplay}
    where $\hat{e}_k \eqdef g^{k} - \nabla f(x^{k})$.
\end{lemma}

\begin{proof}
    {\color{navy}Let $d^k\eqdef x^k-x^{k-1}$ and
    $\hat x^k=u^kx^k+(1-u^k)x^{k-1}$, where
    $u^k\sim\mathcal U(0,1)$. The estimator update gives, for the error term $\hat e_k = g^k - \nabla f(x^k)$,
        \begin{align}
            \hat e_k &= (1-\eta)\hat e_{k-1} + \eta e_k + (1-\eta)\hat R_k + (1-\eta)\hat S_k,
            \label{ch5:2f56c833-5913-40e8-a00c-20c7b1a4b32b-2_chapter_5}
            \shortintertext{where}
            e_k &\eqdef \nabla f(x^k;\xi^k) - \nabla f(x^k), \notag\\
            \hat R_k &\eqdef \nabla^2 f(\hat x^k) d^k - [\nabla f(x^k) - \nabla f(x^{k-1})],
            \label{ch5:14872232-13b0-4801-9fa1-bccc0be89a0b_chapter_5}\\
            \hat S_k &\eqdef [\nabla^2 f(\hat x^k;\hat\xi^k) - \nabla^2 f(\hat x^k)] d^k.
            \label{ch5:891801fc-22ef-4562-81bb-3b5040b0f520_chapter_5}
        \end{align}
    Define
    \begin{thesischfivefitdisplay}
        \mathcal F_k\eqdef
        \sigma\!\left(g^0,
        (\xi^1,\hat\xi^1,u^1),\ldots,
        (\xi^k,\hat\xi^k,u^k)\right).
    \end{thesischfivefitdisplay}
    The fresh variables $\xi^k,\hat\xi^k,u^k$ are mutually independent
    conditionally on $\mathcal F_{k-1}$. Hence $e_k,\hat R_k,\hat S_k$
    are $\mathcal F_k$-measurable martingale differences. For $\hat R_k$,
    the fundamental theorem of calculus gives
        \begin{align*}
            \Exp{\nabla^2 f(\hat x^k)d^k\mid\mathcal F_{k-1}}
            &=\int_0^1\nabla^2 f(x^{k-1}+s d^k)d^k\,\mathrm ds \\
            &=\nabla f(x^k)-\nabla f(x^{k-1}).
            \thesischfivenumberthis\label{ch5:ff082357-0807-4b61-b6c0-12bd29b8b854_chapter_5}
        \end{align*}}

        Then, unrolling the recursion~\eqref{ch5:2f56c833-5913-40e8-a00c-20c7b1a4b32b-2_chapter_5} gives
            \begin{align*}
                \hat{e}_k & = (1 - \eta)^k \hat{e}_0 + \eta \sum_{j = 0}^{k - 1} (1 - \eta)^{k - j - 1} e_{j + 1} \\
                &\quad + \sum_{j = 0}^{k - 1} (1 - \eta)^{k  - j} \hat{R}_{j + 1} + \sum_{j = 0}^{k - 1} (1 - \eta)^{k  - j} \hat{S}_{j + 1},
            \end{align*}
        and taking the norm followed by the total expectation yields
            \begin{align*}
                \thesischfiveE{\norm{\hat{e}_k}} & \le (1 - \eta)^k \thesischfiveE{\norm{\hat{e}_0}} + \eta \thesischfiveE{\norm{\sum_{j = 0}^{k - 1} (1 - \eta)^{k - j - 1} e_{j + 1}}} \\
                &\quad + (1 - \eta) \thesischfiveE{\norm{\sum_{j = 0}^{k - 1} (1 - \eta)^{k - j - 1} \hat{R}_{j + 1}}} \\
                &\quad + (1 - \eta) \thesischfiveE{\norm{\sum_{j = 0}^{k - 1} (1 - \eta)^{k - j - 1} \hat{S}_{j + 1}}}.
            \end{align*}
        We now need to upper bound the last three terms of the previous inequality. For the first term, exactly as in~\eqref{ch5:cea5881d-8688-41f3-8083-df58a5981198_chapter_5}, we have
            \begin{thesischfivefitalignatnum}{2}
                \thesischfiveE{\norm{\sum_{j = 0}^{k - 1} (1 - \eta)^{k - j - 1} e_{j + 1}}} & \le 2 \sigma_1 \eta^{-1/p}. \label{ch5:cea5881d-8688-41f3-8083-df58a5981198-2_chapter_5}
            \end{thesischfivefitalignatnum}
        Then, for the second term, Jensen's inequality (Lemma~\ref{ch5:lem:jensen-inequality_chapter_5}) and Lemma~\ref{ch5:appdx-technical-lem:app-von-Bahr-and-Essen_chapter_5} give
            \begin{thesischfivefitalignatnum}{2}
                \thesischfiveE{\norm{\sum_{j = 0}^{k - 1} (1 - \eta)^{k - j - 1} \hat{R}_{j + 1}}} & \le \left( \thesischfiveE{\norm{\sum_{j = 0}^{k - 1} (1 - \eta)^{k - j - 1} \hat{R}_{j + 1}}^2} \right)^{1/2} \\
                & \le \left( 2 \sum_{j = 0}^{k - 1} (1 - \eta)^{2 (k - j - 1)} \thesischfiveE{\norm{\hat{R}_{j + 1}}^2} \right)^{1/2}, \label{ch5:d98da1f4-086e-4004-a809-2eb471da42a6-3_chapter_5}
            \end{thesischfivefitalignatnum}
        and the second moment of $\hat{R}_k$ can be bounded in two different ways. First, by using~Assumption~\ref{ch5:ass:L-lipschitz-gradients_chapter_5} we have
            \begin{thesischfivefitalignatnum}{2}
                \thesischfiveE{\sqnorm{\hat{R}_k}} \thesischfiveoversetlab{\eqref{ch5:14872232-13b0-4801-9fa1-bccc0be89a0b_chapter_5}}&{=} \thesischfiveE{\sqnorm{\nabla^2 f(\hat{x}^k) d^k - [\nabla f(x^{k}) - \nabla f(x^{k-1})]}} \\
                & \le 2 \thesischfiveE{\sqnorm{\nabla f(x^{k}) - \nabla f(x^{k-1})}} + 2 \thesischfiveE{\sqnorm{\nabla^2 f(\hat{x}^k) d^k}} \\
                & \le 2 L_1^2 \thesischfiveE{\sqnorm{d^k}} + 2 \thesischfiveE{\normop{\nabla^2 f(\hat{x}^k)}^2 \sqnorm{d^k}} \\
                & \le 4 L_1^2 \thesischfiveE{\sqnorm{d^k}} \\
                & \le 4 \gamma^2 L_1^2, \label{ch5:46d7ba91-ef00-4ec2-9b1e-d222e9968289a_chapter_5}
            \end{thesischfivefitalignatnum}
        where the first inequality follows from Lemma~\ref{ch5:lem:jensen-form-1_chapter_5}, the second from Assumption~\ref{ch5:ass:L-lipschitz-gradients_chapter_5}, the third from Lemma~\ref{ch5:appdx-lem:lipschitz-gradients-implies-bounded-hessian_chapter_5}, and the last from $\norm{d^k} \le \gamma$. Second, using the connection between $\nabla^2 f(\hat{x}^k) d^k$ and $\nabla f(x^{k}) - \nabla f(x^{k-1})$ as displayed in~\eqref{ch5:ff082357-0807-4b61-b6c0-12bd29b8b854_chapter_5} we also have
            \begin{thesischfivefitalignatnum}{2}
                \thesischfiveE{\sqnorm{\hat{R}_k}} \thesischfiveoversetlab{\eqref{ch5:14872232-13b0-4801-9fa1-bccc0be89a0b_chapter_5}}&{=} \thesischfiveE{\sqnorm{\nabla^2 f(\hat{x}^k) d^k - [\nabla f(x^{k}) - \nabla f(x^{k-1})]}} \\
                & = \thesischfiveE{\sqnorm{\int_0^1 \left( \nabla^2 f(x^{k-1} + s d^k) - \nabla^2 f(\hat{x}^k) \right) d^k \odif{s}}} \\
                & \le \thesischfiveE{\left( \int_0^1 \norm{\left( \nabla^2 f(x^{k-1} + s d^k) - \nabla^2 f(\hat{x}^k) \right) d^k} \odif{s} \right)^2} \\
                & \le \thesischfiveE{\left( \int_0^1 \normop{\nabla^2 f(x^{k-1} + s d^k) - \nabla^2 f(\hat{x}^k)} \cdot \norm{d^k} \odif{s} \right)^2} \\
                & \le \gamma^2 \, \thesischfiveE{\left( \int_0^1 \normop{\nabla^2 f(x^{k-1} + s d^k) - \nabla^2 f(\hat{x}^k)} \odif{s} \right)^2} \\
                & \le \gamma^2 L_2^2 \, \thesischfiveE{\left( \int_0^1 \norm{s (x^{k} - \hat{x}^k) + (1 - s) (x^{k-1} - \hat{x}^k)} \odif{s} \right)^2} \\
                & \le \gamma^2 L_2^2 \, \thesischfiveE{\left( \int_0^1 s \norm{x^{k} - \hat{x}^k} \odif{s} + \int_0^1 (1 - s) \norm{x^{k-1} - \hat{x}^k} \odif{s} \right)^2} \\
                & = \frac{\gamma^2 L_2^2}{4} \, \thesischfiveE{\left( \norm{x^{k} - \hat{x}^k} + \norm{x^{k-1} - \hat{x}^k} \right)^2} \\
                & = \frac{\gamma^2 L_2^2}{4} \, \thesischfiveE{\left( (1 - {\color{navy}u^k}) \norm{d^k} + {\color{navy}u^k} \norm{d^k} \right)^2} \\
                & \le \frac{\gamma^4 L_2^2}{4}, \label{ch5:46d7ba91-ef00-4ec2-9b1e-d222e9968289b_chapter_5}
            \end{thesischfivefitalignatnum}
        where in the first two inequalities we use the triangle inequality and the definition of $\normop{\cdot}$, in the third and the last $\norm{d^k} \le \gamma$, in the fourth Assumption~\ref{ch5:ass:L-lipschitz-hessians_chapter_5}, and in the fifth again the triangle inequality. {\color{navy}In the second-to-last equality we use $\hat{x}^k \eqdef u^k x^{k} + (1 - u^k)x^{k-1}$. The bounds~\eqref{ch5:46d7ba91-ef00-4ec2-9b1e-d222e9968289a_chapter_5} and~\eqref{ch5:46d7ba91-ef00-4ec2-9b1e-d222e9968289b_chapter_5} also hold pointwise, without taking expectation.}

        Thus, using~\eqref{ch5:46d7ba91-ef00-4ec2-9b1e-d222e9968289a_chapter_5} and~\eqref{ch5:46d7ba91-ef00-4ec2-9b1e-d222e9968289b_chapter_5}, which give $\thesischfiveE{\sqnorm{\hat{R}_{j+1}}} \le 4 \gamma^2 \min\{L_1, \gamma L_2\}^2$, and $\sqrt{2} \le 2$ we obtain
             \begin{alignat}{2}
                \thesischfiveE{\norm{\sum_{j = 0}^{k - 1} (1 - \eta)^{k - j - 1} \hat{R}_{j + 1}}} \thesischfiveoversetlab{\eqref{ch5:d98da1f4-086e-4004-a809-2eb471da42a6-3_chapter_5}}&{\le} 4 \gamma \min\thesischfiveens{L_1, \gamma L_2} \left( \sum_{j = 0}^{k - 1} (1 - \eta)^{2 (k - j - 1)} \right)^{1/2} \notag\\
                & \le 4 \gamma \min\thesischfiveens{L_1, \gamma L_2} \eta^{-1/2}. \label{ch5:91e3e64b-54a6-46f3-96a5-a9b4f24dc6ee-3_chapter_5}
            \end{alignat}
        
        Finally, for the last term, we can write
            \begin{thesischfivefitalignatnum}{2}
                \thesischfiveE{\norm{\sum_{j = 0}^{k - 1} (1 - \eta)^{k  - j - 1} \hat{S}_{j + 1}}} & \le \left( \thesischfiveE{\norm{\sum_{j = 0}^{k - 1} (1 - \eta)^{k  - j - 1} \hat{S}_{j + 1}}^q} \right)^{1/q} \\
                & \le \left( 2 \sum_{j = 0}^{k - 1} (1 - \eta)^{q (k - j - 1)} \thesischfiveE{\norm{\hat{S}_{j + 1}}^q} \right)^{1/q}, \label{ch5:d98da1f4-086e-4004-a809-2eb471da42a6-2_chapter_5}
            \end{thesischfivefitalignatnum}
        by Jensen's inequality and Lemma~\ref{ch5:appdx-technical-lem:app-von-Bahr-and-Essen_chapter_5}, and, using~Assumption~\ref{ch5:ass:q-bounded-central-moment-hessian_chapter_5} we obtain
            \begin{thesischfivefitalignatnum}{2}
                \thesischfiveE{\norm{\hat{S}_{j + 1}}^q} \thesischfiveoversetlab{\eqref{ch5:891801fc-22ef-4562-81bb-3b5040b0f520_chapter_5}}&{=} {\color{navy}\thesischfiveE{\norm{[\nabla^2 f(\hat{x}^{j+1};\hat{\xi}^{j+1}) - \nabla^2 f(\hat{x}^{j+1})] d^{j+1}}^q}} \\
                & \le {\color{navy}\thesischfiveE{\normop{\nabla^2 f(\hat{x}^{j+1};\hat{\xi}^{j+1}) - \nabla^2 f(\hat{x}^{j+1})}^q \norm{d^{j+1}}^q}} \\
                & \le \gamma^q \sigma_2^q, \label{ch5:f79e21c1-04c0-497e-81e9-c10afd7f04f9a_chapter_5}
            \end{thesischfivefitalignatnum}
        where the last inequality uses $\norm{d^{j+1}} \le \gamma$ and Assumption~\ref{ch5:ass:q-bounded-central-moment-hessian_chapter_5} (conditionally on $\hat{x}^{j+1}$).
        Thus, using~\eqref{ch5:f79e21c1-04c0-497e-81e9-c10afd7f04f9a_chapter_5} and the fact that $2^{1/q} \le 2$ %
        we obtain
             \begin{thesischfivefitalignatnum}{2}
                \thesischfiveE{\norm{\sum_{j = 0}^{k - 1} (1 - \eta)^{k  - j - 1} \hat{S}_{j + 1}}} \thesischfiveoversetlab{\eqref{ch5:d98da1f4-086e-4004-a809-2eb471da42a6-2_chapter_5}}&{\le} 2 \gamma \sigma_2 \left( \sum_{j = 0}^{k - 1} (1 - \eta)^{q (k - j - 1)} \right)^{1/q} \\
                & \le 2 \gamma\sigma_2 \eta^{-1/q}. \label{ch5:91e3e64b-54a6-46f3-96a5-a9b4f24dc6ee-2_chapter_5}
            \end{thesischfivefitalignatnum}

        Then, combining the bounds~\eqref{ch5:cea5881d-8688-41f3-8083-df58a5981198-2_chapter_5},~\eqref{ch5:91e3e64b-54a6-46f3-96a5-a9b4f24dc6ee-3_chapter_5} and~\eqref{ch5:91e3e64b-54a6-46f3-96a5-a9b4f24dc6ee-2_chapter_5} we have
            \begin{thesischfivefitdisplay} \thesischfiveE{\norm{\hat{e}_k}} \le (1 - \eta)^k \thesischfiveE{\norm{\hat{e}_0}} + 2 \sigma_1 \eta^{(p - 1)/p} + 2 \gamma \sigma_2 \eta^{-1/q} + 4 \gamma \min\thesischfiveens{L_1, \gamma L_2} \eta^{-1/2}, \end{thesischfivefitdisplay}
        {\color{navy}which completes the proof of the lemma.}
\end{proof}

\subsection{Proof of~Theorem~\ref{ch5:thm:nsgd-mvr-hess-convergence-analysis_chapter_5}}\label{ch5:appdx-subsec:proof-nsgd-mvr-hess-convergence-analysis_chapter_5}

{\color{navy}Combining Lemmas~\ref{ch5:appdx-lem:unrolling-descent-lemma-nsgd-mvr-hess_chapter_5} and~\ref{ch5:appdx-lem:bounding-error-term-nsgd-mvr-hess_chapter_5} yields the in-expectation convergence analysis of Algorithm~\ref{ch5:algo:nsgd-mvr-hess_chapter_5}.}

\begin{theorem}[Restatement of Theorem~\ref{ch5:thm:nsgd-mvr-hess-convergence-analysis_chapter_5}]
    Under~Assumption~\ref{ch5:ass:lower-boundedness_chapter_5}, Assumption~\ref{ch5:ass:L-lipschitz-gradients_chapter_5}, Assumption~\ref{ch5:ass:L-lipschitz-hessians_chapter_5}, Assumption~\ref{ch5:ass:p-bounded-central-moment-gradient_chapter_5}, and Assumption~\ref{ch5:ass:q-bounded-central-moment-hessian_chapter_5}, {\color{navy}choose, for a sufficiently large universal constant $C>0$,}
        \begin{thesischfivefitalignat}{2}
            \color{navy}K=\Bigg\lceil C\Bigg(&\color{navy}1+\left(\frac{\sigma_1}{\thesischfiveeps}\right)^{p/(p-1)}
            +\frac{\Delta}{\thesischfiveeps}\left(\frac{L_1+\sigma_2}{\thesischfiveeps}
            +\sqrt{\frac{L_2}{\thesischfiveeps}}\right)\\
            &\color{navy}+\frac{\Delta\sigma_2}{\thesischfiveeps^2}
            \left(\frac{\sigma_1}{\thesischfiveeps}\right)^{p/(q(p-1))}
            +\frac{L_2^{1/2}\Delta\sigma_1^{1/4}}{\thesischfiveeps^{7/4}}
            \left(\frac{\sigma_1}{\thesischfiveeps}\right)^{1/(4(p-1))}\Bigg)\Bigg\rceil .
        \end{thesischfivefitalignat}
        {\color{navy}Let the initial gradient estimate $g^0$ be given by}
        \begin{thesischfivefitdisplay} g^0 = \frac{1}{B_{\textnormal{init}}} \sum_{j = 0}^{B_{\textnormal{init}} - 1} \nabla f\left( x^0; \xi_{j}^{0} \right), \end{thesischfivefitdisplay}
    {\color{navy}where}
        \begin{thesischfivefitdisplay}
            {\color{navy}B_{\textnormal{init}}
            =\left\lceil\max\thesischfiveens{1,
            \left(\frac{2\sigma_1}{\thesischfiveeps}\right)^{p/(p-1)}}
            \right\rceil.}
        \end{thesischfivefitdisplay}
    {\color{navy}We use the stepsize}
        \begin{thesischfivefitdisplay}
            {\color{navy}\gamma=\min\thesischfiveens{
            \sqrt{\frac{\Delta}{L_1K}},
            \sqrt{\frac{\Delta\eta^{1/q}}{\sigma_2K}},
            \sqrt[3]{\frac{\Delta\eta^{1/2}}{L_2K}}},}
        \end{thesischfivefitdisplay}
    {\color{navy}and the momentum parameter
    $\eta=\min\thesischfiveens{1,\eta_{\textnormal{eff}}}$, where}
        \begin{align*}
            \eta_{\textnormal{eff}} = \max\Biggl\{ & \left( \frac{\thesischfiveeps}{\sigma_1 K} \right)^{p/(2p - 1)}, \left( \frac{\Delta \sigma_2}{\sigma_1^2 K} \right)^{pq/(p (2 q + 1) - 2 q)}, \\
            & \left( \frac{L_2^{1/2} \Delta}{\sigma_1^{3/2} K} \right)^{(4 p)/(7 p - 6)}\Biggr\}. \thesischfivenumberthis\label{ch5:9fe7782d-3762-4699-bd42-02bdf270abe0_chapter_5}
        \end{align*}
    {\color{navy}When $\sigma_1=0$, we use the deterministic convention $\eta=1$ and $B_{\textnormal{init}}=1$.}
    {\color{navy}A stepsize term with zero denominator is interpreted as $+\infty$ and omitted from the minimum.}
    {\color{navy}For $x_{\mathrm{out}}\sim\mathcal U\{x^0,\ldots,x^{K-1}\}$, Algorithm~\ref{ch5:algo:nsgd-mvr-hess_chapter_5} guarantees $\Exp{\norm{\nabla f(x_{\mathrm{out}})}}\le\thesischfiveeps$ with total sample complexity}
        \begin{thesischfivefitalignat}{2}
            \cO\biggl( & {\color{navy}1+} \left( \frac{\sigma_1}{\thesischfiveeps} \right)^{p/(p - 1)} + \frac{\Delta}{\thesischfiveeps} \left( \frac{L_1 + \sigma_2}{\thesischfiveeps} + \sqrt{\frac{L_2}{\thesischfiveeps}} \right) \\
            & + \frac{\Delta \sigma_2}{\thesischfiveeps^2} \left( \frac{\sigma_1}{\thesischfiveeps} \right)^{p/(q (p - 1))} + {\color{navy}\frac{L_2^{1/2} \Delta \sigma_1^{1/4}}{\thesischfiveeps^{7/4}} \left( \frac{\sigma_1}{\thesischfiveeps} \right)^{1/(4 (p - 1))}} \biggr).
        \end{thesischfivefitalignat}
\end{theorem}

\begin{proof}
    According to~Lemma~\ref{ch5:appdx-lem:unrolling-descent-lemma-nsgd-mvr-hess_chapter_5} we have
        \begin{thesischfivefitdisplay} \frac{1}{K} \sum_{k=0}^{K - 1} \norm{\nabla f(x^{k})} \le \frac{\Delta}{\gamma K} + \frac{2}{K} \sum_{k=0}^{K - 1} \norm{\hat{e}_k} + \frac{\gamma L_1}{2}, \end{thesischfivefitdisplay}
    and using~Lemma~\ref{ch5:appdx-lem:bounding-error-term-nsgd-mvr-hess_chapter_5} this yields
        \begin{thesischfivefitalignat}{2}
            \MoveEqLeft[1] \frac{1}{K} \sum_{k=0}^{K - 1} \thesischfiveE{\norm{\hat{e}_k}} \\
            & \le \frac{1}{K} \sum_{k=0}^{K - 1} (1 - \eta)^k \thesischfiveE{\norm{\hat{e}_0}} + 2 \sigma_1 \eta^{(p - 1)/p} + 2 \gamma \sigma_2 \eta^{-1/q} + 4 \gamma^2 L_2 \eta^{-1/2} \\
            & \le \frac{\thesischfiveE{\norm{\hat{e}_0}}}{\eta K} + 2 \sigma_1 \eta^{(p - 1)/p} + 2 \gamma \sigma_2 \eta^{-1/q} + 4 \gamma^2 L_2 \eta^{-1/2},
        \end{thesischfivefitalignat}
    where in the last term we drop the $\min\thesischfiveens{\cdots}$ to only keep the $\gamma^2 L_2$ term. Hence we obtain the bound
        \begin{thesischfivefitalignat}{2}
            \MoveEqLeft[1] \frac{1}{K} \sum_{k=0}^{K - 1} \thesischfiveE{\norm{\nabla f(x^{k})}} \\
            & \le \frac{\Delta}{\gamma K} + \frac{2 \thesischfiveE{\norm{\hat{e}_0}}}{\eta K} + 4 \sigma_1 \eta^{(p - 1)/p} + 4 \gamma \sigma_2 \eta^{-1/q} + 8 \gamma^2 L_2 \eta^{-1/2} + \frac{\gamma L_1}{2}.
        \end{thesischfivefitalignat}
    Now, using our choice of the stepsize $\gamma$ we have
        \begin{thesischfivefitalignatnum}{2}
            \MoveEqLeft[1] \frac{1}{K} \sum_{k=0}^{K - 1} \thesischfiveE{\norm{\nabla f(x^{k})}} \\
            & = \cO\biggl( \sqrt{\frac{L_1 \Delta}{K}} + \eta^{-1/(2 q)} \sqrt{\frac{\Delta \sigma_2}{K}} + \eta^{-1/6} \left( \frac{L_2^{1/2} \Delta}{K} \right)^{2/3} \\
            &\qquad + \frac{\thesischfiveE{\norm{\hat{e}_0}}}{\eta K} + \sigma_1 \eta^{(p - 1)/p} \biggr). \label{ch5:36a6b109-a2ff-445f-937c-7040717cb89d-2_chapter_5}
        \end{thesischfivefitalignatnum}
    {\color{navy}Now, by our choice of $g^0$, as in the preceding section, we have}
        \begin{thesischfivefitalignatnum}{2}
            \thesischfiveE{\norm{\hat{e}_0}} \thesischfiveoversetlab{\eqref{ch5:6a9db682-cb1e-43bb-878f-e4b60ebd3183_chapter_5}}&{\le}  \frac{2 \sigma_1}{B_{\textnormal{init}}^{(p - 1)/p}}, \label{ch5:c83db8de-7f7d-4c84-82b2-a4e56a709c5f_chapter_5}
        \end{thesischfivefitalignatnum}
    and since {\color{navy}$B_{\textnormal{init}} \ge (2\sigma_1/\thesischfiveeps)^{p/(p - 1)}$}, we have {\color{navy}$\thesischfiveE{\norm{\hat{e}_0}}\le\thesischfiveeps$}. Using also $\eta = \min\{1, \eta_{\textnormal{eff}}\}$ as in the proof of Theorem~\ref{ch5:thm:nsgd-mvr-convergence-analysis_chapter_5}, this gives
        \begin{thesischfivefitalignatnum}{2}
            \MoveEqLeft[1] \frac{1}{K} \sum_{k=0}^{K - 1} \thesischfiveE{\norm{\nabla f(x^{k})}} \\
            \thesischfiveoversetlab{\eqref{ch5:36a6b109-a2ff-445f-937c-7040717cb89d-2_chapter_5}}&{=} \cO\biggl( \sqrt{\frac{L_1 \Delta}{K}} + \eta^{-1/(2 q)} \sqrt{\frac{\Delta \sigma_2}{K}} + \eta^{-1/6} \left( \frac{L_2^{1/2} \Delta}{K} \right)^{2/3} \\
            &\qquad + \frac{\thesischfiveE{\norm{\hat{e}_0}}}{\eta K} + \sigma_1 \eta^{(p - 1)/p} \biggr) \\
            \thesischfiveoversetlab{\eqref{ch5:c83db8de-7f7d-4c84-82b2-a4e56a709c5f_chapter_5}}&{=} \cO\biggl( \sqrt{\frac{L_1 \Delta}{K}} + \eta^{-1/(2 q)} \sqrt{\frac{\Delta \sigma_2}{K}} + \eta^{-1/6} \left( \frac{L_2^{1/2} \Delta}{K} \right)^{2/3} \\
            &\qquad + \frac{\thesischfiveeps}{\eta K} + \sigma_1 \eta^{(p - 1)/p} \biggr) \\
            & = \cO\biggl( \sqrt{\frac{L_1 \Delta}{K}} + \sqrt{\frac{\Delta \sigma_2}{K}} + \eta_{\textnormal{eff}}^{-1/(2 q)} \sqrt{\frac{\Delta \sigma_2}{K}} \\
            &\qquad + \left( \frac{L_2^{1/2} \Delta}{K} \right)^{2/3} + \eta_{\textnormal{eff}}^{-1/6} \left( \frac{L_2^{1/2} \Delta}{K} \right)^{2/3} \\
            &\qquad + \frac{\thesischfiveeps}{K} + \frac{\thesischfiveeps}{K} \eta_{\textnormal{eff}}^{-1} + \sigma_1 \eta_{\textnormal{eff}}^{(p - 1)/p} \biggr) \\
            & = \cO\biggl( \frac{\thesischfiveeps}{K} + \sqrt{\frac{L_1 \Delta}{K}} + \sqrt{\frac{\Delta \sigma_2}{K}} + \left( \frac{L_2^{1/2} \Delta}{K} \right)^{2/3} \\
            &\qquad + \sigma_1 \left( \frac{\Delta \sigma_2}{\sigma_1^2 K} \right)^{q (p - 1)/(p (2 q + 1) - 2 q)} \\
            &\qquad + \sigma_1 \left( \frac{\thesischfiveeps}{\sigma_1 K} \right)^{(p - 1)/(2 p - 1)} + \sigma_1 \left( \frac{L_2^{1/2} \Delta}{\sigma_1^{3/2} K} \right)^{4 (p - 1)/(7 p - 6)} \biggr), \label{ch5:0f0a2e3d-d04e-4e64-948f-d38d2619b1fb-2_chapter_5}
        \end{thesischfivefitalignatnum}
    where in the last equality we use the choice of $\eta_{\textnormal{eff}}$ from~\eqref{ch5:9fe7782d-3762-4699-bd42-02bdf270abe0_chapter_5}; in particular,
        \begin{align*}
            \frac{\thesischfiveeps}{K} \eta_{\textnormal{eff}}^{-1} & \le \sigma_1 \left( \frac{\thesischfiveeps}{\sigma_1 K} \right)^{(p - 1)/(2 p - 1)}, \\
            \eta_{\textnormal{eff}}^{-1/(2 q)} \sqrt{\frac{\Delta \sigma_2}{K}} & \le \sigma_1 \left( \frac{\Delta \sigma_2}{\sigma_1^2 K} \right)^{q (p - 1)/(p (2 q + 1) - 2 q)}, \\
            \eta_{\textnormal{eff}}^{-1/6} \left( \frac{L_2^{1/2} \Delta}{K} \right)^{2/3} & \le \sigma_1 \left( \frac{L_2^{1/2} \Delta}{\sigma_1^{3/2} K} \right)^{4 (p - 1)/(7 p - 6)},
        \end{align*}
    and
        \begin{align*}
            \sigma_1 \eta_{\textnormal{eff}}^{(p - 1)/p} & \le \sigma_1 \left( \frac{\Delta \sigma_2}{\sigma_1^2 K} \right)^{q (p - 1)/(p (2 q + 1) - 2 q)} \\
            &\quad + \sigma_1 \left( \frac{\thesischfiveeps}{\sigma_1 K} \right)^{(p - 1)/(2 p - 1)} + \sigma_1 \left( \frac{L_2^{1/2} \Delta}{\sigma_1^{3/2} K} \right)^{4 (p - 1)/(7 p - 6)}.
        \end{align*}

    Finally, since $\thesischfiveeps/K \le \thesischfiveeps$, from the bound~\eqref{ch5:0f0a2e3d-d04e-4e64-948f-d38d2619b1fb-2_chapter_5} we deduce that the sample complexity of~\ref{ch5:algo:nsgd-mvr-hess_chapter_5} is exactly
         \begin{thesischfivefitalignatnum}{2}
            \cO\biggl( & {\color{navy}1+} \left( \frac{\sigma_1}{\thesischfiveeps} \right)^{p/(p - 1)} + \frac{\Delta}{\thesischfiveeps} \left( \frac{L_1 + \sigma_2}{\thesischfiveeps} + \sqrt{\frac{L_2}{\thesischfiveeps}} \right) \\
            & + \frac{\Delta \sigma_2}{\thesischfiveeps^2} \left( \frac{\sigma_1}{\thesischfiveeps} \right)^{p/(q (p - 1))} + {\color{navy}\frac{L_2^{1/2} \Delta \sigma_1^{1/4}}{\thesischfiveeps^{7/4}} \left( \frac{\sigma_1}{\thesischfiveeps} \right)^{1/(4 (p - 1))}} \biggr), \label{ch5:b5be16f7-9e74-48d0-af57-e31254eaedf7_chapter_5}
        \end{thesischfivefitalignatnum}
    as claimed.
\end{proof}

\begin{remark}
    In particular, when $\sigma_2 = 0$ we observe that the sample complexity~\eqref{ch5:b5be16f7-9e74-48d0-af57-e31254eaedf7_chapter_5} does not depend on the exponent $q$ anymore, as in the lower bound we derived in~Theorem~\ref{ch5:thm:lower-bound-sofo-bounded-central-moments_chapter_5}.
\end{remark}

\newpage
\section{Additional Proofs for Section~\ref{ch5:sec:clipped-nsgd-mvr_chapter_5}}\label{ch5:appdx:proof-clipped-nsgd-mvr_chapter_5}

\subsection{Preliminary Lemmas}

{\color{navy}We now begin the high-probability convergence analysis of Algorithm~\ref{ch5:algo:clipped-nsgd-mvr_chapter_5}, following the proof strategy of~\textcite{sadiev2025second}.}

\subsubsection{Some Descent Lemma}

\begin{lemma}[A Descent Lemma (for~\ref{ch5:algo:clipped-nsgd-mvr_chapter_5})]\label{ch5:appdx-lem:high-probability-analysis-descent-lemma_chapter_5}
    Under~Assumption~\ref{ch5:ass:lower-boundedness_chapter_5} and Assumption~\ref{ch5:ass:mean-squared-smoothness_chapter_5}, for any choice of stepsize $\gamma > 0$ and any choice of momentum parameter $\eta \in \thesischfiveintof{0}{1}$,~\ref{ch5:algo:clipped-nsgd-mvr_chapter_5} generates iterates $\{x^{k}\}_{k \in \thesischfiveInt{0}{K}}$ which satisfy almost surely (a.s.) the inequality
        \begin{thesischfivefitalignat}{2}
            \MoveEqLeft[1] \gamma \sum_{k=0}^{K - 1} \norm{\nabla f(x^{k})} + \Delta_K \\
            & \le \Delta_0 + {\color{navy}2 \gamma \eta \sum_{k=1}^{K-1} \norm{\sum_{j = 1}^k (1 - \eta)^{k - j} \theta_j}} \\
            &\quad + {\color{navy}2 \gamma (1 - \eta) \sum_{k=1}^{K-1} \norm{\sum_{j = 1}^k (1 - \eta)^{k - j} \omega_j}} + \frac{2 \gamma \norm{\hat{e}_0}}{\eta} + \frac{\gamma^2 \bar{L} K}{2},
        \end{thesischfivefitalignat}
    where for any $j \in \thesischfiveInt{1}{K - 1}$, the vectors $\theta_j$ and $\omega_j$ are defined as
        \begin{alignat*}{2}
            \theta_j \eqdef* \thesischfiveclip\left( \nabla f(x^{j};\xi^j), \thesischfivelbd_2 \right) - \nabla f(x^{j}), \thesischfivenumberthis\label{ch5:appdx-lem-proof:def-theta_chapter_5} \\
            \omega_j \eqdef* \thesischfiveclip\left( \nabla f(x^{j};\xi^j) - \nabla f(x^{j-1};\xi^j) , \thesischfivelbd_1 \right) - \left( \nabla f(x^{j}) - \nabla f(x^{j-1}) \right), \thesischfivenumberthis\label{ch5:appdx-lem-proof:def-omega_chapter_5}
        \end{alignat*}
    and for any $k \in \thesischfiveInt{0}{K - 1}$ we let $\Delta_k \eqdef f(x^{k}) - f^{\inf}$.
\end{lemma}

\begin{proof}
    First of all, let us observe that~Lemma~\ref{ch5:appdx-lem:descent-lemma-nsgd-mvr_chapter_5} still holds for~\ref{ch5:algo:clipped-nsgd-mvr_chapter_5} (and can be proved analogously since~\ref{ch5:algo:nsgd-mvr_chapter_5} and \ref{ch5:algo:clipped-nsgd-mvr_chapter_5} have the same gradient update rule). Moreover, as the momentum term does not play any role in the proof of~Lemma~\ref{ch5:appdx-lem:descent-lemma-nsgd-mvr_chapter_5}, the iterates $\{x^{k}\}_{k \in \thesischfiveInt{0}{K}}$ of~\ref{ch5:algo:clipped-nsgd-mvr_chapter_5} satisfy almost surely (a.s.)
        \begin{thesischfivefitdisplay} \gamma \sum_{k=0}^{K - 1} \norm{\nabla f(x^{k})} + \Delta_K \le \Delta_0 + 2 \gamma \sum_{k=0}^{K - 1} \norm{\hat{e}_k} + \frac{\gamma^2 \bar{L} K}{2}. \thesischfivenumberthis\label{ch5:ef85feaf-c75f-4a8a-9afc-c33134870e38_chapter_5} \end{thesischfivefitdisplay}

   Next, we bound $\norm{\hat{e}_k}$ in an analogous way as we did in~Lemma~\ref{ch5:appdx-lem:bounding-error-term-nsgd-mvr_chapter_5}. According to the update rule for the momentum term in~\ref{ch5:algo:clipped-nsgd-mvr_chapter_5}, we have
        \begin{thesischfivefitalignatnum}{1}
            \hat{e}_k
            & = g^{k} - \nabla f(x^{k}) \\
            & = (1 - \eta) \left(
                g^{k-1}
                + \thesischfiveclip\left(
                    \nabla f(x^{k};\xi^{k})-\nabla f(x^{k-1};\xi^{k}),
                    \thesischfivelbd_1
                \right)
            \right) \\
            & \quad
                + \eta\,\thesischfiveclip\left(
                    \nabla f(x^{k};\xi^{k}),\thesischfivelbd_2
                \right)
                - \nabla f(x^{k}) \\
            & = (1 - \eta)\left(g^{k-1}-\nabla f(x^{k-1})\right) \\
            & \quad
                + (1 - \eta)\Bigl(
                    \thesischfiveclip\left(
                        \nabla f(x^{k};\xi^{k})-\nabla f(x^{k-1};\xi^{k}),
                        \thesischfivelbd_1
                    \right) \\
            & \hspace{9em}
                    - \left[\nabla f(x^{k})-\nabla f(x^{k-1})\right]
                \Bigr) \\
            & \quad
                + \eta\left(
                    \thesischfiveclip\left(
                        \nabla f(x^{k};\xi^{k}),\thesischfivelbd_2
                    \right)
                    - \nabla f(x^{k})
                \right) \\
            \thesischfiveoversetlab{\eqref{ch5:appdx-lem-proof:def-theta_chapter_5}+\eqref{ch5:appdx-lem-proof:def-omega_chapter_5}}&{=} (1 - \eta) \hat{e}_{k - 1} + \eta \theta_k + (1 - \eta) \omega_k \\
            & =
                (1 - \eta)^k \hat{e}_0
                + \eta \sum_{j = 1}^k (1 - \eta)^{k - j} \theta_j \\
            & \quad
                + (1 - \eta)\sum_{j = 1}^k (1 - \eta)^{k - j}\omega_j,
                \label{ch5:6f49529c-cbca-40f2-af32-aea3df8c50f2_chapter_5}
        \end{thesischfivefitalignatnum}
    where in the last equality we unroll the recursion. Then, taking the norm and applying the triangle inequality in this series of equalities, we obtain
        \begin{alignat}{2}
            \norm{\hat{e}_k} & \le (1 - \eta)^k \norm{\hat{e}_0} + \eta \norm{\sum_{j = 1}^k (1 - \eta)^{k - j} \theta_j} \notag\\
            & \quad + (1 - \eta) \norm{\sum_{j = 1}^k (1 - \eta)^{k - j} \omega_j}, \label{ch5:cc2cc25c-1b54-471f-ab22-55d547ad29a6_chapter_5}%
        \end{alignat}
    and plugging~\eqref{ch5:cc2cc25c-1b54-471f-ab22-55d547ad29a6_chapter_5} in~\eqref{ch5:ef85feaf-c75f-4a8a-9afc-c33134870e38_chapter_5} leads to
        \begin{thesischfivefitalignat}{1}
            \gamma \sum_{k=0}^{K - 1} \norm{\nabla f(x^{k})} + \Delta_K & \le \Delta_0 + 2 \gamma \sum_{k=0}^{K - 1} \norm{\hat{e}_k} + \frac{\gamma^2 \bar{L} K}{2} \\
            & \le \Delta_0
                + 2 \gamma \sum_{k=0}^{K - 1}(1-\eta)^k\norm{\hat{e}_0} \\
            & \quad
                + 2 \gamma\eta\sum_{k=0}^{K - 1}
                    \norm{\sum_{j=1}^k(1-\eta)^{k-j}\theta_j} \\
            & \quad
                + 2 \gamma(1-\eta)\sum_{k=0}^{K - 1}
                    \norm{\sum_{j=1}^k(1-\eta)^{k-j}\omega_j}
                + \frac{\gamma^2\bar{L}K}{2} \\
            & \le
                \Delta_0+\frac{2\gamma\norm{\hat{e}_0}}{\eta}
                +\frac{\gamma^2\bar{L}K}{2} \\
            & \quad
                + 2\gamma\eta\sum_{k=0}^{K - 1}
                    \norm{\sum_{j=1}^k(1-\eta)^{k-j}\theta_j} \\
            & \quad
                + 2\gamma(1-\eta)\sum_{k=0}^{K - 1}
                    \norm{\sum_{j=1}^k(1-\eta)^{k-j}\omega_j},
        \end{thesischfivefitalignat}
    where in the last inequality we use
        \begin{thesischfivefitdisplay} 2 \gamma \sum_{k=0}^{K - 1} (1 - \eta)^k \norm{\hat{e}_0} \le 2 \gamma \norm{\hat{e}_0} \sum_{k \ge 0} (1 - \eta)^k = \frac{2 \gamma \norm{\hat{e}_0}}{\eta}. \end{thesischfivefitdisplay}

    This proves the desired result.
\end{proof}

\begin{remark}
    In particular, if we assume $g^0 = 0$ in~\ref{ch5:algo:clipped-nsgd-mvr_chapter_5} then $x^0 = x^1$ and {\color{navy}$\Delta_1 = \Delta_0$}. Moreover, we have
        \begin{thesischfivefitdisplay} \norm{\hat{e}_0} = \norm{g^0 - \nabla f(x^0)} = \norm{\nabla f(x^1)} \le \sqrt{2 \bar{L} \Delta_0}, \end{thesischfivefitdisplay}
    as by~Assumption~\ref{ch5:ass:mean-squared-smoothness_chapter_5} we know that $f$ has $\bar{L}$--Lipschitz continuous gradients.
\end{remark}

\begin{corollary}[Another Descent Lemma]\label{ch5:appdx-lem:high-probability-analysis-descent-lemma-2_chapter_5}
    Under~Assumption~\ref{ch5:ass:lower-boundedness_chapter_5} and Assumption~\ref{ch5:ass:L-lipschitz-gradients_chapter_5}, for any choice of stepsize $\gamma > 0$ and any choice of momentum parameter $\eta \in \thesischfiveintof{0}{1}$,~\ref{ch5:algo:clipped-nsgd-mvr_chapter_5} generates iterates $\{x^{k}\}_{k \in \thesischfiveInt{0}{K}}$ which satisfy almost surely (a.s.) the inequality
        \begin{thesischfivefitalignat}{2}
            \MoveEqLeft[1] \gamma \sum_{k=0}^{K - 1} \norm{\nabla f(x^{k})} + \Delta_K \\
            & \le \Delta_0 + {\color{navy}2 \gamma \eta \sum_{k=1}^{K-1} \norm{\sum_{j = 1}^k (1 - \eta)^{k - j} \theta_j}} \\
            &\quad + {\color{navy}2 \gamma (1 - \eta) \sum_{k=1}^{K-1} \norm{\sum_{j = 1}^k (1 - \eta)^{k - j} \omega_j}} + \frac{2 \gamma \norm{\hat{e}_0}}{\eta} + \frac{\gamma^2 L_1 K}{2},
        \end{thesischfivefitalignat}
    where for any $j \in \thesischfiveInt{1}{K - 1}$, the vectors $\theta_j$ and $\omega_j$ are defined in~\eqref{ch5:appdx-lem-proof:def-theta_chapter_5} and~\eqref{ch5:appdx-lem-proof:def-omega_chapter_5} and for any $k \in \thesischfiveInt{0}{K - 1}$ we let $\Delta_k \eqdef f(x^{k}) - f^{\inf}$.    
\end{corollary}

\begin{proof}
    The proof follows the exact same steps as in the previous lemma (Lemma~\ref{ch5:appdx-lem:high-probability-analysis-descent-lemma_chapter_5}), with the exception that the function $f$ has $L_1$--Lipschitz continuous gradients (instead of $\bar{L}$).
\end{proof}

\begin{lemma}[A Descent Lemma (for~\ref{ch5:algo:clipped-nsgd-mvr-hess_chapter_5})]\label{ch5:appdx-lem:high-probability-analysis-descent-lemma-hess_chapter_5}
    Under~Assumption~\ref{ch5:ass:lower-boundedness_chapter_5} and Assumption~\ref{ch5:ass:L-lipschitz-gradients_chapter_5}, for any choice of stepsize $\gamma > 0$ and any choice of momentum parameter $\eta \in \thesischfiveintof{0}{1}$,~\ref{ch5:algo:clipped-nsgd-mvr-hess_chapter_5} generates iterates $\{x^{k}\}_{k \in \thesischfiveInt{0}{K}}$ which satisfy almost surely (a.s.) the inequality
        \begin{thesischfivefitalignatnum}{2}
            \MoveEqLeft[1] \gamma \sum_{k=0}^{K - 1} \norm{\nabla f(x^{k})} + \Delta_K \\
            & \le \Delta_0 + {\color{navy}2 \gamma \eta \sum_{k=1}^{K-1} \norm{\sum_{j = 1}^k (1 - \eta)^{k - j} \theta_j}} \\
            &\quad + {\color{navy}2 \gamma (1 - \eta) \sum_{k=1}^{K-1} \norm{\sum_{j = 1}^k (1 - \eta)^{k - j} \omega_j}} + \frac{2 \gamma \norm{\hat{e}_0}}{\eta} + \frac{\gamma^2 L_1 K}{2}, \label{ch5:5a8226aa-756e-43a3-92e4-1096c9fe1a58_chapter_5}
        \end{thesischfivefitalignatnum}
    where for any $j \in \thesischfiveInt{1}{K - 1}$, the vectors $\theta_j$ and $\omega_j$ are defined as
        \begin{alignat*}{2}
            \theta_j \eqdef* \thesischfiveclip\left( \nabla f(x^{j};\xi^j), \thesischfivelbd_2 \right) - \nabla f(x^{j}), \thesischfivenumberthis\label{ch5:appdx-lem-proof:def-theta-hess_chapter_5} \\
            \omega_j \eqdef* \thesischfiveclip\left( \nabla^2 f(\hat{x}^j;\hat{\xi}^j)(x^{j} - x^{j-1}), \thesischfivelbd_1 \right) - \left( \nabla f(x^{j}) - \nabla f(x^{j-1}) \right), \thesischfivenumberthis\label{ch5:appdx-lem-proof:def-omega-hess_chapter_5}
        \end{alignat*}
    and for any $k \in \thesischfiveInt{0}{K - 1}$ we let $\Delta_k \eqdef f(x^{k}) - f^{\inf}$.   
\end{lemma}

\begin{proof}
    The proof is very similar to the proof of~Lemma~\ref{ch5:appdx-lem:high-probability-analysis-descent-lemma_chapter_5}, except that $f$ has $L_1$--Lipschitz continuous gradients (instead of $\bar{L}$) and the computation~\eqref{ch5:6f49529c-cbca-40f2-af32-aea3df8c50f2_chapter_5} changes slightly because it uses the Hessian term {\color{navy}$\nabla^2 f(\hat{x}^j;\hat{\xi}^j)(x^{j} - x^{j-1})$} instead of the gradient difference {\color{navy}$\nabla f(x^{j};\xi^j) - \nabla f(x^{j-1};\xi^j)$}. We thus have
        \begin{thesischfivefitalignat}{2}
            \hat{e}_k & = g^{k} - \nabla f(x^{k}) \\
            & = (1 - \eta) \left( g^{k-1} + {\color{navy}\thesischfiveclip(\nabla^2 f(\hat{x}^k;\hat{\xi}^{k}) (x^{k} - x^{k-1}), \thesischfivelbd_1)} \right) \\
            &\quad + \eta \, \thesischfiveclip(\nabla f(x^{k};\xi^{k}), \thesischfivelbd_2) - \nabla f(x^{k}) \\
            & = (1 - \eta) \left( g^{k-1} - \nabla f(x^{k-1}) \right) \\
            &\quad + (1 - \eta) \Bigl( {\color{navy}\thesischfiveclip(\nabla^2 f(\hat{x}^k;\hat{\xi}^{k}) (x^{k} - x^{k-1}), \thesischfivelbd_1)} \\
            &\hspace{9em} - \left[ \nabla f(x^{k}) - \nabla f(x^{k-1}) \right] \Bigr) \\
            &\quad + \eta \left( \thesischfiveclip(\nabla f(x^{k};\xi^{k}), \thesischfivelbd_2) - \nabla f(x^{k}) \right) \\
            \thesischfiveoversetlab{\eqref{ch5:appdx-lem-proof:def-theta-hess_chapter_5}+\eqref{ch5:appdx-lem-proof:def-omega-hess_chapter_5}}&{=} (1 - \eta) \hat{e}_{k - 1} + \eta \theta_k + (1 - \eta) \omega_k \\
            & = (1 - \eta)^k \hat{e}_0 + \eta \sum_{j = 1}^k (1 - \eta)^{k - j} \theta_j \\
            &\quad + (1 - \eta) \sum_{j = 1}^k (1 - \eta)^{k - j} \omega_j,
        \end{thesischfivefitalignat}
    where in the last equality we unroll the recursion. The next steps are exactly as in the proof of~Lemma~\ref{ch5:appdx-lem:high-probability-analysis-descent-lemma_chapter_5} and the claimed result~\eqref{ch5:5a8226aa-756e-43a3-92e4-1096c9fe1a58_chapter_5} follows.
\end{proof}

\subsubsection{High-Probability Analysis}

From the previous~Lemma~\ref{ch5:appdx-lem:high-probability-analysis-descent-lemma_chapter_5}, we need to bound the two new terms:
    \begin{thesischfivefitdisplay} \norm{\sum_{j = 1}^t (1 - \eta)^{t - j} \theta_j} \quad \text{and} \quad \norm{\sum_{j = 1}^t (1 - \eta)^{t - j} \omega_j}, \end{thesischfivefitdisplay}
where $t \in \thesischfiveInt{0}{K - 1}$. To do so, we use the same strategy as in previous works~\parencite{10.5555/3495724.3496985,pmlr-v202-sadiev23a,pmlr-v195-liu23c,sadiev2025second} and we introduce the biased and unbiased parts of $\theta_j$ and $\omega_j$, that is, for any $j \in [K - 1]$ we let $\theta_j = \theta_j^b + \theta_j^u$ with
    \begin{thesischfivefitdisplay}
        {\color{navy}\theta_j^b \eqdef
        \thesischfiveExpSub{\xi^{j}}{\thesischfiveclip\left(
        \nabla f(x^{j};\xi^{j}),\thesischfivelbd_2\right)}
        -\nabla f(x^{j}),}
        \thesischfivenumberthis\label{ch5:647ac97a-ff0a-48a4-a4ed-13d555d50876_chapter_5}
    \end{thesischfivefitdisplay}
    \begin{thesischfivefitdisplay}
        {\color{navy}\theta_j^u \eqdef
        \thesischfiveclip\left(\nabla f(x^{j};\xi^{j}),\thesischfivelbd_2\right)
        -\thesischfiveExpSub{\xi^{j}}{\thesischfiveclip\left(
        \nabla f(x^{j};\xi^{j}),\thesischfivelbd_2\right)},}
        \thesischfivenumberthis\label{ch5:0ea28dd7-2e3c-42c1-af0a-754531b8eac2_chapter_5}
    \end{thesischfivefitdisplay}
and $\omega_j = \omega_j^b + \omega_j^u$ where
    \begin{thesischfivefitdisplay}
        {\color{navy}\omega_j^b \eqdef
        \thesischfiveExpSub{\xi^{j}}{\thesischfiveclip\left(
        \nabla f(x^{j};\xi^j)-\nabla f(x^{j-1};\xi^j),\thesischfivelbd_1\right)}
        -\left(\nabla f(x^{j})-\nabla f(x^{j-1})\right),}
        \thesischfivenumberthis\label{ch5:9715675c-c075-4d13-a6e5-978530f178fd_chapter_5}
    \end{thesischfivefitdisplay}
    \begin{thesischfivefitalignatnum}{2}
        \color{navy}\omega_j^u & \color{navy}\eqdef \thesischfiveclip\left(\nabla f(x^{j};\xi^j)-\nabla f(x^{j-1};\xi^j),\thesischfivelbd_1\right) \\
        &\quad \color{navy}- \thesischfiveExpSub{\xi^{j}}{\thesischfiveclip\left(\nabla f(x^{j};\xi^j)-\nabla f(x^{j-1};\xi^j),\thesischfivelbd_1\right)}, \label{ch5:5bf5f4fb-a644-41c4-a212-9c4930a0075d_chapter_5}
    \end{thesischfivefitalignatnum}
for~Lemma~\ref{ch5:appdx-lem:high-probability-analysis-descent-lemma_chapter_5} and~Lemma~\ref{ch5:appdx-lem:high-probability-analysis-descent-lemma-2_chapter_5}. For~Lemma~\ref{ch5:appdx-lem:high-probability-analysis-descent-lemma-hess_chapter_5} (hessian clipping) we let $\omega_j = \omega_j^r + \omega_j^b + \omega_j^u$ where
    \begin{thesischfivefitdisplay}
        {\color{navy}\omega_j^r \eqdef
        \nabla^2 f(\hat{x}^j)(x^{j}-x^{j-1})
        -\left(\nabla f(x^{j})-\nabla f(x^{j-1})\right),}
        \thesischfivenumberthis\label{ch5:9715675c-c075-4d13-a6e5-978530f178fe-hess_chapter_5}
    \end{thesischfivefitdisplay}
    \begin{thesischfivefitdisplay}
        {\color{navy}\omega_j^b \eqdef
        \thesischfiveExpSub{\hat{\xi}^{j}}{\thesischfiveclip\left(
        \nabla^2 f(\hat{x}^j;\hat{\xi}^j)(x^{j}-x^{j-1}),\thesischfivelbd_1\right)}
        -\nabla^2 f(\hat{x}^j)(x^{j}-x^{j-1}),}
        \thesischfivenumberthis\label{ch5:9715675c-c075-4d13-a6e5-978530f178fd-hess_chapter_5}
    \end{thesischfivefitdisplay}
    \begin{thesischfivefitalignatnum}{2}
        \color{navy}\omega_j^u & \color{navy}\eqdef \thesischfiveclip\left(\nabla^2 f(\hat{x}^j;\hat{\xi}^j)(x^{j}-x^{j-1}),\thesischfivelbd_1\right) \\
        &\quad \color{navy}- \thesischfiveExpSub{\hat{\xi}^{j}}{\thesischfiveclip\left(\nabla^2 f(\hat{x}^j;\hat{\xi}^j)(x^{j}-x^{j-1}),\thesischfivelbd_1\right)}, \label{ch5:5bf5f4fb-a644-41c4-a212-9c4930a0075d-hess_chapter_5}
    \end{thesischfivefitalignatnum}

\begin{lemma}\label{ch5:appdx-lem:high-probability-analysis-bound-1_chapter_5}
    Under~Assumption~\ref{ch5:ass:p-bounded-central-moment-gradient_chapter_5}, for any $t \in \thesischfiveInt{0}{K - 1}$ and any $\beta' \in \thesischfiveintof{0}{1/2}$, if the clipping threshold satisfies
        \begin{thesischfivefitdisplay} \thesischfivelbd_2 \ge \max\thesischfiveens{2 \norm{\nabla f(x^{j})}, \sigma_1 \eta^{-1/p}}, \end{thesischfivefitdisplay}
    for all $j \in [t]$ then, with probability at least $1 - 2 \beta'$, we have
        \begin{thesischfivefitdisplay} \norm{\sum_{j = 1}^t (1 - \eta)^{t - j} \theta_j} \le 22 \thesischfivelbd_2 \log(2/\beta'). \end{thesischfivefitdisplay}
\end{lemma}

\begin{proof}
    First of all, using~\eqref{ch5:647ac97a-ff0a-48a4-a4ed-13d555d50876_chapter_5} and~\eqref{ch5:0ea28dd7-2e3c-42c1-af0a-754531b8eac2_chapter_5} it follows
        \begin{thesischfivefitdisplay} \norm{\sum_{j = 1}^t (1 - \eta)^{t - j} \theta_j} \le \underbrace{\norm{\sum_{j = 1}^t (1 - \eta)^{t - j} \theta_j^b}}_{\Circled{4}} + \norm{\sum_{j = 1}^t (1 - \eta)^{t - j} \theta_j^u}, \end{thesischfivefitdisplay}
    and now, we need to bound both terms above. For the second term, we use~Lemma~\ref{ch5:appdx-technical-lem:cutkosky-mehta_chapter_5} with exponent $2$ to obtain
        \begin{thesischfivefitalignatnum}{2}
            \MoveEqLeft[1] \norm{\sum_{j = 1}^t (1 - \eta)^{t - j} \theta_j^u} \\
            & \le \abs{\sum_{j = 1}^t V_j^t} + \sqrt{\max_{j \in [t]} \sqnorm{(1 - \eta)^{t - j} \theta_j^u} + \sum_{j = 1}^t \sqnorm{(1 - \eta)^{t - j} \theta_j^u}} \\
            & \le \abs{\sum_{j = 1}^t V_j^t} + \sqrt{2 \sum_{j = 1}^t \sqnorm{(1 - \eta)^{t - j} \theta_j^u}} \\
            & = \underbrace{\abs{\sum_{j = 1}^t V_j^t}}_{\Circled{1}} + \sqrt{2 \Bigl( \underbrace{\sum_{j = 1}^t Y_j^t}_{\Circled{2}} + \underbrace{\sum_{j = 1}^t \thesischfiveExpSub{\xi^j}{\sqnorm{(1 - \eta)^{t - j} \theta_j^u}}}_{\Circled{3}} \Bigr)}, \label{ch5:21eaf977-44d0-4b4b-8b63-5af59e5cc089_chapter_5}
        \end{thesischfivefitalignatnum}
    where in the first inequality we use Lemma~\ref{ch5:appdx-technical-lem:cutkosky-mehta_chapter_5} with the sequence defined, using the partial sums $S_{j-1}^t \eqdef \sum_{i = 1}^{j - 1} (1 - \eta)^{t - i} \theta_i^u$, by
        \begin{thesischfivefitdisplay}
            V_j^t \eqdef
            \begin{cases}
                \thesischfivesign\left( \sum_{i = 1}^{j - 1} V_i^t \right) \thesischfiveps{S_{j-1}^t}{(1 - \eta)^{t - j} \theta_j^u} / \norm{S_{j-1}^t}, & \text{if } S_{j-1}^t \neq 0; \\
                0, & \text{if } S_{j-1}^t = 0.
            \end{cases}
        \end{thesischfivefitdisplay}
    while in the equality, for any $j \in [t]$ we let
        \begin{thesischfivefitdisplay} Y_j^t \eqdef \sqnorm{(1 - \eta)^{t - j} \theta_j^u} - \thesischfiveExpSub{\xi^j}{\sqnorm{(1 - \eta)^{t - j} \theta_j^u}}. \end{thesischfivefitdisplay}

    We now bound all terms \Circled{1}, \Circled{2}, \Circled{3} and \Circled{4}.

    \paragraph{Bounding \Circled{1}.} The sequence $V_1^t, \ldots, V_t^t$ is a martingale difference sequence since, by definition of $\theta_j^u$ and $V_j^t$, for all $j \in [t]$ we have $\thesischfiveExpCond{V_j^t}{V_{j - 1}^t, \ldots, V_1^t} = 0$ since $\thesischfiveExpSub{\xi^j}{\theta_j^u} = 0$. Moreover, by~Lemma~\ref{ch5:appdx-technical-lem:cutkosky-mehta_chapter_5} we also have
        \begin{thesischfivefitdisplay} \abs{V_j^t} \le \norm{(1 - \eta)^{t - j} \theta_j^u} \le \norm{\theta_j^u} \le 2 \thesischfivelbd_2, \thesischfivenumberthis\label{ch5:77cc40ad-0220-4787-ab09-7c07a2d6dfa7_chapter_5} \end{thesischfivefitdisplay}
    where the last inequality follows from Lemma~\ref{ch5:appdx-technical-lem:sadiev-generalization_chapter_5}, and let $c_2 \eqdef 2 \thesischfivelbd_2$ be the upper bound on the random variables $V_1^t, \ldots, V_t^t$. Additionally, if we denote by $s_j^2 \eqdef \thesischfiveExpCond{\left( V_j^t \right)^2}{V_{j - 1}^t, \ldots, V_1^t}$ then by~Lemma~\ref{ch5:appdx-technical-lem:freedman-inequality_chapter_5} for any real number $b_2 > 0$ and any $G_2> 0$, as $0 < \beta' \le 1/2$ then $\log(2/\beta') \ge 1$ and we have
        \begin{thesischfivefitalignatnum}{2}
            \MoveEqLeft[1] \thesischfiveProba{\abs{\sum_{j = 1}^t V_j^t} > b_2 \quad \text{and} \quad \sum_{j = 1}^t s_j^2 \le G_2 \log(2/\beta')} \\
            & \le 2 \exp\left( -\frac{b_2^2}{2 G_2 \log(2/\beta') + 2 b_2 c_2/3} \right) = \beta', \label{ch5:158c1570-485a-49ca-9ec3-32f47e8d8e8b_chapter_5}
        \end{thesischfivefitalignatnum}
    where the last equality holds provided $b_2 \eqdef \bigl( c_2/3 + \sqrt{c_2^2/9 + 2 G_2} \bigr) \log(2/\beta') > 0$. We now need to define the constant $G_2$. To do so, we need to bound the sum $s_1^2 + \cdots + s_t^2$, this gives
        \begin{thesischfivefitalignatnum}{2}
            \sum_{j = 1}^t s_j^2 & = \sum_{j = 1}^t \thesischfiveExpCond{\left( V_j^t \right)^2}{V_{j - 1}^t, \ldots, V_1^t} \\
            \thesischfiveoversetlab{\eqref{ch5:77cc40ad-0220-4787-ab09-7c07a2d6dfa7_chapter_5}}&{\le} \sum_{j = 1}^t \thesischfiveExpSub{\xi^j}{\sqnorm{(1 - \eta)^{t - j} \theta_j^u}} \\
            & = \sum_{j = 1}^t ( 1 - \eta)^{2 (t - j)} \, \thesischfiveExpSub{\xi^j}{\sqnorm{\theta_j^u}} \\
            & \le 18 \thesischfivelbd_2^{2 - p} \sigma_1^p \sum_{j = 1}^t (1 - \eta)^{2 (t - j)} \\
            & \le \frac{18 \thesischfivelbd_2^{2 - p} \sigma_1^p}{1 - (1 - \eta)^2} \\
            & \le \frac{18 \thesischfivelbd_2^{2 - p} \sigma_1^p}{\eta}, \label{ch5:74b4e2bb-4754-4dd3-ae90-37ef40362f93_chapter_5}
        \end{thesischfivefitalignatnum}
    where in the application of~Lemma~\ref{ch5:appdx-technical-lem:sadiev-generalization_chapter_5} we require $\thesischfivelbd_2 \ge 2 \max_{j \in [t]} \norm{\nabla f(x^{j})}$. Setting $G_2 \eqdef 18 \thesischfivelbd_2^{2 - p} \sigma_1^p/\eta > 0$ gives the desired bound~\eqref{ch5:158c1570-485a-49ca-9ec3-32f47e8d8e8b_chapter_5}.

    \paragraph{Bounding \Circled{2}.} As in the previous paragraph, the sequence $Y_1^t, \ldots, Y_t^t$ is also a martingale difference sequence, since $\thesischfiveExpSub{\xi^j}{Y_j^t} = 0$ for all $j \in [t]$ by definition of $Y_j^t$. Moreover, according to~Lemma~\ref{ch5:appdx-technical-lem:sadiev-generalization_chapter_5}, we also have, as we did in~\eqref{ch5:77cc40ad-0220-4787-ab09-7c07a2d6dfa7_chapter_5}, for any $j \in [t]$
        \begin{thesischfivefitdisplay} \abs{Y_j^t} \le \sqnorm{(1 - \eta)^{t - j} \theta_j^u} + \thesischfiveExpSub{\xi^j}{\sqnorm{(1 - \eta)^{t - j} \theta_j^u}} \le 4 \thesischfivelbd_2^2 + 4 \thesischfivelbd_2^2 = 8 \thesischfivelbd_2^2. \thesischfivenumberthis\label{ch5:4923f3a4-39e5-43f3-b260-7deeb46a8c67_chapter_5} \end{thesischfivefitdisplay}
    Hence we define $\widetilde{c}_2 \eqdef 8 \thesischfivelbd_2^2$. Now, denoting the conditional variance of $Y_j^t$ as $\widetilde{s}_j^2 \eqdef \thesischfiveExpCond{\left( Y_j^t \right)^2}{Y_{j - 1}^t, \ldots, Y_1^t} = \thesischfiveExpSub{\xi^j}{\left( Y_j^t \right)^2}$ we have the bound
        \begin{thesischfivefitalignatnum}{2}
            \widetilde{s}_j^2 \thesischfiveoversetlab{\eqref{ch5:4923f3a4-39e5-43f3-b260-7deeb46a8c67_chapter_5}}&{\le} 8 \thesischfivelbd_2^2 \, \thesischfiveExpSub{\xi^j}{\abs{Y_j^t}} \le 16 \thesischfivelbd_2^2 \, \thesischfiveExpSub{\xi^j}{\sqnorm{(1 - \eta)^{t - j} \theta_j^u}}, \label{ch5:a26983bb-60b7-49ab-ac0e-a1d026f1c51e_chapter_5}
        \end{thesischfivefitalignatnum}
    for all $j \in [t]$. Hence, 
        \begin{thesischfivefitalignat}{2}
            \sum_{j = 1}^t \widetilde{s}_j^2 \thesischfiveoversetlab{\eqref{ch5:a26983bb-60b7-49ab-ac0e-a1d026f1c51e_chapter_5}}&{\le} 16 \thesischfivelbd_2^2 \sum_{j = 1}^t \thesischfiveExpSub{\xi^j}{\sqnorm{(1 - \eta)^{t - j} \theta_j^u}} \\
            \thesischfiveoversetlab{\eqref{ch5:74b4e2bb-4754-4dd3-ae90-37ef40362f93_chapter_5}}&{\le} 16 \thesischfivelbd_2^2 \cdot \frac{18 \thesischfivelbd_2^{2 - p} \sigma_1^p}{\eta} = \frac{288 \thesischfivelbd_2^{4 - p} \sigma_1^p}{\eta},
        \end{thesischfivefitalignat}
    (where we require $\thesischfivelbd_2 \ge 2 \max_{j \in [t]} \norm{\nabla f(x^{j})}$). Next, if we define $\widetilde{G}_2 \eqdef 288 \thesischfivelbd_2^{4 - p} \sigma_1^p/\eta$ then, applying~Lemma~\ref{ch5:appdx-technical-lem:freedman-inequality_chapter_5} we obtain for any real number $\widetilde{b}_2$
        \begin{thesischfivefitalignatnum}{2}
            \MoveEqLeft[1] \thesischfiveProba{\abs{\sum_{j = 1}^t Y_j^t} > \widetilde{b}_2 \quad \text{and} \quad \sum_{j = 1}^t \widetilde{s}_j^2 \le \widetilde{G}_2 \log(2/\beta')} \\
            & \le 2 \exp\left( -\frac{\widetilde{b}_2^2}{2 \widetilde{G}_2 \log(2/\beta') + 2 \widetilde{b}_2 \widetilde{c}_2/3} \right) = \beta', \label{ch5:04635ad8-7080-4c15-a1de-f0ada3a80f78_chapter_5}
        \end{thesischfivefitalignatnum}
    where the last identity holds if we set $\widetilde{b}_2 \eqdef \bigl( \widetilde{c}_2/3 + \sqrt{\widetilde{c}_2^2/9 + 2 \widetilde{G}_2} \bigr) \log(2/\beta')$. This establishes the desired bound.

    \paragraph{Bounding \Circled{3}.} {\color{navy}As in the two preceding paragraphs, we have}
        \begin{thesischfivefitdisplay} \Circled{3} \eqdef \sum_{j = 1}^t \thesischfiveExpSub{\xi^j}{\sqnorm{(1 - \eta)^{t - j} \theta_j^u}} \thesischfiveoversetlab{\eqref{ch5:74b4e2bb-4754-4dd3-ae90-37ef40362f93_chapter_5}}{\le} \frac{18 \thesischfivelbd_2^{2 - p} \sigma_1^p}{\eta}, \end{thesischfivefitdisplay}
    as desired.

    \paragraph{Bounding \Circled{4}.} If we assume that $\thesischfivelbd_2 \ge 2 \norm{\nabla f(x^{j})}$ for all $j \in [t]$ then, with probability one we have
        \begin{thesischfivefitalignat}{2}
            \Circled{4} \eqdef* \norm{\sum_{j = 1}^t (1 - \eta)^{t - j} \theta_j^b} \\
            & \le \sum_{j = 1}^t (1 - \eta)^{t - j} \norm{\theta_j^b} \\
            & \le 4 \thesischfivelbd_2^{1 - p} \sigma_1^p \sum_{j = 1}^t (1 - \eta)^{t - j} \\
            & \le \frac{4 \thesischfivelbd_2^{1 - p} \sigma_1^p}{\eta},
        \end{thesischfivefitalignat}
    where in the second inequality we use~Lemma~\ref{ch5:appdx-technical-lem:sadiev-generalization_chapter_5}, more precisely, for any $j \in [t]$,
        \begin{thesischfivefitdisplay} \norm{\theta_j^b} \le 2^p \thesischfivelbd_2^{1 - p} \sigma_1^p \le 4 \thesischfivelbd_2^{1 - p} \sigma_1^p. \end{thesischfivefitdisplay}

    \paragraph{Summing up all bounds \Circled{1}, \Circled{2}, \Circled{3} and \Circled{4}.} To sum up, we introduce the event $E_{\Circled{1}, t}$ as follows
        \begin{thesischfivefitdisplay} E_{\Circled{1}, t} \eqdef \thesischfiveens{\abs{\sum_{j = 1}^t V_j^t} \le b_2 \quad \text{or} \quad \sum_{j = 1}^t s_j^2 > G_2 \log(2/\beta')}, \end{thesischfivefitdisplay}
    where $c_2$, $G_2$ and $b_2$ are defined above, and we can bound $b_2$ as
        \begin{thesischfivefitalignatnum}{2}
            b_2 & = \left( \frac{c_2}{3} + \sqrt{\frac{c_2^2}{9} + 2 G_2} \right) \log(2/\beta') \\
            & \le \left( \frac{2 c_2}{3} + \sqrt{2 G_2} \right) \log(2/\beta') \\
            & = \left( \frac{4 \thesischfivelbd_2}{3} + \sqrt{\frac{36 \thesischfivelbd_2^{2 - p} \sigma_1^p}{\eta}} \right) \log(2/\beta') \\
            & = \thesischfivelbd_2 \left( \frac{4}{3} + 6 \sqrt{\frac{1}{\eta} \left(\frac{\sigma_1}{\thesischfivelbd_2} \right)^p} \right) \log(2/\beta') \\
            & \le \frac{22 \thesischfivelbd_2}{3} \log(2/\beta'), \label{ch5:35c58053-3ad9-4b62-ae23-2a6f658c08da_chapter_5}
        \end{thesischfivefitalignatnum}
    where the last inequality holds provided $\thesischfivelbd_2 \ge \sigma_1 \eta^{-1/p}$. On the other hand, we also define the event $E_{\Circled{2}, t}$ as follows
        \begin{thesischfivefitdisplay} E_{\Circled{2}, t} \eqdef \thesischfiveens{\abs{\sum_{j = 1}^t Y_j^t} \le \widetilde{b}_2 \quad \text{or} \quad \sum_{j = 1}^t \widetilde{s}_j^2 > \widetilde{G}_2 \log(2/\beta')}, \end{thesischfivefitdisplay}
    where $\widetilde{c}_2$, $\widetilde{G}_2$ and $\widetilde{b}_2$ are defined above, and we can bound $\widetilde{b}_2$ as
        \begin{thesischfivefitalignatnum}{2}
            \widetilde{b}_2 & = \left( \frac{\widetilde{c}_2}{3} + \sqrt{\frac{\widetilde{c}_2^2}{9} + 2 \widetilde{G}_2} \right) \log(2/\beta') \\
            & \le \left( \frac{2 \widetilde{c}_2}{3} + \sqrt{2 \widetilde{G}_2} \right) \log(2/\beta') \\
            & = \left( \frac{16 \thesischfivelbd_2^2}{3} + \sqrt{\frac{576 \thesischfivelbd_2^{4 - p} \sigma_1^p}{\eta}} \right) \log(2/\beta') \\
            & = \thesischfivelbd_2^2 \left( \frac{16}{3} + 24 \sqrt{\frac{1}{\eta} \left(\frac{\sigma_1}{\thesischfivelbd_2} \right)^p} \right) \log(2/\beta') \\
            & \le \frac{88 \thesischfivelbd_2^2}{3} \log(2/\beta'), \label{ch5:c33e8593-ede2-4c3e-b76d-7e7ba9b43443_chapter_5}
        \end{thesischfivefitalignatnum}
    where the last inequality holds given $\thesischfivelbd_2 \ge \sigma_1 \eta^{-1/p}$.

    Next, if $\thesischfivelbd_2 \ge \max\thesischfiveens{2 \norm{\nabla f(x^{j})}, \sigma_1 \eta^{-1/p}}$ for all $j \in [t]$, then we proved in~\eqref{ch5:158c1570-485a-49ca-9ec3-32f47e8d8e8b_chapter_5} and~\eqref{ch5:04635ad8-7080-4c15-a1de-f0ada3a80f78_chapter_5} that $\thesischfiveProba{E_{\Circled{1}, t}} \ge 1 - \beta'$ and $\thesischfiveProba{E_{\Circled{2}, t}} \ge 1 - \beta'$ hence, by the union bound inequality we have
        \begin{thesischfivefitdisplay} \thesischfiveProba{E_{\Circled{1}, t} \cap E_{\Circled{2}, t}} \ge 1 - 2 \beta', \end{thesischfivefitdisplay}
    and on the event $E_{\Circled{1}, t} \cap E_{\Circled{2}, t}$ we obtain the inequality
        \begin{thesischfivefitalignat}{2}
            \norm{\sum_{j = 1}^t (1 - \eta)^{t - j} \theta_j} & \le \norm{\sum_{j = 1}^t (1 - \eta)^{t - j} \theta_j^b} + \norm{\sum_{j = 1}^t (1 - \eta)^{t - j} \theta_j^u} \\
            \thesischfiveoversetlab{\eqref{ch5:21eaf977-44d0-4b4b-8b63-5af59e5cc089_chapter_5}}&{\le} \Circled{1} + \sqrt{2 \cdot \Circled{2} + 2 \cdot \Circled{3}} + \Circled{4} \\
            \thesischfiveoversetlab{\eqref{ch5:35c58053-3ad9-4b62-ae23-2a6f658c08da_chapter_5}+\eqref{ch5:c33e8593-ede2-4c3e-b76d-7e7ba9b43443_chapter_5}}&{\le} \frac{22 \thesischfivelbd_2}{3} \log(2/\beta') \\
            &\qquad + \sqrt{2 \cdot \frac{88 \thesischfivelbd_2^2}{3} \log(2/\beta') + 2 \cdot \frac{18 \thesischfivelbd_2^{2 - p} \sigma_1^p}{\eta}} \\
            &\qquad + \frac{4 \thesischfivelbd_2^{1 - p} \sigma_1^p}{\eta} \\
            & \le \thesischfivelbd_2 \left( \frac{22}{3} + \sqrt{\frac{176}{3} + \frac{36}{\eta} \left( \frac{\sigma_1}{\thesischfivelbd_2} \right)^p} + \frac{4}{\eta} \left( \frac{\sigma_1}{\thesischfivelbd_2} \right)^p \right) \log(2/\beta') \\
            & \le \thesischfivelbd_2 \left( \frac{22}{3} + \sqrt{\frac{176}{3} + 36} + 4 \right) \log(2/\beta') \\
            & \le 22 \thesischfivelbd_2 \log(2/\beta'),
        \end{thesischfivefitalignat}
    where in the fourth inequality we use the fact that $\log(2/\beta') \ge 1$ and in the fifth we use $\thesischfivelbd_2 \ge \sigma_1 \eta^{-1/p}$. This proves the lemma.
\end{proof}

\begin{lemma}\label{ch5:appdx-lem:high-probability-analysis-bound-2_chapter_5}
    Under~Assumption~\ref{ch5:ass:mean-squared-smoothness_chapter_5}, for any $t \in \thesischfiveInt{0}{K - 1}$ and any $\beta'' \in \thesischfiveintof{0}{1/2}$, if the clipping threshold satisfies
        \begin{thesischfivefitdisplay} \thesischfivelbd_1 \ge \max\thesischfiveens{2 \gamma \bar{L}, \gamma \bar{L} \eta^{-1/q}}, \end{thesischfivefitdisplay}
    for all $j \in [t]$ then, with probability at least $1 - 2 \beta''$, we have
        \begin{thesischfivefitdisplay} \norm{\sum_{j = 1}^t (1 - \eta)^{t - j} \omega_j} \le 46 \thesischfivelbd_1 \log(2/\beta''). \end{thesischfivefitdisplay}
\end{lemma}

\begin{proof}
    First, using~\eqref{ch5:9715675c-c075-4d13-a6e5-978530f178fd_chapter_5} and~\eqref{ch5:5bf5f4fb-a644-41c4-a212-9c4930a0075d_chapter_5} we have
        \begin{thesischfivefitdisplay} \norm{\sum_{j = 1}^t (1 - \eta)^{t - j} \omega_j} \le \underbrace{\norm{\sum_{j = 1}^t (1 - \eta)^{t - j} \omega_j^b}}_{\Circled{8}} + \norm{\sum_{j = 1}^t (1 - \eta)^{t - j} \omega_j^u}, \thesischfivenumberthis\label{ch5:9221fb23-db5b-497f-9062-1c71b7f12882_chapter_5} \end{thesischfivefitdisplay}
    and, as before, we need to bound both terms above. For the second term, we use~Lemma~\ref{ch5:appdx-technical-lem:cutkosky-mehta_chapter_5} with exponent $2$ to obtain
        \begin{thesischfivefitalignatnum}{2}
            \MoveEqLeft[1] \norm{\sum_{j = 1}^t (1 - \eta)^{t - j} \omega_j^u} \\
            & \le \abs{\sum_{j = 1}^t W_j^t} + \sqrt{\max_{j \in [t]} \sqnorm{(1 - \eta)^{t - j} \omega_j^u} + \sum_{j = 1}^t \sqnorm{(1 - \eta)^{t - j} \omega_j^u}} \\
            & \le \abs{\sum_{j = 1}^t W_j^t} + \sqrt{2 \sum_{j = 1}^t \sqnorm{(1 - \eta)^{t - j} \omega_j^u}} \\
            & = \underbrace{\abs{\sum_{j = 1}^t W_j^t}}_{\Circled{5}} + \sqrt{2 \Bigl( \underbrace{\sum_{j = 1}^t Z_j^t}_{\Circled{6}} + \underbrace{\sum_{j = 1}^t \thesischfiveExpSub{\xi^j}{\sqnorm{(1 - \eta)^{t - j} \omega_j^u}}}_{\Circled{7}} \Bigr)}, \label{ch5:21eaf977-44d0-4b4b-8b63-5af59e5cc089-2_chapter_5}
        \end{thesischfivefitalignatnum}
    where in the first inequality we use Lemma~\ref{ch5:appdx-technical-lem:cutkosky-mehta_chapter_5} with the sequence defined, using the partial sums $S_{j-1}^t \eqdef \sum_{i = 1}^{j - 1} (1 - \eta)^{t - i} \omega_i^u$, by
        \begin{thesischfivefitdisplay}
            W_j^t \eqdef
            \begin{cases}
                \thesischfivesign\left( \sum_{i = 1}^{j - 1} W_i^t \right) \thesischfiveps{S_{j-1}^t}{(1 - \eta)^{t - j} \omega_j^u} / \norm{S_{j-1}^t}, & \text{if } S_{j-1}^t \neq 0; \\
                0, & \text{if } S_{j-1}^t = 0.
            \end{cases}
        \end{thesischfivefitdisplay}
    while in the equality, for any $j \in [t]$ we let
        \begin{thesischfivefitdisplay} Z_j^t \eqdef \sqnorm{(1 - \eta)^{t - j} \omega_j^u} - \thesischfiveExpSub{\xi^j}{\sqnorm{(1 - \eta)^{t - j} \omega_j^u}}. \end{thesischfivefitdisplay}

    We now bound all terms \Circled{5}, \Circled{6}, \Circled{7} and \Circled{8}.

    \paragraph{Bounding \Circled{5}.} The sequence $W_1^t, \ldots, W_t^t$ is a martingale difference sequence since, by definition of $\omega_j^u$ and $W_j^t$, for all $j \in [t]$ we have $\thesischfiveExpCond{W_j^t}{W_{j - 1}^t, \ldots, W_1^t} = 0$ since $\thesischfiveExpSub{\xi^j}{\omega_j^u} = 0$. Moreover, by~Lemma~\ref{ch5:appdx-technical-lem:cutkosky-mehta_chapter_5} we also have
        \begin{thesischfivefitdisplay} \abs{W_j^t} \le \norm{(1 - \eta)^{t - j} \omega_j^u} \le \norm{\omega_j^u} \le 2 \thesischfivelbd_1, \thesischfivenumberthis\label{ch5:77cc40ad-0220-4787-ab09-7c07a2d6dfa7-2_chapter_5} \end{thesischfivefitdisplay}
    where the last inequality follows from Lemma~\ref{ch5:appdx-technical-lem:sadiev-generalization_chapter_5}, and let $c_1 \eqdef 2 \thesischfivelbd_1$ be the upper bound on the random variables $W_1^t, \ldots, W_t^t$. Additionally, if we denote by $s_j^2 \eqdef \thesischfiveExpCond{\left( W_j^t \right)^2}{W_{j - 1}^t, \ldots, W_1^t}$ then by~Lemma~\ref{ch5:appdx-technical-lem:freedman-inequality_chapter_5} for any real number $b_1 > 0$ and any $G_1 > 0$, as $0 < \beta'' \le 1/2$ then $\log(2/\beta'') \ge 1$ and we have
        \begin{thesischfivefitalignatnum}{2}
            \MoveEqLeft[1] \thesischfiveProba{\abs{\sum_{j = 1}^t W_j^t} > b_1 \quad \text{and} \quad \sum_{j = 1}^t s_j^2 \le G_1 \log(2/\beta'')} \\
            & \le 2 \exp\left( -\frac{b_1^2}{2 G_1 \log(2/\beta'') + 2 b_1 c_1/3} \right) = \beta'', \label{ch5:158c1570-485a-49ca-9ec3-32f47e8d8e8b-2_chapter_5}
        \end{thesischfivefitalignatnum}
    where the last equality holds provided $b_1 \eqdef \bigl( c_1/3 + \sqrt{c_1^2/9 + 2 G_1} \bigr) \log(2/\beta'') > 0$. We now need to define the constant $G_1$. To do so, we need to bound the sum $s_1^2 + \cdots + s_t^2$, this gives
        \begin{thesischfivefitalignatnum}{2}
            \sum_{j = 1}^t s_j^2 & = \sum_{j = 1}^t \thesischfiveExpCond{\left( W_j^t \right)^2}{W_{j - 1}^t, \ldots, W_1^t} \\
            \thesischfiveoversetlab{\eqref{ch5:77cc40ad-0220-4787-ab09-7c07a2d6dfa7-2_chapter_5}}&{\le} \sum_{j = 1}^t \thesischfiveExpSub{\xi^j}{\sqnorm{(1 - \eta)^{t - j} \omega_j^u}} \\
            & = \sum_{j = 1}^t ( 1 - \eta)^{2 (t - j)} \, \thesischfiveExpSub{\xi^j}{\sqnorm{\omega_j^u}} \\
            & \le 72 \thesischfivelbd_1^{2 - q} \gamma^q \bar{L}^q \sum_{j = 1}^t (1 - \eta)^{2 (t - j)} \\
            & \le \frac{72 \thesischfivelbd_1^{2 - q} \gamma^q \bar{L}^q}{1 - (1 - \eta)^2} \\
            & \le \frac{72 \thesischfivelbd_1^{2 - q} \gamma^q \bar{L}^q}{\eta}, \label{ch5:74b4e2bb-4754-4dd3-ae90-37ef40362f93-2_chapter_5}
        \end{thesischfivefitalignatnum}
    where in the application of~Lemma~\ref{ch5:appdx-technical-lem:sadiev-generalization_chapter_5} we use the bound we proved earlier in~Lemma~\ref{ch5:appdx-lem:bounding-error-term-nsgd-mvr_chapter_5} (see more precisely at~\eqref{ch5:f79e21c1-04c0-497e-81e9-c10afd7f04f9_chapter_5}), that is,
        \begin{thesischfivefitdisplay} \thesischfiveExpSub{\xi^j}{\norm{\left[ \nabla f(x^{j}) - \nabla f(x^{j-1}) \right] - \left[ \nabla f(x^{j};\xi^j) - \nabla f(x^{j-1};\xi^j) \right]}^q} \le 2^q \gamma^q \bar{L}^q. \end{thesischfivefitdisplay}
    {\color{navy}Moreover, applying~Lemma~\ref{ch5:appdx-technical-lem:sadiev-generalization_chapter_5} requires}
        \begin{thesischfivefitdisplay}
            {\color{navy}\thesischfivelbd_1
            \ge 2\max_{j\in[t]}
            \norm{\nabla f(x^j)-\nabla f(x^{j-1})}.}
        \end{thesischfivefitdisplay}
    {\color{navy}Since}
        \begin{thesischfivefitdisplay} \norm{\nabla f(x^{j}) - \nabla f(x^{j-1})} \le \bar{L} \norm{x^{j} - x^{j-1}} \le \gamma \bar{L}, \end{thesischfivefitdisplay}
    by Assumption~\ref{ch5:ass:mean-squared-smoothness_chapter_5} (see~\eqref{ch5:94828cd2-cf02-492b-a6ec-40b3012f66c1_chapter_5}), it is enough to have $\thesischfivelbd_1 \ge 2 \gamma \bar{L}$. Setting $G_1 \eqdef 72 \thesischfivelbd_1^{2 - q} \gamma^q \bar{L}^q/\eta > 0$ gives the desired bound~\eqref{ch5:158c1570-485a-49ca-9ec3-32f47e8d8e8b-2_chapter_5}.

    \paragraph{Bounding \Circled{6}.} As in the previous paragraph, the sequence $Z_1^t, \ldots, Z_t^t$ is also a martingale difference sequence, since $\thesischfiveExpSub{\xi^j}{Z_j^t} = 0$ for all $j \in [t]$ by definition of $Z_j^t$. Moreover, according to~Lemma~\ref{ch5:appdx-technical-lem:sadiev-generalization_chapter_5}, we also have, as we did in~\eqref{ch5:77cc40ad-0220-4787-ab09-7c07a2d6dfa7-2_chapter_5}, for any $j \in [t]$
        \begin{thesischfivefitdisplay} \abs{Z_j^t} \le \sqnorm{(1 - \eta)^{t - j} \omega_j^u} + \thesischfiveExpSub{\xi^j}{\sqnorm{(1 - \eta)^{t - j} \omega_j^u}} \le 4 \thesischfivelbd_1^2 + 4 \thesischfivelbd_1^2 = 8 \thesischfivelbd_1^2. \thesischfivenumberthis\label{ch5:4923f3a4-39e5-43f3-b260-7deeb46a8c67-2_chapter_5} \end{thesischfivefitdisplay}
    Hence we define $\widetilde{c}_1 \eqdef 8 \thesischfivelbd_1^2$. Now, denoting the conditional variance of $Z_j^t$ as $\widetilde{s}_j^2 \eqdef \thesischfiveExpCond{\left( Z_j^t \right)^2}{Z_{j - 1}^t, \ldots, Z_1^t} = \thesischfiveExpSub{\xi^j}{\left( Z_j^t \right)^2}$ we have the bound
        \begin{thesischfivefitalignatnum}{2}
            \widetilde{s}_j^2 \thesischfiveoversetlab{\eqref{ch5:4923f3a4-39e5-43f3-b260-7deeb46a8c67-2_chapter_5}}&{\le} 8 \thesischfivelbd_1^2 \, \thesischfiveExpSub{\xi^j}{\abs{Z_j^t}} \le 16 \thesischfivelbd_1^2 \, \thesischfiveExpSub{\xi^j}{\sqnorm{(1 - \eta)^{t - j} \omega_j^u}}, \label{ch5:a26983bb-60b7-49ab-ac0e-a1d026f1c51e-2_chapter_5}
        \end{thesischfivefitalignatnum}
    for all $j \in [t]$. Hence, 
        \begin{thesischfivefitalignat}{2}
            \sum_{j = 1}^t \widetilde{s}_j^2 \thesischfiveoversetlab{\eqref{ch5:a26983bb-60b7-49ab-ac0e-a1d026f1c51e-2_chapter_5}}&{\le} 16 \thesischfivelbd_1^2 \sum_{j = 1}^t \thesischfiveExpSub{\xi^j}{\sqnorm{(1 - \eta)^{t - j} \omega_j^u}} \\
            \thesischfiveoversetlab{\eqref{ch5:74b4e2bb-4754-4dd3-ae90-37ef40362f93-2_chapter_5}}&{\le} 16 \thesischfivelbd_1^2 \cdot \frac{72 \thesischfivelbd_1^{2 - q} \gamma^q \bar{L}^q}{\eta} = \frac{1152 \thesischfivelbd_1^{4 - q} \gamma^q \bar{L}^q}{\eta},
        \end{thesischfivefitalignat}
    (where we require $\thesischfivelbd_1 \ge 2 \gamma \bar{L}$). Next, if we define $\widetilde{G}_1 \eqdef 1152 \thesischfivelbd_1^{4 - q} \gamma^q \bar{L}^q/\eta$ then, applying~Lemma~\ref{ch5:appdx-technical-lem:freedman-inequality_chapter_5} we obtain for any real number $\widetilde{b}_1$
        \begin{thesischfivefitalignatnum}{2}
            \MoveEqLeft[1] \thesischfiveProba{\abs{\sum_{j = 1}^t Z_j^t} > \widetilde{b}_1 \quad \text{and} \quad \sum_{j = 1}^t \widetilde{s}_j^2 \le \widetilde{G}_1 \log(2/\beta'')} \\
            & \le 2 \exp\left( -\frac{\widetilde{b}_1^2}{2 \widetilde{G}_1 \log(2/\beta'') + 2 \widetilde{b}_1 \widetilde{c}_1/3} \right) = \beta'', \label{ch5:04635ad8-7080-4c15-a1de-f0ada3a80f78-2_chapter_5}
        \end{thesischfivefitalignatnum}
    where the last identity holds if we set $\widetilde{b}_1 \eqdef \bigl( \widetilde{c}_1/3 + \sqrt{\widetilde{c}_1^2/9 + 2 \widetilde{G}_1} \bigr) \log(2/\beta'')$. This establishes the desired bound.

    \paragraph{Bounding \Circled{7}.} {\color{navy}As in the two preceding paragraphs, we have}
        \begin{thesischfivefitdisplay} \Circled{7} \eqdef \sum_{j = 1}^t \thesischfiveExpSub{\xi^j}{\sqnorm{(1 - \eta)^{t - j} \omega_j^u}} \thesischfiveoversetlab{\eqref{ch5:74b4e2bb-4754-4dd3-ae90-37ef40362f93-2_chapter_5}}{\le} \frac{72 \thesischfivelbd_1^{2 - q} \gamma^q \bar{L}^q}{\eta}, \end{thesischfivefitdisplay}
    as desired.

    \paragraph{Bounding \Circled{8}.} If we assume that $\thesischfivelbd_1 \ge 2 \gamma \bar{L}$ then, with probability one we have
        \begin{thesischfivefitalignat}{2}
            \Circled{8} \eqdef* \norm{\sum_{j = 1}^t (1 - \eta)^{t - j} \omega_j^b} \\
            & \le \sum_{j = 1}^t (1 - \eta)^{t - j} \norm{\omega_j^b} \\
            & \le 16 \thesischfivelbd_1^{1 - q} \gamma^q \bar{L}^q \sum_{j = 1}^t (1 - \eta)^{t - j} \\
            & \le \frac{16 \thesischfivelbd_1^{1 - q} \gamma^q \bar{L}^q}{\eta},
        \end{thesischfivefitalignat}
    where in the second inequality we use~Lemma~\ref{ch5:appdx-technical-lem:sadiev-generalization_chapter_5}, more precisely, for any $j \in [t]$,
        \begin{thesischfivefitdisplay} \norm{\omega_j^b} \le 4^q \thesischfivelbd_1^{1 - q} \gamma^q \bar{L}^q \le 16 \thesischfivelbd_1^{1 - q} \gamma^q \bar{L}^q, \end{thesischfivefitdisplay}
    since, by~Lemma~\ref{ch5:appdx-lem:bounding-error-term-nsgd-mvr_chapter_5} (more precisely by~\eqref{ch5:f79e21c1-04c0-497e-81e9-c10afd7f04f9_chapter_5}),
        \begin{thesischfivefitdisplay} \thesischfiveE{\norm{\left[ \nabla f(x^{j}) - \nabla f(x^{j-1}) \right] - \left[ \nabla f(x^{j};\xi^{j}) - \nabla f(x^{j-1};\xi^{j}) \right]}^q} \le 2^q \gamma^q \bar{L}^q. \end{thesischfivefitdisplay}

    \paragraph{Summing up all bounds \Circled{5}, \Circled{6}, \Circled{7} and \Circled{8}.} To sum up, we introduce the event $E_{\Circled{5}, t}$ as follows
        \begin{thesischfivefitdisplay} E_{\Circled{5}, t} \eqdef \thesischfiveens{\abs{\sum_{j = 1}^t W_j^t} \le b_1 \quad \text{or} \quad \sum_{j = 1}^t s_j^2 > G_1 \log(2/\beta'')}, \thesischfivenumberthis\label{ch5:1ded530f-71dc-46eb-9fdb-5308c908cb80_chapter_5} \end{thesischfivefitdisplay}
    where $c_1$, $G_1$ and $b_1$ are defined above, and we can bound $b_1$ as
        \begin{thesischfivefitalignatnum}{2}
            b_1 & = \left( \frac{c_1}{3} + \sqrt{\frac{c_1^2}{9} + 2 G_1} \right) \log(2/\beta'') \\
            & \le \left( \frac{2 c_1}{3} + \sqrt{2 G_1} \right) \log(2/\beta'') \\
            & = \left( \frac{4 \thesischfivelbd_1}{3} + \sqrt{\frac{144 \thesischfivelbd_1^{2 - q} \gamma^q \bar{L}^q}{\eta}} \right) \log(2/\beta'') \\
            & = \thesischfivelbd_1 \left( \frac{4}{3} + 12 \sqrt{\frac{1}{\eta} \left(\frac{\gamma \bar{L}}{\thesischfivelbd_1} \right)^q} \right) \log(2/\beta'') \\
            & \le \frac{40 \thesischfivelbd_1}{3} \log(2/\beta''), \label{ch5:35c58053-3ad9-4b62-ae23-2a6f658c08da-2_chapter_5}
        \end{thesischfivefitalignatnum}
    where the last inequality holds provided $\thesischfivelbd_1 \ge \gamma \bar{L} \eta^{-1/q}$. On the other hand, we also define the event $E_{\Circled{6}, t}$ as follows
        \begin{thesischfivefitdisplay} E_{\Circled{6}, t} \eqdef \thesischfiveens{\abs{\sum_{j = 1}^t Z_j^t} \le \widetilde{b}_1 \quad \text{or} \quad \sum_{j = 1}^t \widetilde{s}_j^2 > \widetilde{G}_1 \log(2/\beta'')}, \thesischfivenumberthis\label{ch5:92e2f28f-a246-45f3-bdea-ba28563d3825_chapter_5} \end{thesischfivefitdisplay}
    where $\widetilde{c}_1$, $\widetilde{G}_1$ and $\widetilde{b}_1$ are defined above, and we can bound $\widetilde{b}_1$ as
        \begin{thesischfivefitalignatnum}{2}
            \widetilde{b}_1 & = \left( \frac{\widetilde{c}_1}{3} + \sqrt{\frac{\widetilde{c}_1^2}{9} + 2 \widetilde{G}_1} \right) \log(2/\beta'') \\
            & \le \left( \frac{2 \widetilde{c}_1}{3} + \sqrt{2 \widetilde{G}_1} \right) \log(2/\beta'') \\
            & = \left( \frac{16 \thesischfivelbd_1^2}{3} + \sqrt{\frac{2304 \thesischfivelbd_1^{4 - q} \gamma^q \bar{L}^q}{\eta}} \right) \log(2/\beta'') \\
            & = \thesischfivelbd_1^2 \left( \frac{16}{3} + 48 \sqrt{\frac{1}{\eta} \left(\frac{\gamma \bar{L}}{\thesischfivelbd_1} \right)^q} \right) \log(2/\beta'') \\
            & \le \frac{160 \thesischfivelbd_1^2}{3} \log(2/\beta''), \label{ch5:c33e8593-ede2-4c3e-b76d-7e7ba9b43443-2_chapter_5}
        \end{thesischfivefitalignatnum}
    where the last inequality holds given $\thesischfivelbd_1 \ge \gamma \bar{L} \eta^{-1/q}$.

    Next, if $\thesischfivelbd_1 \ge \max\thesischfiveens{2 \gamma \bar{L}, \gamma \bar{L} \eta^{-1/q}}$ for all $j \in [t]$, then we proved in~\eqref{ch5:158c1570-485a-49ca-9ec3-32f47e8d8e8b-2_chapter_5} and~\eqref{ch5:04635ad8-7080-4c15-a1de-f0ada3a80f78-2_chapter_5} that $\thesischfiveProba{E_{\Circled{5}, t}} \ge 1 - \beta''$ and $\thesischfiveProba{E_{\Circled{6}, t}} \ge 1 - \beta''$ hence, by the union bound inequality we have
        \begin{thesischfivefitdisplay} \thesischfiveProba{E_{\Circled{5}, t} \cap E_{\Circled{6}, t}} \ge 1 - 2 \beta'', \end{thesischfivefitdisplay}
    and on the event $E_{\Circled{5}, t} \cap E_{\Circled{6}, t}$ we obtain the inequality
        \begin{thesischfivefitalignat}{2}
            \norm{\sum_{j = 1}^t (1 - \eta)^{t - j} \omega_j} & \le \norm{\sum_{j = 1}^t (1 - \eta)^{t - j} \omega_j^b} + \norm{\sum_{j = 1}^t (1 - \eta)^{t - j} \omega_j^u} \\
            \thesischfiveoversetlab{\eqref{ch5:21eaf977-44d0-4b4b-8b63-5af59e5cc089-2_chapter_5}}&{\le} \Circled{5} + \sqrt{2 \cdot \Circled{6} + 2 \cdot \Circled{7}} + \Circled{8} \\
            \thesischfiveoversetlab{\eqref{ch5:35c58053-3ad9-4b62-ae23-2a6f658c08da-2_chapter_5}+\eqref{ch5:c33e8593-ede2-4c3e-b76d-7e7ba9b43443-2_chapter_5}}&{\le} \frac{40 \thesischfivelbd_1}{3} \log(2/\beta'') \\
            &\qquad + \sqrt{2 \cdot \frac{160 \thesischfivelbd_1^2}{3} \log(2/\beta'') + 2 \cdot \frac{72 \thesischfivelbd_1^{2 - q} \gamma^q \bar{L}^q}{\eta}} \\
            &\qquad + \frac{16 \thesischfivelbd_1^{1 - q} \gamma^q \bar{L}^q}{\eta} \\
            & \le \thesischfivelbd_1 \left( \frac{40}{3} + \sqrt{\frac{320}{3} + \frac{144}{\eta} \left( \frac{\gamma \bar{L}}{\thesischfivelbd_1} \right)^q} + \frac{16}{\eta} \left( \frac{\gamma \bar{L}}{\thesischfivelbd_1} \right)^q \right) \log(2/\beta'') \\
            & \le \thesischfivelbd_1 \left( \frac{40}{3} + \sqrt{\frac{320}{3} + 144} + 16 \right) \log(2/\beta'') \\
            & \le 46 \thesischfivelbd_1 \log(2/\beta''),
        \end{thesischfivefitalignat}
    where in the fourth inequality we use the fact that $\log(2/\beta'') \ge 1$ and in the fifth we use $\thesischfivelbd_1 \ge \gamma \bar{L} \eta^{-1/q}$. This proves the lemma.
\end{proof}

\begin{lemma}\label{ch5:appdx-lem:high-probability-analysis-bound-3_chapter_5}
    Under~Assumption~\ref{ch5:ass:L-lipschitz-gradients_chapter_5} and Assumption~\ref{ch5:ass:mean-squared-smoothness-2_chapter_5}, for any $\beta'' \in \thesischfiveintof{0}{1/2}$ and any $t \in \thesischfiveInt{0}{K - 1}$, if the clipping threshold satisfies
        \begin{thesischfivefitdisplay} \thesischfivelbd_1 \ge \max\thesischfiveens{2 \gamma L_1, \gamma {\color{navy}\delta} \eta^{-1/q}}, \end{thesischfivefitdisplay}
    for all $j \in [t]$ then, with probability at least $1 - 2 \beta''$, we have
        \begin{thesischfivefitdisplay} \norm{\sum_{j = 1}^t (1 - \eta)^{t - j} \omega_j} \le 22 \thesischfivelbd_1 \log(2/\beta''). \end{thesischfivefitdisplay}
\end{lemma}

\begin{proof}
    First, note that the computations~\eqref{ch5:9221fb23-db5b-497f-9062-1c71b7f12882_chapter_5} and~\eqref{ch5:21eaf977-44d0-4b4b-8b63-5af59e5cc089-2_chapter_5} are identical and we now need to bound all terms \Circled{5}, \Circled{6}, \Circled{7} and \Circled{8} where
        \begin{align*}
            & \Circled{5} \eqdef \abs{\sum_{j = 1}^t W_j^t}, \quad \Circled{6} \eqdef \sum_{j = 1}^t Z_j^t, \\
            & \Circled{7} \eqdef \sum_{j = 1}^t \thesischfiveExpSub{\xi^j}{\sqnorm{(1 - \eta)^{t - j} \omega_j^u}} \quad \text{and} \quad \Circled{8} \eqdef \norm{\sum_{j = 1}^t (1 - \eta)^{t - j} \omega_j^b}.
        \end{align*}

    \paragraph{Bounding \Circled{5}.} As in the previous lemma, we choose $c_1 \eqdef 2 \thesischfivelbd_1$ and $b_1 \eqdef \bigl( c_1/3 + \sqrt{c_1^2/9 + 2 G_1} \bigr) \log(2/\beta'') > 0$. Then, to select $G_1 > 0$ we need to bound the sum $s_1^2 + \cdots + s_t^2$, i.e.,
        \begin{thesischfivefitalignatnum}{2}
            \sum_{j = 1}^t s_j^2 \thesischfiveoversetlab{\eqref{ch5:74b4e2bb-4754-4dd3-ae90-37ef40362f93-2_chapter_5}}&{\le} \sum_{j = 1}^t ( 1 - \eta)^{2 (t - j)} \, \thesischfiveExpSub{\xi^j}{\sqnorm{\omega_j^u}} \\
            & \le 18 \thesischfivelbd_1^{2 - q} \gamma^q {\color{navy}\delta}^q \sum_{j = 1}^t (1 - \eta)^{2 (t - j)} \\
            & \le \frac{18 \thesischfivelbd_1^{2 - q} \gamma^q {\color{navy}\delta}^q}{\eta}, \label{ch5:74b4e2bb-4754-4dd3-ae90-37ef40362f93-3_chapter_5}
        \end{thesischfivefitalignatnum}
    where ${\color{navy}\delta} \ge 0$ is the parameter in~Assumption~\ref{ch5:ass:mean-squared-smoothness-2_chapter_5}. In the application of~Lemma~\ref{ch5:appdx-technical-lem:sadiev-generalization_chapter_5} we use~Assumption~\ref{ch5:ass:mean-squared-smoothness-2_chapter_5}
        \begin{thesischfivefitalignat}{2}
            \MoveEqLeft[1] \thesischfiveExpSub{\xi^j}{\norm{\left[ \nabla f(x^{j}) - \nabla f(x^{j-1}) \right] - \left[ \nabla f(x^{j};\xi^j) - \nabla f(x^{j-1};\xi^j) \right]}^q} \\
            & \le {\color{navy}\delta}^q \norm{x^{j} - x^{j-1}}^q \le \gamma^q {\color{navy}\delta}^q.
        \end{thesischfivefitalignat}
    {\color{navy}Moreover, applying~Lemma~\ref{ch5:appdx-technical-lem:sadiev-generalization_chapter_5} requires}
        \begin{thesischfivefitdisplay}
            {\color{navy}\thesischfivelbd_1
            \ge 2\max_{j\in[t]}
            \norm{\nabla f(x^j)-\nabla f(x^{j-1})}.}
        \end{thesischfivefitdisplay}
    {\color{navy}Since}
        \begin{thesischfivefitdisplay} \norm{\nabla f(x^{j}) - \nabla f(x^{j-1})} \le L_1 \norm{x^{j} - x^{j-1}} \le \gamma L_1, \end{thesischfivefitdisplay}
    by Assumption~\ref{ch5:ass:L-lipschitz-gradients_chapter_5}, it is enough to have $\thesischfivelbd_1 \ge 2 \gamma L_1$. Setting $G_1 \eqdef 18 \thesischfivelbd_1^{2 - q} \gamma^q {\color{navy}\delta}^q/\eta > 0$ gives the desired bound.

    \paragraph{Bounding \Circled{6}.} Similarly to the previous lemma, we set $\widetilde{c}_1 \eqdef 8 \thesischfivelbd_1^2$ and $\widetilde{b}_1 \eqdef \bigl( \widetilde{c}_1/3 + \sqrt{\widetilde{c}_1^2/9 + 2 \widetilde{G}_1} \bigr) \log(2/\beta'')$. For the choice $\widetilde{G}_1$, we have the bound 
        \begin{thesischfivefitalignat}{2}
            \sum_{j = 1}^t \widetilde{s}_j^2 \thesischfiveoversetlab{\eqref{ch5:a26983bb-60b7-49ab-ac0e-a1d026f1c51e-2_chapter_5}}&{\le} 16 \thesischfivelbd_1^2 \sum_{j = 1}^t \thesischfiveExpSub{\xi^j}{\sqnorm{(1 - \eta)^{t - j} \omega_j^u}} \\
            \thesischfiveoversetlab{\eqref{ch5:74b4e2bb-4754-4dd3-ae90-37ef40362f93-3_chapter_5}}&{\le} 16 \thesischfivelbd_1^2 \cdot \frac{18 \thesischfivelbd_1^{2 - q} \gamma^q {\color{navy}\delta}^q}{\eta} = \frac{288 \thesischfivelbd_1^{4 - q} \gamma^q {\color{navy}\delta}^q}{\eta},
        \end{thesischfivefitalignat}
    (where we require $\thesischfivelbd_1 \ge 2 \gamma L_1$). Hence, if we let $\widetilde{G}_1 \eqdef 288 \thesischfivelbd_1^{4 - q} \gamma^q {\color{navy}\delta}^q/\eta$ this establishes the desired bound.

    \paragraph{Bounding \Circled{7}.} {\color{navy}As in the two preceding paragraphs, we have}
        \begin{thesischfivefitdisplay} \Circled{7} \eqdef \sum_{j = 1}^t \thesischfiveExpSub{\xi^j}{\sqnorm{(1 - \eta)^{t - j} \omega_j^u}} \thesischfiveoversetlab{\eqref{ch5:74b4e2bb-4754-4dd3-ae90-37ef40362f93-3_chapter_5}}{\le} \frac{18 \thesischfivelbd_1^{2 - q} \gamma^q {\color{navy}\delta}^q}{\eta}, \end{thesischfivefitdisplay}
    as desired.

    \paragraph{Bounding \Circled{8}.} For the last bound, if we assume that $\thesischfivelbd_1 \ge 2 \gamma L_1$ then, with probability one we have
        \begin{thesischfivefitalignat}{2}
            \Circled{8} \eqdef* \norm{\sum_{j = 1}^t (1 - \eta)^{t - j} \omega_j^b} \\
            & \le \sum_{j = 1}^t (1 - \eta)^{t - j} \norm{\omega_j^b} \\
            & \le 4 \thesischfivelbd_1^{1 - q} \gamma^q {\color{navy}\delta}^q \sum_{j = 1}^t (1 - \eta)^{t - j} \\
            & \le \frac{4 \thesischfivelbd_1^{1 - q} \gamma^q {\color{navy}\delta}^q}{\eta},
        \end{thesischfivefitalignat}
    where in the second inequality we use~Lemma~\ref{ch5:appdx-technical-lem:sadiev-generalization_chapter_5}, more precisely, for any $j \in [t]$,
        \begin{thesischfivefitdisplay} \norm{\omega_j^b} \le 2^q \thesischfivelbd_1^{1 - q} \gamma^q {\color{navy}\delta}^q \le 4 \thesischfivelbd_1^{1 - q} \gamma^q {\color{navy}\delta}^q, \end{thesischfivefitdisplay}
    since, by~Assumption~\ref{ch5:ass:mean-squared-smoothness-2_chapter_5},
        \begin{thesischfivefitdisplay} \thesischfiveE{\norm{\left[ \nabla f(x^{j}) - \nabla f(x^{j-1}) \right] - \left[ \nabla f(x^{j};\xi^{j}) - \nabla f(x^{j-1};\xi^{j}) \right]}^q} \le \gamma^q {\color{navy}\delta}^q. \end{thesischfivefitdisplay}

    \paragraph{Summing up all bounds \Circled{5}, \Circled{6}, \Circled{7} and \Circled{8}.} We introduce the same events $E_{\Circled{5}, t}$ and $E_{\Circled{6}, t}$ as in~\eqref{ch5:1ded530f-71dc-46eb-9fdb-5308c908cb80_chapter_5} and~\eqref{ch5:92e2f28f-a246-45f3-bdea-ba28563d3825_chapter_5} respectively. Then, by our choice of $c_1$, $b_1$ and $G_1$ we have the bound
        \begin{thesischfivefitalignatnum}{2}
            b_1 & = \left( \frac{c_1}{3} + \sqrt{\frac{c_1^2}{9} + 2 G_1} \right) \log(2/\beta'') \\
            & \le \left( \frac{2 c_1}{3} + \sqrt{2 G_1} \right) \log(2/\beta'') \\
            & = \left( \frac{4 \thesischfivelbd_1}{3} + \sqrt{\frac{36 \thesischfivelbd_1^{2 - q} \gamma^q {\color{navy}\delta}^q}{\eta}} \right) \log(2/\beta'') \\
            & = \thesischfivelbd_1 \left( \frac{4}{3} + 6 \sqrt{\frac{1}{\eta} \left(\frac{\gamma {\color{navy}\delta}}{\thesischfivelbd_1} \right)^q} \right) \log(2/\beta'') \\
            & \le \frac{22 \thesischfivelbd_1}{3} \log(2/\beta''), \label{ch5:35c58053-3ad9-4b62-ae23-2a6f658c08da-3_chapter_5}
        \end{thesischfivefitalignatnum}
    where the last inequality holds provided $\thesischfivelbd_1 \ge \gamma {\color{navy}\delta} \eta^{-1/q}$. Moreover, by our choice of $\widetilde{c}_1$, $\widetilde{b}_1$ and $\widetilde{G}_1$ we also have the bound
        \begin{thesischfivefitalignatnum}{2}
            \widetilde{b}_1 & = \left( \frac{\widetilde{c}_1}{3} + \sqrt{\frac{\widetilde{c}_1^2}{9} + 2 \widetilde{G}_1} \right) \log(2/\beta'') \\
            & \le \left( \frac{2 \widetilde{c}_1}{3} + \sqrt{2 \widetilde{G}_1} \right) \log(2/\beta'') \\
            & = \left( \frac{16 \thesischfivelbd_1^2}{3} + \sqrt{\frac{576 \thesischfivelbd_1^{4 - q} \gamma^q {\color{navy}\delta}^q}{\eta}} \right) \log(2/\beta'') \\
            & = \thesischfivelbd_1^2 \left( \frac{16}{3} + 24 \sqrt{\frac{1}{\eta} \left(\frac{\gamma {\color{navy}\delta}}{\thesischfivelbd_1} \right)^q} \right) \log(2/\beta'') \\
            & \le \frac{88 \thesischfivelbd_1^2}{3} \log(2/\beta''), \label{ch5:c33e8593-ede2-4c3e-b76d-7e7ba9b43443-3_chapter_5}
        \end{thesischfivefitalignatnum}
    where the last inequality holds given $\thesischfivelbd_1 \ge \gamma {\color{navy}\delta} \eta^{-1/q}$.

    Next, given $\thesischfivelbd_1 \ge \max\thesischfiveens{2 \gamma L_1, \gamma {\color{navy}\delta} \eta^{-1/q}}$, as in the previous lemma, we have $\thesischfiveProba{E_{\Circled{5}, t}} \ge 1 - \beta''$ and $\thesischfiveProba{E_{\Circled{6}, t}} \ge 1 - \beta''$ hence, by the union bound inequality we have
        \begin{thesischfivefitdisplay} \thesischfiveProba{E_{\Circled{5}, t} \cap E_{\Circled{6}, t}} \ge 1 - 2 \beta'', \end{thesischfivefitdisplay}
    and on the event $E_{\Circled{5}, t} \cap E_{\Circled{6}, t}$ we obtain the inequality
        \begin{thesischfivefitalignat}{2}
            \norm{\sum_{j = 1}^t (1 - \eta)^{t - j} \omega_j} & \le \norm{\sum_{j = 1}^t (1 - \eta)^{t - j} \omega_j^b} + \norm{\sum_{j = 1}^t (1 - \eta)^{t - j} \omega_j^u} \\
            \thesischfiveoversetlab{\eqref{ch5:21eaf977-44d0-4b4b-8b63-5af59e5cc089-2_chapter_5}}&{\le} \Circled{5} + \sqrt{2 \cdot \Circled{6} + 2 \cdot \Circled{7}} + \Circled{8} \\
            \thesischfiveoversetlab{\eqref{ch5:35c58053-3ad9-4b62-ae23-2a6f658c08da-3_chapter_5}+\eqref{ch5:c33e8593-ede2-4c3e-b76d-7e7ba9b43443-3_chapter_5}}&{\le} \frac{22 \thesischfivelbd_1}{3} \log(2/\beta'') \\
            &\qquad + \sqrt{2 \cdot \frac{88 \thesischfivelbd_1^2}{3} \log(2/\beta'') + 2 \cdot \frac{18 \thesischfivelbd_1^{2 - q} \gamma^q {\color{navy}\delta}^q}{\eta}} \\
            &\qquad + \frac{4 \thesischfivelbd_1^{1 - q} \gamma^q {\color{navy}\delta}^q}{\eta} \\
            & \le \thesischfivelbd_1 \left( \frac{22}{3} + \sqrt{\frac{176}{3} + \frac{36}{\eta} \left( \frac{\gamma {\color{navy}\delta}}{\thesischfivelbd_1} \right)^q} + \frac{4}{\eta} \left( \frac{\gamma {\color{navy}\delta}}{\thesischfivelbd_1} \right)^q \right) \log(2/\beta'') \\
            & \le \thesischfivelbd_1 \left( \frac{22}{3} + \sqrt{\frac{176}{3} + 36} + 4 \right) \log(2/\beta'') \\
            & \le 22 \thesischfivelbd_1 \log(2/\beta''),
        \end{thesischfivefitalignat}
    where in the fourth inequality we use the fact that $\log(2/\beta'') \ge 1$ and in the fifth we use $\thesischfivelbd_1 \ge \gamma {\color{navy}\delta} \eta^{-1/q}$. This proves the lemma.
\end{proof}

\begin{lemma}\label{ch5:appdx-lem:high-probability-analysis-bound-hess_chapter_5}
    Under~Assumption~\ref{ch5:ass:L-lipschitz-gradients_chapter_5}, Assumption~\ref{ch5:ass:L-lipschitz-hessians_chapter_5}, and Assumption~\ref{ch5:ass:q-bounded-central-moment-hessian_chapter_5}, for any $\beta'' \in \thesischfiveintof{0}{1/2}$ and any $t \in \thesischfiveInt{0}{K - 1}$, if the clipping threshold satisfies
        \begin{thesischfivefitdisplay} \thesischfivelbd_1 \ge \max\thesischfiveens{2 \gamma L_1, \gamma \sigma_2 \eta^{-1/q}}, \end{thesischfivefitdisplay}
    for all $j \in [t]$ then, with probability at least $1 - 2 \beta''$, we have
        \begin{thesischfivefitdisplay} \norm{\sum_{j = 1}^t (1 - \eta)^{t - j} \omega_j} \le \frac{2 \min\thesischfiveens{\gamma L_1, \gamma^2 L_2}}{\eta} + 22 \thesischfivelbd_1 \log(2/\beta''). \end{thesischfivefitdisplay}
\end{lemma}

\begin{proof}
    First, using~\eqref{ch5:9715675c-c075-4d13-a6e5-978530f178fe-hess_chapter_5},~\eqref{ch5:9715675c-c075-4d13-a6e5-978530f178fd-hess_chapter_5} and~\eqref{ch5:5bf5f4fb-a644-41c4-a212-9c4930a0075d-hess_chapter_5} we have
        \begin{thesischfivefitalignatnum}{2}
            \MoveEqLeft[1] \norm{\sum_{j = 1}^t (1 - \eta)^{t - j} \omega_j} \\
            & \le \underbrace{\norm{\sum_{j = 1}^t (1 - \eta)^{t - j} \omega_j^r}}_{\Circled{9}} + \underbrace{\norm{\sum_{j = 1}^t (1 - \eta)^{t - j} \omega_j^b}}_{\Circled{8}} + \norm{\sum_{j = 1}^t (1 - \eta)^{t - j} \omega_j^u}, \label{ch5:9221fb23-db5b-497f-9062-1c71b7f12882-hess_chapter_5}
        \end{thesischfivefitalignatnum}
    As in the previous lemmas, the computation~\eqref{ch5:21eaf977-44d0-4b4b-8b63-5af59e5cc089-2_chapter_5} is identical and we now need to bound all terms \Circled{5}, \Circled{6}, \Circled{7}, \Circled{8} and \Circled{9} where
        \begin{thesischfivefitdisplay} \Circled{5} \eqdef \abs{\sum_{j = 1}^t W_j^t}, \quad \Circled{6} \eqdef \sum_{j = 1}^t Z_j^t \quad \text{and} \quad \Circled{7} \eqdef \sum_{j = 1}^t \thesischfiveExpSub{\hat{\xi}^j}{\sqnorm{(1 - \eta)^{t - j} \omega_j^u}}, \end{thesischfivefitdisplay}
    and $W_1^t, \ldots, W_t^t$ is a martingale difference sequence since, by definition of $\omega_j^u$ and $W_j^t$, for all $j \in [t]$ we have $\thesischfiveExpCond{W_j^t}{W_{j - 1}^t, \ldots, W_1^t} = 0$ since 
        \begin{thesischfivefitdisplay} \thesischfiveExpSub{u^j, \hat{\xi}^j}{\omega_j^u} = \thesischfiveExpSub{u^j}{\thesischfiveExpSubCond{\hat{\xi}^j}{\omega_j^u}{u^j}} = 0. \end{thesischfivefitdisplay}
    The same argument applies to the sequence $Z_1^t, \ldots, Z_t^t$.

    \paragraph{Bounding \Circled{5}.} As already done before, we choose $c_1 \eqdef 2 \thesischfivelbd_1$ and $b_1 \eqdef \bigl( c_1/3 + \sqrt{c_1^2/9 + 2 G_1} \bigr) \log(2/\beta'') > 0$. Then, to select $G_1 > 0$ we need to bound the sum $s_1^2 + \cdots + s_t^2$, i.e.,
        \begin{thesischfivefitalignatnum}{2}
            \sum_{j = 1}^t s_j^2 \thesischfiveoversetlab{\eqref{ch5:74b4e2bb-4754-4dd3-ae90-37ef40362f93-2_chapter_5}}&{\le} \sum_{j = 1}^t ( 1 - \eta)^{2 (t - j)} \, \thesischfiveExpSub{\hat{\xi}^j}{\sqnorm{\omega_j^u}} \\
            & \le 18 \thesischfivelbd_1^{2 - q} \gamma^q \sigma_2^q \sum_{j = 1}^t (1 - \eta)^{2 (t - j)} \\
            & \le \frac{18 \thesischfivelbd_1^{2 - q} \gamma^q \sigma_2^q}{\eta}. \label{ch5:74b4e2bb-4754-4dd3-ae90-37ef40362f93-hess_chapter_5}
        \end{thesischfivefitalignatnum}
    In the application of~Lemma~\ref{ch5:appdx-technical-lem:sadiev-generalization_chapter_5} we use~Assumption~\ref{ch5:ass:q-bounded-central-moment-hessian_chapter_5}
        \begin{thesischfivefitalignat}{2}
            \MoveEqLeft[1] \color{navy}\thesischfiveExpSub{\hat{\xi}^j}{\norm{\nabla^2 f(\hat{x}^j;\hat{\xi}^j)(x^{j} - x^{j-1}) - \nabla^2 f(\hat{x}^j) (x^{j} - x^{j-1})}^q} \\
            & \color{navy}\le \thesischfiveExpSub{\hat{\xi}^j}{\normop{\nabla^2 f(\hat{x}^j;\hat{\xi}^j) - \nabla^2 f(\hat{x}^j)}^q \cdot \norm{x^{j} - x^{j-1}}^q} \le \gamma^q \sigma_2^q,
        \end{thesischfivefitalignat}
    as
        \begin{thesischfivefitdisplay} {\color{navy}\thesischfiveExpSub{\hat{\xi}^j}{\nabla^2 f(\hat{x}^j;\hat{\xi}^j)(x^{j} - x^{j-1})} = \nabla^2 f(\hat{x}^j)(x^{j} - x^{j-1}).} \end{thesischfivefitdisplay}
    {\color{navy}Moreover, applying~Lemma~\ref{ch5:appdx-technical-lem:sadiev-generalization_chapter_5} requires}
        \begin{thesischfivefitdisplay}
            {\color{navy}\thesischfivelbd_1
            \ge 2\max_{j\in[t]}
            \norm{\nabla^2 f(\hat{x}^j)(x^j-x^{j-1})}.}
        \end{thesischfivefitdisplay}
    {\color{navy}Since}
        \begin{thesischfivefitdisplay} \norm{\nabla^2 f(\hat{x}^j) (x^{j} - x^{j-1})} \le \normop{\nabla^2 f(\hat{x}^j)} \cdot \norm{x^{j} - x^{j-1}} \le L_1 \norm{x^{j} - x^{j-1}} \le \gamma L_1, \end{thesischfivefitdisplay}
    by Assumption~\ref{ch5:ass:L-lipschitz-gradients_chapter_5} and Lemma~\ref{ch5:appdx-lem:lipschitz-gradients-implies-bounded-hessian_chapter_5}, it is enough to have $\thesischfivelbd_1 \ge 2 \gamma L_1$. Setting $G_1 \eqdef 18 \thesischfivelbd_1^{2 - q} \gamma^q \sigma_2^q/\eta > 0$ gives the desired bound.

    \paragraph{Bounding \Circled{6}.} Similarly to the previous lemma, we set $\widetilde{c}_1 \eqdef 8 \thesischfivelbd_1^2$ and $\widetilde{b}_1 \eqdef \bigl( \widetilde{c}_1/3 + \sqrt{\widetilde{c}_1^2/9 + 2 \widetilde{G}_1} \bigr) \log(2/\beta'')$. For the choice $\widetilde{G}_1$, we have the bound 
        \begin{thesischfivefitalignat}{2}
            \sum_{j = 1}^t \widetilde{s}_j^2 \thesischfiveoversetlab{\eqref{ch5:a26983bb-60b7-49ab-ac0e-a1d026f1c51e-2_chapter_5}}&{\le} 16 \thesischfivelbd_1^2 \sum_{j = 1}^t \thesischfiveExpSub{u^j, \hat{\xi}^j}{\sqnorm{(1 - \eta)^{t - j} \omega_j^u}} \\
            \thesischfiveoversetlab{\eqref{ch5:74b4e2bb-4754-4dd3-ae90-37ef40362f93-hess_chapter_5}}&{\le} 16 \thesischfivelbd_1^2 \cdot \frac{18 \thesischfivelbd_1^{2 - q} \gamma^q \sigma_2^q}{\eta} = \frac{288 \thesischfivelbd_1^{4 - q} \gamma^q \sigma_2^q}{\eta},
        \end{thesischfivefitalignat}
    (where we require $\thesischfivelbd_1 \ge 2 \gamma L_1$). Hence, if we let $\widetilde{G}_1 \eqdef 288 \thesischfivelbd_1^{4 - q} \gamma^q \sigma_2^q/\eta$ this establishes the desired bound.

    \paragraph{Bounding \Circled{7}.} {\color{navy}As in the two preceding paragraphs, we have}
        \begin{thesischfivefitdisplay} \Circled{7} \eqdef \sum_{j = 1}^t \thesischfiveExpSub{\hat{\xi}^j}{\sqnorm{(1 - \eta)^{t - j} \omega_j^u}} \thesischfiveoversetlab{\eqref{ch5:74b4e2bb-4754-4dd3-ae90-37ef40362f93-hess_chapter_5}}{\le} \frac{18 \thesischfivelbd_1^{2 - q} \gamma^q \sigma_2^q}{\eta}, \end{thesischfivefitdisplay}
    as desired.

    \paragraph{Bounding \Circled{8}.} For this bound, if we assume that $\thesischfivelbd_1 \ge 2 \gamma L_1$ then, with probability one we have
        \begin{thesischfivefitalignat}{2}
            \Circled{8} \eqdef* \norm{\sum_{j = 1}^t (1 - \eta)^{t - j} \omega_j^b} \\
            & \le \sum_{j = 1}^t (1 - \eta)^{t - j} \norm{\omega_j^b} \\
            & \le 4 \thesischfivelbd_1^{1 - q} \gamma^q \sigma_2^q \sum_{j = 1}^t (1 - \eta)^{t - j} \\
            & \le \frac{4 \thesischfivelbd_1^{1 - q} \gamma^q \sigma_2^q}{\eta},
        \end{thesischfivefitalignat}
    where in the second inequality we use~Lemma~\ref{ch5:appdx-technical-lem:sadiev-generalization_chapter_5} (with $\thesischfiveExpSubCond{\hat{\xi}^j}{\cdot}{u^j}$), more precisely, for any $j \in [t]$,
        \begin{thesischfivefitdisplay} \norm{\omega_j^b} \le 2^q \thesischfivelbd_1^{1 - q} \gamma^q \sigma_2^q \le 4 \thesischfivelbd_1^{1 - q} \gamma^q \sigma_2^q, \end{thesischfivefitdisplay}
    since 
        \begin{thesischfivefitalignat}{2}
            \MoveEqLeft[1] {\color{navy}\thesischfiveE{\norm{\nabla^2 f(\hat{x}^j;\hat{\xi}^j)(x^{j} - x^{j-1}) - \nabla^2 f(\hat{x}^j) (x^{j} - x^{j-1})}^q}} \\
            & \le {\color{navy}\thesischfiveE{\normop{\nabla^2 f(\hat{x}^j;\hat{\xi}^j) - \nabla^2 f(\hat{x}^j)}^q \cdot \norm{x^{j} - x^{j-1}}^q}} \\
            & \le \gamma^q \, {\color{navy}\thesischfiveE{\thesischfiveExpCond{\normop{\nabla^2 f(\hat{x}^j;\hat{\xi}^j) - \nabla^2 f(\hat{x}^j)}^q}{\hat{x}^j}}} \\
            & \le \gamma^q \sigma_2^q,
        \end{thesischfivefitalignat}
    where in the second inequality we use the tower property (Lemma~\ref{ch5:lem:tower-property_chapter_5}) together with $\norm{x^{j} - x^{j-1}} \le \gamma$, and in the last one we use Assumption~\ref{ch5:ass:q-bounded-central-moment-hessian_chapter_5}.

    \paragraph{Bounding \Circled{9}.} For the last bound, we use the triangle inequality, this gives
        {\color{navy}
        \begin{thesischfivefitalignatnum}{2}
            \Circled{9} \eqdef* \norm{\sum_{j = 1}^{t} (1 - \eta)^{t-j} \omega_j^r} \\
            & \le \sum_{j = 1}^t (1 - \eta)^{t - j} \norm{\omega_j^r} \\
            & = \sum_{j = 1}^t (1 - \eta)^{t - j} \norm{\nabla^2 f(\hat{x}^j)(x^{j} - x^{j-1}) - \left[ \nabla f(x^{j}) - \nabla f(x^{j-1}) \right]} \\
            & \le \min\thesischfiveens{2 \gamma L_1, \frac{\gamma^2 L_2}{2}} \sum_{j = 1}^t (1 - \eta)^{t - j} \\
            & \le \frac{2 \min\thesischfiveens{\gamma L_1, \gamma^2 L_2}}{\eta}, \label{ch5:2752ee44-9282-4441-80e7-2703274875e6_chapter_5}
        \end{thesischfivefitalignatnum}
        }
    {\color{navy}In the second inequality we use the bounds~\eqref{ch5:46d7ba91-ef00-4ec2-9b1e-d222e9968289a_chapter_5} and~\eqref{ch5:46d7ba91-ef00-4ec2-9b1e-d222e9968289b_chapter_5}, both of which hold with probability one.}

    \paragraph{Summing up all bounds \Circled{5}, \Circled{6}, \Circled{7}, \Circled{8} and \Circled{9}.} We introduce the same events $E_{\Circled{5}, t}$ and $E_{\Circled{6}, t}$ as in~\eqref{ch5:1ded530f-71dc-46eb-9fdb-5308c908cb80_chapter_5} and~\eqref{ch5:92e2f28f-a246-45f3-bdea-ba28563d3825_chapter_5} respectively. Then, by our choice of $c_1$, $b_1$ and $G_1$ we have the bound
        \begin{thesischfivefitalignatnum}{2}
            b_1 & = \left( \frac{c_1}{3} + \sqrt{\frac{c_1^2}{9} + 2 G_1} \right) \log(2/\beta'') \\
            & \le \left( \frac{2 c_1}{3} + \sqrt{2 G_1} \right) \log(2/\beta'') \\
            & = \left( \frac{4 \thesischfivelbd_1}{3} + \sqrt{\frac{36 \thesischfivelbd_1^{2 - q} \gamma^q \sigma_2^q}{\eta}} \right) \log(2/\beta'') \\
            & = \thesischfivelbd_1 \left( \frac{4}{3} + 6 \sqrt{\frac{1}{\eta} \left(\frac{\gamma \sigma_2}{\thesischfivelbd_1} \right)^q} \right) \log(2/\beta'') \\
            & \le \frac{22 \thesischfivelbd_1}{3} \log(2/\beta''), \label{ch5:35c58053-3ad9-4b62-ae23-2a6f658c08da-hess_chapter_5}
        \end{thesischfivefitalignatnum}
    where the last inequality holds provided {\color{navy}$\thesischfivelbd_1 \ge \gamma \sigma_2 \eta^{-1/q}$}. Moreover, by our choice of $\widetilde{c}_1$, $\widetilde{b}_1$ and $\widetilde{G}_1$ we also have the bound
        \begin{thesischfivefitalignatnum}{2}
            \widetilde{b}_1 & = \left( \frac{\widetilde{c}_1}{3} + \sqrt{\frac{\widetilde{c}_1^2}{9} + 2 \widetilde{G}_1} \right) \log(2/\beta'') \\
            & \le \left( \frac{2 \widetilde{c}_1}{3} + \sqrt{2 \widetilde{G}_1} \right) \log(2/\beta'') \\
            & = \left( \frac{16 \thesischfivelbd_1^2}{3} + \sqrt{\frac{576 \thesischfivelbd_1^{4 - q} \gamma^q \sigma_2^q}{\eta}} \right) \log(2/\beta'') \\
            & = \thesischfivelbd_1^2 \left( \frac{16}{3} + 24 \sqrt{\frac{1}{\eta} \left(\frac{\gamma \sigma_2}{\thesischfivelbd_1} \right)^q} \right) \log(2/\beta'') \\
            & \le \frac{88 \thesischfivelbd_1^2}{3} \log(2/\beta''), \label{ch5:c33e8593-ede2-4c3e-b76d-7e7ba9b43443-hess_chapter_5}
        \end{thesischfivefitalignatnum}
    where the last inequality holds given $\thesischfivelbd_1 \ge \gamma \sigma_2 \eta^{-1/q}$.

    Next, given $\thesischfivelbd_1 \ge \max\thesischfiveens{2 \gamma L_1, \gamma \sigma_2 \eta^{-1/q}}$, as in the previous lemma, we have $\thesischfiveProba{E_{\Circled{5}, t}} \ge 1 - \beta''$ and $\thesischfiveProba{E_{\Circled{6}, t}} \ge 1 - \beta''$ hence, by the union bound inequality we have
        \begin{thesischfivefitdisplay} \thesischfiveProba{E_{\Circled{5}, t} \cap E_{\Circled{6}, t}} \ge 1 - 2 \beta'', \end{thesischfivefitdisplay}
    and on the event $E_{\Circled{5}, t} \cap E_{\Circled{6}, t}$ we obtain the inequality
        \begin{thesischfivefitalignat}{2}
            \MoveEqLeft[1] \norm{\sum_{j = 1}^t (1 - \eta)^{t - j} \omega_j} \\
            & \le \norm{\sum_{j = 1}^t (1 - \eta)^{t - j} \omega_j^r} + \norm{\sum_{j = 1}^t (1 - \eta)^{t - j} \omega_j^b} + \norm{\sum_{j = 1}^t (1 - \eta)^{t - j} \omega_j^u} \\
            \thesischfiveoversetlab{\eqref{ch5:21eaf977-44d0-4b4b-8b63-5af59e5cc089-2_chapter_5}+\eqref{ch5:9221fb23-db5b-497f-9062-1c71b7f12882-hess_chapter_5}}&{\le} \Circled{9} + \Circled{5} + \sqrt{2 \cdot \Circled{6} + 2 \cdot \Circled{7}} + \Circled{8} \\
            \thesischfiveoversetlab{\eqref{ch5:2752ee44-9282-4441-80e7-2703274875e6_chapter_5}+\eqref{ch5:35c58053-3ad9-4b62-ae23-2a6f658c08da-hess_chapter_5}+\eqref{ch5:c33e8593-ede2-4c3e-b76d-7e7ba9b43443-hess_chapter_5}}&{\le} \frac{2 \min\thesischfiveens{\gamma L_1, \gamma^2 L_2}}{\eta} + \frac{22 \thesischfivelbd_1}{3} \log(2/\beta'') \\
            &\qquad + \sqrt{2 \cdot \frac{88 \thesischfivelbd_1^2}{3} \log(2/\beta'') + 2 \cdot \frac{18 \thesischfivelbd_1^{2 - q} \gamma^q \sigma_2^q}{\eta}} \\
            &\qquad + \frac{4 \thesischfivelbd_1^{1 - q} \gamma^q \sigma_2^q}{\eta} \\
            & \le \frac{2 \min\thesischfiveens{\gamma L_1, \gamma^2 L_2}}{\eta} \\
            &\qquad + \thesischfivelbd_1 \left( \frac{22}{3} + \sqrt{\frac{176}{3} + \frac{36}{\eta} \left( \frac{\gamma \sigma_2}{\thesischfivelbd_1} \right)^q} + \frac{4}{\eta} \left( \frac{\gamma \sigma_2}{\thesischfivelbd_1} \right)^q \right) \log(2/\beta'') \\
            & \le \frac{2 \min\thesischfiveens{\gamma L_1, \gamma^2 L_2}}{\eta} + \thesischfivelbd_1 \left( \frac{22}{3} + \sqrt{\frac{176}{3} + 36} + 4 \right) \log(2/\beta'') \\
            & \le \frac{2 \min\thesischfiveens{\gamma L_1, \gamma^2 L_2}}{\eta} + 22 \thesischfivelbd_1 \log(2/\beta''),
        \end{thesischfivefitalignat}
    where in the fourth inequality we use the fact that $\log(2/\beta'') \ge 1$ and in the fifth we use $\thesischfivelbd_1 \ge \gamma \sigma_2 \eta^{-1/q}$. This proves the lemma.
\end{proof}

\subsection{Proof of~Theorem~\ref{ch5:thm:clipped-nsgd-mvr-convergence-analysis_chapter_5}}

{\color{navy}Lemmas~\ref{ch5:appdx-lem:high-probability-analysis-bound-1_chapter_5} and~\ref{ch5:appdx-lem:high-probability-analysis-bound-2_chapter_5} now yield the high-probability convergence guarantee for \algname{\thesischfiveclip-NSGD-MVR} (Algorithm~\ref{ch5:algo:clipped-nsgd-mvr_chapter_5}).}

\begin{theorem}[Restatement of Theorem~\ref{ch5:thm:clipped-nsgd-mvr-convergence-analysis_chapter_5}]
    {\color{navy}Under Assumptions~\ref{ch5:ass:lower-boundedness_chapter_5}, \ref{ch5:ass:p-bounded-central-moment-gradient_chapter_5}, and~\ref{ch5:ass:mean-squared-smoothness_chapter_5}, let $K\ge1$ and $\beta\in\thesischfiveintof{0}{1}$ satisfy $\log(8K/\beta)\ge1$. Set $\Delta_0\eqdef f(x^0)-f^{\inf}$ and $g^0=0$ in Algorithm~\ref{ch5:algo:clipped-nsgd-mvr_chapter_5}. Run the algorithm with momentum parameter and clipping thresholds}
        {\color{navy}
        \begin{align*}
            \eta&=\max\{K^{-p/(2p-1)},K^{-pq/(p(2q+1)-2q)}\}, \\
            \thesischfivelbd_1&=2\gamma\bar L\eta^{-1/q}, \quad \thesischfivelbd_2=\max\{4\sqrt{\bar L\Delta_0},\sigma_1\eta^{-1/p}\},
        \end{align*}}
    {\color{navy}and stepsize}
        {\color{navy}
        \begin{thesischfivefitalignat}{2}
            \gamma = \min\biggl\{ & \sqrt{\frac{\Delta_0}{2\bar L K}},\ \frac{\eta}{8}\sqrt{\frac{\Delta_0}{2\bar L}},\ \frac{1}{704\eta K \log(8K/\beta)} \sqrt{\frac{\Delta_0}{\bar L}}, \\
            & \frac{\Delta_0}{176\sigma_1\eta^{(p-1)/p}K \log(8K/\beta)},\ \sqrt{\frac{\Delta_0\eta^{1/q}}{736\bar L K \log(8K/\beta)}} \biggr\}.
        \end{thesischfivefitalignat}}
        {\color{navy}Any term whose denominator contains a zero value of $\bar L$ or $\sigma_1$ is interpreted as $+\infty$ and hence does not constrain the minimum.}
    
    {\color{navy}Then, with probability at least $1 - \beta$, the iterates generated by~\ref{ch5:algo:clipped-nsgd-mvr_chapter_5} satisfy}
        \begin{thesischfivefitdisplay} \frac{1}{K} \sum_{k=0}^{K - 1} \norm{\nabla f(x^{k})} \le \frac{2 \Delta_0}{\gamma K}, \end{thesischfivefitdisplay}
    and, by our choice of parameters, the norm of the gradients converges at the rate
        \begin{thesischfivefitdisplay} \frac{1}{K} \sum_{k=0}^{K - 1} \norm{\nabla f(x^{k})} = \cO\left( \frac{\sqrt{\bar{L} \Delta_0} + \sigma_1}{K^{r}} \log \frac{K}{\beta} \right), \end{thesischfivefitdisplay}
    with high probability, where $r \eqdef \min\{(p - 1)/(2 p - 1),\, q (p - 1)/(p (2 q + 1) - 2 q)\}$.
\end{theorem}

\begin{proof}
    Let us remind from the proofs of~Lemma~\ref{ch5:appdx-lem:high-probability-analysis-bound-1_chapter_5} and Lemma~\ref{ch5:appdx-lem:high-probability-analysis-bound-2_chapter_5} that we defined
        \begin{thesischfivefitalignat}{1}
            &E_{\Circled{1},t}\eqdef
            \thesischfiveens{\abs{\sum_{j=1}^tV_j^t}\le b_2
            \ \text{or}\ \sum_{j=1}^ts_j^2>G_2\log(2/\beta')},
            \\
            &E_{\Circled{2},t}\eqdef
            \thesischfiveens{\abs{\sum_{j=1}^tY_j^t}\le\widetilde b_2
            \ \text{or}\ \sum_{j=1}^t\widetilde{s}_j^2
            >\widetilde G_2\log(2/\beta')},
            \\
            &E_{\Circled{5},t}\eqdef
            \thesischfiveens{\abs{\sum_{j=1}^tW_j^t}\le b_1
            \ \text{or}\ \sum_{j=1}^ts_j^2>G_1\log(2/\beta'')},
            \\
            &E_{\Circled{6},t}\eqdef
            \thesischfiveens{\abs{\sum_{j=1}^tZ_j^t}\le\widetilde b_1
            \ \text{or}\ \sum_{j=1}^t\widetilde{s}_j^2
            >\widetilde G_1\log(2/\beta'')}.
        \end{thesischfivefitalignat}
    where
        \begin{thesischfivefitdisplay} c_1 \eqdef 2 \thesischfivelbd_1, \quad \widetilde{c}_1 \eqdef 8 \thesischfivelbd_1^2, \quad c_2 \eqdef 2 \thesischfivelbd_2, \quad \widetilde{c}_2 \eqdef 8 \thesischfivelbd_2^2, \end{thesischfivefitdisplay}
        \begin{align*}
            G_1 & \eqdef \frac{72 \thesischfivelbd_1^{2 - q} \gamma^q \bar{L}^q}{\eta}, & \widetilde{G}_1 & \eqdef \frac{1152 \thesischfivelbd_1^{4 - q} \gamma^q \bar{L}^q}{\eta}, \\
            G_2 & \eqdef \frac{18 \thesischfivelbd_2^{2 - p} \sigma_1^p}{\eta}, & \widetilde{G}_2 & \eqdef \frac{288 \thesischfivelbd_2^{4 - p} \sigma_1^p}{\eta}, \\
            b_1 & \le \frac{40 \thesischfivelbd_1}{3} \log \frac{2}{\beta''}, & \widetilde{b}_1 & \le  \frac{160 \thesischfivelbd_1^2}{3} \log \frac{2}{\beta''}, \\
            b_2 & \le \frac{22 \thesischfivelbd_2}{3} \log \frac{2}{\beta'}, & \widetilde{b}_2 & \le \frac{88 \thesischfivelbd_2^2}{3} \log \frac{2}{\beta'},
        \end{align*}
    and we have shown the following bounds:
        \begin{gather*}
            \thesischfiveProba{E_{\Circled{1}, t}} \ge 1 - \beta', \quad \thesischfiveProba{E_{\Circled{2}, t}} \ge 1 - \beta', \\
            \thesischfiveProba{E_{\Circled{5}, t}} \ge 1 - \beta'', \quad \thesischfiveProba{E_{\Circled{6}, t}} \ge 1 - \beta'', \thesischfivenumberthis\label{ch5:4ead95f0-2b11-4c81-9020-65644456a5a3_chapter_5}
        \end{gather*}
    valid for all $t \in \thesischfiveInt{0}{K - 1}$.

    {\color{navy}For $\tau\in\thesischfiveInt{0}{K-1}$, define
    \begin{thesischfivefitdisplay}
        E_\tau\eqdef\bigcap_{t=0}^{\tau}
        \left\{\gamma\sum_{j=0}^{t}\norm{\nabla f(x^j)}
        +\Delta_{t+1}\le2\Delta_0\right\},
    \end{thesischfivefitdisplay}
    let
    \begin{thesischfivefitdisplay}
        C_\tau\eqdef
        \bigcap_{t=1}^{\tau}
        \left(E_{\Circled{1},t}\cap E_{\Circled{2},t}
        \cap E_{\Circled{5},t}\cap E_{\Circled{6},t}\right),
        \qquad C_0\eqdef\Omega,
    \end{thesischfivefitdisplay}
    and set $G_\tau\eqdef E_\tau\cap C_\tau$.

    The concentration lemmas must be used with stopping because their
    threshold hypotheses are guaranteed only while $E_{\tau-1}$ holds.
    For example, replace each martingale increment $X_j$ by
    $X_j\mathbf1_{E_{j-1}}$, with $E_{-1}\eqdef\Omega$.
    The indicator is predictable, so the stopped increments remain
    martingale differences; after the first exit they vanish, and before
    exit all moment and clipping bounds used in
    Lemmas~\ref{ch5:appdx-lem:high-probability-analysis-bound-1_chapter_5}
    and~\ref{ch5:appdx-lem:high-probability-analysis-bound-2_chapter_5}
    hold. Consequently, for every $\tau\ge1$,
    \begin{thesischfivefitdisplay}
    \begin{aligned}
        \Prob{G_{\tau-1}\cap\overline{E_{\Circled{1},\tau}}}
        &\le\beta',
        &\Prob{G_{\tau-1}\cap\overline{E_{\Circled{2},\tau}}}
        &\le\beta',\\
        \Prob{G_{\tau-1}\cap\overline{E_{\Circled{5},\tau}}}
        &\le\beta'',
        &\Prob{G_{\tau-1}\cap\overline{E_{\Circled{6},\tau}}}
        &\le\beta''.
    \end{aligned}
    \thesischfivenumberthis\label{ch5:092ed50d-cb6a-420c-b2dc-26ff706c2ac0_chapter_5}
    \end{thesischfivefitdisplay}

    For the base case, $g^0=0$ and the zero-safe update gives
    $x^1=x^0$. Moreover,
    $\norm{\nabla f(x^0)}\le\sqrt{2\bar L\Delta_0}$ and the stated
    stepsize satisfies
    $\gamma\le\frac{\eta}{8}\sqrt{\Delta_0/(2\bar L)}$.
    Therefore
    $\gamma\norm{\nabla f(x^0)}+\Delta_1\le(1+\eta/8)\Delta_0
    \le2\Delta_0$, so $\Prob{G_0}=1$.

    For the inductive step, assume
    $\Prob{G_{\tau-1}}\ge1-\frac{(\tau-1)\beta}{K}$ and take
    $\beta'=\beta''=\frac{\beta}{4K}$. The stopped bounds imply
    \begin{thesischfivefitdisplay}
        \Prob{G_{\tau-1}\cap E_{\Circled{1},\tau}
        \cap E_{\Circled{2},\tau}\cap E_{\Circled{5},\tau}
        \cap E_{\Circled{6},\tau}}
        \ge1-\frac{\tau\beta}{K}.
        \thesischfivenumberthis\label{ch5:37c9a99d-499b-4b0e-ab10-001e0f246a9e_chapter_5}
    \end{thesischfivefitdisplay}}
    Now, given $\tau \in \thesischfiveInt{1}{K - 1}$, on the event $E_{\tau - 1}$ we have
        \begin{thesischfivefitdisplay} \Delta_t \le \gamma \sum_{j = 0}^{t - 1} \norm{\nabla f(x^{j})} + \Delta_t \le 2 \Delta_0 \end{thesischfivefitdisplay}
    for all integer $1 \le t \le \tau$ and therefore, for all $j \in \thesischfiveInt{0}{\tau}$, we have
        \begin{thesischfivefitdisplay} \norm{\nabla f(x^{j})} \le \sqrt{2 \bar{L} \Delta_j} \le 2 \sqrt{\bar{L} \Delta_0}, \end{thesischfivefitdisplay}
    hence, by our choice of clipping threshold $\thesischfivelbd_2 = \max\{4 \sqrt{\bar{L} \Delta_0}, \sigma_1 \eta^{-1/p}\}$ we have $\thesischfivelbd_2 \ge 2 \max_{j \in [\tau]} \norm{\nabla f(x^{j})}$. Additionally, with the clipping level $\thesischfivelbd_1 = 2 \gamma \bar{L} \eta^{-1/q}$, on the event {\color{navy}$G_{\tau - 1} \cap E_{\Circled{1}, \tau} \cap E_{\Circled{2}, \tau} \cap E_{\Circled{5}, \tau} \cap E_{\Circled{6}, \tau}$} we have
        \begin{thesischfivefitalignat}{2}
            \MoveEqLeft[1] \gamma \sum_{k = 0}^{\tau} \norm{\nabla f(x^{k})} + \Delta_{\tau + 1} \\
            & \le \Delta_0 + \frac{2 \gamma \sqrt{2 \bar{L} \Delta_0}}{\eta} + \frac{\gamma^2 \bar{L} (\tau + 1)}{2} \\
            &\quad + 2 \gamma \eta \sum_{t = 1}^{\tau} \norm{\sum_{j = 1}^t (1 - \eta)^{t - j} \theta_j} + 2 \gamma (1 - \eta) \sum_{t = 1}^{\tau} \norm{\sum_{j = 1}^t (1 - \eta)^{t - j} \omega_j} \\
            & \le \Delta_0 + \frac{2 \gamma \sqrt{2 \bar{L} \Delta_0}}{\eta} + \frac{\gamma^2 \bar{L} (\tau + 1)}{2} \\
            &\quad + 44 \gamma \eta \tau \thesischfivelbd_2 \log(8K/\beta) + 92 \gamma (1 - \eta) \tau \thesischfivelbd_1 \log(8K/\beta),
        \end{thesischfivefitalignat}
    where the first inequality follows from Lemma~\ref{ch5:appdx-lem:high-probability-analysis-descent-lemma_chapter_5} with $\norm{\hat{e}_0} = \norm{\nabla f(x^0)} \le \sqrt{2 \bar{L} \Delta_0}$, and the second from Lemmas~\ref{ch5:appdx-lem:high-probability-analysis-bound-1_chapter_5} and~\ref{ch5:appdx-lem:high-probability-analysis-bound-2_chapter_5} with $2/\beta' = 2/\beta'' = 8K/\beta$. As $0 \le \tau \le K - 1$, by our choice of stepsize
        \begin{thesischfivefitalignatnum}{2}
            \gamma = \min\biggl\{ & \sqrt{\frac{\Delta_0}{2\bar L K}},\ \frac{\eta}{8}\sqrt{\frac{\Delta_0}{2\bar L}},\ \frac{1}{704\eta K \log(8K/\beta)} \sqrt{\frac{\Delta_0}{\bar L}}, \\
            & \frac{\Delta_0}{176\sigma_1\eta^{(p-1)/p}K \log(8K/\beta)},\ \sqrt{\frac{\Delta_0\eta^{1/q}}{736\bar L K \log(8K/\beta)}} \biggr\}, \label{ch5:787efe83-7ee7-4154-be0c-7c0b97faf778_chapter_5}
        \end{thesischfivefitalignatnum}
    each of the last four terms is at most $\Delta_0/4$ (for the term with $\thesischfivelbd_2$ we use both entries of the maximum), and we have
        \begin{thesischfivefitdisplay} \gamma \sum_{k = 0}^{\tau} \norm{\nabla f(x^{k})} + \Delta_{\tau + 1} \thesischfiveoversetlab{\eqref{ch5:787efe83-7ee7-4154-be0c-7c0b97faf778_chapter_5}}{\le} \Delta_0 + \frac{\Delta_0}{4} + \frac{\Delta_0}{4} + \frac{\Delta_0}{4} + \frac{\Delta_0}{4} = 2 \Delta_0. \end{thesischfivefitdisplay}
    {\color{navy}Thus
    \begin{thesischfivefitdisplay}
        G_{\tau-1}\cap E_{\Circled{1},\tau}
        \cap E_{\Circled{2},\tau}\cap E_{\Circled{5},\tau}
        \cap E_{\Circled{6},\tau}\subseteq G_\tau.
    \end{thesischfivefitdisplay}
    Combining this inclusion with~\eqref{ch5:37c9a99d-499b-4b0e-ab10-001e0f246a9e_chapter_5}
    yields $\Prob{G_\tau}\ge1-\tau\beta/K$, completing the induction.}
    
    \paragraph{Bounding $\frac{1}{K} \sum\limits_{k = 0}^{K - 1} \norm{\nabla f(x^{k})}$ with high probability.} {\color{navy}Taking $\tau=K-1$ gives
        \begin{thesischfivefitdisplay} \thesischfiveProba{E_{K-1}} \ge \thesischfiveProba{G_{K-1}} \ge 1-\frac{(K-1)\beta}{K}\ge1-\beta. \end{thesischfivefitdisplay}}
    hence, with probability at least $1 - \beta$ we have
        \begin{thesischfivefitdisplay} \gamma \sum_{k=0}^{K - 1} \norm{\nabla f(x^{k})} + \Delta_K \le 2 \Delta_0, \end{thesischfivefitdisplay}
    which implies the bound
        \begin{thesischfivefitdisplay} \frac{1}{K} \sum_{k=0}^{K - 1} \norm{\nabla f(x^{k})} \le \frac{2 \Delta_0}{\gamma K}, \end{thesischfivefitdisplay}
    with probability at least $1 - \beta$. By our choice of stepsize~\eqref{ch5:787efe83-7ee7-4154-be0c-7c0b97faf778_chapter_5} we obtain
        \begin{thesischfivefitalignat}{2}
            \MoveEqLeft[1] \frac1K\sum_{k=0}^{K-1}\norm{\nabla f(x^k)} \le\frac{2\Delta_0}{\gamma K} \\
            & = \cO\biggl(\max\biggl\{ \sqrt{\frac{\bar L\Delta_0}{K}},\ \frac{\sqrt{\bar L\Delta_0}}{\eta K},\ \eta\sqrt{\bar L\Delta_0}\log \frac{K}{\beta}, \\
            &\qquad \sigma_1\eta^{(p-1)/p}\log \frac{K}{\beta},\ \sqrt{\frac{\bar L\Delta_0\log(K/\beta)}{K\eta^{1/q}}} \biggr\}\biggr),
        \end{thesischfivefitalignat}
    and choosing $\eta = \max\{K^{-p/(2 p - 1)}, K^{-p q/(p (2 q + 1) - 2 q)}\}$ we have
        \begin{align*}
            \eta^{(p - 1)/p} & = \max\thesischfiveens{K^{-(p - 1)/(2 p - 1)}, K^{-q (p - 1)/(p (2 q + 1) - 2 q)}}, \\
            K^{1/2} \eta^{1/(2q)} & \ge K^{\frac{1}{2} \left( 1 - p/(p (2 q + 1) - 2 q) \right)} = K^{q (p - 1)/(p (2 q + 1) - 2 q)}, \\
            \eta K & \ge K^{1 - p/(2 p - 1)} = K^{(p - 1)/(2 p - 1)}, \\
            \eta & \le \max\thesischfiveens{K^{-(p - 1)/(2 p - 1)}, K^{-q (p - 1)/(p (2 q + 1) - 2 q)}}.
        \end{align*}
    Thus, assuming $\log(8K/\beta) \ge 1$, we get
        \begin{thesischfivefitalignat}{2}
            \MoveEqLeft[1] \frac{1}{K} \sum_{k=0}^{K - 1} \norm{\nabla f(x^{k})} \\
            & = \cO\biggl( \max\biggl\{ \sqrt{\frac{\bar{L} \Delta_0}{K}},\ \frac{\sqrt{\bar{L} \Delta_0}}{K^{(p - 1)/(2 p - 1)}} \log \frac{K}{\beta},\ \frac{\sqrt{\bar{L} \Delta_0}}{K^{q (p - 1)/(p (2 q + 1) - 2 q)}} \log \frac{K}{\beta}, \\
            &\qquad \frac{\sigma_1}{K^{(p - 1)/(2 p - 1)}} \log \frac{K}{\beta},\ \frac{\sigma_1}{K^{q (p - 1)/(p (2 q + 1) - 2 q)}} \log \frac{K}{\beta},\ \frac{\sqrt{\bar{L} \Delta_0 \log(K/\beta)}}{K^{q (p - 1)/(p (2 q + 1) - 2 q)}} \biggr\} \biggr) \\
            & = \cO\left( \frac{\sqrt{\bar{L} \Delta_0} + \sigma_1}{K^{r}} \log \frac{K}{\beta} \right),
        \end{thesischfivefitalignat}
    where $r \eqdef \min\{(p - 1)/(2 p - 1), q (p - 1)/(p (2 q + 1) - 2 q)\}$, since $(p - 1)/(2 p - 1) \le 1/2$ and $\sqrt{\log(K/\beta)} \le \log(K/\beta)$ up to constants.
    
    This concludes the proof of the theorem.
\end{proof}

\subsection{Proof of~Theorem~\ref{ch5:thm:clipped-nsgd-mvr-convergence-analysis-2_chapter_5}}

\begin{theorem}\label{ch5:thm:clipped-nsgd-mvr-convergence-analysis-2_chapter_5}
    {\color{navy}Under Assumptions~\ref{ch5:ass:lower-boundedness_chapter_5}, \ref{ch5:ass:L-lipschitz-gradients_chapter_5}, \ref{ch5:ass:p-bounded-central-moment-gradient_chapter_5}, and~\ref{ch5:ass:mean-squared-smoothness-2_chapter_5}, let $K\ge1$ and $\beta\in\thesischfiveintof{0}{1}$ satisfy $\log(8K/\beta)\ge1$. Set $\Delta_0\eqdef f(x^0)-f^{\inf}$ and $g^0=0$ in Algorithm~\ref{ch5:algo:clipped-nsgd-mvr_chapter_5}. Run the algorithm with momentum parameter and clipping thresholds}
        {\color{navy}
        \begin{align*}
            \eta&=\max\{K^{-p/(2p-1)},K^{-pq/(p(2q+1)-2q)}\}, \\
            \thesischfivelbd_1&=\max\{2\gamma L_1,\gamma\delta\eta^{-1/q}\}, \quad \thesischfivelbd_2=\max\{4\sqrt{L_1\Delta_0},\sigma_1\eta^{-1/p}\},
        \end{align*}}
    {\color{navy}and stepsize}
        {\color{navy}
        \begin{thesischfivefitalignat}{2}
            \gamma = \min\biggl\{ & \sqrt{\frac{\Delta_0}{2L_1K}},\ \frac{\eta}{8}\sqrt{\frac{\Delta_0}{2L_1}},\ \frac{1}{704\eta K \log(8K/\beta)} \sqrt{\frac{\Delta_0}{L_1}}, \\
            & \frac{\Delta_0}{176\sigma_1\eta^{(p-1)/p}K \log(8K/\beta)},\ \sqrt{\frac{\Delta_0}{352L_1K \log(8K/\beta)}}, \\
            & \sqrt{\frac{\Delta_0\eta^{1/q}}{176\delta K \log(8K/\beta)}} \biggr\}.
        \end{thesischfivefitalignat}}
    {\color{navy}Any term whose denominator contains a zero value of $L_1$, $\sigma_1$, or $\delta$ is interpreted as $+\infty$ and hence does not constrain the minimum.}
    
    {\color{navy}Then, with probability at least $1 - \beta$, the iterates generated by~\ref{ch5:algo:clipped-nsgd-mvr_chapter_5} satisfy}
        \begin{thesischfivefitdisplay} \frac{1}{K} \sum_{k=0}^{K - 1} \norm{\nabla f(x^{k})} \le \frac{2 \Delta_0}{\gamma K}, \end{thesischfivefitdisplay}
    and, by our choice of parameters, the norm of the gradients converges at the rate
        \begin{thesischfivefitalignat}{2}
            \MoveEqLeft[1] \frac{1}{K} \sum_{k=0}^{K - 1} \norm{\nabla f(x^{k})} \\
            & = \cO\left( \frac{\sqrt{L_1 \Delta_0} + \sigma_1}{K^{r}} \log \frac{K}{\beta} + \frac{\sqrt{{\color{navy}\delta} \Delta_0}}{K^{r_q}} \sqrt{\log \frac{K}{\beta}} \right),
        \end{thesischfivefitalignat}
    with high probability, where $r_q \eqdef q (p - 1)/(p (2 q + 1) - 2 q)$ and $r \eqdef \min\{(p - 1)/(2 p - 1), r_q\}$.
\end{theorem}

\begin{proof}
       From the proofs of~Lemma~\ref{ch5:appdx-lem:high-probability-analysis-bound-1_chapter_5} and Lemma~\ref{ch5:appdx-lem:high-probability-analysis-bound-3_chapter_5} we have
        \begin{thesischfivefitalignat}{1}
            &E_{\Circled{1},t}\eqdef
            \thesischfiveens{\abs{\sum_{j=1}^tV_j^t}\le b_2
            \ \text{or}\ \sum_{j=1}^ts_j^2>G_2\log(2/\beta')},
            \\
            &E_{\Circled{2},t}\eqdef
            \thesischfiveens{\abs{\sum_{j=1}^tY_j^t}\le\widetilde b_2
            \ \text{or}\ \sum_{j=1}^t\widetilde{s}_j^2
            >\widetilde G_2\log(2/\beta')},
            \\
            &E_{\Circled{5},t}\eqdef
            \thesischfiveens{\abs{\sum_{j=1}^tW_j^t}\le b_1
            \ \text{or}\ \sum_{j=1}^ts_j^2>G_1\log(2/\beta'')},
            \\
            &E_{\Circled{6},t}\eqdef
            \thesischfiveens{\abs{\sum_{j=1}^tZ_j^t}\le\widetilde b_1
            \ \text{or}\ \sum_{j=1}^t\widetilde{s}_j^2
            >\widetilde G_1\log(2/\beta'')}.
        \end{thesischfivefitalignat}
    where
        \begin{thesischfivefitdisplay} c_1 \eqdef 2 \thesischfivelbd_1, \quad \widetilde{c}_1 \eqdef 8 \thesischfivelbd_1^2, \quad c_2 \eqdef 2 \thesischfivelbd_2, \quad \widetilde{c}_2 \eqdef 8 \thesischfivelbd_2^2, \end{thesischfivefitdisplay}
        \begin{align*}
            G_1 & \eqdef \frac{18 \thesischfivelbd_1^{2 - q} \gamma^q {\color{navy}\delta}^q}{\eta}, & \widetilde{G}_1 & \eqdef \frac{288 \thesischfivelbd_1^{4 - q} \gamma^q {\color{navy}\delta}^q}{\eta}, \\
            G_2 & \eqdef \frac{18 \thesischfivelbd_2^{2 - p} \sigma_1^p}{\eta}, & \widetilde{G}_2 & \eqdef \frac{288 \thesischfivelbd_2^{4 - p} \sigma_1^p}{\eta}, \\
            b_1 & \le \frac{22 \thesischfivelbd_1}{3} \log \frac{2}{\beta''}, & \widetilde{b}_1 & \le  \frac{88 \thesischfivelbd_1^2}{3} \log \frac{2}{\beta''}, \\
            b_2 & \le \frac{22 \thesischfivelbd_2}{3} \log \frac{2}{\beta'}, & \widetilde{b}_2 & \le \frac{88 \thesischfivelbd_2^2}{3} \log \frac{2}{\beta'},
        \end{align*}
    and we have shown the following bounds:
        \begin{gather*}
            \thesischfiveProba{E_{\Circled{1}, t}} \ge 1 - \beta', \quad \thesischfiveProba{E_{\Circled{2}, t}} \ge 1 - \beta', \\
            \thesischfiveProba{E_{\Circled{5}, t}} \ge 1 - \beta'', \quad \thesischfiveProba{E_{\Circled{6}, t}} \ge 1 - \beta'', \thesischfivenumberthis\label{ch5:4ead95f0-2b11-4c81-9020-65644456a5a3-2_chapter_5}
        \end{gather*}
    valid for all $t \in \thesischfiveInt{0}{K - 1}$.

    {\color{navy}We use exactly the stopped-process construction
    $G_\tau=E_\tau\cap C_\tau$ from the preceding proof, now with
    Lemma~\ref{ch5:appdx-lem:high-probability-analysis-bound-3_chapter_5}
    in place of Lemma~\ref{ch5:appdx-lem:high-probability-analysis-bound-2_chapter_5}.
    Thus the four conditional intersection-failure probabilities are bounded
    by $\beta',\beta',\beta'',\beta''$, respectively, while
    $G_{\tau-1}$ holds. With
    $\beta'=\beta''=\beta/(4K)$, the same union-bound induction gives
    $\Prob{G_\tau}\ge1-\tau\beta/K$.}
    Then, given $\tau \in \thesischfiveInt{1}{K - 1}$, on the event $E_{\tau - 1}$ we have
        \begin{thesischfivefitdisplay} \Delta_t \le \gamma \sum_{j = 0}^{t - 1} \norm{\nabla f(x^{j})} + \Delta_t \le 2 \Delta_0 \end{thesischfivefitdisplay}
    for all integer $1 \le t \le \tau$ and therefore, for all $j \in \thesischfiveInt{0}{\tau}$, we have
        \begin{thesischfivefitdisplay} {\color{navy}\norm{\nabla f(x^{j})} \le \sqrt{2 L_1 \Delta_j} \le 2 \sqrt{L_1 \Delta_0},} \end{thesischfivefitdisplay}
    hence, by our choice of clipping threshold $\thesischfivelbd_2 = \max\{4 \sqrt{L_1 \Delta_0}, \sigma_1 \eta^{-1/p}\}$ we have $\thesischfivelbd_2 \ge 2 \max_{j \in [\tau]} \norm{\nabla f(x^{j})}$. Additionally, with the clipping level $\thesischfivelbd_1 = \max\{2 \gamma L_1, \gamma {\color{navy}\delta} \eta^{-1/q}\}$, on the event {\color{navy}$G_{\tau - 1} \cap E_{\Circled{1}, \tau} \cap E_{\Circled{2}, \tau} \cap E_{\Circled{5}, \tau} \cap E_{\Circled{6}, \tau}$} we have
        \begin{thesischfivefitalignatnum}{2}
            \MoveEqLeft[1] \gamma \sum_{k = 0}^{\tau} \norm{\nabla f(x^{k})} + \Delta_{\tau + 1} \\
            & \le \Delta_0 + \frac{2 \gamma \sqrt{2 L_1 \Delta_0}}{\eta} + \frac{\gamma^2 L_1 (\tau + 1)}{2} \\
            &\quad + 2 \gamma \eta \sum_{t = 1}^{\tau} \norm{\sum_{j = 1}^t (1 - \eta)^{t - j} \theta_j} + 2 \gamma (1 - \eta) \sum_{t = 1}^{\tau} \norm{\sum_{j = 1}^t (1 - \eta)^{t - j} \omega_j} \\
            & \le \Delta_0 + \frac{2 \gamma \sqrt{2 L_1 \Delta_0}}{\eta} + \frac{\gamma^2 L_1 (\tau + 1)}{2} \\
            &\quad + 44 \gamma \eta \tau \thesischfivelbd_2 \log(8K/\beta) + 44 \gamma (1 - \eta) \tau \thesischfivelbd_1 \log(8K/\beta), \label{ch5:0f3f34c1-906a-43c0-b66b-9472d1bce764_chapter_5}
        \end{thesischfivefitalignatnum}
    where the first inequality follows from Corollary~\ref{ch5:appdx-lem:high-probability-analysis-descent-lemma-2_chapter_5} and the second from Lemmas~\ref{ch5:appdx-lem:high-probability-analysis-bound-1_chapter_5} and~\ref{ch5:appdx-lem:high-probability-analysis-bound-3_chapter_5}. As $0 \le \tau \le K - 1$, by our choice of stepsize
        \begin{thesischfivefitalignatnum}{2}
            \gamma = \min\biggl\{ & \sqrt{\frac{\Delta_0}{2L_1K}},\ \frac{\eta}{8}\sqrt{\frac{\Delta_0}{2L_1}},\ \frac{1}{704\eta K \log(8K/\beta)} \sqrt{\frac{\Delta_0}{L_1}}, \\
            & \frac{\Delta_0}{176\sigma_1\eta^{(p-1)/p}K \log(8K/\beta)},\ \sqrt{\frac{\Delta_0}{352L_1K \log(8K/\beta)}}, \\
            & \sqrt{\frac{\Delta_0\eta^{1/q}}{176{\color{navy}\delta} K \log(8K/\beta)}} \biggr\}, \label{ch5:787efe83-7ee7-4154-be0c-7c0b97faf778-2_chapter_5}
        \end{thesischfivefitalignatnum}
    each of the last four terms is at most $\Delta_0/4$ (for the terms with $\thesischfivelbd_1$ and $\thesischfivelbd_2$ we use both entries of the maximum), and we have
        \begin{thesischfivefitdisplay} \gamma \sum_{k = 0}^{\tau} \norm{\nabla f(x^{k})} + \Delta_{\tau + 1} \thesischfiveoversetlab{\eqref{ch5:787efe83-7ee7-4154-be0c-7c0b97faf778-2_chapter_5}}{\le} \Delta_0 + \frac{\Delta_0}{4} + \frac{\Delta_0}{4} + \frac{\Delta_0}{4} + \frac{\Delta_0}{4} = 2 \Delta_0. \end{thesischfivefitdisplay}
    {\color{navy}This deterministic implication closes the stopped-event induction exactly as in the preceding proof.}

    Then, as before, with probability {\color{navy}at least} $1 - \beta$, we have
        \begin{thesischfivefitdisplay} \gamma \sum_{k=0}^{K - 1} \norm{\nabla f(x^{k})} + \Delta_K \le 2 \Delta_0, \end{thesischfivefitdisplay}
    which implies the bound
        \begin{thesischfivefitdisplay} \frac{1}{K} \sum_{k=0}^{K - 1} \norm{\nabla f(x^{k})} \le \frac{2 \Delta_0}{\gamma K}, \end{thesischfivefitdisplay}
    with probability at least $1 - \beta$. By our choice of stepsize~\eqref{ch5:787efe83-7ee7-4154-be0c-7c0b97faf778-2_chapter_5} we obtain
        \begin{thesischfivefitalignat}{2}
            \MoveEqLeft[1] \frac1K\sum_{k=0}^{K-1}\norm{\nabla f(x^k)} \le\frac{2\Delta_0}{\gamma K} \\
            & = \cO\biggl(\max\biggl\{ \sqrt{\frac{L_1\Delta_0}{K}},\ \frac{\sqrt{L_1\Delta_0}}{\eta K},\ \eta\sqrt{L_1\Delta_0}\log \frac{K}{\beta}, \\
            &\qquad \sigma_1\eta^{(p-1)/p}\log \frac{K}{\beta},\ \sqrt{\frac{L_1\Delta_0\log(K/\beta)}{K}},\ \sqrt{\frac{{\color{navy}\delta}\Delta_0\log(K/\beta)}{K\eta^{1/q}}} \biggr\}\biggr) \\
            & = \cO\biggl(\max\biggl\{ \frac{\sqrt{L_1\Delta_0}}{\eta K},\ \eta\sqrt{L_1\Delta_0}\log \frac{K}{\beta},\ \sigma_1\eta^{(p-1)/p}\log \frac{K}{\beta}, \\
            &\qquad \sqrt{\frac{L_1\Delta_0\log(K/\beta)}{K}},\ \sqrt{\frac{{\color{navy}\delta}\Delta_0\log(K/\beta)}{K\eta^{1/q}}} \biggr\}\biggr),
        \end{thesischfivefitalignat}
    as we assumed $\log(8K/\beta) \ge 1$. Now, with $r_q \eqdef q (p - 1)/(p (2 q + 1) - 2 q)$ and our choice of $\eta$, we have
        \begin{align*}
            \eta^{(p - 1)/p} & = \max\thesischfiveens{K^{-(p - 1)/(2 p - 1)}, K^{-r_q}}, \\
            K^{1/2} \eta^{1/(2q)} & \ge K^{\frac{1}{2} \left( 1 - p/(p (2 q + 1) - 2 q) \right)} = K^{r_q}, \\
            \eta K & \ge K^{1 - p/(2 p - 1)} = K^{(p - 1)/(2 p - 1)}, \\
            \eta & \le \max\thesischfiveens{K^{-(p - 1)/(2 p - 1)}, K^{-r_q}}.
        \end{align*}
    Thus, given $\log(8K/\beta) \ge 1$, we get
        \begin{thesischfivefitalignat}{2}
            \MoveEqLeft[1] \frac{1}{K}\sum_{k=0}^{K-1}\norm{\nabla f(x^k)} \\
            & = \cO\biggl(\max\biggl\{ \frac{\sqrt{L_1\Delta_0}}{K^{(p-1)/(2p-1)}}\log \frac{K}{\beta},\ \frac{\sqrt{L_1\Delta_0}}{K^{r_q}}\log \frac{K}{\beta}, \\
            &\qquad \frac{\sigma_1}{K^{(p-1)/(2p-1)}}\log \frac{K}{\beta},\ \frac{\sigma_1}{K^{r_q}}\log \frac{K}{\beta},\ \frac{\sqrt{{\color{navy}\delta}\Delta_0\log(K/\beta)}}{K^{r_q}} \biggr\}\biggr) \\
            & = \cO\left( \frac{\sqrt{L_1\Delta_0}+\sigma_1}{K^{r}}\log \frac{K}{\beta} + \frac{\sqrt{{\color{navy}\delta}\Delta_0}}{K^{r_q}} \sqrt{\log \frac{K}{\beta}} \right),
        \end{thesischfivefitalignat}
    where $r \eqdef \min\{(p - 1)/(2 p - 1), r_q\}$, since $(p - 1)/(2 p - 1) \le 1/2$ and $\sqrt{\log(K/\beta)} \lesssim \log(K/\beta)$.
    
    This concludes the proof of the theorem.
\end{proof}

\subsection{Proof of~Theorem~\ref{ch5:thm:clipped-nsgd-mvr-convergence-analysis-hess_chapter_5}}

\begin{algorithm}[htbp]
\caption{{\color{navy}\algname{D-Clip-NSGDHess}} (Double-Clipped Normalized SGD with Hessian-corrected Momentum)}
\label{ch5:algo:clipped-nsgd-mvr-hess_chapter_5}
\begin{algorithmic}[1]
\Require $x^0 \in \R^d$, the starting point;              $K > 0$, the number of iterations;              $g^0 \in \R^d$, an initial vector;              $\gamma > 0$, the stepsize;              $\eta \in \thesischfiveintof{0}{1}$, the momentum parameter for \texttt{Hess};              $\lambda_1, \lambda_2 >0$, the clipping thresholds;
\State ${\color{navy}x^1 = x^0 - \gamma \unitdir{g^0}}$
\For{$k = 1, 2, \ldots, K - 1$}
\State {\color{navy}Sample $u^k \sim \mathcal{U}\left(\thesischfiveintff{0}{1}\right)$}
\State {\color{navy}$\hat{x}^k = u^k x^{k} + (1 - u^k) x^{k-1}$}
\State \textit{Apply \algname{NSGDHess}; here $\xi^{k}, \hat{\xi}^{k} \sim \thesischfivecD$ are independent.}
\State ${\color{navy}g^{k} = (1 - \eta) \left( g^{k-1}
    + \thesischfiveclip\left(\nabla^2 f(\hat{x}^k;\hat{\xi}^{k})
    (x^{k}-x^{k-1}),\lambda_1\right) \right)}$
\Statex \hspace{\algorithmicindent}${\color{navy}\phantom{g^{k}}
    + \eta \thesischfiveclip\left(\nabla f(x^{k};\xi^{k}),\lambda_2\right)}$
\State \textit{Do one descent step.}
\State ${\color{navy}x^{k+1} = x^{k} - \gamma \unitdir{g^{k}}}$
\EndFor
\Ensure ${\color{navy}x_{\mathrm{out}}\sim\mathcal U\{x^0,\ldots,x^{K-1}\}}$
\end{algorithmic}
\end{algorithm}%
\begin{theorem}\label{ch5:thm:clipped-nsgd-mvr-convergence-analysis-hess_chapter_5}
    {\color{navy}Under Assumptions~\ref{ch5:ass:lower-boundedness_chapter_5}, \ref{ch5:ass:L-lipschitz-gradients_chapter_5}, \ref{ch5:ass:p-bounded-central-moment-gradient_chapter_5}, \ref{ch5:ass:L-lipschitz-hessians_chapter_5}, and~\ref{ch5:ass:q-bounded-central-moment-hessian_chapter_5}, let $K\ge1$ and $\beta\in\thesischfiveintof{0}{1}$ satisfy $\log(8K/\beta)\ge1$. Set $\Delta_0\eqdef f(x^0)-f^{\inf}$ and $g^0=0$ in Algorithm~\ref{ch5:algo:clipped-nsgd-mvr-hess_chapter_5}. Run the algorithm with momentum parameter and clipping thresholds}
        {\color{navy}
        \begin{align*}
            \eta&=\max\{K^{-p/(2p-1)},K^{-pq/(p(2q+1)-2q)}\}, \\
            \thesischfivelbd_1&=\max\{2\gamma L_1,\gamma\sigma_2\eta^{-1/q}\}, \quad \thesischfivelbd_2=\max\{4\sqrt{L_1\Delta_0},\sigma_1\eta^{-1/p}\},
        \end{align*}}
    {\color{navy}and stepsize}
        {\color{navy}
        \begin{thesischfivefitalignat}{2}
            \gamma = \min\biggl\{ & \sqrt{\frac{2\Delta_0}{5L_1K}},\ \frac{\eta}{10}\sqrt{\frac{\Delta_0}{2L_1}},\ \sqrt[3]{\frac{\Delta_0\eta}{20L_2K}}, \\
            & \frac{1}{880\eta K \log(8K/\beta)}\sqrt{\frac{\Delta_0}{L_1}},\ \frac{\Delta_0}{220\sigma_1\eta^{(p-1)/p}K \log(8K/\beta)}, \\
            & \sqrt{\frac{\Delta_0}{440L_1K \log(8K/\beta)}},\ \sqrt{\frac{\Delta_0\eta^{1/q}}{220\sigma_2K \log(8K/\beta)}} \biggr\}.
        \end{thesischfivefitalignat}}
    {\color{navy}Any term whose denominator contains a zero value of $L_1$, $L_2$, $\sigma_1$, or $\sigma_2$ is interpreted as $+\infty$ and hence does not constrain the minimum.}
    
    {\color{navy}Then, with probability at least $1 - \beta$, the iterates generated by~\ref{ch5:algo:clipped-nsgd-mvr-hess_chapter_5} satisfy}
        \begin{thesischfivefitdisplay} \frac{1}{K} \sum_{k=0}^{K - 1} \norm{\nabla f(x^{k})} \le \frac{2 \Delta_0}{\gamma K}, \end{thesischfivefitdisplay}
    and, by our choice of parameters, the norm of the gradients converges at the rate
        \begin{thesischfivefitalignat}{2}
            \MoveEqLeft[1] \frac{1}{K} \sum_{k=0}^{K - 1} \norm{\nabla f(x^{k})} \\
            & = \cO\left( \frac{( L_2^{1/2} \Delta_0 )^{2/3}}{K^{2 (p - 1)/(4 p - 3)}} + \frac{\sqrt{L_1 \Delta_0} + \sigma_1}{K^{r}} \log \frac{K}{\beta} + \frac{\sqrt{\sigma_2 \Delta_0}}{K^{r_q}} \sqrt{\log \frac{K}{\beta}} \right),
        \end{thesischfivefitalignat}
    with high probability, where $r_q \eqdef q (p - 1)/(p (2 q + 1) - 2 q)$ and $r \eqdef \min\{(p - 1)/(2 p - 1), r_q\}$.
\end{theorem}

\begin{proof}
    {\color{navy}We use the same stopped-process induction as in the two preceding proofs, with the Hessian concentration event supplied by Lemma~\ref{ch5:appdx-lem:high-probability-analysis-bound-hess_chapter_5}. This is essential because its clipping hypotheses, too, are guaranteed on the preceding good event.} The proof otherwise follows from the proof of~Theorem~\ref{ch5:thm:clipped-nsgd-mvr-convergence-analysis-2_chapter_5} with $\delta$ replaced by $\sigma_2$. Additionally, by~Lemma~\ref{ch5:appdx-lem:high-probability-analysis-bound-hess_chapter_5} we need to bound an extra term in~\eqref{ch5:0f3f34c1-906a-43c0-b66b-9472d1bce764_chapter_5}, i.e., we have
        \begin{thesischfivefitalignat}{2}
            \MoveEqLeft[1] \gamma \sum_{k = 0}^{\tau} \norm{\nabla f(x^{k})} + \Delta_{\tau + 1} \\
            & \le \Delta_0 + \frac{2 \gamma \sqrt{2 L_1 \Delta_0}}{\eta} + \frac{\gamma^2 L_1 (\tau + 1)}{2} \\
            &\quad + 2 \gamma \eta \sum_{t = 1}^{\tau} \norm{\sum_{j = 1}^t (1 - \eta)^{t - j} \theta_j} + 2 \gamma (1 - \eta) \sum_{t = 1}^{\tau} \norm{\sum_{j = 1}^t (1 - \eta)^{t - j} \omega_j} \\
            & \le \Delta_0 + \frac{2 \gamma \sqrt{2 L_1 \Delta_0}}{\eta} + \frac{\gamma^2 L_1 (\tau + 1)}{2} + \frac{4 (1 - \eta) \gamma^3 L_2 (\tau + 1)}{\eta} \\
            &\quad + 44 \gamma \eta \tau \thesischfivelbd_2 \log(8K/\beta) + 44 \gamma (1 - \eta) \tau \thesischfivelbd_1 \log(8K/\beta),
        \end{thesischfivefitalignat}
    where the first inequality follows from Lemma~{\color{navy}\ref{ch5:appdx-lem:high-probability-analysis-descent-lemma-hess_chapter_5}} and the second from Lemmas~{\color{navy}\ref{ch5:appdx-lem:high-probability-analysis-bound-1_chapter_5} and~\ref{ch5:appdx-lem:high-probability-analysis-bound-hess_chapter_5}}; in the last inequality we drop the $\min\thesischfiveens{\cdots}$ term and simply keep the $\gamma^3 L_2$ term. Additionally, as we have $0 \le \tau \le K - 1$, choosing
        \begin{thesischfivefitalignatnum}{2}
            \gamma = \min\biggl\{ & \sqrt{\frac{2\Delta_0}{5L_1K}},\ \frac{\eta}{10}\sqrt{\frac{\Delta_0}{2L_1}},\ \sqrt[3]{\frac{\Delta_0\eta}{20L_2K}}, \\
            & \frac{1}{880\eta K \log(8K/\beta)}\sqrt{\frac{\Delta_0}{L_1}},\ \frac{\Delta_0}{220\sigma_1\eta^{(p-1)/p}K \log(8K/\beta)}, \\
            & \sqrt{\frac{\Delta_0}{440L_1K \log(8K/\beta)}},\ \sqrt{\frac{\Delta_0\eta^{1/q}}{220\sigma_2K \log(8K/\beta)}} \biggr\}, \label{ch5:787efe83-7ee7-4154-be0c-7c0b97faf778-hess_chapter_5}
        \end{thesischfivefitalignatnum}
    makes each of the last five terms at most $\Delta_0/5$, which gives
        \begin{thesischfivefitdisplay} \gamma \sum_{k = 0}^{\tau} \norm{\nabla f(x^{k})} + \Delta_{\tau + 1} \thesischfiveoversetlab{\eqref{ch5:787efe83-7ee7-4154-be0c-7c0b97faf778-hess_chapter_5}}{\le} \Delta_0 + \frac{\Delta_0}{5} + \frac{\Delta_0}{5} + \frac{\Delta_0}{5} + \frac{\Delta_0}{5} + \frac{\Delta_0}{5} = 2 \Delta_0. \end{thesischfivefitdisplay}
    Then, as before, with probability {\color{navy}at least} $1 - \beta$, we have
        \begin{thesischfivefitdisplay} \gamma \sum_{k=0}^{K - 1} \norm{\nabla f(x^{k})} + \Delta_K \le 2 \Delta_0, \end{thesischfivefitdisplay}
    which implies the bound
        \begin{thesischfivefitdisplay} \frac{1}{K} \sum_{k=0}^{K - 1} \norm{\nabla f(x^{k})} \le \frac{2 \Delta_0}{\gamma K}, \end{thesischfivefitdisplay}
    with probability at least $1 - \beta$. By our choice of stepsize~\eqref{ch5:787efe83-7ee7-4154-be0c-7c0b97faf778-hess_chapter_5} we obtain
        \begin{thesischfivefitalignat}{2}
            \MoveEqLeft[1] \frac{1}{K} \sum_{k=0}^{K - 1} \norm{\nabla f(x^{k})} \le \frac{2 \Delta_0}{\gamma K} \\
            & = \cO\biggl( \max\biggl\{\sqrt{\frac{L_1 \Delta_0}{K}},\ \frac{\sqrt{L_1 \Delta_0}}{\eta K},\ \eta^{-1/3} \left( \frac{L_2^{1/2} \Delta_0}{K} \right)^{2/3},\ \eta \sqrt{L_1 \Delta_0} \log \frac{K}{\beta}, \\
            &\qquad \sigma_1 \eta^{(p - 1)/p} \log \frac{K}{\beta},\ \sqrt{\frac{L_1 \Delta_0 \log(K/\beta)}{K}},\ \sqrt{\frac{\sigma_2 \Delta_0 \log(K/\beta)}{K \eta^{1/q}}}\biggr\} \biggr) \\
            & = \cO\biggl( \max\biggl\{\frac{\sqrt{L_1 \Delta_0}}{\eta K},\ \eta^{-1/3} \left( \frac{L_2^{1/2} \Delta_0}{K} \right)^{2/3},\ \eta \sqrt{L_1 \Delta_0} \log \frac{K}{\beta}, \\
            &\qquad \sigma_1 \eta^{(p - 1)/p} \log \frac{K}{\beta},\ \sqrt{\frac{L_1 \Delta_0 \log(K/\beta)}{K}},\ \sqrt{\frac{\sigma_2 \Delta_0 \log(K/\beta)}{K \eta^{1/q}}}\biggr\} \biggr),
        \end{thesischfivefitalignat}
    as we assumed $\log(8K/\beta) \ge 1$. Now, choosing 
        \begin{align*}
            \eta & = \max\thesischfiveens{K^{-p/(2 p - 1)}, K^{-p q/(p (2 q + 1) - 2 q)}, K^{-2 p/(4 p - 3)}} \\
            & = \max\thesischfiveens{K^{-p/(2 p - 1)}, K^{-p q/(p (2 q + 1) - 2 q)}},
        \end{align*}
    since $p/(2 p - 1) \le 2 p/(4 p - 3)$, we have, with $r_q \eqdef q (p - 1)/(p (2 q + 1) - 2 q)$,
        \begin{align*}
            \eta^{(p - 1)/p} & = \max\thesischfiveens{K^{-(p - 1)/(2 p - 1)}, K^{-r_q}}, \\
            K^{1/2} \eta^{1/(2q)} & \ge K^{\frac{1}{2} \left( 1 - p/(p (2 q + 1) - 2 q) \right)} = K^{r_q}, \\
            \eta K & \ge K^{1 - p/(2 p - 1)} = K^{(p - 1)/(2 p - 1)}, \\
            {\color{navy}K^{2/3} \eta^{1/3}} & {\color{navy}\ge K^{2/3 - 2 p/(3 (4 p - 3))} = K^{2 (p - 1)/(4 p - 3)}}, \\
            \eta & \le \max\thesischfiveens{K^{-(p - 1)/(2 p - 1)}, K^{-r_q}}.
        \end{align*}
    Thus, given $\log(8K/\beta) \ge 1$, we get
        \begin{thesischfivefitalignat}{2}
            \MoveEqLeft[1] \frac{1}{K}\sum_{k=0}^{K-1}\norm{\nabla f(x^k)} \\
            & = \cO\biggl(\max\biggl\{ \frac{(L_2^{1/2}\Delta_0)^{2/3}}{K^{2(p-1)/(4p-3)}},\ \frac{\sqrt{L_1\Delta_0}}{K^{(p-1)/(2p-1)}}\log \frac{K}{\beta},\ \frac{\sqrt{L_1\Delta_0}}{K^{r_q}}\log \frac{K}{\beta}, \\
            &\qquad \frac{\sigma_1}{K^{(p-1)/(2p-1)}}\log \frac{K}{\beta},\ \frac{\sigma_1}{K^{r_q}}\log \frac{K}{\beta},\ \frac{\sqrt{\sigma_2\Delta_0\log(K/\beta)}}{K^{r_q}} \biggr\}\biggr) \\
            & = \cO\left( \frac{(L_2^{1/2}\Delta_0)^{2/3}}{K^{2(p-1)/(4p-3)}} + \frac{\sqrt{L_1\Delta_0}+\sigma_1}{K^{r}}\log \frac{K}{\beta} + \frac{\sqrt{\sigma_2\Delta_0}}{K^{r_q}}\sqrt{\log \frac{K}{\beta}} \right),
        \end{thesischfivefitalignat}
    where $r \eqdef \min\{(p - 1)/(2 p - 1), r_q\}$, since $(p - 1)/(2 p - 1) \le 1/2$ and $\sqrt{\log(K/\beta)} \lesssim \log(K/\beta)$.
    
    This concludes the proof of the theorem.
\end{proof}

\newpage

\section{Technical Lemmas}

{\color{navy}We next collect the technical lemmas used throughout this chapter.}

\begin{lemma}[{\textcite[Lemma~10]{hübler2025gradientclippingnormalizationheavy}}]\label{ch5:appdx-technical-lem:app-von-Bahr-and-Essen_chapter_5}
    Let $p \in \thesischfiveintff{1}{2}$, and $X_1, \ldots, X_n \in \R^d$ be a martingale difference sequence (MDS), i.e., $\thesischfiveExpCond{X_j}{X_{j-1}, \ldots, X_1} = 0$ a.s. for all $j = 1, \ldots, n$ satisfying 
	\begin{align*}
	   \thesischfiveE{\norm{X_j}^p } < \infty, \qquad \text{for all } j = 1, \ldots, n. 
	\end{align*}
	Define $S_n \eqdef \sum\limits_{j=1}^n X_j$, then 
	\begin{align*}
		\thesischfiveE{\norm{S_n}^p} \leq 2 \sum_{j = 1}^n \thesischfiveE{\norm{X_j}^p}.
	\end{align*}
\end{lemma}

The following lemma will be useful to derive high-probability convergence bounds. The constant $2^{2 - p}$ appearing in~\eqref{ch5:2425acdd-8fb2-4df1-ae5c-908f6e744926_chapter_5} is a consequence of~\textcite[Lemma~1]{he2025complexitynormalizedstochasticfirstorder} which we recall in~Lemma~\ref{ch5:lem:euclidean-norm-expansion-generalized_chapter_5}.
\begin{lemma}[{\textcite[Lemma~10]{NEURIPS2021_26901deb}}]\label{ch5:appdx-technical-lem:cutkosky-mehta_chapter_5}
    Let $X_1, \ldots, X_n \in \R^d$ be a sequence of random vectors. Let us define the following sequence $w_1, \ldots, w_n$ recursively by
        \begin{enumerate}
            \item $w_0 = 0$,
            \item if $\sum\limits_{i = 1}^{j - 1} X_i \neq 0$, then we set
                \begin{thesischfivefitdisplay} w_j = \thesischfivesign\left( \sum\limits_{i = 1}^{j - 1} w_i \right) \frac{\thesischfiveps{\sum\limits_{i = 1}^{j - 1} X_i}{X_j}}{\norm{\sum\limits_{i = 1}^{j - 1} X_i}}, \end{thesischfivefitdisplay}
            \item if $\sum\limits_{i = 1}^{j - 1} X_i = 0$, we set $w_j = 0$.
        \end{enumerate}
    Then, $\abs{w_j} \le \norm{X_j}$ for all $j \in \thesischfiveInt{1}{n}$ and
        \begin{thesischfivefitdisplay} \norm{\sum_{i = 1}^n X_i} \le \abs{\sum_{i = 1}^n w_i} + \left( \max_{i \in [n]} \norm{X_i}^p + 2^{2 - p} \sum_{i = 1}^n \norm{X_i}^p  \right)^{\frac{1}{p}}. \thesischfivenumberthis\label{ch5:2425acdd-8fb2-4df1-ae5c-908f6e744926_chapter_5} \end{thesischfivefitdisplay}
\end{lemma}

\begin{lemma}[Freedman's Inequality]\label{ch5:appdx-technical-lem:freedman-inequality_chapter_5}
    Let $X_1, \ldots, X_n \in \R^d$ be a martingale difference sequence (MDS), i.e., $\thesischfiveExpCond{X_j}{X_{j-1}, \ldots, X_1} = 0$ a.s. for $j = 1, \ldots, n$. Assume there exists a constant $c > 0$ such that $\abs{X_j} \le c$ a.s. for $j = 1, \ldots, n$ and define $\sigma_j^2 \eqdef \thesischfiveExpCond{X_j^2}{X_{j - 1}, \ldots, X_1}$. Then, for all $b > 0$ and all $G > 0$ we have
        \begin{thesischfivefitalignat}{2}
            \MoveEqLeft[1] \thesischfiveProba{\abs{\sum_{j = 1}^n X_j} > b \quad \text{and} \quad \sum_{j = 1}^n \sigma_j^2 \le G} \\
            & \le 2 \exp\left( -\frac{b^2}{2 G + \frac{2 b c}{3}} \right).
        \end{thesischfivefitalignat}
\end{lemma}

The next lemma is a generalization of~\textcite[Lemma~5.1]{pmlr-v202-sadiev23a}.
\begin{lemma}\label{ch5:appdx-technical-lem:sadiev-generalization_chapter_5}
    Let $\thesischfivelbd > 0$ be a scalar, $X \in \R^d$ be a random vector and let $\widetilde{X} = \thesischfiveclip(X, \thesischfivelbd)$ then
        \begin{thesischfivefitdisplay} \norm{\widetilde{X} - \thesischfiveE{\widetilde{X}}} \le 2 \thesischfivelbd. \thesischfivenumberthis\label{ch5:3f712cbe-7d84-446c-a934-ad43a82e08dd-1_chapter_5} \end{thesischfivefitdisplay}
    Moreover, assume that $\thesischfiveE{X} = x \in \R^d$, $\thesischfiveE{\norm{X - x}^p} \le \sigma_1^p$ and $\norm{x} \le \thesischfivelbd/2$ for some $\sigma_1 \ge 0$ and $p \in \thesischfiveintof{1}{2}$. Then, for any $q \in \thesischfiveintff{p}{2}$ we have
        \begin{alignat*}{2}
            \norm{\thesischfiveE{\widetilde{X}} - x} & \le \frac{2^p \sigma_1^p}{\thesischfivelbd^{p - 1}}, \thesischfivenumberthis\label{ch5:3f712cbe-7d84-446c-a934-ad43a82e08dd-2_chapter_5} \\
            \thesischfiveE{\norm{\widetilde{X} - x}^q} & \le 2 \cdot 3^q \thesischfivelbd^{q - p} \sigma_1^p, \thesischfivenumberthis\label{ch5:3f712cbe-7d84-446c-a934-ad43a82e08dd-3_chapter_5} \\
            \thesischfiveE{\norm{\widetilde{X} - \thesischfiveE{\widetilde{X}}}^q} & \le 2^{3 - q} \cdot 3^q \thesischfivelbd^{q - p} \sigma_1^p. \thesischfivenumberthis\label{ch5:3f712cbe-7d84-446c-a934-ad43a82e08dd-4_chapter_5}
        \end{alignat*}
\end{lemma}

\begin{proof}
    The proof follows the same strategy as in~\textcite{pmlr-v202-sadiev23a}; inequalities~\eqref{ch5:3f712cbe-7d84-446c-a934-ad43a82e08dd-1_chapter_5} and~\eqref{ch5:3f712cbe-7d84-446c-a934-ad43a82e08dd-2_chapter_5} have already been established in~\textcite{pmlr-v202-sadiev23a}.

    As in~\textcite{pmlr-v202-sadiev23a}, {\color{navy}let us define the indicator variables
        \begin{align*}
            I & = \mathbb{I}\,\thesischfiveens{\norm{X} > \thesischfivelbd} = \begin{cases}
            1, & \text{if $\norm{X} > \thesischfivelbd$;} \\
            0, & \text{otherwise;}
        \end{cases} \\
            J & = \mathbb{I}\,\thesischfiveens{\norm{X - x} > \thesischfivelbd/2} = \begin{cases}
            1, & \text{if $\norm{X - x} > \thesischfivelbd/2$;} \\
            0, & \text{otherwise.}
        \end{cases}
        \end{align*}}
    then, observe that
        \begin{align*}
            \color{navy}\thesischfiveclip(X, \thesischfivelbd) & \color{navy}= \min\thesischfiveens{1, \frac{\thesischfivelbd}{\norm{X}}}X = \frac{\thesischfivelbd X}{\norm{X}} I + X (1 - I) \\
            & \color{navy}= X + I \left( \frac{\thesischfivelbd}{\norm{X}} - 1 \right) X,
        \end{align*}
    and {\color{navy}$I \le J$} since $\norm{X} > \thesischfivelbd$ implies 
        \begin{thesischfivefitdisplay} \norm{X - x} \ge \norm{X} - \norm{x} > \thesischfivelbd - \norm{x} \ge \frac{\thesischfivelbd}{2}. \end{thesischfivefitdisplay}

    Moreover, thanks to Markov's inequality, we have
        \begin{thesischfivefitalignatnum}{2}
            {\color{navy}\thesischfiveE{J}} & {\color{navy}= \thesischfiveProba{\norm{X - x} > \thesischfivelbd/2}} \\
            & = \thesischfiveProba{\norm{X - x}^p > \left( \frac{\thesischfivelbd}{2} \right)^p} \\
            & \le \left( \frac{2}{\thesischfivelbd} \right)^p \thesischfiveE{\norm{X - x}^p} \\
            & \le \left( \frac{2}{\thesischfivelbd} \right)^p \sigma_1^p, \label{ch5:a8e6cfe7-3338-4dff-9a19-431faa006da8_chapter_5}
        \end{thesischfivefitalignatnum}
    where the first inequality is Markov's inequality (Lemma~\ref{ch5:lem:markov-inequality_chapter_5}).

    \paragraph{Proof of~\eqref{ch5:3f712cbe-7d84-446c-a934-ad43a82e08dd-3_chapter_5}:} using $\norm{x} \le \thesischfivelbd/2$ and $\norm{\widetilde{X} - x} \le \norm{\widetilde{X}} + \norm{x} \le 3 \thesischfivelbd/2$ we have for any $q \in \thesischfiveintff{p}{2}$
        \begin{thesischfivefitalignat}{2}
            \thesischfiveE{\norm{\widetilde{X} - x}^q} & = \thesischfiveE{\norm{\widetilde{X} - x}^p \cdot \norm{\widetilde{X} - x}^{q - p}} \\
            & \le \left( \frac{3 \thesischfivelbd}{2} \right)^{q - p} \thesischfiveE{\norm{\widetilde{X} - x}^p I + \norm{\widetilde{X} - x}^p (1 - I)} \\
            & = \left( \frac{3 \thesischfivelbd}{2} \right)^{q - p} \thesischfiveE{\norm{\frac{\thesischfivelbd X}{\norm{X}} - x}^p I + \norm{X - x}^p (1 - I)} \\
            & \le \left( \frac{3 \thesischfivelbd}{2} \right)^{q - p} \thesischfiveE{\left( \norm{\frac{\thesischfivelbd X}{\norm{X}}} + \norm{x} \right)^p I + \norm{X - x}^p (1 - I)} \\
            & \le \left( \frac{3 \thesischfivelbd}{2} \right)^{q - p} \thesischfiveE{\left( \frac{3 \thesischfivelbd}{2} \right)^p I + \norm{X - x}^p (1 - I)} \\
            & \le \left( \frac{3 \thesischfivelbd}{2} \right)^q \thesischfiveE{I} + \left( \frac{3 \thesischfivelbd}{2} \right)^{q - p} \sigma_1^p \\
            \thesischfiveoversetlab{\eqref{ch5:a8e6cfe7-3338-4dff-9a19-431faa006da8_chapter_5}}&{\le} \left[ \left( \frac{3 \thesischfivelbd}{2} \right)^q \left( \frac{2}{\thesischfivelbd} \right)^p + \left( \frac{3 \thesischfivelbd}{2} \right)^{q - p} \right] \sigma_1^p \\
            & \le 2 \cdot 3^q \thesischfivelbd^{q - p} \sigma_1^p,
        \end{thesischfivefitalignat}
    where in the second equality we use the definition of $\thesischfiveclip$ and {\color{navy}$I$}, in the fourth inequality we use {\color{navy}$1 - I \le 1$} and $\thesischfiveE{\norm{X - x}^p} \le \sigma_1^p$, and in the fifth {\color{navy}$I \le J$}.
    
    \paragraph{Proof of~\eqref{ch5:3f712cbe-7d84-446c-a934-ad43a82e08dd-4_chapter_5}:} according to~Lemma~\ref{ch5:lem:variance-decomposition-any-exponent_chapter_5} with $c = x$, we have
        \begin{thesischfivefitalignat}{2}
            {\color{navy}\thesischfiveE{\norm{\widetilde{X} - \thesischfiveE{\widetilde{X}}}^q}} & \le 2^{2 - q} \thesischfiveE{\norm{\widetilde{X} - x}^q} \\
            \thesischfiveoversetlab{\eqref{ch5:3f712cbe-7d84-446c-a934-ad43a82e08dd-3_chapter_5}}&{\le} 2^{3 - q} \cdot 3^q \thesischfivelbd^{q - p} \sigma_1^p,
        \end{thesischfivefitalignat}
    as claimed.
\end{proof}

\newpage

\section{Useful Identities and Inequalities}

For any vectors $x, y \in \R^d$, we have
\begin{thesischfivefitalignatnum}{2}
    2 \thesischfiveps{x}{y} & = \sqnorm{x} + \sqnorm{y} - \sqnorm{x - y}. \label{ch5:eq:identity-1_chapter_5}
\end{thesischfivefitalignatnum}

\begin{lemma}[Variance Decomposition]\label{ch5:lem:variance-decomposition_chapter_5}
    For any random vector $X \in \R^d$ and any non-random vector $c \in \R^d$ we have
        \begin{thesischfivefitdisplay} \thesischfiveE{\sqnorm{X - c}} = \thesischfiveE{\sqnorm{X - \thesischfiveE{X}}} + \sqnorm{\thesischfiveE{X} - c}. \end{thesischfivefitdisplay}
\end{lemma}

\begin{remark}
    In particular, this implies that for any random vector $X \in \R^d$ and any non-random vector $c \in \R^d$ we have
        \begin{thesischfivefitdisplay} \thesischfiveE{\sqnorm{X - \thesischfiveE{X}}} \le \thesischfiveE{\sqnorm{X - c}}. \end{thesischfivefitdisplay}
    In~Lemma~\ref{ch5:lem:variance-decomposition-any-exponent_chapter_5} we extend this inequality to any exponent $\eta > 1$ instead of just $2$.
\end{remark}

\begin{lemma}[Tower Property of the Expectation]\label{ch5:lem:tower-property_chapter_5}
    For any random variables $X \in \R^d$ and $Y_1, \ldots, Y_n$ we have
        \begin{thesischfivefitdisplay} \thesischfiveE{\thesischfiveExpCond{X}{Y_1, \ldots, Y_n}} = \thesischfiveE{X}. \end{thesischfivefitdisplay}
\end{lemma}

\begin{lemma}[Cauchy Schwarz's Inequality]\label{ch5:lem:cauchy-schwarz_chapter_5}
    For any vectors $a, b \in \R^d$ we have
        \begin{thesischfivefitdisplay} \thesischfiveps{a}{b} \le \abs{\thesischfiveps{a}{b}} \le \norm{a} \norm{b}. \end{thesischfivefitdisplay}
\end{lemma}

\begin{lemma}[Young's Inequality (Norm Form)]\label{ch5:lem:youg-inequality_chapter_5}
    For any vectors $a, b \in \R^d$ and any scalar $\eta > 0$ we have
        \begin{thesischfivefitdisplay} {\color{navy}\sqnorm{a + b} \le (1 + \eta) \sqnorm{a} + \left( 1 + \frac{1}{\eta} \right) \sqnorm{b}.} \end{thesischfivefitdisplay}
\end{lemma}

\begin{lemma}[Young's Inequality (Inner Product Form)]\label{ch5:lem:youg-inequality-inner-product_chapter_5}
    For any vectors $a, b \in \R^d$ and any scalar $\eta > 0$ we have
        \begin{thesischfivefitdisplay} {\color{navy}2\thesischfiveps{a}{b} \le 2 \abs{\thesischfiveps{a}{b}} \le \eta \sqnorm{a} + \frac{1}{\eta} \sqnorm{b}.} \thesischfivenumberthis\label{ch5:55c90e49-73cb-4448-8f5b-fb2c30a7c439_chapter_5} \end{thesischfivefitdisplay}
\end{lemma}

\begin{proof}
  It's enough to prove inequality~\eqref{ch5:55c90e49-73cb-4448-8f5b-fb2c30a7c439_chapter_5} when $d = 1$. Hence, consider $a, b \in \R$, we have given $\eta > 0$
    \begin{thesischfivefitdisplay} 2ab \le 2\abs{ab} = 2\abs{a} \cdot \abs{b} = 2\abs{\sqrt{\eta} \, a} \cdot \abs{\frac{b}{\sqrt{\eta}}} \thesischfiveoversetrel{rel:c481fc03-dcbc-4a9b-b12e-1b8282d05b5d}{\le} \eta \abs{a}^2 + \frac{1}{\eta} \abs{b}^2 = \eta \, a^2 + \frac{b^2}{\eta}, \end{thesischfivefitdisplay}
  where in~\thesischfiverelref{rel:c481fc03-dcbc-4a9b-b12e-1b8282d05b5d} we use the arithmetic-geometric inequality in $n = 2$ variables $\left( \sqrt{\eta} \abs{a}, \frac{1}{\sqrt{\eta}} \abs{b} \right)$.
\end{proof}

\begin{lemma}\label{ch5:lem:norm-power-alpha-inequality_chapter_5}
    For any vectors $a, b \in \R^d$ and any scalar $\eta \ge 1$ we have
        \begin{thesischfivefitdisplay} \norm{a + b}^{\eta} \le 2^{\eta - 1} \left( \norm{a}^{\eta} + \norm{b}^{\eta} \right). \thesischfivenumberthis\label{ch5:e4d36f60-854e-4737-9f49-0ac946fb5b07_chapter_5} \end{thesischfivefitdisplay}
\end{lemma}

\begin{proof}
    Note that the function $x \mapsto \norm{x}^{\eta}$ is convex over $\R^d$ since $\norm{\cdot}$ is convex (as a norm) and $t \mapsto t^{\eta}$ is convex since $\eta \ge 1$. Hence, by the Jensen's inequality (Lemma~\ref{ch5:lem:jensen-inequality_chapter_5}), we have
        \begin{thesischfivefitdisplay} \norm{\frac{a + b}{2}}^{\eta} \le \frac{1}{2} \left( \norm{a}^{\eta} + \norm{b}^{\eta} \right), \end{thesischfivefitdisplay}
    i.e.,
        \begin{thesischfivefitdisplay} \norm{a + b}^{\eta} \le 2^{\eta - 1} \left( \norm{a}^{\eta} + \norm{b}^{\eta} \right), \end{thesischfivefitdisplay}
    as claimed.
\end{proof}

\begin{remark}
    Similar inequalities to~\eqref{ch5:e4d36f60-854e-4737-9f49-0ac946fb5b07_chapter_5} appear in~\textcite{he2025complexitynormalizedstochasticfirstorder}. In~Lemma~\ref{ch5:lem:euclidean-norm-expansion-generalized_chapter_5} we state one of them that we use in our proofs.
\end{remark}

\begin{lemma}[{\textcite[Lemma~1]{he2025complexitynormalizedstochasticfirstorder}}]\label{ch5:lem:euclidean-norm-expansion-generalized_chapter_5}
     For any vectors $a, b \in \R^d$, $a \neq \boldsymbol{0}$, and any scalar $\eta\in \thesischfiveintof{1}{2}$ we have\footnote{Note that, when $a = 0$ we have $\norm{a}^{\eta - 2} a = 0$ since $\eta > 1$ thus $\eta - 2 > -1$.}
        \begin{thesischfivefitdisplay} \norm{a + b}^{\eta} \le \norm{a}^{\eta} + \eta \norm{a}^{\eta - 2} \thesischfiveps{a}{b} + 2^{2 - \eta} \norm{b}^{\eta}. \end{thesischfivefitdisplay}
\end{lemma}

\begin{lemma}\label{ch5:lem:variance-decomposition-any-exponent_chapter_5}
    For any scalar $\eta \in \thesischfiveintof{1}{2}$, any random vector $X \in \R^d$ and any non-random vector $c \in \R^d$ we have
        \begin{thesischfivefitdisplay} \thesischfiveE{\norm{X - \thesischfiveE{X}}^{\eta}} \le 2^{2 - \eta} \thesischfiveE{\norm{X - c}^{\eta}}. \end{thesischfivefitdisplay}
\end{lemma}

\begin{proof}
    {\color{navy}If $c=\thesischfiveE{X}$, the claimed inequality is immediate. Assume henceforth that $c\ne\thesischfiveE{X}$.} By~Lemma~\ref{ch5:lem:euclidean-norm-expansion-generalized_chapter_5} we have
        \begin{thesischfivefitalignat}{2}
            \MoveEqLeft[1] \thesischfiveE{\norm{X - \thesischfiveE{X}}^{\eta}} \\
            & = \thesischfiveE{\norm{\left( X - c \right) + \left( c - \thesischfiveE{X} \right)}^{\eta}} \\
            & \le \norm{c - \thesischfiveE{X}}^{\eta} + \eta \norm{c - \thesischfiveE{X}}^{\eta - 2} \thesischfiveE{\thesischfiveps{c - \thesischfiveE{X}}{X - c}} \\
            &\quad + 2^{2 - \eta} \thesischfiveE{\norm{X - c}^{\eta}} \\
            & = \norm{c - \thesischfiveE{X}}^{\eta} - \eta \norm{c - \thesischfiveE{X}}^{\eta} + 2^{2 - \eta} \thesischfiveE{\norm{X - c}^{\eta}} \\
            & \le 2^{2 - \eta} \thesischfiveE{\norm{X - c}^{\eta}},
        \end{thesischfivefitalignat}
    where the first inequality is Lemma~\ref{ch5:lem:euclidean-norm-expansion-generalized_chapter_5}, the second equality uses $\thesischfiveE{X - c} = \thesischfiveE{X} - c$, and the last inequality uses $\eta > 1$.
\end{proof}

\begin{lemma}[Jensen's Inequality]\label{ch5:lem:jensen-inequality_chapter_5}
    Let $f \colon \R^d \to \R$ be a convex function then
    \begin{enumerate}
        \item (Probabilistic Form) for any random vector $X \in \R^d$ we have
            $$\thesischfiveE{f(X)} \ge f\left(\thesischfiveE{X}\right).$$

        \item (Deterministic Form) for any vectors $v_1, \ldots, v_n \in \R^d$ and scalars $\thesischfivelbd_1, \ldots, \thesischfivelbd_n \in \R_+$ we have
            $$\sum_{i = 1}^n \thesischfivelbd_i f(v_i) \ge f\left( \sum_{i = 1}^n \thesischfivelbd_i v_i \right),$$
        provided $\thesischfivelbd_i \ge 0$ for all $i \in [n]$ and $\sum\limits_{i = 1}^n \thesischfivelbd_i = 1$.
    \end{enumerate}
\end{lemma}

\begin{lemma}\label{ch5:lem:jensen-form-1_chapter_5}
    For any vectors $v_1, \ldots, v_n \in \R^d$ we have
        $$\sqnorm{\sum_{i = 1}^n v_i} \leq n \sum_{i = 1}^n \sqnorm{v_i}.$$
\end{lemma}

\begin{proof}
    The function $\sqnorm{\cdot} \colon \R^d \to \R$ is $\mu$-strongly convex with $\mu = 2$ so is convex thus applying Jensen's inequality~\ref{ch5:lem:jensen-inequality_chapter_5} with $\thesischfivelbd_1 = \cdots = \thesischfivelbd_n = \frac{1}{n}$ gives
        $$\sqnorm{\sum_{i = 1}^n \frac{v_i}{n}} \le \frac{1}{n} \sum_{i = 1}^n \sqnorm{v_i},$$
    and multiplying both sides by $n^2$ gives the desired inequality.
\end{proof}

\begin{lemma}[Markov's Inequality]\label{ch5:lem:markov-inequality_chapter_5}
    For any non-negative random variable $X$ and any scalar $a > 0$, we have
        \begin{thesischfivefitdisplay} \thesischfiveProba{X \ge a} \le \frac{\thesischfiveE{X}}{a}. \end{thesischfivefitdisplay}
\end{lemma}

\begin{lemma}[Bernoulli's Inequality]\label{ch5:lem:bernoulli-inequality_chapter_5}
    For any real number $x \ge -1$ and $r \in \thesischfiveens{0} \cup \thesischfiveintfo{1}{+\infty}$ we have
        \begin{thesischfivefitdisplay} (1 + x)^r \ge 1 + r x. \end{thesischfivefitdisplay}
\end{lemma}

\begin{lemma}[{Absolute $p^{\textnormal{th}}$ Central Moment of Normal Distribution~\textcite[eq.~(18)]{winkelbauer2014momentsabsolutemomentsnormal}}]\label{ch5:appdx-lem:absolute-central-moment-gaussian_chapter_5}
    For any real number $p > -1$, $\sigma_1 > 0$ and $\mu \in \R$, if $X \sim \mathcal{N}(\mu, \sigma_1^2)$ then
        \begin{thesischfivefitdisplay} \thesischfiveE{\abs{X - \mu}^p} = 2^{\frac{p}{2}} \sigma_1^p \cdot \frac{\Gamma\left( \frac{p + 1}{2} \right)}{\sqrt{\pi}}, \end{thesischfivefitdisplay}
    where $\Gamma$ is the Euler's gamma function (see~\textcite[{(5.2.1)}]{NIST:DLMF} for formal definition).
\end{lemma}

\begin{definition}[Operator Norm]\label{ch5:appdx-def:operator-norm_chapter_5}
    Given a $d \times d$ matrix $A \in \R^{d \times d}$, the operator norm of $A$ is the norm
        \begin{thesischfivefitdisplay} \normop{A} \eqdef \sup_{y \in \R^d, \norm{y} = 1} \norm{A y}, \end{thesischfivefitdisplay}
    where $\norm{\cdot}$ is the standard euclidean norm in $\R^d$.
\end{definition}

\begin{lemma}[{\color{navy}Lipschitz Gradients Imply a Bounded Hessian}]\label{ch5:appdx-lem:lipschitz-gradients-implies-bounded-hessian_chapter_5}
    Given $f \colon \R^d \to \R$ a twice continuously differentiable function over $\R^d$ such that its gradients are $L_1$--Lipschitz continuous for some constant $L_1 \ge 0$ then, for all $x \in \R^d$ we have
        \begin{thesischfivefitdisplay} \normop{\nabla^2 f(x)} \le L_1. \end{thesischfivefitdisplay}
\end{lemma}
{\color{navy}The preceding lemma extends to an $m$-times continuously differentiable function, $m\ge2$, whose $(m-1)$st derivative is Lipschitz continuous in the operator norm induced by the Euclidean norm.}

\begin{lemma}\label{ch5:appdx-lem:bound-norm-gradient-by-L-and-Delta_chapter_5}
    Given $f \colon \R^d \to \R$ a continuously differentiable and lower bounded function over $\R^d$ such that its gradients are $L_1$--Lipschitz continuous for some constant $L_1 \ge 0$ then, for all $x \in \R^d$ we have
        \begin{thesischfivefitdisplay} \norm{\nabla f(x)}^2 \le 2 L_1 \left( f(x) - f^{\inf} \right), \end{thesischfivefitdisplay}
    for any $f^{\inf} \le \inf_{x \in \R^d} f(x)$.
\end{lemma}

\begin{lemma}\label{ch5:appdx-lem:bound-taylor-gradient-hessian_chapter_5}
    Given $f \colon \R^d \to \R$ a twice continuously differentiable and lower bounded function over $\R^d$ such that its Hessians are $L_2$--Lipschitz continuous for some constant $L_2 \ge 0$ then, for all $x, y \in \R^d$ we have
        \begin{thesischfivefitdisplay} \norm{\nabla f(x) - \nabla f(y) - \nabla^2 f(y) (x - y)} \le \frac{L_2}{2} \sqnorm{x - y}. \end{thesischfivefitdisplay}
\end{lemma}

\begin{lemma}[{Landau-Kolmogorov's Inequalities on a Finite Interval~\textcite[Theorem~2]{Fabry_1987}}]\label{ch5:appdx-lem:landau-kolmogorov-inequalities_chapter_5}
    {\color{navy}Let $m\ge1$ be an integer, let $a\le b$ be real numbers, and let $f\in\mathcal C^m(\thesischfiveintff{a}{b},\R)$. For each integer $0 \le k \le m$, define $M_k\eqdef\sup_{x\in\thesischfiveintff{a}{b}}\abs{f^{(k)}(x)}$. If $M_0$ and $M_m$ are finite, then, for every $k\in[m]$,}
        \begin{enumerate}
            \item {\color{navy}$M_k$ is finite;}
            \item {\color{navy}the inequality}
                \begin{thesischfivefitdisplay} {\color{navy}M_k \le c_{k, m} M_0^{(m - k)/m} \max\thesischfiveens{\frac{M_m}{2^{m - 1} m!}, 4^m (b - a)^{-m} M_0}^{k/m},} \thesischfivenumberthis\label{ch5:ce6d9b55-0646-4949-b2b3-2311a3f3f37b_chapter_5} \end{thesischfivefitdisplay}
            {\color{navy}holds, where $c_{k,m}>0$ is the universal constant}
                \begin{thesischfivefitdisplay} {\color{navy}c_{k,m}\eqdef\frac{m2^k k!}{m+k}\binom{m+k}{m-k}.} \end{thesischfivefitdisplay}
        \end{enumerate}
\end{lemma}

\begin{remark}
    {\color{navy}The constant $c_{k,m}$ in~\eqref{ch5:ce6d9b55-0646-4949-b2b3-2311a3f3f37b_chapter_5} is independent of $f$, whereas the quantities $\thesischfiveens{M_k}_{k\in\thesischfiveInt{0}{m}}$ depend on $f$.}
\end{remark}

\begin{remark}
    Other related inequalities can be found in~\textcite{Mitrinović1991,Chen_1993}.
\end{remark}

\begin{lemma}\label{ch5:appdx-lem:lower-bounding-F-steepest-descent-direction_chapter_5}
    Let $f \colon \R^d \to \R$ be a continuously differentiable function over $\R^d$ with $L$--Lipschitz gradients for some $L > 0$ and $x^{0} \in \R^d$. Assume that $\norm{\nabla f(x^{0})} > 0$ then for all $t \in \thesischfiveintff{0}{2}$ we have
        \begin{thesischfivefitdisplay} f(x^{0} + t v) \le f(x^{0}), \end{thesischfivefitdisplay}
    where $v = - \nabla f(x^{0})/L$.
\end{lemma}

\begin{proof}
    Using the fact that the function $f$ has $L$--Lipschitz gradients then we know $f$ is $L$--smooth \parencite{nesterov2018lectures}; hence, for all $t \in \thesischfiveintff{0}{2}$ we have
        \begin{thesischfivefitalignat}{2}
            D_f(x^{0} + t v, x^{0}) & \eqdef f(x^{0} + t v) - f(x^{0}) - \thesischfiveps{\nabla f(x^{0})}{(x^{0} + t v) - x^{0}} \\
            & \le \frac{L}{2} \sqnorm{(x^{0} + t v) - x^{0}} = \frac{L t^2}{2} \sqnorm{v},
        \end{thesischfivefitalignat}
    where $D_f(x, y)$ denotes the Bregman divergence of $f$ at $x, y \in \R^d$. Rewriting the above inequality using the choice of $v$ gives
        \begin{thesischfivefitalignat}{2}
            f(x^{0} + t v) & \le f(x^{0}) + t \thesischfiveps{\nabla f(x^{0})}{v} + \frac{L t^2}{2} \sqnorm{v} \\
            & = f(x^{0}) - \frac{t}{L} \sqnorm{\nabla f(x^{0})} + \frac{t^2}{2 L} \sqnorm{\nabla f(x^{0})} \\
            & = f(x^{0}) - \frac{t}{L} \sqnorm{\nabla f(x^{0})} \left( 1 - \frac{t}{2} \right) \\
            & \le f(x^{0}),
        \end{thesischfivefitalignat}
    since $0 \le t \le 2$ and $L > 0$. {\color{navy}This completes the proof of the lemma.}
\end{proof}

\begin{remark}
    {\color{navy}The preceding lemma also covers $\nabla f(x^0)=0$, in which case $v=0$ and the conclusion is immediate.}
\end{remark}

\begin{lemma}
    {\color{navy}Let $f\colon\R^d\to\R$ be $m+1$ times continuously differentiable for an integer $m\ge1$. Assume that}
        \begin{enumerate}
            \item {\color{navy}for every $k\in[m]$, the function $\nabla^k f$ is $L_k$--Lipschitz continuous for some $L_k>0$;}
            \item the function $f$ is lower bounded over $\R^d$ and we denote $f^{\inf} \eqdef \inf_{x \in \R^d} f(x)$,
        \end{enumerate}
    {\color{navy}Then, for every $k\in[m]$, the constant $c_k\eqdef(k+1)^2$ satisfies, for all $x^0\in\R^d$,}
        \begin{thesischfivefitalignatnum}{2}
            \MoveEqLeft[1] \norm{\nabla f(x^{0})} \\
            & \le c_k \left( f(x^{0}) - f^{\inf} \right)^{k/(k + 1)} \\
            &\qquad \cdot \max\thesischfiveens{\frac{L_k^{1/(k + 1)}}{2^{k/(k + 1)} [(k + 1)!]^{1/(k + 1)}}, \frac{2 L_1 \left( f(x^{0}) - f^{\inf} \right)^{1/(k + 1)}}{\norm{\nabla f(x^{0})}}}. \label{ch5:e7a373a5-ff82-4f5f-9ece-a932573b577d_chapter_5}
        \end{thesischfivefitalignatnum}
    {\color{navy}When $\nabla f(x^0)=0$, the fraction containing $\norm{\nabla f(x^0)}$ in the denominator is interpreted as $+\infty$.}

    In particular, if the Lipschitz constant $L_k$ dominates and the gradient of $f$ at $x^{0}$ is large, for instance of the order of $\sqrt{L_1 ( f(x^{0}) - f^{\inf} )}$ (see~Lemma~\ref{ch5:appdx-lem:bound-norm-gradient-by-L-and-Delta_chapter_5}), then~\eqref{ch5:e7a373a5-ff82-4f5f-9ece-a932573b577d_chapter_5} simplifies to
        \begin{thesischfivefitdisplay} \norm{\nabla f(x^{0})} \le \widetilde{c}_k L_k^{1/(k + 1)} \left( f(x^{0}) - f^{\inf} \right)^{k/(k + 1)}, \thesischfivenumberthis\label{ch5:5ea0787a-1b21-4f89-9e4e-032978f6ec42_chapter_5} \end{thesischfivefitdisplay}
    for some universal constant $\widetilde{c}_k$ which depends only on $k$.
\end{lemma}

\begin{remark}
    The bound~\eqref{ch5:e7a373a5-ff82-4f5f-9ece-a932573b577d_chapter_5} is (almost) a generalization to high-order Lipschitz constant of the well-known inequality
        \begin{thesischfivefitdisplay} \norm{\nabla f(x^{0})} \le \sqrt{2 L_1 \left( f(x^{0}) - f^{\inf} \right)}, \end{thesischfivefitdisplay}
    which we recall in~Lemma~\ref{ch5:appdx-lem:bound-norm-gradient-by-L-and-Delta_chapter_5} and which corresponds to the case $k = 1$.

    {\color{navy}In full generality, one cannot obtain~\eqref{ch5:5ea0787a-1b21-4f89-9e4e-032978f6ec42_chapter_5} without an additional dominance assumption.}\footnote{Inequality~\eqref{ch5:5ea0787a-1b21-4f89-9e4e-032978f6ec42_chapter_5} already fails for the case $k = 2$ with quadratic functions. For instance, consider $f \colon x \mapsto \frac{1}{2} \norm{x}^2$ over $\R^d$ and let $x^{0} = (1, 0, \ldots, 0)$ then $\norm{\nabla f(x^{0})} = \norm{x^{0}} = 1$, but $\nabla^2 f(x) = \mathrm{Id}$ for all $x \in \R^d$ hence we can take $L_2 = 0$ since the Hessian of $f$ is constant, but then
        \begin{thesischfivefitdisplay} L_2^{1/3} \left( f(x^{0}) - f^{\inf} \right)^{2/3} = 0, \end{thesischfivefitdisplay}
    and there is no universal constant $c_2$ for which~\eqref{ch5:5ea0787a-1b21-4f89-9e4e-032978f6ec42_chapter_5} can hold.} {\color{navy}Accordingly, the simplified bound should be invoked only after verifying the stated dominance condition.}
\end{remark}

\begin{proof}
    Let $f \colon \R^d \to \R$ be defined as above. Then, if $\norm{\nabla f(x^{0})} = 0$ the inequality~\eqref{ch5:e7a373a5-ff82-4f5f-9ece-a932573b577d_chapter_5} holds since $c_k > 0$ and $f(x^{0}) - f^{\inf} \ge 0$. Now, assume $\norm{\nabla f(x^{0})} > 0$, we define the scalar function $\varphi \colon \thesischfiveintff{0}{2} \to \R$ as
        \begin{thesischfivefitdisplay} \varphi \colon t \mapsto f(x^{0} + t v) - f^{\inf}, \end{thesischfivefitdisplay}
    where $v = -\nabla f(x^{0})/L_1 \neq 0$ and $L_1$ is the Lipschitz constant of $\nabla f$. {\color{navy}Then $\varphi$ is $m+1$ times continuously differentiable on $\thesischfiveintff{0}{2}$, and, for every $t\in\thesischfiveintff{0}{2}$ and $k\in\thesischfiveInt{2}{m+1}$,}
        \begin{thesischfivefitdisplay} \varphi^{(k)}(t) = \nabla^k f(x^{0} + t v) [v, \ldots, v] = \frac{(-1)^k}{L_1^k} \nabla^k f(x^{0} + t v) [g, \ldots, g], \thesischfivenumberthis\label{ch5:441999ae-9106-410c-a85b-fdcd81183237_chapter_5} \end{thesischfivefitdisplay}
    where we let $g \eqdef \nabla f(x^{0})$ and $\nabla^k f(x^{0} + t v) [\cdot, \ldots, \cdot]$ is the $k$--linear form induced by the $k^{\textnormal{th}}$--derivative of $f$. From equality~\eqref{ch5:441999ae-9106-410c-a85b-fdcd81183237_chapter_5} we obtain
        \begin{thesischfivefitalignatnum}{2}
            \abs{\varphi^{(k)}(t)} & = \frac{1}{L_1^k} \abs{\nabla^k f(x^{0} + t v) [g, \ldots, g]} \\
            & \le \frac{1}{L_1^k} \normop{\nabla^k f(x^{0} + t v)} \cdot \norm{g}^k \\
            & \le \left( \frac{\norm{\nabla f(x^{0})}}{L_1} \right)^k L_{k - 1}, \label{ch5:7d81208b-8c62-409a-af2f-97a7bcb40436_chapter_5}
        \end{thesischfivefitalignatnum}
    {\color{navy}In the first inequality we use the Cauchy--Schwarz inequality (Lemma~\ref{ch5:lem:cauchy-schwarz_chapter_5}). In the second, the $L_{k-1}$--Lipschitz continuity of $\nabla^{k-1}f$ implies $\normop{\nabla^k f}\le L_{k-1}$; Lemma~\ref{ch5:appdx-lem:lipschitz-gradients-implies-bounded-hessian_chapter_5} is the case $k=2$.} Hence, if we define $M_k \eqdef \sup_{t \in \thesischfiveintff{0}{2}} \abs{\varphi^{(k)} (t)}$ we have, by~\eqref{ch5:7d81208b-8c62-409a-af2f-97a7bcb40436_chapter_5},
        \begin{thesischfivefitdisplay} M_k \le \left( \frac{\norm{\nabla f(x^{0})}}{L_1} \right)^k L_{k - 1} < +\infty, \thesischfivenumberthis\label{ch5:e0aab2e4-82dd-4824-b2b2-df7ad67d99b1_chapter_5} \end{thesischfivefitdisplay}
    {\color{navy}This quantity is finite for every $k\in\thesischfiveInt{2}{m+1}$. Moreover,}
        \begin{thesischfivefitalignatnum}{2}
            M_1 & = \sup_{t \in \thesischfiveintff{0}{2}} \abs{\varphi'(t)} = \sup_{t \in \thesischfiveintff{0}{2}} \abs{\thesischfiveps{\nabla f(x^{0} + t v)}{v}} \\
            & = \frac{1}{L_1} \sup_{t \in \thesischfiveintff{0}{2}} \abs{\thesischfiveps{\nabla f(x^{0} + t v)}{\nabla f(x^{0})}} \\
            & \ge \frac{1}{L_1} \sqnorm{\nabla f(x^{0})}. \label{ch5:a718a294-777a-4426-9f05-5289df231fe8_chapter_5}
        \end{thesischfivefitalignatnum}
    Additionally,
        \begin{thesischfivefitalignatnum}{2}
            M_0 & \eqdef \sup_{t \in \thesischfiveintff{0}{2}} \abs{\varphi(t)} \\
            & = \sup_{t \in \thesischfiveintff{0}{2}} \left( f(x^{0} + t v) - f^{\inf} \right) \\
            & = f(x^{0}) - f^{\inf} < +\infty, \label{ch5:6c704741-f58f-49c3-be4f-986e40421895_chapter_5}
        \end{thesischfivefitalignatnum}
    where the first equality uses $f \ge f^{\inf}$ and the second uses Lemma~\ref{ch5:appdx-lem:lower-bounding-F-steepest-descent-direction_chapter_5}, since $f$ is continuously differentiable with $L_1$--Lipschitz gradients.

    {\color{navy}Fix any $k\in[m]$. The function $\varphi$ is $k+1$ times continuously differentiable on $\thesischfiveintff{0}{2}$, and the quantities $\thesischfiveens{M_{\ell}}_{\ell\in\thesischfiveInt{0}{k+1}}$ are finite by~\eqref{ch5:e0aab2e4-82dd-4824-b2b2-df7ad67d99b1_chapter_5} and~\eqref{ch5:6c704741-f58f-49c3-be4f-986e40421895_chapter_5}. Therefore, the Landau--Kolmogorov inequality (Lemma~\ref{ch5:appdx-lem:landau-kolmogorov-inequalities_chapter_5}), applied with order $k+1$, gives}
        \begin{thesischfivefitdisplay} M_1 \le c_k \, M_0^{k/(k + 1)} \max\thesischfiveens{\frac{M_{k + 1}}{2^k (k + 1)!}, 2^{k + 1} M_0}^{1/(k + 1)}, \thesischfivenumberthis\label{ch5:80ccea34-41b0-4eb4-9743-22232a23d1f2_chapter_5} \end{thesischfivefitdisplay}
    where $c_k \eqdef c_{1, k + 1} = (k + 1)^2$ is {\color{navy}a universal constant} depending only on $k$ (by~Lemma~\ref{ch5:appdx-lem:landau-kolmogorov-inequalities_chapter_5}). From the inequalities~\eqref{ch5:e0aab2e4-82dd-4824-b2b2-df7ad67d99b1_chapter_5} and~\eqref{ch5:a718a294-777a-4426-9f05-5289df231fe8_chapter_5},~\eqref{ch5:6c704741-f58f-49c3-be4f-986e40421895_chapter_5} along with~\eqref{ch5:80ccea34-41b0-4eb4-9743-22232a23d1f2_chapter_5} we obtain
        \begin{thesischfivefitalignat}{2}
            \MoveEqLeft[1] \frac{1}{L_1} \sqnorm{\nabla f(x^{0})} \\
            & \le c_k \, \left( f(x^{0}) - f^{\inf} \right)^{k/(k + 1)} \\
            &\qquad \cdot \max\thesischfiveens{\left( \left( \frac{\norm{\nabla f(x^{0})}}{L_1} \right)^{k + 1} \frac{L_k}{2^k (k + 1)!} \right)^{1/(k + 1)}, 2 M_0^{1/(k + 1)}} \\
            & = \frac{c_k}{L_1} \left( f(x^{0}) - f^{\inf} \right)^{k/(k + 1)} \norm{\nabla f(x^{0})} \\
            &\qquad \cdot \max\thesischfiveens{\frac{L_k^{1/(k + 1)}}{2^{k/(k + 1)} [(k + 1)!]^{1/(k + 1)}}, \frac{2 L_1 \left( f(x^{0}) - f^{\inf} \right)^{1/(k + 1)}}{\norm{\nabla f(x^{0})}}},
        \end{thesischfivefitalignat}
    which simplifies to 
        \begin{align*}
            \norm{\nabla f(x^{0})} & \le c_k \left( f(x^{0}) - f^{\inf} \right)^{k/(k + 1)} \\
            &\qquad \cdot \max\thesischfiveens{\frac{L_k^{1/(k + 1)}}{2^{k/(k + 1)} [(k + 1)!]^{1/(k + 1)}}, \frac{2 L_1 \left( f(x^{0}) - f^{\inf} \right)^{1/(k + 1)}}{\norm{\nabla f(x^{0})}}},
        \end{align*}
    which leads to the desired inequality.
\end{proof}

\section{\texorpdfstring{{\color{navy}Additional Analysis in the Case $p=q=2$}}{Additional Analysis in the Case p=q=2}}

In this section, which is of independent interest, we study regular, not normalized, \algname{SGD} with Momentum Variance Reduction (\algname{SGD-MVR}). This method is the original one by \textcite{cutkosky2019momentum}, where it is called \algname{STORM}. The method is shown as Algorithm~\ref{ch5:algo:sgd-mvr_chapter_5}. We show that in the case $p=q=2$, it achieves the same optimal complexity as \algname{NSGD-MVR}. %
{\color{navy}That is, normalization is not required in this particular case.}

  \begin{algorithm}[htbp]
\caption{\algname{SGD-MVR} (\algname{SGD} with \texttt{MVR})}
\label{ch5:algo:sgd-mvr_chapter_5}
\begin{algorithmic}[1]
\Require $x^0 \in \R^d$, the starting point;              $K > 0$, the number of iterations;              $g^0 \in \R^d$, an initial vector;              $\gamma > 0$, the stepsize;              $\eta \in \thesischfiveintof{0}{1}$, the momentum parameter for \texttt{MVR};
\State $x^1 = x^0 - \gamma g^0$
\For{$k = 1, 2, \ldots, K - 1$}
\State \textit{Apply \texttt{MVR}.}
\State $g^{k} = (1 - \eta) \left( g^{k-1} + \nabla f\left(x^{k};\xi^{k}\right) - \nabla f\left(x^{k-1};\xi^{k}\right) \right) + \eta \nabla f\left(x^{k};\xi^{k}\right)$
\State \textit{Do one descent step.}
\State $x^{k+1} = x^{k} - \gamma g^{k}$
\EndFor
\Ensure ${\color{navy}x_{\mathrm{out}}\sim\mathcal U\{x^0,\ldots,x^{K-1}\}}$
\end{algorithmic}
\end{algorithm}%

\begin{theorem}\label{ch5:thm:sgd-mvr-convergence-analysis_chapter_5}
    {\color{navy}Under~Assumption~\ref{ch5:ass:lower-boundedness_chapter_5},
    Assumption~\ref{ch5:ass:L-lipschitz-gradients_chapter_5},
    Assumption~\ref{ch5:ass:p-bounded-central-moment-gradient_chapter_5},
    and Assumption~\ref{ch5:ass:mean-squared-smoothness-2_chapter_5},
    with $p=q=2$, let $\Delta\eqdef f(x^0)-f^{\inf}$ and
    $\thesischfiveeps>0$. With the convention that $\eta=1$ when
    $\sigma_1=0$, choose
    \begin{thesischfivefitdisplay}
        \eta=\min\thesischfiveens{1,\frac{\thesischfiveeps^2}{6\sigma_1^2}},
        \qquad
        B_{\mathrm{init}}=\left\lceil\eta^{-1}\right\rceil,
        \qquad
        g^0=\frac1{B_{\mathrm{init}}}
        \sum_{j=0}^{B_{\mathrm{init}}-1}\nabla f(x^0;\xi_j^0),
    \end{thesischfivefitdisplay}
    and
    \begin{thesischfivefitdisplay}
        \gamma=\frac{1}{
        L_1+\frac{\sqrt2\,\delta(1-\eta)}{\sqrt\eta}},
        \qquad
        K\ge
        \left\lceil
        \max\thesischfiveens{
        1,\frac{3\sigma_1^2}{\thesischfiveeps^2},
        \frac{6\Delta}{\gamma\thesischfiveeps^2}}
        \right\rceil .
    \end{thesischfivefitdisplay}
    Then Algorithm~\ref{ch5:algo:sgd-mvr_chapter_5} satisfies
    \begin{thesischfivefitdisplay}
        \Exp{\norm{\nabla f(x_{\mathrm{out}})}^2}
        =\frac1K\sum_{k=0}^{K-1}
        \Exp{\norm{\nabla f(x^k)}^2}
        \le\thesischfiveeps^2.
    \end{thesischfivefitdisplay}
    Consequently,
    $\Exp{\norm{\nabla f(x_{\mathrm{out}})}}\le\thesischfiveeps$, and the
    total number of stochastic-gradient evaluations is
    \begin{thesischfivefitdisplay}
        \cO\left(
        1+\frac{\sigma_1^2}{\thesischfiveeps^2}
        +\frac{L_1\Delta}{\thesischfiveeps^2}
        +\frac{\delta\Delta\sigma_1}{\thesischfiveeps^3}
        \right).
    \end{thesischfivefitdisplay}}
\end{theorem}

\begin{proof}
{\color{navy}
Let $e^k\eqdef g^k-\nabla f(x^k)$ and
$\mathcal F_k\eqdef\sigma(g^0,\xi^1,\ldots,\xi^k)$, with
$\mathcal F_0=\sigma(g^0)$. Smoothness and the update
$x^k=x^{k-1}-\gamma g^{k-1}$ imply, for every $k\ge1$,
\begin{thesischfivefitdisplay}
\begin{aligned}
f(x^k)
&\le f(x^{k-1})
-\frac{\gamma}{2}\norm{\nabla f(x^{k-1})}^2
+\frac{\gamma}{2}\norm{e^{k-1}}^2\\
&\quad+
\left(\frac{L_1}{2}-\frac{1}{2\gamma}\right)
\norm{x^k-x^{k-1}}^2 .
\end{aligned}
\thesischfivenumberthis\label{ch5:p2-descent_chapter_5}
\end{thesischfivefitdisplay}
The exact \algname{MVR} recursion and the conditional unbiasedness of the
stochastic gradients give
\begin{thesischfivefitdisplay}
\begin{aligned}
\thesischfiveExpCond{\norm{e^k}^2}{\mathcal F_{k-1}}
&\le(1-\eta)^2\norm{e^{k-1}}^2+2\eta^2\sigma_1^2\\
&\quad+2(1-\eta)^2\delta^2\norm{x^k-x^{k-1}}^2 .
\end{aligned}
\thesischfivenumberthis\label{ch5:p2-error_chapter_5}
\end{thesischfivefitdisplay}
Define
\begin{thesischfivefitdisplay}
    \mathcal V_k\eqdef f(x^k)-f^{\inf}
    +\frac{\gamma}{2\eta}\norm{e^k}^2 .
\end{thesischfivefitdisplay}
Adding $\gamma/(2\eta)$ times~\eqref{ch5:p2-error_chapter_5}
to~\eqref{ch5:p2-descent_chapter_5}, and using
$\eta+(1-\eta)^2\le1$, yields
\begin{thesischfivefitdisplay}
\begin{aligned}
\thesischfiveExpCond{\mathcal V_k}{\mathcal F_{k-1}}
&\le\mathcal V_{k-1}
-\frac{\gamma}{2}\norm{\nabla f(x^{k-1})}^2
+\gamma\eta\sigma_1^2\\
&\quad+
\left(\frac{L_1}{2}-\frac{1}{2\gamma}
+\frac{\gamma(1-\eta)^2\delta^2}{\eta}\right)
\norm{x^k-x^{k-1}}^2 .
\end{aligned}
\end{thesischfivefitdisplay}
The stated stepsize makes the last coefficient nonpositive: indeed, with
$a=\sqrt2\,\delta(1-\eta)/\sqrt\eta$ and
$\gamma=(L_1+a)^{-1}$, we have
$\gamma a\le1$ and hence
$1/\gamma-L_1=a\ge\gamma a^2
=2\gamma(1-\eta)^2\delta^2/\eta$.
Taking total expectations and summing from $k=1$ to $K$ therefore gives
\begin{thesischfivefitdisplay}
    \frac1K\sum_{k=0}^{K-1}
    \Exp{\norm{\nabla f(x^k)}^2}
    \le
    \frac{2\Delta}{\gamma K}
    +2\eta\sigma_1^2
    +\frac{\Exp{\norm{e^0}^2}}{\eta K}.
    \thesischfivenumberthis\label{ch5:p2-master_chapter_5}
\end{thesischfivefitdisplay}
Because $g^0$ is the average of $B_{\mathrm{init}}$ independent unbiased
samples,
$\Exp{\norm{e^0}^2}\le\sigma_1^2/B_{\mathrm{init}}$.
The choices of $\eta$, $B_{\mathrm{init}}$, and $K$ make the three
    terms on the right-hand side of~\eqref{ch5:p2-master_chapter_5} at most
    $\thesischfiveeps^2/3$ each. Since $x_{\mathrm{out}}$ is uniform on
$\{x^0,\ldots,x^{K-1}\}$ and independent of the algorithm conditional
on its iterates, the claimed squared-gradient bound follows. Jensen's
inequality then gives
    $\Exp{\norm{\nabla f(x_{\mathrm{out}})}}\le\thesischfiveeps$.

Finally,
\begin{thesischfivefitdisplay}
    \frac1\gamma
    =L_1+\frac{\sqrt2\,\delta(1-\eta)}{\sqrt\eta}
    =\cO\left(L_1+\frac{\delta\sigma_1}{\thesischfiveeps}\right),
\end{thesischfivefitdisplay}
and each iteration after initialization uses two stochastic-gradient
evaluations. Thus $B_{\mathrm{init}}+2(K-1)$ has the stated order.}
\end{proof}

\newpage 
\section{Empirical Evaluation}
\label{ch5:sec:experiments_chapter_5}

{\color{navy}To complement the theory and provide intuition about the behavior and stability of the proposed algorithms, we conduct synthetic experiments. These experiments compare convergence trajectories, illustrate the effect of double clipping, and explore how empirical behavior changes with the Pareto tail exponent.}

\subsection{Experimental Setup}
We consider the minimization of an ill-conditioned quadratic objective $f(x) = \frac{1}{2} x^\top A x$, where $A \in \R^{d \times d}$ is a diagonal matrix with $A_{1,1} = 0.01$ and $A_{i,i} = 1$ for $i > 1$, and dimension $d=10$. To simulate heavy-tailed noise, we inject synthetic noise into the exact gradients and Hessian-vector products. {\color{navy}Specifically, the noise vectors are drawn from a symmetric Pareto-like distribution generated via $s \cdot (u^{-1/r} - 1)$, where $u \sim \mathcal{U}(0,1)$, $s \sim \{-1, 1\}$ uniformly, and $r$ is the displayed tail exponent. This distribution has finite moments of every order strictly smaller than $r$, but its $r$-th moment is infinite. Consequently, these experiments illustrate behavior arbitrarily close to the $p$-BCM boundary (take any $p<r$); they do not instantiate the endpoint $p=r$ assumption exactly. The empirical comparisons below should therefore be interpreted as qualitative stability evidence rather than as a numerical verification of the theorem constants or exponents.}

\subsection{The Necessity of Double-Clipping and Algorithmic Stability}

\begin{figure}[htbp]
    \centering
    
    \includegraphics[width=0.78\linewidth]{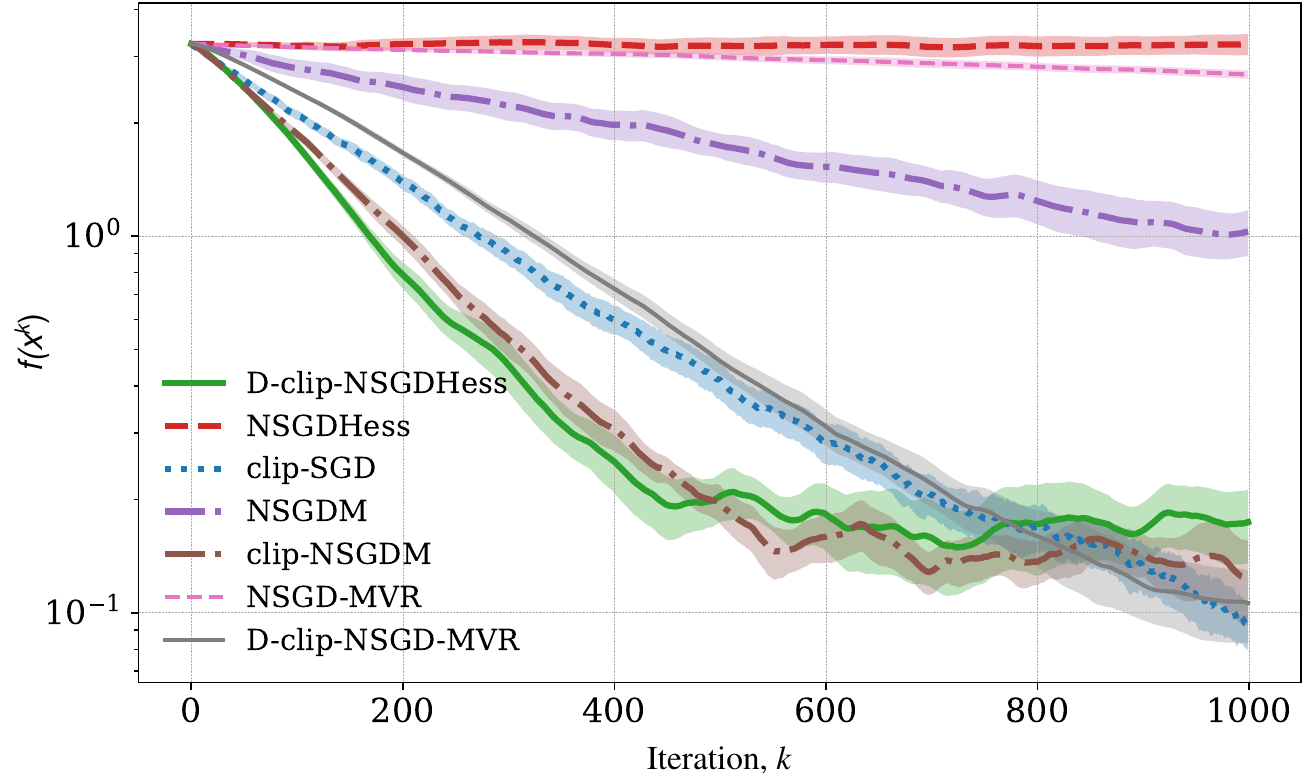} \\[1em]
    
    \includegraphics[width=\linewidth]{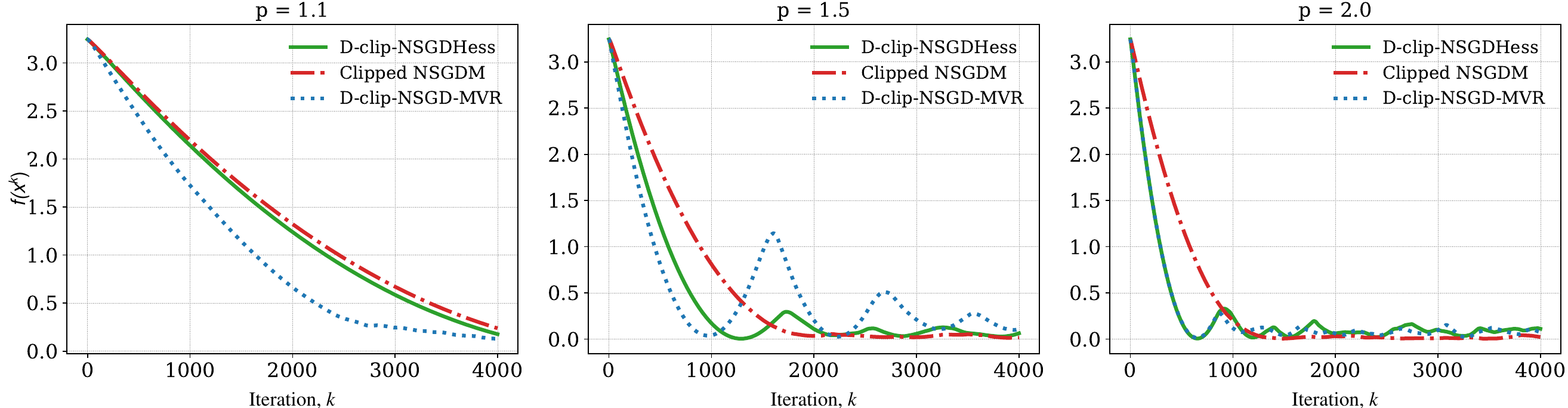}
    
    \caption[Synthetic heavy-tailed convergence and tail-exponent comparisons]{\textbf{Top:} Convergence trajectories of $f(x^{k})$ over $K=1000$ iterations for displayed tail exponent $r=1.1$. Solid/dashed lines represent the mean over 20 independent runs, and shaded regions show the empirical standard deviation across these runs. {\color{navy}Because the underlying variance may be infinite, this shading is descriptive and is not a confidence interval.} \textbf{Bottom:} Algorithm performance across displayed tail exponents $r \in \{1.1, 1.5, 2.0\}$. Hyperparameters use the corresponding theoretical power-law scalings.}
    \label{ch5:fig:convergence_trajectories_chapter_5}
\end{figure}

In our first experiment (Figure~\ref{ch5:fig:convergence_trajectories_chapter_5}, {\color{navy}top}), we fix the displayed tail exponent to {\color{navy}$r=1.1$} (an extreme heavy-tailed regime) and compare our proposed double-clipped methods (\algname{D-\thesischfiveclip-NSGD-MVR} and \algname{D-\thesischfiveclip-NSGDHess}) against standard clipped SGD \parencite{10.5555/3495724.3497014}, clipped NSGD with momentum, \algname{clip-NSGDM}~\parencite{NEURIPS2021_26901deb,pmlr-v195-liu23c}, and their unclipped counterparts. We run $K=1000$ iterations and average the results over 20 independent trials. 

{\color{navy}In these trials,} the unclipped methods \algname{NSGDM}~\parencite{pmlr-v119-cutkosky20b,hübler2025gradientclippingnormalizationheavy} and \algname{NSGD-MVR} fail to converge and exhibit severe instability, eventually diverging. In contrast, all clipped variants decrease the objective. Notably, our proposed \algname{D-\thesischfiveclip-NSGD-MVR} and \algname{D-\thesischfiveclip-NSGDHess} exhibit the fastest and most stable trajectories {\color{navy}among the tested configurations, supporting---but not by itself proving---the usefulness of combining variance reduction or Hessian correction with double clipping.}

Next, we evaluate the performance of these algorithms across a spectrum of heavy-tailed noise distributions by varying the displayed tail exponent {\color{navy}$r \in \{1.1, 1.5, 2.0\}$} over $K=4000$ iterations (Figure~\ref{ch5:fig:convergence_trajectories_chapter_5}, {\color{navy}bottom}). The stepsize $\gamma$, momentum parameter $\eta$, and clipping thresholds for each method use the $K$-dependent power-law exponents suggested by our theorems. {\color{navy}The results are consistent with the predicted qualitative trend:} in the most challenging regime ($r=1.1$), \algname{D-\thesischfiveclip-NSGD-MVR} significantly outperforms the standard \algname{\thesischfiveclip-NSGDM} \parencite{NEURIPS2021_26901deb,pmlr-v195-liu23c} baseline. As $r$ increases toward $2.0$, the gap between the methods narrows, and all algorithms exhibit rapid, stable convergence. {\color{navy}Here $r=2$ remains a borderline infinite-variance Pareto case; it should not be described as the bounded-variance regime.}

\subsection{Theoretical Vs. Empirical Complexity Scaling}

\begin{figure}[htbp]
    \centering
    \begin{minipage}{0.48\textwidth}
        \centering
        \includegraphics[width=\linewidth]{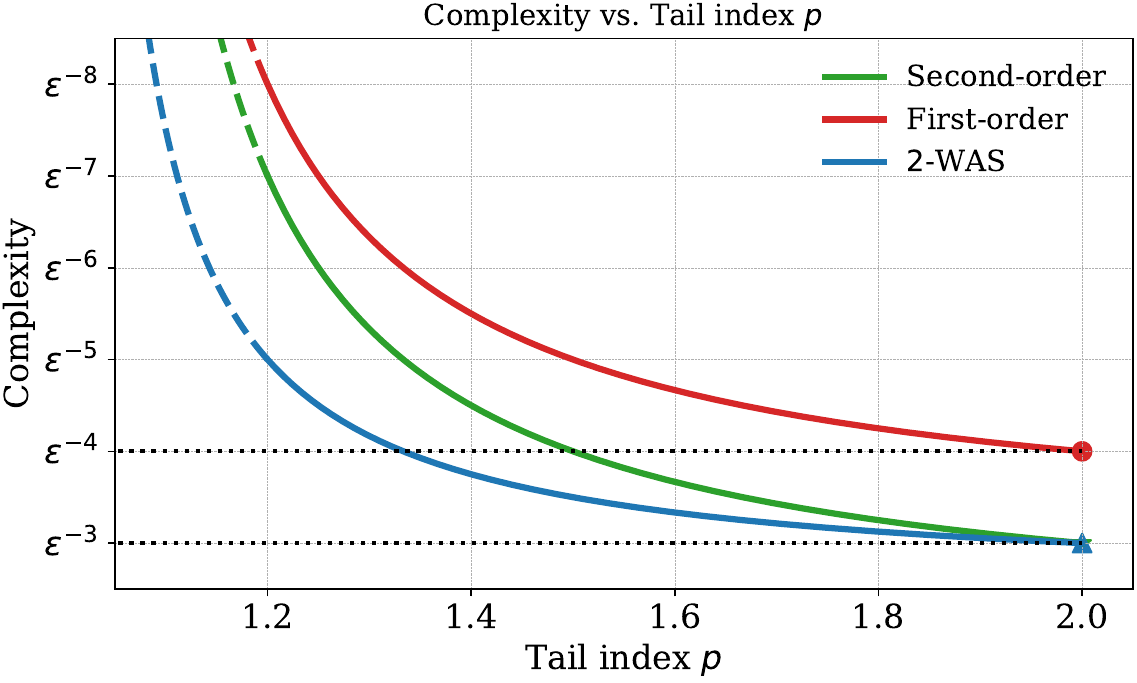} %
    \end{minipage}\hfill
    \begin{minipage}{0.48\textwidth}
        \centering
        \includegraphics[width=\linewidth]{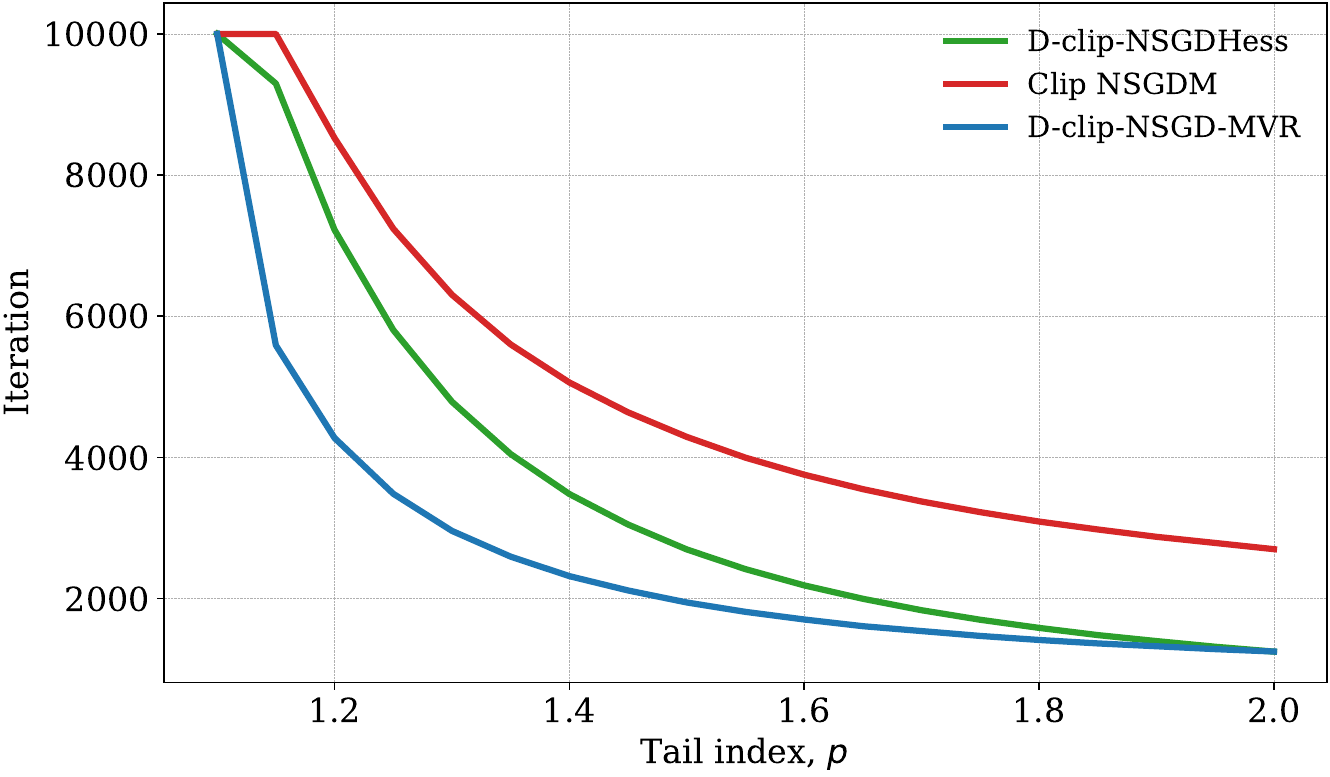} %
    \end{minipage}
    \caption{\textbf{Left:} Theoretical complexity $\mathcal{O}(\varepsilon^{-c})$ vs. moment index $p$. The y-axis represents the exponent $c$ of $\varepsilon$. \textbf{Right:} Empirical iterations required to reach a target suboptimality $f(x^{k}) < 1.2$ across different displayed Pareto exponents $r$. {\color{navy}The empirical trend is qualitatively consistent with the theoretical scaling and suggests that \algname{D-clip-NSGD-MVR} ($2$-WAS) is robust in the tested heavy-tailed regimes.}}
    \label{ch5:fig:complexity_vs_p_chapter_5}
\end{figure}

To further illustrate the robustness of our analysis, we investigate how the heavy-tail index $p$ explicitly impacts sample complexity. We compare three distinct regimes: standard first-order optimization with \algname{\thesischfiveclip NSGDM}~\parencite{NEURIPS2021_26901deb,pmlr-v195-liu23c}, second-order optimization with \algname{D-\thesischfiveclip-NSGDHess}, and first-order optimization with variance reduction under the 2-WAS assumption using \algname{D-\thesischfiveclip-NSGD-MVR}.

\textbf{Theoretical Scaling:} In the left panel of Figure~\ref{ch5:fig:complexity_vs_p_chapter_5}, we plot the theoretical complexity exponents $c$ for the bound $\mathcal{O}(\varepsilon^{-c})$ as a function of the tail index $p \in (1, 2]$. As expected, in the bounded variance regime ($p=2$), the standard first-order method requires $\mathcal{O}(\varepsilon^{-4})$, while both the second-order method and the 2-WAS variance-reduced method achieve the optimal $\mathcal{O}(\varepsilon^{-3})$ complexity. However, as the noise becomes increasingly heavy-tailed ($p \to 1$), the required complexity approaches infinity across all methods. Notably, the theoretical curves reveal that under the 2-WAS assumption, our variance-reduced method (\algname{D-clip-NSGD-MVR}) achieves a strictly better polynomial dependence on $\varepsilon$ than even the second-order method for all $p < 2$.

{\color{navy}\textbf{Empirical comparison:} We track the number of iterations required for each algorithm to reach the fixed suboptimality threshold $f(x^{k}) < 1.2$ as the displayed Pareto exponent $r$ varies from $1.1$ to $2.0$ (Figure~\ref{ch5:fig:complexity_vs_p_chapter_5}, right). In these trials, \algname{\thesischfiveclip-NSGDM} requires the most iterations, while \algname{D-\thesischfiveclip-NSGD-MVR} requires fewer iterations than \algname{D-\thesischfiveclip-NSGDHess} for the displayed exponents below $2$. This finite synthetic comparison is qualitatively consistent with the ordering suggested by the theoretical exponents, but it does not establish the exact asymptotic rates or a statistical crossover law.}
\endgroup
  \endgroup

\chapter{Supplementary Material and Proofs for Chapter 6}
  \begingroup
  \thesischaptersetup{6}%
  \stepEightAAppendixColors{6}%
  
\section{Lemmas}
	In this section, we introduce useful lemmas  for our analysis. 
	Lemmas~\ref{ch6:lemma:LzeroLoneSmooth_chapter_6} and~\ref{ch6:lemma:bound_avg_norm_grad_multinodeef21_chapter_6} introduce inequalities by generalized smoothness, while Lemmas~\ref{ch6:lemma:descent_ineq_generalized_hubler_chapter_6} and~\ref{ch6:lemma:conv_v2_chapter_6} present the descent inequality and convergence rate, respectively, when the normalized gradient descent update is applied.
	
	\begin{lemma}\label{ch6:lemma:LzeroLoneSmooth_chapter_6}
		Let each $f_i(x)$ be generalized smooth with  parameters $L_0,L_1>0$, and lower bounded by $f^{\inf}_{i}$, and  let $f(x)=\frac{1}{n}\sum_{i=1}^n f_i(x)$. Then, for any $x,y\in\R^d$
		{\color{darkred}
		\begin{equation}
						\norm{\nabla f_i(x)-\nabla f_i(y)}
			\leq (L_0+L_1\norm{\nabla f_i(y)})
			\exp\!\left(L_1\norm{x-y}\right)\norm{x-y}
			\label{ch6:eqn:LzeroLoneSmooth_ineq_0_chapter_6}
		\end{equation}
		\begin{equation}
			\begin{aligned}
			f_i(y)&\leq f_i(x)+\inner{\nabla f_i(x)}{y-x}\\
			&\quad+\frac{L_0+L_1\norm{\nabla f_i(x)}}{2}
			\exp\!\left(L_1\norm{x-y}\right)\sqnorm{y-x}
			\end{aligned}
			\label{ch6:eqn:LzeroLoneSmooth_ineq_1_chapter_6}
		\end{equation}
		\begin{equation}
			\frac{\sqnorm{\nabla f_i(x)}}{4(L_0+L_1\norm{\nabla f_i(x)})}
			\leq f_i(x)-f_i^{\inf}
			\label{ch6:eqn:LzeroLoneSmooth_ineq_2_chapter_6}
		\end{equation}
		\begin{equation}
			\begin{aligned}
			f(y)&\leq f(x)+\inner{\nabla f(x)}{y-x}\\
			&\quad+\frac{L_0+\frac{L_1}{n}\sum_{i=1}^n\norm{\nabla f_i(x)}}{2}
			\exp\!\left(L_1\norm{x-y}\right)\sqnorm{y-x}
			\end{aligned}
			\label{ch6:eqn:LzeroLoneSmooth_ineq_3_chapter_6}
		\end{equation}}
	\end{lemma}
	\begin{proof}
		The first and second statements are derived in \textcite[Proposition 3.2]{chen2023generalized}. Next, the third inequality follows from \textcite[Lemma 2.2]{gorbunov2024methods}. Finally, averaging \eqref{ch6:eqn:LzeroLoneSmooth_ineq_1_chapter_6} for $i=1,\ldots,n$ and taking into account that $f(x) = \frac{1}{n}\sum_{i=1}^n f_i(x)$, we get \eqref{ch6:eqn:LzeroLoneSmooth_ineq_3_chapter_6}. 
	\end{proof}

	\begin{lemma}\label{ch6:lemma:bound_avg_norm_grad_multinodeef21_chapter_6}
		{\color{darkred}Let $f_i(x)$ be generalized smooth with parameters $L_0>0$ and $L_1>0$, lower bounded by $f^{\inf}_{i}$, and let $f(x)$ be lower bounded by $f^{\inf}$. Then, for any $x\in\R^d$,}
		\begin{align}
			\frac{1}{n}\sum_{i=1}^n \norm{\nabla f_i(x)} &\leq 8L_1 (f(x)-f^{\inf}) \notag\\
			&\quad + \frac{8L_1}{n}\sum_{i=1}^n (f^{\inf}-f^{\inf}_{i}) + L_0/L_1. \label{ch6:eq:sum_of_grad_norms_bound_chapter_6}
		\end{align}
	\end{lemma}
	\begin{proof}
		By the $(L_0,L_1)$-smoothness of $f_i(x)$, 
		\begin{equation*}
			4(f_i(x)-f^{\inf}_{i}) \overset{\eqref{ch6:eqn:LzeroLoneSmooth_ineq_2_chapter_6}}{\geq} \frac{\sqnorm{\nabla f_i(x)}}{L_0+L_1\norm{\nabla f_i(x)}} \geq \begin{cases} 
				\frac{\sqnorm{\nabla f_i(x)}}{2L_0} \quad \text{if } \norm{\nabla f_i(x)} \leq \frac{L_0}{L_1} \\ 
				\frac{\norm{\nabla f_i(x)}}{2L_1} \quad \text{otherwise}.
			\end{cases}
		\end{equation*}  
		This condition implies
		\begin{align*}
			\norm{\nabla f_i(x)} 
			& \leq \max( 8L_1 (f_i(x)-f^{\inf}_{i}), L_0/L_1) \\
			& \leq  8L_1 (f_i(x)-f^{\inf}_{i}) + L_0/L_1 \\
			& \leq  8L_1 (f_i(x)-f^{\inf}) + 8L_1(f^{\inf}-f^{\inf}_{i})  + L_0/L_1.
		\end{align*}
		Finally, by the fact that $f(x)=\frac{1}{n}\sum_{i=1}^n f_i(x)$, 
		\begin{equation*}
			\frac{1}{n}\sum_{i=1}^n \norm{\nabla f_i(x)}  \leq 8L_1 (f(x)-f^{\inf}) + \frac{8L_1}{n}\sum_{i=1}^n (f^{\inf}-f^{\inf}_{i}) + L_0/L_1.
		\end{equation*}    
	\end{proof}
	
	\begin{lemma}\label{ch6:lemma:descent_ineq_generalized_hubler_chapter_6}
		Let $f(x) = \frac{1}{n}\sum_{i=1}^n f_i(x)$, where each $f_i(x)$ is  generalized smooth with parameters $L_0,L_1 >0$. 
		{\color{darkred}Let $x^{k+1}=x^k - \frac{\gamma_k}{\norm{v^k}} v^k$ for $\gamma_k>0$, with the quotient defined as zero when $v^k=0$, as in Algorithm~\ref{ch6:alg:normalized_ef21_chapter_6}. Then,}
		\begin{align*}
			f(x^{k+1}) 
			& \leq f(x^k) - \gamma_k \norm{\nabla f(x^k)} + 2\gamma_k \norm{\nabla f(x^k)-v^k}  \\ 
			&\quad+\frac{\gamma_k^2}{2}\exp\left(\gamma_k L_1\right) \left(L_0 + \frac{L_1}{n}\sum_{i=1}^n \norm{\nabla f_i(x^k)}\right). 
		\end{align*}   
	\end{lemma}
	\begin{proof}
			If $v^k=0$, then $x^{k+1}=x^k$ and $\norm{\nabla f(x^k)-v^k}=\norm{\nabla f(x^k)}$, so the claim holds trivially. Let $v^k\neq0$. Let each $f_i(x)$ be  generalized smooth with $L_0,L_1>0$, and $f(x)=\frac{1}{n}\sum_{i=1}^n f_i(x)$. By \eqref{ch6:eqn:LzeroLoneSmooth_ineq_3_chapter_6} of Lemma~\ref{ch6:lemma:LzeroLoneSmooth_chapter_6}, and by the fact that $x^{k+1}=x^k - \frac{\gamma_k}{\norm{v^k}} v^k$ for $\gamma_k>0$,
		\begin{align*}
			f(x^{k+1}) 
			& \leq f(x^k) - \frac{\gamma_k}{\norm{v^k}} \langle \nabla f(x^k) , v^k \rangle \\
			&\quad + \frac{\gamma_k^2}{2}\exp(\gamma_k L_1) \left( L_0 + \frac{L_1}{n}\sum_{i=1}^n \norm{\nabla f_i(x^k)} \right) \\
			& =  f(x^k) - \frac{\gamma_k}{\norm{v^k}} \langle \nabla f(x^k) - v^k , v^k \rangle - \gamma_k \norm{v^k} \\
			&\quad +  \frac{\gamma_k^2}{2}\exp(\gamma_k L_1) \left( L_0 + \frac{L_1}{n}\sum_{i=1}^n \norm{\nabla f_i(x^k)} \right) \\
			& \leq f(x^k) + \gamma_k \norm{\nabla f(x^k)-v^k}- \gamma_k \norm{v^k} \\
			&\quad +  \frac{\gamma_k^2}{2}\exp(\gamma_k L_1) \left( L_0 + \frac{L_1}{n}\sum_{i=1}^n \norm{\nabla f_i(x^k)} \right),
		\end{align*}
		where we reach the last inequality by Cauchy-Schwarz inequality. 
		Next, since 
		\begin{equation*}
			- \norm{v^k} \overset{\text{triangle ineq.}}{\leq} - \norm{\nabla f(x^k)} + \norm{\nabla f(x^k)-v^k},
		\end{equation*}
		we get  
		\begin{align*}
			f(x^{k+1}) 
			& \leq f(x^k) - \gamma_k \norm{\nabla f(x^k)} + 2\gamma_k \norm{\nabla f(x^k)-v^k} \\
			&\quad + \frac{\gamma_k^2}{2}\exp(\gamma_k L_1) \left( L_0 + \frac{L_1}{n}\sum_{i=1}^n \norm{\nabla f_i(x^k)} \right).
		\end{align*}

	\end{proof}

	\begin{lemma}\label{ch6:lemma:conv_v2_chapter_6}
		Let $\{V^k\}_{k\geq 0}, \{W^k\}_{k\geq 0}$ be non-negative sequences satisfying
		\begin{equation*}
			V^{k+1} \leq (1+b_1\exp(L_1\gamma)\gamma^2)V^k -b_2 \gamma W^k + b_3 \exp(L_1\gamma) \gamma^2,
		\end{equation*}
		for $\gamma,b_1,b_2,b_3>0$. Then, 
		\begin{equation*}
			\underset{k=0,1,\ldots,K}{\min} W^k
			\leq  \frac{V^0 \exp(b_1\exp(L_1\gamma)\gamma^2 (K+1))}{b_2\gamma (K+1)} + \frac{b_3}{b_2}\exp(L_1\gamma)\gamma.
		\end{equation*} 
	\end{lemma}
	\begin{proof}
		Define $\beta_k = \frac{\beta_{k-1}}{1+b_1\exp(L_1\gamma)\gamma^2}$ for $k=0,1,\ldots$ and $\beta_{-1}=1$. Then, we can show that  $\beta_k=\frac{1}{(1+b_1\exp(L_1\gamma)\gamma^2)^{k+1}}$ for $k=0,1,\ldots$, and that 
		\begin{align*}
			\beta_k V^{k+1} 
			& \leq (1+b_1\exp(L_1\gamma)\gamma^2)\beta_kV^k -b_2 \gamma \beta_k W^k \\
			&\quad + b_3 \exp(L_1\gamma) \gamma^2\beta_k \\ 
			& = \beta_{k-1 } V^k - b_2\gamma \beta_k W^k + b_3 \exp(L_1\gamma)\gamma^2 \beta_k.
		\end{align*}
		Therefore, 
		\begin{align*}
			\underset{k=0,1,\ldots,K}{\min} W^k
			& \leq \frac{1}{\sum_{k=0}^K \beta_k} \sum_{k=0}^K \beta_k W^k \\ 
			& \leq \frac{\sum_{k=0}^K (\beta_{k-1}V^k - \beta_k V^{k+1})}{b_2\gamma \sum_{k=0}^K \beta_k} + \frac{b_3}{b_2}\exp(L_1\gamma)\gamma  \\
			& = {\color{darkred}\frac{\beta_{-1}V^0 - \beta_K V^{K+1}}{b_2\gamma \sum_{k=0}^K \beta_k}} + \frac{b_3}{b_2}\exp(L_1\gamma)\gamma.
		\end{align*}
		{\color{darkred}By the fact that $\beta_{-1}=1$, $\beta_K >0$, and $V^{K+1} \geq 0$,}
		\begin{equation*}
			\underset{k=0,1,\ldots,K}{\min} W^k
			\leq \frac{V^0}{b_2\gamma \sum_{k=0}^K \beta_k} + \frac{b_3}{b_2}\exp(L_1\gamma)\gamma.
		\end{equation*}
		Next, since
		\begin{align*}
\sum_{k=0}^K \beta_k & \geq  (K+1) \underset{k=0,1,\ldots,K}{\min} \beta_k \\
& = \frac{K+1}{(1+b_1\exp(L_1\gamma)\gamma^2)^{K+1}},
    \end{align*}
		we have 
		\begin{align*}
			\MoveEqLeft[1] \underset{k=0,1,\ldots,K}{\min} W^k \\
			& \leq \frac{V^0 (1+b_1\exp(L_1\gamma)\gamma^2)^{K+1}}{b_2\gamma (K+1)} + \frac{b_3}{b_2}\exp(L_1\gamma)\gamma  \\
			& \overset{1+x \leq \exp(x)}{\leq} \frac{V^0 \exp(b_1\exp(L_1\gamma)\gamma^2 (K+1))}{b_2\gamma (K+1)} + \frac{b_3}{b_2}\exp(L_1\gamma)\gamma.
		\end{align*} 
		
	\end{proof}

	\clearpage

	\section{Convergence Proof for \normEFname{} (Theorem~\ref{ch6:thm:ef21_chapter_6})}
	
	{\color{darkred}In this section, we derive the convergence rate of \normEFname{}. We start with the following technical lemma.}
	\begin{lemma}\label{ch6:lemma:ef21_bound_for_diff_chapter_6}
		Let Assumptions~\ref{ch6:assum:LzeroLoneSmooth_chapter_6} and \ref{ch6:assum:contractive_comp_chapter_6} hold. Then, the iterates $\{x^k\}$ generated by \normEFname{} (Algorithm~\ref{ch6:alg:normalized_ef21_chapter_6}) satisfy
			\begin{align}
				\Exp{\norm{\nabla f_i(x^{k+1}) - g_i^{k+1}}} 
				& \leq   \sqrt{1-\alpha}\Exp{\norm{\nabla f_i(x^k) - g_i^k}} \notag \\
				&\quad + \sqrt{1-\alpha}\exp(L_1\gamma_k)\gamma_k(L_0 + L_1\Exp{\norm{\nabla f_i(x^k)}}). \label{ch6:eqn:bound_variance_multief21_di_2_chapter_6}
			\end{align}
		
	\end{lemma}
	\begin{proof}
		{\color{darkred}Condition on the history $\gH_{k+1}$ immediately before the compression at iteration $k+1$; in particular, $x^{k+1}$ and $g_i^k$ are then fixed. By the update of $g_i^{k+1}$ in Algorithm~\ref{ch6:alg:normalized_ef21_chapter_6},}
		    
			\begin{align*}
				\MoveEqLeft {\color{darkred}\thesischsixExpCond{\norm{\nabla f_i(x^{k+1}) - g_i^{k+1}}}{\gH_{k+1}}} \\
				& {\color{darkred}{}= \thesischsixExpCond{ \norm{\nabla f_i(x^{k+1}) - g_i^k - \thesischsixcC^{k+1}(\nabla f_i(x^{k+1}) - g_i^k )} }{\gH_{k+1}}} \\
				& {\color{darkred}{}\leq \sqrt{ \thesischsixExpCond{ \sqnorm{\nabla f_i(x^{k+1}) - g_i^k - \thesischsixcC^{k+1}(\nabla f_i(x^{k+1}) - g_i^k )} }{\gH_{k+1}}}}.
		\end{align*}
		{\color{darkred}The inequality is Jensen's inequality for the concave square-root function. By the $\alpha$-contractive property in~\eqref{ch6:eqn:contractive_comp_chapter_6}, the $\gH_{k+1}$-measurability of $\nabla f_i(x^{k+1})-g_i^k$, and the triangle inequality,}
		
			\begin{align*}
			\MoveEqLeft[1] {\color{darkred}\thesischsixExpCond{\norm{\nabla f_i(x^{k+1}) - g_i^{k+1}}}{\gH_{k+1}}} \\
			& \leq  {\color{darkred}\sqrt{(1-\alpha) \thesischsixExpCond{ \sqnorm{\nabla f_i(x^{k+1}) - g_i^k} }{\gH_{k+1}}  }} \\
				& =  \sqrt{1-\alpha}  \norm{\nabla f_i(x^{k+1}) - g_i^k} \\
				& \leq      \sqrt{1-\alpha}\norm{\nabla f_i(x^k) - g_i^k} \\
				&\quad + \sqrt{1-\alpha}\norm{\nabla f_i(x^{k+1})-\nabla f_i(x^k)}.  
			\end{align*}

		By the generalized smoothness of $f_i(x)$ in \eqref{ch6:eqn:LzeroLoneSmooth_chapter_6}, and by the fact that $x^{k+1}=x^k-\gamma_k \frac{g^k}{\norm{g^k}}$,
		
			\begin{align*}
			\MoveEqLeft[1] {\color{darkred}\thesischsixExpCond{\norm{\nabla f_i(x^{k+1}) - g_i^{k+1}}}{\gH_{k+1}}}
\\
			& \leq  \sqrt{1-\alpha}\norm{\nabla f_i(x^k) - g_i^k} \\
				&\quad + \sqrt{1-\alpha}(L_0 + L_1\norm{\nabla f_i(x^k)})\exp(L_1\gamma_k)\gamma_k.
			\end{align*}

		{\color{darkred}Taking total expectation and using the tower property,}
			\begin{equation*}
			{\color{darkred}\Exp{\norm{\nabla f_i(x^{k+1}) - g_i^{k+1}}} =  \Exp{ \thesischsixExpCond{\norm{\nabla f_i(x^{k+1}) - g_i^{k+1}}}{\gH_{k+1}} },}
			\end{equation*}
			we have 
			\begin{align}
				\MoveEqLeft[1] \Exp{\norm{\nabla f_i(x^{k+1}) - g_i^{k+1}}} 
\notag \\
				& \leq   \sqrt{1-\alpha}\Exp{\norm{\nabla f_i(x^k) - g_i^k}} \notag \\
				&\quad + \sqrt{1-\alpha}\exp(L_1\gamma_k)\gamma_k(L_0 + L_1\Exp{\norm{\nabla f_i(x^k)}}). \notag
			\end{align}
		
		This concludes the proof.
	\end{proof}

	Next, we present the following descent lemma for \normEFname{}. 
	
	\begin{lemma}\label{ch6:lemma:di_multinodeef21_chapter_6}
		Let Assumptions~\ref{ch6:assum:lowerbound_whole_f_chapter_6}-\ref{ch6:assum:contractive_comp_chapter_6} hold. 
		Then, the iterates $\{x^k\}$ generated by \normEFname{} (Algorithm~\ref{ch6:alg:normalized_ef21_chapter_6}) with nonincreasing step sizes $\gamma_{k+1}\leq\gamma_k$ satisfy 
		{\color{darkred}
		\begin{align*}
			\Exp{V^{k+1}} 
			& \leq  \Exp{V^k}  +  c_1 \exp(L_1\gamma_k)\gamma_k^2\frac{1}{n}\sum_{i=1}^n \Exp{\norm{\nabla f_i(x^k)}} \\
			&\quad -   \gamma_k  \Exp{ \norm{\nabla f(x^k)} }  + c_0\exp(L_1\gamma_k)\gamma_k^2,
		\end{align*}
		}
		where $V^k \eqdef f(x^k)-f^{\inf} + \frac{2\gamma_k}{1-\sqrt{1-\alpha}} \frac{1}{n}\sum_{i=1}^n \norm{\nabla f_i(x^k)-g_i^k}$, and  $c_i = \frac{L_i}{2} + 2\frac{\sqrt{1-\alpha} L_i}{1-\sqrt{1-\alpha}}$ for $i=0,1$.  
	\end{lemma}
	\begin{proof}

		{\color{darkred}For brevity, define
		\[
			A_k\eqdef\frac{2\gamma_k}{1-\sqrt{1-\alpha}},
			\qquad
			V^k\eqdef f(x^k)-f^{\inf}
			+A_k\frac{1}{n}\sum_{i=1}^n
			\norm{\nabla f_i(x^k)-g_i^k}.
		\]
		Lemma~\ref{ch6:lemma:descent_ineq_generalized_hubler_chapter_6} then gives}
		\begin{align*}
			\Exp{V^{k+1}} 
			& \leq \Exp{f(x^k)-f^{\inf}} - \gamma_k \Exp{\norm{\nabla f(x^k)}}  \\
			&\quad + \exp(L_1\gamma_k)\gamma_k^2 \frac{L_1}{2n}\sum_{i=1}^n \Exp{\norm{\nabla f_i(x^k)}}+ \exp(L_1\gamma_k)\gamma_k^2\frac{L_0}{2} \\
			&\quad + 2\gamma_k \Exp{\norm{\nabla f(x^k)- g^k}} + A_{k+1} \frac{1}{n}\sum_{i=1}^n \Exp{\norm{\nabla f_i(x^{k+1}) - g_i^{k+1}}}.
		\end{align*}
		Identities $\nabla f(x^k) = \frac{1}{n}\sum_{i=1}^n \nabla f_i(x^k)$ and $g^k = \frac{1}{n}\sum_{i=1}^n g_i^k$ and the triangle inequality imply
		\begin{align*}
			\Exp{V^{k+1}} 
			& \leq \Exp{f(x^k)-f^{\inf}} - \gamma_k \Exp{\norm{\nabla f(x^k)}} \\
			&\quad + \exp(L_1\gamma_k)\gamma_k^2 \frac{L_1}{2n}\sum_{i=1}^n \Exp{\norm{\nabla f_i(x^k)}}+ \exp(L_1\gamma_k)\gamma_k^2\frac{L_0}{2} \\
			& \quad + 2\gamma_k \frac{1}{n}\sum_{i=1}^n \Exp{\norm{\nabla f_i(x^{k}) - g_i^{k}}} \\
			&\quad +  A_{k+1} \frac{1}{n}\sum_{i=1}^n \Exp{\norm{\nabla f_i(x^{k+1}) - g_i^{k+1}}}.
		\end{align*}
		Next, we apply \eqref{ch6:eqn:bound_variance_multief21_di_2_chapter_6}:
		\begin{align*}
			\Exp{V^{k+1}}  & \leq \Exp{f(x^k) - f^{\inf}} - \gamma_k \Exp{\norm{\nabla f(x^k)}} \\
			&\quad + \left( \frac{\gamma_k^2 }{2} + A_{k+1} \sqrt{1-\alpha}  \gamma_k \right)\exp(L_1\gamma_k)  L_0   \\
			&\quad + \left( \frac{\gamma_k^2}{2} + A_{k+1} \sqrt{1-\alpha} \gamma_k \right) \exp(L_1\gamma_k) L_1 \frac{1}{n}\sum_{i=1}^n \Exp{\norm{\nabla f_i(x^k)}} \\
			&\quad  +  \left( 2\gamma_k + A_{k+1} \sqrt{1-\alpha} \right) \frac{1}{n}\sum_{i=1}^n \Exp{\norm{\nabla f_i(x^k) - g_i^k}}.
		\end{align*}
		
		Since $A_k = \frac{2\gamma_k}{1-\sqrt{1-\alpha}}$ and $\gamma_{k+1} \leq \gamma_k$,  
		\begin{equation*}
			2\gamma_k + A_{k+1} \sqrt{1-\alpha} \leq 2\gamma_k + A_k\sqrt{1-\alpha} = A_k. 
		\end{equation*}
		Therefore, 
		\begin{align*}
			\Exp{V^{k+1}} 
			& \leq  \Exp{V^k}  +  c_1 \exp(L_1\gamma_k)\gamma_k^2\frac{1}{n}\sum_{i=1}^n \Exp{\norm{\nabla f_i(x^k)}} \\
			&\quad -   \gamma_k  \Exp{ \norm{\nabla f(x^k)} }  + c_0\exp(L_1\gamma_k)\gamma_k^2,
		\end{align*}
		where $c_i = \frac{L_i}{2} + 2\frac{\sqrt{1-\alpha} L_i}{1-\sqrt{1-\alpha}}$ for $i=0,1$.
	\end{proof}

	\subsection{Proof of Theorem~\ref{ch6:thm:ef21_chapter_6}}
	{\color{darkred}Set $\gamma_k\equiv\gamma_K=\gamma_0/\sqrt{K+1}$ for $k=0,\ldots,K$. We are now ready to prove Theorem~\ref{ch6:thm:ef21_chapter_6}.}
	From Lemma~\ref{ch6:lemma:di_multinodeef21_chapter_6} and~\ref{ch6:lemma:bound_avg_norm_grad_multinodeef21_chapter_6}, and by the fact that $c_1 L_0/L_1 = c_0$, we have
	\begin{align*}
		\Exp{V^{k+1}} 
		& \leq  \Exp{V^k}  + 8c_1 L_1 \exp(L_1\gamma_k) \gamma_k^2 \Exp{f(x^k)-f^{\inf}} \\
		&\quad -  \gamma_k  \Exp{ \norm{\nabla f(x^k)} }  + B\exp(L_1\gamma_k)\gamma_k^2,
	\end{align*}
	where $B = 2c_0 +  \frac{8c_1L_1}{n}\sum_{i=1}^n (f^{\inf}-f^{\inf}_{i})$. Using the fact that $f(x^k)-f^{\inf} \leq V^k$, we derive
	\begin{equation*}
		\Exp{V^{k+1}} 
		\leq  (1+8c_1 L_1\exp(L_1\gamma_k)\gamma_k^2)\Exp{V^k}  -  \gamma_k  \Exp{ \norm{\nabla f(x^k)} }  + B\exp(L_1\gamma_k)\gamma_k^2.
	\end{equation*}
	{\color{darkred}Applying Lemma~\ref{ch6:lemma:conv_v2_chapter_6} with the deterministic sequences $\widehat V^k = \Exp{V^k}$ and $W^k = \Exp{\norm{\nabla f(x^k)}}$, with $b_1 = 8c_1L_1$, $b_2 = 1$, and $b_3 = B$, we get}
	\begin{equation*}
		\underset{k=0,1,\ldots,K}{\min} W^k
		\leq  {\color{darkred}\frac{\Exp{V^0} \exp(b_1\exp(L_1\gamma_K)\gamma_K^2 (K+1))}{b_2\gamma_K (K+1)} + \frac{b_3}{b_2}\exp(L_1\gamma_K)\gamma_K.}
	\end{equation*} 
	
	{\color{darkred}Finally, since $\gamma_K^2(K+1)=\gamma_0^2$,}
	\begin{equation*}
		\underset{k=0,1,\ldots,K}{\min} W^k
		\leq  {\color{darkred}\frac{\Exp{V^0} \exp(b_1 \exp(L_1\gamma_K)\gamma_0^2)}{b_2\gamma_0 \sqrt{K+1}} + \frac{b_3}{b_2}\frac{\gamma_0 \exp(L_1\gamma_K)}{\sqrt{K+1}}.}
	\end{equation*}

\section{Discussion on Theorem~\ref{ch6:thm:ef21_chapter_6}}\label{ch6:sec:diss_thm_ef21_chapter_6}
{\color{darkred}
The proof of Theorem~\ref{ch6:thm:ef21_chapter_6} assumes $L_1>0$.
We therefore treat the classical $L$-smooth case directly instead of
substituting $L_1=0$ into an expression containing $L_0/L_1$.
Let $r\eqdef\sqrt{1-\alpha}$ and initialize
$g_i^{-1}=\nabla f_i(x^0)$, so that $g_i^0=\nabla f_i(x^0)$ almost
surely. For
\[
	b\eqdef\frac{L}{2}+\frac{2rL}{1-r},
	\qquad
	V^k\eqdef f(x^k)-f^{\inf}
	+\frac{2\gamma}{1-r}\frac1n\sum_{i=1}^n
	\norm{\nabla f_i(x^k)-g_i^k},
\]
the classical smoothness descent inequality and compressor recursion give
\[
	\Exp{V^{k+1}}
	\leq\Exp{V^k}-\gamma\Exp{\norm{\nabla f(x^k)}}+b\gamma^2.
\]
Set $\gamma=\gamma_0/\sqrt{K+1}$. Telescoping yields
\[
	\min_{0\leq k\leq K}\Exp{\norm{\nabla f(x^k)}}
	\leq\frac{1}{\sqrt{K+1}}
	\left(\frac{f(x^0)-f^{\inf}}{\gamma_0}+b\gamma_0\right).
\]
With $\gamma_0=\sqrt{(f(x^0)-f^{\inf})/b}$, this becomes
\begin{align*}
	\min_{0\leq k\leq K}\Exp{\norm{\nabla f(x^k)}}
	&\leq
	\sqrt{\frac{2L(1+3r)(1+r)}{\alpha}}\,
	\sqrt{\frac{f(x^0)-f^{\inf}}{K+1}}\\
	&\leq
	4\sqrt{\frac{L}{\alpha}}\,
	\sqrt{\frac{f(x^0)-f^{\inf}}{K+1}}.
\end{align*}
By comparison, Theorem~1 of \textcite{richtarik2021ef21}, with common
smoothness constant $L$, gives at most
$2\sqrt{L/\alpha}\sqrt{(f(x^0)-f^{\inf})/(K+1)}$ after converting its
squared-gradient guarantee by Jensen's inequality. Thus, the displayed
constant for normalized \algname{EF21} is at most twice that of
\algname{EF21}; both bounds have the same $\cO(K^{-1/2})$ rate.}

		\section{Convergence of \normEFname{} for a Single-node Case}
	In this section, we provide the convergence of \normEFname{} for a single-node case. 
	{\color{darkred}In particular, the algorithm enjoys $\cO(1/K)$ convergence up to the residual
	\[
		\frac{c_0\exp(L_1\gamma)\gamma}{1-c_1\exp(L_1\gamma)\gamma}.
	\]
	In contrast to Theorem~\ref{ch6:thm:ef21_chapter_6} for the multi-node case, the next result uses a fixed stepsize $\gamma=1/(\kappa c_1)$ with $\kappa\geq3$.}
	
	\begin{theorem}\label{ch6:thm:normalized_ef21_singlenode_chapter_6}
		{\color{darkred}Let Assumptions~\ref{ch6:assum:lowerbound_whole_f_chapter_6}--\ref{ch6:assum:contractive_comp_chapter_6} hold with $L_1>0$. Then, the iterates $\{x^k\}$ generated by \normEFname{} (Algorithm~\ref{ch6:alg:normalized_ef21_chapter_6}) with $n=1$, $\gamma_k\equiv\gamma=1/(\kappa c_1)$, and $\kappa\geq3$ satisfy}
		\begin{equation*}
			\underset{k=0,1,\ldots,K}{\min} \Exp{\norm{\nabla f(x^k)}} \leq \frac{\Exp{V^0}-\Exp{V^{K+1}}}{\gamma(1-c_1\exp(L_1\gamma)\gamma)(K+1)} + {\color{darkred}\frac{c_0\exp(L_1\gamma)\gamma}{1-c_1\exp(L_1\gamma)\gamma}},
		\end{equation*}
		where $V^k = f(x^k)-f^{\inf} + \frac{2\gamma}{1-\sqrt{1-\alpha}} \norm{\nabla f(x^k)- g^k}$, and $c_i = \frac{L_i}{2} + 2\frac{\sqrt{1-\alpha} L_i}{1-\sqrt{1-\alpha}}$ for $i=0,1$.  
	\end{theorem}
	\begin{proof}
		In the single-node case, Lemma~\ref{ch6:lemma:ef21_bound_for_diff_chapter_6} implies
		\begin{align}
			\Exp{\norm{\nabla f(x^{k+1}) - g^{k+1} }} 
			& \leq   \sqrt{1-\alpha}\Exp{\norm{\nabla f(x^k) - g^k }} \notag \\
			&\quad + \sqrt{1-\alpha}\exp(L_1\gamma_k)\gamma_k(L_0 + L_1\Exp{\norm{\nabla f(x^k)}}). \label{ch6:eqn:bound_variance_ef21_di_2_chapter_6}
		\end{align}
		Next, for brevity, let $A_k = \frac{2\gamma_k}{1-\sqrt{1-\alpha}}$. Then, we have $V^k \eqdef f(x^k)-f^{\inf} + A_k \norm{\nabla f(x^k)-g^k}$, and from Lemma~\ref{ch6:lemma:descent_ineq_generalized_hubler_chapter_6}, we derive
		\begin{align*}
			\Exp{V^{k+1}} 
			& \leq \Exp{f(x^k) - f^{\inf}} \\
			&\quad - \left( \gamma_k  - \frac{\gamma_k^2 L_1}{2}\exp(L_1\gamma_k) \right)\Exp{\norm{\nabla f(x^k)}}+ \frac{\gamma_k^2 L_0}{2}\exp(L_1\gamma_k)  \\
			&\quad +  2\gamma_k \Exp{\norm{\nabla f(x^k) - g^k}} + A_{k+1} \Exp{\norm{\nabla f(x^{k+1}) - g^{k+1}}}  \\
			& \overset{\eqref{ch6:eqn:bound_variance_ef21_di_2_chapter_6}}{\leq} \Exp{f(x^k) - f^{\inf}} \\
			&\quad +  \left( 2\gamma_k + A_{k+1} \sqrt{1-\alpha} \right) \Exp{\norm{\nabla f(x^k) - g^k}} \\
			&\quad - \biggl( \gamma_k  - \frac{\gamma_k^2 L_1}{2}\exp(L_1\gamma_k)\\
			&\qquad\quad -  A_{k+1} \sqrt{1-\alpha}L_1 \gamma_k \exp(L_1\gamma_k) \biggr)\Exp{\norm{\nabla f(x^k)}} \\
			&\quad + \frac{\gamma_k^2 L_0}{2} \exp(L_1\gamma_k) + A_{k+1} \sqrt{1-\alpha} L_0 \gamma_k \exp(L_1\gamma_k). 
		\end{align*}
		
		Since $A_k = \frac{2\gamma_k}{1-\sqrt{1-\alpha}}$ and $\gamma_{k+1} \leq \gamma_k$,  
		\begin{equation*}
			2\gamma_k + A_{k+1} \sqrt{1-\alpha} \leq 2\gamma_k + A_k\sqrt{1-\alpha} = A_k. 
		\end{equation*}
		Therefore, 
		\begin{equation*}
			\Exp{V^{k+1}} 
			\leq  \Exp{V^k}  - \left( \gamma_k  - c_1\exp(L_1\gamma_k)\gamma_k^2 \right)\Exp{ \norm{\nabla f(x^k)} } + c_0\exp(L_1\gamma_k)\gamma_k^2,
		\end{equation*}
		where $c_i = \frac{L_i}{2}  + 2\frac{\sqrt{1-\alpha} L_i}{1-\sqrt{1-\alpha}}$ for $i=0,1$.

		{\color{darkred}Finally, take $\gamma_k=\gamma=1/(\kappa c_1)$ for $\kappa\geq3$. Since $L_1/c_1\leq2$ and $\exp(2/\kappa)/\kappa$ is decreasing in $\kappa$,}
		\[
			{\color{darkred}
			c_1\exp(L_1\gamma)\gamma
			=\frac{\exp(L_1/(\kappa c_1))}{\kappa}
			\leq\frac{\exp(2/\kappa)}{\kappa}
			\leq\frac{\exp(2/3)}{3}
			\leq0.7<1.}
		\]
		{\color{darkred}Consequently,}
		\begin{align*}
			\Exp{V^{k+1}} 
			& \leq   \Exp{V^k} - \gamma \left( 1  - c_1 \exp(L_1\gamma)\gamma \right)\Exp{ \norm{\nabla f(x^k)} } + {\color{darkred}c_0\exp(L_1\gamma)\gamma^2}.
		\end{align*}
		Rearranging the terms, we derive
		\begin{align*}
			\underset{k= 0,1,\ldots,K}{\min} \Exp{\norm{\nabla f(x^k)}} 
			& \leq \frac{1}{K+1} \sum_{k=0}^K \Exp{\norm{\nabla f(x^k)}} \\
			& \leq \frac{\Exp{V^0}-\Exp{V^{K+1}}}{\gamma(1-c_1\exp(L_1\gamma)\gamma)(K+1)} \\
			&\quad + {\color{darkred}\frac{c_0\exp(L_1\gamma)\gamma}{1-c_1\exp(L_1\gamma)\gamma}}.
		\end{align*}
		Noticing that $V^k \geq 0$, we complete the proof. 
	\end{proof}

	\newpage 
	\section{Convergence of \normEFname{SGDM}  (Theorem~\ref{ch6:thm:ef21_sgdm_chapter_6})}

	In this section, we derive the convergence rate results of \normEFname{SGDM}. We first introduce auxiliary lemmas in Section~\ref{ch6:app:sec:lemma_ef21_sgdm_chapter_6}, and later prove the convergence theorem (Theorem~\ref{ch6:thm:ef21_sgdm_chapter_6}) in Section~\ref{ch6:app:sec:ef21_sgdm_chapter_6}.

	\subsection{Auxiliary Lemmas}\label{ch6:app:sec:lemma_ef21_sgdm_chapter_6}

	Now, we provide useful lemmas for analyzing \normEFname{SGDM}. First, Lemma~\ref{ch6:lem:descent_lemma_ef21_sgdm_chapter_6} shows the descent inequality of the normalized gradient descent update under Assumption~\ref{ch6:assum:LzeroLoneSmooth_chapter_6} (generalized smoothness of $f_i$). 
	Second, Lemmas~\ref{ch6:lem:compression_error_ef21_sgdm_chapter_6} and \ref{ch6:lem:momentum_error_ef21_sgdm_chapter_6} provide the upper-bound of the Euclidean distance between $v_i^{k}$ and $g_i^k$, and of the Euclidean distance between $v_i^k$ and $\nabla f_i(x^k)$, respectively.

	\begin{lemma}
		\label{ch6:lem:descent_lemma_ef21_sgdm_chapter_6}
		Consider the iterates $\{x^k\}$ generated by  Algorithm~\ref{ch6:alg:normalized_ef21_sgdm_chapter_6}. If 
		{\color{darkred}Assumption~\ref{ch6:assum:LzeroLoneSmooth_chapter_6} holds, then for any $\gamma_k>0$,}
		\begin{align}
			f(x^{k+1}) & \leq f(x^k) -\gamma_k \norm{\nabla f(x^k)} +2\gamma_k\norm{\nabla f(x^k) -v^k} \notag \\
			&\quad +2\gamma_k \norm{v^k -g^k} \notag\\
			&\quad + L_0\gamma_k^2\exp(\gamma_k L_1) + 4L_1^2\gamma_k^2\exp(\gamma_k L_1)\left(f(x^k) - f^{\inf}\right) \notag\\
			&\quad + \frac{4L_1^2\gamma_k^2\exp(\gamma_k L_1)}{n}\sum_{i=1}^n\left(f^{\inf} - f_i^{\inf}\right). \notag
		\end{align}
	\end{lemma}
	\begin{proof}
		Applying the triangle inequality in Lemma~\ref{ch6:lemma:descent_ineq_generalized_hubler_chapter_6}, i.e., $\norm{\nabla f(x^k) -g^k} \leq \norm{\nabla f(x^k) -v^k} + \norm{v^k - g^k}$, we get
			\begin{align*}
				f(x^{k+1}) & \leq f(x^k) -\gamma_k \norm{\nabla f(x^k)} +2\gamma_k\norm{\nabla f(x^k) -v^k} \\
				&\quad +2\gamma_k \norm{v^k -g^k} \\
				&\quad+ \frac{\gamma^2_k}{2}\exp\left(\gamma_kL_1\right) \left(L_0+ \frac{L_1}{n}\sum^n_{i=1}\norm{\nabla f_i (x^k)}\right) \\
				& \overset{\eqref{ch6:eq:sum_of_grad_norms_bound_chapter_6}}{\leq} f(x^k) -\gamma_k \norm{\nabla f(x^k)} +2\gamma_k\norm{\nabla f(x^k) -v^k} +2\gamma_k \norm{v^k -g^k} \\
				&\quad + L_0\gamma_k^2\exp(\gamma_k L_1) + 4L_1^2\gamma_k^2\exp(\gamma_k L_1)\left(f(x^k) - f^{\inf}\right)\\
				&\quad + \frac{4L_1^2\gamma_k^2\exp(\gamma_k L_1)}{n}\sum_{i=1}^n\left(f^{\inf} - f_i^{\inf}\right), 
			\end{align*}
			which concludes the proof.
	\end{proof}
	
	\begin{lemma}
		\label{ch6:lem:compression_error_ef21_sgdm_chapter_6}
		Consider the iterates $\{x^k\}$ generated by Algorithm~\ref{ch6:alg:normalized_ef21_sgdm_chapter_6}. If Assumptions~\ref{ch6:assum:LzeroLoneSmooth_chapter_6},~\ref{ch6:assum:contractive_comp_chapter_6}, and~\ref{ch6:assum:bounded_variance_chapter_6} hold, then for $\gamma_k>0, \eta_k \in [0,1]$, and $k \geq 0$,
		\begin{align*}
			\MoveEqLeft \frac{1}{n}\sum^n_{i=1}\Exp{\norm{v^{k+1}_i - g^{k+1}_i}}\\
			& \leq \frac{\sqrt{1-\alpha}}{n}\sum^n_{i=1}\Exp{\norm{v^{k}_i - g^{k}_i}} \\
			&\quad + \frac{\sqrt{1-\alpha}\eta_{k+1}}{n}\sum^n_{i=1}\Exp{\norm{v^{k}_i -\nabla f_i(x^k)}} \\
			&\quad +8L_1^2\sqrt{1-\alpha}\eta_{k+1}\gamma_k\exp\left(\gamma_k L_1\right) \Exp{f(x^k) - f^{\inf}} \\
			&\quad +\frac{8L_1^2\sqrt{1-\alpha}\eta_{k+1}\gamma_k\exp\left(\gamma_k L_1\right)}{n} \sum_{i=1}^n (f^{\inf}-f^{\inf}_{i}) \\
			&\quad + 2L_0\sqrt{1-\alpha}\eta_{k+1}\gamma_k\exp\left(\gamma_k L_1\right) + \sqrt{1-\alpha} \eta_{k+1} \sigma.
		\end{align*}
	\end{lemma}
	\begin{proof}
		{\color{darkred}Condition on the history $\gH_{k+1}$ immediately before the compression at iteration $k+1$, so that $v_i^{k+1},x^{k+1}$, and $g_i^k$ are fixed. Using concavity of the square-root function and the update of $g_i^{k+1}$ in Algorithm~\ref{ch6:alg:normalized_ef21_sgdm_chapter_6}, we have}
		\begin{align*}
			\thesischsixExpCond{\norm{v^{k+1}_i- g^{k+1}_i}}{\gH_{k+1}} & \leq \sqrt{\thesischsixExpCond{\sqnorm{v^{k+1}_i- g^{k+1}_i}}{\gH_{k+1}}}\\
			& ={\color{darkred}\sqrt{\thesischsixExpCond{\sqnorm{v^{k+1}_i- g^{k}_i - \gC^{k+1}\left(v^{k+1}_i- g^{k}_i\right) }}{\gH_{k+1}}}}\\
			& \overset{\eqref{ch6:eqn:contractive_comp_chapter_6}}{\leq}\sqrt{\thesischsixExpCond{(1-\alpha)\sqnorm{v^{k+1}_i- g^{k}_i  }}{\gH_{k+1}}}.
		\end{align*}
		{\color{darkred}Because $v_i^{k+1}$ and $g_i^k$ are fixed conditional on $\gH_{k+1}$, the triangle inequality gives}
		\begin{align*}
			\thesischsixExpCond{\norm{v^{k+1}_i- g^{k+1}_i}}{\gH_{k+1}}   %
			& \leq\sqrt{1-\alpha}\norm{v^{k}_i- g^{k}_i} + \sqrt{1-\alpha}\norm{v^{k+1}_i- v^{k}_i}\\
			& =\sqrt{1-\alpha}\norm{v^{k}_i- g^{k}_i} \\
			&\quad + {\color{darkred}\sqrt{1-\alpha}\eta_{k+1}\norm{\nabla f_i(x^{k+1};\xi^{k+1}_i) - v^k_i}.}
		\end{align*}
		Here, the equality comes from the definition of $v_i^{k+1}$ in Algorithm~\ref{ch6:alg:normalized_ef21_sgdm_chapter_6}. Next, by the triangle inequality,
		\begin{align*}
			\MoveEqLeft[1] \thesischsixExpCond{\norm{v^{k+1}_i- g^{k+1}_i}}{\gH_{k+1}} \\
			& \leq\sqrt{1-\alpha}\norm{v^{k}_i- g^{k}_i} + \sqrt{1-\alpha} \eta_{k+1} \|v^k_i -\nabla f_i(x^k)\|\\
			&\quad + \sqrt{1-\alpha}\eta_{k+1}\norm{\nabla f_i(x^{k}) - \nabla f_i(x^{k+1})}\\
			&\quad + \sqrt{1-\alpha}\eta_{k+1}\norm{\nabla f_i(x^{k+1};\xi^{k+1}_i) - \nabla f_i(x^{k+1})}\\
			& \overset{\eqref{ch6:eqn:LzeroLoneSmooth_ineq_0_chapter_6}}{\leq} \sqrt{1-\alpha}\norm{v^{k}_i- g^{k}_i} + \sqrt{1-\alpha} \eta_{k+1} \|v^k_i -\nabla f_i(x^k)\|\\
			&\quad + \sqrt{1-\alpha}\eta_{k+1}\left(L_0 +L_1\norm{\nabla f_i(x^k)}\right)\exp\left(L_1\norm{x^{k+1}-x^k}\right)\norm{x^{k+1}-x^k}\\
			&\quad {\color{darkred}+ \sqrt{1-\alpha}\eta_{k+1}\norm{\nabla f_i(x^{k+1};\xi^{k+1}_i) - \nabla f_i(x^{k+1})}.}
		\end{align*}
		Next, using $\norm{x^{k+1} -x^k}\leq\gamma_k$ (the step is $-\gamma_k g^k/\norm{g^k}$, or zero when $g^k=0$), and taking the expectation, we obtain
		\begin{align*}
			\Exp{\norm{v^{k+1}_i- g^{k+1}_i}} & \leq \sqrt{1-\alpha}\Exp{\norm{v^{k}_i- g^{k}_i}} \\
			&\quad + \sqrt{1-\alpha}\eta_{k+1}\Exp{\norm{v^k_i - \nabla f_i(x^k)}}\\
			&\quad +\sqrt{1-\alpha}\eta_{k+1}\gamma_{k}\exp\left(\gamma_k L_1\right) \left( L_0 +  L_1\Exp{\norm{\nabla f_i(x^k)}}\right)\\
			&\quad + \sqrt{1-\alpha}\eta_{k+1} \Exp{\norm{\nabla f_i(x^{k+1};\xi^{k+1}_i) - \nabla f_i(x^{k+1})}}.
		\end{align*}
		Finally, since 
		\begin{align*}
			\MoveEqLeft[1] \Exp{\norm{\nabla f_i(x^{k+1};\xi^{k+1}_i) - \nabla f_i(x^{k+1})}} \\
			& \leq \sqrt{ \Exp{\sqnorm{\nabla f_i(x^{k+1};\xi^{k+1}_i) - \nabla f_i(x^{k+1})}}} \\
			& \overset{\eqref{ch6:eqn:bounded_variance_chapter_6}}{\leq} \sigma, 
		\end{align*}
			we derive
			\begin{equation*}
			\begin{aligned}
				&\frac{1}{n}\sum_{i=1}^n
				\Exp{\norm{v_i^{k+1}-g_i^{k+1}}}\\[-0.2em]
				&\leq
				\frac{\sqrt{1-\alpha}}{n}\sum_{i=1}^n
				\Exp{\norm{v_i^k-g_i^k}}\\
				&\quad+
				\frac{\sqrt{1-\alpha}\,\eta_{k+1}}{n}\sum_{i=1}^n
				\Exp{\norm{v_i^k-\nabla f_i(x^k)}}\\
				&\quad+
				\sqrt{1-\alpha}\,\eta_{k+1}\gamma_k
				{\exp(\gamma_kL_1)}
				\left(
				L_0+\frac{L_1}{n}\sum_{i=1}^n
				\Exp{\norm{\nabla f_i(x^k)}}
				\right)\\
				&\quad+\sqrt{1-\alpha}\,\eta_{k+1}\sigma\\
				&\overset{\eqref{ch6:eq:sum_of_grad_norms_bound_chapter_6}}{\leq}
				\frac{\sqrt{1-\alpha}}{n}\sum_{i=1}^n
				\Exp{\norm{v_i^k-g_i^k}}\\
				&\quad+
				\frac{\sqrt{1-\alpha}\,\eta_{k+1}}{n}\sum_{i=1}^n
				\Exp{\norm{v_i^k-\nabla f_i(x^k)}}\\
				&\quad+
				8L_1^2\sqrt{1-\alpha}\,\eta_{k+1}\gamma_k
				{\exp(\gamma_kL_1)}\Exp{f(x^k)-f^{\inf}}\\
				&\quad+
				\frac{8L_1^2\sqrt{1-\alpha}\,\eta_{k+1}\gamma_k
				{\exp(\gamma_kL_1)}}{n}
				\sum_{i=1}^n\left(f^{\inf}-f_i^{\inf}\right)\\
				&\quad+
				2L_0\sqrt{1-\alpha}\,\eta_{k+1}\gamma_k
				{\exp(\gamma_kL_1)}
				+\sqrt{1-\alpha}\,\eta_{k+1}\sigma.
			\end{aligned}
			\end{equation*}
				This concludes the proof.
	\end{proof}
	
	\begin{lemma}
		\label{ch6:lem:momentum_error_ef21_sgdm_chapter_6}
		{\color{darkred}Consider the iterates $\{x^k\}$ generated by Algorithm~\ref{ch6:alg:normalized_ef21_sgdm_chapter_6}. If Assumptions~\ref{ch6:assum:LzeroLoneSmooth_chapter_6} and~\ref{ch6:assum:bounded_variance_chapter_6} hold, then, for any $\gamma_k \equiv \gamma >0$, $\eta_k \equiv \eta\in(0,1]$, and $k\geq 0$,}
		\begin{align}
			\Exp{\norm{v^{k}-\nabla f(x^{k})}} & \leq (1-\eta)^{k}\Exp{\norm{v^0 - \nabla f(x^0) }} +\frac{\sqrt{\eta}\sigma}{\sqrt{n}} \notag \\
			&\quad + \frac{2L_0 \gamma \exp\left(\gamma L_1\right)}{\eta} \notag\\
			&\quad + 8L_1^2 \gamma \exp\left(\gamma L_1\right)\sum^{k-1}_{t=0}(1-\eta)^{k-t} \Exp{f(x^t) - f^{\inf}}\notag\\
			&\quad +  \frac{8L_1^2 \gamma \exp\left(\gamma L_1\right)}{\eta n}\sum_{i=1}^n\left(f^{\inf} - f_i^{\inf}\right). \label{ch6:eq:vk_nablafxk_bound_chapter_6}
		\end{align}
		In addition,  for any $k\geq 0$,
		\begin{equation}
			{\color{darkred}
			\begin{aligned}
				&\frac{1}{n}\sum_{i=1}^n
				\Exp{\norm{v_i^{k+1}-\nabla f_i(x^{k+1})}}\\
				&\leq \frac{1-\eta}{n}\sum_{i=1}^n
				\Exp{\norm{v_i^k-\nabla f_i(x^k)}}+\eta\sigma
				+2L_0\gamma{\exp(\gamma L_1)}\\
				&\quad+8L_1^2\gamma{\exp(\gamma L_1)}
				\Exp{f(x^k)-f^{\inf}}\\
				&\quad+\frac{8L_1^2\gamma{\exp(\gamma L_1)}}{n}
				\sum_{i=1}^n(f^{\inf}-f_i^{\inf}).
			\end{aligned}}
			\label{ch6:eq:sum_vik_nablafixk_bound_chapter_6}
		\end{equation}
	\end{lemma}
	\begin{proof}
		We prove the result using the arguments similar to those given in the proof of Theorem 1 from \textcite{cutkosky2020momentum}.
			From the definition of $v_i^{k+1}$, we have the following recursion for any $k \geq 0$:
			\begin{align*}
				v^{k+1}_i & = (1-\eta)v^k_i +\eta\nabla f_i(x^{k+1};\xi^{k+1}_i) \\
				& =\nabla f_i(x^{k+1}) + (1-\eta)(v^k_i -\nabla f_i(x^k)) \\
				&\quad + (1-\eta)(\nabla f_i(x^k)-\nabla f_i(x^{k+1}))\\
				&\quad + \eta(\nabla f_i(x^{k+1};\xi^{k+1}_i) - \nabla f_i(x^{k+1})).
			\end{align*}
			Next, from the recursion of  $v^{k+1}_i$, we obtain the following recursion for $k \geq 0$:  
			\begin{align}
				H^{k+1}_i & = (1-\eta)H^{k}_i + (1-\eta)G^k_i + \eta U^{k+1}_i, \label{ch6:eq:H_i^k_recursion_chapter_6}
			\end{align}    
			where
			{\color{darkred}
			\begin{align*}
				U_i^{k+1}&\eqdef\nabla f_i(x^{k+1};\xi_i^{k+1})
				-\nabla f_i(x^{k+1}),&
				G_i^k&\eqdef\nabla f_i(x^k)-\nabla f_i(x^{k+1}),\\
				H_i^k&\eqdef v_i^k-\nabla f_i(x^k),&
				U^{k+1}&\eqdef\frac{1}{n}\sum_{i=1}^nU_i^{k+1},\\
				G^k&\eqdef\frac{1}{n}\sum_{i=1}^nG_i^k,&
				H^k&\eqdef\frac{1}{n}\sum_{i=1}^nH_i^k.
			\end{align*}}
			Unrolling the recursion for $H_i^k$, we derive
			\begin{align*}
				H^{k+1}_i 
				=& (1-\eta)^{k+1} H^0_i \\
				&\quad + \sum^{k}_{t=0}(1-\eta)^{k-t+1}G^t_i + \eta\sum^{k}_{t=0}(1-\eta)^{k-t}U^{t+1}_i.
			\end{align*}
			Averaging the above identity, we get
			\begin{align*}
				H^{k+1} & = (1-\eta)^{k+1} H^0 + \sum^{k}_{t=0}(1-\eta)^{k-t+1}G^t \\
				&\quad + \eta\sum^{k}_{t=0}(1-\eta)^{k-t}U^{t+1}.
			\end{align*}
			Next, taking the Euclidean norm, using the triangle inequality, and then taking the expectation, we obtain
			\begin{align}
				\Exp{\norm{H^{k+1}}}  & \leq (1-\eta)^{k+1} \Exp{\norm{H^0}} + \underbrace{\sum^{k}_{t=0}(1-\eta)^{k-t+1}\Exp{\norm{G^t}}}_{=: \gA_1} \notag \\
				&\quad + \eta\underbrace{\Exp{\norm{\sum^{k}_{t=0}(1-\eta)^{k-t}U^{t+1}}}}_{=:\gA_2}.\label{ch6:eqn:ef21sgdm:exp_norm_H_k_chapter_6}
			\end{align}
			To bound $\Exp{\norm{H^{k+1}}}$, we need to bound the expectation of the last two terms. First, we bound term $\gA_1$. Using the fact that $\norm{G^t} \leq \frac{1}{n}\sum_{i=1}^n \norm{G_i^t}$, and the definition of $G_i^t$, we obtain
			\begin{align*}
				\gA_1 
				& \leq \frac{1}{n}\sum^n_{i=1}\sum^{k}_{t=0}(1-\eta)^{k-t+1}\Exp{\norm{\nabla f_i(x^t)-\nabla f_i(x^{t+1})}}\\
				& \overset{\eqref{ch6:eqn:LzeroLoneSmooth_ineq_0_chapter_6}}{\leq} \frac{1}{n}\sum^n_{i=1}\sum^{k}_{t=0}(1-\eta)^{k-t+1}\Exp{L_0\exp\left(L_1\norm{x^{t+1}-x^t}\right)\norm{x^{t+1}-x^t}}\\
				&\quad + \frac{1}{n}\sum^n_{i=1}\sum^{k}_{t=0}(1-\eta)^{k-t+1}\Exp{L_1\norm{\nabla f_i(x^t)}\exp\left(L_1\norm{x^{t+1}-x^t}\right)\norm{x^{t+1}-x^t}}\\
				& \leq \sum^{k}_{t=0}(1-\eta)^{k-t+1} \gamma \exp(\gamma L_1) L_0 \\
				&\quad + \frac{L_1}{n}\sum^n_{i=1}\sum^{k}_{t=0}(1-\eta)^{k-t+1} \gamma \exp\left(\gamma L_1\right) \Exp{\norm{\nabla f_i(x^t)}}\\
				& \overset{\eqref{ch6:eq:sum_of_grad_norms_bound_chapter_6}}{\leq} 2L_0 \gamma \exp\left(\gamma L_1\right)\sum^{k}_{t=0}(1-\eta)^{k-t+1} \\
				&\quad + 8L_1^2 \gamma \exp\left(\gamma L_1\right)\sum^{k}_{t=0}(1-\eta)^{k-t+1} \Exp{f(x^t) - f^{\inf}}\\
				&\quad +  \frac{8L_1^2 \gamma \exp\left(\gamma L_1\right)}{n}\sum_{i=1}^n\sum^{k}_{t=0}(1-\eta)^{k-t+1} \left(f^{\inf} - f_i^{\inf}\right)\\
				& \leq 2L_0 \gamma \exp\left(\gamma L_1\right)\sum^{\infty}_{t=0}(1-\eta)^{t} \\
				&\quad + 8L_1^2 \gamma \exp\left(\gamma L_1\right)\sum^{k}_{t=0}(1-\eta)^{k-t+1} \Exp{f(x^t) - f^{\inf}}\\
				&\quad +  \frac{8L_1^2 \gamma \exp\left(\gamma L_1\right)}{n}\sum_{i=1}^n\left(f^{\inf} - f_i^{\inf}\right) \sum^{\infty}_{t=0}(1-\eta)^{t}\\
				& = \frac{2L_0 \gamma \exp\left(\gamma L_1\right)}{\eta} \\
				&\quad + 8L_1^2 \gamma \exp\left(\gamma L_1\right)\sum^{k}_{t=0}(1-\eta)^{k-t+1} \Exp{f(x^t) - f^{\inf}}\\
				&\quad +  \frac{8L_1^2 \gamma \exp\left(\gamma L_1\right)}{\eta n}\sum_{i=1}^n\left(f^{\inf} - f_i^{\inf}\right).
			\end{align*}
			Next, we bound term $\gA_2$. Jensen's inequality and the tower property of the conditional expectation imply
			\begin{align*}
				\gA_2 & \leq \sqrt{\Exp{\sqnorm{\sum^{k}_{t=0}(1-\eta)^{k-t}U^{t+1}}}} \\
				& =  \sqrt{\sum^{k}_{t=0}(1-\eta)^{2(k-t)}\Exp{\sqnorm{U^{t+1}}}}.
			\end{align*}
		{\color{darkred}Moreover, the conditional independence and zero-mean properties in Assumption~\ref{ch6:assum:bounded_variance_chapter_6} make all cross-time and cross-client inner products vanish, so}
			\begin{align*}    				\gA_2 & \leq \sqrt{\sum^{k}_{t=0}\frac{(1-\eta)^{2(k-t)}}{n^2}\sum_{i=1}^n\Exp{\sqnorm{U_i^{t+1}}}} \\
    				& \overset{\eqref{ch6:eqn:bounded_variance_chapter_6}}{\leq} \sqrt{\sum^{k}_{t=0}(1-\eta)^{2(k-t)}\frac{\sigma^2}{n}}\\
				& \leq\frac{\sigma}{\sqrt{n}}\sqrt{\sum^{\infty}_{t=0}(1-\eta)^{2t}} = \frac{\sigma}{\sqrt{n\eta(2-\eta)}} \overset{\eta\in(0,1]}{\leq} \frac{\sigma}{\sqrt{n\eta}}.
			\end{align*}
			Therefore, plugging the derived upper-bounds for $\gA_1$, and for $\gA_2$ into \eqref{ch6:eqn:ef21sgdm:exp_norm_H_k_chapter_6}, we obtain 
			\begin{align*}
				\Exp{\norm{H^{k+1}}} & \leq (1-\eta)^{k+1} \Exp{\norm{H^0}} + \frac{2L_0 \gamma \exp\left(\gamma L_1\right)}{\eta}\\
				&\quad + 8L_1^2 \gamma \exp\left(\gamma L_1\right)\sum^{k}_{t=0}(1-\eta)^{k-t+1} \Exp{f(x^t) - f^{\inf}}\\
				&\quad +  \frac{8L_1^2 \gamma \exp\left(\gamma L_1\right)}{\eta n}\sum_{i=1}^n\left(f^{\inf} - f_i^{\inf}\right) + \frac{\sqrt{\eta}\sigma}{\sqrt{n}},
			\end{align*}
			which is equivalent to \eqref{ch6:eq:vk_nablafxk_bound_chapter_6}.
			
			To derive \eqref{ch6:eq:sum_vik_nablafixk_bound_chapter_6}, we make a step back to the recursion from \eqref{ch6:eq:H_i^k_recursion_chapter_6}, which implies
			\begin{align}
				\frac{1}{n}\sum_{i=1}^n \Exp{\norm{H_i^{k+1}}} & \leq \frac{1-\eta}{n}\sum_{i=1}^n \Exp{\norm{H_i^{k}}} + \underbrace{\frac{1-\eta}{n}\sum_{i=1}^n \Exp{\norm{G_i^{k}}}}_{=: \thesischsixcB_1} \notag\\
				&\quad + \underbrace{\frac{\eta}{n}\sum_{i=1}^n \Exp{\norm{U_i^{k+1}}}}_{=: \thesischsixcB_2}. \label{ch6:eq:bound_for_sum_tecnical_ineq_chapter_6}
			\end{align}
			Next, we derive the upper bounds for $\thesischsixcB_1$ and $\thesischsixcB_2$. For $\thesischsixcB_1$, we have
			\begin{align*}
				\thesischsixcB_1 & = \frac{1-\eta}{n}\sum_{i=1}^n \Exp{\norm{\nabla f_i(x^k) - \nabla f_i(x^{k+1})}}\\
				& \overset{\eqref{ch6:eqn:LzeroLoneSmooth_ineq_0_chapter_6}}{\leq} \frac{1-\eta}{n}\sum_{i=1}^n\Exp{(L_0 + L_1\norm{\nabla f_i(x^k)})\exp\left(L_1\norm{x^k - x^{k+1}}\right)\norm{x^k - x^{k+1}}}\\
				& \leq (1-\eta)L_0\gamma\exp(\gamma L_1) \\
				&\quad +  \frac{(1-\eta)L_1\gamma\exp(\gamma L_1)}{n}\sum_{i=1}^n\Exp{\norm{\nabla f_i(x^k)}}\\
				& \overset{\eqref{ch6:eq:sum_of_grad_norms_bound_chapter_6}}{\leq} 2(1-\eta)L_0\gamma\exp(\gamma L_1) + 8(1-\eta)L_1^2\gamma\exp(\gamma L_1)\Exp{f(x^k) - f^{\inf}}\\
				&\quad + \frac{8(1-\eta)L_1^2\gamma\exp(\gamma L_1)}{n}\sum_{i=1}^n\left(f^{\inf} - f_i^{\inf}\right),
			\end{align*}
			and for $\thesischsixcB_2$, we obtain
			\begin{align*}
				\thesischsixcB_2 & = \frac{\eta}{n}\sum_{i=1}^n \Exp{\norm{\nabla f_i(x^{k+1};\xi_i^{k+1}) - \nabla f_i(x^{k+1})}} \overset{\eqref{ch6:eqn:bounded_variance_chapter_6}}{\leq} \eta\sigma.
			\end{align*}
			Plugging the derived upper bounds for $\thesischsixcB_1$ and $\thesischsixcB_2$ into \eqref{ch6:eq:bound_for_sum_tecnical_ineq_chapter_6} and using $1-\eta \leq 1$, we get
			\begin{align*}
				\frac{1}{n}\sum_{i=1}^n \Exp{\norm{H_i^{k+1}}} & \leq \frac{1-\eta}{n}\sum_{i=1}^n \Exp{\norm{H_i^{k}}} + 2L_0\gamma\exp(\gamma L_1)\\
				&\quad + 8L_1^2\gamma\exp(\gamma L_1)\Exp{f(x^k) - f^{\inf}}\\
				&\quad + \frac{8L_1^2\gamma\exp(\gamma L_1)}{n}\sum_{i=1}^n\left(f^{\inf} - f_i^{\inf}\right) +  \eta\sigma,
			\end{align*}
			which is equivalent to \eqref{ch6:eq:sum_vik_nablafixk_bound_chapter_6}.

	\end{proof}

	\subsection{Proof of Theorem~\ref{ch6:thm:ef21_sgdm_chapter_6}}\label{ch6:app:sec:ef21_sgdm_chapter_6}
	
	Now, we are ready to prove Theorem~\ref{ch6:thm:ef21_sgdm_chapter_6}. For convenience, we introduce new notation: 
		{\color{darkred}
		\begin{align*}
			\delta^k&\eqdef\Exp{f(x^k)-f^{\inf}},&
			A_k&\eqdef\frac{1}{n}\sum_{i=1}^n
			\Exp{\norm{v_i^k-g_i^k}},\\
			B_k & \eqdef\Exp{\norm{v^k-\nabla f(x^k)}},&
			C_k& \\
			& \eqdef\frac{1}{n}\sum_{i=1}^n
			\Exp{\norm{v_i^k-\nabla f_i(x^k)}},\\
			\delta^{\inf}&\eqdef\frac{1}{n}\sum_{i=1}^n
			(f^{\inf}-f_i^{\inf}).&&
		\end{align*}}
		Using the new notation and noticing that $\Exp{\|v^k - g^k\|} \leq A_k$, we rewrite the results of Lemmas~\ref{ch6:lem:descent_lemma_ef21_sgdm_chapter_6}, \ref{ch6:lem:compression_error_ef21_sgdm_chapter_6}, and \ref{ch6:lem:momentum_error_ef21_sgdm_chapter_6} as
		\begin{align*}
			\delta^{k+1} & \leq \left(1 + 4L_1^2\gamma^2\exp(L_1\gamma)\right)\delta^k + 2\gamma A_k + 2\gamma B_k - \gamma\Exp{\norm{\nabla f(x^k)}}\\
			&\quad+ \gamma^2 \exp(L_1\gamma)\left(L_0 + 4L_1^2 \delta^{\inf}\right),\\
			A_{k+1} & \leq \sqrt{1 - \alpha} A_k + \eta\sqrt{1-\alpha}C_k + 8L_1^2\sqrt{1-\alpha}\eta\gamma\exp\left(\gamma L_1\right) \delta^k\\
			&\quad+ 2\sqrt{1-\alpha} \eta\gamma \exp\left(\gamma L_1\right)\left(L_0 + 4L_1^2 \delta^{\inf}\right) + \sqrt{1-\alpha}\eta\sigma,\\
			B_k & \leq (1-\eta)^kB_0 + \frac{\sqrt{\eta}\sigma}{\sqrt{n}} + \frac{2 \gamma \exp\left(L_1\gamma\right)}{\eta}\left(L_0 + 4L_1^2 \delta^{\inf}\right)\\
			&\quad+ 8L_1^2\gamma \exp\left(L_1\gamma\right)\sum\limits_{t=0}^{k-1}(1-\eta)^{k-t}\delta^t,\\
			C_{k+1} & \leq (1-\eta)C_k + 8L_1^2\gamma \exp\left(L_1\gamma\right)\delta^k + \eta\sigma \\
			&\quad + 2\gamma\exp\left(\gamma L_1\right)\left(L_0 + 4L_1^2 \delta^{\inf}\right).
		\end{align*}
		Moreover, since $\gamma = \frac{\gamma_0}{(K+1)^{3/4}}$ with $\gamma_0 \leq \frac{1}{2L_1}$, we have $\exp(L_1\gamma) \leq \exp(L_1\gamma_0) \leq 2$ and the above inequalities can be further simplified as
		{\color{darkred}
		\begin{align}
			\delta^{k+1}
			& \leq \left(1 + 8L_1^2\gamma^2\right)\delta^k
			+ 2\gamma A_k \notag \\
			&\quad + 2\gamma B_k
			- \gamma\Exp{\norm{\nabla f(x^k)}} \notag\\
			&\quad + 2\gamma^2 \left(L_0 + 4L_1^2 \delta^{\inf}\right),
			\label{ch6:eq:ef21_sgdm_main_ineq1_chapter_6}\\
			A_{k+1}
			&\leq \sqrt{1 - \alpha} A_k
			+ \eta\sqrt{1-\alpha}C_k
			+ 16L_1^2\sqrt{1-\alpha}\eta\gamma \delta^k \notag\\
			&\quad + 4\sqrt{1-\alpha}\eta\gamma
			\left(L_0 + 4L_1^2 \delta^{\inf}\right)
			+ \sqrt{1-\alpha}\eta\sigma,
			\label{ch6:eq:ef21_sgdm_main_ineq2_chapter_6}\\
			B_k
			&\leq (1-\eta)^kB_0 + \frac{\sqrt{\eta}\sigma}{\sqrt{n}}
			+ \frac{4\gamma}{\eta}\left(L_0 + 4L_1^2 \delta^{\inf}\right) \notag\\
			&\quad + 16L_1^2\gamma
			\sum_{t=0}^{k-1}(1-\eta)^{k-t}\delta^t,
			\label{ch6:eq:ef21_sgdm_main_ineq3_chapter_6}\\
			C_{k+1}
			& \leq (1-\eta)C_k \notag \\
			&\quad + 16L_1^2\gamma \delta^k
			+ \eta\sigma
			+ 4\gamma\left(L_0 + 4L_1^2 \delta^{\inf}\right). \label{ch6:eq:ef21_sgdm_main_ineq4_chapter_6}
		\end{align}
		}
		Next, we introduce the Lyapunov function $V_k$ defined for any $k \geq 0$ as
		\begin{equation*}
			V_{k} = \delta^k + aA_k + cC_k,
		\end{equation*}
		where $a \eqdef \frac{2\gamma}{1 - \sqrt{1-\alpha}}$ and $c\eqdef a\sqrt{1-\alpha}$. Then, using \eqref{ch6:eq:ef21_sgdm_main_ineq1_chapter_6}, \eqref{ch6:eq:ef21_sgdm_main_ineq2_chapter_6}, \eqref{ch6:eq:ef21_sgdm_main_ineq4_chapter_6}, we get
		\begin{align*}
			V_{k+1}
			& \leq \left(1 + 8L_1^2\gamma^2\right)\delta^k
			+ 2\gamma A_k \\
			&\quad + 2\gamma B_k
			- \gamma\Exp{\norm{\nabla f(x^k)}}\\
			&\quad+ 2\gamma^2 \left(L_0 + 4L_1^2 \delta^{\inf}\right)\\
			&\quad+ a\Bigl(\sqrt{1-\alpha}A_k
			+ \eta\sqrt{1-\alpha}C_k
			+ 16L_1^2\sqrt{1-\alpha}\eta\gamma\delta^k\Bigr)\\
			&\quad+ a\Bigl(4\sqrt{1-\alpha}\eta\gamma
			\left(L_0 + 4L_1^2 \delta^{\inf}\right)
			+ \sqrt{1-\alpha}\eta\sigma\Bigr)\\
			&\quad+ c\Bigl((1-\eta)C_k + 16L_1^2\gamma\delta^k
			+ \eta\sigma
			+ 4\gamma\left(L_0 + 4L_1^2 \delta^{\inf}\right)\Bigr).
		\end{align*}
		To proceed, we rearrange the terms and then use $c=a\sqrt{1-\alpha}$ and $\eta\leq1$:
		\begin{align*}
			V_{k+1}
			& \leq \Bigl(1 + 8L_1^2\gamma^2
			+ 16aL_1^2\sqrt{1-\alpha}\eta\gamma
			+ 16cL_1^2\gamma\Bigr)\delta^k\\
			&\quad+ \left(\frac{2\gamma}{a}+\sqrt{1-\alpha}\right)aA_k
			+ \left(\frac{a\eta\sqrt{1-\alpha}}{c}+1-\eta\right)cC_k\\
			&\quad+ 2\gamma B_k - \gamma\Exp{\norm{\nabla f(x^k)}}\\
			&\quad+ \Bigl(2\gamma^2+4a\sqrt{1-\alpha}\eta\gamma+4c\gamma\Bigr)
			\left(L_0+4L_1^2\delta^{\inf}\right)\\
			&\quad+ \eta\left(a\sqrt{1-\alpha}+c\right)\sigma\\
			& \leq
			\Bigl(1+8L_1^2\gamma^2
			+32aL_1^2\sqrt{1-\alpha}\gamma\Bigr)\delta^k\\
			&\quad+ \left(\frac{2\gamma}{a}+\sqrt{1-\alpha}\right)aA_k
			+ cC_k + 2\gamma B_k
			- \gamma\Exp{\norm{\nabla f(x^k)}}\\
			&\quad+ \Bigl(2\gamma^2+8a\sqrt{1-\alpha}\gamma\Bigr)
			\left(L_0+4L_1^2\delta^{\inf}\right)
			+2\eta a\sqrt{1-\alpha}\sigma.
		\end{align*}
		Since $a = \frac{2\gamma}{1 - \sqrt{1-\alpha}}$, we have $\frac{2\gamma}{a} + \sqrt{1-\alpha} = 1$ and
		\begin{align*}
			V_{k+1}
			& \leq \left(1 + 8L_1^2\gamma^2
			+ \frac{64L_1^2\gamma^2\sqrt{1-\alpha}}{1-\sqrt{1-\alpha}}\right)\delta^k\\
			&\quad+ aA_k+cC_k+2\gamma B_k
			-\gamma\Exp{\norm{\nabla f(x^k)}}\\
			&\quad+ \left(2\gamma^2
			+\frac{16\gamma^2\sqrt{1-\alpha}}{1-\sqrt{1-\alpha}}\right)
			\left(L_0+4L_1^2\delta^{\inf}\right)\\
			&\quad+ \frac{4\gamma\eta\sqrt{1-\alpha}}{1-\sqrt{1-\alpha}}\,\sigma\\
			& \leq \left(1+\frac{64L_1^2\gamma^2}{1-\sqrt{1-\alpha}}\right)V_k
			+2\gamma B_k-\gamma\Exp{\norm{\nabla f(x^k)}}\\
			&\quad+ \frac{16\gamma^2}{1-\sqrt{1-\alpha}}
			\left(L_0+4L_1^2\delta^{\inf}\right)
			+\frac{4\gamma\eta\sqrt{1-\alpha}}{1-\sqrt{1-\alpha}}\,\sigma.
		\end{align*}
		Next, we bound $B_k$ using \eqref{ch6:eq:ef21_sgdm_main_ineq3_chapter_6} and $\delta^k \leq V_k$:
		\begin{align*}
			V_{k+1}
			& \leq \left(1+\frac{64L_1^2\gamma^2}{1-\sqrt{1-\alpha}}\right)V_k
			+32L_1^2\gamma^2\sum_{t=0}^{k-1}(1-\eta)^{k-t}V_t\\
			&\quad+2\gamma(1-\eta)^kB_0
			-\gamma\Exp{\norm{\nabla f(x^k)}}\\
			&\quad+\left(\frac{16\gamma^2}{1-\sqrt{1-\alpha}}
			+\frac{8\gamma^2}{\eta}\right)
			\left(L_0+4L_1^2\delta^{\inf}\right)\\
			&\quad+\left(\frac{4\gamma\eta\sqrt{1-\alpha}}{1-\sqrt{1-\alpha}}
			+\frac{2\gamma\sqrt{\eta}}{\sqrt n}\right)\sigma.
		\end{align*}
		Summing up the above inequality with weights $\beta_k \eqdef \left(1 + \frac{64L_1^2\gamma^2}{1-\sqrt{1-\alpha}} + \frac{32L_1^2\gamma^2}{\eta}\right)^{-(k+1)}$ for $k = 0,\ldots,K$ and denoting $S_K \eqdef \sum_{k=0}^K\beta_k$ and $\beta_{-1}\eqdef 1$, we get
		\begin{align*}
			\sum_{k=0}^K\beta_k V_{k+1}
			& \leq \sum_{k=0}^K
			\left(1+\frac{64L_1^2\gamma^2}{1-\sqrt{1-\alpha}}\right)\beta_kV_k\\
			&\quad+32L_1^2\gamma^2\sum_{k=0}^K\beta_k
			\sum_{t=0}^{k-1}(1-\eta)^{k-t}V_t\\
			&\quad+2\gamma B_0\sum_{k=0}^K(1-\eta)^k\beta_k
			-\gamma\sum_{k=0}^K\beta_k\Exp{\norm{\nabla f(x^k)}}\\
			&\quad+S_K\left(\frac{16\gamma^2}{1-\sqrt{1-\alpha}}
			+\frac{8\gamma^2}{\eta}\right)
			\times\left(L_0+4L_1^2\delta^{\inf}\right)\\
			&\quad+S_K\left(\frac{4\gamma\eta\sqrt{1-\alpha}}{1-\sqrt{1-\alpha}}
			+\frac{2\gamma\sqrt{\eta}}{\sqrt n}\right)\sigma.
		\end{align*}
		By definition of $\beta_k$, we have $\beta_k \leq \beta_{k-1}$ and, in particular, $\beta_k \leq 1$ for all $k\geq 0$. Using these inequalities, we continue the derivation as follows:
		\begin{align*}
			\sum_{k=0}^K\beta_k V_{k+1}
			& \leq \sum_{k=0}^K
			\left(1+\frac{64L_1^2\gamma^2}{1-\sqrt{1-\alpha}}\right)\beta_kV_k\\
			&\quad+32L_1^2\gamma^2\sum_{k=0}^K\sum_{t=0}^{k-1}
			(1-\eta)^{k-t}\beta_tV_t\\
			&\quad+2\gamma B_0\sum_{k=0}^K(1-\eta)^k
			-\gamma\sum_{k=0}^K\beta_k\Exp{\norm{\nabla f(x^k)}}\\
			&\quad+S_K\left(\frac{16\gamma^2}{1-\sqrt{1-\alpha}}
			+\frac{8\gamma^2}{\eta}\right)
			\times\left(L_0+4L_1^2\delta^{\inf}\right)\\
			&\quad+S_K\left(\frac{4\gamma\eta\sqrt{1-\alpha}}{1-\sqrt{1-\alpha}}
			+\frac{2\gamma\sqrt{\eta}}{\sqrt n}\right)\sigma\\
			& \leq \sum_{k=0}^K
			\left(1+\frac{64L_1^2\gamma^2}{1-\sqrt{1-\alpha}}\right)\beta_kV_k\\
			&\quad+32L_1^2\gamma^2
			\left(\sum_{t=0}^{\infty}(1-\eta)^t\right)
			\left(\sum_{k=0}^K\beta_kV_k\right)\\
			&\quad+2\gamma B_0\sum_{k=0}^{\infty}(1-\eta)^k
			-\gamma S_K\min_{k=0,\ldots,K}\Exp{\norm{\nabla f(x^k)}}\\
			&\quad+S_K\left(\frac{16\gamma^2}{1-\sqrt{1-\alpha}}
			+\frac{8\gamma^2}{\eta}\right)
			\times\left(L_0+4L_1^2\delta^{\inf}\right)\\
			&\quad+S_K\left(\frac{4\gamma\eta\sqrt{1-\alpha}}{1-\sqrt{1-\alpha}}
			+\frac{2\gamma\sqrt{\eta}}{\sqrt n}\right)\sigma\\
			& = \sum_{k=0}^K
			\underbrace{\left(1+\frac{64L_1^2\gamma^2}{1-\sqrt{1-\alpha}}
			+\frac{32L_1^2\gamma^2}{\eta}\right)\beta_k}_{=\beta_{k-1}}V_k\\
			&\quad+\frac{2\gamma B_0}{\eta}
			-\gamma S_K\min_{k=0,\ldots,K}\Exp{\norm{\nabla f(x^k)}}\\
			&\quad+S_K\left(\frac{16\gamma^2}{1-\sqrt{1-\alpha}}
			+\frac{8\gamma^2}{\eta}\right)
			\times\left(L_0+4L_1^2\delta^{\inf}\right)\\
			&\quad+S_K\left(\frac{4\gamma\eta\sqrt{1-\alpha}}{1-\sqrt{1-\alpha}}
			+\frac{2\gamma\sqrt{\eta}}{\sqrt n}\right)\sigma.
		\end{align*}
		Rearranging the terms and dividing both sides of the above inequality by $\gamma S_K$, we obtain
		\begin{equation}
			{\color{darkred}
						\begin{aligned}
				\min_{k=0,\ldots,K}\Exp{\norm{\nabla f(x^k)}}
				&\leq\frac{1}{\gamma S_K}\sum_{k=0}^K
				(\beta_{k-1}V_k-\beta_kV_{k+1})
				+\frac{2B_0}{\eta S_K}\\
				&\quad+\left(\frac{16\gamma}{1-\sqrt{1-\alpha}}
				+\frac{8\gamma}{\eta}\right)
				\left(L_0+4L_1^2\delta^{\inf}\right)\\
				&\quad+\left(\frac{4\eta\sqrt{1-\alpha}}{1-\sqrt{1-\alpha}}
				+\frac{2\sqrt{\eta}}{\sqrt n}\right)\sigma\\
				&\leq\frac{V_0}{\gamma S_K}+\frac{2B_0}{\eta S_K}\\
				&\quad+\left(\frac{16\gamma}{1-\sqrt{1-\alpha}}
				+\frac{8\gamma}{\eta}\right)
				\left(L_0+4L_1^2\delta^{\inf}\right)\\
				&\quad+\left(\frac{4\eta\sqrt{1-\alpha}}{1-\sqrt{1-\alpha}}
				+\frac{2\sqrt{\eta}}{\sqrt n}\right)\sigma.
			\end{aligned}}
			\label{ch6:eq:pre_final_bound_EF21_SGDM_chapter_6}
		\end{equation}
		where in the last inequality we use $V_{K+1} \geq 0$ and $\beta_{-1} = 1$. Next, we estimate $S_K$:
		\begin{align}
			S_K & = \sum_{k=0}^K \beta_k \geq (K+1)\beta_K \notag \\
			& = \frac{K+1}{\left(1 + \frac{64L_1^2\gamma^2}{1-\sqrt{1-\alpha}} + \frac{32L_1^2\gamma^2}{\eta}\right)^{K+1}} \notag\\
			& \geq \frac{K+1}{\exp\left(\frac{64L_1^2\gamma^2 (K+1)}{1-\sqrt{1-\alpha}} + \frac{32L_1^2\gamma^2(K+1)}{\eta}\right)}. \label{ch6:eq:S_K_lower_bound_chapter_6}
		\end{align}
		{\color{darkred}Set
		\begin{align*}
			\eta&=\frac{1}{(K+1)^{1/2}},
			\qquad
			\gamma=\frac{\gamma_0}{(K+1)^{3/4}},\\
			\gamma_0&\leq\frac{1}{16L_1}
			\min\!\left\{(K+1)^{1/2}(1-\sqrt{1-\alpha}),1\right\}.
		\end{align*}
		Then
		\[
			\frac{32L_1^2\gamma^2(K+1)}{\eta}\leq\frac14,
			\qquad
			\frac{64L_1^2\gamma^2(K+1)}{1-\sqrt{1-\alpha}}\leq\frac14.
		\]
		Hence, \eqref{ch6:eq:S_K_lower_bound_chapter_6} yields
		$S_K\geq(K+1)/\exp(1/2)\geq(K+1)/2$. Substituting these
		parameter values into
		\eqref{ch6:eq:pre_final_bound_EF21_SGDM_chapter_6} gives}
		{\color{darkred}
		\begin{align*}
			\MoveEqLeft[1] \min\limits_{k=0,\ldots,K} \Exp{\norm{\nabla f(x^k)}} \\
			& \leq \frac{2V_0}{\gamma_0 (K+1)^{1/4}} + \frac{4B_0}{(K+1)^{1/2}}\notag\\
			&\quad + \left(\frac{16\gamma_0}{(1-\sqrt{1-\alpha})(K+1)^{3/4}} + \frac{8\gamma_0}{(K+1)^{1/4}}\right)\left(L_0 + 4L_1^2 \delta^{\inf}\right)\notag\\
			&\quad  + \left(\frac{4\sqrt{1-\alpha}}{(1-\sqrt{1-\alpha})(K+1)^{1/2}} + \frac{2}{\sqrt{n}(K+1)^{1/4}}\right)\sigma. \notag
		\end{align*}}
		{\color{darkred}
		For convenience, define $C_\alpha \eqdef 1-\sqrt{1-\alpha}$. By the definition of $V_0$,
		\begin{equation*}
			\frac{2V_0}{\gamma_0 (K+1)^{1/4}}
			=
			\frac{2\delta^0}{\gamma_0 (K+1)^{1/4}}
			+\frac{4A_0}{C_\alpha(K+1)}
			+\frac{4(1-C_\alpha)C_0}{C_\alpha(K+1)}.
		\end{equation*}
		At initialization, $v_i^0=\nabla f_i(x^0;\xi_i^0)$,
		$g_i^0=\thesischsixcC^0(v_i^0)$, and the samples
		$\{\xi_i^0\}_{i=1}^n$ are independent. Hence,
		\begin{align*}
			A_0 & \leq \frac{\sqrt{1-\alpha}}{n}\sum_{i=1}^n\Exp{\norm{v_i^0}} \\
			& \leq \sqrt{1-\alpha}\left(
				\frac{1}{n}\sum_{i=1}^n\norm{\nabla f_i(x^0)}+\sigma
			\right)\\
			&\leq \sqrt{1-\alpha}\left(
				8L_1(\delta^0+\delta^{\inf})+\frac{L_0}{L_1}+\sigma
			\right),\\
			C_0
			&=\frac{1}{n}\sum_{i=1}^n
			\Exp{\norm{\nabla f_i(x^0;\xi_i^0)-\nabla f_i(x^0)}}\leq\sigma,\\
			B_0
			&=\Exp{\norm{\frac{1}{n}\sum_{i=1}^n
			\left(\nabla f_i(x^0;\xi_i^0)-\nabla f_i(x^0)\right)}}\\
			&\leq
			\sqrt{\frac{1}{n^2}\sum_{i=1}^n
			\Exp{\sqnorm{\nabla f_i(x^0;\xi_i^0)-\nabla f_i(x^0)}}}
			\leq\frac{\sigma}{\sqrt n}.
		\end{align*}
		The first estimate uses contractivity, the triangle inequality, the
		bounded-variance assumption, and
		\eqref{ch6:eq:sum_of_grad_norms_bound_chapter_6}. The estimate for
		$B_0$ uses Jensen's inequality and the conditional independence and
		zero-mean properties from Assumption~\ref{ch6:assum:bounded_variance_chapter_6}.
		Substituting these estimates into
		\eqref{ch6:eq:pre_final_bound_EF21_SGDM_chapter_6} yields
		\begin{align*}
			\MoveEqLeft[1]\min_{0\leq k\leq K}\Exp{\norm{\nabla f(x^k)}}\\
			&\leq
			\gO\left(
			\frac{\delta^0/\gamma_0+\sigma/\sqrt n
			+\gamma_0(L_0+L_1^2\delta^{\inf})}{(K+1)^{1/4}}\right)\\
			&\quad+\gO\Bigg(
			\frac{\sqrt{1-\alpha}}{\alpha}\Bigg[
			\frac{L_1(\delta^0+\delta^{\inf})
				+L_0/L_1+\sigma}{K+1}
			+\frac{\sigma}{(K+1)^{1/2}}\\
			&\qquad\qquad
			+\frac{\gamma_0(L_0+L_1^2\delta^{\inf})}{(K+1)^{3/4}}
			\Bigg]\Bigg).
		\end{align*}
		Here we used
		$\frac{1-C_\alpha}{C_\alpha}
		=\frac{\sqrt{1-\alpha}}{1-\sqrt{1-\alpha}}
		\leq\frac{2\sqrt{1-\alpha}}{\alpha}$
		and
		$\frac{1}{C_\alpha}
		=\frac{1+\sqrt{1-\alpha}}{\alpha}\leq\frac{2}{\alpha}$; for the term with $\gamma_0/C_\alpha$ we also used $\frac{1}{C_\alpha}=1+\frac{1-C_\alpha}{C_\alpha}$ and $(K+1)^{-3/4}\leq(K+1)^{-1/4}$.
		This proves the theorem.
		}

	\clearpage

    {\color{black}
    \section{Extension to Strongly Convex and Convex Problems}

    Our current analysis for \normEFname{} and  \normEFname{SGDM}, which are initially developed for minimizing non-convex functions, can be extended to  strongly convex and convex functions.

    \paragraph{Strongly convex problems.}

    We can extend the convergence for \normEFname{} and  \normEFname{SGDM} to minimize strongly convex functions.
    Applying the $\mu$-strong convexity condition of the function $f$, i.e. $\sqnorm{\nabla f(x^k)} \geq 2\mu(f(x^k)-f(x^\star))$, where $x^\star = \arg\min_{x\in\R^d} f(x)$, into the convergence bounds in Theorems~\ref{ch6:thm:ef21_chapter_6} and~\ref{ch6:thm:ef21_sgdm_chapter_6}
    yields the convergence results in $\min_{k=0,1,\ldots,K} \Exp{\sqrt{f(x^k)-f(x^\star)}}.$
    However, these results do not imply the standard exponential convergence typically expected in strongly convex problems. 
    This theoretical gap suggests a need for new analytical techniques, which involves tighter Lyapunov functions or more refined descent inequalities tailored to strongly convex functions.

    \paragraph{Convex problems.}

    We can extend the convergence for minimizing convex functions. 
    This can be achieved by assuming that the iterates $\{x^k\}$ satisfy $\norm{x^k-x^\star}\leq R$ for some $R>0$.
    Hence, the convexity of the function $f$ implies that
    \begin{equation*}
        f(x^k)-f(x^\star) \leq \norm{\nabla f(x^k)} \norm{x^k-x^\star} \leq R \norm{\nabla f(x^k)}. 
    \end{equation*}
    Applying the above inequality to Theorems~\ref{ch6:thm:ef21_chapter_6} and~\ref{ch6:thm:ef21_sgdm_chapter_6} yields the convergence bounds in $\min_{k=0,1,\ldots,K} \Exp{{f(x^k)-f(x^\star)}}$.
    
	}

	\clearpage
	
	\section{Additional Experimental Results}

	In this section, we provide additional results for minimizing nonconvex polynomial functions, and for training the ResNet-20 model over the CIFAR-10 dataset.

	\subsection{Minimization of Nonconvex Polynomial Functions}\label{ch6:app:sub:minimize_nonconvex_poly_chapter_6}
	
	We ran \normEFname{} and \algname{EF21} in a single-node setting ($n=1$) for solving the following problem:
	\begin{align}\label{ch6:eqn:simple_polynomials_chapter_6}
		\min _{x \in \mathbb{R}^d}\biggl\{f(x)\eqdef\underbrace{\sum_{i=1}^d a_i x_i^4}_{=: g(x)}+\underbrace{\lambda \sum_{i=1}^d \frac{x_i^2}{1+x_i^2}}_{=: h(x)}\biggr\},
	\end{align}
	where $a_i > 0$, $i = 1, \ldots, d$, $\lambda > 0$.
	
	Let us show that $f(x)$ is non-convex (for the specific choice of $a_i$) and $(L_0, L_1)$-smooth. First, we  prove that $f(x)$ is non-convex. Indeed,
	\begin{align*}
		\nabla^{2} f(x) & =\nabla^{2} g(x) + \nabla^{2} h(x) \\
		& =12 \operatorname{diag}\left\{a_1 x_{1}^{2}, \ldots, a_d x_{d}^{2}\right\} \\
		&\quad + 2 \lambda \operatorname{diag}\biggl\{\frac{1-3 x_{1}^{2}}{\left(1+x_{1}^{2}\right)^{3}}, \ldots, \frac{1-3 x_{d}^{2}}{\left(1+x_{d}^{2}\right)^{3}}\biggr\},
	\end{align*} 
	{\color{darkred}is not positive semidefinite if, for example, we choose
	$a_i=\frac{\lambda}{24}$ and $x_i=\pm1/\sqrt2$ for every $i$. Indeed,
	each diagonal entry then equals
	$\lambda/4-8\lambda/27=-5\lambda/108<0$.}
	
	Second, we  find $L_0, L_1 > 0$ such that \[\norm{\nabla^{2} f(x)} \leq L_{0}+L_{1}\norm{\nabla f(x)}, \quad \forall x\in\mathbb{R}^d.\] This condition is equivalent to Assumption~\ref{ch6:assum:LzeroLoneSmooth_chapter_6} (generalized smoothness) with $L_0,L_1$ \parencite[Theorem 1]{chen2023generalized}.
	{\color{darkred}Let us keep $a_i=\lambda/24$ for every $i$, fix some $L_1 > 0$, and choose $L_0 = \frac{9 \lambda d^2}{2 L_1^2} + 2 \lambda$. Since $\norm{\nabla^2 h(x)}\leq2\lambda$,}
	\begin{align*}
		\norm{\nabla^{2} f(x)}=\norm{\nabla^{2} g(x)+\nabla^{2} h(x)} & \leq \norm{\nabla^{2} g(x)}+\norm{\nabla^{2} h(x)} \\
		&\leq 12 \sqrt{a_1^2 x_{1}^{4}+\ldots+a_d^2 x_{d}^{4}}+2 \lambda\\
		&\leq 12\left(a_1 x_{1}^{2}+\ldots+a_d x_{d}^{2}\right) + 2\lambda.
	\end{align*}
	Also, notice that
	\begin{align*}
		\MoveEqLeft[1] \norm{\nabla f(x)}=\norm{\nabla g(x)+\nabla h(x)} \\
		&=\sqrt{\left(4a_1x_{1}^{2} + \frac{2 \lambda}{(1+x_1^2)^2}\right)^2 x_1^2 +\ldots+\left(4a_dx_d^{2} + \frac{2 \lambda}{(1+x_d^2)^2}\right)^2 x_d^2}\\
		& \geq 4\sqrt{a_1^2 x_1^6+\ldots+a_d^2 x_d^6} \\
		& \overset{(*)}{\geq} \frac{4}{\sqrt{d}}\left(a_1 \left|x_1\right|^3+\ldots+a_d \left|x_d\right|^{3}\right),
	\end{align*}
	where (*) results from the fact that $\norm{x}_1 \leq \sqrt{d}\norm{x}$ for $x\in\R^d$. 
	Our goal is to show that 
	\begin{equation*}
		12\left(a_1 x_{1}^{2}+\ldots+a_d x_{d}^{2}\right) \leq \tilde{L}_{0}+\frac{4 L_1}{\sqrt{d}}\left(a_1 \left|x_{1}\right|^{3}+\ldots+a_d \left|x_d\right|^{3}\right),
	\end{equation*}
	where $\tilde{L}_{0}=L_{0}-2 \lambda$. To show this, we consider two cases: if $\left|x_{i}\right| \leq \frac{3 \sqrt{d}}{L_{1}}$, and otherwise. 
	
	\begin{enumerate}
		\item {\color{darkred}If $\left|x_i\right|\leq3\sqrt d/L_1$ for every
		$i\in[d]$, then $12a_ix_i^2\leq108a_id/L_1^2$. Therefore,
		\[
			12\sum_{i=1}^d a_ix_i^2
			\leq\frac{108d}{L_1^2}\sum_{i=1}^d a_i
			=\frac{108\lambda d^2}{24L_1^2}
			=\widetilde L_0.
		\]}
		\item {\color{darkred}For every index $j$ with
		$|x_j|>3\sqrt d/L_1$,}
		\[
			{\color{darkred}
			12a_jx_j^2<\frac{4L_1}{\sqrt d}a_j|x_j|^3.}
		\]
		{\color{darkred}The contribution of the remaining indices, for which
		$|x_i|\leq3\sqrt d/L_1$, is at most $\widetilde L_0$.}
	\end{enumerate}
	In conclusion, $f(x)$ is $(L_0, L_1)$-smooth, where $L_1$ is any positive constant and $L_0 = \frac{9 \lambda d^2}{2 L_1^2} + 2 \lambda$.
	
	{\color{darkred}Additionally, on the box $\{x:\lvert x_i\rvert\leq D,\ i\in[d]\}$, $f$ is $L$-smooth with $L=\frac{\lambda\sqrt d\,D^2}{2}+2\lambda$. If $\left|x_i\right| \leq D$ for all $i=1, \ldots, d$, then}
	\begin{align*}
		\left\|\nabla^{2} f(x)\right\| \leq 12 \sqrt{a_1^2 x_{1}^{4}+\ldots+a_d^2 x_{d}^{4}}+2 \lambda \leq \frac{ \lambda \sqrt{d} D^2}{2} + 2 \lambda = L.
	\end{align*}
	{\color{darkred}In the experiments, $D=20$ is used only as a heuristic local estimate based on the scale of the initial point. It is not a verified global smoothness bound: the theoretical \algname{EF21} stepsize based on this $L$ is valid only while all iterates remain in the displayed box.}

	{\color{darkred}In the following experiments, we used a top-$s$ sparsifier with $s=1$ and $\alpha=s/d$, setting $d=4$, $L_1\in\{1,4,8\}$, and $L_0=4$ (adjusting $\lambda$ to maintain a constant $L_0$). The coordinates of $x^0$ were drawn independently as $x_i^0\sim\mathcal N(20,1)$, and $D=20$ was used as the local heuristic described above.}
	For \algname{EF21}, we set $\gamma_k = \frac{1}{L + L \sqrt{\frac{\beta}{\theta}}}$, using $\theta = 1 - \sqrt{1 - \alpha}$ and $\beta = \frac{1 - \alpha}{1 - \sqrt{1 - \alpha}}$, according to Theorem 1 of \textcite{richtarik2021ef21}. 
	For \normEFname{}, we chose $\gamma_k = \frac{1}{2c_1}$ with $c_1 = \frac{L_1}{2} + 2\frac{\sqrt{1 - \alpha} L_1}{1 - \sqrt{1 - \alpha}}$ from Theorem~\ref{ch6:thm:normalized_ef21_singlenode_chapter_6}, and $\gamma_k = \frac{\gamma_0}{\sqrt{K + 1}}$ with $\gamma_0 > 0$, as specified in Theorem~\ref{ch6:thm:ef21_chapter_6} with $n=1$.

	\paragraph{The impact of $\gamma_0$ and $K$ on the convergence of \normEFname{}.}
	First, we investigate the impact of $\gamma_0$ and $K$ on the convergence of \normEFname{}.
	We evaluated $\gamma_0$ from the set $\{0.1, 1, 10\}$, and plotted the histogram representing the number of iterations required to achieve the target accuracy of $\norm{\nabla f(x)}^2 < \varepsilon$ with $\varepsilon = 10^{-4}$, using the stepsize rule $\gamma = \frac{\gamma_0}{\sqrt{K+1}}$. 
	For each $\gamma_0$, we determined $K$ as the minimum number of iterations required to achieve the desired accuracy, found through a grid search with step sizes of 500 for $\gamma_0 \in \{1, 10\}$ and $5000$ for $\gamma_0 = 0.1$.
	\begin{figure}[h]
		\centering
		\includegraphics[width=0.5\textwidth]{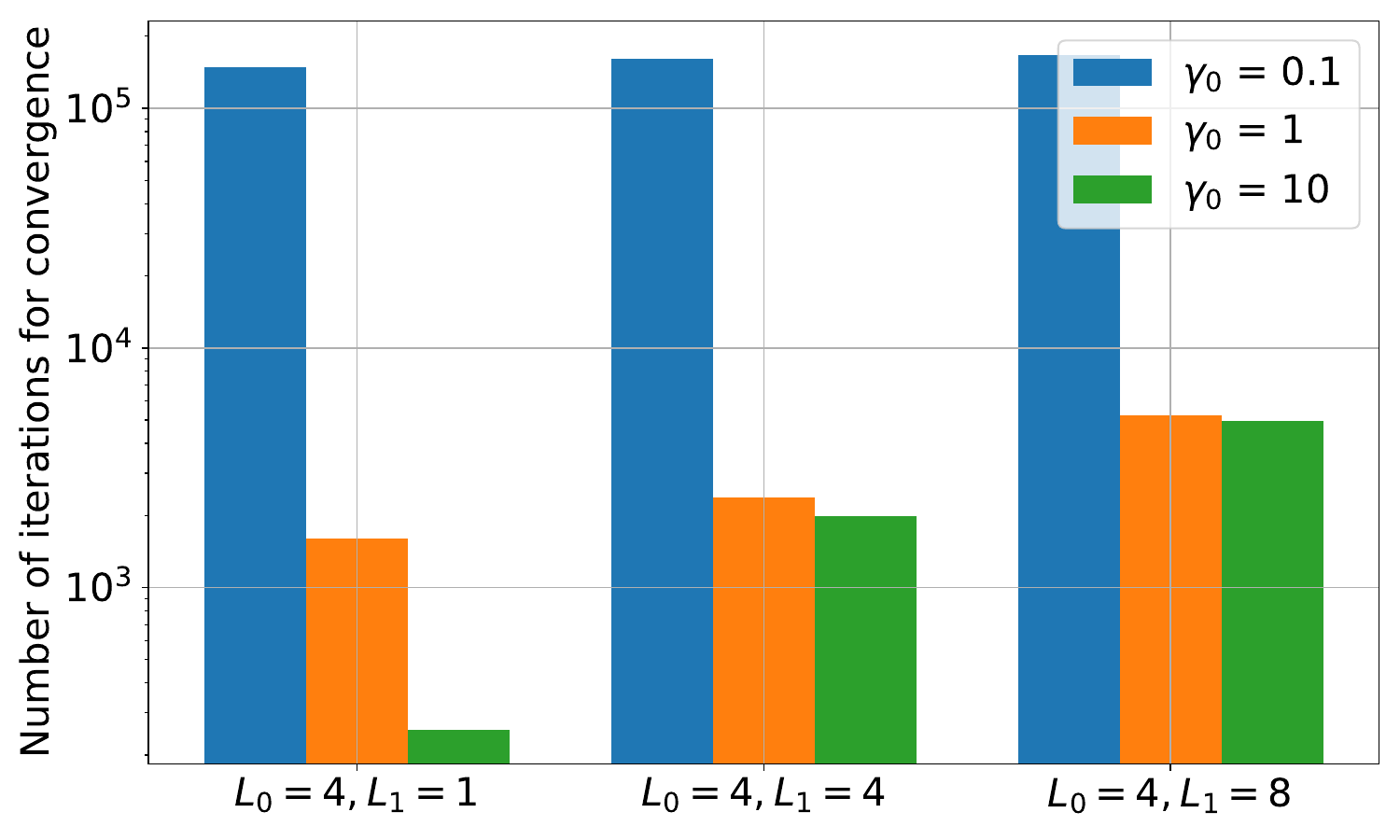}
		\caption{Number of iterations required to achieve the desired accuracy, $\norm{\nabla f(x)}^2 < \varepsilon$, $\varepsilon = 10^{-4}$, using \normEFname{}  with $\gamma = \frac{\gamma_0}{\sqrt{K+1}}$ for different values of $L_0$ and $L_1$.}
		\label{ch6:fig:histogram_chapter_6}
	\end{figure}
	From Figure~\ref{ch6:fig:histogram_chapter_6}, for small values of $\gamma_0$, such as $0.1$,  significantly more iterations are required to reach convergence compared to $\gamma_0$ values of $1$ and $10$, which show similar performance (with the exception of the $L_0=4$, $L_1=1$ case, where  $\gamma_0=10$ converges faster). Based on this observation, we use $\gamma_0 = 1$ in all subsequent experiments and adjust only $K$ to achieve convergence,  identifying the minimum number of iterations needed to reach the target accuracy through a grid search with a step size of $500$. 
	
	\paragraph{Comparisons between \algname{EF21} and \normEFname{}.}
	Next, we evaluate the performance of \algname{EF21} and \normEFname{} for a fixed \( L_0 = 4 \) and varying \( L_1 \) values of \{1, 4, 8\}. From Figure~\ref{ch6:fig:exp1_chapter_6}, \normEFname{}, regardless of the chosen stepsize \(\gamma\), achieves the desired accuracy \(\norm{\nabla f(x)}^2 < \varepsilon\) with \(\varepsilon = 10^{-4}\) faster than \algname{EF21}. Initially, however, \algname{EF21} converges more quickly, likely because \normEFname{} employs normalized gradients, which can be slower at the start due to the large gradients when the initial point is far from the stationary point. Moreover, as \( L_1 \) increases, both methods show slower convergence.

	\subsection{ResNet20 Training over CIFAR-10}

	{\color{darkred}We included additional results for training ResNet20 on CIFAR-10 with \algname{EF21} and \normEFname{}. The settings are the same as in Section~\ref{ch6:subsec:resnet_20_chapter_6}, except that the top-$s$ sparsifier uses $s\in\{0.01d,0.05d\}$. Figures~\ref{ch6:fig:plot_exp_DL_2_chapter_6} and~\ref{ch6:fig:plot_exp_DL_3_chapter_6} show that \normEFname{} attains higher accuracy across these sparsification levels.}

	\begin{figure}[h]
		\centering
		\includegraphics[width=\textwidth,viewport=0 0 2293 831,clip]{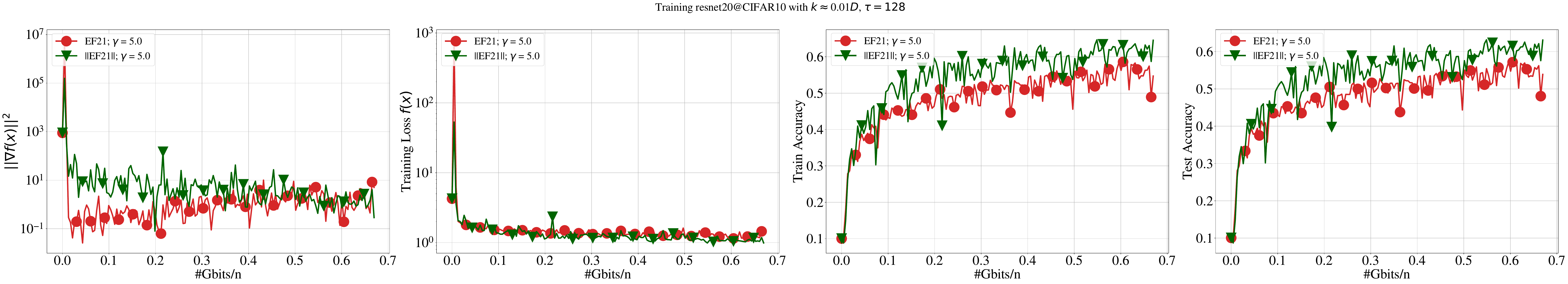}\\[0.5em]
		\includegraphics[width=\textwidth,viewport=2293 0 4586 831,clip]{plot_exp_DL_2.pdf}
		\caption[ResNet20 on CIFAR-10 with a top-$s$ sparsifier, $s=0.01d$]{{\color{darkred}ResNet20 training on CIFAR-10 using \algname{EF21} and \normEFname{} with the same stepsize $\gamma=5$ and a top-$s$ sparsifier with $s=0.01d$.}}
	\label{ch6:fig:plot_exp_DL_2_chapter_6}
	\end{figure}
	
	\begin{figure}[h]
	\centering
	\includegraphics[width=\textwidth,viewport=0 0 2281 831,clip]{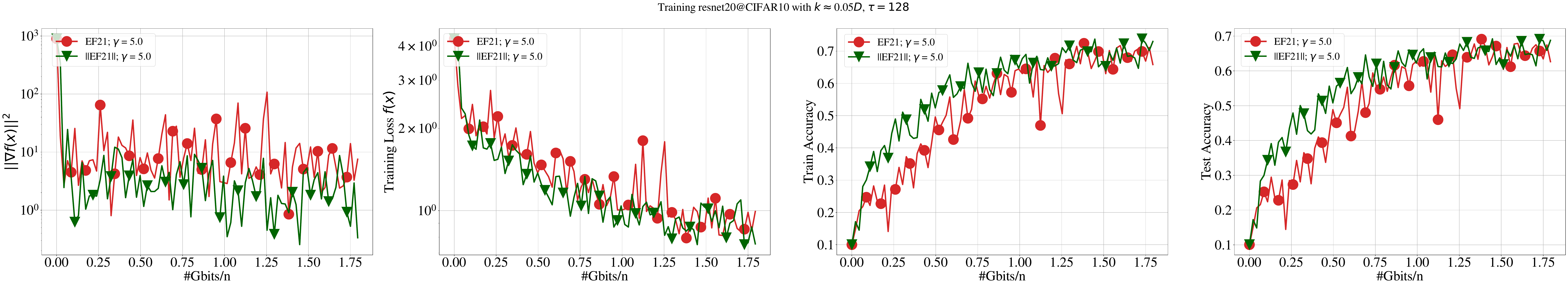}\\[0.5em]
	\includegraphics[width=\textwidth,viewport=2281 0 4563 831,clip]{plot_exp_DL_3.pdf}
	\caption[ResNet20 on CIFAR-10 with a top-$s$ sparsifier, $s=0.05d$]{{\color{darkred}ResNet20 training on CIFAR-10 using \algname{EF21} and \normEFname{} with the same stepsize $\gamma=5$ and a top-$s$ sparsifier with $s=0.05d$.}}
	\label{ch6:fig:plot_exp_DL_3_chapter_6}
\end{figure}

\clearpage

\section{{\color{darkred}Derivation of Smoothness Parameters for Logistic Regression}}\label{ch6:app:smooth_logistic_chapter_6}

In this section, we prove the generalized smoothness parameters $L_0,L_1$ for logistic regression problems with a nonconvex regularizer, which are the following problems 
\begin{align*}
\min _{x \in \mathbb{R}^d}\biggl\{f(x)\eqdef\frac{1}{n}\sum_{i=1}^n f_i(x)\eqdef \frac{1}{n}\sum_{i=1}^n \underbrace{\log(1+\exp(-b_i a_i^T x))}_{=:\tilde f_i(x)} + \underbrace{\lambda \sum_{j=1}^d \frac{x_j^2}{1+x_j^2}}_{=:h(x)}\biggr\},
\end{align*}
where $a_i \in \R^d$ is the $i^{\rm th}$ feature vector of matrix $A$ with its class label $b_i \in \{-1,1\}$, $\lambda>0$. 

First, we can prove that $f(x)$ is $L$-smooth with $L= \frac{1}{4n} \|A\|^2 + 2 \lambda$, and that each $f_i(x)$ is $\hat L_i$-smooth with $\hat L_i=\frac{1}{4} \|a_i\|^2 + 2 \lambda$.

{\color{darkred}Next, we show that each $f_i$ is generalized smooth, whenever
the Hessian exists, with}
\[
	{\color{darkred}
	L_0\eqdef2\lambda+\lambda\sqrt d\max_i\norm{a_i},
	\qquad
	L_1\eqdef\max_i\norm{a_i}.}
\]
By the fact that
\begin{align*}
\nabla \tilde f_i(x) &= - \frac{\exp(-b_i a_i^Tx)}{1+\exp(-b_i a_i^T x)} b_i a_i \quad \text{and}\\
\nabla^2 \tilde f_i(x) &= \frac{\exp(-b_i a_i^Tx)}{(1+\exp(-b_i a_i^T x))^2} b_i^2 a_i a_i^T, 
\end{align*}
we have 
\begin{align}
\norm{\nabla^2 \tilde f_i(x)} 
& \overset{b_i \in \{-1,1\} }{=} \frac{\exp(-b_i a_i^Tx)}{(1+\exp(-b_i a_i^T x))^2} \lambda_{\max}(a_i a_i^T)  \nonumber \\
& = \frac{\exp(-b_i a_i^Tx)}{(1+\exp(-b_i a_i^T x))^2} \sqnorm{a_i}  \nonumber \\
& = \frac{\norm{a_i}}{1+\exp(-b_ia_i^Tx)} \norm{\nabla \tilde f_i(x)} \nonumber \\
& \leq \norm{a_i}\norm{\nabla \tilde f_i(x)} \label{ch6:eq:hession_bound_chapter_6}.
\end{align}
After adding the nonconvex regularizer $h(x)$, we can show the following inequalities:
\begin{align}
\norm{\nabla^2 f_i(x)} 
& \leq \norm{\nabla^2 \tilde f_i(x)}+\norm{\nabla^2 h(x)} \nonumber\\
& \leq \norm{\nabla^2 \tilde f_i(x)} + 2\lambda \label{ch6:eq:norm_ineq_chapter_6},
\end{align}
and
\begin{equation}
{\color{darkred}
\begin{aligned}
\norm{\nabla f_i(x)}
&\geq\norm{\nabla\widetilde f_i(x)}-\norm{\nabla h(x)}\\
&=\norm{\nabla\widetilde f_i(x)}
-\sqrt{\sum_{j=1}^d\left(\frac{2\lambda x_j}{(1+x_j^2)^2}\right)^2}\\
&\geq\norm{\nabla\widetilde f_i(x)}-\lambda\sqrt d.
\end{aligned}}
\label{ch6:eq:triangular_ineq_chapter_6}
\end{equation}

By combining inequalities \eqref{ch6:eq:hession_bound_chapter_6}, \eqref{ch6:eq:norm_ineq_chapter_6}, and \eqref{ch6:eq:triangular_ineq_chapter_6}, we obtain
\begin{align*}
\norm{\nabla^2 f_i(x)}
& \leq \norm{\nabla^2 \tilde f_i(x)} + 2\lambda \\
& \leq \norm{a_i} \norm{\nabla \tilde f_i(x)} + 2\lambda \\ 
& \leq 2\lambda + \lambda \sqrt{d}\|a_i\| + \norm{a_i} \norm{\nabla f_i(x)}. 
\end{align*}
{\color{darkred}Consequently,
\(\norm{\nabla^2 f_i(x)}\leq L_0+L_1\norm{\nabla f_i(x)}\)
holds simultaneously for every client when
\[
L_0=2\lambda+\lambda\sqrt d\,\max_{j\in[n]}\norm{a_j},
\qquad
L_1=\max_{j\in[n]}\norm{a_j}.
\]
For twice differentiable functions, this Hessian bound implies
Assumption~\ref{ch6:assum:LzeroLoneSmooth_chapter_6}; see
\textcite[Theorem 1]{chen2023generalized}.}

	\newpage
  \endgroup

\chapter{Supplementary Material and Proofs for Chapter 7}
  \begingroup
  \thesischaptersetup{7}%
  \stepEightAAppendixColors{7}%
  
\section{Key Inequalities}

In this section, we introduce basic lemmas for facilitating our convergence analysis.

The first lemma, similar to \textcite[Lemma 8]{nguyen2018sgd}, establishes an explicit expression for the iterates $\{e_k\}$ defined by the specific recursion. 

\begin{lemma}\label{ch7:lemma:summing_the_recursion_chapter_7}
Let the sequence $\{e_k\}$ be governed by 
$e_{k+1} = (1-\eta_k) e_k + (1-\eta_k) A_{k+1} + \eta_k B_{k+1}$ for arbitrary vectors $A_k,B_k$. Then,
\begin{align*}
e_{k+1} &= \prod_{\tau=0}^k (1-\eta_\tau) {e}_0 + \sum_{j=0}^k \left(\prod_{\tau=j+1}^k (1-\eta_\tau)\right) (1-\eta_j) {A}_{j+1}\\
&\quad + \sum_{j=0}^k \left(\prod_{\tau=j+1}^k (1-\eta_\tau)\right) \eta_j B_{j+1}.
\end{align*}
\end{lemma}

The next lemma provides the descent inequality for the normalized gradient descent iteration $x^{k+1}=x^k-\gamma_k g^k/\norm{g^k}$.
\begin{lemma}[Descent Lemma]
    \label{ch7:lem:descent_lemma_chapter_7}
    Let~Assumption~\ref{ch7:assum:L_fcn_chapter_7} hold. Then for  the iterates $\{x^k\}$ generated by the following gradient update
    \begin{equation}
    \label{ch7:eq:grad_update_chapter_7}
        x^{k+1} = x^k - \gamma_k \frac{g^k}{\|g^k\|}
    \end{equation}
    with $g^k = (1/n)\sum_{i=1}^n g_i^k$ and $v^k = (1/n)\sum_{i=1}^n v_i^k$ 
    satisfy 
    \begin{equation}
    \label{ch7:eq:descent_lemma_chapter_7}
        \Delta_{k+1} + \gamma_k \|\nabla f(x^k)\| \le \Delta_k  + 2 \gamma_k \norm{v^k- \nabla f(x^k)} +\frac{ 2\gamma_k}{n}\sum_{i=1}^n\norm{g_i^k-v_i^k} + \frac{\gamma_k^2L}{2} ,
    \end{equation}
    where $\Delta_k \eqdef f(x^k) - f^{\inf}$ for any $k \geq 0$. 
\end{lemma}

\begin{proof}
Applying $L$-smoothness of $f(x)$ (Assumption~\ref{ch7:assum:L_fcn_chapter_7}) and~\eqref{ch7:eq:grad_update_chapter_7}, we have 
\begin{align*}
f(x^{k+1}) & \le f(x^k) + \langle \nabla f(x^k), x^{k+1}-x^k \rangle + \frac{L}{2} \|x^{k+1}-x^k\|^2  \\
& = f(x^k) - \gamma_k \langle \nabla f(x^k), g^k/\|g^k\| \rangle + \frac{L}{2} \gamma_k^2  \\
& = f(x^k) - \gamma_k \|g^k\| - \gamma_k \langle \nabla f(x^k) - g^k, g^k/\|g^k\| \rangle + \frac{L}{2} \gamma_k^2  \\
& \le f(x^k) - \gamma_k \|g^k\| + \gamma_k \norm{ \nabla f(x^k) - g^k} + \frac{L}{2} \gamma_k^2 \\
& \le f(x^k) - \gamma_k \|\nabla f(x^k)\| + 2\gamma_k \|\nabla f(x^k) - g^k\| + \frac{L}{2} \gamma_k^2,
\end{align*}
where in the second inequality we used the Cauchy--Schwarz inequality, and in the third the triangle inequality. Finally, by Triangle inequality, we have 
\begin{align*}
    \norm{\nabla f(x^k) - g^k} 
    & \leq  \norm{v^k- \nabla f(x^k)} + \norm{g^k-v^k} \\
    & \leq\norm{v^k- \nabla f(x^k)} + \frac{1}{n}\sum_{i=1}^n\norm{g_i^k-v_i^k}, 
\end{align*}
 by denoting $\Delta_{k+1} := f(x^{k+1}) - f^{\inf}$ for $k \in \{0, \dots, K-1\}$, we complete the proof. 
\end{proof}

Next, the following lemma bounds $(1/n)\sum_{i=1}^n\norm{g_i^k-v_i^k}$
from the iteration $g_i^{k+1}=g_i^k  + \thesischsevencC_i^{k+1}(v_i^{k+1}-g_i^k)$. 
\begin{lemma}\label{ch7:lemma:boundg_and_v_chapter_7}
Let $g_i^{k+1}=g_i^k + \thesischsevencC_i^{k+1}(v_i^{k+1}-g_i^k)$ for $i=1,2,\ldots,n$, where $\thesischsevencC_i^k(\cdot)$ satisfies~Assumption~\ref{ch7:assum:contract_comp_chapter_7}. Then,  
\begin{equation*}
   \Exp{\thesischsevencV_{k+1}}  \leq \sqrt{1-\alpha}    \Exp{\thesischsevencV_{k}} + \frac{\sqrt{1-\alpha}}{n} \sum_{i=1}^n \Exp{\norm{v_i^{k+1}-v_i^k}},
\end{equation*}
where $\thesischsevencV_k = (1/n)\sum_{i=1}^n \norm{g_i^k-v_i^k}$.
\end{lemma}
\begin{proof}
Define $\thesischsevencV_k \eqdef (1/n)\sum_{i=1}^n \norm{g_i^k-v_i^k}$. Then, 
\begin{align*}
    \thesischsevencV_{k+1} 
    & = \frac{1}{n}\sum_{i=1}^n \norm{ v_i^{k+1} - g_i^k -  \mathcal{C}_i^{k+1}(v_i^{k+1} - g_i^k)},
\end{align*}

where we used the definition of $g_i^{k+1}$. Next, by taking the conditional expectation with fixed $\thesischsevencF_{k+1} =\{v_i^{k+1},x^{k+1},g_i^k\}$, and by using Jensen's inequality, 
\begin{align*}
    \thesischsevenExpCond{\thesischsevencV_{k+1}}{\thesischsevencF_{k+1}} 
      & \leq \frac{1}{n}\sum_{i=1}^n\sqrt{ \thesischsevenExpCond{\sqnorm{ v_i^{k+1} - g_i^k -  \mathcal{C}_i^{k+1}(v_i^{k+1} - g_i^k)}}{\thesischsevencF_{k+1}}} \\
      & \leq \frac{\sqrt{1-\alpha}}{n} \sum_{i=1}^n \sqrt{\thesischsevenExpCond{\sqnorm{ v_i^{k+1} - g_i^k}}{\thesischsevencF_{k+1}}} \\
      & =  \frac{\sqrt{1-\alpha}}{n} \sum_{i=1}^n \norm{v_i^{k+1}-g_i^k} \\
      & \leq   \frac{\sqrt{1-\alpha}}{n} \sum_{i=1}^n \norm{v_i^{k+1}-v_i^k} + \frac{\sqrt{1-\alpha}}{n} \sum_{i=1}^n \norm{v_i^k-g_i^k},
\end{align*}
where in the second inequality we used Assumption~\ref{ch7:assum:contract_comp_chapter_7}, and in the third the triangle inequality. Therefore, 
\begin{align*}
   \Exp{\thesischsevencV_{k+1}}  & = \Exp{\thesischsevenExpCond{\thesischsevencV_{k+1}}{\thesischsevencF_{k+1}}} \\
   & \leq   \frac{\sqrt{1-\alpha}}{n} \sum_{i=1}^n \Exp{\norm{v_i^{k+1}-v_i^k}} + \frac{\sqrt{1-\alpha}}{n} \sum_{i=1}^n \Exp{\norm{v_i^k-g_i^k}}. 
\end{align*}

Finally, from the definition of $\thesischsevencV_k$, we obtain the final result. 
\end{proof}

Finally, the following lemma presents an explicit upper bound on the summation involving the tuning parameters $\gamma_k,\eta_k$.
\begin{lemma}[Lemma 15 of \textcite{fatkhullin2023stochastic}]\label{ch7:lemma:sum_prod_chapter_7}
   Let $q\in [0,1), p \geq 0, \gamma_0 >0$, and let $\eta_k = (2/(k+2))^q$ and $\gamma_k = (2/(k+2))^p$ for every integer $k$. Then, for any integer $k$ and $K \geq 1$, 
   \begin{equation*}
       \sum_{k=0}^{K-1} \gamma_k \prod_{\tau=k+1}^{K-1} (1-\eta_\tau) \leq C(p,q) \frac{\gamma_K}{\eta_K},
   \end{equation*}
   where
   \begin{align*}
   C(p,q) &= 2^{p-q} (1-q)^{-1} t_0(p,q) \exp\bigl( 2^q (1-q) t_0^{1-q}\bigr)  + 2^{2p+1-q} (1-q)^{-2},\\
   t_0(p,q) &= \max\left\{ \left(  \frac{p}{(1-q) 2^q}\right)^{1/(1-q)} , \frac{2(p-q)}{(1-q)^2}  \right\}^{1/(1-q)}.
   \end{align*}
\end{lemma}

\newpage 
\section{\normEFname{SGDM}}\label{ch7:sec:EF21_SGDM_chapter_7}

In \normEFname{SGDM}, we update the iterates $\{x^k\}$ according to: 
\begin{equation*}
    x^{k+1} = x^k - \gamma_k \frac{g^k}{\norm{g^k}},
\end{equation*}
where $g^k = (1/n)\sum_{i=1}^n g_i^k$ and $v^k = (1/n)\sum_{i=1}^n v_i^k$ with $g_i^k,v_i^k$ being governed by 
\begin{align*}
    g_i^{k+1} & = g_i^k + \mathcal{C}_i^{k+1}(v_i^{k+1} - g_i^k), \quad \text{and} \\
    v_i^{k+1} & = (1-\eta_k)v_i^k + \eta_k \nabla f_i \left( x^{k+1}  ;\xi_i^{k+1} \right).
\end{align*}

\subsection{Convergence Proof}

We prove the result in the following steps. 

\subsubsection{Deriving the descent inequality} 

From Lemma~\ref{ch7:lem:descent_lemma_chapter_7}, we obtain \eqref{ch7:eq:descent_lemma_chapter_7}:
\begin{equation*}
    \Delta_{k+1} + \gamma_k \|\nabla f(x^k)\| \le \Delta_k  + 2 \gamma_k \norm{v^k- \nabla f(x^k)} +\frac{ 2\gamma_k}{n}\sum_{i=1}^n\norm{g_i^k-v_i^k} + \frac{\gamma_k^2L}{2}.
\end{equation*}

\subsubsection{Error bound I}

Define $\thesischsevencV_k \eqdef (1/n)\sum_{i=1}^n \norm{g_i^k-v_i^k}$. From~Lemma~\ref{ch7:lemma:boundg_and_v_chapter_7}, 
\begin{equation*}
   \Exp{\thesischsevencV_{k+1}}  \leq \sqrt{1-\alpha}    \Exp{\thesischsevencV_{k}} + \frac{\sqrt{1-\alpha}}{n} \sum_{i=1}^n \Exp{\norm{v_i^{k+1}-v_i^k}}.
\end{equation*}

To complete the upper-bound for $ \Exp{\thesischsevencV_{k+1}}$, we bound $\Exp{\norm{v_i^{k+1}-v_i^k}}$ as follows:
\begin{align*}
    \Exp{\norm{v_i^{k+1}-v_i^k}}
    & \leq \eta_k \Exp{\norm{ \nabla f_i(x^{k+1};\xi_i^{k+1}) - v_i^k }} \\
    & \leq \eta_k \Exp{\norm{v_i^k-\nabla f_i(x^k)}} + \eta_k \Exp{\norm{\nabla f_i(x^k) - \nabla f_i(x^{k+1})}} \\
    &\quad + \eta_k \Exp{\norm{\nabla f_i(x^{k+1};\xi_i^{k+1}) - \nabla f_i(x^{k+1})}}.
\end{align*}

Next, by the $L_i$-smoothness of $f_i(\cdot)$, by Jensen's inequality, and by the fact that $\Exp{\sqnorm{\nabla f_i(x^{k+1};\xi_i^{k+1}) - \nabla f_i(x^{k+1})}} \leq \sigma_g^2$, 
\begin{equation*}
    \Exp{\norm{v_i^{k+1}-v_i^k}}
   \leq \eta_k \Exp{\norm{v_i^k-\nabla f_i(x^k)}} + \eta_k L_i \gamma_k + \eta_k \sigma_g.
\end{equation*}
Therefore, 
\begin{align*}
    \Exp{\thesischsevencV_{k+1}}  &\leq \sqrt{1-\alpha} \Exp{\thesischsevencV_k} + \sqrt{1-\alpha} \eta_k \Exp{\thesischsevencU_k}\\
    &\quad + \sqrt{1-\alpha} \eta_k\gamma_k \bar L + \sqrt{1-\alpha} \eta_k \sigma_g,
\end{align*}
where $\bar L = (1/n)\sum_{i=1}^n L_i$.

\subsubsection{Error bound II}

We consider $\norm{v_i^{k+1} - \nabla f_i(x^{k+1})}$. From the definition of $v_i^{k+1}$, 
\begin{align*}
    v_i^{k+1} - \nabla f_i(x^{k+1})
    & = (1-\eta_k)(v_i^k - \nabla f_i(x^k)) + (1-\eta_k)(\nabla f_i(x^{k}) - \nabla f_i(x^{k+1})) \\
    &\quad + \eta_k (\nabla f_i(x^{k+1};\xi_i^{k+1}) - \nabla f_i(x^{k+1})).
\end{align*}

Next, define $\thesischsevencU_k = (1/n)\sum_{i=1}^n \norm{v_i^k - \nabla f_i(x^k)}$. Then, by the triangle inequality, 
\begin{align*}
    \Exp{\thesischsevencU_{k+1}} 
    & \leq \frac{1}{n}\sum_{i=1}^n \Exp{ \norm{v_i^{k+1} - \nabla f_i(x^{k+1})} } \\
    & \leq (1-\eta_k) \Exp{\thesischsevencU_k} + (1-\eta_k) \frac{1}{n}\sum_{i=1}^n \Exp{\norm{ \nabla f_i(x^{k}) - \nabla f_i(x^{k+1} }}  \\
    &\quad + \eta_k \frac{1}{n}\sum_{i=1}^n \Exp{\norm{ \nabla f_i(x^{k+1};\xi_i^{k+1}) - \nabla f_i(x^{k+1}) }}.
\end{align*}

Next, by the $L_i$-smoothness of $f_i(\cdot)$, by Jensen's inequality, and by the fact that $\Exp{\sqnorm{\nabla f_i(x^{k+1};\xi_i^{k+1}) - \nabla f_i(x^{k+1})}} \leq \sigma_g^2$,
\begin{align*}
    \Exp{\thesischsevencU_{k+1}} 
    & \leq (1-\eta_k) \Exp{\thesischsevencU_k} + (1-\eta_k) \frac{1}{n}\sum_{i=1}^n  L_i \Exp{\norm{ x^{k} - x^{k+1} }} + \eta_k \sigma_g \\
    & \leq (1-\eta_k) \Exp{\thesischsevencU_k} + (1-\eta_k) \gamma_k \bar L + \eta_k \sigma_g,
\end{align*}
where in the last inequality we used the update rule for $x^{k+1}$.

\subsubsection{Error Bound III} 
From the recursion of $v_i - \nabla f_i(x^k)$, and by the fact that $v^k = (1/n)\sum_{i=1}^n v_i^k$ and $\nabla f(x) = (1/n)\sum_{i=1}^n \nabla f_i(x)$,
\begin{align*}
    v^{k+1} - \nabla f(x^{k+1})
    & = (1-\eta_k)(v^k - \nabla f(x^k)) + (1-\eta_k)(\nabla f(x^{k}) - \nabla f(x^{k+1})) \\
    &\quad + \eta_k \frac{1}{n}\sum_{i=1}^n [\nabla f_i(x^{k+1};\xi_i^{k+1}) - \nabla f_i(x^{k+1})].
\end{align*}

Next, from~Lemma~\ref{ch7:lemma:summing_the_recursion_chapter_7}, 
\begin{multline*}
 e_{k+1} = \prod_{\tau=0}^k (1-\eta_\tau) {e}_0 + \sum_{j=0}^k \left(\prod_{\tau=j+1}^k (1-\eta_\tau)\right) (1-\eta_j) {A}_{j+1} \\
 + \sum_{j=0}^k \left(\prod_{\tau=j+1}^k (1-\eta_\tau)\right) \eta_j B_{j+1},
\end{multline*}
where
\begin{align*}
    e_k &= v^k - \nabla f(x^k), \qquad A_{k+1} = \nabla f(x^k)-\nabla f(x^{k+1}), \quad\text{and}\\
    B_{k+1} &= \frac{1}{n}\sum_{i=1}^n [\nabla f_i(x^{k+1};\xi_i^{k+1}) - \nabla f_i(x^{k+1})].
\end{align*}

Therefore, from the definition of the Euclidean norm and by the triangle inequality, 
\begin{align*}
    \Exp{\norm{e_{k+1}}} 
    & \leq \prod_{\tau=0}^k (1-\eta_\tau) \Exp{\norm{{e}_0}} + \sum_{j=0}^k \left(\prod_{\tau=j+1}^k (1-\eta_\tau)\right) (1-\eta_j) \Exp{\norm{{A}_{j+1}}} \\
    &\quad + \Exp{\norm{\sum_{j=0}^k \left(\prod_{\tau=j+1}^k (1-\eta_\tau)\right) \eta_j B_{j+1}}}.
\end{align*}

Next, by the $L_i$-smoothness of $f_i(\cdot)$ and by the definition of $x^{k+1}$, we can show that $\norm{A_{j+1}} \leq \bar L \gamma_j$. Therefore, 
\begin{align*}
    \Exp{\norm{e_{k+1}}} 
    & \leq \prod_{\tau=0}^k (1-\eta_\tau) \Exp{\norm{{e}_0}} + \bar L \cdot \sum_{j=0}^k \left(\prod_{\tau=j+1}^k (1-\eta_\tau)\right) (1-\eta_j) \gamma_j  \\
    &\quad + \Exp{\norm{\sum_{j=0}^k \left(\prod_{\tau=j+1}^k (1-\eta_\tau)\right) \eta_j B_{j+1}}}.
\end{align*}

Next, since 
\begin{align*}
    \MoveEqLeft \Exp{\norm{\sum_{j=0}^k \left(\prod_{\tau=j+1}^k (1-\eta_\tau)\right) \eta_j B_{j+1}}} \\
    & \leq \sqrt{  \Exp{\sqnorm{\sum_{j=0}^k \left(\prod_{\tau=j+1}^k (1-\eta_\tau)\right) \eta_j B_{j+1}}} } 
    \leq \frac{\sigma_g}{\sqrt{n}} \sqrt{ \sum_{j=0}^k \left(\prod_{\tau=j+1}^k (1-\eta_\tau)^2\right) \eta_j^2 },
\end{align*}
where the first inequality follows from Jensen's inequality and the second from Assumption~\ref{ch7:assum:stoc_g_h_chapter_7},
we obtain
\begin{align*}
    \Exp{\norm{e_{k+1}}} 
    & \leq \prod_{\tau=0}^k (1-\eta_\tau) \Exp{\norm{{e}_0}} + \bar L \cdot \sum_{j=0}^k \left(\prod_{\tau=j+1}^k (1-\eta_\tau)\right) (1-\eta_j) \gamma_j  \\
    &\quad + \frac{\sigma_g}{\sqrt{n}} \sqrt{ \sum_{j=0}^k \left(\prod_{\tau=j+1}^k (1-\eta_\tau)^2\right) \eta_j^2 }\\
    & \leq \prod_{\tau=0}^k (1-\eta_\tau) \Exp{\norm{{e}_0}} + \bar L \cdot \sum_{j=0}^k \left(\prod_{\tau=j+1}^k (1-\eta_\tau)\right)  \gamma_j  \\
    &\quad + \frac{\sigma_g}{\sqrt{n}} \sqrt{ \sum_{j=0}^k \left(\prod_{\tau=j+1}^k (1-\eta_\tau)\right) \eta_j^2 },
\end{align*}
where in the last inequality we used $1-\eta_j \leq 1$.

If $\eta_k = (2/(k+2))^{1/2}$ and $\gamma_k = \gamma_0 (2/(k+2))^{3/4}$, then we can prove  that $\eta_0=1$, and that
\begin{equation*}
    \Exp{\norm{e_{k+1}}} 
    \leq  \bar L \cdot \sum_{j=0}^k \left(\prod_{\tau=j+1}^k (1-\eta_\tau)\right)  \gamma_j  + \frac{\sigma_g}{\sqrt{n}} \sqrt{ \sum_{j=0}^k \left(\prod_{\tau=j+1}^k (1-\eta_\tau)\right) \eta_j^2 }.
\end{equation*}
In conclusion, from~Lemma~\ref{ch7:lemma:sum_prod_chapter_7}, 
\begin{equation*}
    \Exp{\norm{v^{k+1}  - \nabla f(x^{k+1})}} 
    \leq  \bar L C(3/4, 1/2) \frac{\gamma_{k+1}}{\eta_{k+1}}   + \frac{\sigma_g}{\sqrt{n}} \sqrt{ C(1, 1/2) \eta_{k+1} }.
\end{equation*}

\subsubsection{Bounding Lyapunov function }

Define Lyapunov function as  $V_k = \Delta_{k} + C_{1,k} \thesischsevencV_k + C_{2,k} \thesischsevencU_k$ with
\begin{equation*}
    C_{1,k} = \frac{2\gamma_k}{1-\sqrt{1-\alpha}} \quad\text{and}\quad C_{2,k} = \frac{2\gamma_k \sqrt{1-\alpha}}{1-\sqrt{1-\alpha}}.
\end{equation*}
Then, 
\begin{align*}
    \Exp{V_{k+1}} 
    & = \Exp{\Delta_{k+1} + C_{1,k+1} \thesischsevencV_{k+1} + C_{2,k+1} \thesischsevencU_{k+1}} \\
    & \leq \Exp{\Delta_k} - \gamma_k \Exp{\norm{\nabla f(x^k)}} + 2 \gamma_k \Exp{\norm{v^k- \nabla f(x^k)}} + 2\gamma_k \Exp{\thesischsevencV_k} \\
    &\quad + \frac{\gamma_k^2L}{2} + C_{1,k+1} \Exp{\thesischsevencV_{k+1}} + C_{2,k+1} \Exp{\thesischsevencU_{k+1}},
\end{align*}
where in the last inequality we used \eqref{ch7:eq:descent_lemma_chapter_7}.
Next, by the upper-bounds for $\Exp{\thesischsevencV_{k+1}}$ and  $\Exp{\thesischsevencU_{k+1}}$, 
\begin{align*}
    \MoveEqLeft \Exp{ V_{k+1} }\\
   & \leq \Exp{\Delta_k} - \gamma_k \Exp{\norm{\nabla f(x^k)}} + 2 \gamma_k \Exp{\norm{v^k- \nabla f(x^k)}} + 2\gamma_k \Exp{\thesischsevencV_k} + \frac{\gamma_k^2L}{2} \\
    &\quad + C_{1,k+1} \bigl(\sqrt{1-\alpha} \Exp{\thesischsevencV_k} + \sqrt{1-\alpha} \eta_k \Exp{\thesischsevencU_k}\\
    &\qquad\qquad + \sqrt{1-\alpha} \eta_k\gamma_k \bar L + \sqrt{1-\alpha} \eta_k \sigma_g\bigr) \\
    &\quad + C_{2,k+1} \bigl((1-\eta_k) \Exp{\thesischsevencU_k} + (1-\eta_k) \gamma_k \bar L + \eta_k \sigma_g\bigr).
\end{align*}

Next,  since $\gamma_{k+1} \leq \gamma_k$, we can prove that $C_{1,k+1} \leq C_{1,k}$ and that $C_{2,k+1} \leq C_{2,k}$, and also that 
\begin{align*}
    2\gamma_k + C_{1,k+1}\sqrt{1-\alpha}  & \leq 2\gamma_k + C_{1,k}\sqrt{1-\alpha} = C_{1,k}, \\
    C_{1,k+1} \sqrt{1-\alpha}\eta_k + C_{2,k+1} (1-\eta_k) & \leq C_{1,k} \sqrt{1-\alpha}\eta_k + C_{2,k} (1-\eta_k)= C_{2,k}.
\end{align*}
Therefore, 
\begin{align*}
   \Exp{ V_{k+1} } 
   & \leq \Exp{V_k} - \gamma_k \Exp{\norm{\nabla f(x^k)}} + 2 \gamma_k \Exp{\norm{v^k- \nabla f(x^k)}}  + \frac{\gamma_k^2L}{2} \\
    &\quad + C_{1,k} (\sqrt{1-\alpha} \eta_k\gamma_k \bar L + \sqrt{1-\alpha} \eta_k \sigma_g) \\
    &\quad + C_{2,k} ( (1-\eta_k) \gamma_k \bar L + \eta_k \sigma_g) \\
     & \leq \Exp{V_k} - \gamma_k \Exp{\norm{\nabla f(x^k)}} + 2 \gamma_k \Exp{\norm{v^k- \nabla f(x^k)}}   \\
    &\quad + \gamma_k^2 \left( \frac{2\sqrt{1-\alpha}}{1-\sqrt{1-\alpha}} \bar L + \frac{L}{2} \right) + \eta_k\gamma_k \frac{4\sqrt{1-\alpha}}{1-\sqrt{1-\alpha}} \sigma_g.
\end{align*}

Next, from the upper-bound of $\Exp{\norm{v^k- \nabla f(x^k)}}$, 
\begin{align*}
    \MoveEqLeft \Exp{ V_{k+1} } \\
    & \leq \Exp{V_k} - \gamma_k \Exp{\norm{\nabla f(x^k)}} + 2 \gamma_k \left( \bar L C(3/4, 1/2) \frac{\gamma_{k}}{\eta_{k}}   + \frac{\sigma_g}{\sqrt{n}} \sqrt{ C(1, 1/2) \eta_{k} } \right)  \\
    &\quad + \gamma_k^2 \left( \frac{2\sqrt{1-\alpha}}{1-\sqrt{1-\alpha}} \bar L + \frac{L}{2} \right) + \eta_k\gamma_k \frac{4\sqrt{1-\alpha}}{1-\sqrt{1-\alpha}} \sigma_g \\
    & = \Exp{V_k} - \gamma_k \Exp{\norm{\nabla f(x^k)}} + \frac{\gamma_k^2}{\eta_k} \cdot 2 \bar L C_1 + \gamma_k \sqrt{\eta_k} \cdot \frac{2 \sigma_g C_2}{\sqrt{n}}  \\
    &\quad + \gamma_k^2 \left( \frac{2\sqrt{1-\alpha}}{1-\sqrt{1-\alpha}} \bar L + \frac{L}{2} \right) +  \eta_k\gamma_k \cdot  \frac{4\sqrt{1-\alpha}}{1-\sqrt{1-\alpha}}  \sigma_g,
\end{align*}
where  we denoted $C_1 = C(3/4, 1/2)$, $C_2 = \sqrt{C(1, 1/2)}$.

\subsubsection{Deriving the convergence rate}
Let $\Gamma_K \eqdef \sum_{k=0}^{K-1} \gamma_k$. By re-arranging the terms and by the telescopic series, 
\begin{align*}
    \MoveEqLeft \frac{1}{\Gamma_K} \sum_{k=0}^{K-1} \gamma_k\Exp{\norm{\nabla f(x^k)}}\\
    & \leq \frac{1}{\Gamma_K} \biggl[ \Exp{V_0} -\Exp{V_K} + 2 \bar L C_1 \sum_{k=0}^{K-1} \gamma_k^2 \eta_k^{-1} + \frac{2\sigma_g C_2}{\sqrt{n}} \sum_{k=0}^{K-1} \gamma_k \sqrt{\eta_k} \\
    &\qquad +  \left( \frac{2\sqrt{1-\alpha}}{1-\sqrt{1-\alpha}} \bar L + \frac{L}{2} \right) \sum_{k=0}^{K-1} \gamma_k^2 + \frac{4\sqrt{1-\alpha}}{1-\sqrt{1-\alpha}}  \sigma_g \sum_{k=0}^{K-1} \eta_k\gamma_k \biggr]\\
    & \leq \frac{1}{\Gamma_K} \biggl[ \Exp{V_0} + 2 \bar L C_1 \sum_{k=0}^{K-1} \gamma_k^2\eta_k^{-1} +  \left( \frac{2}{\alpha^2} \bar L + \frac{L}{2} \right) \sum_{k=0}^{K-1} \gamma_k^2 \\
    &\qquad + \frac{2\sigma_g C_2}{\sqrt{n}} \sum_{k=0}^{K-1} \gamma_k \sqrt{\eta_k} +  \frac{4\sigma_g}{\alpha^2} \sum_{k=0}^{K-1} \eta_k\gamma_k \biggr],
\end{align*}
where in the last inequality we used $\sqrt{1-\alpha}/(1-\sqrt{1-\alpha}) \leq 1/\alpha^2$ for any $\alpha \in (0,1]$.
To continue the proof, we need to provide bounds on several sum-type terms. Since we set $\eta_k = (2/(k+2))^{1/2}$ and $\gamma_k = \gamma_0 (2/(k+2))^{3/4}$, we have 
\begin{align*}
    \Gamma_K = \sum_{k=0}^{K-1} \gamma_k &\geq  K \gamma_{K-1} = \gamma_0K  \left(\frac{2}{K+1}\right)^{3/4} \geq \gamma_0 K^{1/4};\\
    \sum_{k=0}^{K-1} \gamma_k^2 \eta_k^{-1}  &= \gamma_0^2\sum_{k=0}^{K-1}  \left( \frac{2}{k+2} \right)^{3/2} \left( \frac{2}{k+2} \right)^{-1/2} = 2\gamma_0^2 \sum_{k=0}^{K-1} \frac{1}{k+2}\\
    &= 2\gamma_0^2 \sum^{K}_{k=1}\frac{1}{1+k} \leq 2\gamma^2_0 \left(\frac{1}{2} + \int^K_{1}\frac{1}{1+s}\,ds\right) \leq 2\gamma^2_0 \log \left(K+1\right); \\
    \sum_{k=0}^{K-1} \gamma_k \sqrt{\eta_k} &= \gamma_0 \sum_{k=0}^{K-1} \left( \frac{2}{k+2} \right)^{3/4} \left( \frac{2}{k+2} \right)^{1/4}\\
    &=   2\gamma_0 \sum_{k=0}^{K-1}\frac{1}{k+2} \leq 2\gamma_0 \log\left(K+1\right);\\
    \sum_{k=0}^{K-1} \gamma_k^2 &= \gamma_0^2 \sum_{k=0}^{K-1} \left( \frac{2}{k+2} \right)^{3/2} = 2^{3/2}\gamma_0^2 \sum_{k=1}^{K} \frac{1}{(1+k)^{3/2}}\\
    &\leq 2^{3/2}\gamma_0^2 \left(\frac{1}{2^{3/2}} + \int^{\infty}_{1}\frac{1}{(1+s)^{3/2}}\,ds\right) = 5\gamma^2_0;\\
    \sum_{k=0}^{K-1} \gamma_k \eta_k &= \gamma_0 \sum_{k=0}^{K-1}  \left( \frac{2}{k+2} \right)^{3/4} \left( \frac{2}{k+2} \right)^{1/2} = 2^{5/4}\gamma_0 \sum^{K}_{k=1} \frac{1}{(1+k)^{5/4}} \\
    &\leq 2^{5/4}\gamma_0 \left(\frac{1}{2^{5/4}} + \int^{\infty}_{1}\frac{1}{(1+s)^{5/4}}\,ds\right) = 9\gamma_0.
\end{align*}
Therefore, denoting $\tilde{x}^K$ as a point randomly chosen from $\{x^0,x^1,\ldots,x^{K-1}\}$ with probability $\gamma_k/\Gamma_K$ for $k = 0,1, \ldots, K-1$, we obtain 
\begin{align*}
     \Exp{\norm{\nabla f(\tilde{x}^K)}} 
    & \leq \frac{\Exp{V_0}}{\gamma_0 K^{1/4}} + 2 \bar L C_1 \frac{2\gamma^2_0 \log \left(K+1\right)}{\gamma_0 K^{1/4}} +  \left( \frac{2}{\alpha^2} \bar L + \frac{L}{2} \right) \frac{5\gamma^2_0}{ \gamma_0 K^{1/4}}\\
    &\quad + \frac{2\sigma_g C_2}{\sqrt{n}} \frac{2\gamma_0 \log \left(K+1\right)}{ \gamma_0 K^{1/4}} +  \frac{4\sigma_g}{\alpha^2}  \frac{9\gamma_0}{ \gamma_0 K^{1/4}}\\
    &= \frac{\Exp{V_0}}{\gamma_0 K^{1/4}} + \frac{\gamma_0}{K^{1/4}}\left(\frac{5L}{2} + 4C_1{\bar L}\log\left(K+1\right) + \frac{10\bar L}{\alpha^2}\right)\\
    &\quad +\frac{4}{K^{1/4}}\left(C_2 \frac{\sigma_g}{\sqrt{n}}\log(K+1) +\frac{9\sigma_g}{\alpha^2} \right)\\
    &= \widetilde{\cO}\left(\frac{1}{K^{1/4}} \left( \frac{\Exp{V_0}}{\gamma_0} + \gamma_0\left(L + \frac{\bar L}{\alpha^2}\right) + \sigma_g\left(\frac{1}{\sqrt{n}} + \frac{1}{\alpha^2}\right) \right) \right).
\end{align*}

\newpage 
\section{\normEFname{IGT}}

In \normEFname{IGT}, we update the iterates $\{x^k\}$ according to: 
\begin{equation*}
    x^{k+1} = x^k - \gamma_k \frac{g^k}{\norm{g^k}},
\end{equation*}
where $g^k = (1/n)\sum_{i=1}^n g_i^k$ and $v^k = (1/n)\sum_{i=1}^n v_i^k$ with $g_i^k,v_i^k$ being governed by 
\begin{align*}
    g_i^{k+1} & = g_i^k + \mathcal{C}_i^{k+1}(v_i^{k+1} - g_i^k), \quad \text{and} \\
    v_i^{k+1} & = (1-\eta_k)v_i^k + \eta_k \nabla f_i \left( x^{k+1} + \frac{1-\eta_k}{\eta_k} (x^{k+1}-x^k) ;\xi_i^{k+1} \right).
\end{align*}
For simplicity, we denote $y^{k+1} = x^{k+1} +\theta_k(x^{k+1} - x^k)$, where $\theta_k = (1-\eta_k)/\eta_k$.

\subsection{Convergence Proof}

We prove the result in the following steps. 

\subsubsection{Deriving the descent inequality} From Lemma~\ref{ch7:lem:descent_lemma_chapter_7}, we obtain \eqref{ch7:eq:descent_lemma_chapter_7}:
\begin{equation*}
    \Delta_{k+1} + \gamma_k \|\nabla f(x^k)\| \le \Delta_k  + 2 \gamma_k \norm{v^k- \nabla f(x^k)} +\frac{ 2\gamma_k}{n}\sum_{i=1}^n\norm{g_i^k-v_i^k} + \frac{\gamma_k^2L}{2}.
\end{equation*}

\subsubsection{Error bound I} 

Define $\thesischsevencV_k \eqdef (1/n)\sum_{i=1}^n \norm{g_i^k-v_i^k}$. From~Lemma~\ref{ch7:lemma:boundg_and_v_chapter_7}, 
\begin{equation*}
   \Exp{\thesischsevencV_{k+1}}  \leq \sqrt{1-\alpha}    \Exp{\thesischsevencV_{k}} + \frac{\sqrt{1-\alpha}}{n} \sum_{i=1}^n \Exp{\norm{v_i^{k+1}-v_i^k}}.
\end{equation*}

To finalize the upper-bound for $\Exp{\thesischsevencV_{k+1}}$, we bound $\Exp{\norm{v_i^{k+1}-v_i^k}}$.
\begin{align*}
    \Exp{\norm{v_i^{k+1}-v_i^k}}
    & = \eta_k \Exp{\norm{ \nabla f_i(y^{k+1};\xi_i^{k+1}) - v_i^k }} \\
    & \leq \eta_k \Exp{\norm{v_i^k-\nabla f_i(x^k)}} + \eta_k \Exp{\norm{\nabla f_i(x^k) - \nabla f_i(y^{k+1})}} \\
    &\quad + \eta_k \Exp{\norm{\nabla f_i(y^{k+1};\xi_i^{k+1}) - \nabla f_i(y^{k+1})}}. 
\end{align*}

Next, by the $L_i$-smoothness of $f_i$ and from the definition of $y^{k+1}$
\begin{align*}
  \Exp{\norm{v_i^{k+1}-v_i^k}}  
    & \leq \eta_k \Exp{\norm{v_i^k-\nabla f_i(x^k)}} + \eta_k L_i \Exp{\norm{ x^k - y^{k+1}}} \\
    &\quad + \eta_k \Exp{\norm{\nabla f_i(y^{k+1};\xi_i^{k+1}) - \nabla f_i(y^{k+1})}} \\
    & \leq \eta_k \Exp{\norm{v_i^k-\nabla f_i(x^k)}} + \eta_k\left(1+ \theta_k \right) L_i \Exp{\norm{ x^k -  x^{k+1}}} \\
    &\quad + \eta_k \Exp{\norm{\nabla f_i(y^{k+1};\xi_i^{k+1}) - \nabla f_i(y^{k+1})}} \\
    & \leq \eta_k \Exp{\norm{v_i^k-\nabla f_i(x^k)}} + \gamma_k L_i \\
    &\quad + \eta_k \Exp{\norm{\nabla f_i(y^{k+1};\xi_i^{k+1}) - \nabla f_i(y^{k+1})}},
\end{align*}
where in the last inequality we used the update rule for $x^{k+1}$ and $\theta_k = (1-\eta_k)/\eta_k$. 
Next, by~Assumption~\ref{ch7:assum:stoc_g_h_chapter_7}, 
\begin{equation*}
  \Exp{\norm{v_i^{k+1}-v_i^k}}
  \leq \eta_k \Exp{\norm{v_i^k-\nabla f_i(x^k)}} +  L_i \gamma_k + \eta_k \sigma_g. 
\end{equation*}

Therefore,
\begin{equation*}
   \Exp{\thesischsevencV_{k+1}}  \leq \sqrt{1-\alpha}    \Exp{\thesischsevencV_{k}} + \sqrt{1-\alpha} \eta_k \Exp{\thesischsevencU_k} + \sqrt{1-\alpha} \bar L \gamma_k  + \sqrt{1-\alpha} \sigma_g \eta_k.
\end{equation*}

\subsubsection{Error bound II}

We consider $\norm{v_i^{k+1} - \nabla f_i(x^{k+1})}$. From the definition of $v_i^{k+1}$, 
\begin{align*}
    v_i^{k+1} - \nabla f_i(x^{k+1})
    & = (1-\eta_k)(v_i^k - \nabla f_i(x^k)) \\
    &\quad + (1-\eta_k)(\nabla f_i(x^{k}) - \nabla f_i(x^{k+1})) \\
    &\quad + \eta_k(\nabla f_i(y^{k+1}) -\nabla f_i(x^{k+1})) \\
    &\quad + \eta_k (\nabla f_i(y^{k+1};\xi_i^{k+1}) - \nabla f_i(y^{k+1}))\\
    & = (1-\eta_k)(v_i^k - \nabla f_i(x^k)) \\
    &\quad + (1-\eta_k)(\nabla f_i(x^{k}) - \nabla f_i(x^{k+1})) \\
    &\quad + \eta_k(\nabla f_i(y^{k+1}) -\nabla f_i(x^{k+1})) \\
    &\quad + \eta_k (\nabla f_i(y^{k+1};\xi_i^{k+1}) - \nabla f_i(y^{k+1}))\\
    &\quad+ (1-\eta_k) \cdot \nabla^2 f_i(x^{k+1})(x^k-x^{k+1}) \\
    &\quad - (1-\eta_k) \cdot \nabla^2 f_i(x^{k+1})(x^k-x^{k+1}). 
\end{align*}

Next, by the fact that $y^{k+1}-x^{k+1} = \theta_k(x^{k+1}-x^k) = ((1-\eta_k)/\eta_k)(x^{k+1}-x^k)$, we can prove that 
\begin{align*}
    (1-\eta_k) \cdot \nabla^2 f_i(x^{k+1})(x^k-x^{k+1}) 
    &= \eta_k \cdot \frac{1-\eta_k}{\eta_k} \nabla^2 f_i(x^{k+1})(x^k-x^{k+1}) \\
    &= -\eta_k \nabla^2 f_i(x^{k+1})(y^{k+1}-x^{k+1}).
\end{align*}
Therefore, 
\begin{align*}
    v_i^{k+1} - \nabla f_i(x^{k+1})
    & = (1-\eta_k)(v_i^k - \nabla f_i(x^k)) \\
    &\quad + (1-\eta_k) Z_{f_i}(x^k,x^{k+1}) + \eta_k Z_{f_i}(y^{k+1},x^{k+1})\\
    &\quad + \eta_k (\nabla f_i(y^{k+1};\xi_i^{k+1}) - \nabla f_i(y^{k+1})),
\end{align*}
where $Z_{f_i}(x,y) :=\nabla f_i(x) - \nabla f_i(y) -\nabla^2 f_i(y)(x-y)$.

Next, define $\thesischsevencU_k = (1/n)\sum_{i=1}^n \norm{v_i^k - \nabla f_i(x^k)}$. Then, by the triangle inequality, 
\begin{align*}
    \Exp{\thesischsevencU_{k+1}} 
    & \leq \frac{1}{n}\sum_{i=1}^n \Exp{ \norm{v_i^{k+1} - \nabla f_i(x^{k+1})} } \\
    & \leq (1-\eta_k) \Exp{\thesischsevencU_k} + \eta_k \frac{1}{n}\sum_{i=1}^n \Exp{\norm{ \nabla f_i(y^{k+1};\xi_i^{k+1}) - \nabla f_i(y^{k+1}) }} \\
    &\quad + (1-\eta_k) \frac{1}{n}\sum_{i=1}^n \Exp{\norm{ Z_{f_i}(x^k,x^{k+1}) }}  + \eta_k  \frac{1}{n}\sum_{i=1}^n \Exp{\norm{ Z_{f_i}(y^{k+1},x^{k+1}) }}.
\end{align*}

Next, by the fact that $\Exp{\sqnorm{\nabla f_i(x^{k+1};\xi_i^{k+1}) - \nabla f_i(x^{k+1})}} \leq \sigma_g^2$,
\begin{align*}
    \Exp{\thesischsevencU_{k+1}} 
    & \leq (1-\eta_k) \Exp{\thesischsevencU_k} + \eta_k \sigma_g \\
    &\quad + (1-\eta_k) \frac{1}{n}\sum_{i=1}^n \Exp{\norm{ Z_{f_i}(x^k,x^{k+1}) }}  + \eta_k  \frac{1}{n}\sum_{i=1}^n \Exp{\norm{ Z_{f_i}(y^{k+1},x^{k+1}) }}.
\end{align*}

By the $L_{h,i}$-Lipschitz continuity of $\nabla^2 f_i(\cdot)$, i.e. $\norm{Z_{f_i}(x,y)} \leq (L_{h,i}/2) \sqnorm{x-y}$,
\begin{align*}
    \Exp{\thesischsevencU_{k+1}} 
    & \leq (1-\eta_k) \Exp{\thesischsevencU_k} + \eta_k \sigma_g \\
    &\quad + (1-\eta_k) \frac{1}{n}\sum_{i=1}^n \frac{L_{h,i}}{2}\Exp{\sqnorm{x^{k} -x^{k+1}}}\\
    &\quad + \eta_k  \frac{1}{n}\sum_{i=1}^n \frac{L_{h,i}}{2}\Exp{\sqnorm{y^{k+1} -x^{k+1}}} \\
    & \leq (1-\eta_k) \Exp{\thesischsevencU_k} + \eta_k \sigma_g \\
    &\quad + (1-\eta_k) \frac{1}{n}\sum_{i=1}^n \frac{L_{h,i}}{2}\Exp{\sqnorm{x^{k} -x^{k+1}}}\\
    &\quad + \eta_k \theta_k^2 \frac{1}{n}\sum_{i=1}^n \frac{L_{h,i}}{2}\Exp{\sqnorm{x^{k+1} -x^{k}}} \\
    & \leq (1-\eta_k) \Exp{\thesischsevencU_k} + \eta_k \sigma_g + (1-\eta_k) \gamma_k^2 \frac{\bar L_{h}}{2}  + \eta_k \theta_k^2 \gamma_k^2  \frac{\bar L_{h}}{2} \\
    & \leq  (1-\eta_k) \Exp{\thesischsevencU_k} + \eta_k \sigma_g +  \frac{1-\eta_k}{\eta_k}\gamma_k^2  \frac{\bar L_{h}}{2},
\end{align*}
where in the second inequality we used the definition of $y^{k+1}$, in the third we denoted $\bar L_h = (1/n)\sum_{i=1}^n L_{h,i}$ and used the update rule for $x^{k+1}$, and in the last one we used $(1-\eta_k) + \eta_k\theta_k^2 = (1-\eta_k)/\eta_k$.

\subsubsection{Error bound III} 

From the recursion of $v_i - \nabla f_i(x^k)$, and by the fact that $v^k = (1/n)\sum_{i=1}^n v_i^k$ and $\nabla f(x) = (1/n)\sum_{i=1}^n \nabla f_i(x)$,
\begin{align*}
    v^{k+1} - \nabla f(x^{k+1})
    & = (1-\eta_k)(v^k - \nabla f(x^k)) \\
    &\quad + (1-\eta_k)   Z_{f}(x^k,x^{k+1}) + \eta_k Z_{f}(y^{k+1},x^{k+1})\\
    &\quad + \eta_k \frac{1}{n}\sum_{i=1}^n (\nabla f_i(y^{k+1};\xi_i^{k+1}) - \nabla f_i(y^{k+1})). 
\end{align*}

Next, from~Lemma~\ref{ch7:lemma:summing_the_recursion_chapter_7}, 
\begin{align*}
 e_{k+1} 
 & = \prod_{\tau=0}^k (1-\eta_\tau) {e}_0 + \sum_{j=0}^k \left(\prod_{\tau=j+1}^k (1-\eta_\tau)\right) (1-\eta_j) {A}_{j+1} \\
 &\quad + \sum_{j=0}^k \left(\prod_{\tau=j+1}^k (1-\eta_\tau)\right) \eta_j B_{j+1}  + \sum_{j=0}^k \left(\prod_{\tau=j+1}^k (1-\eta_\tau)\right) \eta_j C_{j+1},
\end{align*}
where we introduced the following notation: 
\begin{align*}
     e_k &\eqdef v^k - \nabla f(x^k), \qquad  A_{k+1} \eqdef  Z_{f}(x^k,x^{k+1}),\\
     B_{k+1} &\eqdef  Z_{f}(y^{k+1},x^{k+1}), \quad\text{and}\\
     C_{k+1} &\eqdef \frac{1}{n}\sum_{i=1}^n [\nabla f_i(y^{k+1};\xi_i^{k+1}) - \nabla f_i(y^{k+1})].
\end{align*}

Therefore, 
from the definition of the Euclidean norm and by the triangle inequality, 
\begin{align*}
\Exp{\norm{e_{k+1}}} 
    & \leq \prod_{\tau=0}^k (1-\eta_\tau) \Exp{\norm{{e}_0}} + \sum_{j=0}^k \left(\prod_{\tau=j+1}^k (1-\eta_\tau)\right) (1-\eta_j) \Exp{\norm{{A}_{j+1}}} \\
    &\quad + \sum_{j=0}^k \left(\prod_{\tau=j+1}^k (1-\eta_\tau)\right) \eta_j \Exp{\norm{B_{j+1}}}\\
    &\quad + \Exp{\norm{\sum_{j=0}^k \left(\prod_{\tau=j+1}^k (1-\eta_\tau)\right) \eta_j C_{j+1}}}.
\end{align*}

Next, by $L_{h}$-Lipschitz continuity of $\nabla^2 f(\cdot)$,
\begin{align*}
    \norm{A_{j+1}} & \leq \frac{L_{h}}{2}\sqnorm{x^j-x^{j+1}} \leq \frac{L_{h}}{2} \gamma_j^2, \quad \text{and} \\
    \norm{B_{j+1}} & \leq \frac{L_{h}}{2}\sqnorm{y^{j+1}-x^{j+1}} \leq \frac{L_{h}}{2} \theta_j^2\gamma_j^2 = \frac{L_{h}}{2} \frac{(1-\eta_j)^2} {\eta_j^2} \gamma_j^2. 
\end{align*}
Therefore, 
\begin{align*}
    \Exp{\norm{e_{k+1}}} 
    & \leq \prod_{\tau=0}^k (1-\eta_\tau) \Exp{\norm{{e}_0}} + \frac{L_h}{2}\sum_{j=0}^k \left(\prod_{\tau=j+1}^k (1-\eta_\tau)\right) (1-\eta_j) \gamma_j^2 \\
    &\quad + \frac{L_h}{2}\sum_{j=0}^k \left(\prod_{\tau=j+1}^k (1-\eta_\tau)\right) \frac{(1-\eta_j)^2}{\eta_j} \gamma_j^2 \\
    &\quad + \Exp{\norm{\sum_{j=0}^k \left(\prod_{\tau=j+1}^k (1-\eta_\tau)\right) \eta_j C_{j+1}}}.
\end{align*}

Next, since $1-\eta_j + (1-\eta_j)^2/\eta_j = (1-\eta_j)/\eta_j \leq 1/\eta_j$, we obtain 
\begin{align*}
    \Exp{\norm{e_{k+1}}} 
    & \leq \prod_{\tau=0}^k (1-\eta_\tau) \Exp{\norm{{e}_0}} + \frac{L_h}{2}\sum_{j=0}^k \left(\prod_{\tau=j+1}^k (1-\eta_\tau)\right) \frac{1}{\eta_j} \gamma_j^2 \\
    &\quad + \Exp{\norm{\sum_{j=0}^k \left(\prod_{\tau=j+1}^k (1-\eta_\tau)\right) \eta_j C_{j+1}}}.
\end{align*}

Next, since 
\begin{align*}
    \MoveEqLeft \Exp{\norm{\sum_{j=0}^k \left(\prod_{\tau=j+1}^k (1-\eta_\tau)\right) \eta_j C_{j+1}}} \\
    & \leq \sqrt{  \Exp{\sqnorm{\sum_{j=0}^k \left(\prod_{\tau=j+1}^k (1-\eta_\tau)\right) \eta_j C_{j+1}}} } \\
    & \leq \frac{\sigma_g}{\sqrt{n}} \sqrt{ \sum_{j=0}^k \left(\prod_{\tau=j+1}^k (1-\eta_\tau)^2\right) \eta_j^2 },
\end{align*}
where the first inequality follows from Jensen's inequality and the second from Assumption~\ref{ch7:assum:stoc_g_h_chapter_7}, we get
\begin{align*}
    \Exp{\norm{e_{k+1}}} 
    & \leq \prod_{\tau=0}^k (1-\eta_\tau) \Exp{\norm{{e}_0}} + \frac{L_h}{2}\sum_{j=0}^k \left(\prod_{\tau=j+1}^k (1-\eta_\tau)\right) \frac{1}{\eta_j} \gamma_j^2 \\
    &\quad + \frac{\sigma_g}{\sqrt{n}} \sqrt{ \sum_{j=0}^k \left(\prod_{\tau=j+1}^k (1-\eta_\tau)\right) \eta_j^2 }.
\end{align*}

If $\eta_k = (2/(k+2))^{4/7}$ and $\gamma_k = \gamma_0 (2/(k+2))^{5/7}$, then we can prove  that $\eta_0=1$, and that
\begin{equation*}
    \Exp{\norm{e_{k+1}}} 
    \leq \frac{L_h}{2}\sum_{j=0}^k \left(\prod_{\tau=j+1}^k (1-\eta_\tau)\right) \frac{1}{\eta_j} \gamma_j^2 + \frac{\sigma_g}{\sqrt{n}} \sqrt{ \sum_{j=0}^k \left(\prod_{\tau=j+1}^k (1-\eta_\tau)\right) \eta_j^2 }.
\end{equation*}

In conclusion, from~Lemma~\ref{ch7:lemma:sum_prod_chapter_7}, 
\begin{equation*}
    \Exp{\norm{e_{k+1}}} 
    \leq \frac{L_h}{2}C(6/7, 4/7) \left( \frac{\gamma_{k+1}}{\eta_{k+1}}\right)^2  + \frac{\sigma_g}{\sqrt{n}} \sqrt{C(8/7, 4/7) \eta_{k+1} }.
\end{equation*}

\subsubsection{Bounding Lyapunov function}

We define Lyapunov function as  $V_k = \Delta_{k} + C_{1,k} \thesischsevencV_k + C_{2,k} \thesischsevencU_k$ with
\begin{equation*}
    C_{1,k} = \frac{2\gamma_k}{1-\sqrt{1-\alpha}} \quad\text{and}\quad C_{2,k} = \frac{2\gamma_k \sqrt{1-\alpha}}{1-\sqrt{1-\alpha}}.
\end{equation*}
Then, 
\begin{align*}
    \Exp{V_{k+1}} 
    & = \Exp{\Delta_{k+1} + C_{1,k+1} \thesischsevencV_{k+1} + C_{2,k+1} \thesischsevencU_{k+1}} \\
    & \leq \Exp{\Delta_k} - \gamma_k \Exp{\norm{\nabla f(x^k)}} + 2 \gamma_k \Exp{\norm{v^k- \nabla f(x^k)}} + 2\gamma_k \Exp{\thesischsevencV_k} \\
    &\quad + \frac{\gamma_k^2L}{2} + C_{1,k+1} \Exp{\thesischsevencV_{k+1}} + C_{2,k+1} \Exp{\thesischsevencU_{k+1}},
\end{align*}
where in the last inequality we used \eqref{ch7:eq:descent_lemma_chapter_7}.
Next, by the upper-bounds for $\Exp{\thesischsevencV_{k+1}}$ and  $\Exp{\thesischsevencU_{k+1}}$, 
\begin{align*}
    \MoveEqLeft \Exp{ V_{k+1} }\\
   & \leq \Exp{\Delta_k} - \gamma_k \Exp{\norm{\nabla f(x^k)}} + 2 \gamma_k \Exp{\norm{v^k- \nabla f(x^k)}} + 2\gamma_k \Exp{\thesischsevencV_k} + \frac{\gamma_k^2L}{2} \\
    &\quad + C_{1,k+1} \bigl(\sqrt{1-\alpha}    \Exp{\thesischsevencV_{k}} + \sqrt{1-\alpha} \eta_k \Exp{\thesischsevencU_k}\\
    &\qquad\qquad + \sqrt{1-\alpha} \bar L \gamma_k  + \sqrt{1-\alpha} \sigma_g \eta_k\bigr) \\
    &\quad + C_{2,k+1} \left( (1-\eta_k) \Exp{\thesischsevencU_k} + \eta_k \sigma_g +  \frac{1-\eta_k}{\eta_k}\gamma_k^2  \frac{\bar L_{h}}{2} \right).
\end{align*}

Next,  since $\gamma_{k+1} \leq \gamma_k$, we can prove that $C_{1,k+1} \leq C_{1,k}$ and that $C_{2,k+1} \leq C_{2,k}$, and also that 
\begin{align*}
    2\gamma_k + C_{1,k+1}\sqrt{1-\alpha}  & \leq 2\gamma_k + C_{1,k}\sqrt{1-\alpha} = C_{1,k}, \\
    C_{1,k+1} \sqrt{1-\alpha}\eta_k + C_{2,k+1} (1-\eta_k) & \leq C_{1,k} \sqrt{1-\alpha}\eta_k + C_{2,k} (1-\eta_k)= C_{2,k}.
\end{align*}
Therefore, 
\begin{align*}
   \MoveEqLeft \Exp{ V_{k+1} } \\
   & \leq \Exp{V_k} - \gamma_k \Exp{\norm{\nabla f(x^k)}} + 2 \gamma_k \Exp{\norm{v^k- \nabla f(x^k)}}  + \frac{\gamma_k^2L}{2} \\
    &\quad + C_{1,k} (\sqrt{1-\alpha} \bar L \gamma_k  + \sqrt{1-\alpha} \sigma_g \eta_k)\\
    &\quad + C_{2,k} \left( \eta_k \sigma_g +  \frac{1-\eta_k}{\eta_k}\gamma_k^2  \frac{\bar L_{h}}{2} \right) \\
     & \leq \Exp{V_k} - \gamma_k \Exp{\norm{\nabla f(x^k)}} + 2 \gamma_k \Exp{\norm{v^k- \nabla f(x^k)}}   \\
    &\quad + \gamma_k^2  \left( \frac{2\sqrt{1-\alpha}}{1-\sqrt{1-\alpha}} \bar L + \frac{L}{2} \right) + \gamma_k^3 \frac{1-\eta_k}{\eta_k} \frac{\sqrt{1-\alpha}}{1-\sqrt{1-\alpha}} \bar L_h\\
    &\quad + \eta_k\gamma_k \frac{4\sqrt{1-\alpha}}{1-\sqrt{1-\alpha}} \sigma_g.
\end{align*}

Next, from the upper-bound of $\Exp{\norm{v^k- \nabla f(x^k)}}$, 
\begin{align*}
   \MoveEqLeft \Exp{ V_{k+1} }\\
   & \leq \Exp{V_k} - \gamma_k \Exp{\norm{\nabla f(x^k)}}\\
   &\quad + 2 \gamma_k \left( \frac{L_h}{2}C(6/7, 4/7) \left( \frac{\gamma_{k}}{\eta_{k}}\right)^2  + \frac{\sigma_g}{\sqrt{n}} \sqrt{C(8/7, 4/7) \eta_{k} }  \right) \\
   &\quad + \gamma_k^2  \left( \frac{2\sqrt{1-\alpha}}{1-\sqrt{1-\alpha}} \bar L + \frac{L}{2} \right)
   + \gamma_k^3 \frac{1-\eta_k}{\eta_k} \frac{\sqrt{1-\alpha}}{1-\sqrt{1-\alpha}} \bar L_h\\
   &\quad + \eta_k\gamma_k \frac{4\sqrt{1-\alpha}}{1-\sqrt{1-\alpha}} \sigma_g \\
   & \leq \Exp{V_k} - \gamma_k \Exp{\norm{\nabla f(x^k)}} + \frac{\gamma_k^3}{\eta_k^2} \cdot  L_h C(6/7, 4/7) \\
   &\quad + \gamma_k \sqrt{\eta_k} \cdot \frac{2 \sigma_g \sqrt{C(8/7, 4/7)}}{\sqrt{n}}
   + \gamma_k^2  \left( \frac{2\sqrt{1-\alpha}}{1-\sqrt{1-\alpha}} \bar L + \frac{L}{2} \right)\\
   &\quad +  \frac{\gamma_k^3}{\eta_k} \frac{\sqrt{1-\alpha}}{1-\sqrt{1-\alpha}} \bar L_h
   + \eta_k\gamma_k \frac{4\sqrt{1-\alpha}}{1-\sqrt{1-\alpha}} \sigma_g.
\end{align*}

\subsubsection{Deriving the convergence rate}
Let $\Gamma_K \eqdef \sum_{k=0}^{K-1} \gamma_k$. By re-arranging the terms and by the telescopic series, 
\begin{align*}
    \MoveEqLeft \frac{1}{\Gamma_K} \sum_{k=0}^{K-1} \gamma_k\Exp{\norm{\nabla f(x^k)}}\\
    & \leq \frac{1}{\Gamma_K} \biggl[ \Exp{V_0} -\Exp{V_K}
    + \left( \frac{2\sqrt{1-\alpha}}{1-\sqrt{1-\alpha}} \bar L + \frac{L}{2} \right) \sum_{k=0}^{K-1} \gamma_k^2 \\
    &\qquad + L_h C(6/7, 4/7) \sum_{k=0}^{K-1} \gamma_k^3\eta_k^{-2}
    + \frac{2\sigma_g \sqrt{C(8/7, 4/7)}}{\sqrt{n}} \sum_{k=0}^{K-1} \gamma_k \sqrt{\eta_k} \\
    &\qquad + \frac{\sqrt{1-\alpha}}{1-\sqrt{1-\alpha}} \bar L_h \sum_{k=0}^{K-1} \gamma_k^3 \eta_k^{-2}
    + \frac{4\sqrt{1-\alpha}}{1-\sqrt{1-\alpha}}  \sigma_g \sum_{k=0}^{K-1} \eta_k\gamma_k \biggr]\\
    & \leq \frac{1}{\Gamma_K} \biggl[ \Exp{V_0}
    + \left( \frac{L}{2} + \frac{2\bar L }{\alpha^2}  \right) \sum_{k=0}^{K-1} \gamma_k^2
    + \frac{2\sigma_g \sqrt{C(8/7, 4/7)}}{\sqrt{n}} \sum_{k=0}^{K-1} \gamma_k \sqrt{\eta_k}\\
    &\qquad + \left( L_h C(6/7, 4/7) + \frac{\bar L_h}{\alpha^2} \right) \sum_{k=0}^{K-1} \gamma_k^3\eta_k^{-2}
    + \frac{4\sigma_g}{\alpha^2} \sum_{k=0}^{K-1} \eta_k\gamma_k \biggr]\\
    & \leq \frac{1}{\Gamma_K} \biggl[ \Exp{V_0}
    + \left( \frac{L}{2} + \frac{2\bar L }{\alpha^2}  \right) \sum_{k=0}^{K-1} \gamma_k^2
    + \left( L_h C_2 + \frac{\bar L_h}{\alpha^2} \right) \sum_{k=0}^{K-1} \gamma_k^3\eta_k^{-2}\\
    &\qquad + \frac{2C_1\sigma_g }{\sqrt{n}} \sum_{k=0}^{K-1} \gamma_k \sqrt{\eta_k}
    + \frac{4\sigma_g}{\alpha^2} \sum_{k=0}^{K-1} \eta_k\gamma_k \biggr],
\end{align*}
where in the second inequality we used $\sqrt{1-\alpha}/(1-\sqrt{1-\alpha}) \leq 1/\alpha^2$ for any $\alpha \in (0,1]$, and in the last step we denoted $C_1 \eqdef \sqrt{C(8/7, 4/7)}$ and $C_2 \eqdef C(6/7, 4/7)$.
Since $\eta_k = (2/(k+2))^{4/7}$ and $\gamma_k = \gamma_0 (2/(k+2))^{5/7}$, for $K \geq 1$
\begin{align*}
    \Gamma_K = \sum_{k=0}^{K-1} \gamma_k
    &\geq K \gamma_{K-1}
    = \gamma_0K  \left(\frac{2}{K+1}\right)^{5/7}
    \geq \gamma_0 K^{2/7};\\
    \sum_{k=0}^{K-1} \gamma_k^3 \eta_k^{-2} &= \gamma_0^3\sum_{k=0}^{K-1}  \left( \frac{2}{k+2} \right)^{15/7} \left( \frac{2}{k+2} \right)^{-8/7} 
    = 2\gamma_0^3 \sum_{k=0}^{K-1} \frac{1}{k+2}\\
    &= 2\gamma_0^3 \sum^{K}_{k=1}\frac{1}{1+k}
    \leq 2\gamma^3_0 \left(\frac{1}{2} + \int^K_{1}\frac{1}{1+s}\,ds\right)
    \leq 2\gamma^3_0 \log \left(K+1\right); \\
    \sum_{k=0}^{K-1} \gamma_k \sqrt{\eta_k}
    &= \gamma_0 \sum_{k=0}^{K-1} \left( \frac{2}{k+2} \right)^{5/7} \left( \frac{2}{k+2} \right)^{2/7}
\\
    &= 2\gamma_0 \sum_{k=0}^{K-1}\frac{1}{k+2}
    \leq 2\gamma_0 \log\left(K+1\right);\\
    \sum_{k=0}^{K-1} \gamma_k^2 &= \gamma_0^2 \sum_{k=0}^{K-1} \left( \frac{2}{k+2} \right)^{10/7} 
    = 2^{10/7}\gamma_0^2 \sum_{k=1}^{K} \frac{1}{(1+k)^{10/7}}\\
    &\leq 2^{10/7}\gamma_0^2 \left(\frac{1}{2^{10/7}} + \int^{\infty}_{1}\frac{1}{(1+s)^{10/7}}\,ds\right)
    = \frac{17}{3}\gamma^2_0 \leq 6\gamma^2_0;\\
    \sum_{k=0}^{K-1} \gamma_k \eta_k &= \gamma_0 \sum_{k=0}^{K-1}  \left( \frac{2}{k+2} \right)^{5/7} \left( \frac{2}{k+2} \right)^{4/7} 
    = 2^{9/7}\gamma_0 \sum^{K}_{k=1} \frac{1}{(1+k)^{9/7}}\\
    &\leq 2^{9/7}\gamma_0 \left(\frac{1}{2^{9/7}} + \int^{\infty}_{1}\frac{1}{(1+s)^{9/7}}\,ds\right)
    = 8\gamma_0.
\end{align*}

Therefore, denoting $\tilde{x}^K$ as a point randomly chosen from $\{x^0,x^1,\ldots,x^{K-1}\}$ with probability $\gamma_k/\Gamma_K$ for $k = 0,1, \ldots, K-1$, we obtain 
\begin{align*}
    \Exp{\norm{\nabla f(\tilde{x}^K)}}  
    & \leq \frac{\Exp{V_0}}{\gamma_0 K^{2/7}} + \left( \frac{L}{2} + \frac{2\bar L }{\alpha^2}  \right) \frac{ 6\gamma^2_0}{\gamma_0 K^{2/7}}\\
    &\quad + \left( L_h C_2 + \frac{\bar L_h}{\alpha^2} \right)\frac{2\gamma^3_0 \log \left(K+1\right)}{\gamma_0 K^{2/7}}\\
    &\quad + \frac{2C_1\sigma_g }{\sqrt{n}} \frac{2\gamma_0 \log \left(K+1\right)}{ \gamma_0 K^{2/7}} + \frac{4\sigma_g}{\alpha^2} \frac{8\gamma_0}{ \gamma_0 K^{2/7}}\\
    &=\widetilde{\cO}\biggl(\frac{1}{K^{2/7}} \biggl( \frac{V_0}{\gamma_0} + \sigma_g\left(\frac{1}{\sqrt{n}} + \frac{1}{\alpha^2}\right) + \gamma_0\left(L+\frac{\bar L}{\alpha^2}\right)\\
    &\qquad\qquad + \gamma_0^2\left(L_h +\frac{\bar L_h}{\alpha^2}\right) \biggr)\biggr).
\end{align*}

\newpage 
\section{\normEFname{RHM}}
\label{ch7:sec:appendix_EF21_RHM_chapter_7}
In \normEFname{RHM}, we update the iterates $\{x^k\}$ according to: 
\begin{equation*}
    x^{k+1} = x^k - \gamma_k \frac{g^k}{\norm{g^k}},
\end{equation*}
where $g^k = (1/n)\sum_{i=1}^n g_i^k$ and $v^k = (1/n)\sum_{i=1}^n v_i^k$ with $g_i^k$ and $v_i^k$ being governed by 
\begin{align*}
    g_i^{k+1} & = g_i^k + \mathcal{C}_i^{k+1}(v_i^{k+1} - g_i^k), \quad \text{and} \\
    v_i^{k+1} & = (1-\eta_k)(v_i^k + \nabla^2 f_i(\hat x^{k+1};\hat{\xi}_i^{k+1}) (x^{k+1} - x^k) ) + \eta_k \nabla f_i(x^{k+1};\xi_i^{k+1}).
\end{align*}
Here, one scalar $q_k\sim\thesischsevencU(0,1)$ is drawn by the master and shared across all clients, $\hat x^{k+1} = q_k x^{k+1} + (1-q_k) x^k$, and the Hessian sample $\hat{\xi}_i^{k+1}$ is independent of the gradient sample $\xi_i^{k+1}$ as specified in Assumption~\ref{ch7:assum:stoc_g_h_chapter_7}.

\subsection{Convergence Proof}

We prove the result in the following steps.  Throughout this proof,
$\mathcal F_k$ is the history from
Assumption~\ref{ch7:assum:stoc_g_h_chapter_7}, immediately before $q_k$ and
the iteration-$k$ oracle samples are drawn, and
$\mathcal H_k=\sigma(\mathcal F_k,q_k)$.  Thus $x^k$, $x^{k+1}$,
$g_i^k$, and $v_i^k$ are $\mathcal F_k$-measurable, whereas
$\hat x^{k+1}$ is $\mathcal H_k$-measurable.  We first state the one-step
oracle estimates conditionally and then use the tower property to pass to
total expectation.

\subsubsection{Deriving the descent inequality} From Lemma~\ref{ch7:lem:descent_lemma_chapter_7}, we obtain \eqref{ch7:eq:descent_lemma_chapter_7}:
\begin{equation*}
    \Delta_{k+1} + \gamma_k \|\nabla f(x^k)\| \le \Delta_k  + 2 \gamma_k \norm{v^k- \nabla f(x^k)} +\frac{ 2\gamma_k}{n}\sum_{i=1}^n\norm{g_i^k-v_i^k} + \frac{\gamma_k^2L}{2}.
\end{equation*}

\subsubsection{Error bound I}

Define $\thesischsevencV_k \eqdef (1/n)\sum_{i=1}^n \norm{g_i^k-v_i^k}$. From~Lemma~\ref{ch7:lemma:boundg_and_v_chapter_7}, 
\begin{equation*}
   \Exp{\thesischsevencV_{k+1}} \leq \sqrt{1-\alpha}    \Exp{\thesischsevencV_{k}} + \frac{\sqrt{1-\alpha}}{n} \sum_{i=1}^n \Exp{\norm{v_i^{k+1}-v_i^k}}.
\end{equation*}

To complete the upper-bound for $\Exp{\thesischsevencV_{k+1}}$, we bound $\Exp{\norm{v_i^{k+1}-v_i^k}}$. 
From the definition of $v_i^{k+1}$, 

\begin{align*}
\MoveEqLeft v_i^{k+1} - v_i^k\\
& = (1-\eta_k)(v_i^k + \nabla^2 f_i(\hat{x}^{k+1}; \hat{\xi}_i^{k+1})(x^{k+1}-x^k))  + \eta_k \nabla f_i(x^{k+1}; \xi_i^{k+1}) - v_i^k \\
& = \eta_k (\nabla f_i(x^k) - v_i^k)  + \eta_k (\nabla f_i(x^{k+1}) - \nabla f_i(x^k)) \\
&\quad + (1-\eta_k)(\nabla^2 f_i(\hat{x}^{k+1}; \hat{\xi}_i^{k+1})(x^{k+1}-x^k) - \nabla^2 f_i(\hat{x}^{k+1})(x^{k+1}-x^k)) \\
&\quad + (1-\eta_k)(\nabla^2 f_i(\hat{x}^{k+1})(x^{k+1}-x^k))\\
&\quad + \eta_k (\nabla f_i(x^{k+1}; \xi_i^{k+1}) - \nabla f_i(x^{k+1})).
\end{align*}

Next, from the definition of the Euclidean norm,  by the triangle inequality, and by the fact that $f_i$ is $L_i$-smooth,
\begin{align*}
\|v_i^{k+1} - v_i^k\|
&\le \eta_k \|v_i^k - \nabla f_i(x^k)\| + \eta_k \|\nabla f_i(x^{k+1}) - \nabla f_i(x^k)\|  \\
&\quad + (1-\eta_k) \|\nabla^2 f_i(\hat{x}^{k+1}; \hat{\xi}_i^{k+1}) - \nabla^2 f_i(\hat{x}^{k+1})\|_{\text{op}} \|x^{k+1}-x^k\|  \\
&\quad + (1-\eta_k) \|\nabla^2 f_i(\hat{x}^{k+1})\|_{\text{op}} \|x^{k+1}-x^k\|  \\
&\quad + \eta_k \|\nabla f_i(x^{k+1}; \xi_i^{k+1}) - \nabla f_i(x^{k+1})\| \\
&\leq \eta_k \|v_i^k - \nabla f_i(x^k)\| + \eta_k \gamma_k L_i  \\
&\quad + (1-\eta_k) \gamma_k \|\nabla^2 f_i(\hat{x}^{k+1}; \hat{\xi}_i^{k+1}) - \nabla^2 f_i(\hat{x}^{k+1})\|_{\text{op}}  \\
&\quad + (1-\eta_k) \gamma_k \|\nabla^2 f_i(\hat{x}^{k+1})\|_{\text{op}}   \\
&\quad + \eta_k \|\nabla f_i(x^{k+1}; \xi_i^{k+1}) - \nabla f_i(x^{k+1})\|,
\end{align*}
where in the last inequality we used the update rule for $x^{k+1}$.
Conditioning on $\mathcal H_k$, applying
Assumption~\ref{ch7:assum:stoc_g_h_chapter_7}, and using
$\norm{\nabla^2 f_i(x)}_{\mathrm{op}}\leq L_i$, we obtain
\begin{align*}
\thesischsevenExpCond{\|v_i^{k+1}-v_i^k\|}{\mathcal H_k}
&\leq \eta_k\|v_i^k-\nabla f_i(x^k)\|
   +\eta_kL_i\gamma_k+(1-\eta_k)L_i\gamma_k \\
&\quad +(1-\eta_k)\gamma_k\sigma_h+\eta_k\sigma_g .
\end{align*}
Taking total expectation and using the tower property gives
\begin{align*}
\Exp{\|v_i^{k+1}-v_i^k\|}
&\leq \eta_k\Exp{\|v_i^k-\nabla f_i(x^k)\|}
   +\eta_kL_i\gamma_k+(1-\eta_k)L_i\gamma_k \\
&\quad +(1-\eta_k)\gamma_k\sigma_h+\eta_k\sigma_g .
\end{align*}
Therefore,
\begin{multline*}
    \Exp{\thesischsevencV_{k+1}} \leq \sqrt{1-\alpha} \Exp{\thesischsevencV_k}
    + \sqrt{1-\alpha} \eta_k \Exp{\thesischsevencU_k}
    + \sqrt{1-\alpha} \gamma_k \bar L \\
    + \sqrt{1-\alpha} (1-\eta_k) \gamma_k \sigma_h + \sqrt{1-\alpha} \eta_k \sigma_g.
\end{multline*}

\subsubsection{Error bound II} 
Define $\thesischsevencU_k \eqdef (1/n)\sum_{i=1}^n \norm{v_i^k - \nabla f_i(x^k)}$. Then, from the definition of  $v_i^{k+1}$,
    \begin{align*}
         \|v_i^{k+1} - \nabla f_i(x^{k+1})\| 
        &= \|(1-\eta_k)(v_i^k - \nabla f_i(x^k)) + (1-\eta_k)\hat S_{i,k+1} + \eta_k e_{i,k+1}\| \\
        &\le (1-\eta_k)\|v_i^k - \nabla f_i(x^k)\| + (1-\eta_k)\|\hat{S}_{i,k+1}\| + \eta_k \|e_{i,k+1}\|,
    \end{align*}
    where 
    \begin{align*}
        \hat{S}_{i,k+1} &= \nabla^2 f_i(\hat{x}^{k+1}; \hat{\xi}_i^{k+1})(x^{k+1}-x^k) - \nabla f_i(x^{k+1}) + \nabla f_i(x^k),\\
        e_{i,k+1} &= \nabla f_i(x^{k+1}; \xi_i^{k+1}) - \nabla f_i(x^{k+1}).
    \end{align*}

Therefore, conditioning on $\mathcal F_k$, applying conditional Jensen's
inequality, and using
Assumption~\ref{ch7:assum:stoc_g_h_chapter_7}, we have
\begin{align*}
\thesischsevenExpCond{\thesischsevencU_{k+1}}{\mathcal F_k}
&\leq (1-\eta_k)\thesischsevencU_k
 +\frac{1-\eta_k}{n}\sum_{i=1}^n
 \sqrt{\thesischsevenExpCond{\sqnorm{\hat S_{i,k+1}}}{\mathcal F_k}}
 +\eta_k\sigma_g .
\end{align*}

Next, we bound the conditional second moment of $\hat S_{i,k+1}$.  Conditional
on $\mathcal H_k$, the sampled-Hessian error is centered.  Hence its cross
term with the remaining $\mathcal H_k$-measurable term vanishes, and the tower
property yields
\begin{align*}
    \MoveEqLeft \thesischsevenExpCond{\norm{\hat{S}_{i,k+1}}^2}{\mathcal F_k}\\
    &= \thesischsevenExpCond{\left\| (\nabla^2 f_i(\hat{x}^{k+1}; \hat{\xi}_i^{k+1})
       - \nabla^2 f_i(\hat{x}^{k+1})) (x^{k+1}-x^k) \right\|^2}{\mathcal F_k} \\
    &\quad + \thesischsevenExpCond{\left\| \nabla^2 f_i(\hat{x}^{k+1}) (x^{k+1}-x^k)
       - (\nabla f_i(x^{k+1}) - \nabla f_i(x^k)) \right\|^2}{\mathcal F_k} \\
    &\leq \thesischsevenExpCond{\left\| \nabla^2 f_i(\hat{x}^{k+1}; \hat{\xi}_i^{k+1})
       - \nabla^2 f_i(\hat{x}^{k+1}) \right\|_{\rm op}^2
       \left\| x^{k+1}-x^k \right\|^2}{\mathcal F_k}  \\
    &\quad + 2\thesischsevenExpCond{\left\| \nabla^2 f_i(\hat{x}^{k+1}) \right\|^2_{\rm op}
       \left\|x^{k+1}-x^k \right\|^2}{\mathcal F_k} \\
    &\quad + 2\left\| \nabla f_i(x^{k+1})-\nabla f_i(x^k)\right\|^2 \\
    &\leq \gamma_k^2\left(\sigma_h^2+4L_i^2\right).
\end{align*}
Here we used the conditional moment bound, $L_i$-smoothness, and
$\norm{x^{k+1}-x^k}\leq\gamma_k$.  Substitution into the preceding conditional
bound gives
\begin{align*}
\thesischsevenExpCond{\thesischsevencU_{k+1}}{\mathcal F_k}
&\leq (1-\eta_k)\thesischsevencU_k
 +(1-\eta_k)\gamma_k(\sigma_h+2\bar L)+\eta_k\sigma_g .
\end{align*}
Finally, taking total expectation and using the tower property,
\begin{align*}
\Exp{\thesischsevencU_{k+1}}
&\leq (1-\eta_k)\Exp{\thesischsevencU_k}
 +(1-\eta_k)\gamma_k(\sigma_h+2\bar L)+\eta_k\sigma_g .
\end{align*}

\subsubsection{Error bound III} 

Next, we bound $\norm{v^{k+1} - \nabla f(x^{k+1})}$: 
\begin{align*}
\hat{e}_{k+1} 
&= (1-\eta_k)\hat{e}_k + (1-\eta_k)\hat{S}_{k+1} + \eta_k e_{k+1} \\
&= \prod_{\tau=0}^k (1-\eta_\tau) \hat{e}_0 
 + \sum_{j=0}^k \left(\prod_{\tau=j}^k (1-\eta_\tau)\right) \hat{S}_{j+1} \\
&\quad + \sum_{j=0}^k \left(\prod_{\tau=j+1}^k (1-\eta_\tau)\right) \eta_j e_{j+1},
\end{align*}
where 
\begin{align*}
\hat S_{k+1} &= \frac{1}{n}\sum_{i=1}^n \left(\nabla^2 f_i(\hat{x}^{k+1}; \hat{\xi}_i^{k+1})(x^{k+1}-x^k) - (\nabla f_i(x^{k+1}) - \nabla f_i(x^k))\right),\\
e_{k+1} &= \frac{1}{n}\sum_{i=1}^n e_{i,k+1}, \text{ and } e_{i,k+1} = \nabla f_i(x^{k+1}; \xi_i^{k+1}) - \nabla f_i(x^{k+1}).
\end{align*}

If $\eta_0=1$, then by taking the Euclidean norm and the expectation, 
\begin{align*}
\Exp{\norm{\hat{e}_{k+1}}} 
&\leq
 \Exp{\norm{ \sum_{j=0}^k \left(\prod_{\tau=j}^k (1-\eta_\tau)\right)  \hat{S}_{j+1}}} \\
&\quad + \Exp{ \norm{ \sum_{j=0}^k \left(\prod_{\tau=j+1}^k (1-\eta_\tau)\right) \eta_j e_{j+1}} }\\
&\leq
 \left(\Exp\left\|\sum_{j=0}^k \prod_{\tau=j}^k (1-\eta_\tau)
 \hat{S}_{j+1}\right\|^2\right)^{1/2} \\
&\quad + \left(\Exp\left\|\sum_{j=0}^k \prod_{\tau=j+1}^k (1-\eta_\tau)
 \eta_j e_{j+1}\right\|^2\right)^{1/2},
\end{align*}
where the first inequality follows from the triangle inequality and the second from
Jensen's inequality.  Assumption~\ref{ch7:assum:stoc_g_h_chapter_7} gives
\[
\thesischsevenExpCond{e_{j+1}}{\mathcal F_j}=0,
\qquad
\thesischsevenExpCond{\sqnorm{e_{j+1}}}{\mathcal F_j}
\leq \frac{\sigma_g^2}{n}.
\]
Consequently, for $j\neq\ell$, the tower property gives
$\Exp{\inner{e_{j+1}}{e_{\ell+1}}}=0$.  Thus,
\begin{align*}
\Exp{\norm{\hat{e}_{k+1}}} 
&\leq \left(\Exp\left\|\sum_{j=0}^k \prod_{\tau=j}^k (1-\eta_\tau)
 \hat{S}_{j+1}\right\|^2\right)^{1/2} \\
&\quad + \frac{\sigma_g}{\sqrt{n}}
 \left( \sum_{j=0}^{k} \prod_{\tau=j+1}^{k}
 (1-\eta_\tau)^2 \eta_j^2 \right)^{1/2}.
\end{align*}

For the Hessian correction, conditional unbiasedness first gives
\[
\thesischsevenExpCond{
\nabla^2 f_i(\hat x^{j+1};\hat\xi_i^{j+1})
}{\mathcal H_j}
=\nabla^2 f_i(\hat x^{j+1}).
\]
Moreover, because $q_j$ is uniform and independent of $\mathcal F_j$, the
fundamental theorem of calculus along the segment from $x^j$ to $x^{j+1}$
gives
\[
\thesischsevenExpCond{
\nabla^2 f_i(\hat x^{j+1})(x^{j+1}-x^j)
}{\mathcal F_j}
=\nabla f_i(x^{j+1})-\nabla f_i(x^j).
\]
Hence $\thesischsevenExpCond{\hat S_{j+1}}{\mathcal F_j}=0$.  For
$j\neq\ell$, the earlier correction is measurable with respect to the later
history, so the tower property yields
$\Exp{\inner{\hat S_{\ell+1}}{\hat S_{j+1}}}=0$.  Therefore,
\begin{align*}
\Exp{\norm{\hat{e}_{k+1}}} 
&\leq \left( \sum_{j=0}^k \prod_{\tau=j}^k (1-\eta_\tau)^2
 \Exp\left\| \hat{S}_{j+1}\right\|^2\right)^{1/2} \\
&\quad + \frac{\sigma_g}{\sqrt{n}}
 \left( \sum_{j=0}^{k} \prod_{\tau=j+1}^{k}
 (1-\eta_\tau)^2 \eta_j^2 \right)^{1/2}.
\end{align*}
Since 
\begin{align*}
\MoveEqLeft \Exp[\|\hat{S}_{j+1}\|^2]\\
&= \Exp\left[\left\|\frac{1}{n}\sum_{i=1}^n \left(
\nabla^2 f_i(\hat{x}^{j+1}; \hat{\xi}_i^{j+1})(x^{j+1}-x^j) 
- (\nabla f_i(x^{j+1}) - \nabla f_i(x^j))
\right)\right\|^2\right] \\
&\leq \frac{1}{n^2}\sum_{i=1}^n
 \Exp{\left\|
 (\nabla^2 f_i(\hat{x}^{j+1}; \hat{\xi}_i^{j+1})
   - \nabla^2 f_i(\hat{x}^{j+1})) (x^{j+1}-x^j)
 \right\|^2} \\
&\quad + \Exp{\left\| \nabla^2 f(\hat{x}^{j+1}) (x^{j+1}-x^j) - (\nabla f(x^{j+1}) - \nabla f(x^j)) \right\|^2} \\
&\le \frac{1}{n^2}\sum_{i=1}^n
 \Exp{\left\|
 (\nabla^2 f_i(\hat{x}^{j+1}; \hat{\xi}_i^{j+1})
   - \nabla^2 f_i(\hat{x}^{j+1})) (x^{j+1}-x^j)
 \right\|^2} \\
&\quad + 2\Exp\left\| \nabla^2 f(\hat{x}^{j+1})\right\|^2_{\text{op}} \left\|x^{j+1}-x^j \right\|^2 \\
&\quad + 2\Exp\left\| \nabla f(x^{j+1}) - \nabla f(x^j) \right\|^2 \\
&\leq \gamma_j^2 \left(\frac{1}{n^2}\sum_{i=1}^n \sigma_{h}^2 + 4L^2 \right) 
\le 4\left(\frac{\sigma_h^2}{n} + L^2\right) \gamma_j^2 ,
\end{align*}
where in the first inequality we applied the variance-bias decomposition, and in the third we used the smoothness of $f$ and Assumption~\ref{ch7:assum:stoc_g_h_chapter_7}, 
we obtain:
\begin{multline*}
\Exp{\norm{\hat{e}_{k+1}}}
 \leq    4 \left(\sum_{j=1}^k \left(\prod_{\tau=j}^k (1-\eta_\tau)^2\right) \gamma_j^2\right)^{1/2} \left(\frac{\sigma_h}{\sqrt{n}} + L\right) \\
 + \frac{\sigma_g}{\sqrt{n}} \left(\sum_{j=0}^k \left(\prod_{\tau=j+1}^k (1-\eta_\tau)^2\right) \eta_j^2\right)^{1/2}.
\end{multline*}

If $\eta_k = (2/(k+2))^{2/3}$ and $\gamma_k = \gamma_0 (2/(k+2))^{2/3}$, then we can prove  that $\eta_0=1$, and from~Lemma~\ref{ch7:lemma:sum_prod_chapter_7} that
\begin{equation*}
\Exp{\norm{\hat{e}_{k+1}}}  
\leq 4  \left(\frac{\sigma_h}{\sqrt{n}} + L\right)  \sqrt{C(4/3, 2/3) \frac{\gamma_{k+1}^2}{\eta_{k+1}}}  + \frac{\sigma_g}{\sqrt{n}} \sqrt{C(4/3, 2/3) \eta_{k+1}}.
\end{equation*}

\subsubsection{Bounding Lyapunov function}

Define $V_k = \Delta_{k} + C_{1,k} \thesischsevencV_k + C_{2,k} \thesischsevencU_k$ with
\begin{equation*}
    C_{1,k} = \frac{2\gamma_k}{1-\sqrt{1-\alpha}} \quad\text{and}\quad C_{2,k} = \frac{2\gamma_k \sqrt{1-\alpha}}{1-\sqrt{1-\alpha}}.
\end{equation*}
Then, 
\begin{align*}
    \Exp{V_{k+1}} 
    & = \Exp{\Delta_{k+1} + C_{1,k+1} \thesischsevencV_{k+1} + C_{2,k+1} \thesischsevencU_{k+1}} \\
    &\leq \Exp{\Delta_k} - \gamma_k \Exp{\norm{\nabla f(x^k)}} + 2 \gamma_k \Exp{\norm{v^k- \nabla f(x^k)}} + 2\gamma_k \Exp{\thesischsevencV_k} \\
    &\quad + \frac{\gamma_k^2L}{2} + C_{1,k+1} \Exp{\thesischsevencV_{k+1}} + C_{2,k+1} \Exp{\thesischsevencU_{k+1}},
\end{align*}
where in the last inequality we used \eqref{ch7:eq:descent_lemma_chapter_7}. 
Next, by the upper-bounds for $\Exp{\thesischsevencV_{k+1}}$ and  $\Exp{\thesischsevencU_{k+1}}$, 
\begin{align*}
    \MoveEqLeft \Exp{ V_{k+1} }\\
   & \leq \Exp{\Delta_k} - \gamma_k \Exp{\norm{\nabla f(x^k)}} + 2 \gamma_k \Exp{\norm{v^k- \nabla f(x^k)}} + 2\gamma_k \Exp{\thesischsevencV_k} + \frac{\gamma_k^2L}{2} \\
    &\quad + C_{1,k+1} (\sqrt{1-\alpha} \Exp{\thesischsevencV_k} + \sqrt{1-\alpha} \eta_k \Exp{\thesischsevencU_k} ) \\
    &\quad + C_{1,k+1} (\sqrt{1-\alpha} \gamma_k \bar L + \sqrt{1-\alpha} (1-\eta_k) \gamma_k \sigma_h + \sqrt{1-\alpha} \eta_k \sigma_g) \\
    &\quad + C_{2,k+1} ( (1-\eta_k)  \Exp{\thesischsevencU_k} + (1-\eta_k)\gamma_k (\sigma_h + 2\bar L) + \eta_k \sigma_g ).
\end{align*}

Next,  since $\gamma_{k+1} \leq \gamma_k$, we can prove that $C_{1,k+1} \leq C_{1,k}$ and that $C_{2,k+1} \leq C_{2,k}$, and also that 
\begin{align*}
    2\gamma_k + C_{1,k+1}\sqrt{1-\alpha}  & \leq 2\gamma_k + C_{1,k}\sqrt{1-\alpha} = C_{1,k}, \\
    C_{1,k+1} \sqrt{1-\alpha}\eta_k + C_{2,k+1} (1-\eta_k) & \leq C_{1,k} \sqrt{1-\alpha}\eta_k + C_{2,k} (1-\eta_k)= C_{2,k}.
\end{align*}
Therefore, 
\begin{align*}
   \Exp{ V_{k+1} } 
   &\leq \Exp{V_k} - \gamma_k \Exp{\norm{\nabla f(x^k)}}
      + 2 \gamma_k \Exp{\norm{v^k- \nabla f(x^k)}} + \frac{\gamma_k^2L}{2} \\
   &\quad + C_{1,k}\bigl(\sqrt{1-\alpha}\gamma_k\bar L
      + \sqrt{1-\alpha}(1-\eta_k)\gamma_k\sigma_h
      + \sqrt{1-\alpha}\eta_k\sigma_g\bigr) \\
   &\quad + C_{2,k}\bigl((1-\eta_k)\gamma_k(\sigma_h+2\bar L)
      + \eta_k\sigma_g\bigr) \\
   &\leq \Exp{V_k} - \gamma_k \Exp{\norm{\nabla f(x^k)}}
      + 2 \gamma_k \Exp{\norm{v^k- \nabla f(x^k)}} \\
   &\quad + \gamma_k^2\left(
      \frac{2\sqrt{1-\alpha}}{1-\sqrt{1-\alpha}}(2\sigma_h+3\bar L)
      + \frac{L}{2}\right) \\
   &\quad + \eta_k\gamma_k
      \frac{4\sqrt{1-\alpha}}{1-\sqrt{1-\alpha}}\sigma_g.
\end{align*}

Next, from the upper-bound of $\Exp{\norm{v^k- \nabla f(x^k)}}$, 
\begin{align*}
   \Exp{ V_{k+1} }
   & \leq \Exp{V_k} - \gamma_k \Exp{\norm{\nabla f(x^k)}} \\
   &\quad + 2 \gamma_k \left(
      4\left(\frac{\sigma_h}{\sqrt{n}}+L\right)
      \sqrt{C(4/3, 2/3)\frac{\gamma_k^2}{\eta_k}}
      + \frac{\sigma_g}{\sqrt{n}}
      \sqrt{C(4/3, 2/3)\eta_k}\right) \\
   &\quad + \gamma_k^2\left(
      \frac{2\sqrt{1-\alpha}}{1-\sqrt{1-\alpha}}(2\sigma_h+3\bar L)
      + \frac{L}{2}\right) 
   + \eta_k\gamma_k
      \frac{4\sqrt{1-\alpha}}{1-\sqrt{1-\alpha}}\sigma_g \\
   & \leq \Exp{V_k} - \gamma_k \Exp{\norm{\nabla f(x^k)}} 
   + \frac{\gamma_k^2}{\sqrt{\eta_k}}\,
      8\left(\frac{\sigma_h}{\sqrt{n}}+L\right)
      \sqrt{C(4/3, 2/3)} \\
   &\quad + \gamma_k\sqrt{\eta_k}\,
      \frac{2\sigma_g\sqrt{C(4/3, 2/3)}}{\sqrt{n}}
   + \gamma_k^2\left(
      \frac{2\sqrt{1-\alpha}}{1-\sqrt{1-\alpha}}(2\sigma_h+3\bar L)
      + \frac{L}{2}\right) \\
   &\quad + \eta_k\gamma_k
      \frac{4\sqrt{1-\alpha}}{1-\sqrt{1-\alpha}}\sigma_g.
\end{align*}

\subsubsection{Deriving the convergence rate}
Let $\Gamma_K \eqdef \sum_{k=0}^{K-1} \gamma_k$. By re-arranging the terms and by the telescopic series, 
\begin{align*}
    \MoveEqLeft \frac{1}{\Gamma_K} \sum_{k=0}^{K-1} \gamma_k \Exp{\norm{\nabla f(x^k)}}\\
    &\leq \frac{1}{\Gamma_K} \biggl[ \Exp{V_0} -\Exp{V_K}
    + 8 \left(\frac{\sigma_h}{\sqrt{n}} + L\right)
      \sqrt{C(4/3, 2/3)} \sum_{k=0}^{K-1} \frac{\gamma_k^2}{\sqrt{\eta_k}}\\
    &\qquad + \frac{2\sigma_g\sqrt{C(4/3, 2/3)}}{\sqrt{n}} \sum_{k=0}^{K-1}\gamma_k\sqrt{\eta_k} 
    + \frac{4\sqrt{1-\alpha}}{1-\sqrt{1-\alpha}}\sigma_g \sum_{k=0}^{K-1}\eta_k\gamma_k\\
    &\qquad + \left(
      \frac{2\sqrt{1-\alpha}}{1-\sqrt{1-\alpha}}(2\sigma_h+3\bar L)
      +\frac{L}{2}\right) \sum_{k=0}^{K-1}\gamma_k^2 \biggr]\\
    &\leq \frac{1}{\Gamma_K} \biggl[ \Exp{V_0}
    + \frac{2\sigma_g\sqrt{C(4/3, 2/3)}}{\sqrt{n}} \sum_{k=0}^{K-1}\gamma_k\sqrt{\eta_k}
    + \frac{4\sigma_g}{\alpha^2} \sum_{k=0}^{K-1}\eta_k\gamma_k \\
    &\qquad + 8\left(\frac{\sigma_h}{\sqrt{n}}+L\right)
      \sqrt{C(4/3, 2/3)} \sum_{k=0}^{K-1}\frac{\gamma_k^2}{\sqrt{\eta_k}}\\
    &\qquad + \left(\frac{L}{2}
      +\frac{2(2\sigma_h+3\bar L)}{\alpha^2}\right) \sum_{k=0}^{K-1}\gamma_k^2 \biggr]\\
    &\leq \frac{1}{\Gamma_K} \biggl[ \Exp{V_0}
    + \frac{2C_1\sigma_g}{\sqrt{n}} \sum_{k=0}^{K-1}\gamma_k\sqrt{\eta_k}
    + \frac{4\sigma_g}{\alpha^2} \sum_{k=0}^{K-1}\eta_k\gamma_k \\
    &\qquad + 8C_1\left(\frac{\sigma_h}{\sqrt{n}}+L\right) \sum_{k=0}^{K-1}\gamma_k^2\eta_k^{-1/2}
    + \left(\frac{L}{2}
      +\frac{4\sigma_h+6\bar L}{\alpha^2}\right) \sum_{k=0}^{K-1}\gamma_k^2 \biggr],
\end{align*}
where in the second inequality we used $\sqrt{1-\alpha}/(1-\sqrt{1-\alpha}) \leq 1/\alpha^2$ for any $\alpha \in (0,1]$, and in the last step we denote $C_1 = \sqrt{C(4/3, 2/3)}$.
Since $\eta_k = (2/(k+2))^{2/3}$ and $\gamma_k = \gamma_0 (2/(k+2))^{2/3}$, 
we obtain
\begin{align*}
    \Gamma_K = \sum_{k=0}^{K-1} \gamma_k &\geq  K \gamma_{K-1} = \gamma_0K  \left(\frac{2}{K+1}\right)^{2/3} \geq \gamma_0 K^{1/3};\\
    \sum_{k=0}^{K-1} \gamma_k^2 \eta_k^{-1/2}  &= \gamma_0^2\sum_{k=0}^{K-1}  \left( \frac{2}{k+2} \right)^{4/3} \left( \frac{2}{k+2} \right)^{-1/3} = 2\gamma_0^2 \sum_{k=0}^{K-1} \frac{1}{k+2}\\
    &= 2\gamma_0^2 \sum^{K}_{k=1}\frac{1}{1+k}
    \leq 2\gamma^2_0 \left(\frac12+\int^K_{1}\frac{1}{1+s}\,ds\right) \leq 2\gamma^2_0 \log \left(K+1\right); \\
    \sum_{k=0}^{K-1} \gamma_k \sqrt{\eta_k} &= \gamma_0 \sum_{k=0}^{K-1} \left( \frac{2}{k+2} \right)^{2/3} \left( \frac{2}{k+2} \right)^{1/3}\\
    &=   2\gamma_0 \sum_{k=0}^{K-1}\frac{1}{k+2} \leq 2\gamma_0 \log\left(K+1\right);\\
    \sum_{k=0}^{K-1} \gamma_k^2 &= 2^{4/3}\gamma_0^2 \sum_{k=1}^{K} \frac{1}{(1+k)^{4/3}}\\
    &\leq 2^{4/3}\gamma_0^2
    \left(\frac{1}{2^{4/3}}
    +\int^\infty_{1}\frac{1}{(1+s)^{4/3}}\,ds\right)
    = 7\gamma_0^2;\\
    \sum_{k=0}^{K-1} \gamma_k \eta_k &= \frac{1}{\gamma_0}\sum_{k=0}^{K-1} \gamma_k^2 \leq 7\gamma_0.
\end{align*}

Therefore, denoting $\tilde{x}^K$ as a point randomly chosen from $\{x^0,x^1,\ldots,x^{K-1}\}$ with probability $\gamma_k/\Gamma_K$ for $k = 0,1, \ldots, K-1$, we obtain 
\begin{align*}
    \Exp{\norm{\nabla f(\tilde{x}^K)}} 
    &\leq \frac{\Exp{V_0}}{\gamma_0 K^{1/3}}
      + \frac{2C_1\sigma_g}{\sqrt{n}}
        \frac{2\gamma_0\log(K+1)}{\gamma_0K^{1/3}} 
      + \frac{4\sigma_g}{\alpha^2}
        \frac{7\gamma_0}{\gamma_0K^{1/3}} \\
    &\quad + 8C_1\left(\frac{\sigma_h}{\sqrt{n}}+L\right)
        \frac{2\gamma_0^2\log(K+1)}
             {\gamma_0K^{1/3}} \\
    &\quad + \left(\frac{L}{2}
        +\frac{4\sigma_h+6\bar L}{\alpha^2}\right)
        \frac{7\gamma_0^2}{\gamma_0K^{1/3}} \\
    &= \frac{\Exp{V_0}}{\gamma_0 K^{1/3}}
      + \frac{4C_1\sigma_g}{\sqrt{n}}
        \frac{\log(K+1)}{K^{1/3}}
      + \frac{28\sigma_g}{\alpha^2K^{1/3}} \\
    &\quad + 16C_1\gamma_0\left(\frac{\sigma_h}{\sqrt{n}}+L\right)
        \frac{\log(K+1)}{K^{1/3}} \\
    &\quad + \frac{7\gamma_0}{2}\left(
        L+\frac{8\sigma_h}{\alpha^2}
        +\frac{12\bar L}{\alpha^2}\right)
        \frac{1}{K^{1/3}} \\
    &= \widetilde{\cO}\biggl(\frac{1}{K^{1/3}} \biggl(
       \frac{\Exp{V_0}}{\gamma_0}
       +\sigma_g\left(\frac{1}{\sqrt{n}}+\frac{1}{\alpha^2}\right)\\
    &\qquad\qquad +\gamma_0\sigma_h\left(\frac{1}{\sqrt{n}}+\frac{1}{\alpha^2}\right)
       +\gamma_0\left(L+\frac{\bar L}{\alpha^2}\right) \biggr)\biggr).
\end{align*}

\newpage 
\section{\normEFname{HM}}

In \normEFname{HM}, we update the iterates $\{x^k\}$ according to: 
\begin{equation*}
    x^{k+1} = x^k - \gamma_k \frac{g^k}{\norm{g^k}}, 
\end{equation*}
where $g^k = (1/n)\sum_{i=1}^n g_i^k$ and $v^k = (1/n)\sum_{i=1}^n v_i^k$ with $g_i^k,v_i^k$ being governed by 
\begin{align}
    g_i^{k+1} & = g_i^k + \mathcal{C}_i^{k+1}(v_i^{k+1} - g_i^k), \quad \text{and} \notag\\
    v_i^{k+1} & = (1-\eta_k)(v_i^k + \nabla^2 f_i(x^{k+1};\xi_i^{k+1}) (x^{k+1} - x^k) ) + \eta_k \nabla f_i(x^{k+1};\xi_i^{k+1}). \notag
\end{align}

\subsection{Convergence Proof}

We prove the result in the following steps.

\subsubsection{Deriving the descent inequality} From Lemma~\ref{ch7:lem:descent_lemma_chapter_7}, we obtain \eqref{ch7:eq:descent_lemma_chapter_7}:
\begin{equation*}
    \Delta_{k+1} + \gamma_k \|\nabla f(x^k)\| \le \Delta_k  + 2 \gamma_k \norm{v^k- \nabla f(x^k)} +\frac{ 2\gamma_k}{n}\sum_{i=1}^n\norm{g_i^k-v_i^k} + \frac{\gamma_k^2L}{2}.
\end{equation*}

\subsubsection{Error bound I}
Define $\thesischsevencV_k \eqdef (1/n)\sum_{i=1}^n \norm{g_i^k-v_i^k}$. From~Lemma~\ref{ch7:lemma:boundg_and_v_chapter_7}, we have
\begin{align*}
   \Exp{\thesischsevencV_{k+1}} & \leq \sqrt{1-\alpha}    \Exp{\thesischsevencV_{k}} + \frac{\sqrt{1-\alpha}}{n} \sum_{i=1}^n \Exp{\norm{v_i^{k+1}-v_i^k}}.
\end{align*}

To complete the upper-bound for $\Exp{\thesischsevencV_{k+1}}$, we bound $\Exp{\norm{v_i^{k+1}-v_i^k}}$. 
Since
\begin{align*}
v_i^{k+1}-v_i^k
&= - \eta_k (v_i^k - \nabla f_i(x^k)) \\
&+ (1-\eta_k) ( \nabla^2 f_i(x^{k+1};\xi_i^{k+1})  - \nabla^2 f_i(x^{k+1}))(x^{k+1} - x^k) \\
&-\eta_k\left[ \nabla f_i(x^k) - \nabla f_i(x^{k+1}) - \nabla^2 f_i(x^{k+1})(x^k -x^{k+1}) \right] \\
&+  \nabla^2 f_i(x^{k+1})(x^{k+1}-x^k) \\
&+ \eta_k \left(\nabla f_i(x^{k+1};\xi_i^{k+1}) - \nabla f_i(x^{k+1})\right) \\
&= - \eta_k (v_i^k - \nabla f_i(x^k)) \\
&+ (1-\eta_k) S_{i,k+1} \\
&+ \nabla^2 f_i(x^{k+1})(x^{k+1}-x^k) \\
&- \eta_k Z_{f_i}(x^{k}, x^{k+1}) \\
&+ \eta_k e_{i,k+1},
\end{align*}
where we denoted
\begin{align*}
S_{i,k+1} &\eqdef ( \nabla^2 f_i(x^{k+1};\xi_i^{k+1})  - \nabla^2 f_i(x^{k+1}))(x^{k+1} - x^k),\\
Z_{f_i}(x^{k}, x^{k+1}) &\eqdef \nabla f_i(x^k) - \nabla f_i(x^{k+1}) - \nabla^2 f_i(x^{k+1})(x^k -x^{k+1}),\\
e_{i,k+1} &\eqdef \nabla f_i(x^{k+1};\xi_i^{k+1}) - \nabla f_i(x^{k+1}),
\end{align*}
we obtain 
\begin{align*}
    \MoveEqLeft \Exp{\norm{v_i^{k+1}-v_i^k}}\\
    & \leq \eta_k \Exp{\norm{v_i^k - \nabla f_i(x^k)}} + (1-\eta_k) \Exp{\norm{S_{i,k+1}}} \\  
    &\quad + \eta_k \Exp{\norm{Z_{f_i}(x^{k}, x^{k+1})}} + \eta_k\Exp{\norm{ e_{i,k+1}}}\\
    &\quad + \Exp{\norm{\nabla^2 f_i(x^{k+1})(x^{k+1}-x^k)}}\\
    & \leq \eta_k \Exp{\norm{v_i^k - \nabla f_i(x^k)}} + \Exp{\norm{\nabla^2 f_i(x^{k+1})}_{\text{op}}\norm{x^{k+1}-x^k}}\\
    &\quad + (1-\eta_k) \Exp{\norm{\nabla^2 f_i(x^{k+1};\xi_i^{k+1})  - \nabla^2 f_i(x^{k+1})}\norm{x^{k+1}-x^k}} \\  
    &\quad + \frac{\eta_k L_{h,i}}{2} \Exp{\norm{x^{k+1} -x^k}^2} + \eta_k\Exp{\norm{ e_{i,k+1}}}\\
    & \leq \eta_k \Exp{\norm{v_i^k - \nabla f_i(x^k)}} + L_{i}\gamma_k + (1-\eta_k)\gamma_k\sigma_h + \frac{1}{2}\eta_k\gamma_k^2 L_{h,i} + \eta_k \sigma_g
\end{align*}
where in the first inequality we used the triangle inequality, in the second we used Assumption~\ref{ch7:assum:L_Hessian_comp_fcn_chapter_7}, 
and in the third we used Assumption~\ref{ch7:assum:L_comp_fcn_chapter_7}, Assumption~\ref{ch7:assum:stoc_g_h_chapter_7} and the update rule for $x^{k+1}$. 

Therefore, we have 
\begin{align*}
   \Exp{\thesischsevencV_{k+1}}  
   & \leq  \sqrt{1-\alpha}    \Exp{\thesischsevencV_{k}} + \frac{\sqrt{1-\alpha}}{n} \sum_{i=1}^n \Exp{\norm{v_i^{k+1}-v_i^k}} \\
   & \leq  \sqrt{1-\alpha}    \Exp{\thesischsevencV_{k}} + \sqrt{1-\alpha} \eta_k \Exp{\thesischsevencU_k}  \\
   &\quad + \sqrt{1-\alpha}(1-\eta_k) \gamma_k \sigma_h + \frac{\sqrt{1-\alpha}}{2} \gamma_k^2 \bar L_{h}\\
   &\quad + \sqrt{1-\alpha} \gamma_k \bar L + \sqrt{1-\alpha} \eta_k \sigma_g,
\end{align*}
where $\bar L_h = (1/n)\sum_{i=1}^n L_{h,i}$ and $\bar L = (1/n)\sum_{i=1}^n L_{i}$.

\subsubsection{Error bound II} 

Define $\thesischsevencU_k \eqdef (1/n)\sum_{i=1}^n \norm{v_i^k - \nabla f_i(x^k)}$. Then, using the update rule for $x^{k+1}$, we have 
\begin{align*}
    v_i^{k+1} - \nabla f_i(x^{k+1}) & = (1-\eta_k)(v_i^k - \nabla f_i(x^k)) \\
    &\quad + (1-\eta_k)(\nabla^2 f_i(x^{k+1};\xi_i^{k+1}) - \nabla^2 f_i(x^{k+1}))(x^{k+1}-x^k) \\
    &\quad + (1-\eta_k) ( \nabla^2 f_i(x^{k+1})(x^{k+1}-x^k) - [\nabla f_i(x^{k+1}) - \nabla f_i(x^k)] ) \\
    &\quad + \eta_k (\nabla f_i(x^{k+1};\xi_i^{k+1}) - \nabla f_i(x^{k+1}))\\
    & = (1-\eta_k)(v_i^k - \nabla f_i(x^k)) + (1-\eta_k)S_{i,k+1}\\
    &\quad+ (1-\eta_k)Z_{f_i}(x^k, x^{k+1}) + \eta_k e_{i,k+1}.
\end{align*} 
Therefore, using triangle inequality, we obtain
\begin{align*}
    \Exp{\thesischsevencU_{k+1}} & \leq   \frac{1-\eta_k}{n}\sum_{i=1}^n \Exp{\norm{v_i^{k} - \nabla f_i(x^{k})}}   + (1-\eta_k) \frac{1}{n}\sum_{i=1}^n \Exp{\norm{S_{i,k+1}}}  \\
    &\quad + (1-\eta_k) \frac{1}{n}\sum_{i=1}^n \Exp{\norm{ Z_{f_i}(x^k,x^{k+1}) } }   + \eta_k \frac{1}{n}\sum_{i=1}^n \Exp{\norm{e_{i,k+1}}}\\
    & \leq \frac{1-\eta_k}{n}\sum_{i=1}^n \Exp{\norm{v_i^{k} - \nabla f_i(x^{k})}}\\
    &\quad +  \frac{1-\eta_k}{n}\sum_{i=1}^n \Exp{\norm{\nabla^2 f_i(x^{k+1};\xi^{k+1}_i) - \nabla^2 f_i(x^{k+1})}\norm{x^{k+1}-x^k}}  \\
    &\quad + \frac{1-\eta_k}{2n}\sum_{i=1}^n L_{h,i}\Exp{\sqnorm{ x^{k+1}-x^k } } + \frac{\eta_k }{n}\sum_{i=1}^n \Exp{\norm{e_{i,k+1}}},
\end{align*}
where in the second inequality we used Assumption~\ref{ch7:assum:L_Hessian_comp_fcn_chapter_7}.
Next, bounding the second term and the fourth term via Assumption~\ref{ch7:assum:stoc_g_h_chapter_7}, we get
\begin{align*}
    \Exp{\thesischsevencU_{k+1}} & \leq  \frac{1-\eta_k}{n}\sum_{i=1}^n \Exp{\norm{v_i^{k} - \nabla f_i(x^{k})}}  + (1-\eta_k) \sigma_h \Exp{\norm{x^{k+1}-x^k}}  \\
    &\quad + (1-\eta_k) \frac{\bar L_{h}}{2} \Exp{\sqnorm{x^{k+1}-x^k}} + \eta_k \sigma_g\\
    & \leq (1-\eta_k)\Exp{\thesischsevencU_{k}} + (1-\eta_k) \gamma_k\sigma_h  + (1-\eta_k)\gamma_k^2 \bar{L}_h + \eta_k\sigma_g,
\end{align*}
where in the last inequality we used the update rule for $x^{k+1}$.

\subsubsection{Error bound III} 
From the definition of $v^k = (1/n)\sum_{i=1}^n v_i^k$ and denoting $\hat{e}_{k+1} \eqdef v^{k+1} - \nabla f(x^{k+1}) $, we have 
\begin{align*}
    \hat{e}_{k+1}
    & =  \frac{1}{n}\sum_{i=1}^n (v_i^{k+1} - \nabla f_i(x^{k+1}))\\
    & = (1-\eta_k)(v^k - \nabla f(x^k)) \\
    &\quad + (1-\eta_k) \frac{1}{n}\sum_{i=1}^n (\nabla^2 f_i(x^{k+1};\xi_i^{k+1}) - \nabla^2 f_i(x^{k+1}))(x^{k+1}-x^k) \\
    &\quad + (1-\eta_k) \frac{1}{n}\sum_{i=1}^n ( \nabla^2 f_i(x^{k+1})(x^{k+1}-x^k) - [\nabla f_i(x^{k+1}) - \nabla f_i(x^k)] ) \\
    &\quad + \eta_k \frac{1}{n}\sum_{i=1}^n (\nabla f_i(x^{k+1};\xi_i^{k+1}) - \nabla f_i(x^{k+1}))\\
    & = (1-\eta_k)\hat{e}_k + (1-\eta_k)S_{k+1} + (1-\eta_k)Z_{f}(x^k, x^{k+1}) + \eta_k e_{k+1},
\end{align*}
where we used the following notation
\begin{align*}
    S_{k+1} & \eqdef \frac{1}{n}\sum^n_{i=1} S_{i, k+1} = \frac{1}{n}\sum_{i=1}^n (\nabla^2 f_i(x^{k+1};\xi_i^{k+1}) - \nabla^2 f_i(x^{k+1}))(x^{k+1}-x^k);\\
    Z_f(x^k, x^{k+1}) & \eqdef \nabla f(x^k) -\nabla f(x^{k+1}) - \nabla^2 f(x^{k+1})(x^{k} -x^{k+1});\\
    e_{k+1} & \eqdef \frac{1}{n}\sum^n_{i=1} e_{i,k+1} = \frac{1}{n}\sum_{i=1}^n (\nabla f_i(x^{k+1};\xi_i^{k+1}) - \nabla f_i(x^{k+1})).
\end{align*}
Therefore, we obtain
\begin{align*}
    \hat{e}_{k+1} 
    & = \prod_{\tau=0}^k (1-\eta_\tau) \hat{e}_{0}   + \sum_{j=0}^k \left(\prod_{\tau=j}^k (1-\eta_\tau)\right) S_{j+1}\\
    &\quad+ \sum_{j=0}^k \left(\prod_{\tau=j}^k (1-\eta_\tau)\right) Z_f(x^j, x^{j+1})  + \sum_{j=0}^k \left(\prod_{\tau=j+1}^k (1-\eta_\tau)\right) \eta_j e_{j+1}.
\end{align*}
If $\eta_0 = 1$, then 
\begin{align*}
  \hat{e}_{k+1} 
    & = \sum_{j=0}^k \left(\prod_{\tau=j}^k (1-\eta_\tau)\right) S_{j+1} + \sum_{j=0}^k \left(\prod_{\tau=j}^k (1-\eta_\tau)\right) Z_f(x^j, x^{j+1})  \\
    &\quad + \sum_{j=0}^k \left(\prod_{\tau=j+1}^k (1-\eta_\tau)\right) \eta_j e_{j+1}.
\end{align*}

Therefore, from the definition of the Euclidean norm, by the triangle inequality, and by Jensen's inequality. 
\begin{align*}
 \Exp{  \norm{\hat{e}_{k+1}}  }
 &\leq   \sqrt{ \Exp{ \sqnorm{\sum_{j=0}^k \left(\prod_{\tau=j}^k (1-\eta_\tau)\right) {S}_{j+1}}  } } \\
 &\quad + \Exp{ \norm{\sum_{j=0}^k \left(\prod_{\tau=j}^k (1-\eta_\tau)\right) Z_f(x^j, x^{j+1}) } } \\
 &\quad + \sqrt{ \Exp{ \sqnorm{\sum_{j=0}^k \left(\prod_{\tau=j+1}^k (1-\eta_\tau)\right) \eta_j e_{j+1}} } }.
\end{align*}

Next, from Assumption~\ref{ch7:assum:stoc_g_h_chapter_7}, we can prove that $\Exp{S_{k+1}} = 0$, $\Exp{ \inner{S_l}{S_{i}}} = 0$ for $l \neq i$, $\Exp{e_{k+1}}=0$, $\Exp{\inner{e_l}{e_i}}=0$ for $l \neq i$.  Hence, since $\prod_{\tau=j}^k (1-\eta_\tau) \leq \prod_{\tau=j+1}^k (1-\eta_\tau)$,
\begin{align*}
    \Exp{  \norm{\hat{e}_{k+1}}  } & \leq \sqrt{ \sum_{j=0}^k \left(\prod_{\tau=j+1}^k (1-\eta_\tau)^2 \right) \Exp{ \sqnorm{ {S}_{j+1}}  } }\\
    &\quad + \sqrt{ \sum_{j=0}^k \left(\prod_{\tau=j+1}^k (1-\eta_\tau)^2\right) \eta_j^2 \Exp{ \sqnorm{ e_{j+1}} } }\\
    &\quad + \sum_{j=0}^k \left(\prod_{\tau=j+1}^k (1-\eta_\tau)\right)\Exp{ \norm{ Z_f(x^j, x^{j+1}) } }.
\end{align*}

Next, from Assumption~\ref{ch7:assum:stoc_g_h_chapter_7} and from the definition of $x^{k+1}$, we can prove that 
\begin{align*}
    \MoveEqLeft \Exp{\sqnorm{S_{j+1}}} \\
    & \leq \frac{1}{n^2} \sum_{i=1}^n \Exp{\sqnorm{\nabla^2 f_i(x^{j+1};\xi_i^{j+1}) - \nabla^2 f_i(x^{j+1})} \cdot \sqnorm{x^{j+1} - x^j}} \\
    & \leq \frac{\gamma_j^2}{n} \sigma_h^2, \quad\text{and}\\
    \Exp{\sqnorm{e_{j+1}}} & \leq \frac{1}{n^2}\sum_{i=1}^n \Exp{\sqnorm{ \nabla f_i(x^{j+1};\xi_i^{j+1}) - \nabla f_i(x^{j+1})  }} \leq \frac{\sigma_g^2}{n}.
\end{align*}
Therefore, 
\begin{align*}
    \Exp{  \norm{\hat{e}_{k+1}}  } & \leq  \frac{\sigma_h}{\sqrt{n}}\sqrt{ \sum_{j=0}^k \left(\prod_{\tau=j+1}^k (1-\eta_\tau)^2 \right) \gamma_j^2  } + \frac{\sigma_g}{\sqrt{n}}\sqrt{ \sum_{j=0}^k \left(\prod_{\tau=j+1}^k (1-\eta_\tau)^2\right) \eta_j^2 }\\
    &\quad + \sum_{j=0}^k \left(\prod_{\tau=j+1}^k (1-\eta_\tau)\right)\Exp{ \norm{ Z_f(x^j, x^{j+1}) } }.
\end{align*}

Next, by Assumption~\ref{ch7:assum:L_Hessian_fcn_chapter_7},  
\begin{align*}
    \Exp{  \norm{\hat{e}_{k+1}}  } & \leq \frac{\sigma_h}{\sqrt{n}}\sqrt{ \sum_{j=0}^k \left(\prod_{\tau=j+1}^k (1-\eta_\tau)^2 \right) \gamma_j^2  } + \frac{ L_h}{2} \sum_{j=0}^k \left(\prod_{\tau=j+1}^k (1-\eta_\tau)\right)  \gamma_j^2  \\
    &\quad + \frac{\sigma_g}{\sqrt{n}}\sqrt{ \sum_{j=0}^k \left(\prod_{\tau=j+1}^k (1-\eta_\tau)^2\right) \eta_j^2 } 
\end{align*}

If $\eta_k = (2/(k+2))^{2/3}$ and $\gamma_k = \gamma_0 (2/(k+2))^{2/3}$, then we can prove  that $\eta_0=1$, and from~Lemma~\ref{ch7:lemma:sum_prod_chapter_7} that
\begin{align*}
    \Exp{  \norm{v^{k+1} - \nabla f(x^{k+1})}  } 
    & \leq \frac{\sigma_h }{\sqrt{n}} \sqrt{C(4/3,2/3) \frac{\gamma_{k+1}^2}{\eta_{k+1}}} + \frac{C(4/3,2/3)}{2} L_h \frac{\gamma_{k+1}^2}{\eta_{k+1}} \\
    &\quad + \frac{\sigma_g}{\sqrt{n}} \sqrt{C(4/3,2/3) \eta_{k+1}}.
\end{align*}

\subsubsection{Bounding Lyapunov function}
Define our Lyapunov function $V_k = \Delta_{k} + C_{1,k} \thesischsevencV_k + C_{2,k} \thesischsevencU_k$ with
\begin{equation*}
    C_{1,k} = \frac{2\gamma_k}{1-\sqrt{1-\alpha}} \quad\text{and}\quad C_{2,k} = \frac{2\gamma_k \sqrt{1-\alpha}}{1-\sqrt{1-\alpha}}.
\end{equation*}
Then, using \eqref{ch7:eq:descent_lemma_chapter_7}, we have 
\begin{align*}
    \Exp{V_{k+1}} 
    & = \Exp{\Delta_{k+1} + C_{1,k+1} \thesischsevencV_{k+1} + C_{2,k+1} \thesischsevencU_{k+1}} \\
    & \leq \Exp{\Delta_k} - \gamma_k \Exp{\norm{\nabla f(x^k)}} + 2 \gamma_k \Exp{\norm{v^k- \nabla f(x^k)}} + 2\gamma_k \Exp{\thesischsevencV_k} + \frac{\gamma_k^2L}{2} \\
    &\quad + C_{1,k+1} \Exp{\thesischsevencV_{k+1}} + C_{2,k+1} \Exp{\thesischsevencU_{k+1}}.
\end{align*}

Next, by the upper-bounds for $\Exp{\thesischsevencV_{k+1}}$ and  $\Exp{\thesischsevencU_{k+1}}$, 
\begin{align*}
   \MoveEqLeft \Exp{ V_{k+1} }\\
   & \leq \Exp{\Delta_k} - \gamma_k \Exp{\norm{\nabla f(x^k)}} + 2 \gamma_k \Exp{\norm{v^k- \nabla f(x^k)}} + 2\gamma_k \Exp{\thesischsevencV_k} + \frac{\gamma_k^2L}{2} \\
    &\quad + C_{1,k+1} ( \sqrt{1-\alpha}    \Exp{\thesischsevencV_{k}} + \sqrt{1-\alpha} \eta_k \Exp{\thesischsevencU_k})  \\
    &\quad + C_{1,k+1}(\sqrt{1-\alpha}(1-\eta_k) \gamma_k \sigma_h + \frac{\sqrt{1-\alpha}}{2} \gamma_k^2 \bar L_{h}\\
    &\qquad\qquad  + \sqrt{1-\alpha} \gamma_k \bar L + \sqrt{1-\alpha} \eta_k \sigma_g ) \\
    &\quad + C_{2,k+1} ( (1-\eta_k) \Exp{\thesischsevencU_k}  + (1-\eta_k) \gamma_k \sigma_h + (1-\eta_k) \gamma_k^2 \bar L_h   + \eta_k \sigma_g ).
\end{align*}

Next,  since $\gamma_{k+1} \leq \gamma_k$, we can prove that $C_{1,k+1} \leq C_{1,k}$ and that $C_{2,k+1} \leq C_{2,k}$, and also that 
\begin{align*}
    2\gamma_k + C_{1,k+1}\sqrt{1-\alpha}  & \leq 2\gamma_k + C_{1,k}\sqrt{1-\alpha} = C_{1,k} \\
    C_{1,k+1} \sqrt{1-\alpha}\eta_k + C_{2,k+1} (1-\eta_k) & \leq C_{1,k} \sqrt{1-\alpha}\eta_k + C_{2,k} (1-\eta_k)= C_{2,k}
\end{align*}
Therefore, 
\begin{align*}
    \MoveEqLeft \Exp{ V_{k+1} } \\
    & \leq \Exp{V_k} - \gamma_k \Exp{\norm{\nabla f(x^k)}} + 2 \gamma_k \Exp{\norm{v^k- \nabla f(x^k)}}  + \frac{\gamma_k^2L}{2} \\
    &\quad + C_{1,k} (\sqrt{1-\alpha}(1-\eta_k) \gamma_k \sigma_h + \frac{\sqrt{1-\alpha}}{2} \gamma_k^2 \bar L_{h}\\
    &\qquad\qquad  + \sqrt{1-\alpha} \gamma_k \bar L + \sqrt{1-\alpha} \eta_k \sigma_g ) \\
    &\quad + C_{2,k} \left( (1-\eta_k) \gamma_k \sigma_h + (1-\eta_k) \gamma_k^2 \bar L_h   + \eta_k \sigma_g  \right) \\
    & \leq \Exp{V_k} - \gamma_k \Exp{\norm{\nabla f(x^k)}} + 2 \gamma_k \Exp{\norm{v^k- \nabla f(x^k)}}   \\
    &\quad + \gamma_k^2  \left( \frac{4\sqrt{1-\alpha}}{1-\sqrt{1-\alpha}} \sigma_h + \frac{2\sqrt{1-\alpha}}{1-\sqrt{1-\alpha}} \bar L  + \frac{L}{2}\right) \\
    &\quad + \gamma_k^3 \left( \frac{4\sqrt{1-\alpha}}{1-\sqrt{1-\alpha}} \right) \bar L_h + \eta_k\gamma_k \frac{4\sqrt{1-\alpha}}{1-\sqrt{1-\alpha}} \sigma_g.
\end{align*}

Next, from the upper-bound of $\Exp{\norm{v^k- \nabla f(x^k)}}$, 
\begin{align*}
    \MoveEqLeft \Exp{ V_{k+1} } \\
    & \leq \Exp{V_k} - \gamma_k \Exp{\norm{\nabla f(x^k)}} \\
    &\quad + 2 \gamma_k \left( \frac{\sigma_h }{\sqrt{n}} \sqrt{C(4/3,2/3) \frac{\gamma_{k}^2}{\eta_{k}}} + \frac{ L_h C(4/3,2/3)}{2} \frac{\gamma_{k}^2}{\eta_{k}}\right.\\
    &\qquad\qquad \left.  + \frac{\sigma_g}{\sqrt{n}} \sqrt{C(4/3,2/3) \eta_{k}}  \right)  \\
    &\quad + \gamma_k^2  \left( \frac{4\sqrt{1-\alpha}}{1-\sqrt{1-\alpha}} \sigma_h + \frac{2\sqrt{1-\alpha}}{1-\sqrt{1-\alpha}} \bar L  + \frac{L}{2}\right) \\
    &\quad + \gamma_k^3 \left( \frac{4\sqrt{1-\alpha}}{1-\sqrt{1-\alpha}} \right) \bar L_h + \eta_k\gamma_k \frac{4\sqrt{1-\alpha}}{1-\sqrt{1-\alpha}} \sigma_g \\
    & \leq \Exp{V_k} - \gamma_k \Exp{\norm{\nabla f(x^k)}}  + \frac{2 C_1\sigma_h }{\sqrt{n}}\frac{\gamma_k^2}{\sqrt{\eta_k}}   + C^2_1 L_h \frac{\gamma_k^3}{\eta_k}\\
    &\quad + \frac{2C_1\sigma_g }{\sqrt{n}} \gamma_k \sqrt{\eta_k}  \\
    &\quad +   \left( \frac{L}{2} +\frac{4\sigma_h}{\alpha^2}  + \frac{2\bar L}{\alpha^2} \right) \gamma_k^2\\
    &\quad + \frac{4\bar L_h}{\alpha^2}\gamma_k^3  + \frac{4\sigma_g}{\alpha^2} \eta_k\gamma_k.
\end{align*}
where we used the notation $C_1 = \sqrt{ C(4/3,2/3) }$ and used $\sqrt{1-\alpha}/(1-\sqrt{1-\alpha}) \leq 1/\alpha^2$.

\subsubsection{Deriving the convergence rate}
Let $\Gamma_K \eqdef \sum_{k=0}^{K-1} \gamma_k$. By re-arranging the terms and by the telescopic series, 
\begin{align*}
    \MoveEqLeft \frac{1}{\Gamma_K} \sum_{k=0}^{K-1} \gamma_k \Exp{\norm{\nabla f(x^k)}}\\
    & \leq \frac{1}{\Gamma_K} \biggl[ \Exp{V_0} -\Exp{V_K}   + \frac{2 C_1 \sigma_h }{\sqrt{n}} \sum_{k=0}^{K-1} \gamma_k^2\eta_k^{-1/2} \\
    &\qquad + C^2_1 L_h  \sum_{k=0}^{K-1} \gamma_k^3 \eta_k^{-1}  + \frac{2C_1 \sigma_g}{\sqrt{n}} \sum_{k=0}^{K-1} \gamma_k \sqrt{\eta_k} \\
    &\qquad +  \left(\frac{L}{2} + \frac{4\sigma_h}{\alpha^2}  + \frac{2\bar L}{\alpha^2} \right) \sum_{k=0}^{K-1} \gamma_k^2 \\
    &\qquad +     \frac{4\bar L_h}{\alpha^2}   \sum_{k=0}^{K-1} \gamma_k^3 + \frac{4\sigma_g }{\alpha^2}  \sum_{k=0}^{K-1} \eta_k\gamma_k \biggr].
\end{align*}

Since $\eta_k = (2/(k+2))^{2/3}$ and $\gamma_k = \gamma_0 (2/(k+2))^{2/3}$,  for $K \geq 1$ we obtain
\begin{align*}
    \Gamma_K = \sum_{k=0}^{K-1} \gamma_k & \geq  K \gamma_{K-1} = \gamma_0K  \left(\frac{2}{K+1}\right)^{2/3} \geq \gamma_0 K^{1/3};\\
    \sum_{k=0}^{K-1} \gamma_k^2 \eta_k^{-1/2}  & = \gamma_0^2\sum_{k=0}^{K-1}  \left( \frac{2}{k+2} \right)^{4/3} \left( \frac{2}{k+2} \right)^{-1/3 } = 2\gamma_0^2 \sum_{k=0}^{K-1} \frac{1}{k+2}\\
    &= 2\gamma_0^2 \sum^{K}_{k=1}\frac{1}{1+k} \leq 2\gamma^2_0 \left(\frac12+\int^K_{1}\frac{1}{1+s}\,ds\right) \leq 2\gamma^2_0 \log \left(K+1\right); \\
    \sum_{k=0}^{K-1} \gamma_k \sqrt{\eta_k} & = \gamma_0 \sum_{k=0}^{K-1} \left( \frac{2}{k+2} \right)^{2/3} \left( \frac{2}{k+2} \right)^{1/3 }\\
    &=   2\gamma_0 \sum_{k=0}^{K-1}\frac{1}{k+2} \leq 2\gamma_0 \log\left(K+1\right);\\
    \sum_{k=0}^{K-1} \gamma_k^2 & = 2^{4/3}\gamma_0^2 \sum_{k=1}^{K} \frac{1}{(1+k)^{4/3}}\\
    & \leq 2^{4/3}\gamma_0^2
    \left(\frac{1}{2^{4/3}}
    +\int^\infty_{1}\frac{1}{(1+s)^{4/3}}\,ds\right)
    = 7\gamma_0^2;\\
    \sum_{k=0}^{K-1} \gamma_k \eta_k & = \frac{1}{\gamma_0}\sum_{k=0}^{K-1} \gamma_k^2 \leq 7\gamma_0;\qquad\quad \sum_{k=0}^{K-1}\gamma^3_k\eta^{-1}_k = \gamma_0\sum_{k=0}^{K-1}\gamma^2_k \leq 7\gamma^3_0;\\
    \sum_{k=0}^{K-1} \gamma_k^3  & = \gamma^3_0\sum_{k=0}^{K-1} \left( \frac{2}{k+2} \right)^{2} = 4\gamma^3_0 \sum_{k=1}^{K} \frac{1}{(1+k)^2} \\
    & \leq 4\gamma^3_0\left(\frac14+\int_{1}^{\infty} \frac{1}{(1+s)^2}\,ds\right) = 3\gamma^3_0.
    \end{align*}
Therefore, denoting $\tilde{x}^K$ as a point randomly chosen from $\{x^0,x^1,\ldots,x^{K-1}\}$ with probability $\gamma_k/\Gamma_K$ for $k = 0,1, \ldots, K-1$, we obtain 
\begin{align*}
    \Exp{\norm{\nabla f(\tilde{x}^K)}}
    & \leq \frac{\Exp{V_0}}{\gamma_0 K^{1/3}}   + \frac{2 C_1 \sigma_h }{\sqrt{n}} \frac{2\gamma_0^2 \log\left(K+1\right)}{\gamma_0 K^{1/3}} + C^2_1  L_h  \frac{7\gamma^3_0}{\gamma_0 K^{1/3}} \\
    &\quad +     \frac{4\bar L_h}{\alpha^2}   \frac{3\gamma^3_0}{ \gamma_0 K^{1/3}} +  \left(\frac{L}{2} + \frac{4\sigma_h}{\alpha^2}  + \frac{2\bar L}{\alpha^2} \right) \frac{7\gamma_0^2}{ \gamma_0 K^{1/3}} \\
    &\quad + \frac{2C_1 \sigma_g}{\sqrt{n}} \frac{2\gamma_0 \log\left(K+1\right)}{ \gamma_0 K^{1/3}} + \frac{4\sigma_g }{\alpha^2}  \frac{7\gamma_0}{ \gamma_0 K^{1/3}}\\
    & = \widetilde{\cO}\biggl(
    \frac{1}{K^{1/3}}\biggl[
    \frac{V_0}{\gamma_0}
    + \left(\sigma_g + \gamma_0\sigma_h\right)
      \left(\frac{1}{\sqrt{n}} + \frac{1}{\alpha^2}\right)\\
    &\qquad\qquad
    + \gamma_0\left(L + \frac{\bar L}{\alpha^2}\right)
    + \gamma^2_0\left(L_h + \frac{\bar L_h}{\alpha^2}\right)
    \biggr]\biggr).
\end{align*}

\newpage

\section{\normEFname{MVR}}

In \normEFname{MVR}, we update the iterates $\{x^k\}$ according to: 
\begin{equation*}
    x^{k+1} = x^k - \gamma_k \frac{g^k}{\norm{g^k}}, 
\end{equation*}
where $g^k = (1/n)\sum_{i=1}^n g_i^k$ and $v^k = (1/n)\sum_{i=1}^n v_i^k$ with $g_i^k,v_i^k$ being governed by 
\begin{align}
    g_i^{k+1}
    &= g_i^k + \mathcal{C}_i^{k+1}(v_i^{k+1} - g_i^k),
       \quad \text{and} \notag\\
    v_i^{k+1}
    &= (1-\eta_k)(v_i^k
       + \nabla f_i(x^{k+1};\xi_i^{k+1})
       - \nabla f_i(x^k;\xi_i^{k+1})) \notag\\
    &\quad + \eta_k \nabla f_i(x^{k+1};\xi_i^{k+1}).
       \label{ch7:eq:mommetum_update_mvr_chapter_7}
\end{align}

\subsection{Convergence Proof}

We prove the result in the following steps. 

\subsubsection{Deriving the descent inequality} From Lemma~\ref{ch7:lem:descent_lemma_chapter_7}, we obtain \eqref{ch7:eq:descent_lemma_chapter_7}:
\begin{equation*}
    \Delta_{k+1} + \gamma_k \|\nabla f(x^k)\| \le \Delta_k  + 2 \gamma_k \norm{v^k- \nabla f(x^k)} +\frac{ 2\gamma_k}{n}\sum_{i=1}^n\norm{g_i^k-v_i^k} + \frac{\gamma_k^2L}{2}.
\end{equation*}

\subsubsection{Error bound I} 
Define $\thesischsevencV_k \eqdef (1/n)\sum_{i=1}^n \norm{g_i^k-v_i^k}$. From~Lemma~\ref{ch7:lemma:boundg_and_v_chapter_7}, 
\begin{equation*}
   \Exp{\thesischsevencV_{k+1}} \leq \sqrt{1-\alpha}    \Exp{\thesischsevencV_{k}} + \frac{\sqrt{1-\alpha}}{n} \sum_{i=1}^n \Exp{\norm{v_i^{k+1}-v_i^k}}.
\end{equation*}

To complete the upper-bound for $\Exp{\thesischsevencV_{k+1}}$, we bound $\Exp{\norm{v_i^{k+1}-v_i^k}}$. 
Since 
\begin{equation*}
    v_i^{k+1} - v_i^k = - \eta_k (v_i^k - \nabla f_i(x^k)) - \eta_k S_{i,k+1} - (1-\eta_k)S_{i,k+1} + \eta_k e_{i,k+1},
\end{equation*}
where we defined
\begin{equation*}
    S_{i,k+1} \eqdef \nabla f_i(x^{k};\xi_i^{k+1}) - \nabla f_i(x^{k+1};\xi_i^{k+1}) \quad\text{and}\quad e_{i,k+1} \eqdef \nabla f_i(x^{k};\xi_i^{k+1}) - \nabla f_i(x^k),
\end{equation*}
from the definition of the Euclidean norm and by the triangle inequality, we get
\begin{align*}
    \Exp{\norm{v_i^{k+1}-v_i^k}} & \leq \eta_k \Exp{\norm{v_i^k -\nabla f_i(x^k)}} + \Exp{ \norm{S_{i,k+1}} } + \eta_k\Exp{\norm{e_{i,k+1}}} \\
    & \leq \eta_k \Exp{\norm{v_i^k -\nabla f_i(x^k)}} + \sqrt{ \Exp{ \sqnorm{S_{i,k+1}} }  } + \eta_k\sqrt{ \Exp{\sqnorm{e_{i,k+1}}} }, 
\end{align*}
where in the last inequality we used Jensen's inequality. 
Next, by Assumption~\ref{ch7:assum:L_stoc_comp_fcn_chapter_7}, Assumption~\ref{ch7:assum:stoc_g_h_chapter_7}, and by the update rule for $x^{k+1}$,  
\begin{align*}
    \Exp{\norm{v_i^{k+1}-v_i^k}} 
    & \leq \eta_k \Exp{\norm{v_i^k - \nabla f_i(x^k)}} +  \gamma_k L_{\text{ms},i} + \eta_k \sigma_g.
\end{align*}
Plugging the above result into the upper-bound for $\Exp{\thesischsevencV_{k+1}}$, and then defining $\thesischsevencU_k \eqdef (1/n)\sum_{i=1}^n \norm{v_i^k-\nabla f_i(x^k)}$,  we obtain
\begin{equation}\label{ch7:eqn:MVR_V_t_chapter_7}
\begin{aligned}
   \Exp{\thesischsevencV_{k+1}}
   & \leq  \sqrt{1-\alpha}    \Exp{\thesischsevencV_{k}} + \sqrt{1-\alpha} \eta_k \Exp{\thesischsevencU_k} \\
   &\quad + \sqrt{1-\alpha} \gamma_k \bar L_{\text{ms}} + \sqrt{1-\alpha} \eta_k \sigma_g,
\end{aligned}
\end{equation}
where $\bar L_{\text{ms}} = (1/n)\sum_{i=1}^n L_{\text{ms},i}$.

\subsubsection{Error bound II} 
Define $\thesischsevencU_k \eqdef (1/n)\sum_{i=1}^n \norm{v_i^k - \nabla f_i(x^k)}$. Then, we have 
\begin{align*}
    v_i^{k+1} - \nabla f_i(x^{k+1}) 
    & = (1-\eta_k)(v_i^k + \nabla f_i(x^{k+1};\xi_i^{k+1}) - \nabla f_i(x^k;\xi_i^{k+1})) \\
    &\quad + \eta_k \nabla f_i(x^{k+1};\xi_i^{k+1}) - \nabla f_i(x^{k+1}) \\
    & = (1-\eta_k)(v_i^k - \nabla f_i(x^k)) - (1-\eta_k) D_{i,k+1} + \eta_k e_{i,k+1},
\end{align*}
where in the first equality we used the update rule for $v^{k+1}_i$, i.e. \eqref{ch7:eq:mommetum_update_mvr_chapter_7}, and in the second we defined
\begin{align*}
    D_{i,k+1} &\eqdef \nabla f_i(x^{k};\xi_i^{k+1}) - \nabla f_i(x^{k+1};\xi_i^{k+1}) - [\nabla f_i(x^k) - \nabla f_i(x^{k+1})] \quad\text{and}\\
    e_{i,k+1} &\eqdef \nabla f_i(x^{k+1};\xi_i^{k+1}) - \nabla f_i(x^{k+1}).
\end{align*}
Therefore, 
\begin{align*}
    \Exp{\thesischsevencU_{k+1}}
    & \leq  \frac{1-\eta_k}{n}\sum_{i=1}^n \Exp{\norm{v_i^{k} - \nabla f_i(x^{k})}} +\frac{1-\eta_k}{n}\sum_{i=1}^n \Exp{\norm{D_{i,k+1}}}\\
    &\quad  +  \frac{\eta_k}{n}\sum_{i=1}^n \Exp{\norm{e_{i,k+1}}} \\
    & \leq  \frac{1-\eta_k}{n}\sum_{i=1}^n \Exp{\norm{v_i^{k} - \nabla f_i(x^{k})}} + \frac{1-\eta_k}{n}\sum_{i=1}^n \sqrt{\Exp{\sqnorm{D_{i,k+1}}}}\\
    &\quad +  \frac{\eta_k}{n}\sum_{i=1}^n \sqrt{\Exp{\sqnorm{e_{i,k+1}}}} \\
    & \leq  \frac{1-\eta_k}{n}\sum_{i=1}^n \Exp{\norm{v_i^{k} - \nabla f_i(x^{k})}} + \frac{1-\eta_k}{n}\sum_{i=1}^n \sqrt{\Exp{\sqnorm{S_{i,k+1}}}}\\
    &\quad +  \frac{\eta_k}{n}\sum_{i=1}^n \sqrt{ \Exp{\sqnorm{e_{i,k+1}}}},
\end{align*}
where in the first inequality we used the triangle inequality, in the second we used Jensen's inequality, and in the third we used the bias-variance decomposition, more precisely 
$\Exp{\sqnorm{D_{i,k+1}}}\leq \Exp{\sqnorm{S_{i,k+1}}} $.
Next, by Assumption~\ref{ch7:assum:L_stoc_comp_fcn_chapter_7} and Assumption~\ref{ch7:assum:stoc_g_h_chapter_7}, 
\begin{align}
    \Exp{\thesischsevencU_{k+1}}
    & \leq  (1-\eta_k) \Exp{\thesischsevencU_k} + (1-\eta_k) \frac{1}{n}\sum_{i=1}^n L_{\text{ms},i}\sqrt{  \Exp{\sqnorm{ x^k - x^{k+1} } } }   + \eta_k \sigma_g \notag \\
    & \leq (1-\eta_k) \Exp{\thesischsevencU_k} + (1-\eta_k)\gamma_k  \bar L_{\text{ms}}  + \eta_k \sigma_g,  \label{ch7:eqn:MVR_U_t_chapter_7}
\end{align}
where  in the last inequality we used the update rule for $x^{k+1}$ and denoted $\bar L_{\text{ms}} = (1/n)\sum_{i=1}^n L_{\text{ms},i}$.

\subsubsection{Error bound III} 
From the definition of $v^k = (1/n)\sum_{i=1}^n v_i^k$,  we denote $\hat{e}_{k+1}$ as follows 
\begin{align*}
    \hat{e}_{k+1} & \eqdef v^{k+1} - \nabla f(x^{k+1}) \\
    & =  \frac{1}{n}\sum_{i=1}^n (v_i^{k+1} - \nabla f_i(x^{k+1}))\\
    & = (1-\eta_k)(v^k - \nabla f(x^k)) - (1-\eta_k) \frac{1}{n}\sum_{i=1}^n D_{i,k+1}  + \eta_k \frac{1}{n}\sum_{i=1}^n e_{i,k+1},
\end{align*}
where $D_{i,k+1}$ and $e_{i,k+1}$ are defined as in Error bound II.

Next, by applying Lemma~\ref{ch7:lemma:summing_the_recursion_chapter_7}, 
\begin{align*}
    \hat{e}_{k+1} & = \prod_{\tau=0}^k (1-\eta_\tau) \hat{e}_{0} - \sum_{j=0}^k \left(\prod_{\tau=j}^k (1-\eta_\tau)\right) {D}_{j+1}  + \sum_{j=0}^k \left(\prod_{\tau=j+1}^k (1-\eta_\tau)\right) \eta_j e_{j+1},
\end{align*}
where $D_{k+1} = (1/n)\sum_{i=1}^n D_{i,k+1}$ and $e_{k+1} = (1/n)\sum_{i=1}^n e_{i,k+1}$. If $\eta_0 = 1$, then we have 
\begin{align*}
    \hat{e}_{k+1} & = - \sum_{j=0}^k \left(\prod_{\tau=j}^k (1-\eta_\tau)\right) {D}_{j+1} + \sum_{j=0}^k \left(\prod_{\tau=j+1}^k (1-\eta_\tau)\right) \eta_j e_{j+1}. 
\end{align*}
Therefore, 
\begin{align*}
    \MoveEqLeft \Exp{\norm{\hat{e}_{k+1}}}\\
   & \leq \Exp{ \norm{ \sum_{j=0}^k \left(\prod_{\tau=j}^k (1-\eta_\tau)\right) {D}_{j+1}  } }\\
   &\quad  + \Exp{ \norm{ \sum_{j=0}^k \left(\prod_{\tau=j+1}^k (1-\eta_\tau)\right) \eta_j e_{j+1} } } \\
    & \leq \sqrt{ \Exp{ \sqnorm{ \sum_{j=0}^k \left(\prod_{\tau=j}^k (1-\eta_\tau)\right) {D}_{j+1}  } }  } \\
    &\quad+ \sqrt{ \Exp{ \sqnorm{ \sum_{j=0}^k \left(\prod_{\tau=j+1}^k (1-\eta_\tau)\right) \eta_j e_{j+1} } } },
\end{align*} 
where in the last inequality we used Jensen's inequality.

Next, from Assumption~\ref{ch7:assum:stoc_g_h_chapter_7}, we can prove that $\Exp{D_{k+1}} = 0$, $\Exp{ \inner{D_l}{D_{i}}} = 0$ for $l \neq i$, $\Exp{e_{k+1}}=0$, $\Exp{\inner{e_l}{e_i}}=0$ for $l \neq i$.  Hence, since $\prod_{\tau=j}^k (1-\eta_\tau) \leq \prod_{\tau=j+1}^k (1-\eta_\tau)$ and the $D_{i,j+1}$ are independent across $i$ with zero mean,
\begin{align*}
    \MoveEqLeft \Exp{\norm{\hat{e}_{k+1}}}\\
    & \leq \sqrt{ \sum_{j=0}^k \left(\prod_{\tau=j+1}^k (1-\eta_\tau)^2\right)  \Exp{ \sqnorm{  {D}_{j+1}  } }  } \\
    &\quad + \sqrt{ \sum_{j=0}^k \left(\prod_{\tau=j+1}^k (1-\eta_\tau)^2\right) \eta_j^2 \Exp{ \sqnorm{  e_{j+1} } } } \\
    & \leq \sqrt{ \sum_{j=0}^k \left(\prod_{\tau=j+1}^k (1-\eta_\tau)^2\right)  \frac{1}{n^2}\sum_{i=1}^n\Exp{ \sqnorm{  \nabla f_i(x^{j};\xi_i^{j+1}) - \nabla f_i(x^{j+1};\xi_i^{j+1})  } }  } \\
    &\quad + \sqrt{ \sum_{j=0}^k \left(\prod_{\tau=j+1}^k (1-\eta_\tau)^2\right) \eta_j^2 \Exp{ \sqnorm{  e_{j+1} } } }.
\end{align*}

Next, by Assumption~\ref{ch7:assum:L_stoc_comp_fcn_chapter_7} and Assumption~\ref{ch7:assum:stoc_g_h_chapter_7}, and since $\Exp{\sqnorm{e_{j+1}}} \leq \sigma_g^2/n$, 
\begin{align*}
    \MoveEqLeft \Exp{\norm{\hat{e}_{k+1}}}\\
  & \leq  \sqrt{\frac{\bar L_{\text{ms}}^2}{n}}\sqrt{ \sum_{j=0}^k \left(\prod_{\tau=j+1}^k (1-\eta_\tau)^2\right) \cdot  \gamma_j^2  } + \frac{\sigma_g}{\sqrt{n}}\sqrt{ \sum_{j=0}^k \left(\prod_{\tau=j+1}^k (1-\eta_\tau)^2\right) \eta_j^2 }\\
  & \leq  \bar L_{\text{ms}}\sqrt{ \sum_{j=0}^k \left(\prod_{\tau=j+1}^k (1-\eta_\tau)\right) \cdot  \gamma_j^2  } + \frac{\sigma_g}{\sqrt{n}}\sqrt{ \sum_{j=0}^k \left(\prod_{\tau=j+1}^k (1-\eta_\tau)\right) \eta_j^2 },
\end{align*} 
where $\bar L_{\text{ms}}^2 = (1/n)\sum_{i=1}^n L_{\text{ms},i}^2$, and in the last inequality we also used $\sqrt{\bar L_{\text{ms}}^2/n} \leq \bar L_{\text{ms}}$, which holds since $(1/n^2)\sum_{i=1}^n L_{\text{ms},i}^2 \leq \bigl((1/n)\sum_{i=1}^n L_{\text{ms},i}\bigr)^2$.

If $\eta_k = (2/(k+2))^{2/3}$ and $\gamma_k = \gamma_0 (2/(k+2))^{2/3}$, then we can prove  that $\eta_0=1$, and from~Lemma~\ref{ch7:lemma:sum_prod_chapter_7} that 
\begin{equation*}\label{ch7:eqn:MVR_whole_v_err_chapter_7}
    \Exp{\norm{v^{k+1} - \nabla f(x^{k+1}) }} 
  \leq \bar L_{\text{ms}}\sqrt{ C(4/3,2/3) \frac{\gamma_{k+1}^2}{\eta_{k+1}} } + \frac{\sigma_g}{\sqrt{n}} \sqrt{ C(4/3,2/3) \eta_{k+1} }.
\end{equation*}

\subsubsection{Bounding Lyapunov function}
Define our Lyapunov function $V_k = \Delta_{k} + C_{1,k} \thesischsevencV_k + C_{2,k} \thesischsevencU_k$ with
\begin{equation*}
    C_{1,k} = \frac{2\gamma_k}{1-\sqrt{1-\alpha}} \quad\text{and}\quad C_{2,k} = \frac{2\gamma_k \sqrt{1-\alpha}}{1-\sqrt{1-\alpha}}.
\end{equation*}
Then, 
\begin{align*}
    \Exp{V_{k+1}} 
    & = \Exp{\Delta_{k+1} + C_{1,k+1} \thesischsevencV_{k+1} + C_{2,k+1} \thesischsevencU_{k+1}} \\
    & \leq \Exp{\Delta_k} - \gamma_k \Exp{\norm{\nabla f(x^k)}} + 2 \gamma_k \Exp{\norm{v^k- \nabla f(x^k)}} + 2\gamma_k \Exp{\thesischsevencV_k} + \frac{\gamma_k^2L}{2} \\
    &\quad + C_{1,k+1} \Exp{\thesischsevencV_{k+1}} + C_{2,k+1} \Exp{\thesischsevencU_{k+1}},
\end{align*}
where we use \eqref{ch7:eq:descent_lemma_chapter_7}.

Next, by the upper-bounds for $\Exp{\thesischsevencV_{k+1}}$ and  $\Exp{\thesischsevencU_{k+1}}$, 
\begin{align*}
    \MoveEqLeft \Exp{ V_{k+1} }\\
   & \leq \Exp{\Delta_k} - \gamma_k \Exp{\norm{\nabla f(x^k)}} + 2 \gamma_k \Exp{\norm{v^k- \nabla f(x^k)}} + 2\gamma_k \Exp{\thesischsevencV_k} + \frac{\gamma_k^2L}{2} \\
    &\quad + C_{1,k+1} (\sqrt{1-\alpha}    \Exp{\thesischsevencV_{k}} + \sqrt{1-\alpha} \eta_k \Exp{\thesischsevencU_k}\\
    &\qquad\qquad + \sqrt{1-\alpha} \gamma_k \bar L_{\text{ms}} + \sqrt{1-\alpha} \eta_k \sigma_g) \\
    &\quad + C_{2,k+1} ((1-\eta_k) \Exp{\thesischsevencU_k} + (1-\eta_k)\gamma_k  \bar L_{\text{ms}}  + \eta_k \sigma_g ).
\end{align*}

Next,  since $\gamma_{k+1} \leq \gamma_k$, we can prove that $C_{1,k+1} \leq C_{1,k}$ and that $C_{2,k+1} \leq C_{2,k}$, and also that 
\begin{align*}
    2\gamma_k + C_{1,k+1}\sqrt{1-\alpha}  & \leq 2\gamma_k + C_{1,k}\sqrt{1-\alpha} = C_{1,k} \\
    C_{1,k+1} \sqrt{1-\alpha}\eta_k + C_{2,k+1} (1-\eta_k) & \leq C_{1,k} \sqrt{1-\alpha}\eta_k + C_{2,k} (1-\eta_k)= C_{2,k}
\end{align*}
Therefore, 
\begin{align*}
   \Exp{ V_{k+1} } 
   & \leq \Exp{V_k} - \gamma_k \Exp{\norm{\nabla f(x^k)}} + 2 \gamma_k \Exp{\norm{v^k- \nabla f(x^k)}}  + \frac{\gamma_k^2L}{2} \\
    &\quad + C_{1,k} (\sqrt{1-\alpha} \gamma_k \bar L_{\text{ms}} + \sqrt{1-\alpha} \eta_k \sigma_g) \\
    &\quad + C_{2,k} \left( (1-\eta_k)\gamma_k  \bar L_{\text{ms}}  + \eta_k \sigma_g  \right) \\
     & \leq \Exp{V_k} - \gamma_k \Exp{\norm{\nabla f(x^k)}} + 2 \gamma_k \Exp{\norm{v^k- \nabla f(x^k)}}   \\
    &\quad + \gamma_k^2  \left( \frac{4\sqrt{1-\alpha}}{1-\sqrt{1-\alpha}} \bar L_{\text{ms}}  + \frac{L}{2}\right) + \eta_k\gamma_k \frac{4\sqrt{1-\alpha}}{1-\sqrt{1-\alpha}} \sigma_g.
\end{align*}

Next, from the upper-bound of $\Exp{\norm{v^k- \nabla f(x^k)}}$, 
\begin{align*}
   \Exp{ V_{k+1} }
   & \leq \Exp{V_k} - \gamma_k \Exp{\norm{\nabla f(x^k)}} \\
   &\quad + 2 \gamma_k \left( \bar L_{\text{ms}}\sqrt{ C(4/3,2/3) \frac{\gamma_{k}^2}{\eta_{k}} } + \frac{\sigma_g}{\sqrt{n}} \sqrt{ C(4/3,2/3) \eta_k }  \right)  \\
   &\quad + \gamma_k^2  \left( \frac{4\sqrt{1-\alpha}}{1-\sqrt{1-\alpha}} \bar L_{\text{ms}}  + \frac{L}{2}\right) + \eta_k\gamma_k \frac{4\sqrt{1-\alpha}}{1-\sqrt{1-\alpha}} \sigma_g \\
   & \leq \Exp{V_k} - \gamma_k \Exp{\norm{\nabla f(x^k)}} \\
   &\quad + \frac{\gamma_k^2}{\sqrt{\eta_k}} \cdot  2 \bar L_{\text{ms}} \sqrt{C(4/3,2/3) } +  \gamma_k\sqrt{\eta_k} \cdot \frac{2\sigma_g \sqrt{C(4/3, 2/3)}}{\sqrt{n}}  \\
   &\quad + \gamma_k^2  \left( \frac{4\sqrt{1-\alpha}}{1-\sqrt{1-\alpha}} \bar L_{\text{ms}}  + \frac{L}{2}\right) + \eta_k\gamma_k \frac{4\sqrt{1-\alpha}}{1-\sqrt{1-\alpha}} \sigma_g\\
   &\leq \Exp{V_k} - \gamma_k \Exp{\norm{\nabla f(x^k)}} \\
   &\quad +  2 C_1  \bar L_{\text{ms}}\frac{\gamma_k^2}{\sqrt{\eta_k}} +  \frac{2 C_1 \sigma_g }{\sqrt{n}}\gamma_k\sqrt{\eta_k} \\
   &\quad + \gamma_k^2  \left( \frac{L}{2} + \frac{4 {\bar L}_{\text{ms}} }{\alpha^2}   \right) + \eta_k\gamma_k \frac{4\sigma_g}{\alpha^2},
\end{align*} 
where in the last inequality we used $\sqrt{1-\alpha}/(1-\sqrt{1-\alpha}) \leq 1/\alpha^2$ and denoted $C_1 = \sqrt{ C(4/3,2/3) }$.

\subsubsection{Deriving the convergence rate}
Let $\Gamma_K \eqdef \sum_{k=0}^{K-1} \gamma_k$. By re-arranging the terms and by the telescopic series, 
\begin{align*}
    \MoveEqLeft \frac{1}{\Gamma_K} \sum_{k=0}^{K-1} \gamma_k\Exp{\norm{\nabla f(x^k)}}\\
    & \leq \frac{1}{\Gamma_K} \biggl[ \Exp{V_0} -\Exp{V_K} + \frac{2C_1\sigma_g }{\sqrt{n}} \sum_{k=0}^{K-1} \gamma_k \sqrt{\eta_k} +     \frac{4 \sigma_g}{\alpha^2}  \sum_{k=0}^{K-1} \gamma_k \eta_k\\
    &\qquad + 2 C_1 {\bar L_{\text{ms}}} \sum_{k=0}^{K-1} \gamma_k^2 \eta_k^{-1/2}  +  \left(\frac{L}{2} + \frac{4 {\bar L}_{\text{ms}}}{\alpha^2} \right) \sum_{k=0}^{K-1} \gamma_k^2 \biggr].
\end{align*}

Since $\eta_k = (2/(k+2))^{2/3}$ and $\gamma_k = \gamma_0 (2/(k+2))^{2/3}$, and  for $K \geq 1$ we have 
\begin{align*}
    \Gamma_K = \sum_{k=0}^{K-1} \gamma_k & \geq  K \gamma_{K-1} = \gamma_0K  \left(\frac{2}{K+1}\right)^{2/3} \geq \gamma_0 K^{1/3};\\
    \sum_{k=0}^{K-1} \gamma_k^2 \eta_k^{-1/2}  & = \gamma_0^2\sum_{k=0}^{K-1}  \left( \frac{2}{k+2} \right)^{4/3} \left( \frac{2}{k+2} \right)^{-1/3 } = 2\gamma_0^2 \sum_{k=0}^{K-1} \frac{1}{k+2}\\
    &= 2\gamma_0^2 \sum^{K}_{k=1}\frac{1}{1+k}
    \leq 2\gamma^2_0 \left(\frac12+\int^K_{1}\frac{1}{1+s}\,ds\right) \leq 2\gamma^2_0 \log \left(K+1\right); \\
    \sum_{k=0}^{K-1} \gamma_k \sqrt{\eta_k} & = \gamma_0 \sum_{k=0}^{K-1} \left( \frac{2}{k+2} \right)^{2/3} \left( \frac{2}{k+2} \right)^{1/3 }\\
    &=   2\gamma_0 \sum_{k=0}^{K-1}\frac{1}{k+2} \leq 2\gamma_0 \log\left(K+1\right);\\
    \sum_{k=0}^{K-1} \gamma_k^2 & = 2^{4/3}\gamma_0^2 \sum_{k=1}^{K} \frac{1}{(1+k)^{4/3}}\\
    & \leq 2^{4/3}\gamma_0^2
    \left(\frac{1}{2^{4/3}}
    +\int^\infty_{1}\frac{1}{(1+s)^{4/3}}\,ds\right)
    = 7\gamma_0^2;\\
    \sum_{k=0}^{K-1} \gamma_k \eta_k & = \frac{1}{\gamma_0}\sum_{k=0}^{K-1} \gamma_k^2 \leq 7\gamma_0.
\end{align*}

Therefore, denoting $\tilde{x}^K$ as a point randomly chosen from $\{x^0,x^1,\ldots,x^{K-1}\}$ with probability $\gamma_k/\Gamma_K$ for $k = 0,1, \ldots, K-1$, we obtain
\begin{align*}
    \Exp{\norm{\nabla f(\tilde{x}^K)}}
    & \leq \frac{\Exp{V_0}}{\gamma_0 K^{1/3}} + \frac{2C_1\sigma_g }{\sqrt{n}} \frac{2\gamma_0 \log\left(K+1\right)}{ \gamma_0 K^{1/3}} +     \frac{4 \sigma_g}{\alpha^2}  \frac{7\gamma_0}{\gamma_0 K^{1/3}}\\
    &\quad + 2 C_1 {\bar L_{\text{ms}}} \frac{2\gamma_0^2 \log\left(K+1\right)}{\gamma_0 K^{1/3}}  +  \left(\frac{L}{2} + \frac{4 {\bar L}_{\text{ms}}}{\alpha^2} \right) \frac{7\gamma_0^2}{\gamma_0 K^{1/3}} \\
    & = \widetilde{\cO}\left(\frac{1}{K^{1/3}}\left(\frac{V_0}{\gamma_0} + \sigma_g\left(\frac{1}{\sqrt{n}} +\frac{1}{\alpha^2}\right) +\gamma_0\left(L + {\bar L}_{\text{ms}} + \frac{{\bar L}_{\text{ms}}}{\alpha^2}\right)\right)\right).
\end{align*}

\newpage
\section{Deep Learning experiments}
\subsection{Hardware and Datasets.}
To evaluate the performance of the proposed methods in training Deep Neural Networks (DNNs), we utilized the ResNet-18 architecture \parencite{he2016deep}. ResNet-18 is a prominent model for image classification, and its architecture is also frequently adapted for tasks such as feature extraction in image segmentation, object detection, image embedding, and image captioning. Our experiments involved training all layers of the ResNet-18 model, corresponding to an optimization problem with $d = 11,173,962$ parameters.

All implementations were developed in PyTorch \parencite{paszke2019pytorch}, and experiments were conducted on the CIFAR-10 dataset \parencite{krizhevsky2009learning}. Numerical evaluations were performed on a server-grade machine running Ubuntu 18.04 (Linux Kernel v5.4.0). This system was equipped with dual 16-core 3.3 GHz Intel Xeon processors (totaling 32 cores) and four NVIDIA A100 GPUs, each with 40GB of memory.

\subsection{Additional implementation details}
\label{ch7:sec:restored_implementation}
The data partition and Top-$K$ compression setting are described in
Section~\ref{ch7:sec:exp_chapter_7}. The following details supplement the
implementation summary there. All methods use PyTorch automatic differentiation;
Hessian-vector products are computed without materializing the Hessian matrix.

\paragraph{\normEFname{HM}.}
The implementation first computes the stochastic gradient at the new point
while retaining its computational graph with \texttt{create\_graph=True}.
A second call to \texttt{torch.autograd.grad} uses the parameter gradients as
outputs and the displacement vector as \texttt{grad\_outputs}, producing the
Hessian-vector product. The gradient and Hessian-vector product use the same
minibatch. The momentum buffer is then updated with this Hessian correction,
as described in Chapter~\ref{ch:chapter_7}. This requires two automatic
differentiation passes per minibatch, but their measured costs need not be equal.

\paragraph{\normEFname{RHM}.}
The master draws a scalar uniformly from $(0,1)$ and shares the corresponding
interpolation point across workers. Consistently with the RHM convention in
Chapter~\ref{ch:chapter_7}, independent minibatches are used for the stochastic
gradient at the new point and the Hessian-vector product at the interpolation
point. The implementation first computes the new-point gradient, temporarily
swaps the parameters to the interpolation point inside a
\texttt{torch.no\_grad()} context, and performs a forward and backward pass on
the independent minibatch to construct the first-order graph there. An
\texttt{autograd.grad} call then computes the Hessian-vector product with the
displacement vector. Finally, the original parameters are restored and the
randomized Hessian-corrected momentum update is applied. The procedure uses
approximately three backpropagation-equivalents per minibatch.

\paragraph{\normEFname{MVR}.}
The variance-reduced update uses gradients at the previous and new points
computed on the same minibatch. The implementation caches the previous
parameters in the optimizer state, temporarily swaps the live parameters to
that point under \texttt{torch.no\_grad()}, and calls \texttt{autograd.grad}
without populating the model's \texttt{.grad} attributes. It then restores the
new parameters, computes the new-point gradient with
\texttt{loss.backward()}, and passes both gradients to the optimizer step.
This requires two first-order differentiation passes and does not construct
a second-order graph.

\paragraph{\normEFname{IGT}.}
The optimizer caches the previous weights and uses their displacement from
the current weights to form the extrapolated point specified by the IGT
update. It temporarily moves the parameters to this point inside a
\texttt{no\_grad} context, performs one backward pass, restores the original
parameters, and applies the optimizer step using the extrapolated gradient.
Only one backpropagation is required per minibatch.

\paragraph{Other methods.}
The implementations of \algname{EF21-SGD}, \algname{EF21-SGDM},
\normEFname{SGDM}, and \algname{EControl} each use a standard first-order
backward pass followed by an optimizer step. Their differences lie in the
optimizer's state updates and do not require additional backpropagations.

\subsection{Hyperparameter tuning}
\label{ch7:sec:restored_tuning}
The proposed momentum variants were compared with \algname{EF21-SGD}
\parencite{fatkhullin2025ef21}, \algname{EF21-SGDM}
\parencite{fatkhullin2023momentum}, \normEFname{SGDM}
\parencite{khirirat2024error}, and \algname{EControl}
\parencite{gao2024econtrol}. Table~\ref{ch7:tab:restored_search_space}
summarizes the reported search space for each method.

\paragraph{Momentum parameters.}
For \algname{EF21-SGDM} and \algname{EControl}, the momentum parameter was
fixed at $\eta=0.1$, following the recommendations in their original
publications. For \normEFname{HM} and \normEFname{SGDM}, both constant
values in $\{0.01,0.1,0.2\}$ and decreasing schedules were explored.
Following the stronger performance observed with decreasing schedules for
these variants, the search for \normEFname{RHM}, \normEFname{MVR}, and
\normEFname{IGT} focused on decreasing schedules. As a practical adaptation
to prevent the momentum parameter from diminishing too rapidly early in
training, $\eta_e$ was updated once per epoch rather than once per iteration.

\paragraph{Stepsizes.}
For the non-normalized baselines, constant stepsizes were selected from
$\{1.0,0.1,0.05,0.01,0.005\}$. For the normalized methods, both constant
stepsizes in $\{1.0,0.1,0.05,0.01\}$ and decreasing schedules
$\gamma_e=\gamma_0(e+1)^{-p}$ were evaluated, with
$\gamma_0\in\{2.0,1.0,0.5,0.1\}$. The exponents are listed in
Table~\ref{ch7:tab:restored_search_space}. These are experimental tuning
choices: the epoch-wise schedules and selected constant stepsizes should
not be identified with the iteration-wise schedules in the convergence
theorems.

\begin{table}[htbp]
\centering
\caption{Hyperparameter search spaces. Here $e$ is the epoch index,
$C_b=\{1.0,0.1,0.05,0.01,0.005\}$,
$C_n=\{1.0,0.1,0.05,0.01\}$, and every decreasing stepsize uses
$\gamma_0\in\{2.0,1.0,0.5,0.1\}$.}
\label{ch7:tab:restored_search_space}
\begingroup
\small
\setlength{\tabcolsep}{4pt}
\renewcommand{\arraystretch}{1.18}
\begin{tabular}{@{}>{\raggedright\arraybackslash\leavevmode}p{0.26\textwidth}>{\raggedright\arraybackslash}p{0.34\textwidth}>{\raggedright\arraybackslash}p{0.34\textwidth}@{}}
\toprule
\textbf{Method} & \textbf{Stepsize} & \textbf{Momentum} \\
\midrule
\algname{EF21-SGDM} & Constant: $\gamma\in C_b$ & Constant: $\eta=0.1$ \\
\algname{EControl} & Constant: $\gamma\in C_b$ & Constant: $\eta=0.1$ \\
\algname{EF21-SGD} & Constant: $\gamma\in C_b$ & Not applicable \\
\midrule
\normEFname{MVR} & Constant: $\gamma\in C_n$; or $\gamma_e=\gamma_0(e+1)^{-2/3}$ & $\eta_e=(2/(e+2))^{2/3}$ \\
\normEFname{RHM} & Constant: $\gamma\in C_n$; or $\gamma_e=\gamma_0(e+1)^{-2/3}$ & $\eta_e=(2/(e+2))^{2/3}$ \\
\normEFname{SGDM} & Constant: $\gamma\in C_n$; or $\gamma_e=\gamma_0(e+1)^{-0.75}$ & Constant: $\eta\in\{0.01,0.1,0.2\}$; or $\eta_e=(2/(e+2))^{0.5}$ \\
\normEFname{IGT} & Constant: $\gamma\in C_n$; or $\gamma_e=\gamma_0(e+1)^{-5/7}$ & $\eta_e=(2/(e+2))^{4/7}$ \\
\normEFname{HM} & Constant: $\gamma\in C_n$; or $\gamma_e=\gamma_0(e+1)^{-2/3}$ & Constant: $\eta\in\{0.01,0.1,0.2\}$; or $\eta_e=(2/(e+2))^{2/3}$ \\
\bottomrule
\end{tabular}
\endgroup
\end{table}

All methods were trained for $90$ epochs. In this fixed experimental setup,
the per-iteration communication cost is identical across methods, so the
epoch count serves as a proxy for the total bits communicated. The selected
hyperparameters were chosen using the best validation accuracy achieved
and observed stable convergence behavior. Table~\ref{ch7:tab:restored_selected_parameters}
reports these choices; the associated accuracy results are already given
in Table~\ref{ch7:main_table:accuracy_results_chapter_7}.

\begin{table}[htbp]
\centering
\caption{Selected hyperparameters reported for the ResNet-18/CIFAR-10 experiments.}
\label{ch7:tab:restored_selected_parameters}
\begin{threeparttable}
\small
\setlength{\tabcolsep}{4pt}
\renewcommand{\arraystretch}{1.15}
\begin{tabular}{@{}>{\raggedright\arraybackslash\leavevmode}p{0.27\textwidth}>{\raggedright\arraybackslash}p{0.29\textwidth}>{\raggedright\arraybackslash}p{0.38\textwidth}@{}}
\toprule
\textbf{Method} & \textbf{Stepsize} & \textbf{Momentum} \\
\midrule
\algname{EF21-SGDM} & Constant $\gamma=0.1$ & Constant $\eta=0.1$ \\
\algname{EControl} & Constant $\gamma=1.0$ & Constant $\eta=0.1$ \\
\algname{EF21-SGD} & Constant $\gamma=1.0$ & --- \\
\midrule
\normEFname{MVR} & Constant $\gamma=0.1$ & $\eta_e=(2/(e+2))^{2/3}$ \\
\normEFname{RHM} & Constant $\gamma=0.1$ & $\eta_e=(2/(e+2))^{2/3}$ \\
\normEFname{SGDM} & Constant $\gamma=0.1$ & $\eta_e=(2/(e+2))^{0.5}$ \\
\normEFname{IGT}\tnote{a} & Constant $\gamma=0.1$ & $\eta_e=(2/(e+2))^{0.57}$ \\
\normEFname{HM} & Constant $\gamma=0.1$ & $\eta_e=(2/(e+2))^{2/3}$ \\
\bottomrule
\end{tabular}
\begin{tablenotes}[flushleft]
\footnotesize
\item All decreasing momentum schedules have $\eta_0=1.0$.
\item[a] The released runs use the exponent $0.57$, i.e., $4/7\approx0.571$ rounded to two
decimals; Table~\ref{ch7:tab:restored_search_space} lists it as $4/7$.
\end{tablenotes}
\end{threeparttable}
\end{table}
\FloatBarrier

\subsection{GPU and wall-clock timing definitions}
\label{ch7:sec:restored_timing}
The timing columns in Table~\ref{ch7:main_table:accuracy_results_chapter_7}
use the following definitions. Let $e^\star$ denote the first epoch attaining
the maximal validation accuracy. For each epoch $e\in\{1,\ldots,e^\star\}$,
CUDA events $\mathrm{start}_e$ and $\mathrm{end}_e$ are created using
\texttt{torch.cuda.Event(}\allowbreak\texttt{enable\_timing=True)}. The start event is recorded
immediately before the first training minibatch of epoch $e$; the end event
is recorded after all GPU work for that epoch, including training,
validation, and test evaluation. After calling
\texttt{torch.cuda.synchronize()}, the per-epoch device time is
\[
\tau_e=\frac{\texttt{start\_e.elapsed\_time(end\_e)}}{1000}
\quad\text{seconds}.
\]
The reported \emph{GPU Time to Best} is
\[
T^{\mathrm{GPU}}_{\leq e^\star}=\sum_{e=1}^{e^\star}\tau_e,
\]
which includes all measured device time up to and including the best epoch.

\emph{Wall Time to Best} is the corresponding real elapsed time measured
with \texttt{time.time()}, reset immediately after the epoch-0 snapshot.
At $e^\star$, it equals the cumulative \texttt{wall\_seconds} logged by the
code. Runs execute one at a time on a single GPU, and the CUDA events bracket
all GPU work with explicit synchronization. Consequently, wall time closely
matches GPU time in this setup, with comparatively small CPU-only overheads.
These are timings for the simulated distributed experiment, not measurements
of communication latency on a physical multi-device network.

\subsection{Wall-clock performance plots}
\label{ch7:sec:restored_wall_clock}
Figure~\ref{ch7:main_fig:performance_plots_epoch_chapter_7} already presents
the epoch-based comparison. Figures~\ref{ch7:fig:restored_time_full}
and~\ref{ch7:fig:restored_time_truncated} provide the complementary
wall-clock views from the paper. The full-duration plots show the time
required by each method to complete the fixed $90$-epoch budget. The
truncated plots restrict the time axis to the completion time of the
fastest methods, exposing the trade-off between per-epoch accuracy and
computational cost. In particular, the IGT trajectories remain competitive
with HM over this shorter time horizon, whereas the Hessian-based methods
require substantially more time to complete all epochs.

\begin{figure}[p]
\centering
\begin{subfigure}{0.48\textwidth}
\includegraphics[width=\linewidth]{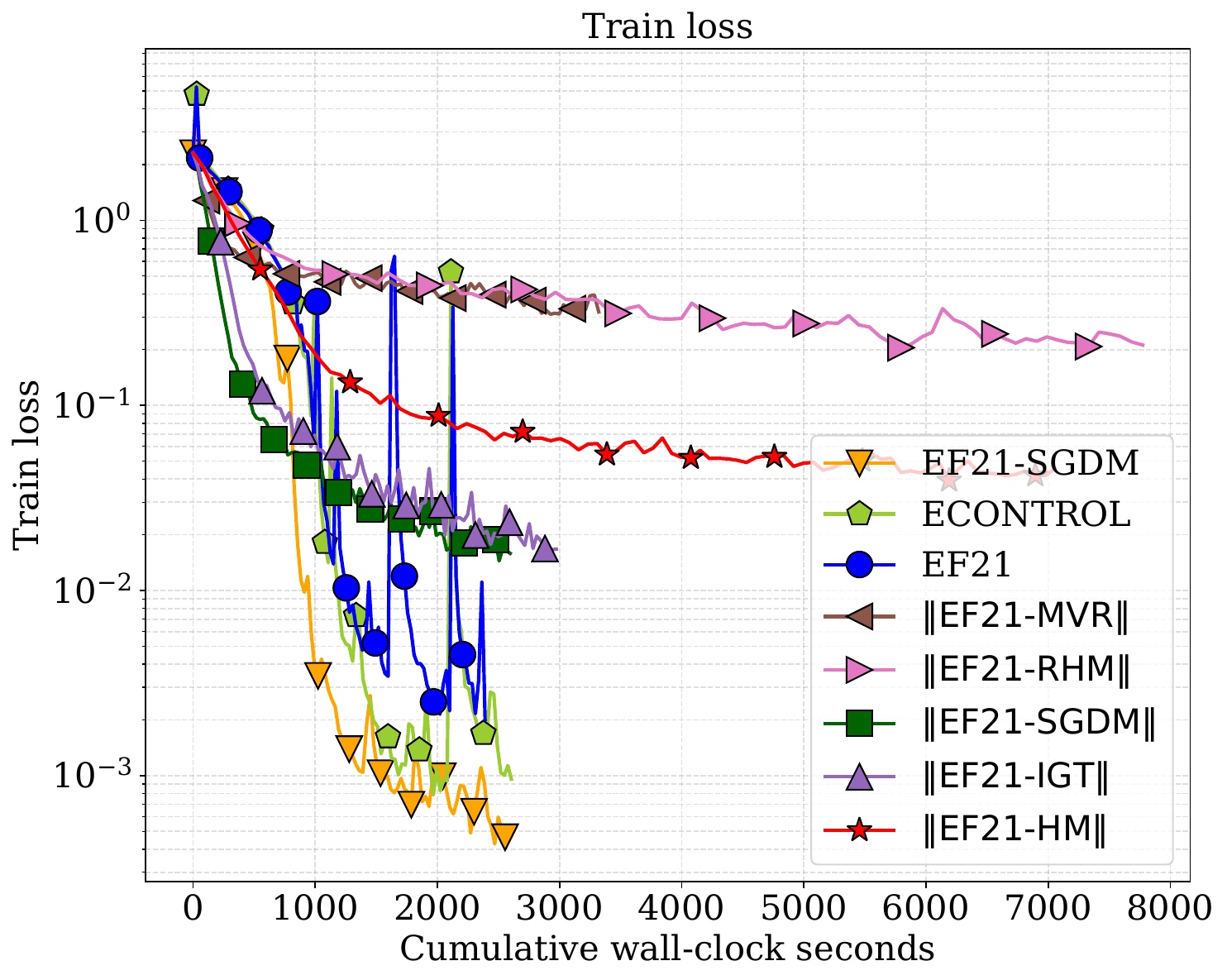}
\caption{Training loss}
\end{subfigure}\hfill
\begin{subfigure}{0.48\textwidth}
\includegraphics[width=\linewidth]{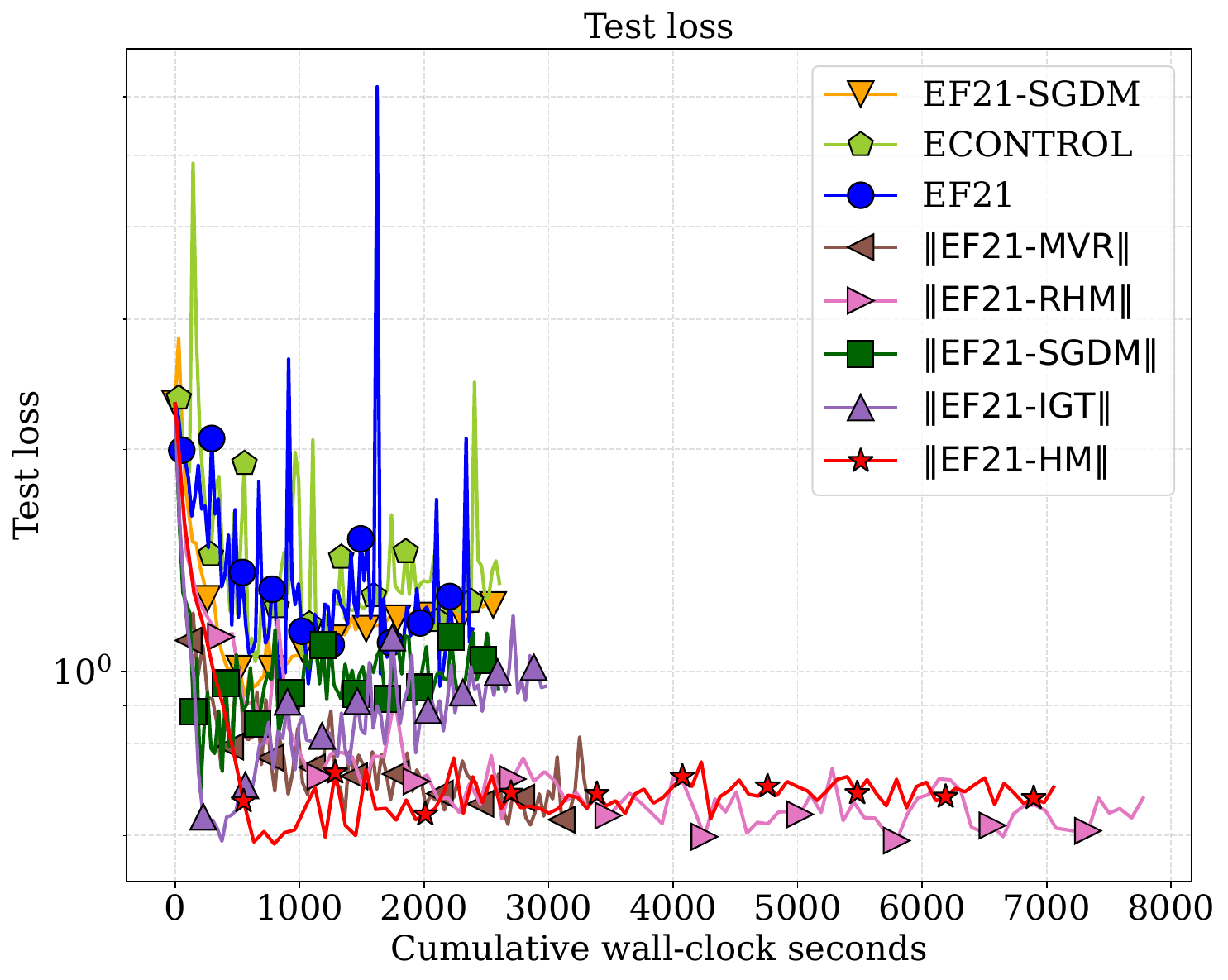}
\caption{Test loss}
\end{subfigure}
\par\medskip
\begin{subfigure}{0.48\textwidth}
\includegraphics[width=\linewidth]{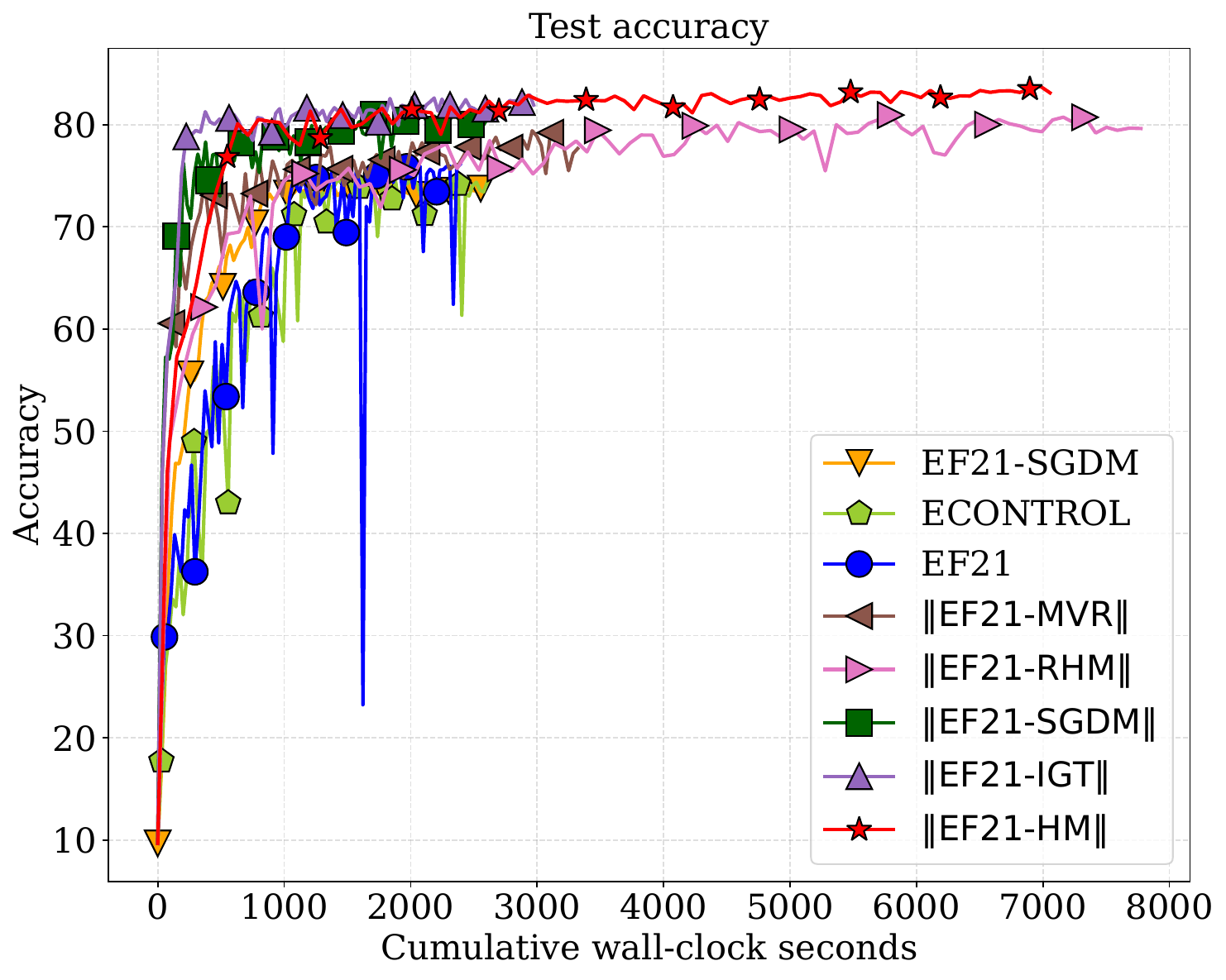}
\caption{Test accuracy}
\end{subfigure}
\caption{ResNet-18/CIFAR-10 performance as a function of cumulative wall-clock
seconds over the full training duration. Methods with higher per-epoch costs
take longer to complete the $90$-epoch training schedule.}
\label{ch7:fig:restored_time_full}
\end{figure}

\begin{figure}[p]
\centering
\begin{subfigure}{0.48\textwidth}
\includegraphics[width=\linewidth]{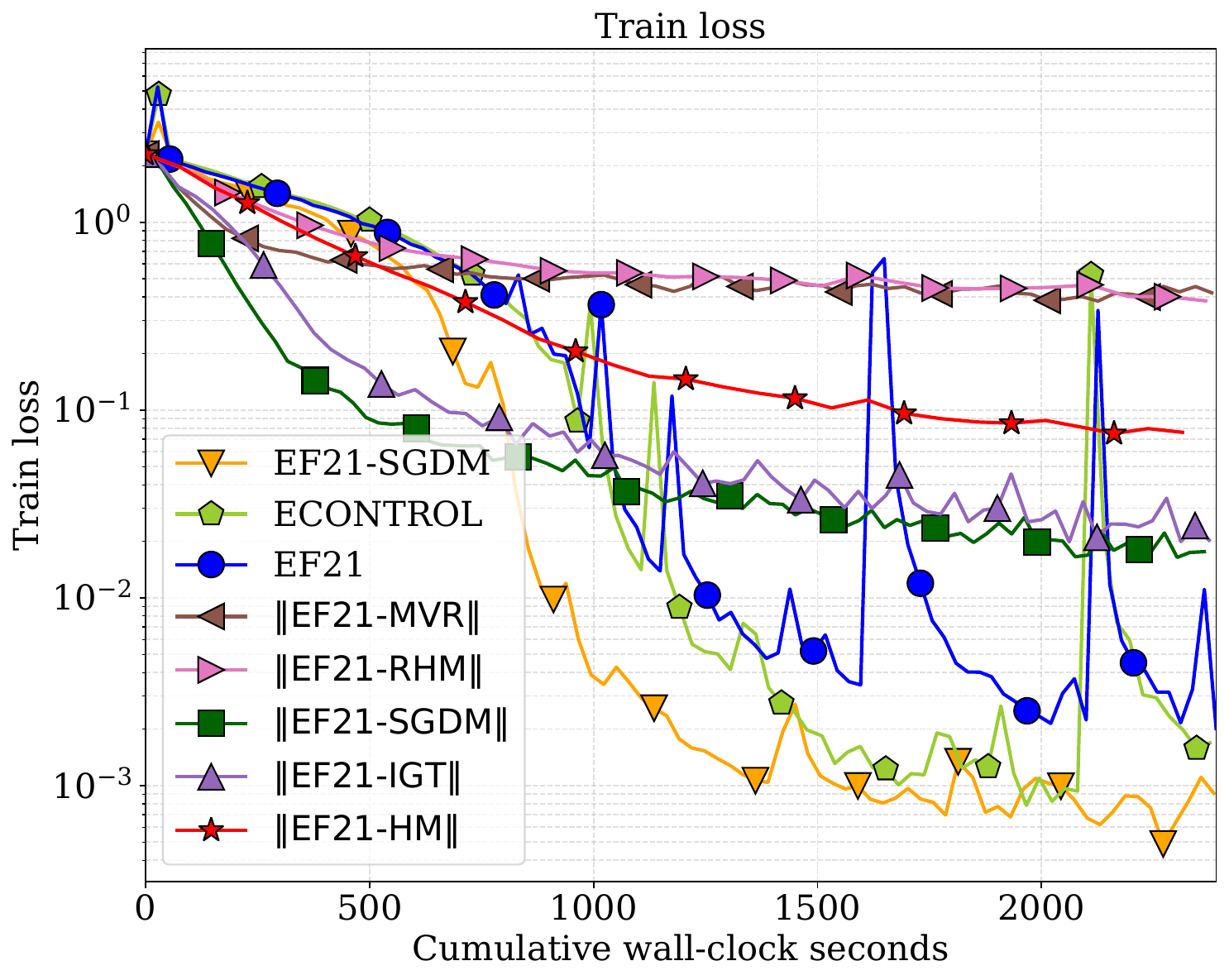}
\caption{Training loss}
\end{subfigure}\hfill
\begin{subfigure}{0.48\textwidth}
\includegraphics[width=\linewidth]{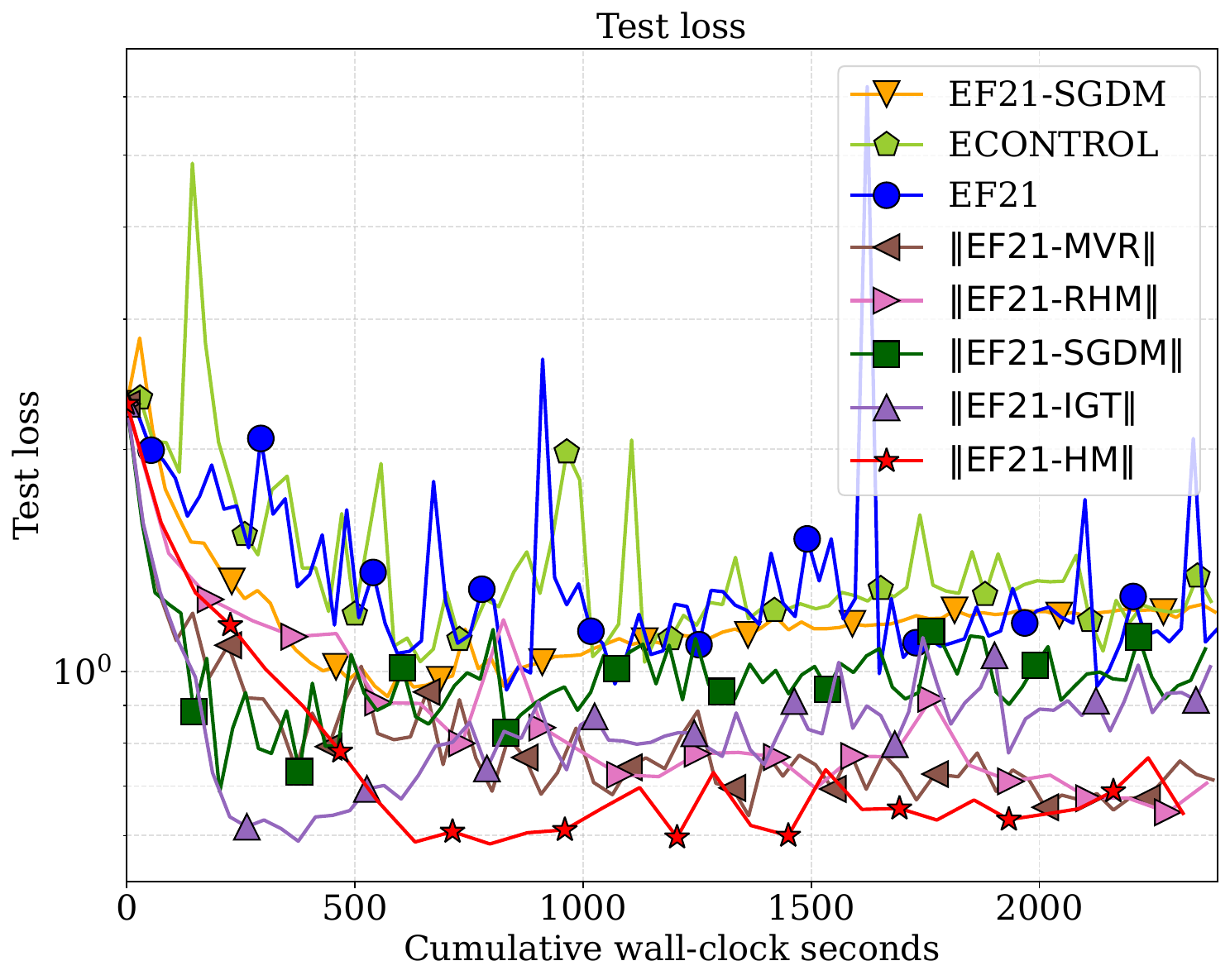}
\caption{Test loss}
\end{subfigure}
\par\medskip
\begin{subfigure}{0.48\textwidth}
\includegraphics[width=\linewidth]{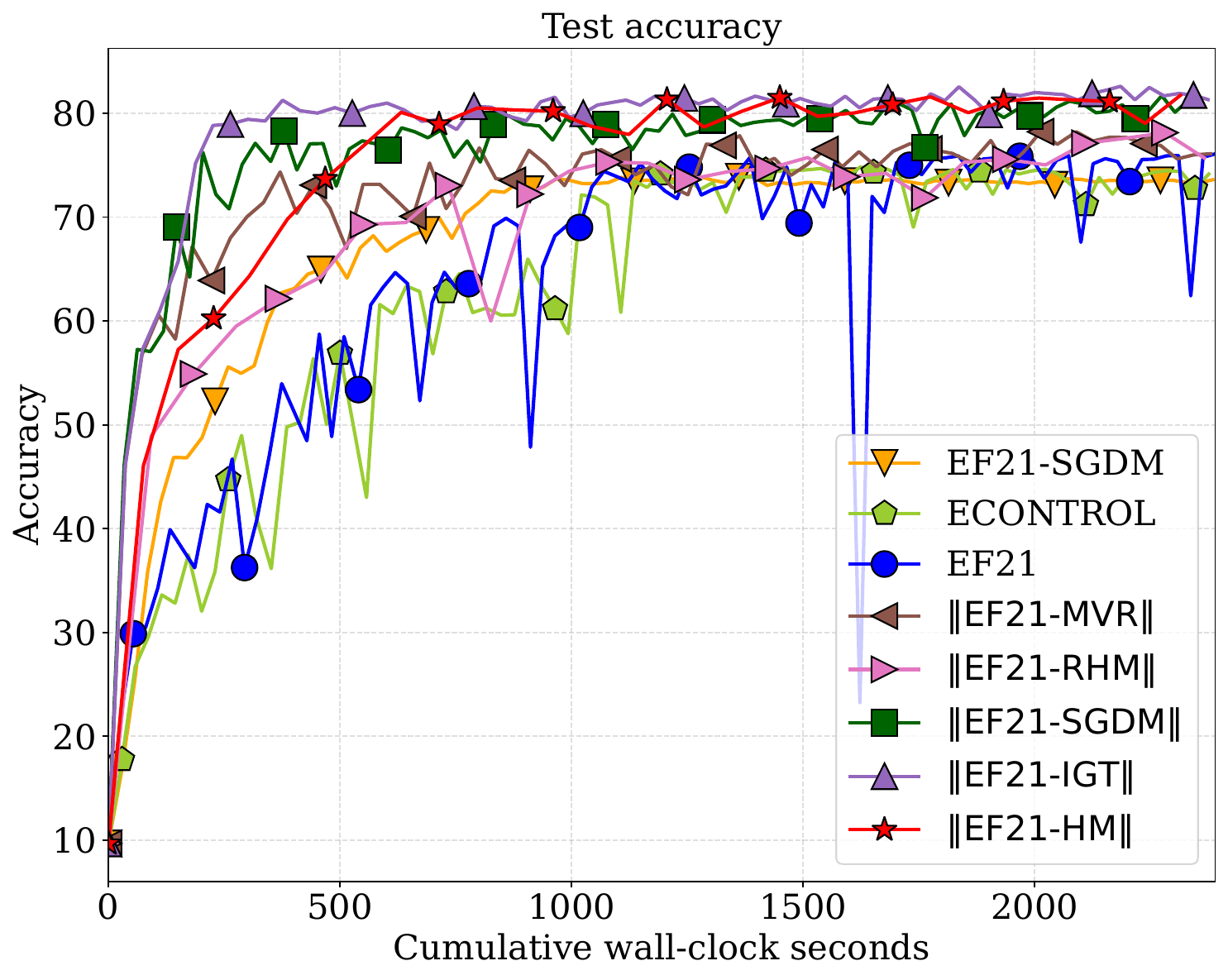}
\caption{Test accuracy}
\end{subfigure}
\caption{Wall-clock comparison with the timeline truncated to the completion
time of the fastest methods. The IGT trajectories are competitive with the
more costly HM method over this time horizon.}
\label{ch7:fig:restored_time_truncated}
\end{figure}
\FloatBarrier

\subsection{Comparison of per-epoch GPU cost}
\label{ch7:sec:restored_gpu_cost}
Table~\ref{ch7:tab:restored_gpu_cost} reports the average GPU time per epoch
and the relative cost normalized by \normEFname{SGDM}, the first-order
baseline. The underlying timings include the GPU work specified in
Section~\ref{ch7:sec:restored_timing}.

\begin{table}[htbp]
\centering
\caption{Average per-epoch GPU runtime and relative cost. The row order
matches Table~\ref{ch7:main_table:accuracy_results_chapter_7}.}
\label{ch7:tab:restored_gpu_cost}
\small
\setlength{\tabcolsep}{5pt}
\renewcommand{\arraystretch}{1.15}
\begin{tabular}{@{}lcc@{}}
\toprule
\textbf{Method} & \textbf{Mean seconds/epoch} & \textbf{Relative cost} \\
\midrule
\algname{EF21-SGDM} & $28.373$ & $0.984\times$ \\
\algname{EControl} & $28.943$ & $1.004\times$ \\
\algname{EF21-SGD} & $26.553$ & $0.921\times$ \\
\midrule
\normEFname{MVR} & $36.892$ & $1.279\times$ \\
\normEFname{RHM} & $86.380$ & $2.996\times$ \\
\normEFname{SGDM} & $28.835$ & $1.000\times$ \\
\normEFname{IGT} & $33.026$ & $1.145\times$ \\
\normEFname{HM} & $78.390$ & $2.719\times$ \\
\bottomrule
\end{tabular}
\end{table}

\paragraph{Hessian-corrected momentum.}
Although HM uses two automatic differentiation passes, these do not have
the same cost as two ordinary first-order backward passes. The first pass
constructs and retains the graph of the gradient computation through
\texttt{create\_graph=True}, increasing computational work and memory use.
The Hessian-vector product then traverses this graph. This explains why a
simple count of backward passes does not predict the measured relative
cost of $2.719\times$. Fixed epoch overheads and implementation-level GPU
optimizations also affect the observed ratio.

\paragraph{Randomized Hessian-corrected momentum.}
RHM combines a gradient evaluation at the new point with graph construction
and a Hessian-vector product at a different point. Its approximately three
backpropagation-equivalents give a measured relative cost of $2.996\times$
in these experiments.

\paragraph{Momentum variance reduction.}
MVR uses two first-order differentiation passes without the expensive
second-order graph construction. The previous-point gradient is obtained
through \texttt{autograd.grad}, without populating parameter
\texttt{.grad} fields, and the same minibatch is reused for both evaluations.
Together with fixed per-epoch overheads, these implementation details
explain why the measured relative cost in
Table~\ref{ch7:tab:restored_gpu_cost} is substantially below a factor of two.
\FloatBarrier

  \endgroup

\chapter{Supplementary Material and Proofs for Chapter 8}
  \begingroup
  \thesischaptersetup{8}%
  \stepEightAAppendixColors{8}%

\newpage 
\newpage 
\section{Notation}

For random variables $u,v$, we use $\Exp{u}$ for the expectation of $u$, and $\thesischeightExpSub{v}{u}$ for the expectation of $u$ with respect to $v$.
For functions $f, g : \R^d \to \R$, we use $\nabla f$ and $\nabla^2 f$ to denote the gradient and Hessian of $f$, respectively. We write $f(x) = \cO(g(x))$ to indicate that there exists a constant $C > 0$ and a value $x^{0} \in \R$ such that $f(x) \leq C \cdot g(x)$ for all $x \geq x^{0}$.
For vectors $x,y\in\R^d$, $\inner{x}{y} = \sum_{i=1}^d x_i y_i$, while $\norm{x}$,
 $\normstar{x}$, $\normtwo{x}$, and $\norm{x}_{\infty}$ denote its arbitrary norm,  associated dual norm,  Euclidean norm, and  $\ell_{\infty}$-norm, respectively. 

For the stochastic proofs, let $\mathcal F_k$ be the sigma-field generated by
the initialization and all random draws through iteration $k-1$, before the
fresh variables at iteration $k$ are sampled. Thus, $x^k$, $m^k$, and the
update $x^{k+1}$ are $\mathcal F_k$-measurable. The shared sample
$\xi^{k+1}$ is independent of $\mathcal F_k$. In Variant~1, $b_{k+1}$ is
also independent of $\mathcal F_k$ and of $\xi^{k+1}$; after drawing it, we
write $\mathcal G_{k+1}\eqdef
\sigma(\mathcal F_k,b_{k+1})$. Consequently, the pointwise oracle conditions
in Assumption~\ref{ch8:assum:Hessian_boundedvariance_v2_chapter_8} imply, for
every $\mathcal F_k$-measurable $x,w$,
\begin{align*}
\Exp{\nabla f(x;\xi^{k+1})\mid\mathcal F_k}
&=\nabla f(x),\\
\Exp{\normtwo{\nabla f(x;\xi^{k+1})-\nabla f(x)}^2\mid\mathcal F_k}
&\leq\sigma_g^2,\\
\Exp{\nabla^2 f(x;\xi^{k+1})\mid\mathcal F_k}
&=\nabla^2 f(x),\\
\Exp{\normtwo{(\nabla^2 f(x;\xi^{k+1})-\nabla^2 f(x))w}^2
      \mid\mathcal F_k}
&\leq\sigma_H^2\normtwo{w}^2.
\end{align*}
The same Hessian statements hold conditionally on $\mathcal G_{k+1}$ for
$\mathcal G_{k+1}$-measurable arguments. Likewise, when the \steponeonechange{mean-squared
generalized smoothness} assumption is invoked, its pointwise bound holds
conditionally on $\mathcal F_k$ for $\mathcal F_k$-measurable arguments.
These conditional forms do not strengthen the oracle model; they follow from
the fresh-sample independence in Algorithm~\ref{ch8:alg:2_order_momentum_chapter_8}.

\section{Useful Lemma from Relaxed Smoothness}

From Assumptions~\ref{ch8:assum:Relaxed_Lipschitz_Grad_chapter_8} and~\ref{ch8:assum:Relaxed_Lipschitz_Hessian_chapter_8}, we obtain the following lemma that is useful for our analysis. 
\begin{lemma}[Proposition 3.2 and Theorem 1 from \textcite{chen2023generalized}]
\label{ch8:lemma:Relaxed_Lipschitz_Grad_chapter_8}
    Let $f$ satisfy Assumption~\ref{ch8:assum:Relaxed_Lipschitz_Grad_chapter_8}. Then, for all $x,y \in \R^d$,
    \begin{align*}
        \normstar{\nabla f(x) -\nabla f(y)} & \leq (L_0 + L_1\normstar{\nabla f(y)})\exp\left(L_1\norm{x-y}\right)\norm{x-y};\\
        D_f(x,y) & \leq \frac{1}{2}(L_0 + L_1\normstar{\nabla f(y)})\exp\left(L_1\norm{x-y}\right)\norm{x-y}^2,
    \end{align*}
    where $D_f(x,y) \eqdef f(x) - f(y) - \la\nabla f(y), x-y\rangle$ denotes the Bregman divergence.
    Furthermore, if $f$ is twice-differentiable, then for any $x \in \R^d$,
    \begin{align*}
        \norm{\nabla^2 f(x)}_{\text{op}} & \eqdef \sup_{u\neq 0}\frac{\normstar{\nabla^2f(x)u}}{\norm{u}} \leq L_0 + L_1 \normstar{\nabla f(x)}.
    \end{align*}
\end{lemma}
\begin{proof}
    The proof follows the same lines as Proposition 3.2 and Theorem 1 in \textcite{chen2023generalized}, adapted to the setting of arbitrary norms. 
\end{proof}
\begin{lemma}\label{ch8:lemma:Relaxed_Grad_and_Hessian_chapter_8}
Let $f$ satisfy  Assumptions~\ref{ch8:assum:Relaxed_Lipschitz_Grad_chapter_8} and~\ref{ch8:assum:Relaxed_Lipschitz_Hessian_chapter_8}. Then, for all $x,y\in\R^d$, 
\begin{align*}
    \normstar{Z_f(x,y)} 
    & \leq \frac{1}{2}(M_0 + M_1 \normstar{\nabla f(y)}) \norm{x-y}^2 \\
    &\quad + \frac{1}{3} M_1 (L_0+L_1 \normstar{\nabla f(y)}) \exp(L_1 \norm{x-y}) \norm{x-y}^3,
\end{align*}
where we denote $Z_f(x,y) \eqdef \nabla f(x) - \nabla f(y) - \nabla^2 f(y)(x-y)$.
\end{lemma}
\begin{proof}
Since $f$ is twice differentiable, with $\hat x \eqdef y + \tau(x-y)$ and $w_\theta(\hat x, y) = \theta\hat x + (1-\theta)y$,
\begin{align*}
  Z_f(x,y) 
  & = \int_{\tau=0}^1 \nabla^2 f(\hat x) (x-y) d\tau  - \nabla^2 f(y)(x-y) \\
  & = \int_{\tau=0}^1 [\nabla^2 f(\hat x)  - \nabla^2 f(y)](x-y) d\tau.  
\end{align*}
Therefore,
\begin{align*}
    \normstar{Z_f(x,y) } 
    & \leq \int_{\tau=0}^1 \normstar{[\nabla^2 f(\hat x)  - \nabla^2 f(y)](x-y)} d\tau. 
\end{align*}

Next, from~Assumption~\ref{ch8:assum:Relaxed_Lipschitz_Hessian_chapter_8} and $\norm{\hat x - y} = \tau\norm{x-y}$, 
\begin{align*}
    \normstar{Z_f(x,y) } 
    & \leq \int_{\tau=0}^1 (M_0 + M_1 \underset{\theta \in [0,1]}{\sup} \normstar{ \nabla f(w_\theta(\hat x,y))}) \norm{\hat x-y}\norm{x-y} d\tau  \\
    & \leq \frac{M_0 \norm{x-y}^2}{2}\\
    &\quad + M_1 \norm{x-y}^2  \int_{\tau=0}^1 \tau  \underset{\theta \in [0,1]}{\sup} \normstar{ \nabla f(w_\theta(\hat x,y)) - \nabla f(y) + \nabla f(y) } d\tau \\
    & \leq \frac{M_0 \norm{x-y}^2}{2} \\
    &\quad + M_1 \norm{x-y}^2  \int_{\tau=0}^1 \tau  \underset{\theta \in [0,1]}{\sup} \normstar{ \nabla f(w_\theta(\hat x,y)) - \nabla f(y) } d\tau \\
    &\quad + M_1 \norm{x-y}^2  \int_{\tau=0}^1 \tau  \normstar{ \nabla f(y) } d\tau \\
    & = \frac{M_0+M_1\normstar{\nabla f(y)}}{2}\norm{x-y}^2 \\
    &\quad+ M_1 \norm{x-y}^2  \int_{\tau=0}^1 \tau  \underset{\theta \in [0,1]}{\sup} \normstar{ \nabla f(w_\theta(\hat x,y)) - \nabla f(y) } d\tau.
\end{align*}

Next, from Lemma~\ref{ch8:lemma:Relaxed_Lipschitz_Grad_chapter_8} and $w_\theta(\hat x,y) - y = \tau\theta(x-y)$,
\begin{align*}
\MoveEqLeft[1] \normstar{ \nabla f(w_\theta(\hat x,y)) - \nabla f(y) } \\
& \leq (L_0+L_1 \normstar{\nabla f(y)}) \exp( \tau\theta L_1 \norm{x-y}) \tau\theta \norm{x-y} \\
& \leq (L_0+L_1 \normstar{\nabla f(y)}) \exp( \tau L_1 \norm{x-y}) \tau \norm{x-y},
\end{align*}
where we used $\theta \leq 1$.
Therefore, plugging this result into the upper bound for
$\normstar{Z_f(x,y)}$ and using $\exp(\tau L_1\norm{x-y}) \leq \exp(L_1\norm{x-y})$ and $\int_0^1\tau^2 d\tau = 1/3$ completes the proof.

\end{proof}

\newpage
\section{Proof of~Theorem~\ref{ch8:thm:V1_relaxed_smoothness_chapter_8}  }

We prove~Theorem~\ref{ch8:thm:V1_relaxed_smoothness_chapter_8} in the following steps.

\paragraph{Step 1) Proving the descent inequality.}

We follow the proof arguments from Lemma D.1. of \textcite{pethick2025training}.
By~Lemma~\ref{ch8:lemma:Relaxed_Lipschitz_Grad_chapter_8}, and by the fact that $ x^{k+1}= x^{k} + \text{lmo}(m^{k})$, 
\begin{align*}
    f(x^{k+1}) 
    & \leq f(x^{k}) + \inner{\nabla f(x^{k})}{\text{lmo}(m^{k})} \\
    &\quad + \frac{L_0 + L_1 \normstar{\nabla f(x^{k})}}{2} \exp(L_1 \norm{x^{k+1}-x^{k}})\norm{\text{lmo}(m^{k})}^2 \\
    & =  f(x^{k}) + \inner{\nabla f(x^{k}) - m^{k}}{\text{lmo}(m^{k})} + \inner{m^{k}}{\text{lmo}(m^{k})} \\
    &\quad + \frac{L_0 + L_1 \normstar{\nabla f(x^{k})}}{2} \exp(L_1 \norm{x^{k+1}-x^{k}})\norm{\text{lmo}(m^{k})}^2.
\end{align*}

By the definition of the dual norm, 
\begin{align*}
    f(x^{k+1}) 
& \leq f(x^{k}) + \normstar{\nabla f(x^{k}) - m^{k}}\norm{\text{lmo}(m^{k})} + \inner{m^{k}}{\text{lmo}(m^{k})} \\
    &\quad + \frac{L_0 + L_1 \normstar{\nabla f(x^{k})}}{2} \exp(L_1 \norm{x^{k+1}-x^{k}})\norm{\text{lmo}(m^{k})}^2.
\end{align*}

By the fact that $\norm{x^{k+1} - x^{k}} = \norm{\text{lmo}(m^{k})} \leq \gamma_k$, 
\begin{align*}
    f(x^{k+1})
    &\quad + \inner{m^{k}}{\text{lmo}(m^{k})} + \frac{L_0 + L_1 \normstar{\nabla f(x^{k})}}{2} \exp(L_1 \gamma_k) \gamma_k^2 \\
    &\leq f(x^{k}) + \gamma_k\normstar{\nabla f(x^{k}) - m^{k}} + \inner{m^{k}}{\text{lmo}(m^{k})}\\
    &\quad + \frac{L_0 + L_1 \normstar{\nabla f(x^{k})}}{2} \exp(L_1 \gamma_k) \gamma_k^2.
\end{align*}

Since for all $u \in \mathcal{X}$
\begin{equation*}
    \normstar{u} = \underset{v: \norm{v} \leq 1}{\max} \inner{u}{v} = \max_{v: \norm{v} \leq \gamma} \inner{u}{\frac{1}{\gamma} v} = - \inner{u}{ \frac{1}{\gamma} \text{lmo}(u)},
\end{equation*}
we obtain
\begin{equation*}
    \inner{m^{k}}{\text{lmo}(m^{k})} = \gamma_k \inner{m^{k}}{ \frac{1}{\gamma_k} \text{lmo}(m^{k})} = - \gamma_k \normstar{m^{k}}.
\end{equation*}
Therefore,
\begin{align*}
    f(x^{k+1})
    &\leq f(x^{k})+\gamma_k\normstar{\nabla f(x^{k})-m^{k}}
      -\gamma_k\normstar{m^{k}} \\
    &\quad
      +\frac{L_0+L_1\normstar{\nabla f(x^{k})}}{2}
       \exp(L_1\gamma_k)\gamma_k^2.
\end{align*}

By the triangle inequality, i.e. $\normstar{a} \geq \normstar{b} - \normstar{a-b}$ for $a,b\in\R^d$,
\begin{align*}
    f(x^{k+1})
    &\leq f(x^{k})+2\gamma_k\normstar{\nabla f(x^{k})-m^{k}}
      -\gamma_k\normstar{\nabla f(x^{k})} \\
    &\quad
      +\frac{L_0+L_1\normstar{\nabla f(x^{k})}}{2}
       \exp(L_1\gamma_k)\gamma_k^2.
\end{align*}

Finally, by summing the inequality over $k=0,1,\ldots,K$ and by re-arranging the terms,
\begin{align}
 \sum_{k=0}^K \gamma_k\varphi_k\normstar{\nabla f(x^{k})}
 &\overset{(a)}{=}
 \sum_{k=0}^K \gamma_k
 \left(1-\exp(L_1\gamma_k)\frac{L_1\gamma_k}{2}\right)
 \normstar{\nabla f(x^{k})} \notag\\
 &\leq \sum_{k=0}^K \left(f(x^{k})-f(x^{k+1})\right) \notag\\
 &\quad+2\sum_{k=0}^K\gamma_k\normstar{\nabla f(x^{k})-m^{k}} \notag\\
 &\quad+\frac{L_0}{2}\sum_{k=0}^K\exp(L_1\gamma_k)\gamma_k^2 \notag\\
 &\overset{(b)}{\leq} (f(x^{0})-f^{\inf}) \notag\\
 &\quad+2\sum_{k=0}^K\gamma_k\normstar{\nabla f(x^{k})-m^{k}} \notag\\
 &\quad+\frac{L_0}{2}\sum_{k=0}^K\exp(L_1\gamma_k)\gamma_k^2 \notag\\
 &\overset{(c)}{=} \Delta \notag\\
 &\quad+2\sum_{k=0}^K\gamma_k\normstar{\nabla f(x^{k})-m^{k}} \notag\\
 &\quad+\frac{L_0}{2}\sum_{k=0}^K\exp(L_1\gamma_k)\gamma_k^2.
 \label{ch8:eq:descent_inequality_chapter_8}
\end{align}
where we reach  $(a)$ by defining $\varphi_k \eqdef\left( 1 - \exp(L_1 \gamma_k) 
 L_1\gamma_k/2 \right) $, $(b)$ by using the fact that $ {f(x) \geq f^{\inf}}$, and $(c)$ by denoting $\Delta \eqdef f(x^{0}) - f^{\inf}$.

\paragraph{Step 2) Bounding the error term.}
Next, we bound $\normstar{r^{k+1}}$, where $r^{k+1} \eqdef m^{k+1}-\nabla f(x^{k+1})$. 
From the definition of $r^{k}$,
\begin{align}
    r^{k+1} 
    & = m^{k+1} - \nabla f(x^{k+1}) \notag \\
    & \overset{(a)}{=} (1-\eta_k) [m^{k} + \nabla^2 f(\hat x^{k+1};\xi^{k+1}) (x^{k+1}-x^{k})  - \nabla f(x^{k+1})]\notag\\
    &\quad+ \eta_k[\nabla f(x^{k+1};\xi^{k+1})-\nabla f(x^{k+1})] \notag \\
    & \overset{(b)}{=} (1-\eta_k) r^{k} + (1-\eta_k)W_{k+1} +  \eta_kV_{k+1}, \label{ch8:eqn:SHARP_e_chapter_8}
\end{align}
where we reach $(a)$ by \eqref{ch8:eqn:second_order_momentum_update_chapter_8}, and $(b)$ by denoting  $W_{k+1} \eqdef  \nabla^2 f(\hat x^{k+1};\xi^{k+1}) (x^{k+1}-x^{k}) - (\nabla f(x^{k+1}) - \nabla f(x^{k}))$ and $V_{k+1} \eqdef\nabla f(x^{k+1};\xi^{k+1})-\nabla f(x^{k+1})$.

By recursively applying this identity, 
\begin{align*}
    r^{k+1} & = \prod_{t=0}^k (1-\eta_t) r^{0} \\
    &\quad +\sum_{t=0}^k \left( \prod_{j=t+1}^k (1-\eta_j) \right) (1-\eta_t)  W_{t+1}\\
    &\quad +  \sum_{t=0}^k \left(\prod_{j=t+1}^k (1-\eta_j) \right) \eta_t V_{t+1}.
\end{align*}

Next, by taking $\normstar{\cdot}$, and by taking the expectation,
\begin{align*}
    \Exp{\normstar{r^{k+1} } } 
    & \overset{(a)}{\leq} \normstar{r^{0}}  \prod_{t=0}^k (1-\eta_t) +\Exp{\normstar{\sum_{t=0}^k \left( \prod_{j=t+1}^k (1-\eta_j)\right) (1-\eta_t)  W_{t+1}}} \\
    &\quad + \Exp{\normstar{\sum_{t=0}^k \left(\prod_{j=t+1}^k (1-\eta_j) \right) \eta_t V_{t+1}}} \\
    & \overset{(b)}{\leq} \normstar{r^{0}}  \prod_{t=0}^k(1-\eta_t) + \Exp{\normstar{\sum_{t=0}^k \left( \prod_{j=t+1}^k (1-\eta_j)\right) (1-\eta_t)  W_{t+1}}} \\
    &\quad + \sqrt{ \Exp{\normstar{\sum_{t=0}^k \left(\prod_{j=t+1}^k (1-\eta_j) \right) \eta_t V_{t+1}}^2} } \\
    & \overset{(c)}{\leq} \normstar{r^{0}}  \prod_{t=0}^k(1-\eta_t) + {\color{red}\bar{\rho}}\underbrace{\Exp{\normtwo{\sum_{t=0}^k \left( \prod_{j=t+1}^k (1-\eta_j)\right) (1-\eta_t)  W_{t+1}}}}_{\eqdef \circledOne} \\
    &\quad + {\color{red}\bar{\rho}}\underbrace{\sqrt{ \Exp{\normtwo{\sum_{t=0}^k \left(\prod_{j=t+1}^k (1-\eta_j) \right) \eta_t V_{t+1}}^2} }}_{\eqdef \circledTwo },
\end{align*}
where $(a)$ follows from the triangle inequality, $(b)$ results from Jensen's inequality applied to the second sum, and $(c)$ obtains from \eqref{ch8:eqn:norm_diff_chapter_8}.

To bound  $\Exp{\normstar{r^{k+1} } }$, we must bound $\circledOne$ ~and $\circledTwo$. First, we bound $\circledTwo$. By~Assumption~\ref{ch8:assum:Hessian_boundedvariance_v2_chapter_8} and the fresh-sample condition,
\[
\Exp{V_{t+1}\mid\mathcal F_t}=0
\quad\text{and}\quad
\Exp{\normtwo{V_{t+1}}^2\mid\mathcal F_t}\leq\sigma_g^2.
\]
For $0\leq i<j$, the vector $V_{i+1}$ is
$\mathcal F_j$-measurable. Hence, by the tower property,
\[
\Exp{\inner{V_{j+1}}{V_{i+1}}}
=\Exp{\inner{\Exp{V_{j+1}\mid\mathcal F_j}}{V_{i+1}}}=0.
\]
Therefore, 
\begin{align*}
    \circledTwo^2 & = \sum_{t=0}^k \left(\prod_{j=t+1}^k (1-\eta_j)^2 \right) \eta_t^2 \Exp{\normtwo{V_{t+1}}^2} \\
    & \leq \sigma_g^2 \sum_{t=0}^k \left(\prod_{j=t+1}^k (1-\eta_j)^2 \right) \eta_t^2.
\end{align*}
Plugging this result into the main inequality yields
\begin{align*}
 \Exp{\normstar{r^{k+1} } } & \leq \normstar{r^{0}}  \prod_{t=0}^k(1-\eta_t) + {\color{red}\bar{\rho}}\,\circledOne \\
 &\quad + {\color{red}\bar{\rho}} \sigma_g \sqrt{\sum_{t=0}^k \left(\prod_{j=t+1}^k (1-\eta_j)^2 \right) \eta_t^2}. 
\end{align*}

Second, we bound $\circledOne$. Conditional first on
$\mathcal G_{t+1}$ and then on $\mathcal F_t$, Assumption~\ref{ch8:assum:Hessian_boundedvariance_v2_chapter_8} and
$b_{t+1}\sim\thesischeightcU(0,1)$ give
\begin{align*}
\Exp{W_{t+1}\mid\mathcal F_t}
&=\int_0^1\nabla^2 f\bigl(bx^{t+1}+(1-b)x^t\bigr)
       (x^{t+1}-x^t)\,db\\
&\quad-\bigl(\nabla f(x^{t+1})-\nabla f(x^t)\bigr)
=0.
\end{align*}
For $0\leq i<j$, $W_{i+1}$ is $\mathcal F_j$-measurable, and therefore
\[
\Exp{\inner{W_{j+1}}{W_{i+1}}}
=\Exp{\inner{\Exp{W_{j+1}\mid\mathcal F_j}}{W_{i+1}}}=0.
\]
Thus, the summands of $\circledOne$ form a martingale difference sequence with respect to $\{\mathcal F_t\}$. To bound $\circledOne$, we first bound the conditional second moments of $W_{t+1}$ given $\mathcal F_t$ in Step~3 and then combine them in Step~4.

\paragraph{Step 3) Bounding $\Exp{\normtwo{ W_{k+1}}^2}$.} 
Next, we bound $\Exp\left[\normtwo{ W_{k+1}}^2\mid\mathcal F_k\right]$. Throughout this step, $\Exp{}$ denotes the conditional expectation given $\mathcal F_k$; recall that $x^k$ and $x^{k+1}$ are $\mathcal F_k$-measurable, while $b_{k+1}$ and $\xi^{k+1}$ are independent of $\mathcal F_k$.
By the fact that $W_{k+1} =  \nabla^2 f(\hat x^{k+1};\xi^{k+1}) (x^{k+1}-x^{k}) - (\nabla f(x^{k+1}) - \nabla f(x^{k}))$, and that $\normtwo{x+y+z}^2 \leq 3 \normtwo{x}^2+3\normtwo{y}^2+3\normtwo{z}^2$ for $x,y,z\in\R^d$,
\begin{align*}
    \Exp{\normtwo{ W_{k+1}}^2} 
    & \leq 3 \Exp{\normtwo{\nabla f(x^{k+1})-\nabla f(x^{k})}^2} \\
    &\quad + 3 \Exp{\normtwo{\nabla^2 f(\hat x^{k+1})(x^{k+1}-x^{k})}^2} \\
    &\quad + 3\Exp{ \normtwo{(\nabla^2 f(\hat x^{k+1};\xi^{k+1}) - \nabla^2 f(\hat x^{k+1}))(x^{k+1}-x^{k})}^2 }.
\end{align*}

Next, by Assumption~\ref{ch8:assum:Hessian_boundedvariance_v2_chapter_8}, 
\begin{align*}
    \Exp{\normtwo{ W_{k+1}}^2} 
    & \leq 3 \Exp{\normtwo{\nabla f(x^{k+1})-\nabla f(x^{k})}^2} \\
    &\quad + 3 \Exp{\normtwo{\nabla^2 f(\hat x^{k+1})(x^{k+1}-x^{k})}^2} \\
    &\quad + 3 \sigma_H^2 \Exp{ \normtwo{x^{k+1}-x^{k}}^2 } \\
    & \overset{(a)}{\leq} \frac{3}{{\color{red} \thesischeightubar \rho }^2} \Exp{\normstar{\nabla f(x^{k+1})-\nabla f(x^{k})}^2} \\
    &\quad +  \frac{3}{{\color{red} \thesischeightubar \rho }^2} \Exp{\normstar{\nabla^2 f(\hat x^{k+1})(x^{k+1}-x^{k})}^2} \\
    &\quad + \frac{3 \sigma_H^2}{{\color{blue}\thesischeightubar\theta}^2} \Exp{ \norm{x^{k+1}-x^{k}}^2 }\\
    & \leq \frac{3}{{\color{red} \thesischeightubar \rho }^2} \Exp{\normstar{\nabla f(x^{k+1})-\nabla f(x^{k})}^2} \\
    &\quad +  \frac{3}{{\color{red} \thesischeightubar \rho}^2} \Exp{ \norm{\nabla^2 f(\hat x^{k+1})}_{\text{op}}^2 \norm{x^{k+1}-x^{k}}^2} \\
    &\quad + \frac{3 \sigma_H^2}{{\color{blue}\thesischeightubar\theta}^2} \Exp{ \norm{x^{k+1}-x^{k}}^2 }
\end{align*}
where in $(a)$ we used \eqref{ch8:eqn:norm_diff_chapter_8}.

Next, applying Lemma~\ref{ch8:lemma:Relaxed_Lipschitz_Grad_chapter_8}, we obtain
\begin{align*}
    \Exp{\normtwo{ W_{k+1}}^2} 
    & \leq \frac{3}{{\color{red} \thesischeightubar \rho }^2} \Exp{ (L_0 + L_1 \normstar{\nabla f(x^{k})})^2 \exp(2 L_1 \norm{x^{k+1}-x^{k}}) \norm{x^{k+1} -x^{k}}^2}  \\
    & \quad +  \frac{3}{{\color{red} \thesischeightubar \rho }^2} \Exp{ (L_0 + L_1 \normstar{\nabla f(\hat x^{k+1})})^2\norm{x^{k+1}-x^{k}}^2} \\
    &\quad + \frac{3 \sigma_H^2}{{\color{blue}\thesischeightubar\theta}^2} \Exp{ \norm{x^{k+1}-x^{k}}^2 }. 
\end{align*}

By the fact that $\norm{x^{k+1}-x^{k}} \leq \gamma_{k}$, 
\begin{align*}
\Exp{\normtwo{ W_{k+1}}^2} 
& \leq \frac{3}{{\color{red} \thesischeightubar \rho }^2}  (L_0 + L_1 \normstar{\nabla f(x^{k})})^2 \exp(2 L_1 \gamma_{k}) \gamma_{k}^2  \\
& \quad +  \frac{3}{{\color{red} \thesischeightubar \rho}^2} \gamma_{k}^2 \Exp{ (L_0 + L_1 \normstar{\nabla f(\hat x^{k+1})})^2} \\
&\quad + \frac{3 \sigma_H^2}{{\color{blue}\thesischeightubar\theta}^2} \gamma_{k}^2. 
\end{align*}

Next, since 
\begin{align*}
    \MoveEqLeft[1] (L_0+L_1\normstar{\nabla f(\hat x^{k+1})})^2 \\
    & \leq (L_0 + L_1\normstar{\nabla f(x^{k})} + L_1\normstar{\nabla f(\hat x^{k+1}) - \nabla f(x^{k})})^2 \\
    & \leq 2 (L_0 + L_1 \normstar{\nabla f(x^{k})})^2 + 2 L_1^2 \normstar{\nabla f(\hat x^{k+1}) - \nabla f(x^{k})}^2 \\
    & \leq  2 (L_0 + L_1 \normstar{\nabla f(x^{k})})^2  \\
    &\quad+ 2 L_1^2 (L_0 + L_1\normstar{\nabla f(x^{k})})^2 \exp(2L_1 \norm{\hat x^{k+1} - x^{k}}) \norm{\hat x^{k+1}-x^{k}}^2,
\end{align*}
we obtain
\begin{align*}
    \Exp{\normtwo{ W_{k+1}}^2}  & \leq  \frac{3}{{\color{red} \thesischeightubar \rho }^2}  (L_0 + L_1 \normstar{\nabla f(x^{k})})^2 (\exp(2 L_1 \gamma_{k}) + 2)\gamma_{k}^2 \\
    &\quad + \frac{3 \sigma_H^2}{{\color{blue}\thesischeightubar\theta}^2} \gamma_{k}^2 \\
    &\quad +  \frac{6}{{\color{red} \thesischeightubar \rho}^2} L_1^2\gamma_{k}^2  (L_0 + L_1 \normstar{\nabla f(x^{k})})^2 \\
        &\qquad\cdot \Exp{ \exp(2L_1 \norm{\hat x^{k+1} - x^{k}})\norm{\hat x^{k+1} -x^{k}}^2}   \\
    & \overset{(a)}{\leq} \frac{3}{{\color{red} \thesischeightubar \rho}^2}  (L_0 + L_1 \normstar{\nabla f(x^{k})})^2 (\exp(2 L_1 \gamma_{k}) + 2)\gamma_{k}^2 \\
    &\quad + \frac{3 \sigma_H^2}{{\color{blue}\thesischeightubar\theta}^2} \gamma_{k}^2 \\
    &\quad +  \frac{6}{{\color{red} \thesischeightubar \rho}^2} L_1^2\gamma_{k}^2  (L_0 + L_1 \normstar{\nabla f(x^{k})})^2 \\
        &\qquad\cdot \Exp{ \exp(2L_1 b_{k+1} \norm{x^{k+1} - x^{k}})b_{k+1}^2\norm{x^{k+1} -x^{k}}^2} ,
\end{align*}
where $(a)$ results from the definition of $\hat x^{k+1}$. 
Next, by the fact that $\norm{x^{k+1}-x^{k}} \leq \gamma_{k}$,
\begin{align*}
    \Exp{\normtwo{ W_{k+1}}^2} & \leq \frac{3}{{\color{red} \thesischeightubar \rho}^2}  (L_0 + L_1 \normstar{\nabla f(x^{k})})^2 (\exp(2 L_1 \gamma_{k}) + 2)\gamma_{k}^2  \\
    & \quad +  \frac{6}{{\color{red} \thesischeightubar \rho}^2} L_1^2\gamma_{k}^2  (L_0 + L_1 \normstar{\nabla f(x^{k})})^2 \\
        &\qquad\cdot \Exp{ \exp(2L_1 b_{k+1} \gamma_{k})b_{k+1}^2\gamma_{k}^2} \\
    &\quad + \frac{3 \sigma_H^2}{{\color{blue}\thesischeightubar\theta}^2} \gamma_{k}^2 \\
    & =  \frac{3}{{\color{red} \thesischeightubar \rho}^2}  (L_0 + L_1 \normstar{\nabla f(x^{k})})^2 (\exp(2 L_1 \gamma_{k}) + 2)\gamma_{k}^2  \\
    &\quad +  \frac{6}{{\color{red} \thesischeightubar \rho}^2} L_1^2\gamma_{k}^2  (L_0 + L_1 \normstar{\nabla f(x^{k})})^2  \gamma_{k}^2 \int_{0}^1 \exp(2L_1 b \gamma_{k})b^2 db  \\
    &\quad + \frac{3 \sigma_H^2}{{\color{blue}\thesischeightubar\theta}^2} \gamma_{k}^2. 
\end{align*}

Since 
\begin{equation*}
    \int \exp(a_k z) z^2 dz = \frac{1}{a_k} z^2 \exp(a_k z) - \frac{2}{a_k}\int z \exp(a_k z) dz,
\end{equation*}
we obtain $\int_{z=0}^1 \exp(a_k z) z^2 dz \leq \frac{1}{a_k} \exp(a_k)$ for $a_k = 2L_1\gamma_k > 0$. 
Therefore,
\begin{align*}
    \Exp{\normtwo{ W_{k+1}}^2}  
& \leq \frac{3}{{\color{red} \thesischeightubar \rho}^2}  (L_0 + L_1 \normstar{\nabla f(x^{k})})^2 (\exp(2 L_1 \gamma_{k}) + 2)\gamma_{k}^2  \\
& \quad +  \frac{6}{{\color{red} \thesischeightubar \rho}^2} L_1^2\gamma_{k}^2  (L_0 + L_1 \normstar{\nabla f(x^{k})})^2  \gamma_{k}^2 \frac{1}{2L_1 \gamma_{k}} \exp(2L_1 \gamma_{k}) \\
&\quad + \frac{3 \sigma_H^2}{{\color{blue}\thesischeightubar\theta}^2} \gamma_{k}^2.
\end{align*}

By re-arranging the terms, and by the fact that $(a+b)^2 \leq 2 a^2 + 2b^2$ for $a,b\in\R$,
\begin{align*}
    \Exp{\normtwo{ W_{k+1}}^2}  
    & \leq  \frac{3}{{\color{red} \thesischeightubar \rho}^2}  (L_0 + L_1 \normstar{\nabla f(x^{k})})^2 \left(\exp(2 L_1 \gamma_{k}) + \gamma_{k}L_1\exp(2 L_1 \gamma_{k}) + 2\right)\gamma_{k}^2  \\
    &\quad + \frac{3 \sigma_H^2}{{\color{blue}\thesischeightubar\theta}^2} \gamma_{k}^2 \\
    & \leq \frac{6}{{\color{red} \thesischeightubar \rho}^2}  L_1^2 \normstar{\nabla f(x^{k})}^2 \left(\exp(2 L_1 \gamma_{k}) + \gamma_{k}L_1\exp(2 L_1 \gamma_{k}) + 2\right)\gamma_{k}^2   \\
    &\quad +  \frac{6}{{\color{red} \thesischeightubar \rho}^2}  L_0^2 \left((\exp(2 L_1 \gamma_{k}) + \gamma_{k}L_1\exp(2 L_1 \gamma_{k}) + 2\right)\gamma_{k}^2 + \frac{3 \sigma_H^2}{{\color{blue}\thesischeightubar\theta}^2} \gamma_{k}^2.
\end{align*}

\paragraph{Step 4) Plugging the bound on $W_{k+1}$ back into the upper-bound for $\Exp{\normstar{r^{k+1}}}$.}
Step 3 gives
\begin{equation*}
\Exp\left[\normtwo{W_{t+1}}^2\mid\mathcal F_t\right] \leq A_{t+1} + B_{t+1},
\end{equation*}
where
\begin{align*}
A_{t+1} &\eqdef \frac{6L_1^2\gamma_{t}^2}{{\color{red} \thesischeightubar \rho}^2}\normstar{\nabla f(x^{t})}^2 \left(\exp(2 L_1 \gamma_{t}) + \gamma_{t}L_1\exp(2 L_1 \gamma_{t}) + 2\right),\\
B_{t+1} & \eqdef \frac{6L_0^2\gamma_{t}^2}{{\color{red} \thesischeightubar \rho}^2} \left(\exp(2 L_1 \gamma_{t}) + \gamma_{t}L_1\exp(2 L_1 \gamma_{t}) + 2\right) \\
&\quad + \frac{3 \sigma_H^2}{{\color{blue}\thesischeightubar\theta}^2} \gamma_{t}^2.
\end{align*}
Here $A_{t+1}$ is $\mathcal F_t$-measurable and $B_{t+1}$ is deterministic. Denote
\begin{equation*}
a_t \eqdef \left( \prod_{j=t+1}^k (1-\eta_j)\right) (1-\eta_t), \qquad S_s \eqdef \sum_{t=0}^{s} a_t W_{t+1}.
\end{equation*}
Since $\Exp\left[W_{s+1}\mid\mathcal F_s\right]=0$ and $S_{s-1}$ is $\mathcal F_s$-measurable, the conditional Jensen inequality gives, for any $\mathcal F_s$-measurable $Q \geq 0$,
\begin{align*}
\MoveEqLeft[1] \Exp\left[\left(\normtwo{S_s}^2 + Q\right)^{1/2}\mid\mathcal F_s\right] \\
&\leq \left(\normtwo{S_{s-1}}^2 + a_s^2\left(A_{s+1}+B_{s+1}\right) + Q\right)^{1/2}\\
&\leq \left(\normtwo{S_{s-1}}^2 + a_s^2B_{s+1} + Q\right)^{1/2} + a_s\sqrt{A_{s+1}},
\end{align*}
where the last step uses $\sqrt{u+v}\leq\sqrt{u}+\sqrt{v}$. Applying this inequality successively for $s = k, k-1, \ldots, 0$ (starting with $Q = 0$ and keeping the deterministic terms $a_s^2B_{s+1}$ in $Q$), we obtain
\begin{align*}
\circledOne = \Exp{\normtwo{S_k}}
\leq \sqrt{\sum_{t=0}^k a_t^2 B_{t+1}} + \sum_{t=0}^k a_t \Exp{\sqrt{A_{t+1}}}.
\end{align*}
Plugging this bound into the upper-bound for $\Exp{\normstar{r^{k+1}}}$, 
\begin{align*}
 \Exp{\normstar{r^{k+1} } } & \leq   \prod_{t=0}^k(1-\eta_t) \Exp{\normstar{r^{0}}} \\
 &\quad + {\color{red}\bar{\rho}}\sqrt{\sum_{t=0}^k \left( \prod_{j=t+1}^k (1-\eta_j)^2\right) (1-\eta_t)^2  B_{t+1}} \\
& \quad + {\color{red}\bar{\rho}}\sum_{t=0}^k \left( \prod_{j=t+1}^k (1-\eta_j)\right) (1-\eta_t)  \Exp{\sqrt{A_{t+1}}} \\
&\quad + {\color{red}\bar{\rho}} \sigma_g  \sqrt{\sum_{t=0}^k \left(\prod_{j=t+1}^k (1-\eta_j)^2 \right) \eta_t^2}.
\end{align*}

\paragraph{Step 5) Deriving the convergence bound under constant tuning parameters.}

If $\gamma_k=\gamma$ and $\eta_k=\eta$, then 
\begin{align*}
 \Exp{\normstar{r^{k+1} } } 
& \leq  (1-\eta)^{k+1}\Exp{\normstar{r^{0}}}  + {\color{red}\bar\rho}\sqrt{ \sum_{t=0}^{k} (1-\eta)^{2(k-t+1)} B }  \\
& \quad + {\color{red}\bar\rho}\sum_{t=0}^{k} (1-\eta)^{k-t+1} \Exp{\sqrt{A_{t+1}}} \\
&\quad + {\color{red}\bar\rho} \sigma_g \sqrt{\sum_{t=0}^{k} (1-\eta)^{2(k-t)} \eta^2},
\end{align*}
where $A_{k+1} = c\gamma^2 \normstar{\nabla f(x^{k})}^2$, $B = \frac{6}{{\color{red} \thesischeightubar \rho}^2}  L_0^2 (\exp(2 L_1 \gamma) +  L_1 \gamma\exp(2 L_1 \gamma) + 2)\gamma^2 + \frac{3 \sigma_H^2}{{\color{blue}\thesischeightubar\theta}^2} \gamma^2$, and $c= \frac{6}{{\color{red} \thesischeightubar \rho}^2}  L_1^2  (\exp(2 L_1 \gamma) + L_1\gamma \exp(2 L_1 \gamma) + 2)$.

Next, since 
\begin{align*}
    \sum_{t=0}^{k} (1-\eta)^{2(k-t+1)} & \leq \sum_{j=0}^\infty ((1-\eta)^2)^j = \frac{1}{1-(1-\eta)^2} \\
    & = \frac{1}{\eta(2-\eta)} \overset{\eta \in [0,1]}{\leq}  \frac{1}{\eta}, \\
    &\text{and}&\\
     \sum_{t=0}^{k} (1-\eta)^{2(k-t)} \eta^2 & \leq  \eta^2 \sum_{j=0}^\infty ((1-\eta)^2)^j = \frac{\eta^2}{1-(1-\eta)^2} \\
     & = \frac{\eta}{2-\eta} \overset{\eta \in [0,1]}{\leq} \eta,
\end{align*}
we obtain
\begin{align*}
    \Exp{\normstar{r^{k+1} } } 
    & \overset{(a)}{\leq}  (1-\eta)^{k+1}\Exp{\normstar{r^{0}}} \\
    &\quad + {\color{red}\bar\rho}\sum_{t=0}^{k} (1-\eta)^{k-t+1} \Exp{\sqrt{A_{t+1}}} \\
    &\quad +  {\color{red}\bar\rho} \frac{\gamma}{\sqrt{\eta}} \hat B \\
    &\quad + {\color{red}\bar\rho} \sqrt{\eta} \sigma_g \\
    & \overset{(b)}{\leq}  (1-\eta)^{k}\Exp{\normstar{r^{0}}} \\
    &\quad + {\color{red}\bar\rho}\hat{c} \gamma \sum_{t=0}^{k} (1-\eta)^{(k-t+1)}  \Exp{\normstar{\nabla f(x^{t})}} \\
    &\quad +  {\color{red}\bar\rho} \frac{\gamma}{\sqrt{\eta}} \hat B \\
    &\quad + {\color{red}\bar\rho} \sqrt{\eta} \sigma_g,
\end{align*}
where we reach $(a)$ by denoting $\hat B = \frac{3}{{\color{red} \thesischeightubar \rho}}  L_0 (\exp( L_1 \gamma) + \sqrt{L_1 \gamma} \exp( L_1 \gamma) + 2) + \frac{2\sigma_H}{{\color{blue}\thesischeightubar\theta}}$, so that $\sqrt{B/\eta} \leq \gamma\hat B/\sqrt{\eta}$, and $(b)$ by using $\eta \in [0,1]$, the definition of $A_{t+1}$, and $\sqrt{c} \leq \hat{c} \eqdef \frac{3}{{\color{red} \thesischeightubar \rho}}  L_1 (\exp( L_1 \gamma) + \sqrt{L_1 \gamma} \exp( L_1 \gamma) + 2)$.

Plugging the above result into the main descent inequality with $\gamma_k=\gamma$ and $\eta_k = \eta$ and denoting $\varphi = \left( 1 - \exp(L_1 \gamma) \frac{L_1\gamma}{2} \right)$, we obtain
\begin{align*}
    \sum_{k=0}^K \gamma\varphi  \Exp{ \normstar{\nabla f(x^{k})} }   
    & \leq  \Delta + 2 \sum_{k=0}^K \gamma \Exp{\normstar{r^{k}}}  + \frac{L_0}{2} \sum_{k=0}^K \exp(L_1 \gamma) \gamma^2 \\
    & \leq \Delta + 2 \gamma \sum_{k=0}^K (1-\eta)^{k}\Exp{\normstar{r^{0}}} + \frac{L_0}{2}  \exp(L_1 \gamma) \gamma^2 (K+1) \\
    &\quad + 2 {\color{red}\bar\rho}\hat{c} \gamma^2 \sum_{k=0}^K \sum_{t=0}^{k-1} (1-\eta)^{(k-t)}  \Exp{ \normstar{\nabla f(x^{t})} } \\
    &\quad + 2 {\color{red}\bar\rho} \frac{\gamma^2}{\sqrt{\eta}} \hat B (K+1) + 2  {\color{red}\bar\rho}  \gamma \sqrt{\eta} \sigma_g (K+1)   .
\end{align*}

By the fact that $\sum_{k=0}^K (1-\eta)^k \leq \sum_{k=0}^\infty (1-\eta)^k = \frac{1}{\eta}$, 
\begin{align*}
    \sum_{k=0}^K \gamma\varphi  \Exp{ \normstar{\nabla f(x^{k})} } & \leq \Delta + \frac{2 \gamma \Exp{\normstar{r^{0}}}}{\eta} \\
    &\quad + 2 {\color{red}\bar\rho}\hat{c} \gamma^2 \sum_{k=0}^K \sum_{t=0}^{k-1} (1-\eta)^{(k-t)}  \Exp{ \normstar{\nabla f(x^{t})} }  \\
    & \quad + 2 {\color{red}\bar\rho} \frac{\gamma^2}{\sqrt{\eta}} \hat B (K+1) + 2  {\color{red}\bar\rho}  \gamma \sqrt{\eta} \sigma_g (K+1) \\
    &\quad + \frac{L_0}{2}  \exp(L_1 \gamma) \gamma^2 (K+1).
\end{align*}

Next, since 
\begin{align*}
    \MoveEqLeft[1] \sum_{k=0}^K \sum_{t=0}^{k-1} (1-\eta)^{k-t} \Exp{\normstar{\nabla f(x^{t})}} \\
    & \leq \sum_{k=0}^K \sum_{t=0}^{k} (1-\eta)^{k-t} \Exp{\normstar{\nabla f(x^{t})}} \\
    & = \sum_{t=0}^K (\sum_{k=t}^K (1-\eta)^{k-t}) \Exp{\normstar{\nabla f(x^{t})}} \\
    & \leq \sum_{t=0}^K (\sum_{k=0}^{\infty} (1-\eta)^{k}) \Exp{\normstar{\nabla f(x^{t})}} \\
    & = \frac{1}{\eta} \sum_{k=0}^K \Exp{\normstar{\nabla f(x^{k})}},
\end{align*}
we obtain
\begin{align*}
    \MoveEqLeft[1] \sum_{k=0}^K \gamma( \varphi - 2 {\color{red}\bar\rho}\hat{c} \frac{\gamma}{\eta}) \Exp{ \normstar{\nabla f(x^{k})} } \\
    & \leq \Delta + \frac{2 \gamma \Exp{\normstar{r^{0}}}}{\eta}     + 2 {\color{red}\bar\rho} \frac{\gamma^2}{\sqrt{\eta}} \hat B (K+1)\\
    &\quad + 2  {\color{red}\bar\rho}  \gamma \sqrt{\eta} \sigma_g (K+1)    + \frac{L_0}{2}  \exp(L_1 \gamma) \gamma^2 (K+1).
\end{align*}

\paragraph{Step 6) Choosing tuning parameters.}

If $\gamma \leq  \frac{\eta}{80 L_1} \left(\frac{{\color{red}\bar\rho}}{{\color{red}\thesischeightubar\rho}}\right)^{-1}$, then $L_1\gamma \leq 1/80$, so $\varphi \geq 0.99$ and $2{\color{red}\bar\rho}\hat c\gamma/\eta \leq 0.24$, i.e., $\varphi - 2{\color{red}\bar\rho}\hat c\gamma/\eta \geq 1/2$; hence 
\begin{align*}
    \frac{\gamma}{2}\sum_{k=0}^K  \Exp{ \normstar{\nabla f(x^{k})} } & \leq \Delta + \frac{2 \gamma \Exp{\normstar{r^{0}}}}{\eta}     + 2 {\color{red}\bar\rho} \frac{\gamma^2}{\sqrt{\eta}} \hat B (K+1) \\
    &\quad + 2  {\color{red}\bar\rho}  \gamma \sqrt{\eta} \sigma_g (K+1)    + \frac{L_0}{2} \exp(L_1 \gamma) \gamma^2 (K+1).
\end{align*}
Therefore, 
\begin{align*}
\min_{k \in \{0,1,\ldots,K\}} \Exp{\normstar{\nabla f(x^{k})}}
& \leq \frac{1}{K+1}\sum_{k=0}^K \Exp{\normstar{\nabla f(x^{k})}}  \\
& \leq \frac{2 \Delta}{\gamma (K+1)} +  \frac{4 \Exp{\normstar{r^{0}}}}{\eta (K+1)}  + 4 {\color{red}\bar\rho} \frac{\gamma}{\sqrt{\eta}} \hat B \\
&\quad + 4 {\color{red}\bar\rho}   \sqrt{\eta} \sigma_g     + L_0 \exp(L_1 \gamma) \gamma.
\end{align*}

If $\gamma = \frac{\hat \gamma}{(K+1)^{2/3}}$ with $\hat \gamma =\frac{1}{80 L_1} \left(\frac{{\color{red}\bar\rho}}{{\color{red}\thesischeightubar\rho}}\right)^{-1}$, and $\eta = \frac{1}{(K+1)^{2/3}}$, then 
\begin{align*}
\MoveEqLeft[1] \min_{k \in \{0,1,\ldots,K\}} \Exp{\normstar{\nabla f(x^{k})}} \\
& \leq \frac{2 \Delta}{\hat \gamma (K+1)^{1/3}} +  \frac{4 \Exp{\normstar{r^{0}}}}{(K+1)^{1/3}} \\
&\quad +  \frac{4 {\color{red}\bar\rho} \hat\gamma \hat B}{(K+1)^{1/3}}  \\
&\quad + 4 {\color{red}\bar\rho}   \frac{1}{(K+1)^{1/3}} \sigma_g     + L_0 \exp(L_1 \hat\gamma) \frac{\hat \gamma}{(K+1)^{2/3}}.
\end{align*}

\newpage

\section{Proof of~Theorem~\ref{ch8:thm:V2_relaxed_smoothness_chapter_8}}

We prove the result in the following steps. 

\paragraph{Step 1) Proving the descent inequality.} 
By following the proof arguments in Step 1) of~Theorem~\ref{ch8:thm:V1_relaxed_smoothness_chapter_8} (see \eqref{ch8:eq:descent_inequality_chapter_8}), we have 
\begin{align*}
 \sum_{k=0}^K \gamma_k\varphi_k  \normstar{\nabla f(x^{k})}    & \leq  \Delta + 2 \sum_{k=0}^K \gamma_k \normstar{\nabla f(x^{k})-m^{k}} \\
 &\quad + \frac{L_0}{2} \sum_{k=0}^K \exp(L_1 \gamma_k) \gamma_k^2,
\end{align*}
where $\varphi_k \eqdef \left( 1 - \exp(L_1 \gamma_k) 
 L_1\gamma_k/2 \right)$ and $\Delta \eqdef f(x^{0}) - f^{\inf}$.

\paragraph{Step 2) Bounding the error term.}
Next, we bound $\normstar{r^{k}}$, where $r^{k} \eqdef  m^{k} - \nabla f(x^{k})$. 
From the definition of $r^{k}$,
\begin{align*}
    r^{k+1} 
    & = m^{k+1} - \nabla f(x^{k+1}) \\
    & = (1-\eta_k)r^{k} + (1-\eta_k)Z_f(x^{k},x^{k+1}) + (1-\eta_k) \zeta^H_{k+1} \\
    &\quad + \eta_k \zeta^{g}_{k+1} ,
\end{align*}
where
\begin{align*}
Z_f(x^{k},x^{k+1}) &= \nabla f(x^{k}) - \nabla f(x^{k+1}) + \nabla^2 f(x^{k+1})(x^{k+1}-x^{k}),\\
\zeta^H_{k+1} &= (\nabla^2 f(x^{k+1};\xi^{k+1})-\nabla^2 f(x^{k+1}))(x^{k+1}-x^{k}),\\
\zeta^{g}_{k+1} &= \nabla f(x^{k+1};\xi^{k+1}) - \nabla f(x^{k+1}).
\end{align*}

By recursively applying the above identity, 
\begin{align*}
    r^{k+1} 
    & = \prod_{t=0}^k (1-\eta_t) r^{0} \\
    &\quad +\sum_{t=0}^k \left( \prod_{j=t+1}^k (1-\eta_j) \right) (1-\eta_t)  Z_f(x^{t},x^{t+1}) \\
    &\quad +\sum_{t=0}^k \left( \prod_{j=t+1}^k (1-\eta_j) \right) (1-\eta_t)  \zeta^H_{t+1}+  \sum_{t=0}^k \left(\prod_{j=t+1}^k (1-\eta_j) \right) \eta_t \zeta^{g}_{t+1}.
\end{align*}
Set
\begin{align*}
\circledThree
&\eqdef
\Exp{\normstar{\sum_{t=0}^k
\left(\prod_{j=t+1}^k(1-\eta_j)\right)
(1-\eta_t)Z_f(x^t,x^{t+1})}},\\
\circledFour
&\eqdef
\Exp{\normstar{\sum_{t=0}^k
\left(\prod_{j=t+1}^k(1-\eta_j)\right)
(1-\eta_t)\zeta^H_{t+1}}},\\
\circledFive
&\eqdef
\Exp{\normstar{\sum_{t=0}^k
\left(\prod_{j=t+1}^k(1-\eta_j)\right)
\eta_t\zeta^g_{t+1}}}.
\end{align*}
Therefore,
\begin{align*}
\Exp{\normstar{r^{k+1}}}
&\leq \Exp{\normstar{\prod_{t=0}^k(1-\eta_t)r^0}}
 +\circledThree+\circledFour+\circledFive.
\end{align*}

Next, we bound $\circledFive$ by using~Assumption~\ref{ch8:assum:Hessian_boundedvariance_v2_chapter_8} and by following the proof arguments for bounding $\circledTwo$ in Step 2) of the convergence proof for~Theorem~\ref{ch8:thm:V1_relaxed_smoothness_chapter_8}. The conditional centering of $\zeta^g_{t+1}$ established above makes all cross-time inner products vanish by the tower property. Then, we have 
\begin{align*}
    \circledFive
    & \leq {\color{red}\bar\rho}   \Exp{\normtwo{\sum_{t=0}^k \left(\prod_{j=t+1}^k (1-\eta_j) \right) \eta_t \zeta^{g}_{t+1}}} \\
    & \leq {\color{red}\bar\rho}  \sqrt{  \Exp{\normtwo{\sum_{t=0}^k \left(\prod_{j=t+1}^k (1-\eta_j) \right) \eta_t \zeta^{g}_{t+1}}^2} } \\
    & \leq {\color{red}\bar\rho} \sigma_g \sqrt{ \sum_{t=0}^k \left(\prod_{j=t+1}^k (1-\eta_j)^2 \right) \eta_t^2}, 
\end{align*}
Therefore, 
\begin{align*}
    \Exp{\normstar{r^{k+1} }}
    & \leq \prod_{t=0}^k (1-\eta_t)\Exp{\normstar{r^{0}}} + ~\circledThree~ +~\circledFour~ \\
    &\quad +   {\color{red}\bar\rho} \sigma_g \sqrt{ \sum_{t=0}^k \left(\prod_{j=t+1}^k (1-\eta_j)^2 \right) \eta_t^2}.
\end{align*}

Next, we bound $\circledFour$. 
\begin{align*}
    \circledFour  
    & \leq {\color{red}\bar\rho} \Exp{\normtwo{\sum_{t=0}^k \left( \prod_{j=t+1}^k (1-\eta_j) \right) (1-\eta_t)  \zeta^H_{t+1} }} \\
    & \leq {\color{red}\bar\rho}  \sqrt{ \Exp{\normtwo{\sum_{t=0}^k \left( \prod_{j=t+1}^k (1-\eta_j) \right) (1-\eta_t)  \zeta^H_{t+1} }^2} }.
\end{align*}
Since $x^{t+1}-x^t$ is $\mathcal F_t$-measurable,
Assumption~\ref{ch8:assum:Hessian_boundedvariance_v2_chapter_8} gives
\begin{align*}
\Exp{\zeta^H_{t+1}\mid\mathcal F_t}&=0,\\
\Exp{\normtwo{\zeta^H_{t+1}}^2\mid\mathcal F_t}
&\leq\sigma_H^2\normtwo{x^{t+1}-x^t}^2.
\end{align*}
Consequently, for $0\leq i<j$,
$\Exp{\inner{\zeta^H_{j+1}}{\zeta^H_{i+1}}}=0$ by the tower property, and
\begin{align*}
    \circledFour
  & \leq {\color{red}\bar\rho}  \sqrt{  \sum_{t=0}^k \left( \prod_{j=t+1}^k (1-\eta_j)^2 \right) (1-\eta_t)^2 \Exp{ \normtwo{  \zeta^H_{t+1} }^2} } \\
   & \leq {\color{red}\bar\rho}  \sqrt{  \sum_{t=0}^k \left( \prod_{j=t+1}^k (1-\eta_j)^2 \right) (1-\eta_t)^2 \sigma_H^2 \normtwo{x^{t+1}-x^{t}}^2 } \\
   & \leq \frac{{\color{red}\bar\rho}}{\color{blue} \thesischeightubar{\theta}}  \sqrt{  \sum_{t=0}^k \left( \prod_{j=t+1}^k (1-\eta_j)^2 \right) (1-\eta_t)^2 \sigma_H^2 \norm{x^{t+1}-x^{t}}^2 } \\
   & \leq \frac{{\color{red}\bar\rho}}{\color{blue} \thesischeightubar{\theta}} \gamma \sigma_H  \sqrt{  \sum_{t=0}^k \left( \prod_{j=t+1}^k (1-\eta_j)^2 \right) (1-\eta_t)^2  }.
\end{align*}
Therefore, 
\begin{align*}
    \Exp{\normstar{r^{k+1} }}
    & \leq \prod_{t=0}^k (1-\eta_t)\Exp{\normstar{ r^{0}}} \\
    &\quad + \frac{{\color{red}\bar\rho}}{\color{blue} \thesischeightubar{\theta}} \gamma \sigma_H  \sqrt{  \sum_{t=0}^k \left( \prod_{j=t+1}^k (1-\eta_j)^2 \right) (1-\eta_t)^2  } \\
    &\quad + ~\circledThree~ +   {\color{red}\bar\rho} \sigma_g \sqrt{ \sum_{t=0}^k \left(\prod_{j=t+1}^k (1-\eta_j)^2 \right) \eta_t^2}.
\end{align*}

Next, we bound $\circledThree$: 
\begin{align*}
    \circledThree & \leq  \sum_{t=0}^k \left( \prod_{j=t+1}^k (1-\eta_j) \right) (1-\eta_t) \Exp{ \normstar{ Z_f(x^{t},x^{t+1}) } }.
\end{align*}

From~Lemma~\ref{ch8:lemma:Relaxed_Grad_and_Hessian_chapter_8},
\begin{align*}
\MoveEqLeft \normstar{Z_f(x^t,x^{t+1})}\\
&\leq \frac{1}{2}\left(M_0+M_1\normstar{\nabla f(x^{t+1})}\right)
       \norm{x^t-x^{t+1}}^2\\
&\quad+\frac{1}{3}M_1
       \left(L_0+L_1\normstar{\nabla f(x^{t+1})}\right)
       \exp\left(L_1\norm{x^t-x^{t+1}}\right)
       \norm{x^t-x^{t+1}}^3\\
&\leq \frac{1}{2}\left(M_0+M_1\normstar{\nabla f(x^{t+1})}\right)\gamma_t^2\\
&\quad+\frac{1}{3}M_1
       \left(L_0+L_1\normstar{\nabla f(x^{t+1})}\right)
       \exp(L_1\gamma_t)\gamma_t^3\\
&=\left(\frac{M_0}{2}+\frac{M_1}{3}L_0\gamma_t\exp(L_1\gamma_t)\right)\gamma_t^2\\
&\quad+\left(\frac{M_1}{2}+\frac{M_1}{3}L_1\gamma_t\exp(L_1\gamma_t)\right)
       \gamma_t^2\normstar{\nabla f(x^{t+1})}.
\end{align*}
Therefore, 
\begin{align*}
\circledThree
&\leq \sum_{t=0}^k
\left(\prod_{j=t+1}^k(1-\eta_j)\right)(1-\eta_t)
\left(\frac{M_0}{2}+\frac{M_1}{3}L_0\gamma_t\exp(L_1\gamma_t)\right)\gamma_t^2\\
&\quad+\sum_{t=0}^k
\left(\prod_{j=t+1}^k(1-\eta_j)\right)(1-\eta_t)
\left(\frac{M_1}{2}+\frac{M_1}{3}L_1\gamma_t\exp(L_1\gamma_t)\right)\\
&\qquad\qquad{}\times\gamma_t^2\Exp{\normstar{\nabla f(x^{t+1})}}.
\end{align*}
Plugging $\circledThree$ into the upper-bound for $\Exp{\normstar{r^{k+1} }}$ yields
\begin{align*}
\Exp{\normstar{r^{k+1}}}
&\leq \prod_{t=0}^k(1-\eta_t)\Exp{\normstar{r^0}}\\
&\quad+\sum_{t=0}^k
\left(\prod_{j=t+1}^k(1-\eta_j)\right)(1-\eta_t)
\left(\frac{M_0}{2}+\frac{M_1}{3}L_0\gamma_t\exp(L_1\gamma_t)\right)\gamma_t^2\\
&\quad+\sum_{t=0}^k
\left(\prod_{j=t+1}^k(1-\eta_j)\right)(1-\eta_t)
\left(\frac{M_1}{2}+\frac{M_1}{3}L_1\gamma_t\exp(L_1\gamma_t)\right)\\
&\qquad\qquad{}\times\gamma_t^2\Exp{\normstar{\nabla f(x^{t+1})}}\\
&\quad+\frac{{\color{red}\bar\rho}}{{\color{blue}\thesischeightubar{\theta}}}
\gamma\sigma_H
\sqrt{\sum_{t=0}^k\left(\prod_{j=t+1}^k(1-\eta_j)^2\right)(1-\eta_t)^2}\\
&\quad+{\color{red}\bar\rho}\sigma_g
\sqrt{\sum_{t=0}^k\left(\prod_{j=t+1}^k(1-\eta_j)^2\right)\eta_t^2}.
\end{align*}

If $\gamma_k = \gamma$ and $\eta_k = \eta$, then 
\begin{align*}
    \Exp{\normstar{r^{k+1} }}
    & \leq (1-\eta)^{k+1} \Exp{\normstar{r^{0}}} \\
    &\quad + \left(\frac{M_0}{2} + \frac{M_1L_0\gamma \exp(L_1\gamma)}{3}\right) \sum_{t=0}^k (1-\eta)^{k-t+1} \gamma^2   \\
    &\quad + \left(\frac{M_1}{2} + \frac{M_1L_1 \exp(L_1 \gamma) \gamma}{3}\right) \sum_{t=0}^k (1-\eta)^{k-t+1}  \gamma^2 \Exp{\normstar{\nabla f(x^{t+1})}} \\
    & \quad + \frac{{\color{red}\bar\rho}}{\color{blue} \thesischeightubar{\theta}} \gamma \sigma_H  \sqrt{  \sum_{t=0}^k (1-\eta)^{2(k-t+1)}  } \\
    &\quad +   {\color{red}\bar\rho} \sigma_g \sqrt{ \sum_{t=0}^k (1-\eta)^{2(k-t)} \eta^2}.
\end{align*}

Next, since 
\begin{align*}
    \sum_{t=0}^{k} (1-\eta)^{2(k-t+1)}
    & \leq \sum_{j=0}^\infty ((1-\eta)^2)^j  = \frac{1}{\eta(2-\eta)}  \overset{\eta \in [0,1]}{\leq} \frac{1}{\eta}, \\
    &\text{and}&\\
     \sum_{t=0}^{k} (1-\eta)^{2(k-t)} \eta^2 & \leq  \eta^2 \sum_{j=0}^\infty ((1-\eta)^2)^j = \frac{\eta^2}{1-(1-\eta)^2} \\
     & = \frac{\eta}{2-\eta} \overset{\eta \in [0,1]}{\leq} \eta,
\end{align*}
we obtain
\begin{align*}
    \Exp{\normstar{r^{k+1} }}
    & \leq (1-\eta)^{k+1} \Exp{\normstar{r^{0}}} + \left(\frac{M_0}{2} + \frac{M_1}{3}L_0\gamma \exp(L_1\gamma)\right) \frac{\gamma^2}{\eta}   \\
    &\quad + \left(\frac{M_1}{2} + \frac{M_1}{3}L_1  \gamma\exp(L_1 \gamma)\right) \sum_{t=0}^k (1-\eta)^{k-t+1}  \gamma^2 \Exp{\normstar{\nabla f(x^{t+1})}} \\
    &\quad + \frac{{\color{red}\bar\rho}}{\color{blue} \thesischeightubar{\theta}} \frac{\gamma}{\sqrt{\eta}} \sigma_H  +   {\color{red}\bar\rho} \sqrt{\eta}\sigma_g.
\end{align*}

Therefore, recalling $\varphi = 1 - \exp(L_1 \gamma) 
 L_1\gamma/2$, we have 
\begin{align*}
\sum_{k=0}^K\gamma\varphi\Exp{\normstar{\nabla f(x^k)}}
&\leq \Delta+\frac{L_0}{2}\sum_{k=0}^K\exp(L_1\gamma)\gamma^2
 +2\gamma\sum_{k=0}^K\Exp{\normstar{r^k}}\\
&\leq \Delta+\frac{L_0}{2}\exp(L_1\gamma)\gamma^2(K+1)\\
&\quad+2\gamma\sum_{k=0}^K(1-\eta)^k\Exp{\normstar{r^0}}\\
&\quad+\gamma\left(M_0+\frac{2M_1}{3}L_0\gamma\exp(L_1\gamma)\right)
\frac{\gamma^2}{\eta}(K+1)\\
&\quad+\gamma\left(M_1+\frac{2M_1}{3}L_1\gamma\exp(L_1\gamma)\right)\\
&\qquad\qquad{}\times
\sum_{k=0}^K\sum_{t=0}^{k-1}(1-\eta)^{k-t}
\gamma^2\Exp{\normstar{\nabla f(x^{t+1})}}\\
&\quad+2\gamma
\left(
\frac{{\color{red}\bar\rho}}{{\color{blue}\thesischeightubar{\theta}}}
\frac{\gamma}{\sqrt{\eta}}\sigma_H
+{\color{red}\bar\rho}\sqrt{\eta}\sigma_g
\right)(K+1).
\end{align*}
where $\Delta = \Exp{f(x^{0})-f^{\inf}}$.

Next, since 
\begin{gather*}
    \sum_{k=0}^K (1-\eta)^k  \leq \frac{1}{\eta} \quad\text{and}\\
    \sum_{k=0}^K \sum_{t=0}^{k-1} (1-\eta)^{k-t} \Exp{\normstar{\nabla f(x^{t+1})}}
    \leq \frac{1}{\eta} \sum_{k=0}^K \Exp{\normstar{\nabla f(x^{k})}},
\end{gather*}
we obtain
\begin{align*}
    \sum_{k=0}^K \gamma\vartheta  \Exp{\normstar{\nabla f(x^{k})}}
    & \leq \Delta  + (K+1)\frac{L_0}{2}  \exp(L_1 \gamma) \gamma^2  + 2  \frac{\gamma}{\eta}  \Exp{\normstar{r^{0}}}\\
    &\quad + \gamma (K+1) \left(M_0 + \frac{2M_1}{3}L_0 \gamma\exp(L_1\gamma)\right) \frac{\gamma^2}{\eta}  \\
    &\quad + 2\gamma\left(\frac{{\color{red}\bar\rho}}{\color{blue} \thesischeightubar{\theta}} \frac{\gamma}{\sqrt{\eta}} \sigma_H  +   {\color{red}\bar\rho} \sqrt{\eta}\sigma_g\right)(K+1),
\end{align*}
where $\vartheta \eqdef \left( 1 - \exp(L_1 \gamma) 
 L_1\gamma/2   -  \left(M_1 + 2M_1/3 \cdot L_1 \gamma \exp(L_1 \gamma) \right) \frac{\gamma^2}{\eta}  \right)$. 
 
 If $\gamma \leq \frac{\eta}{3}\min\left\{ \frac{1}{L_1} , \frac{1}{\sqrt{M_1 } }  \right\} $, then $L_1\gamma \leq 1/3$ and $M_1\gamma^2/\eta \leq \eta/9 \leq 1/9$, so $\vartheta \geq 1 - \exp(1/3)/6 - 1.31/9 \geq 1/2$, and 
\begin{align*}
    \frac{\gamma}{2}\sum_{k=0}^K   \Exp{\normstar{\nabla f(x^{k})}}   
    & \leq \Delta  + (K+1)\frac{L_0}{2}  \exp(L_1 \gamma) \gamma^2  + 2  \frac{\gamma}{\eta}  \Exp{\normstar{r^{0}}} \\
    &\quad+ \gamma (K+1) \left(M_0 + \frac{2M_1}{3}L_0\gamma \exp(L_1\gamma)\right) \frac{\gamma^2}{\eta}  \\
    &\quad + 2\gamma\left(\frac{{\color{red}\bar\rho}}{\color{blue} \thesischeightubar{\theta}} \frac{\gamma}{\sqrt{\eta}} \sigma_H  +   {\color{red}\bar\rho} \sqrt{\eta}\sigma_g\right)(K+1).
\end{align*}
Therefore, 
\begin{align*}
    \frac{1}{K+1}\sum_{k=0}^K   \Exp{\normstar{\nabla f(x^{k})}} 
    & \leq \frac{2\Delta}{\gamma(K+1)}  +L_0 \exp(L_1 \gamma) \gamma  +  \frac{4\Exp{\normstar{r^{0}}}}{\eta (K+1)}  \\
    &\quad + 2 \left(M_0 + \frac{2M_1}{3}L_0\gamma \exp(L_1\gamma)\right) \frac{\gamma^2}{\eta}  \\
    &\quad + 4\left(\frac{{\color{red}\bar\rho}}{\color{blue} \thesischeightubar{\theta}} \frac{\gamma}{\sqrt{\eta}} \sigma_H  +   {\color{red}\bar\rho} \sqrt{\eta}\sigma_g\right).
\end{align*}

If $\gamma = \hat \gamma \eta$ with $\hat \gamma = \frac{1}{3}\min\left\{ \frac{1}{L_1} , \frac{1}{\sqrt{M_1}}  \right\} $ and $\eta = \frac{1}{(K+1)^{2/3}}$, then 
\begin{align*}
    \frac{1}{K+1}\sum_{k=0}^K   \Exp{\normstar{\nabla f(x^{k})}}   
    & \leq \frac{2\Delta}{\hat\gamma(K+1)^{1/3}} \\
    &\quad + \frac{\hat\gamma L_0 \exp(L_1 \hat \gamma)}{(K+1)^{2/3}} + \frac{4 \Exp{\normstar{r^{0}}}}{ (K+1)^{1/3}} \\
    &\quad + 2\left(M_0 + \frac{2M_1}{3} L_0\hat\gamma \exp(L_1\hat \gamma)\right) \frac{\hat \gamma^2}{(K+1)^{2/3}}  \\
    & \quad + 4\frac{{\color{red}\bar\rho}}{\color{blue} \thesischeightubar{\theta}} \frac{\hat \gamma \sigma_H }{(K+1)^{1/3}} \\
    &\quad +   4\frac{{\color{red}\bar\rho} \sigma_g}{(K+1)^{1/3}}.
\end{align*}

Finally, by the fact that
\begin{equation*}
\min_{k \in \{0,1,\ldots,K\}} \Exp{\normstar{\nabla f(x^{k})}}
\leq \frac{1}{K+1}\sum_{k=0}^K \Exp{\normstar{\nabla f(x^{k})}},
\end{equation*}
we obtain the final result.

\newpage
\section{LMO-based Methods with Extrapolated Momentum}\label{ch8:app:IGT_relaxed_smoothness_chapter_8}

In this section, we present the convergence of LMO-based methods in~\eqref{ch8:eqn:LMO_momentum_chapter_8} that leverage extrapolated momentum in~\eqref{ch8:eqn:extrapolated_m_chapter_8}. 

\begin{theorem}
    \label{ch8:thm:IGT_relaxed_smoothness_chapter_8}
    Consider the problem of minimizing $f(x)=\thesischeightExpSub{\xi \sim \thesischeightcD}{f(x;\xi)}$. Let $f$ be twice differentiable, and let   Assumptions~\ref{ch8:assum:Hessian_boundedvariance_v2_chapter_8},~\ref{ch8:assum:lower_bound_chapter_8},~\ref{ch8:assum:Relaxed_Lipschitz_Grad_chapter_8}, and~\ref{ch8:assum:Relaxed_Lipschitz_Hessian_chapter_8} hold. 
    Then, the iterates $\{x^{k}\}$ generated by LMO-based methods in~\eqref{ch8:eqn:LMO_momentum_chapter_8} that leverage extrapolated momentum in~\eqref{ch8:eqn:extrapolated_m_chapter_8} with
    \begin{equation*}
        \eta_k =  \eta = \frac{1}{(K+1)^{4/7}}, \quad \text{and} \quad  \gamma_k= \gamma = \frac{\hat \gamma}{(K+1)^{5/7}},
    \end{equation*}
    where $\hat \gamma = \min\left\{\frac{1}{3L_1}, \frac{1}{3\sqrt{M_1}}\right\}$  satisfy  
    \begin{align*}
        \MoveEqLeft[1] \min_{k \in \{0,1,\ldots,K\}} \Exp{\normstar{\nabla f(x^{k})}} \\
        & \leq  \frac{2(f(x^{0})-f^{\inf})}{\hat\gamma(K+1)^{2/7}} + \frac{\hat\gamma L_0 \exp(L_1 \hat \gamma)}{(K+1)^{5/7}} + \frac{4 \left(\Exp{\normstar{r^{0}}}+ {\color{red}\bar\rho} \sigma_g\right)}{ (K+1)^{2/7}} \\
        &\quad + 2\left(M_0 + \frac{2M_1}{3} L_0\hat \gamma \exp(L_1\hat \gamma)\right) \frac{\hat \gamma^2}{(K+1)^{2/7}}.
    \end{align*}

\end{theorem}

Theorem~\ref{ch8:thm:IGT_relaxed_smoothness_chapter_8} establishes the $\cO(1/K^{2/7})$ convergence of LMO-based methods using extrapolated momentum under relaxed smoothness with respect to the arbitrary norm, which matches the known rate obtained by \textcite[Corollary 5]{kovalev2025understanding} under traditional smoothness with respect to the arbitrary norm. 

\subsection{Proof of Theorem~\ref{ch8:thm:IGT_relaxed_smoothness_chapter_8}}

We prove the result in the following steps. 

\paragraph{Step 1) Proving the descent inequality.} 
By following the proof arguments in Step 1) of~Theorem~\ref{ch8:thm:V1_relaxed_smoothness_chapter_8} (see \eqref{ch8:eq:descent_inequality_chapter_8}), we have 
\begin{align*}
    \sum_{k=0}^K \gamma_k\varphi_k  \normstar{\nabla f(x^{k})}    & \leq  \Delta + 2 \sum_{k=0}^K \gamma_k \normstar{\nabla f(x^{k})-m^{k}} \\
    &\quad + \frac{L_0}{2} \sum_{k=0}^K \exp(L_1 \gamma_k) \gamma_k^2,
\end{align*}
where $\varphi_k \eqdef \left( 1 - \exp(L_1 \gamma_k) 
 L_1\gamma_k/2 \right)$ and $\Delta \eqdef f(x^{0}) - f^{\inf}$.

\paragraph{Step 2) Bounding the error term.}
Next, we bound $\normstar{r^{k}}$, where $r^{k} \eqdef  m^{k} - \nabla f(x^{k})$. 
From the definition of $r^{k}$,
\begin{align*}
    r^{k+1} 
    & = m^{k+1} - \nabla f(x^{k+1}) \\
    & = (1-\eta_k)r^{k} \\
    &\quad +\eta_k\left(\nabla f(y^{k+1}) - \nabla f(x^{k+1}) -\nabla^2 f(x^{k+1})(y^{k+1}-x^{k+1})\right)\\
    &\quad +(1-\eta_k)\left(\nabla f(x^{k}) - \nabla f(x^{k+1}) +\frac{\eta_k}{1-\eta_k}\nabla^2 f(x^{k+1})(y^{k+1}-x^{k+1})\right)\\
    &\quad + \eta_k\left(\nabla f(y^{k+1};\xi^{k+1}) - \nabla f(y^{k+1})\right)\\
    & = (1-\eta_k)r^{k} + (1-\eta_k)Z_f(x^{k},x^{k+1}) + \eta_k Z_f(y^{k+1},x^{k+1}) \\
    &\quad + \eta_k \zeta^{g}_{k+1} ,
\end{align*}
where $Z_f$ is defined in Lemma~\ref{ch8:lemma:Relaxed_Grad_and_Hessian_chapter_8}, $\zeta^{g}_{k+1} = \nabla f(y^{k+1};\xi^{k+1}) - \nabla f(y^{k+1})$, and we used $\eta_k(y^{k+1}-x^{k+1}) = (1-\eta_k)(x^{k+1}-x^k)$. 

By recursively applying the above identity, 
\begin{align*}
    r^{k+1} 
    & = \prod_{t=0}^k (1-\eta_t) r^{0} \\
    &\quad +\sum_{t=0}^k \left( \prod_{j=t+1}^k (1-\eta_j) \right) (1-\eta_t)  Z_f(x^{t},x^{t+1}) \\
    &\quad +\sum_{t=0}^k \left( \prod_{j=t+1}^k (1-\eta_j) \right) \eta_t Z_f(y^{t+1},x^{t+1})+  \sum_{t=0}^k \left(\prod_{j=t+1}^k (1-\eta_j) \right) \eta_t \zeta^{g}_{t+1}.
\end{align*}
Set
\begin{align*}
\circledSix
&\eqdef
\Exp{\normstar{\sum_{t=0}^k
\left(\prod_{j=t+1}^k(1-\eta_j)\right)
(1-\eta_t)Z_f(x^t,x^{t+1})}},\\
\circledSeven
&\eqdef
\Exp{\normstar{\sum_{t=0}^k
\left(\prod_{j=t+1}^k(1-\eta_j)\right)
\eta_t Z_f(y^{t+1},x^{t+1})}},\\
\circledEight
&\eqdef
\Exp{\normstar{\sum_{t=0}^k
\left(\prod_{j=t+1}^k(1-\eta_j)\right)
\eta_t\zeta^g_{t+1}}}.
\end{align*}
Therefore,
\begin{align*}
\Exp{\normstar{r^{k+1}}}
&\leq \Exp{\normstar{\prod_{t=0}^k(1-\eta_t)r^0}}
+\circledSix+\circledSeven+\circledEight.
\end{align*}

Next, we bound $\circledEight$ by using~Assumption~\ref{ch8:assum:Hessian_boundedvariance_v2_chapter_8} and by following the proof arguments for bounding $\circledTwo$ in Step 2) of the convergence proof for~Theorem~\ref{ch8:thm:V1_relaxed_smoothness_chapter_8}. Since $y^{t+1}$ is $\mathcal F_t$-measurable, the conditional centering of $\zeta^g_{t+1}$ makes all cross-time inner products vanish by the tower property. Then, we have 
\begin{align*}
    \circledEight
    & \leq {\color{red}\bar\rho}   \Exp{\normtwo{\sum_{t=0}^k \left(\prod_{j=t+1}^k (1-\eta_j) \right) \eta_t \zeta^{g}_{t+1}}} \\
    & \leq {\color{red}\bar\rho}  \sqrt{  \Exp{\normtwo{\sum_{t=0}^k \left(\prod_{j=t+1}^k (1-\eta_j) \right) \eta_t \zeta^{g}_{t+1}}^2} } \\
    & \leq {\color{red}\bar\rho} \sigma_g \sqrt{ \sum_{t=0}^k \left(\prod_{j=t+1}^k (1-\eta_j)^2 \right) \eta_t^2}, 
\end{align*}
Therefore, 
\begin{align*}
    \Exp{\normstar{r^{k+1} }}
    & \leq \prod_{t=0}^k (1-\eta_t)\Exp{\normstar{r^{0}}} + ~\circledSix~ +~\circledSeven~ \\
    &\quad +   {\color{red}\bar\rho} \sigma_g \sqrt{ \sum_{t=0}^k \left(\prod_{j=t+1}^k (1-\eta_j)^2 \right) \eta_t^2}.
\end{align*}

Next, we bound $\circledSix$ and $\circledSeven$: 
\begin{align*}
    \circledSix & \leq  \sum_{t=0}^k \left( \prod_{j=t+1}^k (1-\eta_j) \right) (1-\eta_t) \Exp{ \normstar{ Z_f(x^{t},x^{t+1}) } };\\
    \circledSeven & \leq  \sum_{t=0}^k \left( \prod_{j=t+1}^k (1-\eta_j) \right) \eta_t \Exp{ \normstar{ Z_f(y^{t+1},x^{t+1}) } }
\end{align*}
From Lemma~\ref{ch8:lemma:Relaxed_Grad_and_Hessian_chapter_8}, $\norm{x^{t+1}-x^t}\leq\gamma_t$ and $\norm{y^{t+1}-x^{t+1}} = \theta_t\norm{x^{t+1}-x^t} \leq \theta_t\gamma_t$ with $\theta_t \eqdef (1-\eta_t)/\eta_t$,
\begin{align*}
\MoveEqLeft \normstar{Z_f(x^t,x^{t+1})}\\
&\leq \frac{1}{2}\left(M_0+M_1\normstar{\nabla f(x^{t+1})}\right)
\norm{x^t-x^{t+1}}^2\\
&\quad+\frac{1}{3}M_1
\left(L_0+L_1\normstar{\nabla f(x^{t+1})}\right)
\exp\left(L_1\norm{x^t-x^{t+1}}\right)
\norm{x^t-x^{t+1}}^3\\
&\leq \frac{1}{2}\left(M_0+M_1\normstar{\nabla f(x^{t+1})}\right)\gamma_t^2\\
&\quad+\frac{1}{3}M_1
\left(L_0+L_1\normstar{\nabla f(x^{t+1})}\right)
\exp(L_1\gamma_t)\gamma_t^3\\
&=\left(\frac{M_0}{2}+\frac{M_1}{3}L_0\gamma_t\exp(L_1\gamma_t)\right)\gamma_t^2\\
&\quad+\left(\frac{M_1}{2}+\frac{M_1}{3}L_1\gamma_t\exp(L_1\gamma_t)\right)
\gamma_t^2\normstar{\nabla f(x^{t+1})};
\end{align*}
and
\begin{align*}
\MoveEqLeft \normstar{Z_f(y^{t+1},x^{t+1})}\\
&\leq \frac{1}{2}\left(M_0+M_1\normstar{\nabla f(x^{t+1})}\right)
\norm{y^{t+1}-x^{t+1}}^2\\
&\quad+\frac{1}{3}M_1
\left(L_0+L_1\normstar{\nabla f(x^{t+1})}\right)
\exp\left(L_1\norm{y^{t+1}-x^{t+1}}\right)
\norm{y^{t+1}-x^{t+1}}^3\\
&\leq \frac{1}{2}\left(M_0+M_1\normstar{\nabla f(x^{t+1})}\right)
\theta_t^2\gamma_t^2\\
&\quad+\frac{1}{3}M_1
\left(L_0+L_1\normstar{\nabla f(x^{t+1})}\right)
\exp(L_1\theta_t\gamma_t)\theta_t^3\gamma_t^3\\
&=\left(\frac{M_0}{2}
+\frac{M_1}{3}L_0\theta_t\gamma_t\exp(L_1\theta_t\gamma_t)\right)
\theta_t^2\gamma_t^2\\
&\quad+\left(\frac{M_1}{2}
+\frac{M_1}{3}L_1\theta_t\gamma_t\exp(L_1\theta_t\gamma_t)\right)
\theta_t^2\gamma_t^2\normstar{\nabla f(x^{t+1})},
\end{align*}
we obtain, using $\theta_t \leq 1/\eta_t$ and $(1-\eta_t) +\eta_t \theta_t^2 = (1-\eta_t)/\eta_t = \theta_t$,
\begin{align*}
\circledSix+\circledSeven
&\leq \sum_{t=0}^k
\left(\prod_{j=t+1}^k(1-\eta_j)\right)
\left(\frac{M_0}{2}
+\frac{M_1}{3}\frac{L_0\gamma_t}{\eta_t}
\exp\left(\frac{L_1\gamma_t}{\eta_t}\right)\right)
\theta_t\gamma_t^2\\
&\quad+\sum_{t=0}^k
\left(\prod_{j=t+1}^k(1-\eta_j)\right)
\left(\frac{M_1}{2}
+\frac{M_1}{3}\frac{L_1\gamma_t}{\eta_t}
\exp\left(\frac{L_1\gamma_t}{\eta_t}\right)\right)\\
&\qquad\qquad{}\times
\theta_t\gamma_t^2\Exp{\normstar{\nabla f(x^{t+1})}}.
\end{align*}

Therefore, 
\begin{align*}
\Exp{\normstar{r^{k+1}}}
&\leq \prod_{t=0}^k(1-\eta_t)\Exp{\normstar{r^0}}\\
&\quad+{\color{red}\bar\rho}\sigma_g
\sqrt{\sum_{t=0}^k
\left(\prod_{j=t+1}^k(1-\eta_j)^2\right)\eta_t^2}\\
&\quad+\sum_{t=0}^k
\left(\prod_{j=t+1}^k(1-\eta_j)\right)
\left(\frac{M_0}{2}
+\frac{M_1}{3}\frac{L_0\gamma_t}{\eta_t}
\exp\left(\frac{L_1\gamma_t}{\eta_t}\right)\right)
\theta_t\gamma_t^2\\
&\quad+\sum_{t=0}^k
\left(\prod_{j=t+1}^k(1-\eta_j)\right)
\left(\frac{M_1}{2}
+\frac{M_1}{3}\frac{L_1\gamma_t}{\eta_t}
\exp\left(\frac{L_1\gamma_t}{\eta_t}\right)\right)\\
&\qquad\qquad{}\times
\theta_t\gamma_t^2\Exp{\normstar{\nabla f(x^{t+1})}}.
\end{align*}

If $\gamma_k = \gamma$ and $\eta_k = \eta$, then 
\begin{align*}
    \Exp{\normstar{r^{k+1} }}
    & \leq (1-\eta)^{k+1} \Exp{\normstar{r^{0}}} \\
    &\quad + \left(\frac{M_0}{2} + \frac{M_1}{3} \cdot\frac{L_0 \gamma}{\eta} \exp\left(\frac{L_1 \gamma}{\eta}\right) \right) \sum_{t=0}^k (1-\eta)^{k-t} \theta\gamma^2   \\
    &\quad + \left(\frac{M_1}{2} + \frac{M_1}{3} \frac{L_1 \gamma}{\eta}\exp\left(\frac{L_1 \gamma}{\eta}\right) \right) \sum_{t=0}^k (1-\eta)^{k-t}  \theta\gamma^2 \Exp{\normstar{\nabla f(x^{t+1})}} \\
    &\quad +   {\color{red}\bar\rho} \sigma_g \sqrt{ \sum_{t=0}^k (1-\eta)^{2(k-t)} \eta^2}.
\end{align*}

Next, since 
\begin{align*}
    \sum_{t=0}^{k} (1-\eta)^{k-t}
    & \leq \sum_{j=0}^\infty (1-\eta)^j  =  \frac{1}{\eta}, \\
    &\text{and}&\\
     \sum_{t=0}^{k} (1-\eta)^{2(k-t)} \eta^2 & \leq  \eta^2 \sum_{j=0}^\infty ((1-\eta)^2)^j = \frac{\eta^2}{1-(1-\eta)^2} \\
     & = \frac{\eta}{2-\eta} \overset{\eta \in [0,1]}{\leq} \eta,
\end{align*}
we obtain
\begin{align*}
    \Exp{\normstar{r^{k+1} }}
    & \leq (1-\eta)^{k+1} \Exp{\normstar{r^{0}}}\\
    &\quad + \left(\frac{M_0}{2} + \frac{M_1}{3}\cdot\frac{L_0 \gamma}{\eta} \exp\left(\frac{L_1\gamma}{\eta}\right)\right) \frac{\gamma^2}{\eta^2} \\
    &\quad +    {\color{red}\bar\rho} \sqrt{\eta}\sigma_g  \\
    &\quad + \left(\frac{M_1}{2} + \frac{M_1}{3}\frac{L_1\gamma}{\eta}\exp\left(\frac{L_1\gamma}{\eta}\right)\right)\\
    &\qquad\cdot \sum_{t=0}^k (1-\eta)^{k-t}  \theta\gamma^2 \Exp{\normstar{\nabla f(x^{t+1})}} .
\end{align*}

Therefore, recalling $\varphi = 1 - \exp(L_1 \gamma) 
 L_1\gamma/2$, we have 
\begin{align*}
\sum_{k=0}^K\gamma\varphi\Exp{\normstar{\nabla f(x^k)}}
&\leq \Delta+\frac{L_0}{2}\sum_{k=0}^K\exp(L_1\gamma)\gamma^2
+2\gamma\sum_{k=0}^K\Exp{\normstar{r^k}}\\
&\leq \Delta+\frac{L_0}{2}\exp(L_1\gamma)\gamma^2(K+1)\\
&\quad+2\gamma\sum_{k=0}^K(1-\eta)^k\Exp{\normstar{r^0}}\\
&\quad+\gamma\left(
M_0+\frac{2M_1}{3}\frac{L_0\gamma}{\eta}
\exp\left(\frac{L_1\gamma}{\eta}\right)
\right)\frac{\gamma^2}{\eta^2}(K+1)\\
&\quad+\gamma\left(
M_1+\frac{2M_1}{3}\frac{L_1\gamma}{\eta}
\exp\left(\frac{L_1\gamma}{\eta}\right)
\right)\\
&\qquad\qquad{}\times
\sum_{k=0}^K\sum_{t=0}^{k-1}(1-\eta)^{k-1-t}
\theta\gamma^2\Exp{\normstar{\nabla f(x^{t+1})}}\\
&\quad+2\gamma{\color{red}\bar\rho}\sqrt{\eta}\sigma_g(K+1).
\end{align*}
where $\Delta = \Exp{f(x^{0})-f^{\inf}}$.

Next, since 
\begin{gather*}
    \sum_{k=0}^K (1-\eta)^k  \leq \frac{1}{\eta} \quad\text{and}\\
    \sum_{k=0}^K \sum_{t=0}^{k-1} (1-\eta)^{k-1-t} \Exp{\normstar{\nabla f(x^{t+1})}}
    \leq \frac{1}{\eta} \sum_{k=0}^K \Exp{\normstar{\nabla f(x^{k})}},
\end{gather*}
we obtain
\begin{align*}
    \sum_{k=0}^K \gamma\vartheta  \Exp{\normstar{\nabla f(x^{k})}}
    & \leq \Delta  + \frac{L_0\gamma^2}{2}  \exp(L_1 \gamma)  (K+1) + 2  \frac{\gamma}{\eta}  \Exp{\normstar{r^{0}}}\\
    &\quad + \left(M_0 + \frac{2M_1}{3}\cdot\frac{L_0 \gamma}{\eta} \exp\left(\frac{L_1\gamma}{\eta}\right)\right) \frac{\gamma^3}{\eta^2}(K+1)  \\
    &\quad +   2\gamma{\color{red}\bar\rho} \sqrt{\eta}\sigma_g(K+1),
\end{align*}
where $\vartheta \eqdef \left( 1 - \exp(L_1 \gamma) 
 \frac{L_1\gamma}{2}   -  \left(M_1 + \frac{2M_1}{3} \cdot\frac{L_1 \gamma}{\eta} \exp\left(\frac{L_1 \gamma}{\eta}\right) \right) \frac{\gamma^2}{\eta^2}  \right)$.

 If $\gamma \leq \min\left\{ \frac{1}{3L_1} , \frac{1}{3\sqrt{M_1}}  \right\} \eta$, then $L_1\gamma \leq L_1\gamma/\eta \leq 1/3$ and $M_1\gamma^2/\eta^2 \leq 1/9$, so $\vartheta \geq 1 - \exp(1/3)/6 - 1.31/9 \geq 1/2$, and 
\begin{align*}
    \frac{\gamma}{2}\sum_{k=0}^K   \Exp{\normstar{\nabla f(x^{k})}}   
    & \leq \Delta + (K+1)\frac{L_0}{2}  \exp(L_1 \gamma) \gamma^2 + 2  \frac{\gamma}{\eta}  \Exp{\normstar{r^{0}}} \\
    &\quad + 2\gamma   {\color{red}\bar\rho} \sqrt{\eta}\sigma_g(K+1)\\
    &\quad+ \left(M_0 + \frac{2M_1}{3}\cdot\frac{L_0 \gamma}{\eta} \exp\left(\frac{L_1\gamma}{\eta}\right)\right) \frac{\gamma^3}{\eta^2}(K+1).
\end{align*}

Therefore, 
\begin{align*}
    \frac{1}{K+1}\sum_{k=0}^K   \Exp{\normstar{\nabla f(x^{k})}} 
    & \leq \frac{2\Delta}{\gamma(K+1)} +L_0 \exp(L_1 \gamma) \gamma \\
    &\quad +  \frac{4\Exp{\normstar{r^{0}}}}{\eta (K+1)} + 4   {\color{red}\bar\rho} \sqrt{\eta}\sigma_g\\
    &\quad + 2 \left(M_0 + \frac{2M_1}{3}\cdot\frac{L_0 \gamma}{\eta} \exp\left(\frac{L_1\gamma}{\eta}\right)\right) \frac{\gamma^2}{\eta^2}.
\end{align*}

If $\gamma = \frac{\hat \gamma}{(K+1)^{5/7}} $ with $\hat \gamma = \min\left\{ \frac{1}{3L_1} , \frac{1}{3\sqrt{M_1}}  \right\} $ and $\eta = \frac{1}{(K+1)^{4/7}}$, then 
\begin{align*}
    \frac{1}{K+1}\sum_{k=0}^K   \Exp{\normstar{\nabla f(x^{k})}}   
    & \leq \frac{2\Delta}{\hat\gamma(K+1)^{2/7}} + \frac{\hat\gamma L_0 \exp(L_1 \hat \gamma)}{(K+1)^{5/7}} \\
    &\quad + \frac{4\Exp{\normstar{r^{0}}}}{ (K+1)^{2/7}} + 4\frac{{\color{red}\bar\rho} \sigma_g}{(K+1)^{2/7}}\\
    &\quad + 2\left(M_0 + \frac{2M_1}{3} L_0\hat \gamma \exp(L_1\hat \gamma)\right) \frac{\hat \gamma^2}{(K+1)^{2/7}}.
\end{align*}

Finally, by the fact that
\begin{equation*}
\min_{k \in \{0,1,\ldots,K\}} \Exp{\normstar{\nabla f(x^{k})}}
\leq \frac{1}{K+1}\sum_{k=0}^K \Exp{\normstar{\nabla f(x^{k})}},
\end{equation*}
and $(K+1)^{-3/7} \leq (K+1)^{-2/7}$, we obtain the final result.

\newpage
\section{LMO-based Methods with Momentum Variance Reduction}\label{ch8:app:MVR_relaxed_smoothness_chapter_8}

This section establishes the convergence properties of the LMO-based algorithms described in~\eqref{ch8:eqn:LMO_momentum_chapter_8}, utilizing the momentum-based variance reduction update given by \eqref{ch8:eq:MARS_chapter_8}.

To facilitate the convergence analysis of MVR, \algname{STORM}, and \algname{MARS}, we use the \steponeonechange{mean-squared generalized smoothness} in Assumption~\ref{ch8:assum:Relaxed_Mean_Squared_Lipschitz_Grad_chapter_8}.

\begin{lemma}[Proposition 5.3 from \textcite{chen2023generalized}]
\label{ch8:lemma:Relaxed_Mean_Squared_Lipschitz_Grad_chapter_8}
    Let $f$ satisfy Assumption~\ref{ch8:assum:Relaxed_Mean_Squared_Lipschitz_Grad_chapter_8}. Then, for all $x,y \in \R^d$,
    \begin{align*}
        \thesischeightExpSub{\xi}{\normstar{\nabla f(x;\xi)-\nabla f(y;\xi)}^2}
        &\leq
        2\left(
        \steponeonechange{\cL_0}^2
        +2\steponeonechange{\cL_1}^2
        \thesischeightExpSub{\xi}{\normstar{\nabla f(y;\xi)}^2}
        \right)\\
        &\quad{}\times
        \exp\left(12\steponeonechange{\cL_1}^2\norm{x-y}^2\right)
        \norm{x-y}^2.
    \end{align*}
\end{lemma}

\begin{proof}
    The proof adapts the arguments of Proposition 5.3 in \textcite{chen2023generalized} to the setting of arbitrary norms.
\end{proof}

Now, we provide the main convergence theorem.

\begin{theorem}
    \label{ch8:thm:MVR_relaxed_smoothness_chapter_8}
    Consider the problem of minimizing $f(x)=\thesischeightExpSub{\xi \sim \thesischeightcD}{f(x;\xi)}$. Suppose Assumptions~\ref{ch8:assum:Hessian_boundedvariance_v2_chapter_8},~\ref{ch8:assum:lower_bound_chapter_8},~\ref{ch8:assum:Relaxed_Lipschitz_Grad_chapter_8}, and~\ref{ch8:assum:Relaxed_Mean_Squared_Lipschitz_Grad_chapter_8} hold.
    Then, the iterates $\{x^{k}\}$ generated by the LMO-based method in~\eqref{ch8:eqn:LMO_momentum_chapter_8} using the momentum variance reduction in~\eqref{ch8:eq:MARS_chapter_8} with  $\rho_k = 1-\eta_k$, and  step sizes
    \begin{equation*}
        \eta_k = \eta = \frac{1}{(K+1)^{2/3}} \quad \text{and} \quad \gamma_k= \gamma = \frac{\hat \gamma}{(K+1)^{2/3}},
    \end{equation*}
    where $\hat \gamma = \frac{1}{3}\min\left\{\frac{1}{L_1}, \frac{1}{11\steponeonechange{\cL_1}}\left(\frac{\color{red}\bar \rho}{\color{red} \thesischeightubar \rho}\right)^{-1}\right\}$, satisfy
    \begin{align*}
        \MoveEqLeft[1] \min_{k \in \{0,1,\ldots,K\}} \Exp{\normstar{\nabla f(x^{k})}} \\
        & \leq \frac{2(f(x^{0}) - f^{\inf})}{\hat\gamma(K+1)^{1/3}} + \frac{4\left( \Exp{ \normstar{r^{0}}} + {\color{red} \bar \rho} \sigma_g\right)}{ (K+1)^{1/3}} +  \frac{\hat\gamma L_0 \exp(L_1 \hat \gamma)}{(K+1)^{2/3}} \\
        &\quad + 8 \frac{\color{red}\bar \rho}{\color{red} \thesischeightubar \rho} \left(\steponeonechange{\cL_0} + 2\steponeonechange{\cL_1} {\color{red} \bar \rho} \sigma_g\right)\exp\left(6 \steponeonechange{\cL_1^2}\hat \gamma^2\right)\frac{\hat \gamma}{(K+1)^{1/3}}.
    \end{align*}
\end{theorem}

\subsection{Proof of Theorem~\ref{ch8:thm:MVR_relaxed_smoothness_chapter_8}}

We prove the result in the following steps. 

\paragraph{Step 1) Proving the descent inequality.} 
By following the proof arguments in Step 1) of~Theorem~\ref{ch8:thm:V1_relaxed_smoothness_chapter_8} (see \eqref{ch8:eq:descent_inequality_chapter_8}), we have 
\begin{align*}
    \sum_{k=0}^K \gamma_k\varphi_k  \normstar{\nabla f(x^{k})}    & \leq  \Delta + 2 \sum_{k=0}^K \gamma_k \normstar{\nabla f(x^{k})-m^{k}} \\
    &\quad + \frac{L_0}{2} \sum_{k=0}^K \exp(L_1 \gamma_k) \gamma_k^2,
\end{align*}
where $\varphi_k \eqdef \left( 1 - \exp(L_1 \gamma_k) 
 L_1\gamma_k/2 \right)$ and $\Delta \eqdef f(x^{0}) - f^{\inf}$.

\paragraph{Step 2) Bounding the error term.}
Next, we bound $\normstar{r^{k}}$, where $r^{k} \eqdef  m^{k} - \nabla f(x^{k})$. 
From the definition of $r^{k}$,
\begin{align*}
r^{k+1}
&=m^{k+1}-\nabla f(x^{k+1})\\
&=(1-\eta_k)\left(m^k-\nabla f(x^k)\right)\\
&\quad+(1-\eta_k)\Bigl[
    \nabla f(x^{k+1};\xi^{k+1})
    -\nabla f(x^k;\xi^{k+1})\\
&\hspace{49mm}
    -\left(\nabla f(x^{k+1})-\nabla f(x^k)\right)
    \Bigr]\\
&\quad+\eta_k\left(
    \nabla f(x^{k+1};\xi^{k+1})
    -\nabla f(x^{k+1})
    \right)\\
&=(1-\eta_k)r^k+(1-\eta_k)\hat S_{k+1}
  +\eta_k\zeta^g_{k+1}.
\end{align*}
For any $t\geq0$, we denote
\begin{align*}
\hat S_{t+1}
&=\left(
\nabla f(x^{t+1};\xi^{t+1})-\nabla f(x^t;\xi^{t+1})
\right)
-\left(\nabla f(x^{t+1})-\nabla f(x^t)\right),\\
\zeta^g_{t+1}
&=\nabla f(x^{t+1};\xi^{t+1})-\nabla f(x^{t+1}).
\end{align*}

By recursively applying the above identity, 
\begin{align*}
r^{k+1}
&=\prod_{t=0}^k(1-\eta_t)r^0\\
&\quad+\sum_{t=0}^k
\left(\prod_{j=t+1}^k(1-\eta_j)\right)
(1-\eta_t)\hat S_{t+1}\\
&\quad+\sum_{t=0}^k
\left(\prod_{j=t+1}^k(1-\eta_j)\right)
\eta_t\zeta^g_{t+1}.
\end{align*}
Set
\begin{align*}
\circledNine
&\eqdef
\Exp{\normstar{\sum_{t=0}^k
\left(\prod_{j=t+1}^k(1-\eta_j)\right)
(1-\eta_t)\hat S_{t+1}}},\\
\circledTen
&\eqdef
\Exp{\normstar{\sum_{t=0}^k
\left(\prod_{j=t+1}^k(1-\eta_j)\right)
\eta_t\zeta^g_{t+1}}}.
\end{align*}
Therefore,
\begin{align*}
\Exp{\normstar{r^{k+1}}}
&\leq
\Exp{\normstar{\prod_{t=0}^k(1-\eta_t)r^0}}
+\circledNine+\circledTen.
\end{align*}

Next, we bound $\circledTen$ by using~Assumption~\ref{ch8:assum:Hessian_boundedvariance_v2_chapter_8} and by following the proof arguments for bounding $\circledTwo$ in Step 2) of the convergence proof for~Theorem~\ref{ch8:thm:V1_relaxed_smoothness_chapter_8}. The conditional centering of $\zeta^g_{t+1}$ makes all cross-time inner products vanish by the tower property. Then, we have 
\begin{align*}
    \circledTen
    & \leq {\color{red}\bar\rho}   \Exp{\normtwo{\sum_{t=0}^k \left(\prod_{j=t+1}^k (1-\eta_j) \right) \eta_t \zeta^{g}_{t+1}}} \\
    & \leq {\color{red}\bar\rho}  \sqrt{  \Exp{\normtwo{\sum_{t=0}^k \left(\prod_{j=t+1}^k (1-\eta_j) \right) \eta_t \zeta^{g}_{t+1}}^2} } \\
    & \leq {\color{red}\bar\rho} \sigma_g \sqrt{ \sum_{t=0}^k \left(\prod_{j=t+1}^k (1-\eta_j)^2 \right) \eta_t^2}. 
\end{align*}
Therefore, 
\begin{align*}
    \Exp{\normstar{r^{k+1} }}
    & \leq \prod_{t=0}^k (1-\eta_t)\Exp{\normstar{r^{0}}} + ~\circledNine~ \\
    &\quad +   {\color{red}\bar\rho} \sigma_g \sqrt{ \sum_{t=0}^k \left(\prod_{j=t+1}^k (1-\eta_j)^2 \right) \eta_t^2}.
\end{align*}

Next, we bound $\circledNine$. For $0\leq\ell\leq k$, define
\begin{align*}
\mathcal S_{\ell}^{(k)}
\eqdef
\sum_{t=0}^{\ell}
\left(\prod_{j=t+1}^k(1-\eta_j)\right)
(1-\eta_t)\hat S_{t+1}.
\end{align*}
The shared sample $\xi^{t+1}$ is fresh conditionally on $\mathcal F_t$, so
\begin{align*}
\Exp{\hat S_{t+1}\mid\mathcal F_t}=0.
\end{align*}
Moreover, the conditional form of
Lemma~\ref{ch8:lemma:Relaxed_Mean_Squared_Lipschitz_Grad_chapter_8} gives, for
any $t\geq0$,
\begin{align*}
\Exp{\norm{\hat S_{t+1}}_2^2\mid\mathcal F_t}
&\leq
\Exp{
{\norm{\nabla f(x^{t+1};\xi^{t+1})
-\nabla f(x^t;\xi^{t+1})}_2^2}
\mid\mathcal F_t}\\
&\leq
\frac{1}{{\color{red}\thesischeightubar{\rho}}^2}
\Exp{
{\normstar{\nabla f(x^{t+1};\xi^{t+1})
-\nabla f(x^t;\xi^{t+1})}^2}
\mid\mathcal F_t}\\
&\leq
\frac{2}{{\color{red}\thesischeightubar{\rho}}^2}
\left(
\steponeonechange{\cL_0}^2
+2\steponeonechange{\cL_1}^2
\Exp{\normstar{\nabla f(x^t;\xi^{t+1})}^2\mid\mathcal F_t}
\right)\\
&\qquad{}\times
\exp\left(12\steponeonechange{\cL_1}^2\norm{x^{t+1}-x^t}^2\right)
\norm{x^{t+1}-x^t}^2\\
&\leq
\frac{2}{{\color{red}\thesischeightubar{\rho}}^2}
\left(
\steponeonechange{\cL_0}^2
+2\steponeonechange{\cL_1}^2
\Exp{\normstar{\nabla f(x^t;\xi^{t+1})}^2\mid\mathcal F_t}
\right)\\
&\qquad{}\times
\exp\left(12\steponeonechange{\cL_1}^2\gamma_t^2\right)\gamma_t^2\\
&\leq
\frac{2}{{\color{red}\thesischeightubar{\rho}}^2}
\left(
\steponeonechange{\cL_0}^2
+4\steponeonechange{\cL_1}^2
\Exp{\normstar{\nabla f(x^t;\xi^{t+1})
-\nabla f(x^t)}^2\mid\mathcal F_t}
\right)\\
&\qquad{}\times
\exp\left(12\steponeonechange{\cL_1}^2\gamma_t^2\right)\gamma_t^2\\
&\quad+
\frac{8}{{\color{red}\thesischeightubar{\rho}}^2}
\gamma_t^2\steponeonechange{\cL_1}^2
\exp\left(12\steponeonechange{\cL_1}^2\gamma_t^2\right)
\normstar{\nabla f(x^t)}^2\\
&\leq A_t^2+B_t^2,
\end{align*}
where the penultimate inequality uses
$\normstar{a+b}^2\leq2\normstar{a}^2+2\normstar{b}^2$, and the last one uses the conditional
gradient-noise bound from Assumption~\ref{ch8:assum:Hessian_boundedvariance_v2_chapter_8} together with $\normstar{\cdot}\leq{\color{red}\bar\rho}\normtwo{\cdot}$. Here
\begin{align*}
A_t^2
&\eqdef
\left(
\frac{2}{{\color{red}\thesischeightubar{\rho}}^2}
\gamma_t^2\steponeonechange{\cL_0}^2
+\frac{8{\color{red}\bar\rho^2}}
{{\color{red}\thesischeightubar{\rho}}^2}
\gamma_t^2\steponeonechange{\cL_1}^2\sigma_g^2
\right)
\exp\left(12\steponeonechange{\cL_1}^2\gamma_t^2\right),\\
B_t^2
&\eqdef
\frac{8}{{\color{red}\thesischeightubar{\rho}}^2}
\gamma_t^2\steponeonechange{\cL_1}^2
\exp\left(12\steponeonechange{\cL_1}^2\gamma_t^2\right)
\normstar{\nabla f(x^t)}^2.
\end{align*}

We now handle the two initial indices directly. For $k=0$, conditional
Jensen's inequality and
$\sqrt{a^2+b^2}\leq a+b$ give
\begin{align*}
\circledNine
&\leq{\color{red}\bar\rho}(1-\eta_0)
\left(A_0+\Exp{B_0}\right).
\end{align*}
For $k=1$, conditioning successively on $\mathcal F_1$ and
$\mathcal F_0$ yields
\begin{align*}
\circledNine
&\leq{\color{red}\bar\rho}
\Bigl(
(1-\eta_0)^2(1-\eta_1)^2A_0^2
+(1-\eta_1)^2A_1^2
\Bigr)^{1/2}\\
&\quad+{\color{red}\bar\rho}
(1-\eta_0)(1-\eta_1)\Exp{B_0}
+{\color{red}\bar\rho}(1-\eta_1)\Exp{B_1}.
\end{align*}

It remains to establish the expansion for $k\geq2$. Define
\begin{align*}
\mathcal Q_k
&\eqdef
\norm{\mathcal S_{k-2}^{(k)}}_2^2
+(1-\eta_{k-1})^2(1-\eta_k)^2
\Exp{\norm{\hat S_k}_2^2\mid\mathcal F_{k-1}}
+(1-\eta_k)^2A_k^2,\\
\widetilde{\mathcal Q}_k
&\eqdef
\norm{\mathcal S_{k-2}^{(k)}}_2^2
+(1-\eta_{k-1})^2(1-\eta_k)^2A_{k-1}^2
+(1-\eta_k)^2A_k^2.
\end{align*}
The conditional centering of $\hat S_{k+1}$ eliminates the cross term, and
conditional Jensen's inequality gives
\begin{align*}
\circledNine
&\leq{\color{red}\bar\rho}
\Exp{\left(
\norm{\mathcal S_{k-1}^{(k)}}_2^2
+(1-\eta_k)^2
\Exp{\norm{\hat S_{k+1}}_2^2\mid\mathcal F_k}
\right)^{1/2}}\\
&\leq{\color{red}\bar\rho}
\Exp{\mathcal Q_k^{1/2}}
+{\color{red}\bar\rho}(1-\eta_k)\Exp{B_k}\\
&\leq{\color{red}\bar\rho}
\Exp{\widetilde{\mathcal Q}_k^{1/2}}
+{\color{red}\bar\rho}(1-\eta_k)\Exp{B_k}\\
&\quad+{\color{red}\bar\rho}
(1-\eta_{k-1})(1-\eta_k)\Exp{B_{k-1}}.
\end{align*}
Iterating this argument for $k\geq2$ yields
\begin{align*}
    \circledNine & \leq {\color{red}\bar\rho}\left( \sum^{k}_{t=0} \left(\prod^k_{j = t}(1-\eta_j)^2\right)A^2_{t} \right)^{1/2} \\
    &\quad + {\color{red}\bar\rho}\sum^{k}_{t=0} \left(\prod^k_{j = t}(1-\eta_j)\right)\Exp{B_{t}}\\
    & \leq \left(2\steponeonechange{\cL_0} +4\steponeonechange{\cL_1}  {\color{red}\bar\rho}\sigma_g \right)\frac{{\color{red}\bar\rho}}{\color{red} \thesischeightubar{\rho}}\left( \sum^{k}_{t=0} \left(\prod^k_{j = t}(1-\eta_j)^2\right)\gamma^2_t\exp\left(12 \steponeonechange{\cL_1^2}\gamma_t^2\right) \right)^{1/2} \\
    &\quad + 4\frac{{\color{red}\bar\rho}}{\color{red} \thesischeightubar{\rho}}\sum^{k}_{t=0} \left(\prod^k_{j = t}(1-\eta_j)\right)  \gamma_t\steponeonechange{\cL_1} \exp\left(6 \steponeonechange{\cL_1^2}\gamma_t^2\right)\Exp{\normstar{\nabla f(x^{t})}}.
\end{align*}
The direct bounds for $k=0$ and $k=1$ are exactly the corresponding
one- and two-term instances of this expansion. Hence the displayed bound
holds for every $k\geq0$, and the theorem's subsequent sums still cover all
indices $0,\ldots,K$.

Therefore,
\begin{align*}
    \Exp{\normstar{r^{k+1} }}
    & \leq \prod_{t=0}^k (1-\eta_t)\Exp{\normstar{r^{0}}} \\
    &\quad +   {\color{red}\bar\rho} \sigma_g \sqrt{ \sum_{t=0}^k \left(\prod_{j=t+1}^k (1-\eta_j)^2 \right) \eta_t^2}\\
    &\quad + \left(2\steponeonechange{\cL_0} +4\steponeonechange{\cL_1}  {\color{red}\bar\rho}\sigma_g \right)\frac{{\color{red}\bar\rho}}{\color{red} \thesischeightubar{\rho}}\left( \sum^{k}_{t=0} \left(\prod^k_{j = t}(1-\eta_j)^2\right)\gamma^2_t\exp\left(12 \steponeonechange{\cL_1^2}\gamma_t^2\right) \right)^{1/2} \\
    &\quad + 4\frac{{\color{red}\bar\rho}}{\color{red} \thesischeightubar{\rho}}\sum^{k}_{t=0} \left(\prod^k_{j = t}(1-\eta_j)\right)  \gamma_t\steponeonechange{\cL_1} \exp\left(6 \steponeonechange{\cL_1^2}\gamma_t^2\right)\Exp{\normstar{\nabla f(x^{t})}}.
\end{align*}

If $\gamma_k = \gamma$ and $\eta_k = \eta$, then 
\begin{align*}
    \Exp{\normstar{r^{k+1} }}
    & \leq (1-\eta)^{k+1}\Exp{\normstar{r^{0}}} +   {\color{red}\bar\rho} \sigma_g \sqrt{ \sum_{t=0}^k  (1-\eta)^{2(k-t)}\eta^2}\\
    &\quad + 2\frac{{\color{red}\bar\rho}}{\color{red} \thesischeightubar{\rho}}\left(\steponeonechange{\cL_0} +2\steponeonechange{\cL_1}  {\color{red}\bar\rho}\sigma_g \right)\gamma\exp\left(6 \steponeonechange{\cL_1^2}\gamma^2\right)\left( \sum^{k}_{t=0} (1-\eta)^{2(k-t+1)} \right)^{1/2} \\
    &\quad + 4\frac{{\color{red}\bar\rho}}{\color{red} \thesischeightubar{\rho}}\steponeonechange{\cL_1} \exp\left(6 \steponeonechange{\cL_1^2}\gamma^2\right)\sum^{k}_{t=0} (1-\eta)^{k-t+1}\gamma\Exp{\normstar{\nabla f(x^{t})}}.
\end{align*}

Next, since 
\begin{align*}
    \sum_{t=0}^{k} (1-\eta)^{2(k-t+1)}
    & \leq \sum_{j=0}^\infty (1-\eta)^{2j}  =  \frac{1}{1-(1-\eta)^2} \leq \frac{1}{\eta} \quad\text{and}\\
     \sum_{t=0}^{k} (1-\eta)^{2(k-t)} \eta^2 & \leq  \eta^2 \sum_{j=0}^\infty ((1-\eta)^2)^j = \frac{\eta^2}{1-(1-\eta)^2} \\
     & = \frac{\eta}{2-\eta} \overset{\eta \in (0,1]}{\leq} \eta,
\end{align*}
we obtain
\begin{align*}
    \Exp{\normstar{r^{k+1} }}
    & \leq (1-\eta)^{k+1}\Exp{\normstar{r^{0}}} +   {\color{red}\bar\rho} \sigma_g \sqrt{\eta} \\
    &\quad + 2\frac{{\color{red}\bar\rho}}{\color{red} \thesischeightubar{\rho}}\left(\steponeonechange{\cL_0} +2\steponeonechange{\cL_1}  {\color{red}\bar\rho}\sigma_g \right)\exp\left(6 \steponeonechange{\cL_1^2}\gamma^2\right) \frac{\gamma}{\sqrt{\eta}} \\
    &\quad + 4\frac{{\color{red}\bar\rho}}{\color{red} \thesischeightubar{\rho}}\steponeonechange{\cL_1} \exp\left(6 \steponeonechange{\cL_1^2}\gamma^2\right)\sum^{k}_{t=0} (1-\eta)^{k-t+1}\gamma\Exp{\normstar{\nabla f(x^{t})}}.
\end{align*}

Therefore, recalling $\varphi = 1 - \exp(L_1 \gamma) 
 L_1\gamma/2$, we have 
\begin{align*}   
    \sum_{k=0}^K  \gamma \varphi \Exp{\normstar{\nabla f(x^{k})}}   & \leq  \Delta + \frac{L_0}{2} \sum_{k=0}^K \exp(L_1 \gamma) \gamma^2 + 2  \gamma \sum_{k=0}^K \Exp{\normstar{r^{k}}}   \\
    & \leq \Delta  + \frac{L_0}{2}  \exp(L_1 \gamma) \gamma^2(K+1)  + 2  \gamma \sum_{k=0}^K (1-\eta)^{k} \Exp{\normstar{r^{0}}}\\
    &\quad+  4\gamma\frac{{\color{red}\bar\rho}}{\color{red} \thesischeightubar{\rho}}\left(\steponeonechange{\cL_0} +2\steponeonechange{\cL_1}  {\color{red}\bar\rho}\sigma_g \right)\exp\left(6 \steponeonechange{\cL_1^2}\gamma^2\right) \frac{\gamma}{\sqrt{\eta}} (K+1) \\
    &\quad + 8\gamma \frac{{\color{red}\bar\rho}}{\color{red} \thesischeightubar{\rho}}\steponeonechange{\cL_1} \exp\left(6 \steponeonechange{\cL_1^2}\gamma^2\right) \sum_{k=0}^K \sum_{t=0}^{k-1} (1-\eta)^{k-t}  \gamma \Exp{\normstar{\nabla f(x^{t})}} \\
    &\quad + 2\gamma{\color{red}\bar\rho} \sqrt{\eta}\sigma_g(K+1),
\end{align*}
where $\Delta = \Exp{f(x^{0})-f^{\inf}}$.

Next, since 
\begin{gather*}
    \sum_{k=0}^K (1-\eta)^k  \leq \frac{1}{\eta} \quad\text{and}\\
    \sum_{k=0}^K \sum_{t=0}^{k-1} (1-\eta)^{k-t} \Exp{\normstar{\nabla f(x^{t})}}
    \leq \frac{1}{\eta} \sum_{k=0}^K \Exp{\normstar{\nabla f(x^{k})}},
\end{gather*}
we obtain
\begin{align*}
    \sum_{k=0}^K \gamma\vartheta  \Exp{\normstar{\nabla f(x^{k})}}
    & \leq \Delta  + \frac{L_0\gamma^2}{2}  \exp(L_1 \gamma)  (K+1) + 2  \frac{\gamma}{\eta} \Exp{ \normstar{r^{0}}}\\
    &\quad + 4\gamma\frac{{\color{red}\bar\rho}}{\color{red} \thesischeightubar{\rho}}\left(\steponeonechange{\cL_0} + 2\steponeonechange{\cL_1}  {\color{red}\bar\rho}\sigma_g \right)\exp\left(6 \steponeonechange{\cL_1^2}\gamma^2\right) \frac{\gamma}{\sqrt{\eta}} (K+1) \\
    &\quad +  2\gamma {\color{red}\bar\rho} \sqrt{\eta}\sigma_g(K+1),
\end{align*}
where $\vartheta \eqdef \left( 1 - \exp(L_1 \gamma) 
 \frac{L_1\gamma}{2}   -  8 \frac{{\color{red}\bar\rho}}{\color{red} \thesischeightubar{\rho}}\steponeonechange{\cL_1} \exp\left(6 \steponeonechange{\cL_1^2}\gamma^2\right)\frac{\gamma}{\eta}  \right)$.

 If $\gamma \leq \min\left\{ \frac{1}{3L_1} , \frac{1}{33\steponeonechange{\cL_1}}\left(\frac{{\color{red}\bar\rho}}{\color{red} \thesischeightubar{\rho}}\right)^{-1}  \right\} \eta$, then $L_1\gamma \leq L_1\gamma/\eta \leq 1/3$, $8({\color{red}\bar\rho}/{\color{red}\thesischeightubar{\rho}})\steponeonechange{\cL_1}\gamma/\eta \leq 8/33$, and $\steponeonechange{\cL_1}\gamma \leq 1/33$ (since ${\color{red}\thesischeightubar{\rho}}\leq{\color{red}\bar\rho}$), so $\vartheta \geq 1 - \exp(1/3)/6 - (8/33)\exp(6/33^2) \geq 1/2$, and 
\begin{align*}
    \frac{\gamma}{2}\sum_{k=0}^K   \Exp{\normstar{\nabla f(x^{k})}}   
    & \leq \Delta + (K+1)\frac{L_0}{2}  \exp(L_1 \gamma) \gamma^2 + 2  \frac{\gamma}{\eta}  \Exp{\normstar{r^{0}}} \\
    &\quad + 2\gamma   {\color{red}\bar\rho} \sqrt{\eta}\sigma_g(K+1)\\
    &\quad+ 4\gamma\frac{{\color{red}\bar\rho}}{\color{red} \thesischeightubar{\rho}}\left(\steponeonechange{\cL_0} + 2\steponeonechange{\cL_1}  {\color{red}\bar\rho}\sigma_g \right)\exp\left(6 \steponeonechange{\cL_1^2}\gamma^2\right) \frac{\gamma}{\sqrt{\eta}} (K+1).
\end{align*}

Therefore, 
\begin{align*}
    \frac{1}{K+1}\sum_{k=0}^K   \Exp{\normstar{\nabla f(x^{k})}} 
    & \leq \frac{2\Delta}{\gamma(K+1)} +L_0 \exp(L_1 \gamma) \gamma \\
    &\quad +  \frac{4\Exp{\normstar{r^{0}}}}{\eta (K+1)} + 4   {\color{red}\bar\rho} \sqrt{\eta}\sigma_g\\
    &\quad +  8\frac{{\color{red}\bar\rho}}{\color{red} \thesischeightubar{\rho}}\left(\steponeonechange{\cL_0} + 2\steponeonechange{\cL_1}  {\color{red}\bar\rho}\sigma_g \right)\exp\left(6 \steponeonechange{\cL_1^2}\gamma^2\right) \frac{\gamma}{\sqrt{\eta}} .
\end{align*}

If $\gamma = \frac{\hat \gamma}{(K+1)^{2/3}} $ with $\hat \gamma = \min\left\{ \frac{1}{3L_1} , \frac{1}{33\steponeonechange{\cL_1}}\left(\frac{{\color{red}\bar\rho}}{\color{red} \thesischeightubar{\rho}}\right)^{-1}  \right\} $ and $\eta = \frac{1}{(K+1)^{2/3}}$, then 
\begin{align*}
    \frac{1}{K+1}\sum_{k=0}^K   \Exp{\normstar{\nabla f(x^{k})}}   
    & \leq \frac{2\Delta}{\hat\gamma(K+1)^{1/3}} + \frac{\hat\gamma L_0 \exp(L_1 \hat \gamma)}{(K+1)^{2/3}} \\
    &\quad + \frac{4 \Exp{\normstar{r^{0}}}}{ (K+1)^{1/3}} + 4\frac{{\color{red}\bar\rho} \sigma_g}{(K+1)^{1/3}}\\
    &\quad +  8\frac{{\color{red}\bar\rho}}{\color{red} \thesischeightubar{\rho}}\left(\steponeonechange{\cL_0} + 2\steponeonechange{\cL_1}  {\color{red}\bar\rho}\sigma_g \right)\exp\left(6 \steponeonechange{\cL_1^2}\hat\gamma^2\right) \frac{\hat\gamma}{(K+1)^{1/3}} .
\end{align*}

Finally, by the fact that
\begin{equation*}
\min_{k \in \{0,1,\ldots,K\}} \Exp{\normstar{\nabla f(x^{k})}}
\leq \frac{1}{K+1}\sum_{k=0}^K \Exp{\normstar{\nabla f(x^{k})}},
\end{equation*}
we obtain the final result.

\newpage
\section{Experiments}\label{ch8:sec:exp_chapter_8}
To ensure a fair and reproducible comparison, the initialization and hyperparameter schedules were standardized across all experiments.

\subsection{Initialization}
For each experiment, a single initial parameter vector, $x^{0}$, was generated by drawing from a normal (Gaussian) distribution, $x^{0} \sim \mathcal{N}(0, 1)$. The random seed was fixed to ensure that this exact same starting point was used for every algorithm evaluated in that experiment. This "far start" initialization is designed to test the robustness and convergence capabilities of the optimizers from a non-trivial region of the parameter space.
\subsection{Stepsize and Momentum Schedules}
The learning rate $\gamma_k$ and the momentum parameter $\eta_k$ are decayed at each iteration $k$. The schedules are chosen based on the theoretical underpinnings of each class of algorithm. Let $\gamma_0$ be a pre-defined initial learning rate.

\begin{itemize}
    \item \textbf{For Polyak Momentum:} This first-order method uses its standard theoretically-backed schedule:
    \begin{align*}
        \eta_k = \frac{1}{\sqrt{k+1}}, \quad
        \gamma_k = \frac{\gamma_0}{(k+1)^{3/4}}.
    \end{align*}

    \item \textbf{For Extrapolated Momentum:} This method uses a distinct schedule designed for its update rule:
    \begin{align*}
        \eta_k = \frac{1}{(k+1)^{4/7}}, \quad
        \gamma_k = \frac{\gamma_0}{(k+1)^{5/7}}.
    \end{align*}

    \item \textbf{For all Second-Order Momentum Variants (SOM-V1, SOM-V2, $\rho$-SOM-V1, and $\rho$-SOM-V2):} These methods share a common schedule for their primary learning rate and momentum parameters:
    \begin{align*}
        \eta_k = \frac{1}{(k+1)^{2/3}}, \quad
        \gamma_k = \frac{\gamma_0}{(k+1)^{2/3}}.
    \end{align*}
\end{itemize}
These schedules ensure that the step sizes and momentum contributions diminish over time, a necessary condition for convergence in stochastic optimization.
\subsection{Nonconvex Logistic Regression}
We evaluate the performance of several stochastic optimization algorithms on a composite non-convex problem. The objective is to benchmark their convergence speed and stability on a logistic regression task augmented with a non-convex regularizer.

The model is a standard logistic regression classifier. The experiment is conducted on the \texttt{splice} dataset from the libsvm library, which contains $1000$ training samples and $60$ features.

The optimization objective is to minimize a composite function $f(x)$, defined as:
\begin{equation*}
    \min_{x \in \mathbb{R}^d} f(x) = \mathcal{L}(x) + R(x)
\end{equation*}
where $x \in \mathbb{R}^{60}$ is the vector of model parameters. The loss term $\mathcal{L}(x)$ is the standard mean logistic loss:
\begin{equation*}
    \mathcal{L}(x) = \frac{1}{N} \sum_{i=1}^{N} \log(1 + \exp(-y_i a_i^T x))
\end{equation*}
The regularization term $R(x)$ is the non-convex Welsch regularizer, defined as:
\begin{equation*}
    R(x) = \lambda \sum_{j=1}^{d} \frac{x_{j}^{2}}{1 + x_{j}^{2}}
\end{equation*}
The regularization hyperparameter $\lambda$ is set to $0.01$.

Figure~\ref{ch8:fig:main_comparison_lr_chapter_8} displays the convergence behavior of the three algorithms over $20000$ iterations. The training loss plot shows that the second-order methods containing Hessian information (Variant 1, Variant 2) outperform the first-order Polyak Momentum. Among these, Variant 2 achieves the fastest convergence and reaches the lowest final loss value.

The gradient norm plot corroborates these findings. While all methods exhibit stochastic oscillations, the overall trend for the second-order methods is a more rapid and consistent decrease in the gradient norm. The Polyak Momentum method converges to a region with a substantially higher gradient norm, indicating a less optimal solution.

\begin{figure}[h!]
    \centering
    \includegraphics[width=0.47\textwidth]{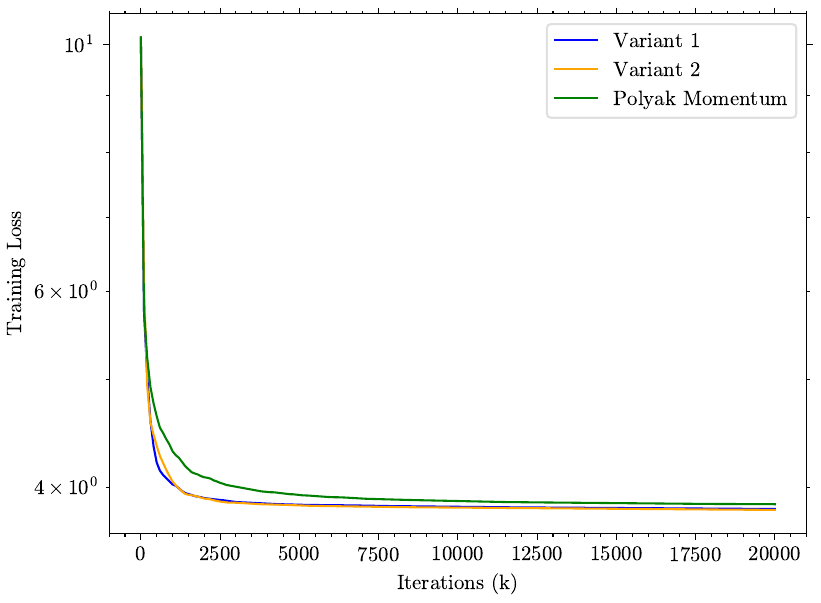}
    \includegraphics[width=0.47\textwidth]{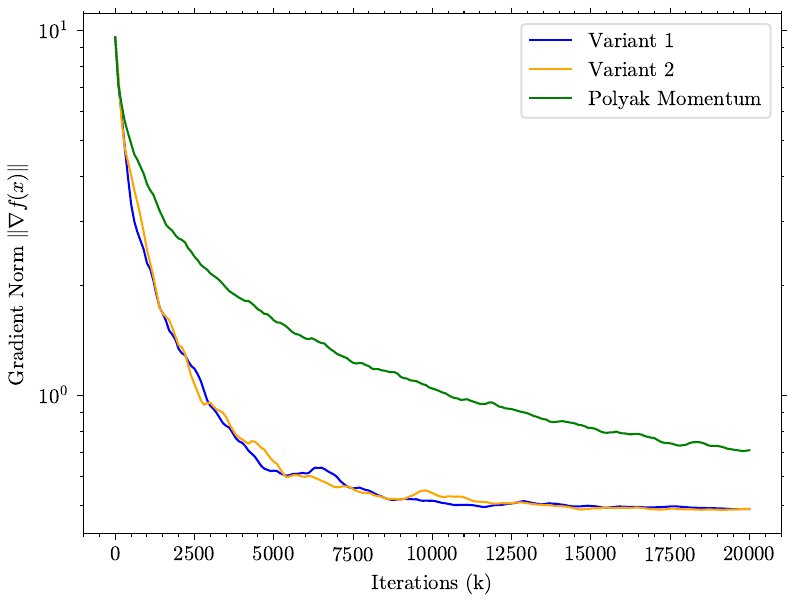}
    \caption{Main algorithm comparison for logistic regression. Training loss vs. iterations. Gradient norm vs. iterations.}
    \label{ch8:fig:main_comparison_lr_chapter_8}
\end{figure}
\subsection{Multi-Layer Perceptron}
We investigate the performance of several stochastic optimization algorithms on a non-convex binary classification problem. Our goal is to compare their convergence properties.

The model is a Multi-Layer Perceptron (MLP) with two hidden layers, implemented in PyTorch. The experiment is conducted on the \texttt{splice} dataset, obtained from the libsvm library. It consists of 1000 training samples, each with 60 features.

The optimization objective is to minimize a composite function $f(x)$, which includes a standard loss term and a non-convex regularizer:
\begin{equation*}
    \min_{x \in \mathbb{R}^d} f(x) = \mathcal{L}(x) + R(x)
\end{equation*}
where $x$ represents the flattened vector of all model parameters.

The loss term, $\mathcal{L}(x)$, is the mean Binary Cross-Entropy with Logits loss, calculated over the entire training dataset.

The regularization term, $R(x)$, is the non-convex Welsch regularizer, chosen to create a more challenging optimization landscape. It is defined as:
\begin{equation*}
    R(x) = \lambda \sum_{i=1}^{d} \frac{x_i^2}{1 + x_i^2}
\end{equation*}
where $\lambda$ is the regularization hyperparameter, set to $0.01$.

We set the theoretical learning rates for the algorithms considered.

The results are presented in two separate comparisons. Figure~\ref{ch8:fig:momentum_comparison_chapter_8} compares four momentum-based algorithms by plotting their running minimum performance. Among these, SOM-V2 achieves the lowest final training loss. Extrapolated Momentum and SOM-V1 perform similarly, while Polyak Momentum converges to a slightly higher loss. The running minimum plot of the gradient norm clarifies the convergence trend for each method by smoothing out stochastic noise.

\begin{figure}[h!]
    \centering
    \includegraphics[width=0.47\textwidth]{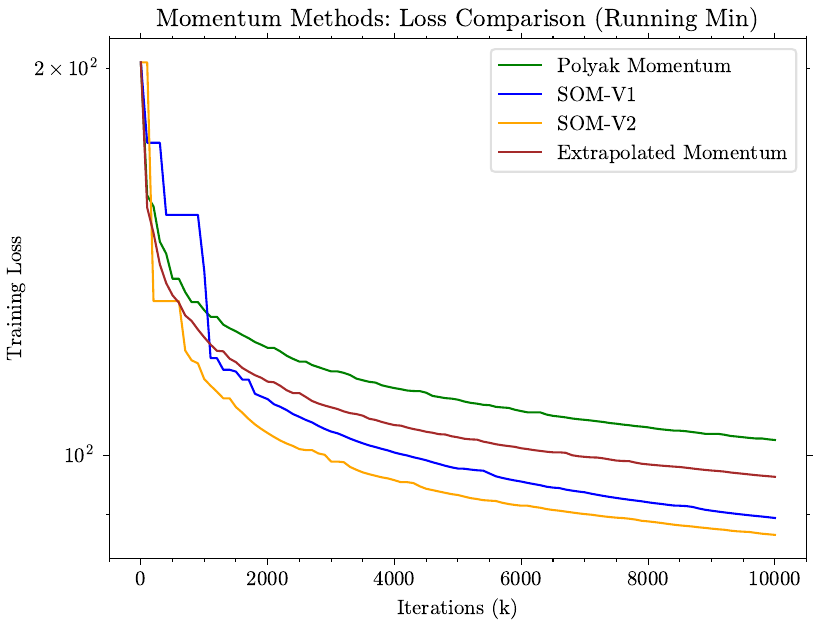}
    \includegraphics[width=0.47\textwidth]{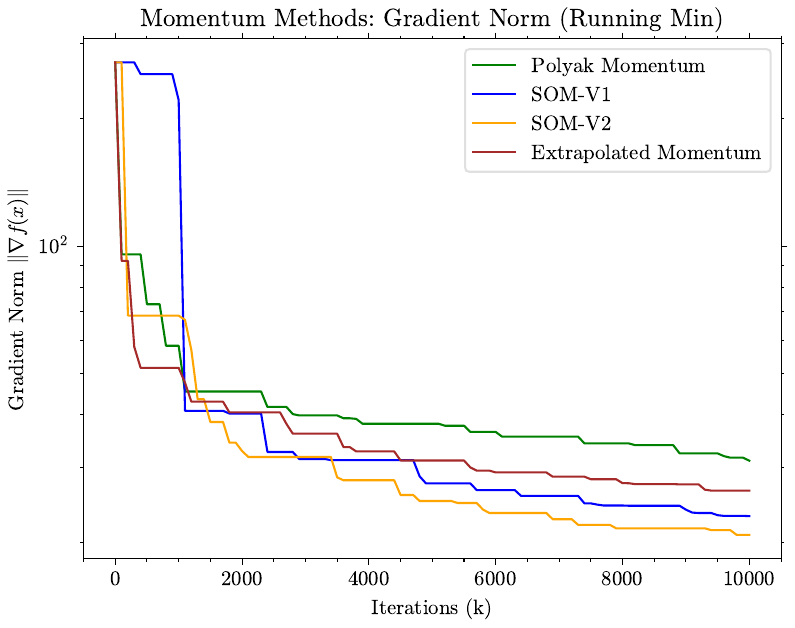}
    \caption{Comparison of momentum-based methods showing the running minimum (best value found so far). Left: Training loss vs. iterations. Right: Gradient norm vs. iterations.}
    \label{ch8:fig:momentum_comparison_chapter_8}
\end{figure}

Figure~\ref{ch8:fig:hvp_comparison_chapter_8} compares the standard second-order momentum methods (SOM-V1 and SOM-V2) with their $\rho$-HVP counterparts. The results clearly show the benefit of the additional Hessian term, as both $\rho$-SOM-V1 and $\rho$-SOM-V2 outperform their respective base variants. $\rho$-SOM-V2 demonstrates the strongest overall performance, converging to the lowest training loss of all tested algorithms. The running minimum gradient norm plot also makes it clear that the $\rho$-SOM methods, particularly $\rho$-SOM-V2, offer a more consistent decrease in gradient magnitude compared to the standard SOM-V1 and SOM-V2 methods.

\begin{figure}[h!]
    \centering
    \includegraphics[width=0.47\textwidth]{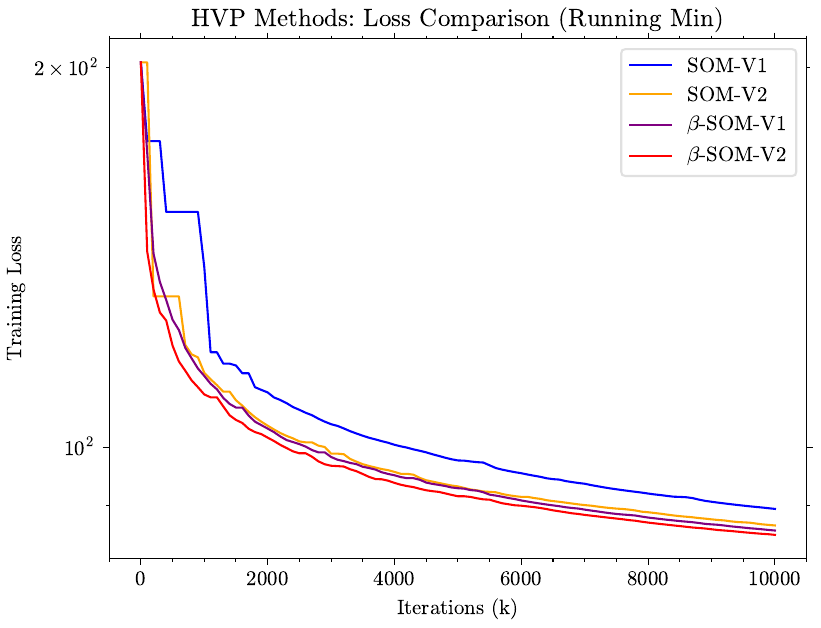}
    \includegraphics[width=0.47\textwidth]{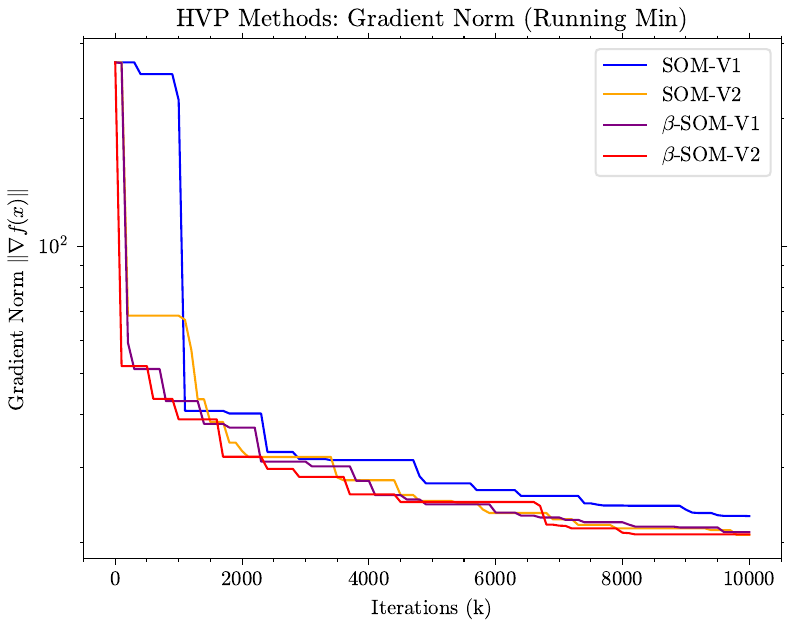}
    \caption{Comparison of HVP-based methods showing the running minimum (best value found so far). Left: Training loss vs. iterations. Right: Gradient norm vs. iterations.}
    \label{ch8:fig:hvp_comparison_chapter_8}
\end{figure}

\subsection{Recurrent Neural Network Training}\label{ch8:sec:LSTMs_chapter_8}
The model is a two-layer LSTM network designed for word-level language modeling. The architecture consists of an embedding layer ($200$ dimensions), a two-layer LSTM core with hidden units of size $200,$ and a final linear decoder layer to produce logits over the vocabulary.

The experiments are conducted on the Penn Treebank (PTB) dataset, a standard and widely-used benchmark for evaluating language models.

Key statistics and preprocessing steps are as follows:
\begin{itemize}
    \item \textbf{Corpus Size:} The dataset is split into training, validation, and testing sets. The training portion, used in our experiment, contains approximately $929,000$ tokens.
    \item \textbf{Vocabulary:} A dictionary is constructed from the unique words present in the training data. An end-of-sentence token, \texttt{<eos>}, is added to each sentence, resulting in a total vocabulary size of $10,000$ unique tokens.
    \item \textbf{Tokenization:} The raw text is tokenized by splitting on whitespace, and each word is converted into its corresponding integer index from the vocabulary.
    \item \textbf{Batching:} For training, the entire sequence of token IDs is reshaped into a fixed number of parallel streams (a batch size of $20$ in our case). The model is then trained on sequential chunks of this data using Truncated Backpropagation Through Time (BPTT) with a sequence length of $35.$
\end{itemize}

The objective is to minimize the standard Cross-Entropy Loss. The parameter update follows the LMO-based rule:
\begin{equation*}
    x^{k+1} = x^{k} + \text{lmo}(m^{k})
\end{equation*}
where the Linear Minimization Oracle, $\text{lmo}(m^{k})$, is defined as:
\begin{equation*}
    \text{lmo}(m^{k}) := \underset{\norm{v} \le \gamma_k}{\argmin} \langle m^{k}, v \rangle.
\end{equation*}

We compare six algorithms: Polyak Momentum, SOM-V1, SOM-V2, our proposed $\rho$-SOM-V1 and $\rho$-SOM-V2 variants, and Extrapolated Momentum. The hyperparameter schedules for $\eta_k$ and $\gamma_k$ are set to their theoretical values for each respective algorithm.
\subsection{The case when $\norm{\cdot}$ is a Euclidean norm}
For this experiment, we use the Euclidean ($\ell_2$) norm. The results of the LMO-based training, presented in Figure~\ref{ch8:fig:lstm_lmo_loss_chapter_8}, reveal a significant divergence in performance between the first-order and more advanced momentum methods.
\begin{figure}[h!]
    \centering
    \includegraphics[width=0.47\textwidth]{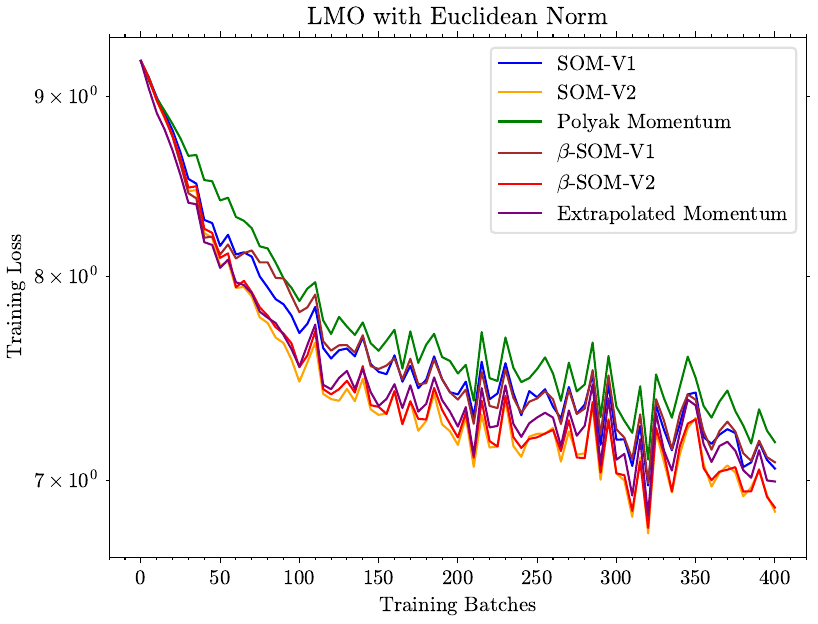}
    \caption{Training loss vs. training batches for the LSTM model using an LMO update with the Euclidean norm. The plot highlights the stronger performance of the other five methods over Polyak Momentum.}
    \label{ch8:fig:lstm_lmo_loss_chapter_8}
\end{figure}
The simple momentum update rule is less efficient at navigating the loss landscape for this task. In contrast, the other five algorithms, which all incorporate either second-order (HVP) or extrapolated gradient information, form a tight cluster of high-performing methods. They converge faster and to a lower loss value than Polyak Momentum. Within this cluster, SOM-V2 and $\rho$-SOM-V2 often achieve a marginally lower loss, suggesting that calculating the curvature information at the current point $x^{k}$ (rather than the interpolated point used by SOM-V1 and $\rho$-SOM-V1) may offer a slight advantage. The strong performance of these five methods indicates that the LMO update rule is a highly effective stabilization technique when paired with momentum strategies that utilize more sophisticated directional information.
\subsection{The case when $\norm{\cdot}$ is $\ell_{\infty}$-norm}
To further investigate the effect of the update geometry, we conducted a parallel experiment where the LMO is constrained by the Infinity norm ($\ell_\infty$). The resulting training loss curves are presented in Figure~\ref{ch8:fig:lstm_linf_loss_chapter_8}.
\begin{figure}[h!]
    \centering
    \includegraphics[width=0.47\textwidth]{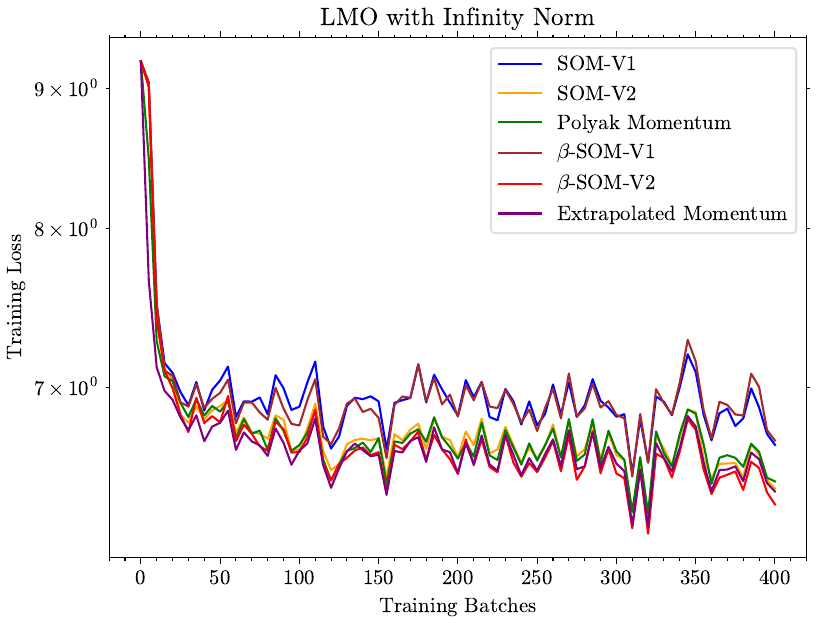}
    \caption{Training loss vs. training batches for the LSTM model using an LMO update with the Infinity norm. The optimization process is successful but visibly less stable than with the Euclidean norm.}
    \label{ch8:fig:lstm_linf_loss_chapter_8}
\end{figure}
The most immediate observation is that while all algorithms still successfully guide the training process toward a lower loss, the convergence path is significantly less stable. The loss curves for all methods are characterized by high-frequency, large-amplitude oscillations.

This instability is an expected consequence of the $\ell_\infty$ LMO. Unlike the smooth, normalized direction vector produced by the $\ell_2$ norm, the $\ell_\infty$ LMO generates a more aggressive update where every component of the step vector is pushed to its maximum value. This results in a more chaotic exploration of the parameter space, leading to the observed volatility in the training loss.

An important consequence of this instability is that the clear performance hierarchy observed in the Euclidean experiment has vanished. The loss curves for all six algorithms are tightly intertwined, and no single method demonstrates a consistent advantage. This suggests that the high variance and non-smooth nature of the $\ell_\infty$ update step dominate the more subtle directional corrections offered by the second-order and extrapolated momentum terms.

While the final loss values appear comparable to those achieved with the $\ell_2$ norm, the convergence path is significantly less stable, making the Euclidean norm the more reliable and predictable choice for this particular language modeling task.

\newpage
  \endgroup

\end{document}